\documentclass[11pt,a4paper]{book}

\usepackage{amsmath,amssymb,amsfonts,amsthm,latexsym,tikz}
\usepackage{graphicx, bm, yfonts, colortbl}
\usepackage[citecolor=blue,urlcolor=blue,hyperindex,breaklinks]{hyperref}
\usepackage{subcaption,float}
\usepackage{MnSymbol}
\usepackage{imakeidx}
\makeindex

\newcommand{\mA}{\mathcal A}
\newcommand {\tc}{\text{c}}
\newcommand {\ttf}{\text{\rm f}}
\newcommand{\ts}{\text{\rm s}}
\newcommand{\trt}{\text{\rm t}}

\newcommand {\bC}{\mathbb{C}}
\newcommand {\bD}{\mathbb{D}}
\newcommand{\bE}{\mathbb E}
\newcommand {\IN}{\mathbb{N}}
\newcommand{\N}{\mathbb{N}}
\newcommand{\bP}{\mathbb P}
\newcommand {\IR}{\mathbb{R}}
\newcommand{\R}{\mathbb{R}}

\newcommand{\re}{\mathbb{R}}
\newcommand{\rek}{\mathbb{R}^{k}}
\newcommand{\red}{\mathbb{R}^{d}}
\newcommand{\T}{\mathbb{T}}
\newcommand {\Z}{{\mathbb Z}}
\newcommand{\cA}{\mathcal A}
\newcommand {\B}{\mathcal B}
\newcommand{\cB}{\mathcal B}
\newcommand{\cC}{\mathcal {C}}
\newcommand{\cD}{\mathcal {D}}
\newcommand {\F}{\mbox{${\mathcal F}$}}
\newcommand{\cF}{{\mathcal F}}
\newcommand{\tf}{\mathcal{F}}
\newcommand{\cG}{\mathcal G}
\newcommand{\cH}{\mathcal H}
\newcommand{\bH}{\mathbb H}
\newcommand{\cI}{\mathcal I}
\newcommand{\cJ}{\mathcal J}
\newcommand{\Jh}{J^{\text{\rm h}}}
\newcommand{\Jw}{J^{\text{\rm w}}}
\newcommand{\Jf}{J^{\text{\rm f}}}
\newcommand{\cK}{{\mathcal K}}
\newcommand{\cKf}{\cK^{\mbox{\scriptsize\rm f}}}
\newcommand{\cKh}{\cK^{\mbox{\scriptsize\rm h}}}
\newcommand{\cKw}{\cK^{\mbox{\scriptsize\rm w}}}
\newcommand{\cL}{\mathcal L}
\newcommand{\cLL}{\underline{\mathcal L}}

\newcommand{\hcF}{{\cLL} {\mathcal F}}
\newcommand{\hcFT}{\cLL_T {\mathcal F}}

\newcommand{\cM}{\mathcal M}

\newcommand {\cO}{\mathcal O}
\newcommand{\cP}{\mathcal P}
\newcommand{\cX}{\mathcal X}
\newcommand {\IS}{\mbox{${\mathcal S}$}}
\newcommand{\hac}{\mathcal{H}}
\newcommand{\cs}{\mathcal{S}}
\newcommand{\csp}{\mathcal{S}^\prime}
\newcommand{\ca}{\mathcal{C}}
\newcommand {\U}{\mathcal U}
\newcommand {\bfone}{\mathbf 1}
\newcommand{\ep}{\varepsilon}
\newcommand{\barh}{\bar{h}}
\newcommand{\laplacef}{\null_\delta G_a}
\newcommand{\half}{\frac{1}{2}}

\newcommand{\alphaf}{\alpha^{\mbox{\scriptsize\rm f}}}
\newcommand{\alphah}{\alpha^{\mbox{\scriptsize\rm h}}}
\newcommand{\alphaw}{\alpha^{\mbox{\scriptsize\rm w}}}

\newcommand{\beq}{\begin{equation}}
\newcommand{\eeq}{\end{equation}}
\newcommand{\beqn}{\begin{equation*}}
\newcommand{\eeqn}{\end{equation*}}

\newtheorem{stat}{Statement}[section]
  \newtheorem{examp}[stat]{Example}
   \newtheorem{assump}[stat]{Assumption}
  \newtheorem{prop}[stat]{Proposition}
  \newtheorem{cor}[stat]{Corollary}
  \newtheorem{thm}[stat]{Theorem}
  \newtheorem{lemma}[stat]{Lemma}
  \newtheorem{remark}[stat]{Remark}
  \newtheorem{def1}[stat]{Definition}

\numberwithin{equation}{section}
\numberwithin{subsection}{section}
\numberwithin{figure}{chapter}

\begin{document}
\frontmatter

\begin{titlepage} 

	\centering 
	
	\scshape 
	
	\vspace*{\baselineskip} 
	
	
	\rule{\textwidth}{1.6pt}\vspace*{-\baselineskip}\vspace*{2pt} 
	\rule{\textwidth}{0.4pt} 
	
	\vspace{0.65\baselineskip} 
	
	{\LARGE Stochastic Partial Differential Equations,\vskip 0.2cm
	Space-time White Noise \vskip 0.3cm
	and Random Fields} 
	
	\vspace{0.75\baselineskip} 
	
	\rule{\textwidth}{0.4pt}\vspace*{-\baselineskip}\vspace{3.2pt} 
	\rule{\textwidth}{1.6pt} 
	
	\vspace{2\baselineskip} 
	
	
	
	\vspace*{3\baselineskip} 
	

	\vspace{0.5\baselineskip} 
	
	{\scshape {\Large Robert C.~Dalang}, EPFL, Switzerland} 
	\vspace{0.5\baselineskip}
	
		{\scshape and}
		\vspace{0.5\baselineskip}
			
	{\scshape{\Large Marta Sanz-Sol\'e,} University of Barcelona} 
	
	\vspace{0.5\baselineskip} 

	\vfill 
	
	
	
	\vspace{0.3\baselineskip} 
	
	$\copyright$ Robert C. Dalang and Marta Sanz-Solé
	
	\vspace{0.4\baselineskip}
	
	
	Version June 26, 2025 


\end{titlepage}

\vfill
\eject

\thispagestyle{myheadings}
\markboth{R.C.~Dalang and M.~Sanz-Sol\'e}{~}

~

\vskip 1 in

{\bf Author information:}
\vskip 24pt

{\scshape
Robert C.~Dalang}

Institut de Mathématiques

\'Ecole Polytechnique Fédérale de Lausanne (EPFL)

CH-1015 Lausanne

Switzerland
\vskip 12pt

robert.dalang@epfl.ch

\vskip 24pt

{\scshape
Marta Sanz-Sol\'e}

Facultat de Matemàtiques i Informàtica

Universitat de Barcelona

Gran Via de les Corts Catalanes 585

08007 Barcelona

Spain
\vskip 12pt

marta.sanz@ub.edu
\vfill

\eject

   \thispagestyle{empty}
     \vspace*{\stretch{1}}
  {\em To my lovely wife Marie-Alix
  
     To my beautiful children Lorraine and Charles
     
      To my wonderful grandchildren Amelia and Sebastien
      \medskip
      
      \hspace*{6cm} R.C.D.}
     \vspace*{\stretch{0.3}}
   
    {\em To the memory of my parents Celestí and Montserrat}
    
     {\em To my son Roger and granddaughter Laia
     \medskip
     
      \hspace*{6cm} M.S.-S.}
    \vfill
    \eject

\tableofcontents
\listoffigures

\chapter{Preface}
\label{ch0}

\pagestyle{myheadings}
\markboth{R.C.~Dalang and M.~Sanz-Sol\'e}{Preface}

This book is an introduction to the theory of stochastic partial differential equations (SPDEs), a field that emerged in the middle of the 1970's motivated by problems inside mathematics and also from other disciplines, such as physics and biology. Since then, the subject has undergone spectacular growth.

The theory of SPDEs can be viewed as an infinite-dimensional extension of the theory of stochastic differential equations (SDEs). Some initial contributions to foundational aspects of SPDEs are due to  K. Itô \cite{ito84}, the main founder of the theory of SDEs. The development of this mathematical area combines methodologies from stochastic analysis and analysis: functional analysis, partial differential equations, semigroup theory and Fourier analysis, among others.

There are several approaches to the theory of SPDEs. These are inspired by existing traditions in the field of partial differential equations and in particular, in evolution systems. The {\em variational approach} was initiated by É. Pardoux  in \cite{pardoux79} and, in greater generality, by N. Krylov and B. Rozovsky \cite{krylov-rozosky-2007} and B. Rozovsky \cite{rozovskii-1990}. An introduction to this approach can be found in the more recent monograph \cite{lr}. The {\em semigroup approach} extends to the random setting the theory of evolution equations in functional spaces defined by differential operators generated by semigroups. The monograph \cite{dz} by G. Da Prato and J. Zabczyk  (first published in 1992) gives a systematic and self-contained presentation of the theory of SPDEs within that approach. In both the variational and in the semigroup approaches, the solutions to the SPDEs are stochastic processes taking values in function spaces (such as Hilbert or  Banach spaces) or in spaces of distributions.
The {\em random field approach} was pioneered by J.B. Walsh \cite{walsh}. In comparison with the preceding approaches, Walsh's setting provides solutions to SPDEs that are random fields, that is, $\re^d$-valued stochastic processes
indexed by several parameters (time, multidimensional space, etc.). This approach is well-suited to the study of sample path space-time properties of the solutions. There is yet another {\em analytical approach}
 by N. Krylov \cite{krylov-1999} (see also \cite{lototsky-rozovsky-2017}). For a large class of SPDEs, using Sobolev embeding theorems, the function space-valued solutions can be realised as random field solutions. The article \cite{dalang-quer} shows connections between
 some of these approaches, which, in various cases, are essentially equivalent.

There are many published documents originating from courses, conferences and other academic activities devoted to specific questions in SPDEs. Sometimes a brief account of the fundamental theory, in one of the existing approaches, is included. Far from being exhaustive, this is a small, hopefully informative, list (in alphabetic order):  \cite{cerrai-2001}, \cite{Khoshnevisan09Mini}, \cite{dz-1996}, \cite{f08}, \cite{holden-1996}, \cite{khosh}, \cite{khoshnevisan-2016}, \cite{leveque}, \cite{pardoux89}, \cite{pardoux2021}, \cite{sanz-sole-2005}, \cite{zab}, $\ldots$, and the books \cite{chow}, \cite{kotolenez-2008} and \cite{gawarecki-mandrekar-2011}.

There are up to now no books that introduce SPDEs via the random field point of view. By writing this volume, our aim is to fill this gap. We are addressing readers with a background in mathematical
sciences and classical stochastic analysis at the graduate level (such as \cite{ks}, \cite{legall-2013} and \cite{ry})
who wish
to learn about the subject and perhaps continue towards research in this field. We assume no prior knowledge on SPDEs (nor even PDEs), and numerous references
 throughout the book
 point the reader to original sources and to supplementary material.

      We have chosen to focus this volume on SPDEs driven by space-time white noise, because there is a large body of literature devoted to this subject and because the technical background required for this situation is minimal. On the other hand, our presentation is designed so that the key steps in the development of the theory of SPDEs for space-time white noise will carry over with limited changes to more general random noises, in particular to spatially homogeneous noise which is white in time. These extensions will be described elsewhere.

The book consists of two blocks: the core matter (Chapters \ref{ch1} to \ref{ch2'}) and the appendices (\ref{app1}, \ref{app2} and \ref{app3}).
Chapter \ref{ch1} introduces the subject, with a discussion of isonormal Gaussian processes and a description of the many facets of space-time white noise, and contains several motivating examples of SPDEs. Chapter \ref{ch2} presents a theory of stochastic integration with respect to space-time white noise and gives fundamental properties of the stochastic integral. Since this integral is defined as a series of Itô integrals, many of its properties can be deduced from the classical Itô theory. We also discuss it relationship with
 Walsh's stochastic integral with respect to martingale measures.

 The SPDEs studied in this book involve a linear partial differential operator and a space-time white noise, possibly multiplied by a non-linear function of the solution, and a possibly nonlinear drift term.

  Chapter \ref{chapter1'} introduces this topic via linear SPDEs (also called SPDEs with {\em additive noise}). We focus the study on the classical examples of the stochastic heat, wave and fractional heat equations, and carry out a detailed analysis of the sample path regularity of the
 random field solutions. In particular, we provide sharp regularity properties of the sample paths, by exhibiting optimal intervals for their H\"older continuity exponents. In Chapter \ref{ch1'-s5}, we formulate and study a general class of SPDEs, in which additive and multiplicative nonlinearities appear (we refer to these equations as {\em nonlinear SPDEs}); in particular, the noise is multiplied by a possibly nonlinear function of the solution  ({\em multiplicative noise}). We prove a general theorem on existence and uniqueness of (strong) random field solutions in a
 framework that covers a wide class of examples. We assume linear growth conditions on the coefficients, but both global and local Lipschitz conditions are discussed. For the nonlinear stochastic heat, wave and fractional heat equations, we discuss the regularity of their sample paths and we exhibit intervals for their H\"older continuity exponents that match those of the linear versions of these equations studied in Chapter \ref{chapter1'}.
  
 Chapter \ref{ch6} discusses various asymptotic properties of the solution to the stochastic heat equation, mostly in the framework of Markov processes in time. Beginning with the Markov and strong Markov properties, we discuss the existence of invariant and reversible measures for linear and nonlinear stochastic heat equations with additive noise and gradient-type drift, as well as properties such as regularity, convergence in law to the invariant distribution, mixing and irreducibility. In the last section of this chapter, we consider asymptotic properties of the solution on a finite interval $[-L, L]$, as $L \to \infty$, and we show that the solution converges to the solution of the stochastic heat equation on $\R$.
In Chapter \ref{ch2'}, we introduce the notion of {\em weak solution in law} for SPDEs, that complements the results of Chapter \ref{ch1'-s5}, and we prove a theorem on existence and uniqueness of solutions in this sense. Then we present a selection of important topics in the theory of SPDEs
that have been the subject of much research over the last thirty years: the Markov field property, asymptotic bounds on moments of solutions that are useful for studying long-time behavior of the solutions, a comparison theorem for the stochastic heat equation, an introduction to potential theory for SPDEs, and a study of SPDEs wth rough initial conditions.

In Appendix \ref{app1}, we first summarize the main results from the theory of stochastic processes and stochastic analysis that are used throughout the book, and then establish a theorem on the existence of versions of processes with values in spaces of Schwartz distributions. We then give an {\em anisotropic} version of Kolmogorov's continuity criterion, for processes indexed by subsets of Euclidean space, that is useful in the study of SPDEs.  Appendix \ref{app2} is devoted to a systematic presentation, along with detailed proofs, of integrability properties of fundamental solutions and Green's functions associated to the classical linear differential operators (heat, fractional heat and wave operators) and upper and lower bounds on their increments in $L^2$- and $L^p$-norm. Many, but not all, of these results are scattered throughout the literature, and we think that having them all together will be useful to many readers. Appendix \ref{app3} is a toolbox section, containing various results from analysis and PDEs, and in particular, a Gronwall-type lemma that is used throughout the book.

 Each chapter is followed by a short ``Notes'' section, which gives historically important references, original sources and points towards other related important contributions.

This book project started many years ago and has been crafted during and in between many working sessions in Barcelona and Lausanne. A second volume is planned.
Throughout the years, we have benefitted from the excellent working conditions provided by our academic institutions--EPFL-École Polytechnique Fédérale de Lausanne, Switzerland, and the University of Barcelona, and for our research, from the financial support of the research councils of our respective countries:
The Swiss National Science Foundation, the Ministerio de Ciencia, Innovación y Universidades, Spain and the Catalan Department of Research and Universities. We express our thanks to these institutions.
\bigskip

\hskip 1cm Lausanne and Barcelona, June 2025.
\bigskip

\hfill {\em Robert C. Dalang} and 
{\em Marta Sanz-Solé}
\vfill

\noindent{\small Thanks to Uijun Kim for reading a preliminary version of this book and for pointing out several typos, to Le Chen for many useful comments on Sections \ref{rdrough}, \ref{rd1+1anderson} and \ref{rd08_26s1}, and to Víctor de la Torre for producing the figures in Chapter \ref{chapter1'} and Appendix \ref{app2}.  }

\mainmatter

\chapter{Basics on noise and SPDEs}
\label{ch1}

\pagestyle{myheadings}
\markboth{R.C.~Dalang and M.~Sanz-Sol\'e}{Basics on noise and SPDEs}

This chapter is devoted to introducing the notion of random noise, to informally describing SPDEs as extensions of PDEs and, through various examples, to giving motivations for the study of SPDEs. We put the focus on Gaussian random fields with an emphasis on isonormal Gaussian processes and (Gaussian) white noise on $\rek$, and the connections between these notions. Because of its importance throughout the book, we give particular attention to  (Gaussian) space-time white noise and describe its many facets, including several series representations. Using  the classical approach in several familiar examples of PDEs, we anticipate the form of the random field solutions to SPDEs which will be developed in later chapters. We close this chapter with a selection of examples of SPDEs that appear in models coming from various disciplines. 


\section{Noise in It\^o stochastic differential equations}
\label{ch1-1.1}
In this section, as a prelude to the study of (space-time) white noise, we examine somewhat informally the noise that drives classical stochastic differential equations on the nonnegative half-line.

One of the most simple stochastic differential equations on $\IR_+$ is
$$
   dX_{t} = \mu X_{t}\, dt + \sigma X_{t} \, dB_{t} ,\qquad  X_{0} = x_0,
$$
where $(B_t,\, t\ge 0)$ is a standard one-dimensional Brownian motion, $\mu\in \re$, $\sigma>0$, and $x_0 \in \IR$ is given. This equation appears for example in the Black-Scholes mathematical model of a financial market containing one risky investment (see \cite{black-scholes1974}).
In principle, this equation could be written
$$
   \frac{dX_{t}}{dt} = \mu\, X_{t} + \sigma X_{t}\, \frac{dB_{t}}{dt},
$$
in which case $\frac{dB_{t}}{dt}$ would be a ``white noise in time." Therefore, we would like to think of white noise in $\IR_{+}$ as the derivative of Brownian motion. However, the meaning of this derivative is not  clear. Indeed, it is well-known that for almost all $\omega \in \Omega$,  the map $B(\omega) : \IR \to \IR$, defined by $t \mapsto B_{t}(\omega)$, is nowhere differentiable on $\re_+$ (here, we set $B_{t} = 0$ if $t  < 0)$.

 Recall that a.s., $B(\omega) \in \mathcal{C}(\IR)$, and by the strong law of large numbers for Brownian motion,
$$
   \lim_{t \to + \infty}  \frac{B_{t}}{t} =0, \qquad\mbox{ a.s.}
$$
Therefore, $B(\omega)$ is a  {\em slowly growing}\index{slowly growing} continuous function, or in other words, $B(\omega) \in \mathcal{S'}(\IR)$ a.s., where $\mathcal{S'}(\IR)$ denotes the space of tempered distributions\index{tempered distribution}\index{distribution!tempered}  (also called Schwartz distributions; see \cite{schwartz}).\index{Schwartz distribution}\index{distribution!Schwartz}

 In particular, $\dot B := \frac{dB}{dt} \in \mathcal{S'}(\IR)$ is well-defined, as the derivative of a Schwartz distribution, by the following property: a.s., for all $\varphi \in \mathcal{S}(\IR)$ (the Schwartz space of $C^{\infty}$ functions with rapid decrease at $\pm\infty$),
\beq
\label{1.1.0}
 \langle\dot B, \varphi \rangle = - \langle B, \dot \varphi \rangle = - \int_{0}^{\infty} B_{t}\, \dot \varphi (t)\, dt.
\eeq

It is therefore possible to view white noise as a tempered distribution. In order to study white noise, a certain amount of discussion of Schwartz spaces and tempered distributions is necessary.
\medskip

\noindent{\em Discrete approximation to white noise}
\smallskip

   It is useful to have a discrete object which provides an approximation of white noise. Recall that a {\em simple random walk}\index{simple random walk}\index{random walk!simple} can be viewed as a discrete version of Brownian motion. We will consider here that the simple random walk is a piecewise constant step function $Z= ( Z_{t} ,\, t \in \IR)$ defined for $t \in \IR$ by
$$
   Z_{t} = \sum^{\infty}_{j=1}\, \xi_{j}\, 1_{\{t \geq j\}},
$$
where the $\xi_{j}$ are independent, identically distributed (i.i.d.) random variables such that
\beqn
P\{\xi_{j} = + 1 \} = P\{\xi_{j} = -1 \}= \frac{1}{2}.
\eeqn

    We can now define a {\em discrete noise}\index{discrete noise}\index{noise!discrete} to be the derivative, in the sense of Schwartz distributions, of the step function $t \mapsto Z_{t}$, which, according to standard facts about Schwartz distributions, is
$$
   \dot Z_t := \frac{dZ_t}{dt} = \sum^{\infty}_{j=1}\, \xi_{j}\, \delta_{j}(t),
$$
where $\delta_{j}(t) = \delta_0(t-j)$ denotes the Dirac delta function. In the classical physical interpretation of the Dirac delta functions as impulses, we see that $\dot Z$ is a sum of independent impulses with random signs.

  We now replace the impulses at integer times by impulses on a finer mesh, with a smaller amplitude: for $n \in \IN^{*}$, set
$$
   Z_{t}^{(n)} = \frac{1}{\sqrt{n}} \sum^{\infty}_{j=1}\, \xi_{j}\, 1_{\left\{t \geq \frac{j}{n}\right\}}.
$$
Donsker's theorem (see for instance \cite{billingsley-1999}) tells us that the sequence $(Z^{(n)},\, n\ge 1)$ converges weakly in the space $\mathcal{C}(\IR)$ to Brownian motion, and so we can expect that
$\dot Z^{(n)}$ converges weakly to $\dot B$, in which case we can view white noise as the cummulative effect of many small independent impulses with random signs.
\bigskip

\noindent{\em Relation with stochastic calculus}
\medskip

According to \eqref{1.1.0}, for every $\varphi \in \mathcal{S}(\IR)$,
\beqn
\langle\dot{B}, \varphi\rangle
  = - \int_{0}^{\infty} B_{t}\, d \varphi (t),
\eeqn
where the right-hand side is a Riemann-Stieltjes integral.
Using integration by parts for It\^o integrals (i.e. the It\^o formula), we see that
\beqn
\int_0^\infty B_t\, d\varphi(t) =   \lim_{t \uparrow + \infty} B_{t}\, \varphi(t)
   - B_{0}\, \varphi (0) - \int_{0}^{\infty} \varphi(t)\, d B_{t}.
\eeqn
Since $\varphi \in \cs(\IR)$ and $B_0 = 0$, we conclude that
\beqn
 \langle\dot{B}, \varphi\rangle =  \int_{0}^{\infty} \varphi{(t)\, d B_{t}}.
\eeqn
In particular, another way to define white noise would be to start from stochastic integrals.
Informally,
$$ \int_{0}^{\infty} \varphi(t)\, \dot{B}(t)\, dt = \langle \dot B, \varphi \rangle = \int_{0}^{\infty} \varphi(t)\, d B_{t},
$$
where the integral on the right-hand side is a Wiener integral (\cite{wiener-1923}). Since the Wiener integral is defined up to a null set, which depends on $\varphi$, some care is needed: see Section \ref{ch1-2.3}.

Up to here, we have discussed white noise on $\IR$. The extension to $\IR^{k}$ will be introduced later on.


\section{Gaussian random fields and white noise}
\label{ch1-2}

In this section, we introduce a class of random fields that plays, in the theory of SPDEs, a role similar to that of Brownian motion for stochastic differential equations.

\subsection{Fundamentals}
\label{ch1-2.0}
We begin this subsection by defining the notion of Gaussian random field.

\begin{def1}
\label{ch1-d1}
 Let $\mathbb{T}$ be an arbitrary set. A family  $G=\left(G(t),\, t\in \mathbb{T}\right)$  of real-valued random variables defined on a probability space $(\Omega, \mathcal{F}, P)$ is  a {\em Gaussian random field}\index{Gaussian!random field}\index{random field!Gaussian} if  for all $r \in \IN^{*}$ and $t_{1}, \dots, t_{r} \in \mathbb{T}$, the random vector $(G(t_{1}), \dots, G(t_{r}))$ is Gaussian.
\end{def1}

The {\em finite-dimensional distributions}\index{finite-dimensional!distributions} \index{distributions!finite-dimensional} of $G$ are given by the family $(\mu_{t_{1}, \dots, t_{r}})$ of probability laws of the Gaussian random vectors $(G(t_{1}), \dots, G(t_{r}))$, that is,
\beq
\label{gauss-ch1}
   \mu_{t_{1}, \dots, t_{r}} (A_{1} \times \dots \times A_{r}) = P\{G(t_{1})\in A_{1}, \dots, G(t_{r}) \in A_{r}\},
\eeq
for all $A_{1}, \dots, A_{r} \in \mathcal{B}_\IR$, $(t_{1}, \dots, t_{r}) \in \mathbb{T}^{r}$, $r \in \IN^{*}.$
\vskip 12pt

The covariance function $(s,t)\mapsto C(s,t) = E[G(s) G(t)] - E[G(s))E(G(t)]$ is obviously {\em symmetric} ($C(s,t) = C(t,s)$).  It defines a
{\em nonnegative definite}\index{nonnegative definite}\index{definite!nonnegative} function on $\mathbb{T}$, that is, for all $r \in \IN^{*}$, for all $x_{1}, \dots, x_{r} \in \IR$ and $t_{1}, \dots t_{r} \in \mathbb{T}$,
\beqn
  \sum^{r}_{i=1} \sum^{r}_{j=1} C(t_{i}, t_{j})\, x_{i}\, x_{j} \geq 0.
\eeqn
Indeed,
\beqn
  \sum^{r}_{i=1} \sum^{r}_{j=1} C(t_{i}, t_{j})\, x_{i}\, x_{j} = {\text{Var}} \left(\sum_{i=1}^r x_i G(t_i) \right) \ge 0.
\eeqn
The following classical lemma discusses the existence of a Gaussian random field with a given covariance function.

\begin{lemma}
\label{ch1.l1}
\noindent (1) Let $G=\left(G(t),\, t\in \mathbb{T}\right)$ be a Gaussian random field. The probability measures  $\mu_{t_{1}, \dots, t_{r}}$ defined in \eqref{gauss-ch1} are entirely determined by the {\em mean function} $m (t) = E[G(t)]$ and the {\em covariance function} $C(s,t)$.

\noindent (2) Given functions $m : \mathbb{T} \to \IR$ and $C: \mathbb{T}^{2} \to \IR$ such that $C(s,t) = C(t,s)$, for all $(s,t) \in \mathbb{T}^{2},$ and $C$ is nonnegative definite, there exists a Gaussian random field $G=\left(G(t),\, t\in \mathbb{T}\right)$ with mean function $m$ and covariance function $C$.
\end{lemma}

\begin{proof}
\noindent (1) Fix $t_{1}, \dots, t_{r}\in \mathbb{T}.$  Let $m_{t_{1}, \ldots, t_{r}}= (m(t_{1}), \ldots, m(t_{r})),$ $C_{t_{1}, \ldots, t_{r}} = (c_{i,j})$, where $c_{i,j} = C(t_i,t_j)=\mbox{Cov} (G(t_{i}), G(t_{j}))$.

Since the column vector $G_{t_{1}, \ldots, t_{r}} = (G(t_1),\dots,G(t_r))$ is Gaussian, it is well-known that its probability law is determined by its mean-vector $m_{t_{1}, \ldots, t_{r}}$ and its $r \times r$ variance-covariance matrix $C_{t_{1}, \ldots, t_{r}}$, which can also be written
$$
   C_{t_{1}, \ldots, t_{r}} = E[(G_{t_{1}, \ldots, t_{r}} - m_{t_{1}, \ldots, t_{r}})(G_{t_{1}, \ldots, t_{r}} - m_{t_{1}, \ldots, t_{r}})^{\intercal}].
$$

In fact, if $\det (C_{t_{1}, \ldots, t_{r}}) \not= 0$, then
$\mu_{t_{1}, \dots t_{r}}$ has a density, and
\begin{align*}
\label{ch1.1}
 &\mu_{t_{1}, \dots t_{r}} (A_{1} \times \dots \times A_{r}) =
   \left({(2 \pi)}^{r}\textrm{det}(C_{t_{1}, \ldots, t_{r}})\right)^{-1/2}
\\
 & \quad \quad\quad\times \int_{A_{1} \times \cdots \times A_r} \exp \left(- \frac{1}{2} {(x-m_{t_{1}, \ldots, t_{r}})}^{\intercal} C_{t_{1}, \ldots, t_{r}}^{-1} (x-m_{t_{1}, \ldots, t_{r}})\right)dx.
\end{align*}

In the general case $\det (C_{t_{1}, \ldots, t_{r}}) \ge 0$, let $O_{t_{1}, \ldots, t_{r}}$ be an orthogonal matrix such that $O_{t_{1}, \ldots, t_{r}}C_{t_{1}, \ldots, t_{r}}O_{t_{1}, \ldots, t_{r}}^{\intercal}= \Lambda_{t_{1}, \ldots, t_{r}}$,
where $\Lambda_{t_{1}, \ldots, t_{r}}$ is the diagonal matrix of (nonnegative) eigenvalues $\lambda_1,\ldots, \lambda_r$ of $C_{t_{1}, \ldots, t_{r}}$. Define a random vector $Y_{t_{1}, \ldots, t_{r}}$ by
$$
   Y_{t_{1}, \ldots, t_{r}} = O_{t_{1}, \ldots, t_{r}} (G_{t_{1}, \ldots, t_{r}} - m_{t_{1}, \ldots, t_{r}}).
$$
Then $E[Y_{t_{1}, \ldots, t_{r}}] = 0$ and
\begin{align*}
   &E[Y_{t_{1}, \ldots, t_{r}} Y_{t_{1}, \ldots, t_{r}}^{\intercal}] \\
    &\qquad = E[O_{t_{1}, \ldots, t_{r}}\, (G_{t_{1}, \ldots, t_{r}}-m_{t_{1}, \ldots, t_{r}})\, (G_{t_{1}, \ldots, t_{r}}-m_{t_{1}, \ldots, t_{r}})^{\intercal}\, O_{t_{1}, \ldots, t_{r}}^{\intercal}] \\
   &\qquad = O_{t_{1}, \ldots, t_{r}}\, E[(G_{t_{1}, \ldots, t_{r}}-m_{t_{1}, \ldots, t_{r}})\, (G_{t_{1}, \ldots, t_{r}}-m_{t_{1}, \ldots, t_{r}})^{\intercal}]\, O_{t_{1}, \ldots, t_{r}}^{\intercal} \\
   &\qquad  = O_{t_{1}, \ldots, t_{r}}\, C_{t_{1}, \ldots, t_{r}}\, O_{t_{1}, \ldots, t_{r}}^{\intercal}\\
   &\qquad = \Lambda_{t_{1}, \ldots, t_{r}}.
\end{align*}
Therefore, the components $(Y_{t_1}, \ldots,Y_{t_r})$ of the random vector $Y_{t_{1}, \ldots, t_{r}}$ are independent random variables, and  $Y_{t_j}$ is ${\rm N}(0,\lambda_j)$ if $\lambda_j > 0$, and $Y_{t_j}=0$ if $\lambda_j= 0$.

Since
\beq
\label{ch1(*1)}
G_{t_{1}, \ldots, t_{r}} = O_{t_{1}, \ldots, t_{r}}^{\intercal}\, Y_{t_{1}, \ldots, t_{r}} + m_{t_{1}, \ldots, t_{r}},
\eeq
one checks by direct calculation that the characteristic function of $G_{t_{1}, \ldots, t_{r}}$ is
\beqn
\varphi_{G_{t_{1}, \ldots, t_{r}}}(z) = \exp\left( i\, z^{\intercal}\, m_{t_{1}, \ldots, t_{r}} - \frac{1}{2}z^{\intercal}\, C_{t_{1}, \ldots, t_{r}}\, z\right), \qquad z\in\re^r,
\eeqn
and this is the Fourier transform of  an $r$-dimensional Gaussian distribution with mean vector $m_{t_{1}, \ldots, t_{r}}$ and covariance matrix $C_{t_{1}, \ldots, t_{r}}$. This Gaussian distribution is supported on the subspace spanned by the rows $j$ of  $O_{t_{1}, \ldots, t_{r}}$ for which $\lambda_j >0$, shifted by $m_{t_{1}, \ldots, t_{r}}$.
\smallskip

 \noindent (2)\ Given the functions $m$ and $C$,
 and a family $(Z_t,\, t \in \T)$ of i.i.d. ${\rm N}(0,1)$ random variables,
 for any $r\ge 1$ and $t_1,\ldots,t_r\in \mathbb{T}$, we construct an $r$-dimensional Gaussian random vector $G_{t_1,\dots,t_r}$
 with mean given by $m_{t_1, \dots, t_r} := (m(t_1), \dots, m(t_r))$ and variance-covariance matrix $C_{t_1, \dots, t_r} :=  (C(t_i, t_j),\, i, j = 1, \dots, r)$
by setting $Y_{t_1,\dots,t_r} = (\lambda_1^{1/2} Z_{t_1},\dots,\lambda_r^{1/2} Z_{t_r})$, where $(\lambda_1,\dots, \lambda_r)$ are the eigenvalues of $C_{t_1,\dots,t_r}$, and then using \eqref{ch1(*1)}.
 Denote by
$\mu_{t_{1}, \dots, t_{r}}$ its probability law. Then
 the claim follows from the Kolmogorov Extension Theorem,\index{Kolmogorov!Extension Theorem}\index{Extension!Theorem, Kolmogorov}\index{theorem!Kolmogorov Extension} since the family of probability measures $(\mu_{t_{1}, \dots, t_{r}})$ satisfies the required consistency conditions (see e.g. \cite[Theorem 36.2, p. 510 and p. 523]{billingsley}), as can be checked using characteristic functions.
\end{proof}

\noindent{\em Examples}
\medskip

We end this section with two fundamental examples.
\medskip

\noindent{\em 1. Brownian motion}
\smallskip

Let $\mathbb{T}=\IR_+$, $m(t)=0$ and $C(s,t) = s\wedge t$. Then
\begin{align*}
 \sum^{r}_{i=1} \sum^{r}_{j=1}\, C(t_{i}, t_{j})\, x_{i}\, x_{j} & =  \sum^{r}_{i=1} \sum^{r}_{j=1}\, (t_{i}\wedge t_{j})\, x_{i}\, x_{j}\\
 & = \sum^{r}_{i=1} \sum^{r}_{j=1}\,  x_{i}\, x_{j} \left(\int_0^\infty 1_{[0,t_i]}(s)\, 1_{[0,t_j]}(s)\, ds\right)\\
 & = \int_0^\infty  \left( \sum^{r}_{i=1}\, x_{i}\, 1_{[0,t_i]}(s)\right)^2 ds \ge 0.
 \end{align*}
Hence, the assumptions of Lemma \ref{ch1.l1} part (2) are satisfied. The Gaussian random field thus defined is a  Brownian motion
$(B_t,\, t\ge 0)$.

 \begin{remark}
\label{ch1-r3}
 The Gaussian process $B = ({B}_{t}, \, t \geq 0)$ obtained by applying Lemma \ref{ch1.l1} does not necessarily have continuous sample paths. However, according to its definition, we have
 \beqn
 E[( B_{t_1} - B_{t_2})^2] = |t_1-t_2|,
 \eeqn
 and therefore, for $p>0$,
 \beqn
 E[\vert B_{t_1} - B_{t_2}\vert^p] = C_p\, |t_1-t_2|^{\frac{p}{2}},
 \eeqn
where $C_p =\left(\tfrac{2^p}{\pi}\right)^{\half}\, \Gamma_E\left( \tfrac{p+1}{2}\right)$ (see Lemma \ref{A3-l1}) and $\Gamma_E$ is the Euler Gamma function (defined in \eqref{Euler-gamma}). Hence, by applying Kolmogorov's continuity criterion\index{Kolmogorov!continuity criterion}\index{criterion!Kolmogorov's continuity}\index{continuity criterion!Kolmogorov's} (see e.g. \cite[Theorem 2.1, p. 26]{ry} or Theorem \ref{ch1'-s7-t2}), for any $\gamma\in\,]0,\frac{1}{2}[$, $B$ has a $\gamma$--H\"older {\em continuous version} (also called a continuous {\em modification}) $({\tilde{B}}_{t},\, t \in \IR_{+}),$ that is, (i) a.s., $t \mapsto {\tilde{B}}_{t}$ is locally $\gamma$-H\"older continuous, and (ii) for all $t \in \IR_{+}$, $P\{{B}_{t} = {\tilde{B}}_{t}\} = 1$. This continuous version is a standard Brownian motion.
\end{remark}
\medskip

\noindent{\em 2. The Brownian sheet on $\IR_+^k$}
\smallskip

Let $\mathbb{T}=\IR_+^k$, and denote by $s=(s_1, \ldots, s_k)$ and $t=(t_1, \ldots, t_k)$ generic points of $\IR_+^k$. Define $m(t)=0$ and
$C(s,t) = \Pi_{i=1}^k(s_i\wedge t_i)$.
By a straightforward extension of the calculation in the preceding example, one can check that the hypotheses of
 Lemma \ref{ch1.l1} part (2) are satisfied. This implies the existence of a Gaussian random field $(W_t,\, t\in \IR_+^k)$, called the {\em Brownian sheet on $\IR_+^k$}.\index{Brownian!sheet}\index{sheet!Brownian} This process satisfies
 \beqn
 E[(W_t - W_s)^2] \le C_J\, |t-s|,
 \eeqn
 for any $s, t \in J$, where $J$ is an arbitrary bounded rectangle of $\IR_+^k$.
 Thus, as for Brownian motion, by applying Kolmogorov's continuity criterion \cite[Theorem 2.1, p. 26]{ry} or Theorem \ref{ch1'-s7-t2}, one obtains the existence of a version $\bar W=(W_t,\, t\in \IR_+^k)$ of this process with locally $\gamma$-H\"older continuous sample paths, for any $\gamma\in\, ]0,\frac{1}{2}[$. We will always use this continuous version. For $k=2$, we will refer to this process simply as the {\em Brownian sheet} (or the {\em Wiener sheet})\index{Wiener!sheet}\index{sheet!Wiener}. According to \cite{walsh}, the Wiener sheet was introduced by Kitagawa in \cite{kitagawa-1951}.

\subsection{Isonormal Gaussian processes}

 An important class of Gaussian random fields are stochastic processes indexed by Hilbert spaces.
 As will be shown later, these appear naturally as stochastic integrals of deterministic processes.
  The next definition gives the precise description.

Let $H$ be a real separable Hilbert space with inner product $\langle \cdot,\cdot\rangle_H$ and norm $\Vert \cdot \Vert_H$.

\begin{def1}
\label{ch1-d3}
A stochastic process $W = (W(h),\, h \in H)$ defined on a complete probability space $(\Omega,\F,P)$ is an {\em isonormal Gaussian process}\index{isonormal Gaussian process}\index{Gaussian!process, isonormal}\index{process!isonormal Gaussian} on $H$ if for all $h\in H$, the random variable $W(h)$ is ${\rm N}(0,\Vert h\Vert_H^2)$, and  $E[W(h)W(g)] = \langle h,g\rangle_H$, for all $h,g \in H$.
\end{def1}

\begin{lemma}
\label{ch1-l2}
 If $(W(h),\, h\in H)$ is an isonormal Gaussian process on $H$, then the mapping $h \mapsto W(h)$, from $H$ into $L^2(\Omega)$, is a linear isometry.
\end{lemma}

\begin{proof} The map $h \mapsto W(h)$ clearly preserves norms, since $\Vert h \Vert_H^2 = E[W(h)^2]$, for all $h \in H$. In order to check that this map is linear, observe that for any $a, b \in \IR$ and $h,g \in H$,
\begin{align*}
   &E[(W(ah+bg) - a W(h) - bW(g))^2] \\
	 &\qquad = \Vert ah+bg\Vert_H^2 + a^2 \Vert h\Vert_H^2 +b^2 \Vert g\Vert_H^2 \\
	    & \qquad\qquad - 2 a \langle ah+bg,h\rangle_H - 2 b \langle ah+bg,g\rangle_H +2 ab \langle h,g\rangle_H = 0.
\end{align*}
\end{proof}
\medskip

The preceding lemma tells us that for an isonormal Gaussian process $(W(h), \, h\in H)$, any linear combination of a finite number of random variables $W(h)$ is also Gaussian. Recall that this property characterizes the Gaussian distribution on finite-dimensional spaces. Therefore an isonormal Gaussian process is indeed a Gaussian random field.
\medskip

\noindent{\em Construction of an isonormal process}
\smallskip

The next proposition gives a way to construct an isonormal Gaussian process on $H$ and provides insight on the structure of this class of processes.

\begin{prop}
\label{ch1-p0}
\begin{enumerate}
\item Let $(e_n, \, n \geq 1)$ be a complete orthonormal system\index{complete!orthonormal system}\index{orthonormal system!complete} (CONS)\index{CONS} in $H$ and let $(\xi_n,\, n \geq 1)$  be a sequence of independent standard Normal random variables defined on $(\Omega,\F,P)$. Then, for any $h\in H$, the series
\beq
\label{ch1-p0.1}
   \sum_{n=1}^\infty\, \langle h,e_n\rangle_H \, \xi_n
 \eeq
 converges in $L^2(\Omega)$ to a random variable which we denote by $W(h)$, and the family $(W(h),\, h\in H)$ thus defined is an isonormal Gaussian process on $H$.

\item Conversely, given an isonormal Gaussian process $(W(h),\, h\in H)$ and a CONS $(e_n, \, n \geq 1)$ in $H$, the sequence $(W(e_n),\, n\ge 1)$ consists of independent standard Normal random variables and
\beq\label{rdch1-p0.1}
 W(h) =  \sum_{n=1}^\infty\, \langle h,e_n\rangle_H\, W(e_n).
\eeq
\end{enumerate}
 \end{prop}
 \begin{proof}
 1. The convergence in $L^2(\Omega)$ of the series \eqref{ch1-p0.1} follows easily from the independence of the random variables  $\xi_n$, since
$$
    \sum_{n=1}^\infty\, |\langle h,e_n\rangle_H^2| = \Vert h\Vert_H^2,
 $$
by Parseval's identity. Moreover, since $W(h)$ is defined as the $L^2(\Omega)$-limit of a sequence of centered Gaussian random variables, it is Gaussian and centered, and by independence of the $\xi_n$ and Parseval's identity, for $h,g \in H$,
\beqn
 E[W(h)W(g)] = \sum_{n=1}^\infty\, \langle h,e_n\rangle_H\, \langle g,e_n\rangle_H = \langle h, g\rangle_H.
\eeqn
Hence, $(W(h),\, h\in H)$ is an isonormal Gaussian process on $H$.
\medskip

2. Since $(e_n, \, n \geq 1)$ is orthonormal, it follows from the definition of an isonormal Gaussian process that the random variables $W(e_n)$, $n\ge 1$, are ${\rm N}(0,1)$ and orthogonal. Because an isonormal Gaussian process is a Gaussian random field, the $W(e_n)$, $n\ge 1$, are independent.
Since $h\mapsto W(h)$ is both linear and continuous (because it is an isometry), and since for any $h\in H$, $h = \sum_{n=1}^\infty\, \langle h,e_n\rangle_H\, e_n$, where the series converges in $H$, \eqref{rdch1-p0.1} follows.
\end{proof}


\subsection{White noise on $\IR^{k}$}
\label{ch1-2.1}

Let $\nu$ be a $\sigma$-finite measure\index{sigma@$\sigma$-finite measure}\index{measure!$\sigma$-finite} on $\IR^{k}$ ($k\ge 1$),  that is, there are compact sets $E_{n} \subset E_{n+1}$ such that $\nu(E_{n}) < + \infty$, for all $n\in \mathbb{N^\ast}$, and $\cup^{\infty}_{n=1}\, E_{n} = \IR^{k}$. We denote by  ${\mathcal{B}}^{f}_{\rek}$ the family  $\{A \in \mathcal{B}_{\IR^{k}} : \nu(A) < + \infty \}$.

\begin{def1}
\label{ch1-d2} A (Gaussian) {\em white noise on $\rek$ based on $\nu$}\index{white noise}\index{noise!white} is a Gaussian random field
\beqn
W=(W(A), \, A \in \mathcal{B}^f_{\IR^{k}}),
\eeqn
 defined on some probability space $(\Omega, \F, P)$, with mean function
\beqn
\mu(A) = E[W(A)] = 0
\eeqn
 and covariance function
 \beqn
 C(A,B) = E[W(A) W(B)] := \nu(A \cap B).
 \eeqn
\end{def1}

The existence of white noise based on $\nu$ follows from Lemma \ref{ch1.l1}. Indeed,
 it suffices to check that the covariance function defined above is nonnegative definite. For this, let $x_1, \ldots, x_r\in\re$ and $A_1, \ldots, A_r\in\mathcal{B}^f_{\IR^{k}}$. Then
\begin{align*}
   \sum_{i=1}^r\sum_{j=1}^r\, x_{i}\, x_{j}\, C(A_{i}, A_{j})& =  \sum_{i=1}^r\sum_{j=1}^r\,  x_{i}\, x_{j} \left(\int_{\IR^{k}} 1_{A_{i}}\, (x) 1_{A_{j}} (x)\, \nu (dx)\right)\\
   & = \int_{\IR^{k}} {\left( \sum_{i=1}^r\, x_{i}\, 1_{A_{i}}(x)\right)}^{2} \nu(dx) \geq 0.
\end{align*}

\begin{remark}
\label{ch1-r1}
 In the case where $\nu$  is Lebesgue measure on $\IR^{k}$, we refer to the white noise based on $\nu$  simply as {\em white noise}.
\end{remark}
\medskip

The proposition below gathers some of the most basic properties of white noise based on $\nu$.

\begin{prop}
\label{ch1-p1a-new}
\begin{enumerate}
 \item Let $A, B \in \mathcal{B}^f_{\IR^{k}}$ be such that $A \cap B = \emptyset$.  Then $W(A)$ and $W(B)$ are independent and
$W(A\cup B) = W(A) + W(B)$.
\item Let $(A_n, \, n\ge 1) \subset \mathcal{B}_{\IR^{k}}$ be a decreasing sequence with $\nu(A_1)<\infty$. Set $A:=\cap_{n\ge 1}A_n$. Then $W(A_n) \to W(A)$
in $L^{2} ( \Omega, \F, P)$.
\item Let $(A_n,\, n\ge 1) \subset \mathcal{B}_{\IR^{k}}$ be increasing. Set $A:=\cup_{n\ge 1}\, A_n$ and assume that $\nu(A)<\infty$. Then $W(A_n) \to W(A)$
in $L^{2} ( \Omega, \F, P)$.
\end{enumerate}
\end{prop}
\begin{proof}
 1. The covariance of $W(A)$ and $W(B)$ is $E(W(A) W(B)) = \nu(A \cap B) = 0$.  Since $(W(A), W(B))$ is Gaussian, this proves the claim about independence.
The claim about additivity follows from the fact that
\begin{align*}
   &E[{(W(A \cup B) - W(A) - W(B))}^{2} ]\\ &=  E[W(A\cup B)^2] + E[W{(A)}^{2}] + E[W{(B)}^{2}] - 2 E[W(A\cup B) W(A)] \\
   &\qquad - 2E[W(A\cup B) W(B)] + 2 E[W(A) W(B)]\\
   & = \nu(A \cup B) + \nu (A) + \nu (B) - 2 \nu (A) - 2\nu (B) \\
   & = 0.
\end{align*}

2. By the additivity property established in 1. applied to the disjoint sets $A$ and $A_n\cap A^c$, we have
\beqn
   E\left[(W(A_n) - W(A))^2\right] = E\left[(W(A_n \setminus A))^2\right] = \nu(A_n \setminus A) \longrightarrow 0
\eeqn
as $ n \to \infty$, since $\cap_n(A_n \setminus A) = \emptyset$.

3. As in the proof of claim 2., since the sets $A_n$ and $A\cap A_n^c$ are disjoint,
\beqn
   E\left[(W(A) - W(A_{n}))^{2}\right] = E\left[(W(A \setminus A_n))^2\right] = \nu(A\setminus A_n) \longrightarrow 0
\eeqn
as $ n \to \infty$, since $\cap_n (A \setminus A_n) = \emptyset$.
\end{proof}

\begin{remark}
\label{ch1-r2}
By (2) and (3) of Proposition \ref{ch1-p1a-new}, the mapping $A \mapsto W(A)$ from  ${\mathcal{B}}^f_{\IR^{k}}$ into $L^{2} (\Omega, \F, P)$ is a $\sigma$-additive vector-valued measure.
 However, for fixed $\omega \in \Omega$, $A \mapsto W(A) (\omega)$ is {\em not} a real-valued signed measure.

Indeed, consider  the case $k=1$ and let $W$ be a white noise on $\IR$ based on the measure $\nu(ds) = 1_{\IR_{+}}(s) ds$. Observe that for any
$t\geq 0$, the functions $1_{]-\infty,t]}$ and  $1_{[0,t]}$ are equal  $\nu$--a.e. By defining
\beq
\label{bm}
B_{t} = W (] {-} \infty, t])= W([0, t]), \quad t\geq 0,
\eeq
we obtain a Brownian motion. Hence, the well-known results on its quadratic variation yield
  $$
   \lim_{n \to \infty}\, \sum_{j =1}^{2^{n}} \left(W \left(\left[\frac{j-1}{2^{n}}, \frac{j}{2^{n}}\right]\right)\right)^{2}
  =  \lim_{n \to \infty}\, \sum_{j =1}^{2^{n}} \left(B_{\frac{j}{2^{n}}}- B_{\frac{j-1}{2^{n}}}\right)^{2}
    =1, \quad\mbox{a.s.}
$$
This implies that
$$
   \lim_{n \to \infty}\, \sum_{j=1}^{2^{n}}\, \left\vert W \left(\left[\frac{j-1}{2^{n}}, \frac{j}{2^{n}}\right]\right) \right\vert = + \infty,  \qquad\mbox{a.s.}
$$
Therefore, if $A\mapsto W(A)(\omega)$ were a signed measure, then it could not be $\sigma$-finite (in fact, the total variation measure of every nonempty open set would be infinite).
\end{remark}

\subsection{Constructing an isonormal process from white noise}
\label{ch1-2.2}

   Let $\nu$ be a $\sigma$-finite measure on $\IR^k$ and $H = L^2(\IR^k, \nu)$. Given an isonormal Gaussian process $W$ on $H$, it is straightforward to define a white noise  based on $\nu$. Indeed, for $A \in \B^f_{\IR^k}$, we set
$ \bar W(A) = W(1_A)$ which obviously satisfies the condition in Definition \ref{ch1-d2} and therefore defines a white noise $\bar W$ based on $\nu$.

Conversely, starting from a white noise $\bar W$ based on $\nu$, we will construct an isonormal Gaussian process
$(W(h),\, h\in H)$, as follows.

For $A\in \B^f_{\IR^{k}}$, set
\beqn
   W (1_{A}):= \bar W(A).
\eeqn
Consider the set of simple functions of the form $h=\sum^{r}_{j=1}c_{j} 1_{A_{j}}$, where
$c_{1}, \dots , c_{r} \in \IR$ and  $A_{1}, \dots , A_{r} \in {\mathcal{B}}^f_{\IR^{k}}$ are pairwise disjoint sets. For $h$ of this form, we define
\beq
\label{def-simple}
   W(h) = W \left(\sum^{r}_{j=1}c_{j} 1_{A_{j}}\right) := \sum^{r}_{j=1} c_{j}\, \bar W(A_{j}).
\eeq
Since the $A_j$ are pairwise disjoint, and by the properties of white noise, we have
\begin{align}
\label{iso-simple}
\left\Vert W  \left(\sum^{r}_{j=1}c_{j} 1_{A_{j}}\right)\right\Vert^{2}_{L^{2} (\Omega)} &= E\left[\left(\sum^{r}_{j=1}c_{j}\, \bar W(A_{j})\right)^{2}\right]\notag\\
  &  =  \sum^{r}_{j=1} c_j^{2}\, E[\bar W(A_{j})^{2}] = \sum^{r}_{j=1} c^{2}_{j}\, \nu (A_{j})\notag\\
  & =  \int_{\IR^{k}} {\left( \sum^{r}_{j=1} c_{j}\, 1_{A_{j}}(x)  \right)}^{2}\nu (dx).
\end{align}

The definition \eqref{def-simple} is legitimate, that is, if $h$ can also be written as  $\sum^{m}_{\ell=1} d_{\ell}\, 1_{B_{\ell}}$, with
$d_{1}, \dots , d_{m} \in \IR$ and $B_{1}, \dots , B_{m} \in {\mathcal{B}}^f_{\IR^{k}}$ pairwise disjoint,  then
\beq
\label{coherence}
\sum^{r}_{j=1} c_{j}\, \bar W(A_{j}) = \sum^{m}_{\ell=1} d_{\ell}\, \bar W(B_{\ell}),\quad {\text{a.s.}}
\eeq
Indeed, taking the second moment of the difference of the two terms of this equality, we obtain
\begin{align*}
&E\left[ \left(\sum^{r}_{j=1}\, c_{j}\, \bar W(A_{j}) - \sum^{m}_{\ell=1}\, d_{\ell}\, \bar W(B_{\ell})\right)^2\right]
= E\left[ \left(\sum^{r}_{j=1}\, c_{j}\, \bar W(A_{j})\right)^2\right] \\
&\qquad \quad+ E\left[ \left(\sum^{m}_{\ell=1}\, d_{\ell}\, \bar W(B_{\ell})\right)^2\right]
- 2 E\left[\sum^{r}_{j=1}\sum^{m}_{\ell=1}\, c_{j}\, d_{\ell}\, \bar W(A_{j})\, \bar W(B_{\ell})\right].
\end{align*}
As in \eqref{iso-simple}, we see that this is equal to
\begin{align*}
&\int_{\rek}\left(\sum^{r}_{j=1}\, c_{j}^2\, 1_{A_{j}} + \sum^{m}_{\ell=1}\, d_{\ell}^2\, 1_{B_{\ell}}\, - 2\sum_{j=1}^r
\sum_{\ell=1}^m\, c_j\, d_\ell\, 1_{ A_j\cap B_\ell} \right) d\nu\\
&\qquad\quad = \int_{\rek}\left(\sum^{r}_{j=1}\, c_{j}\, 1_{A_{j}} - \sum^{m}_{\ell=1}\, d_{\ell}\, 1_{B_{\ell}}\right)^2 d\nu = 0,
\end{align*}
since both sums are equal to $h$,
 proving \eqref{coherence}.

Because of \eqref{iso-simple}, on the set of simple functions $h=\sum^{r}_{j=1}c_{j}\, 1_{A_{j}}$, the mapping $h \mapsto W (h)$ is an isometry from $L^{2}(\IR^{k}, \nu)$ into $L^2 (\Omega)$, and one easily checks that this mapping is linear.
Since the set of simple functions is dense in $L^{2} (\IR^{k}, \nu)$, this isometry admits a unique extension from $L^{2} (\IR^{k}, \nu)$ into $L^2 (\Omega)$, given as follows. For a fixed $h \in L^{2} ( \IR^{k},\, \nu)$, let $(h_{n})$ be a sequence of simple functions such that $\Vert h-h_{n} \Vert_{L^{2} (\IR^{k},\, \nu)} \to 0$. Then
$$
   W (h) := \lim_{n \to \infty} W (h_{n}),
$$
where the limit is in $L^{2} (\Omega, \F, P)$.

The above isometry is known as {\em Wiener's isometry}.\index{Wiener!isometry}\index{isometry!Wiener's} By definition, the random variable $W(h)$ is the Wiener integral of $h$ with respect to the white noise $\bar W$:
\beqn
W (h) = \int_{\IR^{k}} h(x)\, \bar W(dx).
\eeqn
It is easy to prove that $W(h)$ does not depend on the particular sequence of simple functions that approximates $h$.
For the sake of simplicity, we will write $W$ instead of $\bar W$ and use the notation
\beq
\label{wiener0}
W (h) = \int_{\IR^{k}} h(x)\, W(dx).
\eeq
Informally, when $\nu(dx) = dx$ is Lebesgue measure, anticipating Example \ref{ch1-r5}, one sometimes writes
$$
W (h) = \int_{\IR^{k}} h(x) \, \dot W(x)\, dx,
$$
in the same way as one sometimes writes the basic property of the Dirac delta function
$$
\langle \delta_0, h \rangle = \int_{\IR^{k}} h(x) \, \delta_0(x) \, dx = h(0).
$$

\begin{prop}
\label{ch1-p1}
 The family $(W (h),\, h \in L^{2} (\IR^{k}, \nu))$ is an isonormal Gaussian process on $L^{2} (\IR^{k}, \nu)$.
\end{prop}

\begin{proof}
It follows directly from the definition that $W(h)$ is Gaussian with mean zero and
$E [W(h)^{2}] = \Vert h \Vert^{2}_{L^{2}(\IR^{k},\, \nu)}$. 
We now check by using polarisation that
$$
   E[W(h_{1}) W(h_{2})] = \langle h_{1}, h_{2}\rangle_{L^2(\IR^{k},\nu)}.
$$
Indeed,
\begin{align*}
   W(h_{1}) W(h_{2}) &= \frac{1}{4} \left(W (h_{1}) +  W (h_{2})\right)^{2} - \frac{1}{4} \left(W(h_{1}) - W(h_{2})\right)^{2} \\
   &= \frac{1}{4} \left(W(h_{1} + h_{2})\right)^{2} - \frac{1}{4} \left(W(h_{1} - h_{2})\right)^{2},
\end{align*}
so
\begin{align*}
   E[W (h_{1}) W(h_{2})] & = \frac{1}{4}\, {\Vert h_{1}+ h_{2} \Vert}^{2}_{L^{2}(\IR^{k},\, \nu)}- \frac{1}{4}\, {\Vert h_{1}- h_{2} \Vert}^{2}_{L^{2}(\IR^{k},\, \nu)}\\
  &  = \langle h_{1}, h_{2}\rangle_{L^2(\IR^{k},\,\nu)}.
\end{align*}
\end{proof}

\noindent{\bf Examples}
\medskip

\noindent{\em 1. The Wiener integral with respect to Brownian motion}
\smallskip

Let $k=1$ and $W$ be a white noise on $\IR$ based on the measure $\nu(ds) = 1_{\IR_+}(s)\, ds$. Consider the Brownian motion $(B_t,\, t\in \IR_+)$ defined in \eqref{bm}. For the sake of simplicity, we will also denote by $B$ its continuous modification (see Remark \ref{ch1-r3}).
\bigskip

\begin{lemma}
\label{cha1-l3}
Let $(W(h),\,  h\in L^{2} (\IR, \nu))$ be the isonormal Gaussian process given in Proposition \ref{ch1-p1}. Then for all $h\in L^{2} (\IR, \nu)$,
 \beqn
W(h)=\int_{0}^{\infty} h(t)\, d{B}_{t},\qquad a.s.,
\eeqn
where the integral on the right-hand is the classical Wiener integral with respect to Brownian motion.
\end{lemma}
\begin{proof}
The conclusion follows from the following remark. Let $h = 1_{]t_1,t_2]}$, $0\le t_1< t_2$. Then, by definition of the Wiener integral and \eqref{bm},
\beqn
\int_{0}^{\infty} h(t)\, d{B}_{t} = B_{t_2}-B_{t_1} = W(]t_1,t_2]) = W(1_{]t_1,t_2]}) = W(h).
\eeqn
By linearity, this identity extends to step functions $ h(t) = \sum^{r}_{j=1} a_{j}\, 1_{]t_{j-1}, t_{j}]} (t)$ and consequently, to every $h\in  L^{2} (\IR, \nu)$,
by the isometry properties of $h \mapsto W(h)$ and $h \mapsto \int_{0}^{\infty} h(t)\, d{B}_{t}$ from $ L^{2} (\IR, \nu)$ to $L^2(\Omega)$.
 \end{proof}
\medskip

If $h\in \mathcal{S}(\IR)$, then $W(h)$ admits the representation given in the next lemma.
Observe that this lemma makes rigorous the informal discussion in the last part of Section \ref{ch1-1.1}.

\begin{lemma}
\label{cha-l4}
If $\varphi \in \mathcal{S}(\IR)$, then
$$
  W(\varphi) = - \int_{0}^{\infty} {B}_{t}\,\dot{\varphi}(t)\, dt,  \qquad \mbox{a.s.}
$$
\end{lemma}
\begin{proof}
We have seen  in the preceding lemma that
\beqn
   W(\varphi) = \int_{0}^{\infty} \varphi(t)\, d {B}_{t},\qquad a.s.
\eeqn
By applying integration by parts for It\^o integrals, and since $\lim_{r\to\infty}\varphi(r) B_r=0$ a.s., $B_0=0$, and $\varphi$ is of bounded variation, we obtain \begin{align*}
  W(\varphi) &= \lim_{r\to\infty}\left[ {B}_{r}\, \varphi (r) - {B}_{0}\, \varphi(0) - \int_{0}^{r} {B}_{t}\, d \varphi(t) \right]\\
   & = - \int_{0}^{+ \infty}  {B}_{t}\,  \dot \varphi(t)\, dt,\qquad a.s.
\end{align*}
\end{proof}
\medskip

\noindent{\em 2. The Wiener integral with respect to the Brownian sheet}
\smallskip

Let us consider the case $k= 2$ and let $(W(A),\, A \in \mathcal{B}^{f}_{\IR^{2}})$ be a white noise based on the measure $\nu(dx) = 1_{\IR^{2}_{+}}(x)\, dx$. We define a two-parameter Gaussian process in a manner similar to that used in \eqref{bm} to derive the Brownian motion from a white noise based on $1_{\IR_{+}}(x)\, dx$. Indeed, for $(t_{1}, t_{2}) \in \IR^{2},$ set
$$
   W_{t_{1}, t_{2}} = W(]{-}\infty, t_{1}] \times ]{-} \infty, t_{2}]) =
 \left\lbrace \begin{array}{ll}
W([0, t_{1}] \times [0, t_{2}]), & \textrm{if } (t_{1}, t_{2}) \in \IR_{+}^{2},\\
 0, &\textrm{otherwise.}
\end{array}\right.
$$

\begin{prop}
\label{ch1-p2}
 $(W_{t_{1}, t_{2}},\, (t_{1}, t_{2}) \in \IR^{2}_{+})$ is a {\em Brownian sheet}.
 \end{prop}
\begin{proof}
 By definition, this process is clearly Gaussian, $E[W_{t_{1}, t_{2}}] =0$ and moreover, for $(t_1,t_2),\, (s_1,s_2) \in\re^2_+$,
\begin{align*}
   E[W_{t_{1}, t_{2}}W_{s_{1}, s_{2}}] &= E[W([0, t_{1}] \times [0, t_{2}])\, W ([0, s_{1}] \times [0, s_{2}])] \\
	&= \nu (([0, t_{1}] \times [0, t_{2}]) \cap [0, s_{1}] \times [0, s_{2}]) \\
	&= \nu ( [0, t_{1} \wedge s_{1}] \times [0, t_{2} \wedge s_{2}]).
\end{align*}
Thus, $ E[W_{t_{1}, t_{2}}W_{s_{1}, s_{2}} ]=(t_{1} \wedge s_{1}) (t_{2} \wedge s_{2})$.
\end{proof}

For the white noise  on $\IR$ based on the measure $1_{\IR_+}(s)\, ds$ and the corresponding isonormal process $(W(h))$, we proved the formula $W(h) = \int_0^\infty h(t)\, dB_t$, where $(B_t,\, t\ge 0)$ is a standard Brownian motion. We now use an analogous identity in the context of the Brownian sheet to define the stochastic integral of a function $h\in L^2(\IR_+^2,\,\nu)$ with respect to the Brownian sheet. 
Namely, we set
\beq
\label{isobs}
      \int_{\IR^{2}_+} h(t_{1}, t_{2})\, dW_{t_{1}, t_{2}}: = W(h),
\eeq
where the right-hand side refers to the isonormal process on $L^2(\IR^2, \nu)$.

Notice that for a rectangle  $A=\, ]a_{1}, b_{1}] \times\, ]a_{2}, b_{2}]\subset \IR^2_+$,
\begin{align*}
   \int_{\IR^{2}_{+}} 1_{A} (t_{1}, t_{2})\, d W_{t_{1}, t_{2}} &  = W(1_A)  = W( ]a_{1},b_{1}] \times\, ]a_{2},b_{2}])\\
	&= W_{b_{1}, b_{2}} - W_{b_{1}, a_{2}}- W_{a_{1}, b_{2}} + W_{a_{1}, a_{2}},
\end{align*}
so this definition coincides with the one given for instance in \cite{wong-zakai-74}.

\vskip 12pt

\noindent{\em White noise on $\IR^{2}_{+}$ as the second cross-derivative of a Brownian sheet ($k=2$)}
\smallskip

The continuous version of a Brownian sheet $\bar W = (W_{t_{1}, t_{2}})$ satisfies the  property
\beqn
\limsup_{t_1+t_2\, \uparrow\, \infty} \frac{W_{t_1,t_2}}{(t_1+t_2)^2} = 0, \quad a.s.
\eeqn
(see \cite{dh015}). Hence for a.a.~$\omega \in \Omega$, $(t_1,t_2) \mapsto W_{t_1,t_2}(\omega)$ is a slowly growing function, that is, $W_{\cdot, \cdot}(\omega)\in \mathcal{S}^{\prime}(\IR^2)$, a.s. Similar to Lemma \ref{cha-l4}, when $h \in \mathcal{S}(\IR^2)$, the Wiener integral $W(h)$ admits the following representation.

\begin{prop}
\label{rd08_03l1}
Let $\bar W = (W_{t_1, t_2})$ be the continuous version of the Brownian sheet defined in Proposition \ref{ch1-p2}. For $\varphi \in \mathcal{S}(\IR^2)$,
$$
  W(\varphi) = \int_{\R_+^2}  {W}_{t_1, t_2}\,\frac{\partial^2}{\partial t_1 \partial t_2}{\varphi}(t_1, t_2)\, dt_1\, dt_2,  \qquad \mbox{a.s.}
$$
\end{prop}


\begin{proof}
Fix $\varphi\in \mathcal{S}(\IR^2)$. Notice that
\beqn
\int_{\IR^2} 1_{]{-}\infty,t_1]\times ]{-}\infty,t_2]}(s_1,s_2)\,\frac{\partial^2}{\partial t_1 \partial t_2}\varphi(t_1,t_2) dt_1 dt_2
= \varphi(s_1,s_2).
\eeqn
Applying the stochastic Fubini's Theorem \ref{ch1'-tfubini}, 
we see that
\begin{align*}
W(\varphi)& = \int_{\re^2_+}  \left( \int_{\IR^2} 1_{]-\infty,t_1]\times ]-\infty,t_2]}(s_1, s_2)\, \frac{\partial^2}{\partial t_1 \partial t_2}\varphi(t_1,t_2)\, dt_1 dt_2\right) dW_{s_1, s_2}\\
& = \int_{\re^2} W(1_{]-\infty,t_1]\times ]-\infty,t_2]}(\cdot,\cdot))\, \frac{\partial^2}{\partial t_1 \partial t_2}\varphi(t_1,t_2)\, dt_1 dt_2 \\
& = \int_{\re^2_+} W_{t_1,t_2}\, \frac{\partial^2}{\partial t_1 \partial t_2} \varphi(t_1,t_2)\, dt_1 dt_2 ,\ {\text{a.s.}}
\end{align*}
This proves the proposition.
\end{proof}
\medskip

\noindent{\em Extension to all $k \geq 1$}
\smallskip

To a white noise $W$ on $\rek$ based on the measure $\nu(dt)=1_{\rek_+}(t)\, dt$, we associate a $k$-parameter Gaussian process $\bar W=(W_{t_1,\dots,t_k},\, (t_1,\dots,t_k) \in \re_+^k)$\index{k@$k$-parameter Brownian sheet}\index{Brownian!sheet, $k$-parameter}\index{sheet!k@$k$-parameter Brownian} defined by
$$
   W_{t_1,\dots,t_k} = \left\{\begin{array}{ll}
      W(]{-} \infty, t_1] \times \cdots \times \,]{-}\infty, t_k]), & \mbox{if } (t_1,\dots,t_k) \in \rek_+,\\
      0, & \mbox{otherwise.}
   \end{array} \right.
$$
Extending Proposition \ref{ch1-p2} and the considerations in the paragraph above, 
we see that the continuous version $\bar W = (W_{t_1,\dots,t_k},\, (t_1,\dots,t_k) \in \rek_+)$ of this process is a $k$-parameter Brownian sheet, and for a.a.~$\omega$, $(t_1,\dots,t_k) \mapsto W_{t_1,\dots,t_k}(\omega)$ is slowly growing (see \cite{dh015}). We define the Wiener integral of a function $h\in L^2(\IR_+^k,\,\nu)$ with respect to the Brownian sheet $\bar W$ 
by
\beqn
      \int_{\IR^{k}_+} h(t_{1},\dots, t_{k})\, dW_{t_{1},\dots, t_{k}}: = W(h),
\eeqn
where the right-hand side refers to the isonormal process on $L^2(\IR^k, \nu)$.

    Noting that for $\varphi\in \mathcal{S}(\IR^k)$,
 \begin{align*}
 &(-1)^k \int_{\IR^k} 1_{]{-}\infty,t_1]\times \cdots \times  ]{-}\infty,t_k]}(s_1,\dots, s_k)\, \frac{\partial^k}{\partial t_1 \cdots \partial t_k}\varphi(t_1,\dots, t_k)\, dt_1 \cdots dt_k \\
&\qquad\qquad = \varphi(s_1,\dots, s_k),
 \end{align*}
and using the stochastic Fubini's Theorem \ref{ch1'-tfubini} as in the proof of Lemma \ref{rd08_03l1}, we see that
\beq
\label{rd08_04e2}
  W(\varphi) = \int_{\R_+^k}   {W}_{t_1,\dots, t_k}\,\frac{\partial^k}{\partial t_1 \cdots \partial t_k}{\varphi}(t_1,\dots, t_k)\, dt_1 \cdots dt_k,  \qquad \mbox{a.s.}
\eeq

\subsection{Distribution-valued versions}
\label{ch1-2.3}

Let
\beq\label{rde1.2.6}
   W=(W(A),\, A \in \cB^f_{\IR^2})
\eeq
be a white noise as in Definition \ref{ch1-d2} with $k = 2$, based on the measure $\nu(dt) = 1_{\IR_+^2}(t)\, dt$. Let $(W(h),\, h \in L^2(\IR^2, dx))$ be the isonormal Gaussian process associated to $W$ as in Proposition \ref{ch1-p1}. Then for $a_1,a_2 \in \IR$ and $\varphi_1,\varphi_2 \in \mathcal{S}(\IR^2)$,
 \beq
\label{rde1.2.600}
   W(a_{1}\,\varphi_{1}+ a_{2}\, \varphi_{2}) = a_{1}\, W(\varphi_{1}) + a_{2}\,   W(\varphi_{2}), \qquad  a.s.
\eeq
However, the null set implicit in the ``a.s." of \eqref{rde1.2.600} depends on $a_{1}$, $a_{2}$, $\varphi_{1}$ and $\varphi_{2}$, so one cannot deduce from \eqref{rde1.2.600} that for a.a.~$\omega \in \Omega$, $\varphi \mapsto W(\varphi)(\omega)$ belongs to $\mathcal{S'} (\IR^2)$ (even if this map were continuous), and in general, this is not the case.

\begin{def1}
\label{ch1-d4}
(1) A family of random variables $X=(X(\varphi),\, \varphi \in \mathcal{S}(\IR^{k}))$ is called a {\em random linear functional}\index{random linear functional}\index{linear!functional, random}\index{functional!random linear} if,  for all $a_{1}, a_{2} \in \IR$ and $\varphi_{1}, \varphi_{2} \in \mathcal{S}(\IR^{k})$,
$X(a_{1}\, \varphi_{1} + a_{2}\, \varphi_{2}) = a_{1}\,  X (\varphi_{1}) + a_{2}\, X (\varphi_{2})$, a.s.

(2) A process $(\hat{X}(\varphi),\, \varphi \in \mathcal{S}(\IR^{k}))$ is a {\em version with values in $\mathcal{S'} (\rek)$ of $X$}\index{version!with values in $\mathcal{S'}$}\index{distribution-valued!version} if:
\begin{description}
\item{(a)} for all $\varphi \in \mathcal{S}(\IR^{k})$, $\hat{X} (\varphi) = X(\varphi)$  a.s.;
\item{(b)}
 for a.a.~$\omega\in \Omega$, the mapping $\varphi \mapsto \hat{X} (\varphi) (\omega) = \hat{X}(\omega) (\varphi)$ is an element of $\mathcal{S'} (\IR^{k})$ (that is, $\hat{X}$ takes values in $\mathcal{S'} (\IR^{k})$ a.s.)
 \end{description}
\end{def1}

In the special case of \eqref{rde1.2.6}, we can modify $W$ slightly so as to create a version with values in $\mathcal{S'} (\IR^2)$.

\begin{prop}
\label{rd08_04p1}
Let $\bar W = (W_t,\, t\in \IR_+^2)$ be the continuous version of the Brownian sheet in Proposition \ref{ch1-p2}. Let
\beqn
  \dot W := \frac{\partial^2}{\partial t_1 \partial t_2} \bar W,
\eeqn
where the derivative is in the sense of (Schwarz) distributions. Then $\dot W$ is a version with values in $\mathcal{S'} (\IR^{2})$ of the white noise in \eqref{rde1.2.6}.
\end{prop}

\begin{proof}
By definition of the second cross derivative (in $\cs'(\IR^2)$) of the Brownian sheet $\bar W$, for almost all $\omega \in \Omega$ (that we omit in the notation), we have
\begin{align}\label{rde129}
\langle \dot W, \varphi\rangle = \left\langle \frac{\partial^2}{\partial t_1 \partial t_2} \bar W, \varphi\right\rangle & =\left \langle \bar W, \frac{\partial^2}{\partial t_1 \partial t_2}\varphi\right\rangle\\ \nonumber
& = \int_{\IR^2_+} W_{t_1,t_2}\, \frac{\partial^2}{\partial t_1 \partial t_2}\varphi(t_1,t_2)\, dt_1 dt_2 .
\end{align}
From the comments that preceed Proposition \ref{rd08_03l1}, 
for almost all $\omega \in \Omega$, the mapping
\beqn
  \varphi \mapsto \dot W(\varphi) = \int_{0}^{\infty} dt_1 \int_{0}^{\infty} dt_2\,  {W}_{t_1, t_2}\,\frac{\partial^2}{\partial t_1 \partial t_2}{\varphi}(t_1, t_2),
\eeqn
defines an element of $\mathcal{S'} (\IR^2)$, and by Proposition \ref{rd08_03l1}, for all $\varphi \in \mathcal{S} (\IR^2)$,
$$
   W(\varphi)= \dot W(\varphi) \quad a.s.
$$
Therefore, $(\dot{W}(\varphi),\, \varphi \in \mathcal{S}(\IR^{k}))$ is a version of $W$ with values in $\mathcal{S'} (\IR^2)$.
\end{proof}

Proposition \ref{rd08_04p1} states that white noise on $\IR_+^2$ is the second cross-derivative of the Brownian sheet, and clearly, using \eqref{rd08_04e2}, Proposition \ref{rd08_04p1} extends to $k \geq 2$, to give an $\mathcal{S'}(\IR^k)$-valued version $\dot W$ of white noise $W$ on $\IR_+^k$, based on $\nu(dt)=1_{\rek_+}(t) dt$, with values in $\mathcal{S'} (\IR^{k})$:
\beqn
  \dot W := (-1)^k \frac{\partial^k}{\partial t_1 \cdots \partial t_k} \bar W  ,
\eeqn
where $\bar W$ is a $k$-parameter Brownian sheet.

The fact that $\dot W$ is the derivative of some order of a slowly growing continuous function is not surprising, since every deterministic element of $\mathcal{S'}(\IR^k)$ is the derivative of some order of a slowly growing continuous function \cite[Théorème VI, Chapter VII, Section 4, p.~239]{schwartz}. However, in general, it may difficult to find this function, and most often, this is not necessary: we can check that a random linear functional has a version with values in $\mathcal{S'}(\IR^k)$ by using general results on the existence of such versions, such as the following theorem (see Appendix \ref{app1}, Corollary \ref{app1-ch1-c1} for its proof), which is a particular case of \cite[Corollary 4.2, p. 332]{walsh} (see also \cite[Chapter 2]{ito84}).

\begin{thm}
\label{ch1-t1}
 Let $(X(\varphi),\, \varphi \in \mathcal{S}(\IR^{k}))$ be a random linear functional which is continuous in $L^p(\Omega)$, for some $p\ge 1$ (that is, $\varphi_{n} \to \varphi$ in $\mathcal{S}(\IR^{k})$ implies  $X(\varphi_{n}) \to X(\varphi)$ in $L^p(\Omega)$). Then $X$ has a version with values in $\mathcal{S'}(\IR^{k})$.
\end{thm}

The proof of this theorem is given in Appendix \ref{app1}, Corollary \ref{app1-ch1-c1}.
In some situations, it is convenient to use a different family of test functions than $\mathcal{S}(\IR^k)$, and its dual space $\mathcal{S'}(\IR^k)$ of distributions: let $\mathcal{D}(\rek)$\label{rdD} denote the set of $\cC^\infty$ functions defined on $\rek$ with compact support, and equipped with the topology defined (for instance) in  \cite[p.~18]{folland}. Its dual space is denoted $\mathcal{D}^\prime(\rek)$,\label{rdD'} which is larger than $\mathcal{S}^\prime(\rek)$. Since $\mathcal{S}^\prime(\rek)$ is continuously embedded in $\mathcal{D}^\prime(\rek)$, a version with values in $\mathcal{S}^\prime(\rek)$ gives rise to a version with values in $\mathcal{D}^\prime(\rek)$. Further, Theorem \ref{ch1-t1} remains valid with $\mathcal{S}(\rek)$ (respectively $\mathcal{S}^\prime(\rek))$ replaced by $\mathcal{D}(\rek)$(respectively $\mathcal{D}^\prime(\rek))$ (\cite{ito84}, \cite{walsh}).

\begin{examp}
\label{ch1-r5}
Let $A\subset \rek$ be a Borel set, $\nu(dx) = 1_A(x)\, dx$,
 and let $W$ be a white noise on $\rek$ based on $\nu$. Consider the isonormal process $(W(h),\, h\in L^2(\rek, \nu))$ of Proposition \ref{ch1-p1}, from which we obtain the random linear functional $(W(\varphi),\, \varphi \in \mathcal{S}(\rek))$. The convergence $\varphi_{n} \to \varphi$ in $\mathcal{S}(\rek)$ implies $\varphi_{n}\to \varphi$ in $L^{2}(\rek,\, \nu)$. Therefore, by the Wiener isometry,  $W(\varphi_{n})\to W(\varphi)$ in $L^{2} (\Omega)$. By Theorem \ref {ch1-t1}, $(W(\varphi),\, \varphi \in \mathcal{S}(\rek))$ has a version
$(\dot W(\varphi),\, \varphi \in \mathcal{S}(\rek))$ with values in $\mathcal{S'} (\rek)$.
\end{examp}

\noindent{\em Brownian sheet and stochastic convolution}
\smallskip

We return to the setting of Proposition \ref{ch1-p2}.
Using the fact that $1_{]{-}\infty, t_{1}] \times ] {-} \infty, t_{2}]}(s_{1}, s_{2})  = 1_{\IR^{2}_{+}} (t_{1}-s_{1}, t_{2}-s_{2})$, we see that
\begin{align}
W_{t_{1}, t_{2}}& = W(]{-}\infty, t_{1}] \times ]{-}\infty, t_{2}])
= W(1_{]{-}\infty, t_{1}] \times ] {-} \infty, t_{2}]})\notag\\
& = W(1_{\IR^{2}_{+}} (t_{1}- \cdot, t_{2} - \cdot))
= \int_{\IR^2_+} 1_{\IR^{2}_{+}} (t_{1}- s_1, t_{2} - s_2)\, dW_{s_1,s_2}. \label{ch1-3}
\end{align}
The last integral is a particular example of\index{stochastic!convolution}\index{convolution!stochastic} {\em stochastic convolution}---a notion that will appear later in the context of SPDEs.
The terminology is suggested by the formula in the following proposition.

\begin{prop}
\label{rd08_05p1}
Let $\dot W$ be the version with values in $\mathcal{D}^\prime(\IR^2)$ of a white noise $W$ based on $\nu(dx) = 1_{\IR^{2}_{+}}(x)\, dx$, and let $\bar W = (W_{t_{1}, t_{2}})$ be the (continuous version of the) Brownian sheet defined in Proposition \ref{ch1-p2}. Then
\begin{equation}
\label{ch1-20}
1_{\IR^{2}_{+}} \ast  \dot W =   ((u_{1}, u_{2}) \mapsto W_{u_1,u_2}), \qquad \mbox{a.s.}
\end{equation}
(where ``$\ast$" denotes the convolution operation\index{convolution!operation} in $\mathcal{D}^\prime(\IR^2)$).
\end{prop}

\begin{proof}
First, notice that we use $\mathcal{D}^\prime(\R^2)$ here in order that the convolution in \eqref{ch1-20} be well-defined. Indeed, the supports of the distributions $1_{\IR^{2}_{+}}$ and $\dot W$ are both equal to $\re_+^2$, which is not a compact set, and neither of these two distributions has rapid decrease \cite[Chap.VII, §5]{schwartz}. However, because $\re_+^2$ is a {\em cone limited from below}, appealing to \cite[Chapitre VI, \S5]{schwartz} (see also \cite[pp.~304-305]{gasquet}), the convolution $1_{\IR^{2}_{+}} \ast  \dot W$ is a well-defined distribution of $\mathcal{D}^\prime(\re^2)$ such that 
\begin{align*}
\left\langle 1_{\IR^{2}_{+}} \ast  \dot W, \varphi\right\rangle &:= 
 \left\langle \dot W, \psi_K\, ( \tilde1_{\IR^{2}_{+}}\ast \varphi)\right\rangle \\
&= \left\langle \dot W(t_1,t_2), \psi_K(t_1,t_2)\left\langle 1_{\IR^{2}_{+}}(u_1-t_1,u_2-t_2), \varphi(u_1,u_2)\right\rangle\right\rangle,
\end{align*}
where $\tilde 1_{\IR^{2}_{+}}(x) = 1_{\IR^{2}_{+}}(-x)$,
$K = {\rm supp}\, \varphi$, and $\psi_K\in \mathcal{D}(\IR^2)$ is such that $\psi_K(t)=1$ whenever $t=(t_1,t_2)$ belongs to a given neighbourhood of the set
\beqn
\re_+^2 \cap(K-\re_+^2):= \left\{t\in\re_+^2:  t+v\in K {\text{ for some }\  v\in\re_+^2} \right\}.
\eeqn

We now prove \eqref{ch1-20}. Fix $\varphi,\, \psi_K \in \mathcal{D}(\IR^2)$ as above.
 Since $\dot W$ is a version with values in $\mathcal{D}^\prime(\R^2)$ of the white noise $W$, we have
\begin{align*}
&\left\langle \dot W(t_1,t_2), \psi_K(t_1,t_2)\left\langle 1_{\IR^{2}_{+}}(u_1-t_1,u_2-t_2),\varphi(u_1,u_2)\right\rangle\right\rangle\\
&\qquad = \dot W\left(\psi_K(t_1,t_2)\int_{t_1}^\infty du_1  \int_{t_2}^\infty u_2\, \varphi(u_1,u_2)\right)\notag\\
& \qquad  = \int_{\IR^{2}_{+}} \psi_K(t_1,t_2) \left(\int_{t_1} ^\infty du_1  \int_{t_2} ^\infty du_2\, \varphi(u_1,u_2)\right) W(dt_1,dt_2),\quad  {\text{a.s.}}
\end{align*}
where the first integral is a Wiener integral.

Apply the stochastic Fubini's theorem (Theorem \ref{ch1'-tfubini}) with $X:=\re_+^2$, $\mu$ equal to Lebesgue measure, and $G(u_1,u_2,t_1,t_2):= \psi_K(t_1,t_2)\, 1_{\re_+^2}(u_1-t_1,u_2-t_2)\, \varphi(u_1,u_2)$ there, to deduce that the stochastic process
\beqn
Z(u_1,u_2):= \int_{\re_+^2} \psi_K(t_1,t_2)\, 1_{\re_+^2}(u_1-t_1,u_2-t_2)\, W(dt_1,dt_2),
\eeqn
has a jointly measurable version in $(u_1,u_2,\omega)$ and that a.s.,
\begin{align*}
&\int_{\IR^{2}_{+}} \psi_K(t_1,t_2) \left(\int_{t_1} ^\infty du_1  \int_{t_2} ^\infty du_2\, \varphi(u_1,u_2)\right)W(dt_1,dt_2)\\
& \qquad = \int_{\re_+^2} du_1 du_2\ \varphi(u_1,u_2)\\
&\qquad\qquad \times \left(\int_{\re_+^2} \psi_K(t_1,t_2) 1_{\IR^{2}_{+}}(u_1-t_1,u_2-t_2)\,  W(dt_1,dt_2)\right)\\
 & \qquad = \int_{K} du_1 du_2\ \varphi(u_1,u_2)\\
  &\qquad\qquad \times\left(\int_{\re_+^2} \psi_K(t_1,t_2) 1_{\IR^{2}_{+}}(u_1-t_1,u_2-t_2)\,  W(dt_1,dt_2)\right)\\
&\qquad = \int_{K} du_1 du_2\  \varphi(u_1,u_2) \left(\int_{\re_+^2} 1_{\IR^{2}_{+}}(u_1-t_1,u_2-t_2)\, W(dt_1,dt_2)\right),
\end{align*}
where the second equality holds beacause the support of $\varphi$ is $K$, and in the last equality, we have used the fact that the conditions $u\in K$,  $t\in\re_+^2$, $u-t\in \R_+^2$ imply $\psi_K(t_1,t_2)=1$.

For any $(u_1,u_2)\in\re_+^2$,
\beqn
\tilde Z(u_1,u_2):=\int_{\IR^{2}_{+}}1_{\IR^{2}_{+}}(u_1-t_1,u_2-t_2)\, W(dt_1,dt_2) = W_{u_1,u_2} \quad {\text{a.s.}}
\eeqn
From the joint measurability property in $(u_1,u_2,\omega)$ of the stochastic integral $\tilde Z(u_1,u_2)$ (see Theorem \ref{ch1'-l-m1}), along with Fubini's theorem, for a.a.~$\omega$, the set
$\{(u_1,u_2)\in\re_+^2: \tilde Z(u_1,u_2,\omega)\ne W_{u_1,u_2}(\omega)\}$ has zero Lebesgue measure.
This yields
\begin{align*}
 &\int_{\IR^{2}_+}du_1 du_2\,  \varphi(u_1,u_2) \left(\int_{\IR^{2}_{+}}1_{\IR^{2}_{+}}(u_1-t_1,u_2-t_2)\, W(dt_1,dt_2)\right)\\
  &\quad\qquad = \int_{\IR^{2}_+}du_1\, du_2\, \varphi(u_1,u_2)\, W_{u_1,u_2}, \quad {\text{a.s.}}
 \end{align*}
  Thus, for any $\varphi\in \mathcal{D}(\IR^2)$, we have proved that
  \begin{align}
  \label{ch1-31}
 \left\langle 1_{\IR^{2}_{+}} \ast  \dot W, \varphi\right\rangle &:=
 \left\langle \dot W, \psi_K\, (\tilde1_{\IR^{2}_{+}}\ast \varphi )\right\rangle \notag\\
& =  \int_{\IR^{2}_+}du_1 du_2\,  \varphi(u_1,u_2)\, W_{u_1,u_2},\quad {\text{a.s.}},
  \end{align}
where the set of probability zero in which the equality \eqref{ch1-31} may fail depends on $\varphi$.

For any $k\ge 1$, the space $\mathcal{D}(\IR^k)$ is separable (i.e.~admits a countable dense subset) (see e.g.~\cite[Corollary 2, p.~144]{rs} or \cite[Chapter 4]{walsh}). Hence, \eqref{ch1-31} yields \eqref{ch1-20}.
\end{proof}

\section{Space-time white noise}
\label{stwn}

Space-time white noise is a special case of white noise as introduced in Subsection \ref{ch1-2.1}.
If $D \subset \rek$ is a nonempty bounded or unbounded open set, then space-time white noise\index{space-time!white noise}\index{white noise!space-time}\index{noise!space-time white} on $\IR_+ \times D$ is obtained by letting $\nu$ be the measure on $\IR^{1+k}$ defined by $\nu(dt,dx) = 1_{\IR_+}(t)\, 1_D(x)\, dtdx$. Since space-time white noise plays a central role in the theory of SPDEs, and for its further use in this book, we gather in this section the key properties of this object.

Throughout this section, $\B^f_{\IR \times D}$ denotes the set of Borel subsets of $\IR \times D$ with finite Lebesgue measure.

\subsection{Main properties}
\label{stwn-s1}

We recall the definition of space-time white noise and recapitulate its fundamental properties.

\begin{def1}\label{rd1.2.18}
A centered Gaussian random field
\beqn
W=\left(W(A),\, A\in \B^f_{\IR \times D}\right),
\eeqn
 defined on some probability space $(\Omega, \F, P)$, is a {\em space-time white noise on $\IR_+ \times D$}
 if it is a white noise on $\R \times D$ based on the measure
  $\nu(dt,dx) = 1_{\IR_+}(t)\, 1_D(x)\, dtdx$. In particular, for all $A, B \in \B^f_{\IR_+ \times D}$, we have
  \beqn
 E[W(A)\, W(B)] = \vert A \cap B\vert,
 \eeqn
  where $\vert \cdot\vert$ denotes Lebesgue measure on $\IR \times D$.
\end{def1}

\begin{prop}
\label{rd1.2.19}
Let $W=(W(A),\, A\in \B^f_{\IR \times D})$ be a space-time white noise on $\IR_+ \times D$. The following properties hold.

(a) If $A, B \in \mathcal{B}^f_{\re \times D}$ are such that $A \cap B = \emptyset$, then $W(A)$ and $W(B)$ are independent and
$W(A\cup B) = W(A) + W(B)$.

(b) If $(A_n,\, n\ge 1) \subset \mathcal{B}_{\re \times D}$ is a decreasing sequence with $|A_1| <\infty$, then $W(A_n) \to W(A)$
in $L^{2} ( \Omega, \F, P)$, where $A:=\cap_{n\ge 1}A_n$.

 If $(A_n,\, n\ge 1) \subset \mathcal{B}_{\re \times D}$ is an increasing sequence and  $A:=\cup_{n\ge 1}A_n$ satisfies $|A| <\infty$, then $W(A_n) \to W(A)$
in $L^{2} ( \Omega, \F, P)$.

(c) The mapping $A \mapsto W(A)$ from $\B^f_{\IR \times D}$ into $L^2(\Omega, \F, P)$ is a $\sigma$-additive vector-valued measure (but for fixed $\omega \in \Omega$, $A \mapsto W(A)(\omega)$ is {\em not} a real-valued signed measure).

(d) An isonormal Gaussian process $\tilde W = (\tilde W(h),\, h \in L^2(\IR_+ \times D))$ is associated to the space-time white noise $W$, so that $\tilde W(1_A) = W(A)$, $A\in \B^f_{\IR \times D}$. For $h \in L^2(\IR \times D)$, $\tilde W(h)$ is termed the Wiener integral of $h$ with respect to $W$ and is denoted
$$
   \tilde W(h) = \int_{\IR \times D} h(t,x) \, W(dt,dx).
$$
We usually write $W$ instead of $\tilde W$ and we sometimes informally write (as in \eqref{rde1.2.12} below)
$$
   W(h) = \int_{\IR \times D} h(t,x) \, \dot W(t,x)\, dtdx.
$$
In the statements (e) to (g) below, we assume that $D = \rek$.

(e) For $a \in \IR$, define
$$
   I(a) = \left\{\begin{array}{ll}
          [0,a], & \mbox{if } a \geq 0,
         \\
      {[a,0]}, & \mbox{if } a < 0.
      \end{array} \right.
$$
For $(t,x) \in \IR_+ \times \rek$ with $x =(x_1,\dots,x_k)$, define
$$
   W(t,x) = W([0,t] \times I(x_1) \times \cdots \times I(x_k)),
$$
and for $t<0$ and $x\in\rek$, let $W(t,x)=0$. Then $(W(t,x),\, (t,x)\in \IR_+ \times \rek)$ is a $(1+k)$-parameter Brownian sheet on $\IR^{1+k}$ (of which we take the continuous version), and so is the restriction of $(W(t,x))$ to each orthant in $\IR_+ \times \rek$, and these restrictions are independent Brownian sheets.

(f) The random linear functional $(W(\varphi),\, \varphi \in \IS(\IR^{1+k}))$ has a version with values in $\IS'(\IR^{1+k})$, which we denote $\dot W = (\dot W(\varphi),\, \varphi \in \IS(\IR^{1+k}))$. For $\varphi \in \IS(\IR^{1+k})$, we sometimes write $\langle \dot W, \varphi\rangle$ instead of $\dot W(\varphi)$ and
\beq\label{rde1.2.12}
   \langle\dot W, \varphi\rangle = \int_{\IR \times \rek} \varphi(t,x) \, \dot W(t,x)\, dtdx = \int_{\IR \times \rek} \varphi(t,x) \, W(dt,dx).
\eeq
We also refer to $\dot W$ as a space-time white noise.

(g) The $\IS'(\IR^{1+k})$-valued version $\dot W$ of $W$ is the $(1+k)$-fold cross-derivative (in $\IS'(\IR^{1+k})$) of the continuous version of
\beqn
\bar W = ((-1)^{L(x)}W(t,x),\, (t,x)\in \IR^{1+k}),
\eeqn
 where $W(t,x)$ is defined in part (e) and $L(x)$ is the number of negative components of $x \in \rek$. That is, a.s., for all $\varphi \in \IS(\IR^{1+k})$,
$$
   \left\langle \frac{\partial^{1+k}}{\partial t \partial x_1 \cdots \partial x_k} \bar W,\, \varphi \right\rangle = (-1)^{1+k} \langle \dot W,\, \varphi\rangle.
$$
\end{prop}

\begin{proof} Statements (a) and (b) follow from the claims 1--3 of Proposition \ref{ch1-p1a-new}. Statement (c) is explained in Remark \ref{ch1-r2}. Statement (d) is explained in Section \ref{ch1-2.2} (in particular, in Proposition \ref{ch1-p1}). Statement (e) is a straightforward extension of Proposition \ref{ch1-p2} to $1+k$ parameters. Statement (f) is discussed in Example \ref{ch1-r5}. Property (g) is an extension to $1+k$ parameters of \eqref{rde129}.
\end{proof}

\subsection{Fourier series representations}
\label{stwn-s2}

In this subsection, we discuss some Fourier series representations of Brownian motion, of the Brownian sheet, and of space-time white noise on $\R_+ \times [0, 1]^k$, $k \in \N^*$. We will use two CONS  $(e_n,\ n\in \N^*)$ and $(g_n,\ n\in \N)$ of $L^2([0,1])$, defined for $n\in \N^*$ by $e_n(x) = \sqrt 2 \sin(n\pi x)$, $g_n(x) = \sqrt 2 \cos(n\pi x)$,  and $g_0(x) = 1$. We will also use the notation $V:= L^2([0,1])$.
\medskip

\noindent{\em A Fourier-cosine series for white noise on $[0, 1]$}
\medskip

The Fourier-sine series $F_\ts B$ of a standard Brownian motion $B = (B_x,\, x \in [0, 1])$ is defined by
\beqn
     F_\ts B(x) = \sqrt{2}\, \sum_{n=1}^\infty U_n \sin(n \pi x),
\eeqn
with
\beqn
   U_n = \langle B, e_n\rangle_V = \sqrt{2} \int_0^1 B_y\, \sin(n \pi y)\, dy.
\eeqn
Integrating by parts, we find that
\beqn
    U_n = - \frac{\sqrt{2} \cos(n\pi)}{n \pi}\, B_1 + \frac{1}{n \pi}\, X_n,\  {\text{where}}\  X_n = \sqrt{2} \int_0^1 \cos(n \pi y)\, dB_y .
    \eeqn
We notice that $(X_n,\, n\in\N^*)$  is a sequence of independent
${\rm{N}}(0, 1)$ random variables, and are these even independent of $B_1$. 

Thus, for $x\in[0,1]$,
\beq
\label{fourier-series-bis}
 F_\ts B(x) = \sqrt 2 B_1\sum_{n=1}^\infty\frac{-\sqrt 2 \cos(n\pi)}{n\pi} \sin(n\pi x)+  \sqrt 2\sum_{n=1}^\infty\frac{1}{n\pi}X_n\sin(n\pi  x).
 \eeq
Computing the Fourier-sine series $F_\ts f$ of $f(x) = x$, $x \in [0, 1]$, we obtain
\beqn
      F_\ts f (x) = \sqrt{2}\, \sum_{n=1}^\infty\, \frac{- \sqrt{2} \cos(n \pi)}{n\pi}\, \sin(n \pi x).
\eeqn
Since $F_\ts f (x) = f(x)$, for $x\in[0,1[$, by comparing this expression with \eqref{fourier-series-bis}, we see that for $x\in[0,1[$,
\beqn
F_\ts B(x )= B_x = xB_1 + \sqrt 2\sum_{n=1}^\infty\frac{1}{n\pi}X_n \sin(n\pi  x)
\eeqn
and, in fact,
\beq
   \label{rdch1(*1)}
     B_x = x B_1 + \sqrt{2}\, \sum_{n=1}^\infty \, \frac{1}{n\pi} X_n \sin(n \pi x),\quad   x \in[0, 1],   
\eeq
because for $x=1$, the series on the right-hand side vanishes.

Formula \eqref{rdch1(*1)} is the famous {\em Paley-Wiener representation}\index{Paley-Wiener!representation}\index{representation!Paley-Wiener} of standard Brownian motion: for any sequence $(X_n,\, n \in \N)$ of i.i.d.~${\rm{N}}(0, 1)$ random variables, the series in \eqref{rdch1(*1)} converges uniformly on $[0, 1]$, a.s.~and in $L^2(\Omega)$ (see e.g. \cite[p.~166]{khosh-2007}).
The equality ${\rm{Cov}}(B_x, B_y) = x \wedge y$ can be verified using the equalities \eqref{rdeB.5.4} and \eqref{rdeB.5.6}. The series in \eqref{rdch1(*1)} is the Fourier series of a Brownian bridge (see Section \ref{ch5-added-1-s2}).

We can use \eqref{rdch1(*1)} to obtain a Fourier-cosine series representation of ``white noise in time". Indeed, according to Lemma \ref{cha-l4}, this is $\dot B(x) = \frac{dB_x}{dx}$, which should have the Fourier-cosine series expansion obtained by differentiating \eqref{rdch1(*1)} term-by-term. More precisely, setting $X_0:=B_1$,
\beqn
    \dot B(x) = X_0 + \sqrt{2}\, \sum_{n=1}^\infty X_n \cos(n \pi x),\quad     x \in [0, 1].
    \eeqn
We then interpret $\langle \dot B, \varphi \rangle$, where $\varphi \in V := L^2([0, 1])$, as
\beq
\label{ch1(*1c)}
    \langle \dot B, \varphi \rangle = \sum_{n=0}^\infty X_n\, \langle \varphi, g_n \rangle_V. 
    \eeq
In particular,
\beqn
   E\left[\langle  \dot B, \varphi \rangle^2\right] = \sum_{n=0}^\infty\, \langle \varphi, g_n \rangle_V^2 = \Vert \varphi \Vert_V^2
   \eeqn
by Parseval's identity.

Fix $\ep>0$ and set $a_n:= n\pi$. For $\omega$ in the event
$
G_\ep := \left\{\sum_{n=1}^\infty  a_n^{-(1 + 2\ep)}\, X_n^2 < \infty \right\}$,
 which has probability $1$ because
 \beqn
 E\left[\sum_{n=1}^\infty\,  \frac{X_n^2}{a_n^{1 + 2\ep}}\right] = \sum_{n=1}^\infty\,  \frac{1}{(n \pi)^{1 + 2\ep}}  < \infty,
 \eeqn
\eqref{ch1(*1c)} defines an element of $\cs'(\R)$ that belongs also to $H^{-\half - \ep}([0, 1])$, where $H^{-\half - \ep}([0, 1])$ is the dual of \beqn
H^{\half + \ep}([0, 1]) := \left\{ f \in L^2([0, 1]): \sum_{n=0}^\infty\, (1 + n^2)^{\half + \ep}\, \langle f, g_n\rangle_V^2 < \infty \right\}.
\eeqn
Indeed,
\begin{align*}
    \langle \dot B(\omega), \varphi \rangle &= \sum_{n=0}^\infty\, X_n(\omega)\, \langle \varphi, g_n\rangle_V\\
    & = \sum_{n=0}^\infty\, \frac{X_n(\omega)}{(1 + n^2)^{\half (\half + \ep)}}\,  (1 + n^2)^{\half (\half + \ep)}\,  \langle \varphi, g_n \rangle_V \\
     &\leq \left(\sum_{n=0}^\infty\, \frac{X_n^2(\omega)}{(1 + n^2)^{\half + \ep}}\,   \sum_{\ell =0}^\infty\, (1 + \ell^2)^{\half + \ep}\,  \langle \varphi, g_\ell \rangle_V^2\right)^\half < \infty
     \end{align*}
when $\omega \in G_\ep $ and  $\varphi \in H^{\half + \ep}([0, 1])$.
\medskip

\noindent{\em A Fourier-sine series for white noise on $[0, 1]$}
\medskip

The Fourier-cosine series $F_\tc B$ of a standard Brownian motion $B = (B_x,\, x \in [0, 1])$ is defined by
\beq\label{rd02_05e2}
    F_\tc B  (x) = \sum_{n=0}^\infty V_n \, g_n(x),
 \eeq
 with 
 \begin{align*}
 V_0 &= \int_0^1 B_z \, dz = \int_0^1 (1-y)\, dB_y,\\
 V_n &= \sqrt 2 \int_0^1 B_y \cos(n\pi y)\, dy = -\frac{\sqrt 2}{n\pi} \int_0^1 \sin(n\pi y)\,dB_y,\quad n\ge 1.
 \end{align*}
 The random variables $V_0$ and $V_n$, $n\ge 1$, are Gaussian, and using the Wiener isometry, we see that in fact, $V_0$ is ${\rm N}(0, 1/3)$ and $V_n$ is ${\rm N}(0, 1/(n^2\pi^2))$. Moreover, the $V_n$, $n\ge 1$, are independent. However, $V_0$ and $V_n$ are not independent because
 \beq
 \label{cov-v0-vn}
 E[V_0 V_n] = -\frac{\sqrt 2}{n^2 \pi^2}.
 \eeq

Defining $Y_0 = \sqrt{3} V_0$ and $Y_n = n\pi V_n$, $n\ge 1$, we deduce the Fourier-cosine series analoguous to \eqref{fourier-series-bis}, namely
 \beq
 \label{rd08_03e1a}
    F_\tc B (x) = \frac{Y_0}{\sqrt{3}} + \sqrt{2}\, \sum_{n=1}^\infty\, \frac{Y_n}{n \pi}\,  \cos(n \pi x) ,\quad       x \in [0, 1],
\eeq
where $Y_0$, $Y_n$, $n \in \N^*$, are all $\rm{N}(0, 1)$, the $Y_n$, $n \in \N^*$ are independent, but $Y_0$ and $Y_n$ are not independent. Indeed, \eqref{cov-v0-vn} is equivalent to 
\beq
\label {cov-y0-yn}
E[Y_0Y_n] = - \frac{\sqrt 6}{n\pi}\, , \qquad n \geq 1.
\eeq
Formula \eqref{rd08_03e1a} is similar to that proposed in
\cite[Exercise 3.9 p. 326]{walsh}.

From \eqref{cov-y0-yn}, we obtain
\beqn
    \sum_{n=1}^\infty E\left[(Y_0Y_n)^2\right] = \frac{6}{\pi^2} \sum_{n=1}^\infty \frac{1}{n^2} = 1= E[Y_0].
\eeqn
Hence, by Parseval's inequality (in the context of the Hilbert space $L^2(\Omega)$ endowed with the inner product $\langle X,Y\rangle=E(XY)$), we deduce that
\beq\label{rd02_06e2}
Y_0 = \sum_{n=1}^\infty \frac{\sqrt 6}{n\pi}\, Y_n,\quad{\text{a.s.}}
\eeq
Along with \eqref{rd08_03e1a}, this yields 
\beq
\label{rd08_03e1-bis}
   F_\tc B (x) = \sqrt{2}\, \sum_{n=1}^\infty\, \frac{Y_n}{n\pi}\, (1 - \cos(n \pi x)), \quad{\text{a.s.}}
  \eeq

As in \eqref{rdch1(*1)}, one can check the uniform a.s.~and $L^2(\Omega)$-convergence of the series in \eqref{rd08_03e1-bis}. Applying \eqref{rdeB.5.6-quatre}, we see that $E[F_\tc B(x)\, F_\tc B(y)] = x\wedge y$, 
therefore, $(F_\tc B(x),\, x \in [0, 1])$ is a standard Brownian motion. 

Because $x \mapsto B_x$ belongs to $V$ a.s., \eqref{rd02_05e2} implies that
\beqn
     \int_0^1 (B_x - F_\tc B (x))^2\, dx = 0, \qquad a.s.
\eeqn
Toegether with the continuity of sample paths, we conclude that for all $x \in [0, 1]$,
\beq\label{rd02_05e3}
    B_x = F_\tc B (x) =  \sqrt{2}\, \sum_{n=1}^\infty\, \frac{Y_n}{n\pi}\, (1 - \cos(n \pi x)), \quad{\text{a.s.}}
\eeq
This is sometimes called a {\em Fourier-Wiener series representation} of Brownian motion (see \cite[Chap. 16]{Kahane-1985}). We will use this representation in Subsection \ref{rd03_06ss1}, for instance.

Differentiating term-by-term in \eqref{rd02_05e3}, we see that $\dot B(x) = \frac{dB_x}{dx} $ will have the following Fourier-sine series expansion:
\beqn
     \dot B(x) = \sqrt{2}\, \sum_{n=1}^\infty\,  X_n  \sin(n \pi x),\quad     x \in [0, 1].
     \eeqn
In this case, we  interpret $\langle \dot B, \varphi \rangle$, where $\varphi \in V$, as
\beqn
     \langle \dot B, \varphi \rangle = \int_0^1 \varphi(x)\, \dot B(x)\, dx := \sum_{n=1}^\infty\, X_n\, \langle \varphi, e_n \rangle_V.
\eeqn
Comparing with \eqref{ch1(*1c)}, this gives another expression for $\dot{B}$ with values in $\cs'(\R)$.
\medskip

\noindent{\em A Fourier series for white noise on $[0,1]$}
\medskip

Consider the CONS of $L^2([0,1])$ defined by $(e^{2 i n \pi  x}, \, n \in \Z)$, or, equivalenty, 
\beqn
    1, \sqrt 2 \cos(2n \pi x), \sqrt 2 \sin(2 n \pi x),\quad n \in \N^*.
\eeqn
The Fourier series $F B$ of standard Brownian motion $B = (B_x,\, x \in [0, 1])$ is
\beq
\label{rd08_03e4}
    F B(x) =  U_0  + \sqrt{2}\sum_{n=1}^\infty [U_n \cos(2n \pi x) + V_n \sin(2 n \pi x)],\quad     x \in [0, 1],
\eeq
where 
\begin{align*}
    U_0 =  \int_0^1 B_y\, dy = \int_0^1 (1 - y)\, dB_y =: \frac{1}{\sqrt{3}} Z_0,
 \end{align*}
and, for $n \geq 1$, 
\begin{align*}
    U_n = \sqrt{2} \int_0^1 B_y \, \cos(2n  \pi y)\, dy , \qquad    
    V_n = \sqrt{2} \int_0^1 B_y \, \sin(2 n  \pi y)\, dy . 
\end{align*}
Integrating by parts, we see that for $n\geq 1$,
\beqn
    U_n =  -\frac{1}{2 n \pi} X_n, \qquad V_n = \frac{1}{2 n \pi} [-\sqrt{2}\, B_1 + Y_n],
\eeqn
where
\begin{align*}
   X_n = \sqrt{2}  \int_0^1 \sin(2 n \pi y)\, dB_y , \qquad Y_n =  \sqrt{2} \int_0^1 \cos(2 n \pi y)\, dB_y.
\end{align*}
Notice that the $X_n$ and $Y_n$, $n \geq 1$, are i.i.d. $\rm{N}(0, 1)$ random variables, and $Z_0$ is $\rm{N}(0, 1)$ and independent of the $Y_n$ but not of the $X_n$, $n \geq 1$. In fact, for $n\geq 1$,
\beqn
    E[Z_0\, X_n] = \sqrt{6} \int_0^1 (1 - y)\, \sin(2 n \pi y)\, dy = \frac{\sqrt{6}}{2 n \pi}.
\eeqn

Substituting back into \eqref{rd08_03e4}, we obtain
\begin{align}\nonumber
   F B(x) &= \frac{Z_0}{\sqrt{3}} - B_1 \sum_{n=1}^\infty \frac{\sin(2 n \pi x)}{n \pi} \\
      &\qquad\qquad + \sqrt{2} \sum_{n=1}^\infty \frac{1}{2 n \pi}  [- X_n \cos(2 n \pi x) + Y_n \sin(2n\pi x)].
\label{rd02_06e1}
\end{align}
An elementary Fourier series calculation shows that
\beqn
     \half - x = \sum_{n=1}^\infty \frac{\sin(2 n \pi x)}{n \pi},\qquad x \in\, ]0, 1[.
\eeqn
Let $X_0 := B_1$ and $A_0 := \frac{Z_0}{\sqrt{3}} - \frac{B_1}{2}$. By \eqref{rd02_06e1},
\beq\label{rd02_06e3}
    F B(x) = x\, X_0 + A_0 + \sqrt{2} \sum_{n=1}^\infty \frac{1}{2 n \pi}  [- X_n \cos(2 n \pi x) + Y_n \sin(2n\pi x)].
\eeq
We easily check that 
\beqn
     E[A_0^2] = \frac{1}{12}, \qquad E[A_0\, X_n] = \frac{1}{n \pi \sqrt{2}}\qquad \text{and }\quad E[A_0\, Y_n] = 0.
\eeqn
Since
\beqn
  \sum_{n=1}^\infty (E[A_0\, X_n])^2 = \sum_{n=1}^\infty \frac{1}{2 n^2 \pi^2} = \frac{1}{12} = E[A_0^2],
\eeqn
we deduce, as in \eqref{rd02_06e2}, that
\beqn
   A_0 = \sum_{n=1}^\infty \frac{1}{n \pi \sqrt{2}}\, X_n.
\eeqn 
Substituting back into \eqref{rd02_06e3}, we obtain
\beq\label{rd02_06e4}
   F B(x) = x X_0 + \sqrt{2} \sum_{n=1}^\infty \frac{1}{2 n \pi}\, [X_n (1 -\cos(2 n \pi x)) + Y_n \sin(2 n \pi x) ].
\eeq
Since the $X_0$, $X_n$ and $Y_n$, $n\ge 1$, are i.i.d. $\rm{N}(0, 1)$ random variables, we see that for any $x, y\in [0,1]$,
\begin{align*}
E\left[F B(x)\, F B(y)\right] &= xy + \sum_{n=1}^\infty \frac{(1-\cos(2\pi n x))\, (1-\cos(2\pi n y))}{2\pi^2 n^2} \\
& \qquad + \sum_{n=1}^\infty \frac{\sin(2\pi n x)\sin(2\pi n y)}{2\pi^2 n^2} .
\end{align*}
Using the identities \eqref{rdeB.5.6-quatre} and \eqref{rdeB.5.6-quatre-bis} for the first and second series respectively, we obtain
\beqn
E\left[F B(x)\, F B(y)\right] =   xy+\half (x\wedge y) + \half (x\wedge y)(1-2(x\vee y))  = x\wedge y .
\eeqn

   As in \eqref{rdch1(*1)}, one can check the uniform a.s. and $L^2(\Omega)$ convergence of the series in \eqref{rd02_06e4}. Therefore, $(F B(x),\, x \in [0, 1])$ has a version with continuous sample paths that is a standard Brownian motion.
   
   Because $x \mapsto B_x$ belongs to $V$ a.s., 
\beqn
    \int_0^1 (B_x - F B(x))^2\, dx = 0, \qquad\text{a.s.}
\eeqn
Together with the continuity of sample paths, we conclude that for all $x \in [0, 1]$,
\beq\label{fourier-sin-cos}
    B_x = F B(x) = x X_0 + \sqrt{2} \sum_{n=1}^\infty \frac{1}{2 n \pi}\, [X_n (1 -\cos(2 n \pi x)) + Y_n \sin(2 n \pi x) ].
\eeq
This is also sometimes called a {\em Fourier-Wiener series representation} of Brownian motion (see \cite[Chap. 16]{Kahane-1985}).

%
Formula \eqref{fourier-sin-cos} gives a third  Fourier series expansion for the white noise $\dot B(x) = \frac{dB_x}{dx}$ on $[0, 1]$:
\beqn
\dot B_x = X_0 + \sum_{n=1}^\infty \left[X_n \sin(2\pi n x) + Y_n\cos(2\pi n x)\right]
\eeqn
and the expression
\begin{align*}
\langle \dot B,\varphi\rangle &= \int_0^1 \varphi(z)\, \dot B_z\, dz\\
& =  X_0\, \langle\varphi,1\rangle_V + \sum_{n=1}^\infty \left[X_n\, \langle \varphi, \sin(2\pi n \ast)\rangle_V + Y_n\, \langle \varphi, \cos(2\pi n \ast)\rangle_V\right],
\end{align*}
for $\varphi \in V$.
\medskip

\noindent{\em A Fourier-sine series for white noise on $[0, 1]^2$ and on $\IR_+ \times [0, 1]$}
\medskip

   Just as for Brownian motion, one can obtain a series representation for the Brownian sheet on $[0, 1]^2$, namely,
   \beq
   \label{ch1(*1d)}
    W(t, x) = 2 \sum_{m,n=1}^\infty \frac{X_{m, n}}{a_m a_n}\,  (1 - \cos(m \pi t)) \, (1 - \cos(n \pi x)),\quad   t\in [0, 1],\   x \in [0, 1],      
    \eeq
 where the $X_{m, n}$ are i.i.d.~${\rm{N}}(0, 1)$ and $a_n := n\pi$, $n \geq 1$. Indeed, the equality
 \beqn
 {\rm{Cov}}(W(t,x), W(s, y)) = (t \wedge s)\, (x \wedge y)
 \eeqn
  follows again from \eqref{rdeB.5.6-quatre}. For any $T > 0$, the a.s.~and $L^2(\Omega)$ uniform convergence on $[0, T] \times [0, 1]$ can be checked in the same way as for \eqref{rdch1(*1)} (see \cite{walsh-1967}).
  
 Writing $\dot W(t, x) = \frac{\partial^2}{\partial t \partial x} W(t, x)$, we obtain a Fourier-sine series representation of white noise on $[0, 1]^2$:
 \beqn
     \dot W(t, x) = 2 \sum_{m,n=1}^\infty X_{m, n}\,  \sin(m \pi t) \sin(n \pi x),
 \eeqn
 and for $\varphi \in L^2([0, 1]^2)$,
 \beq\label{int-bis}
  \langle \dot W, \varphi \rangle = 2 \sum_{m,n=1}^\infty X_{m, n} \int_0^1 dt \int_0^1 dx\, \varphi(t, x)\, \sin(m \pi t) \sin(n \pi x).
 \eeq

 Performing first the sum over $m$ in \eqref{ch1(*1d)}, then using \eqref{rd02_05e3}, we see that
\beq
\label{ch1(*1a)}
    W(t, x) = \sqrt{2}\, \sum_{n=1}^\infty\, \frac{B^n_t }{a_n}\, (1 - \cos(n \pi x)),\quad   t \in[0,1],\   x \in [0, 1],      
    \eeq
where the $(B^n_t)$ are independent standard Brownian motions; assuming that the $B^n_t$ are defined for $t \in \R_+$, formula \eqref{ch1(*1a)} defines in fact a Brownian sheet on $\R_+ \times [0, 1]$.

Taking the representation \eqref{ch1(*1a)}, the space-time white noise $\dot W(t, x) = \frac{\partial^2 W}{\partial t \partial x}(t, x)$ on $\R_+ \times [0, 1]$ can be interpreted as
\beqn
    \dot W(t, x) = \frac{\partial^2}{\partial t \partial x} W(t, x) = \sqrt{2}\, \sum_{n=1}^\infty\, \frac{dB^n_t }{dt}\, \sin(n \pi x).
    \eeqn
Let $\varphi \in V_2 := L^2([0, T] \times [0, 1])$. Then
\begin{align}
\label{int}
    \langle \dot W, \varphi \rangle &= \int_0^T dt \int_0^1 dx\, \varphi(t, x)\, \dot W(t, x)\notag\\
    & = \sqrt{2} \sum_{n=1}^\infty \int_0^T dB^n_t \int_0^1 dx\, \varphi(t, x)\, \sin(n \pi x)\notag\\
     & =  \sqrt{2}\, \sum_{n=1}^\infty \int_0^T dB^n_t\,  \langle \varphi(t, *), e_n \rangle_{V}.
     \end{align}
On the other hand, using the representation \eqref{int-bis}, we see that
 for $T = 1$ and any $\ep > 0$, since $\sum_{m, n = 1}^\infty \frac{1}{(1 + m^2 + n^2)^{1 + 2\ep}} < \infty$, $\dot W$ belongs a.s.~to $H^{-1 - \ep}([0, 1]^2)$, where $H^{-1 - \ep}([0, 1]^2)$ is the dual of
 \begin{align*}
 &H^{1 + \ep}([0, 1]^2)\\
 &\qquad := \left\{ f \in L^2([0, 1]^2): \sum_{m, n=1}^\infty (1 + m^2 + n^2)^{1 + \ep}\, \langle f, e_{m, n}\rangle_{L^2([0, 1]^2)}^2 < \infty \right\}
 \end{align*}
 and $e_{m, n} (t, x) = e_n(t)\, e_m(x)$.

    Extensions of these formulas to $\R_+ \times [0,1]^k$ are similar and are easily obtained. Analogously, one can also obtain a Fourier-cosine series and a Fourier series representation for space-time white noise on $[0, 1]^2$.

\section{From PDEs to SPDEs}
\label{ch1-s3}

The stochastic partial differential equations studied in this book are obtained from (deterministic) partial differential equations (PDEs) by adding a random forcing term. In this section, we give an elementary introduction to these objects. First, we recall briefly three fundamental partial differential equations motivated by physics, namely the heat equation, the wave equation and the Poisson equation. We will denote by $t$ and $x$ the time and space variables, respectively, and by $D$ a domain in $\IR^k$, that is, $D$ is a non-empty open connected subset of $\rek$.
\medskip

\noindent{\em Heat equation}
\smallskip

The heat equation\index{equation!heat}\index{heat!equation} is also known as the {\em diffusion equation}. It describes the evolution in time of the density $u$ of some quantity such as heat or chemical concentration.
For a given function $f: [0,T] \times D \to \re$, it is the PDE
$$
   \frac{\partial u}{\partial t}(t, x) - \Delta u (t, x) = f(t, x),\qquad (t, x) \in \,]0, T] \times D,
$$
with  initial condition $u(0, x) = u_{0}(x)$, for all $x = (x_1,\dots, x_k) \in D$, where $\Delta$ denotes the Laplacian\index{Laplacian!operator}\index{operator!Laplacian}\label{rdlaplacian} operator in $\IR^k$:
\beqn
    \Delta v(x) = \sum_{j=1}^k \frac{\partial^2 v}{\partial x_j^2} (x). 
\eeqn
If $D$ has a boundary $\partial D$, then one must also impose a boundary condition, such as the value of $u (t, x) $, for all $x \in \partial D $ (this is called a {\em Dirichlet} boundary condition),\index{Dirichlet boundary condition}\index{boundary condition!Dirichlet} or the value of the normal derivative $\frac{\partial u}{\partial \vec{n}}$ at the boundary, for all $x \in \partial D$ (this is called a {\em Neumann} boundary condition).\index{Neumann boundary condition}\index{boundary condition!Neumann}
\medskip

\noindent{\em Wave equation}
\smallskip

The wave equation\index{wave!equation}\index{equation!wave} appears in simplified models for a vibrating medium, namely a string if $k=1$, a membrane if $k=2$, and an elastic solid if $k=3$. It is the second order in time PDE
$$
\frac{\partial^{2} u}{\partial t^{2}}(t, x) - \Delta u (t, x) = f(t, x),\qquad (t, x) \in\, ]0, T] \times D,
$$
where $u(t,x)$ typically represents the displacement of position $x$ at time $t$.

There are two initial conditions: the initial position $u(0, x) = u_{0} (x)$ and the initial velocity
$\frac{\partial u}{\partial t} (0, x) = u_{0} (x)$, $x \in D$. If the domain $D$ has a boundary, then one must also specify a boundary condition as in the case of the heat equation.
\medskip

\noindent{\em Poisson equation}\index{Poisson equation}\index{equation!Poisson}
\smallskip

 This is the PDE defined by
\beqn
   \Delta u (x) = f(x),\   x \in D, \quad
    u(x)  =  u_{0}(x),\   x \in \partial D. \\
\eeqn
For a given electrical charge distribution $f$, solving this equation gives the electric potential $u$. The function $u$ can also represent the density of some quantity, like chemical concentration, at equilibrium.
If $f\equiv 0$, then this equation is called the {\em Laplace equation}\index{Laplace!equation}\index{equation!Laplace}. The solutions of the Laplace equation are called harmonic functions. Notice that in this class of PDEs, there is no time variable.
\medskip

The above equations are of the form
$\mathcal{L}u(t,x) = f(t, x)$, where $\mathcal{L}$ is a linear partial differential operator, with some given initial conditions and, if $D$ has a boundary, with appropriate boundary conditions.
They are generically called {\em inhomogeneous PDEs}. Since the operator $\mathcal{L}$ is  linear, the solution $u$ should be a linear functional of the right-hand side $f$. It turns out that this linear functional often has a very specific form.
\medskip

\noindent{\em Basic principle}
\smallskip

For simplicity, we formulate this principle in the case where the partial differential operator $\mathcal{L}$ is first order in time.

There is an object $\Gamma(t, x; s, y)$, which depends on the PDE and is a real-valued function defined on $\{(t,x;s,y)\in (\re_+ \times \rek)^2: \ s<t\}$, such that the solution to the PDE $\mathcal{L}u(t,x) = f(t,x)$  with initial condition $\psi$ (and appropriate boundary conditions) is
\beq
\label{BP}
   u(t, x) =  \int_D \Gamma(t,x;0,y)\, \psi(y)\, dy +\int^{t}_{0} ds \int_D dy \, \Gamma(t, x; s, y)\, f(s, y),
\eeq
$(t,x)\in\, ]0,\infty[ \times \re^k$.
If $D=\rek$ then $\Gamma$ is called the {\em fundamental solution},\index{fundamental solution}\index{solution!fundamental} and for $D\subset \rek$, it is called the {\em Green's function}.\index{Green's function}\index{function!Green's} In non rigorous terms, one can view the function $\Gamma$ as a description of the inverse of the operator $\mathcal{L}$. More details are given in Chapter \ref{chapter1'}.
\medskip

\noindent{\em Some examples of PDEs and their solutions}
\smallskip

Let us consider some particular cases of heat equations where it is possible to exhibit an explicit form of the solution.
\medskip

\noindent {\em 1. The heat equation on $\IR$ with external forcing}
\smallskip

 This equation is given by
\beq
\label{heatpde}
\begin{cases}
\frac{\partial u }{\partial t}(t, x) - \frac{\partial^{2}u}{\partial x^{2}}(t, x) = f(t,x), &  t > 0,\ x \in \IR,\\
u(0, x) = u_0(x), & x \in \IR,
\end{cases}
\eeq
where the function $f$ represents an external forcing. Notice that the domain $\R$ here has no boundary.

It is well-known that the fundamental solution of the heat equation is
\begin{equation}
\label{fdheat}
   \Gamma(t,x; s,y): = \Gamma(t-s, x-y),
\end{equation}
where
\beq
\label{preheatcauchy-1'}
   \Gamma(s,y) = \frac{1}{\sqrt{4 \pi s}}\, \exp \left(- \frac{ y^{2}}{4 s}\right), \quad (s,y)\in\, ]0,\infty[ \times \IR,
\eeq
(see \cite[p.~46]{evans}).
This is the density of a a one-dimensional mean zero Gaussian random variable with variance $2s$.

Assume that $f$ vanishes and $u_0$ is continuous and bounded. Then
\beqn
u(t, x) = \int_{\IR} dy\ \Gamma(t, x-y)\, u_{0}(y)
\eeqn
is $\mathcal{C^\infty}$ for $x\in \re^n$, $t>0$, and satisfies \eqref{heatpde}; see e.g.~\cite[Theorem on p.~209] {john}. However, this is one out of infinitely many solutions of  \eqref{heatpde}, as illustrated in \cite[p.~211]{john}.

In general, under appropriate conditions on $u_0$ and $f$, a classical solution of \eqref{heatpde}
is given by
\beq
\label{ch1-40}
   u(t, x) = \int_{\IR} dy\ \Gamma(t, x-y)\, u_{0}(y) +
   \int_{0}^{t} ds \int_{\IR} dy\, \Gamma(t-s, x-y)\,  f(s,y),
\eeq
$(t,x)\in\, ]0,\infty[ \times \re$, $u(0,x)=u_0(x)$. In particular, both integrals on the right-hand side should be well defined.
\medskip

\noindent {\em 2. The homogeneous heat equation on $[0, L]$ with vanishing Dirichlet boundary conditions}
\smallskip

This is the PDE
\beq
\label{ch1.500}
\begin{cases}
\frac{\partial u }{\partial t}(t, x) - \frac{\partial^{2}u}{\partial x^{2}}(t, x) = 0, &  t > 0,\ x \in\, ]0, L[,\\
u(0, x) = u_0(x), & x \in\, ]0,L[ ,\\
u(t, 0)= u(t, L) =0, & t > 0.
\end{cases}
\eeq
\medskip

Assuming for example that $u_0\in L^2([0,L])$, a solution to this equation in the sense of \eqref{BP} (with $\psi=u_0$ and $f=0$) exists and for any $t>0$, $u(t,\ast)\in L^2([0,L])$. Here, ``$\ast$'' refers to the space variable in $[0,L]$, a notation that will be used throughout the book. We now show how to find the Green's function.
\medskip

\noindent{\em Calculating the Green's function}
\smallskip

Let $e_{n,L}(x) := \sqrt{\frac{2}{L}}\, \sin\left(\frac{n\pi}{L}x\right)$, $n\in \IN^*$. The family $(e_{n,L},\, n \in \IN^*)$ is a complete orthonormal basis of $L^{2} ([0,L])$ whose elements satisfy the vanishing Dirichlet boundary conditions. Furthermore, $e_{n,L}$ is an eigenvector of the differential operator $\frac{d^{2}}{d x^{2}}$ with associated eigenvalue $-\frac{\pi^2}{L^2}n^2$.

Let $a_n(t)=\langle u(t,\ast), e_{n,L}\rangle$, where $\langle \cdot,\cdot\rangle$ denotes the inner product in $L^2([0,L])$. Since $u(t,x) = \sum_{n=1}^\infty\, \langle u(t,\ast),e_{n,L}\rangle\, e_{n,L}(x)$, we informally have
\beqn
\frac{\partial^{2} u}{\partial x^{2}}  = -\sum_{n=1}^{\infty}\, a_{n} (t)\, \frac{\pi^2}{L^2}\,n^2\, e_{n,L}(x),\qquad
\frac{\partial u}{\partial t} = \sum_{n=1}^{\infty}\, a_{n}^\prime(t)\, e_{n,L}(x).
\eeqn
The equation $\frac{\partial u}{\partial t} - \frac{\partial^{2} u}{\partial x^{2}} =0 $ implies that
$$
   a_{n}' (t) = -\frac{\pi^2}{L^2}\,n^{2}\, a_{n} (t).
$$
On the other hand, the initial condition $u(0,\ast) = u_{0}$ implies that the $a_n(0)$ are the Fourier coefficients of $u_0$:
$$
   a_{n}(0) =
   \langle u_{0}, e_{n,L}\rangle.
$$
Solving the differential equation for  $a_{n}$, for every $n$, gives
\beqn
a_{n} (t) = e^{-\frac{\pi^2}{L^2}n^{2}t}\, \langle u_{0}, e_{n,L}\rangle,
\eeqn
 and so
$$
   u(t, x) = \sum_{n=1}^{\infty}\,  e^{-\frac{\pi^2}{L^2}n^{2} t}\,  \langle u_{0}, e_{n,L}\rangle\, e_{n,L}(x).
   $$
Writing explicitly the inner product and using Fubini's theorem yields
\begin{eqnarray}
\label{explicit}
u(t, x) &=& \sum_{n=1}^{\infty}\, e^{-\frac{\pi^2}{L^2}n^{2}t} e_{n,L} (x) \int^{L}_{0} dy\ u_{0} (y)\, e_{n,L} (y)\notag\\
&=& \int_{0}^{L}dy\ u_{0} (y) \sum^{\infty}_{n=1}\, e^{-\frac{\pi^2}{L^2}n^{2} t}\, e_{n,L} (x)\, e_{n,L} (y),
\end{eqnarray}
from which we see that the Green's function of the heat equation \eqref{ch1.500} on $[0,L]$ with Dirichlet boundary conditions is
\begin{align}
\label{ch1.600}
  \Gamma(t,x;s,y) &:= G_L(t-s;x,y)\notag\\
  & = \sum^{\infty}_{n=1}\, e^{-\frac{\pi^2}{L^2}n^{2} (t-s)}\, e_{n,L} (x)\, e_{n,L} (y),\quad t>s\ge 0, \ x,y\in [0,L].
\end{align}
Notice that
\beqn
G_L(t;x,y)= G_L(t;y,x).
\eeqn
An equivalent expression of $G_L$ is given in Lemma \ref{ch1-equivGD} below.

Recall the version of the Poisson summation formula\index{Poisson summation formula}\index{summation formula!Poisson} given in \cite[37.2.2 Theorem, p.~347]{gasquet}: if $f$ and its derivative $f^\prime$ belong to $L^1(\re)$ and if $\hat f(\xi) = \frac{1}{\sqrt{2\pi}}\int_{\re} f(y)\, e^{-i \xi y}\, dy$ denotes its Fourier transform, then for any $z\in\re$,
\beq
\label{poissonformula-ch1}
\sum^{\infty}_{n=-\infty} f(z + 2nL) = \sqrt{\frac{\pi}{2L^2}}\sum^{\infty}_{n=-\infty} \hat f\left(\frac{\pi n}{L}\right) e^{i\frac{\pi}{L}nz}.
\eeq

\begin{lemma}
\label{ch1-equivGD}
The Green's function in \eqref{ch1.600}
 has the equivalent expression
\begin{align}
\label{ch1.6000-double}
   &G_L(t;x,y)\notag\\
   &\qquad = \frac{1}{\sqrt{4 \pi t}}\notag\\
   &\qquad\qquad\times \sum^{+\infty}_{m=-\infty} \left[ \exp \left(-\frac{{(y-x+2mL)}^{2}}{4t}\right)- \exp \left(-\frac{{(x+y+2mL)}^{2}}{4t}\right) \right]\notag\\
   &\qquad  = \sum^{+\infty}_{m=-\infty} \left[\Gamma(t,y-x+2mL) - \Gamma(t,x+y+2mL)\right],
\end{align}
where $\Gamma$ is the Gaussian density defined in \eqref{preheatcauchy-1'}.
\end{lemma}
\begin{proof}
By the definition of the complex exponential, for all $n\ge 1$,
\begin{align*}
e_{n,L} (x)\, e_{n,L} (y) &= \frac{2}{L} \left[\left(\frac{e^{i\frac{\pi}{L}nx}-e^{-i\frac{\pi}{L}nx}}{2i}\right) \left( \frac{e^{i\frac{\pi}{L}ny}-e^{-i\frac{\pi}{L}ny}}{2i}\right)\right]\\
&= \frac{1}{2L}\left(e^{i\frac{\pi}{L}n(y-x)} + e^{-i\frac{\pi}{L}n(y-x)} - e^{i\frac{\pi}{L}n(x+y)} - e^{-i\frac{\pi}{L}n(x+y)}\right).
\end{align*}
Hence,
\begin{align*}
G_L(t;x,y) &= \frac{1}{2L}\sum^{\infty}_{n=1}\,  e^{-\frac{\pi^2}{L^2}n^{2} t}
\left(e^{i\frac{\pi}{L}n(y-x)} + e^{-i\frac{\pi}{L}n(y-x)} \right)\\
&\qquad -  \frac{1}{2L}\sum^{\infty}_{n=1}\,  e^{-\frac{\pi^2}{L^2}n^{2} t}\left(e^{i\frac{\pi}{L}n(x+y)} + e^{-i\frac{\pi}{L}n(x+y)}\right)\\
& = \frac{1}{2L}\sum^{\infty}_{n=-\infty} e^{-\frac{\pi^2}{L^2}n^{2} t}\, e^{i\frac{\pi}{L}n(y-x)}
-  \frac{1}{2L}\sum^{\infty}_{n=-\infty} e^{-\frac{\pi^2}{L^2}n^{2} t}\, e^{i\frac{\pi}{L}n(x+y)}.
\end{align*}
Apply \eqref{poissonformula-ch1} to $f(z) = \frac{1}{\sqrt{4\pi t}} e^{-\frac{z^2}{4t}}$ and $\hat f(\xi) = \frac{1}{\sqrt{2\pi}}e^{-\xi^2 t}$, to see that
 this is equal to
 \beqn
\frac{1}{\sqrt{4 \pi t}} \sum^{+\infty}_{m=-\infty} \left[ \exp \left(-\frac{{(y-x+2mL)}^{2}}{4t}\right)- \exp \left(-\frac{{(x+y+2mL)}^{2}}{4t}\right) \right],
\eeqn
proving \eqref{ch1.6000-double}.
\end{proof}

Given $v: [0,L]\longrightarrow \IR$, we let $v^o: [-L,L]\longrightarrow \IR$ be the {\em odd extension}\index{odd extension}\index{extension!odd} of $v$, that is,
\beqn
    v^o(x) = \left\{\begin{array}{cl}
        -v(-x), & \text{if } x\in\, ]-L,0[, \\
          v(x), &  \text{if } x\in [0,L],
       \end{array}
       \right.
\eeqn
 and we let $v^{o,p}$ be the $2L$-periodic extension\index{periodic extension}\index{extension!periodic} of $v^o$:
 \beq
\label{vop}
v^{o,p}(x)=
v^o(x-2mL),\   {\text{if}}\  x\in\, ](2m-1)L,(2m+1)L], \ m\in\mathbb{Z}.
\eeq
Lemma \ref{ch1-equivGD} has the following consequence.
\begin{prop}
\label{ch1-equivGD-bis}
Let $G_L(t;x,y)$ be as in Lemma \ref{ch1-equivGD}. Fix $x\in [0,L]$, $t>0$, and suppose that
\beq
\label{equivGD-bis.1}
\int_{-\infty}^\infty \Gamma(t,x-y)\, \vert  v^{o,p}(y)\vert \, dy < \infty.
\eeq
Then
\beqn
\int_0^L G_L(t;x,y)\, v(y)\, dy = \int_{-\infty}^\infty \Gamma(t,x-y)\,  v^{o,p}(y)\, dy.
\eeqn
\end{prop}
\begin{proof}
According to Lemma \ref{ch1-equivGD},
\begin{align*}
&\int_0^L G_L(t;x,y)\, v(y)\, dy\\
&\qquad = \int_0^L \sum_{m\in\mathbb{Z}}\left(\Gamma(t,x-y+2mL) - \Gamma(t,x+y+2mL)\right) v(y)\, dy\\
&\qquad=\sum_{m\in\mathbb{Z}}\int_0^L \Gamma(t,x-y-2mL)\, v(y)\, dy\\
&\qquad \qquad -  \sum_{n\in\mathbb{Z}}\int_0^L \Gamma(t,x+y-2nL)\, v(y)\, dy \\
&\qquad= \sum_{m\in\mathbb{Z}}\int_{2mL}^{(2m+1)L} \Gamma(t,x-z)\, v(z-2mL)\, dz\\
&\qquad \qquad - \sum_{n\in\mathbb{Z}} \int_{(2n-1)L}^{2nL} \Gamma(t,x-z)\, v(2nL-z)\, dz\\
&\qquad= \int_{-\infty}^\infty \Gamma(t,x-z)\, v^{o,p}(z)\, dz.
\end{align*}
We note that permuting the integral and sum is justified by assumption \eqref{equivGD-bis.1}
(first do the calculation with $v(y)$ replaced by $\vert v(y) \vert$ and after the first equality, which becomes an inequality, replace $-\Gamma$ by $\Gamma$). In the sum over $m$, we have used the change of variables $z=y+2mL$, and in the sum over $n$, we used $-z=y-2nL$.
\end{proof}

\noindent {\em 3. The homogeneous heat equation on $[0, L]$ with vanishing Neumann boundary conditions}
\smallskip

This is the PDE
\begin{equation}
\label{1'.14-double}
\begin{cases}
 \frac{\partial u}{\partial t}- \frac{\partial^{2}u}{\partial x^{2}} = 0,\ \qquad\qquad (t, x) \in\, ]0,\infty[ \times\, ]0, L[\, ,\\
 u(0, x) = u_{0}(x), \qquad\qquad  x\in\, ]0, L[\, ,\\
\frac{\partial u}{\partial x}(t, 0) = \frac{\partial u}{\partial x}(t, L) =0,\ \  t\in\,]0,\infty[\, .
 \end{cases}
\end{equation}
\smallskip

The Green's function is
\begin{align}
\label{1'.400-double}
 \Gamma(t,x;s,y):&=G_L(t-s;x,y)\notag\\
  &=\sum^{\infty}_{n=0}\, e^{-\frac{\pi^2}{L^2}n^{2} (t-s)}\, g_{n,L} (x)\, g_{n,L} (y),\quad t>s\ge 0,\ x,y\in[0,L],
\end{align}
where
\beqn
g_{0,L}(x)= \frac{1}{\sqrt{L}}, \quad  g_{n,L}(x)= \sqrt{\frac{2}{L}} \cos\left(\frac{n\pi}{L}x\right), \quad n\ge1.
\eeqn
The sequence $(g_{n,L}, \, n \in \IN)$ is a complete orthonormal basis of $L^2([0,L])$, and each $g_{n,L}$ satisfies the Neumann boundary conditions
\beqn
\frac{\partial g_{n,L}}{\partial x}(t, 0) = \frac{\partial g_{n,L}}{\partial x}(t, L) =0.
\eeqn
Formula \eqref{1'.400-double} can be obtained by the same method as \eqref{ch1.600} for Dirichlet boundary conditions.

As in Lemma \ref{ch1-equivGD}, there is another useful expression for $G_L$.

 \begin{lemma}
 \label{ch1-lequivN}
 The Green's function $G_L(t;x,y)$ defined in \eqref{1'.400-double} has the equivalent expression
\begin{align}
\label{1'.15-double}
&G_L(t;x,y)\notag\\
&\qquad = \frac{1}{\sqrt{4\pi t}}\notag\\
&\qquad\qquad\times\sum_{m=-\infty}^\infty \left[ \exp\left(-\frac{(x-y+2mL)^2}{4t}\right) + \exp\left(-\frac{(x+y+2mL)^2}{4t}\right)\right]\notag\\
&\qquad = \sum_{m=-\infty}^\infty \left[\Gamma(t,x-y-2mL)+\Gamma(t,x+y-2mL)\right],
\end{align}
where $\Gamma$ is defined in \eqref{preheatcauchy-1'}.
\end{lemma}

\begin{proof}
The proof follows the same lines as for \eqref{ch1.6000-double}. First, we notice that
\beqn
g_{0,L}(x)\, g_{0,L}(y) = \frac{1}{L},
\eeqn
and for $n\ge 1$,
\begin{align*}
g_{n,L}(x) g_{n,L}(y) &= \frac{2}{L} \left[\left(\frac{e^{i\frac{\pi}{L}nx}+e^{-i\frac{\pi}{L}nx}}{2}\right) \left(\frac{e^{i\frac{\pi}{L}ny}+e^{-i\frac{\pi}{L}ny}}{2}\right)\right]\\
&= \frac{1}{2L}\left(e^{i\frac{\pi}{L}n(x+y)} + e^{-i\frac{\pi}{L}n(x-y)} + e^{i\frac{\pi}{L}n(x-y)} + e^{-i\frac{\pi}{L}n(x+y)}\right).
\end{align*}
Substituting these expressions in the right-hand side of \eqref{1'.400-double}, we have
\beqn
G_L(t;x,y) = \frac{1}{2L}\left[\sum_{n=-\infty}^\infty e^{-\frac{\pi^2}{L^2}n^2 t}\, e^{i\frac{\pi}{L}n(x-y)}
+\sum_{n=-\infty}^\infty e^{-\frac{\pi^2}{L^2}n^2 t}\, e^{i\frac{\pi}{L}n(x+y)}\right].
\eeqn
From this, \eqref{1'.15-double} follows by applying \eqref{poissonformula-ch1} to $f(z) = \frac{1}{\sqrt{4\pi t}}\, e^{-\frac{z^2}{4t}}$, and $\hat f(\xi) = \frac{1}{\sqrt{2\pi}}\, e^{-\xi^2 t}$.
\end{proof}

Given $v: [0,L] \longrightarrow \IR$, we let $v^e: [-L,L]\longrightarrow \IR$ be the even extension\index{even extension}\index{extension!even} of $v$, that is,
\beqn
   v^e(x) = \left\{\begin{array}{cl}
       v(-x), &\text{if } x\in\, ]-L,0[, \\
       v(x), &\text{if } x\in [0,L],
       \end{array}
      \right.
\eeqn
 and we let $v^{e,p}$ be the $2L$-periodic extension\index{periodic extension}\index{extension!periodic} of $v^e$:
\beq
\label{vep}
v^{e,p}(x)=
v^e(x-2mL),\  {\text{if}}\ x\in\, ](2m-1)L,(2m+1)L], \ m\in\mathbb{Z}.
\eeq
Lemma \ref{ch1-lequivN} has the following consequence.
\begin{prop}
\label{ch1-lequivN-bis} Fix $x\in[0,L]$, $t>0$, and suppose that
\beqn
\int_{-\infty}^\infty \Gamma(t,x-y)\, \vert v^{e,p}(y)\vert\, dy < \infty.
\eeqn
Then
\beqn
\int_0^L G_L(t;x,y)\, v(y)\, dy\ = \int_{-\infty}^\infty \Gamma(t,x-y)\,  v^{e,p}(y)\, dy.
\eeqn
\end{prop}
\begin{proof}
The proof is similar to that of Proposition \ref{ch1-equivGD-bis}  and is omitted.
\end{proof}
\medskip

\noindent{\em Towards linear spde's driven by space-time white noise}
\medskip

As has already been mentioned, we are interested in SPDEs of the type $\mathcal{L}u = \tilde f(t,x)$, where $\tilde f$ stands for a {\em random} external forcing. In the simplest cases, $\tilde f$ is just a noise (such as space-time white noise), and we term this class of equations {\em SPDEs with additive noise},\index{additive noise}\index{noise!additive} and also {\em linear SPDEs}.\index{linear!SPDE}\index{SPDE!linear} We end this section with an introductory example of such an equation.
\smallskip

Let $W = (W(A),\, A\in\mathcal{B}^f_{\IR_+\times \IR })$ be a space-time white noise. According to  Proposition \ref{rd1.2.19}(d) and (e), we can construct the associated isonormal Gaussian process $(W(h),\, h\in L^2(\IR_+\times D))$ and the stochastic process $(W(t,x),\, (t,x)\in \IR_+\times \IR)$.
We will consider the $\mathcal{S}^\prime(\IR^{1+1})$-valued version $\dot W = (\dot W(\varphi),\, \varphi\in \mathcal{S}(\IR^{1+1}))$ of  $W$, as given in
Proposition \ref{rd1.2.19}(f). Recall that, according to the statement (g) of this Proposition,
$\dot W(t,x)$ is the  second cross-derivative $\frac{\partial^2}{\partial t \partial x}\tilde W$, where $\tilde W(t,x) = \left(1_{\{x\ge 0\}}-1_{\{x<0\}}\right)W(t,x)$.

Consider the initial value problem (with Dirichlet boundary conditions)
\begin{equation}
\label{ch1.77}
\begin{cases}
 \frac{\partial u}{\partial t}- \frac{\partial^{2}u}{\partial x^{2}} = \dot W(t, x),\qquad (t, x) \in\, ]0,\infty[\, \times\, ]0, L[\, ,\\
 u(0, x) = u_{0}(x), \qquad\qquad  x\in\, ]0, L[\, , \\
 u(t, 0) = u(t, L) =0, \qquad t\in\, ]0,\infty[\ .
 \end{cases}
\end{equation}

By the {\em basic principle} formulated in \eqref{BP} and \eqref{rde1.2.12}, a solution to \eqref{ch1.77} should be
\beq
\label{giveanumber(*2)}
u(t, x)= \int^{L}_0 G_L(t;x,y)\, u_{0}(y)\, dy + \int_{[0,t]\times [0,L]} G_L(t-s; x,y)\, W(ds,dy),
\eeq
$(t,x)\in\, ]0,\infty[\times [0,L]$, $u(0,x)=u_0(x)$, with $G_L(t;x,y)$ given in \eqref{ch1.600} (or in \eqref{ch1.6000-double}). This statement will be made rigorous in Chapter \ref{chapter1'}.
The first integral is the solution to the (deterministic) homogeneous heat equation
$$
   \frac{\partial v}{\partial t}- \frac{\partial^{2}v}{\partial x^{2}} = 0
$$
with the same initial and boundary conditions as \eqref{ch1.77} (given in \eqref{explicit}), while according to Proposition \ref{rd1.2.19}(d), for a space-time white noise,  the second integral is in fact
$$
  W\left( G_L(t-\cdot; x, \ast)\, 1_{[0,t] \times [0, L]}\right),
$$
the isonormal Gaussian process $W$ evaluated at $G_L(t-\cdot; x, \ast) 1_{[0,t] \times [0, L]}$. Here the notation $``\cdot"$ refers to the time variable in $[0,t]$, while $``\ast"$ refers to the space variable in $[0,L]$. Observe that
this makes sense because $(s,y) \mapsto G_L(t-s, x,y)\, 1_{[0,t] \times [0, L]}(s,y)$ belongs to $L^2(\IR_+ \times[0, L],\, dsdy )$.
Indeed, from the expression \eqref{ch1.600} and Parseval's identity,  we see that
\begin{align}
 \label{27.1-Ch1}
&\Vert G_L(t-\cdot;x,\ast) 1_{[0,t] \times [0, L]}\Vert_{L^2(\re_+ \times[0, L])}^2\notag \\
&\qquad = \int_0^t ds \int_0^L dy\, G_L^2(t-s; x,y)
 = \int_0^t ds \int_0^L dy\, G_L^2(s; x,y) \notag\\
& \qquad = \int_0^t ds \sum_{n=1}^\infty e^{-\frac{2\pi^2 }{L^2}n^2 s}\, e^2_{n,L}(x)
  \le \sum_{n=1}^\infty e_{n,L}^2(x)\, \frac{L^2}{2\pi^2n^2} < \infty.
\end{align}
\bigskip

A similar discussion can be done with the Neumann boundary value problem
\begin{equation}
\label{1'.14neumann}
\begin{cases}
 \frac{\partial u}{\partial t}- \frac{\partial^{2}u}{\partial x^{2}} = \dot W(t,x),\qquad (t, x) \in\, ]0,\infty[\, \times\, ]0, L[,\\
 u(0, x) = u_{0}(x), \qquad\qquad  x\in\, ]0, L[\, ,\\
\frac{\partial u}{\partial x}(t, 0) = \frac{\partial u}{\partial x}(t, L) =0, \quad t\in\, ]0,\infty[\, .
 \end{cases}
\end{equation}
Indeed, by using \eqref{1'.400-double}, a direct integration yields
\begin{align}
\label{cansada}
&\Vert G_L(t-\cdot;x,\ast) 1_{[0,t] \times [0, L]}\Vert_{L^2(\re_+ \times[0, L])}^2
= \int_0^t ds \int_0^L dy\,  G_L^2(s; x,y) \notag\\
&\qquad = \int_0^t ds \sum_{n=0}^\infty e^{-\frac{2\pi^2}{L^2}n^2 s}\, g_{n,L}^2(x)= \frac{t}{L} + \sum_{n=1}^\infty\, g^2_{n,L}(x)\int_0^t ds\, e^{-\frac{2\pi^2}{L^2}n^2 s}\notag\\
&\qquad\le \frac{t}{L}+ \sum_{n=1}^{\infty}\, g^2_{n,L}(x)\, \frac{L^2}{2\pi^2 n^2}< \infty.
\end{align}
\bigskip

Through these elementary examples, we see the important role that the Green's function plays in giving a rigorous meaning to the stochastic integral accounting for the random forcing.

In Chapter \ref{chapter1'}, we will undertake a deeper study of the stochastic heat and wave equations with additive noise.
\medskip

\noindent{\em A connection between SDEs and SPDEs}
\medskip

    The formula \eqref{ch1(*1a)} for the Brownian sheet provides a natural connection between SDEs and SPDEs. Indeed, consider the sequence of linear SDEs
    \beq
    \label{ch1(*2)}
     dY^n_t = - (a_n^2 + a^2)\, Y^n_t\, dt + \sigma\, dB^n_t,\quad   t > 0,\     Y^n_0 = 0,  
     \eeq
$n\ge 1$, with the same $a_n=n\pi$ and independent standard Brownian motions $(B^n_t)$ as in \eqref{ch1(*1a)}, where $a \in \R$ and $\sigma > 0$. Suppose that for each $t \geq 0$, the $Y^n_t$ are the Fourier coefficients of a function $x \mapsto u(t, x)$, $x \in [0, 1]$, so that
\beq
\label{ch1(*2a)}
    u(t, x) = \sqrt{2}\, \sum_{n=1}^\infty\, Y^n_t\, \sin(n \pi x).   
    \eeq
   Does the random field $(u(t, x))$ satisfy a SPDE?
   To answer this question, we can take  derivatives in \eqref{ch1(*2a)} (informally, or in the sense of Schwartz distributions). We begin with the first derivative in time, obtaining
   \begin{align}
   \label{ch1(*3)}
     \frac{\partial}{\partial t} u(t, x) &= \sqrt{2}\, \sum_{n=1}^\infty\, \frac{d Y^n_t}{dt}\, \sin(n \pi x)\notag\\
   &=\sqrt{2}\, \sum_{n=1}^\infty  \left(- (a_n^2 + a^2)\, Y^n_t + \sigma\, \frac{d B^n_t}{dt}\right)\, \sin(n \pi x)\notag\\
    &= \sqrt{2}\, \sum_{n=1}^\infty  - (a_n^2 + a^2)\, Y^n_t\, \sin(n \pi x)+ \sigma\, \frac{\partial}{\partial t} Z(t, x),
    \end{align}
where
\beq
\label{giveanumber(*1a)}
    Z(t, x) = \sqrt{2}\, \sum_{n=1}^\infty B^n_t\, \sin(n \pi x).
    \eeq
 Next, we differentiate twice with respect to $x$ in \eqref{ch1(*2a)} to see that
 \beqn
 \frac{\partial^2}{\partial x^2} u(t,x) = -\sqrt 2\, \sum_{n=1}^\infty\, a_n^2\, Y_t^n\, \sin(n\pi x).
 \eeqn
 Thus,
 \begin{align}
 \label{ch1-(*3)-bis}
  \sqrt 2 \, \sum_{n=1}^\infty - (a_n^2 + a^2)\, Y^n_t \,\sin(n \pi x) &
  = \frac{\partial^2}{\partial x^2} u(t, x) - \sqrt{2}\, a^2\, \sum_{n=1}^\infty\, Y^n_t\, \sin(n \pi x)\notag\\
 & =  \frac{\partial^2}{\partial x^2} u(t, x) - a^2\, u(t, x).
\end{align}
In order to interpret the equality \eqref{ch1-(*3)-bis}, let $\cA$ be  the operator defined on $L^2([0, 1])$ as follows: for $f(x) = \sqrt{2} \, \sum_{n=1}^\infty f_n \sin(n \pi x)$,
\beqn
   \mathcal{A}(f)(x) =  - \sqrt{2} \sum_{n=1}^\infty (a_n^2 + a^2)\, f_n\, \sin(n \pi x).
\eeqn
Then \eqref{ch1-(*3)-bis} tells us that for fixed $t>0$,
\beqn
   \mathcal{A}(u(t,*))(x) = \left(\frac{\partial^2}{\partial x^2} - a^2\right) u(t, x).
\eeqn
Let
\begin{align*}
    W(t, x) :&= \int_0^x  dy\, Z(t, y) = \sqrt{2}\, \sum_{n=1}^\infty\, B^n_t \int_0^x  dy\, \sin(n \pi y)\\
    &= \sqrt{2}\, \sum_{n=1}^\infty B^n_t\, \frac{1 - \cos(n \pi x)}{n \pi}.
    \end{align*}
Then $W(t,x)$ is the Brownian sheet from \eqref{ch1(*1a)} and, by definition of $W$,
\beq
\label{ch1(*5)}
    \frac{\partial}{\partial x} W(t, x) = Z(t, x).     
    \eeq
From \eqref{ch1(*3)}--\eqref{ch1(*5)}, we see that
\beq
\label{ch1(*7)}
     \frac{\partial}{\partial t} u(t, x) = \left(\frac{\partial^2}{\partial x^2} - a^2\right) u(t, x) + \sigma\, \frac{\partial^2}{\partial t\partial x} W(t, x).    
     \eeq
     Noting that $u(0, x) = 0$ and $u(t, 0) = u(t, 1) = 0$, we see that $(u(t, x))$ is the solution of the stochastic heat equation \eqref{ch1(*7)}, with vanishing initial and Dirichlet boundary conditions, driven by $\sigma\, \frac{\partial^2}{\partial t \partial x} W(t, x)$, which is a multiple of space-time white noise according to Proposition \ref{rd1.2.19} (g). We will encounter this random field again in Chapter \ref{ch6}.
     \smallskip

Formulas for the Ornstein-Uhlenbeck processes\index{Ornstein-Uhlenbeck process}\index{process!Ornstein-Uhlenbeck} $(Y^n_t)$ are well-known \cite[Chapter 5, §6, Example 6.8]{ks}: if the $(Y^n_0)$ are random variables that are independent of the $(B^n_s)$, then
\beqn
    Y^n_t = Y^n_0\, e^{-(a^2 + a_n^2) t} + \sigma \int_0^t e^{- (a^2 + a_n^2)(t - s)}\, dB^n_s,
    \eeqn
so, when $Y_0^n=0$, we should have
\beq
\label{giveanumber(*1)}
    u(t, x) = \sqrt 2\, \sigma\, \sum_{n=1}^\infty\, \sin(n \pi x) \int_0^t e^{- (a^2 + a_n^2)(t - s)}\, dB^n_s,
    \eeq
a formula that will be confirmed in Section \ref{ch4-ss2.2-1'}.

Notice that if $a = 0$, $\sigma = \sqrt{2}$ and if the law of $Y_0^n$ is ${\rm{N}}\left(0, \frac{1}{2a_n^2}\right)$, then the law of $Y^n_t$ is also ${\rm{N}}\left(0, \frac{1}{2a_n^2}\right)$, so for each $t \geq 0$, $x \mapsto u(t, x)$ is a Brownian bridge as for the series in \eqref{rdch1(*1)}.

\begin{remark}
\label{ch1-r- *1}
(a)\  If we seek an expression for $u(t, x)$ in \eqref{giveanumber(*1)} that relates directly to the space-time white noise $\frac{\partial^2}{\partial t \partial x} W(t, x)$, we can proceed informally as follows, taking $a = 0$ and $\sigma = 1$ for simplicity. From \eqref{giveanumber(*1a)}, we have
\beqn
      B^n_t = \langle Z(t, *), e_n \rangle = \int_0^1 e_n(y)\, Z(t, y)\, dy,
      \eeqn
where $e_n(y) = \sqrt{2} \sin(n \pi y)$, so
\beqn
    dB^n_s = ds \int_0^1 e_n(y)\, \frac{\partial}{\partial s} Z(s, y)\, dy.
    \eeqn
From \eqref{giveanumber(*1)}, we see that
\begin{align*}
    u(t, x)&= \sum_{n=1}^\infty\, e_n(x) \int_0^t ds\, e^{-a_n^2 (t-s)} \int_0^1 e_n(y)\, \frac{\partial}{\partial s} Z(s, y)\, dy\\
    &= \int_{[0, t] \times [0, 1]} \left(\sum_{n=1}^\infty\, e^{-a_n^2 (t-s)}\, e_n(x)\, e_n(y)\right) \frac{\partial}{\partial s} Z(s, y)\, ds dy.
\end{align*}
By \eqref{ch1(*5)}, we see that
\beqn
   u(t, x) =  \int_{[0, t] \times [0, 1]} G_1(t-s; x, y)\, \frac{\partial^2}{\partial s \partial y} W(s, y)\, ds dy,
   \eeqn
 with $G_1(t-s; x, y)$ given in \eqref{ch1.600}.
Interpreting the integral as in \eqref{rde1.2.12}, this is formula \eqref{giveanumber(*2)} with $u_0 \equiv 0$ and $L = 1$.
\smallskip

(b)\ As in Section \ref{ch1-1.1}, it is traditional to write SDE's as
\beqn
    dY_t = - a\, Y_t\, dt + \sigma\, dB_t
\eeqn
rather than
\beqn
   \frac{dY_t}{dt} = - a\, Y_t + \sigma\, \frac{dB_t}{dt}.
   \eeqn
As mentioned at the beginning of Section \ref{ch1-1.1}, the term $\frac{dB_t}{dt}$ would be a white noise in time. On the other hand, in an SPDE, one usually puts $``\dot W (t,x)"$ on the right-hand side. As we can see in \eqref{ch1(*7)}, when there are partial derivatives with respect to two or more variables, the notation $\dot W$ is convenient.
\end{remark}

\section{Examples of SPDEs}
\label{ch1-1.2}

We describe here some situations from physics, biology and economics which are naturally modelled using SPDEs. The rigorous background will be given in later chapters.
\medskip

\noindent{\em The vibrating string}\index{vibrating string}\index{string!vibrating}
\smallskip

Imagine a guitar that has been left outdoors during a sandstorm. The impacts of the grains of sand will cause the strings to vibrate.
What tune will the guitar play?

   Let $u(t,x)$ be the vertical displacement of that string at position $x\in \IR$ and at time $t\ge 0$. The motion of the string will be described by the stochastic wave equation:
$$
   \rho\,\frac{\partial^{2} u(t,x)}{\partial t^{2}} = \frac{\partial^{2} u(t,x)}{\partial x^{2}} + \dot {F}(t, x),
$$
where $\rho$ is the density of the string per unit length, the second derivative with respect to $t$ represents the acceleration, the second derivative with respect to $x$ is the contribution of elastic forces, and the term $\dot F(t,x)$ represents the forces due to the random impacts of the grains of sand. Assuming that the numbers of impacts in disjoint regions of space-time are independent of each other, a reasonable model for $\dot F(t,x)$ will be space-time white noise. This example was introduced by J.B. Walsh in \cite{walsh}.
\medskip

\noindent{\em String in a random environment}\index{string!in a random environment}
\smallskip

 In \cite{funaki83}, Funaki proposed a model for the motion of an elastic string in a viscous random environment. It is obtained as the limit in law of discrete approximations. Let $W^1, \ldots, W^N$ be a collection of independent $d$-dimensional Brownian motions (or Wiener processes). Consider a system of $N$ particles that move under the influence of three kinds of forces: elastic forces acting between particles of intensity proportional to their distance, an external force $f$ and a random force.
The movement of the $k$-th particle is described by
\beq
\label{ch1-funaki}
dy_k = \left[\frac{\kappa}{2} N^2 (y_{k+1}+y_{k-1} - 2y_k) + f(y_k)\right] dt + \sqrt{N} b(y_k)\, dW_t^k,
\eeq
with given $y_k(0)$, $k=1,\ldots,N$. Here $\kappa$ is the modulus of the elastic forces and $b(y)$, $y\in \IR$, is a function describing the intensities of the random forces. This is inspired by the equation governing the motion of a particle moving in a viscous environment in $\IR^d$ under a forcing field $f(y)$, $y\in \IR^d$, which is
\beqn
\frac{d}{dt} y(t) = f(y(t)), \qquad y(0) \in \IR^d.
\eeqn

Assume that $f$ and $b$ are Lipschitz continuous and that we are given $y_0(t)$ and $y_{N+1}(t)$ for any $t\ge 0$. Then the process $(y_k(t),\, t\ge 0,\, k=1, \ldots N)$ is uniquely determined by \eqref{ch1-funaki}.

Let $x_k^N =( k-1)/(N-1)$, $k=1,\ldots,N$, and fix a continuous function $\varphi : [0,1] \rightarrow \IR^d$. We prescribe the initial conditions
\beqn
y_k(0) = \varphi(x_k^N), \quad  k=1,\ldots,N,
\eeqn
and the boundary conditions
\beqn
y_0(t) = y_1(t), \quad y_{N+1}(t) = y_N(t), \quad  t\ge 0.
\eeqn

We introduce a stochastic process with values in $\mathcal{C}([0,1];\IR^d)$ that will turn out to be a discrete approximation of the moving string for $t\in [0,T]$:
\beqn
X_N(t,x) = y_k(t) + \frac{x- x_k^N}{x_{k+1}^N-x_k^N}\, y_{k+1}(t), \qquad  x\in [x_k^N, x_{k+1}^N] ,\ t\in[0,T],
\eeqn
$k=1,\ldots, N-1$.

The following result is proved in \cite{funaki83}. The sequence of random vectors $(X_N, \, N \in \IN)$, converges weakly in $\mathcal{C}([0,1];\IR^d)$ to a process $(X(t),\, t\in [0,T])$ solution to
\beq
\label{ch1-funaki-2}
dX_t = \left[\frac{d^2}{d x^2}X_t + f(X_t)\right] dt + b(X_t)\, dW_t,
\eeq
with domain of $A= \frac{d^2}{d x^2}$ given by
\beqn
\mathcal{D}(A)= \left\{z\in W^{2,2}([0,1]; \IR^d): \ \frac{d^2 z}{d x^2}\in \mathcal{C}([0,1];\IR^d), \ \frac{dz}{d x}(0)= \frac{dz}{d x}(1)=0 \right\},
\eeqn
where $W^{2,2}([0,1]; \IR^d)$ is the $(2,2)$-Sobolev space of $\IR^d$-valued functions defined on $[0,1]$ (see e.g.~\cite{adams}).
In \eqref{ch1-funaki-2}, $W$ is a {\em cylindrical Wiener process} on $L^2([0,1]; \IR)$, a notion that is defined for instance in \cite{dz}.
\medskip

\noindent{\em Motion of a strand of DNA}
\smallskip

This example is quoted from \cite{Khoshnevisan09Mini}. A DNA molecule can be viewed as a long elastic string, whose length is  infinitely long compared to its diameter. We can describe the position of the string by using a parameterization defined on $\IR_+ \times [0,1]$ with values in $\IR^3$:
$$
   \vec{u}(t,x) = \left(\begin{array}{c}
            u_1(t,x)\\
            u_2(t,x)\\
            u_3(t,x)
            \end{array}\right).
$$
Here, $\vec{u}(t,x)$ is the position at time $t$ of the point labelled $x$ on the string, where $x \in [0,1]$ represents the distance from this point to one extremity of the string if the string were straightened out. The unit of length is chosen so that the entire string has length $1$.

   A DNA molecule typically ``floats" in a fluid, so it is constantly in motion, just as a particle of pollen floating in a fluid moves according to Brownian motion. The motion of the particle can be described by Newton's law of motion, which equates the sum of forces acting on the string with the product of the mass and the acceleration. Let $\mu$ be the mass of the string per unit length. The acceleration at position $x$ along the string, at time $t$, is
$$
   \frac{\partial^2\vec{u}}{\partial{t}^2} (t,x),
$$
and the forces acting on the string are mainly of three kinds: elastic forces $\vec{F_1}$, which include torsion forces, friction due to viscosity of the fluid $\vec{F_2}$, and random impulses $\vec{F_3}$, due to the impacts on the string of the fluid's molecules. Newton's equation of motion can therefore be written
$$
   	\mu\frac{\partial^2\vec{u}}{\partial{t}^2}
   = 
   \vec{F_1}
   - \vec{F_2}
   + \vec{F_3}.
$$

   This is a rather complicated system of three stochastic partial differential equations, and it is not even clear how to write down the torsion forces or the friction term. Elastic forces are generally related to the second derivative in the spatial variable, and the molecular forces are reasonably modelled by a stochastic noise term.

   Assume $\mu = 1$. The simplest $1$-dimensional equation related to this problem, in which one only considers vertical displacement and forgets about torsion, is the following one, in which $u(t,x)$ is now scalar valued:
\begin{equation}\label{DNAs}
   \frac{\partial^2u}{\partial t^2} (t,x)=
   \frac{\partial^2u}{\partial x^2}(t,x)
   -\int^1_0 \,k(x,y)\,u(t,y)\,dy + \dot{F}(t,x).
\end{equation}
Here the first term on the right hand side represents the elastic forces, the second term is a (non-local) friction term, and the third term $\dot{F}(t,y)$ is a Gaussian noise, with spatial correlation $k(\cdot,\cdot)$, that is,
$$
   E(\dot{F}(t,x)\,\dot F(s,y)) = \delta_0(t-s)\, k(x,y),
$$
where $\delta_0$ denotes the Dirac delta function; a rigorous definition of this equation can be found in \cite{zab} (see also \cite{dalang}, \cite{pz2}).
The function $k(\cdot,\cdot)$ is the same in the friction term and in the correlation.

The motion of a DNA strand is of biological interest, since when it moves
 around and two normally distant parts of the string get close enough together, it can happen that a biological event occurs, for instance, an enzyme may be released. Therefore, some biological events are related to the motion of the DNA string. This particular question could be translated as follows.

   Fix $0<x<y<1$ and $\epsilon > 0$. Estimate $P\{ \Vert \vec{u}(t,x) - \vec{u} (t,y) \Vert < \epsilon \}$.
A mathematical idealization of this question could be: is $P\{\exists t > 0: u(t,x) = u (t,y) \}$ positive? That is, do distant points on the string come together at some positive time? An even simpler question, that is already highly non-trivial from a mathematical point of view would be: given $u_{0} \in \IR^{d}$, is $P\{\exists (t,x): \vec{u} (t,x) = u_{0} \}$ positive?

    Some mathematical results for equation (\ref{DNAs}) can be found in \cite{zab}. Some of the biological motivation can be found in \cite{maddocks}.
\medskip

\noindent{\em Multi-dimensional waves}
\smallskip

In the previous example, $x\in\re$. There are also interesting cases where $x$ is multi-dimensional.
We present now two examples that are studied in \cite{leveque}.
 \medskip

\noindent{\em (1) 2-$d$ surface waves}
\smallskip

Consider raindrops falling on the surface of a lake, each raindrop generating a (small) surface wave. Let $u(t,x)$ be the vertical displacement at time $t$ of position $x \in D \subset \IR^{2}$. The evolution of $u$ is given by the following wave equation:
\beqn
   \frac{\partial^{2}u(t,x)}{\partial t^{2}} = \Delta u(t,x) - \frac{\partial u(t,x)}{\partial t} + b(t, x, u(t,x)) + a(t, x, u)\, \dot{F}(t,x).
\eeqn
The term $\dot F(t,x)$ models the impacts of the raindrops, after compensating for the average effect, which may be contained in the term $b$.
\medskip

\noindent{\em (2) Pressure waves}
\smallskip

 A different situation occurs if we consider another effect of the raindrops falling on the surface of the lake. Each raindrop generates a sound, or pressure wave, that moves down into the depth of the water. Let $u(t,x)$ denote the pressure at time $t$ and position $x \in D \subset \IR^{3}$. The evolution of $u$ is given by the following wave equation:
$$
   \frac{\partial^{2} u(t,x)}{\partial t^{2}} = \Delta u (t, x_{1}, x_{2}, x_{3}) + a (t, x, u)\, \dot{F}(t,x_{1}, x_{2})\, \delta_{0} (x_{3})
$$
where the noise $\dot F(t,x_1,x_2)\, \delta_{0}(x_3)$ is concentrated on a lower-dimensional surface (the surface of the sea).
\medskip

\noindent{\em The internal structure of the sun}
\smallskip

This example is quoted from \cite{Khoshnevisan09Mini}. The study of the internal structure of the sun is an active area of research. One important international project was known as Project SOHO \cite{SOHO}. Its objective was to use measurements of the motion of the sun's surface to obtain information about the internal structure of the sun. Indeed, the sun's surface moves in a rather complex manner: at any given time, any point on the surface is typically moving towards or away from the center. There are also waves going around the surface, as well as shock waves propagating through the sun itself, which cause the surface to pulsate.

A question of interest to solar geophysicists is to determine the origin of these shock waves. One school of thought is that they are due to turbulence, but the location and intensities of the shocks are unknown; thus, a probabilistic model may be appropriate.

P. Stark (U.C.~Berkeley) proposed a model that assumes that the sun is a ball of radius $R$, and that the main
source of shocks is located in a spherical zone inside the sun. Assuming that the shocks are randomly located on this sphere, the equation for the pressure variations throughout the sun would be
\begin{equation}\label{sun}
   \frac{\partial^2u}{\partial t^2}(t,x) =
   c^2(x)\,\rho_0(x)\left(\vec\nabla\cdot\left(\frac{1}{\rho_0(x)}\,\vec\nabla
   u\right)-\vec\nabla \cdot \vec F(t,x)\right),
\end{equation}
where $x \in B(0,R)$, the ball centered at the origin with radius $R$, $c^2(x)$ is the speed of wave propagation at position $x$, $\rho_0(x)$ is the density at position $x$ and $\vec F(t,x)$ models the shock that originates at time $t$ and position $x$. The notation ``$\cdot$'' here refers to the Euclidean inner product.

A model for $\vec F$ that corresponds to the description of the situation would be 3-dimensional Gaussian noise concentrated on the sphere $\partial B(0,r)$, where $0 < r < R$. A possible choice of the spatial correlation for the components of $\vec F$ would be
$\delta(t-s) \, f(x\cdot y)$. A problem of interest is to estimate $r$ from the available observations of the sun's surface. Some mathematical results relevant to this problem are developed in \cite{dalangleveque}.
\medskip

\noindent{\em Neural response}
\smallskip

 In \cite{walsh}, Walsh considers synapses, that send impulses of current into a neuron. The axon of a neuron can be identified with a long thin cable, say $[0,L]$. Let $u(t,x)$ be the electrical potential (called {\em action potential}) at position $x$ and time $t$. The evolution of this potential is described by the SPDE
\beq
\label{walsh-neuro}
   \frac{\partial u(t,x)}{\partial t} = \frac{\partial^{2}u(t,x)}{\partial x^{2}} - u(t,x) + g(u(t,x))\, \dot F(t,x),
\eeq
$(t,x)\in\,]0,\infty[\, \times\, ]0,L[$,  with Neuman boundary conditions, where $\dot F(t,x)$ is a Gaussian noise that models the random electrical impulses.

In the case that $g\equiv 0$, the SPDE \eqref{walsh-neuro} is the cable equation. It is an approximation of the Hodgkin and Huxley model for the description of action potentials in neurons (see \cite{hodgkin-huxley}).
\medskip

\noindent{\em Parabolic Anderson model}\index{Anderson model}\index{model!Anderson}
\smallskip

The continuous parabolic Anderson model\index{parabolic Anderson model}\index{Anderson model!parabolic} (see \cite{carmona-molchanov-1994}) is described by the SPDE
\beq
\label{parabolic-anderson}
   \frac{\partial u(t,x)}{\partial t} = \frac{\nu}{2}\, \Delta u(t,x)+ \rho\, u(t,x)\, \dot W(t,x),
   \eeq
$(t,x)\in\,]0,\infty[ \times \re^k$, where $\nu>0$, $\rho\in\re\setminus\{0\}$ and $\dot W(t,x)$ is a space-time white noise (see Definition \ref{rd1.2.18}). The initial data $u(0,\ast)$ is assumed to satisfy certain conditions (in some physically relevant cases, $u(0,\ast)$ may be a nonnegative measure rather than a function). This SPDE arises in connection with the physical phenomenon of {\em strong localization}\index{strong!localization}\index{localisation!strong} in a random potential, stated by P. W. Anderson in \cite{anderson1958} in the following terms: if the disorder is strong enough, then localization of states will occur no matter the dimension of the system. This is also a statement on absence of diffusion of waves (or localization of electrons) in a random medium.
In probabilistic terms,
the continuous parabolic Anderson model can be seen as a limit of a $k$-dimensional random walk  in a random environment.

 Related to the phenomenon of localisation, the notion of {\em intermittency}\index{intermittency}\index{intermittent} expresses the property that a random field takes values close to zero in vast regions of space-time, but develops high peaks on some small ``islands." This was translated into mathematical terms by Zeldovich, Molchanov and several coauthors (\cite{molchanov1994}, \cite{shandarin1989}, \cite{zeldovich1987}, \cite{zeldovich21987}), who formulated this property in terms of Lyapounov exponents (see \eqref{u-Lyapounov}). For the SPDE \eqref{parabolic-anderson} with $k=1$, this property was intensively studied in \cite{carmona-molchanov-1994} and many subsequent papers (see, for instance \cite{carmona2001}, \cite{cranston2007}, \cite{cranston2005}, \cite{gartner2000}, \cite{gartner2007}, \cite{tindel2002}), and later in \cite{conus2012}, \cite{conus-joseph-khoshnevisan-2013}, \cite{chen}, \cite{chendalang2015-2}, among others). For a review, see \cite{konig2016}.

In the mathematical study of \eqref{parabolic-anderson}, a crucial point that remained open for many years is the construction of random field solutions in the physically relevant dimensions $k=2, 3$. The difficulty lies in defining the product $u(t,x) \dot W(t,x)$ in these dimensions. This problem was addressed non-rigorously until the development of the theory of  {\em regularity structures} by M. Hairer in \cite{H-2013}, \cite{H-2014} (see also \cite{g-i-p2015}).


To a certain extent, a similar study was undertaken for the hyperbolic SPDE\index{hyperbolic Anderson model}\index{Anderson model!hyperbolic}
\beqn
\label{hyperbolic-anderson}
   \frac{\partial^2 u}{\partial t^2}(t,x) = \Delta u(t,x) + \rho\, u(t,x))\, \dot W(t,x),
   \eeqn
   with initial conditions $u(0,\ast)=f(\ast)$ and $\frac{\partial}{\partial t}u(0,\ast)= g(\ast)$,
beginning with \cite{dalang-mueller2009}, and continued in particular by R. Balan and coauthors (see e.g. \cite{balan-conus2016}). Some results on the Anderson model are given in Section \ref{rd1+1anderson}.

\medskip

\noindent{\em Population dynamics}
\smallskip

Two popular and important examples of SPDEs originating in
 population dynamics are
 \begin{align}
\frac{\partial u(t,x)}{\partial t} & = \frac{\partial^{2}u(t,x)}{\partial x^{2}} +\sqrt{u(t,x)}\, \dot F(t,x), \label{d-w}\\
\frac{\partial u(t,x)}{\partial t} & = \frac{\partial^{2}u(t,x)}{\partial x^{2}} +\sqrt{u(t,x)(1-u(t,x))}\, \dot F(t,x),\label{f-v}
\end{align}
where $t>0$, $x\in\IR$, and the initial conditions are given.

These are equations satisfied by the densities of some classes of measure-valued processes. Indeed, consider a sequence  of stochastic processes $(X_n, n\ge 1)$ indexed by a time parameter $t\ge0$. For each $n\ge 1$, we define
\begin{equation*}
\mu_n(t) = \sum_{i,j=1}^n q_{i,j}(n,t)\, \delta_{X_j(t)},
\end{equation*}
where $\delta_{x}$ denotes the Dirac measure. The factors $q_{i,j}(n,t)$ may be random; they represent  dynamical interactions or relationships between the particles. The simplest example corresponds to the choice
$q_{i,j}(n,t)= \delta_i^j/n$, with $\delta_i^j$ denoting the Kronecker symbol.\index{Kronecker symbol}\index{symbol!Kronecker}\label{rd05_04i1} In this case,
$\mu_n(t)$ is the empirical measure of the particle system. A fundamental and difficult question is the existence of the limit of $(\mu_n, n\ge 1)$ after an appropriate rescaling. If such a limit exists, then it will define a measure-valued process $(\mu(t),\, t\ge 0)$.  In specific cases, $\mu(t)$ has a density with respect to Lebesgue measure, and this density satisfies a SPDE. Two famous examples are the {\em Dawson-Watanabe} and the {\em Fleming-Viot} processes, whose densities $u(t,x)$ are described by \eqref{d-w} and \eqref{f-v}, respectively (see e.g. \cite{dawson72}, \cite{dawson75}, and \cite{watanabe-1968}, \cite{fleming75}).

Similar to the above example {\em ``String in a random environment''}, this is yet another illustration of the connections between particle systems and SPDEs.
\medskip

\noindent{\em Interest rate models}
\smallskip

In \cite{ct06}, infinite dimensional stochastic analysis is used in the study of interest rate models. In stochastic models for the term structure of interest rates, a central role is played by the Heath-Jarrow-Morton framework (HJM).  Let $P(t,T)$ denote the bond price at time $0\le t \le T$ with maturity date $T>0$ in a market satisfying certain assumptions. Musiela's instantaneous forward rate is defined as
\beqn
f(t,T) = -\frac{\partial \log P(t,T)}{\partial T}.
\eeqn
The HJM evolution equation is the SPDE for $u(t,x):=f(t,t+x)$,\ $t>0$, $x\ge 0$, given by
\beqn
\frac{\partial u(t,x)}{\partial t} = \frac{\partial u(t,x)}{\partial x} + a(t,x) + \sum_{k\ge 1} \sigma_k(t,x)\, \dot W_k(t),
\eeqn
where $(W_k)$, $k\ge 1$ is a sequence of independent standard Brownian motions and $a(t,x)$ is defined by the {\em HJM no-arbitrage condition}.

This equation, and several other abstract SPDEs with more general multiplicative noises, are studied in \cite[Chapter 6]{ct06}
(see also \cite{mt06}).
\medskip

\noindent{\em Nonlinear filtering}
\smallskip

Consider a stochastic evolution system denoted by $(Z_t = (X_t,Y_t),\, t\in [0,T])$, described by the following stochastic differential equations:
\begin{align}
\label{filter}
dX_t & = h(Z_t)\, dt + f(Z_t)\, dW_t + g(Z_t)\, dV_t, \qquad X_0=x_0,\notag\\
dY_t & = B(Z_t)\, dt + dV_t, \qquad\qquad\qquad\qquad\quad\quad Y_0=y_0,
\end{align}
where $W$ and $V$ are independent Brownian motions. The $\rek$-valued process $X$ denotes the {\em state process}, and $Y$ is the {\em observation process}. As the terminology suggests, only the process $Y$ is observed and  $X$ is not observed. The problem of nonlinear filtering\index{nonlinear!filtering}\index{filtering!nonlinear} consists in estimating the random vector $X_t$ from the observed values $Y_{[0,t]}:=(Y_s,\, 0\le s\le t)$.

 The search for an {\em optimal} estimate, in the sense that it minimizes the expected quadratic error, can be formulated as follows: for any real-valued measurable and bounded function $\varphi$, find a random variable $\tilde\Phi_t$ such that
\beqn
E\left[(\tilde\Phi_t - \varphi(X_t))^2\right] = \inf_F E\left[(F(Y_{[0,t]})- \varphi(X_t))^2\right],
\eeqn
where the infimum  is over all measurable functionals $F$ on the space of continuous functions defined on $[0,t]$ such that $F(Y_{[0,t]})$ has a finite second moment. It is a well-known fact that the solution to this problem is given by the conditional expectation of the random variable $\varphi(X_t)$ given the $\sigma$-field $\mathcal{G}_ t = \sigma(Y_s,\, 0\le s\le t)$:
\beqn
\tilde\Phi_t = E\left[\varphi(X_t)|\mathcal{G}_ t\right] = \int_{\IR^k} \varphi(x)\, P_t(dx),
\eeqn
where $P_t(dx)$ denotes the conditional probability law.

Assume that $P_t(dx)$ is absolutely continuous with respect to Lebesgue measure, and denote by $\pi_t(x)$ its density. For theoretical and practical motivations, it is relevant to have a specific description of $P_t$ and $\pi_t$.

Under suitable regularity conditions on the coefficients of the system \eqref{filter}, $\pi_t$ and a normalised version of this density satisfy certain SPDEs.

For example, in the particular case $g\equiv 0$, $h(Z_t)= h(X_t)$, $f(Z_t)= f(X_t)$, $B(Z_t)= B(X_t)$, the density $\pi_t$ (if it exists), satisfies
\beq
\label{ch1.filter}
\frac{\partial \pi}{\partial t} = \mathcal{L}(\pi) + B(\pi)\, dY_t,
\eeq
 where
 \beqn
 \mathcal{L}\psi = \sum_{i,j=1}^d\, \frac{\partial^2((\sigma\sigma^T)_{i,j} \psi)}{\partial x_i \partial x_j} - \sum_{i=1}^d \, \frac{\partial (h_i\psi)}{\partial {x_i}}.
 \eeqn
 This is a particular case of {\em Zakai's equation}\index{Zakai's equation}\index{equation!Zakai's} (see \cite{zakai}).
 We refer to \cite{gyongy2016} for the details of the proof of this fact, and point to \cite{pardoux79}, \cite{pardoux89},
 \cite{cr2011} for some references on this topic.
 \medskip

\noindent{\em Other examples}
\smallskip

We have already seen several examples of SPDEs directly related with famous equations in physics, like the heat, wave and Poisson equations. These share the property of being {\em linear}, meaning here that they are of the form $\mathcal{L}u = f$, where $\mathcal{L}$ is a {\em linear} partial differential operator.
There are important  examples that do not share this property, such as the SPDE \eqref{ch1-bu} below.
In the context of SPDEs with nonlinear partial differential operators, the study of  fundamental questions, such as well-posedness, requires in most the cases the use of forcing noises that are smoother than white noise.

However, there are also some particular examples where the choice of space-time white noise (or other types of white noises) is possible, for example, for the {\em stochastic Burgers' equation}.\index{Burger's equation}\index{equation!Burger's}
The deterministic version of this equation was introduced in \cite{bu48}, in connection with the study of turbulence. The stochastic counterpart is the equation
\beq
\label{ch1-bu}
\frac{\partial u(t,x)}{\partial t} = \frac{\partial^2 u(t,x)}{\partial x^2} - u\frac{\partial u(t,x)}{\partial x} + \sigma(u(t,x)) \dot W(t,x),
\eeq
on $[0,1]$ or on $\IR$, where $\dot W(t,x)$ is a space-time white noise, with suitable initial and boundary conditions.
The first results on  \eqref{ch1-bu} are in the case of $\sigma$ constant (see \cite{bcj-n94}). Existence and uniqueness of an $L^p([0,1])$-valued solution, for a suitable value of $p$, is  proved in \cite{dp95}. A more general setting has been considered in \cite{gy98} where, in particular, these results have been extended and cover examples of {\em stochastic reaction-diffusion equations} (see also \cite{a81}, \cite{f88}, \cite{c2002}). Reaction-diffusion equations describe chemical reactions and are fundamental
in the field of physical chemistry.

With Hairer's theory of {\em regularity structures}\index{regularity structures}\index{structures!regularity} (\cite{H-2014}),  SPDEs driven by space-time white noise and with  singular coefficients can be considered, the paradigmatic example being the  Kardar-Parisi-Zhang (KPZ) equation (\cite{H-2013}). This type of SPDE is beyond the scope of this book.
\bigskip

 \section{Notes on Chapter \ref{ch1}}
 \label{ch1-notes}
The basic theory of Brownian motion and stochastic differential equations is covered in many classical books such as \cite{ks}, \cite{oksendal-1998}, \cite{ry}. Section \ref{ch1-1.1} emphasizes the distributional aspect of the time-derivative of Brownian motion.

 The fundamental notions on Gaussian processes in Section \ref{ch1-2} can be found in classical books such as \cite{grimmett-2001}, \cite{K-1997}. More advanced accounts of the theory are presented in \cite{neveu-1968}, \cite{adler} and \cite{hida-1993}).

The notion of isonormal Gaussian process was introduced by Segal \cite{segal-1954}; this  is a fundamental notion in the theory of abstract Wiener spaces (\cite{gross-1965}) and Wiener calculus (\cite{wiener-1938}, \cite{ito-1951}), and in most of the approaches to Malliavin calculus (\cite{watanabe-1984}, \cite{watanabe-1987}, \cite{nualart-1995}, \cite{nourdin-peccati-2012}).

In signal processing, white noise refers to signals with equal intensity at all frequencies. In more rigorous terms, white noise is a Gaussian random field that can equivalently be viewed as an isonormal Gaussian process on $L^2(D, dx)$, where $D \subset \R^k$ and $dx$ denotes Lebesgue measure. The theory of distribution-valued processes (see for instance \cite{gv}, \cite{ito84}) provides the mathematical framework for the study of such random fields, and
space-time white noise is a fundamental example of such a distribution-valued process that plays an important role in the theory of parabolic and hyperbolic SPDEs. Its main properties are gathered in Section \ref{stwn-s1}.

Random series representations (termed Fourier-Wiener series in \cite[Chapter 16, § 3]{Kahane-1985}) of Brownian motion were already used by Paul Lévy \cite{levy-1965} in the construction of Brownian motion and the study of its sample paths. Wiener's construction of Brownian motion relies on the Fourier-sine series \eqref{rdch1(*1)} (see \cite{wiener-1923}), but the Fourier-cosine series and the Fourier series of Brownian motion are also useful. In Section \ref{stwn-s2}, we use these to give Fourier series representations for white noise on an interval, and for various extensions of Brownian motion, such as the Brownian sheet (in $[0,1]^2$ and in $\re_+\times [0,1]$) and we deduce Fourier series representations for the corresponding white noises.

The literature on partial differential equations is quite vast. However, the background material that is directly relevant to the topics developed in this book  is comparatively standard and can be found, for instance, in \cite{evans}, \cite{folland}, \cite{john}.

Partial differential equations with an external forcing noise can be found in the physics literature already in the 1950's and 1960's. In the mathematical literature, examples of SPDEs driven by space-time white noise and a systematic study of their properties began in the 1960's and 1970's; see for instance \cite{dalets-1967}, \cite{zakai}, \cite{cabanya-1970}, \cite{dawson72}, \cite{pardoux79}. A brief account of the historical development of the subject can be found in \cite{zambotti-2021}.
The stochastic heat equation discussed at the end of Section \ref{ch1-s3} is one of the most studied examples in the theory of SPDEs. For early results on this equation, see \cite{walsh}.

The selection of SPDEs given in Section \ref{ch1-1.2}  aims to highlight the diversity of fields in which SPDEs can be useful: as mathematical models of physical phenomena, in the understanding of biological mechanisms, in cosmology and in the analysis of financial markets, to mention only a few. Additional examples are presented in the monographs \cite{dz} and \cite{lototsky-rozovsky-2017}.


\chapter[Stochastic integrals with respect to space-time white\texorpdfstring{\\}{}noise]{Stochastic integrals with respect to space-time white noise}

\label{ch2}

\pagestyle{myheadings}
\markboth{R.C.~Dalang and M.~Sanz-Sol\'e}{Stochastic integral with respect to space-time white noise}

In this chapter, we develop a theory of stochastic integration that is suitable for integrating random functions of space and time with respect to space-time white noise $W=(W(A),\, A\in \B^f_{\IR \times D})$. Recall that space-time white noise is defined in Section \ref{stwn} and its properties are summarized in Proposition \ref{rd1.2.19}.
We associate to $W$ a sequence of independent standard Brownian motions, leading to a natural definition of the stochastic integral with respect to $W$ as a series of one-dimensional It\^o stochastic integrals with respect to these Brownian motions. This definition coincides with Walsh's integral as developed in \cite{walsh}. We present the main properties of the stochastic integral and some of its extensions, and we give fundamental results of the infinite-dimensional stochastic analysis toolbox: a stochastic Fubini's theorem, a theorem on differentiation under the stochastic integral, and a theorem on joint measurability for stochastic integrals depending on a parameter.

\section{Preliminaries}
\label{ch2new-s1}

Let $W=(W(A),\, A\in \B^f_{\IR_+ \times D})$ be a space-time white noise (Definition \ref{rd1.2.18}) defined on a complete probability space $(\Omega,\F,P)$, where $D\subset\rek$ is a non-empty bounded or unbounded connected open set. Fix $T >0$ and set $H= L^2([0,T] \times D)$.  From the isonormal process on $H$ associated to $W$ as in Proposition \ref{rd1.2.19}~(d), we can define a Gaussian random field $(W_s(\varphi),\, s\in[0,T],\, \varphi \in L^2(D))$, as follows:
\beqn
   W_s(\varphi) = W(1_{[0,s]}(\cdot)\, \varphi(\ast)),
\eeqn
where ``$\cdot$'' and ``$\ast$'' refer to the temporal and spatial variables, respectively.
\begin{lemma}
\label{ch1'-lc}
The process $(W_s(\varphi))$ has the following properties:
\begin{enumerate}
\item for any $\varphi\in L^2(D)$, $(W_s(\varphi),\, s\in[0,T])$ defines a Brownian motion with variance $s\Vert\varphi\Vert^2_{L^2(D)}$,
\item for all $s,t\in[0,T]$, and $\varphi, \psi\in L^2(D)$,
\beqn
E[W_s(\varphi)\, W_t(\psi)] = (s\wedge t)\, \langle\varphi,\psi\rangle_{L^2(D)},
\eeqn
with $s\wedge t= \min(s,t)$.
\end{enumerate}
\end{lemma}
\begin{proof}
Let $\nu(dt, dx) = 1_{[0,T]}(t)\, 1_D(x)\, dt dx$. Since $(W(h),\, h\in L^2([0,T]\times D, \nu))$ is an isonormal process, we see from Definition \ref{ch1-d3} that
\begin{align*}
E\left[W_s(\varphi)\right] & = E\left[W(1_{[0,s]}(\cdot)\, \varphi(\ast))\right] =0,
\end{align*}
and from the isometry property of the isonormal process, for any $\varphi, \psi\in L^2(D)$,
\begin{align*}
E\left[W_s(\varphi)\,W_t(\psi)\right]& = E\left[W(1_{[0,s]}(\cdot)\, \varphi(\ast))\, W(1_{[0,t]}(\cdot)\, \psi(\ast))\right] \\
& = \langle 1_{[0,s]}(\cdot)\, \varphi(\ast), 1_{[0,t]}(\cdot)\, \psi(\ast)\rangle_{H}\\
& = (s\wedge t)\, \langle \varphi,\psi\rangle_{L^2(D)}.
\end{align*}
When $s=t$ and $\varphi=\psi$, this yields
\beqn
E\left[W_s(\varphi)^2\right] = s\, \Vert\varphi\Vert^2_{L^2(D)}.
\eeqn
This completes the proof.
\end{proof}

The process $(W_s(\varphi) = W(1_{[0,s]}(\cdot)\, \varphi(\ast)),\, s\ge 0,\, \varphi\in L^2(D))$ is an example of a {\em cylindrical Wiener process}\index{cylindrical!Wiener process}\index{Wiener!process, cylindrical} \cite{metivier-pellaumail}.

In the sequel, we will work with the continuous versions of the Brownian motions $(W_s(\varphi),\, s\in[0,T])$.


Consider a right-continuous complete filtration $\left({\mathcal{F}}_s,\, s\in[0,T]\right)$ consisting of
sub-$\sigma$--fields of $\F$ satisfying the following conditions:
\begin{description}
\item (i) for all fixed $s\in[0,T]$ and $\varphi\in L^2(D)$, the random variable $W_s(\varphi)$ is $\mathcal{F}_s$--measurable;
\item (ii) for any $s\in[0,T]$, the family $(W_t(\varphi)-W_s(\varphi),\, \varphi\in L^2(D),\, t\in[s,T])$ is independent of  $\mathcal{F}_s$.
\end{description}
An example of filtration satisfying the above conditions is the completed natural filtration associated with $W$, that is, for $s\ge 0$, ${\mathcal{F}}_s$ is the $\sigma$-field generated by the random variables $(W_t(\varphi),\, 0\le t\le s,\, \varphi\in L^2(D))$ and the $P$--null sets.

Throughout this chapter, we will denote by $V$ the Hilbert space $L^2(D)$ endowed with a complete orthonormal basis $(e_j, \, j\ge 1)$.

\begin{lemma}
\label{ch1'-lsi}
\begin{itemize}
\item[(1)] The sequence $(W_s(e_j),\, s\in[0,T],\, j\ge 1)$ consists of independent standard Brownian motions. These are adapted to $(\tf_s)$, and for $s\geq 0$, $(W_{t+s}(e_j) - W_s(e_j),\, t \in [0,T-s],\, j \geq 1)$ is independent of $\tf_s$. Further, for all $\varphi \in V$ and $s \in [0,T]$,
\beq
\label{2.1.2(*1)}
       W_s(\varphi) = \sum_{j=1}^\infty\, \langle \varphi, e_j \rangle_V\,  W_s(e_j),
       \eeq
where the series converges a.s.~and in $L^2(\Omega)$.
\item[(2)] Conversely, given a sequence $(B_s^j,\, s\in[0,T], \, j\ge 1)$ of independent Brownian motions, the series
\beq
\label{ch1'-s4.id}
     \sum_{j=1}^\infty\, B_s^j\, \langle \varphi, e_j\rangle_V =: \tilde W_s(\varphi), \quad  s\in[0,T], \  \varphi\in V,
\eeq
converges a.s.~and in $L^2(\Omega)$, and $(\tilde W_s(\varphi))$ defines a random field satisfying the conclusions of Lemma \ref{ch1'-lc}.

For $h \in H = L^2([0, T] \times D)$, define
\beq
\label{2.1(*J)}
    \tilde W(h) := \sum_{j=1}^\infty \int_0^T \langle h(s, *), e_j\rangle_V\, d B_s^j,
    \eeq
where the series converges a.s.~and in $L^2(\Omega)$. Then $\tilde W =(\tilde W(1_A),\, A \in \B_{[0, T] \times D}^f)$ is a space-time white noise on $[0,T] \times D$ and $(\tilde W(h),\, h \in H)$ is the associated isonormal Gaussian process.
\item[(3)] Let $ W$  be the space-time white noise considered at the beginning of this section. If in part (2), we take $ B_s^j = W_s(e_j)$, then the space-time white noise $\tilde W$  coincides with $W$.
\end{itemize}
\end{lemma}

\begin{proof}
(1) Since $(e_j,\, j\ge 1)$ is an orthonormal basis of $V$, it follows from Lemma \ref{ch1'-lc} that the $(W_s(e_j),\, s \in [0,T])$, $j\geq 1$, are independent standard Brownian motions. The next two properties in (1) follow directly from the conditions (i) and (ii) on $(\tf_s)$.

Finally, notice that for $\varphi \in V$, $\varphi(*) = \sum_{j=1}^\infty \langle \varphi, e_j \rangle_V\,  e_j(*)$, where the series converges in $V$. Using the linear isometry property of an isonormal Gaussian process given in Lemma \ref{ch1-l2}, we have
\beqn
     W_s(\varphi) = W(1_{[0,s]} (\cdot)\, \varphi(*)) = \sum_{j=1}^\infty\, \langle \varphi, e_j \rangle_V\, W(1_{[0,s]} (\cdot)\, e_j(*)),
\eeqn
where the series converges in $L^2(\Omega)$, and this is the right-hand side of \eqref{2.1.2(*1)}. By the Khintchine-Kolmogorov convergence theorem (see e.g.~\cite[Theorem 1, p.~110]{chow-teicher-1978}),
 the series in  \eqref{2.1.2(*1)} also converges a.s.

(2) For fixed $s\in[0,T]$ and $\varphi\in V$, $(X_j:= B_s^j\ \langle\varphi,e_j\rangle_V,\, j\ge 1)$ is a sequence of independent random variables in $L^2(\Omega)$ with $E[X_j] =0$, $j\ge 1$, and
\beqn
\sum_{j=1}^\infty E[X_j^2] = s\, \sum_{j=1}^\infty\, \langle \varphi,e_j\rangle^2_V = s\, \Vert\varphi\Vert_V^2 <\infty.
\eeqn
Hence by the Khintchine-Kolmogorov convergence theorem, the series in \eqref{ch1'-s4.id} converges a.s.~(and in $L^2(\Omega)$, as we already knew).

The convergence in $L^2(\Omega) $ of the series in \eqref{2.1(*J)} is a special case of the calculation in \eqref{revised-1} below, so we do not give details here. The a.s.~convergence is again a consequence of the Khintchine-Kolmogorv convergence theorem. For $h, g \in H$, it is clear that $E(\tilde W(h)) = 0$ and
\begin{align}
\label{2.1(*K)}
 E[\tilde W(h)\, \tilde W(g)] &= E\Bigg[\left(\sum_{j=1}^\infty \int_0^T \langle h(s, *), e_j\rangle_V\, d B_s^j\right)\notag\\
 &\qquad\qquad\times \left(\sum_{k=1}^\infty \int_0^T \langle g(s, *), e_k\rangle_V\, d B_s^j\right)\Bigg]\notag\\
   &= \sum_{j, k = 1}^\infty E\left[\int_0^T ds\, \langle h(s, *), e_j\rangle_V\, \langle g(s, *), e_k\rangle_V\, d\langle B_\cdot^j , B_\cdot^k  \rangle_s\right]\notag\\
   &= \int_0^T ds\, \sum_{j=1}^\infty\, \langle h(s, *), e_j\rangle_V\, \langle g(s, *), e_j\rangle_V\notag\\
   &= \int_0^T \langle h(s, *), g(s, *) \rangle_V\, ds\notag\\
   &= \langle h, g \rangle_H.
   \end{align}
In the third equality, we have used the independence of $B^j$ and $B^k$ for $j \neq k$, and in the last equality, we have applied Parseval's identity. Therefore, $(W(h),\, h \in H)$ is an isonormal Gaussian process on $H$.

   It follows immediately from \eqref{2.1(*K)} that $\tilde W := (\tilde W(1_A),\, A \in \cB_{[0, T] \times D}^f)$ is a space-time white noise on $[0,T] \times D$. Since $h \mapsto \tilde W(h)$ is linear and by definition,
   \beq
   \label{2.1(*L)}
    \tilde W(1_A) = \sum_{j=1}^\infty \int_0^T \langle 1_A(s, *), e_j\rangle_V\, d B_s^j .
    \eeq
If we apply the construction of Section \ref{ch1-2.2}, we see that \eqref{2.1(*L)} extends to functions $h$ which are linear combinations of indicators of disjoints sets in $\cB_{[0, T] \times D}^f$. Using the isometry property \eqref{iso-simple} and the just established isometry property \eqref{2.1(*K)}, we conclude that $(\tilde W(h),\, h \in H)$ is the isonormal Gaussian process associated to $\tilde W$.

 (3) With this choice of the $(B_s^j)$, for $A \in B_{[0, T] \times D}^f$ of the form $A = [0,t] \times F$, where $F \subset D$, by definition of $\tilde W$,
 \begin{align*}
 \tilde W(A) :&= \tilde W(1_A)\\
 & = \sum_{j=1}^\infty \int_0^t \langle 1_F, e_j\rangle_V\, dW_s(e_j)\\
   &=  \sum_{j=1}^\infty  \langle 1_F, e_j\rangle_V\, W_t(e_j)
    = W_t\left(\sum_{j=1}^\infty\,  \langle 1_F, e_j \rangle_V\, e_j\right)\\
    &= W_t(1_F) = W(1_{[0,t]}(\cdot)\, 1_F(*)) = W(A),
    \end{align*}
so $\tilde W$ and $W$ coincide. This completes the proof.
\end{proof}

\section{The stochastic integral}
\label{ch2new-s2}

Let $W=(W(A),\, A\in \B^f_{\IR_+ \times D})$ be a space-time white noise and define $(W_s(\varphi),\, s\in[0,T],\, \varphi \in L^2(D))$ as in Section \ref{ch2new-s1}. Let $(\tf_s,\, s\in [0,T])$ be a filtration satisfying (i) and (ii) of Section \ref{ch2new-s1}.  Fix a complete orthonormal basis $(e_j,\, j\ge 1)$ of $V= L^2(D)$, and let $\left(W_s(e_j),\, s\in[0,T],\, j\ge 1\right)$ be the sequence of independent standard Brownian motions given in Lemma \ref{ch1'-lsi}.

We want to integrate {\em real-valued, jointly measurable, adapted and square-integrable} stochastic processes $G=(G(s,y),\, (s,y)\in [0,T] \times D)$  with respect to the space-time white noise $W$. Given $G$, for $(s,\omega) \in [0,T]\times \Omega$, we denote $G(s,\ast; \omega)$ the partial function $y \mapsto G(s,y; \omega)$ from $D$ into $\IR$. The precise assumptions on $G$ are:
\begin{itemize}
   \item[(1)] $(s,y; \omega) \mapsto G(s,y; \omega)$ from $[0,T]\times D\times \IR$ into $\IR$ is $\cB_{[0,T]}\times \cB_D \times \cF$-measurable;
   \item[(2)] for $s \in [0,T]$, $(y,\omega) \mapsto G(s,y;\omega)$ from $D \times \Omega$ into $\IR$ is $\cB_D \times \cF_s$-measurable;
   \item[(3)] $E\left[ \int_0^T ds\int_D dy\, G^2(s,y)  \right] < \infty$.
\end{itemize}
From (3), we see that for $dsdP$-almost all $(s,\omega)$, the map $G(s,\ast; \omega)$ belongs to the Hilbert space $V=L^2(D)$, so
\beqn
G(s,\ast,\omega) = \sum_{j=1}^\infty\, \langle G(s,\ast,\omega),e_j\rangle_V\, e_j(\ast),
\eeqn
where the series converges in $V$.

   We now define the stochastic integral of $G$ with respect to $W$.

\begin{def1}\label{ch1'-s4.d1}
Let $G=(G(s,y),\, (s,y)\in [0,T] \times D)$ be a jointly measurable, adapted and square-integrable random field, that is, assumptions {\rm (1)--(3)} above are satisfied. The stochastic integral\index{stochastic!integral}\index{integral!stochastic} of $G$ with respect to the space-time white noise $W$ is the random variable
\beq
 \label{ch1'-s4.1}
 \int_0^T \int_D G(s,y)\, W(ds,dy):= \sum_{j=1}^\infty \int_0^T \langle G(s,\ast),e_j\rangle_V\, dW_s(e_j),
 \eeq
 where the series converges in $L^2(\Omega)$.
\end{def1}

In this definition, the integrals on the right-hand side are It\^o integrals with respect to the one-dimensional independent continuous Brownian motions $(W_s(e_j),\, s\in[0,T])$, $j\ge 1$. According to Lemma \ref{ch1'-lsi}~(1), these are well-defined by \cite[Theorem 3.8]{chung-williams}, since $(s,\omega)\mapsto \langle G(s,\ast),e_j\rangle_V$ is $\cB_{[0,T]}\times \cF$-measurable by (1), and for fixed $s \in [0,T]$, $\omega \mapsto \langle G(s,\ast),e_j\rangle_V$ is $\cF_s$-measurable by (2). Further,
$$
   E\left[\int_0^T \langle G(s,\ast),e_j\rangle_V^2\, ds \right] \leq E\left[\int_0^T \Vert G(s,\ast) \Vert_V^2\, ds \right] < \infty
$$
by (3). Moreover, the terms in the series \eqref{ch1'-s4.1} are orthogonal in
$L^2(\Omega)$, because for $j\ne k$, $W_s(e_j)$ and $W_s(e_k)$ are independent Brownian motions.

The convergence in $L^2(\Omega)$ of the series in \eqref{ch1'-s4.1} is ensured by the assumptions on $G$. Indeed, by the isometry property of the It\^o integral with respect to Brownian motion and Fubini's theorem,
 \begin{align}
 \label{revised-1}
  \sum_{j=1}^\infty\, \left\Vert \int_0^T\langle G(s,*),e_j\rangle_V\, dW_s(e_j) \right\Vert_{L^2(\Omega)}^2
  &= \sum_{j=1}^\infty\, E\left[\int_0^T \langle G(s,\ast),e_j\rangle_V^2\, ds\right]\notag\\
 &= E\left[\int_0^T\left(\sum_{j=1}^\infty\, \langle G(s,\ast),e_j\rangle_V^2\right) ds\right]\notag\\
  &= E\left[\int_0^T \Vert G(s,\ast)\Vert^2_V\, ds\right] < \infty,
\end{align}
where we have used Parseval's identity for the third equality.

The stochastic integral in \eqref{ch1'-s4.1} will also be denoted by $(G\cdot W)_T$.
We note that the intuition behind the formula \eqref{ch1'-s4.1} comes from the classical Parseval's identity
\beqn
\int_0^T ds\int_D dy\ f_1(s,y)\, f_2(s,y) = \sum_{j=1}^\infty \int_0^T ds\ \langle f_1(s,\ast), e_j\rangle_V\, \langle f_2(s,\ast), e_j\rangle_V,
\eeqn
valid for functions $f_1, f_2\in L^2([0,T]\times D)$.
\smallskip

For convenience, if $0\leq r \leq t$ and $A \subset [0,T] \times D$ is a Borel set, we will sometimes write
  $ \int_r^t \int_D G(s,y)\, W(ds,dy)$ and  $\int_A G(s,y)\, W(ds,dy)$
instead of
    $\int_0^t \int_D 1_{[r,t]}(s)\, G(s,y)\, W(ds,dy)$ and  $\int_0^T\int_D 1_A(s,y)\, G(s,y)\, W(ds,dy)$,
respectively, which are both well-defined when G satisfies the conditions of Definition \ref{ch1'-s4.d1}.

\begin{prop}
\label{ch1'-s4-p1}
The stochastic integral satisfies the isometry property\index{isometry!property}\index{property!isometry}
\beq
\label{ch1'-s4.2}
E\left[\left(\int_0^T \int_D G(s,y)\, W(ds,dy)\right)^2\right] = E\left[\int_0^T \Vert G(s,\ast)\Vert^2_V \, ds\right].
\eeq
\end{prop}

\begin{proof}
Because the elements of the family $(W_t(e_j),\, 0\le t\le T)$, $j\ge 1$, are mutually independent and centered,
\begin{align*}
   E\left[\left(\int_0^T G(s,y)\, W(ds,dy)\right)^2\right] &= E\left[\left( \sum_{j=1}^\infty \int_0^T \langle G(s,\ast),e_j\rangle_V\, dW_s(e_j)\right)^2\right] \\
  & = \sum_{j=1}^\infty E\left[\left(\int_0^T \langle G(s,\ast),e_j\rangle_V \, dW_s(e_j)\right)^2\right].
\end{align*}
Using \eqref{revised-1} we obtain the result.
\end{proof}

\begin{remark}
\label{ch2(*A)}
 Let $(W(h),\, h \in L^2([0,T] \times D))$ be the isonormal Gaussian process associated with $W$. Notice that if $G$ is a deterministic function $g(\cdot,*) \in L^2([0, T] \times D)$, then according to Definition \ref{ch1'-s4.d1},
 \beqn
      \int_0^T \int_D  g(s,y)\, W(ds,dy) = \sum_{j = 1}^\infty \int_0^T  \langle g(s, *), e_j \rangle_V\, dW_s(e_j) = W(g),
      \eeqn
by Lemma \ref{ch1'-lsi} (2) and (3) applied to the sequence $(W_s(e_j),\,j \geq 1)$.
In particular, Definition \ref{ch1'-s4.d1} is compatible with the construction of the isonormal process $(W(h),\, h \in L^2([0,T] \times D)$ and notations such as \eqref{wiener0}.
\end{remark}

\begin{lemma}
\label{ch1'-lib}
The definition of the stochastic integral in \eqref{ch1'-s4.1} does not depend on the particular orthonormal basis in $V$.
\end{lemma}

\begin{proof}
Consider an orthonormal basis $(v_j,\, j\ge 1)$ in $V$ and write
\beqn
G(s,\ast) = \sum_{k=1}^\infty\, \langle G(s,\ast),v_k\rangle_V\, v_k.
\eeqn
Since $\varphi\mapsto W_s(\varphi)$ is linear, if we assume that the series and integrals can be permuted, then we would have
\begin{align}
\label{ind.0}
&\sum_{j=1}^\infty \int_0^T \langle G(s,\ast), e_j\rangle_V\, dW_s(e_j)\notag\\
&\qquad\qquad = \sum_{j=1}^\infty \int_0^T \Big \langle  \sum_{k=1}^\infty\, \langle G(s,\ast),v_k\rangle_V\, v_k, e_j\Big\rangle_V\, dW_s(e_j)\notag\\
& \qquad\qquad= \int_0^T  \sum_{k=1}^\infty\, \langle G(s,\ast), v_k\rangle_V \left[\sum_{j=1}^\infty\, \langle v_k,e_j\rangle_V\, dW_s(e_j)\right]\notag\\
& \qquad\qquad= \sum_{k=1}^\infty \int_0^T \langle G(s,\ast), v_k\rangle_V \, dW_s(v_k),
 \end{align}
 and this would prove the lemma.

We now check \eqref{ind.0} with some care. Since $v_k = \sum_{j=1}^\infty \langle v_k, e_j\rangle_V\, e_j$, we have
\beq
\label{ind.1}
W_s(v_k) = \sum_{j=1}^\infty\, \langle v_k, e_j\rangle_V\, W_s(e_j),
\eeq
where the series converges a.s. and in $L^2(\Omega)$.

We claim that for all jointly measurable adapted real-valued processes $g\in L^2([0,T]\times \Omega)$ and for all $k\ge 1$,
\beq
\label{ind.2}
\int_0^T g_s\, dW_s(v_k) = \sum_{j=1}^\infty \int_0^T g_s\, \langle v_k, e_j\rangle_V\, dW_s(e_j),
\eeq
where the series converges in $L^2(\Omega)$. Indeed, if $g_s = X\, 1_{]s_0,t_0]}(s)$ with $X$ bounded and $\tf_{s_0}$-measurable, then \eqref{ind.2} is an immediate consequence of \eqref{ind.1}. For general jointly measurable adapted $g=(g_s,\, s\in[0,T]) \in L^2( [0,T]\times \Omega)$, there is a sequence of simple processes $(g_s^n,\, s\in[0,T])$ such that
\beq
\label{ind.3}
\lim_{n\to\infty} E\left[\int_0^T (g_s - g_s^n)^2\, ds \right] = 0
\eeq
(see  \cite[p.~35]{chung-williams} for the notion of simple process).
For $(g_s^n)$, \eqref{ind.2} holds by linearity, and
\beqn
E\left[\left( \int_0^T (g_s - g_s^n)\, dW_s(v_k)\right)^2\right] = E\left[\int_0^T (g_s - g_s^n)^2\, ds \right] \rightarrow 0.
\eeqn
Using the independence of the $(W_s(e_j))$, $j\ge 1$, we have
\begin{align*}
&E\left[\left(\sum_{j=1}^\infty \int_0^T (g_s - g_s^n)\, \langle v_k,e_j\rangle_V\, dW_s(e_j)\right)^2\right]\\
& \qquad\qquad= E\left[\sum_{j=1}^\infty \left(\int_0^T (g_s - g_s^n)^2\, \langle v_k, e_j\rangle_V^2\, ds \right)\right]\\
& \qquad\qquad= \Vert v_k\Vert^2_V\,  E\left[\int_0^T (g_s - g_s^n)^2\, ds\right]
\rightarrow 0,
\end{align*}
by \eqref{ind.3}. We conclude that \eqref{ind.2} holds for $(g_s)$.

We now observe that for fixed $M, N\ge 1$,
\begin{align*}
&\sum_{j=1}^M \int_0^T \left\langle \sum_{k=1}^N\, \langle G(s,\ast), v_k\rangle_V\, v_k, e_j\right\rangle_V\, dW_s(e_j)\\
&\qquad = \sum_{k=1}^N \sum_{j=1}^M  \int_0^T \langle G(s,\ast), v_k\rangle_V\, \langle v_k, e_j\rangle_V\, dW_s(e_j).
\end{align*}
For fixed $N$, by \eqref{ind.2} applied to $g_s:=\langle G(s,\ast),v_k\rangle_V$, the right-hand side converges in $L^2(\Omega)$ as $M\to\infty$ to
\beq
\label{ind.4}
 \sum_{k=1}^N \int_0^T  \langle G(s,\ast), v_k\rangle_V\, dW_s(v_k),
 \eeq
 while the left-hand side converges in $L^2(\Omega)$ to
 \begin{align}
 \label{ind.5}
 &\sum_{j=1}^\infty \int_0^T  \left\langle\sum_{k=1}^N \langle G(s,\ast), v_k\rangle_V\, v_k, e_j\right\rangle_V\, dW_s(e_j)\notag\\
 &\qquad = \int_0^T \int_D \left(\sum_{k=1}^N \langle G(s,\ast), v_k\rangle_V\, v_k(y)\right)\, W(ds,dy).
 \end{align}
 Since
 \beq
 \label{ind.5bis}
 \lim_{N\to\infty} E\left[\int_0^T \left\Vert G(s,\ast) - \sum_{k=1}^N \langle G(s,\ast), v_k\rangle_V\ v_k\right\Vert^2_V\, ds\right] =0,
 \eeq
 we let $N\to\infty$ in \eqref{ind.4} and \eqref{ind.5} to conclude from the isometry property \eqref{ch1'-s4.2} and the equality of \eqref{ind.4} and \eqref{ind.5} that
 \beqn
 \int_0^T \int_D G(s,y)\, W(ds,dy) = \sum_{k=1}^\infty \int_0^T \langle G(s,\ast), v_k\rangle_V\, dW_s(v_k).
 \eeqn
 This completes the proof of \eqref{ind.0}.
 \end{proof}

We note that the stochastic integral with respect to the space-time white noise defined above is the stochastic integral with respect to the standard {\em cylindrical Wiener process} $(W_t(\varphi) = W(1_{[0,t]}(\cdot)\, \varphi(\ast)),\, t\ge 0,\, \varphi\in L^2(D))$  (see \cite[Section 4.3.2]{dz} and \cite{metivier-pellaumail}).
\medskip

\noindent{\em Indefinite integral}
\smallskip

Let $G=(G(s,y),\, (s,y)\in [0,T] \times D)$ be as in Definition \ref{ch1'-s4.d1}.
Let ${\bH}^2$ be the vector space of continuous $(\cF_t)$-martingales on $[0,T]$ which vanish at time $0$, in which indistinguishable processes are identified. Recall that the space ${\bH}^2$\label{rdH2} with the inner product $\langle M, N \rangle := E[M_T N_T]  = E[\langle M, N \rangle_T]$ (where the bracket $\langle\cdot, \cdot \rangle_T$ denotes the quadratic covariation process at time $T$) is a Hilbert space (see for instance  \cite[Sec. 5.1]{legall-2013}).
Consider the sequence $Z^n$ of elements of ${\bH}^2$ defined by
\beq
\label{rd2.2.12}
Z^n=\left(Z_t^n=\sum_{j=1}^n \int_0^t \langle G(s,\ast),e_j\rangle_V\, dW_s(e_j),\  t\in[0,T]\right), \quad n\ge 1.
\eeq
We notice that by the independence of the $(W_t(e_j))$, for each $n\ge 1$, the quadratic variation process of $(Z^n)$
is
\beqn
\langle Z^n\rangle=\left(\langle Z^n\rangle_t=\sum_{j=1}^n \int_0^t \langle G(s,\ast),e_j\rangle_V^2\, ds, \  t\in[0,T]\right).
\eeqn
Since the sequence $(Z^n_T)$ converges in $L^2(\Omega)$ to $(G \cdot W)_T$, we deduce that the sequence $(Z^n)$ converges in ${\bH}^2$. We denote its limit
\beq
\label{(*2)}
    G \cdot W = ((G \cdot W)_t,\,  t \in [0,T])
    \eeq
and call this the {\em indefinite integral process of G with respect to W}.\index{indefinite integral process}\index{integral!process, indefinite} Since it belongs to ${\bH}^2$, it is a continuous $L^2(\Omega)$-bounded $(\F_t)$-martingale on $[0,T]$ which vanishes at time $0$. We note that as $n \to \infty$, $Z^n_t$ converges to $(G \cdot W)_t$ in $L^2(\Omega)$, uniformly in $t \in [0,T]$, and there is a subsequence $(n_k)$ such that a.s., $Z^{n_k}_t$ converges to $(G \cdot W)_t$ uniformly in $t \in [0,T]$ (see \cite[proof of Proposition 5.1]{legall-2013}).
\begin{lemma}
\label{lemma 5}
For each $t \in [0,T]$,
\beqn
    (G \cdot W)_t = \int_0^T \int_D G(s,y)\, 1_{[0, t]}(s)\, W(ds,dy)\quad     {\text a.s}.
\eeqn
\end{lemma}
\begin{proof}
Using the definition of the stochastic integral on the right-hand side, we see that it is the $L^2(\Omega)$-limit of $Z^n_t$.
\end{proof}

Instead of $G \cdot W$, we will often write
\beqn
\left(\int_0^t \int_D G(s,y)\, W(ds,dy),\ t \in [0,T]\right),
\eeqn
 with the understanding that this is indistinguishable from the continuous martingale $G \cdot W$.

\begin{prop}
\label{Proposition 6}
The indefinite integral process in \eqref{(*2)} is an $L^2(\Omega)$-bounded continuous martingale with quadratic variation  process
   \beq
\label{ch1'-s4.3}
\left(\int_0^t \Vert G(s,\ast)\Vert^2_V\, ds, \ t\in [0,T]\right).
\eeq
As a consequence, there is a Burkholder's inequality.\index{Burkholder's inequality}\index{inequality!Burkholder's} More precisely, for any $p>0$, there is a constant $C_p$, depending only on $p$, such that for any stopping time $\tau$,
\beq
\label{ch1'-s4.4}
E\left[\sup_{r\in [0,\tau\wedge T]} \left\vert (G\cdot W)_r\right\vert^p\right] \le C_p\, E\left[\left(\int_0^{\tau\wedge T} \Vert G(s,\ast)\Vert_V^2\, ds\right)^{\frac{p}{2}}\right].
\eeq
\end{prop}
\begin{proof}
We only need to prove the statement concerning quadratic variation.
Let $Z= (Z_t) := ((G \cdot W)_t,\,  t \in [0,T])$ and let $Z^n$ be defined as in \eqref{rd2.2.12}.
Applying the Cauchy-Schwarz inequality and the isometry property of the stochastic integral, we have
\beqn
E\left[\left\vert (Z_t^n)^2 - Z_t^2\right\vert\right] \le \sqrt 2 \left(E\left[\int_0^t \Vert G(s,\ast)\Vert_V^2\, ds\right]\right)^{\frac{1}{2}}
\Vert Z_t^n-Z_t\Vert_{L^2(\Omega)}\ \to 0,
\eeqn
as $n\to \infty$.

The stochastic process
\beqn
\left( \langle Z\rangle_t = \sum_{j=1}^\infty \int_0^t \langle G(s,\ast),e_j\rangle_V^2\, ds= \int_0^t \Vert G(s,\ast)\Vert^2_V\, ds, \quad  t\in[0,T]\right)
\eeqn
is adapted, continuous, increasing, $\langle Z\rangle_0=0$ a.s., and satisfies
\beqn
E\left[\left\vert \langle Z_t^n\rangle_t - \langle Z\rangle_t\right\vert\right]
= E\left[\sum_{j=n+1}^\infty \int_0^t  \langle G(s,\ast),e_j\rangle_V^2\, ds\right] \to 0,
\eeqn
as $n\to \infty$.

Since the $L^1(\Omega)$-limit of a sequence of continuous martingales with respect to some filtration is a continuous martingale (with respect to the same filtration), we deduce that the stochastic process $(Z_t^2-\langle Z\rangle_t,\, t\in[0,T])$ is a  martingale with respect to the filtration $(\tf_ t,\,  t\in [0,T])$. This proves that \eqref{ch1'-s4.3} is the quadratic variation of the indefinite integral process (see e.g. \cite{ks}).

The Burkholder-Davis-Gundy inequality can be found in \cite[Theorem 3.2.8, p. 166]{ks} or in \cite[(4.2) Corollary, p. 161]{ry}.
\end{proof}

\medskip

\begin{remark}\label{ch2-r.burkholder}

The inequality \eqref{ch1'-s4.4} clearly implies
\beq
\label{ch1'-s4.4.mod}
\sup_{r\in [0,t]}E\left[ \left\vert\int_0^r \int_D G(s,y)\, W(ds,dy)\right\vert^p\right] \le \tilde C_p E\left[\left(\int_0^t \Vert G(s,\ast)\Vert_V^2\, ds\right)^{\frac{p}{2}}\right].
\eeq
In this form, the optimal constant for $p=2$ is $\tilde C_2 = 1$ and for $p\geq 1$, one has $\tilde C_p \leq (4p)^{\frac{p}{2}}$. Indeed, this follows from the version of the Burkholder-Davis-Gundy inequality for bounded continuous martingales proved in  \cite[Theorems 1 and  A, pages 354 and 365, respectively]{C-K-1991}. An extension to continuous $L^2$-martingales is proved in \cite[Theorem B.1, p.~97]{khosh}.
\end{remark}
\medskip

\noindent{\em Examples of integrands}
\smallskip

In Chapter \ref{ch1'-s5}, we will frequently encounter integrands of the form
\beq\label{rd2.2.12a}
G(s,y) = \Gamma(t,x;s,y)\, Z(s,y),
\eeq
where $t\in\, ]0,T]$ and $x\in D\subset \rek$ are fixed, $D$ is a bounded or unbounded domain in $\rek$, $0\le s<t\le T$ and $y\in D$. The function $\Gamma$ is usually 	the fundamental solution or the Green's function corresponding to some partial differential operator, and $Z=(Z(s,y),\, (s,y)\in[0,T]\times D)$ is a jointly measurable and adapted random field satisfying
\beq\label{rd2.2.12b}
\sup_{(s,y)\in[0,T]\times D} E \left[ Z^2(s,y)\right] =: C <\infty.
\eeq
 The assumptions on $\Gamma$ (see ($\bf H_\Gamma$) in Section \ref{ch4-section0}) are such that $\Gamma$ is measurable and
\beq\label{rd2.2.12c}
\int_0^t ds\, \sup_{x\in D} \int_D dy\, \Gamma^2(t,x;s,y) <\infty.
\eeq
Then, for $(t,x) \in [0,T]\times D$,
\begin{align*}
&E\left[\int_0^t ds \int_D dy\ (\Gamma(t,x;s,y)\, Z(s,y))^2\right] \\
&\qquad \le\sup_{(s,y)\in[0,T]\times D} E \left[ Z^2(s,y)\right] \left(\int_0^t ds \int_D dy\, \Gamma^2(t,x;s,y)\right)\\
&\qquad = C\int_0^t ds \int_D dy\, \Gamma^2(t,x;s,y) < \infty.
\end{align*}
Thus,
the stochastic integral
\beqn
\int_0^t \int_D \Gamma(t,x;s,y)\, Z(s,y)\, W(ds,dy), \qquad t\in[0,T],
\eeqn
is well-defined according to Definition \ref{ch1'-s4.d1}.

Notice that for fixed $t\in[0,T]$, the process
\beqn
\left(\int_0^r \int_D\Gamma(t,x;s,y)\, Z(s,y)\ W(ds,dy), \  r\in[0,T]\right),
\eeqn
is a martingale. Hence, according to \eqref{ch1'-s4.4.mod},
\beq
\label{burkholder-spdes}
\sup_{r\in[0,t]} \left\Vert \int_0^r \int_D\Gamma(t,x;s,y)\, Z(s,y)\, W(ds,dy)\right\Vert_p^2 \le (\tilde{C}_p)^{\frac{2}{p}} \left\Vert \int_0^t \Vert G(s,\ast)\Vert_V^2\, ds \right\Vert_{\frac{p}{2}}.
\eeq
For $p\ge 2$, this version of Burkholder's inequality is extensively used in the theory of SPDEs.
\medskip

\noindent{\em Local property in $\Omega$ of the stochastic integral}
\smallskip

It is well-known that the It\^o stochastic integral has the {\em local property},\index{local property!in $\Omega$} meaning that on the subset of $\Omega$ where
the integrand vanishes, the stochastic integral also vanishes (see \cite[ Th\'eor\`eme 23, p. 346]{dm2}). Because of the definition \eqref{ch1'-s4.1}, this property directly carries over to the stochastic integral with respect to space-time white noise, as stated below.

\begin{lemma}
\label{ch1'-llp}
Let $G^{(1)}$, $G^{(2)}$, be two random fields satisfying the conditions of Definition \ref{ch1'-s4.d1}. Assume that on some $F\in\tf$, the sample paths of $G_1$ and $G_2$ are the same, that is, for almost all $\omega\in F$, $G^{(1)}(s,y,\omega)=G^{(2)}(s,y,\omega)$, for $dsdy$-almost all $(s,y)\in[0,T]\times D$. Then a.s., for all $t\in[0,T]$,
\beqn
1_F\int_0^t \int_D G^{(1)}(s,y)\, W(ds,dy) = 1_F\int_0^t \int_D G^{(2)}(s,y)\, W(ds,dy).
\eeqn
\end{lemma}
\medskip

\noindent{\em Local property in space of the stochastic integral}\index{local property!in space}
\smallskip

   Suppose that $D_1$ is a sub-domain of $D$ and we want to integrate with respect to $W$ a random field $G = (G(s,y),\, (s,y) \in [0,T] \times D_1)$, which we extend to $D$ by setting $G(s,y) = 0$ for $y \in D \setminus D_1$. Then we can use either of the two following procedures: (1) integrate this extension using formula \eqref{ch1'-s4.1}, or (2) use an orthonormal basis $(v_i)$, $i \geq 1$, of $V_1 = L^2(D_1)$, the restriction of $W$ to
  to $\re_+ \times D_1$ (that is, $(W(A),\, A \in \B^f_{\R_+ \times D_1})$) and the analogue
 of \eqref{ch1'-s4.1} for $D_1$. It turns out that both procedures give the same result, as the next proposition shows.
   \begin{prop}
 \label{prop-local-in-space}
Let $D_1 \subset D$ be a domain and let $(v_i)$, $i \geq 1$, be  an orthonormal basis of $V_1 = L^2(D_1)$. Let $G = (G(s,y),\, (s,y) \in [0,T] \times D_1)$. Suppose that the assumptions (1)--(3) at the beginning of the section are satisfied with $D$ there replaced by $D_1$. We extend $G$ to $[0,T] \times D$ by setting $G(s,y) = 0$ for $s \in [0,T]$ and $y \in D \setminus D_1$. Then assumptions (1)--(3) are satisfied with $D$ and
\beq
\label{(*1)}
    \sum_{i=1}^\infty \int_0^T \langle G(s,*), v_i \rangle_{V_1}\, dW_s(v_i)
    = \sum_{j=1}^\infty \int_0^T \langle G(s,*), e_j \rangle_V\, dW_s(e_j),
\eeq
where both series converge in $L^2(\Omega)$.
\end{prop}
\begin{proof} Assumptions (1) and (2) for $D$ follow from the fact that the extension of $G$ to $D$ is simply $(s,y,\omega) \mapsto G(s,y,\omega) 1_{D_1}(y)$.

    Let $\varphi \in V = L^2(D)$. Note that  $\langle G, \varphi \vert_{D_1} \rangle_{V_1} = \langle G, \varphi \rangle_V$, therefore assumption (3) for $D$ follows from assumption (3) for $D_1$.

    Note also that $W_s(v_i)$, obtained by restricting $W$ to Borel subsets of the open set $\re_+ \times D_1$ and using the procedure of Section \ref{ch1-2.2}, can equivalently be obtained by extending $v_i$ to $D$ by setting $v_i(y) = 0$ for $y \in D \setminus D_1$ and using the original white noise $W$.

  We now prove \eqref{(*1)}. The second moment of the difference of the two sides of \eqref{(*1)} is equal to the sum of three terms:
  \begin{align*}
A_1 &= E\left[\left(\sum_{i=1}^\infty \int_0^T \langle G(s,*), v_i \rangle_{V_1}\, dW_s(v_i)\right)^2\right],\\
     A_2 &=  E\left[\left(\sum_{j=1}^\infty \int_0^T \langle G(s,*), e_j \rangle_V\, dW_s(e_j)\right)^2\right],\\
     A_3 &= -2 E\left[\left(\sum_{i=1}^\infty \int_0^T \langle G(s,*), v_i \rangle_{V_1}\, dW_s(v_i)\right)\right.\\
    &\left.\qquad\qquad\times \left( \sum_{j=1}^\infty \int_0^T \langle G(s,*), e_j \rangle_V\, dW_s(e_j)\right)\right].
     \end{align*}
Notice that
\beqn
    A_1 = E\left[\int_0^T \Vert G(s,*) \Vert_{V_1}^2\, ds\right]
            = E\left[\int_0^T \Vert G(s,*) \Vert_{V}^2\, ds\right]
            = A_2,
  \eeqn
where the first and last equality are due to Proposition \ref{ch1'-s4-p1}. Further, since both series in $A_3$ converge in $L^2(\Omega)$, we can permute the sums and expectation in $A_3$ to obtain
\begin{align*}
   A_3 &= -2 \sum_{i=1}^\infty \sum_{j=1}^\infty E\left[\int_0^T \langle G(s,*), v_i \rangle_{V_1}\, dW_s(v_i) \int_0^T \langle G(s,*), e_j \rangle_V\, dW_s(e_j)\right]\\
        &= -2 \sum_{i=1}^\infty \sum_{j=1}^\infty E\left[\int_0^T \langle G(s,*), v_i \rangle_{V_1}  \langle G(s,*), e_j \rangle_V\,  d \langle W_\cdot(v_i), W_\cdot(e_j)\rangle_s\right].
\end{align*}
From assertion 2. of Lemma \ref{ch1'-lc}, we have $E\left(W_s(v_i) W_s(e_j)\right) = s\, \langle v_i,e_j\rangle_{V_1}$. Consequently, the quadratic covariation $\langle W_\cdot(v_i), W_\cdot(e_j)\rangle_s$ is equal to $s\, \langle v_i, e_j \rangle_{V_1}$. Therefore,
\begin{align*}
   A_3 &= -2\sum_{j=1}^\infty E\left[\int_0^T ds\ \langle G(s,*), e_j\rangle_V \left(\sum_{i=1}^\infty\langle G(s,\ast),v_i\rangle_{V_1}\, \langle v_i, e_j \rangle_{V_1}\right)\right] \\
&= -2 E\left[\int_0^T ds\, \Vert G(s,*) \Vert_{V_1}^2\right] = -2A_1.
\end{align*}
Indeed, by Parseval's identity, the sum over $i$ is equal to  $\langle G(s,\ast), e_j\rangle_{V_1}$, then we bring the series inside the $ds$-integral to see that the remaining sum over $j$ is equal to $\Vert G(s,*)\Vert_{V_1}^2$, because
$\langle G(s,\ast), e_j\rangle_{V}= \langle G(s,\ast), e_j\rangle_{V_1}$.

It follows that $A_1 +A_2 - 2 A_3 = 0$, proving \eqref{(*1)}.
\end{proof}

\noindent{\em Stochastic integral and stopping times}
\smallskip

Let $G=(G(s,y),\, (s,y)\in [0,T] \times D)$ be as in Definition \ref{ch1'-s4.d1}. Consider the continuous version of the indefinite integral process
\beqn
 \left(M_t = \int_0^t\int_D G(s,y)\, W(ds,dy),\quad  t \in [0,T]\right).
\eeqn

\begin{lemma}\label{rdlem2.2.6}
Let $\tau$ be a stopping time $($with respect to the filtration $(\tf_t))$ with values in $[0,T]$. Then
\beq\label{rde2.2.13}
   M_\tau = \int_0^T \int_D 1_{[0,\tau]}(s)\, G(s,y)\, W(ds,dy),\qquad \text{a.s.}
\eeq
\end{lemma}

\begin{proof} We recall that $M_\tau $ denotes the random variable defined by $(M_\tau)(\omega) = M_{\tau(\omega)}(\omega)$, $\omega\in\Omega$. For $n \in \IN$, let $\bD_n = \{k 2^{-n},\, k \in \IN \}$, and let
$
   \tau_n := \inf \{t\in \bD_n: t \geq \tau \} \wedge T.
$
Then $\bD_n \cap [0,T]$ is a finite set and $\tau_n \in (\bD_n \cap [0,T])\cup\{ T\}$ a.s. Further, $(\tau_n,\, n \in \IN)$ is a decreasing sequence of stopping times such that
\beq\label{rde2.2.14}
   \lim_{n\to \infty} \tau_n = \tau\qquad \text{a.s.,}
\eeq
and
\beq \label{rde2.2.15}
    \lim_{n\to \infty} M_{\tau_n} = M_\tau \qquad \text{a.s.}
\eeq
For $t \in (\bD_n \cap [0,T])\cup\{ T\}$ such that $P\{\tau_n = t \} >0$, a.s.~on $\{\tau_n = t \}$,
\begin{align}\nonumber
  M_{\tau_n} &= M_t = \int_0^t \int_D G(s,y)\, W(ds,dy) \\ \nonumber
   &= \sum_{j=1}^\infty \int_0^t \langle G(s,\ast),e_j\rangle_V\, dW_s(e_j) \\ \nonumber
   &= \sum_{j=1}^\infty \int_0^T 1_{[0,\tau_n]}(s)\, \langle G(s,\ast),e_j\rangle_V\, dW_s(e_j) \\
   &=  \int_0^T \int_D  1_{[0,\tau_n]}(s)\, G(s,y)\, W(ds,dy) \qquad\text{a.s.},
\label{rde2.2.16}
\end{align}
where in the next to last identity, we have used Lemma \ref{ch1'-llp}. It follows that
\eqref{rde2.2.16} holds a.s.~(on $\Omega$).
Since
$$
   \lim_{n\to \infty} E\left[\int_0^T ds \int_D dy\, \left(1_{[0,\tau]}(s)- 1_{[0,\tau_n]}(s)\right)^2\, G^2(s,y) \right] = 0
$$
by assumption (3) on $G$, \eqref{rde2.2.14} and dominated convergence, we let $n\to\infty$ in \eqref{rde2.2.16} to conclude from \eqref{rde2.2.15} that \eqref{rde2.2.13} holds.
\end{proof}

In view of Lemma \ref{rdlem2.2.6}, we will sometimes use the notation
\beq
   M_\tau =: \int_0^\tau \int_D G(s,y)\, W(ds,dy).
\eeq

\medskip

\noindent{\em Relation with Walsh's integral}\index{Walsh's integral}\index{integral!Walsh's}
\smallskip

Walsh's theory developed in \cite{walsh} defines in particular the stochastic integral of a predictable square integrable processes $G$
with respect to space-time white noise $W$. We refer to Section \ref{rdsecA.4} in Appendix \ref{app1} for the definition of predictable process. On this class of processes (which is smaller than the class of jointly measurable and adapted square integrable processes) Walsh's integral coincides with that of Definition \ref{ch1'-s4.d1}. Indeed, we prove this claim by checking the equality of both integrals on a class of {\em elementary processes,} the linear combinations of which are dense in the set of  predictable square integrable processes.

Indeed, consider the class of processes of the form
\beqn
\left(G(s,y;\omega) = X(\omega)\, 1_{]a,b]}(s) \, 1_A(y),\ (s,y)\in[0,T]\times D\right),
\eeqn
 where $0\leq a < b \leq T$, $X$ is $\tf_a$-measurable and $A\subset D$ is a bounded Borel set.
For $G$ in this class and $t\in[0,T]$, Walsh's stochastic integral of $(G(s,y))$ with respect to space-time white noise is defined by
\beq
\label{w-e}
  \int_{0}^t \int_D G(s,y)\, W(ds,dy)
  = X \left[W([0,t\wedge b]\times A) - W([0,t\wedge a]\times A)\right].
 \eeq
 On the other hand, according to Definition \ref{ch1'-s4.d1},
 \begin{align*}
 &\int_0^t\int_D G(s,y)\, W(ds,dy)\\
  &\qquad = \sum_{j=1}^\infty \int_0^t \langle G(s,\ast), e_j\rangle_V\, dW_s(e_j)\\
 &\qquad  = \sum_{j=1}^\infty  \int_{t\wedge a}^{t\wedge b} X\, \langle 1_{A}, e_j\rangle_V \, dW_s(e_j)\\
 &\qquad  = \sum_{j=1}^\infty\,  X\, \langle 1_{A}, e_j\rangle_V \left(W_{t\wedge b}(e_j) - W_{t\wedge a}(e_j)\right)\\
 &\qquad  = X\left[W_{t\wedge b}\left(\sum_{j=1}^\infty\, \langle 1_{A}, e_j\rangle_V\, e_j\right) - W_{t\wedge a}\left(\sum_{j=1}^\infty\, \langle 1_{A}, e_j\rangle_V\, e_j\right)\right]\\
&\qquad   = X\, [W_{t\wedge b} (1_A) - W_{t\wedge a} (1_A)]\\
&\qquad  = X \left[W([0,t\wedge b]\times A) - W([0,t\wedge a]\times A)\right].
\end{align*}
Since the last term is equal to the right-hand side of \eqref{w-e}, the claim is proved.
\bigskip

\section{Extensions of the stochastic integral}
\label{ch2'new-s3}

As in the case of the stochastic integral with respect to a finite-dimensional Brownian motion, the stochastic integral introduced in Section \ref{ch2new-s2} can be extended to integrands $G$ that are jointly measurable and adapted processes and satisfy
\beq
\label{weakint}
\int_0^T \Vert G(s,\ast)\Vert_V^2\, ds = \int_0^T \sum_{j=1}^\infty \langle G(s,\ast), e_j\rangle_V^2\, ds < \infty,\ \text{a.s.}
\eeq
This is done by localisation. Indeed, for any integer $N\ge 0$, define
\beq
\label{stop}
\tau_N = \inf\left\{t\in[0,T]: \int_0^t  \Vert G(s,\ast)\Vert_V^2\, ds\ge N\right\}\wedge T.
\eeq
Clearly, $(\tau_N,\, N\ge 1)$ is an increasing sequence of stopping times, and because of the assumption \eqref{weakint}, $\tau_N\uparrow T$, a.s. and even $\lim_{N\to\infty}P\{\tau_N=T\} = 1$.  Then, for any $t\in[0,T]$, we define
\beq
\label{weaksi}
\int_0^t \int_D G(s,y)\, W(ds,dy) = \lim_{N\to\infty}\int_0^t \int_D \left(1_{[0,\tau_N]}(s)\, G(s,y) \right) W(ds,dy).
\eeq
This a.s.~limit is well-defined. Indeed, since $\tau_N$ is a stopping time, the process $\left\{1_{[0,\tau_N]} (s)\, G(s,y),\, (s,y)\in[0,T]\times D\right\}$ is a jointly measurable and adapted process. Moreover,
\beq\label{rde2.3.4}
   \int_0^T 1_{[0,\tau_N]}(s)\, \Vert G(s,\ast)\Vert_V^2\, ds \le N, \qquad \text{a.s.},
\eeq
so for fixed $N\ge 1$, taking expectations on both sides of \eqref{rde2.3.4}, we see that the stochastic integral process $(Z^N_t)$ on the right-hand side of \eqref{weaksi} is a well-defined continuous martingale as in Proposition \ref{Proposition 6}.

The local property of the stochastic integral given in Lemma \ref{ch1'-llp} ensures that, for $1\le N\le M$ and $r\in [0,T]$, on $\{r\le \tau_N\}$, for $t\in[0,r]$,
\beq\label{rde2.3.5}
\int_0^t \int_D 1_{[0,\tau_M]} (s)\, G(s,y)\, W(ds,dy) = \int_0^t \int_D 1_{[0,\tau_N]} (s)\, G(s,y)\, W(ds,dy),
\eeq
a.s.
Since both stochastic integral processes are continuous in $t$, a.s.~on $\{r \leq \tau_N\}$,
\eqref{rde2.3.5} holds for all $t \in [0,r]$. Therefore, the limit in \eqref{weaksi} is stationary on $\{r \leq \tau_N\}$ for $t \in [0,r]$,
hence is stationary on $\{\tau_N = T\}$ for $t \in [0, T]$. It follows that the left-hand
side of \eqref{weaksi} is a well-defined process $(Z_t)$ with continuous sample paths a.s. In addition, a.s.~on $\{r \leq \tau_N\}$, for $t \in [0,r]$, $Z_t = Z^N_t$. Since $\{r \leq \tau_N\}$ increases to $\Omega$ a.s.~as $N\to\infty$,
 $(Z_t,\, t \in [0,T])$ is a continuous local martingale with respect to the filtration $(\tf_t,\, t\in[0,T])$, denoted
 \beqn
\left(\int_0^t\int_D  G(s,y)\, W(ds,dy),\ t\in[0,T]\right),
\eeqn
with quadratic variation process
\beq
\label
{cucuvar}
   \langle M\rangle_t = \int_0^t \Vert G(s,\ast) \Vert_V^2\, ds.
\eeq

\begin{prop}
\label{moreproper}
\noindent 1.\  The stochastic integral defined in \eqref{weaksi} satisfies the analogue of the local property in $\Omega$ stated in Lemma \ref{ch1'-llp} (for the stochastic integral constructed assuming
$E\left(\int_0^T \Vert G(s,\ast)\Vert_V^2\, ds\right) <\infty$).

\noindent 2. \ The stochastic integral defined in \eqref{weaksi} satisfies the analogue of the local property in space stated in Proposition \ref{prop-local-in-space}.

\noindent 3. \ The local martingale $\left(\int_0^t\int_D  G(s,y)\, W(ds,dy),\, t\in[0,T]\right)$ defined above satisfies Burkholder's inequality \eqref{ch1'-s4.4}.

\end{prop}
\begin{proof}\ Because of Lemma \ref{ch1'-llp}, the local property in $\Omega$ holds for the approximating sequence of integrals on the right-hand side of \eqref{weaksi}. Hence, it also holds for the limit of that sequence. With the same argument, applying Lemma \ref{prop-local-in-space}, we obtain the validity of the local property in space.

We refer to \cite[(4.1) Theorem, p.~160]{ry}) for a proof of Burkholder's  inequality in the setting of this proposition.
\end{proof}

\begin{prop}
\label{conv-probab}
Under condition \eqref{weakint}, for all $t \in [0,T]$, a.s.,
\beqn
   \int_0^t \int_D G(s,y)\, W(ds,dy) = \sum_{j=1}^\infty\int_0^t \langle G(s,\ast),e_j\rangle_V\, dW_s(e_j),
\eeqn
where the series converges in probability, uniformly in $t\in[0,T]$.
\end{prop}

\begin{proof}
Fix $\eta>0$ and let $\tau_N$ be as in \eqref{stop}. Since $\lim_{N\to\infty}P\{\tau_N<T\}=0$, there exists $N_0$ such that $P\{\tau_{N_0}<T\} \le \eta$. Let $\int_0^t \int_D G(s,y)\, W(ds,dy)$ be defined as in \eqref{weaksi} and let $\varepsilon>0$ be fixed. Then, for any $M\ge 1$,
\begin{align*}
&P\left\{\sup_{t\in[0,T]}\Big\vert \sum_{j=1}^M \int_0^t \langle G(s,\ast),e_j\rangle_V\, dW_s (e_j) - \int_0^t\int_D G(s,y)\, W(ds,dy)\Big\vert >\varepsilon\right\}\\
&\qquad \le P\Big\{\tau_{N_0}<T\Big\}\\
&\qquad \qquad+ P\Bigg\{\sup_{t\in[0,T]}\Big\vert \sum_{j=1}^M \int_0^t \langle G(s,\ast),e_j\rangle_V\, dW_s(e_j)\\
&\qquad \qquad\qquad- \int_0^t \int_D G(s,y)\, W(ds,dy)\Big\vert >\varepsilon, \tau_{N_0}=T\Bigg\}.
\end{align*}
The first term on the right-hand side is bounded above by $\eta$. By Lemma \ref{ch1'-llp}, on $\{ \tau_{N_0}=T\}$,
$$
    \int_0^t \langle G(s,\ast),e_j\rangle_V\, dW_s(e_j) = \int_0^T 1_{[0,\tau_{N_0} \wedge t]}(s)\, \langle G(s,\ast),e_j\rangle_V\, dW_s(e_j)
$$
and
$$
   \int_0^t \int_D G(s,y)\, W(ds,dy) = \int_0^T \int_D 1_{[0,\tau_{N_0} \wedge t]}(s)\, G(s,y)\, W(ds,dy).
$$
By \eqref{rde2.3.4},
\beq
\label{boundN0}
   E\left[\int_0^T \sum_{j=1}^\infty 1_{[0,\tau_{N_0} \wedge t]}(s)\, \langle G(s,\ast),e_j\rangle_V^2\, ds  \right] \leq N_0,
\eeq
therefore a.s. on $\Omega$,
\begin{align}
\label{rde2.3.6}
 & \int_0^T \int_D 1_{[0,\tau_{N_0} \wedge t]}(s)\, G(s,y)\, W(ds,dy)\notag \\
 &\qquad \qquad \qquad= \sum_{j=1}^\infty \int_0^T 1_{[0,\tau_{N_0} \wedge t]}(s)\, \langle G(s,\ast),e_j\rangle_V \, dW_s(e_j),
\end{align}
where the series converges in ${\bH}^2$ (see Section \ref{ch2new-s2} for the definition of this space). Along a subsequence $(m_k)$, the series converges a.s., uniformly in $t\in[0,T]$. 

We deduce from the Chebychev and the Burkholder inequalities that
\begin{align*}
   &P\left\{\sup_{t\in[0,T]}\left\vert \sum_{j=1}^M \int_0^t \langle G(s,\ast),e_j\rangle_V\, dW_s (e_j) - \int_0^t\int_D G(s,y)\, W(ds,dy)\right\vert >\varepsilon\right\}\\
&\qquad \le \eta + P\left\{\sup_{t\in[0,T]}\left\vert \sum_{j=M+1}^\infty \int_0^t1_{[0,\tau_{N_0} \wedge t]}(s)\, \langle G(s,\ast),e_j\rangle_V \, dW_s(e_j)\right\vert > \varepsilon  \right\} \\
&\qquad \le \eta + \frac{1}{\varepsilon^2} E\left[\sup_{t\in[0,T]}\left( \sum_{j=M+1}^\infty \int_0^t 1_{[0,\tau_{N_0} \wedge t]}(s)\, \langle G(s,\ast),e_j\rangle_V \, dW_s(e_j)\right)^2 \right]\\
&\qquad \le \eta + \frac{1}{\varepsilon^2} E\left[\sum_{j=M+1}^\infty \int_0^T 1_{[0,\tau_{N_0} \wedge t]}(s)\, \langle G(s,\ast),e_j\rangle_V^2 \, ds\right]
\end{align*}
and this converges to $0$ as $M \to \infty$ by \eqref{boundN0}.  This proves the proposition.
\end{proof}

\medskip

The next proposition gives a condition on the process $G$ under which the indefinite stochastic integral is an $L^1(\Omega)$--martingale (rather than a square-integrable martingale as in Proposition \ref{Proposition 6}). Recall the notation $V=L^2(D)$.

\begin{prop}
\label{ch1'-s4-p0}
Let $G$ be a jointly measurable and adapted stochastic process such that
\beq
\label{ch1'-s4.10}
E\left[\left(\int_0^T  \Vert G(s,\ast)\Vert_V^2\, ds\right)^{\half}\right] < \infty.
\eeq
Then
\beqn
\left(\int_0^t \int_D G(s,y) \, W(ds,dy), \ t\in[0,T]\right)
\eeqn
is an $L^1(\Omega)$--martingale, and for $t\in [0,T]$,
\beq
\label{ch1'-s4.11}
\int_0^t \int_D G(s,y) \, W(ds,dy) = \sum_{j=1}^\infty \int_0^t \langle G(s,\ast),e_j\rangle_V\, dW_s(e_j) \qquad\mbox{a.s.},
\eeq
where the series converges in $L^1(\Omega)$, uniformly in $t\in[0,T]$.
\end{prop}

\begin{proof}
Consider the local martingales
\begin{align*}
M_t = \int_0^t \int_D G(s,y)\, W(ds,dy),\qquad
M_t^j = \int_0^t \langle G(s,\ast),e_j\rangle_{V}\, dW_s(e_j).
\end{align*}
By the Burkholder-Davis-Gundy inequality \eqref{ch1'-s4.4} with $p=1$ (see Proposition \ref{moreproper}, claim 3.),
\beqn
E\left[\sup_{t\in[0,T]} |M_t|\right] \le cE\left[\langle M\rangle_T^\half\right] = c E\left[\left(\int_0^t  \Vert G(s,\ast)\Vert_V^2 \, ds \right)^{\half}\right] < \infty,
\eeqn
therefore, $(M_t,\, t\in[0,T])$ is in fact a martingale in $L^1(\Omega)$: see \cite[Propositions 1.8 and 1.1]{chung-williams}. Similarly, for $j\ge 1$,
\begin{align*}
E\left[\sup_{t\in[0,T]} |M_t^j|\right] & \le c E\left[\langle M^j\rangle_T^\half\right] = c E\left[\left(\int_0^t   \langle G(s,\ast), e_j\rangle^2_{V}\, ds\right)^{\frac{1}{2}}\right]\\
& \le c E\left[\left(\int_0^T  \Vert G(s,\ast)\Vert^2_{V} \, ds\right)^{\frac{1}{2}}\right] < \infty,
\end{align*}
so $(M_t^j,\, t\in[0,T])$ is also an $L^1(\Omega)$--martingale.

Since \eqref{ch1'-s4.10} implies \eqref{weakint}, Proposition \ref{conv-probab} shows that
 the series $\sum_{j=1}^\infty M_t^j$ on the right-hand side of \eqref{ch1'-s4.11} converges in probability and is equal to $M_t$ a.s.

It remains to prove that the series $\sum_{j=1}^\infty M_t^j$ converges in $L^1(\Omega)$, uniformly in $t\in[0,T]$. Indeed, for $1\le n\le m$,
\beqn
\left\Vert \sup_{t\in[0,T]}\left(\sum_{j=1}^m M_t^j - \sum_{j=1}^n M_t^j\right)\right\Vert_{L^1(\Omega)} = \left\Vert\, \sup_{t\in[0,T]}\sum_{j=n+1}^m M_t^j\, \right\Vert_{L^1(\Omega)}.
\eeqn
By the Burkholder-Davis-Gundy inequality referred to above,
\beq
\label{ch1'-s4.12}
\left\Vert\, \sup_{t\in[0,T]}\left\vert \sum_{j=n+1}^m M_t^j\right\vert\,\right\Vert_{L^1(\Omega)}
\le cE\left[\Big\langle\sum_{j=n+1}^m M^j\Big\rangle_T^{\half}\right].
\eeq
Because the $(W_s(e_j),\, j\in\mathbb{N})$, are independent,
\beqn
\Big\langle\sum_{j=n+1}^m M^j\Big\rangle_T = \sum_{j=n+1}^m\langle M^j\rangle_T ,
\eeqn
so the right-hand side of \eqref{ch1'-s4.12} is equal to
\beq
\label{ch1'-s4.13}
c E\left[\left(\sum_{j=n+1}^m \int_0^T \langle G(s,\ast),e_j\rangle_{V}^2\, ds\right)^\half\right].
\eeq
This converges to $0$ as $n,m\to\infty$, since
\beq
\label{ch1'-s4.13-a}
E\left[\left(\sum_{j=1}^\infty \int_0^T \langle G(s,\ast),e_j\rangle_{V}^2\, ds\right)^\half\right]
= E\left[\left(\int_0^T \Vert G(s,\ast)\Vert_V^2 \, ds \right)^\half\right] < \infty,
\eeq
by \eqref{ch1'-s4.10}.
Indeed, let $Z = \sum_{j=1}^\infty \int_0^T \langle G(s,*), e_j\rangle_V^2 \, ds$. By \eqref{ch1'-s4.13-a}, $E[Z^\half] < \infty$, so $0 \leq Z < \infty$ a.s. Define $Z_n = \sum_{j=n+1}^\infty \int_0^T \langle G(s,*), e_j\rangle_V^2 \, ds$. Then $0 \leq Z_n \leq Z < \infty$ a.s., and $Z_n \downarrow 0 $ a.s.~as $n \to \infty$ because $Z_n$ is the tail sum of the a.s.~convergent series that defines $Z$. By the dominated convergence theorem, $\lim_{n\to\infty} E[Z_n^\half] = 0$. Since the expression in \eqref{ch1'-s4.13} is bounded above by $E[Z_n^\half]$, we obtain that this expression converges to $0$.

This shows that the right-hand side of \eqref{ch1'-s4.11} converges in $L^1(\Omega)$, uniformly in $t\in[0,T]$. The proposition is proved.
\end{proof}
\medskip

\noindent{\em Weakening the measurability requirements on integrands}
\smallskip

We have defined the stochastic integral with respect to space-time white noise for stochastic processes $G=(G(s,y),\, (s,y)\in [0,T] \times D)$ that are jointly measurable, adapted and satisfy \eqref{weakint}. It is possible to weaken the measurability requirement, as we now explain.

In the classical It\^o theory of stochastic integrals with respect to Brownian motion \cite{chung-williams}, one begins by defining the stochastic integral of predictable processes (see Appendix \ref{app1}, Section \ref{rdsecA.4}) that satisfy \eqref{weakint}. However, observe that if $(X_1(s),\, s \in [0,T])$ and $(X_2(s),\, s \in [0,T])$ are predictable processes such that
$$
   \int_0^T(X_1(s) - X_2(s))^2\, ds = 0\qquad \text{a.s.,}
$$
then they will have the same stochastic integral. It is therefore natural to extend the stochastic integral with respect to Brownian motion to processes that are $\cP^*$-measurable, where $\cP^*$ is the completion of $\cP$ with respect to $dsdP$-null sets. It then turns out that processes $(X(s),\, s \in \IR_+)$ that are $\cB_{\IR_+}\times \cF$-measurable and adapted are in fact $\cP^*$-measurable \cite [Theorem 3.8]{chung-williams}.

Applying these ideas in the context of the stochastic integral with respect to space-time white noise, we see that the assumption ``jointly measurable and adapted" can be replaced by ``$(y,s,\omega) \mapsto G(s,y,\omega)$ is $\cB_D\times \cP^*$-measurable, in Definition \ref{ch1'-s4.d1}, and all the results of Sections \ref{ch2new-s2} and \ref{ch2'new-s3} remain valid.

\section{Stochastic Fubini's theorem}\index{stochastic!Fubini's theorem}\index{Fubini's theorem!stochastic}\index{theorem!stochastic Fubini's}
\label{ch2'new-s4}

Let $(X, \mathcal{X})$ be a measure space and let $\mu$ be a $\sigma$--finite (nonnegative) measure on $X$. We let $W$, $(\cF_s,\, s\in[0,T])$ and $(e_j,\, j\ge 1)$ be as defined at the beginning of Section \ref{ch2new-s2}.

We recall (see Appendix \ref{app1}, Section \ref{rdsecA.4}) that two stochastic processes $(u(s,y),\, (s,y) \in [0,T]\times D)$ and $(v(s,y),\, (s,y) \in [0,T]\times D)$ are {\em indistinguishable} if $P\{u(s,y) = v(s,y),\ {\text{for all}}\ (s,y) \in [0,T]\times D\}=1$,
whereas  $(u(s,y),\, (s,y) \in [0,T]\times D)$
 is a {\em modification} or {\em version} of $(v(s,y),\, (s,y) \in [0,T]\times D)$ if for all $(s,y) \in [0,T]\times D$, we have $v(s,y) = u(s,y)$ a.s.~(where the implied null set may depend on $(s,y)$).

\begin{thm}
\label{ch1'-tfubini}
Let $G: X\times [0,T]\times D\times \Omega \rightarrow \re$ be $\cX\times \cB_{[0,T]} \times \cB_D\times \cF$-measurable and such that, for fixed $s\in [0,T]$, the partial function $(x,y,\omega) \mapsto G(x,s,y,\omega)$ from $X\times D\times \Omega$ into $\re$ is $\cX\times\mathcal{B}_D\times \tf_s$--measurable. Suppose that
\beq
\label{ch1'.f1}
\int_{X} \mu(dx)\, \Vert G(x, \cdot, \ast)\Vert_{L^2([0,T]\times D)} < \infty, \qquad a.s.
\eeq
Then the following statements hold:
\begin{description}
\item{(a)} There exists $X_0\in\mathcal{X}$ with $\mu(X\setminus X_0)=0$ such that for any $x\in X_0$,
$G(x, \cdot, \ast)\in L^2([0,T]\times D)$ a.s.
There is an $\cX\times \cB_{[0,T]}\times \cF_T$-measurable
 map $Z: X \times [0,T] \times \Omega \to \R$ such that, for all $x \in X_0$, $Z(x, \cdot)$ is indistinguishable from the stochastic integral process $G(x,\cdot,\ast) \cdot W$, and its sample paths are continuous. In addition,
\beq
\label{ch1'.f3}
\sup_{t\in[0,T]} \int_ X \mu(dx)\, \vert Z(x,t)\vert < \infty, \quad a.s.
\eeq
\item{(b)} Almost surely,
\begin{align}
\label{ch1'.f4}
\left\Vert \int_X \mu(dx)\ |G(x, \cdot, \ast)|\right\Vert_{L^2([0,T]\times D)}<\infty.
\end{align}

Consequently, for $dsdydP$-almost all $(s, y, \omega)\in[0,T]\times D\times \Omega$, $x\mapsto G(x, s, y, \omega)$ is $\mu$--integrable. Further,
 the stochastic integral process
\beq
\label{ch1'.f5}
\left(\int_0^t \int_D \left(\int_X \mu(dx)\ G(x, s, y)\right) \, W(ds, dy),\ t\in[0,T]\right)
\eeq
is well-defined in the sense of Section \ref{ch2'new-s3}.
\item{(c)} Almost surely, for all $t\in[0,T]$,
\begin{align}
\label{ch1'.f6}
&\int_X \mu(dx)  \left(\int_0^t \int_D G(x,s,y)\, W(ds,dy)\right)\notag\\
&\qquad \qquad=  \int_0^t \int_D \left(\int_X \mu(dx)\, G(x, s, y)\right) \, W(ds, dy),
\end{align}
 where, by definition, the left-hand side is equal to $\int_{X_0}  \mu(dx)\,  Z(x,t)$. A.s., this process has continuous sample paths.

 \end{description}
\end{thm}

\begin{remark}
\label{fubini-name}
The name ``stochastic Fubini's theorem'' refers to the identity \eqref{ch1'.f6}. Part (a) of the statement implies that the integral on the left-hand side of \eqref{ch1'.f6} is well-defined, while part (b) leads to a similar conclusion for the integral of the right-hand side.
\end{remark}

\noindent{\em Proof of Theorem \ref{ch1'-tfubini}.}
We will proceed through several steps.
\smallskip

\noindent{\em Step 1. Some elements of the proof of (a).}
By \eqref{ch1'.f1}, there is a $dP$-null set $F_0$ such that \eqref{ch1'.f1} holds outside of $F_0$. Therefore, for $\omega \not\in F_0$,
there is a $\mu(dx)$-null set $X_1(\omega)$ such that for $x \not\in X_1(\omega)$, $\Vert G(x,\cdot,*,\omega) \Vert_{L^2([0,T] \times D)} < \infty$. Since
\beqn
   \{(x,\omega): \Vert G(x,\cdot,*) \Vert_{L^2([0,T] \times D)} = \infty\} \in \cX \times \F,
   \eeqn
we deduce that this is a $\mu(dx)dP$-null set. Hence, by Fubini's theorem, there is a $\mu(dx)$-null set $X \setminus X_0$, which can be chosen in $\cX$, such that for $x \in X_0$, $G(x,\cdot,*) \in L^2([0,T] \times D)$ a.s.

   Let $Z$ be the jointly measurable function given by Theorem \ref{ch1'-l-m1} (b), with $X$ there replaced by $X_0$. We already know that this function $Z$ satisfies the conclusions of (a) except for \eqref{ch1'.f3},  which will be checked at the end of the proof. For $x \not\in X_0$ and all $(t,\omega)$, we set $Z(x,t,\omega) = 0$.
\smallskip

\noindent{\em Step 2. Proof of (b)}.
For any fixed $(s,y,\omega)\in[0,T]\times D\times \Omega$, the map $x\mapsto G(x,s,y,\omega)$ is measurable, and  by Minkowski's inequality,
\beq
\label{ch1'.f6b}
\left\Vert \int_X \mu(dx)\, |G(x, \cdot, \ast)|\right\Vert_{L^2([0,T]\times D)}\le \int_X \mu(dx)\, \Vert G(x, \cdot, \ast)\Vert_{L^2([0,T]\times D)}.
\eeq
From \eqref{ch1'.f1}, we obtain
\beqn
\left\Vert \int_X \mu(dx)\, |G(x, \cdot, \ast)|\right\Vert_{L^2([0,T]\times D)} <\infty,\quad {\text{a.s.}}
\eeqn
This is property \eqref{ch1'.f4}, which  implies that, for $dsdydP$-almost all $(s, y, \omega)\in[0,T]\times D\times \Omega$, $x\mapsto G(x, s, y, \omega)$ is $\mu$--integrable.

By the deterministic Fubini's theorem, the process
\beqn
\left(\int_X \mu(dx)\, G(x,s,y),\ (s,y)\in[0,T]\times D\right)
\eeqn
 is jointly measurable and adapted, therefore the indefinite integral process \eqref{ch1'.f5} is well-defined (in the sense of Section \ref{ch2'new-s3}).
 \smallskip

 \noindent{\em Step 3. Proof of a localised version of \eqref{ch1'.f6}.}
 We now turn to part (c). First, we will establish a localised version of \eqref{ch1'.f6}: 

Define the increasing sequence of stopping times
\beqn
\tau_n=
\inf\left\{t\in[0,T]: \int_X \mu(dx)\, \Vert G(x, \cdot, \ast)\Vert_{L^2([0,t]\times D)}\ge n\right\}\wedge T,\
n\in \mathbb{N}.
\eeqn
By \eqref{ch1'.f1}, we have $\lim_{n\to\infty}P\{\tau_n = T\} =1$.

Observe that  $t\mapsto \int_X \mu(dx)\, \Vert G(x, \cdot, \ast)\Vert_{L^2([0,t]\times D)}$
is continuous a.s., and the inequality in $\eqref{ch1'.f6b}$ is path-by-path. Hence, setting $t:=\tau_n$ there, we see that a.s.,
\begin{align*}
\left\Vert \int_X \mu(dx)\, G(x, \cdot, \ast)\right\Vert_{L^2([0,\tau_n]\times D)}
 & \le\int_X \mu(dx)\, \Vert G(x,\cdot,\ast)\Vert_{L^2([0,\tau_n] \times D)}\\
 &\le n.
\end{align*}
In particular,
\beq
\label{p65(*1)}
 \int_X \mu(dx)\, E\left( \Vert G(x,\cdot,*) \Vert_{L^2([0,\tau_n] \times D)}\right) \leq n ,
 \eeq
 and
\beqn
E\left[\int_0^{\tau_n} ds \int_D dy \left(\int_X \mu(dx)\, G(x, s, y)\right) ^2\right] \le n^2.
\eeqn
Therefore, using Lemma \ref{rdlem2.2.6}, we can define the indefinite integral process (which is a
square integrable continuous martingale) by
\beqn
M_t:= \int_0^{t} \int_D 1_{\{s \leq \tau_n\}}\left(\int_X \mu(dx)\, G(x, s, y)\right) W(ds,dy), \quad t\in[0,T],
\eeqn
as in Lemma \ref{lemma 5}.
Let
\beq
\label{ch1'.f6a}
Z^N_t= \sum_{j=1}^N \int_0^{t\wedge \tau_n}\Big\langle \int_X \mu(dx)\, G(x, s, \ast), e_j\Big\rangle_{V}\, dW_s(e_j),
\quad t\in[0,T].
\eeq
Then the martingale $(Z^N_t)$ converges in ${\bH}^2$ to $(M_t,\, t\in[0,T])$.
\smallskip

We next prove that a.s., for almost all $s\in[0,T]$,
\beq
\label{innerproducts}
\Big\langle\int_X \mu(dx)\, G(x, s, \ast), e_j\Big\rangle_{V}= \int_X \mu(dx)\, \langle G(x, s, \ast), e_j\rangle_{V},
\eeq
for all $j\ge 1$.
Since the right-hand side is a jointly measurable and adapted process, this will imply that $Z^N$ is indistinguishable from
\beq
\label{(*3)}
   \sum_{j=1}^N \int_0^{t \wedge \tau_n}\left (\int_X \mu(dx)\, \langle G(x,s,*), e_j\rangle_V\right)  dW_s(e_j),
   \eeq
which we denote again by $Z^N$, and preserves the convergence in ${\bH}^2$ to $(M_t,\, t\in[0,T])$.


We now prove \eqref{innerproducts}. Since
\beq
\label{ch1'.f7}
\left\langle \int_X \mu(dx)\, G(x, s, \ast), e_j\right\rangle_{V} = \int_D dy \left(\int_X \mu(dx)\, G(x, s, y)\right)e_j(y),
\eeq
the identity \eqref{innerproducts} will follow by applying
the deterministic Fubini's theorem. Let us check that the assumptions of this theorem are satisfied.
Indeed, by Minkowski's inequality and \eqref{ch1'.f1},
\begin{align*}
   \left\Vert \int_ X \mu(dx)\, \Vert G(x, \cdot, \ast)\Vert_V \right\Vert_{L^2([0,T])} &\leq \int_ X \mu(dx)\, \left\Vert\, \Vert G(x, \cdot, \ast)\Vert_V \right\Vert_{L^2([0,T])}\\
   &= \int_ X \mu(dx)\, \Vert G(x, \cdot, \ast)\Vert_{L^2([0,T]\times D)}\\
   &< \infty\quad {\text{a.s.}}
\end{align*}
Thus, a.s., for almost all $s\in[0,T]$,
\beqn
\int_ X \mu(dx)\, \Vert G(x, s,\ast)\Vert_V < \infty.
\eeqn
By the Cauchy-Schwarz inequality,
\begin{align*}
\int_X \mu(dx) \int_D dy\, |G(x, s, y)|\, |e_j(y)| &\le \int_ X \mu(dx)\, \Vert G(x, s, \ast)\Vert_V\, \Vert e_j\Vert_V\\
&= \int_ X \mu(dx)\, \Vert G(x, s, \ast)\Vert_V,
\end{align*}
because $\Vert e_j\Vert_V=1$. Therefore, a.s., for almost all $s\in[0,T]$,
\beqn
\int_X \mu(dx) \int_D dy\, |G(x, s, y)|\, |e_j(y)| < \infty.
\eeqn
This implies \eqref{innerproducts}.
\smallskip

Continuing towards the proof of a localized version of  \eqref{ch1'.f6}, we write the expression in \eqref{(*3)} as $\sum_{j=1}^N I_j(t)$, where, for $j \geq 1$,
\beqn
   I_j(t) = \int_0^{t \wedge \tau_n}\left (\int_X \mu(dx)\, \langle G(x,s,*), e_j\rangle_V\right)  dW_s(e_j).
\eeqn
Our next goal is to permute the integrals in each term $I_j(t)$ of \eqref{(*3)}. For this, we apply the stochastic Fubini's Theorem for Brownian motion (Theorem \ref{app1-2.t1}) to $g(x,s,\omega) := \langle G(x,s,*), e_j\rangle_V$  and $B_s := W_s(e_j)$ there. The hypotheses of Theorem \ref{app1-2.t1} are satisfied (notice that we can use the same $X_0$ as for $Z$ in part (a) above). We obtain that a.s., for all $t \in [0,T]$,
\beqn
     I_j(t) = \int_X \mu(dx) \int_0^{t \wedge \tau_n} \langle G(x,s,*), e_j\rangle_V\,  dW_s(e_j),
     \eeqn
where the stochastic integral on the right-hand-side is the jointly measurable function $\Psi_j$ with continuous sample paths given by Theorem \ref{app1-2.t1}, evaluated at $(x,t \wedge \tau_n)$. Therefore, by \eqref{(*3)}, $Z^N$ is indistinguishable from
  \beqn
  \int_X \mu(dx)\, \sum_{j=1}^N \int_0^{t \wedge \tau_n} \langle G(x,s,*), e_j\rangle_V\,  dW_s(e_j),
  \eeqn
which is indistinguishable from
 \beqn
 Z^{\prime,N}_t := \int_X \mu(dx)  \int_0^{t} \int_D G_N(x, s, y)\, W(ds, dy),
 \eeqn
where
   \beqn
   G_N(x,s,y) =   \sum_{j=1}^N 1_{[0,\tau_n]}(s)\, \langle G(x,s,*) , e_j\rangle_V\,  e_j(y).
   \eeqn
For fixed $x \in X$, the process $J^N$ defined by
\beqn
     J^N_t(x) = \sum_{j=1}^N \int_0^{t \wedge \tau_n} \langle G(x,s,*), e_j\rangle_V\,  dW_s(e_j),
     \eeqn
where the stochastic integral on the right-hand side is the function $\Psi_j(x,t \wedge \tau_n)$ mentioned above, is indistinguishable from the $L^2(\Omega)$-bounded continuous martingale
\beqn
     \left( \int_0^{t} \int_D G_N(x,s,y)\, W(ds,dy), \  t \in [0,T]\right),
     \eeqn
by the definition of $(G_N(x,\cdot,*) \cdot W)$ (we take the function given by Theorem \ref{ch1'-l-m1} (b)).
  Consider the set
  \beqn
     A= \{(x,\omega): t \mapsto J^N_t(x,\omega)\  {\text{is the same as}}\  t \mapsto (G_N(x,\cdot,*) \cdot W)_t(\omega) \},
     \eeqn
which belongs to $\cX \times \F_T$. For fixed $x \in X$, $\{\omega: (x,\omega) \not\in A\}$ has probability $0$, so the $d\mu dP$-measure of $A^c$ is $0$. Therefore, a.s., for $d\mu$-a.a.~$x \in X$,  $t \mapsto J^N_t(x)$ is the same as $t \mapsto (G_N(x,\cdot,*) \cdot W)_t$. Hence, a.s., for all $t \in [0,T]$,
\beqn
    Z^N_t = Z^{\prime,N}_t = \int_X \mu(dx)\, (G_N(x,\cdot,*) \cdot W)_{t}.
 \eeqn

We now establish two properties of the map $Z$ given in part (a).
\smallskip

\noindent{\em Property (i)}. The map $Z$ given in (a) satisfies the following two properties:
\begin{align*}
E\left[\sup_{t\in [0,T]} \int_X \mu(dx)\, \vert Z(x,t \wedge \tau_n) \vert\right] & \leq n,\\
\sup_{t\in [0,\tau_n]} \int_X \mu(dx)\, \vert Z(x,t) \vert &< \infty\quad  \text{a.s.}
\end{align*}

    Indeed, for every $\omega \in \Omega$, if $t \leq \tau_n(\omega)$, then for all $x \in X_0$, $Z(x,t \wedge \tau_n)(\omega) = Z(x,t,\omega)$, and if $T \geq t > \tau_n(\omega)$, then for all $x \in X_0$, $Z(x,t \wedge \tau_n)(\omega) =Z(x,\tau_n(\omega),\omega)$. It follows that for every $\omega \in \Omega$,
     \beq
     \label{(*A)}
      \sup_{t\in [0,T]} \int_X \mu(dx)\, \vert Z(x,t \wedge \tau_n) (\omega)\vert  = \sup_{t\in [0,\tau_n(\omega)]} \int_X   \mu(dx)\, \vert Z(x,t,\omega) \vert.
      \eeq
Moreover,
   \begin{align*}
 E\left[\sup_{t\in[0,\tau_n]}  \int_X \mu(dx)\, |Z(x, t)|\right] & \le E\left[ \int_X \mu(dx)\, \sup_{t\in[0,\tau_n]} |Z(x, t)|\right] \\
 & = \int_X \mu(dx)\, E\left[\sup_{t\in[0,\tau_n]} |Z(x, t)|\right].
 \end{align*}
 By the Burkholder-Davis-Gundy inequality \eqref{ch1'-s4.4} and \eqref{p65(*1)}, this is bounded above by
 \beqn
  \int_X \mu(dx)\, E\left[\left(\int_0^{\tau_n} ds \int_D dy\,  G^2(x,s,y)\right)^\half\right] \le n.
  \eeqn
  Therefore, the first property holds and
  \beqn
  \sup_{t\in[0,\tau_n]} \int_X \mu(dx)\, |Z(t,x)| < \infty, \qquad a.s.,
  \eeqn
  which is the second property.
\smallskip

\noindent{\em Property (ii)}. The following holds:
\beqn
\lim_{N \to \infty} E\left[\sup_{t \in [0,T]} \left\vert Z^N_t - \int_X \mu(dx)\, Z(x,t \wedge \tau_n)\right\vert \right] = 0.
\eeqn
Indeed, the expectation is equal to
   \begin{align*}
   &E\left[\sup_{t \in [0,T]} \left\vert\int_X \mu(dx)\, (G_N(x,\cdot,*) \cdot W)_{t} - \int_X \mu(dx)\, Z(x,t \wedge \tau_n)\right\vert \right]\\
   &\qquad \leq \int_X \mu(dx)\, E\left[\sup_{t \in [0,T]} \left\vert (G_N(x,\cdot,*) \cdot W)_{t} - Z(x,t \wedge \tau_n) \right\vert\right].
   \end{align*}
For fixed $x \in X$, we can replace $Z(x, \cdot)$ by the indistinguishable process $(G(x,\cdot,*) \cdot W)$, to see that this is equal to
\begin{align}
\label{(*4)}
  & \int_X \mu(dx)\, E\left[\sup_{t \in [0,T]} \left\vert (G_N(x,\cdot,*) \cdot W)_{t} - (G(x,\cdot,*) \cdot W)_{t \wedge \tau_n}\right\vert\right]\notag\\
   &\qquad \qquad= \int_X \mu(dx)\, E\left[\sup_{t \in [0,T]} \left\vert ((G_N(x,\cdot,*) - G(x,\cdot,*)) \cdot W)_{t \wedge \tau_n}\right\vert\right]\notag\\
  &\qquad \qquad\leq  \int_X \mu(dx)\, E\left[\Vert G_N(x,\cdot,*) - G(x,\cdot,*)\Vert_{L^2([0,\tau_n] \times D)}\right],
  \end{align}
  by the Burkholder-Davis-Gundy inequality \eqref{ch1'-s4.4}.
This converges to $0$ as $N \to \infty$ as we now show.

    Since $G_N(x,s,*)$ is the projection of $G(x,s,*) 1_{[0,\tau_n]}(s)$ onto a finite-dimensional subspace of $L^2(D)$, for all $x \in X$ and all $\omega \in \Omega$,
    \beqn
       \Vert G_N(x,\cdot,*) \Vert _{L^2([0,\tau_n] \times D)} \leq \Vert G(x,\cdot,*) \Vert _{L^2([0,\tau_n] \times D)},
       \eeqn
so
\beqn
   \Vert G_N(x,\cdot,*) - G(x,\cdot,*)\Vert_{L^2([0,\tau_n] \times D)} \leq 2 \Vert G(x,\cdot,*) \Vert _{L^2([0,\tau_n] \times D)}.
   \eeqn
Applying Fubini's theorem to the product measure $\mu \times P$, and using \eqref{p65(*1)},
\begin{align*}
    \int_X \mu(dx)\, E\left[\Vert G(x,\cdot,*)\Vert_{L^2([0,\tau_n] \times D)}\right]  &= E\left[\int_X \mu(dx)\, \Vert G(x,\cdot,*)\Vert_{L^2([0,\tau_n] \times D)}\right]\\
    & \leq n < \infty.
    \end{align*}
Therefore, for $d\mu dP$- a.a. $(x,\omega)$,
   $\Vert G(x,\cdot,*)\Vert_{L^2([0,\tau_n] \times D)} < \infty$,
so
   \beqn
   \lim_{N\to\infty}\Vert G_N(x,\cdot,*) - G(x,\cdot,*)\Vert_{L^2([0,\tau_n] \times D)} =0.
   \eeqn
Apply the dominated convergence theorem with $d\mu dP$-measure to conclude that the expression in \eqref{(*4)} converges to $0$ as $N \to \infty$.
This ends the proof of Property (ii).

Property (ii) implies that for some subsequence $(N_k)$, $k\geq 1$, a.s.,
\beqn
  \lim_{N_k \to \infty} \sup_{t \in [0,T]} \left\vert Z^{N_k}_t - \int_X \mu(dx)\, Z(x,t \wedge \tau_n)\right\vert  = 0.
  \eeqn
We conclude that a.s.,
  $ t \mapsto \int_X \mu(dx)\, Z(x,t \wedge \tau_n)$
is continuous. Since $Z^N$ converges to $(M_t,\, t\in[0,T])$ in ${\bH}^2$, we conclude that a.s., for all $ t\in [0,T]$,
 \beqn
 M_t = \int_X \mu(dx)\, Z(x,t \wedge \tau_n).
 \eeqn
This means that the two processes
\beq
\label{(*5)}
    \left(\int_X \mu(dx)\, Z(x,t \wedge \tau_n),\  t \in [0,T]\right)
    \eeq
and
\beq
\label{(*6)}
    \left(\int_0^{t \wedge \tau_n} \int_D \left(\int_X \mu(dx)\, G(x,s,y)\right) W(ds,dy),\  t \in [0,T]\right)
    \eeq
are indistinguishable, proving \eqref{ch1'.f6} with $t$ there replaced by $t \wedge \tau_n$.
This completes the proof of the localized version of \eqref{ch1'.f6}.
\smallskip

\noindent{\em Step 4. Proof of \eqref{ch1'.f6} and \eqref{ch1'.f3}}. 
On the event $\{\tau_n = T\}$, by the local property in $\Omega$ of the stochastic integral (see part 1 of Proposition \ref{moreproper}), we can replace $t \wedge \tau_n$ by $t$ in the upper bound of the integral in \eqref{(*6)}.
And as in the proof of Property (i), on the event $\{\tau_n = T\}$, we can also replace $t \wedge \tau_n$ by $t$ in \eqref{(*5)}. This means that the two processes, restricted to $\{\tau_n = T\}$, are indistinguishable. Since $P(\cup_{n=1}^\infty \{\tau_n = T\}) = 1$, we get \eqref{ch1'.f6}.

   Finally, using \eqref{(*A)} and Property (i), we see that on the event $\{\tau_n = T\}$, \eqref{ch1'.f3} holds.
This completes the proof of Theorem \ref{ch1'-tfubini}.
\qed
\medskip

 We can also obtain Fubini's theorem under a Walsh-type condition (\cite[Theorem 2.6]{walsh}), as follows.

  \begin{cor}
 \label{ch1'-cfubini}
 If, in Theorem \ref{ch1'-tfubini}, we assume that the measure $\mu$ is finite  and condition \eqref{ch1'.f1} is replaced by
  \beq
 \label{ch1'.f13}
 E\left[\int_X \mu(dx)\, \Vert G(x,\cdot,\ast)\Vert^2_{L^2([0,T]\times D)}\right] < \infty,
 \eeq
 then  \eqref{ch1'.f1} holds as do all the conclusions of Theorem \ref{ch1'-tfubini}.
 \end{cor}
 \begin{proof}
 Condition \eqref{ch1'.f13} clearly implies $\int_X \mu(dx) \Vert G(x,\cdot,\ast)\Vert_{L^2([0,T]\times D)}^2 < \infty$, a.s.  Because  $\mu(X)<\infty$, applying
 the Cauchy-Schwarz inequality, we see that
 \begin{align*}
& \int_X \mu(dx)\, \Vert G(x,\cdot,\ast)\Vert_{L^2([0,T]\times D)}\\
&\qquad\qquad  \le
 \left[\mu(X)\right]^\half \left(\int_X \mu(dx)\, \Vert G(x,\cdot,\ast)\Vert_{L^2([0,T]\times D)}^2\right)^\half < \infty\quad  \text{a.s.,}
 \end{align*}
 and \eqref{ch1'.f1} holds.
\end{proof}
\medskip

\begin{remark}
\label{Fubini-det}
Theorem \ref{ch1'-tfubini} also holds in the case where the integrand $G$ is deterministic and does not depend on the variable $t$, and the stochastic integral is with respect to white noise on $\rek$.
\end{remark}


\section{Differentiation under the stochastic integral}
\label{ch2'new-s5}

In this section, we address the question of differentiability of the stochastic integral of a process $G(\lambda,s,y)$, that depends on a parameter $\lambda\in\re$, with respect to this parameter, and prove a formula for the derivative of the integral.
We let $W$, $(\cF_s,\, s\in[0,T])$ and $(e_j, \, j\ge 1)$ be as defined at the beginning of Section \ref{ch2new-s2}.

Let $I\subset \IR$ be a bounded open interval. Recall that a function $f: I\to \IR$ is {\em absolutely continuous} if there is a locally integrable function $g:I\to\IR$ such that, for all $a,b\in I$, $f(b) -f(a) = \int_a^b g(\lambda)\, d\lambda$. The function $g$ is often denoted by $\frac{df}{d\lambda}$. Later in this section, we will refer to $\frac{df}{d\lambda}$ as the ``derivative'' of $f$. An interesting sufficient condition for $f$ to be absolutely continuous is the following \cite[Chapter 3, Section 6, Problem 5]{stein_shak}: $f'(\lambda)$ exists for every $\lambda \in I$, and $f'$ (which is necessarily a Borel function) is locally (Lebesgue) integrable. In this case, for all $a,b\in I$, $f(b) - f(a) = \int_a^b f'(\lambda)\, d\lambda$.

Consider the following set of assumptions:
\smallskip

\noindent{\bf(H)}
\begin{itemize}
\item[(i)]  $G: I\times [0,T]\times D\times \Omega \rightarrow \re$ is $\cB_I\times \cB_{[0,T]} \times \cB_D\times \cF$-measurable and such that, for fixed $s\in [0,T]$, the partial function $(\lambda,y,\omega) \mapsto G(\lambda,s,y,\omega)$ from $I\times D\times \Omega$ into $\re$ is $\cB_I\times\mathcal{B}_D\times \tf_s$--measurable. Furthermore, for all $\lambda\in I$,
\beqn
\Vert G(\lambda,\cdot,\ast)\Vert_{L^2([0,T]\times D)}\,
 < \infty, \qquad \text{a.s.}
\eeqn
\item[(ii)] For $ds dy dP$-almost all $(s, y,\omega)\in[0,T]\times D\times \Omega$, the map
\beqn
\lambda \mapsto  G(\lambda, s, y;\omega)
\eeqn
 is absolutely continuous. Let us denote by $\lambda\mapsto \frac{\partial}{\partial\lambda}G(\lambda,s,y)$ its ``derivative''. This function is well defined  for a.a.~$\lambda$, where the ``a.a.'' depends on $(s,y,\omega)$.

Denote by $I_0(s,y,\omega)\in\cB_I$ the implied full measure set, and define
 \beqn
 \bar G(\lambda,s,y,\omega)=
 \begin{cases}
 \frac{\partial}{\partial\lambda}G(\lambda,s,y,\omega), & {\text{if}}\  \lambda\in I_0(s,y,\omega),\\
 0, & {\text{if}}\   \lambda \in I\setminus I_0(s,y,\omega).
 \end{cases}
 \eeqn
 On the $ds dy dP$-null set of points $(s,y,\omega)\in[0,T]\times D\times \Omega$ where $\lambda \mapsto  G(\lambda, s, y;\omega)$ fails to be absolutely continuous, we set
 $\bar G(\lambda,s,y)=0$ for all $\lambda\in I$.

\item[(iii)] The map $(\lambda,s,y,\omega)\mapsto \bar G(\lambda,s,y,\omega)$ is $\cB_I\times \cB_{[0,T]} \times \cB_D\times \cF$-measurable and such that, for fixed $s\in [0,T]$, the partial function $(\lambda,y,\omega) \mapsto \bar G(\lambda,s,y,\omega)$ from $I\times D\times \Omega$ into $\re$ is $\cB_I\times\mathcal{B}_D\times \tf_s$--measurable. Furthermore,
\beqn
\int_{I} d\lambda\, \Vert \bar G(\lambda, \cdot, \ast)\Vert_{L^2([0,T]\times D)} < \infty,\qquad  a.s.
\eeqn
\end{itemize}

\begin{thm}
\label{ch1'-tdif}
Let $(G(\lambda,s,y),\, (s,y)\in [0,T] \times D)$, $\lambda\in I$,  be a family of stochastic processes.
\begin{description}
\item{(1)} Suppose that assumptions {\bf(H)} above are satisfied. Then the process
\beq\label{rde2.5.1}
\left(\int_0^T \int_D G(\lambda, s, y)\, W(ds, dy),\ \lambda\in I\right)
\eeq
has a version $(H(\lambda),\, \lambda \in I)$ that is jointly measurable in $(\lambda,\omega)$
and such that a.s., $\lambda \mapsto H(\lambda)$ is absolutely continuous.

Further, the process  $\left(\int_0^T \int_D \bar G(\lambda,s,y)\, W(ds,dy),\ \lambda\in I\right)$ has a jointly measurable version in $(\lambda,\omega)$, that we denote by $(K(\lambda),\, \lambda\in I)$,
such that a.s., for a.a.~$\lambda \in I$,
\beq
\label{ch1'.dif1}
\frac{d}{d \lambda} H(\lambda)
= K(\lambda).
\eeq
\item{(2)} In addition to the assumptions {\bf(H)}, we assume:

(iv)\  The process $\left( \int_0^T \int_D \bar G(\lambda,s,y)\, W(ds,dy),\  \lambda\in I\right)$ has a version
$(\bar K(\lambda),\, \lambda \in I)$
with continuous sample paths.
\end{description}
 Then, a.s.,
$\lambda \mapsto H(\lambda)$ from part (1) is  continuously differentiable on $I$, and for all $\lambda \in I$,
\beq\label{rde2.5.4}
  H^\prime(\lambda) = \bar K(\lambda).
\eeq

The equalities \eqref{ch1'.dif1}  and \eqref{rde2.5.4} are informally written
$$
   \frac{d}{d \lambda} \int_0^T \int_D G(\lambda, s, y)\, W(ds, dy) = \int_0^T \int_D \frac{\partial}{\partial\lambda} G(\lambda, s, y)\, W(ds, dy).
$$
\end{thm}
\begin{proof}
The assumption $(\bf H)$(iii) tells us that the map $(\lambda,s,y,\omega)\mapsto \bar G(\lambda,s,y,\omega)$ satisfies the hypotheses of Theorem \ref{ch1'-tfubini} (with $G:=\bar G$ and $\mu$ the Lebesgue measure there).  Hence,
according to the assertion (b) of that theorem, for any $\lambda_1,\lambda_2 \in I$, the stochastic integral
\beqn
\int_0^T\int_D \left(\int_{\lambda_1}^{\lambda_2} d\lambda\, \bar G(\lambda, s, y)\right)  W(ds, dy)
\eeqn
is well-defined in the sense of Section \ref{ch2'new-s3}.

Furthermore, the assertion (a) of Theorem \ref{ch1'-tfubini} (with $G$ there replaced by $\bar G$) implies that the stochastic process
\beqn
\left( \int_0^T \int_D \bar G(\lambda,s,y)\, W(ds,dy),\ \lambda\in I\right)
\eeqn
has a jointly measurable (in $(\lambda,\omega)$) version, denoted by $(K(\lambda),\, \lambda\in I)$, and a.s., the map $\lambda\mapsto K(\lambda)$ is $d\lambda$-integrable.

The assumption $(\bf H)$(i) ensures that the stochastic integral
\beqn
\int_0^T \int_D \left(G(\lambda_2, s, y)-G(\lambda_1, s, y)\right)\, W(ds, dy),
\eeqn
is also well-defined in the sense of Section \ref{ch2'new-s3}, and from assumption
$(\bf H)$(ii), we deduce that
\begin{align}
\label{ch1'.dif2}
&\int_0^T \int_D \left(G(\lambda_2, s, y)-G(\lambda_1, s, y)\right)\, W(ds, dy)\notag\\
& \qquad  \qquad = \int_0^T\int_D \left(\int_{\lambda_1}^{\lambda_2} d\lambda\ \bar G(\lambda, s, y)\right)  W(ds, dy)\qquad \text{a.s.}
\end{align}

We apply the stochastic Fubini's theorem (Theorem \ref{ch1'-tfubini}), for fixed $\lambda_1, \lambda_2\in I$, and obtain
\begin{align}\nonumber
&\int_0^T \int_D G(\lambda_2, s, y)\, W(ds, dy) -  \int_0^T \int_D G(\lambda_1, s, y)\, W(ds, dy) \\
&\qquad  \qquad = \int_{\lambda_1}^{\lambda_2} d\lambda \int_0^T \int_D \bar G(\lambda, s, y)\, W(ds, dy)\quad {\text{a.s.}}
\label{rde2.5.5}
\end{align}

Fix $\lambda_1 \in I$ and for $\lambda \in I$, on the event where $\lambda \mapsto K(\lambda)$ is integrable, define
\beq
\label{def-H}
   H(\lambda) = \int_{\lambda_1}^{\lambda} d\tilde\lambda\, K(\tilde\lambda)
     + \int_0^T \int_D G(\lambda_1, s, y)\, W(ds, dy),
\eeq
and set $H(\lambda) = 0$ on the complement of this event. Then $\lambda \mapsto H(\lambda)$ is absolutely continuous, and the identity \eqref{rde2.5.5}  tells us that $(H(\lambda),\, \lambda \in I)$ is a (jointly measurable in $(\lambda,\omega)$) version of
\beq
\label{version-hlambda}
\left(\int_0^T \int_D G(\lambda, t, x)\, W(dt, dx),\quad \lambda \in I\right).
\eeq

Summarising the above discussion, we see that assuming $\bf{(H)}$, on the event where $\lambda \mapsto K(\lambda)$ is $d\lambda$-integrable on $I$, $\lambda \mapsto H(\lambda)$ is absolutely continuous and for a.a.~$\lambda$, $\frac{dH(\lambda)}{d\lambda} = K(\lambda)$. This proves \eqref{ch1'.dif1}, and completes the proof of (1).

 Assuming (iv), since $(\bar K(\lambda),\, \lambda \in I)$ has continuous sample
paths, it is necessarily jointly measurable (and $d\lambda$-integrable on compact intervals in $I$). Therefore, since for each $\lambda$, $K(\lambda) = \bar K(\lambda)$ a.s., by Fubini's theorem, a.s, for a.a.~$\lambda$, $K(\lambda) = \bar K(\lambda)$. In particular, if we replace $K(\tilde\lambda)$ in \eqref{def-H} by $\bar K(\tilde\lambda)$, we obtain a process $(\bar H(\lambda),\, \lambda\in I)$ that is indistinguishable from $(H(\lambda),\, \lambda\in I)$. This remains a jointly measurable version of the process \eqref{version-hlambda}. By the fundamental theorem of calculus, a.s., $\lambda \mapsto \bar H(\lambda)$ is continuously differentiable on $I$, and a.s., for all $\lambda \in I$, $\bar H'(\lambda) = \bar K(\lambda)$, proving \eqref{rde2.5.4} since $H$ and $\bar H$ are indistinguishable.
\end{proof}

\begin{remark} \label{rdrem2.5.2}
When the assumptions $\bf{(H)}$ are satisfied, a sufficient condition for condition (iv) is the following:
\smallskip

For each compact interval  $J\subset I$, for $dsdydP$-almost all $(s,y,\omega)\in[0,T]\times D\times\Omega$, the function $\lambda \mapsto \bar G(\lambda, s, y)$ is  continuously differentiable and
\beqn
 \sup_{\lambda \in J}\ E\left[\int_0^T ds \int_D dy\, \left(\frac{\partial \bar G}{\partial \lambda}(\lambda, s, y)\right)^2\right] < \infty.
\eeqn
Indeed, set
\beqn
 X(\lambda) = \int_0^T  \int_D \bar G(\lambda, s, y)\, W(ds,dy).
\eeqn
For $\lambda_1, \lambda_2\in J$ with $\lambda_1<\lambda_2$,
\begin{align*}
 &E\left[\left(X(\lambda_2)-X(\lambda_1)\right)^2\right]\\
 &\qquad = E\left[\int_0^T ds \int_D dy \left(\bar G(\lambda_2, s, y) - \bar G(\lambda_1, s, y)\right)^2 \right] \\
 & \qquad = E\left[\int_0^T ds \int_D dy \left(\int_{\lambda_1}^{\lambda_2} d\lambda\, \frac{\partial\bar G}{\partial \lambda}(\lambda, s, y)\right)^2\right]\\
 & \qquad = \int_{\lambda_1}^{\lambda_2} d\lambda \int_{\lambda_1}^{\lambda_2} d\tilde \lambda\,  E\left[\int_0^T ds \int_D dy \, \frac{\partial\bar G}{\partial \lambda}(\lambda, s, y)\,
 \frac{\partial\bar G}{\partial \tilde\lambda}(\tilde\lambda, s, y)\right] \\
 &\qquad  \le  \int_{\lambda_1}^{\lambda_2} d\lambda \int_{\lambda_1}^{\lambda_2} d\tilde \lambda\,  \sup_{\lambda \in J} E\left[\int_0^T ds \int_D dy \, \left(\frac{\partial \bar G}{\partial \lambda}(\lambda, s, y)\right)^2\right]\\
 & \qquad = C (\lambda_2-\lambda_1)^2,
\end{align*}
where we have applied the Cauchy-Schwarz inequality.  By Kolmogorov's continuity criterion (see Theorem \ref{ch1'-s7-t2}), the process $(X(\lambda),\, \lambda\in J)$ has a continuous version $(\bar K(\lambda),\, \lambda\in I)$. Therefore, condition (iv) is satisfied.
\end{remark}

\begin{remark}
\label{dif-det}
Theorem \ref{ch1'-tdif} also holds in the case where the integrand $G$ is deterministic and does not depend on the variable $t$, and the stochastic integration is with respect to white noise on $\rek$. In this case, the set of assumptions is:
\begin{itemize}
\item[($i^\prime$)] The function $(\lambda, y) \mapsto G(\lambda,y)$ from $I\times D$ into $\re$ is $\cB_I\times \cB_D$-measurable, and for all $\lambda\in I$,
\beqn
 \Vert G(\lambda,\ast)\Vert_{L^2(D)}\, < \infty.
\eeqn
\item[($ii^\prime$)] For $dy$-almost all $y\in D$, the mapping $\lambda \mapsto  G(\lambda,y)$ is absolutely continuous.
Let $\lambda\mapsto \frac{\partial}{\partial\lambda} G(\lambda,y)$ denote its ``derivative''. This function is well-defined for a.a.~$\lambda$, where the ``a.a.'' depends on $y$.

Denote by $I_0(y)\in\cB_I$ the implied set of full measure, and define
\beqn
\bar G(\lambda,y)=
\begin{cases}
\frac{\partial}{\partial\lambda} G(\lambda,y), & {\text{if}}\ \lambda\in I_0(y),\\
0, & {\text{if}}\ \lambda\in I\setminus I_0(y).
\end{cases}
\eeqn
On the $dy$-null set of points $y\in D$ where $\lambda \mapsto  G(\lambda,y)$ is not absolutely continuous, we set
 $\bar G(\lambda,y)=0$ for all $\lambda\in I$.

\item[($iii^\prime$)] The mapping $(\lambda,y)\mapsto \bar G(\lambda, y)$
 from $I\times D$ into $\re$ is $\cB_I\times \B_D$-measurable, and such that
\beqn
\int_{I} d\lambda\, \Vert \bar G(\lambda, \ast)\Vert_{L^2(D)} < \infty.
\eeqn
\item[($iv^\prime$)] The mapping $\lambda \mapsto \int_D \bar G(\lambda, y)\, W(dy)$ has a version with continuous sample paths.
\end{itemize}
\end{remark}

\section{Joint measurability of the stochastic integral}
\label{ch2'new-s6}

In this section, we investigate the question of joint measurability of the stochastic integral of a process $G(x,s,y)$ that depends on a parameter $x\in X$, where $(X,\cX)$ is a measure space. We let $T>0$, $W$, $(\cF_s,\ s\in[0,T])$ and $(e_j,\, j\ge 1)$ be as defined at the beginning of Section \ref{ch2new-s2}.

\begin{thm}
\label{ch1'-l-m1}
Let $(X,\cX)$ be a measure space. Consider a function $G:X \times [0,T]\times D\times \Omega \to \IR$ that is $\cX\times\cB_{[0,T]} \times \cB_D\times \cF$--measurable, and for fixed $(x,s) \in X\times [0,T]$, $(y,\omega) \mapsto G(x,s,y,\omega)$ is $\cB_D \times \cF_s$-measurable.
Suppose in addition that for each $x \in X$,
\beq\label{rde2.6.1}
   \int_0^T ds \int_D dy\, G^2(x,s,y) < \infty\qquad a.s.
\eeq
We have the following assertions:

(a) Fix $t \in [0,T]$. There is a function $H_t: X \times \Omega \to \IR$ that is $\cX\times \cF_t$--measurable and such that, for all $x \in X$,
$$
   H_t(x) = \int_0^t \int_D G(x,s,y)\, W(ds,dy), \qquad a.s.
$$
That is, $(H_t(x),\ x\in X)$ is a $\cX\times \cF_t$-measurable version of the process
\beqn
\left(\int_0^t \int_D G(x,s,y)\, W(ds,dy),\ x\in X\right).
\eeqn

(b) There is an $\cX\times \cB_{[0,T]}\times \cF_T$-measurable function $C:X\times [0,T]\times \Omega \to \IR$
such that, for all $x\in X$,  the process $C(x, \cdot)$ is indistinguishable from $G(x,\cdot,*)\cdot W$, and its sample paths are continuous.
\end{thm}
\begin{proof}
(a) By \eqref{rde2.6.1} and Proposition \ref{conv-probab}, for each $x\in \cX$, a.s. for all $t \in [0,T]$,
\beq\label{rde2.6.3}
\int_0^t \int_D G(x,s,y)\, W(ds,dy) = \sum_{j=1}^\infty \int_0^t \langle G(x,s,\ast) , e_j \rangle_{V}\, dW_s(e_j),
\eeq
where the series converges in probability, uniformly in $t\in[0,T].$ Since
$$
   \langle G(x,s,\ast) , e_j \rangle_{V} = \int_D G(x,s,y)  e_j (y)\, dy,
$$
this is an $\cX\times \cB_{[0,T]}\times \cF$-measurable function of $(x,s,\omega)\in X \times [0,T]\times \Omega$, and for fixed $(x,s)$, this is an $\cF_s$--measurable random variable, by hypothesis.
By Theorem \ref{rdthm9.4.1} applied to $Z(x, s, \omega) := \langle G(x,s, *), e_j \rangle_V$, each term in the series on the right-hand side of \eqref{rde2.6.3} is indistinguishable from a process $(x, t, \omega) \mapsto I_{t, j}(x, \omega)$ that is $\cX \times \cB_{[0, T]} \times \cF_T$-measurable and adapted (that is, for fixed $t \in [0, T]$, $(x, \omega) \mapsto I_{t, j}(x, \omega)$ is $\cX \times \cF_t$-measurable).
  The next goal is to partially extend this property to the series
$\sum_{j=1}^\infty I_{t,j}(x,\omega)$.

  For fixed $t\in[0,T]$ and $x \in X$,
$$
   \int_0^t \int_D G(x,s,y)\, W(ds,dy) = \sum_{j=1}^\infty I_{t,j}(x)\qquad \text{a.s.,}
$$
where the series converges in probability. By applying Lemma \ref{rdlem9.4.2a} to the sequence of partial sums $(\sum_{j=1}^n I_{t,j}(x),\, n \geq 1)$, $t$ fixed, there is a function $H_t: X\times \Omega \to \IR$ that is $\cX \times \cF_t$--measurable and such that for all $x \in X$,
$$
  H_t(x) = \sum_{j=1}^\infty I_{t,j}(x) = \int_0^t \int_D G(x,s,y)\, W(ds,dy)\qquad \text{a.s.}
$$
This proves (a).

(b) Considering $\omega \mapsto (t \mapsto I_{t,j}(x, \omega))$ as a random variable with values in $\cC([0, T])$, and \eqref{rde2.6.3} as an equality between $\cC([0, T])$-valued $\cF_T$-measurable random variables, we use Lemma \ref{rdlem9.4.2a} to obtain a function $(x, \omega) \mapsto (t \mapsto C(x, t, \omega))$ that is $\cX \times \cF_T$-measurable with values in $\cC([0, T])$ and such that for all $x \in X$, a.s, $t \mapsto C(x, t)$ is equal to $G(x, \cdot, *) \cdot W$. This is the desired function $C$.

The proof of Theorem \ref{ch1'-l-m1} is complete.
\end{proof}


\medskip

In Chapter \ref{ch1'-s5}, we will need joint measurability properties of stochastic integrals of the form \eqref{rd2.2.12a}. We consider here a slightly more general integrand
\begin{align*}
   U:[0,T] \times D \times [0,T] \times D \times \Omega &\to \IR \\
   (t,x,s,y,\omega) &\mapsto U(t,x,s,y,\omega)
\end{align*}
satisfying the following assumptions:
\begin{itemize}
   \item[(1)] the function $U$ is $\cB_{[0,T]} \times \cB_D \times \cB_{[0,T]} \times \cB_D \times \cF$-measurable;
   \item[(2)] for fixed $(t,x,s) \in [0,T] \times D\times [0,T]$,  $(y,\omega) \mapsto U(t,x,s,y,\omega)$  is \break
   $\cB_D \times \cF_s$-measurable;
   \item[(3)] for all $(t,x) \in [0,T]\times D$,
$$
   \int_0^T ds \int_D dy \, U^2(t,x,s,y)\, 1_{[0,t]}(s)  < \infty \qquad \text{a.s.}
$$
\end{itemize}

In the next proposition, we mention the notion of optional $\sigma$-field $\cO$ on $[0,T]\times \Omega$ associated with the filtration $(\cF_t,\, t \in [0,T])$. We refer the reader to Appendix \ref{app1}, Section \ref{rdsecA.4} for its definition.

\begin{prop}
\label{rdprop2.6.2}
 Under assumptions (1)--(3) above, for each $(t,x) \in [0,T]\times D$, the stochastic integral
\begin{align}
\label{rde2.6.3a}
 I(t,x) &= \int_0^t \int_D U(t,x,s,y)\, W(ds,dy)\notag\\
  &:= \int_0^T \int_D U(t,x,s,y) 1_{[0,t]}(s)\, W(ds,dy)
\end{align}
is well-defined in the sense of Section \ref{ch2'new-s3}. In addition, $I=(I(t,x),\, (t,x) \in [0,T]\times D)$ has a jointly measurable and adapted modification, that is, there is a function $Y: [0,T]\times D\times \Omega \to \IR$ such that:
\begin{itemize}
\item[(a)] $(t,x,\omega)\mapsto Y(t,x,\omega)$ is $\cB_{[0,T]}\times \cB_D \times \cF_T$-measurable;
\item[(b)] for fixed $t\in [0,T]$, $(x,\omega) \mapsto Y(t,x,\omega)$ is $\cB_D\times \cF_t$-measurable;
\item[(c)] for each $(t,x) \in [0,T]\times D$, $Y(t,x) = I(t,x)$ a.s.
\end{itemize}
This modification can in fact be chosen so that $(x,t,\omega)\mapsto Y(t,x,\omega)$ is $\cB_D \times \cO$-measurable.
\end{prop}
\begin{proof} By assumptions (1)--(3), for fixed $(t,x) \in [0,T]\times D$,
$$
(s,y,\omega) \mapsto U(t,x,s,y,\omega)\, 1_{[0,t]}(s)
$$
satisfies the conditions (1)--(2) at the beginning of Section \ref{ch2new-s2}, as well as \eqref{weakint}, therefore the stochastic integral in \eqref{rde2.6.3a} and $I(t,x)$ are well-defined.

According to Theorem \ref{ch1'-l-m1} (a) with $t := T$ and $(X,\cX) := ([0,T] \times D, \cB_{0,T]} \times \cB_D)$, there is a function $\tilde I : [0,T] \times D \times \Omega \to \R$ that is $\cB_{[0,T]} \times \cB_D \times \cF_T$-measurable and such that, for all $(t,x) \in [0,T] \times D$, $\tilde I(t,x) = I(t,x)$ a.s.

For fixed $(s,x)$,  the random variable $\omega \mapsto I(s,x,\omega)$ is $\cF_s$-measurable, so the same is true of $\omega \mapsto \tilde I(s,x,\omega)$, since $\cF_s$ is complete. Therefore, $\tilde I$ satisfies the assumptions of Lemma \ref{rdlem9.4.2} with $(X,\cX) := (D, \B_D)$ and $Z (x,s,\omega) := \tilde I(s,x,\omega)$. Denote by $(x,s,\omega) \mapsto Y(s,x,\omega)$ the $\cB_D \times \cO$-measurable function $K $ given in Lemma \ref{rdlem9.4.2}. By the conclusion of Lemma \ref{rdlem9.4.2} (a), we see that $(x,t,\omega) \mapsto Y(t,x,\omega)$ satisfies the claim (c) of the proposition, and $(x,t,\omega) \mapsto Y(t,x,\omega)$ is $\B_D \times \cO$-measurable. Therefore the assertions (a) and (b) also hold (because for all $t\in [0,T], \cO\vert_{[0,t] \times \Omega} \subset \cB_{[0,t]} \times \cF_t)$.
\end{proof}

We also want to study joint measurability and adaptedness properties of integrals
with respect to Lebesgue measure of stochastic processes $U$ as above.
For this, we consider the following variations on assumptions (2) and (3):
\begin{itemize}
   \item[(2')] for fixed $(t,x,s,y) \in [0,T] \times D\times [0,T]\times D$, the function $\omega \mapsto U(t,x,s,y,\omega)$  is $\cF_s$-measurable;
   \item[(3')] for all $(t,x) \in [0,T]\times D$,
$$
   \int_0^T ds \int_D dy \, \vert U(t,x,s,y)\vert\, 1_{[0,t]}(s)  < \infty \qquad \text{a.s.}
$$
\end{itemize}

\begin{prop}\label{rdprop2.6.4}
Under assumptions {\rm (1), (2') and (3'),} for each $(t,x) \in [0,T]\times D$,
\begin{align}
\label{rde2.6.4}
   J(t,x) &= \int_0^t ds \int_D dy \,  U(t,x,s,y)\notag \\
   &:= \int_0^T ds \int_D dy \,  U(t,x,s,y)\, 1_{[0,t]}(s)
\end{align}
is a well-defined random variable. In addition, $J=(J(t,x),\, (t,x) \in [0,T]\times D)$ has a jointly measurable and adapted modification, that is, there is a function $Y:[0,T]\times D\times \Omega \to \IR$ that satisfies properties (a)--(c) of Proposition \ref{rdprop2.6.2} (with $I$ replaced by $J$ in part (c)). This modification can further be chosen to be optional, that is, $(x,t,\omega)\mapsto Y(t,x,\omega)$ is $\cB_D\times \cO$-measurable.
\end{prop}

\begin{remark} If {\rm (2')} were replaced by ``for fixed $t \in [0,T]$, $(x,s,y,\omega) \mapsto U(t,x;s,y;\omega)$ from $D\times [0,t]\times D\times \Omega$ is $\cB_D\times \cB_{[0,t]}\times \cB_D\times \cF_t$-measurable," then it would be clear that property (b) of Proposition \ref{rdprop2.6.2} is satisfied by $J$. However, {\rm (2')} is not a ``progressive measurability" type of condition.
\end{remark}

\noindent{\em Proof of Proposition \ref{rdprop2.6.4}.} By (1) and (3'), for all $(t,x) \in [0,T]\times D$, the Lebesgue integral in \eqref{rde2.6.4} and $J(t,x,\omega)$ are well-defined a.s. Let $(X,\cX) = ([0,T]\times D\times D, \cB_{[0,T]}\times \cB_D\times \cB_D)$, and for the process $(t,x,y,s,\omega)\mapsto U(t,x,s,y,\omega)1_{[0,t]}(s)$, we denote by  $H(t,x,y,s,\omega)$ the function resulting from Lemma \ref{rdlem9.4.2} (a). According to this lemma, this function is $\cB_{[0,T]}\times \cB_D\times \cB_D \times \cO$-measurable. Define
\beq\label{rde2.6.5}
   \tilde J(t,x,\omega) = \int_0^T ds \int_D dy \,  H(t,x,y,s,\omega).
\eeq
Observe that:

   (i) $(t,x,\omega) \mapsto \tilde J(t,x,\omega)$ is $\cB_{[0,T]}\times \cB_D\times \cF_T$-measurable (by Fubini's theorem);

   (ii) for fixed $t\in [0,T]$, $(x,\omega) \mapsto \tilde J(t,x,\omega)$ is $\cB_D\times \cF_t$-measurable (because $\cO\vert_{[0,t]\times\Omega} \subset \cB_{[0,t]} \times \cF_t$);

   (iii) for fixed $(t,x)\in [0,T]\times D$, $\tilde J(t,x) = J(t,x)$ a.s. Indeed,
by\eqref{rde9.4.3} in Lemma \ref{rdlem9.4.2}, for fixed $(t,x) \in [0,T]\times D$,
$$
 \{(s,y,\omega) \in [0,T]\times D\times\Omega: H(t,x;y,s,\omega) \neq U(t,x,s,y,\omega)1_{[0,t]}(s) \}
$$
is a $dsdydP$-null set in $\cB_{[0,T]}\times \cB_D\times \cF_T$. Therefore, by Fubini's theorem, for a.a.~$\omega\in\Omega$,
$$
   H(t,x,\ast,\cdot,\omega) = U(t,x,\cdot,\ast,\omega) 1_{[0,t]}(\cdot)\qquad dsdy\text{-a.e.},
$$
therefore, by \eqref{rde2.6.5},
$$
  \tilde J(t,x) = \int_0^T ds \int_D dy \,  U(t,x,s,y)\, 1_{[0,t]}(s) = J(t,x)\qquad\text{a.s.}
$$
From (i), (ii) and (iii), we see that the properties (a)--(c) of Proposition \ref{rdprop2.6.2} hold (with $Y$ there replaced by $ \tilde J$ and $I$ there replaced by $J$).

Notice that $(x,t,\omega) \mapsto \tilde J(t,x,\omega)$ satisfies the hypotheses of Lemma \ref{rdlem9.4.2} with $(X,\cX) = (D,\cB_D)$. Denote by $Y$ the $\cB_D\times \cO$-measurable version of $\tilde J$ given by Lemma \ref{rdlem9.4.2} (a). Then $Y$ clearly also satisfies the properties (a)--(c) of Proposition \ref{rdprop2.6.2}.
\hfill $\Box$


\section{The Girsanov theorem for space-time white noise}
\label{ch2'-s3}

Let $W$, $\left({\mathcal{F}}_s,\, s\in [0,T]\right)$ and $(e_j,\, j\ge 1)$ be defined as at the beginning of Section \ref{ch2new-s2}.



Let $\left(h(t,x),\, (t,x)\in[0,T]\times D\right)$ be a jointly measurable and adapted random field such that
\beq
\label{ch2'-s3.2}
\int_0^T dt \int_D dx\ h^2(t,x) < \infty\quad {\text{a.s.}}
\eeq
Define a measure $\tilde P$ on $(\Omega,\mathcal{F})$ that is absolutely continuous with respect to $P$ and with density given by
\beq
\label{ch2'-s3.200}
\frac{d\tilde P}{dP} = \exp\left(-\int_0^T  \int_D \ h(t,x)\, W(dt,dx) - \frac{1}{2}\int_0^T dt \int_D dx\, h^2(t,x)\right),
\eeq
where the stochastic integral is defined in Section \ref{ch2'new-s3}.

The following is a version of {\em Girsanov's theorem}.\index{Girsanov's theorem}\index{theorem!Girsanov's} We refer the reader to \cite[p.~253]{lr} for a similar result.
We recall that in the context of measure theory, two measures are (mutually) equivalent\index{mutually equivalent}\index{equivalent!mutually}\index{measure!equivalent} if they have the same null sets. Clearly, $\tilde P$ and $P$ are mutually equivalent, since the density defined in \eqref{ch2'-s3.200} is strictly positive a.s.

\begin{thm}
\label{ch2'-s3-t1} Assume that $\tilde P$ is a probability measure.
Define a set function
$\tilde W: \mathcal{B}^f_{[0,T]\times D} \to L^2(\Omega, \tf,P)$ by
\beq
\label{ch2'-s3.2a}
\tilde W(A) = W(A) + \int_0^T dt \int_D dx\, 1_A(t,x)\, h(t,x).
\eeq
Then under $\tilde P$, $\left(\tilde W(A),\, A\in\mathcal{B}^f_{[0,T]\times D}\right)$ is a space-time white noise on $[0,T]\times D$ that satisfies conditions (i) and (ii) of Section \ref{ch2new-s1}.
Further, under $P$, the laws of $W$ and $\tilde W$ are mutually equivalent.
\end{thm}
\begin{proof}
By Proposition \ref{conv-probab}, the process
\beq
\label{ch2'-s3.2ac}
Y_t := \int_0^t \int_D h(s,y)\, W(ds,dy) = \sum_{j=1}^\infty \int_0^t\langle h(s, \ast),e_j\rangle_V\, dW_s(e_j) , \ t\in[0,T]
\eeq
(where $\langle \cdot,\cdot\rangle_V$ denotes the inner product in $V=L^2(D)$),
is well-defined. Moreover, it is a local martingale with continuous sample paths a.s., and with  quadratic variation
\beq
\label{ch2'-s3.2aa}
\langle Y\rangle_t = \int_0^t ds \, \Vert h(s, \ast)\Vert^2_{L^2(D)} = \int_0^t ds \int_D dx \, h^2(s,x), \quad t\in[0,T]
\eeq
(see \eqref{cucuvar}).

Set
\beq
\label{ch2'-s3.2b}
D_t = \exp\left(-Y_t-\frac{1}{2}\langle Y\rangle_t\right).
\eeq
Then $(D_t,\ t\in[0,T])$ is a nonnegative local martingale, hence a supermartingale. Since $\tilde P$ is a probability measure, $E[D_T] = E\left[\frac{d\tilde P}{dP}\right] =1$ and therefore, $(D_t,\, t\in[0,T])$ is in fact a martingale.

Let $\tilde P_t$ (respectively, $P_t$) denote the restriction of $\tilde P$ (respectively, $P$) to $\cF_t$.
Observe that $\frac{d\tilde P_t}{dP_t}=D_t$. According to \cite[(1.7) Theorem on p.~329]{ry}, for each $j\in \IN^*$,
\beq
\label{ch2'-s3.3g}
\tilde W_j(t) := W_t(e_j) + \langle W_{\cdot}(e_j),Y\rangle_t,\quad t\in[0,T],
\eeq
is, under $\tilde P$, a continuous local martingale relative to $(\tf_s,\, s\in[0,T])$.

Since the stochastic processes $(W_t(e_j),\ 0\le t\le T)$, $j\ge 1$, are mutually independent,
\beqn
\langle W_{\cdot}(e_j), Y\rangle_t = \int_0^t \langle h(s,\ast),e_j\rangle_V\, ds,\quad {\text{a.s.}}
\eeqn
Therefore,
for any $t\in[0,T]$ and $j\in\N^*$, we have
\begin{align}
\label{ch2'-s3.3g-bis}
\tilde W_j(t)
 = W_t(e_j) + \int_0^t  \langle h(s,\ast), e_j\rangle_V\, ds.
\end{align}

 Next, we prove that for each $n\in \IN^*$, the stochastic process
\beq
\label{ch2'-s3.3a}
\left((\tilde W_1(t), \ldots, \tilde W_n(t)),\ t\in[0,T]\right)
\eeq
is, under $\tilde P$, a standard $n$-dimensional $(\tf_t)$-Brownian motion.

Indeed, by \eqref{ch2'-s3.3g}, $P$-a.s., for all $j, \ell \ge 1$,
\beq
\label{G(*1)}
\langle \tilde W_{j}(\cdot), \tilde W_{\ell}(\cdot)\rangle_t=\langle W_{\cdot}(e_j), W_{\cdot}(e_l)\rangle_t=\delta_{j}^\ell\, t, \quad t\in[0,T].
\eeq
Since $P$ is equivalent to $\tilde P$, this also holds $\tilde P$- a.s. Because under $\tilde P$, \eqref{ch2'-s3.3a} is a continuous $(\tf_t)$-local martingale,  by L\'evy's characterisation of Brownian motion
(\cite[(36), Theorem p.~150]{ry}), it
is a standard $(\tf_t)$-Brownian motion under $\tilde P$.

By the definition of the isonormal process associated to the space-time white noise $W$, and  Remark \ref{ch2(*A)}, we see that, for any $A\in \mathcal{B}^f_{[0,T]\times D}$,
\begin{align*}
&W(A) + \int_0^T dt \int_D dx\, 1_A(t,x)h(t,x)\\
&\quad = W(1_A) + \int_0^T dt\, \langle 1_A(t,*),h(t,*)\rangle_V\\
&\quad = \sum_{j=1}^\infty\int_0^T \langle 1_A(t,*),e_j\rangle_V\left(dW_t(e_j) + \langle h(t,*),e_j\rangle_V\, dt\right)\\
&\quad = \sum_{j=1}^\infty\int_0^T \langle 1_A(t,*),e_j\rangle_V\, d\tilde W_j(t),
\end{align*}
where we have used Parseval's identity and \eqref{ch2'-s3.3g-bis}.
Thus, from \eqref{ch2'-s3.2a}, we obtain
\beq
\label{ch2-G(*2)}
\tilde W(A) = \sum_{j=1}^\infty\int_0^T \langle 1_A(t,*),e_j\rangle_V\, d\tilde W_j(t).
\eeq
It follows that under $\tilde P$, the process $(\tilde W(A),\, A\in \mathcal{B}^f_{[0,T]\times D})$ is Gaussian, with mean zero, and using \eqref{G(*1)},
for any $A, B \in \mathcal{B}^f_{[0,T]\times D}$,
\begin{align*}
E[\tilde W(A)\tilde W(B)] &= \sum_{j=1}^\infty \int_0^T dt\, \langle 1_A(t,*),e_j\rangle_V\, \langle 1_B(t,*),e_j\rangle_V\\
&= \int_0^T dt\, \langle 1_A(t,*), 1_B(t,*)\rangle_V = |A\cap B|.
\end{align*}
Hence, under $\tilde P$, the process $(\tilde W(A),\, A\in \mathcal{B}^f_{[0,T]\times D})$ is a space-time white noise.

We note for future reference that by \eqref{ch2-G(*2)} and Lemma \ref{ch1'-lsi} (2), the isonormal process associated to $\tilde W$ under $\tilde P$ is $(\tilde W(h),\, h \in L^2([0,T] \times D)$, where
\beqn
     \tilde W(h) = \sum_{j=1}^\infty \int_0^T \langle h(s, *), e_j\rangle_V\, d\tilde W_j(t).
     \eeqn
In particular,
\beq
\label{ch2-G(*3)}
   \tilde W_t (e_j) = \tilde W(1_{[0,t]}(\cdot) e_j(*)) = \tilde W_j(t).
   \eeq
   Since for all $n \in \N$, $(\tilde W_j(\cdot),\, 1 \leq j \leq n)$ is an $n$-dimensional $(\cF_t)$-Brownian motion, it also follows from Lemma \ref{ch1'-lsi} (2) that conditions (i) and (ii) of Section \ref{ch2new-s1} are satisfied.

Because $\tilde P$ and $P$ are equivalent, the laws of $\tilde W$ under $\tilde P$ and under $P$ are equivalent, and since the law of $\tilde W$ under $\tilde P$ is the same as the law of $W$ under $P$, the last statement of the theorem follows.
 \end{proof}

\noindent{\em Multidimensional version of Girsanov's theorem}.
\medskip

Assume that $W=(W^1, \ldots, W^d)$ is a $d$-dimensional space-time white noise, that is, the components $W^i$, $1\le i\le d$, are mutually independent space-time white noises that satisfy conditions (i) and (ii) of Section \ref{ch2new-s1}. Let $h= \left(h(t,x),\, (t,x)\in[0,T]\times D\right)$ be an $\red$-valued jointly measurable and adapted random field such that
\beq
\label{unmes(*1)}
\int_0^T dt \int_D dx\, \vert h(t,x)\vert^2 < \infty \quad {\text{a.s.}}
\eeq

Define a measure $\tilde P$ by
\beqn
\frac{d\tilde P}{dP} = \exp\left(-\int_0^T  \int_D \ h(t,x)\cdot W(dt,dx) - \frac{1}{2}\int_0^T dt \int_D dx\,  \vert h(t,x)\vert^2\right),
\eeqn
where, recalling the Euclidean inner product, $h(t,x)\cdot W(dt,dx)$ denotes
\beqn
\sum_{i=1}^d h^i(t,x)\ W^i(dt,dx)
\eeqn
 The proof of Theorem \ref{ch2'-s3-t1} can be easily adapted to obtain the following multidimensional version of Girsanov's theorem.\index{Girsanov's theorem!multidimensional}
\begin{thm}
\label{ch2'-s3-t1-md}
Assume that $\tilde P$ is a probability measure. For each $1\le i\le d$,  define a set function
$\tilde W^i: \mathcal{B}^f_{[0,T]\times D} \to L^2(\Omega, \tf,P)$ by
\beq
\label{ch2'-s3.2abis}
\tilde W^i(A) = W^i(A) + \int_0^T dt \int_D dx\, 1_A(t,x)\, h^i(t,x).
\eeq
Then under $\tilde P$, the process
\beqn
\left(\tilde W(A)=(W^1(A), \ldots, W^d(A)),\ A\in\mathcal{B}^f_{[0,T]\times D}\right)
\eeqn
is a $d$-dimensional space-time white noise that satisfies conditions (i) and (ii) of Section \ref{stwn}. Furthermore, under $P$, the laws of $W$ and $\tilde W$ are equivalent.
\end{thm}

For a fixed integrand process, the stochastic integral with respect to $\tilde W$ can be derived from the stochastic integral with respect to $W$. In the next remark, we identify the correction term. We use the notation $E_P$ to emphasise the probability measure used in the computation of the mathematical expectation.
 \begin{prop}
\label{ch2'-s3-r1}
Assume the hypotheses of Theorem \ref{ch2'-s3-t1}.
 Let
 \beqn
 (G(t,x),\ (t,x)\in[0,T]\times D)
 \eeqn
  be a jointly measurable and adapted stochastic process such that
\beqn
\int_0^T dt \int_D dx\ G^2(t,x) < \infty,\quad \tilde P{\text{-a.s.}}
\eeqn
Then
\begin{align}
\label{ch2'-s3.integral}
\int_0^t \int_D G(s,y)\, \tilde W(ds,dy) &= \int_0^t \int_D G(s,y)\, W(ds,dy)\notag\\
&\qquad + \int_0^t ds \int_D dy\ G(s,y)\, h(s,y).
\end{align}
\end{prop}
\begin{proof}
By definition,
\beqn
\int_0^t \int_D G(s,y)\, \tilde W(ds,dy) = \sum_{j=1}^\infty \int_0^t \langle G(s), e_j\rangle_V\, d\tilde W_s(e_j), \ \quad \tilde P{\text{-a.s.,}}
\eeqn
where the series converges in probability. By \eqref{ch2-G(*3)}, $\tilde W_t(e_j)=\tilde W_j(t)$, with $\tilde W_j(t)$ given in \eqref{ch2'-s3.3g}. Consequently,
\begin{align*}
&\sum_{j=1}^\infty \int_0^t \langle G(s), e_j\rangle_V\, d\tilde W_s(e_j)\\
&\qquad = \sum_{j=1}^\infty \int_0^t \langle G(s), e_j\rangle_V \left[dW_s(e_j)+d\langle W_\cdot(e_j),Y\rangle_s\right]\\
&\qquad = \sum_{j=1}^\infty \int_0^t \langle G(s), e_j\rangle_V\, dW_s(e_j) + \sum_{j=1}^\infty \int_0^t \langle G(s), e_j\rangle_V
\langle h(s),e_j\rangle_V\, ds\\
&\qquad = \int_0^t \int_D G(s,y)\, W(ds,dy) + \int_0^t ds \int_D dy\, G(s,y)\, h(s,y).
\end{align*}
\end{proof}
\bigskip

\noindent{\em Sufficient conditions for $\tilde P$ to be a probability measure}
\medskip

Next, we address the issue of finding sufficient conditions for $\tilde P$ to be a probability measure on $(\Omega, \mathcal{F})$. According to \eqref{ch2'-s3.200}, this will be the case if and only if
\beq
\label{ch2'-s3.4}
E\left[\frac{d\tilde P}{dP}\right] = 1.
\eeq
\begin{prop}
\label{ch2'-s3-p1}
Let $h$ be as in \eqref{unmes(*1)}, set
\beqn
 Y_t = \int_0^t \int_D h(s,y) \cdot W(ds, dy),\quad
   D_t = \exp\left(- Y_t - \half\, \langle Y \rangle_t\right).
   \eeqn
The following are sufficient conditions for \eqref{ch2'-s3.4} to hold.
\begin{description}
\item{(a)} {\em Kazamaki's criterion.}\index{Kazamaki's criterion}\index{criterion!Kazamaki's} If $\left(\exp\left(\frac{1}{2}\,Y_t\right),\, t\in[0,T]\right)$ is a uniformly integrable submartingale, then $\left(D_t,\, t\in[0,T]\right)$ is a uniformly integrable martingale. In particular \eqref{ch2'-s3.4} holds.
\item{(b)}  {\em Novikov's condition.}\index{Novikov's criterion}\index{criterion!Novikov's} If
\beq
\label{ch2'-s3.5}
E\left[\exp\left(\frac{1}{2}\int_0^T dt \int_D dx\, \vert h(t,x)\vert^2\right)\right] < \infty,
\eeq
then $\left(D_t,\, t\in[0,T]\right)$ is a martingale and therefore, \eqref{ch2'-s3.4} holds.
\item{(c)} If there is a partition $0=t_0<t_1<\cdots <t_n=T$ such that
\beqn
E\left[\exp\left(\frac{1}{2}\int_{t_{k-1}}^{t_k} dt \int_D dx\, \vert h(t,x)\vert^2\right)\right] < \infty,\ k=1,\dots,n,
\eeqn
then $(D_t,\, t\in[0,T])$ is a martingale and \eqref{ch2'-s3.4} holds.
\item{(d)} Assume there exist constants $\epsilon>0$ and $C<\infty$ such that
\beqn
\sup_{s\in[0,T]}\exp\left(\epsilon\int_D \vert h(t,x)\vert^2\, dx\right) \le C, \  \text{a.s.}
\eeqn
Then $(D_t,\, t\in[0,T])$ is a martingale and \eqref{ch2'-s3.4} holds.
\end{description}
\end{prop}

\begin{proof}
The proof of (a) and (b) can be found in \cite[(1.14) Proposition p. 331 and (1.16) Corollary, p.~333]{ry}, and the proof of (c) in \cite[5.14, Corollary p.~199]{ks}.

Condition (d) appears in \cite{gp} (see also  \cite[Vol 1. p. 233, Example 3]{ls}).
It is stronger than condition (c). Indeed, if condition (d) holds, then since $x \mapsto \exp(x)$ is increasing, any sufficiently fine partition can be used in condition (c).
\end{proof}

\section{Notes on Chapter \ref{ch2}}
\label{ch2-notes}


The stochastic integral with respect to a cylindrical Brownian motion was developed in particular in \cite{metivier-pellaumail}. Later authors, including \cite{walsh} and \cite{dz}, were aware of this  but preferred other approaches. Many authors used series of one-dimensional It\^o integrals to define stochastic integrals with respect to infinite-dimensional Wiener processes (see for instance  \cite[Definition 3.3.2]{lototsky-rozovsky-2017} and \cite{pardoux2021}). However, Nualart and Quer-Sardanyons \cite{nualartquer} seem to be the first to have observed explicitly that an integral with respect to a certain cylindrical Wiener process provides an extension of the Walsh-Dalang integral introduced in \cite{walsh} and \cite{dalang}. In particular,  the stochastic integral introduced in \cite{dalang} was limited to spatially homogeneous integrands. The presentation of \cite{nualartquer}  was simplified in \cite{dalang-quer}, where it was noticed that the series representation \eqref{ch1'-s4.1} provides a simplification of the presentation of \cite{nualartquer} as well a a unification of the theories presented in \cite{walsh} and \cite{dz}, when the Hilbert space $V = L^2(D)$ is replaced by a more general Hilbert space and the integrands are extended to a Hilbert-space valued setting as introduced in \cite{quer-sanz2004}. Another advantage of Definition \ref{ch1'-s4.d1} is the easy transfer of properties of the classical It\^o integral  to similar properties for the infinite-dimensional integral discussed in this chapter.

In the classical Itô theory of stochastic integrals with respect to continuous martingales, the integrands are often taken to be progressively measurable. However, since we only use Brownian motions as integrators , we have relaxed this condition to ``jointly measurable and adapted'' (see the end of Section \ref{ch2'new-s3}).

A Fubini's theorem for stochastic integrals with respect to martingale measures is given in \cite[p.~296]{walsh}. For stochastic integrals with respect to Hilbert-space-valued Wiener processes, and under weaker assumptions than \cite{walsh}, two versions are available: \cite[Section 4.6, Theorem 4.18]{dz} and \cite{neerven-veraar-2006}, respectively. Theorem \ref{ch1'-tfubini} has the same type of assumptions as these last two references. Its proof relies on Theorem \ref{app1-2.t1}, a Fubini's theorem for stochastic integrals with respect to standard Brownian motion (see \cite[Lemma 2.6]{krylov-2011}  and the more general statement in \cite[Theorem 2.2]{V2012}).

In many contexts, it is useful to differentiate under the stochastic integral. For It\^o integrals, results that justify this operation appear in many references. In Section \ref{ch2'new-s5}, we present analogous statements for stochastic integrals with respect to space-time white noise. Similarly, measurable dependence of stochastic integrals with respect to semimartingales in a important topic, and a systemiatic study of this property was undertaken by Stricker and Yor in \cite{strickeryor}. In Section \ref{ch2'new-s6}, we present some of their results in the context of stochastic integrals with respect to space-time white noise.


In general terms, Girsanov theorems address the issue of absolute continuity of the measure on an abstract Wiener space obtained by a shift transformation of the Wiener measure. When the shift is defined by a smooth deterministic function, this question was studied in \cite{c-m1944} and the result is called the Cameron-Martin theorem. The first result for stochastic shifts appears in \cite{girsanov-1960}. Extensions of these early works can be found in \cite{u-z2000}, \cite[Chapter 4]{nualart-1995} and in many other references.  In the setting of abstract Wiener space, such theorems are in the core of the development of Malliavin calculus (see for instance \cite{malliavin-1978} and \cite{bismut-1984}).
Girsanov-type theorems are fundamental in the study of stochastic differential and stochastic partial differential equations with additive noise, particularly in the study of weak solutions to such equations.  These questions will be developed in Section \ref{ch2'-s1}.


\chapter{Linear SPDEs driven by space-time white noise}
\label{chapter1'}

\pagestyle{myheadings}
\markboth{R.C.~Dalang and M.~Sanz-Sol\'e}{SPDEs driven by space-time white noise}

This chapter initiates the study of SPDEs in an elementary setting. We consider a space-time white noise as a random forcing, and we mostly restrict the spatial dimension to $k=1$.
We start by introducing two notions of solution: random field solutions and weak solutions. Although in this book, we mostly emphasize the former notion, the latter is also widely present in the theory of PDEs and of SPDEs. We will then consider SPDEs with a linear differential operator driven by additive noise.
We study two fundamental examples, namely the stochastic heat and wave equations in several different settings (on the real line, on finite intervals, etc.) and we prove sharp regularity properties of their sample paths.


\section{Notions of solution}
\label{ch1'-s1}

In this section, we introduce some notions of solution that are commonly used in the theory of PDEs, and we will discuss how they can be adapted to the framework of SPDEs driven by a space-time white noise. We consider two cases, namely $D=\rek$ and $D\subset \R^k$, a domain (that is, a non-empty open connected subset) with smooth boundaries. When $D$ is bounded, we denote by $\bar D$ and $\partial D$ the closure and the boundary of $D$, respectively.
\medskip

\noindent{\em PDEs on $\rek$: the classical case}
\smallskip

Let $k\ge 1$ and $\mathcal{L}$ be a linear partial differential operator on $\re_+\times\rek$, possibly with non constant coefficients. For a partial differential equation defined by $\mathcal{L}$, we have outlined in Section \ref{ch1-s3} the notion of {\em fundamental solution}.\index{fundamental solution}\index{solution!fundamental} In the classical case,
which involves smooth functions, we give some illustrations of this notion.

Let $\mathcal{L}= \frac{\partial}{\partial t} + \mA$, where $\mA$ is a partial differential operator in the variable $x$. The {\em fundamental solution} associated to $\mathcal{L}$ is a Borel function $\Gamma(t,x;s,y)$ defined for all $(t,x)$, $(s,y)$,
$0\le s<t\le T$, $x,y\in \rek$, such that $\mathcal{L}\ \Gamma(t,x;s,y)=0$ for $(t,x)\in\ ]s,T]\times \rek$, and
\beq
\label{fs}
\lim_{t\downarrow s}\Gamma(t,x; s,y) = \delta_{0}(x-y)
\eeq
(see e.g. \cite[p.~182]{ez}).

Consider the parabolic Cauchy problem:
\beqn
\begin{cases}
\mathcal{L}u(t,x) & = f(t,x),\qquad  (t,x)\in\, ]0,T]\times \rek, \notag\\
u(0,x) & = \psi(x), \qquad\quad x\in \rek,
\end{cases}
\eeqn
for a given function $f$ and initial condition $\psi$. Under suitable conditions, a solution is given by the formula
\beq
\label{4.2-1'}
u(t,x) = \int_{\rek} \Gamma(t,x; 0,y) \psi(y)\ dy + \int_0^t ds \int_{\rek} dy\ \Gamma(t,x; s,y) f(s,y),
\eeq
$(t,x)\in\, ]0,T]\times \rek$
(see e.g. \cite[pp. 141 and 142]{friedman}, \cite[pp. 205, Theorem VI.13]{ez}).
It is implicitely assumed that the integrals on the right-hand side of \eqref{4.2-1'} are well-defined.
In particular, for all $\varphi \in \cC_0^\infty(]0, \infty[ \times \R^k)$,
\beqn
   \cL\left(\int_0^\cdot  \int_{\R^k}  \Gamma(\cdot, *; s, y)\, \varphi(s,y)\, ds dy\right)(t,x) = \varphi(t,x).
\eeqn

The specific form of the first term on the right-hand side of \eqref{4.2-1'} is due to the fact that
 $\mathcal{L}$ is of order one in $t$, which is the case for instance for the heat equation. If $\mathcal{L}$ is of order $m>1$ in $t$, then one must also specify the values of the partial derivatives of $u$ in $t$ of all orders less than or equal to $m-1$ at time $t=0$. In this case, the first term on the right-hand side of \eqref{4.2-1'} is more complicated. The wave equations discussed in Section \ref{ch1-1.2} provide  examples of PDEs with $m=2$.

Notice that the first term in \eqref{4.2-1'}, which we denote $I_0(t,x)$, involves only the initial value $\psi$, while the second term involves only the function $f$. In fact, $I_0(t,x)$ is the solution to the homogeneous PDE $\mathcal{L} u(t,x) = 0$ with the same initial condition $\psi$,
 and the second term solves $\cL u(t, x) = f(t,x)$ with initial condition $0$.

When the coefficients of the linear operator $\mathcal{L}$ depend neither on $t$ nor on $x$, the fundamental solution is
{\em homogeneous}, that is, $\Gamma(t,x; s,y): = \Gamma(t-s,x-y)$. In this case, \eqref{4.2-1'} is obtained by a convolution operation:
\beq
\label{4.3-1'}
u(t,x) = [\Gamma(t, \cdot) * \psi] (x) + [\Gamma \star f](t,x),
\eeq
where in the first term on the right-hand side, the convolution is in the space variable, while in the second, it is in both variables. Equations of this kind are called {\em autonomous PDEs}.
\medskip

One can also consider differential operators $\mathcal{L}$ on $\rek$ in connection with PDEs without time evolution. An illustrative example is the Poisson equation
\beqn
\Delta u(x) = f(x), \quad x\in \rek.
\eeqn
In this case, the fundamental solution, in the classical sense, is a function $\Gamma(x)$, defined for all $x\in\rek$, such that
\beqn
u(x) = \int_{\rek}\Gamma(x-y)\ f(y)\ dy,\quad x\in\rek,
\eeqn
and  $\mathcal{L}\ \Gamma(x) = \delta_{0}(x)$ for  $x\in\rek$
(see \cite[p. 22]{evans} or \cite[p. 75]{folland} for the expression of $\Gamma$ for $k\ge 2$).
\medskip

\noindent{\em Distribution-valued solutions}
\medskip

A more abstract approach to partial differential equations deals with distribution-valued solutions\index{distribution-valued!solution}\index{solution!distribution-valued} (see e.g. \cite{schwartz}).
In this framework, the {\em fundamental solution}\index{fundamental solution}\index{solution!fundamental} corresponding to a linear partial differential operator  $\mathcal{L}$ on $\rek$ with constant coefficients is a distribution $S$ (on $\rek$) such that
\beqn
\mathcal{L}S = \delta_{0}.
\eeqn
If $\mathcal{L}$ has constant coefficients and is not null, Theorem 10.2.1 in \cite{hormanderII} (or \cite[p. 62, Theorem 1.56]{folland}) states the existence of $S$. This is the Malgrange-Ehrenpreis Theorem.
In the case of PDEs with a time variable $t$ and a spatial variable $x \in \rek$ ($k\ge 1$), the existence of a
 fundamental solution with support in $\re_+\times \rek$ is a more subtle issue (see e.g. \cite[Section 12.5]{hormanderII}).

In some examples, $S$ is a rather smooth function, for instance for the heat operator. In others, it is not even a function, for example if $\mathcal{L}$ is the wave operator in spatial dimension $k\ge 3$.

In this context, one can consider the PDE $\mathcal{L}u= T$, with $T$ a distribution with compact support, for which a distribution solution is given by the convolution\index{convolution!in space}
\beq
\label{4.4-1'}
u = S * T.
\eeq
Indeed, the convolution $S * T$ is well-defined; furthermore,
\beqn
\mathcal{L}(S * T) = \mathcal{L}S * T  = \delta_{0} * T = T
\eeqn
(see \cite[Chapitre VI, \S 3]{schwartz}, \cite[Lesson 32]{gasquet}).

Applying \cite[Chapitre VI, \S 3, Théorème VII]{schwartz}, we have $S * T=T* S$. Hence, $u=T* S$. This expression
shows an analogy between \eqref{4.3-1'} and \eqref{4.4-1'} when $\psi\equiv 0$.
\medskip

\noindent{\em Linear SPDEs on $\rek$ with additive noise}
\medskip

Let $\dot W$ be a space-time white noise on $\re_+\times \rek$ (see the  Definition \ref{rd1.2.18} and Proposition \ref{rd1.2.19}(f)). Consider now the SPDE
\beq
\label{1'.100}
\begin{cases}
\mathcal{L}u(t,x) &= \dot W(t,x), \quad (t,x)\in\ ]0,T]\times \rek,\\
u(0,x) &= \psi(x), \quad\quad\  x\in \rek.
\end{cases}
\eeq
with deterministic initial condition $\psi$. In the sequel, equations of this kind are called {\em linear SPDEs on $\IR^k$ with additive noise}. \index{additive noise}\index{noise!additive} Since $\dot W$ is neither a  smooth function nor a distribution with compact support, this equation
does not quite fall into either of the settings described above. Nevertheless, the above discussion suggests two possibilities for defining a solution to \eqref{1'.100}.

Indeed, assume first that the differential operator $\mathcal{L}$ is such that there is a fundamental solution $\Gamma$, in the classical sense. Fix a finite time horizon $T>0$, and for simplicity, assume that the initial condition $\psi$ vanishes.
Then, by analogy with \eqref{4.2-1'}, for any $(t,x)\in\ ]0,T]\times \rek$, we should put
\beq
\label{4.5-1'}
u(t,x) = \int_0^t \int_{\rek} \Gamma(t,x; s,y)\, W(ds,dy),\quad \text{a.s.}
\eeq

The random field $(u(t,x),\ (t,x)\in[0,T]\times \rek))$, defined by $u(0,x)=0$ if  $x\in\rek$, and \eqref{4.5-1'} if $(t,x)\in\ ]0,T]\times \rek$, will be called the {\em random field solution}\index{random field!solution}\index{solution!random field} of the SPDE \eqref{1'.100}. Notice that according to Proposition \ref{rd1.2.19} (d), the right-hand side of \eqref{4.5-1'} is well-defined if for fixed $(t,x)\in\ ]0,T]\times \rek$, we have $\Gamma(t,x; \cdot,\ast)\in L^2(\ ]0,T]\times \rek)$
(see Definition \ref{ch4-d1-1'} below).
\smallskip

We now consider the framework of distributions.
Suppose that there is a fundamental solution $S \in \mathcal{S}^\prime(\R^{1+k})$  of $\mathcal{L}$, with support in $\R_+ \times \R^k$, such that, for any $\varphi \in \mathcal{S}(\R^{1+k})$,  $\check{S} * \varphi \in \mathcal{S}(\re^{1+k})$. Here, $\check{S}$  denotes {\it reflection}, defined by $\langle\check{S},\varphi\rangle = \langle S,\check{\varphi}\rangle$ with $\check\varphi(t,x) := \varphi(-t,-x)$. Then for any $T\in \csp(\re^{1+k})$, the convolution
$T * S$ is well-defined by the formula
\beq
\label{def-conv-dist}
\langle T * S, \varphi\rangle : = \langle T, \check{S} *\varphi\rangle.
\eeq
Notice that if $S$ is the fundamental solution of $\mathcal{L}$ (i.e., $\cL S = \delta_0$), then $\check{S}$ is the fundamental solution of the adjoint $\mathcal{L}^\star$. Indeed, for $\varphi \in \cs(\re^{1+k})$, elementary properties of distributions show that
\begin{align*}
   \langle \cL^\star{\check{S}},\varphi \rangle = \langle \check{S}, \cL \varphi \rangle = (S \ast \cL\varphi)(0) = (\cL S \ast\varphi)(0) = (\delta_0\ast\varphi)(0) = \varphi(0),
\end{align*}
therefore,
$\cL^\star\check{S} = \delta_0$.

In view of \eqref{4.4-1'} and using \eqref{def-conv-dist}, we define the notion of {\em weak solution}\index{weak!solution}\index{solution!weak} to the SPDE \eqref{1'.100} with vanishing initial condition as a random linear functional $u$ such that
\beq
\label{4.6-1'}
u(\varphi) = \int_{\re^{1+k}}  (\check{S}\ast \varphi)(t,x)\, W(dt,dx).
\eeq
This is well-defined for any fundamental solution $S$ such that $\check{S} \ast \varphi \in L^2(\re^{1+k})$.

The term {\em weak solution} is justified since, by \eqref{rde1.2.12}, when $\check S * \varphi \in \mathcal{S}(\re^{1+k})$, the right-hand side of \eqref{4.6-1'} is a version of $\langle \dot W,\check S\ast\varphi\rangle$, that is, $u$ is a version of $\dot W * S$, which we also denote $u$, and is the tempered distribution-valued solution to \eqref{1'.100} as in \eqref{4.4-1'}.
Therefore, according to \eqref{4.4-1'},
\beqn
\mathcal{L} u = \mathcal{L}(\dot W * S) = \dot W * \mathcal{L} S = \dot W * \delta_0 = \dot W,
\eeqn
where we have used the property $\mathcal{L}S=\delta_0$.

In the case of linear partial differential operators with constant coefficients, both types of solutions--random field and weak--are related: see  Section \ref{chapter3-s-new}.



 \medskip

\noindent{\em Equations on a bounded domain}
\medskip

Let us now consider a bounded domain $D\subset \rek$ with smooth boundaries. As before, $\mathcal{L}$ is a linear partial differential operator of first order in time, possibly with non constant coefficients acting on functions of $(t,x)\in\ ]0,\infty[ \times D$.

Consider the {\em Dirichlet} boundary value problem,
\beq
\left\{\begin{array}{ll}
\mathcal{L}u(t,x) & = f(t,x),\qquad (t,x)\in\ ]0,T]\times D, \\
u(0,x) & = \psi(x), \qquad\quad x\in D,\\
u(t, x) & = \phi(t,x), \qquad\ t\in [0,T], \ x\in\partial D.
\end{array}\right.
\label{bv-det-1'}
\eeq
If the domain $D$ is suitably regular, there exists a {\em Green's function}\index{Green's function}\index{function!Green's} $G$ associated to $\mathcal{L}$ which plays a role similar to the fundamental solution in Cauchy problems. For example, if the function $\phi$ in \eqref{bv-det-1'} vanishes, then the solution to \eqref{bv-det-1'} is given by
\beq
\label{3.1-rep}
u(t,x) = \int_{D} G(t,x; 0,y) \psi(y)\, dy + \int_0^t ds \int_{D} dy\ G(t,x; s,y) f(s,y).
\eeq
If $\phi$ does not vanish, then one adds to $u$ a solution to a homogenous PDE $\cL v(t,x) = 0$ that satisfies the boundary conditions and has vanishing initial conditions. We emphasize that the Green's function depends on the type of boundary conditions (Dirichlet, Neumann, ...), but not on the specific boundary conditions.
We refer the reader to \cite[p. 228]{ez}
 for an extensive presentation.

 Notice that here too, the first term in \eqref{3.1-rep}, which we denote $I_0(t,x)$, involves only the initial value $\psi$, while the second term involves only the function $f$. In fact, $I_0(t,x)$ is the solution to the homogeneous PDE $\mathcal{L} u(t,x) = 0$ with the same initial condition $\psi$ and boundary condition $\phi=0$, and the second term is the solution of $\cL u(t,x) = f(t,x)$ with initial condition $\psi = 0$ and boundary condition $\phi = 0$.

 \begin{remark}
 \label{3.1-remark-bounded-domain}
For a class of parabolic partial differential operators that includes the heat operator (see \cite[p.3]{friedman64}), and for bounded domains $D\subset \rek$, the Green's function $G(t,x;s,y)$ (sometimes also termed ``fundamental solution'') can be defined for all $(t,x), (s,y)$, $0\le s<t\le T$, and $x,y\in \bar D$, where $\bar D$ is the closure of $D$. We refer to \cite[Chapter 1, Theorem 8, p.19]{friedman64} for details. In this case, we have the validity of the analogue of \eqref{3.1-rep} for $(t,x)\in[0,T]\times \bar D$:
 \beq
\label{3.1-representation}
u(t,x) = \int_{D} G(t,x; 0,y) \psi(y)\, dy + \int_0^t ds \int_{D} dy\ G(t,x; s,y) f(s,y).
\eeq
 \end{remark}

Consider an operator $\mathcal{L}$ on a bounded domain $D\subset\rek$. A Dirichlet boundary value problem takes the form
\beq
\label{bv-notime}
\begin{cases}
\mathcal{L}u(x) & = f(x),\qquad x\in D, \\
u(x) & = \phi(x),\qquad x\in\partial D.
\end{cases}
\eeq
In this case, a {\em Green's function} is a mapping $G$ defined on $D\times D$ verifying suitable conditions
(see \cite{GT2001}). The Green's function depends on the choice of boundary conditions: for example, if $\mathcal{L} = \Delta$ then, in the case of Dirichlet boundary conditions, $G$ is obtained by solving the problem
\beqn
\begin{cases}
\Delta_y G(x,y) & = \delta_{x}(y), \qquad\ y\in D,\\
G(x,y) & = 0,\qquad\qquad y\in \partial D,
\end{cases}
\eeqn
for fixed $x\in D$ (see e.g. \cite[p. 35]{evans}.
In addition, if $\phi\equiv 0$, then we have the following representation of the solution of \eqref{bv-notime}:
\beqn
u(x) = \int_D G(x,y) f(y) dy.
\eeqn
If we replace $f$ by white noise on $D$ and we consider the SPDE
\beq
\label{3.1(*1)}
\begin{cases}
\mathcal{L}u(x) & = \dot W(x),\qquad x\in D, \\
u(x) & = 0, \qquad\qquad x\in\partial D,
\end{cases}
\eeq
then the random field solution should be
\beq
\label{poiss}
u(x) = \int_D G(x,y)\, W(dy),
\eeq
provided $G(x,\ast)\in L^2(D)$. We make this statement more precise below.
\medskip

\noindent{\em Random field solutions}
\medskip

After this introductory discussion, we begin our study of random field solutions to linear SPDEs in spatial dimension $k\ge1$ driven by space-time white noise.
\smallskip

Let $\mathcal{L}$ be a partial differential operator on $]0,\infty[ \times \rek$, such as the heat operator $\mathcal{L}= \frac{\partial}{\partial t}- \Delta$, or the wave operator $\mathcal{L} = \frac{\partial^2}{\partial t^2}- \Delta$, on a bounded or unbounded domain $]0,\infty[ \times D$, where $D \subset \rek$. We consider the SPDE
\beq
\label{1'.100-bis}
\mathcal{L}u(t,x) = \dot W(t,x), \quad (t,x)\in\ ]0,\infty[  \times D,
\eeq
with given initial conditions (at $t=0$), and if necessary, also boundary conditions.

Suppose that there is a Borel function $\Gamma (t,x;s,y)$ which is the fundamental solution or the Green's function on $]0,\infty[  \times D$ associated to $\mathcal{L}$ (and the corresponding type of boundary conditions).

\begin{assump}
\label{ch4-a1-1'}
For all $(t,x) \in \IR_+ \times D$,
$$
\re_+\times D\ni(s,y) \mapsto \Gamma(t,x;s,y) 1_{[0,t[}(s)
$$
belongs to $L^2(\re_+\times D)$.
\end{assump}

\begin{def1}
\label{ch4-d1-1'}
 Let $W$ be a space-time white noise on $\re_+\times D$. Under Assumption \ref{ch4-a1-1'}, the random field solution\index{random field!solution}\index{solution!random field} to the SPDE $\mathcal{L}u=\dot W$ on $\IR_+ \times D$,  with the specified initial conditions and boundary conditions, is the random field
\beq
\label{4.11-1'}
   u(t,x) = I_0(t,x) + \int_0^t \int_{D} \Gamma(t,x;s,y) W(ds,dy),
\eeq
$(t,x)\in\ \re_+\times D$, where $I_0(t,x)$ is the solution to the homogeneous PDE $\mathcal{L} u = 0$ with the same initial and boundary conditions, and the stochastic integral is as defined in Section \ref{ch2new-s2}.
\end{def1}

According to \eqref{4.11-1'}, we have in particular $u(0,x)=I_0(0,x)$, since the stochastic integral vanishes at $t=0$.
\medskip

Analogously, suppose that $\mathcal{L}$ is a partial differential operator on $\rek$, such as the Laplacian $\Delta$. Suppose that there is a Borel function $\Gamma(x,y)$ that is the fundamental solution or the Green's function associated to $\mathcal{L}$ on $D$ (and the corresponding type of boundary conditions) satisfying:
\begin{assump}
\label{ch4-a20-1'}
For all $x\in D$,
$$
 D\ni y \mapsto \Gamma(x,y)
$$
belongs to $L^2(D)$.
\end{assump}
The notion of random field solution to \eqref{3.1(*1)} is as follows.
\begin{def1}
\label{ch4-d2-1'}
 Let $W$ be a white noise on $D$. Under Assumption \ref{ch4-a20-1'}, the random field solution to the SPDE $\mathcal{L}u=\dot W$ in $D$, with the specified  boundary conditions, is
\beq
\label{4.11-1'-bis}
   u(x) = I_0(x) + \int_{D} \Gamma(x,y)\, W(dy),
\eeq
where $I_0(x)$ is the solution to the homogeneous PDE $\mathcal{L}u = 0$ with the same boundary conditions.
\end{def1}
\smallskip

In agreement with Remark \ref{3.1-remark-bounded-domain}, when possible, we replace $D$ by $\bar D$ in the above assumptions and definitions.
\medskip

It turns out that, in many interesting cases, Assumptions \ref{ch4-a1-1'} and \ref{ch4-a20-1'} are only satisfied in low dimensions. For example, for the heat (or the wave) equation, the restriction is $k=1$ (see Lemma \ref{app2-l2}). Hence, if the random forcing is a space-time white noise, it makes sense to start the study with $k=1$.
\medskip

In the remainder of this chapter, we discuss some fundamental examples of linear SPDEs. We establish the existence of random field solutions and prove several properties of their sample paths.


\section{The stochastic heat equation on $\IR$}
\label{ch4-ss2.1-1'}

   We consider the SPDE
   \beq
   \label{4.20-1'}
   \begin{cases}
   \frac{\partial}{\partial t} u(t,x)- \frac{\partial^{2}}{\partial x^{2}} u(t,x) = \dot W(t,x), &  (t,x)\in\ ]0,\infty[\times \re,\\
   u(0,x)= u_0(x), & x\in \re,
   \end{cases}
   \eeq
   where $\dot W$ is the space-time white noise on $\re_+\times \R$ given in Definition \ref{rd1.2.18} and Proposition \ref{rd1.2.19}, and $u_0$ is a function from $\re$ into $\re$.

   Observe that \eqref{4.20-1'} is obtained from the PDE \eqref{heatpde} by replacing the external deterministic forcing $f(t,x)$ by the random forcing $\dot W(t,x)$.

   \subsection{Existence of a random field solution}

   The fundamental solution of the heat operator\index{heat!operator}\index{operator!heat}
   \beqn
   \mathcal{L}=  \frac{\partial}{\partial t} - \frac{\partial^{2}}{\partial x^{2}}
   \eeqn
    is
   \beqn
   \Gamma(t,x;s,y) := \Gamma(t-s,x-y),
   \eeqn
   with
   \beq
   \label{heatcauchy-1'}
  \Gamma(s, y)= \frac{1}{\sqrt{4\pi s}} \exp\left(-\frac{y^2}{4s}\right)\, 1_{]0,\infty[}(s), \quad y\in\re.
   \eeq
   This function satisfies
   \beq
   \label{value-at-zero}
   \lim_{s\downarrow 0}\Gamma(s, y) = \delta_0(y),
   \eeq
   in $\cs^\prime(\re)$ (see e.g. \cite[p. 217]{taylor1}).

  Notice that the map $(s,y) \mapsto \Gamma(t-s,x-y)$ belongs to $L^2(\IR_+\times \IR)$. Indeed, using the properties of the Gaussian density, we see that
  \begin{align}
  \label{heatcauchy-11'}
  \int_0^t ds \int_{\IR}dy\, \Gamma^2(t-s,x-y) &=  \int_0^t ds \int_{\IR}dy\, \Gamma^2(s,y)\notag\\
  & = \int_0^t ds \frac{1}{\sqrt{8\pi s}} \int_{\IR}dy\, \frac{1}{\sqrt{2\pi s}}\exp\left(-\frac{y^2}{2s}\right)\notag\\
  &=\int_0^t ds\, \frac{1}{\sqrt{8\pi s}} = \left(\frac{t}{2\pi}\right)^\half.
\end{align}
Hence, Assumption \ref{ch4-a1-1'} is satisfied.

The fundamental solution $\Gamma$ satisfies the semigroup property\index{semigroup}\index{property!semigroup}
\beq
\label{semig-heat}
\int_{\re} dz\  \Gamma(r,z) \Gamma(s,x-z) = \Gamma(s+r,x),
\eeq
for any $s,r\in\ ]0,\infty[$, a property that is checked by direct integration. With the convention $\Gamma(0,x) = \delta_0(x)$ (motivated by \eqref{value-at-zero}), for any $x\in\re$, we see that \eqref{semig-heat} actually holds for any $s,r\in\re_+$.

Assume that the initial value $u_0$ is such that for any $(t,x)\in\ ]0,\infty[ \times \re$, the function $\Gamma(t,x-\ast)u_0(\ast)$ belongs to $L^1(\re)$.
This condition on $u_0$ is equivalent to
\beq
\label{ch1'-v00}
\int_{\re} e^{-ay^2} |u_0(y)|\ dy <\infty, \ {\text{for all}}\  a>0.
\eeq

The solution to the homogeneous PDE $\mathcal{L}u=0$ (with initial condition $u_0$) is
\beq
\label{ch1'-v001}
I_0(t,x) =
\begin{cases}
\int_\re\Gamma(t,x-y) u_0(y)\,dy,& (t,x)\in\ ]0,\infty[ \times \re,\\
u_0(x),& (t,x)\in \{0\}\times\re.
\end{cases}
\eeq
This is well-defined and, on $]0,\infty[\times \re$,
$(t,x) \mapsto I_0(t,x)$ is $\mathcal{C}^\infty$, and for any $T, L >0$, and $0<t_0\le T$, each partial derivative of $I_0(t,x)$ is uniformly bounded on $[t_0,T]\times [-L,L]$: See e.g. \cite[Lemma 2.3.5, p. 27 and Lemma 2.6.13, p. 88]{chen} and for the case $u_0\in \mathcal{S}^\prime(\re)$,  \cite[Proposition 5.1, p 217]{taylor1}.


According to Definition \ref{ch4-d1-1'}, the random field solution to the SPDE \eqref{4.20-1'}
 is the random field $u=(u(t,x),\ (t,x)\in \re_+\times \re)$ given by
\begin{equation}
\label{sollinheat-1'}
   u(t,x) =  I_0(t,x) + \int_0^t \int_{\re} \Gamma(t-s,x-y)\, W(ds,dy).
\end{equation}

The random field  $u=(u(t,x),\ (t,x)\in\re_+\times \re)$ is Gaussian with $E[u(t,x)] = I_0(t,x)$ and
\begin{align*}
{\rm{Var}}\left( u(t,x)\right) &=
E\left[\left(\int_0^t \int_{\re} \Gamma(t-s,x-y)\, W(ds,dy)\right)^2\right]\\
&= \int_0^t ds \int_{\re} dy\ \Gamma^2(t-s,x-y) = \left(\frac{t}{2\pi}\right)^\half,
\end{align*}
by the Wiener isometry \eqref{revised-1} and \eqref{heatcauchy-11'}.

\subsection{H\"older continuity properties of the sample paths}
\label{ch1'-ss3.3.1}

A function $g: \rek\longrightarrow \re$ is {\em locally} H\"older continuous with exponent $\eta\in\ ]0,1]$\index{H\"older continuous!locally}\index{locally H\"older continuous}\index{continuous!locally H\"older} if for any compact set $O\subset \rek$, the constant 
\beqn
\Vert g\Vert_{\mathcal{C}^\eta(O)}:= \sup_{x,y\in O,\, x\ne y} \frac{|g(x)- g(y)|}{|x-y|^\eta}
\eeqn
is finite.

In the case where the property
\beqn
\Vert g\Vert_{\mathcal{C}^\eta(\rek)}:= \sup_{x,y\in \rek,\, x\ne y} \frac{|g(x)- g(y)|}{|x-y|^\eta} <\infty
\eeqn
holds, the function $g$ is H\"older continuous (or {\em globally} H\"older continuous).\index{H\"older continuous!globally}\index{globally H\"older continuous}\index{continuous!globally H\"older}

The set $\mathcal{C}^\eta(O)$\label{rd06_14l1} (respectively, $\mathcal{C}^\eta(\rek)$) of $\eta$-H\"older continuous functions on $O$ (respectively, on $\rek$) is defined as
\beqn
\mathcal{C}^\eta(O)= \left\{g: O\rightarrow \re: \Vert g\Vert_{\mathcal{C}^\eta(O)}<\infty\right\}
\eeqn
(respectively,
$
\mathcal{C}^\eta(\rek)= \{g: \rek\rightarrow \re: \Vert g\Vert_{\mathcal{C}^\eta(\rek)}<\infty\}$).

In the study of sample path properties of SPDEs, we are led to consider functions $g: \re_+\times \rek\longrightarrow \re$ depending on two variables $(t,x)\in \re_+\times \rek$.
We say that $g$ is {\em jointly} locally  H\"older continuous with exponents $(\eta_1,\eta_2)$\index{H\"older continuous!jointly}\index{jointly!H\"older continuous}\index{continuous!globally H\"older} if, for all compact sets $A\subset \re_+$, $B\subset \rek$, the constant 
\beqn
\Vert g\Vert_{\mathcal{C}^{\eta_1,\eta_2}(A\times B)}
:=\sup_{(t,x), (s,y)\in A\times B,\,  (t,x)\ne (s,y)} \frac{\vert g(t,x)- g(s,y)\vert}{|s-t|^{\eta_1}+|x-y|^{\eta_2}}
\eeqn
is finite.

If this property holds with $A\times B$ replaced by  $\re_+\times \rek$, the function $g$ is said to be jointly (globally) H\"older continuous.

The set $\mathcal{C}^{\eta_1,\eta_2}(A\times B)$\label{rdjholder} of $(\eta_1,\eta_2)$-H\"older continuous functions on $A\times B$ is defined as
\beqn
\mathcal{C}^{\eta_1,\eta_2}(A\times B)= \left\{g: A\times B \rightarrow \re: \Vert g\Vert_{\mathcal{C}^{\eta_1,\eta_2}(A\times B)}<\infty\right\}.
\eeqn
Similarly, the set $\mathcal{C}^{\eta_1,\eta_2}(\re_+\times \rek)$ of $(\eta_1,\eta_2)$-H\"older continuous functions on $\re_+\times \rek$ is defined as
\beqn
\mathcal{C}^{\eta_1,\eta_2}(\re_+\times \rek)= \left\{g: \re_+\times \rek \rightarrow \re: \Vert g\Vert_{\mathcal{C}^{\eta_1,\eta_2}(\re_+\times \rek)}<\infty\right\}.
\eeqn

In the next propositions, we study the H\"older continuity properties of the sample paths of the random field  solution $u$ to \eqref{4.20-1'}.

For $x \in \re$, let $v(0,x) = 0$ and for $t > 0$ and $x \in \re$, set
\beq
\label{ch1'.v}
v(t,x) =  \int_0^t \int_{\re} \Gamma(t-s,x-y) W(ds,dy).
\eeq

\begin{remark}
\label{continuity-Dunbounded}
In view of \eqref{sollinheat-1'}, the (H\"older) continuity properties of $u$ are related to those of $I_0$ and $v$ (given in \eqref{ch1'-v001} and \eqref{ch1'.v}, respectively). These can be studied separately. In all cases where $v$ is continuous (respectively H\"older continuous), $u$ will be continuous (respectively H\"older continuous) if $I_0$ is. As has been pointed out above, $I_0$ is even $\mathcal{C}^\infty$ on $]0,\infty[\times \re$ if the initial value $u_0$ satisfies \eqref{ch1'-v00}.
\end{remark}

\begin{prop}
\label{ch1'-p0}
For all $(t,x), (s,y) \in \re_+ \times \re$,
\beq
\label{ch1'.v1}
   E\left[( v(t,x) - v(s,y))^2\right] \leq  \left(\pi^{-\frac{1}{4}}\vert t-s\vert^{\frac{1}{4}} + 2^{-\frac{1}{2}}\vert x-y\vert^{\half}\right)^{2}.
\eeq
Therefore, for any $\alpha \in ]0, \frac{1}{4}[$ and any $\beta \in ]0, \half[$, there exists a version of $v=(v(t,x),\  (t,x)\in[0,\infty[\times \re)$ with locally jointly H\"older continuous sample paths with exponents $(\alpha, \beta)$.
\end{prop}

\begin{remark} The constants $\pi^{-\frac{1}{4}}$ and $2^{-\frac{1}{2}}$ on the right-hand side are best possible. This follows from Theorem \ref{ch1'-s3-t1} below.
\end{remark}

\noindent{\em Proof of Proposition \ref{ch1'-p0}}.
By the isometry property \eqref{ch1'-s4.2},
\begin{align}
\label{ch1'-v3}
&E\left[\vert v(t,x) - v(s,y) \vert^2\right]\notag\\
&\quad=E\left[\left(\int_0^t \int_{\re} \Gamma(t-r,x-z)\, W(dr,dz)\right.\right.\notag\\
&\left.\left.\qquad \qquad \qquad - \int_0^s \int_{\re} \Gamma(s-r,y-z)\, W(dr,dz)\right)^2\right]\notag\\
&\quad= \int_{\IR_+} dr \int_{\IR}  dz \left(\Gamma(t-r,x-z)-\Gamma(s-r,y-z)\right)^2 \notag\\
&\quad \le \left[\pi^{-\frac{1}{4}}\vert t-s\vert^{\frac{1}{4}} + 2^{-\frac{1}{2}}\vert x-y\vert^{\half}\right]^2,
\end{align}
where the last inequality follows from  \eqref{1'.500} in Lemma \ref{ch1'-l0}. This proves \eqref{ch1'.v1}.

The claim on ``local joint H\"older continuity'' of the sample paths follows from the version of Kolmogorov's continuity criterion given in Theorem \ref{app1-3-t1}. Indeed, the process $v$ is Gaussian and \eqref{ch1'.v1} implies in particular \eqref{aap1-3.1}, with $\Delta(t, x; s, y) = \vert t - s \vert^\frac{1}{4} +  \vert x - y \vert^\half$, for all bounded intervals $I\subset \re_+$ and $J\subset\re$.
\qed

\begin{lemma}
\label{ch1'-prep-lemma3.2.2}
Assume $u_0\in \mathcal{C}^\eta(\re)$, for some $\eta\in \, ]0,1]$. Then \eqref{ch1'-v00} holds and
 the function $[0,T]\times \IR\ni (t,x) \longrightarrow I_0(t,x)$ defined by
\beq
\label{doublestar-67}
I_0(t,x)=
\begin{cases}
\int_{\IR} \Gamma(t,x-y) u_0(y)\, dy, & (t,x)\in\, ]0,T]\times \re,\\
u_0(x), & t=0,\ x\in\re,
\end{cases}
\eeq
is globally H\"older continuous, jointly in $(t,x)$, with exponents $(\frac{\eta}{2},\eta)$.
\end{lemma}
\begin{proof}
Since $|u_0(y)|\le |u_0(0)| + |y|^\eta$, we see that \eqref{ch1'-v00} is satisfied.

We consider separately the increments in space and in time of $I_0(t,x)$.
Let $x\in \IR$, $h\in\IR_+$. By changing the spatial variable, for any $t> 0$ we have
\begin{align}
\label{ch1'-initial 2}
I_0(t,x+h)-I_0(t,x)& = \int_{\IR} \left[\Gamma(t,x+h-y)-\Gamma(t,x-y)\right] u_0(y)\, dy\notag\\
& = \int_{\IR} \Gamma(t,z)[u_0(x+h-z)-u_0(x-z)]\, dz.
\end{align}
Since $u_0\in\mathcal{C}^{\eta}(\re)$, we obtain
\beq
\label{ch1'-initial 3}
\sup_{t\in \re_+}\sup_{x\in\re}\left\vert I_0(t,x+h)-I_0(t,x)\right\vert \le \Vert u_0\Vert_{\mathcal{C}^\eta(\re)}h^\eta,
\eeq
because $\int_{\IR} \Gamma(t,z)\ dz = 1$.

Consider the family of convolution operators $\left(\Gamma_t,\, t\ge 0\right)$ defined on measurable real functions $f$ satisfying \eqref{ch1'-v00} by
\beqn
\Gamma_t(f)(x):=(\Gamma_t\ast f)(x)=\int_{\IR} \Gamma(t,x-y) f(y)\, dy, \quad t > 0,
\eeqn
and $\Gamma_0 f = f$  (recall the convention $\Gamma(0,x)=\delta_0(x)$).

Notice that $\Gamma_t(u_0)(x)=I_0(t,x)$.
Then, for $f\in\mathcal{C}^\eta(\re)$ and $t>0$,
\beqn
\Gamma_t(f)(x)-f(x)
= \int_{\IR} \Gamma(t,x-y)[f(y)-f(x)]\, dy,
\eeqn
since $\int_{\IR} \Gamma(t,x-y) dy = 1$. Consequently, for any $t \ge 0$,
\begin{align}\nonumber
\left\vert\Gamma_t(f)(x)-f(x)\right\vert &\le  \int_{\IR} \Gamma(t,x-y)\vert f(y)-f(x)\vert\ dy\\ \nonumber
&\le \Vert f\Vert_{\mathcal{C}^\eta(\re)} \int_{\IR} \Gamma(t,x-y)|x-y|^\eta\ dy\\
& = \Vert f\Vert_{\mathcal{C}^\eta(\re)}\frac{2^{\eta}}{\sqrt \pi}\Gamma_E\left(\frac{\eta+1}{2}\right) t^{\frac{\eta}{2}},
\label{ch1'-initial 4}
\end{align}
by \eqref{A3.4},
 where $\Gamma_E$ denotes the Euler Gamma function.
 In particular, $\Gamma_t(f)$ also satisfies \eqref{ch1'-v00}.

Recall the semigroup property \eqref{semig-heat}, which can be written $\Gamma_s(\Gamma_t(f)) = \Gamma_{s+t}(f)$. Using \eqref{ch1'-initial 4}, we see that for any $0\le s\le t$,
\begin{align}
\label{ch1'-initial 5}
\sup_{x\in\IR}\left\vert I_0(t,x)-I_0(s,x)\right\vert & = \sup_{x\in\IR}\left\vert \Gamma_t(u_0)(x) - \Gamma_s(u_0)(x)\right\vert\notag\\
& =  \sup_{x\in\IR}\left\vert \Gamma_{t-s}\left(\Gamma_s(u_0)\right)(x)-\Gamma_s(u_0)(x)\right\vert\notag\\
& \le \Vert \Gamma_s(u_0)\Vert_{\mathcal{C}^\eta(\re)}\frac{2^{\eta}}{\sqrt \pi}\Gamma_E\left(\frac{\eta+1}{2}\right) |t-s|^{\frac{\eta}{2}}\notag\\
& \le \Vert u_0\Vert_{\mathcal{C}^\eta(\re)}\frac{2^{\eta}}{\sqrt \pi}\Gamma_E\left(\frac{\eta+1}{2}\right) |t-s|^{\frac{\eta}{2}},
\end{align}
where, in the last inequality, we have used  \eqref{ch1'-initial 3}.

From \eqref{ch1'-initial 3} and \eqref{ch1'-initial 5}, we obtain
\beq
\label{ch1'-initial 6}
\left\vert I_0(t,x+h)-I_0(s,x)\right\vert \le C \Vert u_0\Vert_{\mathcal{C}^\eta(\re)}\left(h^\eta+|t-s|^{\frac{\eta}{2}}\right),
\eeq
with $C= \max\left(1,\frac{2^{\eta}}{\sqrt \pi}\Gamma_E\left(\frac{\eta+1}{2}\right)\right)$, and this inequality holds for any $0\le s\le t$ and any $x\in\re$, $h\in\re_+$. This proves the property on global H\"older continuity of $I_0$.
\end{proof}
\medskip

\noindent{\em H\"older continuity of the random field solution}
\medskip


\begin{prop}
\label{ch1'-p1}
 We assume that the initial condition $u_0$ satisfies \eqref{ch1'-v00}. Let $u=(u(t,x), (t,x)\in\R_+\times \re)$ be the random field given in \eqref{sollinheat-1'}.  Fix $T, L >0$.
\begin{enumerate}
 \item {Continuity away from $t=0$.} Fix $0<t_0\le T$. For every $p\in[2,\infty[$, there exists a constant $C:=C(p,t_0,T,L,u_0)>0$ such that, for all $(t,x), (s,y) \in [t_0,T] \times [-L,L]$,
\beq
\label{incheat}
   E\left[\vert u(t,x) - u(s,y) \vert^p\right] \leq C \left(\vert t-s\vert^{\frac{1}{4}} + \vert x-y\vert^{\half}\right)^{p}.
\eeq
Therefore, there exists a version of the process $u$ with locally jointly H\"older continuous sample paths with exponents  $(\alpha,\beta)$, where $\alpha\in\ ]0,\frac{1}{4}[$ and $\beta\in\ ]0,\half[$.

\item {Continuity including $t=0$.} Assume in addition that $u_0\in\mathcal{C}^\eta(\re)$ for some $\eta\in\ ]0,1]$. Then
for every $p\in[2,\infty[$, there exists a constant $C:=C(p,T,L,u_0,\eta)>0$ such that, for all $(t,x), (s,y) \in [0,T]\times [-L,L]$,
\beq
\label{incheat-bis}
   E\left[\vert u(t,x) - u(s,y) \vert^p\right] \leq C \left(\vert t-s\vert^{\frac{1}{4}\wedge\frac{\eta}{2}} + \vert x-y\vert^{\half\wedge \eta}\right)^{p}.
\eeq
Thus, there exists a version of $u$ with jointly H\"older continuous sample paths in $[0,T]\times [-L,L]$. The Hölder exponents are $(\alpha,\beta)$, where
in the time variable $t$ the constraints on $\alpha$ are
\beqn
\alpha\in \left]0,\tfrac{1}{4}\right[,  {\text{ if}}\ \eta\ge \tfrac{1}{2},\qquad \alpha\in  \left]0,\tfrac{\eta}{2}\right],  {\text{ if}}\ \eta < \tfrac{1}{2},
\eeqn
while in the space variable $x$, the constraints on $\beta$ are
\beqn
\beta\in\left]0,\tfrac{1}{2}\right[,  {\text{ if}}\ \eta\ge \tfrac{1}{2},\qquad \beta\in \left]0,\eta\right],  {\text{ if}}\ \eta < \tfrac{1}{2}.
\eeqn
\end{enumerate}
\end{prop}
\begin{proof} 1. Since $u(t,x) = I_0(t,x) + v(t,x)$, by the triangle inequality,
\beq
\label{dec-i+si}
\Vert u(t,x)-u(s,y)\Vert_{L^p(\Omega)} \le |I_0(t,x) - I_0(s,y)| + \Vert v(t,x)-v(s,y)\Vert_{L^p(\Omega)} .
\eeq
The function $I_0$ is $\mathcal{C}^\infty$ on $]0,\infty[\times \re$. Therefore, for any $t_0\in\ ]0,T]$ and $L>0$, there exists a constant $C(t_0,T,L,u_0)>0$ such that, for any $(t,x), (s,y) \in\  ]t_0,T] \times [-L,L]$,
\beq
\label{1'.initial}
\left\vert I_0(t,x) - I_0(s,y)\right\vert \le C(t_0,T,L,u_0) (|t-s| + |x-y|).
\eeq
Since $(v(t,x))$ is a centred Gaussian random field, using \eqref{A3.4.0-bis} we have
\beq
\label{conmens}
\Vert v(t,x) - v(s,y)\Vert _{L^p(\Omega)} = c_p^{1/p} \Vert v(t,x) - v(s,y)\Vert _{L^2(\Omega)},
\eeq
with $c_p = \left(2^p/\pi\right)^\half \Gamma_E\left((p+1)/2\right)$. From  \eqref{dec-i+si}, and using \eqref{1'.initial},
\eqref{conmens} and \eqref{ch1'.v1}, we see that there is a constant $C:=C(p,t_0,T,L,u_0)$ such that
\beqn
E\left[|u(t,x) - u(s,y)|^p\right] \le C\left(|t-s| + |x-y| +  |t-s|^{\frac{1}{4}} + |x-y|^{\frac{1}{2}}\right)^p.
\eeqn
This proves \eqref{incheat}.

The statement on H\"older continuity of the sample paths follows from Kolmogorov's continuity criterion Theorem \ref{app1-3-t1} (or Theorem \ref{ch1'-s7-t2}). This proves claim 1.
\smallskip

2. For the proof of claim 2, we apply \eqref{dec-i+si} and \eqref{conmens}. Using the estimates \eqref{ch1'-initial 6} and \eqref{ch1'.v1}, we obtain
\eqref{incheat-bis}, since on compact sets, increments with the smaller exponents dominate.
The claim about H\"older continuity of the sample paths is obtained  by applying Kolmogorov's continuity criterion.
\end{proof}

\medskip

It is also possible to obtain lower bounds for the second moment of increments of the random field $v(t,x)$. Along with Proposition \ref{ch1'-p0}, this yields the following.

\begin{prop}
\label{ch1'-p3}
 Let $(v(t,x),\, (t,x)\in\R_+\times \re)$ be the random field given in \eqref{ch1'.v}. Fix $0<t_0\le T$.
 Then there exists a constant $C_1:=C_1(t_0)>0$ such that, for all $(t,x), (s,y) \in [t_0,T] \times \IR$ with $|x-y| \le \sqrt{t_0}$,
\begin{align}
\label{lower-st-heat-r}
 C_1 \left(\vert t-s\vert^{\half} + \vert x-y\vert\right)&\le E\left[(v(t,x) - v(s,y))^2\right]\notag\\
 & \leq \left(\pi^{-\frac{1}{4}}\vert t-s\vert^{\half} + 2^{-\frac{1}{2}}\vert x-y\vert\right).
\end{align}
\end{prop}
\begin{proof}
The upper bound is \eqref{ch1'.v1}, and holds for any $(t,x), (s,y) \in \re_+ \times \IR$.

For the lower bound, let
\[
   C_2 = \frac{1}{\sqrt{2\pi}}, \qquad C_3 = c_{0}, \qquad C_4 = \frac{1}{\sqrt{\pi}},
\]
where $c_{0}$ is the constant that appears in \eqref{1'.50} for $C=1$ there. We observe first that the following inequalities hold:
\begin{align}
   E\left[( v(t,x) - v(s,y) )^2\right] &\geq C_2 \vert t-s\vert^{\half}, \quad x,y \in \IR,\ s,t \in \IR_+ \label{rde3.2.22} \\
   E\left[( v(t,x) - v(t,y) )^2\right] &\geq C_3 \vert x-y\vert,\quad x,y \in \IR\text{ with } |x-y| \le \sqrt{t_0},\ t \geq t_0, \label{rde3.2.23}\\
   E\left[( v(t,y) - v(s,y) )^2\right] &\leq C_4 \vert t-s\vert^{\half}, \quad y \in \IR,\ s,t \in \IR_+. \label{rde3.2.24}
\end{align}
Indeed, for \eqref{rde3.2.22}, observe by the isometry property of the stochastic integral that for $t\geq s \geq 0$,
\begin{align*}
E\left[( v(t,x) - v(s,y))^2\right] &=\int_0^t d r\int_{\R} d z\ [\Gamma(t-r,x-z)-\Gamma(s-r,y-z)]^2\\
 &\geq \int_s^t d r\int_{\R} d z\ \Gamma^2(t-r,x-z) \\
 &= \int_0^{t-s} d r\int_{\R} d z\ \Gamma^2(r,x-z)
 =\left(\frac{t-s}{2\pi}\right)^\half,
\end{align*}
by \eqref{heatcauchy-11'}. For \eqref{rde3.2.23}, it suffices to apply \eqref{1'.50}. The upper bound \eqref{rde3.2.24} follows from \eqref{ch1'.v1} (or from \eqref{1'.500}).

We now deduce the lower bound in \eqref{lower-st-heat-r} from \eqref{rde3.2.22}--\eqref{rde3.2.24}, by considering two cases.
\smallskip

\textit{Case 1: $\vert t-s\vert^{\half}\geq\frac{C_3}{4C_4}\vert x-y\vert$}. By \eqref{rde3.2.22} and by the inequality that defines this Case,
\begin{align*}
E\left[( v(t,x) - v(s,y) )^2)\right] \geq C_2 \vert t-s\vert^{\half} \geq \frac{C_2}{2} \left[\vert t-s\vert^{\half} + \frac{C_3}{4C_4}  \vert x-y\vert \right].
\end{align*}

\textit{Case 2: $\vert t-s\vert^{\half}\leq \frac{C_3}{4C_4}\vert x-y\vert$}. Recall that $|x-y|\le \sqrt{t_0}$. Apply the inequality $(a+b)^2\ge \half a^2-b^2$, then \eqref{rde3.2.23} and \eqref{rde3.2.24}, to write
\begin{align*}
E\left[( v(t,x) - v(s,y))^2\right] &\geq \half E\left[(v(t,x)-v(t,y))^2\right]-E\left[(v(t,y)-v(s,y))^2\right] \\
   &\geq \frac{C_3}{2} |x-y| - C_4 \vert t-s\vert^\half \\
   &\geq \frac{C_3}{2} |x-y| - \frac{C_3}{4}|x-y| = \frac{C_3}{4}|x-y|\\
   &\geq \frac{C_3}{4}\left[ \frac12  |x-y| + \frac12 \frac{4C_4}{C_3}  \vert t-s\vert^{\half} \right].
\end{align*}
This completes the proof with $C_1 = \min(\tfrac{C_2}{2}, \tfrac{C_2C_3}{8C_4}, \tfrac{C_3}{8}, \tfrac{C_4}{2})$.
\end{proof}


\begin{remark}
\label{ch1'-r1} For a Gaussian stochastic process $X = (X_\tau,\, \tau\in \mathbb{T})$  indexed by a set $\mathbb{T}$, there is the notion of canonical pseudo-metric,\index{pseudo-metric}\index{canonical!pseudo-metric} defined by
\beqn
\delta(\tau,\bar\tau) = \left(E\left[( X_\tau - X_{\bar\tau})^2\right]\right)^{\half}.
\eeqn
Proposition \ref{ch1'-p3} tells us that the canonical pseudo-metric associated to the random field $v=(v(t,x),\  (t,x)\in\R_+\times \re)$ is locally equivalent (up to multiplicative constants) to
\beqn
\delta((t,x),(s,y))=
\vert t-s\vert^{\frac{1}{4}} + \vert x-y\vert^\half.
\eeqn
This property has a relevant role for instance in the study of hitting probabilities of the random field $v$ (see \cite[Theorem 7.6, p. 188]{xiao} and Section \ref{ch2'-s7} in this book).
\end{remark}
\medskip

\noindent{\em Sharpness of the degree of H\"older continuity}
\medskip

The following result shows that the H\"older exponents obtained in Proposition \ref{ch1'-p1} are sharp.
\begin{prop}
\label{ch1'-p100}
Fix $T>0$, and let $(u(t,x),\, (t,x)\in [0,T]\times \re)$ be the random field defined in \eqref{sollinheat-1'}.
\begin{enumerate}
\item Fix $x\in\re$, $K\subset [0,T]$ a closed interval with positive length, and $\alpha\in\left]\frac{1}{4},1\right]$. Then a.s., the sample paths of the stochastic process $(u(t,x),\, t\in K)$ are not H\"older continuous with exponent $\alpha$.
\item Fix $t\in\,]0,\infty[$, $J\subset\re$ a closed interval with positive length, and $\beta\in\left]\frac{1}{2},1\right]$. Then a.s., the sample paths of the stochastic process $(u(t,x),\, x\in J)$ are not  H\"older continuous with exponent
$\beta.$
\end{enumerate}
\end{prop}

\begin{proof} Since $(t,x) \mapsto I_0(t,x)$ is $\ca^\infty$ on $]0,\infty[\times \IR$, it suffices to prove the proposition for the centred Gaussian process $v$ defined in \eqref{ch1'.v}.

For the proof of Claim 1, we apply the first two equalities in \eqref{ch1'-v3} and \eqref{1'.4-bis} to see that for fixed $x\in \re$ and for $|t-s|$ sufficiently small,
$$
  E\left[( v(t,x) - v(s,x))^2\right] \geq c_0 |t-s|^\half,\qquad \text{with } c_0 = \frac{\sqrt{2}-1}{\sqrt{2\pi}}.
$$
Therefore, condition \eqref{app1-3.5} is satisfied with $\alpha = \tfrac{1}{4}$. Claim 1 now follows from Theorem \ref{app1-3-t2}.

For the proof of Claim 2, we appeal to \eqref{1'.50} to deduce that for fixed $t_0\in]0,\infty[$ and $h_0>0$, for any $t\ge t_0$,  the stochastic process $(v(t,x),\, x\in J)$ satisfies, for $|x-y| \leq \sqrt{t_0}$,
$$
   E\left[( v(t,x) - v(t,y))^2\right] \geq c_{0} |x-y|,
$$
where $c_{0}$ is the constant on the left-hand side of \eqref{1'.50} (with $C=1$ there). Therefore, condition \eqref{app1-3.5} is satisfied with $\alpha = \frac{1}{2}$ for all small enough sub-intervals of $J$. Claim 2 now follows from Theorem \ref{app1-3-t2}.
\end{proof}
\medskip


\subsection{Structure of the solution restricted to lines}
\label{ch1'-ss3.3.2}

A fractional Brownian motion\index{fractional!Brownian motion}\index{Brownian!motion, fractional} $\left(B^H_t, \ t\in \re_+\right)$ with Hurst parameter\index{Hurst parameter}\index{parameter!Hurst} $H\in\ ]0,1[$ is a mean-zero Gaussian process with continuous sample paths and covariance
\beqn
E\left[B^H_s B^H_t\right] = \frac{1}{2} \left(s^{2H} + t^{2H} - |s-t|^{2H}\right), \quad  s,t\in\re_+
\eeqn
(see e.g. \cite{hu2005}).

A two-sided standard Brownian motion\index{two-sided Brownian motion}\index{Brownian!motion, two-sided} $(Z(x),\ x\in \re)$ is a mean zero Gaussian process with continuous sample paths such that $Z(0)=0$ and $E[(Z(x)-Z(y)^2] = |x-y|$. In particular, using the formula for the covariance in \eqref{ix} below, we have
\beqn
    E[Z(x) Z(y)] = \half (\vert x \vert  + \vert y \vert - \vert x - y \vert),\quad   x, y \in \R,
    \eeqn
so $(Z(x),\, x \in \R_+)$ and $(Z(-x),\, x \in \R_+)$ are independent standard Brownian motions.

The next theorem provides additional insight on the trajectories of the random fields  $(u(t,x), \, (t,x)\in\re_+\times \re)$ and $(v(t,x), \, (t,x)\in\re_+\times \re)$ of \eqref{sollinheat-1'} and \eqref{ch1'.v} respectively, when one of its variables, time or space, are fixed.

\begin{thm}
\label{ch1'-s3-t1}
Let $u$ be the solution to \eqref{4.20-1'}, given in \eqref{sollinheat-1'} and let $v$ be as defined in \eqref{ch1'.v}.

\noindent (a) Fix $x\in\re$. There exists a fractional Brownian motion $\left(X_t(x),\,  t\in\re_+\right)$ with Hurst parameter $H=\frac{1}{4}$, and a mean-zero Gaussian process $\left(R_t,\, t\in\re_+\right)$ that is independent of the space-time white noise $\dot W$ (and does not depend on $x$) whose sample paths are a.s. locally H\"older continuous with exponent $\frac{1}{4}-\varepsilon$ on $\re_+$, for all $\varepsilon>0$, and $\mathcal{C}^\infty$ on $]0,\infty[$, such that a.s., for all $t\in\re_+$,
\beq
\label{ch1'-s3-100}
u(t,x) = I_0(t,x) + R_t + \pi^{-\tfrac{1}{4}} X_t(x).
 \eeq

\noindent (b) Fix $t>0$. There is a two-sided standard Brownian motion $\left(B_x(t),\, x\in \re\right)$, and a mean-zero Gaussian process $\left(S_x(t),\, x\in \re\right)$ independent of $\dot W$ with $\mathcal{C}^\infty$ sample paths, such that a.s., for all $x\in\re$,
\beq
 \label{ch1'-s3-101}
 u(t,x) = I_0(t,x) + v(t,0) + S_x(t) +  2^{-\tfrac{1}{2}} B_x(t).
 \eeq
\end{thm}
\begin{remark}\label{rdrem3.2.9}
 It is satisfying that the constants $\pi^{-\frac{1}{4}}$ and $2^{-\frac{1}{2}}$ also appear in the upper bound of Proposition \ref{ch1'-p0}.
\end{remark}
\medskip

\noindent{\em Proof of Theorem \ref{ch1'-s3-t1}}.
(a) Let $\dot w$ be a white noise on $\re$ that is independent of $\dot W$. For $t\in\re_+$, define
\beqn
R_t = \frac{1}{\sqrt{4\pi}} \int_{\re} \frac{1-e^{-tz^2}}{z}\ w(dz).
\eeqn
We will show that
\beq
\label{ch1'-s3-102}
X_t(x) = \pi^{1/4}\left(v(t,x)-R_t\right)
\eeq
is a fractional Brownian motion with $H=\frac{1}{4}$, and that  $R=\left(R_t,\, t\in\re_+\right)$ has the required properties. This will prove \eqref{ch1'-s3-100}. Since we are using continuous versions of each process, it suffices to check that $(X_t(x),\, t\in\re_+)$ has the appropriate covariance.

We observe that $R$ and $\dot W$ are independent, and that for all $t\in\re_+$, $E[R_t] =0$ and
\beqn
E\left[R_t^2\right] = \frac{1}{4 \pi} \int_{\re} \left(\frac{1-e^{-tz^2}}{z}\right)^2 \ dz < \infty.
\eeqn
Further, for $h\ge 0$,
\begin{align}
\label{ch1'-s3-103}
E\left[\left(R_{t+h} - R_t\right)^2\right] & = \frac{1}{4 \pi} \int_{\re} \left(\frac{e^{-tz^2}-e^{-(t+h)z^2}}{z}\right)^2\ dz\notag\\
& = \frac{1}{4 \pi} \int_{\re} e^{-2tz^2}\left(\frac{1-e^{-hz^2}}{z}\right)^2\ dz\notag\\
& = \frac{1}{2 \pi}  \int_{\re} dz \left(1-e^{-hz^2}\right)^2 \int_t^\infty  ds\ e^{-2sz^2}\notag\\
& = \frac{1}{2 \pi}  \int_t^\infty ds \int_{\re} dz \left(e^{-sz^2} - e^{-(s+h)z^2} \right)^2\notag\\
& = \int_t^\infty ds  \int_{\re} dy \left(\Gamma(s,y) - \Gamma(s+h,y)\right)^2,
\end{align}
where we have used Plancherel's theorem in the last equality.

We now observe that
\beq
\label{ch1'-s3-104}
E\left[\left(v(t+h,x) - v(t,x)\right)^2\right] = E\left[A_1(t,h,x)^2\right] + E\left[A_2(t,h,x)^2\right],
\eeq
where
\begin{align*}
A_1(t,h,x)& = \int_0^t \int_{\re}\left(\Gamma(t+h-s,x-y) - \Gamma(t-s,x-y)\right)\, W(ds,dy),\\
A_2(t,h,x)& = \int_t^{t+h} \int_{\re} \Gamma(t+h-s,x-y)\, W(ds,dy).
\end{align*}
The random fields $A_1$ and $A_2$ are independent, and
\beq
\label{ch1'-s3-105}
E\left[A_2(t,h,x)^2\right] =  \int_t^{t+h} ds \int_{\re} dy\,  \Gamma(t+h-s,x-y)^2 = \left(\frac{h}{2\pi}\right)^{\half},
\eeq
by \eqref{heatcauchy-11'} after a change of variable.

 As for $A_1$,
 \begin{align*}
E\left[A_1(t,h,x)^2\right] &= \int_0^t ds \int_{\re} dy \left(\Gamma(t+h-s,x-y) - \Gamma(t-s,x-y)\right)^2\\
& = \int_0^t ds \int_{\re} dy \left( \Gamma(s+h,y)-\Gamma(s,y)\right)^2\\
& = \int_0^\infty ds \int_{\re} dy \left( \Gamma(s+h,y) - \Gamma(s,y)\right)^2 - E\left[\left(R_{t+h} - R_t\right)^2\right],
\end{align*}
where the last equality uses \eqref{ch1'-s3-103}.

By \eqref{1'.4}, we see that
\beqn
E\left[A_1(t,h,x)^2\right] = \left(\frac{h}{2\pi}\right)^{\half}(\sqrt 2-1) - E\left[\left(R_{t+h} - R_t\right)^2\right],
\eeqn
so from \eqref{ch1'-s3-104} and \eqref{ch1'-s3-105},
\beqn
E\left[\left(v(t+h,x) - v(t,x)\right)^2\right] = \left(\frac{h}{\pi}\right)^\half - E\left[\left(R_{t+h} - R_t\right)^2\right].
\eeqn
Moving $E\left[\left(R_{t+h} - R_t\right)^2\right]$ to the left-hand side, using the independence of $(v(t,x))$ and $(R_t)$ and recalling \eqref{ch1'-s3-102}, we see that
\beqn
E\left[\left(X_{t+h}(x) - X_t(x)\right)^2\right] = \sqrt{h}.
\eeqn
Since $X_0(x) = \pi^{1/4}\left(v(0,x)-R_0\right)=0$, we deduce that
\beqn
E\left[X_t(x)^2\right] = E\left[\left(X_t(x)-X_0(x)\right)^2\right] = \sqrt{t},
\eeqn
and therefore
\begin{align}
\label{ix}
E[X_s(x) X_t(x)] & = \frac{1}{2} \Big(E[X_s(x)^2] +E[X_t(x)^2]-E\left[(X_s(x)-X_t(x))^2\right]\Big)\notag\\
& =  \frac{1}{2}\left(\sqrt s + \sqrt t- \sqrt{|t-s|}\right).
\end{align}
Since $(X_t)$ is a Gaussian process with $E[X_t] =0$, $(X_t)$ is a fractional Brownian motion with $H=\frac{1}{4}$.
\smallskip

It remains to check that $(R_t)$ has the requested H\"older continuity and differentiability properties. As in the first equality in \eqref{ch1'-s3-103}, for $s,t\in\re_+$,
\begin{align*}
E\left[\left(R_{t} - R_s\right)^2\right] & = \frac{1}{4 \pi} \int_{\re} \left(\frac{e^{-tz^2}-e^{-sz^2}}{z}\right)^2\ dz\\
& \le  \frac{1}{4 \pi} \int_{\re} \left(\frac{1-e^{-|t-s|z^2}}{z}\right)^2\ dz\\
&= \frac{\sqrt{|t-s|}}{4 \pi}\int_{\re} \left(\frac{1-e^{-y^2}}{y}\right)^2\ dy.
\end{align*}
From Kolmogorov's continuity criterion (see Theorem \ref{app1-3-t1}), $R$ has a version that is locally H\"older continuous with exponent $\frac{1}{4}-\varepsilon$, for all $\varepsilon>0$.

In order to check the differentiability of $t\mapsto R_t$, for $t>0$ and $n\ge 1$, define
\begin{align*}
R_t^{(n)} & = \frac{1}{\sqrt{4 \pi}} \int_{\re} \frac{\partial^n}{\partial t^n} \left(\frac{1-e^{-tz^2}}{z}\right)\ w(dz)\\
& = \frac{(-1)^{n+1}}{\sqrt{4 \pi}} \int_{\re} z^{2n-1} e^{-tz^2}\ w(dz).
\end{align*}
The integrand belongs to $L^2(\re)$ and therefore $(R_t^{(n)},\, t\in\,]0,\infty[)$ is a well-defined mean zero Gaussian process.

By applying Remarks \ref{dif-det} and \ref{rdrem2.5.2} to the deterministic function $t\mapsto z^{2n-1} e^{-tz^2}$, and by setting
$R_t^{(0)}:=R_t$, we infer that for any $n\in\mathbb{N}$, $(R_t^{(n)},\, t\in\, ]0,\infty[)$ has a continuously  differentiable version and
$\frac{d}{dt}R_t^{(n)}= R_t^{(n+1)}$. This completes the proof of part (a).

(b) 
We assume here that the space-time white noise that appears in \eqref{4.20-1'} is defined on $\R \times \R$, rather than only on $\R_+ \times \R$. Fix $t > 0$ and, for $x \in \R$, define
\begin{align}
\label{ch1'-s3-107}
S_x(t) &= 
   \int_{-\infty}^0 \int_\R (\Gamma(t-r,x-z) - \Gamma(t-r,-z))\, W(dr, dz).
\end{align}
This stochastic integral is well-defined because
\begin{align*}
   &\int_{-\infty}^0 dr \int_\R dz\, (\Gamma(t-r,x-z) - \Gamma(t-r,-z))^2 \\
   &\qquad\qquad =  \int_{t}^\infty ds \int_\R dz\, (\Gamma(s,x-z) - \Gamma(s,-z))^2 
    < \infty
\end{align*}
by \eqref{1'.3}.

We will show that
\beqn
B_x(t) = \sqrt{2} \left(v(t,x)-v(t,0)-S_x(t)\right), \qquad  x\in\re,
\eeqn
is a two-sided Brownian motion and that $\left(S_x(t),\,  x\in\re\right)$ has the required properties.

Observe that for $h\in \re$,
\begin{align*}
&E\left[\left(S_{x+h}(t)-S_x(t)\right)^2\right] \\
&\qquad\qquad = \int_{-\infty}^0 dr \int_\R dz\, (\Gamma(t-r,x+h-z) - \Gamma(t-r,x-z))^2 \\
&\qquad\qquad  = \int_t^\infty ds \int_{\re} dy\, \left[\Gamma(s,y+h) - \Gamma(s,y)\right]^2.
\end{align*}

Therefore,
\begin{align*}
&E\left[\left(v(t,x+h) - v(t,x)\right)^2\right] \\
 &\qquad = \int_0^t ds \int_{\re} dy \left(\Gamma(t-s,x+h-y) - \Gamma(t-s,x-y)\right)^2\\
&\qquad  =  \int_0^t ds \int_{\re} dy  \left(\Gamma(s,y+h) - \Gamma(s,y)\right)^2\\
&\qquad  =  \int_0^\infty ds \int_{\re} dy  \left(\Gamma(s,y+h) - \Gamma(s,y)\right)^2
 - E\left[\left(S_{x+h}(t)-S_x(t)\right)^2\right].
 \end{align*}

 By using the independence of $v$ and $(S_x(t))$,
  we obtain
 \begin{align}\nonumber
 E\left[\left(B_{x+h}(t) - B_x(t)\right)^2\right] 
   &= 2 \int_0^\infty ds \int_{\re} dy \left(\Gamma(s,y+h)-\Gamma(s,y)\right)^2\\
 & =|h|,
 \label{rd04_17e1}
 \end{align}
by \eqref{1'.3}. Since $B_0(t)=0$, this shows that $(B_x(t),\,  x\in\re)$ is a standard two-sided Brownian motion and \eqref{ch1'-s3-101} holds.

 It remains to check that $(S_x(t), \, x\in \re)$ has a version with $\mathcal{C}^\infty$ sample paths. However, this follows from \eqref{ch1'-s3-107} and the same arguments as those used to check the differentiability of $(R_t)$, because the domain of integration $\R_- \times \R$ does not contain the two singularities of $(r, z) \mapsto \Gamma(t-r,x-z) - \Gamma(t-r,-z)$, and we can use the bounds in Lemma \ref{rdlemC.2.2}. This proves (b).
\qed
\medskip

\noindent{\em The stationary pinned string}
\medskip

In view of equality \eqref{rd04_17e1}, it is natural to expect that if the initial condition $u_0$ is random and has the law of a two-sided standard Brownian motion multiplied by $2^{-1/2}$, then $u(t, *) - u(t, 0)$ also has this property at all time $t > 0$. This property was used in \cite{m-t2002} and is the object of the next proposition. In Subsection \ref{rd03_06ss1}, we will establish a similar result for the solution to the stochastic heat equation on $[0, L]$ with Neumann boundary conditions.

\begin{prop}
\label{rd04_23p1}
Assume that $u_0 = (u_0(x), \, x \in \R)$ is a stochastic process independent of the space-time white noise $\dot W$ in the SPDE \eqref{4.20-1'}, with the law of a two-sided standard Brownian motion multiplied by $2^{-1/2}$. Let $u$ be the solution to \eqref{4.20-1'} given in \eqref{sollinheat-1'} with this initial condition $u_0$. Then for all $t \geq 0$, the process $(u(t,x) - u(t, 0),\, x \in \R)$  also has the law of a two-sided standard Brownian motion multiplied by $2^{-1/2}$.
\end{prop}

\begin{proof}
In order to create a two-sided standard Brownian motion multiplied by $2^{-1/2}$ that is independent of $\dot W$ on $\R_+ \times \R$, we assume that the space-time white noise $\dot W$ is defined on $\R \times \R$ and for $x \in \R$, we set
\beqn
    u_0(x) := \int_{-\infty}^0 \int_\R (\Gamma(-r,x-z) - \Gamma(-r,-z))\, W(dr, dz).
\eeqn 
Then $u_0(0) = 0$ a.s., and by \eqref{rd04_17e1}, $E[(u_0(x) - u_0(y))^2] = \half\, \vert x - y \vert$, therefore $u_0$ is a standard two-sided Brownian motion multiplied by $2^{-1/2}$. We work with its continuous modification, that we also denote $u_0$.

The solution of \eqref{4.20-1'} is
\beqn
    u(t, x) = I_0(t, x) + \int_0^t \int_\R \Gamma(t-s,x-z) W(ds, dz),
\eeqn
where
\begin{align*}
   I_0(t, x) &= \int_\R \Gamma(t,x-y) u_0(y) \, dy \\
   & = \int_\R dy \,  \Gamma(t,x-y) \int_{-\infty}^0 \int_\R (\Gamma(-r,y-z) - \Gamma(-r,-z))\, W(dr, dz).
\end{align*}
Using the stochastic Fubini's theorem Theorem \ref{ch1'-tfubini}, whose assumptions are easily seen to be satisfied, and the semigroup property \eqref{semig-heat} of $\Gamma$, we see that 
\beqn
    I_0(t, x) = \int_{-\infty}^0 \int_\R (\Gamma(t-r,x-z) - \Gamma(-r,-z))\, W(dr, dz).
\eeqn
Therefore, for $ x \in \R$,
\begin{align*}
   u(t, x) - u(t, 0) &=  \int_{-\infty}^0 \int_\R (\Gamma(t-r,x-z) - \Gamma(t-r,-z))\, W(dr, dz) \\
   &\qquad + \int_0^t \int_\R (\Gamma(t-s,x-z) - \Gamma(t-s,-z)) W(ds, dz),
\end{align*} 
and for $x, y \in \R$,
\begin{align*}
  & u(t, x) - u(t, 0) - (u(t, y) - u(t, 0)) \\
  &\qquad\qquad =  \int_{-\infty}^0 \int_\R (\Gamma(t-r,x-z) - \Gamma(t-r,y-z))\, W(dr, dz) \\
   &\qquad\qquad\qquad + \int_0^t \int_\R (\Gamma(t-s,x-z) - \Gamma(t-s,y-z))\, W(ds, dz).
\end{align*} 
Using the independence properties of space-time white noise, we see that the second moment of this expression is equal to
\begin{align*}
    & \int_{-\infty}^0 \int_\R (\Gamma(t-r,x-z) - \Gamma(t-r,y-z))^2\, dr dz \\
   &\qquad\qquad  + \int_0^t \int_\R (\Gamma(t-s,x-z) - \Gamma(t-s,y-z))^2\, ds dz \\
   &\qquad =  \int_0^\infty \int_\R (\Gamma(s,x-z) - \Gamma(s,y-z))^2\, ds dz \\
   &\qquad = \half\, \vert x - y \vert
\end{align*}
by \eqref{rd04_17e1}. Therefore, $u(t, *) - u(t, 0)$ is a standard two-sided Brownian motion multiplied by $2^{-1/2}$.
\end{proof}


\section{The stochastic heat equation on a bounded real interval}
\label{ch4-ss2.2-1'}

In this section, we study the stochastic heat equation on the interval $[0,L]$, with additive noise and with either Dirichlet or Neumann vanishing boundary conditions.
As we did for the stochastic heat equation on $\re$, we discuss the existence of random field solutions and the regularity of their sample paths.

\subsection{Existence of a random field solution}
\label{ch1'-ss3.4.0}

We will examine the validity of Assumption \ref{ch4-a1-1'} for the Green's functions of the two examples considered in this section.
\medskip


\noindent{\em Dirichlet boundary conditions}
\smallskip

We consider the SPDE \eqref{ch1.77} from Chapter \ref{ch1}:
\begin{equation}
\label{ch1'.HD}
\begin{cases}
 \frac{\partial u}{\partial t}- \frac{\partial^{2}u}{\partial x^{2}} = \dot W(t, x),\qquad (t, x) \in\, ]0,\infty[ \times\, ]0, L[\ ,\\
 u(0, x) = u_{0}(x), \qquad\qquad x \in \, [0,L]\ ,\\
 u(t, 0) = u(t, L) =0, \qquad t \in \, ]0,\infty[\ ,
 \end{cases}
\end{equation}
with $u_0\in L^1([0,L])$.

The associated Green's function, which has been calculated in Chapter \ref{ch1}  (see \eqref{ch1.600}), takes the form $\Gamma(t,x;s,y)=G_L(t-s;x,y)$, where
\beq
\label{ch1'.600}
  G_L(t;x,y)= \sum^{\infty}_{n=1} e^{-\frac{\pi^2}{L^2}n^{2} t}e_{n,L} (x) e_{n,L} (y),  \quad t>0, \ x,y \in[0,L],
\eeq
and $e_{n,L}(x) := \sqrt{\frac{2}{L}} \sin\left(\frac{n\pi}{L}x\right)$, $n\in \IN^*$. For an equivalent expression, see Lemma \ref{ch1-equivGD}.

It is a consequence of \eqref{27.1-Ch1} that
\beq
\label{monyo(*1)}
  \sup_{x \in [0,L]} \Vert G_L(\cdot, x, *) \Vert_{L^2(\R_+ \times [0,L])} < \infty.
  \eeq
In particular, Assumption \ref{ch4-a1-1'} is satisfied for $\Gamma(t,x; s,y) := G_L(t-s; x, y)$.

For its further use, we collect in the next proposition some relevant properties of $G_L(t;x,y)$.
\begin{prop}
\label{ch1'-pPD}
The Green's function $G_L(t;x,y)$ defined in \eqref{ch1'.600} satisfies the following.
\begin{description}
\item {(i)} {\em Semigroup property}. For any $s,t >0$, $x,z\in [0,L]$,
\beq
\label{semig}
\int_0^L dy\ G_L(s;x,y)G_L(t;y,z) = G_L(s+t;x,z).
\eeq
\item{(ii)} {\em Comparison with the heat kernel}. For any $t>0$, $x,y\in[0,L]$,
\beqn
0\le G_L(t;x,y)\le \Gamma(t, x-y),
\eeqn
where $\Gamma(s,y)$ is defined in \eqref{heatcauchy-1'}, and $G_L(t,x,y)>0$ for $x,y\in\, ]0,L[$.
\item {(iii)} {\em Sub-density}. For $t>0$, $x\in[0,L]$,
\beqn
\int_0^L dy\ G_L(t;x,y)<1.
\eeqn
\item {(iv)} {\em Scaling property}. For any $t>0$, $x,y\in[0,L]$,
\beq
\label{scalingD}
G_L(t;x,y) = \frac{1}{L}G_1\left(\frac{t}{L^2};\frac{x}{L},\frac{y}{L}\right).
\eeq
\end{description}
\end{prop}
\begin{proof}
Property (i) follows easily by computing the integral on $[0,L]$ of the left-hand side of \eqref{semig} and using the expression \eqref{ch1'.600} of the Green's function.
As for properties (ii) and (iii), one can give a probabilistic argument. Indeed, from classical connections between Brownian motion and PDEs, it is well-known that
the function
\beqn
]0,\infty[\times [0,L]\times [0,L] \ni (t,x,y) \mapsto G_L(t;x,y)
\eeqn
 is the transition density of a Brownian motion killed upon reaching the boundary points $0, L$
(see, e.g. \cite[Section 2.8]{ks}). This yields property (ii) above.

Property (iii) is a direct consequence of (ii), since $\int_0^L\Gamma(t,x-y)\ dz < 1$. Property (iv) follows immediately from \eqref{ch1'.600}.
 \end{proof}

 \begin{remark}
 \label{rem: parabolic-ub}
 Property (ii) holds in the more general setting of a class of parabolic boundary value problems in which the Laplacian operator is replaced by a uniformly elliptic partial differential operator. For more details, we refer the reader to \cite[Theorem 1.1]{ei70}.
\end{remark}

Property (ii) above along with \eqref{heatcauchy-11'} yield
 \beq
\label{1'.12}
\int_0^t dr \int_0^L dy\ G_L^2(r;x,y)\le \left(\frac{t}{2\pi}\right)^{\half}.
\eeq
This is another way of checking that $\Gamma(t,x;s,y):=G_L(t-s;x,y)$ satisfies Assumption \ref{ch4-a1-1'}.

Consider the homogeneous PDE corresponding to \eqref{ch1'.HD} (that is, with $\dot W(t,x)$ there replaced by $0$), whose solution is the function
\beq
\label{cor1.0}
I_0(t,x)=
\begin{cases}
\int_0^L G_L(t;x,z)\ u_{0}(z)\, dz ,& (t,x)\in\ ]0,\infty[\ \times [0,L],\\
u_0(x), & (t,x)\in \{0\}\times [0,L],
\end{cases}
\eeq
where $G_L(t;x,z)$ is given in \eqref{ch1'.600}. Notice that the integral is well-defined because $u_0 \in L^1([0, L])$.

 We can use
Definition \ref{ch4-d1-1'} to say that the random field solution to \eqref{ch1'.HD} is the random field $u=(u(t,x),\ (t,x)\in\re_+\times [0,L])$
given by
\beq
\label{c}
u(t, x)= I_0(t,x) + \int_0^t \int_0^L G_L(t-s; x,y)\, W(ds,dy).
\eeq
Observe that $(u(t,x),\  (t,x)\in\re_+\times [0,L])$ defines a Gaussian process, and that the values $u(t,0)$ and $u(t,L)$ given by this formula are compatible with the homogeneous Dirichlet boundary conditions.
\smallskip

By Definition \ref{ch1'-s4.d1} of the stochastic integral and \eqref{ch1'.600}, we see that
\beq
\label{series-dirichlet}
u(t,x)= \sum_{n=1}^\infty Y_t(n)\, e_{n,L}(x),\quad (t,x)\in\,[0,\infty[ \times [0,L],
\eeq
with
\beq
\label{OU-coordinates}
Y_t(n)= e^{-\frac{\pi^2}{L^2}n^2 t} \langle u_0,e_{n,L}\rangle_V + \int_0^t e^{-\frac{\pi^2}{L^2}n^2 (t-s)}\, dW_s(e_{n,L}).
\eeq
Notice that $(Y_t(n),\, t\in[0,\infty[)$ is the Ornstein-Uhlenbeck process\index{Ornstein-Uhlenbeck process}\index{process!Ornstein-Uhlenbeck} solution to the stochastic differential equation
\begin{align*}
dY_t(n) &= -\frac{\pi^2}{L^2}\, n^2\, Y_t(n) + dW_t(e_{n,L}),\qquad
Y_0(n)  = \langle u_0,e_{n,L}\rangle_V.
\end{align*}
This representation will play a role in Section \ref{ch5-ss-in-rev}.
\medskip

 \noindent{\em Neumann boundary conditions}
\smallskip

Consider the SPDE on $\re_+\times [0,L]$ with Neumann boundary conditions:
\begin{equation}
\label{1'.14}
\begin{cases}
 \frac{\partial u}{\partial t}- \frac{\partial^{2}u}{\partial x^{2}} = \dot W(t, x),\qquad (t, x) \in\, ]0,\infty[ \times\, ]0, L[,\\
 u(0, x) = u_{0}(x),\qquad\qquad x \in \,[0,L], \\
\frac{\partial u}{\partial x}(t, 0) = \frac{\partial u}{\partial x}(t, L) =0, \quad t \in \,]0,\infty[,
 \end{cases}
\end{equation}
with $u_0\in L^1([0,L])$.

As we have seen in Section \ref{ch1-s3}, the Green's function takes the form
\beqn
\Gamma(t,x;s,y)=G_L(t-s;x,y),
\eeqn
with
\beq
\label{1'.400}
 G_L(t;x,y)= \sum^{\infty}_{n=0} e^{-\frac{\pi^2}{L^2}n^{2} t}g_{n,L} (x) g_{n,L} (y),\quad t>0, \ x,y\in[0,L],
\eeq
where
\beq
\label{basis-neumann}
g_{0,L}(x)= \frac{1}{\sqrt{L}}, \qquad  g_{n,L}(x)= \sqrt{\frac{2}{L}} \cos\left(\frac{n\pi}{L}x\right), \qquad  n\ge1.
\eeq
For an equivalent expression, see Lemma \ref{ch1-lequivN}.

 It is a consequence of \eqref{cansada} that
$\sup_{x\in [0,L]}\Vert G_L(\cdot; x,\ast)\Vert_{L^2([0,T]\times [0,L])}  < \infty$.

\begin{prop}
\label{ch1'-pPN}
The Green's function $G_L(t;x,y)$ satisfies the following properties:
\begin{description}
\item {(a)} {\em Semigroup property}. For any $s,t >0$, $x,z\in [0,L]$,
\beqn
\int_0^L dy\ G_L(s;x,y) G_L(t;y,z) = G_L(s+t;x,z).
\eeqn
\item{(b)} {\em Comparison with the heat kernel}. There exists a constant $c_1(L)$ such that for all $t>0$, and $x,y\in[0,L]$,
\beq
\label{neumann-bad}
0< \Gamma(t,x-y)\le G_L(t;x,y)
\le c_1(L)\left(\frac{1}{\sqrt t}\vee 1\right)\exp\left(-\frac{|x-y|^2}{8t}\right).
\eeq

\item{(c)} {\em Full density}. For any $t>0$, $x\in[0,L]$,
\beq
\label{int1}
\int_0^L G_L(t; x,y)\ dy = 1.
\eeq

\item{(d)} {\em Scaling property}. For any $t>0$, $x, y\in[0,L]$,
\beq
\label{scalingN}
G_L(t;x,y) = \frac{1}{L}G_1\left(\frac{t}{L^2};\frac{x}{L},\frac{y}{L}\right).
\eeq
\end{description}
\end{prop}
\begin{proof}
Property (a) can be easily checked by using the expression \eqref{1'.400}.

The lower bound $G_L(t;x,y)\ge \Gamma(t,x-y)$ in part (b) follows from \eqref{1'.15-double}.
The last upper bound in \eqref {neumann-bad}, follows from \cite[Theorem 3.2.9, p. 90]{davies}.

 The proof of \eqref{int1} follows by integrating term-by-term the expression of $G_L(t;x,y)$ given in \eqref{1'.400}: for $n\ge 1$, the integrals vanish and for $n=0$, the integral is $1$.
 The scaling property follows immediately from \eqref{1'.400}
\end{proof}
\begin{remark}
\label{3.3.9-rem-newmann}
Fix $T > 0$. By \eqref{neumann-bad}, there is a constant $c(T,L)$ such that
 \beq
 \label{compneumann}
 G_L(t;x,y) \le \frac{c(T,L)}{\sqrt t} \exp\left(-\frac{|x-y|^2}{8t}\right)
  \eeq
 for all $t\in\, ]0,T[$, $x,y\in[0,L]$.

In particular,
 \beqn
 \int_0^t dr \int_0^L dy\, G_L^2(r;x,y)\le \tilde c(t,L) t^{\half}
 \eeqn
 and therefore $\Gamma(t,x;s,y):= G_L(t-s;x,y)$ satisfies Assumption \ref{ch4-a1-1'}.
 \end{remark}

 Consider the homogeneous PDE corresponding to \eqref{1'.14} (that is, with $\dot W(t,x)$ there replaced by $0$), whose solution is the function
 \beq
\label{cor1.0-N}
I_0(t,x)=
\begin{cases}
\int_0^L G_L(t;x,z)\ u_{0}(z)\, dz ,& (t,x)\in\ ]0,\infty[\ \times [0,L],\\
u_0(x), & (t,x)\in \{0\}\times [0,L],
\end{cases}
\eeq
where $G_L(t;x,y)$ is given in \eqref{1'.400}. Observe that the integral is well-defined because $u_0 \in L^1([0, L])$.

According to Definition \ref{ch4-d1-1'}, a random field solution to \eqref{1'.14} is the random field $u=(u(t,x),\ (t,x)\in\re_+\times [0,L])$
given by
\beq
\label{1'.N}
u(t, x)= I_0(t,x)+ \int_0^t \int_0^L G_L(t-s; x,y)\, W(ds,dy),
\eeq
with $G_L(t;x,y)$ defined by \eqref{1'.400} (or by \eqref{1'.15-double}).
 Notice that for  $t > 0$, the boundary values $u(t, 0)$ and $u(t, L)$ are given by the right-hand side.

As in the case of Dirichlet boundary conditions, the random field solution \eqref{1'.N} admits another representation, namely
\beq
\label{series-neuman}
u(t,x)= \sum_{n=0}^\infty Z_t(n)\, g_{n,L}(x),\quad (t,x)\in\,]0,\infty[ \times [0,L],
\eeq
with $g_{n,L}$, $n\ge 0$, defined in \eqref{basis-neumann}. For each $n\ge 1$,  $(Z_t(n),\ t\in\,]0,\infty[)$ is the Ornstein-Uhlenbeck process\index{Ornstein-Uhlenbeck process}\index{process!Ornstein-Uhlenbeck}
\beqn
Z_t(n)=  e^{-\frac{\pi^2}{L^2}n^2 t} \langle u_0,g_{n,L}\rangle_V + \int_0^t e^{-\frac{\pi^2}{L^2}n^2 (t-s)}\, dW_s(g_{n,L})
\eeqn
(observe the analogy with \eqref{OU-coordinates}),
while for $n=0$,
\beqn
Z_t(0) = \frac{1}{\sqrt L} \int_0^L u_0(z)\, dz + W_t(g_{0,L}).
\eeqn
This representation is used in Proposition \ref{rem6.1.15-*1}.



\subsection{Relationship with the stochastic heat equation on $\re$}
\label{ch1'-ss3.4.1}

We will see that there is a relationship between solutions to the linear stochastic heat equation
on the real line, considered in Section \ref{ch4-ss2.1-1'}, and of the same equation on a bounded interval. This is due to the relationship between the fundamental solution $\Gamma(t,x-y)$ (see \eqref{heatcauchy-1'}) and the expressions of the Green's functions given in \eqref{ch1.6000-double} and \eqref{1'.15-double} for vanishing Dirichlet and Neumann boundary conditions, respectively.

\begin{lemma}
\label{ch1'-ss2.1-l40}
Let $G_L(t;x,y)$ be either as in \eqref{ch1.6000-double} or in \eqref{1'.15-double}. Fix a closed interval $J\subset\ ]0,L[$. The function
$H_1: \R\times J\times [0,L]\longrightarrow \re$ defined by
\beq
\label{HH-previous}
H_1(t;x,y)
= \begin{cases}
G_L(t;x,y) - \Gamma(t,x-y), & t>0,\\
0,  & t\le 0,
\end{cases}
\eeq
is $\mathcal{C}^\infty$ on $\R \times J\times [0,L]$, and for
$n_1, n_2, n_3 \in \mathbb{N}$, there is $c_{n_1, n_2, n_3, J}<\infty$ such that
\beq
 \label{ch1'-tDL.3}
\sup_{(t,x,y)\in \R\times J\times [0,L]}\left\vert\frac{\partial^{n_1+n_2+n_3}}{\partial t^{n_1}\partial x^{n_2}\partial y^{n_3}} H_1(t;x,y)\right\vert \le c_{n_1, n_2, n_3, J}.
\eeq
\end{lemma}
\begin{proof}
We consider only the case where $G_L(t;x,y)$ is defined in \eqref{ch1.6000-double} (Dirichlet case), since the Neumann case is similar.

From \eqref{ch1.6000-double}, we see that for all $t\in\re$,
\beqn
H_1(t;x,y) = -\tilde H_0(t,x,y) + \tilde H_1(t; x,y),
\eeqn
where
\beqn
\tilde H_0(t,x,y) =
\frac{1}{\sqrt{4 \pi t}}\left[ \exp \left(-\frac{{(x+y)}^{2}}{4t}\right) + \exp \left(-\frac{{(x+y-2L)}^{2}}{4t}\right)\right] 1_{\{t>0\}},
\eeqn
and
 \begin{align*}
 \tilde H_1(t; x,y) &= \frac{1}{\sqrt{4 \pi t}} \left[ \sum_{m\in\mathbb{Z}\setminus \{0\}} \exp \left(-\frac{{(y-x+2mL)}^{2}}{4t}\right)\right.\\
 &\left. \qquad \qquad -
 \sum_{m\in\mathbb{Z}\setminus \{0,-1\}}\exp \left(-\frac{{(y+x+2mL)}^{2}}{4t}\right) \right]1_{\{t>0\}}.
\end{align*}

Assume $J=[\varepsilon_0,L-\varepsilon_0]$ for some $\varepsilon_0>0$. Then for $(x,y)\in J\times [0,L]$,
\beqn
x+y\ge \varepsilon_0, \ {\text{and}}\ x+y-2L\le L-\varepsilon_0+L-2L = -\varepsilon_0.
\eeqn
Therefore, the function $\tilde H_0$ 
is $\mathcal{C}^\infty$ on $\R\times J\times [0,L]$, the same is true for $\tilde H_1(t; x,y)$ 
 and for
$n_1, n_2, n_3 \in \mathbb{N}$, there is $c_{n_1, n_2, n_3, J}<\infty$ such that
\beqn
\sup_{(t,x,y)\in \R \times J\times [0,L]}\left\vert\frac{\partial^{n_1+n_2+n_3}}{\partial t^{n_1}\partial x^{n_2}\partial y^{n_3}} H_1(t;x,y)\right\vert \le c_{n_1, n_2, n_3, J}.
\eeqn
This ends the proof.
\end{proof}

\begin{remark}
\label{abans-ch1'-tDL-*1}
 Observe that if in \eqref{ch1'-tDL.3}, we replace $H_1$ by $\tilde H_1$, then \eqref{ch1'-tDL.3} remains valid even if the supremum is taken over $\R \times \R \times [0, L]$. In certain cases, the decomposition
 \beqn
   G_L(t; x, y) = \Gamma(t, x - y) \pm \Gamma(t, x + y) \pm \Gamma(t, x + y - 2L) + \tilde H_1(t; x, y)
   \eeqn
may be useful (see Lemma \ref{app2-r1}).
\end{remark}

\begin{thm}
\label{ch1'-tDL}
Let $u_D = (u_D(t,x), (t,x)\in\re_+\times [0,L])$ be the random field solution to \eqref{ch1'.HD} with initial condition $u_{D,0}$ and  $u=(u(t,x), (t,x)\in\re_+\times \re)$ be the random field solution to \eqref{4.20-1'} with initial condition $u_0$ satisfying \eqref{ch1'-v00}, where, in both SPDEs, we use the same space-time white noise $(\dot W(t,x), (t,x)\in\re_+\times \re)$. Fix compact intervals $I\subset \ ]0,\infty[$ and $J\subset \ ]0,L[$. Then on $I\times J$, the random field $v_D=u_D-u$ has $\mathcal{C}^\infty$ sample paths, and for all $n_1, n_2\in \mathbb{N}$, there is a constant $c_{n_1, n_2, I, J} < \infty$ such that
\beq
\label{ch1'-tDL.1}
\sup_{(t,x)\in I\times J}E\left[\left(\frac{\partial^{n_1+n_2}}{\partial t^{n_1}\partial x^{n_2}} v_D(t,x)\right)^2\right] \le c_{n_1, n_2, I, J},
\eeq

A similar statement holds with $u_D$ replaced by the random field solution $u_N = (u_N(t,x), (t,x)\in\re_+\times [0,L])$ to \eqref{1'.14} and $v_D$ replaced by $v_N=u_N-u$.
  \end{thm}
 \begin{proof}
 We only consider the case $v_D=u_D-u$ (Dirichlet boundary conditions), since the case with Neumann boundary conditions is similar.

 Let $G_L(t;x,y)$ be as in \eqref{ch1.6000-double} and let
 \begin{align*}
 I_{D,0}(t,x) &= \int_0^L G_L(t;x,y) u_{D,0}(y)\ dy,\\
 I_0(t,x) &= \int_{\re} \Gamma(t,x-y) u_0(y)\ dy.
 \end{align*}
 Since
 $I_{D,0}(\cdot, \ast)$ and $I_0(\cdot, \ast)$ are $\mathcal{C}^\infty$ on $]0, \infty[\times J$, we can and will assume that $u_{D,0}\equiv 0$ and $u_0\equiv 0$. With these vanishing initial conditions, we will in fact show that $v_D$ is $\mathcal{C}^\infty$ on $]0, \infty[\times J$, and for all $T>0$, $n_1, n_2\in \mathbb{N}$, there is a constant
 $c_{n_1, n_2, I, J}$ such that
 \beq
 \label{ch1'-tDL.2}
\sup_{(t,x)\in I\times J}E\left[\left(\frac{\partial^{n_1+n_2}}{\partial t^{n_1}\partial x^{n_2}} v_D(t,x)\right)^2\right] \le c_{n_1, n_2, I, J}.
\eeq
Consider the function $H_1: \R\times J\times [0,L]\longrightarrow \re$ given in \eqref{HH-previous}. Define $H_2: \R\times J\times [0,L]^c \longrightarrow \re$ by
\beqn
H_2(t;x,y)=
\begin{cases}
 \Gamma(t,x-y), & t>0,\\
 0, &t\le 0.
 \end{cases}
 \eeqn
Recall Lemma \ref{ch1'-ss2.1-l40} and, in particular, \eqref{ch1'-tDL.3}. 

The function $H_2$ is also $\cC^\infty$ on $\R\times J\times [0,L]^c$.
Moreover, appealing to Lemma \ref{rdlemC.2.2}, there are constants $0<C=C_{n_1,n_2,n_3,I,J}$ and $0<c=c_{n_1,n_2,n_3,I,J}$ such that for any $(t,x,y)\in ]0,T]\times \re^2$,
\begin{align}
 \label{ch1'-tDL.30-bis}
\left\vert\frac{\partial^{n_1+n_2+n_3}}{\partial t^{n_1}\partial x^{n_2}\partial y^{n_3}} H_2(t;x,y)\right\vert &\le
C\, t^{-(1+2n_1+n_2+n_3)/2} \exp\left(c\ \frac{|x-y|^2}{t}\right)\notag\\
& = C\, t^{-(2n_1+n_2+n_3)/2} \Gamma \left(t/(4c), x-y\right).
\end{align}
This inequality also holds for $t\le 0$, with the convention that the right-hand side is $0$.

Assume that $I=[t_0,t_1]$, $0<t_0<t_1$.
For all  $(t,x)\in \R_+\times [0,L]$,
\begin{align*}
v_D(t,x) &= \int_0^t \int_{[0,L]} H_1(t-s; x,y)\ W(ds,dy)\\
&\qquad -  \int_0^t \int_{[0,L]^c} H_2(t-s; x,y)\ W(ds,dy), \quad a.s.
\end{align*}
However, because $H_i(t-s;x,y)=0$ if $s\ge t$ ($i=1,2$), for all $(t,x)\in I\times J$, we can write
\begin{align*}
v_D(t,x) &= \int_0^{t_1} \int_{[0,L]} H_1(t-s; x,y)\ W(ds,dy)\\
&\qquad -  \int_0^{t_1} \int_{[0,L]^c} H_2(t-s; x,y)\ W(ds,dy), \quad a.s.
\end{align*}
Using the properties of $H_i$,  $i=1,2$,  in particular \eqref{ch1'-tDL.3} and \eqref{ch1'-tDL.30-bis}, we can apply Theorem \ref{ch1'-tdif} (and Remark \ref{rdrem2.5.2}) on differentiation under the stochastic integral to the deterministic family of functions depending on the parameter $(t,x) \in\ ]0,\infty[\times J$,
\beqn
(s,y) \mapsto H_i(t-s; x,y), \qquad  i=1,2,
\eeqn
to differentiate successively several times with respect to $t$ and $x$, and to obtain
\begin{align}
\label{ch1'-tDL.4}
\frac{\partial^{n_1+n_2}}{\partial t^{n_1}\partial x^{n_2}} v_D(t,x) &=\int_0^{t_1} \int_{[0,L]} \frac{\partial^{n_1+n_2}}{\partial t^{n_1}\partial x^{n_2}} H_1(t-s;x,y)\ W(ds,dy)\notag\\
& \qquad - \int_0^{t_1} \int_{[0,L]^c} \frac{\partial^{n_1+n_2}}{\partial t^{n_1}\partial x^{n_2}} H_2(t-s;x,y)\ W(ds,dy),
\end{align}
$(t,x)\in[0,t_1]\times J$.

By the It\^o isometry and \eqref{ch1'-tDL.3},
\begin{align*}
\sup_{(t,x)\in [0,T]\times J} &E\left[\left(\int_0^{t_1} \int_{[0,L]}  \frac{\partial^{n_1+n_2}}{\partial t^{n_1}\partial x^{n_2}} H_1(t-s;x,y)\ W(ds,dy)\right)^2\right]\\
&\qquad  \le c^{(1)}_{n_1,n_2,T,J}.
\end{align*}
Furthermore,
\begin{align}
\label{ch1'-tDL.4-bis}
&\sup_{(t,x)\in I\times J} E\left[\left(\int_0^{t_1} \int_{[0,L]^c}  \frac{\partial^{n_1+n_2}}{\partial t^{n_1}\partial x^{n_2}} H_2(t-s;x,y)\ W(ds,dy)\right)^2\right]\notag\\
&\qquad\qquad= \sup_{(t,x)\in I\times J} \int_0^t ds \int_{[0,L]^c}  dy\, \left(\frac{\partial^{n_1+n_2}}{\partial t^{n_1}\partial x^{n_2}} H_2 (t-s,x,y)\right)^2\notag\\
&\qquad\qquad= \sup_{(t,x)\in I\times J} \int_0^t ds \int_{[0,L]^c}  dy\, \left(\frac{\partial^{n_1+n_2}}{\partial t^{n_1}\partial x^{n_2}} H_2 (s,x,y)\right)^2.
\end{align}
Assuming that $J=[\varepsilon_0, L-\varepsilon_0]$ for some $\varepsilon_0>0$, for any $x\in J$ and $y\in[0,L]^c$, we have $|x-y| \ge \varepsilon_0$. Thus, using \eqref{ch1'-tDL.30-bis} (with $n_3=0)$ and \eqref{rdlemC.2.2-before-2}, yields
\begin{align*}
&\int_{[0,L]^c}  dy\, \left(\frac{\partial^{n_1+n_2}}{\partial t^{n_1}\partial x^{n_2}} H_2 (s,x,y)\right)^2\\
 &\qquad \le C s^{-(2n_1+n_2)} \int_{[0,L]^c}  dy\, \Gamma^2(s/4c, x-y)\\
&\qquad \le C s^{-(2n_1+n_2)} \int_{\varepsilon_0}^\infty dz\, \Gamma^2(s/4c, z)\\
&\qquad \le C s^{-(2n_1+n_2)} \exp\left(- \frac{2c\varepsilon_0^2}{s}\right).
\end{align*}
Using this estimate in \eqref{ch1'-tDL.4-bis}, we obtain
\begin{align*}
&\sup_{(t,x)\in I\times J} E\left[\left(\int_0^{t_1} \int_{[0,L]^c}  \frac{\partial^{n_1+n_2}}{\partial t^{n_1}\partial x^{n_2}} H_2(t-s;x,y)\ W(ds,dy)\right)^2\right]\notag\\
&\qquad \le C\sup_{t\in I} \int_0^t ds\, s^{-(2n_1+n_2)} \exp\left(- \frac{2c\varepsilon_0^2}{s}\right) <\infty,
\end{align*}
because the integrand is bounded and $I$ is a compact set.

Hence, \eqref{ch1'-tDL.2} follows and the theorem is proved.
 \end{proof}
 \bigskip


\subsection[Local H\"older continuity of the sample paths away \texorpdfstring{\\}{} from the space-time boundary]{Local H\"older continuity of the sample paths away from the space-time boundary}
 \label{ch1'-ss3.4.2}

Building on the results of Section \ref{ch1'-ss3.4.1}, we obtain a first statement on local H\"older continuity.
In the next proposition, $u = (u(t,x),\ (t,x)\in\re_+\times [0,L])$ denotes the random field solution to the stochastic heat equation with either Dirichlet or Neumann vanishing boundary conditions.
\begin{prop}
\label{4-p1-1'-weak}
Assume that the initial condition $u_0$ in \eqref{ch1'.HD} and in \eqref{1'.14} is such that $u_0\in L^1([0,L])$.
Fix compact intervals $I\subset\ ]0,\infty[$ and $J\subset\ ]0,L[$.
Then the random field $u = (u(t,x),\, (t,x)\in\re_+\times [0,L])$ satisfies the following property:

Let  $p\in[2,\infty[$; there exists a constant $C>0$ (depending on $p$) such that, for any $(t,x), (s,y)\in I\times J$,
\beq
\label{4-p1-1'-weak.0}
E\left[\left\vert u(t,x) - u(s, y)\right\vert^p\right] \le C \left(|t-s|^{\frac{1}{4}} + |x-y|^{\frac{1}{2}}\right)^p.
\eeq
As a consequence, for $\alpha \in ]0, \frac{1}{4}[$ and $\beta \in ]0, \half[$, there is a version of the random field
$(u(t,x), (t,x) \in\, ]0,\infty[ \times ]0,L[)$
 with locally jointly Hölder continuous sample paths with exponents $(\alpha, \beta)$.
\end{prop}
\begin{proof}
We give the proof in the case of Dirichlet boundary conditions (equation \eqref{ch1'.HD}): we write therefore $u_D$ instead of $u$. For Neumann boundary conditions (equation \eqref{1'.14}), the proof is analogous.

Denote by $v$ the
solution to \eqref{4.20-1'}  (stochastic heat equation on $\re$) with vanishing initial conditions. Let
\beqn
v_D(t,x) = u_D(t,x)-v(t,x), \qquad  (t,x)\in \IR_+ \times [0,L].
\eeqn
For $(t,x), (s,y)\in I\times J$, assuming without loss of generality that $t\le s$ and $x\le y$,
consider the path in $\re_+\times [0,L]$ consisting of two segments, one from $(t,x)$ to $(s,x)$, the other from $(s,x)$ to $(s,y)$. Applying Theorem \ref{ch1'-tDL}, we obtain
\beqn
E\left[\left( v_D(t,x) - v_D(s,y)\right)^2\right] = E\left[\left\vert \int_t^s \frac{\partial v_D}{\partial r}(r,x)\ dr + \int_x^y \frac{\partial v_D}{\partial z}(s,z)\, dz\right\vert^2\right].
\eeqn
By the Cauchy-Schwarz inequality, this is bounded above by
\beqn
 |t-s|\, E\left[\int_t^s \left\vert\frac{\partial v_D}{\partial r}(r,x)\right\vert^2 \, dr \right]
+|x-y|\, E\left[\int_x^y \left\vert\frac{\partial v_D}{\partial z}(s,z)\right\vert^2 dz\right].
\eeqn
Using again Theorem \ref{ch1'-tDL}, we conclude that
\begin{align}
\label{4-p1-1'-weak.1}
E\left[\left( v_D(t,x) - v_D(s,y)\right)^2\right]& \le C
\left(|t-s|^2 + |x-y|^2\right),
\end{align}
where the constant $C$ depends on $I$ and $J$ (see \eqref{ch1'-tDL.3}).

Since $u_D(t,x) =  v_D(t,x) + v(t,x)$, we conclude from \eqref{4-p1-1'-weak.1}  and \eqref{ch1'.v1} that
\beq
\label{4-p1-1'-weak.conclu}
E\left[\left( u_D(t,x) - u_D(s,y)\right)^2\right] \le C\left(|t-s|^{\frac{1}{4}} + |x-y|^{\frac{1}{2}}\right)^2,
\eeq
for any $(t,x), (s,y) \in\ I\times J$, where the constant $C$ depends on $I$ and $J$.

By  \eqref{4-p1-1'-weak.conclu}, the non-centred Gaussian process $u:= u_D$ satisfies the inequality \eqref{aap1-3.1} with ${\bf \Delta}(t,x;s,y) = |t-s|^{\frac{1}{4}} + |x-y|^{\frac{1}{2}}$. By Theorem \ref{app1-3-t1}, \eqref{ch1'-s7.18} holds, which is  \eqref{4-p1-1'-weak.0}, as well as the claim concerning joint local Hölder continuity.
 \end{proof}

\medskip

\noindent{\em Optimality of the H\"older exponents}
\smallskip

As we did in Proposition \ref{ch1'-p100} for the stochastic heat equation on the real line, we now analyse the optimality of the H\"older exponents in Proposition \ref{4-p1-1'-weak}. This question is addressed in the next statement, where $T>0$ is fixed.
\begin{prop}
\label{ch1'.1203}
Let $(u(t,x),\, (t,x)\in[0,T]\times [0,L])$ denote the random field solution to either \eqref{ch1'.HD} or \eqref{1'.14}.
\begin{enumerate}
\item Let $I\subset\ ]0,T]$ be a closed interval with positive length. Fix $x\in\ ]0,L[$. A.s., the sample paths of the process $(u(t,x),\,t\in I)$
 are not H\"older continuous with exponent $\alpha\in\, ]\frac{1}{4},1]$.
\item Let $J\subset\ ]0,L[$ be a closed interval with positive length. Fix $t\in\,]0,T]$. Almost surely, the sample paths of the process
$
(u(t,x),\, x\in J)
$
 are not H\"older continuous with exponent $\beta\in\, ]\frac{1}{2}, 1]$.
\end{enumerate}
\end{prop}
\begin{proof}
Both statements are a consequence of Proposition \ref{ch1'-p100} and Theorem \ref{ch1'-tDL}. Indeed, we only consider the case of Dirichlet boundary conditions, since with the Neumann conditions, the arguments are similar.

Consider the setting of Claim 1 and let $I=[t_0,t_1]\subset\ ]0,T]$ and $x\in\ ]0,L[$ be fixed. Using the notations of Proposition \ref{4-p1-1'-weak}, we have the decomposition
\beqn
u_D = v + v_D.
\eeqn
The function $I\ni t\mapsto v_D(t,x)$ is $\mathcal{C}^\infty$, by Theorem \ref{ch1'-tDL}. Since by Proposition \ref{ch1'-p100}, a.s., $I\ni t\mapsto v(t,x)$ is not H\"older continuous with exponent $\alpha\in\ \left]\frac{1}{4},1\right]$, the same property holds for $u_D$.

The proof of Claim 2. is similar.
\end{proof}

\begin{remark}
\label{ch3new-r1} Because of Theorem \ref{ch1'-tDL}, the conclusions of Theorem \ref{ch1'-s3-t1}, with $u$ replaced by $u_D$ or $u_N$, $\re_+$ replaced by $]0,\infty[$ and $x\in\IR$ replaced by $x\in\, ]0,L[$ apply to the SPDEs \eqref{ch1'.HD} and \eqref{1'.14}.
\end{remark}
\bigskip

\subsection{Global H\"older continuity of the sample paths}
\label{ch1'-ss3.4.3}

Our next objective is to find conditions that imply global H\"older continuity of the sample paths to the solutions to \eqref{ch1'.HD} and \eqref{1'.14}, including at time $t=0$ and at the boundary points $0$ and $L$.
We will address separately the two types of boundary conditions.

\begin{remark}
\label{continuity-Dbounded}
In view of \eqref{c} and \eqref{1'.N}, the (H\"older) continuity properties of $u$ are related to those of $I_0$ and $v$, where
\beq
\label{Io-v}
I_0(t,x) = \int^{L}_0 dy \, G_L(t;x,y) u_{0}(y),\ v(t,x) = \int_0^t \int_0^L G_L(t-s; x,y)\, W(ds,dy).
\eeq
These can be studied separately.
The same comments as in Remark \ref{continuity-Dunbounded} concerning (Hölder) continuity are valid here. In particular, $I_0$ is $C^\infty$ on $]0, \infty[ \times [0,L]$ (\cite[Section 2.3.3]{evans}).
\end{remark}
\medskip

\noindent{\em Dirichlet boundary conditions}
\medskip

For any $\eta\in\, ]0,1]$, let  $\mathcal{C}_0^\eta([0,L])$\label{rd06_14l2} denote the set of functions $f\in\mathcal{C}^\eta([0,L])$ such that $f(0)=f(L)=0$, and define
\beqn
\Vert f\Vert_{\mathcal{C}_0^\eta([0,L])}:= \sup_{0\le x<y\le L}\frac{|f(x)-f(y)|}{|x-y|^\eta} < \infty.
\eeqn
Notice that $\Vert f\Vert_{\mathcal{C}_0^\eta([0,L])} = \Vert f\Vert_{\mathcal{C}^\eta([0,L])}$.

\begin{lemma}
\label{cor1}
We assume that for some $\eta\in\ ]0,1]$, the initial condition $u_0$ of equation \eqref{ch1'.HD} belongs to $\mathcal{C}_0^\eta([0,L])$. Consider the function $(t,x)\mapsto I_0(t,x)$ on $\re_+\times [0,L]$ defined in \eqref{cor1.0}.

Then $\sup_{(t,x) \in \R_+ \times [0,L]} \vert I_0(t,x) \vert \leq \sup_{x \in [0,L]} \vert u_0(t,x) \vert$ and there exists a constant $C>0$ such that, for any $s,t\in\re_+$ and every $x, y \in[0,L]$,
\beq
\label{cor1.1}
\left\vert I_0(t,x) - I_0(s,y)\right\vert\le C \Vert u_0\Vert_{\mathcal{C}_0^\eta([0,L])}\left(|t-s|^{\frac{\eta}{2}} + |x - y|^\eta  \right).
\eeq
Consequently, the mapping $(t,x)\mapsto I_0(t,x)$ is H\"older continuous, jointly in $(t,x)\in\re_+\times[0,L]$ with exponents $(\frac{\eta}{2},\eta)$.
\end{lemma}

\begin{proof}
The uniform bound on $I_0(t,x)$ follows from the formula in \eqref{Io-v} and Proposition \ref{ch1'-pPD} (iii). Let $u_0^{o,p}$ be the odd and $2L$-periodic extension of $u_0$ as defined in \eqref{vop}. According to Proposition \ref{ch1-equivGD-bis},
\beqn
I_0(t,x) = \int_{-\infty}^{+\infty} \Gamma(t,x-z) u_0^{o,p} (z)\, dz.
\eeqn
Since $u_0(0)= u_0(L) = 0$, $u_0^{o,p}$ is continuous on $\IR$, and since $u_0\in \mathcal{C}_0^\eta([0,L])$, Lemma \ref{A3-l2}
implies that $u_0^{o,p} \in \cC^\eta(\R)$ and
$\Vert u_0^{o,p}\Vert_{\mathcal{C}^\eta(\IR)} \le 2 \Vert u_0\Vert_{ \mathcal{C}_0^\eta([0,L])}$. Therefore, the conclusions follows directly from Lemma \ref{ch1'-prep-lemma3.2.2} and \eqref{ch1'-initial 6}.
\end{proof}

\begin{prop}
\label{pre-4-p1-1'}
Let $(v(t,x),\, (t,x)\in\re_+\times [0,L])$ be as in \eqref{Io-v} with $G_L$ given in \eqref{ch1'.600}.
Then $\sup_{(t,x) \in \R_+ \times [0, L]} E(v^2(t,x)) < \infty$ and there exists a constant $C < \infty$ such that for all $(t,x)\in\re_+\times [0,L]$,
\beq
\label{4-p1-1'.1-p2}
E\left[\left(v(t,x) - v(s,y)\right)^2\right] \le C \left(|t-s|^{\frac{1}{4}} + |x-y|^{\half}\right)^2.
\eeq
Therefore, for any $\alpha \in\, ]0, \frac{1}{4}[$ and any $\beta \in\, ]0, \half[$, there exists a version of $v=(v(t,x),\,  (t,x)\in[0,\infty[\times \re)$ with locally jointly H\"older continuous sample paths with exponents $(\alpha, \beta)$.
\end{prop}
\begin{proof}
Let $0\le s<t<\infty$, $x, y\in [0,L]$. To simplify the presentation, we use the convention $G_L(s;x,z)=0$, for all  $s\le0$ and $x,z\in[0,L]$.

The uniform $L^2(\Omega)$-bound on $v(t,x)$ follows from \eqref{monyo(*1)}. By  the It\^o isometry,
\begin{align*}
E\left(\left(v(t,x) - v(s,y)\right)^2\right) &= \int_0^\infty dr \int_0^L dz\  [G_L(t-r;x,z) - G_L(s-r;y,z)]^2 \\
&\le C \left(|t-s|^{\frac{1}{4}} + |x-y|^{\half}\right)^2,
\end{align*}
where the last inequality follows from
\eqref{1'1100} in Lemma \ref{ch1'-l2}, and $C$ does not depend on $s, t, x, y$ (or even $L$). This is \eqref{4-p1-1'.1-p2}.

The claim about H\"older continuity follows from Kolmogorov's continuity criterion Theorem \ref{app1-3-t1}.
\end{proof}


\begin{prop}
\label{4-p1-1'}
Assume that the initial condition $u_0$ in \eqref{ch1'.HD} belongs to $\mathcal{C}^{\eta}_0([0,L])$, $\eta\in\ ]0,1]$. Then the random field solution $u = (u(t,x), \, (t,x)\in\re_+\times [0,L])$ to \eqref{ch1'.HD} given by \eqref{c},  satisfies the following property:

For any $p\in[2,\infty[$, there exists a finite constant $C=C(p,u_0,L)$ such that for every $(t,x), (s,y)\in \re_+\times [0,L]$,
\beq
\label{4-p1-1'.0}
E\left[\left\vert u(t,x) - u(s, y)\right\vert^p\right] \le C \left(|t-s|^{\frac{1}{4}\wedge \frac{\eta}{2}} + |x-y|^{\frac{1}{2}\wedge \eta}\right)^p.
\eeq
As a consequence, there is a version of $(u(t,x), (t,x) \in \re_+\times [0,L])$ with locally jointly Hölder continuous sample paths with exponents $(\alpha, \beta)$; in the time variable $t$, the constraints on $\alpha$ are
\beqn
\alpha\in\left]0,\tfrac{1}{4}\right[,  {\text{ if}}\ \eta\ge \tfrac{1}{2},\qquad
 \alpha\in\left]0,\tfrac{\eta}{2}\right],  {\text{ if}}\  \eta < \tfrac{1}{2},
\eeqn
while in the space variable $x$, the constraints on $\beta$ are
\beqn
\beta\in
\left]0,\tfrac{1}{2}\right[,  {\text{ if}}\ \eta\ge \tfrac{1}{2},\qquad
 \beta\in \left]0,\eta\right],  {\text{ if}}\ \eta < \tfrac{1}{2}.
\eeqn
\end{prop}
\begin{proof}
The random field $(v(t,x),\, (t,x)\in\re_+\times [0,L])$ of Proposition \ref{pre-4-p1-1'} is centred and Gaussian. Thus, \eqref{4-p1-1'.1-p2}
yields
\begin{align}
\label{4-p1-1'.1}
E\left[\left|v(t,x) - v(s,y)\right|^p\right] &
\le c_p \left(E\left[\left|v(t,x) -v(s,y)\right|^2\right]\right)^{\frac{p}{2}}\notag\\
& \le c_p\left(|t-s|^{\frac{1}{4}} + |x-y|^{\half}\right)^p,
\end{align}
with $c_p= \left(2^p/\pi\right)^\half \Gamma_E\left((p+1)/2\right)$ (see \eqref{A3.4.0-bis}).

Notice that by Lemma \ref{cor1} and Proposition \ref{pre-4-p1-1'}, the left-hand side of \eqref{4-p1-1'.0} is uniformly bounded. Therefore, appealing to  \eqref{c} and using \eqref{cor1.1} and \eqref{4-p1-1'.1} we obtain \eqref{4-p1-1'.0}.

The claim about H\"older continuity follows from  Kolmogorov's continuity criterion Theorem \ref{ch1'-s7-t2}.
\end{proof}

The upper bound on increments of $v(t,x)$ in Proposition \ref{4-p1-1'} has a corresponding lower bound, similar to that given in Proposition \ref{ch1'-p3} for the stochastic heat equation on $\R$, which we now establish.

\begin{prop}
\label{eqiv-norms-D}
Let $(v(t,x), (t,x)\in \re_+\times [0,L])$ be the solution to \eqref{ch1'.HD} with $u_0\equiv 0$. Let $0<t_0\le T$, $\alpha\in\, ]0,L/2[$, and define  $I=[t_0,T]$, $J=[\alpha,L-\alpha]$. Then there exist constants $0<C_1\le C_0$ such that, for all $(t,x), (s,y)\in I\times J$,
\beq
\label{eqiv-norms-D.1}
C_1\left(|t-s|^\half + |x-y|\right) \le E\left[(v(t,x) - v(s,y ))^2\right] \le C_0\left(|t-s|^\half + |x-y|\right).
\eeq
The upper bound holds for any $(t,x), (s,y)\in \re_+\times [0,L]$. The constant $C_0$ is universal  and the constant $C_1$ depends on $\alpha$.
\end{prop}
\begin{proof}
The upper bound follows from \eqref{4-p1-1'.1-p2}.
For the lower bound, let
$$
   C_2 = c_0.\qquad C_3 = (1-e^{-2\pi^2 t_0}) \frac{\alpha}{L}, \qquad C_4 = C,
$$
where $c_0$ is defined in Lemma \ref{ch1'-l200}  and $C$ is the constant in \eqref{1'1100}. We observe that we have the three inequalities
\beq
\label{rde3.3.30}
   E\left[( v(t,x) - v(s,y) )^2\right] \geq C_2 \vert t-s\vert^{\half},
   \eeq
 $x,y \in J$, $s,t \in \IR_+$ with $|t-s| < \tfrac{c_0}{c(\alpha)}$, where $c(\alpha)$ appears in Lemma \ref{ch1'-l200},
 \begin{align}
   E\left[( v(t,x) - v(t,y) )^2\right] &\geq C_3 \vert x-y\vert,\quad \ x,y \in J,\ t \geq t_0, \label{rde3.3.31}\\
   E\left[( v(t,y) - v(s,y) )^2\right] &\leq C_4 \vert t-s\vert^{\half}, \quad y \in J,\ s,t \in \IR_+, \label{rde3.3.32}
\end{align}
Indeed, \eqref{rde3.3.30} follows from \eqref{1'.1111} and \eqref{rde3.3.32} from \eqref{1'1100}.
For the proof of \eqref{rde3.3.31}, we fix $t\ge t_0$ and apply
\eqref{1'.dlb0} to deduce that
\begin{align*}
E\left[( v(t,x) - v(t,y) )^2\right] & \ge \left(\frac{1-e^{-2\tfrac{\pi^2}{L} t_0}}{2}\right)
(L - \vert x - y \vert) \frac{1}{L} \vert x - y \vert\\
&\ge\left (1-e^{-2\tfrac{\pi^2}{L} t_0}\right) \frac{\alpha}{L}|x-y|,
\end{align*}
where, in the last inequality, we have used that $\vert x - y \vert \leq L - 2 \alpha$ since $x, y \in J$.

Proceeding exactly as in the proof of the lower bound in Proposition \ref{ch1'-p3}, we see that there is a constant $C_1 >0$ such that for all $s,t \in [t_0,T]$ with $\vert t-s\vert < c_0/c(\alpha)$, for all $x,y \in J$,
\beq\label{rde3.3.33}
C_1\left(|t-s|^\half + |x-y|\right) \le E\left[(v(t,x) - v(s,y ))^2\right].
\eeq

In order to obtain this inequality for all $s,t \in [t_0,T]$, observe that the function $(t,x; s,y) \mapsto E[|v(t,x) - v(s,y )|^2]$ is continuous  (by the upper bound in \eqref{eqiv-norms-D.1}) and positive.  Indeed, by the isometry property, for $s \leq t$,
\begin{align*}
  E[(v(t,x) - v(s,y ))^2] &= \int_0^s dr \int_0^L dz\, (G_L(t-r;x,z) - G_L(s-r;y,z))^2 \\
  &\qquad + \int_s^t dr \int_0^L dz\, G_L^2(t-r;x,z).
\end{align*}
When $s< t$, the second term is positive, and when $s=t$, the first term is positive since the integrand it not identically $0$. Therefore the minimum $m$ of this function on the compact set
$$
 A = \{(t,x;s,y) \in ([t_0,T]\times J)^2: \  \vert t-s\vert \geq c_0/c(\alpha) \}
$$
is positive. Let $M$ be the maximum of $|t-s|^\half + |x-y|$ on $([t_0,T]\times J)^2$. Then for all $(t,x;s,y) \in A$,
$$
   \frac{m}{M} \left(|t-s|^\half + |x-y|\right) \leq E\left[(v(t,x) - v(s,y ))^2\right].
$$
Together with \eqref{rde3.3.33}, this proves \eqref{eqiv-norms-D.1}.
\end{proof}
\medskip


\noindent{\em Neumann boundary conditions}
\smallskip

Recall that, for any $\eta\in\, ]0,1]$, $\mathcal{C}^\eta([0,L])$ denotes the set of functions $f:[0,L]\rightarrow \re$ satisfying
\beqn
\Vert f\Vert_{\mathcal{C}^\eta([0,L])}:= \sup_{0\le x<y\le L}\frac{|f(x)-f(y)\vert}{|x-y|^\eta} < \infty.
\eeqn

\begin{lemma}
\label{cor1-N}
We assume that for some $\eta \in\, ]0, 1]$, the initial condition  $u_0$ of equation \eqref{1'.14} belongs to $\mathcal{C}^\eta([0,L])$. Let
$(t,x)\mapsto I_0(t,x)$ be the function on $\re_+\times [0,L]$ given in \eqref{cor1.0-N}.
There exists a constant $C>0$ such that, for any $s,t\in\re_+$ and every $x, y \in[0,L]$,
\beq
\label{cor1.1-N}
\left\vert I_0(t,x) - I_0(s,y)\right\vert\le C\Vert u_0\Vert_{\mathcal{C}^\eta([0,L])}\left(|t-s|^{\frac{\eta}{2}} + |x - y|^\eta  \right).
\eeq
Consequently, the mapping $(t,x)\mapsto I_0(t,x)$ is H\"older continuous, jointly in $(t,x)\in\re_+\times[0,L]$ with exponents $(\frac{\eta}{2}, \eta)$.
\end{lemma}

\begin{proof}
The proof is similar to that of Lemma \ref{cor1}. We use the even and $2L$-periodic extension $u_0^{e,p}$ of $u_0$ (see \eqref{vep}), instead of $u_0^{o,p}$ there, and Proposition \ref{ch1-lequivN-bis}  instead of Proposition \ref{ch1-equivGD-bis}.
By Lemma \ref{A3-l2}, $u_0^{e, p} \in \mathcal{C}^\eta(\R)$ and in fact, $\Vert u_0^{e,p}\Vert_{\mathcal{C}^\eta(\IR)} = \Vert u_0\Vert_{\mathcal{C}^\eta([0,L])}$.
Therefore, the conclusions follow directly from Lemma \ref{ch1'-prep-lemma3.2.2} and \eqref{ch1'-initial 6}.
\end{proof}
\begin{prop}
\label{pre-4-p1-1'-N}
Let $(v(t,x),\,(t,x)\in\re_+\times [0,L])$ be as in \eqref{Io-v} with $G_L$ given in \eqref{1'.400}.
Fix $T > 0$. There exists a finite constant C = C(T, L) such that  for all $(t,x)\in [0,T]\times [0,L]$,
\beq
\label{dia 5}
E\left[\left(v(t,x) - v(s,y)\right)^2\right] \le C \left(|t-s|^{\frac{1}{4}} + |x-y|^{\half}\right)^2.
\eeq
Therefore, for any $\alpha \in\, ]0, \frac{1}{4}[$ and any $\beta \in\, ]0, \half[$, there exists a version of $v=(v(t,x),\,  (t,x)\in[0,\infty[\times \re)$ with locally jointly H\"older continuous sample paths with exponents $(\alpha, \beta)$.
\end{prop}
\begin{proof}
The proof is similar to that of Proposition \ref{pre-4-p1-1'}. After using the It\^o isometry, and applying \eqref{1'11000} in Lemma \ref{ch1'-l3}, we obtain \eqref{dia 5} (with a constant $C$ that depends on $T$ and $L$).
Finally, the claim about H\"older continuity follows from Kolmogorov's continuity criterion Theorem \ref{app1-3-t1}.
\end{proof}

\begin{prop}
\label{4-p1-1'-N}
Assume that  for some $\eta\in\,]0,1]$, the initial condition $u_0$ in \eqref{1'.14} belongs to $\mathcal{C}^{\eta}([0,L])$. Then the random field solution $u$ to the SPDE \eqref{1'.14} given by \eqref{1'.N} satisfies the following property:

Fix $T>0$. For any $p\in[2,\infty[$, there exists a finite constant $C=C(p,u_0,T,L)$  such that, for any $t,s\in [0,T]$ and every $x,y\in[0,L]$,
\beq
\label{dia5-bis}
E\left[\left\vert u(t,x) - u(s,y)\right\vert^p\right] \le C \left(|t-s|^{\frac{1}{4}\wedge \frac{\eta}{2}} + |x- y|^{\frac{1}{2}\wedge \eta}\right)^p.
\eeq
As a consequence, there is a version of $(u(t,x), (t,x) \in [0, T] \times [0, L])$ with locally jointly Hölder continuous sample paths with exponents $(\alpha, \beta)$. In the time variable $t$, the constraints on $\alpha$ are
\beqn
\alpha\in \left]0,\tfrac{1}{4}\right[,\ {\text{ if}}\ \eta\ge \tfrac{1}{2},\quad
\alpha\in \left]0,\tfrac{\eta}{2}\right],\ {\text{ if}}\ \eta < \tfrac{1}{2},
\eeqn
while in the space variable $x$, the constraints on $\beta$ are
\beqn
\beta\in
\left]0,\tfrac{1}{2}\right[,\  {\text{ if}}\ \eta\ge \tfrac{1}{2},\quad
\beta\in \left]0,\eta\right],\  {\text{ if}}\ \eta < \tfrac{1}{2}.
\eeqn
\end{prop}
\begin{proof}
Let $(v(t,x),\, (t,x)\in\re_+\times [0,L])$ be as in Proposition \ref{pre-4-p1-1'-N}. Since this is a centred Gaussian process, \eqref{dia 5} implies
\beq
\label{4-p1-1'.1-N}
E\left[\left|v(t,x) - v(s,y)\right|^p\right] \le c_p\left(|t-s|^{\frac{1}{4}} + |x-y|^{\half}\right)^p,
\eeq
with $c_p = \left(2^p/\sqrt \pi\right)^\half \Gamma_E \left((p+1)/2\right)$ (see \eqref{A3.4.0-bis}).
From \eqref{1'.N}, consider the decomposition
$u(t,x) = I_0(t,x) + v(t,x)$,
with $I_0(t,x)$ given in \eqref{cor1.0-N}. Then \eqref{dia5-bis} is a consequence of \eqref{cor1.1-N} and \eqref{4-p1-1'.1-N}.

Notice that even though \eqref{cor1.1-N} and \eqref{dia 5} are valid for all $s, t \in \R_+$, \eqref{dia5-bis} is only valid for $s, t \in [0, T]$,
since the left-hand side of \eqref{4-p1-1'.1-N} is unbounded (see \eqref{cansada}) and  it is only on compact sets that increments with the smaller exponents dominate.

The claim about H\"older continuity follows from  Kolmogorov's continuity criterion Theorem \ref{ch1'-s7-t2}.
\end{proof}

As for the heat equation with Dirichlet boundary conditions,
we can obtain upper and lower  bounds for the canonical pseudo-metric associated to the process in \eqref{1'.N} when $u_0\equiv 0$. Here, however, the lower bounds are
valid up to and including the boundary points $0$ and $L$. This will imply that the intervals for possible Hölder exponents obtained in Proposition \ref{4-p1-1'-N} are sharp. The last part of this section addresses these questions.
\begin{prop}
\label{ch1'.pN1}
Fix $0<t_0\le T$ and let $(v(t,x),\, (t,x)\in[0,T]\times [0,L])$ be the random field solution to \eqref{1'.14}, with $u_0\equiv 0$, given by \eqref{1'.N}. There exist constants $c_1>0$ (depending on $t_0$) and a constant $c_2(T,L)$ such that,
for all $(t,x), (s,y)\in [t_0,T]\times [0,L]$,
\beq
\label{1'.N100}
c_1\left(|t-s|^{\half} + |x-y|\right) \le E\left[(v(t,x)-v(s,y))^2\right] \le c_2\left(|t-s|^{\half} + |x-y|\right).
\eeq
The upper bound holds for any $(t,x), (s,y)\in [0,T]\times [0,L]$.
\end{prop}
\begin{proof}
The upper bound follows from \eqref{dia 5}.

For the proof of the lower bound in \eqref{1'.N100}, let
$$
   C_2 = \frac{1}{\sqrt{2\pi}},\qquad C_3 = \frac{1-e^{-2\pi^2 t_0}}{2}, \qquad C_4 = C,
$$
where the constant $C$ is the same as that on the right-hand side of  \eqref{1'11000}. Then we have the following inequalities:
\begin{align}
   E\left[( v(t,x) - v(s,y) )^2\right] &\geq C_2 \vert t-s\vert^{\half}, \quad x,y \in [0,L],\ s,t \in \IR_+, \label{rde3.3.40} \\
   E\left[( v(t,x) - v(t,y) )^2\right] &\geq C_3 \vert x-y\vert,\quad x,y \in [0,L],\ t \geq t_0, \label{rde3.3.41}\\
   E\left[( v(t,y) - v(s,y) )^2\right] &\leq C_4 \vert t-s\vert^{\half}, \quad y \in [0,L],\ s,t \in[0,T]. \label{rde3.3.42}
\end{align}
Indeed, \eqref{rde3.3.40} obviously holds if $s=t$. If $0\le s<t$, then applying  the isometry property (and the convention $G_L(s;y,z)=0$, for any $s<0$, $y,z\in[0,L]$), we have
\begin{align*}
  E\left[(v(t,x) - v(s,y ))^2\right] &= \int_0^t dr \int_0^L dz\, (G_L(t-r;x,z) - G_L(s-r;y,z))^2 \\
  &\geq \int_s^t dr \int_0^L dz\, G_L^2(t-r;x,z)\\
  &\geq C_2 \vert t-s\vert^{\half}
\end{align*}
by \eqref{1'.20}. Hence \eqref{rde3.3.40} holds.
The lower bound \eqref{rde3.3.41} follows from \eqref{1'.2000}. Finally, \eqref{rde3.3.42} follows from \eqref{1'11000}. Proceeding exactly as in the proof of the lower bound in Proposition \ref{eqiv-norms-D}, we see that the lower bound in \eqref{1'.N100} holds (observe that the situation is simpler than in Proposition \ref{eqiv-norms-D}, since there is no constraint on $\vert t-s\vert$ in \eqref{rde3.3.40}, as there is in \eqref{rde3.3.30}).
This ends the proof.
\end{proof}
\medskip

Applying Proposition \ref{ch1'.pN1}, we obtain sharpness of the H\"older exponents of the sample paths of equation \eqref{1'.14}
including  at $x=0$ and $x=L$, whereas, in Proposition \ref{ch1'.1203}, we only obtained this for $x \in ]0, L[$.

\begin{prop}
\label{sharpheatneumann}
Fix $T>0$, $t_0\in\ ]0,T]$ and let
$(u(t,x),\,  (t,x)\in[0,T]\times [0,L])$ be the random field solution to \eqref{1'.14} with vanishing initial conditions.
\begin{enumerate}
\item Fix $x\in[0,L]$ and $\alpha\in\, \left]\frac{1}{4}, 1\right]$. Then a.s., the sample paths of the stochastic process $(u(t,x),\, t\in [t_0,T])$ are not H\"older continuous with exponent $\alpha$.
\item Fix $t\in  [t_0,T]$ and $\beta \in\, \left]\frac{1}{2},1\right]$. Then a.s., the sample paths of the stochastic process $(u(t,x),\, x\in[0,L])$ are not H\"older continuous with exponent $\beta$.
\end{enumerate}
\end{prop}
\begin{proof}
1. Set $I=[t_0,t_1]\subset\ ]0,T]$ and observe (from \eqref{1'.N100}) that if $t_0<s<t<T$, then for any $x\in[0,L]$, we have
\beqn
E\left[(u(t,x)-u(s,x))^2\right] \ge c_1|t-s|^\half.
\eeqn
Hence, 1. follows from Theorem \ref{app1-3-t2} with $\alpha=\tfrac{1}{4}$.

2. Letting $x,y\in[0,L]$ and applying \eqref{1'.N100}, for any $t\in[t_0,T]$, we obtain
\begin{align*}
E\left[(u(t,x)-u(t,y))^2\right] \ge c_1|x-y|.
\end{align*}
As above, we deduce 2. by applying Theorem \ref{app1-3-t2} with $\alpha=\tfrac{1}{2}$.
\end{proof}


\section{The stochastic wave equation}
\label{ch4-ss2.3-1'}

Let $\mathcal{L}= \frac{\partial^2}{\partial t^2}- \frac{\partial^{2}}{\partial x^{2}}$ be the wave operator.\index{operator!wave}\index{wave!operator}
In this section, we will consider three cases of spatial domains $D\subset \re$ namely, $D=\re$, $D=\ ]0,\infty[$ and $D=\ ]0,L[$. As in the case of the stochastic heat equation, we consider the SPDE $\mathcal{L}u = \dot W$ in $\IR_+ \times D$,
where $\dot W$ is a space-time white noise, and we establish existence of random field solutions and determine the regularity of their sample paths.

\subsection{Existence of random field solutions}

The stochastic wave equation on $D$ driven by space-time white noise $\dot W$ is \beq
\label{wave-1'}
  \frac{\partial^2 u}{\partial t^2}(t,x)- \frac{\partial^{2} u}{\partial x^{2}}(t,x) = \dot W(t,x),\qquad  t>0,\ x \in D,
  \eeq
  with initial conditions
  \beq
  \label{3.4(*2)}
  u(0,x)= f(x),\qquad
  \frac{\partial}{\partial t}u(0,x) = g(x),\qquad x\in D,
  \eeq
and, when $\partial D \neq \emptyset$, the vanishing Dirichlet boundary conditions
\beq
\label{3.4(*3)}
     u(t,x) = 0, \quad    t > 0,\  x \in \partial D.
\eeq

In the cases that we are considering, there is no boundary condition when $D = \re$, and $\partial D = \{0\}$ (respectively $\partial D = \{0, L\}$) when $D =\, ]0,\infty[ $ (respectively $D =\, ]0,L[$).

 Let $\Gamma(t,x; s,y)$ be the fundamental solution or the Green's function associated to $\mathcal{L} = \frac{\partial^2}{\partial t^2}- \frac{\partial^2}{\partial x^2}$ in $D$ (and the boundary conditions if $\partial D \neq \emptyset$). We are going to check below that Assumption \ref{ch4-a1-1'} is satisfied in each of the three cases under consideration.

 Let $I_0(t,x)$ be the solution to the homogeneous wave equation
 \beqn
  \frac{\partial^2 u}{\partial t^2}(t,x)- \frac{\partial^{2} u}{\partial x^{2}}(t,x) = 0,\qquad  t>0,\ x \in D,
  \eeqn
  with the same initial conditions \eqref{3.4(*2)}, and, when $\partial D \neq \emptyset$, the same Dirichlet boundary conditions\eqref{3.4(*3)}. The random field solution to \eqref{wave-1'}--\eqref{3.4(*3)} is, according to Definition \ref{ch4-d1-1'},
  \beq
  \label{anotherone}
     u(t,x) = I_0(t,x) + \int_0^t \int_D  \Gamma(t,x; s,y) \, W(ds, dy),
     \eeq
$(t,x)\in \re_+\times D$.

We now present the explicit formulas for $\Gamma$ and $I_0$.
\medskip

\noindent{\em Stochastic wave equation on $\re$}
\medskip

\begin{figure}[h]
		\centering
\begin{tikzpicture}[scale=1.5]
\draw [<->] (3,-1.5) -- (0,-1.5) -- (0,2.5);
\draw (0,0) --(0,-1.5);
\node[below] at (3,-1.5) {$s$};
\node[above] at (0,2.5) {$y$};
\draw [fill=orange](0,-1)--(0,2)--(1.5,0.5)--(0,-1);
\node[right] at (1.5,0.5) {$(t,x)$};
\node[left] at (0,2) {$x+t$};
\node[left] at (0,-1) {$x-t$};
\end{tikzpicture}
        \caption{The region $D(t,x)$ \label{fig1-chapter3}}
    \end{figure}
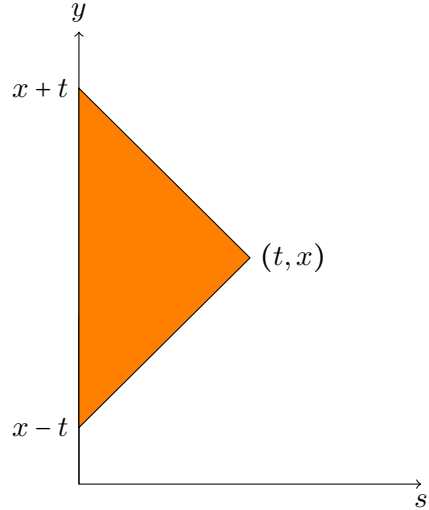

The fundamental solution associated to the wave operator $\mathcal{L}$ on $\re$ is $\Gamma(t,x; s,y): = \Gamma(t-s,x-y)$, where
\beq
\label{wfs}
    \Gamma(t,x) = \half 1_{\IR_+}(t)\, 1_{[-t,t ]}(x), \qquad (t,x)\in \IR^2,
\eeq
(see \cite[Chapter 1, \S 7]{treves}). Defining
\beq
\label{cone}
D(t,x)= \{(s,y)\in[0,t]\times \re: |x-y|\le t-s\}, \qquad (t,x)\in \re_+\times \re
\eeq
(see Figure \ref{fig1-chapter3}), we have
\beq
\label{gamma-cone}
\Gamma\left(t-s,x-y\right) = \half 1_{D(t,x)}(s,y).
\eeq

Clearly, the mapping $(s,y) \mapsto \Gamma(t-s,x-y)$ belongs to $L^2(\IR_+\times \IR)$ and
\begin{align}
\label{gamma-cone-norm}
\Vert \Gamma(t-\cdot,x-\ast)\Vert_{L^2(\IR_+\times \IR)}^2 &= \frac{1}{4} \int_{\re_+} ds
\int_{\re} dy\ 1_{D(t,x)}(s,y)\notag\\
&=\frac{1}{4} \int_0^t ds \int_{x-(t-s)}^{x+(t-s)} dy
= \frac{t^2}{4}.
\end{align}

The function $I_0(t,x)$ is given by d'Alembert's formula\index{d'Alembert's formula}\index{formula!d'Alembert's}
\begin{align}\nonumber
I_0(t,x)& = \left[\frac{d}{dt}\Gamma(t)\ast f + \Gamma(t) \ast g\right](x)\\
& = \frac{1}{2}[f(x+t) + f(x-t)] + \frac{1}{2}\int_{x-t}^{x+t} g(y)\, dy
\label{p93.1}
\end{align}
(see  \cite[Chapter 2, p. 67]{evans}, \cite[Chapter 2, p. 36]{strauss}).
If, for instance,  $f$ is a continuous function and $g\in L^1_{\text{loc}}(\re)$, then the function $(t,x) \mapsto I_0(t,x)$ from $\re_+\times \re$ into $\re$ is well-defined and continuous.
\medskip


\noindent{\em Stochastic wave equation on $\re_+$}
\medskip

\begin{figure}[h]
		\centering
\begin{tikzpicture}[scale=2]
\draw [<->] (2.5,0) -- (0,0) -- (0,2.5);
\draw (0,0) --(0,-1.5);
\node[below] at (2.5,0) {$s$};
\node[above] at (0,2.5) {$y$};
\draw [fill=orange](0,1)--(1,0) -- (1.5,0.5)--(0,2)--(0,1);
\draw [dashed] (0,-1)--(1,0);
\node[right] at (1.5,0.5) {$(t,x)$};
\node[left] at (0,2) {$x+t$};
\node[left] at (0,1) {$t-x$};
\node[left] at (0,-1) {$x-t$};
\end{tikzpicture}
\caption{The region $E(t,x)$ when $t>x$ \label{fig3.2}}
    \end{figure}
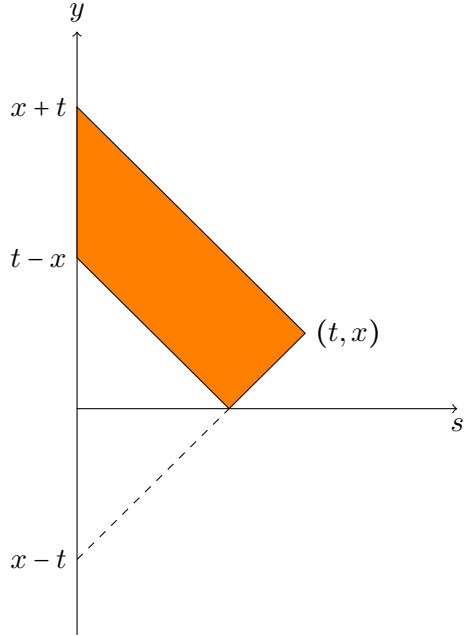

The Green's function associated to the wave operator $\mathcal{L}$ on $\re_+$ with Dirichlet boundary conditions is
\beq
\label{wavehalfline}
\Gamma(t,x;s,y) : = G(t-s;x,y) = \frac{1}{2} 1_{\{|x-(t-s)|\le y\le x+t-s\}} = \frac{1}{2} 1_{E(t,x)}(s,y),
\eeq
which is well-defined for all $0\le s\le t$ and $x,y\ge 0$, where
\beq
\label{Ewave}
E(t,x) = \{(s,y)\in[0,t]\times \re_+:\,  |x-t+s|\le y\le x+t-s\}.
\eeq
Notice that if $t\le x$, then $E(t,x)=D(t,x)$, where $D(t,x)$ has been defined in \eqref{cone}. If $t>x$, then $E(t,x)$ is the shadowed region in Figure \ref{fig3.2}.

The expression of $\Gamma$ can be found using the {\em reflection method} (also called {\em method of images}). Clearly, $G(t-\cdot,x,\ast)\in L^2([0,t]\times\re)$.

According to \cite[Chapter 3, p. 62]{strauss} (see also \cite[pg. 69]{evans}), the function $I_0(t,x)$ is given by
\beq
\label{p98.2}
I_0(t,x) =
\begin{cases}
\half [f(x+t) + f(|x-t|)] + \half \int_{|x-t|}^{x+t} g(z)\ dz, & x\ge t\ge 0,\\
\half [f(x+t) - f(|x-t|)] + \half \int_{|x-t|}^{x+t} g(z)\ dz, & 0\le x < t.
\end{cases}
\eeq
We will assume that $f$ is a continuous function, $f(0)=0$ and $g\in L^1_{\rm{loc}}(\IR_+)$, in which case $I_0(t,x)$ is well-defined and continuous. Notice that these formulas are compatible with the boundary condition $I_0(t, 0) = 0$.
\medskip

\noindent{\em Stochastic wave equation on a finite interval}
\smallskip

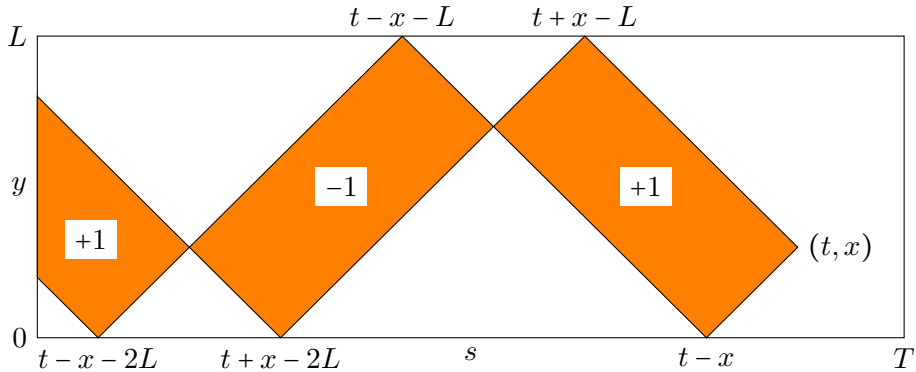
\begin{figure}[h]
\centering
\begin{tikzpicture}[scale=2]
\draw (0,0) -- (5.7,0) -- (5.7,2)--(0,2)--(0,0);
\node[left] at (0,0) {$0$};
\node[below] at (5.7,0) {$T$};
\node[left] at (0,2) {$L$};
\draw [fill=orange](3,1.4)--(3.6,2)--(5.,0.6)--(4.4,0)--(3,1.4);
\draw [fill=orange](3,1.4)--(1.6,0)--(1,0.6)--(2.4,2)--(3,1.4);
\draw [fill=orange](1,0.6)--(0.4,0)--(0,.4)--(0,1.6)--(1,0.6);
\node[right] at (5.,0.6) {$(t,x)$};
\node[above] at (3.6,2) {$t+x-L$};
\node[above] at (2.4,2) {$t-x-L$};
\node[below] at (4.4,0) {$t-x$};
\node[below] at (1.6,0) {$t+x-2L$};
\node[below] at (0.4,0) {$t-x-2L$};
\node [fill=white]at(4,1) {$+1$};
\node [fill=white]at(2,1) {$-1$};
\node [fill=white]at(0.35,0.65) {$+1$};
\node[below] at (2.85,0) {$s$};
\node[left] at (0,1) {$y$};
\end{tikzpicture}
\caption{The values of $2G_L(t-s;x,y)$ according to \eqref{wave-bc1-100'}}\label{fig3.3}
\end{figure}

Fix a bounded interval $[0,L]$. The Green's function for the operator  $\mathcal{L}$ on $]0,L[$ with vanishing Dirichlet boundary conditions is $\Gamma(t,x;s,y):=G_L(t-s; x,y)$, where
\beq
\label{wave-bc1-1'}
G_L(t-s; x,y) = \frac{1}{2} \sum_{m=-\infty}^\infty
 \left[ 1_{\{|x-2mL-y|\le t-s\}} - 1_{\{|x-2mL+y|\le t-s\}} \right],
 \eeq
 which is well-defined for $t\ge s\ge 0$, $x,y\in [0,L].$
This can be found
using the reflection method (see e.g. \cite[p. 235]{taylor1}).
In particular,
 \beq
 \label{wave-bc1-100'}
  G_L(t-s; x, y) = \half \left[1_{F_1(t,x)}(s,y) - 1_{F_2(t,x)}(s,y)\right],
 \eeq
 where, for $i=1,2$, $F_i(t,x)$ is a finite union of (possibly truncated) disjoint open rectangles and $F_1(t,x)\cap F_2(t,x)=\emptyset$: see Figure \ref{fig3.3}.
 \medskip

\medskip

 Alternatively, using the Fourier series expansion in terms of the CONS $e_{n,L}(x) = \sqrt{\frac{2}{L}}\sin(\frac{n\pi x}{L})$, $n\in\IN^*$, of $L^2([0,L])$, one has
 \beq
 \label{wave-bc2-1'}
G_L(t; x,y) = \sum_{m=1}^\infty
 \frac{2}{\pi m} \sin\left(\frac{m\pi x}{L}\right)\sin\left(\frac{m\pi y}{L}\right)\sin\left(\frac{m\pi t}{L}\right)
\eeq
(see e.g. \cite[p. 94, Equation (3.2.25)]{duffy}).

For any $(t,x)\in\ ]0,T[\times [0,L]$, we have $G_L(t-\cdot; x,\ast)\in L^2([0,t]\times [0,L])$.
Indeed, using the expression \eqref{wave-bc2-1'}, we see that
\beq
\label{ltwo}
\Vert G_L(t-s; x, \ast)\Vert_{L^2([0,L])}^2  = L \sum_{m=1}^\infty  \frac{2}{\pi^2 m^2}\sin^2\left(\frac{m\pi x}{L}\right) \sin^2\left(\frac{m\pi (t-s)}{L}\right).
\eeq
Bounding by $1$ the factors with $\sin^2$, we obtain
\beqn
\sup_{r\ge 0}\sup_{x\in[0,L]}\Vert G_L(r; x, \ast)\Vert_{L^2([0,L])}^2\le  \frac{2L}{\pi^2 } \sum_{m=1}^\infty \frac{1}{m^2}<\infty.
\eeqn
This implies $G_L(t-\cdot; x,\ast)\in L^2([0,t]\times [0,L])$.


Assume that the functions $f$ and $g$ belong to $L^2([0,L])$. Using \eqref{wave-bc2-1'}, one can check that the function $I_0(t,x)$ is given by
\beq
\label{i0-bis}
I_0(t,x) = \sum_{m=1}^\infty \left[A_m \cos\left(\frac{m\pi t}{L}\right) + B_m \sin \left(\frac{m\pi t}{L}\right)\right] \sin\left(\frac{m\pi x}{L}\right),
\eeq
where
\begin{align*}
A_m = \frac{2}{L}\int_0^L \sin\left(\frac{m\pi y}{L}\right)f(y)\, dy,\quad
B_m =  \frac{2}{m\pi}\int_0^L \sin\left(\frac{m\pi y}{L}\right)g(y)\, dy.
\end{align*}
Notice that these formulas are compatible with the boundary conditions $I_0(t, 0) = 0 = I_0(t, L)$.
 \medskip

An alternate expression for $I_0(t,x)$ is given in \cite[Chapter 3, p. 65, equation (5)]{strauss}. Take the odd periodic extensions
of the initial conditions $f$ and $g$ as in \eqref{vop}, that is,
\beqn
\phi^{o,p}=
\begin{cases}
\phi(x), & x\in [0,L],\\
-\phi(-x), &  x\in\, ]-L,0[,\\
\phi(x - 2kL), &  x \in\, ](2k-1) L, (2k+1)L], \ k \in \mathbb{Z},
\end{cases}
\eeqn
where $\phi$ stands for either function $f$ or $g$.
Then
\beq
\label{i0-tris}
I_0(t,x) = \half f^{o,p}(x+t) + \half f^{o,p}(x-t) + \half \int_{x-t}^{x+t} g^{o,p}(r)\, dr.
\eeq

 It follows, for example, that if $f$ is continuous on $[0,L]$, $f(0)=f(L)=0$, and $g\in L^1([0,L])$, then $I_0(t,x)$ is well-defined and is continuous on $\re_+\times [0,L]$.

\subsection{H\"older continuity properties of the sample paths}
We start by studying the regularity of the solution to the homogeneous wave equation. Notice that in order for $I_0(t,x)$ to be continuous on $[0,T] \times D$, the given data on the space-time boundary of $[0,T] \times D$ should be continuous, and this implies in particular the compatibility condition $f(y) = 0$ on $\partial D$.
\medskip

\noindent {\em Regularity of the function $(t,x)\mapsto I_0(t,x)$}
\smallskip

\begin{lemma}
\label{ch1'-ss2.3-hi}
Let $\gamma\in\, ]0,1]$. For the three forms of stochastic wave equations discussed above, we assume that
$f \in \cC^\gamma(D \cup \partial D)$ satisfies the compatibility condition $f(x)=0$ on $\partial D$ and $g$ is continuous on $D \cup \partial D$.

Fix $T>0$. Then the function
\beqn
(t,x) \longrightarrow I_0(t,x)
\eeqn
defined on $[0,T] \times(D \cup \partial D)$
is jointly locally  H\"older continuous in $(t,x)$ with exponent $\gamma$.
\end{lemma}
\begin{proof}
1. {\em Case $D=\re$}. Let $x,y \in \IR$, $t,s \in[0,T]$; since $f\in\mathcal{C}^{\gamma}(\IR)$,
\begin{align}
\label{p94.1}
&|f(x+t)+f(x-t)-f(y+s)-f(y-s)|\notag\\
&\quad\quad \quad  \le |f(x+t) - f(y+s)| + |f(x-t )- f(y-s)|\notag\\
&\quad \quad \quad \le \Vert f\Vert_{\mathcal{C}^{\gamma}(\re)} \left(|x-y+t-s|^\gamma + |x-y+s-t|^\gamma\right)\notag\\
&\quad\quad \quad \le 2\Vert f\Vert_{\mathcal{C}^{\gamma}(\re)}\left(|x-y|^{\gamma}+|t-s|^{\gamma} \right).
\end{align}
For any $b\ge 0$, set $F(b)=\int_0^b g(y)\ dy$. On any interval $]{-}M,M[$, the function
 $b\mapsto F(b)$ is differentiable with derivative bounded by
$C_M:=\sup_{x\in\ [-M,M]} |g(x)|$.
Fix $0\le t\le s\le T$ and $x, y\in \ [-M,M]$. Then
\beqn
\int_{y-s}^{y+s} g(z)\, dz - \int_{x-t}^{x+t} g(z)\, dz
 = F(y+s) - F(y-s)- F(x+t) + F(x-t).
  \eeqn
 Hence
\beq
\label{p94.2}
\left\vert\int_{y-s}^{y+s} g(y)\, dy - \int_{x-t}^{x+t} g(y)\, dy\right\vert \le 2C_{M+T} \left(|x-y| + |t-s|\right).
\eeq
Using the expression \eqref{p93.1} together with \eqref{p94.1} and \eqref{p94.2}, we obtain the assertion.
\smallskip

2. {\em Case $D=\, ]0,\infty[$.} The proof is nearly identical to that of the previous case, except that $x-t$ and $y-s$ there are replaced respectively by $|x-t|$, $|y-s|$. Using the reverse triangle inequality $||z|-|w||\le |z-w|$, the bounds \eqref{p94.1} and \eqref{p94.2} remain valid and therefore, also the conclusion.
\smallskip

3. {\em Case $D=\, ]0,L[$}. By \eqref{i0-tris}, we can proceed as in Case 1, with $f$ and $g$ replaced respectively by $f^{o,p}$ and $g^{o,p}$. By Part 1. of Lemma \ref{A3-l2},
\beqn
\Vert f^{o,p} \Vert_{\cC^\gamma(\R)} \leq 2 \Vert f \Vert_{\cC_0^\gamma([0,L])},
\eeqn
 and for $M > 2L$,
\beqn
\sup_{x\in[-M,M]} \vert g^{o,p}(x) \vert = \sup_{x\in[0,L]} \vert g(x) \vert.
\eeqn
 Therefore, we obtain the result as in Case 1.

The proof of the lemma is complete.
\end{proof}
\bigskip

The next theorem summarizes the results on regularity of the sample paths of the random field solutions to the wave equations studied in this section.  We use the notation
\beq
\Gamma(t,x;r,z)=
\begin{cases}
\label{defgammas}
\half\, 1_{D(t,x)}(r,z), & \text{if} \  D=\IR,\\
\half\, 1_{E(t,x)}(r,z), & \text{if} \  D = \, ]0,\infty[,\\
G_L(t-s;x,y), & \text{if} \  D = \,]0,L[,
\end{cases}
\eeq
where $G_L(t-s;x,y)$ is defined in \eqref{wave-bc1-1'}.

\begin{thm}
\label{ch4-pwave-1'}
Fix $T>0$. For the three forms of the stochastic wave equation considered in this section $(D=\IR$, $D = \, ]0,\infty[$, $D=\,]0,L[\,)$, set
\beq
\label{vw}
v(t,x) = \int_0^t \int_D \Gamma(t,x;r,z)\, W(dr,dz), \quad\quad (t,x)\in[0,T]\times D,
\eeq
with $\Gamma(t,x;r,z)$ given in \eqref{defgammas}.
\begin{enumerate}
\item There exists a constant $C$ such that for any $(t,x), (s,y)\in [0,T]\times (D\cup\partial D)$,
\beq
\label{mmm}
E\left[\left(v(t,x)-v(s,y)\right)^2\right] \le C \left(|t-s|^\half + |x-y|^\half\right)^2.
\eeq
If $D=\IR$ or $D = \, ]0,\infty[$, the constant $C$ depends only on $T$. If $D=\,]0,L[$, $C$ depends on $T$ and $L$.
\item For any $p>0$ and any $(t,x), (s,y)\in [0,T]\times (D\cup\partial D)$,
\beq
\label{mmm-p}
E\left[\left|v(t,x)-v(s,y)\right|^p\right] \le C_p \left(|t-s|^\half + |x-y|^\half\right)^p,
\eeq
with $C_p = C^{\frac{p}{2}}\left(\frac{2^ p}{\pi}\right)^\half\ \Gamma_E\left(\frac{p+1}{2}\right)$ and the constant $C$ is that in \eqref{mmm}.

Consequently, $(v(t,x), \, (t,x)\in \re_+\times D)$ has a version with locally  H\"older continuous sample paths, jointly in $(t,x)$, with exponent $\eta\in\, ]0,\half[$.
\item Let $\gamma\in\, ]0,1]$. Assume that the initial conditions $f$ and $g$ satisfy the hypotheses of Lemma \ref{ch1'-ss2.3-hi}. Fix a compact interval $D_0\subset (D\cup \partial D)$. Let $u$ be the random field solution to the stochastic wave equation on $D$ driven by space-time white noise. Then for any  $p\in[2,\infty[$ and any $(t,x), (s,y)\in[0,T]\times D_0$, there exists a constant $C_{p,T,D_0,\gamma}$ such that
\beq
\label{mmm-p-i}
E\left[\left |u(t,x)-u(s,y)\right|^p\right] \le C_{p,T,D_0,\gamma} \left(|t-s|^{\gamma\wedge \half} + |x-y|^{\gamma\wedge \half}\right)^p.
\eeq
Consequently,
$(u(t,x),\, (t,x)\in \re_+\times D)$ has a version with locally H\"older continuous sample paths, jointly in $(t,x)$. If $\gamma\in [\half,1]$ (respectively, $\gamma\in\, ]0,\half[$), then the common H\"older exponent is $\eta\in\, ]0,\half[$ (respectively, $\eta\in\, ]0,\gamma]$).
\end{enumerate}
\end{thm}
\begin{proof}
1.
By the It\^o isometry
\beq
\label{isokey}
E\left[\left(v(t,x)-v(s,y)\right)^2\right] =  \int_0^T dr \int_D dz\ \left(\Gamma(t,x;r,z) - \Gamma(s,y;r,z)\right)^2.
\eeq
If $D = \re$ (respectively $D =\, ]0, \infty[\,$), we use \eqref{1'.w1-bis} in Lemma \ref{P4:G-1'} (respectively \eqref{rdt3-bis}  in Lemma \ref{rdprop_wave_R_ub}) to see that the right hand-side of \eqref{isokey} is bounded above by $\tfrac{T}{2}(|t-s| + |x-y|)$. This yields \eqref{mmm} with $C := \frac{T}{2}$.

Let $D=\, ]0,L[$. From \eqref{rdeB.5.5a} in Lemma \ref{app2-5-l-w1}, we obtain that the right-hand side of \eqref{isokey} is bounded above (up to a multiplicative constant $C=C(T,L)$) by $(|t-s| + |x-y|)$, which implies \eqref{mmm} with $C:=C(T,L)$.
\smallskip

2. Since $(v(t,x))$ is a  centred Gaussian process, the $L^p$-estimate \eqref{mmm-p} follows from \eqref{A3.4.0-bis} in Lemma \ref{A3-l1} and \eqref{mmm}.
The claim about H\"older continuity is a consequence of Theorem \ref{app1-3-t1}.
\smallskip

3. The inequality \eqref{mmm-p-i} follows from Lemma \ref{ch1'-ss2.3-hi} and \eqref{mmm-p} above. As above, the claim about H\"older continuity is a consequence of Theorem \ref{app1-3-t1}.
\end{proof}
  As for the stochastic heat equation, we can obtain upper and lower bounds for the canonical pseudo-metric associated with the process in \eqref{vw}.
  \begin{prop}
  \label{Prop-Sec 3.4-(*1)}
  Let $(v(t,x),\, (t,x) \in [0, T] \times (D \cup \partial D))$ be the random field solution to \eqref{wave-1'} with vanishing initial and boundary conditions, given by \eqref{vw}. Fix $0 < t_0 < T$ and let $J \subset D$ be a compact interval. There are constants $c_1 > 0$ and $c_2 < \infty$ such that, for all $(t,x), (s,y) \in [t_0, T] \times J$,
  \beqn
     c_1 (\vert t - s \vert + \vert x - y \vert) \leq E\left[(v(t,x) - v(s,y))^2\right] \leq c_2 (\vert t - s \vert + \vert x - y \vert).
\eeqn
\end{prop}

\begin{proof}
The upper bound is \eqref{mmm}. Concerning the lower bound, for each of the three cases of wave equations considered in this section, $D = \R$, $D = ]0, \infty[$ and $D = ]0, L[$, we use \eqref{isokey} and Lemmas \ref{ch1'.pW1-lemma}, \ref{rdprop_wave_R+_lb}, and \ref{app2-5-l-w2}, respectively, to see that there is a constant $c_1$ such that for all $(t,x), (s,y) \in [t_0, T] \times J$,
\beqn
   E[(v(t,x) - v(s,y))^2] \geq c_1 (\vert t - s \vert + \vert x - y \vert),
   \eeqn
which is the desired lower bound.
   \end{proof}

   Finally, we check that the constraints on the Hölder exponents given in Theorem 3.4.2 are sharp.

\begin{prop}
\label{ch4-pwave-1'-opt}
Let $v(t,x)$ be defined in \eqref{vw}.
\begin{enumerate}
\item Fix $x\in D$,  $K\subset\, ]0,\infty[$ a closed interval of positive length, and $\eta\in\, ]\half,1]$. Then a.s., the sample paths of the process $(v(t,x),\, t\in K)$ are not H\"older continuous with exponent $\eta$.
\item Fix $t>0$, $J\subset D$ a closed interval with positive length, and $\eta\in\, ]\half,1]$ Then a.s., the sample paths of the stochastic process $(v(t,x),\, x\in J)$ are not H\"older continuous with exponent $\eta$.
\end{enumerate}
\end{prop}
\begin{proof}
The proofs of the two statements are similar. They rely on Theorem \ref{app1-3-t2} applied to the stochastic processes
$(v_1(t) = v(t,x),\, t\in K)$ and  $(v_2(x) = v(t,x),\, x\in J)$, respectively. For the three cases of wave equations considered in this section, $D=\re$, $D=\, ]0,\infty[$, and $D=\, ]0,L[$, we see that the assumptions of Theorem \ref{app1-3-t2}
 are satisfied with $\alpha = \half$, thanks to Proposition \ref{Prop-Sec 3.4-(*1)}.
\end{proof}

\section{Stochastic heat equation with a fractional Laplacian}
\label{ch3-sec3.5}

For $a > 0$, the fractional Laplacian $(-\Delta)^{a/2}$\index{fractional!Laplacian}\index{Laplacian!fractional}\index{operator!fractional Laplacian}\label{rdflaplacian} of an integrable function $f: \R^k \to \R$  is defined by means of its Fourier transform
\beq\label{rd02_08e1}
     \cF( (-\Delta)^{a/2} f)(\xi) = \vert \xi \vert^a \cF f(\xi),\quad    \xi \in \R^k,
 \eeq
with
   $\cF f (\xi) = \int_{\R^k} e^{-i \xi \cdot x} f(x)\, dx$,
and ``$\cdot$'' denotes the Euclidean inner product. This is a pseudo-differential operator, with Fourier multiplier $\vert \xi \vert^a$. For $a = 2$, this is simply the opposite of the ordinary Laplacian $\Delta f(x) = \sum_{i=1}^k \frac{\partial^2 f}{\partial x_i^2}(x)$.

    For $a \in\, ]0, 2[$, in the notation of \eqref{ch1'-ss7.3.1}, $(-\Delta)^{a/2} f = \null_xD_\delta^a f$ with $\delta = 0$. If $f\in \cC_0^\infty(\R^k)$, then $(-\Delta)^{a/2} f$ is a function, otherwise, it may only belong to $\cs^\prime(\R^k)$.
    For $a/2 = n$, $n \in \N^\ast$, $(-\Delta)^{a/2} f = (-1)^n \Delta^n f$ is obtained by iterating the Laplacian $n$ times.
   For $a/2 = n + s/2$, where $n \in \N^\ast$ and $s \in\, ]0, 2[$,
   \beqn
   (-\Delta)^{a/2} = (-\Delta)^{s/2}  \circ (-\Delta)^n =  (-\Delta)^n \circ  (-\Delta)^{s/2},
   \eeqn
    so one can give a formula for $(-\Delta)^{a/2} u$ (see Lemma \ref{ch1'-ss7.3-l0}). Further discussion of these fractional differential operators is deferred to Section \ref{ch1'-ss7.3}.

   In this section, we consider the SPDE 
 \beq
   \label{frac-linear}
      \begin{cases}
        \frac{\partial}{\partial t}u(t,x) + (- \Delta)^{a/2}\, u(t,x) = \dot W(t,x),&      (t,x) \in\, ]0, \infty[ \times \R^k,\\
        u(0, x) = u_0(x),&                 x \in \R^k,
       \end{cases}
\eeq
where $\dot W$ is space-time white noise on $\R_+ \times \R^k$ and $u_0 \in L^2(\R^k)$. This is an SPDE as in 
\eqref{ch1'-s5.0} with the partial differential operator $\cL = \frac{\partial}{\partial t} + (- \Delta)^{a/2}$, $\sigma\equiv 1$, $b\equiv 0$ and  $D = \R^k$. We assume that $a > k \geq 1$ because using a result of \cite[Theorem 11]{dalang}, it can be shown that for $a \in\, ]0, k]$, there is no random field solution to \eqref{ch1'-s5.0}.  

  The fundamental solution associated to $\cL$ on $\R^k$ is
   $\Gamma(t,x ; s,y) := G_a(t-s, x-y)$,
where
\beq\label{fs-frac}
    G_a(s,y) = \frac{1}{(2\pi)^k} \int_{\R^k} d\xi\, \exp(i \xi \cdot y - s \vert \xi \vert^a) 1_{]0, \infty[}(s),
\eeq
as can be found using the method outlined at the end of Section \ref{ch1'-ss7.3}. For each $s > 0$, this defines a bounded continuous real-valued function of $x$ since $\exp(- s \vert \xi \vert^a)$ is integrable over $\R^k$ and is an even function of each coordinate of $\xi$.

   Notice that $(s,y) \mapsto G_a(t-s, x-y)$ belongs to $L^2(\R_+ \times \R^k)$. Indeed, using Plancherel's theorem, we see that
   \begin{align}
   \label{ch3-sec3.5(*1)}
   \int_0^t ds \int_{\R^k} dy\, G_a^2(s, x - y) &= \int_0^t ds\, \frac{1}{(2\pi)^{k}} \int_{\R^k} d \xi\, \vert \exp(i \xi \cdot x - s \vert \xi \vert^a) \vert^2\notag\\
    &= \frac{1}{(2\pi)^{k}} \int_0^t ds \int_{\R^k} d \xi\, \exp(- 2 s \vert \xi \vert^a)\notag\\
    &=  \frac{1}{(2\pi)^{k}} \int_{\R^k} d \xi\,  \frac{1- \exp(- 2 t \vert \xi \vert^a)}{2 t \vert \xi \vert^a}\notag\\
   &\leq \tilde c_{a,k,t} \int_{\R^k} d \xi\,  \frac{1}{1 + \vert \xi \vert^a}
   < \infty   
   \end{align}
since $a > k$. Hence Assumption \ref{ch4-a1-1'} is satisfied.

   Let $u_0 \in L^2(\R^k)$. The solution to the  homogeneous PDE $\cL u = 0$ (the {\em fractional heat equation})\index{equation!fractional heat}\index{fractional!heat equation}\index{heat!equation, fractional} with initial condition $u_0$ is
   \beqn
    I_0(t,x) = \begin{cases}
      \int_{\R^k} dy\, G_a(t, x - y) u_0(y),&   (t, x) \in\, ]0, \infty[ \times \R^k,\\
               u_0(x),&             (t,x) \in \{0\}\times \R^k.
               \end{cases}
               \eeqn
This is well-defined because for fixed $(t,x)$ with $t > 0$, $y \mapsto G_a(t, x - y)$ belongs to $L^2(\R^k)$, as can be checked via its Fourier transform (see the computations in \eqref{ch3-sec3.5(*1)}).

According to Definition \ref{ch4-d1-1'}, the random field solution to the SPDE \eqref{frac-linear} is given by
\beq
\label{ch3-sec3.5(*1aa)}
   u(t,x) = I_0(t,x) + \int_0^t \int_{\R^k} G_a(t-r, x-y)\, W(dr, dy).   
   \eeq
The random field $u = (u(t,x),\, (t,x) \in \R_+ \times \R^k)$ is Gaussian with $E(u(t,x)) = I_0(t,x)$ and finite variance (by \eqref{ch3-sec3.5(*1)}).
\medskip

\noindent{\em Hölder continuity of the sample paths}
\medskip

We consider first the homogeneous equation $(u_0 \equiv 0)$, with solution $v = (v(t,x),\ (t,x) \in \R_+ \times \R^k)$ defined for $x \in \R^k$ by $v(0, x) = 0$ and for $t > 0$ and $x \in \R^k$, by
\beq
\label{ch3-sec3.5(*1a)}
   v(t,x) = \int_0^t \int_{\R^k} G_a(t-r, x-y)\, W(dr, dy).    
   \eeq

   \begin{prop}
   \label{ch3-sec3.5-p1}
   Fix $a > k \geq 1$, $T > 0$ and $L > 0$. Let $(v(t,x),\, (t,x) \in \R_+ \times \R^k)$ be the random field given in \eqref{ch3-sec3.5(*1aa)}. There is $C = C_{a,k,L} < \infty$ such that, for all $(t,x), (s,y) \in\, ]0, T] \times [-L, L]^k$,
   \begin{align}
   \label{ch3-sec3.5(*2)}
  & E\left[(v(t,x) - v(s, y))^2\right]\notag\\
  &\qquad \leq C\left[\vert t - s \vert^{\half - \frac{k}{2a}} + \vert x - y \vert^{\frac{a-k}{2} \wedge 1}
   \left(1 +  1_{\{a = 2 + k\}} \log\left(\frac{2L}{\vert x-y \vert}\right)\right)\right]^2.     
   \end{align}
Therefore, for any $\alpha \in\, ]0, \frac{a-k}{2a}[$ and any $\beta \in\, ]0, \frac{a-k}{2} \wedge 1[$, there exists a version of $v = (v(t,x),\, (t,x) \in\, ]0, \infty[ \times \R^k)$ with locally jointly Hölder continuous sample paths with exponents $(\alpha, \beta)$.
\end{prop}
\begin{proof}
Fix $s, t \geq 0$ and $x, y \in \R^k$. By the Itô isometry,
\beqn
   E\left[(v(t,x) - v(s, y))^2\right] =  \int_0^t dr \int_{\R^k} dz\ \left(G_a(t-r, x-z) - G_a(s-r, y-z)\right)^2.
   \eeqn
By Lemma \ref{ch3-sec3.5-lA} (c) in Appendix \ref{app2}, the right-hand side is, up to a multiplicative constant $C$, bounded above by
\beqn
    \vert t - s \vert^{1-k/a} + \vert x-y \vert^{(a - k) \wedge 2} \left(1 +  1_{\{a = 2 + k\}} \log\left(\frac{2L}{\vert x-y \vert}\right)\right),
    \eeqn
which is equivalent to \eqref{ch3-sec3.5(*2)}.
\end{proof}

\begin{remark}
\label{ch3-sec3.5-r1a}
Write $v_a$ instead of $v$ in order to emphasize the dependence of v on the parameter $a$.
Proposition \ref{ch3-sec3.5-p1} implies that for $t > 0$ and $a > k + 2$, $x \mapsto v_a(t,x)$ is almost Lipschitz continuous. In fact, as the parameter $a$ increases, so does the smoothness of $x \mapsto v_a(t,x)$. Indeed, let  $a > k + 2n$, $n\in\N^\ast$.  Then, for integers $ i_1, \dots, i_k$ with $i_1 + \cdots + i_k = n$, let
\beqn
   G_a^{i_1,\dots, i_k}(t,x) = \frac{\partial^n}{\partial x_1^{i_1} \cdots \partial x_k^{i_k}} G_a(t,x)
   \eeqn
be a weak derivative of $G_a$ in the sense of Sobolev spaces (\cite[Chapter 5]{evans}. Then
\beqn
     \cF( G_a^{i_1,\dots, i_k}(t,*))(\xi) = i^n \xi_1^{i_1} \cdots \xi_k^{i_k} \F G_a(t,*)(\xi),
     \eeqn
therefore
\begin{align*}
   \int_0^t dr \int_{\R^k} d\xi\, \vert  \F( G_a^{i_1,\dots, i_k}(t-r,*))(\xi) \vert^2
      &\leq  \int_0^t dr \int_{\R^k} d\xi\, \vert \xi\vert^{2n} \exp(- 2 (t-r) \vert \xi \vert^a) \notag\\
   &= \int_{\R^k} d\xi\, \vert \xi\vert^{2n}\, \frac{1 - \exp(- 2 t \vert \xi \vert^a)}{2 \vert \xi \vert^a},
   \end{align*}
and this is finite provided $2n - a + k < 0$, that is, $a > k + 2n$. In this case, $G_a^{i_1,\dots, i_k}(t - \cdot,x - *) \in L^2(\R_+ \times \R^k)$, so it is possible to check that
\beqn
     \frac{\partial^n}{\partial x_1^{i_1} \cdots \partial x_k^{i_k}} v_a(t,x) = \int_0^t dr \int_{\R^k} dy\, G_a^{i_1,\dots, i_k}(t - r,x - y)\, W(dr, dy),
     \eeqn
     in the sense of weak derivatives.

     \end{remark}

 \begin{remark}
\label{ch3-sec3.5-r1b}
Since $(t, \xi) \mapsto \vert \xi \vert^p \F G_a(t, *)(\xi)$ belongs to $L^2([t_1, t_2] \times \R^k)$ for all $p \geq 0$ and $t_2 > t_1 > 0$, one can check
 that the function $(t,x) \mapsto I_0(t,x)$ is $C^\infty(]0, \infty[ \times \R^k)$ (but the regularity at $t$ near $0$ depends on the regularity of the initial condition). Therefore, the conclusions of Proposition \ref{ch3-sec3.5-p1} are also valid for the random field $u = (u(t,x),\, (t,x) \in \R_+ \times \R^k)$ defined in \eqref{ch3-sec3.5(*1aa)}.
\end{remark}
\medskip

\noindent{\em Sharpness of the degree of Hölder continuity}
\medskip

We establish the following lower bounds on the second moment of increments of the random field $v(t,x)$.

   \begin{prop}
   \label{ch3-sec3.5-p2}
    Let $(v(t,x),\, (t,x) \in \R_+ \times \R^k)$ be the random field given in \eqref{ch3-sec3.5(*1a)}.
    \smallskip

    (a)\  Fix $0 < t_0 \leq T$ and $a \in\, ]k, k+2[$. There is a constant $c = c_{a,k,t_0,T} > 0$ such that, for all $(t,x), (s,y) \in [t_0, T] \times \R^k$ with $\vert x - y \vert \leq 1$,
    \beq
      \label{ch3-sec3.5(*2a)}
    E\left[(v(t,x) - v(s,y))^2\right] \geq  c \left(\vert t - s \vert^{1 - \frac{k}{a}} + \vert x - y \vert^{a-k}\right).
    \eeq

     (b)\  Fix $0 < t_0 \leq T$ and $a > k$. There is  a constant $c = c_{a,k} > 0$ such that, for all $s, t \in [t_0, T]$ and $x, y \in \R^k$,
  \beq
      \label{ch3-sec3.5(*2b)}
     E\left[(v(t,x) - v(s,y))^2\right] \geq  c \vert t - s \vert^{1 - \frac{k}{a}}.   
     \eeq
     \end{prop}
\begin{proof}
In this proof, we shall use the estimates in Lemmas \ref{ch3-sec3.5-lA} and \ref{ch3-sec3.5-lB} in Appendix \ref{app2}.

(a)\  Fix $0 < t_0 \leq T$ and $a \in\, ]k, k+2[$. We proceed exactly as in the proof of the lower bound in Proposition \ref{ch1'-p3}. Let $c_2 = c_4$ (respectively $c_3$) be the constant $C_3$ that appears in \eqref{ch3-sec3.5(*5)} (respectively $c$ that appears in \eqref{ch3-sec3.5(*10)}), and note that we have the three inequalities
\begin{align}
    E\left[(v(t,x) - v(s,y))^2\right] &\geq c_2\, \vert t - s \vert^{1 - \frac{k}{a}},\quad   x, y \in \R^k,\ s, t \in \R_+,        \label{ch3-sec3.5(*3a)}\\
      E\left[(v(t,x) - v(t,y))^2\right] &\geq c_3\, \vert x - y \vert^{a-k},\ \ x, y \in \R^k,\ \vert x - y \vert \leq 1,\ t \geq t_0,      \label{ch3-sec3.5(*3b)}\\
       E\left[(v(t,y) - v(s,y))^2\right] & \leq c_4\, \vert t - s \vert^{1 - \frac{k}{a}},\quad   y \in \R^k,\ s, t \in \R_+.        \label{ch3-sec3.5(*3c)}
\end{align}
Indeed, \eqref{ch3-sec3.5(*3a)} follows from the fact that for $t \geq s$,
\beqn
    E\left[(v(t,x) - v(s,y))^2\right] \geq \int_s^t dr \int_\R dy \, G_a^2(t-r, x-y) = C_3 (t-s)^{1-k/a}
    \eeqn
by \eqref{ch3-sec3.5(*5)}, the inequality \eqref{ch3-sec3.5(*3b)} follows from Lemma \ref{ch3-sec3.5-lB}, and \eqref{ch3-sec3.5(*3c)} follows from \eqref{ch3-sec3.5(*4)}.

As in the proof of Proposition \ref{ch1'-p3}, we now distinguish the two cases:
\beqn
    \vert t - s \vert^{1 - \frac{k}{a}} \geq \frac{c_3}{4 c_4} \vert x - y \vert^{a-k}\   {\text{and}}\   \vert t - s \vert^{1 - \frac{k}{a}} \leq \frac{c_3}{4 c_4} \vert x - y \vert^{a-k},
    \eeqn
and follow the remaining steps there to obtain \eqref{ch3-sec3.5(*2a)}. This completes the proof of (a).

   (b)\  The inequality \eqref{ch3-sec3.5(*2b)} follows from \eqref{ch3-sec3.5(*3a)}, which is valid for all $a > k$ by \eqref{ch3-sec3.5(*5)}.
   \end{proof}

The next proposition shows that the constraints on Hölder exponents obtained in Proposition \ref{ch3-sec3.5-p1} are sharp.
\begin{prop}
   \label{ch3-sec3.5-p3}
  Let  $(v(t,x),\ (t,x) \in \R_+ \times \R^k)$ be the random field given in \eqref{ch3-sec3.5(*1a)}.
    \smallskip

   (a)\  Fix $a \in\, ]k, k+2[$, $t > 0$, $i\in \{1,\dots, k\}$, $J \subset \R$ a closed interval with positive length, $x_1,\dots, x_{i-1}, x_{i+1},\dots, x_k \in \R$ and $\beta\in\, ]\frac{a-k}{2}, 1]$. Then a.s., the sample paths of the process
    $x_i \mapsto v(t,(x_1,\dots, x_i,\dots, x_k)),\ x_i \in J$,
are not Hölder continuous with exponent $\beta$.

   (b)\ Fix $a > k$, $x \in \R^k$, $K\subset\, ]0, \infty[$ a closed interval with positive length, and $\beta \in\, ]\half - \frac{k}{2a}, 1]$. Then a.s., the sample paths of the process $(v(t,x),\ t \in K)$ are not Hölder continuous with exponent $\beta$.
   \end{prop}
\begin{proof}
For (a), notice that by Proposition \ref{ch3-sec3.5-p2} (a), condition \eqref{app1-3.5} in Theorem \ref{app1-3-t2} is satisfied with $\alpha = \frac{a-k}{2}$. For (b), by Proposition \ref{ch3-sec3.5-p2} (b), condition \eqref{app1-3.5} is satisfied with $\alpha = \half - \frac{k}{2a}$. Therefore, the two conclusions follow from Theorem \ref{app1-3-t2}.
\end{proof}

    \section{Relationship between random field and weak solutions}
  \label{chapter3-s-new}

In this section, we give conditions under which a random field solution to an SPDE in $\R^k$ is also a weak solution.

Let $D\subset \rek$ and let $\cL$ be a partial differential operator on $\R_+ \times D$.  We assume that there is a fundamental solution (or a Green's function) $\Gamma(t,x; s,y)$ with support in $\R_+ \times D$ (for vanishing boundary conditions) that satisfies Assumption \ref{ch4-a1-1'}. Recall that by definition, for all  $\varphi \in \cC_0^\infty(]0, \infty[ \times D)$,
\beq
\label{chapter3-s-new(*0a)}
        \cL[\Gamma(\varphi)] = \varphi,   
        \eeq
        where
        \beqn
        [\Gamma(\varphi)](t,x) = \int_0^t ds\int_D dy\, \Gamma(t,x;s,y) \varphi(s,y) .
 \eeqn
 Consider the SPDE
\beq
\label{chapter3-s-new(*0)}
     \cL u = \dot W\   {\text{ in}}\  ]0, \infty[ \times D,   
     \eeq
with vanishing boundary conditions, and vanishing initial conditions.

\begin{def1}
\label{chapter3-s-new*0} A random linear functional 
\[
u=(u(\varphi),\, \varphi\in\cC_0^\infty(]0, \infty[ \times D))
\]
is a weak solution\index{weak!solution}\index{solution!weak} to \eqref{chapter3-s-new(*0)} if, for all  $\varphi \in \cC_0^\infty(]0, \infty[ \times D)$,
\beq
\label{sort(*1a)}
     \int_0^\infty ds \int_D dy\, ([\Gamma^* (\varphi)](s,y))^2 < \infty
     \eeq
and
\beqn
   u(\varphi) = \int_0^\infty \int_D\,  [\Gamma^*(\varphi)](s,y)\, W(ds, dy),
   \eeqn
where
\beqn
       [\Gamma^* (\varphi)](s,y) = \int_s^\infty dt  \int_D dx\, \varphi(t,x)\, \Gamma(t,x; s,y).
\eeqn
\end{def1}

Next, we formulate conditions on $\cL$ and $\Gamma$ that will play a role in Proposition \ref{chapter3-s-new*2}.

 \begin{assump}
 \label{chapter3-s-new*L}
\begin{description}
\item{(i)} For all  $\psi \in \cC_0^\infty(]0, \infty[ \times D)$,
\beq
      \int_0^\infty dt \int_D dx\,  \vert \psi(t,x) \vert \left(\int_0^t ds \int_D dy\ \Gamma^2(t,x; s,y)\right)^\half < \infty.
      \eeq
\item{(ii)} For all  $\psi \in \cC_0^\infty(]0, \infty[ \times D)$, 
    $\Gamma(\cL \psi) = \psi$.
\item{(iii)} Let
     \beqn
    \cD = \cC_0^\infty(]0, \infty[ \times D) \cup \{\Gamma^*(\varphi): \varphi \in \cC_0^\infty(]0, \infty[ \times D) \}.
    \eeqn
The operator $\cL$ has an adjoint $\cL^\star$ such that, for all $\psi \in \cD$ and $\varphi \in \cC_0^\infty(]0, \infty[ \times D)$, $\cL^\star \psi \in L^1(]0, \infty[ \times D)$ and
   \beqn
   \int_0^\infty dt \int_D dx\, \cL^\star \psi(t,x)\, \varphi(t,x) = \int_0^\infty dt \int_D dx\,  \psi(t,x)\, \cL \varphi(t,x).
   \eeqn
\item{(iv)}  For all  $\varphi,  \psi \in \cC_0^\infty(]0, \infty[ \times D)$,
\beqn
    \int_0^\infty dt \int_D dx \int_0^t ds \int_D dy\, \vert \varphi(t,x)\,  \Gamma(t,x; s,y)\, \cL \psi(s,y) \vert < \infty.
    \eeqn
\end{description}    
\end{assump}

\begin{prop}
\label{chapter3-s-new*2}
     (a) Under Assumptions \ref{ch4-a1-1'} and \ref{chapter3-s-new*L} (i),
 a strong solution $u$ to \eqref{chapter3-s-new(*0)} (in the sense of Definition \ref{ch4-d1-1'}) is also a weak solution to \eqref{chapter3-s-new(*0)} (in the sense of \eqref{4.6-1'} and Definition \ref{chapter3-s-new*0});

  (b) Suppose that for all $\psi \in \cC_0^\infty(]0, \infty[ \times D)$, \eqref{sort(*1a)} holds. Let $u$ be a random linear functional defined in particular on $\cD$ such that for all $\varphi \in \cD$,
   \beq
     \label{chapter3-s-new(1*a)}
   u(\cL^\star \varphi) = W(\varphi). 
   \eeq
Under Assumptions \ref{chapter3-s-new*L} (ii), (iii) and (iv), $u$ is a weak solution of \eqref{chapter3-s-new(*0)} (in the sense of Definition \ref{chapter3-s-new*0}).
\end{prop}
\begin{proof}

 (a) Note that by Assumption \ref{chapter3-s-new*L} (i), the hypotheses of the stochastic Fubini's theorem \ref{ch1'-tfubini}, with the measure $\vert \varphi(t,x) \vert dtdx$ and the function $(t,x, s,y) \mapsto \Gamma(t,x; s,y)$, are satisfied. Therefore
\beqn
   u(\varphi) = \langle u, \varphi \rangle = \int_{\R_+} dt \int_D dx \, \varphi(t,x)\left(\int_0^t \int_D \Gamma(t,x; s,y)\, W(ds, dy)\right)
   \eeqn
is well-defined.

By Assumption \ref{chapter3-s-new*L} (i), we can apply the stochastic Fubini's theorem \ref{ch1'-tfubini} to the positive and negative parts of $dtdx \varphi(t,x)$  to see that
\beqn
       \Vert \Gamma^* (\varphi)(\cdot,*) \Vert_{L^2([0, T] \times D)} < \infty
       \eeqn
and
\begin{align*}
    \langle u, \varphi \rangle &= \int_0^\infty \int_D \left(\int_s^\infty dt \int_D dx\, \varphi(t,x)\, \Gamma(t,x; s,y)\right)W(ds, dy) \\
   & = \int_0^\infty \int_D  [\Gamma^*(\varphi)](s,y)\, W(ds, dy).
   \end{align*}
Hence, the conditions of Definition \ref {chapter3-s-new*0} are satisfied. This proves (a).
\smallskip

    (b) We first check that  for $\psi \in \cC_0^\infty(]0, \infty[ \times D)$,
    \beq
    \label{chapter3-s-new(*1)}
     \cL^\star[\Gamma^*(\psi)] = \psi.    
     \eeq
Indeed, for $\varphi \in \cC_0^\infty(]0, \infty[ \times D)$, by Assumption \ref{chapter3-s-new*L} (iii),
\begin{align*}
     &\int_0^\infty ds \int_D dy\, \cL^\star \Gamma^*(\psi)(s,y)\, \varphi(s,y)
     = \int_0^\infty ds \int_D dy\,  \Gamma^*(\psi)(s,y)\, \cL \varphi(s,y)\\
    &\qquad= \int_0^\infty ds \int_D dy\, \left(\int_s^\infty dt  \int_D dx\, \psi(t,x)\, \Gamma(t,x; s,y)\right) \cL \varphi(s,y).
    \end{align*}
By Assumption \ref{chapter3-s-new*L} (iv), we can apply Fubini's theorem to see that this is equal to
\begin{align*}
     &\int_0^\infty dt \int_D dx\, \psi(t,x) \left(\int_0^t ds \int_D dy\, \Gamma(t,x; s,y)\, \cL \varphi(s,y)\right)\\
  &\qquad  = \int_0^\infty dt \int_D dx\, \psi(t,x)\, \Gamma(\cL \varphi)(t,x)\\
  &\qquad  = \int_0^\infty dt \int_D dx\, \psi(t,x)\, \varphi(t,x),
  \end{align*}
where the last equality follows from Assumption \ref{chapter3-s-new*L} (ii). Therefore, \eqref{chapter3-s-new(*1)} holds.

    Now suppose that $u(\cL^\star \varphi) = W(\varphi)$, for all $\varphi \in \cD$. Fix $\psi \in \cC_0^\infty(]0, \infty[ \times D)$ and let $\varphi = \Gamma^*(\psi)$.
By \eqref{sort(*1a)}, $\varphi \in \cD \cap L^2([0, T] \times D)$ and by \eqref{chapter3-s-new(*1)}, $\cL^\star \varphi = \psi$, so by  \eqref{chapter3-s-new(1*a)},
   $u(\cL^\star \varphi) = W(\varphi)$,
that is
     $u(\psi) = W( \Gamma^\ast(\psi))$,
so the conditions of Definition \ref{chapter3-s-new*0} are satisfied.
   \end{proof}

   \begin{remark}
   \label{chapter3-s-new-r*1}
   Assumption \ref{chapter3-s-new*L}
 is satisfied in particular for the heat and wave operators in spatial dimension $1$ considered in Sections \ref{ch4-ss2.2-1'} and \ref{ch4-ss2.3-1'}.
   \end{remark}

\section{Notes on Chapter \ref{chapter1'}}
\label{ch1'-notes}

In the classical theory of SPDEs,
 there are three main notions of solution: the {\em variational solution} \cite{pardoux79}, \cite{rozovskii-1990}, the {\em random field solution} \cite{walsh}, and the {\em mild solution} \cite{dz}. They all stem from the theory of PDEs. In this book, we  take the  random field approach.

We put the focus on three fundamental examples: the stochastic heat and wave equations
 in spatial dimension $1$, and the stochastic heat equation in $\R^k$ with a fractional Laplacian. This
 choice is  motivated by our decision to give priority to simplicity at this early stage in the introduction to the theory, and also by the fact that most SDPEs that are needed to describe physical phenomena are extensions of these three. In addition, as mentioned in \cite{zambotti-2021} (see also the references therein), starting in the 1950s, these SPDEs were already proposed by several physicists who sought to describe various phenomena arising in random evolutions. For instance, the linear stochastic heat equation appears in \cite{edwards-wilkinson-1982} as a simplified model for the evolution of the profile of a growing interface, a question that was later considered with more completeness and depth in \cite{k-p-z-1986}, and gave rise to the famous nonlinear KPZ equation.
 
 With the background provided in this chapter, linear SPDEs with space-time white noise, but defined by many other partial differential operators, can be handled, including, for instance, the fractional wave operator (see \cite{dalang-mueller2003}), the damped heat or wave operators or strictly parabolic operators (\cite{dalang}, \cite{ss-v2003}). In the case of boundary value problems, one can also consider mixed boundary conditions. 
 
In many problems in physics, the stochastic Poisson equation \eqref{poiss} is as important as the stochastic heat and wave equations. However, when using the classical non-anticipative stochastic calculus, it is not possible to give a sound meaning to the random field solution to the equation when the noise term is multiplied by a nonlinear function of the solution. For this reason, and also because our main interest is on evolution equations, we do not develop this example in detail. However, using the anticipative Skorohod integral \cite{skorohod-1975}, the solution to non-linear versions of the Poisson equation can be rigorously defined: see  \cite{rozovskii-2009}.

Just as for PDEs, regularity properties of the sample paths of random field solutions lie in the core of the theory of SPDEs, for example in questions of well-posedness. They are also crucial in finding optimal numerical schemes and in the analysis of fractal properties, such as the Hausdorff dimension of the range of the sample paths. Although the most interesting setting is that of non-linear SPDEs (see Chapter \ref{ch1'-s5}), obtaining sharp results in the linear case provides the critical exponents to be targeted in the corresponding non-linear equations. Also, the simplicity of the linear case allows for a smooth introduction into the methodology (see \cite{walsh} for some early work). For stochastic wave equations, our approach seems to avoid some of the difficulties discussed in \cite[Section 4.3.1]{lototsky-rozovsky-2017}.

The random field solutions to the SPDEs considered in this chapter are Gaussian processes. The proofs of the regularity results rely therefore on the theory of Gaussian processes, but to a large extent, also on precise upper and lower bounds on $L^2$-norms of increments in time and in space of the fundamental solutions (or the Green's functions) corresponding to the partial differential operators that define the SPDEs. In fact, the same bounds will be used again in Chapter \ref{ch1'-s5} in the study of regularity of solutions to non-linear SPDEs. The needed analytic results have been scattered throughout the literature in numerous articles.  For this reason, we give in Appendix \ref{app2} a systematic compilation of sharp estimates, along with their proofs.

 Section \ref{ch1'-ss3.3.2} gives a different proof of \cite[Theorem 3.3]{khosh}, which builds on \cite{Lei-Nualart-2009} for part (a) and \cite{Foondun-Khoshnevisan-Mahboubi-2015} for part (b). Proposition \ref{rd04_23p1} is adapted from \cite{m-t2002}.


\chapter{Non-linear SPDEs driven by space-time white noise}
\label{ch1'-s5}

\pagestyle{myheadings}
\markboth{R.C.~Dalang and M.~Sanz-Sol\'e}{Non-linear SPDEs}

This chapter is devoted to the study of non-linear SPDEs defined by a linear partial differential operator $\mathcal{L}$  on $\re_+\times \rek$ and driven by space-time white noise. We begin by defining the notion of solution and the basic assumptions that we use in this chapter. Then we consider the case where the nonlinearities in the equation are globally Lipschitz continuous functions (see Section \ref{ch4-section1} for the precise hypotheses) and prove a theorem on existence and uniqueness of random field solutions.
We also formulate sufficient conditions on the fundamental solution associated to $\cL$ that ensure the H\"older continuity of the sample paths of the solution. In Section \ref{ch1'-s6}, we  apply these results to examples of SPDEs in spatial dimension $1$, namely the stochastic heat and wave equations and a fractional stochastic heat equation. In Section \ref{ch1'-s60}, we present a theorem that concerns approximation of the SPDEs studied in Section \ref{ch4-section1} by a sequence of SPDEs driven by finite-dimensional projections of the noise. Finally, in Section \ref{ch1'-s7}, we study the existence and uniqueness of solutions assuming that the nonlinearities are locally  Lipschitz functions which may depend on the entire past of the solution.


\section{Formulation and basic definitions}
\label{ch4-section0}

Throughout this chapter, $D\subset \rek$ is a bounded or unbounded domain with smooth boundary, and $T>0$. Let $W$ be a space-time white noise as defined in
Proposition \ref{rd1.2.19} and let $(\cF_s,\  s\in[0,T])$ be a filtration as described at the beginning of Section \ref{ch2new-s2}.
A non-linear SPDE is an equation of the form
\begin{align}
\label{ch1'-s5.0}
\mathcal{L}u(t,x) = \sigma(t,x,u(t,x)) \dot W(t,x) + b(t,x,u(t,x)), \quad (t,x)\in\, ]0,T]\times D,
\end{align}
with given initial conditions, together with boundary conditions if $D$ has boundaries. As in Section \ref{ch1'-s1}, $\mathcal{L}$ is a linear partial differential operator on $\re_+\times \rek$.

We interpret equation \eqref{ch1'-s5.0} as the integral equation
\begin{align}
\label{ch1'-s5.1}
u(t,x) &= I_0(t,x) + \int_0^t \int_D \Gamma(t,x;s,y) \sigma(s,y,u(s,y))\, W(ds,dy)\notag\\
&\qquad +  \int_0^t \int_D \Gamma(t,x;s,y) b(s,y,u(s,y))\, ds dy,
\end{align}
$(t,x)\in[0,T]\times D$. Here $\Gamma$ denotes the fundamental solution or the Green's function associated to the operator $\mathcal{L}$ (and the boundary conditions, if present), as has been introduced in Section \ref{ch1'-s1}. The term $I_0(t,x)$ is the solution to the homogeneous PDE $\mathcal{L} I_0 =0$ (with the same initial and boundary conditions as for \eqref{ch1'-s5.0}). For example, if $D\subset \rek$, $\mathcal{L}$ is of first order in $t$, and $u(0,x):= u_0(x)$, then
\beq\label{rde4.0.3}
I_0(t,x) = \int_{D} \Gamma(t,x;0,y) u_0(y)\, dy, \qquad\text{for } t > 0,
\eeq
and $I_0(0, x) = u_0(x)$. For the stochastic integral term, we use the theory developed in Chapter \ref{ch2}.

\begin{def1}
\label{def4.1.1} A {\em random field solution}\index{solution!random field}\index{random field!solution} to \eqref{ch1'-s5.0} (or to \eqref{ch1'-s5.1}) is a jointly measurable and adapted real-valued random field
\beqn
u=\left(u(t,x),\, (t,x)\in[0,T]\times D\right)
\eeqn
 such that, for all $(t,x)\in[0,T]\times D$,
 the two integrals in \eqref{ch1'-s5.1} are well-defined and  \eqref{ch1'-s5.1} holds a.s.
\end{def1}

Throughout this chapter, we will consider the following assumptions on the function $\Gamma$.
\medskip

\noindent ($\bf H_\Gamma$) {\em Assumptions on the fundamental solution/Green's function}\label{rdHGamma}
\begin{description}
\item{(i)} The mapping $(t,x;s,y)\mapsto \Gamma(t,x;s,y)$ from $\{(t,x;s,y)\in[0,T]\times D\times [0,T]\times D: 0\le s<t\le T\}$ into $\IR$ is jointly measurable.
\item{(ii)} There is a Borel function $H: [0,T]\times D^2 \longrightarrow \IR_+$ such that
\beqn
\vert\Gamma(t,x;s,y)\vert \le H(t-s,x,y), \qquad 0\le s<t\le T, \quad x,y\in D.
\eeqn
\item{(iiia)} If in \eqref{ch1'-s5.0}$\  \sigma \not\equiv 0$, then
\beqn
\int_0^T ds \ \sup_{x\in D}\int_D dy\, H^2(s,x,y) < \infty.
\eeqn
\item{(iiib)} If in \eqref{ch1'-s5.0}$\  b \not\equiv 0$, then
\beqn
\int_0^T ds \ \sup_{x\in D}\int_D dy\, H(s,x,y) < \infty.
\eeqn
\end{description}
Notice that if $D$ is bounded, then condition (iiib) follows from (iiia).

\section{Non-linear SPDEs with globally Lipschitz coefficients}
\label{ch4-section1}

Throughout this section, we assume the following.
\medskip

\noindent ($\bf H_I$)\label{rdHI} {\em Assumptions on the initial conditions}
\smallskip

\quad The function $(t,x) \mapsto I_0(t,x)$ is Borel and bounded over $[0,T]\times D$.
\bigskip

\noindent ($\bf H_L$)\label{rdHL} {\em Assumptions on the coefficients $\sigma$ and $b$}
\begin{description}
\item{(iv)} {\em Measurability and adaptedness.} The functions $\sigma$ and $b$  are defined on $[0,T]\times D\times \IR\times \Omega$ with values in $\IR$ and are jointly measurable, that is, $\mathcal{B}_{[0,T]}\times \mathcal{B}_D \times\mathcal{B}_{\IR}\times \tf$-measurable. These two functions are also adapted to $(\tf_s,\, s\in[0,T])$, that is, for fixed $s \in [0,T]$, $(y,z,\omega)\mapsto \sigma(s,y,z;\omega)$ and $(y,z,\omega)\mapsto b(s,y,z;\omega)$ are $\cB_D\times \cB_{\IR}\times \cF_s$-measurable. \item{(v)} {\em Global Lipschitz condition.}
There exists $C:=C(T,D)\in \IR_+$ such that for all $(s,y,\omega)\in [0,T]\times D\times \Omega$ and $z_1,z_2\in \IR$,
\beqn
\vert \sigma(s,y,z_1;\omega)-\sigma(s,y,z_2;\omega)\vert + \vert b(s,y,z_1;\omega)-b(s,y,z_2;\omega)\vert  \le C|z_1-z_2|.
\eeqn
\item{(vi)} {\em Uniform linear growth.} There exists a constant $\bar c :=\bar c(T,D)\in \IR_+$ such that for all $(s,y,\omega)\in [0,T]\times D \times \Omega$ and all $z\in \IR$,
\beqn
\vert \sigma(s,y,z;\omega)\vert + \vert b(s,y,z;\omega)\vert \le \bar c \left(1+|z|\right).
\eeqn
\end{description}

Notice that, in the case where the functions $\sigma$ and $b$ do not depend on the first two parameters $s\in[0,T]$, $y\in D$, the assumption (v) implies (vi).

If $\sigma(\cdot, *, z)$ is a constant function of $z$ but $b(\cdot, *, z)$ is not, then we refer to \eqref{ch1'-s5.0} as a {\em nonlinear SPDE with additive noise,}\index{additive noise}\index{SPDE!nonlinear}\index{noise!additive}\index{nonlinear!SPDE} and if $\sigma(\cdot, *, z)$ is not a constant function of $z$, then we refer to \eqref{ch1'-s5.0} as a {\em nonlinear SPDE with multiplicative noise.}\index{multiplicative noise}\index{noise!multiplicative}

\subsection{Existence and uniqueness of solutions}
\label{ch1'-ss5.1}

The next statement is a theorem on existence and uniqueness of solutions to a class of SPDEs with globally Lipschitz coefficients.

\begin{thm}
\label{ch1'-s5.t1}
Under $(\bf H_\Gamma)$, $(\bf H_I)$ and  $(\bf H_L)$, there exists a random field solution $$u=\left(u(t,x),\, (t,x)\in[0,T]\times D\right)$$ to \eqref{ch1'-s5.0}. In addition, for any $p>0$,
\beq
\label{ch1'-s5.2}
\sup_{(t,x)\in[0,T]\times D} E\left[\vert u(t,x)\vert^p\right] <\infty,
\eeq
and the solution $u$ is unique (in the sense of versions) among random field solutions that satisfy \eqref{ch1'-s5.2} with $p=2$.
\end{thm}

\begin{remark}\label{ch1'-rmk1}
(a) The uniqueness statement is different from that for SDEs, in which no condition on moments is required.

(b) According to Remark \ref{after-rdprop2.6.2}, any jointly measurable and adapted random field has an optional version, so this is the case for the solution $u$ to \eqref{ch1'-s5.0}.

(c) In the proof of the theorem, we will take $p\ge 2$. Since $L^p$-norms increase with $p$, this suffices to have \eqref{ch1'-s5.2} for any $p>0$.

(d) Extensions of this theorem to a wider class of initial conditions is possible in certain cases: see Section \ref{rdrough}.
\end{remark}
\noindent{\em Proof of Theorem \ref{ch1'-s5.t1}.}\  We will apply a fixed point argument based on the Picard iteration scheme
\begin{align}
\label{rde4.1.2}
 u^0(t,x) &= I_0(t,x),\notag\\
 u^{n+1}(t,x) & = I_0(t,x) + \int_0^t \int_D \Gamma(t,x;s,y) \sigma(s,y,u^n(s,y))\, W(ds,dy)\notag\\
 &\qquad+ \int_0^t \int_D \Gamma(t,x;s,y) b(s,y,u^n(s,y))\, ds dy, \qquad n\ge 0.
\end{align}
 \noindent{\em Step 1.} We will prove below by induction on $n$ that for each $n\ge 0$, the process
 \beqn
 u^n = \left(u^n(t,x),\, (t,x)\in [0,T] \times D\right)
 \eeqn
 is well-defined, jointly measurable (that is, has a jointly measurable version) and adapted, and satisfies
 \beq
 \label{ch1'-s5.2a}
\sup_{(t,x)\in [0,T] \times D}E\left[\vert u^n(t,x)\vert^p\right] <\infty,
 \eeq
 for any $p\ge 2$, hence for any $p>0$.

First, we explain why these properties of $u^n$ imply that the stochastic integral in \eqref{rde4.1.2} is well-defined according to Definition \ref{ch1'-s4.d1}. Indeed, assuming these properties of $u^n$, we see that the map
\beq\label{rde4.1.3}
   (s,y,\omega) \mapsto \sigma(s,y,u^n(s,y;\omega);\omega)
\eeq
from $([0,T]\times D\times \Omega, \cB_{[0,T]}\times \cB_D \times \cF)$ into $(\IR,\cB_\IR)$ is measurable, since it is the composition of the maps $(s,y,\omega) \mapsto (s,y,u^n(s,y;\omega);\omega)$ from $([0,T]\times D\times \Omega, \cB_{[0,T]}\times \cB_D \times \cF)$ to $([0,T]\times D\times\IR\times \Omega, \cB_{[0,T]}\times \cB_D\times \cB_\IR \times \cF)$ with the map $(s,y,z,\omega) \mapsto \sigma(s,y,z;\omega)$ from $([0,T]\times D\times\IR\times \Omega, \cB_{[0,T]}\times \cB_D\times \cB_\IR \times \cF)$ to $(\IR,\cB_\IR)$.

The map in \eqref{rde4.1.3} is also adapted, since for fixed $s\in [0,T]$, the map $(y,\omega) \mapsto \sigma(s,y,u^n(s,y;\omega);\omega)$ is the composition of the map $(y,\omega) \mapsto\break (y,u^n(s,y;\omega);\omega)$ from $(D\times \Omega,\cB_D\times \cF_s)$ into $(D\times \IR\times \Omega, \cB_D\times \cB_\IR\times \cF_s)$ with the map $(y,z,\omega)\mapsto \sigma(s,y,z;\omega)$ from $(D\times \IR\times \Omega, \cB_D \times \cB_\IR \times \cF_s)$ into $(\IR,\cB_\IR)$.

We now discuss the square-integrability of the integrand in \eqref{rde4.1.2}. Notice that this integrand is of the form \eqref{rd2.2.12a}, with $Z(s,y):= \sigma(s,y,u^n(s,y))$. We check condition \eqref{rd2.2.12b}. By the uniform linear growth of $\sigma$,
\begin{align*}
 &  \sup_{(s,y)\in [0,T] \times D} E[\sigma^2(s,y,u^n(s,y))] \leq \bar c^2 \sup_{(s,y)\in [0,T] \times D} E\left[(1+\vert u^n(s,y)\vert)^2\right]\\
   &\qquad\qquad \qquad\leq 2\ \bar c^2 \left(1+ \sup_{(s,y)\in [0,T] \times D} E[(u^n(s,y))^2]\right)
   <\infty
\end{align*}
by \eqref{ch1'-s5.2a}. Further, condition \eqref{rd2.2.12c} is satisfied since
$$
   \int_0^t ds \sup_{x\in D} \int_D dy\, \Gamma^2(t,x;s,y) \leq
   \int_0^t ds \sup_{x\in D} \int_D dy\, H^2(t-s,x,y) < \infty
$$
by Assumption (iiia). These considerations show that the stochastic integral in \eqref{rde4.1.2} is well-defined.

In a similar way, using Assumption (iiib), one checks that the pathwise (Lebesgue) integral in \eqref{rde4.1.2} is also well-defined.

In order to start our induction, let $n=0$ and fix $p\ge 2$. By Assumption ($\bf H_I$), $u^0 = I_0$ satisfies the properties described at the beginning of this Step 1.

Assume now that for some $n\ge 0$, the process
\beqn
\left(u^n(t,x),\, (t,x)\in[0,T]\times D\right)
\eeqn
 is well-defined, jointly measurable and adapted, and for any $p\ge 2$,
\beqn
 \sup_{(t,x)\in[0,T]\times D}E\left[\vert u^n(t,x)\vert^p\right] < \infty.
\eeqn
According to what we have just established, $u^{n+1}= (u^{n+1}(t,x),\, (t,x) \in [0,T]\times D)$ given in \eqref{rde4.1.2} is well-defined. We want to show that $u^{n+1}$ is jointly measurable, adapted, and \eqref{ch1'-s5.2a} is satisfied with $n$ there replaced by $n+1$.

Define
\begin{align*}
\cI^n(t,x) &:= \int_0^t \int_D \Gamma(t,x;s,y) \sigma(s,y,u^n(s,y))\, W(ds,dy),\\
\cJ^n(t,x) &:=  \int_0^t ds \int_D dy\, \Gamma(t,x;s,y) b(s,y,u^n(s,y)).
\end{align*}
The existence of a jointly measurable and adapted modification of $\cI^n(t,x)$ follows from Assumption (${\bf H_\Gamma}$) (i), the properties of the map in \eqref{rde4.1.3} and Proposition \ref{rdprop2.6.2}.  The existence of a jointly measurable and adapted modification of $\cJ^n(t,x)$ follows  similarly from Proposition \ref{rdprop2.6.4}. These modifications are now used in \eqref{rde4.1.2} to define a jointly measurable and adapted version of $u^{n+1}$, which we again denote by $u^{n+1}$ and which therefore satisfies the measurability conditions listed at the beginning of Step 1.

We now check \eqref{ch1'-s5.2a} for $n+1$. Set
\beq
\label{j}
J_1(s) = \sup_{x\in D} \int_D H^2(s,x,y)\, dy,\qquad
J_2(s) = \sup_{x\in D} \int_D H(s,x,y)\, dy.
\eeq
By (iiia) and (iiib) in Assumption (${\bf H_\Gamma}$),
\beq\label{rdja}
\int_0^T  J_1(s)\, ds < \infty \qquad \text{and}\qquad
\int_0^T  J_2(s)\, ds < \infty.
\eeq
Fix $p\ge 2$. By Burkholder's inequality (see \eqref{ch1'-s4.4}), then H\"older's inequality, and Assumptions (ii) and (vi), we see that
\begin{align}
\label{ch1'-s5.6a}
E\left[\vert \cI^n(t,x)\vert^p\right]
  & \le C_p E\left[\left\vert\int_0^t ds\int_D dy\left(\Gamma(t,x;s,y) \sigma(s,y,u^n(s,y))\right)^2\right\vert^{\frac{p}{2}}\right] \notag\\
& \le \tilde C_p \left(\int_0^t ds\, J_1(t-s)\right)^{\frac{p}{2}-1} \int_0^t ds \int_D dy\, H^2(t-s,x,y)\notag\\
& \qquad \times \left(1+E[\vert u^n(s,y)\vert^p]\right)\notag\\
 & \le \tilde C_p
 \left(\int_0^t ds\, J_1(s)\right)^{\frac{p}{2}-1}\notag\\
 &\quad\times \int_0^t ds \sup_{x\in  D}\left(1+E[\vert u^n(s,x)\vert^p]\right) J_1(t-s).
\end{align}

From H\"older's inequality and Assumptions (ii) and (vi), we also see that
\begin{align}
 E\left[\vert \cJ^n (t,x)\vert^p\right] &\le \left(\int_0^t ds \int_D dy\, \vert \Gamma(t,x;s,y)\vert\right)^{p-1}\notag\\
&\quad\times \int_0^t ds \int_D dy\, \vert \Gamma(t,x;s,y)\vert E\left[\vert b(s,y,u^n(s,y))\vert^p\right] \notag\\
 &\leq \bar c_p \left(\int_0^t ds\, J_2(s)\right)^{p-1}\notag\\
 &\quad\times \int_0^t ds \sup_{x\in  D}\left(1+E[\vert u^n(s,x)\vert^p]\right) J_2(t-s),
 \label{ch1'-s5.6b}
\end{align}
and, adding \eqref{ch1'-s5.6a} and \eqref{ch1'-s5.6b}, using \eqref{rdja} and Assumption ($\bf H_I$), we obtain
\beq
\label{ch1'-s5.7}
 \sup_{(t,x)\in [0,T] \times D} E\left[\vert u^{n+1}(t,x) \vert^p \right]
 \le \bar C_p \left(1+ \sup_{(t,x)\in [0,T] \times D}E\left[\vert  u^n(t,x)\vert^p\right]\right) <\infty.
\eeq
This proves \eqref{ch1'-s5.2a} for $n+1$ and completes the induction and Step 1.
\medskip

\noindent{\em Step 2.} We now show that the sequence of processes
\beqn
\left(u^n(t,x),\, (t,x)\in [0,T] \times D,\  n\ge 0\right),
\eeqn
converges in $L^p(\Omega)$ uniformly in $(t,x)\in [0,T] \times D$ to a process
\beqn
\left(u(t,x),\, (t,x)\in [0,T] \times D\right)
\eeqn
 that satisfies \eqref{ch1'-s5.2} and has a jointly measurable and adapted version.

Indeed, set
\beq
\label{rd08_20e8}
M_n(t)= \sup_{y\in D} E\left[\left\vert u^{n+1}(t,y)-u^n(t,y)\right\vert^p\right], \quad n\ge 0.
\eeq
With arguments similar to those used to deduce \eqref{ch1'-s5.6a} and \eqref{ch1'-s5.6b}, but applying the Lipschitz continuity of $\sigma$ and $b$ (Assumption (v)), instead of the properties of linear growth, we obtain
\beqn
M_n(t) \le C_p\int_0^t ds \, M_{n-1}(s)\left(J_1(t-s)+J_2(t-s)\right).
\eeqn

Consider the sequence of functions defined for $t\in[0,T]$ by  $f_n(t):=M_n(t)$ and let $J(t) = J_1(t)+J_2(t)$, $t\in[0,T]$. From \eqref{rdja}, we see that \eqref{A3.00} and the assumptions of the Gronwall-type
Lemma \ref{A3-l01} hold with $z_0=0$ and $z\equiv 0$. Furthermore, because of \eqref{ch1'-s5.2a}, we have
\beq
\label{ch1'-s5.7-bis}
\sup_{s\in[0,T]} M_0(s) \le C_p\sup_{s\in[0,T]}\left( E[|u^1(s,y)|^p] + E[|u^0(s,y)|^p)\right] < \infty.
\eeq
Using \eqref{rdA3.04bis} and \eqref{ch1'-s5.7-bis}, we obtain
\beq\label{rde4.1.10}
\sum_{n=0}^\infty\, \sup_{(t,x)\in [0,T] \times D}\left\Vert u^{n+1}(t,x)-u^n(t,x)\right\Vert_{L^p(\Omega)}
 < \infty.
\eeq
This implies that the sequence $\left(u^n(t,x),\, (t,x)\in [0,T] \times D\right)$, $n\ge 0$, converges in $L^p(\Omega)$, uniformly in $(t,x)\in [0,T] \times D$. That is, there exists $\left(u(t,x),\, (t,x)\in [0,T] \times D\right)$ such that
\beq\label{rde4.1.11}
   \lim_{n\to\infty}\, \sup_{(t,x)\in [0,T] \times D} \left\Vert u^{n}(t,x)-u(t,x)\right\Vert_{L^p(\Omega)} = 0.
\eeq
In fact,
$$
   u(t,x) = I_0(t,x) + \sum_{n=0}^\infty \left(u^{n+1}(t,x)-u^n(t,x) \right),
$$
where the series converges in $L^p(\Omega)$, uniformly in $(t,x)\in [0,T] \times D$. By the boundedness assumption $(\bf{H_I})$ on $I_0$ and \eqref{rde4.1.10}, $\left(u(t,x),\, (t,x)\in [0,T] \times D\right)$ satisfies \eqref{ch1'-s5.2}.

For each $(t,x)$, $u^n(t,x)$ converges to $u(t,x)$ in probability (in fact, in $L^p(\Omega)$), so $(u(t,x))$ has a jointly measurable version by Lemma \ref{rdlem9.4.2a}, which we again denote $(u(t,x))$.

For each $(t,x)$, $u(t,x)$ is $\cF_t$-measurable.
 Applying Lemma \ref{rdlem9.4.2} (a) with $(X,\cX) = (D, \B_D)$, there is a $\cB_D \times \cO$-measurable function $(x,t,\omega) \mapsto \bar u(t,x,\omega)$ such that, for all $(t,x) \in [0,T]\times D, u(t,x) = \bar u(t,x)$ a.s.  This modification $\bar u$ of $u$ is jointly measurable and adapted (in the sense of Definition \ref{ch1'-s4.d1}) since for all $t \in [0,T]$, $\cO\vert_{[0,t] \times \Omega} \subset \cB_{[0,t]} \times \F_t $. In the sequel, we use this modification and denote it $u$ instead of $\bar u$.


We have already seen that \eqref{ch1'-s5.2} is satisfied.
 Therefore, the stochastic and deterministic integrals in \eqref{ch1'-s5.1} are well-defined.

\smallskip

\noindent{\em Step 3.} We show that the stochastic process $\left(u(t,x),\, (t,x)\in [0,T] \times D\right)$ satisfies equation \eqref{ch1'-s5.1}. Indeed, let
\begin{align*}
\mathcal {I}(t,x) &= \int_0^t \int_D \Gamma(t,x;s,y) \sigma(s,y,u(s,y))\, W(ds,dy),\\
\mathcal{J}(t,x) &=  \int_0^t ds \int_D dy\, \Gamma(t,x;s,y) b(s,y,u(s,y)).
\end{align*}
Proceeding as in the proof of  \eqref{ch1'-s5.6a} and \eqref{ch1'-s5.6b}, but using the Lipschitz continuity assumption (v) instead of the linear growth (vi), we obtain
\begin{align}
\label{ch1'-s5.10}
E\left[\left\vert\mathcal {I}^n(t,x) -\mathcal {I}(t,x)\right\vert^p\right] & \le C_p \left(\int_0^t ds\, J_1(t-s)\right)^{\frac{p}{2}-1}\notag\\
&\quad\times\int_0^t ds \ J_1(t-s) \sup_{y\in D} E\left[\left\vert u^n(s,y)-u(s,y)\right\vert^p\right] \notag\\
&\le C_p \left(\int_0^t ds\, J_1(t-s)\right)^{\frac{p}{2}}\notag\\
&\quad\times \sup_{(s,x)\in [0,t] \times D}E\left[\left\vert u^n(s,x)-u(s,x)\right\vert^p\right],
\end{align}
and
\begin{align}
\label{ch1'-s5.11}
E\left[\left\vert\mathcal {J}^n(t,x) -\mathcal {J}(t,x)\right\vert^p\right] &\le \bar C_p\left(\int_0^t ds\, J_2(t-s)\right)^{p-1}\notag\\
&\quad \times\int_0^t ds \, J_2(t-s) \sup_{y\in D} E\left[\left\vert u^n(s,y)-u(s,y)\right\vert^p\right] \notag\\
&\le \bar C_p \left(\int_0^t ds\, J_2(t-s)\right)^p\notag\\
&\quad\times \sup_{(s,x)\in [0,t] \times D}E\left[\left\vert u^n(s,x)-u(s,x)\right\vert^p\right].
\end{align}
We have seen in \eqref{rde4.1.11} that the last right-hand sides of \eqref{ch1'-s5.10} and \eqref{ch1'-s5.11} converge to $0$ as $n\to\infty$.

With the notation introduced in Step 1, we have that
\beqn
u^{n+1}(t,x) = I_0(t,x) + \mathcal {I}^n(t,x) + \mathcal {J}^n(t,x).
\eeqn

The left-hand side converges to $u(t,x)$ in $L^p(\Omega)$, uniformly in $(t,x)\in [0,T] \times D$, while from \eqref{ch1'-s5.10} and \eqref{ch1'-s5.11}, the right-hand side converges to $I_0(t,x) + \mathcal {I}(t,x) + \mathcal {J}(t,x)$. Therefore, for each $(t,x)\in [0,T] \times D$,
\beq\label{red4.1.14}
   u(t,x) = I_0(t,x) + \cI(t,x) + \cJ(t,x) \qquad \text{a.s.},
\eeq
that is, equation \eqref{ch1'-s5.1} holds a.s.
\smallskip

\noindent{\em Step 4: Uniqueness}. Let
\beqn
\left(u(t,x), (t,x)\in [0,T] \times D\right),\qquad \left(\bar u(t,x), (t,x)\in [0,T] \times D\right),
\eeqn
 be two random field solutions to \eqref{ch1'-s5.1} satisfying \eqref{ch1'-s5.2}
with $p=2$. Using the same arguments as in \eqref{ch1'-s5.10}, \eqref{ch1'-s5.11}, we obtain
\begin{align*}
\sup_{x\in D} E\left[\left(u(t,x)-\bar u(t,x)\right)^2\right] & \le C_2 \int_0^t ds\ (J_1(t-s)+J_2(t-s)) \\
&\quad\qquad\times\sup_{y\in D} E\left[\left\vert u(s,y)-\bar u(s,y)\right\vert^2\right].
\end{align*}
Because condition \eqref{ch1'-s5.2} is assumed for both $u$ and $\bar u$, we can apply \eqref{A3.04bis} in Lemma \ref{A3-l01} to the constant sequence
\beqn
f(t) = f_n(t) := \sup_{x\in D} E\left[\left( u(t,x)-\bar u(t,x)\right)^2\right],
\eeqn
 with $J(t):= J_1(t)+J_2(t)$, $z_0 = 0$ and $z\equiv 0$, to get
 \beqn
 \sup_{x\in D} E\left[\left( u(t,x)-\bar u(t,x)\right)^2\right] = 0, \quad {\text {for all}}\ t\in[0,T],
 \eeqn
 therefore $\bar u$ is a version of $u$.
\qed

\begin{remark}
\label{4.2.1-remark-bounded-domain}
Considering Remark \ref{3.1-remark-bounded-domain}, when $\Gamma(t,x; s,y)$ is defined for $x \in \bar D$ and $y \in D$ (or $\bar D$), if in assumptions $({\bf H_\Gamma})$, $({\bf H_I})$, we replace  $x,y\in D$ by $x\in \bar D$ and $y\in D$, then Theorem \ref{ch1'-s5.t1} remains valid with $D$ replaced by $\bar D$.
\end{remark}


\subsection{Regularity of the sample paths}
\label{ch1'-ss5.2}

We begin this section by stating sufficient conditions for $L^p(\Omega)$-continuity of random fields defined by stochastic or deterministic integrals. These will later on be applied to the integral terms on the right-hand side of \eqref{ch1'-s5.1} and will imply the H\"older-continuity of the trajectories of the solution to \eqref{ch1'-s5.1}.

Given a random field $Z=(Z(t,x), (t,x)\in[0,T]\times D)$, where $D\subset \rek$ is a bounded or unbounded domain, define\label{rdTinftyp}
\beq
\label{bzp}
\Vert Z\Vert_{T,\infty,p} = \sup_{(t,x)\in[0,T]\times D}\Vert Z(t,x)\Vert_{L^p(\Omega)},\quad p\in [1,\infty[\ .
\eeq

In the next lemma, for $i=1,2$, the notation ${\bf\Delta}_i(t,x;s,y)$\label{rdDeltai1}  refers to an arbitrary nonnegative function defined for $(t,x), (s,y)\in \re_+\times D$.

\begin{lemma}
\label{ch1'-ss5.2-l1}
Consider $\Gamma(t,x;s,y)$ satisfying assumptions $({\bf H_\Gamma})$, and set
$\Gamma(t,x;s,y)=0$ if $t\le s$.
Assume that the random field  $Z$  is jointly measurable and adapted and satisfies $\Vert Z\Vert_{T,\infty,p}<\infty$ for some $p\ge2$.
Fix sets $I \subset [0,T]$ and $\tilde D \subset D$.

\noindent (a) Let
\beqn
u_1(t,x) = \int_0^t \int_D \Gamma(t,x;r,z) Z(r,z) W(dr,dz).
\eeqn
Suppose that
for all $(t,x), (s,y)\in I\times \tilde D$,
\beq
\label{ch1'-s5.12}
\int_0^T dr \int_D dz \left(\Gamma(t,x;r,z)-\Gamma(s,y;r,z)\right)^2 \le C^2\, [{\bf\Delta}_1(t,x;s,y)]^2,
\eeq
for some constant $C>0$.

Then for all $(t,x), (s,y)\in I\times \tilde D$,
\beq
\label{ch1'-s5.13}
\Vert u_1(t,x)-u_1(s,y)\Vert_{L^p(\Omega)}\le \tilde C_p^{1/p} C\Vert Z\Vert_{T,\infty,p}\, {\bf\Delta}_1(t,x;s,y),
\eeq
where $\tilde C_p$ is the constant of Burkholder's inequality \eqref{ch1'-s4.4.mod} and $C$ is the constant  in \eqref{ch1'-s5.12}.

\noindent (b) Let
\beqn
u_2(t,x) = \int_0^t dr \int_D dz \ \Gamma(t,x;r,z) Z(r,z).
\eeqn
Suppose that
for all $(t,x), (s,y)\in I\times \tilde D$,
\beq
\label{ch1'-s5.14}
\int_0^T dr \int_D dz\, \left\vert\Gamma(t,x;r,z)-\Gamma(s,y;r,z)\right\vert \le c\, {\bf\Delta}_2(t,x;s,y),
\eeq
for some constant $c>0$.

Then
for all $(t,x), (s,y)\in I\times \tilde D$,
\beq
\label{ch1'-s5.15}
\Vert u_2(t,x)-u_2(s,y)\Vert_{L^p(\Omega)}\le c\, \Vert Z\Vert_{T,\infty,p}\, {\bf \Delta}_2(t,x;s,y).
\eeq
\end{lemma}

\begin{remark}
\label{ch1'-ss5.2-r1}
If $D$ is bounded and \eqref{ch1'-s5.12}  holds, then by the Cauchy-Schwarz inequality, \eqref{ch1'-s5.14} also holds with ${\bf \Delta_2}= {\bf \Delta_1}$.
\end{remark}

\noindent{\em Proof of Lemma \ref{ch1'-ss5.2-l1}.}
Using the comments that follow \eqref{rd2.2.12a}, we first observe that the conditions on $\Gamma $ and $Z$ imply that the stochastic integral $u_1(t,x)$ is well-defined according to Definition \ref{ch1'-s4.d1}. Without loss of generality, we assume that $0\le s\le t\le T$.

(a) For $(t,x), (s,y)\in I\times \tilde D$, we write
\beq
\label{ch1'-s5.16}
u_1(t,x)-u_1(s,y)  =\int_0^t \int_D \left(\Gamma(t,x;r,z)-\Gamma(s,y;r,z)\right) Z(r,z)\, W(dr,dz).
\eeq
Using Burkholder's inequality, and then H\"older's inequality, we obtain
\begin{align}
\label{computations-1}
&\Vert u_1(t,x)-u_1(s,y)\Vert_{L^p(\Omega)}^p = E\left[\left\vert u_1(t,x)-u_1(s,y)\right\vert^p\right] \notag\\
&\qquad\le \tilde C_p E\left[\left(\int_0^t dr \int_D dz\left(\Gamma(t,x;r,z)-\Gamma(s,y;r,z)\right)^2 Z^2(r,z)\right)^{\frac{p}{2}}\right]\notag\\
&\qquad\le \tilde C_p \left( \int_0^t dr \int_D dz\left(\Gamma(t,x;r,z)-\Gamma(s,y;r,z)\right)^2\right)^{\frac{p}{2}-1}\notag\\
&\qquad\qquad\times \int_0^t dr \int_D dz\left(\Gamma(t,x;r,z)-\Gamma(s,y;r,z)\right)^2\ E[\vert Z(r,z)\vert^p] \notag\\
&\qquad\le \tilde C_p  \Vert Z\Vert_{T,\infty,p}^p\left[C\, {\bf \Delta}_1(t,x;s,y)\right]^p,
\end{align}
where we have used \eqref{ch1'-s5.12}. This implies \eqref{ch1'-s5.13}.
\medskip

(b) For $(t,x), (s,y)\in I\times \tilde D$, using Minkowski's inequality, we obtain
\begin{align*}
&\Vert u_2(t,x)-u_2(s,y)\Vert_{L^p(\Omega)}\\
&\qquad\le\int_0^t dr \int_0^L dz\ |\Gamma(t,x;r,z)-\Gamma(s,y;r,z)|\, \Vert Z(r,z)\Vert_{L^p(\Omega)}\\
&\qquad\le \Vert Z\Vert_{T,\infty,p}\int_0^t dr \int_0^L dz\,  |\Gamma(t,x;r,z)-\Gamma(s,y;r,z)|\\
&\qquad \le c\,\Vert Z\Vert_{T,\infty,p}\, {\bf \Delta}_2(t,x;s,y),
\end{align*}
where we have used \eqref{ch1'-s5.14}. This implies \eqref{ch1'-s5.15}.
\hfill $\Box$
\bigskip

Recall the decomposition  \eqref{red4.1.14}.
Each of the terms $I_0(t,x)$, $\mathcal{I}(t,x)$ and $\mathcal{J}(t,x)$ there will contribute to the increments of moments of $u$, as we shall now see.

\begin{prop}
\label{ch1'-ss5.2-prelude}
The assumptions are as in Theorem \ref{ch1'-s5.t1}. In addition, we suppose that there are sets $I\subset [0,T]$ and $\tilde D \subset D$ such that
$\Gamma(t,x;s,y)$ satisfies \eqref{ch1'-s5.12} and \eqref{ch1'-s5.14} of Lemma \ref{ch1'-ss5.2-l1}.
Then for any $p\ge 2$, there is a constant $0\le c_p<\infty$ such that, for all $(t,x), (s,y)\in I\times \tilde D$,
\begin{align}\label{rde4.1.21}
    \Vert \cI(t,x)-\cI(s,y)\Vert_{L^p(\Omega)} &\le c_p\, {\bf\Delta}_1(t,x;s,y),\\
    \Vert \cJ(t,x)-\cJ(s,y)\Vert_{L^p(\Omega)} &\le c_p\, {\bf\Delta}_2(t,x;s,y).
\label{rde4.1.22}
\end{align}
Therefore,
\begin{align}\nonumber
\Vert u(t,x)-u(s,y)\Vert_{L^p(\Omega)} &\le \vert I_0(t,x) - I_0(s,y)\vert\\
&\qquad  + c_p   \left[{\bf\Delta}_1(t,x;s,y) + {\bf\Delta}_2(t,x;s,y)\right].
\label{ch1'-ss5.2-t1.prelude}
\end{align}

\end{prop}
\begin{proof}
Let
\beqn
Z_1(s,y):= \sigma(s,y,u(s,y)), \qquad (s,y)\in[0,T]\times D.
\eeqn
From Assumption $({\bf H_L})$ (vi) on $\sigma$, for any $p\ge 2$, there is a constant $c < \infty$ such that
\beq
\label{ch1'-s5.17}
\Vert Z_1\Vert_{T,\infty,p} \le c\left(1+\sup_{(t,x)\in[0,T]\times D}\Vert u(t,x)\Vert_{L^p(\Omega)}\right).
\eeq
Hence, by \eqref{ch1'-s5.2}, $\Vert Z_1\Vert_{T,\infty,p}<\infty$. By hypothesis \eqref{ch1'-s5.12}, \eqref{ch1'-s5.13} holds, and therefore
$
\cI(t,x)
$
satisfies, for $(t,x),(s,y) \in I \times \tilde D$,
\beq
\label{u1}
\Vert \cI(t,x)-\cI(s,y)\Vert_{L^p(\Omega)}\le c_p\, {\bf\Delta}_1(t,x;s,y),
\eeq
for some positive constant $c_p$.

Analogously, the process
\beqn
Z_2(s,y):= b(s,y,u(s,y)), \qquad (s,y)\in[0,T]\times D,
\eeqn
satisfies $\Vert Z_2\Vert_{T,\infty,p}<\infty$. Consequently, by hypothesis \eqref{ch1'-s5.14}, \eqref{ch1'-s5.15} holds and therefore,
for $(t,x),(s,y) \in I \times \tilde D$,
\beq
\label{u2}
\Vert \cJ(t,x)-\cJ(s,y)\Vert_{L^p(\Omega)}\le \tilde c_p \, {\bf\Delta}_2(t,x;s,y).
\eeq
Finally, we obtain \eqref{ch1'-ss5.2-t1.prelude} from \eqref{red4.1.14} by adding together \eqref{u1} and \eqref{u2}.
\end{proof}


\begin{remark}
\label{continuity-nonlinear}
In view of \eqref{red4.1.14}, the (H\"older-) continuity properties of $u$ are related to those of $I_0$, $\mathcal I$ and $\mathcal J$. These can often be studied separately. In all cases where  $\mathcal I$ and $\mathcal J$ are continuous (respectively H\"older continuous), $u$ will be continuous (respectively H\"older continuous) if $I_0$ is.
\end{remark}

Consider the particular case\label{rdDeltai2}
\beq
\label{choice-rho}
{\bf\Delta}_1(t,x;s,y) = \vert t-s \vert^{\alpha_1}+ \vert x-y \vert^{\alpha_2}, \quad {\bf\Delta}_2(t,x;s,y) = \vert t-s \vert^{\beta_1}+ \vert x-y \vert^{\beta_2},
\eeq
$\alpha_1,\alpha_2,\beta_1,\beta_2 \in\,]0,1]$.
The discussion above yields
the H\"older continuity of the sample paths of the solution of \eqref{ch1'-s5.1}, as the following theorem shows.

\begin{thm}
\label{ch1'-ss5.2-t1}
Consider the hypotheses of Theorem \ref{ch1'-s5.t1}. Suppose in addition that there are sets $I\subset [0,T]$ and $\tilde D \subset D$ such that $\Gamma(t,x;s,y)$ satisfies \eqref{ch1'-s5.12} and \eqref{ch1'-s5.14} of Lemma \ref{ch1'-ss5.2-l1} with ${\bf\Delta}_1$, ${\bf\Delta}_2$ given in \eqref{choice-rho}. Moreover, assume that the function $(t,x)\mapsto I_0(t,x)$ is H\"older continuous, jointly in $(t,x)\in I\times \tilde D$, with exponents $\eta_1, \eta_2\in\, ]0,1]$, respectively. Then the random field solution of \eqref{ch1'-s5.1} satisfies the following:
\smallskip

For any $p\ge 2$, there is a constant $0\le c_p<\infty$ such that, for all $(t,x), (s,y)\in I\times \tilde D$,
\beq
\label{ch1'-ss5.2-t1.1}
\Vert u(t,x)-u(s,y)\Vert_{L^p(\Omega)} \le c_p\left(|t-s|^{\eta_1\wedge\alpha_1\wedge\beta_1}+|x-y|^{\eta_2\wedge\alpha_2\wedge\beta_2}\right).
\eeq
Consequently, there is a version of $(u(t,x),\, (t,x) \in I \times\tilde D)$, that is locally Hölder continuous, jointly in $(t,x)$, with exponents $\gamma_1\in\, ]0,\eta_1\wedge\alpha_1\wedge\beta_1[$, $\gamma_2\in\, ]0,\eta_2\wedge\alpha_2\wedge\beta_2[$, respectively.
If $\tilde D$ is bounded, then this version extends continuously to $\bar I \times D_1$,  where $D_1$ is the closure of $\tilde D$, and is Hölder continuous on $\bar I \times D_1$.

\end{thm}
\begin{proof}
Because of the assumptions on $I_0$, by writing \eqref{ch1'-ss5.2-t1.prelude} for $(t,x),(s,y) \in I\times \tilde D$ with ${\bf\Delta}_1$, ${\bf\Delta}_2$ given in \eqref{choice-rho}, we have
\begin{align*}
\Vert u(t,x)-u(s,y)\Vert_{L^p(\Omega)}
& \le c_p \left[\left(|t-s|^{\eta_1} + |x-y|^{\eta_2}\right)\right.\\
&\left.\qquad + \left(|t-s|^{\alpha_1} + |x-y|^{\alpha_2}\right) + \left(|t-s|^{\beta_1}+ |x-y|^{\beta_2}\right)\right].
\end{align*}
By \eqref{ch1'-s5.2}, the left-hand side of this inequality is bounded, therefore, we obtain \eqref{ch1'-ss5.2-t1.1} (whether or not $\tilde D$ is bounded).

The statements on H\"older continuity are a consequence of Kolmogorov's continuity criterion Theorem \ref{ch1'-s7-t2}.
\end{proof}

\section{Examples of non-linear SPDEs}
\label{ch1'-s6}

In this section, we apply the results of Section \ref{ch4-section1} to selected examples of SPDEs.

\subsection{Stochastic heat equation in spatial dimension $1$}
\label{ch1'-ss6.1}

Let $\mathcal{L}= \frac{\partial}{\partial t} - \frac{\partial^2}{\partial x^2}$ be the heat operator. Extending the results of Sections \ref{ch4-ss2.1-1'} and \ref{ch4-ss2.2-1'}, we can now consider the nonlinear stochastic heat equation on $D=\IR$ or on a bounded interval $D=[0,L]$ and, in the latter case, we consider either homogeneous Dirichlet or Neumann boundary conditions.
Each situation has its own fundamental solution/Green's function $\Gamma$, and the equation is
 \begin{align*}
u(t,x) &= I_0(t,x) + \int_0^t \int_D \Gamma(t,x;s,y) \sigma(s,y,u(s,y))\, W(ds,dy)\notag\\
&\qquad +  \int_0^t \int_D \Gamma(t,x;s,y) b(s,y,u(s,y))\, ds dy,
\end{align*}
  with this $\Gamma$.
First, we start by proving that in each one of these cases, the function $\Gamma$ satisfies the Assumptions ($\bf H_\Gamma$) of Theorem \ref{ch1'-s5.t1}.
\medskip

\noindent{\em Stochastic heat equation on $\re$}
\smallskip

The fundamental solution $\Gamma(t,x;s,y)$ is given by
\beq
\label{ch1'-s6.1}
\Gamma(t,x;s,y):=\Gamma(t-s,x-y)= \frac{1}{\sqrt{4\pi(t-s)}} \exp\left(-\frac{(x-y)^2}{4(t-s)}\right)1_{[0,t[}(s).
\eeq
For $r,x,y \in \IR$, define
\beq
\label{ch1'-s6.1(*1)}
H(r,x,y) = \Gamma(r, x-y) = \frac{1}{\sqrt{4\pi r}} \exp\left(-\frac{(x-y)^2}{4r}\right)1_{]0,\infty[}(r).
\eeq
Then  $0<\Gamma(t,x;s,y)= H(t-s,x,y)$ for $0 \leq s < t$. Clearly, Assumptions (i) and (ii) of ($\bf H_\Gamma$) are satisfied. Moreover, with a change of variable,
\begin{align*}
\int_0^t ds \sup_{x\in\re}\int_{\re} dy\, H^2(s,x,y) & = \int_0^t ds \left[\sup_{x\in\re}\int_{\re} dy \frac{1}{4\pi s} \exp\left(-\frac{(x-y)^2}{2s}\right)\right]\\
& = \int_0^t ds \int_\re dz  \frac{1}{4\pi s}\exp\left(-\frac{z^2}{2s}\right)\\
& = \left(\frac{t}{2\pi}\right)^{\frac{1}{2}},
\end{align*}
(see \eqref{heatcauchy-11'}). This proves Assumption (iiia) of $({\bf H_\Gamma})$, for any $T>0$.

Assumption (iiib) of $({\bf H_\Gamma})$ follows from the fact that, for any $s>0$,  $y\mapsto H(s,x,y)$ is a Gaussian density on $\re$ with mean $x$ and variance $2s$ and therefore,
\beqn
\sup_{x\in\re}\int_{\re} dy\, H(s,x,y)=1.
\eeqn

We notice that since $(x,y)\mapsto H(s,x,y)$ is symmetric, $H(s,x,y)$ also satisfies
\beq
\label{buniclaw}
\sup_{s\in[0,T]}\sup_{y\in \re}\int_\re dx \ H(s,x,y) < \infty,\quad {\text{and}}\quad \sup_{x\in \re} \int_0^T ds \int_\re dy\ H(s,x,y) < \infty.
\eeq
\smallskip

For its further use in Section \ref{ch1'-s7}, we notice that 
\beq
\label{buniclawbis}
\sup_{x\in \re} \int_0^T ds \int_\re dy \ H^\gamma (s,x,y) < \infty,
\eeq
for each $\gamma \in\ ]1,3[$ (see Assumption (iii) of the set of hypotheses ${\bf(h_\Gamma)}$ in Section \ref{ch1'-s7}), by Lemma \ref{app2-l2} with $k=1$.
\medskip

\noindent{\em Stochastic heat equation on a bounded interval with Dirichlet boundary conditions}
\smallskip

 In this example, $D=\ ]0,L[$,
 \beq
 \label{ch4(*1)}
\Gamma(t,x;s,y)= G_L(t-s;x,y),
\eeq
 where $G_L(r;x,y)$ is defined in \eqref{ch1'.600} (with the equivalent expression \eqref{ch1.6000-double}).
By Proposition \ref{ch1'-pPD} (ii), since the Green's function $G_L$ satisfies
\beq
\label{ch1'-s6.1(*2)}
0\le G_L(t;x,y)\le \Gamma(t,x-y),\quad t>0, x, y\in[0,L],
\eeq
 where $\Gamma(t,x-y)$ is the heat kernel defined in
\eqref{ch1'-s6.1}, the computations in the study of the previous example show that $G_L$ satisfies Assumptions ($\bf H_\Gamma$) of Theorem \ref{ch1'-s5.t1}, as well as \eqref{buniclaw} and \eqref{buniclawbis}, for any $T>0$, even with $D$ replaced by $\bar D=[0,L]$. Notice that
with $\Gamma$ defined in \eqref{ch4(*1)}, the right-hand side of \eqref{ch1'-s5.1} vanishes for $x = 0$ and $x = L$, so the vanishing Dirichlet boundary conditions $u(t, 0) = u(t, L) = 0$ are satisfied.
\medskip

\noindent{\em Stochastic heat equation on a bounded interval with Neumann boundary conditions}
\smallskip

In this example, $D=\, ]0,L[$ and
\beq
\label{ch4(*2)}
\Gamma(t,x;s,y) = G_L(t-s;x,y),
\eeq
 with $G_L(r;x,y)$ defined in \eqref{1'.400} (with the equivalent expression  \eqref{1'.15-double}). As in the case of Dirichlet boundary conditions, the Green's function $G_L(t-s;x,y)$ is bounded above by a multiple of a Gaussian density (see \eqref{compneumann}). Therefore the Assumptions ($\bf H_\Gamma$) of Theorem \ref{ch1'-s5.t1}, as well as \eqref{buniclaw} and \eqref{buniclawbis} are satisfied, for any $T>0$, even with $D$ replaced by $\bar D= [0,L]$. The boundary values $u(t,0)$ and $u(t,L)$ are given by the right-hand side of \eqref{ch1'-s5.1} (when $\Gamma$ is defined by \eqref{ch4(*2)}).
\medskip

In the three cases just discussed, the function $I_0$ is given by \eqref{ch1'-v001} (respectively \eqref{cor1.0}, \eqref{cor1.0-N}) for some function $u_0$, that we assume to be bounded so that assumption $({\bf H_I})$ holds.
These considerations along with Remark \ref{4.2.1-remark-bounded-domain}, yield the following.

\begin{thm}
\label{recap}
For the three forms discussed above of the nonlinear stochastic heat equation in spatial dimension $1$ driven by space-time white noise and initial condition $u_0$, under assumption $(\bf H_L)$, the conclusions of Theorem \ref{ch1'-s5.t1} apply, with $\Gamma$ there replaced by the expression \eqref{ch1'-s6.1}, \eqref{ch4(*1)} or \eqref{ch4(*2)}.
That is,
 there exists a random field solution $$u=\left(u(t,x),\, (t,x)\in[0,T]\times (D\cup\partial D)\right)$$ to \eqref{ch1'-s5.1}. This solution satisfies
 \beqn
\sup_{(t,x)\in[0,T]\times (D\cup\partial D)} E\left[\vert u(t,x)\vert^p\right] <\infty,
\eeqn
for any $p>0$,
and the solution $u$ is unique (in the sense of versions) among random field solutions that satisfy this property with $p=2$.
\end{thm}

We address next the question of regularity of the solution to the stochastic heat equation on $D$ with initial condition
\beqn
u(0,x)= u_0(x), \qquad x\in D,
\eeqn
(and boundary conditions if $D=\, ]0,L[$\ ).

Let us begin by considering the function $I_0(t,x)$, which is the solution to the homogeneous PDE $\mathcal{L} I_0 = 0$. Recall that we assume ($\bf H_I$). In the case $D=\IR$, this condition is stronger than \eqref{ch1'-v00}. In particular, when $t=0$, since $I_0(0,x) = u_0(x)$ by definition, it implies that $u_0$ is bounded.

 It is well-known (see e.g. \cite[Theorem 12, Section 5, Chapter 3, p.75]{friedman64}) that $(t,x) \mapsto I_0(t,x)$ is $\cC^\infty$ on $]0,T]\times D$. For $t=0$, $I_0(t,\ast)=u_0(\ast)$ and therefore, $x\mapsto I_0(0,x)$ is continuous in $D$ if and only if $u_0$ is continuous in $D$.

Moreover, we have seen in Chapter \ref{chapter1'} that in the three cases:
\begin{description}
\item{(i)} $D=\IR$ and $u_0\in\mathcal{C}^\eta(\re)$ for some $\eta\in\, ]0,1]$;
\item{(ii)} $D=\, ]0,L[$, with homogeneous Dirichlet boundary conditions,
and $u_0 \in \cC^\eta_0([0,L])$ for some $\eta \in\, ]0, 1]$;
\item{(iii)} $D=\, ]0,L[$, with homogeneous Neumann boundary conditions, and $u_0 \in \cC^\eta([0,L])$ for some $\eta \in\, ]0, 1]$;
\end{description}
\noindent  the function
\beq\label{rde4.2.3}
[0,T]\times D\ni (t,x) \longrightarrow I_0(t,x)=\int_{D} \Gamma(t,x;0,y) u_0(y) dy
\eeq
is H\"older continuous, jointly in $(t,x)$, with exponents $(\frac{\eta}{2},\eta)$ (see \eqref{ch1'-initial 6}, \eqref{cor1.1} and \eqref{cor1.1-N}, respectively).

For the other two terms in the decomposition \eqref{red4.1.14}, the time $t=0$ does not play a special role with regard to sample path regularity. This will be seen in the proofs of the next two statements.

\begin{prop}
\label{reg-si}
The setting and assumptions are the same as in Theorem \ref{recap}. Let
\beqn
\mathcal {I}(t,x)= \int_0^t \int_D \Gamma(t,x;s,y) \sigma(s,y,u(s,y))\, W(ds,dy).
\eeqn
 Then for any $p\geq 2$, there is a constant $C_p$ such that for all $(t,x), (s,y) \in [0,T]\times (D\cup\partial D)$,
\beq
\label{rde4.2.4}
   \Vert \cI(t,x) - \cI(s,y) \Vert_{L^p(\Omega)} \leq C_p \left(\vert t - s \vert^{\frac{1}{4}} + \vert x-y \vert^{\frac{1}{2}} \right).
\eeq
\end{prop}

\begin{proof} We claim that for all $(t,x), (s,y) \in [0,T]\times (D\cup\partial D)$,
\beq\label{rde4.2.5}
   \int_0^Tdr \int_D dz\, (\Gamma(t,x;r,z) - \Gamma(s,y;r,z))^2 \leq C_p^2 \left(|t-s|^{\frac{1}{4}}+|x-y|^{\frac{1}{2}}\right)^2.
\eeq
Indeed, if $D=\IR$, then this follows from \eqref{1'.500} of Lemma \ref{ch1'-l0}. If $D=\,]0,L[$, then under Dirichlet  (resp.~Neumann) boundary conditions, this follows from \eqref{1'1100} in Lemma \ref{ch1'-l2} (resp.~\eqref{1'11000} in Lemma \ref{ch1'-l3}).

Let ${\bf\Delta}_1$ be as in \eqref{choice-rho} with $\alpha_1 = \tfrac{1}{4}$, $\alpha_2 = \tfrac{1}{2}$. Then condition \eqref{ch1'-s5.12} of Lemma \ref{ch1'-ss5.2-l1} is satisfied (with $I=[0,T]$ and $\tilde D = D\cup\partial D$). From \eqref{rde4.1.21} in Proposition \ref{ch1'-ss5.2-prelude}, we deduce that \eqref{rde4.2.4} holds.
\end{proof}

\begin{prop}
\label{reg-li}
The setting and assumptions are the same as in Theorem \ref{recap}. Let
\beqn
\mathcal {J}(t,x)= \int_0^t ds \int_D dy\ \Gamma(t,x;s,y) b(s,y,u(s,y)) .
\eeqn
 Then for any $p\geq 2$, there is a constant $C_p$ such that for all $(t,x), (s,y) \in [0,T]\times (D\cup\partial D)$,
\beq
\label{rde4.2.4-bis}
\Vert \cJ(t,x) - \cJ(s,y) \Vert_{L^p(\Omega)} \leq C_p \left(\vert t - s \vert^{\beta_1} + \vert x-y \vert^{\beta_2} \right).
\eeq
When $D=\re$, the values of the exponents are $\beta_1\in\, ]0,1[$ and $\beta_2=1$. When $D=\, ]0,L[$ with either Dirichlet or Neumann boundary conditions, $\beta_1=\tfrac{1}{2}$ and $\beta_2\in\,]0,1[$.
\end{prop}
\begin{proof}
Consider first the case $D=\re$.
Write $\cJ(t,x) - \cJ(s,y) = \cJ(t,x) - \cJ(t,y) + \cJ(t,y) - \cJ(s,y)$ and use the triangle inequality to split the left-hand side of \eqref{rde4.2.4-bis} into the sum of two terms. To the first term, apply Lemma  \ref{app2-l1} with $k=1$ and $p=1$, observing that $\varphi_p(h) = |h|$ in \eqref{app2.1a}. Then apply to the second term the estimate \eqref{app2.20} in Lemma \ref{app2-l3}.   We deduce that there is $c < \infty$ such that for all $(t,x), (s,y) \in [0,T] \times \IR$ with $0\leq s < t$,
\begin{align*}
&\int_0^T dr \int_{\IR} dz\, |\Gamma(t,x;r,z)-\Gamma(s,y;r,z)|\\
&\qquad = \int_0^T dr \int_{\IR} dz\, |\Gamma(t-r,x-z)-\Gamma(s-r,y-z)|\\
& \qquad \le c\left(|x-y|+|t-s|\log\left(\frac{s}{|t-s|}\right)1_{|t-s|< s}+|t-s|\right).
\end{align*}

Let ${\bf\Delta}_2$ be as in \eqref{choice-rho} with $\beta_1 \in\, ]0,1[$, $\beta_2 = 1$.
 Then condition \eqref{ch1'-s5.14} of Lemma \ref{ch1'-ss5.2-l1} is satisfied with $I = [0, T]$, $\tilde D = \R$, $\Delta_2 = \vert t - s \vert^{\beta_1} + \vert x - y \vert^{\beta_2}$. By \eqref{rde4.1.22} in
Proposition \ref{ch1'-ss5.2-prelude}, we obtain \eqref{rde4.2.4-bis}.

Consider now the case $D=\, ]0,L[$. Since $D$ is bounded, we deduce from Remark \ref{ch1'-ss5.2-r1} and \eqref{rde4.2.5} that the estimate
\eqref{ch1'-s5.14} holds with
\beqn
{\bf\Delta}_2 (t,x;s,y)= {\bf\Delta}_1(t,x;s,y)= |t-s|^{\frac{1}{4}} + |x-y|^{\frac{1}{2}}.
\eeqn
 However, this estimate can be improved.
 Indeed, note that for $h \geq 0$,
 $ \int_t^{t+h} dr \int_0^L dz\, G_L(t+h-r, y, z) \leq h$
by Proposition \ref{ch1'-pPD} (iii) (respectively \eqref{int1}) in the case of Dirichlet (respectively Neumann) boundary conditions. Splitting the left-hand side of \eqref{rde4.2.4-bis} into the sum of two terms as in the case $D = \R$, we appeal to \eqref{nonameD} and \eqref{nonameD-bis} (respectively, Remark \ref{app2-r111}) to find that for all $(t,x), (s,y) \in [0,T] \times [0, L]$,
   \beqn
 \int_0^T ds \int_D dz\, |\Gamma(t,x;r,z) - \Gamma(s,y;r,z)| \le c\, {\bf \Delta}_2(t,x;s,y),
 \eeqn
 where ${\bf \Delta}_2(t,x;s,y) :=  \vert t - s \vert^{\beta_1} + \vert x - y \vert^{\beta_2}$ with $\beta_ 1 = \half$, $\beta_2 \in\, ]0, 1[$. Therefore, condition \eqref{ch1'-s5.14} is satisfied with $I = [0, T]$, $\tilde D = [0, L]$. Using \eqref{rde4.1.22} in Proposition \ref{ch1'-ss5.2-prelude}, we obtain \eqref{rde4.2.4-bis}.
\end{proof}
\begin{thm}
\label{cor-recap}
With the same setting and hypotheses as in Theorem \ref{recap}, the random field solution $(u(t,x),\, (t,x)\in[0,T]\times D)$ to \eqref{ch1'-s5.1} satisfies the following.

(a)
Fix compact intervals $I\subset\ ]0,T]$ and $J\subset (D\cup\partial D)$.
Then for any $p\in[2,\infty[$, there exists a constant $C>0$ (depending on $p$) such that, for any $(t,x),\, (s,y)\in I\times J$,
\beq
\label{4-p1-1'-weak.0D}
E\left[\left\vert u(t,x) - u(s, y)\right\vert^p\right] \le C \left(|t-s|^{\frac{1}{4}} + |x-y|^{\frac{1}{2}}\right)^p.
\eeq
Hence, $(u(t,x),\ (t,x)\in I\times J)$ has a version with jointly Hölder continuous sample paths with exponents $\gamma_1\in\left]0,\frac{1}{4}\right[$ in the time variable $t$, and  $\gamma_2\in\left]0,\frac{1}{2}\right[$ in the spatial variable $x$.

(b) Consider each one of the instances {\rm(i)}, {\rm(ii)} and {\rm(iii)} above relative to the initial condition $u_0$.
 There is $C=C_{p,\eta} < \infty$ such that, for all $(t,x),\, (s,y) \in [0,T]\times (D \cup \partial D)$,
\beq
\label{ch1'-s6.3'}
E\left[\left\vert u(t,x)-u(s,y)\right\vert^p\right] \le C^p \left(|t-s|^{\frac{1}{4}\wedge\frac{\eta}{2}} + |x-y|^{\frac{1}{2}\wedge \eta}\right)^p.
\eeq
Hence, $(u(t,x),\, (t,x)\in[0,T]\times (D \cup \partial D))$ has a version with locally H\"older continuous sample paths.

In the time variable $t$, the H\"older exponent is any
$$
\alpha\in \left]0,\tfrac{1}{4}\right[ \text{ if } \eta\ge \tfrac{1}{2},\qquad
\alpha\in \left]0,\tfrac{\eta}{2}\right]  \text{ if } \eta < \tfrac{1}{2},
$$
while in the space variable $x$, the H\"older exponent is any
$$
   \beta\in \left]0,\tfrac{1}{2}\right[ \text{ if }\eta\ge \tfrac{1}{2},\qquad
   \beta\in \left]0,\eta\right] \text{ if }\eta < \tfrac{1}{2}.
$$
\end{thm}
\begin{proof}
We recall that $(t,x) \mapsto I_0(t,x)$ is $\cC^\infty$ on $]0,T]\times D$ and thus, jointly locally Lipschitz continuous. From \eqref{rde4.2.4} and \eqref{rde4.2.4-bis}, we obtain \eqref{4-p1-1'-weak.0D}. The claim about Hölder continuity follows from
Kolmogorov's continuity criterion Theorem \ref{ch1'-s7-t2}. This completes the proof of (a).

(b)
Recall that by Assumption $(\bf{H_I})$, $I_0$ and $u_0$ are bounded even in the case $D = \R$.
From \eqref{rde4.2.4} and \eqref{rde4.2.4-bis}, we see that for all $(t,x), (s,y) \in [0,T]\times D$,
\eqref{ch1'-ss5.2-t1.prelude} holds with ${\bf\Delta}_1(t,x;s,y) = {\bf\Delta}_2(t,x;s,y) = |t-s|^{\frac{1}{4}} + |x-y|^{\frac{1}{2}}$.
By the H\"older-continuity property of $(t,x) \mapsto I_0(t,x)$ mentioned after \eqref{rde4.2.3},
and since, by \eqref{ch1'-s5.2}, the left-hand side of \eqref{ch1'-s6.3'} is bounded, we have, for all $(t,x), (s,y) \in [0, T] \times (D \cup \partial D)$,
$$
   \Vert u(t,x) - u(s,y) \Vert_{L^p(\Omega)} \leq C_{p,\eta,T} \left(\vert t-s\vert^{\frac{1}{4}\wedge \frac{\eta}{2}} + \vert x-y\vert^{\frac{1}{2}\wedge \eta} \right),
$$
even if $D=\re$. The statement about sample path H\"older continuity follows again from Kolmogorov's continuity criterion Theorem \ref{ch1'-s7-t2}.
\end{proof}

\subsection{Stochastic wave equation in spatial dimension 1}
\label{ch1'-ss6.2}
In Section \ref{ch4-ss2.3-1'}, we studied linear stochastic wave equations. These are defined by the partial differential operator $\mathcal{L}=\frac{\partial^2}{\partial t^2}- \frac{\partial^2}{\partial x^2}$. In this section, we extend the analysis to the non-linear setting
 \begin{align*}
u(t,x) &= I_0(t,x) + \int_0^t \int_D \Gamma(t,x;s,y) \sigma(s,y,u(s,y))\, W(ds,dy)\notag\\
&\qquad +  \int_0^t \int_D \Gamma(t,x;s,y) b(s,y,u(s,y))\, ds dy,
\end{align*}
where $\Gamma$ is the fundamental solution (or the Green's function) associated to the wave operator $\mathcal{L}$.
First, we prove that   $\Gamma$ satisfies the Assumptions (${\bf H_\Gamma}$) of Theorem \ref{ch1'-s5.t1}. For this,
we consider three cases: $D=\IR$, then $D=\, ]0,\infty[$  and $D= \, ]0,L[$ with homogeneous Dirichlet boundary conditions. We state the theorem on existence and uniqueness of random field solutions and, finally, we present the regularity properties of the sample paths.
\medskip

\noindent{\em Stochastic wave equation on $\re$}
\smallskip

 Let $D=\IR$. The fundamental solution is $\Gamma(t,x;s,y)= \Gamma(t-s,x-y)=\frac{1}{2} 1_{\{|x-y|\le t-s\}}$ (see \eqref{gamma-cone}). Define $H(r,x,y)=  \frac{1}{2} 1_{\{|x-y|\le r\}}$. With this choice of $H$, (i) and (ii) of Assumptions (${\bf H_\Gamma}$) are satisfied. Furthermore,
\beqn
\int_0^T ds\ \sup_{x\in\re}\int_{\re} dy\, H^2(s,x,y) = \frac{1}{4} \int_0^T ds\ \sup_{x\in\re}\int_{x-s}^{x+s} dy = \frac{T^2}{4},
\eeqn
 and since $H^2(r,x,y) = \half H(r,x,y)$,
 \beqn
\int_0^T ds\  \sup_{x\in\re}\int_{\re} dy\, H(s,x,y) =
\frac{T^2}{2},
\eeqn
 which shows that the conditions (iiia) and (iiib) of (${\bf H_\Gamma}$) are satisfied as well, for any $T>0$.
 \medskip

 \noindent{\em Stochastic wave equation on $\IR_+$}
 \smallskip

 Let $D=\, ]0,\infty[$. According to \eqref{wavehalfline}, the Green's function is given by
 \beqn
 \Gamma(t,x;s,y):=G(t-s;x,y) = \half 1_{\{|x-(t-s)|\le y\le x+t-s\}}.
 \eeqn
 Let $H(r,x,y)=  \half 1_{\{|x-r|\le y\le x+r\}}$. Clearly, (i) and (ii) of Assumptions (${\bf H_\Gamma}$) hold. Furthermore, since
$\sup_{x\in\IR_+} \int_{\IR_+} 1_{\{|x-r|\le y\le x+r\}}\ dy = 2r$, we see that the conditions (iiia) and (iib) of (${\bf H_\Gamma}$) are satisfied for any $T>0$, even with $D$ replaced by $\bar D$.
\medskip

\noindent{\em Stochastic wave equation on a finite interval}
\medskip

Let $D=\, ]0,L[$. Since we are considering the case of Dirichlet boundary conditions,
we see  from the expression \eqref{wave-bc2-1'}
that $\Gamma(t, x ; s, y)$ is equal to $G_L(t-s; x, y)$, where
\beqn
G_L(r; x,y) = \sum_{m=1}^\infty
 \frac{2}{\pi m} \sin\left(\frac{m\pi x}{L}\right)\sin\left(\frac{m\pi y}{L}\right)\sin\left(\frac{m\pi r}{L}\right).
\eeqn
The series converges for all $r\ge 0$, $x\in[0,L]$ and $y\in[0,L]$, as can be checked using the Fourier series of the expression \eqref{wave-bc1-100'} and the Dirichlet-Jordan convergence test (see e.g. \cite[p. 57]{zygmund1952}).

In addition, the function $H(r, x, y) = \vert G_L(r, x, y)\vert$ satisfies the condition (i) and (ii) of Assumptions (${\bf H_\Gamma}$) even with $D$ replaced by $\bar D$.

Moreover,  we have
 \beqn
\sup_{x\in[0,L]} \int_0^L dy\ H^2(r,x,y) \le \frac{4}{\pi^2} \sum_{m=1}^\infty \frac{1}{m^2} < \infty,
 \eeqn
proving the condition (iiia) for any $T>0$. Since $D$ is bounded, this also implies the validity of  (iiib) for any $T>0$ as well.
\bigskip

In the case $D = \R$ (respectively $D = \R_+$, $D =\, ]0, L[$), the function $I_0$ is given by \eqref{p93.1} (respectively \eqref{p98.2}, \eqref{i0-tris}) for some functions $f$ and $g$. We assume that $f$ is bounded and continuous and that $g \in L^1(D)$, so that assumption $({\bf H_I})$ holds.

For the three choices of $\Gamma$ just discussed, we refer to \eqref{ch1'-s5.1} as the {\em nonlinear stochastic wave equation in spatial dimension $1$ driven by space-time white noise.}  The considerations above yield the following.

\begin{thm}
\label{recapwave}
For the three forms discussed above of the nonlinear stochastic {\em wave} equation in spatial dimension $1$ driven by space-time white noise and initial conditions $f$ and $g$, under assumption
 $({\bf H_L})$, the conclusions of Theorem \ref{ch1'-s5.t1} apply with $D$ there replaced by $\bar D$.
\end{thm}

We now study the sample path regularity of the random field solution $(u(t,x),\, (t,x)\in [0,T]\times \bar D)$ to the nonlinear stochastic wave equation  given in Theorem \ref{recapwave}.

Recall the decomposition \eqref{red4.1.14}:
\beqn
   u(t,x)= I_0(t,x) + \mathcal{I}(t,x) + \mathcal{J}(t,x).
\eeqn
The next proposition discuss the regularity properties of the terms  $\mathcal{I}$ and $\mathcal{J}$.
\begin{prop}
\label{wave-both-integrals}
The setting and hypotheses are those of Theorem \ref{recapwave}. Define
\beq
\label{deltas}
{\bf \Delta}_1 (t,x;s,y)= |t-s|^\half  + |x-y|^\half,\quad {\bf \Delta}_2 (t,x;s,y)= {\bf \Delta}_1 (t,x;s,y)^2,
\eeq
$(t,x),\, (s,y)\in\re_+\times (D\cup\partial D)$.

Then for any $p\ge 2$ there is a constant $C_p$ such that for all $(t,x),\, (s,y)\in[0,T]\times (D\cup\partial D)$,
\beq
\label{rde4.2.4-wave}
\Vert \mathcal{I}(t,x) - \mathcal{I}(s,y)\Vert_{L^p(\Omega)} \le C_p {\bf \Delta}_1 (t,x;s,y),
\eeq
and
\beq
\label{rde4.2.4-wave-bis}
\Vert \mathcal{J}(t,x) - \mathcal{J}(s,y)\Vert_{L^p(\Omega)} \le C_p {\bf \Delta}_2 (t,x;s,y).
\eeq
\end{prop}
\begin{proof}
We will apply Lemma \ref{ch1'-ss5.2-l1} with $Z(r,z) :=\sigma(r,z,u(r,z))$ and $Z(r,z) :=b(r,z,u(r,z))$ to establish
\eqref{rde4.2.4-wave} and \eqref{rde4.2.4-wave-bis}, respectively.

Observe that by condition (${\bf H_L}$)(vi) and Theorem \ref{recapwave}, we have in both instances $\Vert Z\Vert_{T,\infty,p} < \infty$. Hence, it remains to check that for all $(t,x),\, (s,y)\in[0,T]\times (D\cup\partial D)$,
\beq
\label{rde4.2.5-wave}
\int_0^T dr \int_D dz\ (\Gamma(t,x;r,z) - \Gamma(s,y;r,z))^2 \le c {{\bf \Delta}_1} (t,x;s,y)^2,
\eeq
and
\beq
\label{rde4.2.5-wave-bis}
\int_0^T dr \int_D dz\ |\Gamma(t,x;r,z) - \Gamma(s,y;r,z)| \le \bar c {{\bf \Delta}_2} (t,x;s,y),
\eeq
for some constants $c, \bar c>0$.

Indeed, when $D=\re$, apply \eqref {1'.w1-bis} in Lemma \ref{P4:G-1'}; if  $D=\, ]0,\infty[$, use Lemma \ref{rdprop_wave_R_ub}; finally, for $D=\, ]0,L[$, apply \eqref{rdeB.5.5a} of Lemma \ref{app2-5-l-w1}. In this way, we deduce that the left-hand side of \eqref{rde4.2.5-wave} is bounded from above by a constant times $(|t-s| + |x-y|)$ and consequently,
\eqref{rde4.2.5-wave} holds. The conclusion \eqref{ch1'-s5.13} of Lemma \ref{ch1'-ss5.2-l1}(a) yields \eqref{rde4.2.4-wave}.

For the proof of \eqref{rde4.2.5-wave-bis}, we first consider the two cases $D=\re$ and $D=\, ]0,\infty[$, and recall that the fundamental solution and the Green's function are
\beqn
\Gamma(t,x;r,z) = \half 1_{D(t,x)}(r,z), \qquad  \Gamma(t,x;r,z) = \half 1_{E(t,x)}(r,z),
\eeqn
respectively, where $D(t,x)$ is defined in \eqref{cone} and $E(t,x)$ in \eqref{Ewave}. We deduce that
\beqn
|\Gamma(t,x;r,z)-\Gamma(s,y;r,z)| = \half\times
\begin{cases}
1_{D(t,x)\Delta D(s,y)}(r,z),\ {\text {if}}\  & D=\re,\\
1_{E(t,x)\Delta E(s,y)}(r,z),\ {\text {if}}\   & D=\, ]0,\infty[,
\end{cases}
\eeqn
where $A\Delta B$ denotes the symmetric difference of two sets $A$ and $B$. Clearly,
\beqn
|\Gamma(t,x;r,z)-\Gamma(s,y;r,z)|  = 2|\Gamma(t,x;r,z)-\Gamma(s,y;r,z)|^2.
\eeqn
Therefore, using \eqref{rde4.2.5-wave}, we see that \eqref{rde4.2.5-wave-bis} holds with
\beqn
{\bf\Delta}_2(t,x;s,y) = {\bf\Delta}_1(t,x;s,y)^2.
\eeqn
 Then \eqref{rde4.2.4-wave-bis} follows from \eqref{ch1'-s5.15} of Lemma
 \ref{ch1'-ss5.2-l1}(b).

 Finally, consider the case $D=\, ]0,L[$. From Remark \ref{ch1'-ss5.2-r1}, we see that we could take ${\bf\Delta}_2(t,x;s,y)={\bf\Delta}_1(t,x;s,y)=|t-s|^\half + |x-y|^\half$. Nevertheless, using the expression \eqref{wave-bc1-100'} of the Green's function, we see that
 \beqn
    \vert \Gamma(t,x; r,z) - \Gamma(s,y; r,z) \vert \leq 2 (\Gamma(t,x; r,z) - \Gamma(s,y; r,z))^2,
    \eeqn
since the left-hand side takes values in $\left\{0, \half, 1 \right\}$. Arguing as above, we see that \eqref{rde4.2.5-wave-bis} holds with
${\bf \Delta}_2(t,x;s,y) = {\bf \Delta}_1(t,x;s,y)^2$, which is a sharper bound.

The conclusion \eqref{ch1'-s5.14} of Lemma \ref{ch1'-ss5.2-l1}(b) yields \eqref{rde4.2.4-wave-bis} and ends the proof of the Proposition.
\end{proof}

\begin{thm}
\label{cor-recapwave}
The setting and hypotheses are those of Theorem \ref{recapwave}.
Assume also that the initial conditiions $f$ and $g$ satisfy the same assumptions as in Lemma \ref{ch1'-ss2.3-hi}.

Fix compact intervals $I\subset\, ]0,\infty[$ and $J\subset (D\cup \partial D)$.
Then for any $p\in[2,\infty[$, there exists a constant $C=C(p,I,J)>0$  such that, for any $(t,x),\, (s,y)\in I\times J$,
\beq
\label{4-p1-1'-weak.0D-bis}
E\left[\left\vert u(t,x) - u(s, y)\right\vert^p\right] \le C \left(|t-s|^{\frac{1}{2}\wedge \gamma} + |x-y|^{\frac{1}{2}\wedge \gamma}\right)^p.
\eeq
As a consequence, $(u(t,x),\, (t,x) \in I \times J)$ has a version with jointly Hölder continuous sample paths. The H\"older exponents in time and in space coincide.  Denoting by $\alpha$ their common value, we have
\beqn
\alpha\in\, ]0,\tfrac{1}{2}[\ {\text{if}}\quad \gamma\ge \tfrac{1}{2},\qquad \alpha\in\, ]0, \gamma] \ {\text{if}}\quad \gamma < \tfrac{1}{2}.
\eeqn
\end{thm}
\begin{proof}
We check that the hypotheses of Theorem \ref{ch1'-ss5.2-t1} are satisfied.
First, from Lemma \ref{ch1'-ss2.3-hi}, we deduce that the hypotheses on the function $I_0(t,x)$ are satisfied with $\eta_1=\eta_2=\gamma$. The conditions on the stochastic integral $\cI(t,x)$ and the pathwise integral $\cJ(t,x)$ are ensured by Proposition \ref{wave-both-integrals} with $\alpha_i = \half$, $\beta_i=1$, $i=1,2$.

Using the decomposition \eqref{red4.1.14}, the above considerations prove \eqref{4-p1-1'-weak.0D-bis}. The statement about the H\"older continuity of the sample paths follows from Kolmogorov's continuity criterion (Theorem \ref{ch1'-s7-t2}).
\end{proof}


\subsection{Fractional stochastic heat equation}
\label{ch1'-ss7.3}

Let $a\in\, ]0,2]$, $|\delta|\le \min(a,2-a)$. According to \cite[Equation (2.2)]{mlg2001}, we define the Riesz-Feller fractional derivative\index{Riesz-Feller fractional derivative}\index{fractional!derivative, Riesz-Feller}
$\null_xD_\delta^a$ of an integrable function
$f:\re\longrightarrow\re$ by means of its Fourier transform
\beq
\label{ch1'-ss7.3.1}
\tf\left(\null_xD_\delta^a f\right)(\xi) = \null_\delta\psi_a(\xi)\tf f(\xi),\quad \xi\in\re,
\eeq
where
\beqn
\null_\delta\psi_a(\xi)=-|\xi|^a \exp\left(-i\pi\delta\  {\rm sgn}(\xi)/2\right),
\eeqn
and
\beq
\label{defi-fourier}
\tf f(\xi) = \int_{\re} e^{-i\xi x} f(x)\ dx.
\eeq\index{Fourier!transform}\index{transform!Fourier}
Observe that $\null_xD_\delta^a$ defines a pseudo-differential operator with Fourier multiplier $\null_\delta\psi_a(\xi)$.
For $a=2$ (and therefore $\delta=0$), $\null_xD_\delta^a= \frac{d^2}{dx^2}$. If $f$ is $\mathcal{C}^2$ and the second derivative $f^{\prime\prime}$ is integrable, then $\null_xD_\delta^a f$ will be a function, otherwise it may only belong to $\mathcal{S}^\prime(\re)$.

In this section, we consider the SPDE \eqref{ch1'-s5.0} with the partial differential operator $\mathcal{L} = \frac{\partial}{\partial t} - \null_xD_\delta^a$, with $a\in\, ]1,2[$, $|\delta|\le 2-a$ and  $D=\re$. Note that for  $a\in\,]0,1]$, there is no random field solution to \eqref{ch1'-s5.0} (see \cite[p. 361]{chendalang2015}).

From \eqref{ch1'-ss7.3.1} and \eqref{rd02_08e1}, we see that if $\delta=0$ and $a>0$, then the operator 
$\null_xD_\delta^a$ is the opposite of the fractional Laplacian $(-\Delta)^{a/2}$. Hence, 
comparing with Section \ref{ch3-sec3.5}, $\delta$ need not be equal to $0$ here. However $k=1$ and $a\in]1,2[$, while there, $k\ge 1$ and $a>k$. In addition, in this section, we will consider below the nonlinear SPDE \eqref{ch1'-ss7.3.9}.

For $a\in\, ]1,2[$, the Riesz-Feller fractional derivative can be also defined by
\beq
\label{ch1'-ss7.3.2}
\null_x D_\delta^a f(x)= \int_{\re} [f(x+z)-f(x)-zf^\prime(x)]\ \nu_a(dz),
\eeq
(for functions $f$ for which the integral is well-defined), where $\nu_a$ is the mesure
\beq
\label{ch1'-ss7.3.levy}
\nu_a(dz)= c^+_a \frac{dz}{z^{1+a}}\ 1_{\{z>0\}} + c^-_a \frac{dz}{(-z)^{1+a}}\ 1_{\{z<0\}},
\eeq
with
\beqn
c_a^{\pm}= \frac{\Gamma_E(1+a)}{\pi}\sin\left((a\pm \delta)\frac{\pi}{2}\right)
\eeqn
(here, $\Gamma_E$ is the Euler Gamma function: see \eqref{Euler-gamma}).
A proof of this is given in the next lemma.
\begin{lemma}
\label{ch1'-ss7.3-l0}
Let $a\in\, ]1,2[$ and $\vert \delta\vert\le 2-a$.
For any $\mathcal{C}^2$ function
$f:\re\longrightarrow \re$ with compact support, the following formulas hold.
\begin{enumerate}
\item For $\xi\in\re$,
\begin{align}
 \label{ch1'-ss7.3.1prebis}
 &\tf \left(\int_{\re} [f(\ast+z)-f(\ast)-zf^\prime(\ast)]\ \nu_a(dz)\right)(\xi)\notag\\
 &\qquad \qquad \qquad = \tf f(\xi) \int_{\re}\left(e^{i\xi z}-1-i\xi z \right)\ \nu_a(dz) .
 \end{align}
\item For $\xi\in\re$,
  \beq
 \label{ch1'-ss7.3.1bis}
 \int_{\re} \left(e^{i\xi z}-1-i\xi z \right)\ \nu_a(dz) = \null_\delta\psi_a(\xi).
 \eeq
 \end{enumerate}
 Consequently, for $x\in\re$,
 \beqn
\null_x D_\delta^a f(x)= \int_{\re} [f(x+z)-f(x)-zf^\prime(x)]\ \nu_a(dz),
\eeqn
where $\null_x D_\delta^a$ is defined in \eqref{ch1'-ss7.3.1}.
\end{lemma}
\begin{proof}
Since $f$ is $\mathcal{C}^2$ with compact support, the Fourier transform can be calculated using
the duality $\langle \tf \varphi , \psi \rangle_{L^2(\R)} = \langle \varphi , \tf \psi \rangle_{L^2(\R)}$ and Fubini's theorem, and then
the identity \eqref{ch1'-ss7.3.1prebis} follows from the elementary properties $\tf f(\ast + z)(\xi) = e^{i \xi z} \tf f(\xi)$ and $\tf f'(\ast)(\xi) = i \xi \F f(\xi)$ of the Fourier transform.

For the proof of \eqref{ch1'-ss7.3.1bis}, we will make use of the following properties of the Euler Gamma function (see \eqref{Euler-gamma}):
\beq
\label{prodgam}
\Gamma_E(a)\Gamma_E(1-a) = \frac{\pi}{\sin{\pi a}}, \ a\notin\mathbb{Z},
\eeq
 (see \cite[p.7]{podlubny1999}) and
\beq
\label{f1}
\int_0^\infty \frac{e^{-qz}-1+qz}{z^{1+a}}\, dz = q^a \Gamma_E(-a), \ q\in\mathbb{C},\  {\text{Re}}(q)\ge 0,\ 1< a< 2,
\eeq
that is proved in Lemma \ref{A3-l0}.
Define
\begin{align*}
I_1 &= \int_0^\infty  \left(e^{i\xi z}-1-i\xi z \right)\, \frac{dz}{z^{1+a}},\\
I_2 &= \int_{-\infty}^0  \left(e^{i\xi z}-1-i\xi z \right)\, \frac{dz}{(-z)^{1+a}}.
\end{align*}
From \eqref{f1}, we have
\beqn
I_1 = (-i\xi)^a \Gamma_E(-a), \ \, I_2 = (i\xi)^a \Gamma_E(-a).
\eeqn
Use the polar decompositions
\beqn
-i\xi = |\xi| e^{- i\frac{\pi}{2}{\rm{sgn}}(\xi)},\quad
i\xi = |\xi| e^{i\frac{\pi}{2}{\rm{sgn}}(\xi)},
\eeqn
to see that
\beqn
(-i\xi)^a = |\xi|^a e^{-i\frac{a\pi}{2}{\rm{sgn}}(\xi)},\quad
(i\xi)^a = |\xi|^a e^{i\frac{a\pi}{2}{\rm{sgn}}(\xi)}.
\eeqn
Thus,
\begin{align*}
 &\int_{\re} \left(e^{i\xi z}-1-i\xi z \right)\, \nu_a(dz) = c_a^+I_1+c_a^-I_2\\
 &\qquad=\frac{\Gamma_E(1+a)\Gamma_E(-a)}{\pi}|\xi|^a \\
 &\qquad\qquad\times \left\{\sin\frac{(a+\delta)\pi}{2}\left[\cos\left(\frac{a\pi}{2}{\rm{sgn}}(\xi)\right)-i\sin\left(\frac{a\pi}{2}{\rm{sgn}}(\xi)\right)\right]\right.\\
 &\left.\qquad\quad + \sin\frac{(a-\delta)\pi}{2}\left[\cos\left(\frac{a\pi}{2}{\rm{sgn}}(\xi)\right)+i\sin\left(\frac{a\pi}{2}{\rm{sgn}}(\xi)\right)\right]\right\}.
  \end{align*}
  Since $1+a\notin\mathbb{Z}$, from \eqref{prodgam} and the formula $\sin(x+\pi)=-\sin(x)$, we have
  \beqn
  \frac{\Gamma_E(1+a)\Gamma_E(-a)}{\pi}=\frac{1}{\sin(\pi(1+a))} = -\frac{1}{\sin(\pi a)}.
  \eeqn
  Hence, by rearranging terms, we obtain
  \begin{align*}
  &\int_{\re} \left(e^{i\xi z}-1-i\xi z \right)\, \nu_a(dz) = -\frac{1}{\sin(\pi a)}|\xi|^a \\
  &\quad \times\left\{\cos\left(\frac{a\pi}{2}{\rm{sgn}}(\xi)\right)\left[\sin\frac{(a+\delta)\pi}{2}+\sin\frac{(a-\delta)\pi}{2}\right]\right.\\
  &\left.\qquad + i \sin\left(\frac{a\pi}{2}{\rm{sgn}}(\xi)\right)\left[ \sin\frac{(a-\delta)\pi}{2}-\sin\frac{(a+\delta)\pi}{2}\right]\right\}\\
  & \qquad = -\frac{1}{\sin(\pi a)}|\xi|^a\left\{2\cos\left(\frac{a\pi}{2}{\rm{sgn}}(\xi)\right)\sin\left(\frac{a\pi}{2}\right)\cos\left(\frac{\delta\pi}{2}\right)\right.\\
  &\left.\quad\qquad \qquad \qquad \qquad \quad-2i \sin\left(\frac{a\pi}{2}{\rm{sgn}}(\xi)\right)\cos\left(\frac{a\pi}{2}\right)\sin\left(\frac{\delta\pi}{2}\right)\right\},
  \end{align*}
  where in the second equality we have applied the formula
  \beqn
  \sin x \pm \sin y = 2\sin\left(\frac{x\pm y}{2}\right) \cos\left(\frac{x\mp y}{2}\right).
  \eeqn
  Observe that, by the formula $\sin(2x)= 2\sin x\cos(\pm x)$,
  \begin{align*}
 2 \cos\left(\frac{a\pi}{2}\rm{sgn}(\xi)\right)\sin\left(\frac{a\pi}{2}\right)\cos\left(\frac{\delta\pi}{2}\right)
  = \sin(a\pi)\cos\left(\frac{\delta\pi}{2}\right).
  \end{align*}
  Analogously,
   \begin{align*}
   &2\sin\left(\frac{a\pi}{2}{\rm{sgn}}(\xi)\right)\cos\left(\frac{a\pi}{2}\right)\sin\left(\frac{\delta\pi}{2}\right)\\
   &\qquad\qquad=
   2 \sin\left(\frac{a\pi}{2}\right)
    \cos\left(\frac{a\pi}{2}\right)\sin\left(\frac{\delta\pi}{2}{\rm{sgn}}(\xi)\right)\\
   &\qquad\qquad= \sin\left( a\pi\right)\sin\left(\frac{\delta\pi}{2}{\rm{sgn}}(\xi)\right).
   \end{align*}
   This yields
   \beqn
   \int_{\re} \left(e^{i\xi z}-1-i\xi z \right)\, \nu_a(dz) = -|\xi|^a \left[\cos\left(\frac{\delta\pi}{2}{\rm{sgn}}(\xi)\right)-i\sin\left(\frac{\delta\pi}{2}{\rm{sgn}}(\xi)\right)\right],
   \eeqn
   which is \eqref{ch1'-ss7.3.1bis}.

   From \eqref{ch1'-ss7.3.1prebis}, \eqref{ch1'-ss7.3.1bis}, we see that  $\int_{\re} [f(\ast+z)-f(\ast)-zf^\prime(\ast)]\, \nu_a(dz)$ and
   $\null_x D_\delta^a f(\ast)$ have the same Fourier transform and therefore these two continuous functions are equal.
      \end{proof}

There is a useful probabilistic interpretation of the operator $\null_xD_\delta^a$.
Indeed, let $X=(X_t,\, t\ge 0)$ be a strictly $a$-stable L\'evy process (where ``strictly'' refers to the fact that the process is centered: see \cite[Chapter 3]{sato2013}) with L\'evy measure $\nu_a$ given in \eqref{ch1'-ss7.3.levy}.
Following \cite[Theorem 31.5, p. 208]{sato2013}, we deduce from \eqref{ch1'-ss7.3.2} that $\null_xD_\delta^a$ is the infinitesimal generator of $X$.

By the L\'evy-Khintchine formula, the characteristic function of the law of $X_t$ is
\beqn
\exp\left[-t\int_{\re}\left(e^{i\xi z}-1-i\xi z\right) \nu_a(dz)\right].
\eeqn
where no ``truncation function'' nor additional drift term appears because of the centering.
By \eqref{ch1'-ss7.3.1bis} this is equal to $\exp\left(-t\null_\delta\psi_a(\xi)\right)$.

Denote by $\null_\delta G_a(t,x)$ the fundamental solution of the operator $\mathcal{L} = \frac{\partial}{\partial t} - \null_xD_\delta^a$.
By the classical approach for PDEs on $\rek$ (that is, using the Fourier transform methods described in Section \ref{AppC-New-0}), we find that
\beq
\label{fs-fractional}
\null_\delta G_a(t,x) = \tf^{-1}\left[\exp(t\ \null_\delta\psi_a(\ast))\right](x)\,1_{]0,\infty[}(t),
\eeq
thus, by the definition of $\null_\delta\psi_a$, for $t>0$,
\beq
\label{ch1'-ss7.3.3}
\null_\delta G_a(t,x) = \frac{1}{2\pi}\int_{\re} d\xi\ \exp\left[i\xi x-t|\xi|^a
e^{-i\pi\delta\  {\rm sgn}(\xi)/2}\right].
\eeq
This provides a formula for the density of the centered $a$-stable random variable $X_t$ with L\'evy measure $t\nu_a$, as in \cite[p.17]{zolotarev1986}.
In particular, $\null_\delta G_a(t,x)$ is real-valued and nonnegative.
From the above considerations, we obviously have
\begin{align}
&\int_\re\null_\delta G_a(t,x)\ dx = 1, \ {\text{for all}}\  t>0,\label{ch1'-ss7.3.4}\\
&\null_\delta G_a(s+t,\ast) = \null_\delta G_a(s,\ast)\ast \null_\delta G_a(t,\ast),\ {\text{for all}}\  s,t>0.\label{ch1'-ss7.3.5}
\end{align}
Moreover, from \eqref{ch1'-ss7.3.3}, with the change of variable $\bar\xi = t^{\frac{1}{a}}\xi$, we obtain
\beq
\label{ch1'-ss7.3.6}
\null_\delta G_a(t,x) = t^{-\frac{1}{a}}\null_\delta G_a(1, t^{-\frac{1}{a}}\, x), \ {\text{for all}}\quad  t>0,\ x \in \R.
\eeq
\medskip

\begin{lemma}
\label{check-HGamma}
For any $x,y\in\re$ and $0\le s<t$, let
\beq
\label{choicegamma}
\Gamma(t,x;s,y):=\laplacef(t-s,x-y).
\eeq
Then for any $T>0$, $\Gamma$ satisfies the assumptions $({\bf H_\Gamma})$ of Section \ref{ch4-section0}.
\end{lemma}
\begin{proof}
Condition (i) of $({\bf H_\Gamma})$ is immediate, and
 clearly, condition (ii) on $\Gamma$ holds with
 $H(s,x,y)= \null_\delta G_a(s,x-y)$.

 Next, we prove that the assumptions (iiia) and (iiib) are also satisfied (with $D=\re$). Indeed, using \eqref{ch1'-ss7.3.5},
since $a>1$,
 \begin{align*}
\int_0^T ds \sup_{x\in\re} \int_\re dy\,  \left\vert\null_\delta G_a(s,x-y)\right\vert^2
&=\int_0^T ds\, _\delta G_a(2s, 0)\\
&= \int_0^T ds\,  (2s)^{-\frac{1}{a}}  \null_\delta G_a(1, 0)\\
&= \frac{a}{a-1}\, 2^{-\frac{1}{a}}\,  T^{\frac{a-1}{a}}\, \null_\delta G_a(1, 0) < \infty,
\end{align*}
since $a>1$. This proves (iiia).

Further,
\beqn
\int_0^T ds\, \sup_{x\in\re} \int_\re dy\  \null_\delta G_a(s,x-y) = \int_0^T ds \int_\re dz \  \null_\delta G_a(s,z),
\eeqn
since the value of the $dy$-integral does not depend on $x$. Thus, using \eqref{ch1'-ss7.3.4},
we deduce  (iiib).
\end{proof}

Using Lemma \ref{check-HGamma} and Theorem \ref{ch1'-s5.t1}, we obtain immediately the following result on existence and uniqueness of solutions for a {\em non-linear fractional stochastic heat equation.}\index{fractional!heat equation}\index{heat!equation, fractional}\index{equation!fractional heat}


\begin{thm}
\label{ch1'-ss7.3-t1}
Consider the SPDE
\beq
\label{ch1'-ss7.3.9}
\left(\frac{\partial}{\partial t} - \null_xD_\delta^a\right)u(t,x)
= \sigma(t,x,u(t,x)) \dot W(t,s) + b(t,x,u(t,x)),
\eeq
$(t,x)\in\ ]0,T]\times \re$, with $u(0,x)= u_0(x)$, where $a\in\, ]1,2]$, $|\delta|\le 2-a$.
Assume that $I_0(t,x)=\int_\re dy\ \null_\delta G_a(t,x-y)u_0(y)$ satisfies the assumption $({\bf H_I})$ of Section \ref{ch4-section1}, and that the functions $\sigma$ and $b$ satisfy the assumptions $({\bf H_L})$. Then
there exists a jointly measurable and adapted process $\left(u(t,x),\, (t,x)\in[0,T]\times \re\right)$ such that for all $(t,x)\in[0,T]\times \re$,
\begin{align}
\label{ch1'-ss7.3.9bis}
u(t,x) & = I_0(t,x) + \int_0^t \int_{\re} \laplacef(t-s,x-y)\sigma(s,y,u(s,y))\, W(ds,dy)\notag\\
& \qquad +  \int_0^t ds \int_{\re} dy\ \laplacef(t-s,x-y) b(s,y,u(s,y)),\ {\text{a.s.}}
\end{align}
In addition, for any $p>0$,
\beq
\label{ch1'-ss7.3.10}
\sup_{(t,x)\in[0,T]\times \re} E\left[\vert u(t,x)\vert^p\right] <\infty,
\eeq
and the solution $u$ is unique (in the sense of versions) among random fields that satisfy \eqref{ch1'-ss7.3.10} with $p=2$.
\end{thm}
For coefficients $\sigma$, $b$ not depending on $\omega$, a similar result has been obtained in \cite[Theorem 1]{dd2005} (see also \cite[Theorem 3.1] {chendalang2015} for properties of the solution). Observe that because of \eqref{ch1'-ss7.3.4}, the assumption ($\bf H_I$) on $I_0$ holds if the initial condition $u_0$ is a bounded Borel function.
\medskip

 \noindent{\em Regularity of the sample paths}
 \medskip

 Next, we would like to study the regularity of the sample paths of the solution to the SPDE  \eqref{ch1'-ss7.3.9bis}. Instead of applying
 Theorem \ref{ch1'-ss5.2-t1}, as we did for example in Section \ref{ch1'-ss6.1} for the stochastic heat equation on $\re$, we will use a slightly different method.
 This is to avoid determining for which ${\bf\Delta}_2$ the fundamental solution $\Gamma(t,x;s,y):=\laplacef(t-s,x-y)$ of  \eqref{ch1'-ss7.3.9} satisfies the estimate \eqref{ch1'-s5.14}.

 By proceeding in this way, we introduce an alternative approach to the study of the H\"older continuity of sample paths of an SPDE, which in principle can be applied in instances where its fundamental solution $\Gamma$ is homogeneous, that is, $\Gamma(t,x;s,y)= \tilde \Gamma(t-s,x-y)$, for some $\tilde\Gamma$.

\begin{thm}
\label{ch1'-ss7.3-t2}
The assumptions are as in Theorem \ref{ch1'-ss7.3-t1}. In addition to the hypothesis $({\bf H_L})$-(v) relative to the coefficient $b$, assume  that there are $\delta_1, \delta_2\in\, ]0,1]$ and $0<c<\infty$  such that, for all $\omega\in\Omega$, $s,t\in[0,T]$, $x,y,z\in \re$,
\beq
\label{lipb}
|b(t,x,z;\omega) - b(s,y,z;\omega)| \le c(1+|z|)\left(|t-s|^{\delta_1} + |x-y|^{\delta_2}\right).
\eeq
 Suppose also that the function $(t,x)\mapsto I_0(t,x)$ is H\"older continuous, jointly in $(t,x)$ with exponents $\eta_1,\eta_2\in\, ]0,1]$, respectively. Then the random field solution to  \eqref{ch1'-ss7.3.9bis} satisfies the following:

Fix $M>0$. Then for any any $p\ge 2$, there exists a finite and positive constant $C=C_{M,p,T}$ such that, for all $s,t\in[0,T]$ and all $x, y\in [-M, M]$,
\beq
\label{ch1'-ss7.3.100}
\Vert u(t,x)-u(s,y)\Vert_{L^p(\Omega)}  \le C\left(|t-s|^{\left(\frac{a-1}{2a}\right)\wedge \delta_1\wedge \eta_1} + |x-y|^{\left(\frac{a-1}{2}\right)\wedge\delta_2\wedge\eta_2}\right).
\eeq
Hence, $(u(t,x),\, (t,x)\in[0,T]\times \re)$ has a version with locally H\"older continuous sample paths, jointly in $(t,x)$, with exponents $\gamma_1\in\, ]0,\left(\frac{a-1}{2a}\right)\wedge \delta_1\wedge \eta_1[$, $\gamma_2\in\, ]0,\left(\frac{a-1}{2}\right)\wedge\delta_2\wedge\eta_2[$, respectively.
\end{thm}

\begin{proof}
For any $(t,x)\in[0,T]\times \re$, set
\beqn
Z(t,x) = g(t,x,u(t,x)),
\eeqn
where $g$ stands for either $\sigma$ or $b$. Since $\sigma$ and $b$ have at most linear growth (see hypothesis $({\bf H_L})$-(vi)), the property \eqref{ch1'-ss7.3.10}
yields
\beq
\label{ch1'-ss7.3.101}
\sup_{(t,x)\in[0,T]\times \re}\Vert Z(t,x)\Vert_{L^p(\Omega)} < \infty.
\eeq
Fix $0\le s\le t\le T$, $x,y\in\re$ and set $v(t,x):=u(t,x)-I_0(t,x)$. Then
\beq
\label{tipa(*1)}
E\left[|v(t,x) - v(s,y)|^p\right] \le 2^{p-1}\left(T_1(t,x;s,y) + T_2(t,x;s,y)\right),
\eeq
with
\begin{align*}
T_1&(t,x;s,y) \\
&= E\Big[\Big(\int_0^T \int_\re \left[\laplacef(t-r,x-z)-\laplacef(s-r,y-z)\right]\\
&\qquad\qquad\qquad\times
\sigma(r,z,u(r,z))\ W(dr,dz)\Big)^p\Big],\\
 T_2&(t,x;s,y)\\
 &=E\Big[\Big( \int_0^T dr \int_\re dz \left[\laplacef(t-r,x-z)-\laplacef(s-r,y-z)\right]\\
&\qquad\qquad\qquad\times b(r,z,u(r,z))\Big)^p\Big].
 \end{align*}
From Proposition \ref{ch1'-ss7.3-p1}  in Appendix \ref{app2}, we see that $\Gamma(t,x;r,z):=\laplacef(t-s,x-z)$ satisfies the condition \eqref{ch1'-s5.12} of Lemma \ref{ch1'-ss5.2-l1} with $\Delta_1(t,x;r,z) = |t-s|^{\frac{a-1}{2a}} + |x-z|^{\frac{a-1}{2}}$.
Thus, from the conclusion (a) of this Lemma, along with \eqref{ch1'-ss7.3.101}, we see that
\beq
\label{ch1'-ss7.3.102}
T_1(t,x;s,y) \le C \left( |t-s|^{\frac{a-1}{2a}} + |x-y|^{\frac{a-1}{2}}\right)^p.
\eeq
Separating into three $dr dz$-integrals, changing variables in time and space and regrouping, we have
\begin{align*}
&\int_0^T dr \int_\re dz \left[\laplacef(t-r,x-z)-\laplacef(s-r,y-z)\right] b(r,z,u(r,z))\\
&\quad = \int_s^t dr \int_\re dz\ \laplacef(r,z) b(t-r,x-z,u(t-r,x-z))\\
&\qquad + \int_0^s dr \int_\re dz\ \laplacef(r,z) [b(t-r,x-z,u(t-r,x-z)) \\
&\qquad \qquad\qquad - b(s-r,y-z,u(s-r,y-z))].
\end{align*}
Apply H\"older's inequality (alternatively, Minkowski's inequality) and use \eqref{ch1'-ss7.3.101}, \eqref{ch1'-ss7.3.4}, to obtain
\beqn
E\left[\left\vert \int_s^t dr \int_\re dz\ \laplacef(r,z) b(t-r,x-z,u(t-r,x-z))\right\vert^p\right] \le C(t-s)^p.
\eeqn
As for the second term in the array above, we first apply H\"older's inequality with respect to the measure on $[0,s]\times \re$ with density $\laplacef$ to obtain
\begin{align}
\label{intermed}
&E\Big[\Big\vert \int_0^s dr \int_\re dz\ \laplacef(r,z)\notag\\
&\qquad\times [b(t-r,x-z,u(t-r,x-z)) - b(s-r,y-z,u(s-r,y-z))]\Big\vert^p\Big] \notag\\
&\quad \le \left( \int_0^s dr \int_\re dz\ \laplacef(r,z)\right)^{p-1}\int_0^s dr \int_\re dz\ \laplacef(r,z)\notag\\
&\qquad \times E\left[\left\vert b(t-r,x-z,u(t-r,x-z)) - b(s-r,y-z,u(s-r,y-z))\right\vert^p\right].
\end{align}
By the triangle inequality,
\begin{align*}
&\left\vert b(t-r,x-z,u(t-r,x-z)) - b(s-r,y-z,u(s-r,y-z))\right\vert\\
&\qquad\le \left\vert b(t-r,x-z,u(t-r,x-z)) - b(t-r,x-z,u(s-r,y-z))\right\vert\\
&\qquad \quad  +  \left\vert b(t-r,x-z,u(s-r,y-z))- b(s-r,y-z,u(s-r,y-z))\right\vert\\
&\qquad \le C\vert u(t-r, x-z) - u(s-r, y-z)\vert\\
&\qquad \quad  +  c(1+|u(s-r,y-z)|)\left(|t-s|^{\delta_1}+|x-y|^{\delta_2}\right),
\end{align*}
where, in the last inequality, we have used the hypotheses $({\bf H_L})$(v) and \eqref{lipb}.
Applying a change in the $r$ and $z$ variables, we see that the last term in the array \eqref{intermed} is bounded from above by
\begin{align*}
&\le C_pT^{p-1} \int_0^s dr \int_\re dz\ \laplacef(r,z)\\
&\qquad \times  E\Big[\left\vert u(t-r,x-z)-u(s-r,y-z)\right\vert^p\\
&\qquad \qquad  + (1+|u(s-r,y-z)|^p)\left(|t-s|^{\delta_1}+|x-y|^{\delta_2}\right)^p\Big] \\
&\le C_pT^{p-1} \int_0^s dr \int_\re dz\ \laplacef(s-r,y-z) \\
&\qquad \times E\left[\left\vert u(t-s+r,x-y+z)-u(r,z)\right\vert^p +  \left(|t-s|^{\delta_1}+|x-y|^{\delta_2}\right)^p\right] \\
&\le C\, \Big(\left(|t-s|^{\delta_1}+|x-y|^{\delta_2}\right)^p \\
&\qquad \qquad+  \int_0^s dr \ \sup_{z\in \re} E\left[\left\vert u(t-s+r,x-y+z)-u(r,z)\right\vert^p\right]\Big),
\end{align*}
where, in the second inequality, we have used \eqref{ch1'-ss7.3.10}.

Set $h:=t-s$, $\bar h:=x-y$. Notice that  $0 \leq h \leq T$ and $\vert \bar h \vert \leq 2M$. By \eqref{prodgam} and (${\bf H_I}$), the left-hand side of \eqref{tipa(*1)} is bounded, therefore the bounds on $T_1$ and $T_2$ obtained above show that
\begin{align*}
 &E\left[\left\vert v(h+s,\bar h+y)-v(s,y)\right\vert^p\right] \le C\, \Big(\left(h^{\frac{a-1}{2a}\wedge \delta_1}+|\bar h|^{\frac{a-1}{2}\wedge \delta_2}\right)^p \\
 & \qquad +\int_0^s dr \ \sup_{z\in \re} E\left[\left\vert u(h+r,\bar h+z)-u(r,z)\right\vert^p\right]\Big).
 \end{align*}
 Consequently, by the H\"older continuity property of $I_0$, this implies
 \begin{align*}
 &E\left[\left\vert u(h+s,\bar h+y)-u(s,y)\right\vert^p\right] \le C\, \Big(\left(h^{\eta_1\wedge\delta_1\wedge\left(\frac{a-1}{2a}\right)}+|\bar h|^{\eta_2\wedge\delta_1\wedge\left(\frac{a-1}{2}\right)}\right)^p \\
 &\qquad +\int_0^s dr \ \sup_{z\in \re} E\left[\left\vert u(h+r,\bar h+z)-u(r,z)\right\vert^p\right]\Big).
 \end{align*}

 Apply Gronwall's Lemma \ref{lemC.1.1}, more precisely \eqref{rdeC.1.5}, to the real-valued function
 \beqn
 s\mapsto \sup_{y\in \re}E\left[\left\vert u(h+s,\bar h+y)-u(s,y)\right\vert^p\right]
 \eeqn
 to conclude that for $h\in [0,T]$, $\bar h\in\re$ and $s\in[0,T-h]$,
 \begin{align}
 \label{final}
 &\sup_{y\in\re}E\left[\left\vert u(h+s,\bar h+y)-u(s,y)\right\vert^p\right] \notag\\
 &\qquad\qquad \qquad\qquad\le C  \left(h^{\eta_1\wedge\delta_1\wedge\left(\frac{a-1}{2a}\right)}+|\bar h|^{\eta_2\wedge\delta_2\wedge\left(\frac{a-1}{2}\right)}\right)^p.
 \end{align}
 This implies \eqref{ch1'-ss7.3.100}. The last claim is a consequence of Kolmogorov's continuity criterion Theorem \ref{ch1'-s7-t2}.
 \end{proof}
  \bigskip

\section{Approximation by finite-dimensional projections}
\label{ch1'-s60}

According to the theory of stochastic integration developed in Chapter \ref{ch2}, the stochastic integral in the formulation \eqref{ch1'-s5.1} of the SPDE \eqref{ch1'-s5.0} is
\begin{align}
\label{ch1'-s60.1}
&\int_0^t \int_D \Gamma(t,x;s,y) \sigma(s,y,u(s,y))\, W(ds,dy)\notag\\
&\qquad = \sum_{j=1}^\infty \int_0^t \langle \Gamma(t,x; s,\ast) \sigma(s,\ast,u(s,\ast)), e_j\rangle_V\, dW_s(e_j),
\end{align}
where $V=L^2(D)$, $(e_j,\, j\ge 1)$ is a CONS of $V$ and $(W_s(e_j),\, s\in[0,T],\, j\ge 1)$ is the sequence of independent standard Brownian motions given in Lemma \ref{ch1'-lsi}.
Recall that in this context, the symbol ``$\ast$'' refers to the spatial variable in $D$.

In \eqref{ch1'-s60.1}, the integrator is the cylindrical Wiener process given in Lemma \ref{ch1'-lc} and, according to Lemma \ref{ch1'-lsi} (1), this process admits the representation
\beqn
W_s(\varphi) = \sum_{j=1}^\infty \langle \varphi, e_j\rangle_V W_s(e_j), \quad  s\in[0,T], \quad  \varphi\in V.
\eeqn

For any $n \geq 1$, let $V_n$ be the subspace of $V$ spanned by $(e_j,\, 1 \leq j \leq n)$ and let $\Pi_{V_n}$\label{rdorthproj} denote the orthogonal projection from $V$ onto $V_n$. Define
\beq
\label{proj-noise}
W_s^n(\varphi) = \sum_{j=1}^n \langle \varphi, e_j\rangle_V W_s(e_j),  \quad  s\in[0,T], \quad  \varphi\in V_n.
\eeq
 Applying Lemma \ref{ch1'-lsi} (2) to the finite sequence $(W_s(e_j),\, s \in [0, T],\, 1 \leq j \leq n)$, this defines an isonormal Gaussian process  $W^n = (W^n(h),\, h \in L^2([0, T], V_n)) $. It is a ``finite-dimensional projection of the noise''\index{finite-dimensional!projections}\index{projections!finite-dimensional} $(W_s(\varphi))$.

  Given a process $G = (G(s,y),\, (s,y) \in [0, T] \times D)$ as in Section \ref{ch2new-s2}, we define a stochastic integral with respect to $W^n$ by
  \begin{align*}
    \int_0^t \int_D  G(s,y)  W^n(ds, dy) :&= \int_0^t \int_D \Pi_{V_n}(G(s,*))(y)\,  W(ds, dy)\\
    &= \sum_{j=1}^n \int_0^t \langle G(s, *), e_j \rangle_V\,  dW_s(e_j).
    \end{align*}

Consider the equation
\begin{align}
\label{ch1'-s60.2}
\bar u_n(t,x) &= I_0(t,x) + \int_0^t \int_D \Gamma(t,x;s,y) \sigma(s,y,\bar u_n(s,y))\, W^n(ds,dy)\notag\\
&\qquad +  \int_0^t \int_D \Gamma(t,x;s,y) b(s,y,\bar u_n(s,y))\, ds dy,
\end{align}
$(t,x)\in[0,T]\times D$, $n\ge 1$.

The purpose of this section is to establish a convergence result of the sequence $(\bar u_n(t,x),\, n\ge 1)$ to $(u(t,x))$ in a sense made precise in Theorem \ref{ch1'-s60-t2} below.  We call the sequence $(\bar u_n,\, n \geq 1)$ the {\em approximation of $u$ by finite-dimensional projections.}
Throughout the section, we assume the hypotheses $(\bf H_\Gamma)$, $(\bf H_I)$ and  $(\bf H_L)$ (introduced in Sections \ref{ch4-section0} and \ref{ch4-section1}).

The convergence result relies on Theorem \ref{ch1'-s60-t1} below  on existence and uniqueness of solutions to \eqref{ch1'-s60.2}. Its proof is a straightforward adaptation of that of Theorem \ref{ch1'-s5.t1}, mainly based on the following remark:
Because $\Pi_{V_n}$ is a contraction operator and by Burkholder's
 inequality (or the isometry property), for $p\ge 2$, $L^p(\Omega)$-norms of stochastic integrals with respect to $W^n$ are bounded by the same expressions as when the integrator is $W$.

\begin{thm}
\label{ch1'-s60-t1}
Fix $n\ge 1$. Under $(\bf H_\Gamma)$, $(\bf H_I)$ and  $(\bf H_L)$, there exists a random field solution
$$\bar u_n=\left(\bar u_n(t,x),\, (t,x)\in[0,T]\times D\right)$$
to \eqref{ch1'-s60.2}. In addition, for any $p>0$,
\beq
\label{ch1'-s60.3}
\sup_{n\ge 1}\, \sup_{(t,x)\in[0,T]\times D} E\left[\vert \bar u_n(t,x)\vert^p\right] < \infty,
\eeq
and the solution $\bar u_n$ is unique (in the sense of versions) among random field solutions that satisfy \eqref{ch1'-s5.2} with $p=2$.
\end{thm}

We now give the statement on approximation of solutions by finite-dimensional projections. Recall that $\bar D$ denotes the closure of $D$ in the Euclidean topology. 

\begin{thm}
\label{ch1'-s60-t2}
Suppose that the domain $D$ is bounded and the
 hypotheses $({\bf H_\Gamma})$, $({\bf H_I})$ and $(\bf H_L)$ of Theorem 4.4.1 hold for $x \in \bar D$. In addition,
we assume that the map $(t,x;r,z)\mapsto \Gamma(t,x;r,z)$ from $\{(t,x;r,z)\in[0,T]\times \bar D \times[0,T]\times D: 0\le r<t\le T\}$ satisfies \eqref{ch1'-s5.12}, that is,
\beq
\label{4.4(*1)}
\int_0^T dr \int_D dz \left(\Gamma(t,x;r,z)-\Gamma(s,y;r,z)\right)^2 \le C^2\, {\bf \Delta}_1(t,x;s,y)^2,
\eeq
where ${\bf \Delta}_1(t,x;s,y)$ is a metric on $[0,T]\times \bar D$
with the same open sets as the Euclidean metric.
Consider the random field solutions to the SPDEs
\eqref{ch1'-s5.1} and \eqref{ch1'-s60.2},
\beqn
 u=\left(u(t,x),\, (t,x)\in[0,T]\times \bar D\right)\quad {\text{and}}\quad \bar u_n=\left(\bar u_n(t,x),\, (t,x)\in[0,T]\times \bar D\right),
 \eeqn
 respectively. Then for any $p>0$,
\beq
\label{ch1'-s60.4}
\lim_{n\to\infty}\, \sup_{(t,x)\in[0,T]\times \bar D} E\left[\left\vert u(t,x) - \bar u_n(t,x)\right\vert^p\right] = 0.
\eeq
\end{thm}
\begin{proof}
From the expressions \eqref{ch1'-s5.1} and \eqref{ch1'-s60.2}, we have
\beqn
u(t,x) - \bar u_n(t,x) = \mathcal{I}_n(t,x) +  \mathcal{R}_n(t,x) + \mathcal{J}_n(t,x),
\eeqn
with
\begin{align}
\label{ch1'-s60.4-bis}
\mathcal{I}_n(t,x) & = \sum_{j=1}^n \int_0^t \langle \Gamma(t,x; s,\ast)[\sigma(s,\ast,u(s,\ast))- \sigma(s,\ast,\bar u_n(s,\ast))],e_j\rangle_V\notag\\
&\qquad\qquad\times dW_s(e_j),\notag\\
\mathcal{R}_n(t,x) & = \sum_{j=n+1}^\infty\int_0^t \langle \Gamma(t,x; s,\ast)\sigma(s,\ast,u(s,\ast)),e_j\rangle_V\, dW_s(e_j),\notag\\
\mathcal{J}_n(t,x) & = \int_0^t \int_D \Gamma(t,x;s,y) [b(s,y,u(s,y))- b(s,y,\bar u_n(s,y))]\, ds dy.
\end{align}

Let $p\ge 2$ and $(t,x)\in[0,T]\times \bar D$. With the same approach used for instance to obtain \eqref{ch1'-s5.6b}, with $b(s,y,u^n(s,y))$ there replaced by $b(s,y,u(s,y))- b(s,y,\bar u_n(s,y))$, and using the Lipschitz property $({\bf H_L})$(v) instead of $({\bf H_L})$(vi), we have
\begin{align}
\label{(*J1)}
E[|\mathcal{J}_n(t,x)|^p] &\le C_p\left(\int_0^t ds\,J_2(s)\right)^{p-1}\notag\\
&\qquad\qquad \times \int_0^t ds\, J_2(t-s) \sup_{y\in D} E[|u(s,y)-\bar u_n(s,y)|^p],
\end{align}
where $J_2$ is defined in \eqref{j}.

Using the same arguments as in \eqref{ch1'-s5.6a}, a similar estimate for the term $\mathcal{I}_n(t,x)$ follows. More precisely,
\begin{align}
\label{(*J2)}
E[|\mathcal{I}_n(t,x)|^p] &\le \tilde C_p\left(\int_0^t ds\, J_1(s)\right)^{\frac{p}{2}-1}\notag\\
&\qquad \qquad\times \int_0^t ds\, J_1(t-s) \sup_{y\in D} E[|u(s,y)-\bar u_n(s,y)|^p]
\end{align}
(see \eqref{j} for the definition of $J_1$).

Let $Z(s, y) = \sigma(s, y, u(s, y))$. By hypothesis $({\bf H_L})$(vi) and \eqref{ch1'-s5.2}, the property \eqref{ch1'-s5.2} also holds with $u$ replaced by $Z$. Notice that
\beqn
    \mathcal{R}_n(t,x) = \int_0^t \int_D \Pi_{V_n^\perp}(\Gamma(t,x; s,*)Z(s,*))(y)\, W(ds, dy),
    \eeqn
where $\Pi_{V_n^\perp}$\label{rdorthperp} denotes the orthogonal projection from $V$ onto the orthogonal complement of $V_n$ (which is the subspace spanned by $(e_j,\, j \geq n+1))$.

To simplify the notation, set
\begin{align}
\label{def-ks}
 \mathcal{K}^{t,x}(s) &= \sum_{j=1}^\infty\langle \Gamma(t,x; s,\ast)\sigma(s,\ast, u(s,\ast)),e_j\rangle_V^2\notag\\
&\qquad = \Vert \Gamma(t,x; s,\ast)\sigma(s,\ast,u(s,\ast))\Vert_V^2,\notag\\
\mathcal{K}_n^{t,x}(s) &= \sum_{j=n+1}^\infty\langle \Gamma(t,x; s,\ast)\sigma(s,\ast, u(s,\ast)),e_j\rangle_V^2\notag\\
&\qquad = \left\Vert \Pi_{V_n^\perp}(\Gamma(t,x;s,\ast)\sigma(s,\ast,u(s,\ast)))\right\Vert_V^2.
\end{align}
Notice  that 
\beqn
E[|\mathcal{R}_n(t,x)|^2] = E\left[\int_0^t ds\, \mathcal{K}_n^{t,x}(s)\right].
\eeqn
Clearly, a.s., for all $s\in[0,T]$,
 \beq
 \label{sup-bound}
 \sup_{n\ge 0} \mathcal{K}_n^{t,x}(s) = \mathcal{K}^{t,x}(s).
 \eeq

 We observe that
\beq
\label{boundcn-2}
E\left[\left(\int_0^t ds\ \mathcal{K}_n^{t,x}(s)\right)^{\frac{p}{2}}\right]
\le E\left[\left(\int_0^t ds\ \mathcal{K}^{t,x}(s)\right)^{\frac{p}{2}}\right]<\infty.
\eeq
Indeed, the first inequality is clear by \eqref{sup-bound}, and the second integral is finite
because by definition,
\begin{align*}
&E\left[\left(\int_0^t ds\ \mathcal{K}^{t,x}(s)\right)^{\frac{p}{2}}\right]\\
&\qquad\qquad= E\left[\left(\int_0^t ds\ \int_D dy\ \left(\Gamma(t,x; s,y)\sigma(s,y, u(s,y))\right)^2\right)^{\frac{p}{2}}\right],
 \end{align*}
and we can argue as in \eqref{ch1'-s5.6a} (with $u^n$ there replaced by $u$), and apply \eqref{ch1'-s5.2}.

Therefore, the series that defines $\mathcal{K}^{t,x}(s)$ converges $ds dP$-a.e., and this implies
  \beq
 \label{conv-decreasing}
 \mathcal{K}_n^{t,x}(s)\downarrow_{n\to\infty} 0,\quad ds dP-{\text{a.e.}}
 \eeq
 We prove next that
 \beq
\label{non-uniform-conv}
\lim_{n\to\infty}E\left[|\mathcal{R}_n(t,x)|^p\right] =0,
\eeq
for any fixed $(t,x)\in [0,T]\times \bar D$.
Indeed, by  Burkholder's inequality (see Proposition \ref{Proposition 6}), we have
\beq
\label{boundcn}
E\left[|\mathcal{R}_n(t,x)|^p\right] \le \tilde c_p E\left[\left(\int_0^t ds\ \mathcal{K}_n^{t,x}(s)\right)^{\frac{p}{2}}\right].
\eeq
Using \eqref{conv-decreasing}, \eqref{sup-bound}, \eqref{boundcn-2}, we obtain
\eqref{non-uniform-conv} by dominated convergence.

Set
\beqn
\Psi_n(t,x) = E\left[\left(\int_0^t ds\ \mathcal{K}_n^{t,x}(s)\right)^{\frac{p}{2}}\right],\quad n\ge 1.
\eeqn
We prove in Lemma \ref{ch1'-s60-l1} below that the decreasing sequence $(\Psi_n,\, n\ge 1)$ consists of continuous functions in the compact metric space $([0,T]\times \bar D, {\bf \Delta}_1)$.
Thus, by Dini's theorem (see e.g. \cite[Theorem 7.13, p. 150]{rudin1}),
we have
\beq
\label{ch1'-s60.5}
\sup_{(t,x)\in[0,T]\times \bar D}\Psi_n(t,x) \downarrow _{n\to\infty}0.
\eeq

 Let $\Phi_n(t) = \sup_{x\in \bar D}E\left[\left\vert u(t,x) - \bar u_n(t,x)\right\vert^p\right]$. By \eqref{(*J1)}, \eqref{(*J2)} and \eqref{boundcn}, we have proved that for $t\in[0,T],$
 \beqn
 \Phi_n(t) \le \tilde c_p \int_0^t ds\ [J_1(t-s) + J_2(t-s)] \Phi_n(s) + \sup_{(t,x)\in[0,T]\times \bar D}\Psi_n(t,x).
 \eeqn
 By \eqref{ch1'-s60.5}, for any $\varepsilon>0$, there exists $n_0\ge 1$ such that for all $n\ge n_0$,
 \beqn
  \Phi_n(t) \le \tilde c_p \int_0^t ds\ [J_1(t-s) + J_2(t-s)] \Phi_n(s) + \varepsilon.
  \eeqn
  Fix $n\ge n_0$, and apply the inequality \eqref{A3.04bis} in Gronwall's Lemma \ref{A3-l01} to $f(t):=\Phi_n(t)$, $z_0 = 0$ and $z(t)=\varepsilon$ there.  We deduce that
  $\sup_{t\in[0,T]} \Phi_n(t) \le C \varepsilon$ (for some constant $C\ge 0$ that does not depend on $n$). Since $\varepsilon>0$ is arbitary, we obtain
  \eqref{ch1'-s60.4}.
\end{proof}

The following result has been used in the proof of Theorem \ref{ch1'-s60-t2}.

\begin{lemma}
\label{ch1'-s60-l1}
Assume that the hypotheses of Theorem \ref{ch1'-s60-t2} are satisfied.
Then
\beqn
[0,T]\times \bar D \ni (t,x)\longrightarrow \left\Vert\int_0^t dr\, \mathcal{K}_n^{t,x}(r)\right\Vert_{L^p(\Omega)}
\eeqn
is continuous with respect to ${\bf \Delta}_1$.
\end{lemma}
\begin{proof}
 Let $U$ denote the vector space of processes $G = (G(s,y), (s,y) \in [0,T] \times D)$ as in Section \ref{ch2new-s2}  that satisfy $\Vert G \Vert_U < \infty$, where
 \beqn
       \Vert G \Vert_U := \Vert \Vert G(\cdot, *)\Vert_{L^2([0, T] \times D)}  \Vert_{L^p(\Omega)}.
       \eeqn
Notice that
    $(t,x) \mapsto \Gamma(t,x; \cdot,*)Z(\cdot,*)$, where $Z(\cdot,\ast) = \sigma(\cdot,\ast, u(\cdot,\ast))$,
defines a (uniformly) continuous function from $[0,T] \times \bar D$ into $U$, such that
\beq
\label{4.4(*2)}
    \Vert \Gamma(t,x; \cdot,*)Z(\cdot,*) - \Gamma(s,y; \cdot,*)Z(\cdot,*) \Vert_U  \leq \tilde C_p {\bf\Delta}(t, x; s, y)).
    \eeq
Indeed, this follows from \eqref{4.4(*1)} and the calculations in \eqref{computations-1}.
   The operator $\Pi_{V_n^\perp}$ is a contraction (hence a continuous function) from $L^2(D)$ into $V_n^\perp \subset L^2(D)$, and we can view it as a contraction from $L^2([0, T] \times D)$ into itself (which does not depend on the $t$-coordinate), therefore,
   $(t,x) \mapsto \Pi_{V_n^\perp}(\Gamma(t,x; \cdot,*)Z(\cdot,*))$
is also a continuous function from from $[0,T] \times \bar D$ into $U$ and
\beq
\label{4.4(*3)}
     \Vert \Pi_{V_n^\perp}(\Gamma(t,x; \cdot,*)Z(\cdot,*) - \Gamma(s,y; \cdot,*)Z(\cdot,*)) \Vert_U  \leq \tilde C_p {\bf\Delta}(t, x; s, y).
     \eeq
    It follows that
     $ (t,x) \mapsto \Vert \Pi_{V_n^\perp}(\Gamma(t,x; \cdot,*)Z(\cdot,*))\Vert_U$
is a continuous function from $[0,T] \times \bar D$ into $\R$. Writing the $\Vert \cdot \Vert_U$ more explicitly, this means that
\begin{align*}
   (t,x) &\mapsto \left(E\left[\left \Vert \Pi_{V_n^\perp}(\Gamma(t,x; \cdot,*)Z(\cdot,*))\right\Vert_{L^2([0, T] \times D}^p\right]\right)^{1/p}\\
&\qquad= \left(E\left[\left(\int_0^t \mathcal{K}_n^{t,x}(r) dr\right)^{p/2}\right]\right)^{1/p}
\end{align*}
is continuous, which is the conclusion of the lemma.
   \end{proof}
\medskip

\noindent{\em Examples}
\medskip

Let $D=\,]0,L[$ and assume $({\bf H_L})$ and that the condition $({\bf H_I})$ holds on $[0,L]$. Then the nonlinear stochastic heat equation on $D$ with homogeneous Dirichlet or Neumann boundary conditions, and the nonlinear stochastic wave equation on $D$ with homogeneous Dirichlet boundary conditions, satisfy the conclusions of Theorem \ref{ch1'-s60-t2} with $\bar D=[0,L]$.

Indeed, for the stochastic heat equation, we apply Theorem \ref{recap}  and use the fact that, by \eqref{rde4.2.5}, we have
\beqn
{\bf \Delta}_1(t,x;s,y) = |t-s|^{\frac{1}{4}}+|x-y|^{\half}.
\eeqn

For the the stochastic wave equation, we apply Theorem \ref{recapwave} and \eqref{rde4.2.5-wave}, which tell us that we can take
\beqn
{\bf \Delta}_1(t,x;s,y) = |t-s|^{\half}+|x-y|^{\half}.
\eeqn


\section[Non-linear SPDEs on bounded domains with locally \texorpdfstring{\\}{}Lipschitz coefficients]
{Non-linear SPDEs on bounded domains with locally Lipschitz coefficients}
\label{ch1'-s7}

In this section, we assume that $D$ is a bounded domain of $\rek$
with smooth boundary. We fix $T > 0$ and consider $W$ and $(\cF_s)$ as in Section \ref{ch4-section0}. We will discuss SPDEs with coefficients more general than in \eqref{ch1'-s5.0}, formally written as
\begin{align}
\label{ch1'-s5.0-bis}
\mathcal{L}u(t,x) = \sigma(t,x,u) \dot W(t,x) + b(t,x,u), \quad (t,x)\in\, ]0,T]\times D,
\end{align}
with given initial conditions and boundary conditions, where $\mathcal{L}$ is a linear partial differential operator. First, we prove a theorem on existence and uniqueness of solutions to \eqref{ch1'-s5.0-bis} when the coefficients depend on the past of the solution and satisfy a global Lipschitz condition (Theorem \ref{ch1'-s7-t1}). Then we relax the assumptions on the coefficients, assuming a local Lipschitz (and linear growth) hypothesis and prove a result on uniqueness among continuous solutions (Theorem \ref{ch1'-s7-t3}) and a theorem on global existence (Theorem \ref{ch1'-s7-t4}).

\subsection[Assumptions on the Green's function and \texorpdfstring{\\}{}$L^p$-bounds on increments]{Assumptions on the Green's function and $L^p$-bounds on increments}
\label{ch1'-s7-new}

Throughout this section, the notation ${\bf \Delta}_i(t,x;s,y)$, $i= 3,4$, refers to an arbitrary nonnegative function defined for $(t,x), (s,y) \in \R_+ \times D$.

We start by introducing the assumptions on the Green's function of the partial differential operator $\mathcal{L}$.
\medskip

\noindent $({\bf h_\Gamma})$\label{rdhgamma} {\em Assumptions on the Green's function}
\begin{description}
\item{(i)} The mapping $(t,x;s,y)\mapsto \Gamma(t,x;s,y)$ from $\{(t,x;s,y)\in[0,T]\times D\times [0,T]\times D: 0\le s<t\le T\}$ into $\IR$ is jointly measurable.
\item{(ii)} There is a Borel function $H: [0,T]\times D^2 \longrightarrow \IR_+$ such that
\beqn
\vert\Gamma(t,x;s,y)\vert \le H(t-s,x,y), \qquad 0\le s<t\le T, \quad x,y\in D.
\eeqn
\item{(iii)} There is $\varepsilon_0>0$ such that if $\gamma: = 2+\varepsilon_0$, then
\beq
\label{C1}
\sup_{x\in D}\int_0^T ds\ \int_D dy\, H^\gamma(s,x,y) < \infty.
\eeq
\item{(iv)} Let $\varepsilon_0$ and $\gamma$ be as in (iii). There exists $c_{T,\gamma}<\infty$ such that for all $(t,x),\, (s,y)\in [0,T]\times D$,
\beq
\label{C2}
\int_0^T dr \int_D dz\ |\Gamma(t,x;r,z)-\Gamma(s,y;r,z)|^\gamma \le c_{T,\gamma} ({\bf \Delta}_3(t,x;s,y))^\gamma.
\eeq
\item{(v)} Let $\varepsilon_0$ be as in (iii) and $\mu=1+\varepsilon_0$. There exists $c_{T,\mu}<\infty$ such that for all $(t,x),\, (s,y)\in [0,T]\times D$,
\beq
\label{C2-bis}
\int_0^T dr \int_D dz\ |\Gamma(t,x;r,z)-\Gamma(s,y;r,z)|^\mu \le c_{T,\mu}\ ({\bf \Delta}_4(t,x;s,y))^\mu.
\eeq
\end{description}
 \begin{remark}
\label{ch1'-s7-r1}
We are assuming that $D$ is bounded. Therefore, {\rm(iii)} impliesthat \eqref{C1} also holds for any positive exponent $\tilde\gamma < \gamma$.
Observe the relation between the condition {\rm(iii)} above and {\rm(iiia)} and {\rm(iiib)} of the set of conditions (${\bf H_\Gamma}$) in Section \ref{ch4-section0}. Also, conditions {\rm(iv)} and {\rm(v)} can be compared with the hypotheses \eqref{ch1'-s5.12} and \eqref{ch1'-s5.14}. Condition {\rm (iv)} implies {\rm(v)} with ${\bf \Delta}_4 = {\bf \Delta}_3$.

Later on, we will consider the case ${\bf \Delta}_3(t,x;s,y)=|t-s|^{\alpha_1} + |x-y|^{\alpha_2}$ and
${\bf \Delta}_4(t,x;s,y)=|t-s|^{\beta_1} + |x-y|^{\beta_2}$ for some exponents $\alpha_1, \alpha_2, \beta_1, \beta_2 \in\ ]0,1]$
(as in \eqref{choice-rho}).
\end{remark}
\medskip

We now state some technical lemmas that will be used later on in this section.

For a jointly measurable and adapted process $Z=(Z(t,x),\, (t,x)\in[0,T]\times D)$, we define\label{rdnormtinfty}
\beqn
\Vert Z \Vert_{t,\infty} = \sup_{(s,x)\in[0,t]\times D} |Z(s,x)|.
\eeqn


\begin{lemma}
\label{ch1'-s7-l1}
We assume that $D$ is bounded and that the function $\Gamma(t,x;s,y)$ satisfies the Assumptions $(\bf {h_\Gamma})$ $(i)-(iv)$. Fix $T>0$ and let $Z=(Z(t,x),\, (t,x)\in[0,T]\times D)$ be a jointly measurable and adapted process such that for all $(t,x)\in[0,T]\times D$,
\beq
\label{ch1'-s7-l1.111}
E\left[\int_0^t dr \int_D dz\,  \left(\Gamma(t,x;r,z) Z(r,z)\right)^2\right] < \infty.
\eeq
For any $(t,x)\in[0,T]\times D$, set
\begin{align*}
A(t,x) &= \int_0^t \int_D \Gamma(t,x;r,z) Z(r,z)\, W(dr,dz).
\end{align*}
Let $p_0>2$ be such that $\frac{2p_0}{p_0-2} < 2+\varepsilon_0$, where $\varepsilon_0>0$ is given in hypothesis $({\bf h_\Gamma}) (iii)$.
Then for any $p\ge p_0$, we have the following:
\begin{enumerate}
\item There exists a constant $C_{p,T,D}<\infty$ such that, for any $t\in[0,T]$,
\begin{align}
\label{ch1'-s7.l1-1}
\sup_{x\in D}E\left[\vert A(t,x)\vert^p\right]
&\le C_{p,T, D}\int_0^t dr\, \sup_{z\in D} E[|Z(r,z)|^p].
\end{align}
Consequently,
\beq
\label{ch1'-s7.l1-1bis}
\sup_{x\in D}E\left[\vert A(t,x)\vert^p\right]
\le C_{p,T,D}\int_0^t dr\,  E\left[\Vert Z\Vert^p_{r,\infty}\right].
\eeq
\item There exists a constant $C_{p,T,D}<\infty$ such that, for any $(t,x),\, (s,y)\in[0,T]\times D$ with $0\le s\le t\le T$,
\begin{align}
\label{ch1-s7.l1-2}
E\left[\vert A(t,x)-A(s,y)\vert^p\right]
& \le C_{p,T,D}\  {\bf \Delta}_3(t,x;s,y)^p \int_0^t dr \sup_{z\in D}E[|Z(r,z)|^p],
\end{align}
which in turn implies
\beq
\label{ch1-s7.l1-2bis}
E\left[\vert A(t,x)-A(s,y)\vert^p\right]
\le  C_{p,T,D}\  {\bf \Delta}_3(t,x;s,y)^p \int_0^t dr\,  E\left[\Vert Z\Vert^p_{r,\infty}\right].
\eeq
\item Consider the particular case where
\beqn
{\bf \Delta}_3(t,x;s,y)=|t-s|^{\alpha_1}+|x-y|^{\alpha_2},
\eeqn
 with $\alpha_1, \alpha_2\in\ ]0,1]$. Suppose that $\int_0^T dr\, E[\Vert Z\Vert_{r,\infty}^p] < \infty$. Then
  $(A(t,x),\, (t,x) \in [0, T] \times D)$ has a Hölder continuous version (extended to $\bar D$) $(\tilde A(t,x),\, (t,x) \in [0, T] \times \bar D)$  with exponents $(\gamma_1, \gamma_2)$, where $\gamma_1 \in\, ]0, \alpha_1[$ and $\gamma_2 \in\, ]0, \alpha_2[$. Further, for any $p>p_0\vee \left(\frac{1}{\alpha_1}+\frac{k}{\alpha_2}\right)$ and for any $t\in[0,T]$,
\beq
\label{ch1'-s7.l1-2000}
E\left[\Vert \tilde A\Vert_{t,\infty}^p\right] \le \tilde C_{p,T,D,\alpha_1,\alpha_2} \int_0^t dr\,  E\left(\Vert Z\Vert^p_{r,\infty}\right),
\eeq
with $\tilde C_{p,T,D,\alpha_1,\alpha_2}<\infty$.
 \end{enumerate}
\end{lemma}

\begin{proof}
Notice that the hypotheses of the lemma imply that the stochastic integral process $(A(t,x),\ (t,x)\in[0,T]\times D)$ is well-defined in the sense of Definition \ref{ch1'-s4.d1}.

1. Fix $p\ge p_0$. By Burkholder's inequality \eqref{ch1'-s4.4.mod},
\beqn
E\left([\vert A(t,x)\vert^p\right] \le C_p\ E\left[\left(\int_0^t dr \int_D dz\ \Gamma^2(t,x;r,z) Z^2(r,z)\right)^{\frac{p}{2}}\right].
\eeqn
Apply H\"older's inequality with exponents $\frac{p}{2}$ and $\tilde p:=\frac{p}{p-2}$ (whose inverses sum to $1$) to the $drdz$-integral. We deduce that the last expression is bounded above by
\begin{align*}
&C_p E\left[\left[\left( \int_0^t dr \int_D dz\  \Gamma^{2\tilde p}(t,x;r,z)\right)^{\frac{1}{\tilde p}} \left( \int_0^t dr \int_D dz\ |Z(r,z)|^p\right)^{\frac{2}{p}}\right]^{\frac{p}{2}}\right]\notag\\
 & = C_p\left(\int_0^t dr \int_D dz\ \Gamma^{2\tilde p}(t,x;r,z)\right)^{\frac{p}{2\tilde p}}  \int_0^t dr \int_D dz\ E\left[|Z(r,z)|^p\right] \notag\\
 &\quad \le C_p \sup_{x\in D}\left( \int_0^t dr \int_D dz\  H^{2\tilde p}(r,x,y)\right)^{\frac{p}{2\tilde p}} \int_0^t dr \int_D dz\ E\left[|Z(r,z)|^p\right].
 \end{align*}
 Since $p \geq p_0$ and $\frac{2 p_0}{p_0 - 2} < 2 + \ep_0 =: \gamma$ by assumption, we have $2\tilde p\in[2,2+\epsilon_0[$.
 Therefore, using \eqref{C1} and the fact that $D$ is bounded, we obtain
 \begin{align*}
 E\left[\vert A(t,x)\vert^p\right] &\le C_{p,T,D}\int_0^t dr \int_D dz\  E[|Z(r,z)|^p] \\
 &\le C_{p,T,D} \int_0^t dr\ \sup_{z\in D} E[|Z(r,z)|^p].
 \end{align*}
Hence, we have proved \eqref{ch1'-s7.l1-1} which in turn implies \eqref{ch1'-s7.l1-1bis}.
\medskip

2.  We prove \eqref{ch1-s7.l1-2}, with the same approach as in part 1 above, based on Burkholder's and H\"older's inequality, but \eqref{C2} instead of \eqref{C1}. By doing so, for $0\le s\le t\le T$, we obtain
\begin{align*}
E\left[\vert A(t,x)-A(s,y)\vert^p\right] &\le C_p  \left(\int_0^t dr \int_D dz\ \vert\Gamma(t,x;r,z)-\Gamma(s,y;r,z)\vert^{2\tilde p}\right)^{\frac{p}{2\tilde p}}\notag\\
&\quad \qquad\qquad\times \int_0^t dr \int_D dz\ E\, [|Z(r,z)|^p] \\
&\le C_p  \left(\int_0^t dr \int_D dz\ \vert\Gamma(t,x;r,z)-\Gamma(s,y;r,z)\vert^{\gamma}\right)^{\frac{p}{\gamma}}\notag\\
&\quad \qquad\qquad\times \int_0^t dr \int_D dz\, E\, [|Z(r,z)|^p].
\end{align*}
Because of
\eqref{C2}, this is bounded by
\beqn
C_{p,T}\  {\bf \Delta}_3(t,x;s,y)^p \int_0^t dr \int_D dz\, E[|Z(r,z)|^p].
\eeqn
Since $D$ is bounded, this implies \eqref{ch1-s7.l1-2}, which in turn implies \eqref{ch1-s7.l1-2bis}.
\medskip

3. For ${\bf \Delta}_3(t,x;s,y)$ as in Claim 3, the inequality \eqref{ch1-s7.l1-2bis} reads
\beq
\label{Aquatretres}
E\left[\vert A(t,x)-A(s,y)\vert^p\right] \le C_{p,T,D} \left(|t-s|^{\alpha_1}+|x-y|^{\alpha_2}\right)^p \int_0^t dr\,  E\left[\Vert Z\Vert^p_{r,\infty}\right].
\eeq
Apply this with $t$ replaced by $r$, for all $0 \leq s \leq r \leq t$. Then bounding the integral from $0$ to $r$ (that appears on the right-hand side) by the integral from $0$ to $t$, we see that for all $0 \leq s \leq r \leq t$,
\beq
\label{Aquatretres(*A)}
E\left[\vert A(r,x)-A(s,y)\vert^p\right] \le C_{p,T,D} \left(|r-s|^{\alpha_1}+|x-y|^{\alpha_2}\right)^p \int_0^t d\rho\,  E\left[\Vert Z\Vert^p_{\rho,\infty}\right].
\eeq
The existence of the process $(\tilde A(t,x),\, (t,x)\in[0,T]\times \bar D)$ follows from Kolmogorov's continuity criterion (Theorem \ref{ch1'-s7-t2}). Indeed,
\eqref{Aquatretres(*A)} is assumption \eqref{ch1'-s7.18} with $u(t,x):=A(t,x)$ and $K:= C_{p,T,D}\ \int_0^T dr\,  E\left[\Vert Z\Vert^p_{r,\infty}\right]$ there.

Fix $t\in[0,T]$ and let $\alpha\in\ \left]\frac{1}{p}\left(\frac{1}{\alpha_1}+\frac{k}{\alpha_2}\right), 1\right[$. Using again \eqref{Aquatretres(*A)} and applying Theorem \ref{ch1'-s7-t2} then \eqref{ch1'-s7.20bis}, but with $\tilde u(s,x):=\tilde A(s,x)$, $I=[0,t]$ and
\beqn
K=C_{p,T,D}\int_0^t dr\, E\left[\Vert Z\Vert_{r,\infty}^p\right],
\eeqn
yields
\begin{align*}
E\left[\Vert \tilde A\Vert^p_{t,\infty}\right]&\le  c_2(I,D,\alpha,p,Q)t^{p \alpha_0 (\alpha - Q/p)}C_{p,T,D}\int_0^t dr\, E\left[\Vert Z\Vert_{r,\infty}^p\right].
\end{align*}
 Indeed, observe that, since $\tilde A(0,x)=0$, the constant $C_1$ on the right-hand side of \eqref{ch1'-s7.20bis} can be set to $0$.
This implies \eqref{ch1'-s7.l1-2000}.
\end{proof}


\begin{lemma}
\label{ch1'-s7-l2}
We assume that the function $\Gamma(t,x;s,y)$ satisfies the Assumptions $(\bf h_\Gamma)$ $(i)-(iii)$ and $(v)$.
Let $Z=(Z(t,x), (t,x)\in[0,T]\times D)$ be a jointly measurable and adapted process such that for all $(t,x)\in[0,T]\times D$,
\beq
\label{ch1'-s7-l2.111}
E\left[\int_0^t dr \int_D dz\  |\Gamma(t,x;r,z) Z(r,z)|\right] < \infty.
\eeq
For any $(t,x)\in[0,T]\times D$, set
\begin{align*}
B(t,x) & = \int_0^t dr \int_D dz\ \Gamma(t,x;r,z) Z(r,z).
\end{align*}
Let $p_0>1$ be such that $\tfrac{p_0}{p_0-1}<1+\varepsilon_0$. Then for any $p\ge  p_0$, we have the following:
\begin{enumerate}
\item There exists a constant $C_{p,T,D}<\infty$ such that, for any $t\in[0,T]$,
\begin{align}
\label{ch1'-s7.l1-3}
\sup_{x\in D}E\left[\vert B(t,x)\vert^p\right] &\le C_{p,T,D}\int_0^t dr \sup_{z\in D}E[|Z(r,z)|^p].
\end{align}
In particular,
\beq
\label{ch1'-s7.l1-3bis}
\sup_{x\in D}E\left[\vert B(t,x)\vert^p\right]
\le C_{p,T,D}\int_0^t dr E\left[\Vert Z\Vert^p_{r,\infty}\right].
\eeq
\item There exists a constant $C_{p,T,D}<\infty$ such that, for any $(t,x), (s,y)\in[0,T]\times D$ with $0\le s\le t\le T$,
\begin{align}
\label{ch1-s7.l1-4}
E\left[\vert B(t,x)-B(s,y)\vert^p\right]
& \le C_{p,T,D}\  {\bf \Delta}_4(t,x;s,y)^p \int_0^t dr\, \sup_{z\in D} E[|Z(r,z)|^p].
\end{align}
In particular,
\beq
\label{ch1-s7.l1-4bis}
E\left[\vert B(t,x)-B(s,y)\vert^p\right]
\le C_{p,T,D}\ {\bf \Delta}_4(t,x;s,y)^p \int_0^t dr\, E\left[\Vert Z\Vert^p_{r,\infty}\right].
\eeq

\item  Consider the case where ${\bf \Delta}_4(t,x;s,y)=|t-s|^{\beta_1}+|x-y|^{\beta_2}$, with $\beta_1, \beta_2\in\ ]0,1]$. Assume that $\int_0^T dr\, E\left[\Vert Z\Vert^p_{r,\infty}\right] < \infty$. Then
$(B(t,x),\, (t,x) \in [0, T] \times D)$ has a Hölder continuous version (extended to $\bar D$) $(\tilde B(t,x),\, (t,x) \in [0, T] \times \bar D)$ with exponents $(\eta_1, \eta_2)$, where $\eta_1 \in\, ]0, \beta_1[$ and $\eta_2 \in\, ]0, \beta_2[$. Further, for any $p>p_0\vee \left(\frac{1}{\beta_1}+\frac{k}{\beta_2}\right)$ and for any $t\in[0,T]$,
\beq
\label{ch1'-s7.l1-201}
E\left[\Vert \tilde B\Vert_{t,\infty}^p\right] \le \tilde C(p,T,D,\beta_1,\beta_2) \int_0^t dr\,  E\left[\Vert Z\Vert^p_{r,\infty}\right],
\eeq
where $\tilde C(T,p,D,\beta_1,\beta_2)<\infty$.
\end{enumerate}
\end{lemma}

\begin{proof}
Fix $p\ge p_0$ and apply H\"older's inequality with exponents $p$ and $q=\frac{p}{p-1}$ to obtain
\beqn
E\left[\vert B(t,x)\vert^p\right] \le C_p\left(\int_0^t dr \int_D dz\ |\Gamma(t,x;r,z)|^q\right)^{\frac{p}{q}}
\int_0^t dr \int_D dz\ E[|Z(r,z)|^p].
\eeqn
Since $p \geq p_0$ and $\frac{p_0}{p_0 - 1} < 1 + \ep_0$ by assumption, we have $q\in [1,1+\varepsilon_0[$. As in Lemma \ref{ch1'-s7-l1}, since $D$ is bounded, \eqref{ch1'-s7.l1-3} follows from \eqref{C1}.

The estimate \eqref{ch1-s7.l1-4} is proved using similar arguments, by replacing the expression $\Gamma(t,x;r,z)$ by $\Gamma(t,x;r,z)-\Gamma(s,y;r,z)$, and applying \eqref{C2-bis}.

Applying Kolmogorov's continuity criterion in the same way as in the proof of Lemma \ref{ch1'-s7-l1}, we obtain Claim 3.
\end{proof}


\subsection[A theorem on existence and uniqueness under\texorpdfstring{\\}{}global Lipschitz conditions]{A theorem on existence and uniqueness under global Lipschitz conditions}
\label{ch1'-ss7.1}

In this section, we consider the integral equation
\begin{align}
\label{ch1'-s7.1}
u(t,x)& = I_0(t,x) + \int_0^t \int_D \Gamma(t,x;r,z) \sigma(r,z,u(\cdot,\ast))\, W(dr,dz)\notag\\
& \qquad + \int_0^t dr \int_D dz \, \Gamma(t,x;r,z) b(r,z,u(\cdot,\ast)),
\end{align}
$(t,x)\in[0,T]\times D$, where for the sake of simplicity, we have omitted the dependence on $\omega$ in $\sigma$ and $b$ (see Assumption $({\bf h_{L}})$ ibelow), and the function $(t,x) \mapsto I_0(t,x)$ is the solution to the homogeneous PDE $\mathcal{L} u = 0$, with the given initial conditions and boundary conditions.

In \eqref{ch1'-s7.1}, the stochastic integral is in the sense of Definition \ref{ch1'-s4.d1} and the second integral is a Lebesgue integral.
\smallskip

\begin{def1}
\label{ch1'-s7-d1}
A  random field solution to \eqref{ch1'-s5.0-bis}  is an adapted process $u=\left(u(t,x),\, (t,x)\in[0,T]\times \bar D\right)$ with continuous sample paths such that, for all $(t,x)\in[0,T]\times D$, the two integrals in \eqref{ch1'-s7.1} are well-defined and  for all $(t,x)\in [0,T]\times D$, \eqref{ch1'-s7.1} holds a.s.
\end{def1}

This section is devoted to proving a theorem on existence and uniqueness of solutions to equation \eqref{ch1'-s5.0-bis}.
The coefficients $\sigma$ and $b$ are functions
\beqn
\sigma, b: [0,T]\times D\times \mathcal{C}([0,T]\times \bar D)\times \Omega \longrightarrow \re
\eeqn
satisfying a global Lipschitz condition, to be specified below, and $T>0$ is fixed.
In comparison with Theorem \ref{ch1'-s5.t1}, these coefficients are more general since they may depend on the past of the trajectories.

We recall that $\mathcal{B}_{\mathcal{C}([0,T]\times D)}$ denotes the $\sigma$-field generated by the open sets in the topology of the uniform convergence of functions in $\mathcal{C}([0,T]\times D)$.
\medskip

We introduce the following assumptions.
 \smallskip

\noindent $({\bf h_{L}})$\label{rdhL} {\em Assumptions on the coefficients $\sigma$ and $b$}
\begin{description}
\item{(vi)} {\em Measurability and adaptedness.} The functions $\sigma$ and $b$  are jointly measurable, that is, $\mathcal{B}_{[0,T]}\times \mathcal{B}_D\times\mathcal{B}_{\mathcal{C}([0,T]\times \bar D)}\times \tf$-measurable. These two functions are also adapted to $(\tf_s,\, s\in[0,T])$, that is, for fixed $s \in [0,T]$, $(y,v,\omega)\mapsto \sigma(s,y,v,\omega)$ and $(y,v,\omega)\mapsto b(s,y,v,\omega)$ are $\cB_D\times \cB_{\mathcal{C}([0,T]\times \bar D)}\times \cF_s$-measurable.
\item{(vii)} {\em Non-anticipating property.}
For any $v\in\mathcal{C}([0,T]\times \bar D)$ and $t>0$, we define $v^t\in\mathcal{C}([0,T]\times \bar D)$ by $v^t(s,y)= v(t\wedge s,y)$. Then for any $v\in \mathcal{C}([0,T]\times \bar D)$ and $(t,x,\omega)\in [0,T]\times  D\times \Omega$,
\beqn
\sigma(t,x,v,\omega)= \sigma(t,x,v^t,\omega), \qquad b(t,x,v,\omega)=b(t,x,v^t,\omega).
\eeqn
\item{(viii-global)} {\em Global Lipschitz condition.}
There exists a constant $c_1(T)\in\IR_+$ such that, for all $(t,x)\in[0,T]\times D$, $v, \bar v\in \mathcal{C}([0,T]\times \bar D)$ and $\omega\in \Omega$,
\begin{align*}
\vert\sigma(t,x,v,\omega)-\sigma(t,x,\bar v,\omega)\vert &+ \vert b(t,x,v,\omega)-b(t,x,\bar v,\omega)\vert \\
&\le c_1(T)\Vert v-\bar v\Vert_{t,\infty}.
\end{align*}
\item{(viii-local)} {\em Local Lipschitz condition.}
For any $M>0$, there exists a constant $c_1(T,M)\in\IR_+$ such that for all $(t,x)\in[0,T]\times D$,
$v, \bar v\in \mathcal{C}([0,T]\times \bar D)$ satisfying $\Vert v\Vert_{T,\infty}\le M$ and $\Vert \bar v\Vert_{T,\infty}\le M$, and $\omega\in\Omega$,
\begin{align*}
\vert\sigma(t,x,v,\omega)-\sigma(t,x,\bar v,\omega)\vert &+ \vert b(t, x, v, \omega)-b(t, x, \bar v, \omega)\vert\\
&\le c_1(T,M)\Vert v-\bar v\Vert_{t,\infty}.
\end{align*}
\item{(ix)} {\em Uniform linear growth.} There exists a constant $c_2(T)$ such that, for all $(t,x)\in[0,T]\times D$, $v\in \mathcal{C}([0,T]\times \bar D)$ and $\omega\in\Omega$,
\beqn
\vert\sigma(t,x,v,\omega)\vert  + \vert b(t,x,v,\omega)\vert \le c_2(T)\left(1+\Vert v\Vert_{t,\infty}\right).
\eeqn
\end{description}

\noindent ${\bf (h_I)}$\label{rdhI} {\em Assumption  on the initial conditions}
\medskip

The real-valued function $(t,x)\mapsto I_0(t,x)$ defined on $[0,T]\times \bar D$ is continuous, jointly in $(t,x)$.
\medskip

We can now present the main result of this section in which ${\bf \Delta}_3(t,x;s,y)$ and ${\bf \Delta}_4(t,x;s,y)$ of hypothesis $({\bf h_\Gamma})$ are given by
\beq
\label{afegir(*1)}
{\bf \Delta}_3 (t,x;s,y)= |t-s|^{\alpha_1} + |x-y|^{\alpha_2}, \quad  {\bf \Delta}_4 (t,x;s,y)= |t-s|^{\beta_1} + |x-y|^{\beta_2},
\eeq
with $\alpha_1, \alpha_2, \beta_1, \beta_2\in\ ]0,1]$,  so that the conclusions of parts 3. of Lemmas \ref{ch1'-s7-l1} and \ref{ch1'-s7-l2} apply.

\begin{thm}
\label{ch1'-s7-t1}
Let ${\bf \Delta}_3$ and ${\bf \Delta}_4$ be as in \eqref{afegir(*1)}.
We assume the hypotheses ${(\bf h_\Gamma})$, ${(\bf h_I})$ and the conditions {\rm (vi), (vii), (viii-global)} and {\rm (ix)} of Assumptions $({\bf h_L})$.

   (a) There exists a random field solution
\beqn
(u(t,x),\, (t,x)\in[0,T]\times \bar D)
\eeqn
to \eqref{ch1'-s5.0-bis} (in the sense of Definition \ref{ch1'-s7-d1}). Furthermore, for any $p\ge 2$,
\beq
\label{ch1'-s7.3}
E\left[\Vert u\Vert^p_{T,\infty}\right] < \infty,
\eeq
and the solution $u$ is unique (up to indistinguishability) among random fields that satisfy \eqref{ch1'-s7.3}
with $p=2$.
\smallskip

   (b) Assume that the function $(t,x)\mapsto I_0(t,x)$ is H\"older continuous jointly in $(t,x)$ with exponents $(\eta_1, \eta_2)$. Then $u$ satisfies the following property:
for any $p\ge 2$, there is a constant $0\le C_p<\infty$ such that, for all $(t,x),\, (s,y)\in[0,T]\times D$,
\beq
\label{ch1'-s7.2}
\Vert u(t,x)-u(s,y)\Vert_{L^p(\Omega)} \le C_p\left(|t-s|^{\eta_1\wedge\alpha_1\wedge\beta_1}+|x-y|^{\eta_2\wedge\alpha_2\wedge\beta_2}\right).
\eeq
Therefore, $(u(t,x),\, (t,x)\in[0,T]\times \bar D)$
has a H\"older continuous version, jointly in $(t,x)$, with exponents $\delta_1\in\, ]0,\eta_1\wedge\alpha_1\wedge\beta_1[$, $\delta_2\in\, ]0,\eta_2\wedge\alpha_2\wedge\beta_2[$, respectively.
\end{thm}
\smallskip

\begin{remark}
\label{4.5.2-rem1} The condition \eqref{ch1'-s7.3} yields a stronger conclusion than \eqref{ch1'-s5.2} in Theorem \ref{ch1'-s5.t1}. Indeed,
\beqn
E\left[\Vert u\Vert^p_{T,\infty}\right] = E\left[\sup_{(s,x)\in [0,T]\times D} |u(s,x)|^p\right],
\eeqn
and we see that in comparison with \eqref{ch1'-s5.2}, here the supremum is inside the expectation.
\end{remark}

\noindent{\em Proof of Theorem \ref{ch1'-s7-t1}}. (a) We begin with the proof of existence and uniqueness of the solution, using the same approach as in the proof of Theorem \ref{ch1'-s5.t1}.

 Define the Picard iteration scheme: for $(t,x)\in[0,T]\times D$,
\begin{align}
\label{P-lL}
 u^0(t,x) &= I_0(t,x),\notag\\
 u^{n+1}(t,x) & = I_0(t,x) + \int_0^t \int_D \Gamma(t,x;s,y) \sigma(s,y,u^n)\, W(ds,dy)\notag\\
 &\qquad + \int_0^t \int_D \Gamma(t,x;s,y) b(s,y,u^n)\, ds dy, \quad n\ge 0.
 \end{align}
 \noindent{\em Step 1.} We prove by induction that, for each $n\ge 0$, the process
 \beqn
 u^n = \left(u^n(t,x),\,  (t,x)\in [0,T] \times \bar D\right)
 \eeqn
 is well-defined, adapted and continuous (meaning that it has a version with continuous sample paths, extended to $\bar D$, which is again denoted by $u^n$, and we will always use this version)
and satisfies
 \beq
 \label{ch1'-s7.4}
E\left[\Vert u^n\Vert^p_{T,\infty}\right] < \infty,
 \eeq
 for any $p\ge 2$.

 Let us first see that these properties imply that the integrals on the right-hand side of \eqref{P-lL} are well-defined.
 Indeed,  the mapping $\omega\mapsto u^n(\cdot,\ast, \omega)$ from $(\Omega;\tf)$ into
 $(\mathcal{C}([0,T]\times \bar D);\mathcal{B}_{\mathcal{C}([0,T]\times \bar D)})$, is measurable because $u^n$ is a continuous version of the random field, which exists by parts 3. of Lemmas \ref{ch1'-s7-l1} and \ref{ch1'-s7-l2}.
 Furthermore, the mapping $(s,y,\omega)\mapsto \sigma(s,y,u^n(\cdot, \ast,\omega),\omega)$ is
 $\mathcal{B}_{[0,T]}\times \mathcal{B}_D\times \tf$-measurable, because it is the composition of the map
 $(s,y,\omega)\mapsto (s,y,u^n(\cdot, \ast, \omega),\omega)$ from $([0,T]\times D\times \Omega; \mathcal{B}_{[0,T]}\times \mathcal{B}_D\times \tf)$ into $([0,T]\times D\times \mathcal{C}([0,T]\times \bar D)\times \Omega; \mathcal{B}_{[0,T]}\times \mathcal{B}_D\times\mathcal{B}_{\mathcal{C}([0,T]\times \bar D)}\times \tf)$
 and $(s,y,v,\omega)\mapsto \sigma(s,y,v,\omega)$ from $([0,T]\times D\times \mathcal{C}([0,T]\times \bar D)\times \Omega;
 \mathcal{B}_{[0,T]}\times \mathcal{B}_D\times\mathcal{B}_{\mathcal{C}([0,T]\times \bar D)}\times \tf)$
 into $(\re, \mathcal{B}_\re)$.

The mapping $(s,y,\omega)\mapsto \sigma(s,t,u^n(\cdot,\ast,\omega),\omega)$ is adapted. Indeed, for fixed $s\in[0,T]$,
 $ (y,\omega)\mapsto \sigma(s,y,u^n(\cdot,\ast, \omega),\omega) = \sigma(s,y,(u^n)^s(\cdot,\ast, \omega),\omega)$
is the composition of the map $(y,\omega)\mapsto (y, (u^n)^s(\cdot,v, \omega),\omega)$ from $(D\times \Omega;  \mathcal{B}_D\times \tf_s)$ into $(D\times\mathcal{C}(\IR_+\times D)\times \Omega;\mathcal{B}_D\times \mathcal{B}_{\mathcal{C}(\re_+\times D)}\times \tf_s )$ and the map $(y,v,\omega)\mapsto \sigma(s,y,v,\omega)$ from $(D\times \mathcal{C}(\IR_+\times D)\times \Omega; \mathcal{B}_D\times\mathcal{B}_{\mathcal{C}(\re_+\times D)}\times \tf_s)$ into $(\re, \mathcal{B}_\re)$.

The integrand in the stochastic integral in \eqref{P-lL} is of the form \eqref{rd2.2.12a} with $Z(s,y)= \sigma(s,y,u^n)$. Let us check condition \eqref{rd2.2.12b}.
As a consequence of (ix) in (${\bf h_L}$), we have
\beq
\label{mesdelmateix}
\sup_{(s,y)\in[0,T]\times D} E\left[\left(\sigma(s,y,u^n)\right)^2\right]  \le c_2(T)\ \left(1+E\left[\Vert u^n\Vert^2_{T,\infty}\right]\right) < \infty,
\eeq
by \eqref{ch1'-s7.4}. Therefore,
\begin{align}
\label{lip-1}
&E\left[\int_0^t ds \int_D dy\ (\Gamma(t,x;s,y)\ \sigma(s,y,u^n))^2\right] \notag\\
&\qquad \le \sup_{(s,y)\in[0,T]\times D} E\left[(\sigma(s,y,u^n))^2\right] \int_0^t ds \int_D dy\ \Gamma^2(t,x;s,y) < \infty,
\end{align}
by (iii) in assumption (${\bf h_\Gamma}$), since $D$ is bounded. These considerations show that the stochastic integral in
\eqref{P-lL} is well-defined.

In a similar way, we can check that the deterministic integral in \eqref{P-lL} is well-defined.
 \smallskip

 We note that
 \beqn
 \Vert u^0\Vert_{T,\infty} = \Vert I_0\Vert_{T,\infty} < \infty,
 \eeqn
 by Assumption (${\bf h_I}$), since $\bar D$ is bounded. Therefore, $u^0$ satisfies \eqref{ch1'-s7.4} and the other properties described at the beginning of this Step 1.

 Assume that for some $n\ge 0$, the process
 \beqn
 \left(u^n(t,x),\, (t,x)\in[0,T]\times D\right)
 \eeqn
 is well-defined, continuous and adapted, and \eqref{ch1'-s7.4}
 holds.
 According to what we have just established, $u^{n+1}=(u^{n+1}(t,x),\, (t,x)\in[0,T]\times D)$ given in \eqref{P-lL} is well-defined. We want to show that $u^{n+1}$ is continuous (and extends to $\bar D$), adapted and \eqref{ch1'-s7.4} is satisfied with $n$ replaced by
 $n+1$.

 Define
 \begin{align*}
\mathcal {I}^n(t,x) &= \int_0^t \int_D \Gamma(t,x;s,y) \sigma(s,y,u^n)\, W(ds,dy),\\
\mathcal{J}^n(t,x) &=  \int_0^t ds \int_D dy\,  \Gamma(t,x;s,y) b(s,y,u^n).
\end{align*}
Let $Z^n(s,y) = \sigma(s,y,u^n)$. By assumption (ix) of (${\bf h_L}$),
\beqn
E\left[\Vert Z^n\Vert_{T,\infty}^p\right] \le c_2(T) \left(1+E\left[\Vert u^n\Vert_{T,\infty}^p\right]\right) < \infty,
\eeqn
by \eqref{ch1'-s7.4} (the induction hypothesis). Therefore, by part 3. of Lemma \ref{ch1'-s7-l1}, $(\mathcal{I}^n(t,x))$ has a continuous version (that extends to $\bar D$), which we again denote $(\mathcal{I}^n(t,x))$, and for $p\ge 2$ large enough,
\beqn
E\left[\Vert\mathcal{I}^n\Vert_{T,\infty}^p\right] \le \tilde c\  E\left[\Vert Z^n\Vert_{T,\infty}^p\right] < \infty,
\eeqn
where the constant $\tilde c$ depends on $p, T$ and $D$.

Similarly, letting  $Z^n(s,y) = b(s,y,u^n)$ and using Lemma \ref{ch1'-s7-l2} instead of Lemma \ref{ch1'-s7-l1}, we see that
$(\mathcal{J}^n(t,x))$ has a continuous version (which extends to $\bar D$), which we again denote $(\mathcal{J}^n(t,x))$, and for $p\ge 2$ large enough,
\beqn
E\left[\Vert\mathcal{J}^n\Vert_{T,\infty}^p\right] \le \tilde c\  E\left[\Vert Z^n\Vert_{T,\infty}^p\right] < \infty.
\eeqn
It follows that $u^{n+1}$ defined in \eqref{P-lL} is well-defined, continuous (and extends continuously to $\bar D$) and \eqref{ch1'-s7.4} holds with $u^n$ replaced by $u^{n+1}$.

For fixed $(t,x)\in [0,T]\times D$, $\mathcal{I}^n(t,x)$ is $\cF_t$-measurable by definition of the stochastic integral, therefore for $t\in[0,T]$, $(x,\omega) \mapsto \mathcal{I}^n(t,x,\omega)$ is $\cB_D \times \cF_t$-measurable since $\mathcal{I}^n$ is continuous. This implies that $\mathcal{I}^n$ is adapted, and the same is true of $\mathcal{J}^n$. Therefore, $u^{n+1}$ is adapted.
\medskip

\noindent{\em Step 2.} We now show that the sequence $(u^n,\, n\ge 0)$ of Picard iterations converges to a stochastic process $u=(u(t,x),\, (t,x)\in[0,T]\times \bar D)$, that is,
\beq
\label{ch1'-s7.146}
\lim_{n\to\infty} E\left[\Vert u^n-u\Vert^p_{T,\infty}\right] = 0.
\eeq
Indeed, for $(t,x)\in[0,T]\times D$, consider the difference of two consecutive Picard iterations,
\beqn
u^{n+1}(t,x)-u^n(t,x) = \left[\mathcal{I}^{n}(t,x)-\mathcal{I}^{n-1}(t,x)\right] + \left[\mathcal{J}^{n}(t,x)-\mathcal{J}^{n-1}(t,x)\right] .
\eeqn
The term $\mathcal{I}^{n}(t,x)-\mathcal{I}^{n-1}(t,x)$ is equal to the stochastic integral
\beq
\label{ch1'-s7.147}
\int_0^t \int_D \Gamma(t,x;r,z) Z(r,z)\, W(dr,dz),
\eeq
with $Z(r,z) = \sigma(r,z,u^n)-\sigma(r,z,u^{n-1})$. From (viii-global) in Assumption (${\bf h_L}$), we see that
\beq
\label{ch1'-s7.148}
|Z(r,x)| \le c_1(T)\Vert u^n-u^{n-1}\Vert_{r,\infty}.
\eeq
Hence,
\begin{align*}
E\left[\int_0^t dr \int_D dz\, \left(\Gamma(t,x;r,z) Z(r,z)\right)^2 \right]
&\le C(T) E\left[\Vert u^n-u^{n-1}\Vert^2_{T,\infty}\right] \\
&\qquad \times \int_0^tdr \int_D dz\, H^2(r,x,z),
\end{align*}
and, from \eqref{ch1'-s7.4}  and \eqref{C1}, the right-hand side is finite. We can therefore apply Lemma \ref{ch1'-s7-l1}
to $A(t,x):= \mathcal{I}^{n}(t,x)-\mathcal{I}^{n-1}(t,x)$ to obtain (see \eqref{ch1'-s7.l1-2000}) for any $t\in[0,T]$ and $p$ large enough,
\begin{align}
\label{is}
E\left[\Vert\mathcal{I}^{n}-\mathcal{I}^{n-1}\Vert_{t,\infty}^p\right]&
\le C(p,T,D)\int_0^t dr \, E\left[\Vert Z\Vert^p_{r,\infty}\right] \notag\\
 &\le \tilde C(p,T,D)\int_0^t dr \, E\left[\Vert u^n-u^{n-1}\Vert^p_{r,\infty}\right],
\end{align}
where we have used \eqref{ch1'-s7.148}.

With the same approach, relying on Lemma \ref{ch1'-s7-l2}, we find that for $t\in[0,T]$ and $p$ large enough,
\beq
\label{jotes}
E\left[\Vert\mathcal{J}^{n}-\mathcal{J}^{n-1}\Vert_{t,\infty}^p\right] \le \tilde C(p,T,D)\int_0^t dr \, E\left[\Vert u^n-u^{n-1}\Vert^p_{r,\infty}\right].
\eeq

Set
\beqn
M_p^n(t) = E\left[\Vert u^n-u^{n-1}\Vert^p_{t,\infty}\right], \qquad t\in[0,T].
\eeqn
 From \eqref{is}, \eqref{jotes}, we obtain
\beqn
M_p^n(t) \le c_1\int_0^t dr\, M_p^{n-1}(r),
\eeqn
and
with the Gronwall-type Lemma \ref{A3-l01}(b), we deduce that
\beq
\label{lip-2}
\sum_{n=0}^\infty \left(E\left[\Vert u^n-u^{n-1}\Vert^p_{T,\infty}\right]\right)^{\frac{1}{p}}<\infty.
\eeq
This implies that there exists a random field $u=(u(t,x),\, (t,x)\in[0,T]\times D)$ such that
\beq
\label{lip-5}
\lim_{n\to\infty}E\left[\Vert u^n-u\Vert^p_{T,\infty}\right] =0.
\eeq
Passing to a subsequence, again denoted $(u^n)$, we see that a.s., $(t,x)\mapsto u^n(t,x)$ converges to $(t,x)\mapsto u(t,x)$ uniformly on $[0,T]\times D$, therefore $u$ has
 uniformly continuous sample paths, extends continuously to $[0, T] \times \bar D$, and is adapted. In fact, for $(t,x) \in [0, T] \times \bar D$,
\beqn
u(t,x) = I_0(t,x) + \sum_{n=0}^\infty \left(u^{n+1}(t,x) - u^n(t,x)\right),
\eeqn
so by \eqref{lip-2}, $u$ satisfies \eqref{ch1'-s7.3}.
\medskip

\noindent{\em Step 3.} We show that the process $u$ satisfies equation \eqref{ch1'-s7.1}. Define
\begin{align*}
&Z_1(s,y) = \sigma(s,y,u), \quad\ \ Z_2(s,y) = b(s,y,u)\\
&Z_1^n(s,y) = \sigma(s,y,u^n), \quad Z_2^n(s,y) = b(s,y,u^n),
\end{align*}
and for $(t,x)\in[0,T]\times D$, let
\begin{align}
\label{lip-5-bis}
\mathcal{I}(t,x) & = \int_0^t \int_D \Gamma(t,x;r,z) \sigma(r,z,u)\,W(dr,dz),\notag\\
\mathcal{J}(t,x) & =  \int_0^t dr \int_D dz\, \Gamma(t,x;r,z) b(r,z,u).
\end{align}
By \eqref{ch1'-s7.3} and assumption $({\bf H_L})$(ix), $\mathcal{I}$ and $\mathcal{J}$ are well-defined, and by Lemmas \ref{ch1'-s7-l1} and \ref{ch1'-s7-l2}, they have Hölder continuous versions (extended to $\bar D$). Then, for $(t,x) \in [0, T] \times D$,
\beqn
\mathcal{I}^n(t,x) - \mathcal{I}(t,x) = \int_0^t \int_D  \Gamma(t,x;r,z)\left(Z_1^n(s,y) - Z_1(s,y)\right)\, W(ds,dy),
\eeqn
so by Lemma \ref{ch1'-s7-l1} part 3., for $p$ large enough,
\begin{align}
\label{lip-3}
E\left[\Vert \mathcal{I}^n-\mathcal{I}\Vert^p_{T,\infty}\right] & \le \tilde C\int_0^t dr\, E\left[\Vert Z_1^n-Z_1\Vert^p_{T,\infty}\right] \notag\\
&\le C E\left[\Vert u^n - u\Vert^p_{T,\infty}\right],
\end{align}
where we have used (viii-global) of Assumption (${\bf h_L}$).

Similarly, with Lemma \ref{ch1'-s7-l2} part 3., for $p$ large enough, we obtain
\beq
\label{lip-4}
E\left[\Vert \mathcal{J}^n-\mathcal{J}\Vert^p_{r,\infty}\right] \le C E\left[\Vert u^n - u\Vert^p_{r,\infty}\right].
\eeq
With the notation introduced in Step 1, for $(t,x)\in[0,T]\times D$, we have
\beqn
u^{n+1}(t,x) = I_0(t,x) + \mathcal{I}^n(t,x) +\mathcal{J}^n(t,x).
\eeqn
The left-hand side converges to $u(t,x)$ in $L^p(\Omega)$, while from \eqref{lip-3},  \eqref{lip-4}, and \eqref{lip-5}, the right-hand side converges to $I_0(t,x) +\mathcal{I}(t,x) +\mathcal{J}(t,x)$. Therefore, for each $(t,x)\in[0,T]\times D$,
\beq
\label{lip-4-bis}
u(t,x) = I_0(t,x) +\mathcal{I}(t,x) +\mathcal{J}(t,x)\quad  a.s.,
\eeq
that is, equation \eqref{ch1'-s7.1} holds (and $u$ is a solution to \eqref{ch1'-s5.0-bis}).
\medskip

\noindent{\em Step 4. Uniqueness.}
Let $(u(t,x),\, (t,x)\in[0,T]\times D)$ and  $(\bar u(t,x),\, (t,x)\in[0,T]\times D)$ be two adapted processes with continuous sample paths satisfying \eqref{ch1'-s7.3} with $p=2$.
By the same arguments used to obtain \eqref{is} and \eqref{jotes} with $u^n$, $u^{n-1}$ there replaced by $u$, $\bar u$, respectively, we obtain
\beqn
E\left[\Vert u - \bar u\Vert_{t,\infty}^2\right] \le c(T,D) \int_0^t dr\,  E\left[\Vert u-\bar u\Vert_{r,\infty}^2\right].
\eeqn
Applying the classical version of Gronwall's Lemma (Lemma \ref{lemC.1.1}) with $z\equiv 0$ there to the function $[0,T]\ni t\mapsto f(t):= E\left[\Vert u - \bar u\Vert_{t,\infty}^2\right]$
yields
\beqn
E\left[\Vert u - \bar u\Vert_{t,\infty}^2\right] = 0,
\eeqn
for all $t\in[0,T]$. Since $u$ and $\bar u$ have continuous sample paths, this implies that $u$ and $\bar u$ are indistinguishable, and completes the proof of (a).
\smallskip

   (b) We now address the question of H\"older continuity of the sample paths of the solution. Consider the identity \eqref{lip-4-bis} with $\mathcal{I}(t,x)$ and $\mathcal{J}(t,x)$ given in \eqref{lip-5-bis}.
Set $Z(r,z): =\sigma(r,z,u)$.
We have
\beqn
|Z(r,z)|\le C(T,p)\left(1+\Vert u\Vert_{r,\infty}\right), \ r\in[0,T],
\eeqn
because of (ix) in (${\bf h_L}$).

Consider points $(t,x),\, (s,y)\in[0,T]\times D$ with $0\le s\le t\le T$, and $p\in[2,\infty[$ large enough.
Then, by applying Lemma \ref{ch1'-s7-l1} (see \eqref{ch1-s7.l1-2bis}), we have
\begin{align}
\label{ch1'-s7.152}
E\left[|\mathcal{I}(t,x)- \mathcal{I}(s,y)|^p\right] &\le C(p,T,D) \left(|t-s|^{\alpha_1} + |x-y|^{\alpha_2}\right)^p \notag\\
&\quad\times \int_0^t dr\, \left(1+E\left[\Vert u\Vert_{r,\infty}^p\right]\right)\notag\\
&\le C(p,T,D) \left(|t-s|^{\alpha_1} + |x-y|^{\alpha_2}\right)^p,
\end{align}
where, in the last inequality, we have used \eqref{ch1'-s7.3}.

In a similar way, by applying Lemma \ref{ch1'-s7-l2} to $Z(r,z)= b(r,z,u)$, we obtain
\beq
\label{ch1'-s7.153}
E\left[|\mathcal{J}(t,x)- \mathcal{J}(s,y)|^p\right] \le C(p,T,D)\left(|t-s|^{\beta_1} + |x-y|^{\beta_2}\right)^p.
\eeq
With \eqref{ch1'-s7.152}, \eqref{ch1'-s7.153}, along with the H\"older continuity assumption on $I_0$ and since $D$ is bounded, we obtain \eqref{ch1'-s7.2}.

The claim about H\"older continuity follows from Theorem \ref{ch1'-s7-t2}. Indeed, we have just proved that the process
$(v(t,x):=\mathcal{I}(t,x)+\mathcal{J}(t,x), (t,x)\in[0,T]\times D)$ satisfies the condition \eqref{ch1'-s7.18} of that theorem with $I=[0,T]$. Hence, there is a version of this process, still denoted $(v(t,x))$ but extended to $\bar D$, that satisfies \eqref{ch1'-s7.19}, for all $p$ large enough and $\alpha$ arbitrarily close to $1$. Together with the H\"older continuity assumption on $I_0$, this
ends the proof of the theorem.
\qed
\smallskip

\subsection{Uniqueness among continuous solutions}
\label{ch1'-ss7.2}

In this section, we address the question of uniqueness of solutions assuming that the coefficients of the SPDE \eqref{ch1'-s5.0-bis} are locally Lipschitz continuous. More specifically, we assume the condition $({\bf h_L})$ (viii-local), and we prove that if there exist solutions to \eqref{ch1'-s5.0-bis}, in the sense given in Definition \ref{ch1'-s7-d1}, then they must be indistinguishable.

\begin{thm}
\label{ch1'-s7-t3}
Consider solutions $u^{(1)}$ and $u^{(2)}$, respectively, to SPDEs as in  \eqref{ch1'-s5.0-bis} or \eqref{ch1'-s7.1}
 with the same initial and boundary conditions, and, respectively, coefficients $\sigma^{(1)}$, $b^{(1)}$ and $\sigma^{(2)}$, $b^{(2)}$. Let ${\bf \Delta}_3$ and ${\bf \Delta}_4$ be as in \eqref{afegir(*1)}.
Suppose that for both equations, the assumptions $({\bf h_\Gamma})$, $({\bf h_I})$ and all the assumptions in $({\bf h_L})$ except possibly (viii-global) are satisfied. Let $\tau$ be a stopping time with $\tau\leq T$ a.s.
Assume that there is $M > 0$ such that for all $t \in [0, T]$, on the event $\{t \leq \tau\}$,
\beqn
\Vert u^{(1)}\Vert_{t,\infty} + \Vert u^{(2)}\Vert_{t,\infty}  \le M
\eeqn
and
\beq
\label{newdisplay}
 \sigma^{(1)}(r,z,u) = \sigma^{(2)}(r,z,u),\quad b^{(1)}(r,z,u) = b^{(2)}(r,z,u),
 \eeq
  for $0\le r\leq t$, $z \in  D$, $u \in \mathcal{C}([0,t] \times \bar D)$. Then on the event $\{t \leq \tau\}$, $u^{(1)}(s,x) = u^{(2)}(s,x)$ for all $0 \le s \le t$ and $x\in \bar D$.
\end{thm}
\begin{proof}
Define
\beqn
\bar u(t,x) = \left(u^{(1)}(t,x) - u^{(2)}(t,x)\right)1_{\{t\le \tau\}}, \quad (t,x)\in[0,T]\times \bar D.
\eeqn
Clearly, $\Vert \bar u \Vert_{t,\infty}\le M$. Further, by the local property in $\Omega$ of the stochastic integral (see Lemma \ref{ch1'-llp}), for all $(t,x)\in[0,T]\times D$,
\begin{align}
\label{ch1'-s7.156}
\bar u(t,x)&= 1_{\{t\le \tau\}}\int_0^{t} \int_D \Gamma\left(t,x;r,z\right) 1_{\{r\le \tau\}}Z_1(r,z)\, W(dr,dz)\notag\\
& \qquad + 1_{\{t\le \tau\}}\int_0^{t} dr \int_D dz\, \Gamma\left(t,x;r,z\right) 1_{\{r\le \tau\}}Z_2(r,z),
\end{align}
where
\begin{align*}
Z_1(r,z)&= \sigma^{(1)}(r,z,u^{(1)})-\sigma^{(2)}(r,z,u^{(2)}),\\
Z_2(r,z)&= b^{(1)}(r,z,u^{(1)})-b^{(2)}(r,z,u^{(2)}).
\end{align*}
Observe that, by (viii-local) and \eqref{newdisplay}, on $\{t \leq \tau\}$ for $0 \leq r \leq t$ and $i = 1, 2$,
\beq
\label{z}
|Z_i(r,z)|\le C(T,M)\Vert u^{(1)}-u^{(2)}\Vert_{r,\infty}= C(T,M)\Vert \bar u\Vert_{r,\infty} \le C(T,M)M.
\eeq

The stochastic integral in \eqref{ch1'-s7.156} is a well-defined random variable in $L^2(\Omega)$. Indeed, for $t\in[0,T]$
and because of \eqref{C1} with $\gamma=2$ and \eqref{z},
\begin{align*}
&E\left[\int_0^{t} dr \int_D dz \, |\Gamma(t,x;r,z) 1_{\{r\le \tau\}}Z_1(r,z)|^2\right] \\
 &\qquad \le C^2(T,M)M^2 \int_0^{t} dr \int_D dz\,  \Gamma^2(t,x;r,z)\\
  &\qquad \le C^2(T,M)M^2 \int_0^T ds \int_D dz\, H^2(r,x,z) <\infty.
\end{align*}
 With similar arguments, we see that the last integral in \eqref{ch1'-s7.156} is also well-defined.

Apply Claim 3 of Lemmas \ref{ch1'-s7-l1} and \ref{ch1'-s7-l2} to the processes
\begin{align*}
A(t,x)&= \int_0^{t} \int_D \Gamma(t,x;r,z) 1_{\{r\le \tau\}}Z_1(r,z)\, W(dr,dz),\\
B(t,x)&= \int_0^{t} dr \int_D dz\, \Gamma(t,x;r,z) 1_{\{r\le \tau\}}Z_2(r,z),
\end{align*}
respectively, to obtain
\begin{align*}
&E\left[\Vert A\Vert_{t,\infty}^p\right] + E\left[\Vert B\Vert_{t,\infty}^p\right] \\
&\qquad \quad \le C(p,T,D) \int_0^{t} dr\,
E \left[\Vert 1_{\{\cdot\le \tau\}}Z_1\Vert_{r,\infty}^p + \Vert 1_{\{\cdot\le \tau\}}Z_2\Vert_{r,\infty}^p\right],
\end{align*}
for any $p>1$ large enough. This implies
\begin{align*}
E\left[\Vert\bar u\Vert_{t,\infty}^p\right]
& \le C(p,T,D) \int_0^t  dr\,  E \left[\Vert 1_{\{\cdot\le \tau\}}Z_1\Vert_{r,\infty}^p + \Vert 1_{\{\cdot\le \tau\}}Z_2\Vert_{r,\infty}^p\right] \\
& \le C(p,T,D,M) \int_0^t dr\, E\left[\Vert \bar u\Vert_{r,\infty}^p\right],
\end{align*}
where we have used \eqref{z}.
We notice that the constants in the two displays above also depend on $\alpha_i$ and $\beta_i$, $i=1,2$, but this is not relevant in the arguments.
Observe that
the integral on the right-hand side of the last display above is finite and even bounded by $TC^p(T,M)M^p$. Therefore, the classical version of Gronwall's Lemma \ref{lemC.1.1} applies and we conclude that for $t\in[0,T]$,
\beqn
E \left[\Vert \bar u\Vert_{t,\infty}^p\right] =0.
\eeqn
In particular, for $t \in [0,T]$,
\begin{align*}
& \sup_{s \in [0, t \wedge \tau] \times \bar D} \vert u^{(1)}(s,x) - u^{(2)}(s,x) \vert
= \sup_{s \in [0,t \wedge \tau] \times \bar D} \vert \bar u(s,x) \vert \\
&\qquad\qquad \qquad\qquad\qquad= \sup_{s \in [0,t] \times D} \vert \bar u(s,x) \vert
 = \Vert \bar u \Vert_{t,\infty}  = 0 \quad {\text{ a.s.}}
 \end{align*}
This ends the proof of the theorem.
\end{proof}

\subsection{Global existence with locally Lipschitz coefficients}
\label{ch1'-ss8.3}

In this section, we extend Theorem \ref{ch1'-s7-t1} to the situation where the coefficients are locally Lipschitz continuous functions.

We start with some preliminaries.
For $N\in\N$, define $\psi_N: \R  \longrightarrow \R$ by
\beqn
\psi_N(\rho) =
\begin{cases}\ \rho, & {\text{if}}\ |\rho|\le N,\\
\ N, &  {\text{if}}\ \rho>N,\\
-N, &  {\text{if}}\ \rho < -N.
\end{cases}
\eeqn
Then for any $u: [0,T]\times D \longrightarrow \R$, let
\beqn
\Psi_N(u)(t,x) = \psi_N(u(t,x))=
\begin{cases}
u(t,x), &  {\text{if}}\ |u(t,x)|\le N,\\
\ N, &  {\text{if}}\ u(t,x) >N,\\
-N, &  {\text{if}}\ u(t,x)<-N,
\end{cases}
\eeqn
$(t,x)\in[0,T]\times D$. Notice that $\Vert \Psi_N(u)\Vert_{T,\infty}\le N$.

Given a function
$g: [0,T]\times D \times \mathcal{C}([0,T]\times \bar D) \longrightarrow \R$ define
\beq
\label{g-l}
g_N(r,z,u) = g(r,z,\Psi_N(u)).
\eeq
If $g$ satisfies the condition $({\bf h_L})$(viii-local), then for any $N\in\N$, $g_N$ satisfies $({\bf h_L})$(viii-global). Indeed, for any $r\in[0,T]$, $z\in D$, $u,v \in \mathcal{C}([0,T]\times \bar D)$,
\begin{align*}
|g_N(r,z,u) - g_N(r,z,v)| & = |g(r,z,\Psi_N(u)) - g(r,z,\Psi_N(v))|\\
& \le C(N,T)\Vert \Psi_N(u) - \Psi_N(v)\Vert_{r,\infty}\\
&= C(N,T)\sup_{(s,z)\in\ [0,r]\times D} |\psi_N(u(s,z))- \psi_N(v(s,z))|.
\end{align*}
Considering all possible cases regarding the values of $|u(s,z)|$ and $|v(s,z)|$ compared with $N$, it is  easy to see that
\begin{align*}
\sup_{(s,z)\in [0,r]\times D} |\psi_N(u(s,z))- \psi_N(v(s,z))|& \le \sup_{(s,z)\in [0,r]\times D} |u(s,z) - v(s,z)|\\
&=\Vert u-v\Vert_{r,\infty}.
\end{align*}
Thus, $g_N: [0,T]\times D \times \mathcal{C}([0,T]\times \bar D) \longrightarrow \R$ satisfies the condition $({\bf h_L})$(viii-global).

\begin{thm}
\label{ch1'-s7-t4}
Let ${\bf \Delta}_3$ and ${\bf \Delta}_4$ be as in \eqref{afegir(*1)}.
We assume $({\bf h_\Gamma})$, $({\bf h_I})$ and  $({\bf h_L})$ without (viii-global). Then the SPDE \eqref{ch1'-s5.0-bis} has a  solution in the sense of Definition \ref{ch1'-s7-d1} such that for all $p>0$,
\beq
\label{ch1'-s7.159}
E\left[\Vert u\Vert_{T,\infty}^p\right] < \infty.
\eeq
Moreover, this solution is unique (up to indistinguishability).

If $(t,x) \mapsto I_0(t,x)$ is Hölder continuous jointly in $(t,x)$ with exponents $(\eta_1, \eta_2)$, then $u$ is Hölder continuous, jointly in $(t,x)$, with exponents $(\delta_1, \delta_2)$, where $\delta_i \in\, ]0, \eta_i \wedge \alpha_i \wedge \beta_i [$, $i = 1, 2$.
\end{thm}
\begin{proof}
For any $N\in\mathbb{N}$, define $\sigma_N(r,z,u)$ and $b_N(r,z,u)$ by the formula \eqref{g-l}, replacing $g$ there by $\sigma$ and $b$, respectively.
As noted above, $\sigma_N$ and $b_N$ satisfy $({\bf h_L})$(viii-global). For these coefficients, all the assumptions of Theorem \ref{ch1'-s7-t1} are satisfied. Hence, there is a solution $(u_N(t,x),\, (t,x)\in[0,T]\times \bar D)$ with continuous sample paths to the equation
\begin{align}
\label{ch1'-s7.160}
u_N(t,x) = I_0(t,x) &+ \int_0^t \int_D \Gamma(t,x;r,z) \sigma_N(r,z,u_N)\, W(dr,dz)\notag\\
&+ \int_0^t dr \int_D dz\, \Gamma(t,x;r,z) b_N(r,z,u_N),
\end{align}
and this solution is unique (up to indistinguishability).

Let
\beqn
\tau_N = \inf\{t\ge 0: \Vert u_N\Vert_{t,\infty}\ge N\}\wedge T,\quad \rho_N = \tau_N\wedge \tau_{N+1}.
\eeqn
On the event $\{t\le \rho_N\}$, $\Vert u_N\Vert_{t,\infty}+\Vert u_{N+1}\Vert_{t,\infty}\le 2N+1$, and
the coefficients of $u_N$ and $u_{N+1}$ coincide. Hence, by Theorem \ref{ch1'-s7-t3}, on that event, $u_N(s,z) = u_{N+1}(s,z)$ for all $s\in[0,t]$ and $z\in \bar D$.

Observe that on the event $\{\rho_N = \tau_{N+1}<\tau_N\}$, we have
\beqn
N>\Vert u_{N}\Vert_{\tau_{N+1},\infty} =\Vert u_{N+1}\Vert_{\tau_{N+1},\infty} =N+1,
\eeqn
which is a contradiction. Thus, $P\{\rho_N = \tau_{N+1}<\tau_N\}=0$. We deduce $\rho_N = \tau_N \le \tau_{N+1}$ a.s., that is, $(\tau_N)_N$ is an increasing sequence of stopping times and on $\{t\le \tau_N\}$,  $u_{N}(s,z) = u_{N+1}(s,z)$ for all $s\in[0,t]$ and $z\in \bar D$.

Observe also that on $\{t\le \tau_N\}$, for $r\le t$, $\sigma_N(r,z,u_N)= \sigma(r,z,u_N)$ and $b_N(r,z,u_N)= b(r,z,u_N)$.
Thus, by the local property of the stochastic integral, on $\{t\le \tau_N\}$,
\begin{align}
\label{ch1'-s7.160-bis}
u_N(t,x)& = I_0(t,x) + \int_0^{t} \int_D \Gamma(t,x;r,z) \sigma(r,z,u_N)\, W(dr,dz)\notag\\
& \qquad + \int_0^{t} dr \int_D dz\, \Gamma(t,x;r,z) b(r,z,u_N).
\end{align}
Apply \eqref{ch1'-s7.l1-2000} in Lemma \ref{ch1'-s7-l1} and \eqref{ch1'-s7.l1-201} in Lemma \ref{ch1'-s7-l2} to
\begin{align*}
A(t,x):= & \int_0^{t} \int_D \Gamma(t,x;r,z) \sigma(r,z,u_N)1_{\{r\le \tau_N\}}\, W(dr,dz),\\
B(t,x):= & \int_0^{t}dr  \int_D dz\, \Gamma(t,x;r,z) b(r,z,u_N)1_{\{r\le \tau_N\}},
\end{align*}
respectively,
and note that, by (ix) in $({\bf h_L})$, for $0\le r\le t$,
\begin{align*}
\Vert\sigma(\cdot,\ast,u_N)1_{\{\cdot\le \tau_N\}}\Vert_{r,\infty}+\Vert b(\cdot,\ast,u_N)1_{\{\cdot\le \tau_N\}}\Vert_{r,\infty}
\le C_2(T)
\left(1+\Vert u_N1_{\{\cdot\le \tau_N\}}\Vert_{r,\infty}\right).
\end{align*}
We deduce that for $p>1$ sufficiently large,
\begin{align*}
E\left[\Vert u_N 1_{\{\cdot\le \tau_N\}}\Vert^p_{t,\infty}\right]
&\le C(p,T,D)\left( \Vert I_0\Vert_{T,\infty}^p + \int_0^{t} dr \left( 1+ E\left[\Vert u_N 1_{\{\cdot \le \tau_N\}}\Vert_{r,\infty}^p\right]\right)\right) \\
&\le C(p,T,D,I_0) + C(p,T,D) \int_0^t dr\, E\left[\Vert u_N1_{\{\cdot \le \tau_N\}}\Vert_{r,\infty}^p\right].
\end{align*}
In fact the constants in this display also depend on $\alpha_i$ and $\beta_i$, $i=1,2$, but this is not relevant in the arguments.
Since by the definition of $\tau_N$, $\Vert u_N1_{\{\cdot\le \tau_N\}}\Vert_{r,\infty}\le N<\infty$, we conclude from Gronwall's Lemma \ref{lemC.1.1}, applied to
$f(t):=E\left[\Vert u_N1_{\{\cdot\le \tau_N\}}\Vert^p_{t,\infty}\right]$, that
\beq
\label{ch1'-s7.161}
 E\left[\Vert u_N 1_{\{\cdot\le \tau_N\}}\Vert^p_{T,\infty}\right] \le \tilde C(p,T,D,I_0),
\eeq
where the right-hand side does not depend on $N$.

As a consequence of \eqref{ch1'-s7.161}, we obtain that $\tau_N \uparrow T$ a.s., and even that
\beq
\label{ch1'-s7.162}
P\left\{\lim_{N\to\infty}\tau_N=T\right\}=1.
\eeq
 Indeed, we have already proved that the sequence $(\tau_N)_N$ is increasing. 
By the definition of $\tau_N$
and by continuity of $u_N$, for $N > \sup_{x \in D} I_0(0, x)$,
on $\{\tau_N<T\}$, we have
\beqn
\sup_{(t,x)\in[0,T]\times D}\left(\vert u_N(t,x)|^p1_{\{t\le \tau_N\}}\right)=N^p.
\eeqn
Along with
\eqref{ch1'-s7.161}, this yields
\begin{align*}
N^p P\left\{\tau_N<T\right\}&= E\left[1_{\{\tau_N<T\}} \left\Vert u_N 1_{\{\cdot\le \tau_N\}}\right\Vert^p_{T,\infty}\right] \\
&\le \tilde C(p,T,D,I_0).
\end{align*}
Consequently, $\lim_{N\to\infty} P\left\{\tau_N<T\right\}=0$, which proves \eqref{ch1'-s7.162}.

We now define $(u(t,x),\, (t,x)\in[0,T]\times \bar D)$ by setting
\beq
\label{ch1'-s7.163}
u(t,x) = u_N(t,x) \ {\text{on}} \ \{t\le \tau_N\}.
\eeq
As observed near the beginning of the proof, we see that on $\{t\le \tau_N\}$, $u(s,x) = u_M(s,x)$ for all $0\le s\le t$, $x\in \bar D$, and $M\ge N$. Since $\{t\le \tau_N\}\uparrow \Omega$ a.s., $u(t,x)$ is a well-defined random variable.

The process $u$ defined in \eqref{ch1'-s7.163} satisfies
\eqref{ch1'-s7.159}. Indeed, set
\beqn
X_N = \sup_{(t,x)\in[0,T]\times D}\left(1_{\{t\le \tau_N\}}\left\vert u(t,x)\right\vert^p\right).
\eeqn
By \eqref{ch1'-s7.161} and \eqref{ch1'-s7.163}, $E\left[X_N\right] \le \tilde C(p,T,D,I_0)$.
 From the monotone convergence theorem, we deduce that $E[X] \le  \tilde C(p,T,D.I_0)$, where
\beqn
X=\lim_{N\to\infty} X_N = \sup_{(t,x)\in[0,T]\times D}|u(t,x)|^p,
\eeqn
and the limit is in the sense of almost sure convergence.
This establishes \eqref{ch1'-s7.159}.

Finally, we check that $u$ defined in \eqref{ch1'-s7.163} satisfies \eqref{ch1'-s7.1}.
Indeed, from \eqref{ch1'-s7.160-bis} and \eqref{ch1'-s7.163}, we see that on $\{t\le \tau_N\}$, for $x\in D$,
\begin{align*}
u(t,x)&=I_0(t,x) + \int_0^t \int_D \Gamma(t,x;r,z) \sigma(r,z,u)\, W(dr,dz)\\
&+ \int_0^t dr \int_D dz\, \Gamma(t,x;r,z)b(r,z,u),
\end{align*}
and we note that because of \eqref{ch1'-s7.159} and (ix) in $({\bf h_L})$, the stochastic and the pathwise integrals in the above equation are well-defined random variables in $L^2(\Omega)$.
Since $\{t\le \tau_N\}\uparrow \Omega$ a.s., $u$ solves \eqref{ch1'-s7.1}.

Uniqueness follows from the definition of $(u(t,x),\, (t,x)\in[0,T]\times \bar D)$ (see \eqref{ch1'-s7.163}) and Theorem \ref{ch1'-s7-t3}.

 By Theorem \ref{ch1'-s7-t1}, if $I_0$ is Hölder continuous with exponents $(\eta_1, \eta_2)$, then each $u_n$ is Hölder continuous with exponents $(\delta_1, \delta_2)$, where $\delta_i \in\,  ]0, \eta_i \wedge \alpha_i \wedge \beta_i [$, $i = 1, 2$.  Because of \eqref{ch1'-s7.163} and \eqref{ch1'-s7.162}, this property is inherited by $u$.

The proof of the theorem is complete.
\end{proof}
\subsection{Examples}
\label{ch4-sec5.5-ex}
In this subsection, we give two examples in which Theorem \ref{ch1'-s7-t4} applies. These are the stochastic heat and wave equations on the interval $D =\, ]0, L[$, with suitable initial conditions, vanishing boundary conditons and coefficients $\sigma$ and $b$ as in Section \ref{ch1'-ss7.1}.  
\medskip

\noindent{\em Stochastic heat equation on $]0, L[$}
\medskip

We consider the solution $u = (u(t,x),\, (t,x) \in [0, T] \times D)$ to the stochastic heat equation on $D =\, ]0, L[$ with vanishing Dirichlet (resp. Neumann) boundary conditions and initial condition $u_0$, that is, $u$ is a solution of \eqref{ch1'-s5.0-bis} (and \eqref{ch1'-s7.1}) with $\cL = \frac{\partial}{\partial t} - \frac{\partial^2}{\partial x^2}$, $\Gamma(t,x; s,y) = G_L(t-s; x, y)$ given by \eqref{ch1'.600} (respectively \eqref{1'.400}), coefficients $\sigma$ and $b$ as in Section \ref{ch1'-ss7.1} and $I_0$ given by  \eqref{cor1.0}(respectively \eqref{cor1.0-N}) for some function $u_0$. We assume that $u_0 \in \cC([0, L])$ with $u_0(0) = u_0(L) = 0$ (respectively $u_0 \in \cC([0, L])$ with no other condition), so that assumption $({\bf h_I})$ holds.

\begin{thm}
\label{ch4-sec5.5-ex-Theorem 1}
 Assume $({\bf h_L})$ without (viii-global).

   (a)\ The stochastic heat equation on $]0, L[$, with the two types of boundary conditions and initial condition $u_0$ just mentioned, has a solution in the sense of Definition \ref{ch1'-s7-d1} such that for all $p > 0$,
   \beqn
    E\left[\sup_{(t, x) \in [0, T] \times [0, L]} \vert u(t, x) \vert^p\right] < \infty,
    \eeqn
and this solution is unique up to indistinguishability.

   (b)\ If, in addition, $u_0 \in \cC_0^\eta([0, L])$ (respectively $\cC^\eta([0, L]))$ for some $\eta \in ]0, 1]$, then u is Hölder continuous, jointly in $(t,x)$, with exponents $(\delta_1, \delta_2)$, where $\delta_1 \in\, ]0, \frac{1}{4} \wedge \frac{\eta}{2}[$ and $\delta_2 \in\, ]0, \frac{1}{2} \wedge \eta [$.
   \end{thm}

\begin{proof}
We only consider the case of vanishing Dirichlet boundary conditions, since the other case is similar.

  (a)\ Let $\ep_0 \in [0, 1[$ and set $\gamma = 2 + \ep_0$. Let
    $\alpha_1 = \half (\frac{3}{\gamma} - 1), \alpha_2 = \frac{3}{\gamma} - 1$.
Notice that for $\ep_0 \in\, ]0, 1[$, we have $\gamma \in\, ]2, 3[, \alpha_1 \in\, ]0, \frac{1}{4}[$ and $\alpha_2 \in\, ]0, \half[$, and $\alpha_1$ and $\alpha_2$ are near $\frac{1}{4}$ and $\half$ when $ \ep_0$ is near $0$.
By Theorem \ref{ch1'-s7-t4}, it suffices to check assumption $({\bf h_\Gamma})$, with
\beq
\label{ch4-sec5.5-ex-Theorem 1(*a)}
   \Delta_3(t, x; s, y) = \vert t - s \vert^{\alpha_1} + \vert x - y \vert^{\alpha_2},\quad   \Delta_4(t, x; s, y) = \Delta_3(t, x; s, y).   
   \eeq
    In Section \ref{ch1'-ss6.1}, we have already checked assumptions $({\bf h_\Gamma})$(i) and (ii), with $H(t-s, x, y) = \Gamma(t-s, x - y)$ and $\Gamma$ is the heat kernel defined in \eqref{ch1'-s6.1}. By Lemma \ref{app2-l2} with $k = 1$, we see that $({\bf h_\Gamma})$(iii) is satisfied  with $\gamma$ as above  since $\ep_0 \in [0, 1[$. For $({\bf h_\Gamma})$(iv) , we use Lemma \ref{app2-r1} to see that
    \beq
    \label{ch4-sec5.5-ex-Theorem 1(*b)}
   \int_0^T dr \int_0^L dz\, \vert G_L(t - r; x, z) - G_L(s - r; y, z) \vert^\gamma \leq C_{\gamma, T} (\vert t - s \vert^{\half (3 - \gamma)} + \vert x - y \vert^{3 - \gamma}). 
   \eeq
Because $\half (3 - \gamma) = \alpha_1 \gamma$ and $3 - \gamma = \alpha_2 \gamma$, we see that \eqref{ch4-sec5.5-ex-Theorem 1(*b)} is bounded above by $c_{\gamma, T}  (\Delta_3(t, x; s, y))^\gamma$, with our choice of $\Delta_3$ in \eqref{ch4-sec5.5-ex-Theorem 1(*a)} and any $\ep_0 \in [0, 1[$, so $({\bf h_\Gamma})$(iv) holds. As mentioned in Remark \ref{ch1'-s7-r1}, this implies $({\bf h_\Gamma})$(v) with the same $\ep_0$ and $\Delta_4 = \Delta_3$.

   (b)\  If $u_0 \in \cC_0^\eta([0, L])$ (respectively $\cC^\eta([0, L])$) for some $\eta \in\, ]0, 1]$, then $I_0$ is Hölder continuous, jointly in $(t,x)$, with exponents $(\eta/2, \eta)$ by Lemma \ref{cor1} (respectively Lemma \ref{cor1-N}), so we apply Theorem \ref{ch1'-s7-t4} to conclude that $u$ is Hölder continuous, jointly in $(t,x)$, with exponents $(\delta_1, \delta_2)$, where $\delta_i \in\, ]0, \eta_i \wedge \alpha_i[$,  $i = 1, 2$. Since we can choose $\ep_0 > 0$ arbitrarily close to $0$, $\alpha_1$ can be taken arbitrarily close to $\frac{1}{4}$ and $\alpha_2$ to $\half$, so we obtain the conclusion.
    \end{proof}
   \medskip

   \noindent{\em Stochastic wave equation on $]0, L[$}
\medskip

We consider the solution $u = (u(t,x),\, (t,x) \in [0, T] \times D)$ to the stochastic wave equation on $D =\, ]0, L[$ with vanishing Dirichlet boundary conditions and initial conditions $f$ and $g$, $f \in \cC([0, L])$, with $f(0) = f(L) = 0$ and $g \in L^1([0, L])$, that is, $u$ is a solution of \eqref{ch1'-s5.0-bis} (and \eqref{ch1'-s7.1}) with $\cL = \frac{\partial^2}{\partial t^2} - \frac{\partial^2}{\partial x^2}$, $\Gamma(t,x; s,y) = G_L(t-s; x, y)$ given by \eqref{wave-bc1-1'} (and \eqref{wave-bc1-100'}), coefficients $\sigma$ and $b$ as in Section \ref{ch1'-ss7.1} and with $I_0$ given by \eqref{i0-tris}, so that assumption $({\bf h_I})$ holds.

\begin{thm}
\label{ch4-sec5.5-ex-Theorem 2}
Assume $({\bf h_L})$ without (viii-global).

   (a)\ The stochastic wave equation on $]0, L[$, with vanishing Dirichlet boundary conditions and initial conditions $f$ and $g$ as above, has a solution in the sense of Definition \ref{ch1'-s7-d1} such that for all $p > 0$,
   \beqn
    E\left[\sup_{(t, x) \in [0, T] \times [0, L]} \vert u(t, x) \vert^p\right] < \infty,
    \eeqn
and this solution is unique up to indistinguishability.

   (b)\ If, in addition, $f \in \cC_0^\eta([0, L])$ for some $\eta \in\, ]0, 1]$ and $g \in \cC([0, L])$, then $u$ is Hölder continuous, jointly in $(t,x)$, with exponent $\delta$, for any $\delta \in\, ]0, \frac{1}{2} \wedge \eta [$.
\end{thm}
\begin{proof}
(a)\  By Theorem \ref{ch1'-s7-t4}, it suffices to check assumption $({\bf h_\Gamma})$, with
\beq
\label{eqdeltas0}
   \Delta_3(t, x; s, y) = \vert t - s \vert^{\frac{1}{\gamma}} + \vert x - y \vert^{\frac{1}{\gamma}},\quad   \Delta_4(t, x; s, y) = \Delta_3(t, x; s, y)
   \eeq  
and  $\gamma  = 2 + \ep_0$  with $\ep_0 > 0$.
   In Section \ref{ch1'-ss6.2}, we have already checked  assumptions $({\bf h_\Gamma})$(i) and (ii) with $H(t-s, x, y) = \vert G_L(t-s, x, y)\vert$. For $({\bf h_\Gamma})$(iii), we note that for any $\ep_0 > 0$ and $\gamma = 2 + \ep_0$, $H^\gamma(t-s, x, y) \leq H^2(t-s, x, y)$ because $H(t-s, x, y) \in \left\{0, \half\right\}$ according to \eqref{wave-bc1-100'} (see also Figure \ref{fig3.3}). Therefore, $({\bf h_\Gamma)}$(iii) follows from the fact that $H$ satisfies $({\bf H_\Gamma})$(iiia), as mentioned prior to Theorem \ref{recapwave}. For $({\bf h_\Gamma})$(iv), we observe that for any $\ep_0 > 0$ and $\gamma = 2 + \ep_0$,
   \beq
   \label{eqdeltas1}
    \vert G_L(t-r; x, z) - G_L(s-r; x, z) \vert^\gamma \leq \vert G_L(t-r; x, z) - G_L(s-r; x, z) \vert^2,   
    \eeq
since $\vert G_L(t-r; x, z) - G_L(s-r; x, z) \vert$  takes values in $\{0, \half, 1 \}$. By \eqref{rdeB.5.5a} in Lemma \ref{app2-5-l-w1}, we see that $({\bf h_\Gamma})$(iv) is satisfied with $\Delta_3$ as in \eqref{eqdeltas0}. As mentioned in Remark \ref{ch1'-s7-r1}, this implies $({\bf h_\Gamma})$(v) with the same $\ep_0$ and $\Delta_4 = \Delta_3$.

    (b)\, By Lemma \ref{ch1'-ss2.3-hi}, if $u_0 \in \cC_0^\eta([0, L])$ and $g \in \cC([0, L])$, then $I_0$ is Hölder continuous, jointly in $(t,x)$, with exponent $\eta$. Since we can choose $\ep_0 > 0$ arbitrarily small, hence $\gamma > 2$ arbitrarily close to $2$, the claimed Hölder continuity of $u$ follows from Theorem \ref{ch1'-s7-t4}.
    This completes the proof of Theorem \ref{ch4-sec5.5-ex-Theorem 2}.
\end{proof}


\section{Notes on Chapter \ref{ch1'-s5}}
\label{ch4-notes}

Existence of solutions to SDEs and SPDEs is most often proved via a Picard iteration scheme. 
An early example can be found  in \cite[Theorem 3.2]{walsh}, where a nonlinear stochastic heat equation driven by space-time white noise is studied. In this chapter, we consider a general setting suitable to the study of a large class of SPDEs and we also use Picard iterations to establish the main theorem (Theorem \ref{ch1'-s5.t1}) on existence and uniqueness of random field solutions. The conditions on the partial differential operator that defines the SPDE are gathered in hypothesis $({\bf H_\Gamma})$ and are expressed in terms of the corresponding fundamental solution or the Green's function. A similar strategy is used in \cite{dalang} in the context of non-linear SPDEs driven by a noise that is white in time and spatially homogeneous (see also \cite{sanz-sole-2005}).

The regularity of the sample paths of the random field solutions is obtained via the Kolmogorov continuity theorem (see Appendix \ref{app1-3}). However, as in the study of existence and uniqueness of solutions,
the strategy in Section \ref{ch1'-ss5.2} is to highlight the role that the regularity of the fundamental solution/Green's function plays in the regularity properties of the solutions.

There are many results on the stochastic heat and wave equations in the SPDE literature. Going back to the origins, we mention \cite{walsh} for examples of heat equations, and \cite{cabanya-1970}, \cite{cn88} and \cite{walsh} for wave equations. Section \ref{ch1'-ss7.3} on a fractional stochastic heat equation in based on \cite{chendalang2015} and extends the study initiated in \cite{dd2005}.
The approximation result of Section \ref{ch1'-s60} is based on \cite[Lemma 2.1]{donati-martin-pardoux-1993}.
Section \ref{ch1'-s7} expands on \cite[Section 3]{donati-martin-pardoux-1993}, where the spatial domain $D\subset \re$ is assumed to be bounded. The recent work \cite{f-k-n-2024} deals with the case $D=\re$.

Integral equations such as \eqref{ch1'-s5.1} are a natural way of formulating SPDEs obtained from PDEs by adding an external forcing noise (in the form of space-time white noise), and in this case, $\Gamma$ is the fundamental solution or the Green's function of a partial differential operator $\cL$.
However, the results presented here are valid for generic {\em kernels} $\Gamma$ as long as they satisfy suitable assumptions, such as $(\bf{H_\Gamma})$, \eqref{ch1'-s5.12}, \eqref{ch1'-s5.14} or $(\bf{h_\Gamma})$.
Certain $d$-dimensional random fields meant to describe the dynamics of stochastic turbulence, risk management, tumour growth, etc., known as {\em ambit fields}\index{ambit field}\index{field!ambit} and introduced by O.~E.~Barndorff-Nielsen and coauthors (see e.g.~\cite{bn2007}), take a form close to \eqref{ch1'-s5.1}. For more information about this topic and its connections with SPDEs, we refer to \cite{bn2011} and \cite{bn2018}.


\chapter[Asymptotic behaviors in the stochastic heat equation]
{Asymptotic behaviors in the stochastic heat equation}
\label{ch6}
\pagestyle{myheadings}
\markboth{R.C.~Dalang and M.~Sanz-Sol\'e}{Long-time behavior}

We devote this chapter to the study of long-time behavior of solutions
to SPDEs, building on ideas that originate in deterministic dynamical systems and in the theory of Markov processes. The main examples are linear and nonlinear stochastic heat equations, as presented in Chapters \ref{chapter1'} and \ref{ch1'-s5},
respectively. In Section \ref{ch5-added-1}, we prove that for the three instances
of stochastic heat equations (on a bounded interval of $\R$ with either Dirichlet
or Neumann vanishing boundary conditions, and on $\R$), the random field solution
possesses the Markov and strong Markov properties. Then, for the linear stochastic heat equation
with vanishing Dirichlet boundary conditions, we identify an invariant probability measure, establish that it is a limit measure, and we prove a property of recurrence. We also discuss the case of Neumann boundary conditions, where there is an invariant measure but it is not finite, and an SPDE on $\rek$ with a fractional Laplacian. In Section \ref{ch5-ss-in-rev}, we recall the notion of reversible
measure and its connections with invariant measures. After a discussion of the finite-dimensional case, we identify reversible probability measures for a class of linear stochastic heat equations and for a class of nonlinear stochastic heat equations with additive noise and gradient-type nonlinear
drift. We show that these reversible measures can be identified with the conditional distribution of the solution $(X(x),\, x \in [0, 1])$ to a related SDE, given that $X(1) = 0$ (that is, a bridge process). In Section \ref{new-5.5} we prove that the law of the solution to the nonlinear stochastic heat equation with additive noise and gradient-type drift converges (under certain conditions) as $t \to \infty$ to the reversible measure, which is the unique invariant measure, and we show in Section \ref{new-5.5.2} that this process is irreducible. Finally, in Section \ref{new-5.5-candil} we study the approximation of the solution to a nonlinear stochastic heat equation on $\re$ by the solutions of nonlinear stochastic heat equations on $[-L,L]$ with vanishing Dirichlet boundary conditions, when $L\to\infty$.
\section{Markov and strong Markov properties}
\label{ch5-added-1}

In this section, we consider  a variation on the nonlinear stochastic heat equation (see \eqref{operator-calL} and \eqref{nbr-(*Mb)}) and we study several basic properties, such as the classical Markov and strong Markov properties, the Feller property, the existence of invariant measures and limit measures, as well as the property of recurrence. We begin by defining these notions in an abstract context and then we show that these properties hold for the solution to \eqref{nbr-(*Mb)}.

\subsection{Infinite-dimensional Markov and strong Markov processes}
\label{ch5-added-1-s1}

We start by recalling some definitions  (see e.g. \cite[Chap. III]{ry}) and fixing the setting. Let $\mathcal{S}$ be a complete separable metric space (the {\em state space}) and let $\mathcal{B}_{\mathcal{S}}$ denote the $\sigma$-field of Borel sets of $\mathcal{S}$. 

\begin{def1}\label{rd02_19d1}
A {\em transition probability}\index{transition!probability}\index{probability!transition} $P_\trt$ on $(\cs, \B_\cs)$ is a map $(g,A)\mapsto P_\trt(g,A)$ from $\cs \times {\mathcal{B}}_{\cs}$ into $\re$ such that:
\begin{description}
\item{(a)} for each $g\in\cs$, the set function $A\mapsto P_\trt(g,A)$ is a probability measure on $(\cs, \cB_{\cs})$;
\item{(b)} for each $A\in\cB_{\cs}$, the function $g\mapsto P_\trt(g,A)$ from $\cs$ into $[0,1]$ is ${\mathcal{B}}_{\cs}$-measurable.
\end{description}
\end{def1}
The set function $A\mapsto P_\trt(g,A)$ can be extended to the set of $\B_{\cs}$-measurable and bounded functions $f: \cs\longrightarrow \re$ by defining
\beq
\label{ch5-added-1-1-0}
P_\trt f(g) = \int_{\cs} P_\trt(g, d\bar g)\, f(\bar g) .
\eeq
Clearly, $P_\trt 1_A(g)=P_\trt(g,A)$. Moreover, via the 
Monotone Class Theorem \cite[Chap. I, §2, Theorem (2.2)]{ry}, the property (b) above extends to the following:

For any  $\B_{\cs}$-measurable and bounded function $f: \cs\longrightarrow \re$, the mapping $g\mapsto P_\trt f(g)$ from $\cs$ into $\re$ is
$\cB_{\cs}$-measurable.
\smallskip

\begin{def1}\label{rd02_19d2}
A (time homogeneous) {\em transition function}\index{transition!function}\index{function!transition} on $(\cs, \B_\cs)$ is a family $(P_t,\, t \in \R_+)$ of transition probabilities on $(\cs, \B_\cs)$ such that for all $s, t \in \R_+$, $g \in \cs$ and $A \in \B_\cs$,
\beq
\label{ch5-added-1-1-1}
P_{s+t}(g,A) = \int_{\cs} P_s(g,d\bar g)\, P_t(\bar g,A).
\eeq
Equality \eqref{ch5-added-1-1-1} is known as the {\em Chapman-Kolmogorov equation.}\index{Chapman-Kolmogorov equation}\index{equation!Chapman-Kolmogorov}
\end{def1}

The equality \eqref{ch5-added-1-1-1} implies that $(P_t,\, t \in \R_+)$ forms a {\em semigroup:} $P_{s + t} = P_s \circ P_t$ . Indeed, for all $g \in \cs$ and bounded measurable functions $f: \cs \to \R$, we show that $P_{s + t}f(g) = (P_s(P_t f))(g)$: by \eqref{ch5-added-1-1-1} and \eqref{ch5-added-1-1-0}, for $g \in \cs$ and $A \in \B_\cs$,
\beqn
  P_{s+t} 1_A(g) = \int_\cs  P_s(g, d\bar g)\, P_t 1_A (\bar g).
\eeqn
By linearity and the usual extension arguments, we can replace $1_A$ by any bounded and measurable function $f: \cs \to \R$, to obtain
\begin{align*}
    P_{s+t} f(g) &=  \int_\cs  P_s(g, d\bar g)\, P_t f (\bar g) = (P_s(P_t f))(g)
\end{align*}
by \eqref{ch5-added-1-1-0}.
\bigskip

The following properties of transitions functions are important in the theory of Markov processes.

\begin{def1}\label{rd02_21d1}
Let $(P_t,\, t \in \R_+)$ be a transition function on $(\cs, \B_\cs)$.

 (a) The transition function $(P_t)$ has the {\em weak Feller property}\index{Feller property!weak}\index{weak!Feller property} if for all bounded and continuous functions $f: \cs \to \R$, the mapping $g \mapsto P_t f(g)$ is continuous. 

   (b) The transition function $(P_t)$ is {\em stochastically continuous}\index{stochastically continuous}\index{continuous!stochastically} if for all bounded and continuous $f: \cs \to \R$ and $g \in \cs$,
\beqn
   \lim_{t \downarrow 0} P_t f(g) = f(g).
\eeqn
A stochastically continuous weak Feller process is called a {\em Feller process.}
\end{def1}

\noindent{\em Markov processes}
\medskip

On a probability space $(\Omega, \cF, P)$, we consider a filtration $(\cF_t,\, t\in\re_+)$ and an adapted stochastic process $X=(X_t,\, t\in\re_+)$ consisting of $\cs$-valued random variables. 

\begin{def1}
\label{rd02_19d3}
The process $X$ is a (time-homogeneous) {\em Markov process}\index{Markov!process}\index{process!Markov}  with respect to $(\cF_t)$ with transition function $(P_t)$ if for any bounded and $\B_\cs$-measurable function $f: \cs \to \R$ and $s, t \in \R_+$,
\beq
\label{ch5-added-1-1}
E[f(X_{s+t}) \mid \cF_s] = P_t f(X_s),\qquad  P\text{-a.s.}
\eeq
The law of $X_0$ is called the {\em initial distribution} of $X$.
\end{def1}

It is well-known that a process $X$ is a Markov process with respect to its natural filtration with transition function $(P_t)$ and initial distribution $\mu$ if and only if for all $n \in \N^*$, $0 = t_0 < t_1 < \cdots < t_n$ and bounded measurable functions $f_0, f_1,\dots, f_n$ from $\cs$ to $\R$,
\begin{align*}
   E\left[\prod_{i=0}^n f_i(X_{t_i})\right] &= \int_\cs \mu(dg_0) f_0(g_0) \int_\cs P_{t_1- t_0}(g_0, dg_1)\, f_1(g_1)\\
       & \qquad\qquad \times \cdots \times \int_\cs P_{t_n - t_{n-1}}(g_{n-1}, dg_n)\, f_n(g_n)
\end{align*}
(see for instance \cite[Chap. III, §1, Prop. (1.4)]{ry}).
\bigskip

\noindent{\em The canonical probability space for continuous Markov processes}
\bigskip

Let $\tilde\Omega= \cC(\re_+,\cs)$ be the space of continuous functions $\tilde\omega$ from $\re_+$ to $\cs$, endowed with the topology of uniform convergence on compact sets of $\re_+$ and the Borel $\sigma$-field $\cB_{\cC(\re_+,\cs)}$. For $t \in \R_+$, we define an $\cs$-valued random variable $\tilde X_t$ on $\tilde\Omega$ by $\tilde X_t(\tilde\omega) = \tilde\omega(t)$, for all $\tilde\omega \in \tilde\Omega$.

For $t \in \R_+$, define the following $\sigma$-fields on $\cC(\re_+,\cs)$:
\beq\label{rd02_21e1}
    \tilde\cF_t^0 = \sigma(\tilde X(s),\, 0\le s\le t), \qquad \tilde\cF_\infty^0 = \bigvee_{t\in\re_+}  \tilde\cF_t^0.
\eeq
The filtration $(\tilde\cF_t^0 ,\, t \in \R_+)$ is the {\em canonical filtration.} The $\sigma$-field $\tilde\cF_\infty^0$ is equal to the Borel $\sigma$-field $\cB_{\cC(\re_+,\cs)}$ (see for instance \cite[Chap. 2, §8]{billingsley-1999}). 

Let $X$ be a continuous Markov process on $(\Omega, \cF, P)$ with respect to $(\cF_t)$ with state space $\mathcal{S}$, with transition function $(P_t)$ and with initial distribution $\mu$, as in Definition \ref{rd02_19d3}. The canonical probability space and canonical process associated to $X$ are $(\tilde\Omega,\tilde\cF_\infty^0, \tilde\bP^\mu)$ and $\tilde X = (\tilde X_t)$, where $\tilde\Omega$, $\tilde\cF_\infty^0$ and $(\tilde X_t)$ are as above, and $\tilde\bP^\mu$ is the law of $X$ on $\tilde\Omega$, that is $\tilde\bP^\mu$ is the probability measure defined by
\beqn
    \tilde\bP^\mu(F ) = P\{ X \in F\} = P\{\omega \in \Omega: (t \mapsto X_t(\omega)) \in F\}, \qquad F \in \cB_{\cC(\re_+,\cs)}.
\eeqn

With the same transition function $(P_t)$ as above, for any probability measure $\nu$ on $(\cs, \B_\cs)$, we can define a projective family of probability measures $P^\nu_{t_1,\dots,t_n}$ by setting
\begin{align}\nonumber
   & P^\nu_{t_1,\dots,t_n}(A_0 \times A_1 \times \cdots \times A_n) \\ \nonumber
       &\qquad =  \int_{\cs} \nu(dg_0)\, 1_{A_0}(g_0) \int_{\cs} P_{t_1- t_0}(g_0, dg_1)\, 1_{A_1}(g_1)\\
       & \qquad\qquad \times \cdots \times \int_{\cs} P_{t_n - t_{n-1}}(g_{n-1}, dg_n)\, 1_{A_n}(g_n),
\label{rd02_19e1}
\end{align}
where $n \in \N^*$, $A_i \in \B_\cs$ for $i= 0,\dots, n$ and $0 = t_0 < t_1 < \cdots < t_n$. Then we apply the Kolmogorov extension theorem \cite[Appendix II]{billingsley-1999}. The resulting probability measure induces a probability measure $\tilde\bP^\nu$ on $\tilde\Omega= \cC(\re_+,\cs)$ because $X$ has continuous sample paths.
If the support of $\nu$ consists of a single point $g \in \cs$, then we write $\tilde\bP^g$ instead of $\tilde\bP^\nu$. 

For any probability measure $\nu$ on $(\cs, \B_\cs)$, under $\tilde\bP^\nu$, the canonical process $\tilde X$ is a Markov process with respect to $(\tilde{\cF}_t)$ with transition function $(P_t)$ and with initial distribution $\nu$ (see \cite[Chap. III, §1, Theorem (1.5)]{ry}. As a consequence of \eqref{ch5-added-1-1} (for $\tilde X$,  $(\tilde{\cF}_t)$ and $\tilde\bP^g$) with $s=0$, we have
\beq\label{rd03_01e1}
    P_t f(g) = \tilde\bE^g[f(\tilde X_t],  \qquad\text{for all } g \in \cs, 
\eeq
where $\tilde \bE^g$ denotes expectation with respect to $\tilde\bP^g$. The measures $\tilde\bP^g$ and $\tilde\bP^\nu$ are related by the formula
\beq\label{rd02_21e2}
    \tilde\bP^\nu(F) = \int_\cs \tilde\bP^g(F)\, \nu(dg), \qquad F \in \tilde{\cF}_\infty^0
\eeq
(see \cite[Chap. III, §1, Prop. (1.6)]{ry}).
\bigskip

\noindent{\em Strong Markov property}
\medskip

The strong Markov property is an extension of \eqref{ch5-added-1-1} in which the deterministic time $s$ is replaced by a stopping time $\tau$. 

Let $\tilde \cF^\nu_\infty$ be the completion of $\tilde\cF^0_\infty$ with respect to $\tilde\bP^\nu$ and let $(\tilde \cF^\nu_t)$ be the filtration obtained by adding to each $\sigma$-field $\tilde \cF^0_t$ all the the $\tilde\bP^\nu$-null sets in $\tilde \cF^\nu_\infty$. Finally, we set
\beq\label{rd02_21e3}
\tilde{\cF}_t= \bigcap_\nu\,  \tilde{\cF}_t^\nu , \qquad \tilde{\cF}_\infty=\bigvee_{\nu}\, \tilde{\cF}_t.
\eeq
It is shown in \cite[Chap. III, § 2, Prop. (2.10)]{ry} that the filtrations $(\tilde \cF^\nu_t)$ and $(\tilde{\tf}_t,\, t\in \re_+)$ are right-continuous.

For any stopping time $\tau$ with respect to the filtration $(\tilde{\cF}_t,\, t\in\re_+)$, we define the {\em $\sigma$-field of events prior to $\tau$} by
\beqn
\tilde{\cF}_\tau = \{A\in \tilde{\cF}_\infty: A\cap\{\tau\le t\}\in\tilde{\cF}_t,\ {\text{for all}} \ t \in\R_+\}.
\eeqn
The following statement is proved in \cite[Chap. III, §3, Theorem (3.1)]{ry}.

\begin{thm}\label{rd02_19t1}
Let $X$ be a continuous Markov process with state space $\cs$ and transition function $(P_t)$ on $(\cs, \B_\cs)$.
If $(P_t)$ has the weak Feller property, then the canonical process $\tilde X$ on $\tilde \Omega$ associated with $(P_t)$ has the {\em strong Markov property,}\index{strong!Markov property}\index{Markov property!strong} that is, for any initial distribution $\nu$ on $\cs$, any $\tilde \bP^\nu$-a.s.-finite stopping time $\tau$ and any bounded random variable $Z$ on $\tilde \Omega$,
\beqn
     \tilde\bE^\nu[Z \circ \theta_\tau \mid \tilde \cF_\tau] = \tilde\bE^{\tilde X_\tau}[Z], \qquad \tilde\bP^\nu-a.s.,
\eeqn
where for $t\in\re_+$, $\theta_t$ is the shift operator from $\tilde\Omega$ into $\tilde\Omega$ defined by
\beqn
\theta_t (\tilde\omega)(s,x) = \tilde\omega(s+t,x),\quad \tilde\omega\in\tilde\Omega,\ s\in\re_+,\ x\in [0,L].
\eeqn
\end{thm}

\smallskip

\subsection{Markov and strong Markov property: the nonlinear case}

In the remainder of this subsection, we will consider the solution to a nonlinear stochastic heat equation and describe how it can be viewed as a (strong) Markov process.

We consider the nonlinear SPDE \eqref{nbr-(*Mb)} below with vanishing Dirichlet boundary conditions, and we study the Markov and strong Markov properties of the random field solution when viewed as a stochastic process indexed by $\R_+$ taking values in an infinite-dimensional space. With a few modifications, the results that we will establish also hold for the same SPDE with vanishing Neumann boundary conditions and for the linear stochastic heat equation on $\re$, as mentioned in Remark \ref{ch5-added-r1} below.

First, we give some preliminaries for the study of the Markov property for a stochastic process that arises from an SPDE.
On a complete probability space $(\Omega, \tf, P)$, we consider a space-time white noise $W$ along with a right-continuous complete filtration $(\tf_t,\, t\in\re_+)$, as at the beginnings of Sections \ref{ch4-section0} and \ref{ch1'-s7}.  Let $\cL$ be the partial differential operator 
\beq
\label{operator-calL}
\cL = \frac{\partial}{\partial t} - \frac{1}{b^2}\left(\frac{\partial^2}{\partial x^2}-a^2\right), \qquad a,b\in \re,\ b\ne 0.
\eeq
   Let $\beta, \sigma: \R \to \R$ be two deterministic Lipschitz continuous functions and let $u(t,x)$ be the solution of the (time-homogeneous) SPDE
\beq
\label{nbr-(*Mb)}
      \cL u(t,x) = \sigma(u(t, x))\, \dot W(t,x) + \beta(u(t,x)), \qquad   t> 0,\   x \in\, ]0, L[,   
       \eeq
subject to vanishing Dirichlet boundary conditions and some suitable initial condition $u_0$. 
For the sake of simplicity, we take $L=1$ in the remainder of this section. 

Let
\beqn
\mathbb{D}=\{f\in\cC([0,1]): f(0)=f(1)=0\}.
\eeqn
Endowed with the distance corresponding to the supremum norm, $\mathbb{D}$ is a complete separable metric space that will play the role of the state space $\cs$.


The Green's function associated to $\cL$ is  $\Gamma(t,x;s,y):=G_{a,b}(t-s;x,y)$, with
\beq
\label{s5.3-new(*1)green}
   G_{a,b}(t;x,y) = e^{-\frac{a^2}{b^2}t} \, G\left(\frac{t}{b^2}; x,y\right),
   \eeq
where $G$ is $G_L$ defined in \eqref{ch1'.600} with $L = 1$.
It is easy to  check that $G_{a,b}$ satisfies the  assumptions (${\bf H_\Gamma}$) of Section \ref{ch4-section0} and (${\bf h_\Gamma}$) of Section \ref{ch1'-s7-new}.
Hence, according to \eqref{ch1'-s5.1}, we rewrite \eqref{nbr-(*Mb)} as the integral equation
\begin{align}
\label{nbr-(*Mb-bis)}
u_{u_0}(t,x)& = \int_0^1 G_{a,b}(t; x,y) u_0(y)\, dy\notag\\
&\quad+ \int_0^t \int_0^1 G_{a,b}(t-s;x,y) \sigma(u_{u_0}(s,y))\, W(ds,dy)\notag\\
&\quad+ \int_0^t ds \int_0^1 dy\, G_{a,b}(t-s;x,y) \beta(u_{u_0}(s,y)).
\end{align}
Here and in the sequel, we write $u_{u_0}(t, x)$ instead of $u(t,x)$ in order to emphasize the initial condition.

Assume that $u_0$ is such that $(t,x)\mapsto \int_0^1 G_{a,b}(t; x,y) u_0(y)\, dy$ is Borel and bounded over $\R_+ \times [0,1]$
(for example, this is the case if $u_0$ is continuous). Then, by Theorem \ref{ch1'-s5.t1}, there exists a random field solution
\beqn
u_{u_0} = \left(u_{u_0}(t,x),\ (t,x)\in \R_+ \times[0,1]\right)
\eeqn
to \eqref{nbr-(*Mb-bis)}.
\medskip

\noindent{\em A state space for $u_{u_0}$}
\medskip

Assume now that $u_0\in \bD$. From the random field $(u_{u_0}(t,x),\, (t,x)\in\re_+\times [0, 1])$, we obtain the stochastic process $u=(u_{u_0}(t,\ast),\, t\in\re_+)$ consisting of  $\mathbb{D}$-valued random variables.
In fact, setting $I_0(t,x) =  \int_0^1 G_{a,b}(t; x,y) u_0(y)\, dy$, we see from Proposition \ref{ch1'-pPD} (iii) that
$$\sup_{(t,x)\in\re_+\times [0,1]} |I_0(t,x)|\le \sup_{x\in[0,1]}|u_0(x)|,$$ and since $u_0\in \mathbb{D}$, we have $I_0(t,\ast)\in \mathbb{D}$ for any $t>0$. As for the two integrals on the right-hand side of \eqref{nbr-(*Mb-bis)}, we appeal to the Hölder continuity results proved in Section \ref{ch1'-s6} to conclude that for fixed $t > 0$, they define an element of $\mathbb{D}$.
Hence,
 for any $u_0\in \mathbb{D}$ and $A\in\mathcal{B}_{\mathbb{D}}$, we can define
\beq
\label{ch5-added-1-2}
P_t(u_0,A)=
P\{u_{u_0}(t,\ast)\in A\}, \qquad t \in \R_+.
\eeq
Clearly, the mapping $A \mapsto P_t(u_0,A)$ is a probability measure on $\mathcal{B}_{\mathbb{D}}$. 

As in \eqref{ch5-added-1-1-0}, for any $t\in\re_+$, $u_0 \in \bD$ and for any $\cB_{\mathbb{D}}$-measurable and bounded function $f: \bD \to \R$, we set
\beq
\label{ch5-added-1-2-2}
P_tf(u_0) = E[f(u_{u_0}(t,\ast))].
\eeq
In particular, for all $A \in \B_\bD$, $P_t 1_A(u_0) = P_t(u_0, A)$.
\bigskip

\noindent{\em Markov and weak Feller properties of $(P_t)$}
\medskip

The next lemma provides properties of the mapping $u_0\mapsto P_tf(u_0)$.

\begin{lemma}
\label{nbr-2ab}
Fix $t \in \R_+$. For any bounded $\cB_\bD$-measurable function $f : \bD \to \R$, consider the mapping $u_0 \mapsto P_t f(u_0)$. We have the following properties:
\begin{description}
\item{1.} If f is bounded and continuous, then this mapping is bounded and continuous.
 \item{2.} If f is bounded and $\cB_\bD$-measurable, then this mapping is bounded and measurable. In particular, the function $(u_0, A) \mapsto P_t(u_0, A)$ from $\bD \times \cB_\bD$ to $\R$ is a transition probability on $(\bD, \cB_\bD)$.
   \end{description}
   \end{lemma}
\begin{proof}
Fix $T > 0$. First, we prove that if $u_0^n \to u_0$ in $\bD$, then
\beqn
u_{u_0^n} (t, *) \to u_{u_0}(t, *)\quad {\text{in}}\ L^p(\Omega; \cC([0, T] \times [0, 1]).
\eeqn
To simplify the notation, set $v_n(t,x) = u_{u_0^n}(t, x) - u_{u_0}(t, x)$ and for $p \geq 2$, let
\beqn
    f_n(t) = E\left[\sup_{s \in [0, t]} \Vert v_n(t, *)\Vert_\bD^p\right].
    \eeqn
   Proceeding as in the proof of Theorem \ref{ch1'-s7-t1} (see in particular the estimates on the terms $\mathcal{I}^n(t,x)$ and $\mathcal{J}^n(t,x)$ there), we obtain
\beqn
    f_n(t) \leq c \Vert u_0^n - u_0 \Vert_\bD^p + c \int_0^t f_n(r)\, dr.
    \eeqn
Thus, by Gronwall's lemma \ref{lemC.1.1}, we deduce that
\beqn
    f_n(t) \leq c \, \Vert u_0^n - u_0 \Vert_\bD^p,
\eeqn
   for $t \in [0, T]$.
Therefore, $f_n(T) \to 0$, which implies that $u_{u_0^n} \to u_0$ in $L^p(\Omega; \cC([0, T] \times [0, 1])$. In particular, $u_{u_0^n}$ converges weakly to $u_0$. Consequently, for any  bounded and continuous function $f$,
       \beqn
    P_t f(u_0^n) = E\left[f\left(u_{u_0^n}(t, *)\right)\right] \stackrel{n\to\infty}{\longrightarrow} E\left[f\left(u_{u_0}(t, *)\right)\right] = P_t f(u_0),
    \eeqn
which proves 1.

For the proof of 2., we fix $f: \mathbb{D} \to \re$ bounded and $\cB_{\mathbb{D}}$-measurable. We want to show that the mapping $u_0 \mapsto E\left[f(u_{u_0}(t,*))\right]$ is bounded and $\cB_{\mathbb{D}}$-measurable. For this, consider first an open set
$A\subset \mathbb{D}$, and the sequence of bounded continuous functions from $\mathbb{D}$ into $\re$ defined by  $f_n(g) = 1\wedge (n \rho(g,A^c))$, $n\in{\mathbb N}^*$, where $\rho$ stands for the distance in $\mathbb{D}$ derived from the supremum norm. Since
\beqn
\lim_{n\to\infty} f_n(g) = 1_{A}(g),
\eeqn
pointwise in $g\in \mathbb{D}$, appealing to the dominated convergence theorem, we deduce that
\beqn
\lim_{n\to\infty}P_t(f_n(u_0)) = P_t(u_0, A).
\eeqn
Since by part 1., the mapping $u_0\mapsto P_t(f_n(u_0))$ is bounded and continuous, we obtain that $u_0\mapsto P_t(u_0,A)$ is $\mathcal{B}_{\mathbb{D}}$-measurable. This property extends to any $A\in \mathcal{B}_{\mathbb{D}}$, 
and then to any  bounded and $\cB_{\mathbb{D}}$-measurable function, by the usual arguments based on the Monotone Class Theorem \cite[Chap. I, §2, Theorem (2.2)]{ry}.

It follows that the map $(u_0,A)\mapsto P_t(u_0,A)$ from $\mathbb{D}\times \cB_{\mathbb{D}}$ into $[0,1]$ defined in \eqref{ch5-added-1-2} satisfies the conditions (a) and (b) stated in Definition \ref{rd02_19d1}. Indeed, the former has been discussed just after \eqref{ch5-added-1-2}, while the latter has just been proved. Therefore, $(u_0,A)\mapsto P_t(u_0,A)$ is a transition probability on $(\bD, \cB_\bD)$.
 \end{proof}

\begin{remark}
\label{ch5-added-l0-rem}
 For $m \in \N^\ast$, let $\mathbb{D}^m$ be equipped with the product topology and its Borel $\sigma$-field $\cB_{\mathbb{D}^m}$. With the same arguments as in the proof of Lemma \ref{nbr-2ab}, we can show the following:
\begin{description}
  \item{(i)} if $f: \mathbb{D}^m \to \re$ is bounded and continuous, then for all $t_1, \ldots, t_m \in \re_+$, $u_0 \mapsto E[f(u_{u_0}(t_1, \ast), \ldots, u_{u_0}(t_m, \ast))]$ is bounded and continuous;
 \item{(ii)} if $f: \mathbb{D}^m \to \re$ is bounded and $\cB_{\mathbb{D}^m}$-measurable, then for all $t_1, \ldots, t_m \in \re_+$, $u_0 \mapsto E[f(u_{u_0}(t_1, \ast), \ldots, u_{u_0}(t_m, \ast))]$ is bounded and measurable.
 \end{description}
\end{remark}

The next theorem concerns the Markov property of $u_{u_0}$.

\begin{thm}
\label{nbr-2aa} The family $(P_t,\, t \in \R_+)$ defined in \eqref{ch5-added-1-2-2} is a transition function with the weak Feller property and for $u_0 \in \bD$, the process $u_{u_0}:=(u_{u_0}(t, *),\, t\in\re_+)$ that satisfies \eqref{nbr-(*Mb-bis)} is a continuous Markov process on $\bD$ with respect to $(\cF_t)$ with transition function $(P_t)$ and initial distribution $\epsilon_{u_0}$, where $\epsilon_{u_0}(A) = 1_A(u_0)$ for $A \in \cB_\bD$. 
\end{thm}

\begin{proof}
The continuity of the process $u_{u_0}$ is a consequence of the H\"older continuity results established in Section \ref{ch1'-s6}. We will check that the process $u_{u_0}$ satisfies
\beq
\label{ch5-added-1-5}
E\left[f(u_{u_0}(s+t,\ast))|\cF_s\right] = P_tf(u_{u_0}(s,\ast)),
\eeq
for any $s, t \in \R_+$ and for any bounded  $\mathcal{B}_{\bD}$-measurable function $f$. This will imply that $(P_t)$ is a transition function, since in this case,
\begin{align*}
P_{s+t} f(u_0) &= E[f(u_{u_0}(s+t,\ast))] = E\left[E[f(u_{u_0}(s+t,\ast))|\cF_s\right]]\\
&=E\left[P_tf(u_{u_0}(s,\ast))\right] =  (P_s(P_tf))(u_0)
\end{align*}
(all the equalities are trivial except for the third one, which follows from \eqref{ch5-added-1-5}).

The weak Feller property is a consequence of Lemma \ref{nbr-2ab} 1. Finally, we establish \eqref{ch5-added-1-5}. To simplify the notation, we will write $G$ instead of $G_{a,b}$ for the Green's function given in \eqref{s5.3-new(*1)green}.

Using \eqref{nbr-(*Mb-bis)}, we see that
\begin{align}
\label{june-1}
u_{u_0}(s+t,\ast) &= A_{s,t}(*) \notag\\
&\quad +  \int_s^{s+t} \int_0^1 G(s+t-r;*,y) \sigma(u_{u_0}(r,y))\, W(dr,dy)\notag\\
 &\quad + \int_s^{s+t} dr \int_0^1 dy\ G(s+t-r;*,y) \beta(u_{u_0}(r,y)),
\end{align}
with
\begin{align*}
A_{s,t}(*)&= \int_0^1 G(s+t;*,y)u_0(y)\ dy\notag\\
 &\quad + \int_0^s \int_0^1 G(s+t-r;*,y) \sigma(u_{u_0}(r,y))\, W(dr,dy)\notag\\
 &\quad + \int_0^s dr \int_0^1 dy\, G(s+t-r;*,y) \beta(u_{u_0}(r,y)).
 \end{align*}

 Applying the semigroup property of $G$ derived from \eqref{semig}, we can write
 \begin{align*}
A_{s,t}(*) &= \int_0^1 dy\, u_0(y )\left(\int_0^1 dz\, G(t;*,z)G(s;z,y)\right)\\
&\qquad +  \int_0^s\int_0^1 W(dr,dy)\, \sigma(u_{u_0}(r,y))\left(\int_0^1 dz\, G(t;*,z)G(s-r;z,y)\right)\\
&\qquad +  \int_0^s\ dr \int_0^1 dy\, \beta(u_{u_0}(r,y))\left(\int_0^1 dz\, G(t;*,z)G(s-r;z,y)\right).
\end{align*}
The first integral on the right-hand side is deterministic and the third one is a pathwise Lebesgue integral. For these terms, we apply the classical version of Fubini's theorem (we leave  to the reader the verification of the required assumptions). For the stochastic integral, we apply the stochastic Fubini's Theorem \ref{ch1'-tfubini}. This yields
 \begin{align}
 \label{june-2}
A_{s,t}(*) &=
\int_0^1 dz\, G(t;*,z)\left(\int_0^1 G(s;z,y) u_0(y)\, dy\right.\notag\\
&\qquad \quad \left.+ \int_0^s \int_0^1 G(s-r;z,y)\, \sigma(u_{u_0}(r,y))\, W(dr,dy)\right.\notag\\
&\qquad \quad \left.+ \int_0^s dr\int_0^1 dy\, G(s-r;z,y)\, \beta(u_{u_0}(r,y))\right)\notag\\
&\qquad = \int_0^1 G(t;*,z) u_{u_0}(s,z)\, dz.
\end{align}
For the sake of completeness, we argue that the application of Theorem \ref{ch1'-tfubini} is legitimate. Indeed,
for fixed$(s, t, x_0) \in \R_+ \times \R_+ \times [0,1]$, we want to check that
\begin{align}
\label{ch5-added-1-7}
&\int_0^s \int_0^1 \left(\int_0^1 dz\, G(t; x_0,z) G(s-r; z,y)\right) \sigma(u_{u_0}(r,y))\, W(dr,dy)\notag\\
&\qquad =\int_0^1 dz\, G(t;x_0,z) \int_0^s \int_0^1 G(s-r;z,y)\sigma(u_{u_0}(r,y))\, W(dr,dy).
\end{align}
In Theorem \ref{ch1'-tfubini}, take $X:=[0,1]$, replace $x,T$ and $s$ by $z,s$ and $r$, respectively, let $\mu(dx)$ be Lebesgue measure on $[0,1]$ and $G: [0,1]\times [0,s]\times [0,1]\times \Omega\longrightarrow \re$ be given by
\beqn
G(z,r,y,\omega) = G(t; x_0, z) G(s-r; z, y) \sigma(u_{u_0}(r,y)).
\eeqn
Condition \eqref{ch1'.f1} can be checked using \eqref{1'.12} and \eqref{ch1'-s5.2}, as follows:
\begin{align*}
&E\left[\int_0^1 dz\, G(t;x_0,z) \left(\int_0^s dr \int_0^1 dy\, G^2(s-r; z,y) \sigma^2(u_{u_0}(r,y))\right)^{\half}\right]\\
&\qquad \le c  \left(1+\sup_{(r,y)\in[0,T]\times [0,1]} E [u_{u_0}^2(r,y)] \right)\ s^{\frac{1}{4}}\int_0^1 dz\, G(t; x_0,z)\\
&\qquad \le C  s^{\frac{1}{4}} < \infty.
\end{align*}
Hence, \eqref{ch5-added-1-7} holds.

Thus, from \eqref{june-1} and \eqref{june-2}, for $t \geq 0$, we obtain,
\begin{align}
   u_{u_0}(s+t, *) &= \int_0^1 G(t, *, z) u_{u_0}(s, z)\, dz \notag\\
   & \quad + \int_s^{s+t} \int_0^1 G(s + t - r, *, y) \sigma(u_{u_0}(r, y))\, W(dr, dy) \notag\\
   &\quad +   \int_s^{s+t} dr \int_0^1 dy\, G(s+ t - r, *, y) \beta(u_{u_0}(r,y)).
 \label{nbr-(*4aaa)}
 \end{align}
 
 Let $\bar W$ be the space-time white noise obtained in Lemma \ref{ch1'-lsi} (2) from the sequence $(B_r^j := W_{s+r}(e_j) - W_s(e_j),\, r \in [0, t],\, j \geq 1)$ of independent standard Brownian motions. Let $\cG_r = \sigma(B^j_\rho,\, \rho \in [0, r] )$, completed with $P$-null sets. Then for all $r \geq 0$, $\cG_r$ is independent of $\cF_s$. We consider the SPDE
 \begin{align}
    \bar u(t, *) &= \int_0^1 dz\, G(t, *, z) u_{u_0}(s, z)\notag\\
    &\qquad + \int_0^t \int_0^1 G(t - r, *, y) \sigma(\bar u(r, y))\, \bar W(dr, dy) \notag\\
    & \qquad +   \int_0^{t} dr \int_0^1 dy\, G(t - r, *, y) \beta(\bar u(r,y)),  
 \label{nbr-(*4aab)}
    \end{align}
   which is \eqref{nbr-(*Mb-bis)} with $W$  replaced by $\bar W$ and initial condition $u_{u_0}(s,*)$. Let $\bar u_{w}$ be the solution of \eqref{nbr-(*4aab)} with initial condition $w \in \bD$. In particular, by \eqref{ch5-added-1-2-2}, for all $w \in \bD$ and bounded and continuous $ f: \bD \to \R$,
    \beq
    \label{nbr-(*4aac)}
      P_t f(w) = E[f(u_{w}(t, *))] = E[f(\bar u_{w}(t, *))].     
      \eeq
    Using the Picard iteration scheme \eqref{rde4.1.2} to solve \eqref{nbr-(*4aab)} with initial condition $u_{u_0}(s, *)$, we find that  $\bar u_{w}(t, *)$ is $\cF_s \vee \cG_t$-measurable.
    In addition, there is a Borel function $\Phi: \bD \times \cC([0, t])^{\N^*}$ such that, for all $w \in \bD$,
\beq
\label{nbr-(*4aad)}
    \bar u_{w}(t, *) = \Phi(w, (B_r^j,\, r \in [0, t],\, j \geq 1))     
    \eeq
and
\beqn
    \bar u (t, *) = \Phi(u_{u_0}(s, *), (B_r^j,\, r \in [0, t],\, j \geq 1)).
    \eeqn
The sum of the last two terms in \eqref{nbr-(*4aab)} is equal to
\begin{align*}
   & \int_0^t \int_0^1 G(s + t - r, *, y) \sigma(\bar u(r-s, y))\, W(dr, dy)\\
   &\qquad\qquad + \int_s^{s+t} dr \int_0^1 dy\, G(s + t - r, *, y) \beta(\bar u(r-s,y)).
   \end{align*}
Set $v(r, *) =  u_{u_0}(s+r, *)$. By \eqref{nbr-(*4aaa)}, for $t \geq 0$,
\begin{align*}
    v(t, *) &= \int_0^1 dz\, G(t, *, z) u_{u_0}(s, z)\\
    &\qquad + \int_s^{s+t} \int_0^1 G(s + t - r, *, y) \sigma(v(r-s, y))\, W(dr, dy) \\
     &\qquad +   \int_s^{s+t} dr \int_0^1 dy\, G(s+ t - r, *, y) \beta(v(r-s,y))\\
     &=  \int_0^1 dz\, G(t, *, z) u_{u_0}(s, z)  + \int_0^{t} \int_0^1 G(t - r, *, y) \sigma(v(r, y))\, \bar W(dr, dy) \\
     &\qquad +   \int_0^{t} dr \int_0^1 dy\, G(t - r, *, y) \beta(v(r,y)),
     \end{align*}
that is, $(v(t, *))$ is also a solution of \eqref{nbr-(*4aab)} with initial condition $u_{u_0}(s, *)$. By the uniqueness statement in Theorem \ref{ch1'-s7-t1}, we have  a.s., for all $t \geq 0$,
\begin{align}\nonumber
     u_{u_0}(s+t, *) &= v(t, *) = \bar u(s, *)(t, *) \\
     &= \Phi(u_{u_0}(s, *), (B_r^j, r \in [0, t],\, j \geq 1)).
\label{june-3}
\end{align}
By Lemma \ref{ch5-added-l1} and \eqref{nbr-(*4aad)},
\begin{align*}
    E[ f( u_{u_0}(s+t, *)) \mid \cF_s] &=  E[ f(\Phi(u_{u_0}(s, *), (B_r^j,\, r \in [0, t],\, j \geq 1))) \mid \cF_s] \\
    & = E[f(\Phi(w, (B_r^j, r \in [0, t],\, j \geq 1)))]\vert_{w = u_{u_0}(s, *)}\\
    & = E[f(\bar u_{w}(t, *))]\vert_{w = u_{u_0}(s, *)}\\
    & = P_t f(u_{u_0}(s, *))
    \end{align*}
by \eqref{nbr-(*4aac)}. This establishes the Markov property \eqref{ch5-added-1-5} (or \eqref{ch5-added-1-1}).
 \end{proof}
 \medskip

\begin{remark}
\label{ch5-added-c1}
The family of operators $(P_t,\, t\in\re_+)$ defined in \eqref{ch5-added-1-2-2}  is stochastically continuous, 
\index{operator!stochastically continuous} in the sense of Definition \ref{rd02_21d1}.  
Indeed, this follows from the continuity of the mapping $t\mapsto u_{u_0}(t,\ast)$, the formula \eqref{ch5-added-1-2-2} and dominated convergence. We note in passing that the transition function of a Markov process with continuous sample paths is necessarily stochastically continuous.
\end{remark}
\bigskip

\noindent{\em Canonical space and strong Markov property}
\medskip

In the last part of this section, we discuss the canonical space and the strong Markov property for the solution $u=(u_{u_0}(t,\ast),\, t\in\re_+)$ to the SPDE \eqref{nbr-(*Mb)} (or \eqref{nbr-(*Mb-bis)}). 

The canonical probability space associated to the process $u_{u_0}$, $u_0 \in \bD$, is obtained as in the discussion between Definitions \ref{rd02_19d3} and \ref{rd02_21d1}, by setting $\cs := \bD$ and $\tilde\Omega= \cC(\re_+,\bD)$, the space of continuous functions $\tilde\omega$ from $\re_+$ to $\mathbb{D}$ endowed with the topology of uniform convergence on compact sets of $\re_+$ and the Borel $\sigma$-field $\cB_{\cC(\re_+,\mathbb{D})}$.
Since the process
\beqn
(u_{u_0}(t,x),\, (t,x)\in \re_+\times [0,1])
\eeqn
 is jointly continuous and $u_{u_0}(t,\ast)\in\mathbb{D}$, the set $(\tilde\Omega, \cB_{\tilde\Omega})$ is the space of sample paths of $u_{u_0}$. For $t\in\re_+$, we define a $\mathbb{D}$-valued random variable
$\tilde u(t,\ast)$ on $\tilde \Omega$ by
 $\tilde u(t,\ast)(\tilde\omega): =\tilde\omega(t)(\ast)$. The $\sigma$-fields $\tilde\cF_t^0$, $\tilde\cF_\infty^0$, $\tilde\cF_t$ and $\tilde\cF_\infty$ are defined as in \eqref{rd02_21e1} and \eqref{rd02_21e3}, and given $u_0 \in \bD$, the probability measure $\tilde \bP^{u_0}$ on $\cC(\re_+,\bD)$ is the law of $u_{u_0}$ on $\cC(\re_+,\bD)$. For a probability measure $\nu$ on $\bD$, the probability measure $\tilde \bP^\nu$ on $\cC(\re_+,\bD)$ can be obtained as in \eqref{rd02_21e2}. If $U_0$ is an $\cF_0$-measurable $\bD$-valued random variable defined on $(\Omega, \cF, P)$ with law $\nu$, then we can solve the SPDE \eqref{nbr-(*Mb)} with the random initial condition $U_0$, and then $\tilde \bP^\nu$ is the law of the solution $u_{U_0}$ on $\cC(\re_+,\bD)$.
 \medskip

The next theorem establishes the strong Markov property of the canonical process $\tilde u$ on $\cC(\re_+,\bD)$ under the probability measures $\tilde \bP^\nu$.

\begin{thm}\label{ch5-added-t2}
The canonical process $\tilde u$ on $\cC(\re_+,\bD)$ associated with the transition function $(P_t)$ defined in \eqref{ch5-added-1-2-2} has the strong Markov property, that is, for any initial distribution $\nu$ on $\bD$, any $\tilde \bP^\nu$-a.s.-finite stopping time $\tau$ and any bounded random variable $Y$ on $\cC(\re_+,\bD)$,
\beq
\label{ch5-added-1-12}
\tilde{\mathbb{E}}^{\nu}\left[Y\circ \theta_\tau | \tilde{\cF}_\tau\right] =\tilde{\mathbb{E}}^{\tilde u (\tau,\ast)}[Y], \quad \tilde{\mathbb{P}}^{\nu}-{\text{a.s.}},
\eeq
where for $t\in\re_+$, $\theta_t$ is the shift operator from $\cC(\re_+,\bD)$ into $\cC(\re_+,\bD)$ defined by
\beqn
\theta_t (\tilde\omega)(s,x) = \tilde\omega(s+t)(x),\quad \tilde\omega\in \cC(\re_+,\bD),\ s\in\re_+,\ x\in [0,1].
\eeqn
\end{thm}
\begin{proof}
By Theorem \ref{nbr-2aa}, $(P_t)$ is the transition function of the continuous Markov process $u_{u_0}$ and $(P_t)$ has the weak Feller property. By Theorem \ref{rd02_19t1}, $\tilde u$ has the strong Markov property \eqref{ch5-added-1-12}.
\end{proof}

\begin{remark}
\label{ch5-added-r1}
For the linear stochastic heat equation on $[0,L]$ with vanishing Neumann boundary conditions, and the linear stochastic heat equation on $\re$, we can also establish statements similar to Theorems \ref{nbr-2aa} and \ref{ch5-added-t2}.

In the first instance, we can  take $\cs =\cC([0,L])$, $u_0\in\cC([0,L])$, and appeal to Proposition \ref{ch1'-pPN}
(c) and the results of Section \ref{ch1'-s6}, to see that the required ingredients are available.

In the second case, we take $\cs$ to be the space of continuous functions on $\re$ that satisfy  \eqref{ch1'-v00}, endowed with the topology of uniform convergence on compact sets. Assuming that $u_0\in \cC(\re)$ and satisfies \eqref{ch1'-v00},
then it is not difficult to check that for all $t>0$,
 $I_0(t,\ast)$ also satisfies  \eqref{ch1'-v00}.
Further, using computations from Chapter \ref{ch1'-s5}, we can indeed check that
$u_{u_0}(s,\ast)$ also satisfies \eqref{ch1'-v00} a.s.

\end{remark}


\section{Invariant and limit measures: the linear case}
\label{ch5-added-1-s2}


In this section, we examine the weak limit as $t\to\infty$ of the law of the solution to 
the stochastic heat and related equations. The existence and identification of this limit relies on the notion of invariant measure. We begin by recalling this notion and we will show that the invariant probability
measure 
for the linear stochastic heat equation with Dirichlet boundary conditions is the law of a Brownian bridge. This law will also be
the weak limit alluded to before (generalizations of these results will be given in Subsection \ref{new-5.4.2} and Section \ref{new-5.5}). As a by-product, we prove a recurrence property of the 
solution. Then we discuss the situation for the case of Neumann boundary conditions: in this case, the solution has no limit measure, but its oscillations around a specific Brownian motion do have a limit probability measure that we identify. The solution also has an invariant measure, which, however, has infinite total mass. Then we consider the linear stochastic heat equation with a linear drift term and Dirichlet boundary conditions. Finally, we discuss the stochastic heat equation with a fractional Laplacian.

\begin{def1}
\label{invariant-measure}
Let $\left(P_t,\, t\in\re_+ \right)$ be a transition function as in Definition \ref{rd02_19d1}. For a measure \index{measure} $\mu$ on $\cB_{\cs}$, we define a probability measure $P^*_t\mu$ on $\cB_\cs$ by
\beqn
P^*_t\mu(A) := \int_{\cs} P_t(g,A)\, \mu(dg), \qquad A \in cB_\cs  
\eeqn
(the notation $P^*_t$ refers to the dual of $P_t$ defined in \eqref{ch5-added-1-1-0} when acting on the set of bounded continuous functions $f:\cs\longrightarrow \re$). We say that $\mu$ is {\em invariant}\index{invariant measure}\index{measure!invariant} with respect to $(P_t)$  if for all  $t>0$,
\beq
\label{ch5-added-1-s2-1}
P^*_t \mu = \mu .
\eeq
\end{def1}

In agreement with \eqref{rd03_01e1} and \eqref{rd02_21e2}, if $(\tilde X_t)$ is the canonical process associated to a Markov process with transition function $(P_t)$, then $P^*_t\mu (A) = \tilde \bP^\mu(A)$. 

\subsection{Dirichlet boundary conditions}\label{rd03_07ss1}
 On a complete probability space $(\Omega, \tf, P)$, we consider a space-time white noise $\dot W$ along with a right-continuous complete filtration $(\tf_s,\, s\in\re_+)$, as in Section \ref{ch2new-s1}, and the SPDE discussed in Section \ref{ch4-ss2.2-1'}:
\begin{equation}
\label{ch1'.HD-ch5-markov}
\begin{cases}
 \frac{\partial u}{\partial t}- \frac{\partial^{2}u}{\partial x^{2}} = \dot W(t, x),\qquad (t, x) \in\, ]0,\infty[ \times\, ]0, L[\ ,\\
 u(0, x) = u_{0}(x), \qquad\qquad x \in  [0,L]\ ,\\
 u(t, 0) = u(t, L) =0, \qquad t \in \, ]0,\infty[\ ,
 \end{cases}
\end{equation}
with $u_0\in \mathbb{D}$. This is a special case of \eqref{nbr-(*Mb)}. In particular, the conclusions of Theorem \ref{ch5-added-t2} hold.
\smallskip

According to \eqref{nbr-(*Mb-bis)} (or \eqref{c}),
 the random field solution $(u_{u_0}(t, x))$ to \eqref{ch1'.HD-ch5-markov} with initial condition $u_0$ is
\beq
\label{ch5-added-1-3a}
u_{u_0}(t, x) = I_0(t,x) + v(t,x),
\eeq
where, for $t \in \R_+$ and $x\in [0,L]$, 
\begin{align}\label{rd03_06e1}
   I_0(t,x) &= \int^{L}_0 dy \, G_L(t;x,y) u_{0}(y),\\ \nonumber
    v(t,x) &=  \int_0^t \int_0^L G_L(t-s; x,y)\, W(ds,dy),
\end{align}
and
\beq
\label{ch5-added-1-4}
  G_L(t;x,y)= \sum^{\infty}_{n=1} e^{-\frac{\pi^2}{L^2}n^{2} t}e_{n,L} (x) e_{n,L} (y),  \quad t>0, \ x,y \in[0,L],
\eeq
$e_{n,L}(x) := \sqrt{\frac{2}{L}} \sin\left(\frac{n\pi}{L}x\right)$, $n\in \IN^*$.

Throughout this section, $B_{0,1}=(B_{0,1}(x),\, x\in[0,1])$ will denote a standard Brownian bridge,\index{Brownian!bridge} that is, a continuous Gaussian process with mean zero and covariance
$E[B_{0,1}(x)B_{0,1}(y)] = x\wedge y - xy$. This process can be characterized as a standard Brownian motion $B=(B(x),\, x\in[0,1])$ conditioned to take the value zero at $x=1$ (see e.g. \cite[35.2]{bass-2011}). Clearly, the sample paths of $B_{0,1}$ belong to $\mathbb{D}$ a.s.
Without loss of generality, we may and will assume that $B_{0,1}$ is defined on $(\Omega, \tf, P)$ and is  $\tf_0$-measurable.
For simplicity, we will also suppose in this subsection that $L=1$.

\begin{prop}
\label{ch5-added-1-s2-p1}
Let $U_0(x) = 2^{-1/2}B_{0,1}(x)$ and $(u_{U_0}(t,x),\, (t,x)\in\re_+\times[0,1])$ be the solution to
the linear stochastic heat equation \eqref{ch1'.HD-ch5-markov} with Dirichlet boundary conditions
(with $L=1$ there) and $u_0$ replaced by the random initial condition $U_0$.
Then the law $\mu$ of $U_0$ is an invariant probability measure with respect to the Markov transition function $(P_t)$, where $P_t$ is defined in \eqref{ch5-added-1-2-2}:
\beq
\label{ch5-added-1-s2-2}
    P_t f(u_0) = E[f(u_{u_0}(t, *))],
\eeq
for $ t \in \R_+$ and $f : \bD \to \R$ bounded and $\cB_\bD$-measurable.
\end{prop}

\begin{proof}
Fix $t > 0$. By \eqref{ch5-added-1-3a} with $L=1$,
\beq
\label{5.2.4(*A)}
    u_{U_0}(t,x) = \int^{1}_0 dy \, G_L(t;x,y) U_{0}(y) + \int_0^t \int_0^1 G_L(t-s; x,y)\, W(ds,dy).
    \eeq
This implies that the law  of $u_{U_0}(t,\cdot)$ on $\mathbb{D}$ is $P^*_t\mu$.
Indeed, $U_0$ and $W$ are independent because $U_0$ is $\tf_0$-measurable. Therefore, for a bounded Borel function $f:\mathbb{D} \to \re$, we can use Lemma \ref{ch5-added-l1} to write
\begin{align*}
E\left[f(u_{U_0}(t,\ast)\right] &= E\left[E\left[f\left(u_{U_0}(t,\ast)\right) \mid U_0\right]\right]\\
&=\int_{\mathbb{D}} E\left[f(u_{u_0}(t,\ast))\right] \mu(d u_0)\\
& = \int_{\mathbb{D}} P_tf(u_0)\, \mu(d u_0).
\end{align*}
Hence, we will prove that for any $t>0$, the law of $(u_{U_0}(t,x),\, x\in[0,1])$ is the same as that of $(U_0(x),\, x\in[0,1])$. Since both of these processes are mean-zero Gaussian processes, it suffices to show that they have the same covariance function.

Because of \eqref{5.2.4(*A)} and \eqref{ch5-added-1-4}, we can write 
\beqn
    u_{U_0}(t,x) = \tilde I_0(t,x)+v(t,x)
\eeqn 
with
\begin{align*}
\tilde I_0(t,x)&= \sum_{n=1}^\infty e_{n,1}(x) e^{-\pi^2 n^2 t} \int_0^1 e_{n,1}(z) U_0(z)\, dz,\\ 
v(t,x)&=\sum_{n=1}^\infty e_{n,1}(x) \int_0^t \int_0^1 e^{-\pi^2 n^2(t-r)} e_{n,1}(z)\, W(dr,dz).
\end{align*}
Observe that for fixed $x\in[0,1]$, the first series converges a.s. and the second one converges a.s. and in $L^2(\Omega)$. Since $U_0(*)$ and $W$ are Gaussian and independent, $(\tilde I_0(t,x),\, (t,x)\in\re_+\times [0,1])$ and $(v(t,x),\, (t,x)\in\re_+\times [0,1])$ are Gaussian and independent.

Set $C_n=\int_0^1 e_{n,1}(z) U_0(z)\, dz$. Then
\begin{align*}
E[C_n C_m] &= \int_0^1 dz \int_0^1  dw\, e_{n,1}(z) e_{m,1}(w) E[U_0(z) U_0(w)] \\
&=\frac{1}{2}\int_0^1dz\int_0^1 dw\, e_{n,1}(z) e_{m,1}(w) (z\wedge w-zw).
\end{align*}
Using \eqref{rdeB.5.4} and \eqref{rdeB.5.6} in Lemma \ref{lemB.5.2}, we see that
\beq
\label{ch5-added-1-s2-2-bis}
z\wedge w-zw = \sum_{n=1}^\infty \frac{e_{n,1}(z)e_{n,1}(w)}{\pi^2 n^2}.
\eeq
This yields
\beq
\label{ch5-added-1-s2-3}
E[C_n C_m] =\delta_n^m\frac{1}{2\pi^2n^2},
\eeq
where $\delta_n^m$ denotes the Kronecker symbol.\index{Kronecker symbol}\index{symbol!Kronecker}

Set $A_t^n:=\int_0^t \int_0^1 e^{-\pi^2 n^2(t-r)} e_{n,1}(z)\, W(dr,dz)$. By the Itô isometry and the orthonormality of the sequence $(e_{n,1},\, n\in \N^*)$, for any $m, n \in\N^*$,
\beqn
E[A_t^m A_t^n] = \delta_n^m \int_0^t e^{-2\pi^2 n^2(t-r)}\, dr = \delta_n^m\, \frac{1-e^{-2\pi^2 n^2t}}{2\pi^2 n^2}.
\eeqn
 This implies that
\beq
\label{ch5-added-1-s2-4}
E[v(t,x)v(t,y)] = \sum_{n=1}^\infty e_{n,1}(x)e_{n,1}(y)\, \frac{1-e^{-2\pi^2 n^2t}}{2\pi^2 n^2}.
\eeq

From the identities \eqref{ch5-added-1-s2-2-bis}--\eqref{ch5-added-1-s2-4}, we deduce that
\begin{align}
\label{ch5-added-1-s2-4-bis}
E\left[u_{U_0}(t,x)u_{U_0}(t,y)\right] &=E\left[\tilde I_0(t,x) \tilde I_0(t,y)\right] +E\left[v(t,x)v(t,y)\right]\notag\\
&=\sum_{n=1}^\infty e_{n,1}(x)e_{n,1}(y)\left(\frac{e^{-2\pi^2n^2 t}}{2\pi^2 n^2} + \frac{1-e^{-2\pi^2n^2 t}}{2\pi^2 n^2}\right)\notag\\
&=\sum_{n=1}^\infty \frac{e_{n,1}(x)e_{n,1}(y)}{2\pi^2n^2} = \frac{x\wedge y-xy}{2},
\end{align}
which is $E[U_0(x) U_0(y)]$. Therefore, the laws of $(u_{U_0}(t,x),\, x\in[0,1])$ and $(U_0(x),\, x\in[0,1])$ are the same.
\end{proof}

\noindent{\em Digression}
\smallskip

The fact that up to a multiplicative constant, the law of the standard Brownian bridge is the invariant measure of the linear stochastic heat equation with Dirichlet boundary conditions can be intuitively derived from the following observations.

Define
\beqn
u^n(t) = \left\langle u_{u_0}(t,\ast), e_{n,1}\right\rangle_{L^2([0,1])},\qquad u_0^n =\langle u_0,e_{n,1}\rangle_{L^2([0,1])},
\eeqn
and use \eqref{series-dirichlet}and \eqref{OU-coordinates} with $L = 1$ to see that
\beq
\label{ch5-added-1-s2-5}
u^n(t) = e^{-\pi^2 n^2t}u_0^n+\int_0^t e^{-\pi^2n^2(t-s)}\, W^n(ds),\quad n\in\N^*,
\eeq
where $W^n(t) = W_t(e_{n,1})$  
and therefore $(W^n(t),\, t\in\re_+, \,n\in\N^*)$ is a sequence of independent Brownian motions.

The process $(u^n(t),\, t\in\re_+)$ given in \eqref{ch5-added-1-s2-5} is an Ornstein-Uhlenbeck process, that is, the solution to the linear stochastic differential equation
\beq
\label{ch5-added-1-s2-6}
u^n(t) =u_0^n - \pi^2n^2 \int_0^t u^n(s)\, ds + W^n(t).
\eeq
It is a well-known fact that  the unique invariant measure of $(u^n(t),\, t\in\re_+)$ is 
${\rm N}\left(0,(2\pi^2n^2)^{-1}\right)$ (for a proof, see Lemma \ref{nbr-0b}). 

According to \eqref{series-dirichlet} with $L = 1$ and $u(t, x) = u_{u_0}(t, x)$,
\beqn
   u_{u_0}(t, x) = \sum_{n=1}^\infty 
                                u^n(t)\, e_{n, 1}(x).
\eeqn
Therefore, we should expect
that the invariant probability measure $\mu$ with respect to the transition function $(P_t)$ given in \eqref{ch5-added-1-s2-2} will be the law of the random function 
\beqn
   V_0(*) = \sum_{n=1}^\infty \frac{Z_n}{\sqrt{2}\, \pi n}\, e_{n,1}(*),
\eeqn
where $(Z_n,\, n\in\N^*)$ is a sequence of i.i.d.~${\rm N}(0,1)$ random variables. Appealing to the identity \eqref{ch5-added-1-s2-2-bis}, we observe that $V_0$ has the law of $2^{-\half} B_{0,1}$, where $B_{0,1}$ is a Brownian bridge.
\bigskip

\noindent{\em Weak convergence to the invariant measure}
\smallskip

\begin{thm}
\label{ch5-added-1-s2-t1}
Let $u_0\in\mathbb{D}$ and  $(u_{u_0}(t,x),\, (t,x)\in\re_+\times[0,1])$ be the solution to the linear
stochastic heat equation \eqref{ch1'.HD-ch5-markov} with vanishing Dirichlet boundary conditions (with $L=1$ there),  and let $P_t$ be as defined in \eqref{ch5-added-1-s2-2}. The (invariant) probability measure $\mu$ defined in Proposition
\ref{ch5-added-1-s2-p1}, that is, the law of the process $2^{-1/2} B_{0,1}$, is the weak limit as $t\to\infty$ of the law of $u_{u_0}(t,\ast)$, that is, for any bounded and continuous function $f:\mathbb{D}\longrightarrow \re$ and $u_0\in\mathbb{D}$,
\beqn
\lim_{t\to\infty} P_tf(u_0) = \mu(f) := \int_{\mathbb{D}} f(g)\, \mu(dg).
\eeqn
\end{thm}

\begin{proof}
Fix $u_0\in\mathbb{D}$. In this proof, in order to simplify the notation, we will write $u(t,x)$ instead of $u_{u_0}(t,x)$. Recall that $\mathbb{D}$ is a complete separable metric space (for the distance associated with the supremum norm). Hence, according to \cite[Theorem 8.1]{billingsley-1999} (with $\mathcal{C}([0,L])$ there replaced by $\mathbb{D}$),
two facts have to be proved:
\begin{enumerate}
\item Convergence of the finite-dimensional distributions. For any $k\in\N^*$ and any $x_1,\ldots,x_k\in [0,1]$, as $t \to \infty$,
\beq
\label{ch5-added-1-s2-7}
   \left(u(t,x_1), \ldots, u(t,x_k)\right)  \stackrel{w}\longrightarrow 2^{-1/2}\left(B_{0,1}(x_1),\ldots, B_{0,1}(x_k)\right).
\eeq
\item  Equicontinuity. For any $\eta>0$,
\beq
\label{ch5-added-1-s2-8}
\lim_{\varepsilon\to 0} \limsup_{t\to\infty} P\left\{\sup_{|x-y|\le \varepsilon} |u(t,x)-u(t,y)|\ge \eta\right\} = 0.
\eeq
\end{enumerate}
Indeed, since $u(t,0) = u(t,L)=0$ for $t>0$, it is obvious that for any $\delta>0$, there exists $a\ge 0$ such that $P\{|u(t,0)|>a\} + P\{|u(t,L)|>a\} =0 \le \delta$. Along with the equicontinuity property, this is equivalent to the tightness of the family of laws of $u(t,\ast)$, $t>0$ (see \cite[Theorem 8.2]{billingsley-1999} (adapted to $\mathbb{D}$ instead of $\mathcal{C}([0,L])$ there).

In order to prove 1., notice that the random vector $\left(u(t,x_1), \ldots, u(t,x_k)\right)$ is Gaussian.
Recalling the definition of $I_0(t,x)$ in \eqref{rd03_06e1} (with $L = 1$ there), and the computations in the proof of Proposition \ref{ch5-added-1-s2-p1}, for any $j \in \{1,\ldots,k\}$, we have
\begin{align*}
   \vert E\left[u(t,x_j)\right] \vert &= \vert I_0(t,x_j) \vert 
    = \left\vert \sum_{n=1}^\infty e_{n,1}(x_j) e^{-\pi^2 n^2 t} \int_0^1 e_{n,1}(z)u_0(z)\, dz \right\vert\\
    &\leq 2 \sum_{n=1}^\infty e^{-\pi^2 n^2 t} \int_0^1 \vert u_0(z)\vert \, dz \\ 
    & \longrightarrow 0 \qquad \text{as } t \to \infty.
\end{align*}
Moreover, for any $\ell \in \{ 1,\dots, k\}$, letting $x = x_j$ and $y = x_\ell$, we have
\begin{align*}
&E\left[(u(t,x)-I_0(t,x))(u(t,y)-I_0(t,x_y))\right] \\
&\qquad\qquad\qquad = \half \sum_{n=1}^\infty e_{n,1}(x)e_{n,1}(y)\, \frac{1-e^{-2\pi^2 n^2t}}{\pi^2 n^2}.
\end{align*}
As $t\to\infty$, the series converges to the series on the right-hand side of \eqref{ch5-added-1-s2-2-bis}.
Applying \eqref{ch5-added-1-s2-2-bis}, we see that
\beqn
\lim_{t\to \infty}E\left[(u(t,x)-I_0(t,x)) (u(t,y)-I_0(t,y))\right]= \frac{x\wedge y-xy}{2}.
\eeqn
This proves Fact 1. above.

We now turn to the proof of Fact 2. From the definition of $I_0(t,x)$, we see that for $t>0$,
\begin{align*}
|I_0(t,x)- I_0(t,y)|&\le  \sum_{n=1}^\infty e^{-\pi^2 n^2 t} \left(\int_0^1 e_{n,1}(z) u_0(z)\, dz\right) |e_{n,1}(x)-e_{n,1}(y)|\notag\\
&\le C\, |x-y| \,  \sum_{n=1}^\infty n e^{-\pi^2 n^2 t} \notag\\
&\le C_t\, |x-y|,
\end{align*}
where $C_t = C \sum_{n=1}^\infty n e^{-\pi^2 n^2 t}$ (which tends to $\infty$ as $t\to 0$). Moreover, by applying \eqref{1'.10}, we have
\beqn
E[(v(t,x) - v(t,y))^2] \le \half\, |x-y|.
\eeqn
Altogether, for any $x,y\in[0,1]$, there is $\tilde C<\infty$ such that
\beqn
\sup_{t\geq 1}\Vert u(t,x)-u(t,y)\Vert_{L^2(\Omega)} \le \tilde C\, |x-y|^{1/2}.
\eeqn
Apply Theorem \ref{app1-3-t1} (see in particular \eqref{ch1'-s7.20}) to deduce that for $p>2$ and $\alpha\in\, ]2/p,1[$, there is $C>0$ such that for all $t\ge 1$,
\beq
\label{ch5-added-1-s2-9}
E\left[\sup_{|x-y|<\varepsilon}\left\vert u(t,x)-u(t,y)\right\vert^p\right] \le C\,  \varepsilon^{\frac{\alpha p-2}{2}}.
\eeq
Applying Chebyshev's inequality,  we obtain \eqref{ch5-added-1-s2-8} from \eqref{ch5-added-1-s2-9}, because $\alpha p-2>0$.
\end{proof}
\medskip

\noindent{\em Recurrence}
\medskip

For a Markov process $X = (X_t)$ as in Definition \ref{rd02_21d1}, following \cite[Section 3.4, (3.4.5)]{dz-1996}, we say that {\em $X$ is recurrent\index{recurrent}\index{Markov!process, recurrent} with respect to a set $A\in\cB_{\cs}$} if
\beq
\label{ch5-added-1-s2-10}
P\{X_t\in A\ {\text {for an unbounded set of}}\ t \in \R_+\} = 1.
\eeq

Here, we consider the Markov process $u_{u_0}$ of Theorem \ref{nbr-2aa}, in the particular case of equation \eqref{ch1'.HD-ch5-markov} with $L=1$, and we prove the following.

\begin{prop}
\label{ch5-added-1-s2-p2}
Let $(P_t, \, t\in\re_+)$ be the transition function defined in \eqref{ch5-added-1-s2-2} and let $\mu$ be the invariant measure given in Proposition \ref{ch5-added-1-s2-p1}. For $R>0$, let $N_R(0)\subset \mathbb{D}$ be the open ball centred at $0$ with radius $R$. Then the Markov process $u_{u_0}$ that solves \eqref{ch1'.HD-ch5-markov} with $L=1$ is recurrent with respect to $N_R(0)$ and with respect to $N_R(0)^c$.
\end{prop}
\begin{proof}
We first prove the statement on $N_R(0)$. 
According to \cite[Corollary 3.4.6]{dz-1996}, $A:= N_R(0)$ is recurrent if
\beq
\label{ch5-added-1-s2-11}
\lim_{t\to\infty} P_t(u_0,A) = \mu(A)>0.
\eeq
This property is a consequence of the general results of Theorem \ref{nbr-t27} (c) and Proposition \ref{nbr-p21} below, but we give here a direct proof that only uses results discussed up to now.

The distribution function of the supremum of the absolute values of a  Brownian bridge $B_{0,1}$ has the following expression (see e.g. \cite[(9.39)]{billingsley-1999}):
\beq
\label{ch5-added-1-s2-12}
F(R):=P\left\{\sup_{x\in[0,1]}|B_{0,1}(x)|\le R\right\} = 1+2\sum_{n=1}^\infty(-1)^n e^{-2n^2R^2},\quad R>0.
\eeq
In particular, the law of the random variable $\sup_{x\in[0,1]}|B_{0,1}(x)|$ is absolutely continuous.

We now show  that $F(R)\in\, ]0,1[$ for any $R>0$. Indeed, the inequality $F(R)<1$ follows from the
explicit form of the alternating series in \eqref{ch5-added-1-s2-12}.
In order to show that $F(R) > 0$, we use the property $B_{0,1}(x) \stackrel{\rm d}{=}B_x-xB_1$, where $(B_x,\, x\in[0,1])$ is a standard Brownian motion. By the triangle inequality,
\begin{align*}
F(R) &\ge P\left\{ \sup_{x\in[0,1]} |B_x| + |B_1| \le R\right\} \ge P\left\{\sup_{x\in[0,1]} |B_x| \le R/2\right\}\\
& = H(R/2),
\end{align*}
where
\beqn
H(r) = \frac{4}{\pi}\sum_{n=0}^\infty\frac{(-1)^n}{2n+1} \exp\left(-\frac{(2n+1)^2 \pi^2}{8r^2}\right), \quad r>0,
\eeqn
is the (continuous) distribution function of $\sup_{x\in[0,1]}|B_x|$ (see e.g. \cite[Section 7.3, Exercise 8, p. 223]{chung-1974}). Clearly, $H(r)>0$ for any $r>0$. Therefore,
\beqn
\mu(N_R(0)) = F(\sqrt 2R) \ge H(2^{-\half}R) >0.
\eeqn

According to Theorem \ref{ch5-added-1-s2-t1}, $u_{u_0}(t,\ast)\stackrel{w}\longrightarrow B_{0,1}$ as $t\to\infty$. By the {\em portmanteau theorem}
\cite[Theorem 30.2]{bass-2011}, and since
\beqn
\mu(\partial N_R(0))=P\left\{\sup_{x\in[0,1]} |B_{0,1}(x)|= \sqrt 2 R\right\} = 0,
\eeqn
we obtain
\beq
\label{ch5-added-1-s2-13}
\lim_{t\to\infty}P_t\left(u_0,N_R(0)\right)=\lim_{t\to\infty}P\left\{u_{u_0}(t,\ast))\in N_R(0)\right\} =\mu(N_R(0))>0.
\eeq
This yields \eqref{ch5-added-1-s2-11} for $A = N_R(0)$.

For $A = N_R(0)^c$,
we note that
\beqn
\mu(N_R(0)^c) = P\left\{\sup_{x \in [0,1]} \vert B_{0,1}(x)\vert \geq \sqrt{2} R\right\} = 1 - F(\sqrt{2} R) > 0,
\eeqn
 and
$\partial N_R(0)^c = \partial N_R(0)$ has $\mu$-measure $0$. We conclude that \eqref{ch5-added-1-s2-11} holds for $A = N_R(0)^c$, completing the proof of the proposition.
\end{proof}
\medskip

\subsection{Neumann boundary conditions}\label{rd03_06ss1}

 In the case of vanishing Neumann boundary conditions, the random field solution $(v_{v_0}(t, x))$ to the stochastic heat equation \eqref{1'.14} with initial condition $v_0$  is given by an expression similar to \eqref{ch5-added-1-3a}, with $u_0$ there replaced by $v_0$ and
\beqn
  G_L(t;x,y)= \sum^{\infty}_{n=0} e^{-\frac{\pi^2}{L^2}n^{2} t}g_{n,L} (x) g_{n,L} (y),  \qquad t>0, \quad x,y \in[0,L],
\eeqn
$g_{0,L}=\frac{1}{\sqrt L}$ and $g_{n,L}(x) := \sqrt{\frac{2}{L}} \cos\left(\frac{n\pi}{L}x\right)$, $n\in \IN^*$, $L > 0$.

Let $(v(t,x))$ be the solution to \eqref{1'.14} when $v_0 \equiv 0$. Computing the second moment of the stochastic integral, we find
\beqn
E\left [v^2(t,x) \right] = \frac{t}{L} + L\sum_{n=1}^\infty \frac{1-e^{2\frac{\pi^2}{L^2}n^2 t}}{\pi^2n^2} \cos^2\left(\frac{n\pi}{L}x\right).
\eeqn
When $t\to\infty$, this expression converges to $\infty$ because of the term $t/L$. Since the process $v_{v_0}$ is Gaussian, it does not converge in law as $t\to\infty$. However, $E\left[v^2(t,x) \right]-\frac{t}{L}$ does converge as $t\to\infty$, to
\beqn
    L\sum_{n=1}^\infty \frac{\cos^2(n \pi x / L)}{\pi^2 n^2},
    \eeqn
which is finite (and, more precisely, bounded by $L/6$).

In fact, we have the following result.
\begin{prop}
   \label{rem6.1.15-*1}
As $t \to \infty$, the $\mathcal{C}([0,L])$-valued stochastic process $v(t, *) - v(t, 0)$ converges in law to the process $ 2^{-\half} B_*$, where $B = (B_x,\, x \in [0, L])$ is a standard Brownian motion. 
 \end{prop}
      
  \begin{proof}
      Note that
    \beqn
    v(t, x)  = W_t(g_{0, L}) + \sum_{n=1}^\infty g_{n, L}(x) \int_0^t \exp\left(-\frac{\pi^2}{L^2}n^2 (t-s)\right) dW_s(g_{n, L}),
    \eeqn
so setting $x = 0$,
\beqn
    v(t, 0) = W_t(g_{0, L}) + \sqrt{\frac{2}{L}} \sum_{n=1}^\infty  \int_0^t \exp\left(-\frac{\pi^2}{L^2}n^2 (t-s)\right) dW_s(g_{n, L}).
 \eeqn
 Therefore,
 \begin{align*}
    &v(t, x) - v(t, 0)\\
    &\qquad = \sqrt{\frac{2}{L}} \sum_{n=1}^\infty \left(\cos\left(\frac{n\pi}{L}x\right)-1 \right) \int_0^t \exp\left(-\frac{\pi^2}{L^2} n^2 (t-s)\right) dW_s(g_{n, L}).
    \end{align*}

For fixed $t>0$ and $x\in\, ]0,L]$,  $v(t, x) - v(t, 0)$ is a mean-zero Gaussian random variable with variance
  \begin{align*}
   \sigma^2_{v(t,x)-v(t,0)}& =\frac{2}{L} \sum_{n=1}^\infty \left(\cos\left(\frac{n\pi}{L}x\right) - 1\right)^2 \int_0^t\exp\left(-\frac{2\pi^2}{L^2}n^2 (t-s)\right) ds \\
   &= \frac{2}{L} \sum_{n=1}^\infty \left(\cos\left(\frac{n\pi}{L}x\right)-1 \right)^2 \frac{L^2}{2\pi^2 n^2} \left(1- \exp\left(-\frac{2 \pi^2}{L^2} n^2 t\right)\right).
 \end{align*}
Using \eqref{rdeB.5.5}, we see that
  \begin{equation}\label{rd04_23e1}
  \lim_{t\to \infty}  \sigma^2_{v(t,x)-v(t,0)} = 
 L \sum_{n=1}^\infty\, \frac{\left(1 - \cos\left(\frac{n\pi}{L} x\right)\right)^2}{\pi^2 n^2} = \frac{x}{2}\, .
  \end{equation}
With similar computations, we can evaluate the covariance of
$\sqrt{2} (v(t, x) - v(t, 0))$ and $\sqrt{2} (v(t, y) - v(t, 0))$ and see that its limit as $t\to \infty$ is
\beqn
   2L \sum_{n=1}^\infty\, \frac{\left(1 - \cos(n \pi x / L) (1 - \cos(n \pi y / L)\right)}{\pi^2 n^2} = x \wedge y
   \eeqn
by \eqref{rdeB.5.6-quatre}. Therefore, this limiting covariance is that of a standard Brownian motion $(B_x,\, x \in [0, L])$.

 
  Since $(v(t, x) - v(t, 0))$ is a Gaussian process, the calculations above show that the finite-dimensional distributions of $(v(t, x) - v(t, 0))$ converge to those of $( 2^{-1/2}\, B_x)$, where $(B_x)$ is a standard Brownian motion.

    In order to show the claimed weak convergence, it remains to check that $\left(v(t, x) - v(t, 0) \right)$ satisfies the equicontinuity property \eqref{ch5-added-1-s2-8}. As in the proof of Theorem \ref{ch5-added-1-s2-t1}, but using \eqref{1'.17} instead of \eqref{1'.10}, we obtain \eqref{ch5-added-1-s2-9} for $v$ instead of $u$. This completes the proof of Proposition \ref{rem6.1.15-*1}.   
\end{proof}

\noindent{\em Another stationary pinned string}
\medskip
    
Just as in Proposition \ref{rd04_23p1}, if the initial condition $u_0$ in \eqref{1'.14} is random, independent of $\dot W$, and is a standard Brownian motion multiplied by $2^{-1/2}$, then for all $t > 0$, the solution $u$ has the property that the law of $u(t,*) - u(t, 0)$ is that of a standard Brownian motion multiplied by $2^{-1/2}$, as we show in the next proposition.

\begin{prop} 
Assume that $u_0 = (u_0(x), \, x \in [0, L])$ is a stochastic process  independent of the space-time white noise $\dot W$ in the SPDE \eqref{1'.14}, with the law of a standard Brownian motion multiplied by $2^{-1/2}$. Let $u$ be the solution to \eqref{1'.14} given in \eqref{1'.N} with this initial condition $u_0$. Then for all $t \geq 0$, the process $(u(t, x) - u(t, 0),\, x \in [0, L])$ also has the law of a standard Brownian motion multiplied by $2^{-1/2}$.
\end{prop}

\begin{proof}
In order to create a standard Brownian motion multiplied by $2^{-1/2}$ that is independent of $\dot W$ on $\R_+ \times [0, L]$, we assume that the space-time white noise is defined on $\R \times [0, L]$ and, for $x \in [0, L]$, we set
\begin{align*}
   u_0(x) &= \int_{-\infty}^0 \int_0^L (G_L(-r; x, z) - G_L(-r; 0, z))\, W(dr, dz). 
\end{align*}
Then $u_0(0) = 0$ a.s., and by \eqref{1'.17}, $E[(u_0(x) - u_0(y))^2] = \half\, \vert x - y \vert$, therefore, $u_0$ is a standard Brownian motion multiplied by $2^{-1/2}$. We work with its continuous modification, that we also denote $u_0$.

The solution to \eqref{1'.14} is
\beqn
     u(t, x) = I_0(t, x) +  \int_0^t \int_0^L G_L(t-r; x,z)\, W(dr, dz),
\eeqn
where
\begin{align*}
    I_0(t, x) &= \int_0^L G_L(t; x, y) u_0(y)\, dy \\
     &= \int_0^L dy\, G_L(t; x, y) \int_{-\infty}^0 \int_0^L (G_L(-r; y, z) - G_L(-r; 0, z))\, W(dr, dz) .
\end{align*}
Using the stochastic Fubini’s theorem Theorem \ref{ch1'-tfubini}, whose assumptions are easily seen to be satisfied, the semigroup property of $G_L$ in Proposition \ref{ch1'-pPN} (a) and the full density property \eqref{int1}, we see that
\beqn
   I_0(t, x) = \int_{-\infty}^0 \int_0^L (G_L(t-r; x, z) - G_L(-r; 0, z))\, W(dr, dz).
\eeqn
Therefore, for $x \in [0, L]$,
\begin{align*}
  u(t, x) - u(t, 0) &=  \int_{-\infty}^0 \int_0^L (G_L(t-r; x, z) - G_L(-r; 0, z))\, W(dr, dz) \\
      &\qquad + \int_0^t \int_0^L (G_L(t-s; x,z) - G_L(t-s; 0,z)) \, W(ds, dz),
\end{align*}
and for $x, y \in [0, L]$, 
\begin{align*}
  & u(t, x) - u(t, 0) - (u(t, y) - u(t, 0)) \\
  &\qquad\qquad =  \int_{-\infty}^0 \int_0^L (G_L(t-r; x,z) - G_L(t-r; y,z))\, W(dr, dz) \\
   &\qquad\qquad\qquad + \int_0^t \int_0^L (G_L(t-s; x,z) - G_L(t-s; y,z))\, W(ds, dz).
\end{align*} 
Using the independence properties of space-time white noise, we see that the second moment of this expression is equal to
\begin{align*}
    & \int_{-\infty}^0 \int_\R (G_L(t-r; x, z) - G_L(t-r; y, z))^2\, dr dz \\
   &\qquad\qquad  + \int_0^t \int_\R (G_L(t-s; x,z) - G_L(t-s; y,z))^2\, ds dz \\
   &\qquad =  \int_0^\infty \int_0^L (G_L(s; x, z) - G_L(s; y, z))^2\, ds dz \\
   &\qquad = \half\, \vert x - y \vert
\end{align*}
by \eqref{1'.17}. Therefore, $u(t, *) - u(t, 0)$ is a standard Brownian motion multiplied by $2^{-1/2}$.
\end{proof}

\medskip

\noindent{\em An invariant measure for $v_{v_0}$}
\medskip

   In the next theorem, we show that in the case of vanishing Neumann boundary conditions, there exists an invariant measure for $v_{v_0}$ which, however, is not a finite measure. Before giving the precise statement, we introduce some notation. Let $(B_x,\, x\in [0, L])$ be a standard Brownian motion (independent of $W$) and for any $a\in \re$,  define $V_L(a, x) = L^{-\half} a + 2^{-\half} B_x$, $x \in [0, L]$.  We denote $v_{V_L(a, *)}$ the solution to the linear stochastic heat equation on $[0, L]$ with vanishing Neumann boundary conditions and random initial condition $V_L(a, *)$.

\begin{thm}
\label{rem6.1.15-*2}
Let $\cC([0, L])$ be the set of continuous functions from $[0, L]$ to $\R$ endowed with the uniform norm. Let $\nu_L$ be the measure on $(\cC([0, L]), \B_{\cC([0, L])})$ defined by
\beqn
     \nu_L(A) = \int_\R da\, P\{V_L(a, *) \in A \},  \qquad  A \in \B_{\cC([0, L])}.
\eeqn
Let $(P_t,\, t \in \R_+)$  be the  transition function associated with $v$ as in (\ref{ch5-added-1-2-2}) (with $u_{u_0}$ replaced by $v_{v_0}$). Then $\nu_L$ is an invariant measure with respect to $(P_t)$, that is, for any $t > 0$ and bounded continuous function $f: \cC([0, L]) \to \R$, 
we have
\beq
\label{rem6.1.15-(*1)}
    \int_\R da\, E\left[f(v_{V_L(a, *)}(t, *))\right] = \int_\R da\, E\left[f(V_L(a, *))\right].        
\eeq
\end{thm}
\vskip 12pt

\begin{remark}
\label{neumann-invariant-interpretation}
The measure $\nu_L$ is the measure on $\cC([0, L])$ induced by a Brownian motion on $[0, L]$ with speed $\half$ and initial measure equal to a multiple of Lebesgue measure on $\R$.
\end{remark}

 \noindent{\em Proof of Theorem \ref{rem6.1.15-*2}}.\
By definition,
\beqn
    v_{V_L(a, *)}(t, x) = I_0(a, t, x) + v(t, x),
    \eeqn
where
\begin{align*}
     I_0(a, t, x) &= \int_0^L dz\, G_L(t; x, z) V_L(a, z)\\
     & = \frac{1}{\sqrt L} \int_0^L V_L(a, z)\, dz + \sum_{n=1}^\infty g_{n, L}(x)\, e^{-\frac{\pi^2}{L^2} n^2 t} \int_0^L g_{n, L}(z)\, V_L(a, z)\, dz
     \end{align*}
and
\begin{align*}
      v(t, x) &= \int_0^t \int_0^L G_L(t-r; x, z)\, W(dr, dz) \\
      &= \frac{1}{\sqrt L} W([0, t] \times [0, L]) \\
        &\qquad\qquad    + \sum_{n=1}^\infty g_{n, L}(x) \int_0^t \int_0^L e^{-\frac{\pi^2}{L^2} n^2 (t-r)} g_{n, L}(z)\, W(dr, dz).
      \end{align*}

      In order to prove \eqref{rem6.1.15-(*1)}, it suffices to take
$f: \cC([0, L]) \to \R$ of the form
    $ f(w) = f_0(w(0)) f_1(w(x_1)) \cdots f_m(w(x_m)),$
where $m \geq 1$, $0 < x_1 < \cdots < x_m \leq L$, and $f_i: \R \to \R$ is of the form $f_i(x) = 1_{A_i}(x)$, where $A_i \in \cB_\R$,  $i = 0,\dots, m$.
     With this choice of $f$, the left-hand side of \eqref{rem6.1.15-(*1)} is equal to
     \begin{align}
     \label{rem6.1.15-(*2)}
      &\int_\R da\, E\left[f_0(v_{V_L(a, *)}(t, 0))\,  f_1(v_{V_L(a, *)}(t, x_1))\, \cdots\, f_m(v_{V_L(a, *)}(t, x_m))\right]\notag\\
     &\qquad =  E\Big[\int_\R da\, f_0( I_0(a, t, 0) + v(t, 0))\notag\\
     &\qquad\qquad\qquad\times f_1( I_0(a, t, x_1) + v(t, x_1))\, \cdots\, f_m( I_0(a, t, x_m) + v(t, x_m))\Big],       
     \end{align}
and the right-hand side of \eqref{rem6.1.15-(*1)} is equal to
\beq
\label{rem6.1.15-(*3)}
     \int_\R da\, E\left[f_0(a)\, f_1(a + 2^{-\half} B_{x_1})\, \cdots\, f_m(a + 2^{-\half} B_{x_m})\right].    
     \eeq
For any $n\ge 1$, let
\begin{align*}
      C_n & = 2^{-\half}\, e^{-\frac{\pi^2}{L^2} n^2 t} \int_0^L g_{n, L}(z)\, B_z\, dz, \\
      D_n & = \int_0^t \int_0^L e^{-\frac{\pi^2}{L^2} n^2 (t-r)} g_{n, L}(z)\, W(dr, dz).
\end{align*}
Using integration by parts, we see that
\beqn
     C_n = -\frac{\sqrt{L}}{n \pi} e^{-\frac{\pi^2}{L^2} n^2 t} \int_0^L \sin\left(\frac{n \pi}{L} z\right)\, dB_z\, .
\eeqn
Therefore, the families $(C_n)$ and $(D_n)$ are independent sequences of independent centred Gaussian random variables. In particular, using the Wiener isometry, we see that 
the laws of $C_n$ and $D_n$ are respectively
\beqn
   {\rm{N}}\left(0, \frac{L^2}{2 n^2 \pi^2}\, e^{- \frac{2\pi^2}{L^2} n^2 t}\right)\qquad \text{and}\qquad
   {\rm {N}}\left(0, \frac{L^2}{2 n^2 \pi^2}\left(1 - e^{- \frac{2\pi^2}{L^2} n^2 t}\right)\right).
\eeqn
This implies that
\beq
\label{rem6.1.15-(*4)}
   \sqrt{\frac{2}{L}} \, (C_n + D_n)\stackrel{d}{=} N\left(0, \frac{L}{n^2 \pi^2}\right), \quad n\ge 1,      
   \eeq
and these random variables are independent. This fact will be used later on in the proof.
\smallskip
      
Since $\int_0^L g_{n, L}(z)\, dz = 0$ when $n \geq 1$, we have
\begin{align*}
   &I_0(a, t, 0) + v(t, 0) \\
   &\qquad = a + \frac{2^{-\half}}{\sqrt L} \int_0^L B_z\, dz + \frac{1}{\sqrt L} W([0, t] \times [0, L]) \\
   &\qquad \qquad
   + \sum_{n=1}^\infty \sqrt{\frac{2}{L}} e^{-\frac{\pi^2}{L^2} n^2 t}\int_0^L g_{n, L}(z) \left(a + 2^{-\half} B_z\right) dz \\
   &\qquad \qquad
    +\sum_{n=1}^\infty \sqrt{\frac{2}{L}} \int_0^t \int_0^L e^{-\frac{\pi^2}{L^2} n^2 (t-r)} g_{n, L}(z)\, W(dr, dz)\\
     &\qquad= a + \frac{2^{-\half}}{\sqrt L} \int_0^L B_z\, dz + \frac{1}{\sqrt L} W([0, t] \times [0, L]) + \sum_{n=1}^\infty \sqrt{\frac{2}{L}} \,(C_n + D_n),
     \end{align*}
     while
    \begin{align*}
     &I_0(a, t, x_i) + v(t, x_i) \\
     &\qquad= a  + \frac{2^{-\half}}{\sqrt L} \int_0^L B_z\, dz +\frac{1}{\sqrt L} W([0,t]\times [0,L])\\
     &\qquad\qquad+ \sum_{n=1}^\infty  g_{n, L}(x_i) e^{-\frac{\pi^2}{L^2} n^2 t}\int_0^L g_{n,L}(z)\left(a+2^{-\half B_z}\right) dz \\
     &\qquad\qquad +\sum_{n=1}^\infty  g_{n, L}(x_i) \int_0^t \int_0^L  e^{-\frac{\pi^2}{L^2} n^2 (t-r)} g_{n, L}(z)\, W(dr, dz)\\
     & \qquad = a + \frac{2^{-\half}}{\sqrt L} \int_0^L B_z\, dz + \frac{1}{\sqrt L} W([0, t] \times [0, L])\\
     & \qquad\qquad + \sum_{n=1}^\infty \cos\left(\frac{n \pi}{L} x_i\right) \sqrt{\frac{2}{L}} \, (C_n + D_n).
     \end{align*}
We now go back to \eqref{rem6.1.15-(*2)} and do the (random) change of variables $a \to \alpha$, where
\beqn
     \alpha = a + \frac{2^{-\half}}{\sqrt L} \int_0^L B_z\, dz + \frac{1}{\sqrt L} W([0, t] \times [0, L]) + \sum_{n=1}^\infty \sqrt{\frac{2}{L}} \, (C_n + D_n),
     \eeqn
to see that \eqref{rem6.1.15-(*2)} is equal to
\begin{align}
\label{rem6.1.15-(*3a)}
     &\ E\left[\int_\R d\alpha\, f_0(\alpha) f_1\left(\alpha + \sum_{n=1}^\infty \left(\cos\left(n \pi \frac{x_1}{L} \right)  - 1\right) \sqrt{\frac{2}{L}}\,  (C_n + D_n)\right)\right.\notag\\
     &\left. \qquad \times\, \cdots \times f_m\left(\alpha + \sum_{n=1}^\infty
     \left(\cos\left(n \pi \frac{x_m}{L} \right)  - 1\right) \sqrt{\frac{2}{L}} \, (C_n + D_n)\right)\right].     
 \end{align}
  
    Define
 \beqn
  Y(x) := \sum_{n=1}^\infty \left(\cos\left(n \pi \frac{x}{L}\right) - 1\right) \sqrt{\frac{2}{L}}\, (C_n + D_n).
\eeqn
By \eqref{rem6.1.15-(*4)} and the Fourier-Wiener series representation \eqref{rd08_03e1-bis}, $x \mapsto Y(x)$ has the same law as $x \mapsto 2^{-\half} Z_x$, where $(Z_x)$ is a standard Brownian motion. Therefore, \eqref{rem6.1.15-(*3a)} (and hence, \eqref{rem6.1.15-(*2)}) is equal to
\beqn
    \int_\R d\alpha\, f_0(\alpha)\, E\left[f_1(\alpha + 2^{-\half} Z_{x_1}) \, \cdots\, f_m(\alpha + 2^{-\half} Z_{x_m})\right],
    \eeqn
which is the expression \eqref{rem6.1.15-(*3)}. This proves \eqref{rem6.1.15-(*1)} and completes the proof of the theorem.
  \qed
\vskip 12pt

\subsection{Dirichlet boundary conditions with a linear drift}
\medskip

We will see in Section \ref{new-5.4.3} that the fact that the limit law of the solution to \eqref{ch1'.HD-ch5-markov} exists and is the law of a multiple of the Brownian bridge is a special case of general results of \cite{funaki83} and \cite{H-S-V-W-2005}. Here, we extend the  results 
of Subsection \ref{rd03_07ss1} to linear SPDEs with the operator $\cL$ defined in  \eqref{operator-calL}. In particular, we identify the limit laws for these SPDEs. 


Let $\cL$ be the partial differential operator given in \eqref{operator-calL}. For $a, b \in \R$ with $b \neq 0$, we let $\mA$ be the partial differential operator
\beq
\label{def-tilde-L}
    \mA = \frac{1}{b^2}\left(\frac{\partial^2}{\partial x^2} - a^2\right),
\eeq
so that $\cL=\frac{\partial}{\partial t}- \mA$.  

Consider the SPDE
\beq
\label{s5.3-new(*1)}
     \cL v(t,x) = \sqrt{2}\, \dot W(t,x),\quad     t > 0,\ x \in\, ]0, 1[, 
\eeq
subject to vanishing Dirichlet boundary conditions, and with a vanishing initial condition.

\begin{thm}
\label{s5.3-new-t*1}
Let $(v(t,x),\, (t, x) \in \R_+ \times [0, 1])$ be the random field solution of \eqref{s5.3-new(*1)}. As $t \to \infty$, the law of $v(t,*)$ converges to a probability measure $\tilde\mu$ on $(\bD, \cB_{\bD})$, which can be described in the following way:
\smallskip

 \begin{quote} On $(\bD, \cB_{\bD}, \tilde\mu)$, the coordinate process $(\pi_x,\, x \in [0,1])$ is Gaussian with mean $0$, finite variance and covariance operator $\tilde C = - \mA^{-1}$,
that is, $\tilde C(-\mA f) = f$  for $f \in H_0^2([0,1])$ and $-\mA \tilde C(g) = g$ for $g \in L^2([0,1])$. The covariance kernel $C(x,y)$ is given in \eqref{s5.3-new(*2a)} 
below.
\end{quote}
\smallskip
\end{thm}


\begin{proof} 
We recall that the {\em covariance kernel}\index{covariance!kernel}\index{kernel!covariance} of $(\pi_x,\, x \in [0,1])$ is
\beqn
    C(x,y) = \int_\bD  \pi_x(\omega) \pi_y(\omega)\,  \tilde\mu(d\omega),
    \eeqn
where $\pi_x(\omega) = \omega(x)$ is the coordinate process on $\bD$, and its {\em covariance operator}\index{covariance!operator}\index{operator!covariance} is
\beqn
   \tilde C(f)(x) = \int_0^1 C(x,y) f(y)\, dy,\quad    f \in \cC_0^\infty(]0,1[).
   \eeqn

From the expression \eqref{s5.3-new(*1)green} of
the Green's function of the operator $\cL$, we obtain
\beqn
   v(t,x) = \sqrt{2} \int_0^t \int_0^1 e^{-\frac{a^2}{b^2}(t-s)} G\left(\frac{t-s}{b^2}; x, z\right) W(ds, dz),
   \eeqn
and the covariance of $v(t,x)$ and $v(t,y)$ is
\begin{align*}
   E[v(t,x) v(t,y)] &= 2  \int_0^t ds \int_0^1 dz\, e^{-2\frac{a^2}{b^2}(t-s)} G\left(\frac{t-s}{b^2}; x, z\right )G\left(\frac{t-s}{b^2}; y, z\right)\\
   &= 2  \int_0^t ds\, e^{-2\frac{a^2}{b^2}(t-s)} G\left(2\frac{t-s}{b^2}; x, y\right).
   \end{align*}
Therefore, $\lim_{t\to\infty} E[v(t,x) v(t,y)] = C(x,y)$, where
\beqn
     C(x,y) = 2  \int_0^\infty ds\, e^{-2\frac{a^2}{b^2} s} G\left(\frac{2s}{b^2}; x, y\right)
      = b^2 \int_0^\infty dr\, e^{- a^2 r} G(r; x, y).
      \eeqn
This expression is often called the resolvent of the Brownian motion absorbed at the boundaries, or the $a$-potential density when the exponential is evaluated at $-a r$. For $a = 0$, this is the $0$-resolvent. Using formula \eqref{ch5-added-1-4} with $L=1$, we see that
\begin{align}
\label{s5.3-new(*2a)}
  C(x,y)&= b^2 \sum_{n=1}^\infty \int_0^\infty dr\, e^{-(\pi^2 n^2 + a^2) r}  e_{n,1}(x) e_{n,1}(y)\notag\\
     &= b^2 \sum_{n=1}^\infty\, \frac{e_{n,1}(x) e_{n,1}(y)}{\pi^2 n^2 + a^2}.    
     \end{align}
    Let $U = (U(x),\, x \in [0, 1])$ be a Gaussian process with mean zero and covariance kernel $E[U(x) U(y)] =  C(x,y)$. Then $E[(U(x) - U(y))^2] = C(x,x) - 2 C(x,y) + C(y,y)$, and this is bounded by the same expression for the covariance \eqref{ch5-added-1-s2-2-bis}, which is equal to the variance of an increment of the Brownian bridge, hence it is bounded above by $\vert x - y \vert$. By Kolmogorov's continuity criterion Theorem \ref{app1-3-t1},
  $U$ has a Hölder-continuous version and we let $\tilde\mu$ be the law on $(\bD, \B_{\bD})$ of $U$. Then \eqref{s5.3-new(*2a)} shows that the finite-dimensional distributions of $v(t, *)$ converge to those of $\tilde\mu$, and the equicontinuity property can be checked as in the proof of Theorem \ref{ch5-added-1-s2-t1}.

We recall that $H_0^2([0,1])$ is the set of $f \in L^2([0,1])$ such that if $f(x) = \sum_{n=1}^\infty f_n\, e_{n,1}(x)$, then $\sum_{n=1}^\infty (1 + n^2)^2 f_n^2 < \infty$.   

The operator $\mA$ is well-defined on $H_0^2([0,1])$, and for $f \in H_0^2([0,1])$,  if $f(x) = \sum_{n=1}^\infty f_n\,  e_{n,1}(x)$, then
   \beqn
    \mA  f(x) = -\frac{1}{b^2} \sum_{n=1}^\infty f_n (\pi^2 n^2 + a^2)\, e_{n,1}(x),
    \eeqn
so $\mA  f \in L^2([0,1]$. For $g \in L^2([0,1])$,  if $g(x) = \sum_{n=1}^\infty g_n\, e_{n,1}(x)$, then
\beqn
     \tilde C(g)(x) = b^2 \sum_{n=1}^\infty \frac{g_n}{\pi^2 n^2 + a^2}\, e_{n,1}(x),
     \eeqn
so $\tilde C(g) \in H_0^2([0,1])$. Clearly,  $\tilde C(- \mA  f) = f$ and $\mA  \tilde C(g) = g.$ This completes the proof of the theorem.
\end{proof}
\bigskip

\subsection{The stochastic heat equation with a fractional Laplacian}\label{rd03_06ss2}

We end this section by studying the limit distribution for another linear SPDE, this time on $\rek$ and with a fractional Laplacian.

For $a>0$, let $H^a(\R^k)$ be the set of functions $f\in L^2(\R^k)$ such that
 \beqn
 \int_{\rek} \vert \cF f(\xi) \vert^2 (1 + \vert \xi \vert^a)^2\, d \xi < \infty.
 \eeqn
This is equal to the space $\cH^{a,2}$ of Bessel potentials (for which the requirement is
     $\int_{\rek} \vert \F f(\xi) \vert^2 (1 + \vert \xi \vert^2)^a\, d \xi < \infty$).

Let $\mA _{\ttf} = -(-\Delta)^{a/2} - b^2$ and $\cL = \frac{\partial}{\partial t} -  \mA _{\ttf}$.
We consider, in spatial dimension $k$, the SPDE
\beq
\label{s5.3-new-bis(*1)}
     \cL u(t,x)  = \sqrt{2}\, \dot W(t,x),\quad       t > 0,\quad x \in \rek,     
     \eeq
with vanishing initial condition.


Let $G_a$ be the fundamental solution associated to $-(-\Delta)^{a/2}$ on $\rek$ given in \eqref{fs-frac}. Then the fundamental solution associated to $\mA _{\ttf}$ is
\beq
\label{s5.3-new-bis(*2)}
  G(t,x) := e^{-b^2t} G_a(t,x)    
  \eeq
and for $a > k$, as in \eqref{ch3-sec3.5(*1aa)}, the random field solution to \eqref{s5.3-new-bis(*1)} is
\beqn
    u(t,x) = \sqrt{2} \int_0^t \int_{\rek} e^{-b^2(t-s)} G_a(t-s, x-z)\, W(ds, dz).
    \eeqn

    \begin{prop}
  \label{s5.3-new-bis-p*1}  
  Suppose that $a > k$ and $b \neq 0$. Let $\dot W$ be space-time white noise on $\R_+ \times \rek$. As $t \to \infty$, the law of $u(t,*)$ converges weakly on compact sets to a centred Gaussian probability measure $\tilde\mu_{\ttf}$ on $\left(\cC(\rek), \cB_{\cC(\rek)}\right)$. The covariance operator of $\tilde\mu_{\ttf}$ is $\tilde C = - \mA _{\ttf}^{-1} : L^2(\rek) \to H^a(\rek)$, where $\mA _{\ttf} : H^a(\rek) \to L^2(\rek)$, and its covariance kernel is
  \beqn
    C_{a,b} (x,y) = G_{a,b}(x-y) = \frac{1}{b} G_{a,1}(b(x-y)),
   \eeqn
where $\cF^{-1} G_{a,1}(\xi) = (1 +\vert \xi \vert^a)^{-1}$ and $\cF^{-1} G_{a,b} = (b^2 +\vert \xi \vert^a)^{-1}$.
\end{prop}

\begin{proof}
The covariance of $u(t,x)$ and $u(t,y)$ is
\begin{align*}
  E[u(t,x) u(t,y)] &= 2 \int_0^t  ds \int_{\rek} dz\, e^{- 2 b^2(t-s)} G_a(t-s, x-z) G_a(t-s, y-z)\\
    &= 2 \int_0^t  ds\, e^{- 2 b^2(t-s)}  G_a(2(t-s), x-y) \\
    &=: C_t(x-y).
    \end{align*}
Using Fubini's theorem, the (inverse) Fourier transform 
of $C_t(x)$ is
 $  2 \int_0^t  ds e^{- 2 b^2 s} e^{- 2 s \vert \xi \vert^a}$.
As $t \to \infty$, this converges to
\beq
\label{noacabamai(*2)}
      \F^{-1} C(\xi) := 2 \int_0^\infty e^{- 2 b^2 s}  e^{- 2 s \vert \xi \vert^a}\, ds
         = \frac{1}{b^2 +\vert \xi \vert^a},
\eeq
which is finite. Since $a > k$ and $b \neq 0$, this function is integrable over $\rek$. Therefore, the limiting covariance kernel $C(x, y) = \lim_{t\to\infty} C_t(x,y)$ is
\beqn
  C(x, y) = C(x-y) = 2 \int_0^\infty  e^{- 2 b^2 s}  G_a(2 s, x-y) \, ds
      = \F \left(\frac{1}{b^2 +\vert \cdot \vert^a}\right)(x-y),
      \eeqn
$C(0) < \infty$ and $C$ is uniformly continuous over $\rek$.

Further,
\beqn
      C(x) = G_{a,b}(x) = \frac{1}{b} G_{a,1}(b x),\quad   x \in \rek,
      \eeqn
where $\cF G_{a,b} = (b^2 +\vert \xi \vert^a)^{-1}$.
   The covariance operator of $\mu_{\ttf}$ is $ \tilde C(f)(x) = \int_{\rek} C(x,y) f(y)\, dy$,
so $\F(\tilde C(f))(\xi) = \cF f(\xi) (b^2 +\vert \xi \vert^a)^{-1}$
and $\tilde C$ is the inverse of $-\mA _{\ttf}$, since the Fourier multiplier of $-\mA _{\ttf}$ is $\vert \xi \vert^a + b^2$.

    Let $(U(x),\, x \in \rek)$ be the centred Gaussian process with covariance kernel $ C$. Then
    \begin{align*}
    E[(U(x) - U(y))^2] &= \int_{\rek} \frac{2\left(1 - e^{-i \xi (x-y)}\right)}{b^2 +\vert \xi \vert^a}\, d \xi\\
    & = \int_{\rek} \frac{2(1 - \cos(\xi (x-y)))}{b^2 +\vert \xi \vert^a}\, d \xi.
    \end{align*}
Proceeding as in the proof of Lemma \ref{ch3-sec3.5-lA} (a), we see that this is bounded above by $C \vert x - y \vert^{(a-k) \wedge 2}$
 (with a logarithmic correction if $a = 2 + k$). It follows from Kolmogorov's continuity criterion Theorem \ref{app1-3-t1} that $U$ has a locally Hölder-continuous version  with exponent $\beta \in\, ]0, 1\wedge (a-k)/2[$, and the same statements are true for $x \mapsto u(t,x)$ (with the constant $C$ not depending on $t$).
Let $\tilde\mu_{\ttf}$ be the law of $U$ on $(\cC(\rek), \B_{\cC(\rek)})$. By \eqref{noacabamai(*2)}, the finite-dimensional distributions of $u(t,*)$ converge to those of $U$ as $t \to \infty$. We obtain the property of equicontinuity on compacts subsets of $\rek$ as in the proof of Theorem \ref{ch5-added-1-s2-t1}.
This establishes the weak convergence on compact subsets of $\rek$ of the law of $u(t, *)$ to $\tilde\mu_{\ttf}$.
\end{proof}

\begin{remark}
\label{s5.3-new-bis-r*1a} 
(a) If $k = 1$, $a = 2$ (heat equation) and $b = 1$, then
\beqn
G_{2,1}(x) = \half e^{- \vert x \vert}\quad {\text{and}}\quad C_{2,1} (x,y) = \half e^{- \vert x - y \vert},
\eeqn
 so $(U(x),\, x\in\re)$ is a strictly stationary Ornstein-Uhlenbeck process with mean $0$, variance $1/2$ and parameter $1$.

   (b) If $k = 1$, then for $1 < a < 2 $
  and $b = 1$, the covariance of $\tilde\mu$ is that of the Gaussian process $(U(x),\, x \in \R)$, where
    $U(x) = \int_\R H_{a,1}(x-y)\, W(dy)$,
 $W$ is white noise on $\R$ and $\cF H_{a,1}(\xi) = (1 +\vert \xi \vert^a)^{-\half}$. This is a real-valued process because
\beqn
   E[U(x)^2] = \int_\R H_{a,1}^2(x-y)\, dy = \int_\R d\xi\, \vert \F H_{a,1}(\xi)\vert^2  = \int_\R d\xi\, \frac{1}{1 +\vert \xi \vert^a} < \infty.
   \eeqn
   \end{remark}

\section{Reversible measures}
   \label{ch5-ss-in-rev}


We begin this section by defining the notion of {\em reversible measure} associated to a Markov process and describing its relationship with the invariant measure. In Section \ref{new-5.4.1}, we give two finite-dimensional examples namely, the Ornstein-Uhlenbeck process on $\R^n$, $n \geq 1$ (see Lemma \ref{nbr-0b} and Proposition \ref{nbr-p-0c}) and, in Proposition \ref{nbr-1}, a diffusion process on $\R^n$ with a gradient-type drift.
In Section \ref{new-5.4.2}, we exhibit reversible measures for the infinite-dimensional analogues of these processes: the solution to a linear stochastic heat equation (Proposition \ref{nbr-4}) and to a nonlinear stochastic heat equation with additive gradient-type nonlinearity (Theorem \ref{nbr-t1}). Finally, in Section \ref{new-5.4.3}, for a stochastic heat equation with a gradient-type drift term and additive space--time white noise, we identify the reversible measure with a bridge-type measure.
\bigskip

\noindent{\em Definition and first properties}
\medskip

We adopt the setting of  Section \ref{ch5-added-1-s1} and consider a continuous Markov process $(Y_t)$ with transition function $(P_t)$.
First, we note that the defining property  \eqref{ch5-added-1-s2-1} of an invariant measure\index{invariant measure}\index{measure!invariant} $\mu$ for $(Y_t)$ (or $(P_t)$) can be formulated equivalently as follows:
\beq
\label{nbr-(*a2)}
    \int_\cs P_t \psi(g) \, \mu(dg) =  \int_\cs \psi(g)\,  \mu(dg), 
\eeq
for all bounded and continuous functions $\psi: \cs \to \re$.

Now, we give the notion of reversible measure.
\begin{def1}
\label{nbr-0aa} A measure $\nu$ on $\cB_S$ is {\em reversible}\index{reversible measure}\index{measure!reversible} for $(P_t)$ $($or $(Y_t))$ if,
 for all bounded continuous functions $\varphi, \psi: \cs \to \re$  and all $t \geq 0$,
\beq
\label{nbr-(*a3)}
    \int_{\cs}  \varphi(g)\, P_t\psi(g) \, \nu(dg) =  \int_{\cs} P_t\varphi(g) \,  \psi(g)\, \nu(dg).
   \eeq
\end{def1}
 We notice that in the above definition, it is equivalent to require \eqref{nbr-(*a3)} for all bounded measurable functions.
 
 It will be convenient to identify the Markov process $(Y_t)$ with its canonical version $(\tilde Y_t)$, so that we can make use of the probability measures $\tilde \bP^g$ and $\tilde \bP^\nu$ defined in Subsection \ref{ch5-added-1-s1}. For simplicity of notation, we remove the tildes and we write $Y_t$, $P^g$ and $P^\nu$ instead of $\tilde Y_t$, $\tilde \bP^g$ and $\tilde \bP^\nu$, respectively.

The next lemma (together with Remark \ref{nbr-r-reversible} below) explains the choice of the term ``reversible" and shows that a reversible measure is also an invariant measure.

\begin{lemma}
\label{nbr-0a}

(a) If $\nu$ is reversible for $(Y_t)$, then $\nu$ is an invariant measure for $(Y_t)$.

(b) The law $\nu$ of $Y_0$ is reversible for $(Y_t)$ if and only if for all $t \geq 0$, the laws of $(Y_0, Y_t)$ and of $(Y_t, Y_0)$ under $P^\nu$ are the same, that is, for all bounded and continuous functions $\varphi, \psi: \cs \to \R$,
\beqn
   E^\nu[\varphi(Y_0)\, \psi(Y_t)] = E^\nu[\varphi(Y_t)\, \psi(Y_0)].
\eeqn
Equivalently,
\beq\label{rd03_10e1}
     \int_{\cs}  \varphi(g)\,  E^g[\psi(Y_t)]\, \nu(dg) =  \int_{\cs} E^g[\varphi(Y_t)]\,  \psi(g)\, \nu(dg).
\eeq
\end{lemma}

\begin{proof} (a) If we set $\varphi \equiv 1$ in \eqref{nbr-(*a3)}, then we obtain \eqref{nbr-(*a2)}. This proves (a).

(b) Property \eqref{nbr-(*a3)} is equivalent to
\beq
\label{nbr-(*a3a0)}
     E^{\nu} [\varphi(Y_0) E^{Y_0}[\psi(Y_t)]] = E^{\nu} [E^{Y_0}[\varphi(Y_t)] \psi(Y_0)] ,   
     \eeq
 where $E^{\nu}$ (respectively $E^{Y_0}$) denotes the mathematical expectation with respect to the probability measure $P^\nu$ (respectively $P^{Y_0}$).
     
     By the Markov property, the left-hand side of \eqref{nbr-(*a3a0)} is equal to
\beqn
    E^\nu[\varphi(Y_0)  E^\nu[\psi(Y_t) \mid Y_0]] = E^\nu[\varphi(Y_0)  \psi(Y_t)]
    \eeqn
and in the same way, the right-hand side is equal to
$
     E^\nu[\varphi(Y_t) \psi(Y_0)].
$
  Therefore, \eqref{nbr-(*a3)} implies that the laws of $(Y_0, Y_t)$ and of $(Y_t, Y_0)$ under $P^\nu$ are the same. Conversely, if $(Y_0, Y_t)$ and $(Y_t, Y_0)$ have the same law under $P^\nu$, then $E^\nu[\varphi(Y_0)\psi(Y_t)] = E^\nu[\varphi(Y_t)\psi(Y_0)]$, so \eqref{nbr-(*a3a0)} and \eqref{nbr-(*a3)} hold. This proves (b).
\end{proof}

\begin{remark}
\label{nbr-r-reversible}
The term {\em reversible} is justified by the following (easily checked) fact: If the law $\nu$ of $Y_0$ is reversible for $(Y_t,\, t\in[0,T])$, then under $P^\nu$, the laws of the processes  $(Y_t,\, t\in[0,T])$ and $(Y_{T-t},\, t\in[0,T])$ are the same.
\end{remark}

\subsection{The finite-dimensional case $\cs = \R^n$}
\label{new-5.4.1}

Let $Y = (Y_t,\, t \in \R_+)$ be the solution of an $n$-dimensional system of SDEs
\beq
\label{nochnoch(*1)}
    dY_t = b(Y_t)\, dt + \sigma(Y_t)\, dB_t,    
    \eeq
subject to a given initial condition $Y_0$, where $(B_t = (B^1_t,\dots, B^n_t))$ is an $n$-dimensional standard Brownian motion, $b = (b_1,\dots, b_n) : \R^n \to \R^n$ and $\sigma = (\sigma_{i, j}) : \R^n \to \R^{n \times n}$ are Lipschitz continuous functions.  We write $Y(y) = (Y_t(y),\, t \in \R_+)$ to indicate that we are considering the solution to \eqref{nochnoch(*1)} with initial condition $Y_0 = y\in\re^m$. The family of processes $Y(y)$, $y\in \R^n$, is defined on a filtered probability space $(\Omega, \F, P, (\cF_t))$ on which $(B_t)$ is an $(\cF_t)$-standard Brownian motion.

   For a bounded Borel function $f: \R^n \to \R$, define
\beq
\label{nochnoch(*2)}
       P_t f(y) = E[f(Y_t(y))],      
\eeq
and denote $P_t(y, A) := P_t 1_A (y)$.
Then $(P_t)$ is a transition function that has the strong Markov property  \cite[Theorem 5.4.20, p.322]{ks}. If $\mu$ is the law of $Y_0$ and if $Y_0$ is independent of $(B_t)$, then for all $t \in \R_+$, the law of $Y_t(Y_0)$ is $\mu$ if and only if $\mu$ is an invariant measure for $(P_t)$.
\vskip 12pt

\noindent{\em Reversibility for the Ornstein-Uhlenbeck processes on $\R^n$}
\vskip 2mm

   We consider here the ($n$-dimensional) Ornstein-Uhlenbeck processes on $\R^n$, $n \geq1$, and exhibit their reversible measures. We begin with the case $n=1$.
\vskip 12pt

\begin{lemma}
\label{nbr-0b} Let $(Y_t)$ be an Ornstein-Uhlenbeck process on $\R$, that is, $(Y_t)$ satisfies the SDE
\beqn
   dY_t = - a Y_t\, dt + \sqrt{2}\, dB_t,\quad t>0,
   \eeqn
where $(B_t)$ is a standard Brownian motion and $a > 0$. Then $\nu_0 = {\rm N}(0,a^{-1})$ is a reversible (and invariant) probability measure for $(Y_t)$.
\end{lemma}

\begin{proof} If $Y_0 \equiv y$, then the law of $Y_t$ is ${\rm N}(ye^{- a t}$, $a^{-1}(1 - e^{-2 a t}))$ (see e.g. \cite[p. 358, Example 6.8]{ks}). Therefore, if $Y_0$ is ${\rm N}(0,a^{-1})$  and is independent of $(B_t)$, then the law of $Y_t$ is that of
    $e^{- a t} Y_0 + \sqrt{a^{-1}(1 - e^{-2 a t})}\, Z$,
where $Z$ is ${\rm N}(0,1)$ and is independent of $Y_0$. This law is ${\rm N}(0, a^{-1})$, therefore the joint law of $(Y_0, Y_t)$ under $\nu_0$ is Normal, with mean zero. Since under $\nu_0$, $Y_0$ and $Y_t$ have the same variance (and covariance $e^{- a t} / a)$, $(Y_0, Y_t)$ and $(Y_t, Y_0)$ have the same variance-covariance matrix, so this joint law is also the joint law of $(Y_t, Y_0)$ under $\nu_0$. The conclusion follows from Lemma \ref{nbr-0a}.
   \end{proof}

We now consider the case $n \geq 1$.

\begin{prop}
\label{nbr-p-0c}
Let $B_t = (B_t^1,\dots, B_t^n)$, where the $(B_t^j)$ are independent (real-valued) standard Brownian motions, $j = 1, \dots, n$. Let $a = (a_1,\dots, a_n)$, with $a_j > 0$,  $j = 1, \dots, n$, and let $A(n)$ be the $n \times n$ diagonal matrix ${\rm diag}(a_1,\dots, a_n)$. Let $(Y_t)$ satisfy the system of SDEs
\beq
\label{nbr-(*6)}
   dY_t = - A(n)\, Y_t\, dt +  \sqrt{2}\, dB_t,\quad t>0,    
   \eeq
or equivalently,
\beqn
   d Y_t^j = - a_j\, Y_t^j\, dt + \sqrt{2}\, d B_t^j,\quad t>0,\      j = 1, \dots, n.
   \eeqn
Then for all bounded continuous functions $\varphi, \psi : \R^n \to \R$,  \eqref{nbr-(*a3)} holds with the product measure
\beq
\label{nbr-(*7a)}
   \mu_n(dx) = \nu^1(dx_1)\, \cdots\, \nu^n(dx_n),   
   \eeq
where $\nu^j$ is the ${\rm N}(0, a_j^{-1})$ probability measure. Hence $\mu_n$ is a reversible (and invariant) probability measure for $(Y_t)$.
\end{prop}

\begin{remark}
\label{after-nbr-p-0c}
If $a_j = j^2 \pi^2$, then $\mu_n$ is the law of the projection onto the subspace spanned by $(e_{1, 1},\dots, e_{n, 1})$ of the Brownian bridge, where the $e_{j, 1}$ are defined in the line that follows \eqref{ch5-added-1-4}. 
\end{remark}

\noindent{\em Proof of Proposition \ref{nbr-p-0c}.}\ 
Using Lemma \ref{nbr-0b}, one checks directly that \eqref{nbr-(*a3)} holds for bounded continuous functions of the form $\psi(x_1,\dots, x_n) = \psi_1(x_1)\, \cdots\, \psi_n(x_n)$, and this extends by linearity and density to all bounded continuous $\varphi$, $\psi$. The conclusion follows.
\qed
\bigskip

 \noindent{\em Reversibility for diffusion processes on $\R^n$ with gradient-type drift}
 \medskip

When a linear (Gaussian) diffusion has a reversible probability measure, then this is often also the case when the diffusion is modified by introducing a suitable gradient-type drift. We illustrate this statement in Proposition \ref{nbr-1}
 below. First, we give some auxiliary results, in which we denote by $C^*$ the transpose of a matrix $C$.

\begin{lemma}
\label{nbr-l-12} 
Fix $T > 0 $. Let  $A $ and  $C $ be  $n\times n $ matrices of real numbers and let $f: \R^n \to \R^n $ be a Lipschitz continuous function. We assume that $C$ is invertible. Given an  $\R^n$-valued standard Brownian motion  $B = (B_t = (B^1_t,\dots, B^n_t),\, t \in \R_+) $ on a filtered probability space  $(\Omega, \cF, P, (\cF_t))$, consider the two systems of SDEs
\beq
\label{nochnochnoch(*1)}
     dX_t = (A X_t + f(X_t))\, dt  + C dB_t, \quad                X_0 = x_0,     
     \eeq
and
\beq
\label{nochnochnoch(*2)}
     dZ_t = A Z_t dt  + C dB_t,  \quad                 Z_0 = x_0,     
     \eeq
$t \in\, ]0, T]$. Let $\mu_X$ (respectively $\mu_Z$) be the law under $P$ of $X = (X_t,\, t \in [0, T])$ (respectively $Z = (Z_t,\, t \in [0, T]))$ on $\cC^n := \cC([0, T], \R^n)$. Then
\begin{align}
\label{nochnochnoch(*3)}
     \frac{d\mu_X}{d\mu_Z}(Z) &= \exp\left( \int_0^T f(Z_t)^* (C C^*)^{-1}\, dZ_t\right.\notag\\
     &\left.\qquad\qquad - \half \int_0^T f(Z_t)^* (C C^*)^{-1} (f(Z_t) + 2 A Z_t)\, dt\right).           
\end{align}
\end{lemma}

\begin{proof}  Let $\tilde P$ be the probability measure such that
\begin{align*}
     \frac{d\tilde P}{d P} &= \exp\left(- \int_0^T f(X_t)^* (C^*)^{-1} \, dB_t\right.\\
      &\left.\qquad\qquad - \half \int_0^T f(X_t)^* (C C^*)^{-1} f(X_t)\, dt\right)
     \end{align*}
     (see \cite[Chapter 3, Section 5, Corollary 5.16, p.~200]{ks}).
     
By Girsanov's theorem for SDEs, if we set
   $ \tilde B_t = B_t +  C^{-1} \int_0^t f(X(s))\, ds$,
then under $\tilde P$, the process $\tilde B = (\tilde B_t,\, t \in [0, T])$ is a standard Brownian motion. In addition,
\beqn
    dX_t = A X_t\, dt + C\, d\tilde B_t, \quad     X(0) = x_0.
    \eeqn
Therefore, the law  under $\tilde P$ of  $X$ is identical to the law $\mu_Z$ under $P$ of  $Z$. In particular, for a bounded and continuous function  $\psi : \cC^n \to \R$,
\begin{align*}
    &\int_{\cC^n}  \psi(w)\, \mu_X(dw) = E_P[\psi(X)] = E_{\tilde P}\left[\psi(X)\, \frac{dP}{d\tilde P}\right] \\
     &\qquad = E_{\tilde P}\left[\psi(X) \exp\left( \int_0^T f(X_t)^* (C^*)^{-1}\, dB_t\right.\right.\\
     &\qquad  \left.\left.\qquad\qquad + \half \int_0^T f(X_t)^* (C C^*)^{-1} f(X_t)\, dt)\right)\right] \\
   &\qquad = E_{\tilde P}\left[\psi(X)  \exp\left(\int_0^T f(X_t)^* (C^*)^{-1}  (d\tilde B_t - C^{-1} f(X_t)\, dt)\right.\right.\\
   &\qquad  \left.\left.\qquad\qquad + \half \int_0^T f(X_t)^* (C C^*)^{-1} f(X_t)\, dt\right)\right] \\
  &\qquad  = E_{\tilde P}\left[\psi(X) \exp\left(\int_0^T  f(X_t)^* (C^*)^{-1} \,  (C^{-1}\, dX_t - C^{-1}  A X_t\, dt)\right.\right.\\
    &\qquad \left.\left.\qquad\qquad-\half \int_0^T f(X_t)^* (C C^*)^{-1} f(X_t)\, dt\right)\right].
  \end{align*}
As in the proof of Proposition \ref{ch2'-s3.4.2-p2}, we can check that the exponential is a Borel function of the sample paths of $X$. Since the law under $\tilde P$ of $X$ is the same as the law under $P$ of $Z$, this is equal to
\begin{align*}
    & E_{P}\left[\psi(Z) \exp\left( \int_0^T  f(Z_t)^* (C^*)^{-1} (C^{-1}\, dZ_t - C^{-1}  A Z_t\, dt)\right.\right.\\
    &\left.\left. \qquad\qquad\qquad\qquad - \half \int_0^T f(Z_t)^* (C C^*)^{-1} f(Z_t)\, dt \right)\right] \\
     &\qquad= E_{P}[\psi(Z) F(Z)] = \int_{\cC^n}  \psi(w) F(w)\, \mu_Z(dw),
     \end{align*}
where
\begin{align*}
    F(Z) & = \exp\left( \int_0^T  f(Z_t)^* (C C^*)^{-1}\, dZ_t \right.\\
     &\qquad\qquad \left. - \half \int_0^T f(Z_t)^* (C C^*)^{-1} (f(Z_t) + 2 A Z_t)\, dt\right).
    \end{align*}
This implies \eqref{nochnochnoch(*3)} and completes the proof.
    \end{proof}

In the next lemma, we compute the term $\frac{d\mu_X}{d\mu_Z}(Z)$ in the particular instance where the function $f$ is a gradient.

For two $n \times n$ matrices $A$ and $B$, define
$$A : B = \sum_{i, j =1}^n A_{i, j} B_{i, j} = {\rm{Tr}}(AB^*),$$
where ${\rm{Tr}}(A)$ is the trace of the matrix $A$. This operation is an inner product, sometimes called the {\em Frobenius inner product}.\index{Frobenius inner product}\index{inner product!Frobenius}

\begin{lemma}
\label{nbr-l-13}
Let $A$, $C$, $(X_t)$ and $(Z_t)$ be as in \eqref{nochnochnoch(*1)} and \eqref{nochnochnoch(*2)}. Suppose that $f(z) = C C^* \nabla \tilde f(z)$, with $\tilde f \in \cC^2(\R^n, \R)$. Define
\beq
\label{nochnochnoch(*4)}
    \Phi(z) = - \half (C C^* : D^2 \tilde f(z) + (\nabla \tilde f(z))^*\, C C^*\, \nabla \tilde f(z) + 2 (\nabla \tilde f(z))^* A z),   
    \eeq
where $ D^2 \tilde f$  is the Hessian matrix of $ \tilde f$, that is, 
\beqn
  D^2 \tilde f(z) =  \left((D^2 \tilde f)_{i, j}(z) = \frac{\partial^2 \tilde f(z)}{\partial x_i \partial x_j}\right).
\eeqn
Let $ \mu_X$ and $\mu_Z$ be as in Lemma \ref{nbr-l-12}. Then
\beqn
      \frac{d\mu_X}{d\mu_Z}(w) = \exp\left(\tilde f(w(T)) - \tilde f(w(0)) + \int_0^T \Phi(w(t))\, dt\right),  
\eeqn
 $ w \in \cC([0, T], \R^n)$.
\end{lemma}
      
\begin{proof}
 Using the It\^o formula, we have
 \beqn
     \tilde f(Z(T)) - \tilde f(Z(0)) =  \int_0^T \nabla \tilde f(Z_t) \cdot dZ_t + \half \int_0^T [C C^* : D^2 \tilde f(Z_t)]\, dt.
     \eeqn
Since $f(Z_t)^* (C C^*)^{-1} = (\nabla \tilde f(z))^*$, \eqref{nochnochnoch(*3)} can be written
\begin{align*}
   \frac{d\mu_X}{d\mu_Z}(Z) &= \exp\left(\tilde f(Z(T)) - \tilde f(Z(0)) - \half \int_0^T [C C^* : D^2 \tilde f(Z_t)]\, dt\right.\\
   &\left.\qquad\qquad - \half \int_0^T (\nabla \tilde f(Z_t))^*\,  (C C^*\, \nabla \tilde f (Z_t) + 2 A Z_t)\, dt\right).
   \end{align*}
Using the definition of $\Phi$, this becomes
\beqn
    \frac{d\mu_X}{d\mu_Z}(Z) = \exp\left(\tilde f(Z(T)) - \tilde f(Z(0)) + \int_0^T \Phi(Z_t)\, dt\right).
    \eeqn
This completes the proof.
\end{proof}

We now address the issue of reversibility for a diffusion process with gradient-type drift.\index{gradient-type drift}\index{drift!gradient-type}

\begin{prop}
    \label{nbr-1}
Let $b: \R^n \to \R$ be a $\cC^1$-function bounded from above with Lipschitz partial derivatives $\frac{\partial b}{\partial x_k}$, $k = 1, \dots, n$.
Let $\lambda = (\lambda_1,\dots, \lambda_n)$ with $\lambda_k > 0$ for $k = 1, \dots, n$, and set $A(\lambda) = {\rm diag}(\lambda_1,\dots, \lambda_n)$. Consider the system of SDEs
\beq
\label{eqforx}
   dX_t = - A(\lambda) X_t\, dt + \nabla b(X_t)\, dt + \sqrt{2}\, dB_t,\quad t>0,
   \eeq
where $(B_t = (B_t^1,\dots,B_t^n))$ is an $n$-dimensional standard Brownian motion. Equivalently,
\beqn
   d X_t^k = - \lambda_k X_t^k\, dt + \frac{\partial b}{\partial x_k}(X_t^1, \dots, X_t^n)\, dt + \sqrt{2}\, d B_t^k,\quad t>0,\     k = 1, \dots, n.
   \eeqn
Let $\mu_n$ be the probability measure  on $\R^n$ defined in \eqref{nbr-(*7a)} with $a_k = \lambda_k$, $k=1,\,\ldots,\ n$, that is, $\mu_n$ is the law of
\beq
\label{nbr-(*7aa)-vector}
       \left(\frac{Z_1}{\sqrt{\lambda_1}}, \dots, \frac{Z_n}{\sqrt{\lambda_n}}\right),  
       \eeq
where the $(Z_k)$ are i.i.d. ${\rm N}(0,1)$ random variables, and let $\nu_n$ be the probability measure on $\R^n$ given by
\beq
\label{nbr-(*7)}
     \nu_n(dx) = \frac{1}{M}\, e^{b(x)}\, \mu_n(dx),    
\eeq
where $M= \int_{\re^n}e^{b(x)}\, \mu_n(dx)$. Then $\nu_n$ is a reversible measure for $(X_t)$.
\end{prop}

The presence of $\nabla b$ in \ref{eqforx} is a restriction. For instance, even when $n=1$ and any function $\beta: \R \to \R$ is of the form $\beta(x) = b'(x)$, where $b: \R \to \R$ is any antiderivative of $\beta$, this function $b$ might not satisfy the hypotheses of the theorem, in particular, it may not be bounded above.

\medskip

\noindent{\em Proof of Proposition \ref{nbr-1}.}\
Suppose first that $b \in \cC^2(\R^n)$. We check that \eqref{nbr-(*a3)} holds. Indeed, let 
\beqn
   c = \frac{1}{M} \left(\frac{\det (A(\lambda))}{(2\pi)^{n}}\right)^\half. 
\eeqn
Then, for bounded continuous functions $f, g: \R^n \to \R$, 
\beq\label{rd03_09e1}
    \int_{\R^n} \nu_n(dx)\, f(x) P_t g(x) = c \int_{\R^n} dx\, e^{-\half \vert \sqrt{A(\lambda)}\, x \vert^2} e^{b(x)} f(x)  E^x[g(X_t)].
\eeq
Let $(X_t)$ be as in \eqref{nochnochnoch(*1)}, with $A:=-A(\lambda)$, $f=\nabla b$ and $C = {\rm{diag}}\left(\sqrt{2},\dots, \sqrt{2}\right)$ there (hence $\tilde f = \half b$), and  let
$(Z_t)$ be defined in \eqref{nochnochnoch(*2)}. By Lemma \ref{nbr-l-13}, the right-hand side of \eqref{rd03_09e1} is equal to
\begin{align*}
     &c \int_{\R^n} dx\, e^{-\half \vert \sqrt{A(\lambda)}\, x \vert^2} e^{b(x)} \\
     &\qquad\qquad \times E^x\left[f(x) g(Z_t) \exp\left(\half (b(Z_t) - b(Z_0)) + \int_0^t \Phi(Z_s)\, ds\right)\right],
     \end{align*}
where $\Phi$ is defined in \eqref{nochnochnoch(*4)}. Therefore, moving $e^{b(x)}$ inside the expectation, the previous display becomes
\begin{align*}
     &c \int_{\R^n} dx\, e^{-\half \vert \sqrt{A(\lambda)}\, x \vert^2} \\
     &\qquad\qquad \times E^x\left[f(Z_0) g(Z_t) \exp\left(\half (b(Z_t) + b(Z_0))
      + \int_0^t \Phi(Z_s)\, ds\right)\right].
      \end{align*}
Because $(Z_t)$ is reversible under $\mu_n$ by Proposition \ref{nbr-p-0c}, this is equal to
\begin{align*}
    & c \int_{\R^n} dx\, e^{-\half \vert \sqrt{A(\lambda)}\, x \vert^2}\\
    &\qquad \quad \times E^x\left[f(Z_t) g(Z_0) \exp\left(\half (b(Z_0) + b(Z_t)) + \int_0^t \Phi(Z_{t - s})\, ds\right)\right] \\
    &\quad = c \int_{\R^n} dx\, e^{-\half \vert \sqrt{A(\lambda)}\, x \vert^2} e^{b(x)}\\
     &\qquad \quad \times E^x\left[f(Z_t) g(Z_0) \exp\left(\half (b(Z_t) - b(Z_0)) + \int_0^t \Phi(Z_{s})\, ds\right)\right]\\
   &\quad = c \int_{\R^n} dx\, e^{-\half \vert \sqrt{A(\lambda)}\, x \vert^2} e^{b(x)} E^x[f(X_t)] \, g(x) \\
   &\quad  = \int_{\R^n} \nu_n(dx)\,P_t f(x) \, g(x) ,
   \end{align*}
where, in the second equality, we have again used Lemma  \ref{nbr-l-13}.

    In order to handle the case where the function $b$ is only $\cC^1$, bounded above and  with Lipschitz partial derivatives $\frac{\partial b}{\partial x_k}$, $k = 1,\dots, n$, we use a standard approximation argument, as follows.

 Let $\varphi_\ell(x) = \ell^n \varphi(\ell x)$, $\ell\ge 1$, where $\varphi \geq 0$, $\varphi$ is even, $\cC^\infty$ with compact support and $\int_{\R^n} \varphi(x)\, dx = 1$. Let $b_\ell = b * \varphi_\ell$. The functions $b_\ell$ are $\cC^\infty$, and one easily checks the following properties:
 \begin{description}
 \item{(i)} The functions $\nabla b_\ell$ are uniformly Lipschitz continuous,\index{uniformly Lipschitz continuous}\index{Lipschitz!continuous, uniformly} that is, there exists a constant $K<\infty$ such that for all $x,y\in\re^n$ and $\ell\ge 1$, we have
 \beqn
 \left\vert \nabla b_\ell(x) - \nabla b_\ell(x)\right\vert\le K|x-y|.
 \eeqn
 \item{(ii)}  For fixed $x \in \R^n$, $b_\ell (x) \to  b(x)$ and $\nabla b_\ell (x) \to \nabla b(x)$ as $\ell\to\infty$.
 \end{description}

  For any fixed $\ell\ge 1$, let $(X_t^{(\ell)})$ be the solution of the equation \eqref{eqforx} where $\frac{\partial b}{\partial x_k}$ there is replaced by $\frac{\partial b_\ell}{\partial x_k}$, $k = 1,\dots, n$. Using Lemma \ref{aproxode}, we see that
  \beqn
  \lim_{\ell\to\infty}\left\vert X_t^{(\ell)}- X_t\right\vert =0, \quad {\text{uniformly in}}\ t\in[0,T],\ {\text{a.s.}}
  \eeqn
  By the first part of the proof,
    \begin{align}
    \label{nochnoch(*5)}
   &\int_{\R^n} dx\, e^{-\half \vert \sqrt{A(\lambda)}\, x \vert^2} e^{b_\ell(x)} f(x) E^x[g(X^{(\ell)}_t)]\notag\\
   &\qquad\qquad= \int_{\R^n} dx\, e^{-\half \vert \sqrt{A(\lambda)}\, x \vert^2} e^{b_\ell(x)} E^x[f(X^{(\ell)}_t)]\,  g(x).    
   \end{align}
For fixed $x$, $E^x[f(x) g(X^{(\ell)}_t)] \to E^x[f(x) g(X_t)]$ by dominated convergence, and $b_\ell(x) \to b(x)$. Since the $b$ and $b_\ell$ are uniformly bounded above, we can apply dominated convergence to \eqref{nochnoch(*5)}  and conclude that
 \begin{align}
    &\int_{\R^n} dx\, e^{-\half \vert \sqrt{A(\lambda)}\, x \vert^2} e^{b(x)} f(x) E^x[g(X_t)]\notag\\
    & \qquad\qquad= \int_{\R^n} dx\, e^{-\half \vert \sqrt{A(\lambda)}\, x \vert^2} e^{b(x)} E^x[f(X_t)] \, g(x),
   \end{align}
that is, $\nu_n$ is a reversible measure for $(X_t)$. This completes the proof.
 \qed
   \medskip


\subsection{The linear case}
\label{new-5.4.2}

After this brief interlude in finite dimensions, we return to the linear stochastic heat equation \eqref{s5.3-new(*1)}, that is,
\beqn
\cL v(t,x) = \sqrt{2}\, \dot W(t,x), \quad t>0,\ x\in\,]0,1[,
\eeqn
with  $\cL = \frac{\partial}{\partial t} - \frac{1}{b^2}\left(\frac{\partial ^2}{\partial x^2} - a^2 \right)$ and Green's function $G_{a,b}(t;x,y)$ given in \eqref{s5.3-new(*1)green},
subject to vanishing Dirichlet boundary conditions and initial condition $v_0 \in \bD$.

 Let $V = L^2([0, 1])$ endowed with a CONS $(e_n,\, n \geq 1)$, denoted $(e_{n, L})$ (with $L = 1$) in \eqref{ch5-added-1-4}, and
denote by $\Pi_n$  the projection operator onto the linear subspace of $V$ spanned by $\{e_1, \dots,e_n\}$ (in the notation of Section \ref{ch1'-s60}, $\Pi_n$ is $\Pi_{V_n}$). We have seen in Theorem \ref{s5.3-new-t*1} that the law of the solution to \eqref{s5.3-new(*1)} has a limit measure,  which is then necessarily an invariant measure. We show here that this measure is in fact reversible.

\begin{prop}
\label{nbr-4}
For $j \in \N^*$, let $a_j = \pi^2j^2 + a^2$, and set
\beq
\label{nbr-(*10)}
       Y = b \sum_{j=1}^\infty \frac{Z_j}{\sqrt{a_j}}\, e_j  ,   
       \eeq
where the $(Z_j)$ are i.i.d.~${\rm N}(0,1)$ random variables. Then:
\begin{description}
 \item{(a)} The law of Y is $\tilde\mu$ given in Theorem \ref{s5.3-new-t*1}, and $\Pi_{n}(Y)$ converges to $Y$ a.s. uniformly on $[0, 1]$, and in $L^2(\Omega; V)$.
 \item{(b)} Let $v = (v(t,x))$ be the solution of the linear SPDE \eqref{s5.3-new(*1)}. Then $\tilde\mu$ is a reversible (and invariant) probability measure for $v$.
 \end{description}
 \end{prop}
 
Before proving this proposition, we introduce some notations.
\bigskip

\noindent{\em Transformations from $\R^n$ to $\bD$ and back}
\medskip

Define $T_n:\re^n\rightarrow \bD$ by
\beq
\label{nbr-(*7aa1)}
     T_n(y)(*) = \sum_{j=1}^n y_j \, e_j(*),\quad y\in \re^n,    
     \eeq
and  $\Phi_n: \bD \rightarrow \re^n$ by
\beq
\label{nbr-(*7aa2)}
      \Phi_n(w) = (\langle w, e_1 \rangle_V,\dots, \langle w, e_n\rangle_V), \quad w \in \bD.   
      \eeq
so that
\beq
\label{nbr-(*7aa3)}
       T_n(\Phi_n(w)) = \Pi_n w,\quad       \Phi_n(T_n(y)) = y.       
\eeq
\medskip

\noindent{\em Proof of Proposition \ref{nbr-4}.}
(a) A direct calculation of the mean and covariance kernel of $Y$ leads to formula \eqref{s5.3-new(*2a)}. Because $a_j \geq (\pi j)^2$, the uniform convergence is part of the Paley-Wiener construction of Brownian motion 
(see \cite[Chapter 1, Problem 3]{ito-mckean} or \cite[Chapter 9, p. 168]{khosh-2007} and Section \ref{stwn-s2} of this book). The convergence in $L^2(\Omega; V)$ follows from the convergence of the series in \eqref{s5.3-new(*2a)}.

   (b) We are going to check that for bounded continuous functions $\varphi, \psi : \bD \to \R$, \eqref{nbr-(*a3)} (in fact \eqref{rd03_10e1}) holds (with $\cs$ there replaced by $\bD$), that is,
   \beq
   \label{nbr-(*a3)-new}
   \int_\bD \varphi(w) E^w[\psi(v(t,*))]\, \tilde\mu(dw) =  \int_\bD E^w[\varphi(v(t, *)] \psi(w)\, \tilde\mu(dw).
   \eeq

 For $j \in \N^*$, consider the Ornstein-Uhlenbeck process $Y_t(j)$ satisfying
\beqn
     d Y_t(j) = - a_j Y_t(j)\, dt + \sqrt{2}\, dW_t(e_j),\quad    Y_0(j) = \langle v_0, e_j \rangle_V.
     \eeqn
As in the proof of Proposition \ref{ch5-added-1-s2-p1},  for all $(t,x)$,
\beq
\label{nbr-(*12e)}
    v(t, x) = \sum_{j=0}^\infty Y_t(j) e_j(x),     
    \eeq
where the series converges a.s. and in $L^2(\Omega)$. Let
\beqn
   v_n(t, *) = \Pi_n v(t, *) = \sum_{j=1}^n Y_t(j) e_j(*).
   \eeqn
Note that $\tilde\mu$-a.s., for fixed $t$, $v_n(t, *) \to v(t, *)$ uniformly on $[0, 1]$. Since the terms of this sum are independent processes with values in $\cC([0, 1])$, by the Itô-Nisio theorem \cite{ito-nisio-1968} (see also \cite{walsh-1967}, which applies to processes indexed by $[0, T] \times [0, 1]$), the series in \eqref{nbr-(*12e)} converges a.s., uniformly in  $t \in [0, T]$ and $x \in [0, 1]$. We note in passing that the Itô-Nisio theorem theorem is also a variation on the Paley-Wiener construction.

Let $(Y^n_t = (Y_t(1),\dots, Y_t(n)))$. From \eqref{nbr-(*12e)} and  \eqref{nbr-(*7aa2)}, we see that 
\beq\label{nbr-(*7aab4)}
    Y^n_t = \Phi_n(v(t, *))\quad \text{and}\quad    T_n(Y^n_t) = \Pi_n v(t, *) = v_n(t, *).    
\eeq
Let $\mu_n$ be the probability measure on $\R^n$ defined in 
\eqref{nbr-(*7a)}. By Proposition \ref{nbr-p-0c}, this is a reversible measure for $(Y^n_t)$. Therefore, for bounded continuous functions $f, g: \R^n \to \R$,
\beq\label{nbr-(*7aa)}
   \int_{\R^n} f(y) E^y[g(Y^n_t)]\, \mu_n(dy) =  \int_{\R^n} E^y[f(Y^n_t)] g(y)\, \mu_n(dy).   
   \eeq

Let $\varphi, \psi : \bD \to \R$ be bounded and continuous. Consider $f, g: \R^n \to \R$ defined by
\beqn
   f(y) = \varphi(T_n(y)), \qquad g(y) = \psi(T_n(y)),
\eeqn
so that
\beqn
   f(Y^n_t) = \varphi(\Pi_n v(t, *)) \quad\text{and} \quad f(\Phi_n(w)) = \varphi(\Pi_n(w))
 \eeqn
by \eqref{nbr-(*7aab4)} and \eqref{nbr-(*7aa3)}. 

In order to emphasize the dependence of the solution $v$ on the initial condition $v_0$, we write $v_{v_0}$ instead of $v$. Let $(P_t)$ be the transition function of $(v_{v_0}(t, *))$ and let  $(P_{n, t})$ denote the transition function of $(Y^n_t)$. By part (a), the law of $\Phi_n(w)$ under $\tilde\mu$ is $\mu_n$, and the law of $v_n(t, *)$ under $P_t(v_0, \cdot)$ is the same as the law of $T_n(Y^n_t)$ under $P_{n, t}(\Phi_n(v_0), \cdot)$, therefore,
\begin{align*}
    P_{n, t} f (\Phi_n(v_0)) &= E^{\Phi_n(v_0)}[f(Y^n_t)] = E^{\Phi_n(v_0)}[\varphi(T_n(Y^n_t))]\\
    & = E^{v_0}[\varphi(\Pi_n v(t, *))] = (P_t (\varphi \circ \Pi_n))(w).
 \end{align*}
In other words, for $w \in \bD$,
\beqn
    P_{n, t} f (\Phi_n(w)) = (P_t (\varphi \circ \Pi_n)) (w).
    \eeqn
Therefore, by \eqref{nbr-(*7aa)},
\beqn
       \int_{\bD} f(\Phi_n(w))\, E^{\Phi_n(w)}[g(Y^n_t)]\, \tilde\mu (dw)
            =  \int_{\bD} E^{\Phi_n(w)}[f(Y^n_t)]\, g(\Phi_n(w))\, \tilde\mu (dw).
            \eeqn
 Equivalently,
 \begin{align*}
   &\int_{\bD} \varphi(\Pi_n(w))\, E^{\Phi_n(w)}[\psi(T_n(Y^n_t))]\, \tilde\mu (dw)\\
       &\qquad\quad=  \int_{\bD} E^{\Phi_n(w)}[\varphi(T_n(Y^n_t))]\, \psi(\Pi_n(w))\, \tilde\mu (dw).
        \end{align*}
By \eqref{nbr-(*7aab4)}, this is equivalent to
\begin{align}
\label{nbr-(*9ac)}
      &\int_{\bD} \varphi(\Pi_n(w))\, E^{w}[\psi(\Pi_n v(t, *))]\, \tilde\mu (dw)\notag\\
        &\qquad\quad=  \int_{\bD} E^{w}[\varphi(\Pi_n v(t, *))]\, \psi(\Pi_n(w))\,\tilde\mu (dw).      
\end{align}

By part (a), under $\tilde\mu$, as $n \to \infty$, $\Pi_n(w) \to w$ in $\bD$ for $\tilde\mu$-a.a. $w$. By dominated convergence,
\beqn
\label{nbr-(*9ad)}
    E^{w}[\psi(\Pi_n v(t, *))] \to E^{w}[\psi(v(t, *))],    
    \eeqn
and similarly for $\varphi$.

From \eqref{nbr-(*9ac)} and dominated convergence, we obtain
\beqn
   \int_{\bD} \varphi(w)\, E^{w}[\psi( v(t, *))]\, \tilde\mu(dw)
        =  \int_{\bD} E^{w}[\varphi(v(t, *))]\, \psi(w)\, \tilde\mu(dw),
        \eeqn
which is the reversibility property \eqref{nbr-(*a3)-new} of $\tilde\mu$ for the process $v$.
\qed
 
\subsection{The case of a gradient-type drift}\label{rd03_14ss1}

Next, we identify the reversible (and invariant) measure for the nonlinear SPDE that we now describe. Let $\cL$ be the partial differential operator given in \eqref{operator-calL} and let $\beta: \R \to \R$ be a $\cC^1$-function bounded from above and with Lipschitz continuous derivative. Consider the SPDE
\beq
\label{nbr-(*1aa)}
       \cL u(t,x) = \beta'(u(t,x)) + \sqrt{2}\, \dot W(t,x),\quad   t > 0,\   x \in\, ]0, 1[,   
       \eeq
subject to vanishing Dirichlet boundary conditions and some suitable initial condition $u_0 \in \bD$. 
   By Theorem \ref{nbr-2aa}, $(u(t, *))$ is a Markov process on $\bD$ with transition function given in \eqref{ch5-added-1-2-2}.

For $w \in \bD$, define
 \beq
 \label{nbr-(*1)}
    B(w) = \int_0^1 \beta(w(x))\, dx,     
    \eeq
so that $B: \bD \to \R$. The next theorem is an extension of Proposition \ref{nbr-4} (b) to the nonlinear SPDE \eqref{nbr-(*1aa)}.

\begin{thm}
\label{nbr-t1}
Let $(u(t, x))$ be the random field solution to \eqref{nbr-(*1aa)}.
Let  $\tilde\mu$ be the law on  $(\bD, \cB_\bD)$ described in Theorem \ref{s5.3-new-t*1} and let $\mu$ be the probability measure on $(\bD, \cB_\bD)$ defined by
\beq
\label{nbr-(*3)}
     \mu(dw) =  \frac{1}{M}\,  e^{B(w)}\,  \tilde\mu (dw),   
     \eeq
where
\beq
\label{nbr-(*4)}
    M = \int_\bD  e^{B(w)}\, \tilde\mu(dw).     
    \eeq
Then $\mu$ is a reversible (and invariant) probability measure for $(u(t,*))$.
\end{thm}

\begin{remark}
\label{nbr-r1a}
(a) With a little additional effort, one can show that the statement remains true if $\beta'$ is only assumed to be locally Lipschitz continuous (see \cite[Example 2 p. 247]{zabczyk-1989}).

    (b) We will see in Lemma \ref{nbr-l0ab} below that the Frechet derivative $B'(w)$ of $B$ at $w$ is the linear functional associated to an element $D(w)$ of $V$, that is, $B'(w)(h) = \langle D(w), h \rangle_V$, and we will see  in the proof of Lemma \ref{nbr-l4a} that $D: \bD \to V$ is Lipschitz continuous: $\Vert D(w_1) - D(w_2) \Vert_V \leq C \Vert w_1 - w_2 \Vert_\bD$. With the fact that $B$ is bounded above, these are the only properties needed to obtain the conclusion of Theorem \ref{nbr-t1} for the SPDE
    \beqn
     \cL u(t,x) = D(u(t, *))(x) + \sqrt{2}\, \dot W(t,x),\quad   t > 0,\    x \in\, ]0, 1[.
     \eeqn
\end{remark}
\medskip

\noindent{\em Preliminaries}
\medskip

The proof of the Theorem \ref{nbr-t1} requires some preliminary work that we now develop.
Recall the transformations $T_n$ and $\Phi_n$ defined in \eqref{nbr-(*7aa1)} and \eqref{nbr-(*7aa2)}, with the properties given in \eqref{nbr-(*7aa3)}, and let $B$ be as defined in \eqref{nbr-(*1)}. For $y = (y_1,\dots, y_n) \in \R^n$, define
\beq
\label{nbr-(*12a)}
      B_n(y) = B(T_n(y)) = \int_0^1 \beta\left(\sum_{j = 1}^n y_j e_j(x)\right) dx,     
      \eeq
and for $w \in \bD$, define $D: \bD \to V$ by
\beqn
    D(w)(x) =\beta'(w(x)) =  \sum_{j = 1}^\infty \langle \beta'(w(*)), e_j \rangle_V\,  e_j(x).
    \eeqn

    Recall also that  $\beta: \R \to \R$ is a $\cC^1$-function bounded from above, with Lipschitz continuous derivative. We define the mapping $\tilde \beta': V \to V$ (associated to $\beta'$) by
\beq
\label{nbr-betatilde}
\tilde \beta'(w)(x) = \beta'(w(x)),\quad x\in[0,1],
\eeq
and we set
\beq
\label{nbr-compo}
   \Pi_n(\beta'(\Pi_n(w)))(x) := (\Pi_n \circ \tilde \beta' \circ \Pi_n)(w)(x).
\eeq

 \begin{lemma}
 \label{nbr-l0ab}
The following properties hold:
\begin{description}
\item{(a)} For  $j=1,\dots, n$,
\beqn
\frac{\partial B_n}{\partial y_j} (y) = \langle \beta'(T_n(y)), e_j\rangle_V,
\eeqn
and
\begin{align*}
      \nabla B_n(y) &= (\langle \beta'(T_n(y)), e_1 \rangle_V,\dots, \langle \beta'(T_n(y)), e_n \rangle_V)\\
       &= \Phi_n(\Pi_n(\beta'(T_n(y)))).
      \end{align*}
          \item{(b)} The linear functional $h \mapsto \langle D(w), h \rangle_V$ (from $\bD$ to $\re$) associated to $D$ is equal to $B'(w)$, the Frechet derivative of $B: \bD \to \re$.
 \item{(c)} Let $\Pi_n(\beta'(\Pi_n(w)))$ be defined as in \eqref{nbr-compo}. Then
 \beqn
 \Pi_n(D(\Pi_n(w))) = \Pi_n(\beta'(\Pi_n(w))) = T_n(\nabla B_n(\Phi_n(w))).
 \eeqn
 \end{description}
\end{lemma}

\begin{proof}
(a) The first statement is due to the fact that
\beqn
    \frac{\partial B_n}{\partial y_j}(y) = \int_0^1 \beta'\left(\sum_{j = 1}^n y_j e_j(x)\right) e_j(x)\, dx = \langle \beta'(T_n(y)), e_j \rangle_V.
    \eeqn
The second statement is a rewrite of the first.
\smallskip

    (b) Let $D: \cC(\bD, \R) \to \cC^*(\bD, \R)$ be the functional defined by 
 \beqn
    D(w)(h) = \int_0^1 \beta'(w(x)) h(x) dx,\qquad h \in \bD. 
 \eeqn
 We want to check that
    \beq
    \label{nbr-(*D1)}
     \lim_{\Vert h \Vert_\bD \to 0} \frac{B(w+h) - B(w) - D(w)(h)}{\Vert h \Vert_\bD} = 0.   
     \eeq
Indeed,
\begin{align*}
&B(w+h) - B(w) - D(w)(h)\\
   &\qquad =\int_0^1 \left(\beta(w(x) + h(x)) - \beta(w(x)) - \beta'(w(x)) h(x)\right) dx\\
    &\qquad= \int_0^1 dx \left(\int_{0}^{1} dy  \left(\beta'(w(x) + y h(x)) - \beta'(w(x))\right) h(x)\right),
    \end{align*}
    therefore, by the Lipschitz property of $\beta'$,
    \begin{align*}
    &\vert B(w+h) - B(w) - D(w)(h)\vert\\
     &\qquad \leq C \left(\int_0^1 h^2(x)\,dx\right) \left(\int_{0}^{1} y\, dy\right) = \frac{C}{2} \Vert h \Vert_V^2\le  \frac{C}{2} \Vert h \Vert_{\bD}^2 .
    \end{align*}
Thus, \eqref{nbr-(*D1)} holds.
\smallskip

    (c) By (a) and \eqref{nbr-(*7aa3)}, $\nabla B_n(\Phi_n(w)) =  \Phi_n(\Pi_n(\beta'(\Pi_n(w))))$, and applying $T_n$ to both sides and using again \eqref{nbr-(*7aa3)}, we obtain $T_n(\nabla B_n(\Phi_n(w))) = \Pi_n(\beta'(\Pi_n(w)))$, because $\Pi_n(\Pi_n(w)) = \Pi_n(w)$.   
\end{proof}

\medskip

\noindent{\em A finite-dimensional approximation}
\medskip

Recall that $a_j = \pi^2 j^2 + a^2$. For a sequence $(Z_j, \ j\ge 1)$ of  i.i.d. ${\rm N}(0,1)$ random variables, let $\tilde\mu$ be the law on $\mathbb{D}$ of
\beq
\label{nbr-(*3aa)}
   b\,\sum_{j=1}^\infty \frac{Z_j}{\sqrt{a_j}}\, e_j,  
   \eeq
 which we considered in Proposition \ref{nbr-4}, and for $n \geq 1$, let $\mu_n$ be the law on $\R^n$ of
\beq
\label{nbr-(*3aa)-bis}
 b\,\left(\frac{Z_1}{\sqrt{a_1}},\dots, \frac{Z_n}{\sqrt{a_n}}\right).
   \eeq
As a consequence of Proposition \ref{nbr-1}, we have the following.

\begin{cor}
\label{nbr-c1a}
Let $(Y_t)$ be the solution of
\beqn
    dY_t =  - b^{-2}A(n) Y_t\, dt + \nabla B_n(Y_t)\, dt + \sqrt{2}\, dB_t,
    \eeqn
    where $A(n) = {\rm{diag}}(a_1,\ldots,a_n)$, and let $\nu_n$ the probability measure on $\R^n$ defined by
\begin{align}
\label{nbr-(*17ac0)}
      \nu_n(dy) &= \frac{1}{M_n}\, e^{B_n(y)}\, \mu_n(dy)\notag\\
              &=  \frac{1}{M_n} \exp\left(\int_0^1 \beta\left(\sum_{j = 1}^n y_j e_j(x)\right) dx\right) \mu_n(dy),   
              \end{align}
where $M_n= \int_{\re^n} e^{B_n(y)}\, \mu_n(dy)$. Then $\nu_n$ is a reversible probability measure for $(Y_t)$.
\end{cor}

\begin{proof}
 By Proposition \ref{nbr-1}, it suffices to check that for $j = 1,\dots, n$, the function $y\mapsto \frac{\partial B_n}{\partial y_j} (y)$ is Lipschitz. For this, let $c$ be the Lipschitz constant of $\beta'$. By Lemma \ref{nbr-l0ab} (a), for $x, y \in \R^n$,
 \begin{align*}
   & \left\vert \frac{\partial B_n}{\partial y_j} (x) - \frac{\partial B_n}{\partial y_j} (y) \right\vert
    = \left\vert \langle \beta'(T_n(x)) - \beta'(T_n(y)), e_j \rangle_V\right\vert \\
    &\qquad\qquad \leq \Vert \beta'(T_n(x)) - \beta'(T_n(y)) \Vert_V
    \leq c\, \Vert T_n(x) - T_n(y) \Vert_V \\
    &\qquad\qquad = c\,  \Vert T_n(x - y) \Vert_V = c\,  \vert x - y \vert,
   \end{align*}
by definition of $T_n$ in \eqref{nbr-(*7aa1)}.
 \end{proof}

\begin{remark}
\label{nbr-r1ab}
The law of $T_n(\cdot)$ (on $\bD$) under $\mu_n$ is the same as the law of $\Pi_n(\cdot)$  under $\tilde\mu$, and the law of $\Phi_n(\cdot)$ (on $\re^n$) under $\tilde\mu$ is $\mu_n$.  By \eqref{nbr-(*12a)},  $B_n(y) = B(T_n(y))$, therefore,
\begin{align*}
     M_n &= \int_{\R^n} e^{B_n(y)}\, \mu_n(dy) = \int_{\R^n} e^{B(T_n(y))}\, \mu_n(dy)\\
     & = \int_\bD e^{B(\Pi_n(w))}\, \tilde\mu(dw).
     \end{align*}
For $n \geq 1$, let $\gamma_n$ be the probability measure on $\bD$ defined by
\beq\label{rd03_11e1}
     \gamma_n(dw) = \frac{1}{M_n}\, e^{B(\Pi_n(w))}\,  \tilde\mu(dw).
 \eeq
Then for $\varphi: \bD \to\R$ bounded and continuous,
\begin{align*}
    \int_{\R^n}  \varphi(T_n(y))\, \nu_n(dy) &= \int_{\R^n}  \varphi(T_n(y))\,\frac{1}{M_n}\, e^{B(T_n(y))}\, \mu_n(dy)\\
    & = \int_\bD \varphi(\Pi_n(w))\, \frac{1}{M_n}\, e^{B(\Pi_n(w))}\, \tilde\mu(dw)\\
    & =  \int_\bD \varphi(\Pi_n(w))\, \gamma_n(dw),
    \end{align*}
that is, the law of $T_n(\cdot)$ under $\nu_n$ is the same as the law of $\Pi_n(\cdot)$ under $\gamma_n$, and similarly, the law of $\Phi_n(\cdot)$ under $\gamma_n$ is $\nu_n$.
\end{remark}

Recall the definition of $\tilde \beta': V \to V$ (associated to $\beta'$) given in
\eqref{nbr-betatilde}.
A straightforward computation shows that $\tilde\beta'$ inherits the Lipschitz continuity property of $\beta'$,
with the same Lipschitz constant. Recall also the definition of $\Pi_n(\beta'(\Pi_n(w)))$ given in \eqref{nbr-compo}, which is used in the next lemma.

\begin{lemma}
\label{nbr-l4a}
Let $\cL$ be the partial differential operator given in \eqref{operator-calL}.
Let $\beta: \R \to \R$ be as above, $B$ as defined in \eqref{nbr-(*1)}, and $W^n$ as defined in \eqref{proj-noise}. Let $(u_n(t,x))$ be the solution of
\beq
\label{nbr-(*11)}
    \cL u_n(t,x) =  \Pi_n(\beta'(\Pi_n(u_n(t,*))))(x) + \sqrt{2}\, \dot W^n(t,x),    
    \eeq
with vanishing Dirichlet boundary conditions and initial condition $\Pi_n(v_0)$. Let $\tilde\mu$ be the law defined just before Corollary \ref{nbr-c1a}. 
Then for all bounded and continuous functions $\varphi, \psi: \bD \to \R$,
\begin{align}
\label{nbr-(*14ab)}
              &\int_{\bD} \varphi(\Pi_n(w))\, E^{\Pi_n(w)}[\psi(u_n(t,*))]\,  e^{B(\Pi_n(w))}\,  \tilde\mu (dw)\notag\\
              &\qquad  =  \int_{\bD} E^{\Pi_n(w)}[\varphi(u_n(t, *))]\, \psi(\Pi_n(w))\, e^{B(\Pi_n(w))}\, \tilde\mu (dw).    
\end{align}
\end{lemma}

\begin{proof}
First, we notice that for all $n\ge 1$, the mapping
\beqn
w \mapsto   \Pi_n(\beta'(\Pi_n(w)))(x) = (\Pi_n\circ\tilde\beta'\circ \Pi_n)(w)(x)
\eeqn
from $\mathbb{D}$ into $\re$,
  is Lipschitz continuous and the Lipschitz constant does not depend on $n$ or $x$, that is, there exists a constant $C$ such that, for any $n \in \N^*$, $w_1, w_2\in \bD$ and $x \in [0, 1]$,
\beq
\label{campnou(*1)}
   \vert (\Pi_n\circ\tilde\beta'\circ \Pi_n)(w_1)(x) - (\Pi_n\circ\tilde\beta'\circ \Pi_n)(w_2)(x) \vert \leq C \Vert w_1 - w_2 \Vert_\bD.
   \eeq
Indeed, applying \eqref{nbr-betatilde}, we see that the left-hand side is bounded above as follows:
\begin{align}
\label{nbr-L1}
    &\vert \Pi_n(\tilde\beta'\circ\Pi_n(w_1) - \tilde\beta'\circ\Pi_n(w_2))(x) \vert \notag\\
    &\qquad=  \left\vert \sum_{j=1}^n\, \langle \tilde\beta'(\Pi_n(w_1)) - \tilde\beta'(\Pi_n(w_2)), e_j \rangle_V e_j(x)\right \vert\notag\\
    &\qquad \leq C\, \Vert \tilde\beta'(\Pi_n(w_1)) - \tilde\beta'(\Pi_n(w_2)) \Vert_V.
    \end{align}
Because $\tilde\beta'$ is Lipschitz, up to a multiplicative constant, this is bounded above by
\beq
\label{nbr-L2}
    \Vert \Pi_n(w_1) - \Pi_n(w_2) \Vert_V \leq  \Vert w_1 - w_2 \Vert_V \leq  \Vert w_1 - w_2 \Vert_\bD.
    \eeq
Therefore \eqref{campnou(*1)} holds. By Theorem \ref{ch1'-s7-t1}, the SPDE \eqref{nbr-(*11)} has a unique global solution.

We recall that the Green's function of $\cL$ is
\beq
\label{green-repe}
G_{a,b}(t;x,y)= \sum_{j=1}^\infty e^{-\frac{\pi^2 j^2 + a^2}{b^2} t} e_j(x)e_j(y) = \sum_{j=1}^\infty e^{-\alpha_j t} e_j(x)e_j(y),
\eeq
where $\alpha_j =  b^{-2} (\pi^2 j^2 + a^2)$ (see \eqref{s5.3-new(*1)green}).

Let $(Y_t(j))$ be the solution to the stochastic differential equation
\beqn
dY_t(j) = -\alpha_jY_t(j)\, dt + \sqrt 2\, dW_t(e_j),\quad Y_0(j) =0,\ j=1\ldots,n.
\eeqn
More explicitely,
\beq
    \label{nbr-(*16a)}
    Y_t(j) =  \sqrt{2} \int_0^t e^{- \alpha_j (t-s)}\, dW_t(e_j),   
\eeq
and the $(Y_t(j))$ are independent.

Define
\beq
\label{nbr-(*16ab)}
    Z_n(t, x) = \sqrt{2}\, \sum_{j=1}^n  e_j(x) \int_0^t  e^{- \alpha_j (t-s)}\, dW_t(e_j) = \sum_{j=1}^n Y_t(j)\, e_j(x).     
    \eeq
Recalling the setting and results from Section \ref{ch1'-s60}, we have 
\begin{align}
\label{nbr-(*14)}
    u_n(t,x) &=  \langle G_{a,b}(t; x, *), \Pi_n(v_0) \rangle_V\notag\\
    &\qquad + \int_0^t ds\, \left\langle G_{a,b}(t-s; x, *), \Pi_n(\beta'(\Pi_n(u_n(s,*))))(*) \right\rangle_V\notag\\
    &\qquad + \sqrt 2\int_0^t \int_0^1 G_{a,b}(t-s; x, y)\, W^n(ds, dy)\notag\\
  &=  \sum_{j=1}^n\, \langle v_0, e_j \rangle_V\,  e^{- \alpha_j t} e_j(x)\notag\\
  &\qquad + \int_0^t ds\, \sum_{j=1}^n e^{- \alpha_j (t-s)} \left\langle \beta'(u_n(s,*)), e_j \right\rangle_V\,  e_j(x)  + Z_n(t, x),     
  \end{align}
  where we have used \eqref{green-repe}.

    Write $u_n(t,x) =  \sum_{j=1}^n  A_t(j)\, e_j(x) = T_n(A_t)$, or, equivalently,
    \beq
    \label{nbr-(*17a02)}
   A_t = \Phi_n(u_n(t, *)), \quad   A_t =  (A_t(1),\dots, A_t(n)), 
   \eeq
where, by \eqref{nbr-(*16ab)},
\beqn
   A_t(j) =  e^{- \alpha_j t} \langle v_0, e_j \rangle_V  + \int_0^t ds\, e^{- \alpha_j (t-s)} \left\langle \beta'(T_n(A_s), e_j \right\rangle_V + Y_t(j).
   \eeqn
In particular, for any $j=1,\ldots, n$, $(A_t(j))$ satisfies
\begin{align*}
   dA_t(j) &= - \alpha_j \left(e^{- \alpha_j t} \langle v_0, e_j \rangle_V  + \int_0^t ds\, e^{- \alpha_j (t-s)} \left\langle \beta'(T_n(A_s), e_j \right\rangle_V\right)\, dt\\
   &\qquad  + \left\langle \beta'(T_n(A_t), e_j \right\rangle_V\, dt  -  \alpha_j Y_t(j)\, dt + \sqrt{2}\, dW_t(e_j)\\
     & = - \alpha_j\, A_t(j)\, dt + \frac{\partial B_n}{\partial \alpha_j} (A_t)\, dt + \sqrt{2}\, dW_t(e_j),
     \end{align*}
where we have used Lemma \ref{nbr-l0ab} (a). Equivalently,
\beqn
   dA_t = - A(n) A_t\, dt + \nabla B_n(A_t)\, dt + \sqrt{2}\, dW_t^n,
   \eeqn
where $A(n) = {\rm{diag}}(\alpha_1,\ldots,\alpha_n)$ and $W_t^n = (W_t(e_1),\dots, W_t(e_n))$.

    Let $\mu_n$ be the law on $\R^n$ given in \eqref{nbr-(*3aa)-bis}. 
By Proposition \ref{nbr-p-0c}, this is a reversible measure for the process $(C_t)$, where
\beqn
d C_t = - A(n) C_t\, dt + \sqrt{2}\, dW_t^n.
\eeqn

     Let  $\nu_n$ be the law on $\R^n$ defined in \eqref{nbr-(*17ac0)}.
By Corollary \ref{nbr-c1a}, this is a reversible measure for $(A_t)$. In particular, for bounded continuous functions $f, g: \R^n \to \R$,
\beq
\label{nbr-(*14aa)}
    \int_{\R^n} f(y)\, E^y[g(A_t)]\, \nu_n(dy) =  \int_{\R^n} E^y[f(A_t)]\, g(y)\, \nu_n(dy).   
    \eeq
Let $\varphi, \psi: \mathbb{D} \to \R$ be bounded and continuous. Choose
\beqn
\label{nbr-(*14aa)-bis}
       f(y) = \varphi(T_n(y)),\quad      g(y) = \psi(T_n(y)), \quad y\in\re^n.
       \eeqn
By \eqref{nbr-(*14aa)},
\beqn
      \int_{\R^n} \varphi(T_n(y))\, E^y[\psi(T_n(A_t))]\, \nu_n(dy) =  \int_{\R^n} E^y[\varphi(T_n(A_t))]\, \psi(T_n(y))\, \nu_n(dy).
      \eeqn
As in Remark \ref{nbr-r1ab}, we can write the above equality as
\begin{align*}
     &\int_{\bD} \varphi(\Pi_n(w))\, E^{\Phi_n(w)}[\psi(T_n(A_t))]\, \gamma_n(dw)\\
     &\qquad\qquad = \int_{\bD} E^{\Phi_n(w)}[\varphi(T_n(A_t))]\, \psi(\Pi_n(w))\, \gamma_n(dw) ,
 \end{align*}
where $\gamma_n$ is defined in \eqref{rd03_11e1}. By \eqref{nbr-(*17a02)}, 
\beqn
    E^y[\psi(T_n(A_t)] = E^{T_n(y)}[\psi(u_n(t, *))], 
\eeqn
and $T_n(\Phi_n(w)) = \Pi_n(w)$ by \eqref{nbr-(*7aa3)}, therefore
\begin{align*}
    &\int_{\bD} \varphi(\Pi_n(w))\, E^{\Pi_n(w)}[\psi(u_n(t, *))]\, \gamma_n(dw)\\
    &\qquad\qquad = \int_{\bD} E^{\Pi_n(w)}[\varphi(u_n(t, *))]\, \psi(\Pi_n(w))\,  \gamma_n(dw),
    \end{align*}
that is,
\begin{align*}
    &\int_{\bD} \varphi(\Pi_n(w))\, E^{\Pi_n(w)}[\psi(u_n(t, *))] \, \frac{1}{M_n} e^{B(\Pi_n(w))}\,  \tilde\mu (dw)\\
    &\qquad =  \int_{\bD} E^{\Pi_n(w)}[\varphi(u_n(t, *))]\, \psi(\Pi_n(w))\, \frac{1}{M_n} e^{B(\Pi_n(w))} \, \tilde\mu (dw).
    \end{align*}
Simplifying the common fraction, we get \eqref{nbr-(*14ab)}. The proof of the lemma is complete.
\end{proof}

\begin{remark}
\label{local-lipschitz-beta} The above lemma remains true assuming that $\beta'$ is locally Lipschitz continuous. Indeed, by Theorem \ref{ch1'-s7-t4}, \eqref{nbr-(*11)} still has a unique global solution.
\end{remark}
\medskip

\noindent{\em Convergence of the finite-dimensional approximations}
\medskip

\begin{lemma}
\label{nbr-*3a} 
Consider a sequence $(w_n,\, n\ge 1)$ of functions in $\bD$ and assume that $w_n \to w$ in $\bD$. Let
$u_{n,\Pi_n(w_n)}$, $n\ge 1$, be the solution of \eqref{nbr-(*11)} with initial condition $\Pi_n(w_n)$ and $u_{w}$ be the solution to \eqref{nbr-(*1aa)} with initial condition $w$. Then, as $n\to \infty$,
$\tilde\mu$-a.s.,
 \beqn
     u_{n,\Pi_n(w_n)}(t, *) \longrightarrow u_w(t, *)\quad {\text in}\ \bD,
     \eeqn
uniformly for $t \in [0, T]$.
\end{lemma}

\begin{proof}
Let $Y_t(j)$ be defined as in \eqref{nbr-(*16a)}. We have also defined 
\beqn
    Z_n(t, x) = \sum_{j=1}^n Y_t(j)\, e_j(x) \quad\text{and}\quad v(t,x) = \sum_{j=1}^\infty Y_t(j)\, e_j(x)
\eeqn
in \eqref{nbr-(*16ab)}and  \eqref{nbr-(*12e)}, respectively.  The random field $(v(t, x))$ is the solution of the linear SPDE
\beqn
       \cL v(t,x) = \sqrt{2}\, \dot W(t,x),\quad   t > 0,\   x \in\, ]0, 1[,
\eeqn
subject to vanishing Dirichlet boundary conditions, and vanishing initial condition.

We noted in the proof of Proposition \ref{nbr-4} that in fact,
\beq
\label{nbr-(*17a)}
     \sup_{t \in [0, T]} \Vert Z_{n}(t, *) - v(t, *) \Vert_\bD \to 0,    
     \eeq
that is, $Z_n(t,x) \to v(t,x)$, uniformly in $x\in[0, 1]$ and uniformly in $t \in [0, T]$. In addition, $\Pi_n(w_n) \to w$ in $V$, because
\begin{align}
\label{nbr-(*17b)}
    \Vert \Pi_n(w_n) - w \Vert_V &\leq \Vert \Pi_n(w_n - w) \Vert_V + \Vert \Pi_n(w) - w \Vert_V\notag\\
    & \leq \Vert w_n - w \Vert_V + \Vert \Pi_n(w) - w \Vert_V,     
    \end{align}
and both terms converge to $0$ as $n \to \infty$.

    Let $v_n(t,x) = u_{n,\Pi_n(w_n)}(t, x) - u_w(t, x)$
and
\beqn
   f_n(t) = \sup_{(s,x) \in [0, t]\times[0, 1]} \vert v_n(s,x) \vert, \quad
   g_n(t) = \sup_{(s,x) \in [0, t]\times[0, 1]}\vert Z_n(t, x) -  v(t,x) \vert.
   \eeqn
Define $F_n,\, F: \bD \to V$ by
\beq
\label{nbr-(*14abc)}
    F_n(w) = (\Pi_n\circ\tilde\beta'\circ\Pi_n)(w),\quad      F(w) = \tilde\beta'(w).     
\eeq
By \eqref{campnou(*1)}, we see that the sequence $(F_n,\, n\ge1)$ consists of uniformly Lipschitz continuous functions, that is, there is $K < \infty$ such that, for all $w_1, w_2 \in \bD$,
\beq
\label{nbr-(*14abbb)}
    \Vert F_n(w_1) - F_n(w_2) \Vert_{V} \leq K \Vert w_1 - w_2 \Vert_{\bD},\quad  {\text {for all}}\ n\ge 1.    
    \eeq

In addition, for all $w \in \bD$,
\beq
\label{nbr-(*14abe)}
   \Vert F_n(w) - F(w) \Vert_{V} \longrightarrow 0\quad  {\text as}\  n \to \infty.    
   \eeq
Indeed,
\begin{align*}
   \Vert F_n(w) - F(w) \Vert_{V}  &\leq  \Vert \Pi_n(\tilde\beta'(\Pi_n(w))) - \Pi_n(\tilde\beta'(w)) \Vert_{V} \\
   &\qquad +  \Vert \Pi_n(\tilde\beta'(w)) - \tilde\beta'(w) \Vert_{V}.
\end{align*}
The first term on the right-hand side is equal to
\beqn     \Vert \Pi_n(\tilde\beta'(\Pi_n(w)) - \tilde\beta'(w)) \Vert_{V}  
    \leq \Vert \tilde\beta'(\Pi_n(w)) - \tilde\beta'(w) \Vert_{V}   
      \leq K \Vert \Pi_n(w) - w \Vert_{V} ,
\eeqn
where in the second inequality, we have applied the Lipschitz property of $\tilde \beta'$. It follows that
\beqn
   \Vert F_n(w) - F(w) \Vert_{V}  \leq K \Vert \Pi_n(w) - w \Vert_{V} + \Vert \Pi_n(\tilde\beta'(w)) - \tilde\beta'(w) \Vert_{V},
\eeqn
and both terms converge to $0$ by definition of $\Pi_n$, proving \eqref{nbr-(*14abe)}.

    Notice that $u_{n,\Pi_n(w)}(t, x)$ solves \eqref{nbr-(*15a)} with
     $g_n = \Pi_n(w_n)$, $\zeta_n = Z_n,$ and $F_n(w)$ as given in  \eqref{nbr-(*14abc)},
    and $u_w(t, x)$ solves \eqref{nbr-(*15b)} with $g = w$, $\zeta = v$ and $F(w)$ as given in \eqref{nbr-(*14abc)}. By \eqref{nbr-(*17a)}, \eqref{nbr-(*17b)}, \eqref{nbr-(*14abbb)} and \eqref{nbr-(*14abe)}, with these choices, the assumptions of Proposition \ref{nbr-*3aa} are satisfied. Therefore,
    \beqn
   \lim_{n \to \infty} \sup_{(t, x) \in [0, T] \times [0, 1]} \vert  v_n(t,x) \vert = 0,
   \eeqn
which is the desired conclusion.
  \end{proof}
  \smallskip

\noindent{\em Proof of Theorem \ref{nbr-t1}}.\  Let $\varphi, \psi: \bD \to \R$ be bounded and continuous functions. First, we prove that for $\tilde\mu$-a.a. $w \in \bD$,
\beq
\label{nbr-(*14ae)}
     e^{B(\Pi_{n}(w))} \to e^{B(w)},   
     \eeq
 and
 \beq
 \label{nbr-(*14ad)}
    E^{\Pi_n(w)}[\varphi(u_n(t, *))] \to E^{w}[\varphi(u(t, *))].    
    \eeq
Indeed, by Proposition \ref{nbr-4} (a), for $\tilde\mu$-a.a. $w$, $\Pi_{n}(w) \to w$ in $\bD$. Fix such a $w \in \bD$ and assume that $w([0, 1]) \subset \, ]-K_0, K_0[$. Because $\beta$ is uniformly continuous on $[-K_0, K_0]$, $B(\Pi_{n}(w)) \to B(w)$, so we also have \eqref{nbr-(*14ae)}.

    In order to establish  \eqref{nbr-(*14ad)}, let $w_n = \Pi_{n}(w)$, so that $w_n \to w$ in $\bD$ for $\tilde\mu$-a.a. $w$. For $w \in \bD$ such that $\Pi_{n}(w) \to w$ in $\bD$, by Lemma \ref{nbr-*3a},
    \beqn
      u_{n,w_n}(t, *) \longrightarrow u_w(t, *)\quad {\text {in}}\ \bD,\quad     \tilde\mu\text {-a.s.}
      \eeqn
Therefore, since $\varphi$ is bounded and continuous, by dominated convergence, for $\tilde\mu$-a.a. $w$,
\beqn
    E[\varphi(u_{n,w_n}(t, *))] \longrightarrow E[\varphi(u_w(t, *))].
    \eeqn
This is equivalent to
\beqn
     E^{\Pi_n(w)}[\varphi(u_n(t, *))] \longrightarrow E^{w}[\varphi(u(t, *))],\quad {\text {for}}\ \tilde\mu\text {-a.a.}\ w \in \bD,
     \eeqn
which is \eqref{nbr-(*14ad)}.

    By Lemma \ref{nbr-l4a},
\begin{align}
\label{nbr-(*17abl)}
& \int_{\bD} \varphi(\Pi_n(w))\, E^{\Pi_n(w)}[\psi(u_n(t,*))]\,  e^{B(\Pi_n(w))}\,  \tilde\mu (dw)\notag\\
      &\qquad\qquad=  \int_{\bD} E^{\Pi_n(w)}[\varphi(u_n(t, *))]\, \psi(\Pi_n(w))\, e^{B(\Pi_n(w))}\,  \tilde\mu (dw).    
\end{align}
Using  \eqref{nbr-(*14ad)},  \eqref{nbr-(*14ae)} and dominated convergence to pass to the limit in  \eqref{nbr-(*17abl)}, we obtain
\beqn
     \int_{\bD} \varphi(w)\, E^{w}[\psi(u(t, *))]\, e^{B(w)}\, \tilde\mu (dw) =  \int_{\bD} E^{w}[\varphi(u(t, *))]\, \psi(w) e^{B(w)}\,  \tilde\mu (dw),
     \eeqn
which establishes the reversibility of $(u(t, *))$ under $\mu$.
 \qed
 \bigskip

 \noindent{\em An example}
 \medskip
 
 Consider the partial differential operator $\cL = \frac{\partial}{\partial t} - \frac{\partial^2}{\partial x^2}$ and the two SPDEs
 \beq
 \label{linear-h}
 \cL v(t,x) = \sqrt 2\, \dot W(t,x)
 \eeq
 and for $a \neq 0$,
 \beq
 \label{nonlinear-h}
  \cL u(t,x) = -a^2 u(t,x)+ \sqrt 2\, \dot W(t,x),
  \eeq
 $t>0$, $x\in\, ]0,1[$, both with vanishing Dirichlet boundary conditions and vanishing initial condition.
  
  Let  $\tilde \mu$ be the invariant measure for $(v(t,x))$ given in Proposition \ref{nbr-4} and $\mu$ be the invariant measure for $(u(t,x))$ given in Theorem \ref{nbr-t1}
  (with $b=1$ in both instances and $\beta^\prime(y) = -a^2y$). According to \eqref{nbr-(*3)} and \eqref{nbr-(*4)},
  \beqn
   \frac{d \mu}{d \tilde\mu}(w) = \frac{1}{M} \exp\left(- \frac{a^2}{2}\int_0^1 w^2(x)\, dx\right) = \frac{1}{M} \exp\left(- \frac{a^2}{2} \Vert w \Vert_{V}^2\right),
\eeqn
where
\beqn
    M = \int_\bD \exp\left(- \frac{a^2}{2} \Vert w \Vert_{V}^2 \right) \tilde\mu(dw).
\eeqn
Because the exponential is strictly positive $\tilde\mu$-a.s., $\mu$ and $\tilde\mu$ are in fact mutually equivalent.



 \subsection{Reversible measures as bridge measures}
\label{new-5.4.3}
In Proposition \ref{ch5-added-1-s2-p1}, we saw that the invariant measure for the solution to the linear stochastic heat equation \eqref{ch1'.HD-ch5-markov} on $[0, 1]$, with vanishing Dirichlet boundary conditions, is the law of a Brownian bridge, that is, the law of a Brownian motion starting from $0$ and conditioned to take the value $0$ at ``time" $1$. In \cite{H-S-V-W-2005}, it was shown that this relation between invariant measures for various parabolic SPDEs and the law of ``bridge" (or conditioned) diffusions holds for many variants of the stochastic heat equation, including the linear SPDE \eqref{s5.3-new(*1)}, and non-linear variants of this SPDE such as \eqref{nbr-(*1aa)}. We present here some instances of these results. In Section \ref{new-5.5}, we will also establish the analog of Theorem \ref{ch5-added-1-s2-t1}, concerning the fact that the invariant measures are limit laws of the solutions to these SPDEs.

\begin{thm}
\label{split-th 5.3.16}
Let $(v(t,x))$ be the solution of the linear SPDE \eqref{s5.3-new(*1)} and let $\tilde\mu$ be the reversible probability measure on $(\bD, \cB_{\bD})$ defined in Proposition \ref{nbr-4}. Let $(X(x),\, x \in [0, 1])$ be the solution to the SDE
\beq
\label{s5.3-new(*2)}
        dX(x) = a X(x)\, dx + b\, dB(x),\quad  X(0) = 0,      
        \eeq
where $(B(x),\, x \in [0, 1])$ is a standard Brownian motion. Then $\tilde\mu$ is the conditional law on $\bD$ of the process $(X(x),\, x \in [0, 1])$ given that $X(1) = 0$.
\end{thm}
\begin{proof}
 The solution of the SDE \eqref{s5.3-new(*2)} 
     is  $X(x) = b \int_0^x e^{a(x-z)}\, dB(z)$ and its unconditional covariance is
     \begin{align}
     \label{s5.3-new(*3)}
   C_0(x, y) &= E[X(x) X(y)] = b^2 \int_0^{x \wedge y} e^{a(x-z)}e^{a(y-z)}\, dz\notag\\
   & = b^2 e^{a(x+y)}\, \frac{1 - e^{-2a(x \wedge y)}}{2a},\quad    {\text{if}}\  a \neq 0,    
   \end{align}
and
\beq
\label{s5.3-new(*4)}
   C_0(x, y) = E[X(x) X(y)] = b^2 (x \wedge y),\quad  {\text{if}}\ a = 0.   
   \eeq
The conditional covariance of $X$ given $X(1)$ is
\beq
\label{s5.3-new(*5)}
    C_1(x,y) := C_0(x, y) - C_0(x, 1) C_0(1,1)^{-1} C_0(1,y).  
    \eeq
Indeed, recalling that the conditional law of $X_1$ given $X_2 = x_2$, where $(X_1, X_2)$ is ${\rm{N}}((m_1, m_2), C)$ and $C=(c_{i,j})$ is the $2 \times 2$ matrix with $c_{i, j} = {\rm Cov}(X_i, X_j)$,  is
\beq
\label{s5.3-new(*4a)}
   {\rm{N}}\left(m_1 + \frac{{\rm {\rm {\rm Cov}}}(X_1, X_2)}{{\rm Cov}(X_2, X_2)} (x_2 - m_2),  {\rm Cov}(X_1, X_1) - \frac{{\rm Cov}(X_1, X_2)^2}{{\rm Cov}(X_2,X_2)}\right),     
   \eeq
we find that
\beqn
   E[X(x) \mid X(1) = x_1] = E[X(x)] + C_0(x, 1) C_0(1, 1)^{-1} (x_1 - E[X(1)]),
   \eeqn
and by applying \eqref{s5.3-new(*4a)} to $X(x)$, $X(y)$ and $X(x) + X(y)$, that ${\rm Cov}((X(x), X(y)) \mid X(1))$ is given by
\eqref{s5.3-new(*5)} and does not depend on $X(1)$.

   We now identify the covariance  $C_1 $ with  the covariance kernel $C$ of Theorem \ref{s5.3-new-t*1}. For  $a = 0$, this is done in \eqref{ch5-added-1-s2-3}. For $a \neq 0$, using \eqref{s5.3-new(*3)} and \eqref{s5.3-new(*5)}, we see that
\begin{align*}
    C_1(x,y) &= b^2 e^{a(x+y)}\, \frac{1 - e^{-2a(x \wedge y)}}{2a} \\
     &\qquad  - b^2 e^{a(x+1)}\, \frac{1 - e^{-2a x}}{2a}\, e^{-2a}\, \frac{2a}{1 - e^{-2a}}\, e^{a(y+1)}\, \frac{1 - e^{-2a y}}{2a}.
\end{align*}
After simplification, we find that
\begin{align*}
    C_1(x,y) &= \frac{b^2}{a} (e^{a(x+y- x \wedge y)} \sinh(a(x \wedge y)) - \frac{e^a}{\sinh(a)}\, \sinh(ax) \sinh(ay))\\
  &=    \begin{cases}  \sinh(ax) \frac{b^2}{a} \big(e^{ay} -  \frac{e^a}{\sinh(a)} \sinh(ay)\big) &     {\text{if}}\quad  x \leq y,\\
                 \frac{b^2}{a} \sinh(ay) \big(e^{ax} - \frac{e^a}{\sinh(a)} \sinh(ax)\big)&      {\text{if}}\quad x > y.
                 \end{cases}
                \end{align*}
Computing the Fourier series of $C_1(*, y)$ ($y$ fixed), we obtain the same expression as \eqref{s5.3-new(*2a)}, therefore $C_1(x,y) = C(x,y)$.
\end{proof}
\medskip

\noindent{\em Bridge measures in the case of gradient-type drifts}
\medskip

Let $a \in \R$, $b \neq 0$ and let $B = (B(x),\, x \in [0, 1])$ be a standard Brownian motion on a filtered probability space $(\Omega, \cF, P, (\cF_t))$. Let $f(x) = b^2 \tilde f'(x)$, where $\tilde f: \R \to \R$ is a $\cC^1$-function bounded from above and with a Lipschitz continuous derivative. Consider the nonlinear SDE
 \beq
 \label{nbr-(000)}
     dX(x) = \left(a X(x) + f(X(x))\right) dx + b\, dB(x),\quad  x\in\, ]0,1[,
     \eeq
with $X(0)=0$.  Our next objective is to produce an SPDE for which the reversible measure $\mu$ of Theorem \ref{nbr-t1} is the conditional law of $X$ given that $X(1)=0$. This will be achieved in Theorem \ref{nbr-t-14} below (see Remark \ref{s5.3-new-r(*3)} for some motivation). In the next two lemmas, we give the preliminaries that are needed to reach this goal.

Let $\cC_0$ be the space of continuous functions $w: [0, 1] \to \R$ such that $w(0) = 0$, endowed with the uniform norm. Since $\bD$ is a closed subset of $\cC_0$, any probability measure $\nu$ on $\bD$ can be considered to be a probability measure on $\cC_0$ such that $\nu(\cC_0 \setminus \bD) = 0$.



We present now a result similar to Theorem \ref{split-th 5.3.16} for the solution $(X(x))$
to
\beq
\label{nbr-(*29a)}
    dX(x) = \left(a X(x) + b^2 \tilde f'(X(x))\right) dx + b\, dB(x),\quad     X(0) = 0,     
\eeq
$x\in\, ]0,1]$,
which is  \eqref{nochnochnoch(*1)} with $f(x): = b^2 \tilde f'(x)$ and $C=b$ (notice that, in general, $(X(x))$ is not a Gaussian process).
It turns out that the drift function $\tilde f'$ that appears in the SDE \eqref{nbr-(*29a)} will in general not be the same as the drift function that will appear in the associated SPDE, which will be
\beq
\label{nbr-(*29)}
\cL u(t,x) = \Phi'(u(t,x)) + \sqrt{2}\, \dot W(t,x), \quad t>0,\ x\in\, ]0,1[,   
  \eeq
with $\Phi$ given in \eqref{nochnochnoch(*4)}, with $n = 1$, $A = a$, $C = b$ and the same $\tilde f$ as in \eqref{nbr-(*29a)},
 with vanishing Dirichlet boundary conditions and some initial condition $u_0 \in \mathbb{D}$.

\begin{thm}
\label{nbr-t-14}
Let $X=(X(x))$ be defined in \eqref{nbr-(*29a)} and let $\nu_X$ be the conditional law of $X$ given that $X(1) = 0$. Assume that $\tilde f\in \cC^2(\re)$, $\tilde f$ is bounded from above and ${\tilde f}^\prime$ is Lipschitz continuous. Then $\nu_X$ is equal to the reversible probability measure $\mu$ for $(u(t, *))$ defined in Theorem \ref{nbr-t1}.
\end{thm}
\begin{proof}
Let $v$ be the solution to the linear SPDE \eqref{s5.3-new(*1)} and let $Z$ be defined in \eqref{nochnochnoch(*2)} (with $n=1$, $A=a$ and $C=b$). Let $\mu_X$ and $\mu_Z$ be as in Lemma \ref{nbr-l-12} (with $f(x) = b^2 \tilde f^{\prime}(x)$), and let $\nu_Z$ be the conditional law of $Z$ given $Z(1) = 0$. By Theorem \ref{split-th 5.3.16}, $\nu_Z$ (denoted there $\tilde\mu$) is a reversible measure for $v$.

    Let $B(w) = \int_0^1 \Phi(w(x))\, dx$. By Theorem \ref{nbr-t1}, the measure $\nu$ on $\mathbb{D}$ defined by
\beq
    \label{nbr-(*30a)}
      \nu(dw) = \frac{1}{M}\,  e^{B(w)}\,  \nu_Z(dw),      
\eeq
where $M = \int_\bD \exp(B(w))\, \nu_Z(dw)$ is a normalisation constant (see \eqref{nbr-(*4)}), is a reversible (and invariant) probability measure for $(u(t,*))$. We now check that $\nu = \nu_X$.

     Consider the function $\psi: \cC_0 \times \R \to \R$ defined by
     \beq
     \label{nbr-(*31)}
    \psi(w, y) = \exp\left(\tilde f(y) + \int_0^1 \Phi(w(x))\, dx\right).      
    \eeq
By Lemma \ref{nbr-l-13}, the density of $\mu_X$ with respect to $\mu_Z$ is
\beq
\label{nbr-(*30)}
     \varphi(w) = c \exp\left(\tilde f(w(1)) + \int_0^1 \Phi(w(x))\, dx\right) =: \psi(w, w(1)),    
     \eeq
where $c= \exp(-\tilde f(0))$.

    Let $P_0$ be the law on $\cC_0 \times \R$ of $(X, X(1))$ under $\mu_X$, $Q_0$ be the law on $\cC_0 \times \R$ of $(Z, Z(1))$ under $\mu_Z$. For $A_1 \in \B_{\cC_0}$ and $A_2 \in \cB_\R$, let $C = \{w \in \cC_0:\ w(1) \in A_2\}$. Then by \eqref{nbr-(*30)},
    \begin{align*}
   P_0(A_1 \times A_2) &= \mu_X(A_1 \cap C)   = \int_{\cC_0} \mu_Z(dw)\, 1_{A_1 \cap C}(w)\, \varphi(w)\\
   & = \int_{\cC_0} \mu_Z(dw)\, 1_{A_1 \times A_2}(w, w(1))\, \psi(w, w(1))\\
   & = \int_{\cC_0 \times \R} Q_0(dw, dy)\, 1_{A_1  \times A_2}(w, y)\, \psi(w, y),
   \end{align*}
therefore,
\beq
\label{nbr-(*32)}
   \frac{d P_0}{d Q_0}(w, y) = \psi(w, y).      
   \eeq
  The conditional law under $Q_0$ of $Z$ given $Z(1) = 0$ is $\nu_Z$, denoted $\tilde\mu$ in Theorem \ref{split-th 5.3.16}. For $y \in \R$, consider the linear function $z_y: [0,1] \to \R$ defined by 
 \beqn
    z_y(x) = {\rm Cov}(Z(x), Z(1))\, {\rm Cov}^{-1}(Z(1), Z(1))\, y. 
\eeqn
Because $Z$ is a Gaussian process, the formulas \eqref{s5.3-new(*5)} and \eqref{s5.3-new(*4a)} imply that the conditional law under $Q_0$ of $Z$ given $Z(1) = y$ is $\mu_y$, where, for $A \in \B_{\cC_0}$,
  \beq
  \label{afegit-(*1)}
         \mu_y(A) = \mu_0(-z_y + A).,     
         \eeq
and $\mu_0$ is $\mu_{y:= 0}$, that is, $\mu_y$ with $y = 0$. We now apply Lemma \ref{nbr-l-11} with 
\beqn
   \Omega = \cC_0 \times \R, \quad \cF = \cB_{\cC_0} \times \cB_\R, \quad P = P_0,\quad Q = Q_0, \quad \cs = \cC_0, \quad U = \R,
\eeqn
\beqn
     X(w, y) := w, \qquad Y(w, y) := y, \qquad \varphi(x, y) := \psi(w, y),
\eeqn
where $\psi$ is defined in \eqref{nbr-(*31)}. By \eqref{nbr-(*32)}, $P_0$ is absolutely continuous with respect to $Q_0$. By \eqref{afegit-(*1)}, the conditional law $Q_y$ under $Q_0$ of $X$ given $Y = y$ is $\mu_y$, and $(y, A) \mapsto \mu_y(A) = \tilde\mu(-z_y + A)$  is a regular conditional probability distribution.
   By Lemma \ref{nbr-l-11}, the conditional law under $P_0 $ of  $X $ given  $Y = y $ is
   \beqn
     P_y(dw) = \frac{1}{c(y)}\, \psi(w, y)\, \mu_y(dw),
     \eeqn
where  $c(y) = \int_{\cC_0}  \psi(w, y)\, \mu_y(dw) $. Setting  $y = 0 $, and recalling that  $\nu_X $ is, by definition, the conditional law of $X$ (under $P$) given that  $X(1) = 0 $, and that  $\tilde\mu = \nu_Z $,
we have $\nu_X = P_{y := 0}$ (that is, $P_y$ for $y=0$), therefore,
for  $A \in \cB_{\cC_0} $,
\beq
\label{nbr-(*33)}
    \nu_X(A) = \int_A \psi_1(w)\, \nu_Z(dw),   
    \eeq
where
\beqn
\psi_1(w) = \frac{1}{c(0)}\, \psi(w, 0) = \frac{1}{c(0)}\, \exp\left(\int_0^1 \Phi(w(x)) dx\right) = \frac{1}{c(0)}\, \exp(B(w)),
\eeqn
and
 \begin{align*}
 c(0) &= \int_{\cC_0} \psi(w, 0)\, \nu_Z(dw) = \int_{\cC_0} \exp\left(\int_0^1 \Phi(w(x))\, dx\right) \nu_Z(dw)\\
 & = \int_\bD  \exp(B(w))\,\nu_Z(dw) = M.
 \end{align*}
 From \eqref{nbr-(*30a)} and \eqref{nbr-(*33)}, we conclude that $\nu_X=\nu$.
    \end{proof}

\begin{remark}
\label{s5.3-new-r(*3)} Hairer et al.~\cite{H-S-V-W-2007} work with systems of SDE's and SPDEs, and consider other types of conditioning, as well as slight variations on $\cL$ above, both for linear and nonlinear SPDEs with additive noise and gradient-type drift. These results are used in many applications to simulate the distributions of conditioned SDEs.
\end{remark}


\section{Regularity and convergence in law}
\label{new-5.5}

In this section, we delve further into the study of the nonlinear stochastic heat equation with additive noise and gradient-type drift considered  in \eqref{nbr-(*1aa)}: we prove that the reversible measure is the unique invariant measure and that the law of the solution $u(t, *)$ converges to this invariant measure as $t \to \infty$  (see Theorem \ref{nbr-t27} parts (d) and (c), respectively).
This extends the result of Theorem \ref{s5.3-new-t*1} for linear stochastic heat equations. We use a method that goes back to Doob \cite{doob-1948} and that relies on regularity properties of the Markov transition function (see Theorem \ref{nbr-20b}).

We note that the existence of a reversible (and invariant) measure for the solution of
\eqref{nbr-(*1aa)} has already been established in Theorem \ref{nbr-t1}. 

We go back to the setting of Theorem \ref{nbr-t1}, that we now recall.
Let $\beta: \R \to \R$ be a function of class $\cC^1$ that is bounded from above and with Lipschitz continuous derivative $\beta'$. For $w \in \bD$, as in \eqref{nbr-(*1)}, we define
\beqn
    B(w) = \int_0^1 \beta(w(x))\, dx,     
    \eeqn
so that $B: \bD \to \R$.

Let $\cL$ be the partial differential operator given in \eqref{operator-calL}. Consider the SPDE \eqref{nbr-(*1aa)}, that is,
\beqn
       \cL u(t,x) = \beta'(u(t,x)) + \sqrt{2}\, \dot W(t,x),\quad   t > 0,\    x \in\, ]0, 1[,   
       \eeqn
subject to vanishing Dirichlet boundary conditions, and some suitable initial condition $u_0 \in \bD$. 

Notice that if $\beta(z) = -(a^2/2) z^2$ (and $a \neq 0$), that is, $\beta'(z) = - a^2\, z$, then the convergence to the invariant measure is a consequence of Theorem \ref{s5.3-new-t*1}. Indeed, in this case, \eqref{nbr-(*1aa)} is a linear SPDE (and a special case of \eqref{s5.3-new(*1)}) and Theorem \ref{s5.3-new-t*1} applies.

Recall that from Theorem \ref{nbr-2aa} that $(u(t, *))$ is a Markov process with transition function $(P_t)$ defined in \eqref{ch5-added-1-2-2}. 

\bigskip

\noindent{\em Some notions from the theory of Markov processes}
\medskip

In order to exploit the properties of $(u(t, *))$ as a Markov process, we now recall some notions from the theory of such processes. The notation and terminology are as in Section \ref{ch5-added-1-s1}, with the simplifications introduced just after Definition \ref{nbr-0aa}.

\begin{def1}
\label{nbr-20a}
Let $(P_t)$ be the transition function associated with a Markov process $(X_t)$ on a complete separable metric space $\cs$.
\begin{description}
\item{(a)} For $t_0 \in \R_+^*$,  $(P_t)$ is {\em $t_0$-regular}\index{regular!$t_0$-}\index{t@$t_0$-regular} if for all $t > t_0$,  the transition probabilities $P_t(g, \cdot)$, $g \in \cs$, are all mutually equivalent.
\item{(b)} If $\mu$ is an invariant measure for $(P_t)$, then $\mu$ is {\em strongly mixing}\index{strongly!mixing}\index{mixing!strongly} if for all bounded Borel functions $f,\, g \in\cB_{\cs}$,
\beqn
\lim_{s \to \infty} E^\mu[f(X_t)\, g(X_{t + s})] =  E^\mu[f(X_t)]\, E^\mu[g(X_t)].
\eeqn
\end{description}
\end{def1}

The next theorem is quoted from \cite[Theorem 4.2.1]{dz-1996} (originating from \cite{doob-1948}).
\begin{thm}
\label{nbr-20b}
 Let $(P_t)$ be a stochastically continuous transition function on $\cs$ and $\mu$ an invariant probability measure for $(P_t)$. If $(P_t)$ is $t_0$-regular for some $t_0 > 0$, then:
 \begin{description}
 \item{(a)}\ $\mu$ is strongly mixing and for all $g \in \cs$ and $A \in \cB_{\cs}$,
 \beqn
          \lim_{t \to \infty} P_t(g, A) = \mu(A);
          \eeqn
\item{(b)}\  $\mu $ is the unique invariant measure for  $(P_t)$
and $\mu$ is equivalent to all the probability measures $P_t(g, \cdot)$, $g \in \cs$, $t > t_0$.
\end{description}
\end{thm}

We refer to \cite{dz-1996} for the proof of this theorem. 

\subsection{Regularity: the linear case}

We will now establish $t_0$-regularity for the linear version \eqref{s5.3-new(*1)} of \eqref{nbr-(*1aa)} (in which $\beta' \equiv 0$), then we will extend this to more general functions $\beta$. Because the transition function $(P_t)$ defined in \eqref{ch5-added-1-2-2} is stochastically continuous by Remark \ref{ch5-added-c1}, we will then use Theorem \ref{nbr-20b} to obtain in particular the desired convergence result (which is stronger than weak convergence).

Let us consider the SPDE \eqref{s5.3-new(*1)} with initial condition $v_0\in\bD$, that is,
\beqn
    \cL v(t,x) = \sqrt{2}\, \dot W(t,x),\quad   t > 0,\   x \in\, ]0, 1[,     
    \eeqn
 and $v(0, *) = v_0(*)$, and denote by $v_{v_0} = (v_{v_0}(t,x),\, (t, x) \in \R_+ \times [0,1])$ its random field solution. The transition function associated to the Markov process $v_{v_0}$ is defined by
 \beq
 \label{nbr-semigroup}
 P_t f(v_0) = E[f(v_{v_0}(t, *)],\quad t>0,
 \eeq
 for $\cB_\bD$-measurable and bounded functions $f: \bD \to \R$.

\begin{thm}
\label{nbr-25}
The transition function $(P_t)$ defined in \eqref{nbr-semigroup} is $t_0$-regular for all $t_0 > 0$.
\end{thm}
\begin{proof}
Given two arbitrary functions $v_0, \tilde v_0\in \bD$, we will show that for all $t>0$, the laws of $v(t,*):=v_{v_0}(t, *)$ and  $\tilde v(t,*):=v_{\tilde v_0}(t, *)$ are equivalent.

Let $V = L^2([0, 1])$ endowed with the CONS defined in the line following \eqref{ch5-added-1-4} with $L=1$, which, for simplicity, we write $e_n(x)$ instead of $e_{n,1}(x)$ there ($n\ge 1$).
Define the transition operator $G_{a,b;t}: V \to V$ by
\beq
\label{noch(*1)}
G_{a,b;t}(w)(x) = \int_0^1 dy\, G_{a,b}(t; x, y) w(y),\quad t>0,
\eeq
where $G_{a,b}$ is the Green's function corresponding to $\cL$ given in \eqref{s5.3-new(*1)green}. By \eqref{s5.3-new(*1)green} and \eqref{ch5-added-1-4},
\beqn
   G_{a,b}(t; x, y)  = \sum_{n=1}^\infty e^{- (a_n/b^2)  t } \, e_n(x) e_n(y),
\eeqn
where $a_n = \pi^2 n^2+a^2$.

Fix $t \in \R_+$ and let
\beqn
   Z(x) = \sqrt{2}  \int_0^t \int_0^1 G_{a,b}(t - s; x, z)\, W(ds, dz).
       \eeqn
 Then, for any $w\in\bD$, we have
 \beqn
    v_w(t,x) = G_{a,b;t}(w)(x) + Z(x).
 \eeqn

Set  $\lambda_n = \frac{b^2}{a_n}\left(1 - e^{- (2a_n/b^2)\,  t}\right)$.
We can write
\beq
\label{nbr-(*47a)}
     Z(x) = \sum_{n=1}^\infty A^{(n)}\, e_n(x)     
     \eeq
where
\beqn
     A^{(n)}= \sqrt{2} \int_0^t  e^{- (a_n/b^2)  (t - s)}\, dW_s(e_n).
     \eeqn
In particular, $A^{(n)}$ is ${\rm{N}}(0, \lambda_n)$ and the $A^{(n)}$ are independent. By \eqref{nbr-(*47a)}, $Z(*)$ takes values in $V$ (and even in $\bD$), and it defines a Gaussian random vector on $V$ with mean $0$, covariance kernel
  \begin{align*}
    C(x, y) &= E[Z(x) Z(y)] = 2 \int_0^t G_{a,b}(2(t - s); x, y)\, ds\\
    & = 2 \int_0^t G_{a,b}(2s, x, y)\, ds  = 2 \sum_{n=1}^\infty \int_0^t ds \, \exp\left(- \frac{a_n}{b^2} \, 2s \right) e_n(x) e_n(y) \\
    &= 2  \sum_{n=1}^\infty \lambda_n\, e_n(x)\, e_n(y) ,
    \end{align*}
and covariance operator
\beq
\label{nbr-(*50)}
     C(f)(x) = \sum_{n=1}^\infty \lambda_n\, \langle f, e_n \rangle_V \,e_n(x) , \qquad  f \in V.   
 \eeq

For $w\in\bD$ and $n\geq 1$, let
\beqn
  w^{(n)} = e^{- \frac{a_n}{b^2} t}\, \langle w, e_n \rangle_V.
    \eeqn
Then
\beqn
  G_{a,b:t}(w)(x) =  \sum_{n=1}^\infty w^{(n)}\,  e_n(x),
      \eeqn
and, along with \eqref{nbr-(*47a)}, this yields
\beqn
     v_w(t, x) = \sum_{n=1}^\infty \left(w^{(n)} + A^{(n)}\right)\, e_n(x), \qquad w \in \bD.
\eeqn

   For $n \geq 1$, set 
\beq\label{nbr-(*47b)}
     v_0^{(n)} = e^{- (a_n/b^2) t}\,  \langle v_0, e_n \rangle_V,\qquad \tilde v_0^{(n)} = e^{- (a_n/b^2) t}\,  \langle \tilde v_0, e_n \rangle_V
\eeq
and let $\mu_n$ (respectively $\nu_n$) be the law on $\R$ of 
\beqn
     X^{(n)} = v_0^{(n)} + A^{(n)}\quad \text{(respectively, } Y^{(n)} = \tilde v_0^{(n)} + A^{(n)}). 
\eeqn
The density of $\mu_n$ with respect to $\nu_n$ is the ratio of their densities with respect to Lebesgue measure, which can be written   \beq
   \label{nbr-(*50a)}
    \frac{d\mu_n}{d\nu_n}(x) = \exp\left(\frac{v_0^{(n)}- \tilde v_0^{(n)}}{\lambda_n} \left(x - \tilde v_0^{(n)}\right) - \half \frac{\left(v_0^{(n)}- \tilde v_0^{(n)}\right)^2}{\lambda_n}\right),\quad      x \in \R.       
    \eeq
Therefore,
\begin{align}
\label{nbr-(*50b)}
    E\left[h\left(X^{(n)}\right)\right] &= \int_\R h(x)\, \mu_n(dx) = \int_\R h(x) \, \frac{d\mu_n}{d\nu_n}(x)\, \nu_n(dx)\notag\\
    & = E\left[h\left(Y^{(n)}\right)\, \frac{d\mu_n}{d\nu_n}\left(Y^{(n)}\right)\right].     
    \end{align}
Let $\mu$ (respectively $\nu$) be the law on $\bD $ of  $v(t, *) $ (respectively $\tilde v(t, *)$), and let  $\varphi: V \to \R $ be a bounded and continuous function. For  $N \in \N $, define $\Phi_N: \R^N \to \R$ by
\beqn
 \Phi_N(x_1,\dots, x_N) = \varphi\left(\sum_{n=1}^N x_n e_n\right).
 \eeqn
  By independence of the $X^{(n)}$,
  \begin{align*}
    E\left[\varphi\left(\sum_{n=1}^N X^{(n)} e_n\right)\right] &= E\left[\Phi_N\left(X^{(1)},\dots, X^{(N)}\right)\right] \\
    & = \int_{\R^N} \Phi_N(x_1,\dots, x_N)\, \mu_1(dx_1)\, \cdots\, \mu_N(dx_N).
    \end{align*}
By \eqref{nbr-(*50a)} and \eqref{nbr-(*50b)}, this is equal to
\begin{align*}
      & \int_{\R^N} \Phi_N(x_1,\dots, x_N)\\
      &\qquad\times \exp\left(\sum_{n=1}^N \left(\frac{v_0^{(n)}- \tilde v_0^{(n)}}{\lambda_n} \left(x_n - \tilde v_0^{(n)}\right) - \half \frac{\left(v_0^{(n)}- \tilde v_0^{(n)}\right)^2}{\lambda_n}\right)\right)\\
      &\qquad\quad\times   \nu_1(dx_1)\, \cdots\, \nu_N(dx_N),
       \end{align*}
so
\begin{align}
\label{nbr-(*50c)}
     &E\left[\varphi\left(\sum_{n=1}^N X^{(n)} e_n\right)\right]\notag\\
     &\quad = E\left[\varphi\left(\sum_{n=1}^N Y^{(n)} e_n\right)\right. \notag\\
      &\left.\qquad \quad\times\exp\left(\sum_{n=1}^N \left(\frac{v_0^{(n)}- \tilde v_0^{(n)}}{\lambda_n} \left(Y^{(n)} - \tilde v_0^{(n)}\right) - \half \frac{\left(v_0^{(n)}- \tilde v_0^{(n)}\right)^2}{\lambda_n}\right)\right)\right].
     \end{align}
From their definition in \eqref{nbr-(*47b)}, observe that
\beq
\label{nbr-(*50d)}
    \sum_{n=1}^\infty \frac{\left(v_0^{(n)}- \tilde v_0^{(n)}\right)^2}{\lambda_n} < \infty,     
    \eeq
and
\beqn
    \sum_{n=1}^\infty \frac{v_0^{(n)}- \tilde v_0^{(n)}}{\lambda_n}\left (Y^{(n)} - \tilde v_0^{(n)}\right) = \sum_{n=1}^\infty \frac{v_0^{(n)}- \tilde v_0^{(n)}}{\lambda_n} \, A^{(n)}
    \eeqn
converges in $V$ (and in $\bD$, as in Lemma \ref{ch1'-lsi}) because
\beqn
    \sum_{n=1}^\infty \frac{\left(v_0^{(n)}- \tilde v_0^{(n)}\right)^2}{\lambda_n^2}\, E\left[\left(A^{(n)}\right)^2\right] = \sum_{n=1}^\infty \frac{\left(v_0^{(n)}- \tilde v_0^{(n)}\right)^2}{\lambda_n} < \infty.
    \eeqn
Since the function $\varphi$ is continuous, we let $N \to \infty$ in \eqref{nbr-(*50c)} to see that the random variable on the right-hand side converges a.s. to
\beqn
     \varphi(\tilde v(t, *))\, \exp\left(\sum_{n=1}^\infty \left(\frac{v_0^{(n)}- \tilde v_0^{(n)}}{\lambda_n} \left(Y^{(n)} - \tilde v_0^{(n)}\right) - \half \frac{\left(v_0^{(n)}- \tilde v_0^{(n)}\right)^2}{\lambda_n}\right)\right).
     \eeqn
In order to conclude that
\begin{align}
\label{nbr-(*50e)}
     &E[\varphi(v(t, *)] \notag\\
     &\ = E\left[\varphi(\tilde v(t, *)) \exp\left(\sum_{n=1}^\infty\left (\frac{v_0^{(n)}- \tilde v_0^{(n)}}{\lambda_n} \left(Y^{(n)} - \tilde v_0^{(n)}\right) - \half \frac{\left(v_0^{(n)}- \tilde v_0^{(n)}\right)^2}{\lambda_n}\right)\right)\right],      
     \end{align}
we use the fact that $\varphi$ is bounded and the elementary fact that $\psi_N \to \psi$ in $L^2(\Omega)$ implies $\psi_N^2 \to \psi^2$ in $L^1(\Omega)$, and apply this to
\begin{align*}
        &\psi_N\left(\sum_{n=1}^\infty Y^{(n)} e_n\right)\\
        &\qquad = \exp\left(\half\sum_{n=1}^N \left(\frac{v_0^{(n)}- \tilde v_0^{(n)}}{\lambda_n} \left(Y^{(n)} - \tilde v_0^{(n)}\right) - \half \frac{\left(v_0^{(n)}- \tilde v_0^{(n)} \right)^2}{\lambda_n}\right)\right)
\end{align*}
and
\begin{align*}
        &\psi\left(\sum_{n=1}^\infty Y^{(n)} e_n\right)\\
        &\qquad = \exp\left(\half\sum_{n=1}^\infty \left(\frac{v_0^{(n)}- \tilde v_0^{(n)}}{\lambda_n} \left(Y^{(n)} - \tilde v_0^{(n)}\right) - \half \frac{\left(v_0^{(n)}- \tilde v_0^{(n)}\right)^2}{\lambda_n}\right)\right).
\end{align*}

In order to check the $L^2(\Omega)$-convergence, we first observe that
\begin{align*}
    &E\left[\left(\psi_{N+m}\left(\sum_{n=1}^\infty Y^{(n)} e_n\right) - \psi_N\left(\sum_{n=1}^\infty Y^{(n)} e_n\right)\right)^2\right]\\
   &\quad = E\left[ \psi_N^2(\tilde v(t, *))
   \left(\exp\left(\half\sum_{n=N+1}^{N+m} \left[\frac{v_0^{(n)}- \tilde v_0^{(n)}}{\lambda_n} \left(Y^{(n)} - \tilde v_0^{(n)}\right)
    \right. \right. \right.\right.\\
    &\left.\left.\left.\left.\qquad\qquad \qquad\qquad\qquad
     - \half \frac{\left(v_0^{(n)}- \tilde v_0^{(n)}\right)^2}{\lambda_n}\right]\right)- 1\right)^2\right].
    \end{align*}
By independence and the fact that $\psi_N^2(\tilde v(t, *))$ is a product of density functions, this is equal to the expectation of the second factor, which is equal to
\begin{align*}
      &2 - 2 \prod_{n = N+1}^{N+m} E\left[\exp\left(\half \frac{v_0^{(n)}- \tilde v_0^{(n)}}{\lambda_n} \left(Y^{(n)} - \tilde v_0^{(n)} \right) - \frac{1}{4}
      \frac{\left(v_0^{(n)}- \tilde v_0^{(n)}\right)^2}{\lambda_n}\right)\right] \\
      &\qquad = 2 - 2 \exp\left(-\sum_{n=N+1}^{N+m} \frac{\left(v_0^{(n)}- \tilde v_0^{(n)}\right)^2}{8 \lambda_n}\right)
\end{align*}
(where we have used the moment generating function of a ${\rm N}(0, 1)$ random variable). By \eqref{nbr-(*50d)}, the sum in the exponential converges to $0$, so this last right-hand side tends to $0$ as $N,\, m \to \infty$. This implies that
  $\psi_N(\sum_{n=1}^\infty Y^{(n)} e_n)$ converges in $L^2(\Omega)$ to a random variable $\Psi \in L^2(\Omega)$. By the elementary fact recalled above,  $\psi_N^2(\sum_{n=1}^\infty Y^{(n)} e_n)$ converges in $L^1(\Omega)$ to $\Psi^2$, that is,
\begin{align*}
   \Psi^2 &= \lim_{N \to \infty} \psi_N^2\left(\sum_{n=1}^\infty Y^{(n)} e_n\right)\\
   & = \exp\left(\sum_{n=1}^\infty \left(\frac{v_0^{(n)}- \tilde v_0^{(n)}}{\lambda_n} \left(Y^{(n)} - \tilde v_0^{(n)}\right) - \half \frac{\left(v_0^{(n)}- \tilde v_0^{(n)}\right)^2}{\lambda_n}\right)\right)
   \end{align*}
in $L^1(\Omega)$, which proves \eqref{nbr-(*50e)}.

   From \eqref{nbr-(*50e)}, we conclude that the law of $v(t, *)$ is absolutely continuous with respect to the law of $\tilde v(t, *)$.  Since we can exchange the roles of $v$ and $\tilde v$, the laws of $v(t, *)$ and $\tilde v(t, *)$ are in fact equivalent.
      \end{proof}

\begin{remark}
\label{nbr-r-25a}
The equivalence of the laws of $v(t, *)$ and $\tilde v(t, *)$ is a special case of the Feldman-Hajek theorem \cite[Theorems 2.23 and 2.25]{dz} for Gaussian measures on Hilbert spaces. The only property of the sequence $(v_0^{(n)} - \tilde v_0^{(n)},\, n \geq 1)$ that is used is  \eqref{nbr-(*50d)}. The proof of the $L^1(\Omega)$-convergence is a special case of a classical result on the density of a countably infinite product of measures (see, e.g. \cite[Proposition 2.21]{dz}).
\end{remark}

 \subsection{Regularity and convergence in law: the case of gradient-type drift}
 
 We begin with a property concerning equivalence of laws of solutions to a linear SPDE. The proof uses a result from Subsection \ref{ch2'-s3.4.2} in Chapter \ref{ch2'}. 

\begin{prop}
\label{nbr-p22}
 Let $u$ be the solution of \eqref{nbr-(*Mb-bis)}
 with $\sigma \equiv \sqrt{2}$, with vanishing Dirichlet boundary conditions and initial condition $u_0$. Consider a jointly measurable, adapted and bounded process $h = (h(t,x))$. Let $v(t,x))$ be the random field solution of
 \beq
 \label{nbr-(*43)}
        \cL v(t,x) = h(t,x) + \sqrt{2}\, \dot W(t,x),\quad   t > 0, \  x \in\, ]0, 1[,   
        \eeq
with the same boundary conditions and initial conditions as for \eqref{nbr-(*Mb-bis)}. Then the laws of $(u(t, x),\, (t, x) \in [0, T] \times [0, 1])$ and $(v(t,x),\, (t, x) \in [0, T] \times [0, 1])$ are mutually equivalent.
\end{prop}
\begin{proof}
Let $\tilde v$ denote a jointly measurable and adapted version of the random field defined in \eqref{ch2'-s3.4.1.4}, with $\Gamma(t, x; s, y) := G_{a, b}(t - s, x, y)$. Because $h$ is bounded and $\beta$ is Lipschitz, Proposition \ref{ch2'-s3.4.2-t1} implies that the laws of $u$ (respectively $v$) and $\tilde v$ are mutually equivalent. This completes the proof.
\end{proof}

\noindent{\em Convergence to the reversible distribution}
\medskip

With these preliminaries, we now state the main result concerning convergence of the law of the solution of \eqref{nbr-(*1aa)} to the reversible distribution.
\begin{thm}
\label{nbr-t27}
 Let $u$ be the solution of \eqref{nbr-(*1aa)} with vanishing Dirichlet boundary conditions and initial condition $u_0 \in \bD$. Let  $(P_t)$ be the transition function of $u$ defined in \eqref{ch5-added-1-2-2} and let $\mu$ be the reversible probability measure given in Theorem \ref{nbr-t1}.
 Then:
 \begin{description}
    \item{(a)}\ $(P_t)$ is $t_0$-regular, for all $t_0 > 0$;
    \item{(b)}\ $\mu$ is strongly mixing and for all $A \in \cB_\bD$,
    \beq
    \label{nbr-(*52)}
         \lim_{t \to \infty} P_t(u_0, A) = \mu(A).   
         \eeq
    \item{(c)}\ $\mu$ is the unique invariant probability measure for $u$, and is equivalent to all the probability measures $P_t(u_0, \cdot)$, $u_0 \in \bD$, $t > 0$.
    \end{description}
    \end{thm}
    
\begin{proof}
In order to prove that $(P_t)$ is $t_0$-regular, we fix $t > 0$, $u_0, \bar u_0 \in \bD$ and we will show that $P_t(u_0, \cdot)$ and $P_t(\bar u_0, \cdot)$ are mutually equivalent.

  For this, let $u_{u_0}$ (respectively $u_{\bar u_0}$) be the random field solution
 of  \eqref{nbr-(*1aa)} with initial condition $u_0$ (respectively $\bar u_0$). Let $v_{u_0}$ (respectively $v_{\bar u_0}$) solve the SPDE
 \beqn
 \cL v(t,x) = \sqrt 2\, \dot W(t,x),\quad t>0,\ x\in\, ]0,1[,
 \eeqn
  with initial condition $u_0$ (respectively $\bar u_0$). By Proposition \ref{nbr-p22} with $h\equiv 0$, the laws of $u_{u_0}$ and $v_{u_0}$ are mutually equivalent, and the laws of $u_{\bar u_0}$ and $v_{\bar u_0}$ are also mutually equivalent. By Theorem \ref{nbr-25}, the laws of $v_{u_0}$ and $v_{\bar u_0}$ are mutually equivalent. Therefore, the same is true for $u_{u_0}$ and $u_{\bar u_0}$.

 Finally, since $(P_t)$ is stochastically continuous by Remark \ref{ch5-added-c1}, Theorem \ref{nbr-20b} implies that $\mu$ is strongly mixing, \eqref{nbr-(*52)} holds, $\mu$ is the unique invariant measure, and $\mu$ is equivalent to all the $P_t(u_0, \cdot)$, $t > 0$, $u_0 \in \bD$.
 \end{proof}

 \begin{remark}
 Because the invariant probability measure $\mu$ associated with $(P_t)$ in Theorem \ref{nbr-t27} is a limit measure and is the unique invariant measure, the semigroup $(P_t)$ is in fact {\em ergodic} (see \cite[Definition 1.9]{liggett}).
 \end{remark}

 \section{The irreducibility property (additive noise)}
 \label{new-5.5.2}

 This section addresses the property of
 irreducibility of the transition function $(P_t)$ associated to the stochastic heat equation
  \eqref{nbr-(*Mb)} (with $L=1$) defined in \eqref{ch5-added-1-2-2}, with additive noise ($\sigma$ constant) and $\beta$ not necessarily a gradient: see Theorem \ref{nbr-p21} below. This result applies, of course, to the SPDE \eqref{nbr-(*1aa)}.

We begin by recalling the notion of irreducibility and some preliminary results.
\begin{def1}
\label{nbr-20a-irred}
Let $(P_t(g, A))$ be the  transition function associated with a Markov process $(X_t)$ on a complete separable metric space $\cs$. We say that
$(P_t)$ is {\em irreducible at time $t_0 > 0$}\index{irreducible} if $P_{t_0}(g, A) > 0$ for all $g \in \cs$ and open sets $A \subset \cs$.
\end{def1}

Informally, in the context of the solution $(u(t, x))$ to \eqref{nbr-(*Mb)}, $(P_t)$ in \eqref{ch5-added-1-2-2} will be irreducible at time $t_0 > 0$ if, for any $u_0, w \in \bD$, starting from $u_0$, the solution will be ``close to $w$ at time $t_0$" with positive probability. In order to establish this, we will fix $s_0 \in\, ]0, t_0[$ and examine the position $u(s_0, *) \in \bD$ of the solution at time $s_0$. Then we construct a deterministic PDE (see \eqref{nbr-(*53)} below) that, in the absence of noise, would bring the solution ``close to $w$ at time $t_0$." Since the noise is additive, if its contribution to the SPDE is sufficiently small, then $u(t_0, *)$ will also be ``close to $w$."

We now introduce the deterministic PDE that we will need. Consider $\hat w \in \bD$ and $\tilde w \in \bD \cap \cC^2([0, 1])$. Let $0 \leq s_0 < t_0$ and let $\cL$ be as defined in \eqref{operator-calL}. We want to find a source term $h(t,x)$ with the property that the solution to the PDE
\beq
\label{nbr-(*53)}
    \cL U(t,x) = h(t,x),\quad t\in\, ]s_0,t_0[,\ x\in\,]0,1[,       
    \eeq
subject to vanishing Dirichlet boundary conditions and starting condition $U(s_0, *) = \hat w$, satisfies $U(t_0, *) = \tilde w$. That is, $(U(t, *),\, t \in [s_0, t_0])$ is a ``bridge" from $\hat w$ to $\tilde w$ that solves \eqref{nbr-(*53)}. The next lemma gives us such a function $h$. Recall that $G_{a,b;t}(w)(z)$ is defined in \eqref{noch(*1)}.

\begin{lemma}
\label{nbr-l23}
Let $0 \leq s_0 < t_0$. Define
\beqn
   h(t, x) = 1_{[s_0, t_0]}(t) f(t,\hat w,\tilde w)(x),
\eeqn
where $f: [s_0, t_0] \times \bD \times \cC^2([0, 1]) \to \bD$ is given by
\beqn
      f(s,\hat w,\tilde w)(z) = \frac{1}{t_0 - s_0}\, G_{a,b;s - s_0}(\tilde w - \hat w)(z) - \mA  \tilde w(z),\quad       s \in [s_0, t_0],\   z \in [0, 1],
      \eeqn
      with $\mA $ defined in \eqref{def-tilde-L}. 
Let $U(\hat w,\tilde w,\cdot, *)$ be the solution of \eqref{nbr-(*53)} on $[s_0, t_0]$, starting at time $s_0$ from $U(\hat w, \tilde w,s_0, *) = \hat w$, with this choice of $h$. Then $U(\hat w, \tilde w, t_0, *) = \tilde w$.
\end{lemma}

\begin{proof}
Notice that
\beq
\label{nbr-(*53a)}
    U(\hat w,\tilde w, t, x) = G_{a,b;t - s_0}(\hat w)(x) + \int_{s_0}^{t} ds\, G_{a,b;t - s}(f(s,\hat w,\tilde w))(x).         
    \eeq
Clearly, $f = f_1 + f_2$, where
\beqn
     f_1(s,\hat w,\tilde w)(z) = \frac{1}{t_0 - s_0}  G_{a,b;s - s_0}(\tilde w - \hat w)(z), \quad
     f_2(s,\hat w,\tilde w)(z) =   - \mA \tilde w(z),
    \eeqn
and since $\cL$ is a linear operator, $U = U_0 + U_1 + U_2$, where
\begin{align*}
      \cL U_0(t,x) &= 0, \qquad\qquad \qquad\;      U_0(s_0, *) = \hat w,\\
     \cL U_1(t,x) &= f_1(t, \hat w, \tilde w)(x),\quad      U_1(s_0, *) = 0,\\
     \cL U_2(t,x) &= f_2(t, \hat w, \tilde w)(x) ,\quad   U_2(s_0, *) = 0.   
   \end{align*}
Therefore, for $s_0 \leq t \leq t_0$,
\begin{align}
\label{nbr-(*54a)}
     U_0(\hat w, t, x) &= G_{a,b;t - s_0}(\hat w)(x),\notag\\
     U_1(\hat w, t, x) &= \int_{s_0}^{t} ds \int_0^1 dz\, G_{a,b}(t-s; x, z) f_1(s, \hat w, \tilde w)(z)\notag\\
     & = \frac{1}{t_0 - s_0} \int_{s_0}^{t} ds \int_0^1 dz\, G_{a,b}(t-s; x, z)\notag\\
     &\qquad \times \int_0^1 d\xi\, G_{a,b}(s-s_0,z, \xi)  (\tilde w(\xi) - \hat w(\xi))\notag\\
     &= \frac{1}{t_0 - s_0} \int_{s_0}^{t} ds \int_0^1 d\xi\, G_{a,b}(t - s_0; x, \xi) (\tilde w(\xi) - \hat w(\xi))\notag\\
     & = \frac{t - s_0}{t_0 - s_0}\, G_{a,b; t - s_0}(\tilde w - \hat w)(x),   
     \end{align}
and
\begin{align*}
    U_2(\hat w, t, x) &= \int_{s_0}^{t} ds \int_0^1 dz\, G_{a,b}(t-s; x, z) f_2(s, \hat w,\tilde w)(z)\\
    & = - \int_{s_0}^{t} ds \int_0^1 dz\, G_{a,b}(t-s; x, z)  \mA \tilde w(x) =  \tilde w(x) - G_{a,b; t - s_0}(\tilde w)(x),
    \end{align*}
where we have used the fact that $(t, x) \mapsto \tilde w(x)$ solves the PDE $\cL u(t, x) = - \mA \tilde w(x)$, with initial condition $u(s_0, *) = \tilde w$. Adding these three formulas, we find that for all $\hat w \in \bD$ and $x \in [0, 1]$,
    $U(\hat w, t_0, x) = \tilde w(x)$,
as claimed.
 \end{proof}

 \begin{remark}
 \label{nbr-r24}
  $U_0(\hat w, t, *)$ gives the evolution of the initial condition  $\hat w $ under  $\cL $, $U_1 $ is obtained by injecting at each time $t$ a fixed multiple of what the solution of the homogeneous PDE with initial condition  $\tilde w - \hat w $ would be at time $t$.  $U_2(\hat w, t, *) $ produces the difference of  $\tilde w $ and its evolution under  $\cL$ by injecting the source term  $- \tilde\cL \tilde w$ at all times  $t$. The superposition of these three evolutions moves the initial condition $\hat w$ at time  $s_0$ to the value  $\tilde w$ at time  $t_0$. This can be achieved for all initial conditions, smooth values at time  $t_0$, and arbitrarily short time intervals $]s_0, t_0]$. During a short time, with high probability, the contribution of $\dot W $ to the evolution of  $u$ will be quite small. This idea provides the basis for the proof of Proposition \ref{nbr-p21}.
\end{remark}
\smallskip

We can now prove the irreducibility of the transition function associated to the SPDE \eqref{nbr-(*Mb-bis)} in the case of additive noise.

\begin{thm}
\label{nbr-p21}
 The transition function $(P_t)$ of the solution $u$ to \eqref{nbr-(*Mb-bis)} with $\sigma\equiv 1$, defined in \eqref{ch5-added-1-2-2}, is irreducible at all times $t_0 > 0$.
 \end{thm}
 
 \begin{proof}
 Fix $t_0 > 0$, $r > 0$ and $u_0, w \in \bD$. We shall choose a bounded adapted process $h$ such that the solution $v$ to \eqref{nbr-(*43)} satisfies
      $P\{\Vert v(t_0, *) - w \Vert_\bD < r \} > 0$.
Proposition \ref{nbr-p22} will then imply that the same strict inequality holds with $v$ replaced by $u$, and therefore $P_t$ is irreducible.

     Fix $\tilde w \in \bD \cap \cC^2([0, 1])$ such that
     \beq
     \label{nbr-(*44a)}
      \Vert \tilde w  - w \Vert_\bD < \frac{r}{2}.     
      \eeq
Consider the solution $(\eta(t, x))$ to the linear SPDE
\beq
\label{nbr-(*44)}
      \cL\, \eta(t, x) = \sqrt{2}\, \dot W(t,x),\quad   t > 0, \  x \in\, ]0, 1[,   
      \eeq
with the same boundary conditions as in \eqref{nbr-(*43)} and initial condition $u_0\in \bD$. Then
\beq
\label{nbr- (*45)}
     \eta(t, x) = G_{a,b; t}(u_0)(x) +\int_0^t \int_0^1 G_{a,b}(t-s; x, y)\, W(ds, dy).      
     \eeq
Fix $s_0 \in [0, t_0[ $ and  $R > 0$, both to be chosen below. Define a Lipschitz truncation function  $\ell: \R_+ \to [0, 1]$ by  $\ell(\rho)=1$ on $[0,R]$,  $\ell$ is linear on  $[R, 2R] $ and  $\ell(\rho)=0$ on $[2R,\infty[$.  For $\hat w \in \bD$, define  $f(s,\hat w,\tilde w)(z) $ as in Lemma \ref{nbr-l23} and set
    \beq
    \label{noch(*2)}
     h(t, x) = 1_{[s_0, t_0[}(t)\,  f(t,\eta(s_0, *), \tilde w)(x) \,\ell(\Vert \eta(s_0, *) \Vert_\bD),
     \eeq
so that  $h: [s_0, t_0] \times [0, 1] \to \R$ is bounded.

Let  $(v(t, x)) $ be the solution of \eqref{nbr-(*43)} on $[0,t_0]$ with this choice of  $h$ and set
\beq
\label{nbr-(*45a)}
        \tilde U_1(\hat w,\tilde w)(x) = G_{a,b; t_0 - s_0}(\hat w)(x) + \int_{s_0}^{t_0} ds\, G_{a,b; t_0 - s}(f(s,\hat w,\tilde w))(x),   
        \eeq
and for $r \in [s_0, t_0]$, set
\beqn
         \tilde U_2(r, x) = \int_{r}^{t_0}\int_0^1 G_{a,b}(t_0 - s; x, y)\, W(ds, dy).
         \eeqn
By \eqref{nbr-(*4aaa)} 
 (with $s = s_0$, $t = t_0 - s_0$),
 \beqn
 \label{nbr-(*46)}
     v(t_0, *) = \tilde U_1(\eta(s_0, *), \tilde w)(*) + \tilde U_2(s_0, *).     
     \eeqn
Let $U(\hat w, \tilde w, \cdot, *)$ be the solution of \eqref{nbr-(*53)} given in Lemma \ref{nbr-l23}. Comparing \eqref{nbr-(*45a)} and \eqref{nbr-(*53a)}, notice that if $\Vert \hat w \Vert_\bD \leq R$,  then $\tilde U_1(\hat w,\tilde w)(*) = U(\hat w, \tilde w, t_0, *) = \tilde w$. In particular, if  at time $s_0$, $\Vert \eta(s_0, *) \Vert_\bD \leq R$, then $\tilde U_1(\eta(s_0, *), \tilde w) = \tilde w$.

    It remains to choose $s_0$ and $R$ so that with high probability, $\Vert \eta(s_0, *) \Vert_\bD \leq R$ and $\Vert \tilde U_2(s_0, *) \Vert_\bD$  is small.

        Since $(s_0, x) \mapsto \tilde U_2(s_0, x)$ has a continuous version, and $\tilde U_2(t_0, x) \equiv 0$, we can choose $s_0 \in\, ]0, t_0[$ sufficiently close to $t_0$ so that $P\{ \Vert \tilde U_2(s_0, *) \Vert_\bD \geq \frac{r}{2} \} \leq \frac{1}{4}$. Since $(t, x) \mapsto \eta(t, x)$ has a continuous version, $\lim_{R \to \infty} P\{\Vert \eta(s_0, *) \Vert_\bD \le R \} = 1$, so we can choose $R < \infty$ such that $P\{\Vert \eta(s_0, *) \Vert_\bD \leq R \} \geq \frac{3}{4}$. In this case, by \eqref{noch(*2)} and \eqref{nbr-(*44a)},
        \begin{align*}
   P\{\Vert v(t_0, *) - w \Vert_\bD < r \} &\geq P\left\{\tilde U_1(\eta(s_0, *), \tilde w) = \tilde w,\, \Vert \tilde U_2(s_0,*) \Vert_\bD < \frac{r}{2} \right\}\\
     &\geq P\left\{ \Vert \eta(s_0, *) \Vert_\bD \leq R,\, \Vert\tilde U_2(s_0,*) \Vert_\bD < \frac{r}{2}\right \}\\
     &\geq P\left\{ \Vert \eta(s_0, *) \Vert_\bD \leq R \right\} - P\left\{\Vert \tilde U_2(s_0,*) \Vert_\bD \geq \frac{r}{2}\right \}\\
     &\geq \frac{3}{4} - \frac{1}{4} = \half
      > 0,
      \end{align*}
as was to be proved.
\end{proof}

\begin{remark}
\label{nbr-r21a}
Proposition \ref{nbr-p21} remains true for the stochastic heat equation with multiplicative noise (see \cite[Theorem 7.3.1]{dz-1996}).
\end{remark}

\section
{Large interval asymptotics (multiplicative noise)}
\label{new-5.5-candil}

In this section, we are interested in comparing, as $L \to \infty$, the solution to the stochastic heat equation on $[-L, L]$ with vanishing Dirichlet boundary conditions, and the solution to the same equation on $\R$, when both equations have the same space-time white noise, the same coefficients and the same initial condition (see Section \ref{ch1'-s6} for the setting). In particular, let $\sigma(t,x,z)$ and $b(t, x, z)$ be two functions satisfying $({\bf H_L})$ with $D = \R$, $v_{0}$ be a function satisfying $({\bf H_I})$ with $D = \R$ (see Section \ref{ch4-section1} for the formulation of these assumptions), and $W = (W(A),\, A\in\mathcal{B}^f_{\IR_+\times \IR })$ be a space-time white noise. For other types of boundary conditions, see Remark \ref{rd01_01r1}.
\bigskip

\noindent{\em Stating the problem}
\medskip

 Let $\cL = \frac{\partial}{\partial t} - \frac{\partial^{2}}{\partial x^{2}}$ and consider the SPDEs
\begin{equation}
\label{rd08_09e1}
\begin{cases}
 \cL\, u_L(t, x) = b(t, x, u_L(t, x)) + \sigma(t,x,u_L(t, x))\, \dot W(t, x),\\
  ~  \qquad\qquad\qquad\qquad\qquad\qquad(t, x) \in\, ]0,T]\, \times\, ]-L, L[\, ,\\
 u_L(0, x) = v_{0}(x), \qquad\qquad\quad \, x\in\, ]-L, L[\, , \\
 u_L(t, 0) = u_L(t, L) =0, \qquad \ \,t\in\, ]0,T] ,
 \end{cases}
\end{equation}
and
\begin{equation}
\label{rd08_09e2}
\begin{cases}
\cL\, u(t, x) = b(t, x, u(t, x)) + \sigma(t,x,u(t, x))\, \dot W(t, x), \\
   ~  \qquad\qquad\qquad\qquad\qquad\qquad (t, x) \in\, ]0,T] \times \R ,\\
 u(0, x) = v_{0}(x), \qquad\qquad\quad\ \, x\in \R . \\
 \end{cases}
\end{equation}

Let $\Gamma_L(t; x, y)$ be the Green's function associated to the heat equation on $[-L, L]$ with Dirichlet boundary conditions and let $\Gamma(t, x)$ be the fundamental solution of the heat operator on $\R$ given in \eqref{heatcauchy-1'}. It is straightforward to check that $\Gamma_L$ is related to the Green's function $G_{2L}$ in \eqref{ch1.600} as follows:
$\Gamma_L(t; x, y) = G_{2L}(t; x+L, y+ L)$. Therefore, by \eqref{ch1.600} and \eqref{ch1.6000-double},
\begin{align}
\nonumber
\Gamma_L(t; x, y) &= \frac{1}{L} \sum_{n=1}^\infty e^{-\frac{\pi^2}{4 L^2}n^2 t} \sin\left(\tfrac{n \pi}{2L}(x+L) \right) \, \sin\left(\tfrac{n \pi}{2L}(y+L) \right) \\
   &= \sum_{m \in \Z} \left[\Gamma(t,y-x+4mL) - \Gamma(t,x+y+(4m + 2)L)\right].
\label{rd08_09e5}
\end{align}
It turns out that for large $L$, $\Gamma_L(t; x, y)$ and $\Gamma(t, x-y)$ are not so different (see Lemma \ref{rd08_10l3}). This will lead to the Theorem \ref{rd08_09t1} below.

In the sequel, we set $D_L = [-L, L]$ and recall that the random field solutions to \eqref{rd08_09e1} (respectively \eqref{rd08_09e2}) are given by
\begin{align}
\label{ul}
u_L(t,x) = I_{0,L}(t,x) &+ \int_0^t \int_{D_L} \Gamma_L(t-s; x,y)\, \sigma(s,y,u_L(s,y))\, W(ds,dy)\notag\\
&+  \int_0^t ds \int_{D_L} dy\, \Gamma_L(t-s; x,y)\, b(s,y,u_L(s,y))
\end{align}
and
\begin{align}
\label{u}
    u(t, x) =  I_0(t,x) &+ \int_0^t \int_\R \Gamma(t-s, x-y)\, \sigma(s,y,u(s,y))\, W(ds,dy)\notag\\
&+  \int_0^t \int_\R \Gamma(t-s, x-y)\, b(s,y,u(s,y))\, ds dy,
\end{align}
respectively, where
\beqn
I_{0,L}(t,x) = \int_{-L}^L\Gamma_L(t;x,y)\, v_0(y)\, dy ,\quad  I_0(t,x)=\int_\re\Gamma(t,x-y)\, v_0(y)\,dy.
\eeqn
Putting $\Vert v_0\Vert_\infty= \sup_{x \in \re}|v_0(x)|$, we have
\beqn
\Vert I_0\Vert_\infty:=\sup_{(t,x) \in [0,T]\times \re} |I_0(t,x)| \le \Vert v_0\Vert_\infty,
\eeqn
and
\beqn
\sup_{L >  0}\Vert I_{0,L}\Vert_\infty:=\sup_{L > 0}\, \sup_{(t,x) \in [0,T] \times [-L,L]}\, |I_{0,L}(t,x)| \le \Vert v_0 \Vert_\infty.
\eeqn
Indeed, $\int_{\re}\Gamma(t,x-y)\,dy = 1$ and, by Proposition \ref{ch1'-pPD} (ii), 
\beqn
\sup_{L > 0}\Vert I_{0,L}\Vert_\infty \le \sup_{(t,x) \in [0,T]\times \re} \int_{\re} \Gamma(t,x-y)\,  |v_0(y)|\,dy.
\eeqn
\medskip

\noindent{\em The limit theorem}
\medskip

   Throughout this section, we assume that $v_0$ is Borel and bounded, that is, $\Vert v_0\Vert_\infty<\infty$, which implies that
\beq
\label{bI0}
\Vert I_0\Vert_\infty <\infty\ {\text{and}}\ \sup_{L > 0} \Vert I_{0,L} \Vert_\infty < \infty.
\eeq

\begin{thm}
\label{rd08_09t1}
Fix $T > 0$ and $p \geq 1$. Let $v_0: \R \to \R$ be Borel and bounded.
Assume that the functions $\sigma$ and $b$ in \eqref{rd08_09e1} (and \eqref{rd08_09e2}) satisfy $({\bf H_L})$ with $D = \R$
and let $K$ be the constant that appears in the global Lipschitz condition.
Let $(u_L(t, x),\, (t, x) \in [0, T] \times [-L, L])$ (respectively $(u(t, x),\, (t, x) \in [0, T] \times \R)$) be the random field solution of \eqref{rd08_09e1} (respectively \eqref{rd08_09e2}). Then there is a constant $c = c(T, K, \Vert v_0 \Vert_{\infty}, p) < \infty$ such that, for all $L > 0$ and $(t, x) \in [0, T] \times [-L, L]$,
\beq
\label{rd08_10e4}
  \Vert u(t, x) - u_L(t, x) \Vert_{L^p(\Omega)} \leq c \left(\exp\left(-\tfrac{(L-x)^2}{8t} \right) + \exp\left(-\tfrac{(L+x)^2}{8t} \right) \right).
\eeq
In particular, for all $(t, x) \in [0, T] \times \R$,
\beq
\label{rd08_13e1}
   \lim_{L \to \infty} u_L(t, x) =  u(t, x) \quad\text{ in }  L^p(\Omega).
\eeq
\end{thm}

\begin{remark}
\label{rem-to-rd08_09t1}
If $v_0: \re\times \Omega\rightarrow \re$ is a measurable random variable that is independent of the noise $\dot W$ and such that $\sup_{x\in\re}\Vert v_0(x)\Vert_{L^p(\Omega)}<\infty$, then the conclusion of
Theorem \ref{rd08_09t1} still holds $($with $\Vert v_0 \Vert_{\infty}$ there replaced by $\sup_{x\in\re}\Vert v_0(x)\Vert_{L^p(\Omega)})$.
\end{remark}

For the proof of Theorem \ref{rd08_09t1}, we will need the following lemma.

\begin{lemma}
\label{rd08_09l1}
The assumptions are the same as in Theorem \ref{rd08_09t1}. For all $p > 0$,
\beq
\label{rd08_09e3}
   \sup_{(t , x) \in [0, T] \times \R} E[\vert u(t, x)\vert^p] < \infty
\eeq
and
\beq
\label{rd08_09e4}
    \sup_{L > 0}\, \sup_{(t , x) \in [0, T] \times D_L} E[\vert u_L(t, x)\vert^p] < \infty.
\eeq
\end{lemma}

\begin{proof}
Property \eqref{rd08_09e3} is a direct consequence of Theorem \ref{recap}. For \eqref{rd08_09e4}, we fix $p \geq 2$. Using \eqref{ul}, we see that
\begin{align*}
   E\left[\vert u_L(t, x)\vert^p\right] &\leq 3^{p-1}\left(\vert I_{0,L}(t,x)\vert^p
   + E\left[\vert \cI_L(t, x)\vert^p + \vert \cJ_L(t, x)\vert^p\right]\right),
\end{align*}
where
\begin{align*}
     \cI_L(t, x) &= \int_0^t \int_{D_L} \Gamma_L(t-s; x,y)\, \sigma(s,y,u_L(s,y))\, W(ds,dy), \\
     \cJ_L(t, x)&= \int_0^t ds \int_{D_L} dy\, \Gamma_L(t-s; x,y)\, b(s,y,u_L(s,y)).
\end{align*}
Proceeding as in the proof of Theorem \ref{ch1'-s5.t1}, that is, using Burkholder's inequality (see \eqref{ch1'-s4.4}) then Hölder's inequality for $\cI_L$, and Hölder's inequality for $\cJ_L$, we see that
\begin{align*}
   &E[\vert u_L(t, x)\vert^p] \\
    &\qquad\leq C_p\, \left(\vphantom{ \left[\int_0^t ds \, J_{1,L}(s) \right]^{\frac{p}{2} - 1}  } \Vert v_0\Vert_\infty^p \right.\\
    &\qquad\qquad +\left[\int_0^t ds \, J_{1,L}(s) \right]^{\frac{p}{2} - 1} \int_0^t ds \, J_{1,L}(t-s)
        \left(1 + \sup_{x \in D_L} E[\vert u_L(s, x)\vert^p]\right) \\
    &\qquad\qquad+ \left. \left[\int_0^t ds \, J_{2,L}(s) \right]^{p - 1} \int_0^t ds \, J_{2,L}(t-s) 
        \left(1 + \sup_{x \in D_L} E[\vert u_L(s, x)\vert^p]\right)\right),
\end{align*}
where $C_p$ depends on $p$ but not on $L$, and
\beqn
J_{1,L}(s) = \sup_{x\in D_L} \int_{D_L} \Gamma_L^2(s; x,y)\, dy,\qquad
J_{2,L}(s) = \sup_{x\in D_L} \int_{D_L} \Gamma_L(s; x,y)\, dy.
\eeqn
By Proposition \ref{ch1'-pPD} (ii), $J_{1,L}(s) + J_{2,L}(s) \leq J(s)$, where
\beqn
    J(s) = \sup_{x\in \R}  \int_\R  \left(\Gamma^2(s, x-y) + \Gamma(s, x-y)\right)\, dy = \frac{1}{\sqrt{8 \pi s}}+1
\eeqn
by \eqref{heatcauchy-11'}. Let $g_L(t) = \sup_{x\in D_L}  E[\vert u_L(t, x)\vert^p]$. Then
\beqn
  g_L(t) \leq \tilde C_{p,T} \left(\Vert v_0 \Vert_\infty^p + \int_0^T ds \, (1 + f_L(s))\, J(t-s) \right).
\eeqn
By the Gronwall-type Lemma \ref{A3-l01} (c), we see that there is $C < \infty$ such that for all $L > 0$ and all $t \in [0, T]$, $g_L(t) \leq C(1+ \Vert v_0\Vert_\infty^p)$. Since the right-hand side does not depend on $t$ or $L$, \eqref{rd08_09e4} is proved.
\end{proof}

\noindent{\em Proof of Theorem \ref{rd08_09t1}.}\
For $L > 0$, let $$H_L(t; x, y) = \Gamma(t, x-y) - \Gamma_L(t; x, y).$$ 
By Proposition \ref{ch1'-pPD} (ii), we have $H_L \geq 0$. Using \eqref{u} and \eqref{ul}, we write
\beq
\label{rd08_09e6}
    u(t, x) - u_L(t, x) = \sum_{i=0}^3 A_i(t, x),
\eeq
where
\begin{align*}
    A_0(t, x) &= \int_{D_L^c} dy\,  \Gamma(t,x-y)\, v_0(y), \\
    A_1(t, x) &= \int_{D_L} dy\, H_L(t; x, y)\, v_0(y) \\
     &\qquad + \int_0^t \int_{D_L} H_L(t-s; x, y)\\
     &\qquad\qquad\qquad \times    [\sigma(s,y,u_L(s,y))\, W(ds,dy) 
         + b(s,y,u_L(s,y))\, ds dy\, ],\\
          A_2(t, x) &=  \int_0^t \int_{D_L^c} \Gamma(t-s, x-y)\, [\sigma(s,y,u(s,y))\, W(ds,dy)  \\
      &\qquad\qquad\qquad\qquad  + b(s,y,u(s,y))\, ds dy\, ], \\
    A_3(t, x) &= \int_0^t \int_{D_L} \Gamma(t-s, x-y)\\
     &\qquad\qquad\qquad \times    [(\sigma(s,y,u(s,y)) - \sigma(s,y,u_L(s,y))) \, W(ds,dy) \\
     &\qquad\qquad\qquad\qquad  + (b(s,y,u(s,y)) - b(s,y,u_L(s,y)))\, ds dy\, ].
\end{align*}

Let $f_L(t, x) = \Vert u(t, x) - u_L(t, x) \Vert_{L^p(\Omega)}$. By the triangle inequality,
\beq
\label{rd08_10e0}
    f_L(t, x) \leq\sum_{i=0}^3 \Vert A_i(t, x) \Vert_{L^p(\Omega)}.
\eeq

Define
\beq
\label{aela}
    a_L(t, x) = \exp\left(-\tfrac{(L-x)^2}{8t} \right) + \exp\left(-\tfrac{(L+x)^2}{8t} \right).
\eeq
 By \eqref{rdlemC.2.2-before-1} and since $v_0$ is bounded,
\beq
\label{bound0}
  \Vert A_0(t, x) \Vert_{L^p(\Omega)} =\vert A_0(t, x)\vert \leq  C_0\,a_L(t, x),
\eeq

Proceeding as in the proof of Theorem \ref{ch1'-s5.t1} and using \eqref{rd08_09e4}, we see that
\begin{align}
 \label{rd08_10e1a}
\Vert A_1(t, x) \Vert_{L^p(\Omega)} &\leq \tilde C_1 \left(\int_{D_L} dy\,  H_L(t; x, y)  + 
 \left[\int_0^t ds \int_{D_L} dy \,  H_L^2(t-s; x, y)\right]^{\frac{1}{2}}\right.\notag\\
 &\left.\qquad\qquad+ \int_0^t ds \int_{D_L} dy \,  H_L(t-s; x, y) \right)\notag \\
 & \leq \tilde C_{1, T} \left(\exp\left(- \frac{(L - x)^2}{4 t} \right) + \exp\left(- \frac{(L + x)^2}{4 t} \right) \right) \notag \\
  & \leq C_1\, a_L(t, x)
\end{align}
by Lemma \ref{rd08_10l3}.

Proceeding again as in the proof of Theorem \ref{ch1'-s5.t1} but using \eqref{rd08_09e3} instead of \eqref{rd08_09e4}, we see that
\begin{align}
 \label{rd08_10e2}
     \Vert A_2(t, x) \Vert_{L^p(\Omega)}  
     &\leq  \tilde C_2 \left(\left[\int_0^t ds \int_{D_L^c} dy\, \Gamma^2(t-s, x-y)\right]^{\frac{1}{2}}\right.\notag\\
      &\left.\qquad\qquad + \int_0^t ds \int_{D_L^c} dy\, \Gamma(t-s, x-y)\right)\notag  \\
      & \leq \tilde C_{2, T} \left(\exp\left(- \frac{(L - x)^2}{4 t} \right) + \exp\left(- \frac{(L + x)^2}{4 t} \right) \right) \notag \\
      &\leq  C_2\, a_L(t, x)
\end{align}
by Lemma \ref{rdlemC.2.2-before}. Putting together \eqref{bound0}--\eqref{rd08_10e2}, we obtain
\beq
\label{zero-one-three}
\sum_{i=0}^2 \Vert A_i(t, x) \Vert_{L^p(\Omega)} \le C a_L(t,x).
\eeq

Next, we study $A_3(t,x)$. Define
\begin{align*}
A_{3,1}(t,x)&= \int_0^t\int_{D_L} \Gamma(t-s,x-y)[\sigma(u(s,y))-\sigma(u_L(s,y))]\, W(ds,dy),\\
A_{3,2}(t,x)&= \int_0^t ds \int_{D_L} dy\,\Gamma(t-s,x-y)[b(u(s,y))-b(u_L(s,y))],
\end{align*}
so that $A_3(t,x) = A_{3,1}(t,x)+A_{3,2}(t,x)$. Using Burkholder's inequality and the Lipschitz property in $(\bf{H_L})$, we obtain
\beqn
\Vert  A_{3,1}(t,x)\Vert_{L^p(\Omega)}^2 \le C\left\Vert  \int_0^t ds \int_{D_L} dy\,\Gamma^2(t-s,x-y)\,|u(s,y)-u_L(s,y)|^2\right\Vert_{L^{p/2}(\Omega)}.
\eeqn
For $(t,x)\in [0,T]\times D_L$, $(s,y)\mapsto \Gamma^2 (t-s,x-y)$ is the density of a finite measure on $[0,T]\times D_L$, therefore we can apply Minkowski's inequality to bound from above the last right-hand side by
\beqn
   C \int_0^t ds \int_{D_L} dy\,\Gamma^2(t-s,x-y) \Vert u(s,y)-u_L(s,y)\Vert_{L^p(\Omega)}^2.
\eeqn
It follows that
\beqn
\label{a3.1}
\Vert  A_{3,1}(t,x)\Vert_{L^p(\Omega)}^2 \le C\int_0^t ds \int_{D_L} dy\,\Gamma^2(t-s,x-y) f_L^2(s,y).
\eeqn
Applying Minkowski's inequality to $A_{3,2}(t,x)$, we see that
\beqn
\label{a3.3}
\Vert  A_{3,2}(t,x)\Vert_{L^p(\Omega)} \le C\int_0^t ds \int_{D_L} dy\,\Gamma(t-s,x-y) f_L(s,y).
\eeqn
Writing $f_L(s,y)  = f_L(s,y) \cdot 1$ and applying the Cauchy-Schwarz inequality to this integral, we obtain
\beq
\label{a3}
\Vert  A_{3}(t,x)\Vert_{L^p(\Omega)}^2\le \bar C\int_0^t ds \int_{D_L} dy\,\left[\Gamma^2(t-s,x-y)+ \Gamma(t-s,x-y)\right] f_L^2(s,y).
\eeq

Notice that
\beqn
   \Gamma(r, z) + \Gamma^2(r, z) = \Gamma(r, z) + 
     \frac{1}{4 \pi r} \exp\left(- \frac{z^2}{2 r} \right) \leq \Gamma(r, z)\left(1 + \frac{1}{\sqrt{4 \pi r}}\right) , 
\eeqn
and for $r \in [0, T]$,
\beqn
\Gamma(r, z) \left(1 + \frac{1}{\sqrt{4 \pi r}}\right) \leq  \frac{c}{\sqrt r} \Gamma(r, z),
\eeqn
where $c = \frac{\sqrt{4 \pi T} + 1}{\sqrt{4 \pi}}$. Let
\beqn
J(r,z):= \frac{1}{\sqrt r}\Gamma(r, z).
\eeqn
Then \eqref{a3} implies that
\beq
\label{rd08_12e7}
    \Vert A_3(t, x) \Vert_{L^p(\Omega)}^2 \leq C_2\int_0^t ds \int_{D_L} dy \, J(t - s, x - y) \, f_L^2(s, y).
\eeq
For its further use, we note that the function $(r,z)\rightarrow J(r,z)$ from $\re_+\times \re$ into $\re$ satisfies \eqref{rd08_10e6}. 

Using \eqref{rd08_10e0}, \eqref{zero-one-three} and \eqref{rd08_12e7}, we see that for $(t, x) \in [0, T] \times D_L$,
\beq
\label{rd08_10e3}
   f_L^2(t, x) \leq  c \, a_L^2(t, x) + \tilde c \int_0^t ds \int_{\R} dy \, J(t - s, x - y)\, 1_{D_L} (y)\, f_L^2(s, y).
\eeq
By the space-time Gronwall's Lemma \ref{rd08_10l2} (d)  (with $z_0\equiv 0$ there), we conclude that for $(t, x) \in [0, T] \times D_L$,
\beq
\label{rd08_13e2}
   f_L^2(t, x) \leq \tilde C \left[ a_L^2(t, x) +  (\cK \star a^2_L) (t, x)\right],
\eeq
where, for $(t, x) \in [0, T] \times \R$,
\beqn
   \cK(t, x) = \sum_{\ell = 1}^\infty \, J^{\star \ell}(t, x)
\eeqn
and
\beqn
(\cK \star a_L^2) (t, x) = \int_0^t ds \int_{\R} dy \, \cK(t-s, x -y)\, a_L^2(s, y).
\eeqn

Defining $\Gamma_\nu(t,x)$ as in \eqref{fsnu}, we see that $J(t,x)=4\sqrt \pi\, \Gamma_{4}^2(t,x)$, and $\Gamma_2(t,x) = \Gamma(t,x)$. Thus, the kernel $\cK$ is given by the expression in \eqref{rd08_12e8}, that is,
\beq
\label{boundfork} 
\cK(t,x) = B(t)\, \Gamma_2(t,x) =B(t) \, \Gamma(t,x),
\eeq
where
\beqn
    B(t) = \left(\frac{4 \sqrt{\pi}}{\sqrt t}  + 16\pi e^{\pi t} \Phi\left(\sqrt{2\pi t}\right)\right)
\eeqn
and $\Phi(\cdot)$ is the standard Normal distribution function. 

The last part of the proof consists in showing that up to a positive multiplicative constant, $(\cK \star a_L^2) (t, x)$ is bounded by $a_L^2(t, x)$.

Indeed, $(\cK \star a_L^2) (t, x)$ is the sum of two terms $I_\pm(t,x)$, where
\beqn
   I_\pm(t,x) = \int_0^t ds \int_{\R} dy \, \cK(t-s, x -y) \exp\left(-\tfrac{(L \pm y)^2}{4s} \right). 
\eeqn
Using \eqref{boundfork}, we see that
\begin{align*}
    I_\pm(t,x) &= \int_0^t ds \int_{\R} dy \, \cK(t-s, x -y)\, \sqrt{4 \pi s}\, \Gamma\left(s, L \pm y\right) \\
      & \leq C \int_0^t ds  \left(\tfrac{1}{\sqrt{t-s}} + e^{\pi\,(t-s)} \Phi\left(\sqrt{2 \pi\, (t-s)} \right) \right) \, \sqrt{s} \\
      &\qquad\qquad \times \int_{\R} dy \, \Gamma(t-s, x-y) \, \Gamma\left(s, L \pm y\right)\\
      &\leq C_T \int_0^t ds  \left(\tfrac{1}{\sqrt{s}} +1 \right) \, \sqrt{t-s} \\
      &\qquad\qquad \times \int_{\R} dy \, \Gamma(s, x-y) \, \Gamma\left(t-s, L \pm y\right).
   \end{align*}  
Next, we apply the semigroup property \eqref{semig-heat} to deduce that
\beq
\int_{\R} dy \, \Gamma(s, x-y) \, \Gamma\left(t-s, L \pm y\right) = \Gamma\left(t,x\pm L \right).
\eeq
Notice that 
\beqn
\Gamma\left(t, x\pm L\right) = \frac{1}{\sqrt{4\pi t}}\, e^{-\frac{(x\pm L)^2}{4 t}}  . 
\eeqn
Thus,
\beqn
 I_\pm(t,x) \le C_T\,  e^{-\frac{(x\pm L)^2}{4t}}\int_0^t ds \left(\frac{1}{\sqrt s}+1\right) \sqrt{\frac{t-s}{t}}\le \tilde C_T\,  e^{-\frac{(x\pm L)^2}{4t}},
 \eeqn
and consequently, for some constant $c_T$, for all $(t, x) \in [0, T] \times D_L$,
\beqn
(\cK \star a_L^2) (t, x)\le c_T\,  a_L^2(t, x).
\eeqn
From the inequality \eqref{rd08_13e2}, we conclude that $f_L(t,x)$ is bounded above by a constant times $a_L(t,x)$, that is,
\eqref{rd08_10e4} is satisfied.
\qed

\begin{remark}\label{rd01_01r1} The case of vanishing Neumann boundary conditions, as well as certain mixed boundary conditions, has been considered in \cite[Chapter 2]{candil2022}; the result for Neumann boundary conditions is slightly different (see also \cite{CDS}).
\end{remark}


\section{Notes on Chapter \ref{ch6}}
\label{notes-ch6}
Our interest in the material of Section \ref{ch5-added-1} originates with \cite[Chapter 3]{Pu-thesis}. Since the theory of Markov processes is so fundamental, we decided to include in Section \ref{ch5-added-1-s1} an introduction to the main concepts of this theory, presented in the context of SPDEs, using \cite{ry} as our main reference. 

The study of invariant measures for parabolic SPDEs with gradient-type drift goes back at least to \cite{marcus-1974}. A more complete study of the existence and the identification of the reversible measures for semilinear stochastic heat equations with gradient-type drift and Dirichlet or Neumann boundary conditions \eqref{nbr-(*1aa)} (including the linear equations \eqref{ch1'.HD-ch5-markov} and \eqref{s5.3-new(*1)}) goes back to T. Funaki \cite[Theorem 4.1]{funaki83}, who used a different method than the one we present. Funaki's results were extended by various authors, including J. Zabczyk \cite{zabczyk-1989}, which also contains additional references. Further extensions are contained in Hairer et al \cite{H-S-V-W-2005} and \cite{H-S-V-W-2007}.

The explicit calculations in Propositions \ref{ch5-added-1-s2-p1} and \ref{ch5-added-1-s2-p2}, and in Theorem \ref{ch5-added-1-s2-t1} are taken from \cite{Pu-thesis}. Proposition \ref{s5.3-new-bis-p*1} gives a related example of convergence to the invariant measure for an SPDE on $\rek$. The results of Propositions \ref{ch5-added-1-s2-p1} and Proposition \ref{rem6.1.15-*2}, as well as Theorem \ref{s5.3-new-t*1}, are also contained in Funaki's paper, which deals in fact with systems of SPDEs. However, we were inspired by the articles of Hairer et al \cite{H-S-V-W-2005} and \cite{H-S-V-W-2007}, who developed the results of Theorems \ref{s5.3-new-t*1}, \ref{nbr-t1}, \ref{nbr-t-14} and \ref{nbr-t27} (and extended the results of \cite{zabczyk-1989}) with the objective of sampling numerically from the distribution of conditioned diffusions. In particular, \cite{H-S-V-W-2007} contains a general form of the results presented in Sections \ref{ch5-ss-in-rev} (reversibility and bridge measures) and \ref{new-5.5} (convergence to the reversible measure).

The study of reversible measures in the finite-dimensional case of Section \ref{new-5.4.1} is well-known, and seems to originate with \cite{kolmogorov-1937}. Nevertheless, it is not available in very many standard references (for instance, \cite[Chapter 5, Section 6, Problem 6.15 p. 361]{ks} discusses the case of a bounded gradient, which, however, does not cover the SDE in Proposition \ref{nbr-1}; a more complete discussion in the case of compact state spaces is given in \cite[Chapter 5, Theorem 5.6]{ikeda-watanabe-1981}), so we have included it for the sake of completeness.

Proposition \ref{nbr-4} can be deduced from \cite[Theorem 1]{zabczyk-1989}. An abstract form of Theorem \ref{nbr-t1} is quoted in \cite[Chapter 8, Theorem 8.6.3]{dz-1996}, with a sketch of proof and a reference to \cite[Theorem 2]{zabczyk-1989}. We have particularized the proof of \cite{zabczyk-1989} to our context, preserving the main ideas with fewer technical prerequisites. A generalization of \cite[Theorem 2]{zabczyk-1989} is given in \cite[Theorem 3.2]{H-S-V-W-2007}.

Theorem \ref{split-th 5.3.16} is a special case of \cite[Theorem 3.4]{H-S-V-W-2005}, in which the ``bridge distribution" can be given explicitly. Theorem \ref{nbr-t-14} is a particular case of \cite[Theorem 3.2]{H-S-V-W-2007}. In order to have a more complete picture of the invariant measures for stochastic heat equations with gradient-type drift, beyond the existence and reversibility properties, we decided to include the main results on uniqueness of and convergence to the invariant measure.
As mentioned in Remark \ref{nbr-r-25a}, Proposition \ref{nbr-25} is a particular case of
the Feldman-Hajek theorem \cite[Theorems 2.23 and 2.25]{dz};
 its proof contains the main ideas for the proof of the Feldman-Hajek theorem, but in a simpler context. In order to avoid entering into too many specifics of Markov processes, we quote Doob's theorem
 (Theorem \ref{nbr-20b}, taken from \cite[Theorem 4.2.1]{dz-1996}, without proof.
 Proposition \ref{nbr-p22} is \cite[Lemma 7.3.2]{dz-1996}.
 The conclusions of Theorem \ref{nbr-t27} can also be obtained (via the strong Feller property) from results in \cite{dz-1996}, including \cite[Theorems 4.2.1, 7.1.1 and 7.3.1]{dz-1996}.
 
 Theorem \ref{nbr-p21} particularises \cite[Theorem 7.3.1]{dz-1996} to the SPDE \eqref{nbr-(*Mb-bis)} in the case of additive noise. We hope that particularising these abstract results to a simpler setting will make them accessible to a wider audience. 
 
 Finally, the results of Section \ref{new-5.5-candil} have been established in \cite[Theorem 2.5]{candil2022}. 



\chapter[Selected results on SPDEs driven by space-time \texorpdfstring{\\}{}white noise]{Selected results on SPDEs driven by space-time white noise}
\label{ch2'}
\pagestyle{myheadings}
\markboth{R.C.~Dalang and M.~Sanz-Sol\'e}{Selected results on SPDEs driven by space-time white noise}

We devote this chapter to a selection of topics on SPDEs 
which we find  to be of particular interest. The first section discusses the notion of {\em weak solution in law}, in contrast with that of {\em random field solution} that has been studied in Chapters \ref{chapter1'} and \ref{ch1'-s5}. In Section \ref{rd03_17s1}, we use Girsanov's theorem and the results of Section \ref{ch2'-s1} to compare the laws of the solutions to the stochastic heat equations on $[0, L]$ and on $\R$, and we establish a germ-field Markov property of their solutions. 
Section \ref{ch2'-section5.2} is devoted to finding explicit exponential $L^p$-bounds on moments of the random field solutions of various SPDEs, providing information on the behaviour of these moments as $t \to \infty$. 
In Section \ref{ch5-added-0}, we prove a comparison theorem for the stochastic heat equation, first on a bounded interval and then on the real line. In Section \ref{ch2'-s7}, we introduce some elements of probabilistic potential theory focused on the notion of polarity of points. In Section \ref{rdrough}, we discuss SPDEs with rough initial conditions. Finally, in Section \ref{rd1+1anderson}, we discuss moment Lyapounov exponents for various SPDEs and we give a formula for the second moment Lyapounov exponent. In the case of the stochastic heat and wave equations, this formula can be evaluated explicitly, while for the fractional stochastic heat equation, it provides good upper and lower bounds.

\section{Weak solutions in law to SPDEs}
\label{ch2'-s1}

In the classical theory of PDEs, a weak solution refers to a generalized function that, when evaluated on Schwartz test functions, satisfies the equation. This notion has been used in the theory of SPDEs since its inception (see e.g. \cite[p. 313]{walsh}). However, in this section, the term {\em weak} does not refer to this analytical notion but to a probabilistic notion that we call {\em weak solution in law} and that is made precise in Definition \ref{ch2'-s3.4.1-d1} below. For stochastic differential equations, the analogous notion has been used and applied extensively (see for instance \cite[Definition 1.2]{ikeda-watanabe-1981}).

In this section, we consider SPDEs such as those considered in Chapter \ref{ch1'-s5}, mostly with additive noise, and we study their weak solutions in law (in the sense of Definition \ref{ch2'-s3.4.1-d1}).
The question of existence is addressed in Section \ref{ch2'-ss1.1}. Uniqueness in law is discussed in Section \ref{ch2'-s3.4.2}.

\subsection{Main definition}
\label{ch2'-s1-ss1}

Let us begin with some preliminaries.
Throughout this section, we will consider the SPDE \eqref{ch1'-s5.0} with $k=1$ and its integral formulation \eqref{ch1'-s5.1}. For the sake of completeness, we recall the setting.

Let $D\subset \re$ be a bounded or unbounded domain with smooth boundary, $T>0$ and let $W$ be a space-time white noise as defined in Proposition \ref{rd1.2.19}. We consider the integral formulation of \eqref{ch1'-s5.0}:
\begin{align}
\label{ch2'-s3.4.1.1}
u(t,x)& = I_0(t,x) + \int_0^t \int_D \Gamma(t,x;s,y) \sigma(s,y,u(s,y))\, W(ds,dy)\notag\\
& \qquad+  \int_0^t \int_D \Gamma(t,x;s,y) b(s,y,u(s,y))\, ds dy,
\end{align}
for $(t,x)\in[0,T]\times D$, where $\Gamma$ denotes the fundamental solution or the Green's function relative to the operator $\mathcal{L}$ (and the boundary conditions, if present). The term $I_0(t,x)$ is the solution to the homogeneous PDE $\mathcal{L} I_0 (t, x) =0$ (with the same initial and boundary conditions as for \eqref{ch1'-s5.0}). In the sequel, we assume that the functions $\sigma$ and $b$ do not depend on $\omega \in \Omega$.

We now give the notion of {\em weak solution in law}.\index{weak!solution in law}\index{law!weak solution in}\index{solution!in law, weak}
\begin{def1}
\label{ch2'-s3.4.1-d1}
Given the initial condition $u_0$ and the functions $\sigma$ and $b$ from $[0, T] \times D \times \R$ to $\R$, a {\em weak solution in law} to \eqref{ch2'-s3.4.1.1}
is a triplet $(\Theta, W, u)$, where $\Theta$ is a stochastic basis $(\Omega,\tf, P, (\tf_t, \, t \in [0, T]))$  carrying a space-time white noise $W$, and
\beqn
u=(u(t,x),\, (t,x)\in[0,T]\times D)
\eeqn
 is a jointly measurable and adapted random field on $\Theta$ such that for each $(t,x)\in[0,T]\times D$, the equation \eqref{ch2'-s3.4.1.1} holds a.s. (It is implicitly assumed that  in \eqref{ch2'-s3.4.1.1}, the stochastic integral is well-defined in the sense of Section \ref{ch2'new-s3} and the deterministic integral is well-defined as a Lebesgue integral).
\end{def1}

Clearly, a random field solution of \eqref{ch2'-s3.4.1.1} (in the sense of Definition \ref{ch4-d2-1'} or \ref{def4.1.1}) provides a weak solution in law.

Consider $(\Omega, \tf, P)$, $(\tf_t,\, t \in [0, T])$ and $W$ as in Section \ref{ch2new-s2}. Let $(h(t,x),\, (t,x)\in[0,T]\times D)$ be a jointly measurable and adapted random field with sample paths in $L^2([0,T]\times D)$ a.s. Define a measure  $\tilde P$ on $(\Omega,\mathcal{F})$  by
\beqn
\frac{d\tilde P}{dP} = \exp\left(-\int_0^T  \int_D \ h(t,x)\, W(dt,dx) - \frac{1}{2}\int_0^T dt \int_D dx\, h^2(t,x)\right).
\eeqn
An immediate consequence of Girsanov's Theorem \ref{ch2'-s3-t1} is the following.

\begin{thm}
\label{ch2'-s3.4.1-t0}
Suppose that $E\left[\frac{d\tilde P}{dP}\right] = 1$, so that $\tilde P$
 is a probability measure. Assume that $\Gamma$ satisfies the conditions $({\bf H_\Gamma})(i), (ii)$ and $(iiia)$ of Section \ref{ch4-section0}. Let $u$ be defined for $(t,x)\in[0,T]\times D$ by
  \begin{align}
\label{ch2'-s3.4.1.30}
u(t,x) &= I_0(t,x)  + \int_0^t \int_D \Gamma(t,x;s,y)\, W(ds,dy)  \notag \\
&\qquad\quad +  \int_0^t \int_D \Gamma(t,x;s,y) h(s,y)\, ds dy,
\end{align}
 and let $\tilde W$ be the set function as defined in \eqref{ch2'-s3.2a}. Then the joint law of $(u, \tilde W)$ under $\tilde P$ is equal to the joint law under $P$ of $(v, W)$, where $v$ is the random field defined for $(t,x)\in[0,T]\times D$ by
\beq
\label{ch2'-s3.4.1.4}
v(t,x) = I_0(t,x) + \int_0^t \int_D \Gamma(t,x;s,y)\, W(ds,dy).
\eeq
Under $P$, the laws of $u$ and $v$ are mutually equivalent. Further, $u$ and $v$ have jointly measurable versions that are adapted to $(\tf_t)$.
\end{thm}

The notion of law of a stochastic process is recalled in Section  \ref{rdsecA.4}.
\medskip

\noindent{\em Proof of Theorem \ref{ch2'-s3.4.1-t0}.}
The assumptions on $\Gamma$ and $h$ imply that the integrals in \eqref{ch2'-s3.4.1.30} are well-defined.
The random fields $u$ and $v$ have jointly measurable and adapted versions by the assumptions on $\Gamma$ and Propositions \ref{rdprop2.6.2} and \ref{rdprop2.6.4}.

By Theorem \ref{ch2'-s3-t1}, the set function $\tilde W$ given in \eqref{ch2'-s3.2a} is a space-time white noise under $\tilde P$, and since $u$ solves \eqref{ch2'-s3.4.1.30}, it follows from \eqref{ch2'-s3.integral} that
\beqn
u(t,x) = I_0(t,x) + \int_0^t \int_D \Gamma(t,x;s,y)\, \tilde W(ds,dy).
\eeqn
Therefore, under $\tilde P$, $(u,\tilde W)$ is a Gaussian random field whose law is determined by its mean and covariance functions, which are identical to those of $(v,W)$ under $P$.

Since the probability measures $\tilde P$ and $P$ are mutually equivalent, the laws of $u$ under $\tilde P$ and under $P$ are mutually equivalent.
Therefore under $P$, the laws of $u$ and $v$ are mutually equivalent.
\qed

\subsection{Existence of weak solutions in law}
\label{ch2'-ss1.1}
The next theorem provides sufficient conditions for the existence of a weak solution in law to \eqref{ch2'-s3.4.1.1} in the particular case where $\sigma\equiv 1$.

\begin{thm}
\label{ch2'-s3.4.1-t1}
Assume the following three conditions:
\begin{enumerate}
\item The function $(t,x) \mapsto I_0(t,x)$ satisfies condition $({\bf H_I})$ of Section \ref{ch4-section1}.
\item The function $\Gamma$ satisfies the conditions $({\bf H_\Gamma})(i), (ii)$ and $(iiia)$ of Section \ref{ch4-section0}.
\item The function $b: [0,T]\times D\times \re \to \re$ is $\mathcal{B}_{[0,T]}\times \mathcal{B}_{D}\times \mathcal{B}_{\re}$-measurable. Furthermore, there exists a function $b_0\in L^2(D)$ such that for all $(t,x,z)\in [0,T]\times D\times \IR$,
\beq
\label{ch2'-s3.4.1.2}
|b(t,x,z)| \le b_0(x)(1+|z|).
\eeq
\end{enumerate}
Then the equation
\begin{align}
\label{ch2'-s3.4.1.3}
u(t,x) &= I_0(t,x) + \int_0^t \int_D \Gamma(t,x;s,y)\, W(ds,dy)\notag\\
 &\qquad\quad+ \int_0^t \int_D \Gamma(t,x;s,y) b(s,y,u(s,y))\, ds dy,
\end{align}
$(t,x)\in[0,T]\times D$, admits a weak solution in law.
\end{thm}

\begin{remark}
\label{ch5-r1}
 (a)\ Examples of functions $\Gamma (t,x;s,y)$ satisfying condition 2.~are given in
Section \ref{ch1'-s6}. These include: (i)  the fundamental solution of the heat equation, of a fractional heat equation and of the wave equation on $\re$; (ii) the Green's function of the heat equation on an interval $[0,L]$ (both with vanishing Dirichlet and Neumann boundary conditions), and (iii)  the wave equation on $\re_+$ and on a bounded interval $[0,L]$ with vanishing Dirichlet boundary conditions.

(b)\ Comparing with the assumptions stated in Section \ref{ch4-section1}, we observe that for the function $b$, we assume neither Lipschitz continuity nor any other kind of regularity. When the domain $D$ is bounded, condition \eqref{ch2'-s3.4.1.2} is weaker than $({\bf H_L})$ (vi). However, when $D$ is unbounded, these two conditions are unrelated.
\end{remark}

\noindent{\em Proof of Theorem \ref{ch2'-s3.4.1-t1}.}
Let $W$ be a space-time white noise defined on some stochastic basis $(\Omega,{\tf}, P, ({\tf}_t,\, t \in [0, T]))$. Let
$v=(v(t,x),\, (t,x)\in[0,T]\times D)$ be
the jointly measurable and adapted version of the process defined in  \eqref{ch2'-s3.4.1.4}.
By Assumptions 1.~and 2.~above, we have
\begin{align}
\label{defC}
&\sup_{(t,x)\in[0,T]\times D}E\left[v^2(t,x)\right] \notag\\
&\qquad \le 2 \sup_{(t,x)\in[0,T]\times D}\left(I_0^2(t,x) + \int_0^t ds \int_D dy\,  \Gamma^2(t,x;s,y)\right)\notag\\
&\qquad \le 2 \sup_{(t,x)\in[0,T]\times D}\left(I_0^2(t,x)+ \int_0^t ds \sup_{x\in D} \int_D dy\,  H^2(s,x,y) \right)\notag\\
 &\qquad \le c,
 \end{align}
for some constant $c < \infty$.

We now check that for any finite constant $C \geq  0$, there is $\ep > 0$ such that
\beq
\label{5.1.3(*1)}
    \sup_{s \in [0,T]} E\left[\exp\left(C \int_s^{s+\ep} dt \int_D dx\, b^2(t,x,v(t,x))\right)\right] < \infty.   
    \eeq
Taking $C=\half$, this will imply the validity of condition (c) of Proposition \ref{ch2'-s3-p1}.

We will in fact prove that for any finite constant $C \geq  0$,  there is $\ep > 0$ such that
\beq
\label{5.1.3(*2)}
   \sup_{s \in [0,T]} E\left[\exp\left(C \int_s^{s+\ep} dt \int_D dx\, b_0^2(x) (1+v^2(t,x))\right)\right] < \infty.   
   \eeq
By \eqref{ch2'-s3.4.1.2}, this implies \eqref{5.1.3(*1)}. We note that given $C > 0$ and $\ep > 0$ satisfying \eqref{5.1.3(*2)}, this property remains valid if we replace $C$ by any nonnegative $C' < C$ and $\ep$ by any nonnegative $\ep' < \ep$.

To prove \eqref{5.1.3(*2)}, we observe that for all $t\in[0,T]$,
\begin{align*}
&E\left[\exp\left(C\int_s^{s+\ep} dt \int_D dx\, b_0^2(x) \left(1+v^2(t,x)\right)\right)\right] \\
&\qquad\qquad =  \exp\left(C\int_s^{s+\ep} dt \int_D dx\,  b_0^2(x) \right)\\
&\qquad\qquad\qquad \times E\left[\exp\left(C\int_s^{s+\ep} dt  \int_D dx\, v^2(t,x)b_0^2(x)\right)\right].
\end{align*}
The first factor on the right-hand side of the last inequality is bounded over $s \in [0, T]$ because $b_0\in L^2(D)$. As for the second factor, we apply Jensen's inequality to the convex function $x\mapsto \exp(x)$ and the finite measure
$dt \mu(dx)$ on $[s,s+\ep]\times D$, where
$\mu(dx) = b_0^2(x) dx$,  to obtain the upper bound
\beq
\label{jensen-trick}
\frac{1}{\ep K_0} \int_s^{s+\ep} dt \int_D \mu(dx)\, E\left(\exp\left(C\varepsilon K_0\, v^2(t,x)\right)\right),
\eeq
where $K_0 =\Vert b_0\Vert_{L^2((D)}^2$.

Let $m_{t,x}$ and $\sigma_{t,x}^2$ be respectively the mean and variance of $v(t,x)$, and let $Z$ be a ${\rm N}(0,1)$ random variable. Then the expectation in \eqref{jensen-trick} can be written
\beqn
    E\left[\exp\left(C\ep K_0 (m_{t,x} + \sigma_{t,x}  Z)^2\right)\right]
    \leq E\left[\exp\left(2 C\ep K_0 (m_{t,x}^2 + \sigma_{t,x}^2  Z^2)\right)\right],
    \eeqn
and this is clearly an increasing function of $m_{t,x}^2$ and  $\sigma_{t,x}^2$. Since these are uniformly bounded by \eqref{defC}, we can replace them by their uniform bound, say $C_0$. For $\ep$ small enough, more precisely, for $\ep $ satisfying $2 C\ep K_0 C_0 < 1$, the expectation is a finite number $C_\ep$.  Since $C_\ep$ is also an upper bound for the expectation in  \eqref{jensen-trick}, we see that \eqref{5.1.3(*2)} holds. Therefore, \eqref{5.1.3(*1)} is proved.

   We note that \eqref{5.1.3(*1)} implies that a.s., $(t,x) \mapsto h(t,x) := b(t,x,v(t,x))$ belongs to $L^2([0,T] \times D)$ (and even $h \in L^2([0,T] \times D \times \Omega)$), since $z \leq \exp(z)$ for all $z \in \R$, so \eqref{ch2'-s3.2} is satisfied.

Define a measure $\tilde P$ on $(\Omega,\tf_T)$ by
\begin{align*}
\frac{d\tilde P}{dP} = &\exp\left(\int_0^T \int_D b(t,x,v(t,x))\, W(dt,dx)\right. \\
&\left. \qquad\qquad- \frac{1}{2}\int_0^T dt \int_D dx\, b^2(t,x,v(t,x))\right).
\end{align*}
Appealing to Proposition \ref{ch2'-s3-p1} (c), we deduce that $E\left[\frac{d\tilde P}{dP}\right]=1$ and $\tilde P$ is a probability measure on $(\Omega,\tf_T)$. In addition, by Theorem \ref{ch2'-s3-t1}, the set function
\beqn
\tilde W(A) = W(A) - \int_0^T dt \int_D dx\, 1_A(t,x) b(t,x,v(t,x))
 \eeqn
 is a space-time white noise under $\tilde P$. Moreover, by \eqref{ch2'-s3.4.1.4},
\begin{align*}
v(t,x) &= I_0(t,x) + \int_0^t ds \int_D dy\,\Gamma(t,x;s,y) b(s,y,v(s,y))\\
&\qquad\quad + \int_0^t ds \int_Ddy\, \Gamma(t,x;s,y)\left[W(ds,dy) - b(s,y,v(s,y))\right]\\
&= I_0(t,x) + \int_0^t ds \int_D dy\,\Gamma(t,x;s,y)b(s,y,v(s,y))\\
&\qquad\quad + \int_0^t \int_D \Gamma(t,x;s,y)\, \tilde W(ds,dy),
\end{align*}
where in the last equality, we have used \eqref{ch2'-s3.integral}. According to Definition \ref{ch2'-s3.4.1-d1}, $(\Theta, \tilde W, v)$ is a weak solution in law of \eqref{ch2'-s3.4.1.3}, where $\Theta = (\Omega,{\tf},\tilde P, ({\tf}_t,\, t\in [0, T]))$.
\qed

\subsection{Uniqueness in law}
\label{ch2'-s3.4.2}
In the previous section, we addressed the issue of existence of weak solutions in law to the SPDE \eqref{ch2'-s3.4.1.30}.  Here,  we examine the issue of uniqueness of such solutions.
\medskip

\noindent{\em Preliminaries}
\smallskip

The next proposition will be used to verify condition (c) of Proposition \ref{ch2'-s3-p1}.

\begin{prop}
\label{ch2'-s3.4.2-p1}
Under the hypotheses of Theorem \ref{ch2'-s3.4.1-t1}, any weak solution in law $u$ of \eqref{ch2'-s3.4.1.3} has the following property:

For some $\ep > 0$,
\beq
\label{5.1.5(*3)}
    \sup_{s \in [0,T]} E\left[\exp\left(\half \int_s^{s+\ep} dt\ \int_D dx\, b^2(t,x,u(t,x)) \right)\right] < \infty.  
    \eeq
\end{prop}

\begin{proof}
Let $(\Theta, W, u)$ be as in Definition \ref{ch2'-s3.4.1-d1}, with \eqref{ch2'-s3.4.1.1} replaced by \eqref{ch2'-s3.4.1.3}. Let $v$ be the jointly measurable and adapted version of the random field $v$ defined in \eqref{ch2'-s3.4.1.4} and set
 \beqn
 \bar u(t,x) = u(t,x) - v(t,x), \qquad (t,x)\in[0,T]\times D.
 \eeqn
  By \eqref{ch2'-s3.4.1.3}, for fixed $(t,x)$,
\begin{align}
\label{ch2'-s3.4.2.3}
\bar u(t,x)& = \int_0^t ds \int_D dy \, \Gamma(t,x;s,y) b\left(s,y,u(s,y)\right)\notag\\
& = \int_0^t ds \int_D dy \, \Gamma(t,x;s,y) b\left(s,y,\bar u(s,y)+v(s,y)\right).
\end{align}
By \eqref{ch2'-s3.4.1.2}, 
\beq
\label{5.1.5(*4)}
    \vert\bar u(t,x)\vert \leq \int_0^t ds \int_D dy\,  |\Gamma(t,x;s,y)| \, b_0(y) (1 +|\bar u(s,y)|+ |v(s,y)|).     
\eeq
Applying the Cauchy-Schwarz inequality to the right-hand side yields
\begin{align*}
 \bar u^2(t,x) 
 &\le C \left(\int_0^t ds \int_D dy\ H^2(t-s,x,y)\right)\\
 &\qquad \times \left(\int_0^t ds \int_D dy\ b_0^2(y)\left[1+\bar u^2(s,y) + v^2(s,y)\right]\right)\\
 &\le C_T \left(1+\int_0^t ds \int_D dy\ b_0^2(y)\, \bar u^2(s,y) + \int_0^T ds \int_D dy\ b_0^2(y)\, v^2(s,y)\right),
 \end{align*}
 where we have used $({\bf H_\Gamma})$(ii) and replaced $t$ by $T$ in the last integral. Notice that the last right-hand side does not depend on $x$. Multiply the left- and right-hand sides by $b_0^2(x)$, integrate over $D$ and use the fact that $b_0 \in L^2(D)$ to see that
 \begin{align}
 \label{5.1.5(*5)}
 &\int_D dx\ b_0^2(x)\, \bar u^2(t,x)\notag\\
 &\quad \le C \left(1+\int_0^t ds \int_D dy\ b_0^2(y)\, \bar u^2(s,y) + \int_0^T ds \int_D dy\ b_0^2(y)\, v^2(s,y)\right).
 \end{align}
 Apply Gronwall's Lemma (Lemma \ref{lemC.1.1}) to the function
 \beqn
 g(t) := \int_D dx\ b_0^2(x)\, \bar u^2(t,x)
 \eeqn
 to obtain
\beq
\label{ch2'-s3.4.2.4}
\int_D dx\ b_0^2(x)\, \bar u^2(t,x) \le C_1\left(1 + \int_0^T ds \int_D dy\  b_0^2(y)\, v^2(s,y)\right) \exp\left(C_1T\right),
\eeq
where, for simplicity, we have removed the dependence of the constants on the specific parameters.

   Fix $C > 0$ and $\ep > 0$ such that \eqref{5.1.3(*2)} holds. By taking $\ep$ slightly smaller, we can assume that $\ep = T/n$ for some integer $n$. Let $0 = t_0 < t_1 < \cdots < t_n = T$ be a partition with $t_k - t_{k-1} = T/n = \ep$ for all
   $k$. We are going to show that for any finite constant $C' > 0$, there is $\ep' > 0$ such that
   \beq
   \label{5.1.5(*6)}
   \sup_{r \in [0,T]}   E\left[\exp\left(C' \int_r^{r+\ep'} dt \int_D dx\ b_0^2(x)\, \bar u^2(t,x)\right)\right] < \infty.
   \eeq   
For this, we integrate both sides of \eqref{ch2'-s3.4.2.4} with respect to $t$, from $r$ to $r + \ep'$, to obtain
    \beqn
   \int_r^{r+\ep'} dt \int_D dx\ b_0^2(x)\, \bar u^2(t,x) \leq  \bar C_1 \ep' + \bar C_1 \ep' \int_0^T ds \int_D dy\ b_0^2(y) v^2(s,y),
   \eeqn
   where $\bar C_1= C_1 \exp\left(C_1T\right)$. The first term on the right-hand side plays no role, and we write the second term as
\beqn
   \bar C_1 \ep' \sum_{k=1}^n \int_{t_{k-1}}^{t_k} ds \int_D dy\ b_0^2(y) v^2(s,y),
   \eeqn
which no longer depends on $r$. We now replace the integral in \eqref{5.1.5(*6)} by this expression, and bound the left-hand side of \eqref{5.1.5(*6)} by
\begin{align*}
    &C_2\, E\left[\exp\left(C' \bar C_1 \ep' \sum_{k=1}^n \int_{t_{k-1}}^{t_k} ds \int_D dy\ b_0^2(y) v^2(s,y)\right)\right] \\
   &\qquad\qquad =   C_2\, E\left[\prod_{k=1}^n \exp\left(C' \bar C_1 \ep' \int_{t_{k-1}}^{t_k} ds \int_D dy\ b_0^2(y) v^2(s,y)\right)\right].
   \end{align*}
We apply Holder's inequality with $n$ exponents equal to $n$, to bound this above by
\begin{align*}
       &C_2  \prod_{k=1}^n \left(E\left[\exp(C' \bar C_1 \ep' n \int_{t_{k-1}}^{t_k} ds \int_D dy\  b_0^2(y) v^2(s,y)\right]\right)^{1/n}\\
   &\qquad\qquad \leq C_2\sup_{r \in [0,T]} E\left[\exp\left(C' \bar C_1 \ep' n \int_{r}^{r + \ep} ds \int_D dy\  b_0^2(y) v^2(s,y)\right)\right].
   \end{align*}
Choose $\ep' > 0$ small enough so that $C' \bar C_1 \ep' n < C$. By \eqref{5.1.3(*2)}, this last right-hand side is finite, proving \eqref{5.1.5(*6)}.

   We now conclude the proof of \eqref{5.1.5(*3)}. Fix $C = C' > 0$. Take $\ep = \ep' > 0$ small enough so that both \eqref{5.1.3(*2)} and \eqref{5.1.5(*6)} hold with $4C$ instead of $C$.  We will show that
   \beq
   \label{5.1.5(*7)}
    \sup_{r \in [0,T]} E\left[\exp\left(C \int_r^{r+\ep} dt \int_D dx\ b_0^2(x) (1 + u^2(t,x))\right)\right] < \infty,    
\eeq
which will imply \eqref{5.1.5(*3)} because of \eqref{ch2'-s3.4.1.2} (by taking $C=\half)$.

In order to prove \eqref{5.1.5(*7)}, observe that since $u(t,x) = \bar u(t,x) + v(t,x)$, the left-hand side  of \eqref{5.1.5(*7)} is bounded above by
   \beq
   \label{5.1.5(*8)}
  \bar C \sup_{r \in [0,T]} E\left[\exp\left(2C \int_r^{r+\ep} dt \int_D dx\ b_0^2(x) (\bar u^2(t,x) + v^2(t,x))\right)\right].
  \eeq      
Applying the Cauchy-Schwarz inequality, this is in turn bounded above (up to a multiplicative constant) by
\begin{align*}
&\qquad \left(E\left[\exp\left(4C \int_r^{r+\ep} dt \int_D dx\  b_0^2(x) \bar u^2(t,x)\right)\right]\right)^\half\\
&\qquad \qquad \times \left(E\left[\exp\left( 4C \int_r^{r+\ep} dt \int_D dx\ b_0^2(x) v^2(t,x)\right)\right]\right)^\half.
\end{align*}
By our choice of $\ep$, we can use \eqref{5.1.5(*6)}  and \eqref{5.1.3(*2)} to conclude that \eqref{5.1.5(*8)} is finite. This establishes \eqref{5.1.5(*7)}.

The proof of Proposition \ref{ch2'-s3.4.2-p1} is complete.
\end{proof}

\begin{prop}
\label{ch2'-s3.4.2-t1}
Under the assumptions of Theorem \ref{ch2'-s3.4.1-t1}, let $u$ be a weak solution in law of \eqref{ch2'-s3.4.1.3} defined on $(\Theta, W)$ as in Definition \ref{ch2'-s3.4.1-d1}. Consider the measure $\tilde P$ defined by
\begin{align*}
\frac{d\tilde P}{dP} &= \exp\left(-\int_0^T \int_D b(t,x,u(t,x))\, W(dt,dx)\right.\\
 &\left.\qquad\qquad\qquad\qquad- \frac{1}{2} \int_0^T dt \int_D dx\, b^2(t,x,u(t,x))\right).
\end{align*}
Then $E\left[\frac{d\tilde P}{dP}\right]=1$, that is, $\tilde P$ is a probability measure on $(\Omega,\mathcal{F})$, and under $\tilde P$,
the set function $\tilde W$ defined by
\beqn
    \tilde W(A) = W(A) + \int_0^T ds \int_D dy\, 1_A(s, y) b(s, y, u(s, y)),\quad   A \in \cB_{[0, T] \times D}^f ,
    \eeqn
is a space-time white noise such that for all $(t, x) \in [0, T] \times D$,
\beq
\label{5.1.6(*0)}
     u(t,x) = I_0(t,x) + \int_0^t \int_D \Gamma(t,x; s, y)\, \tilde W(ds, dy)\quad  \tilde P-{\text{a.s}}.
     \eeq
 Moreover, under $P$, the laws of $u$ and $v$ are mutually equivalent, where v is a jointly measurable and adapted version of the process defined in \eqref{ch2'-s3.4.1.4}.
\end{prop}
\begin{proof}
By Proposition \ref{ch2'-s3.4.2-p1}, condition (c) of Proposition \ref{ch2'-s3-p1} is satisfied with $h(t,x)=b(t,x,u(t,x))$, and this implies  $E\left[\frac{d\tilde P}{dP}\right] =1$. By Girsanov's Theorem \ref{ch2'-s3-t1} and Proposition \ref{ch2'-s3-r1}, since $u$ is a weak solution in law of \eqref{ch2'-s3.4.1.3},
\begin{align*}
u(t,x) &= I_0(t,x) + \int_0^t \int_D \Gamma(t,x;s,y)\left[ W(ds,dy) + b(s,y,u(s,y))\, dsdy\right]\\
&= I_0(t,x) + \int_0^t \int_D \Gamma(t,x;s,y)\, \tilde W(ds,dy),
\end{align*}
where, under $\tilde P$, $\tilde W$ is a space-time white noise on $\Theta$. Therefore, under $\tilde P$, the random field $(u(t,x),\, (t,x)\in[0,T]\times D)$ is Gaussian and
the joint law of $(u, \tilde W)$ is also Gaussian and is determined by its mean and covariance functions, which are identical to those of $(v, W)$ under $P$.

The proof of the last statement follows in the same way as in the proof of Theorem \ref{ch2'-s3.4.1-t0}.
\end{proof}

\noindent{\em Proving the uniqueness in law}
\medskip

For the proof of uniqueness in law, we need the additional assumption (${\bf H_w}$) below. In the last part of this section, we will discuss conditions for its validity.
\medskip

\noindent  (${\bf H_w}$)\  For $i=1,2$, let $u_i$  be a weak solution in law to \eqref{ch2'-s3.4.1.3} corresponding to $(\Theta_i, W_i)$ and a probability measure $P_i$. Let
\begin{align}\nonumber
X_i &:= -\int_0^T \int_D b(t, x, u_i(t,x))\, W_i(dt, dx) \\
      &\qquad \qquad - \half\int_0^T dt \int_Ddx\, b^2(t,x,u_i(t,x)).
\label{Hw}
\end{align}
Define $\tilde P_i$ as in Proposition \ref{ch2'-s3.4.2-t1}, that is, 
\beq
\label{pisi}
\frac{d\tilde P_i}{dP_i} = \exp(X_i).
\eeq
 Then the joint law under $\tilde P_1$ of $(u_1, W_1, \exp(-X_1))$  is the same as the joint law under $\tilde P_2$ of $(u_2, W_2, \exp(-X_2))$.
\medskip

\begin{thm}
\label{ch2'-s3.4.2-t2}
Suppose that  the assumptions of Theorem \ref{ch2'-s3.4.1-t1}  are satisfied and that $({\bf H_w})$ holds. Then all weak solutions in law to \eqref{ch2'-s3.4.1.3} have the same probability distribution, that is, the joint law of $(u,W)$ is unique.
\end{thm}
\begin{proof}
 For $i=1,2$, let $\tilde P_i$ be defined as in \eqref{pisi}.
We have seen in Proposition \ref{ch2'-s3.4.2-t1} that $\tilde P_i$ is a probability measure,
and that the law under $\tilde P_i$ of $u_i$ is the same as the law under $P_i$ of $v_i$, where $v_i$ is defined in \eqref{ch2'-s3.4.1.4} with $W$ replaced by $W_i$. Since $v_1$ and $v_2$ are mean zero Gaussian processes with the same covariance function, they have the same law, and therefore the law under $\tilde P_1$ of $u_1$ is the same as the law under $\tilde P_2$ of $u_2$.

Fix $n\ge 1$, let $(t_j,x_j)\in[0,T]\times D$ and $A_j \in \mathcal{B}^f_{[0,T]\times D}$, $j=1,\ldots,n$. Then for any Borel set $B$ of $\re^{2n}$,
\begin{align}
\label{ch2'-s3.4.2.5}
&P_1\left\{\left(u_1(t_1,x_1),\ldots,u_1(t_n,x_n),W_1(A_1),\ldots, W_1(A_n)\right)\in B\right\}\notag\\
&\qquad\qquad = E_{\tilde P_1}\left[1_B(\left(u_1(t_1,x_1),\ldots,u_1(t_n,x_n),W_1(A_1),\ldots, W_1(A_n)\right))\frac{d P_1}{d\tilde P_1}\right].
\end{align}
Because $\frac{d P_i}{d \tilde P_i} = \exp(- X_i)$, it follows from  Hypothesis (${\bf H_w}$) that
$\left(u_1, W_1, \frac{dP_1}{d\tilde P_1}\right)$ under $\tilde P_1$ and $\left(u_2, W_2, \frac{dP_2}{d\tilde P_2}\right)$ under $\tilde P_2$ have the same (joint) law, therefore the right-hand side of
\eqref{ch2'-s3.4.2.5} is equal to
\begin{align*}
&E_{\tilde P_2}\left[1_B(\left(u_2(t_1,x_1),\ldots,u_2(t_n,x_n),W_2(A_1),\ldots W_2(A_n)\right)) \frac{d P_2}{d\tilde P_2}\right]\\
& \qquad\qquad = P_2\left\{\left(u_2(t_1,x_1),\ldots,u_2(t_n,x_n),W_2(A_1),\ldots W_2(A_n)\right)\in B\right\}.
\end{align*}
The proof of the proposition is complete.
\end{proof}
\begin{remark}
\label{5.1.7(*1)}
If $u_i$ has continuous sample paths a.s., one can view $u_i$ as a random variable with values in $\cC([0,T] \times D)$, equipped with the metric of uniform convergence on compact sets and its Borel $\sigma$-field $\cB_{\cC([0,T] \times D)}$. In this case, the term {\em law of $u_i$} refers to the probability measure $Q_{u_i}$ on $\cB_{\cC([0,T] \times D)}$ defined by $Q_{u_i}(A) := P\{u_i \in A\}$, $A \in \cB_{\cC([0,T] \times D)}$. Because $\cC([0,T] \times D)$ (with the metric of uniform convergence on compact sets) is a complete separable metric space, the probability measure $Q_{u_i}$ is determined by the finite-dimensional distributions of $u_i$. These facts are recalled in Section \ref{rdsecA.4}.
\end{remark}
\smallskip

\noindent{\em A sufficient condition for hypothesis $({\bf H_w})$}
\medskip

We continue this section by providing a sufficient condition for hypothesis (${\bf H_w}$) to hold. The main point is to identify the random variables $W_i(A_j)$ and $\frac{dP_i}{d\tilde P_i}$ as concrete measurable functions of $u_i$ that do not depend on $i \in \{1,2\}$.

\begin{prop}
\label{ch2'-s3.4.2-p2}
Suppose that the hypotheses of Theorem \ref{ch2'-s3.4.1-t1} are satisfied, that $\Gamma$ satisfies \eqref{ch1'-s5.12} with ${\bf \Delta}_1$ of the same anisotropic form as ${\bf \Delta}$ in \eqref{norm-aniso}, and there is a partial differential operator $\cL^\star$ such that, for all weak solutions in law to \eqref{ch2'-s3.4.1.3} and for all $\psi \in \cC_0^\infty(]0, T[ \times D)$,
\begin{align}
\label{ch2'-s3.4.2.6}
   & \int_0^T dt \int _D dx\, \cL^\star \psi(t,x) (u(t,x) - I_0(t,x))\notag\\
   &\qquad\qquad  = \int_0^T dt \int_D dx\, b(t, x, u(t, x)) \psi(t, x) + \int_0^T \int_D \psi(s,y)\, W(ds, dy),    
    \end{align}
    $P$-a.s.
Then Hypothesis $({\bf H_w})$ holds.
\end{prop}

\begin{remark}
\label{after-ch2'-s3.4.2-p2}
The partial differential operator $\cL^\star$ will usually be the adjoint of the operator $\cL$ (the partial differential operator of which $\Gamma$ is the fundamental solution).
\end{remark}


The proof of Proposition \ref{ch2'-s3.4.2-p2} relies on some technical results which are stated and proved after Proposition \ref{5.1.3-ex-wave*2}. We assume them for the moment and proceed with the proof of the proposition.
\bigskip

\noindent{\em Proof of Proposition \ref{ch2'-s3.4.2-p2}}.
For $i = 1,2$, let $u_i$  be a weak solution in law to \eqref{ch2'-s3.4.1.3}
corresponding to $(\Theta_i, W_i)$ and the probability measure $P_i$. Let $X_i$ be as in \eqref{Hw}.

Define $\tilde P_i$ by
   $\frac{d\tilde P_i}{dP_i} = \exp\left(X_i\right) $.
By Proposition \ref{ch2'-s3.4.2-t1}, $\tilde P_i$ is a probability measure and the set function $\tilde W_i$ defined by
\beqn
    \tilde W_i(A) = W_i(A) + \int_0^T ds \int_D dy\, 1_A(s, y) b(s, y, u_i(s, y)),\quad   A \in \B_{[0, T] \times D}^f,
    \eeqn
is, under $\tilde P_i $, a space-time white noise such that
\beqn
    u_i(t,x) = I_0(t,x) + \int_0^T \int_D  \Gamma(t, x; s, y)\, \tilde W_i(ds, dy),\qquad   \tilde P_i{\text{-a.s.}}
    \eeqn
In the definition of $X_i$ in \eqref{Hw}, by Corollary \ref{Lemma *2} below, we replace $u_i$ by a version $\hat u_i$ such that $\hat u_i - I_0$ has continuous sample paths, and this does not change the value of $X_i$. Observe that $\frac{d P_i}{d \tilde P_i} = \exp(- X_i)$ and
\begin{align*}
   - X_i &= \int_0^T \int_D b(t, x, \hat u_i(t,x))\,  W(dt, dx) + \half \int_0^T dt \int_D dx\, b^2(t, x,\hat u_i(t, x))\\
   & = Z_{i,1} - Z_{i,2},\quad  \tilde P_i-{\text{a.s.}},
   \end{align*}
where
\beqn
    Z_{i,1} = \int_0^T \int_D b(t, x, \hat u_i(t,x))\,  \tilde W(dt, dx),\quad
    Z_{i,2} = \half \int_0^T dt \int_D dx\, b^2(t, x, \hat u_i(t, x))
    \eeqn
and we have used Proposition \ref{ch2'-s3-r1}

Let $\theta = (\theta(t,x),\, (t,x) \in [0, T] \times D)$ be the coordinate process on $\cC$. For $w \in \cC$ and $\hat W$ as in Lemma \ref{Facts C1-C2-D} (ii) below, define
   \beqn
     Y_1 := \int_0^T \int_D b(t,x, I_0(t,x) + \theta(t,x)(w))\, \hat W(dt, dx)
     \eeqn
and
\beqn
    Y_2 := \half \int_0^T dt \int_D dx\, b^2(t, x, I_0(t,x) + \theta(t,x)(w)).
    \eeqn

   Let $n \geq 1$ and $A_\ell \in \cB_{[0, T] \times D}^f$, $\ell = 1, \dots, n$. Using Lemmas \ref{Facts C1-C2-D} (iii), \ref{Facts E1-E2} (2) and \ref{Fact A0}, we see that there is a Borel random variable $\Phi : \cC \to \re^n \times \re$ such that
   \beqn
    \Phi = \left(\hat W(A_1) - \Phi_{A_1}, \dots, \hat W(A_n) - \Phi_{A_n}, \exp(Y_1 - Y_2)\right),\qquad  Q_v{\text{-a.s.}},
    \eeqn
and
\beqn
    \Phi(\hat u_i - I_0) = \left(W_i(A_1),\dots, W_i(A_n), \exp(Z_{i,1}-Z_{i,2})\right),\qquad    \tilde P_i{\text{-a.s.}}
    \eeqn
In particular, the joint law under $\tilde P_i$ of $(\hat u_i - I_0, W_i, \exp(-X_i))$ does not depend on $i \in \{ 1,2\}$,
 where the law of $v_i := \hat u_i - I_0$ refers to its law on $\cC$. Since $u_i$ is a version of $v_i + I_0$, this implies that the law of $(u_i, W_i, \exp(- X_i)$) (in the sense of finite-dimensional distributions) does not depend on $i = 1,2$. This establishes condition $({\bf H_w})$ and completes the proof of the proposition.
\qed
\medskip

As a consequence of Proposition \ref{ch2'-s3.4.2-t1}, we have the following regularity property of the solution to \eqref{ch2'-s3.4.1.3}.

\begin{cor}
\label{Lemma *2} Under the hypotheses of Proposition \ref{ch2'-s3.4.2-p2}, let $(\Theta, W, u)$ be a weak solution in law to \eqref{ch2'-s3.4.1.3}. Then $u$ has a version, again denoted $u$, such that $u - I_0$ has continuous sample paths on $[0, T] \times D$.
\end{cor}

\begin{proof} Let $\tilde P$ be the probability measure defined in Proposition \ref{ch2'-s3.4.2-t1} and let $\tilde W$ be the set function defined there. By Proposition \ref{ch2'-s3.4.2-t1}, under $\tilde P$, $u$ has the same law as $v$ in \eqref{ch2'-s3.4.1.4} under $P$, and $u - I_0$ can be written as in \eqref{5.1.6(*0)}. Because of the assumption that $\Gamma$ satisfies \eqref{ch1'-s5.12}, we can use Lemma \ref{ch1'-ss5.2-l1} (a) and Theorem \ref{app1-3-t1} to see that for each bounded sub-domain $\tilde D$ of $D$, $u - I_0$ has a version $\hat u$ with $\tilde P$-a.s.-continuous sample paths on $[0, T] \times \tilde D$. Writing $D$ as a countable union of bounded sub-domains and because $\tilde P$ and $P$ are mutually equivalent, $\hat u$ also has  $P$-a.s.-continuous sample paths on $[0, T] \times D$. The desired version of $u$ is $I_0 + \hat u$.

We note in passing that if the initial condition $u_0$ that determines $I_0(t, x)$ in \eqref{ch2'-s3.4.1.3} is not continuous, then the sample paths of $u$ alone will not be continuous.
\end{proof}

\medskip

We now give some examples of SPDEs that correspond to partial differential operators $\cL$ that satisfy condition \eqref{ch2'-s3.4.2.6}
of Proposition \ref{ch2'-s3.4.2-p2} or, equivalently by Lemma \ref{Fact C} (a), condition \eqref{5.1.19-prime} below.
\medskip

\noindent{\em Example 1: Stochastic heat equation}
\medskip

Consider the setting of Section \ref{ch1'-ss6.1}, where $\cL=\frac{\partial}{\partial t}- \frac{\partial^2}{\partial x^2}$ is the heat operator and $\cL^\star = - \frac{\partial}{\partial t} - \frac{\partial^2}{\partial x^2}$ is its adjoint.
The domain $D$ is either $\re$ or a bounded interval $[0,L]$. In the latter case, we consider either homogeneous Dirichlet or Neumann boundary conditions, and we denote by $\Gamma$ the fundamental solution (or the Green's function) corresponding to $\cL$,
given by \eqref{ch1'-s6.1}, \eqref{ch4(*1)} or \eqref{ch4(*2)}. The function $I_0$ is given by \eqref{ch1'-v001}, \eqref{cor1.0} or \eqref{cor1.0-N}, for some Borel function $u_0$, that we assume to be bounded so that assumption $({\bf H_I})$ holds. With these choices, we consider the weak solution in law $u$ to \eqref{ch2'-s3.4.1.3}, in which the function $b$ is assumed to satisfy condition 3. of Theorem \ref{ch2'-s3.4.1-t1}.

\begin{prop}
\label{5.1.3-ex-heat*1}
For the three considered forms of the stochastic heat equation, the assumptions of Proposition \ref{ch2'-s3.4.2-p2} hold, as well as the (uniqueness) conclusion of Theorem \ref{ch2'-s3.4.2-t2}.
\end{prop}
\begin{proof}
We have seen in Section \ref{ch1'-ss6.1} that $\Gamma$ satisfies assumption $({\bf H_\Gamma})$, and in \eqref{rde4.2.5} that \eqref{rde4.2.5-wave} is satisfied with ${\bf \Delta_1}(t,x; s, y) = \vert t - x \vert^{\frac{1}{4}} + \vert x - y\vert^\half$. We now check \eqref{ch2'-s3.4.2.6}, or rather, the equivalent condition \eqref{5.1.19-prime}.

Let $\psi\in \cC_0^\infty(]0,T[\times D)$. On the left-hand side of \eqref{5.1.19-prime}, we replace $u(t, x)$ by the expression on the right-hand side of \eqref{5.1.6(*0)},
and then we use the stochastic Fubini's theorem Theorem \ref{ch1'-tfubini}, with $X=[0,T]\times D$, $\mu$ equal to Lebesgue measure and $G$ there defined by
 \beqn
 G((r,x),s,y):= 1_{[s,T]}(r)  \cL^\star \psi(r, x) \Gamma(r,x;s,y),
 \eeqn 
 to obtain
    \begin{align}
  \label{Example(*H1)}
   &\int_0^T dr \int_D dx\, \cL^\star \psi(r, x) \int_0^r \int_D  \Gamma(r,x;s,y) \tilde W(ds, dy)\notag\\
    &\qquad = \int_0^T  \int_D \tilde W(ds, dy) \int_s^T dr \int_D dx\, \cL^\star \psi(r, x) \Gamma(r,x;s,y),\quad \tilde P{\text{-a.s.}}
    \end{align}

 We will argue that
  \begin{align}
  \label{truc}
  & \int_s^T dr \int_D dx\, \cL^\star \psi(r, x)\, \Gamma(r,x;s,y) 
      = \psi(s,y).
 \end{align}
 Plugging back this equality into \eqref{Example(*H1)}, we obtain the right-hand side of \eqref{5.1.19-prime}, and this proves \eqref{ch2'-s3.4.2.6}.

 Equality \eqref{truc} is a generic equality, which however requires verification due to issues of integrability and boundary conditions.
Let $\varphi \in \cC_0^\infty(]0, T[ \times D)$. Then
\beq\label{rd03_25e1}
    \int_0^T dt \int_D dx \int_0^t ds \int_D dy\, \vert \cL^\star\psi(t,x)  \Gamma(t,x; s,y) \varphi(s,y) \vert < \infty,
\eeq
so we can apply Fubini's theorem to see that
\begin{align}
\label{coronilla(*1)rd}
 & \int_0^T ds \int_D dy \left(\int_s^T dt \int_D dx\, \cL^\star \psi(t,x) \Gamma(t,x; s,y)\right) \varphi(s,y)\notag\\
 &\qquad =  \int_0^\infty dt \int_D dx\, \cL^\star\psi(t,x) \left(\int_0^t ds \int_D dy\, \Gamma(t,x; s,y) \varphi(s,y)\right).
  \end{align}   
Because $\psi$ has compact support, and vanishes near $\partial D$ when $D = ]0, L[$, and also near $0$ and $T$, we can integrate by parts (twice in $x$, and once in $t$), with no contribution of boundary terms, to see that \eqref{coronilla(*1)rd} is equal to
\begin{align}\nonumber
    &\int_0^\infty dt \int_D dx\,  \psi(t,x)\, \cL \left(\int_0^t ds \int_D dy\, \Gamma(t,x; s,y) \varphi(s,y)\right)\\
    &\qquad =  \int_0^\infty dt \int_D dx\, \psi(t,x) \varphi(t,x),
\label{rd03_25e2}
    \end{align}
by definition of $\Gamma$. This proves \eqref{truc} and completes the proof of the proposition.
\end{proof}

 \medskip


\noindent{\em Example 2: Stochastic wave equation}
\medskip

We consider the setting of Section \ref{ch1'-ss6.2}, where $\cL=\frac{\partial^2}{\partial t^2}- \frac{\partial^2}{\partial x^2}$ is the wave operator and
\beqn
   \cL^\star = \frac{\partial^2}{\partial t^2} - \frac{\partial^2}{\partial x^2}=\cL
\eeqn
 is its adjoint.
The domain $D$ is either $\re$, $]0,\infty[$ or $]0,L[$. In the last two cases, we consider  homogeneous Dirichlet boundary conditions, and we denote by $\Gamma$ the fundamental solution (or the Green's function) corresponding to $\cL$
given by \eqref{wfs}, \eqref{wavehalfline} or \eqref{wave-bc1-100'}. The function $I_0$ is given by \eqref{p93.1}, \eqref{p98.2}, or \eqref{i0-tris}, for some initial conditions $f$ and $g$. We assume that $f$ is bounded and continuous and that $g \in L^1(D)$, so that assumption $({\bf H_I})$ holds. With these choices, we consider the weak solution in law $u$ to \eqref{ch2'-s3.4.1.3}, in which the function $b$ is assumed to satisfy condition 3.~of Theorem \ref{ch2'-s3.4.1-t1}.

\begin{prop}
\label{5.1.3-ex-wave*2}
For the three considered forms of the stochastic wave equation, the assumptions of Proposition \ref{ch2'-s3.4.2-p2} hold, as well as the (uniqueness) conclusion of Theorem \ref{ch2'-s3.4.2-t2}.
\end{prop}
\begin{proof}
We have seen in Section \ref{ch1'-ss6.2} that $\Gamma$ satisfies assumption $({\bf H_\Gamma})$, and in \eqref{rde4.2.5-wave} that \eqref{ch1'-s5.12} is satisfied with ${\bf \Delta_1}(t,x; s, y) = \vert t - x \vert^\half + \vert x - y\vert^\half$. We now check \eqref{ch2'-s3.4.2.6}, or rather, the equivalent condition \ref{5.1.19-prime}.

As in the proof of Proposition \ref{5.1.3-ex-heat*1},
using \eqref{5.1.6(*0)} and then the stochastic Fubini's Theorem \ref{ch1'-tfubini}, for all  $\psi \in \cC_0^\infty(]0, \infty[ \times D)$, the left-hand side of \eqref{5.1.19-prime} is equal to
\begin{align}
 \label{Example(*H2)}
   &\int_0^T dt \int_D dx\, \cL^\star \psi(t, x) \int_0^r \int_D   \Gamma(t,x;s,y)\, \tilde W(ds, dy)\notag\\
     &\qquad = \int_0^T  \int_D \tilde W(ds, dy) \int_s^T dt \int_D dx\, \cL^\star \psi(t, x) \Gamma(t,x;s,y).     
     \end{align}

 We now check that
 \beq
\label{intermedi}
\int_s^T dt \int_D dx\, \cL^\star \psi(t, x) \Gamma(t,x;s,y) = \psi(s,y),
\eeq
which is the same property as in \eqref{truc}.  Notice that \eqref{rd03_25e1} also holds here and therefore, we can apply Fubini's theorem as in \eqref{coronilla(*1)rd} and then integration by parts (twice in $x$ and in $t$) as in \eqref{rd03_25e2}. This establishes \eqref{intermedi}.



   From \eqref{Example(*H2)} and \eqref{intermedi}, we conclude that \eqref{5.1.19-prime} holds, hence also \eqref{ch2'-s3.4.2.6}.
\end{proof}

\begin{remark}
We can deduce from these two examples that \eqref{ch2'-s3.4.2.6} will hold whenever $\cL$ has an adjoint $\cL^*$ such that \eqref{rd03_25e1} holds and we can use Fubini's theorem as in \eqref{coronilla(*1)rd} and then integrate by parts as in \eqref{rd03_25e2}. 
\end{remark}
\medskip

\noindent{\em Deterministic and stochastic integrals as Borel functions on path space}
\medskip

The next four lemmas were used in the proof of Proposition \ref{ch2'-s3.4.2-p2} and are proved under the hypotheses of that proposition. They show how to express relevant deterministic and stochastic integrals as Borel functions of their arguments.
Let $(\Theta, W, u)$ be a weak solution in law to \eqref{ch2'-s3.4.1.3} and let $\tilde P$ be the probability measure defined in Proposition \ref{ch2'-s3.4.2-t1}.

To simplify the notation, we set $\cC = \cC([0, T] \times D)$, and we let $\theta$ be the identity function on $\cC$ (that is, $\theta(w) = w$), and we let $(\theta(t,x), (t,x) \in [0, T] \times D)$ be the coordinate process $\theta(t,x)(w) = w(t,x)$.

On some probability space $(\Omega^\prime, \tf^\prime, P^\prime)$, let $W^\prime$ be a space-time white noise and let $v = (v(t,x))$ be the random field  such that for all $(t,x) \in [0, T] \times D$,
\beq
\label{Prop5.1.7(*A1)}
      v(t,x) = \int_0^t \int_D \Gamma(t,x; s,y)\, W^\prime(ds, dy).
      \eeq
      Under our assumptions, $v$ has a version with $P^\prime$-a.s.~continuous sample paths, again denoted $v$. Let $Q_v$ be the law on $\cC$ of $\omega \mapsto ((t,x) \mapsto v(t,x; \omega))$ and let $\Theta_0$ be the probability space $\Theta_0 := (\cC, \cB_{\cC}, Q_v)$.

\medskip

\noindent{\em Lebesgue integrals as Borel functions on path-space}
\medskip

\begin{lemma}
\label{Fact A0}
Let $b(t,x,z)$ be jointly measurable and satisfy \eqref{ch2'-s3.4.1.2}.

(a)  Let $\Phi: \cC \to \re$ be defined by
\beqn
    \Phi(w) = \int_0^T dt \int_D dx\, b^2(t, x, I_0(t,x) + w(t,x)).
    \eeqn
Then this is a Borel function with the following property:
 Let $(\Theta, W, u)$ be a weak solution to \eqref{ch2'-s3.4.1.3} such that $u - I_0$ has continuous sample paths. Define
\beqn
   Y = \int_0^T dt \int_D dx\, b^2(t,x, u(t,x)).
   \eeqn
   Then
   $Y = \Phi(u - I_0)$, $\tilde P$-a.s.

(b) Let $A \in \cB_{[0, T] \times D}^f$. Then
\beqn
F_A := \left\{w \in \cC : \int_0^T dt \int_D dx\, 1_A(t,x)\, \vert b(t, x, I_0(t,x) + w(t,x))\vert < \infty \right\}
\eeqn
 belongs to $\cB_{\cC}$ and $Q_v(F_A) = 1$. Let $\Phi_A : \cC \to \R$ be  the function defined by
 \beq
 \label{5.1.10(*1)}
    \Phi_A(w) = \int_0^T dt \int_D dx\, 1_A(t,x)\, b(t, x, I_0(t,x) + w(t,x))   
    \eeq
for $w \in F_A$ and  $\Phi_A(w) = 0$ for $w \in F_A^c$. Then $\Phi_A$ is Borel, and for $u$ as in part (a),
\beqn
\Phi_A(u - I_0)  = \int_0^T dt \int_D dx\, 1_A(t,x)\,  b(t, x, u(t,x)) =: Y_A,\qquad  \tilde P{\text{-a.s.}}
\eeqn
\end{lemma}

\begin{proof}
(a)\ By \eqref{5.1.5(*3)}, the map $(t,x; \omega) \mapsto b(t,x, u(t,x; \omega))$ belongs to $L^2([0, T] \times D \times \Omega, dt dx dP)$.
Define $f: \re_+ \times D \times \cC \to \re$ by
\beqn
f(t,x,w) = b^2(t,x, I_0(t,x) + w(t,x)).
\eeqn
 Since $g: \re_+ \times D \times \re \times \cC \to \re$ defined by $g(t,x,z,w) = (t,x,z+w(t,x))$ is jointly continuous, it is Borel, and $\tilde g(t,x,w) := g(t,x,I_0(t,x),w)$ and $f(t,x,w) = b^2(\tilde g(t,x,w))$ are therefore also Borel. By Fubini's theorem, $\Phi$ is a Borel function of $w$, and by definition, $f(t,x, u - I_0) = b^2(t,x, u(t,x))$, so $Y = \Phi(u - I_0)$, $\tilde P$-a.s.

(b)\ In the same way as above, we see that \beqn
(t, x, w) \mapsto 1_A(t,x)\, b(t, x, I_0(t,x) + w(t,x))
\eeqn
 is a Borel function, so $F_A \in \cB_{\cC}$ by Fubini's theorem. Notice that 
 \beqn
     1_{F_A}(u - I_0) = 1 \quad\text{on}\quad \left\{\int_0^T dt \int_D dx\, 1_A(t,x)\, \vert b(t, x, u(t,x))\vert < \infty \right\}, 
\eeqn
and this event has $\tilde P$-probability $1$ by Proposition \ref{ch2'-s3.4.2-p1} and the Cauchy-Schwarz inequality. Because the law under $Q_v$ of $\theta$ is the same as the law under $\tilde P$ of $u - I_0$, we deduce that $Q_v(F_A) = 1$, so $\Phi_A(w)$ is defined by \eqref{5.1.10(*1)} for $Q_v$-a.a.~$w \in \cC$ and $\Phi_A(u - I_0)  = Y_A$, $\tilde P$-a.s.
\end{proof}
\medskip

\noindent{\em Cylindrical Brownian motion as a Borel function on path space}
\medskip

In the sequel, $(e_j,\, j\ge 1)$ denotes a complete orthonormal basis of $L^2(D)$.

\begin{lemma}
\label{Fact C}
 Let $u$, $W$ and $\tilde W$ be as in Proposition \ref{ch2'-s3.4.2-t1}, such that $u - I_0$ has continuous sample paths. We have the following:

 (a) For all $\psi\in\cC_0^\infty(]0,T[\times D)$,
 \begin{align}\nonumber
    &\int_0^T ds \int_D dx\, \cL^\star \psi(s,y)(u(s,y)-I_0(s,y)) \\
     &\qquad\qquad = \int_0^T \int_D \psi(s,y)\, \tilde W(ds,dy),\qquad \tilde P{\text{-a.s.}}
\label{5.1.19-prime}
\end{align}
This property is equivalent to \eqref{ch2'-s3.4.2.6}.
\smallskip

 (b) For all $j \geq 1$, there is a Borel function $\Phi_j : \cC \to \cC([0, T])$ such that $\tilde P$-a.s., for all $ t\in [0, T]$, $\Phi_j(u - I_0)(t) = \tilde W_t(e_j)$.
 \smallskip

   (c) On the probability space $\Theta_0 = (\cC, \cB_{\cC}, Q_v)$, $(\Phi_j,\, j \geq 1)$ is a sequence of independent standard Brownian motions.
\end{lemma}
\begin{proof}
(a) By the definition of $\tilde W$ and Proposition \ref{ch2'-s3-r1}, 
this equality follows from \eqref{ch2'-s3.4.2.6}.
\smallskip

(b) The proof consists of two steps.
\smallskip

\noindent{\em Step 1.}\
Fix $j \geq 1$ and $t > 0$. We will show that there is a Borel function $\Phi_j(t): \cC \to \re$ such that $\Phi_j(t)(u - I_0) = \tilde W_t(e_j)$,   $\tilde P$-a.s.

Indeed, let $\varphi_{n}: [0, T] \times D \to \re$ be a sequence of $C_0^\infty$-functions with support in $]0, t[ \times D$ that converge in $L^2([0, T] \times D)$ to $1_{[0,t]}(\cdot) e_j(*)$. Define $\Psi_{n}: \cC \to \re$ by
\beqn
   \Psi_{n}(w) = \int_0^T ds \int_D dx\, \cL^\star \varphi_{n}(s,x)\, w(s,x) .
   \eeqn
This is a continuous, hence Borel, function of $w$, and by \eqref{5.1.19-prime},
\beqn
   \Psi_{n}(u - I_0) = \int_0^T \int_D \varphi_{n}(s,y)\, \tilde W(ds, dy),\quad  \tilde P{\text{-a.s.}}
   \eeqn

    As $n \to \infty$, $\Psi_{n}(u - I_0) \longrightarrow \tilde W_t(e_j)$ in $L^2(\Omega, \tilde P)$. Therefore, viewing $\Psi_{n}$ as a random variable on $(\cC, \cB_{\cC}, Q_v)$, $\Psi_{n}$ converges in $L^2(\cC, Q_v)$ to a random variable with the same law as $\tilde W_t(e_j)$, and along a subsequence $(n_i)$, this convergence is $Q_v$-a.s. Define $\Psi : \cC \to \re$ by $\Psi(w) = \limsup_{i \to \infty} \Psi_{n_i}(w)$. This is a Borel function of $w$. Under $Q_v,$ the sequence $(\Psi_{n})$ has the same law as $(\Psi_{n}(u - I_0))$ under $\tilde P$. Therefore,
    \begin{align*}
   \Psi(u - I_0) &= \limsup_{i \to \infty} \Psi_{n_i}(u - I_0)
    = \lim_{i \to \infty} \Psi_{n_i}(u - I_0) = \tilde W_t(e_j),\quad  \tilde P{\text{-a.s.}},
   \end{align*}
that is, $\Psi$ is the function $\Phi_j(t)$ in the statement of Step 1.
\smallskip

Before proceeding to the next step, we introduce some technical elements.
Let $\bD$ be the set of dyadic numbers in $[0, T]$. We endow the sets of continuous functions $\cC_1 := \cC(\bD, \re)$ and $\cC_2 := \cC([0,T], \re)$ with the topologies induced by their supremum norms, and with their Borel $\sigma$-fields.

 The following property holds:
 Let $\U_1:= \U(\bD) \subset \cC_1$ (resp.~$\U_2 := \U([0, T]) = \cC_2$) be the set of uniformly continuous functions from $\bD$ (resp.~$[0,T]$) to $\re$, each with the supremum norm. For $g \in \U_1$, let $I(g)$ be the (unique) continuous extension of $g$ to $\U_2$. Then $I: \U_1 \to \U_2$ is a bijective isometry. In particular,  $\U_1$ is a complete separable metric space, and $\U_1$ belongs to the $\sigma$-field $\cG$ on $\re^{\bD}$ generated by the coordinate functions.

Indeed, with the supremum norm, it is clear that $I$ is a bijective isometry, since an element $g$ of $\U_2$ determines and is determined by its restriction to $\bD$, which is $I^{-1}(g)$. Since $\U_1$ can be written $\cap_n \cup_m \cap_{r\in \bD} \cap _{s \in \bD} A(n, m, r, s)$, where
\beqn
   A(n, m, r, s) = \{g \in \cC_1: \vert g(r) - g(s) \vert < 1/n\}\quad  {\text{if}}\quad  \vert r - s \vert < 1/m,
   \eeqn
and $A(n, m, r, s) = \emptyset$ otherwise, $\U_1$ belongs to $\cG$. This ends the proof of the property that $I: \U_1 \to \U_2$ is a bijective isometry and $\U_1 \in \cG$.
\medskip

\noindent{\em Step 2.} \
For $ r \in \bD$ and $j\geq 1$, let $\Phi_j(r): \cC \to \re$  be the function given in Step 1. Define $ \Psi_j : \cC \to \R^\bD$ by $\Psi_j(w) = (\Phi_j(r)(w),\, r \in \bD)$. Then
$ \Psi_j(u - I_0) = (\tilde W_r(e_j),\, r \in \bD)$  $\tilde P$-a.s. 
In particular, $ Q_v$-a.s., $ \Psi_j$  is a uniformly continuous function from $\bD$  to $\re$, and there is a Borel function $\Phi_j: \cC\to \U_2$ such that $\tilde P$-a.s., $\Phi_j(u - I_0) = \tilde W_\cdot(e_j)$.
\medskip

Indeed, by Step 1, the law under $Q_v$ of the process $X_j := (\Phi_j(r),\, r \in \bD)$ is the same as the law under $\tilde P$ of $(\Phi_j(r)(u - I_0),\, r \in \bD)$, and 
\beqn
\left(\Phi_j(r)(u - I_0),\, r \in \bD\right) = (\tilde W_r(e_j),\, r \in \bD),\qquad  \tilde P{\text{-a.s}}.
\eeqn
Since this is the restriction of a Brownian motion on $[0, T]$ to the countable dense set $\bD$, it is $\tilde P$-a.s.~uniformly continuous on $\bD$. Equivalently, $X_j \in \U_1$, $Q_v$-a.s. Clearly, $\Psi_j: \cC \to \R^\bD$ defined by $\Psi_j(w) = (\Phi_j(r)(w),\, r \in \bD)$ is measurable, and
the inverse image $F_j$ of $\R^{\bD} \setminus \U_1$ under $\Psi_j $ is a $Q_v$-null set and a Borel subset of $\cC$. For $w \in \cC$, define $\Phi_j(w) = I(\Psi_j(w)) 1_{F_j^c}(w)$. 
Then $\Phi_j$ is a Borel function from $\cC$ into $\U_2$, and
\begin{align*}
   \Phi_j(u - I_0) &= I(\Psi_j(u - I_0)) 1_{F_j^c}(u - I_0) = I(\Psi_j(u - I_0)) \\
   &= I(\Phi_j(r)(u-I_0),\, r \in \bD)
    = I(\tilde W_r(e_j),\, r \in \bD)\\
    & = \tilde W_\cdot(e_j),\quad    \tilde P-{\text{a.s.}}
   \end{align*}
This ends the proof of Step 2 and of part (b) of the lemma.
\medskip

For part (c), note that by Step 2, $(\Phi_j(u - I_0),\, j \geq 1)$ is, under $\tilde P$, a sequence of independent standard Brownian motions. Therefore, the same is true of $(\Phi_j,\, j \geq 1)$ under $Q_v$.
\end{proof}
\smallskip

\noindent{\em Wiener integrals as Borel functions on path-space}
\medskip

\begin{lemma}
\label{Facts C1-C2-D}
Fix $j\ge 1$ and let $\Phi_j$ be as in Lemma \ref{Fact C}(b). 
Recall that $\Theta_0 = (\cC, \cB_\cC, Q_v)$.

(i)\, Let $h \in L^2([0, T])$, and on $\Theta_0$, define
\beqn
   X_j = \int_0^T h(t)\, \Phi_j(dt).
   \eeqn
Then there is a Borel function $\Psi_j: \cC \to \re$ such that $X_j = \Psi_j$, $Q_v$-a.s., and
\beqn
    \Psi_j(u - I_0) = \int_0^T h(t)\, d\tilde W_t(e_j),\qquad  \tilde P{\text{-a.s.}}
    \eeqn

    (ii)\, Fix $g \in L^2([0, T] \times D)$. Let $\hat W$  be the space-time white noise on $\Theta_0$  associated to $ (\Phi_j,\, j \geq 1)$  as in Lemma \ref {ch1'-lsi}(2).  On $\Theta_0$, define
    \beqn
    Y = \int_0^T \int_D g(t,x)\, \hat W(dt, dx).
    \eeqn Then there is a Borel function $ \Phi: \cC \to \re$  such that $ Y = \Phi$, $Q_v$-a.s., and
\beq
    \Phi(u - I_0) = \int_0^T \int_D g(t, x)\, \tilde W(dt, dx),\qquad  \tilde P{\text{-a.s.}}
    \eeq

    (iii)\, For  $n \geq 1$ and $A_\ell \in \cB_{[0, T] \times D}^f $, $\ell= 1, \dots, n $,
    let $\Phi_{A_\ell}$ be as defined in Lemma \ref{Fact A0}(b). There is a Borel-measurable random variable  $\Phi_n $, defined on  $\Theta_0 $ and with values in $\R^n$, such that
    \beqn
     \Phi_n = (\hat W(A_1) - \Phi_{A_1}, \dots, \hat W(A_n) - \Phi_{A_n}), \qquad   Q_v{\text{-a.s.}},
   \eeqn
 and
 \beqn
    \Phi_n(u - I_0) = (W(A_1),\dots, W(A_n)),\qquad   \tilde P{\text{-a.s.}}
 \eeqn
\end{lemma}

\begin{proof}
(i)\,  Suppose that $h$ is a simple function: $h(t) = h_0 1_{]t_1, t_2]}(t)$, with $h_0 \in \re$ and $0 \leq t_1 < t_2$. Then
\beqn
   X_j(w) = \left(\int_0^T h(t)\, \Phi_j(dt)\right)(w) = h_0\, (\Phi_j(w)(t_2) - \Phi_j(w)(t_1)),
   \eeqn
and the right-hand side is a Borel function $\Psi_j$ of $w$ from $\cC$ into $\re$. Further, by Lemma \ref{Fact C} (b),
\begin{align*}
     h_0\, (\Phi_j(u - I_0)(t_2) - \Phi_j(u - I_0)(t_1)) &= h_0\, (\tilde W_{t_2}(e_j) - \tilde W_{t_1}(e_j))\\
     & =  \int_0^T h(t)\, d\tilde W_t(e_j),\qquad   \tilde P{\text{-a.s.}}
     \end{align*}
     This proves (i) for simple functions and by linearity, (i) also holds for linear combinations of simple functions.

    Let $h \in L^2([0, T])$ and let $(h_n)$ be a sequence of linear combinations of simple functions such that $h_n \to h$ in $L^2([0, T])$. Let $\Psi^n_j: \cC \to \re$ be the Borel function associated to $h_n$ as in part (i). Then for all $w \in \cC$,
\beqn
    \Psi^n_j(w) = \left(\int_0^T h_n(t)\, \Phi_j(dt)\right)(w)
\eeqn
and
\beqn
    \Psi^n_j(u - I_0) = \int_0^T h_n(t)\, d\tilde W_t(e_j),\qquad   \tilde P\text{-a.s.}
\eeqn
Since $\Psi^n_j$ converges in $L^2(\Theta_0, Q_v)$ to $X_j = \int_0^T h(t)\, \Phi_j(dt)$ as $n\to \infty$, there is a subsequence $(n_i(j))$ such that $\Psi^{n_i(j)}_j$ converges $Q_v$-a.s.~to $X_j$ (this subsequence depends on $j$, $Q_v$ and the $h_n$, but on nothing else).

Define $\Psi_j : \cC \to \re$ by $\Psi_j(w) = \limsup_{i \to \infty} \Psi^{n_i(j)}_j(w)$. Then $\Psi_j$ is Borel, and
\beqn
   \Psi_j = \limsup_{i \to \infty} \Psi^{n_i(j)}_j = \lim_{i \to \infty} \Psi^{n_i(j)}_j = X_j, \qquad  Q_v{\text{-a.s.}}
   \eeqn
Further,
\begin{align}
\label{Prop5.1.7(*3)}
   \Psi_j(u - I_0) &= \limsup_{i \to \infty} \Psi^{n_i(j)}(u - I_0)
   =  \lim_{i \to \infty} \Psi^{n_i(j)}_j (u - I_0)\notag\\
   & = \lim_{i \to \infty} \int_0^T h_{n_i(j)}(t)\, d\tilde W_t(e_j), \qquad  \tilde P{\text{-a.s.}}    
\end{align}
Indeed, the third equality holds because the law under $\tilde P$ of $u - I_0$ is equal to the law under $Q_v$ of $\theta$. Since the $\tilde P$-a.s.-limit in \eqref{Prop5.1.7(*3)} is the same as the $L^2(\Omega, \tilde P)$-limit, we see that
    \beqn
    \Psi_j(u - I_0) = \int_0^T h(t)\, d\tilde W_t(e_j),\qquad    \tilde P{\text{-a.s.}}
    \eeqn
Therefore, $\Psi_j$ has the properties claimed in statement (i).

We now prove (ii).
By definition,
  $Y = \lim_{n \to \infty} Y_n$  in $L^2(\Theta_0, Q_v)$ and $Q_v$-a.s.~by Lemma \ref{ch1'-lsi}(2),
where
      \beqn
      Y_n = \sum_{j=1}^n \int_0^T \langle g(t, *), e_j \rangle_V\, \Phi_j(dt).
      \eeqn
By part (i), there is a Borel function $\Psi_{n}: \cC \to \re$ such that
  $ \Psi_n = Y_{n}$,  $Q_v$-a.s., and
\beq
\label{Prop5.1.7(*4a)}
   \Psi_n(u - I_0) = \sum_{j=1}^{n} \int_0^T \langle g(t, *), e_j \rangle_V\, d\tilde W_t(e_j),\qquad   \tilde P{\text{-a.s.}}  
   \eeq
Define $\Psi := \limsup_{n\to \infty} \Psi_n$. Then $\Psi : \cC \to \re$ is Borel and
\beqn
   \Psi = \limsup_{n \to \infty} \Psi_n = \limsup_{n \to \infty} Y_{n} = \lim_{n\to \infty} Y_{n} = Y,\qquad   Q_v{\text{-a.s.}}
   \eeqn
Further,
   \begin{align}
\label{Prop5.1.7(*4)}
   \Psi(u - I_0) &= \limsup_{n \to \infty} \Psi_n(u - I_0) =  \lim_{n \to \infty}  \Psi_n(u - I_0)\notag\\
   & = \lim_{n \to \infty}  \sum_{j=1}^{n} \int_0^T \langle g(t, *), e_j \rangle_V\, d\tilde W_t(e_j),\qquad   \tilde P{\text{-a.s.}}
   \end{align}  
Indeed, the second equality holds because the law under $\tilde P$ of $u - I_0$ is the same as the law under $Q_v$ of $\theta(w) = w$, and the last equality holds by \eqref{Prop5.1.7(*4a)}.
   Since the a.s.-limit in \eqref{Prop5.1.7(*4)} is the same as the  $L^2(\Omega, \tilde P) $-limit, we see that
   \beqn
   \Psi(u - I_0) = \int_0^T \int_D g(t, x)\, \tilde W(dt, dx),\qquad  \tilde P{\text{-a.s.}}
   \eeqn
This ends the proof of (ii).

Finally, we prove (iii). For $\ell = 1,\dots, n$, let $\Phi_{A_\ell} + \tilde \Phi_\ell : \cC \to \re$ be the Borel function given in (ii) for $g := 1_{A_\ell}$. Then
$\tilde \Phi_\ell = \hat W(A_\ell) - \Phi_{A_\ell}$,  $Q_v$-a.s., and
\beqn
\tilde \Phi_\ell(u - I_0) = \tilde W(A_\ell) - \int_0^T dt \int_D dx\, 1_{A_\ell}(t,x) b(t, x, u(t,x)) = W(A_\ell), \quad  \tilde P{\text{-a.s.}},
\eeqn
where we have used Proposition \ref{ch2'-s3-r1}.
It suffices therefore to set $\Phi_n(w) = (\tilde \Phi_1(w),\dots, \tilde\Phi_n(w))$.

This completes the proof of the lemma.
\end{proof}

\medskip

\noindent{\em Stochastic integrals as Borel functions on path-space}
\medskip

\begin{lemma}
\label{Facts E1-E2}
Fix $j\ge 1$ and let $\Phi_j$ be as in Lemma \ref{Fact C}(b). Let
$\Theta_0$ be the probability space given in Lemma \ref{Fact C}(c). 
\smallskip

(1)\, Let $b : [0, T] \times [0,L] \times \re \to \re$ be a Borel function. Let $G_j: [0,T] \times \cC$ be defined by 
\beqn
   G_j(t, w) = \int_D dx\, b(t, x, I_0(t,x) + w(t,x))\, e_j(x)
\eeqn 
and set $X_j := \int_0^T G_j(t, w)\, \Phi_j(dt)$. Then there is a Borel-measurable random variable $\Psi_j$ defined on $\Theta_0$ such that:\smallskip

   (1.a) $X_j = \Psi_j$,  $Q_v$-{\text{a.s.}}\smallskip

\noindent and\smallskip

   (1.b) $\Psi_j(u - I_0) = \int_0^T K_j(t)\, d\tilde W_t(e_j)$,  $\tilde P$-{\text{a.s.}},\smallskip
   
\noindent where 
\beqn
    K_j(t) = G_j(t, u - I_0) = \int_D dx\, b(t, x, u(t,x)) e_j(x).
\eeqn

(2)\, Let $G(t,x)(w) = b(t,x, I_0(t,x) + \theta(t,x)(w))$ and let $\hat W$ be as in Lemma \ref{Facts C1-C2-D}(ii). There is a Borel-measurable random variable $\Psi$ defined on $\Theta_0$ such that:\smallskip

    (2.a) $\Psi = X$,   $Q_v$-{\text{a.s.}}, where $X := \int_0^T \int_D G(t,x)\, \hat W(dt, dx)$,\smallskip
    
\noindent 
and \smallskip

    (2.b) $\Psi(u - I_0) = \int_0^T \int_D b(t,x, u(t,x))\, \tilde W(dt, dy)$,  $\tilde P$-{\text{a.s.}}
\end{lemma}

\begin{proof}
Let $\cG_t^0$ be the $\sigma$-field on $\cC$ generated by $(\theta(s, *),\, s \leq t)$. We complete $\cG_t^0$ with $Q_v$-null sets and make this into a complete and right-continuous filtration $(\cG_t)$. Let $\cG$ be the smallest $\sigma$-field generated by all the $\cG_t$.

(1)\, Notice that $(t, w) \mapsto G_j(t, w)$ is progressively measurable relative to the (uncompleted) filtration $(\cG_t^0)$, because $(t,w) \mapsto w(t, \ast)$ is continuous, hence progressively measurable relative to the filtration $(\cG_t^0)$, therefore $(t,x,w) \mapsto  b(t, x, I_0(t,x) + w(t,x))$ is also progressively measurable.

   For $\ell = 1,\dots, \lfloor T 2^n\rfloor$ (where $\lfloor T 2^n\rfloor $ denotes the integer part of $T 2^n$), define
   \beqn
   G_j^n(t, w) := 2^n \int_{(\ell-1)2^{-n}}^{\ell 2^{-n}} G_j(r, w)\, dr,\qquad   {\text{for}}\ t \in [ \ell 2^{-n} , (\ell+1) 2^{-n}[.
   \eeqn
Then $(t, w) \mapsto  G_j^n(t, w)$ is progressively measurable relative to the filtration $(\cG_t^0)$. In particular, $(G_j^n(t),\, t \in [0, T])$ is an elementary process that is jointly measurable and adapted to $(\cG_t^0)$, and
\beqn
   X_j^n := \int_0^T G_j^n(t, w)\, \Phi_j(dt) = \sum_{\ell = 0}^{\lfloor T 2^n\rfloor -1} G_j^n(\ell 2^{-n}, w) (\Phi_j((\ell + 1) 2^{-n}) - \Phi_j(\ell 2^{-n})).
   \eeqn
Clearly, $\Psi_j^n := X_j^n$ is a Borel function from $\cC$ into $\re$. Further,
\begin{align}
\label{Prop5.1.7(*5)}
   \Psi_j^n(u - I_0) &=  \sum_{\ell = 0}^{\lfloor T 2^n\rfloor -1} G_j^n(\ell 2^{-n}, u - I_0) (\Phi_j((\ell + 1) 2^{-n}, u - I_0) - \Phi_j(\ell 2^{-n}, u - I_0)) \notag\\
   &= \int_0^T K_j^n(t)\, d\tilde W_t(e_j),\qquad  \tilde P{\text{-a.s.}},    
   \end{align}
where
\beqn
   K_j^n(t) = G_j^n(t, u - I_0) = 2^n \int_{(\ell-1)2^{-n}}^{\ell 2^{-n}}\, dr \int_0^L dx\, b(r, x, u(r,x)) e_j(x)
   \eeqn
    if $\ell 2^{-n} \leq t < (\ell+1) 2^{-n}$.
  This proves (1) for $G_j$ there replaced by $G_j^n$.

  Next, we will apply the following fact:  On any filtered probability space $(\Omega, \F, P)$, let $f \in L^2([0,T] \times \Omega, dtdP)$ be jointly measurable and adapted. For $n \geq 1$, define
\beqn
   f_n(t) := 2^n \int_{(\ell-1)2^{-n}}^{\ell 2^{-n}} f(r)\, dr\quad   {\text{if}}\quad  \ell 2^{-n} \leq t < (\ell+1) 2^{-n}.
   \eeqn
Then $(f_n)$ is a sequence of simple functions and $\lim_{n\to\infty} f_n = f$ in $L^2([0,T] \times \Omega, dtdP)$.

Indeed, let
\beqn
    g_n(t) := 2^n \int_{\ell 2^{-n}}^{(\ell +1) 2^{-n}} f(r)\, dr\quad   {\text{if}}\quad  \ell  2^{-n} \leq t < (\ell +1) 2^{-n},\quad n\ge 1.
    \eeqn
By the $L^2$-martingale convergence theorem, $\lim_{n\to\infty }g_n = f$ in $L^2([0,T] \times \Omega, dt dP)$. A  simple calculation that shows that $\lim_{n\to\infty} ( f_n - g_n) = 0$ in $L^2([0,T] \times \Omega, dtdP)$. Hence the fact is proved.

 As a consequence of this fact, $G_j^n \to G_j$ in $L^2([0, T] \times \cC, dt dQ_v)$ as $n\to\infty$. Therefore, $X_j^n \to X_j$ in $L^2(\cC, Q_v)$, and along a subsequence $(n_i(j))$, $X_j^{n_i(j)} \to X_j$,    $Q_v$-{\text{a.s.}} Define $\Psi_j := \limsup_{i \to \infty} \Psi_j^{n_i(j)}$. Then $\Psi_j$ is a Borel function, and $Q_v$-{\text{a.s.}},
    \beqn
     \Psi_j := \limsup_{i \to \infty} X_j^{n_i(j)} = \lim_{i \to \infty} X_j^{n_i(j)} = X_j.
     \eeqn
This proves (1.a).

For (1.b), we note that
\begin{align*}
     \Psi_j(u - I_0) &= \limsup_{i \to \infty} \Psi_j^{n_i(j)}(u - I_0) \\
     &= \lim_{i \to \infty} \Psi_j^{n_i(j)}(u - I_0) = \lim_{i \to \infty} \int_0^T K_j^{n_i(j)}(t)\, d\tilde W_t(e_j),\qquad   \tilde P{\text{-a.s.}}
     \end{align*}
Indeed, the $\limsup$ is a limit because $u - I_0$ has the same law under $\tilde P$ as $\theta$ under $Q_v$, and the third equality holds by \eqref{Prop5.1.7(*5)}. Since $G_j^n \to G_j$ in $L^2([0, T] \times \cC, dt dQ_v)$, we have $K_j^n \to K_j$ in $L^2([0, T] \times \Omega, dt d\tilde P)$, therefore,
\beqn
   \int_0^T K_j^{n_i(j)}(t)\, d\tilde W_t(e_j) \longrightarrow \int_0^T K_j(t)\, d\tilde W_t(e_j)
   \eeqn
      in $L^2(\Omega, \tilde P)$,
and since the $L^2$-limit is the same as the a.s.-limit,
\beqn
   \Psi_j(u - I_0)  =  \int_0^T K_j(t)\, d\tilde W_t(e_j),\qquad    \tilde P{\text{-a.s.}}
\eeqn

We now prove (2).
Recall that $X = \sum_{j=1}^\infty \int_0^T \langle G(t,*), e_j \rangle_V\, \Phi_j(dt)$
and the series converges in $L^2(\cC, Q_v)$. Notice that
\beq
\label{Prop5.1.7(*6)}
   \langle G(t,*), e_j \rangle_V = G_j(t),   
   \eeq
where $G_j$ is defined in part (1).

By part (1), there is a Borel-measurable random variable $\Psi_j$ defined on $\Theta_0$ that satisfies (1.a) and (1.b) of (1). By Lemma \ref{ch1'-lsi} (2),
\beqn
    (\hat W_s(e_j),\, s \in [0, T],\, j \geq 1) =(\Phi_j(s),\, s \in [0, T],\, j \geq 1),\qquad  Q_v{\text{-a.s.}}
    \eeqn
Therefore, as checked just after Definition \ref{ch1'-s4.d1}, as $n\to\infty$,
\beqn
   X_n := \sum_{j=1}^n \int_0^T \langle G(t,x), e_j \rangle_V\, \Phi_j(dt)
   \eeqn
converges in $L^2(\cC, Q_v)$ to
   $X = \int_0^T \int_D G(t,x)\, \hat W(dt, dx)$.
Along a subsequence $(n_i)$, $X_{n_i} \to X$,  $Q_v$-a.s. Define the Borel function $Z^i : \cC \to \re$ by $Z^i = \sum_{j=1}^{n_i} \Psi_j$ and $\Psi = \limsup_{i \to \infty} Z^i$. Then $\Psi$ is Borel and $\Psi = X$  $Q_v$-a.s., establishing (2.a).

  From \eqref{Prop5.1.7(*6)}, it follows that
  \beq
\label{Prop5.1.7(*7)}
    \langle G(t,*)(u - I_0), e_j \rangle_V = G_j(t, u - I_0) = K_j(t),\qquad  \tilde P{\text{-a.s.}},   
    \eeq
where $K_j$ is defined in (1.b), and
\beqn
   \Psi_j(u - I_0) = \int_0^T K_j(t)\, d\tilde W_t(e_j),\qquad  \tilde P{\text{-a.s.}}
   \eeqn
Further,
\begin{align}
    \Psi(u - I_0) &=  \limsup_{i \to \infty} Z^i(u - I_0) = \lim_{i \to \infty} Z^i(u - I_0) 
     =  \lim_{i \to \infty} \sum_{j=1}^{n_i} \Psi_j(u - I_0)\notag\\
    & = \lim_{i \to \infty} \sum_{j=1}^{n_i} \int_0^T \langle G(t,*)(u - I_0), e_j \rangle_V\, d\tilde W_t(e_j),\qquad  \tilde P{\text{-a.s.}}     
 \label{Prop5.1.7(*8)}
\end{align}
Indeed, the second equality holds because the law under $\tilde P$ of $u- I_0$ is equal to the law under $Q_v$ of $\theta$ (recall that $\theta(w) = w$), and the fourth equality is due to \eqref{Prop5.1.7(*7)} and (1.b).
Note that 
\begin{align*}
    & \sum_{j=1}^{n_i} \int_0^T \langle G(t,*)(u - I_0), e_j \rangle_V\, d\tilde W_t(e_j) \\
  &\qquad \qquad \longrightarrow
  \int_0^T \int_0^L  G(t, x)(u - I_0)\, \tilde  W(dt, dx) 
\end{align*}
in $L^2(\Omega, \tilde P)$ as $i \to \infty$, and this limit must be the same as the a.s.-limit in \eqref{Prop5.1.7(*8)}. Since $G(t,x)(u - I_0) = b(t,x, u(t,x)$),  this show that
\beqn
    \Psi(u - I_0) = \int_0^T \int_0^L  b(t,x, u(t,x))\, \tilde W(dt, dx),\qquad   \tilde P{\text{-a.s.}},
    \eeqn
and (2.2) is proved.
\end{proof}





\section{Two properties of the stochastic heat equation (additive noise)}\label{rd03_17s1}

In this section, we consider two properties of the stochastic heat equation: in Subsection \ref{ch2'-s3.4.3}, we prove the equivalence of the laws of the solutions on $D=\re$ and  $D=\, ]0,L[$ when these solutions are restricted to a bounded space-time rectangle away from the space-time boundaries, and in Subsection \ref{ch2'-s3.4.4}, we study a (space-time) Markov field property of these solutions.

A common feature in the approach to these two questions is the use of Girsanov's theorem (Section \ref{ch2'-s3}).

\subsection{Local equivalence in law of solutions on $[0, L]$ and on $\R$}
\label{ch2'-s3.4.3}

In this section, we first consider a linear stochastic heat equation on $\re$, as in Section \ref{ch4-ss2.1-1'} (see \eqref{4.20-1'}),
along with a linear stochastic heat equation on the interval $[0,L]$, as in Section \ref{ch4-ss2.2-1'}, with either vanishing Dirichlet or Neumann boundary conditions (see \eqref{ch1'.HD} and \eqref{1'.14}). Each equation is driven by a space-time white noise. 
We first show that the laws of these two solutions are mutually equivalent (after restricting to a closed rectangle in $]0,T] \times ]0,L[\,$: see Theorem \ref{ch1'-s3.4.3-t1}). Then we extend this result to nonlinear stochastic heat equations with additive noise (see Theorem \ref{ch5-t-5.1.11}). The linear case was considered in \cite{m-t2002}.
\medskip

\noindent{\em The linear case}
\medskip

We will denote by $u$ the solution to \eqref{4.20-1'} with initial condition  $u(0,\ast)=u_0(\ast)$, and assume that $u_0: \re \to \re$ satisfies \eqref{ch1'-v00}. The notation $u_L$ will refer to the solution to \eqref{ch1'.HD} (or  \eqref{1'.14}) with initial condition $u_L(0,\ast)=u_{0, L}(\ast)$, where $u_{0, L}: [0,L] \to \re$ is a bounded Borel function. Recall  that for $t \in [0, T]$, $x \in \R$ and $z \in [0, L]$,
\begin{align}
 u(t,x)&=I_{0}(t,x)  + \int_0^t \int_{\re} \Gamma(t-s,x-y)\, W(ds,dy), \label{ch2'-s3.4.3.1}\\
 u_L(t,z)&=  I_{0,L}(t,z) + \int_0^t \int_0^L G_L(t-s;z,y)\, W(ds,dy), \label{ch2'-s3.4.3.2}
 \end{align}
  where
   \begin{align}
  \label{in-cond}
  I_{0}(t,x) & = \int_{\re} dy\  \Gamma(t,x-y)\, u_0(y),\notag\\
  I_{0,L}(t,z) & =  \int_0^L dy\ G_L(t;z,y)\, u_{0, L}(y).
  \end{align}
  The function $\Gamma$ is defined in \eqref{heatcauchy-1'}, and $G_L$ is given in \eqref{ch1'.600} in the case of vanishing Dirichlet boundary conditions (respectively \eqref{1'.400} in the case of vanishing Neumann boundary conditions).

  \begin{thm}
 \label{ch1'-s3.4.3-t1}
 Fix $t_0\in\, ]0,T]$ and $\varepsilon\in\, ]0,L/3[$. Then the laws of the processes $(u_L(t,x),\, (t,x)\in[t_0,T[\times[\varepsilon,L-\varepsilon])$ and
 $(u(t,x),\, (t,x)\in[t_0,T[\times[\varepsilon, L-\varepsilon])$ are mutually equivalent.
 \end{thm}
 
 \begin{proof}
 We can and will assume that the same space-time white noise $\dot W$ is used in the equations for $u$ and $u_L$, since the law of the solution does not depend on the specific choice of the space-time white noise. It suffices to prove the theorem for $\ep > 0$ small enough. For $(t,x) \in [0,T] \times [0,L]$, define $\tilde u(t,x) = u(t,x) - u_L(t,x)$. Let $I = [t_0,T]$, $J = [\ep, L - \ep]$. Recall from Theorem \ref{ch1'-tDL} that on $I \times J$, $\tilde u = u - u_L$ has $\cC^\infty$ sample paths a.s. and satisfies \eqref{ch1'-tDL.1}.

    Fix $\psi : [0,T] \times \re \to \re$ such that $\psi \in \cC^\infty$, ${\rm supp}\ \psi \subset I_{\ep/2} \times J_{3\ep/4}$, and $\psi\vert_{I_\ep \times J_\ep} \equiv 1$, where $I_\ep = [t_0 + \ep, T]$  and $J_\ep = [2 \ep, L - 2\ep]$.
   For $(t, x) \in [0, T] \times [0,L]$, define $b(t,x) = \psi(t,x)\, \tilde u(t,x)$ and
      $\tilde v(t,x) = u_L(t,x) + b(t,x)$.
Notice that $\tilde v(t,x) = u(t,x)$ for $(t,x) \in I_\ep \times J_\ep$, and $b$ is $\cC^\infty$ on $[0,T] \times [0,L]$  a.s.

For each $\cC^\infty$ sample path of $(b(t,x))$, define
\beq
\label{ch2'-s3.4.3.2-bis}
 h(t,x) = \left(\frac{\partial}{\partial t}- \frac{\partial^2}{\partial x^2}\right) b(t,x).
 \eeq
Notice that $h$ is $\cC^\infty$ a.s.; moreover, a.s., for $(t,x) \in [0,T] \times [0,L]$,
\beq
\label{5.1.3(*1)-bis}
   b(t,x) = \int_0^t \int_0^L  G_L(t-s;x,y) h(s,y)\,  dsdy.   
   \eeq
Indeed, because of the choice of the support of $\psi$, the function $b$ clearly satisfies the vanishing initial condition, the vanishing Dirichlet (resp. Neumann) boundary conditions, and, by definition, $b$ satisfies the deterministic heat equation
\eqref{ch2'-s3.4.3.2-bis}.
 The same properties are true for the right-hand side of \eqref{5.1.3(*1)-bis}. Hence, by uniqueness of the solution to the heat equation with smooth inhomogeneity and given initial and Dirichlet (resp. Neumann) boundary conditions, \eqref{5.1.3(*1)-bis} holds.

   By \eqref{5.1.3(*1)-bis}, a.s., for $(t,x) \in [0,T] \times [0,L]$,
   \begin{align}
   \label{5.1.3(*2)-bis}
   \tilde v(t,x) &= u_L(t,x) + \int_0^t \int_0^L  G_L (t-s;x,y) h(s,y)\,  dsdy\notag\\
                      & = I_{0,L}(t,x) + \int_0^t \int_0^L G_L(t-s;x,y)\, W(ds,dy)\notag\\
                      &\qquad \qquad+ \int_0^t \int_0^L  G_L (t-s;x,y) h(s,y)\,  dsdy.      \end{align}
Since $\tilde v(t,x) = u(t,x)$ for $(t,x) \in I_\ep \times J_\ep,$ $t_0 > 0$, and $\ep$ can be taken arbitrarily small, the conclusion of the theorem will follow from \eqref{5.1.3(*2)-bis} and Theorem \ref{ch2'-s3.4.1-t0} provided that the assumptions of that theorem are satisfied.

 Define a measure $\tilde P$ by
 \beqn
 \frac{d\tilde P}{dP} = \exp\left(-\int_0^T \int_0^L h(t,x)\, W(dt,dx)-\frac{1}{2}\int_0^Tdt \int_0^L dx\, h^2(t,x)\right).
 \eeqn
 We need to check that $E\left[\frac{d\tilde P}{dP}\right] =1$. According to condition (c) of Proposition \ref{ch2'-s3-p1}, it suffices to verify that there is $\varepsilon>0$ such that
 \beq
 \label{ch1'-s3.4.3.7}
 \sup_{s\in[0,T]}E\left[\exp\left(\half \int_s^{s+\ep} dt \int_0^L dx\, h^2(t,x)\right)\right] < \infty.
 \eeq
  In order to establish \eqref{ch1'-s3.4.3.7}, notice that by definition, $h$ vanishes outside of $I_{\ep/2} \times J_{3\ep/4}$, and inside this rectangle, $h$ is a linear combination of derivatives of orders $0$, $1$ and $2$ of $\psi$ and $\tilde u$. Those of $\psi$ are deterministic  and uniformly bounded over $[0,T] \times [0,L]$. Those of $\tilde u$ satisfy \eqref{ch1'-tDL.1}. Therefore,
  \beqn
    \sup_{(t,x) \in [0,T] \times [0,L]} E[h^2(t,x)] < \infty.
    \eeqn
In addition, inside $I_{\ep/2} \times J_{3\ep/4}$, according to \eqref{ch1'-tDL.4}, each partial derivative of $\tilde u$ is a difference of two stochastic integrals of deterministic integrands with respect to space-time white noise, to which must be added the partial derivatives of the contributions of $I_0(t,x)$ and $I_{0,L}(t,x)$ (see \eqref{in-cond}) that play no role here.  Therefore, $h$ is a Gaussian process. Proceeding as for \eqref{5.1.3(*2)}, we deduce that for any finite $C > 0$, there is $\ep > 0$ such that
\beqn
   \sup_{s \in [0,T]} E\left[\exp\left(C \int_s^{s+\ep} \int_0^L h^2(t,x)\, dt dx\right)\right] < \infty.
   \eeqn
Taking $C = \half,$ we obtain \eqref{ch1'-s3.4.3.7}, and
 Theorem \ref{ch1'-s3.4.3-t1} is proved.
 \end{proof}
 \medskip
 
 \noindent{\em The nonlinear case with additive noise}
 \medskip

 For $i=1,2$, let
 \begin{align}
\label{ch2'-s3.4.1.3-two}
u_i(t,x) &= I_{i,0}(t,x) + \int_0^t \int_{D_i} \Gamma_i(t,x;s,y)\, W(ds,dy)\notag\\
 &\qquad\qquad + \int_0^t \int_{D_i} \Gamma_i(t,x;s,y) b_i(s,y,u_i(s,y))\, ds dy,
\end{align}
$(t,x)\in[0,T]\times D_i$, where $D_1 = \re$, $D_2 = [0,L]$, and the function $\Gamma_i$ is the fundamental solution (or the Green's function) of the heat equation on $D_i$ (with vanishing Dirichlet or Neumann boundary conditions if $ i= 2$).

We assume that $I_{i,0}$ is a bounded Borel function on $[0,T] \times D_i$, and that $b_i : [0,T] \times D_i \times \re \to \re$  satisfies condition 3.~of Theorem \ref{ch2'-s3.4.1-t1}, with $D$ replaced by $D_i$. Then, according to Theorem \ref{ch2'-s3.4.1-t1}, the SPDE \eqref{ch2'-s3.4.1.3-two} has a weak solution in law.

The following theorem is an extension of Theorem \ref{ch1'-s3.4.3-t1}.
 \begin{thm} 
 \label{ch5-t-5.1.11}
   Fix  $t_0 \in\, ]0,T]$ and $\ep \in\, ]0, L/3[$. Then the laws of the processes
      $(u_i(t,x),\, (t,x) \in [t_0, T] \times [\ep,  L - \ep])$, $i=1, 2$,
given in \eqref{ch2'-s3.4.1.3-two} are mutually equivalent.
\end{thm}
\begin{proof}
 Notice that for $i=1, 2$, conditions 1.~and 2.~of Theorem \ref{ch2'-s3.4.1-t1} are satisfied (see Remark \ref{ch5-r1}). Let $u$ be defined by \eqref{ch2'-s3.4.3.1} and $u_L$ by \eqref{ch2'-s3.4.3.2}
 with $I_0$ and $I_{0,L}$ replaced respectively by $I_{1,0}$ and $I_{2,0}$. By Proposition \ref{ch2'-s3.4.2-t1}, the laws of $u_1$ and $u$ are mutually equivalent (on $D_{1,T} := [0,T] \times \re$), and the laws of $u_2$ and $u_L$ are also mutually equivalent (on $D_{2,T} := [0,T] \times [0,L]$). By Theorem \ref{ch1'-s3.4.3-t1}, the laws of $u$ and $u_L$ restricted to $[t_0, T] \times [\ep, L - \ep]$ are mutually equivalent. Since this is a subset of $D_{1,T}$ and $D_{2,T}$, the theorem is proved.
\end{proof}


\subsection{The Markov field property }
\label{ch2'-s3.4.4}

In this section, we study the Markov field property of the weak solution in
law $u = (u(t,x),\, (t,x) \in [0,T] \times [0,L])$, in the sense of Definition \ref{ch2'-s3.4.1-d1}, to the nonlinear stochastic heat equation with additive space-time white noise and vanishing Dirichlet boundary conditions:
\beq
\label{ch2'-s3.4.4.1}
\begin{cases}
\frac{\partial u}{\partial t}(t,x)-\frac{\partial^2u}{\partial x^2}(t,x) = b(t,x,u(t,x)) + \dot W(t,x), & (t,x)\in\,]0,T[\times ]0,L[,\\
u(0,x)= u_0(x), & x\in[0,L],\\
u(t,0)=u(t,L)=0, & 0<t\le T,
\end{cases}
\eeq
where $u_0 : [0,L] \to \R$ is a bounded Borel function.
According to Remark \ref{ch5-r1}, the associated Green's function $G_L$ satisfies condition 2.~of Theorem \ref{ch2'-s3.4.1-t1}. We will require of $u_0$ and $b$ that the other two assumptions of this theorem be satisfied, so that, in particular, the conclusion of Proposition \ref{ch2'-s3.4.2-t1} holds.

We want to show that $u$ satisfies the so-called {\em germ-field Markov property}, which we now define. For any set $A\subset [0,T]\times[0,L]$, and any random field $\xi = \left(\xi(t,x),\, (t,x)\in[0,T]\times[0,L] \right)$, we define the following $\sigma$-fields:
\begin{itemize}
\item {\em the sharp field}: \qquad $\tf^\xi(A) = \sigma\left(\xi(t,x),\, (t,x)\in A\right)$,
\item {\em the germ field}: \qquad\  $\mathcal{G}^\xi(A) = \bigcap_{\mathcal{O}\supset A,\, \mathcal{O}\,  \text{open}} \ \tf^\xi(\mathcal{O})$,
\end{itemize}
completed with the $\sigma$-field generated by $P$-null sets.\index{sharp field}\index{field!sharp}\index{germ-field}\index{field!germ-}
\begin{def1}
\label{ch2'-s3.4.4-d1}
The random field $\xi$ has the {\em germ-field Markov property}\index{Markov property!germ-field}\index{germ-field!Markov property} (respectively the {\em sharp Markov property})\index{Markov property!sharp}
if for any open set $A\subset [0,T]\times[0,L]$, the $\sigma$-fields
$\mathcal{F}^\xi(A)$ and $\mathcal{F}^\xi(A^c)$ are conditionally independent given $\mathcal{G}^\xi(\partial A)$ (respectively $\cF^\xi (\partial A)$), where $\bar A$, $A^c$ and $\partial A$ are, respectively, the closure, the complement and the boundary of $A$.
\end{def1}

\begin{remark}
\label{ch2'-s3.4.4-r0}
It is possible to show that: (1) the germ-field Markov property implies that the property ``$\tf^\xi(A)$ is conditionally independent of $\mathcal{F}^\xi(A^c)$ given  $\mathcal{G}^\xi(\partial A)$" holds for all Borel sets $A\subset [0,T]\times[0,L]$; 
and (2)  if the property  ``$\tf^\xi(A)$ is conditionally independent of $\mathcal{F}^\xi(A^c)$ given  $\mathcal{G}^\xi(\partial A)$" holds for open sets $A$ with smooth boundary, then $\xi$ has the germ-field Markov property  (see \cite[p. 20]{nualartpardoux1}).
\end{remark}
\medskip

\noindent{\em The linear case}
\medskip

We begin by studying the germ-field Markov property of the weak solution in law
$v = (v(t,x),\, (t,x) \in [0,T] \times [0,L])$ to \eqref{ch2'-s3.4.4.1}, assuming $b \equiv 0$, that is,
\begin{align}
\nonumber
    v(t,x) &= \int_0^L dy \, G_L(t;x,y) u_0(y) \\
        &\qquad \qquad  + \int_0^t \int_0^L G_L(t-s;x,y)\, W(ds,dy),
\label{ch2'-s3.4.4.2}
\end{align}
with $G_L$ defined in \eqref{ch1'.600}.
By Proposition \ref{4-p1-1'-weak}, $v$ has a continuous version on $]0,T[ \times ]0,L[\ $, and by Proposition \ref{ch2'-s3.4.2-t1}, the same is true of $u$.
Further, if $u_0 \equiv 0$, then by Proposition \ref{4-p1-1'}, $v$ even has a continuous version on $[0,T] \times [0,L]$, and $(t,x) \mapsto v(t,x) \in L^2(\Omega)$ is continuous.  We will use these continuous versions of $u$ and $v$.


\begin{prop}
\label{ch2'-s3.4.4-p1}
The random field $v$ has the germ-field Markov property.
\end{prop}
Before proving this proposition, we introduce some preliminary results.
Since the first integral on the right-hand side of \eqref{ch2'-s3.4.4.2} is deterministic, it does not affect the definition of the relevant $\sigma$-fields, so we will assume up to the end of the proof of Proposition \ref{ch2'-s3.4.4-p1} that $u_0 \equiv 0$. In this case, $v$ is a centred Gaussian random field, for which the following notion is useful.
\medskip

\noindent{\em Reproducing Kernel Hilbert Space}
\medskip

Let $H$ denote the closed linear Gaussian subspace of $L^2(\Omega,\cF,P)$ spanned by the random variables $\left(v(t,x),\, (t,x)\in[0,T]\times[0,L] \right)$. Clearly, $H$ is the closure in $L^2(\Omega)$ of the vector space of finite linear combinations of these random variables.

The space $H$ admits an equivalent description, as the following lemma shows. This will be used in the proof of Lemma \ref{ch2'-s3.4.4-l2} below.
\begin{lemma}
\label{ch2'-s3.4.4-l1}
Let $\tilde H:=\{W(\varphi),\, \varphi\in L^2([0,T]\times[0,L])\}$. Then $\tilde H = H$.
\end{lemma}
\begin{proof}
Since $(s,y)\mapsto G_L(t-s;x,y) 1_{]0,t[}(s)$ belongs to $L^2([0,T]\times[0,L)]$, it is clear that $H\subset\tilde H$.

For the converse inclusion, we use Lemma \ref{5.1.16a} (a) below with $u_0\equiv 0$ to see that for any $\varphi \in \cC_0^\infty(]0,T[ \times ]0,L[)$,
 \beqn
 W(\varphi) = \int_0^T \int_0^L \varphi(t,x)\, W(dt,dx)
 \eeqn
belongs to $H$, since the left-hand side of \eqref{5.1.4(*1)} is a Riemann integral, hence the $L^2(\Omega)$-limit of linear combinations of the $v(t,x)$.  Since $\cC_0^\infty(]0,T[ \times ]0,L[)$ is dense in $L^2([0,T] \times [0,L])$ (see \cite[Theorem 4.12, p. 57]{brezis}), the proof is complete.
 \end{proof}

\begin{def1}
\label{ch2'-s3.4.4-d2}
The {\em Reproducing Kernel Hilbert Space}\index{Reproducing Kernel Hilbert Space}\index{Hilbert Space!Reproducing Kernel}\index{kernel!Hilbert Space, Reproducing} {\em (RKHS)}\index{RKHS} of $v$ is the vector space $\hac$ of functions $f: [0,T]\times[0,L]\rightarrow \IR$ of the form
\beqn
f(t,x) = E\left[X v(t,x)\right], \quad  (t,x)\in[0,T]\times[0,L],
\eeqn
where $X\in H$. This space is equipped with the inner product $\langle f, g\rangle_\hac = E[XY]$ if $g \in \cH$ is given by $g(t,x)=E\left[Y v(t,x)\right]$.
\end{def1}

We notice that since the random field $v$ is $L^2(\Omega)$-continuous, every $f\in\hac$ is continuous on $[0,T]\times [0,L]$.

\begin{lemma}
\label{ch2'-s3.4.4-l2}
$f\in\hac$ if and only if there exists $\varphi\in L^2([0,T]\times[0,L])$ such that for all $(t,x) \in [0, T] \times [0, L]$,
\beq
\label{ch2'-s3.4.4.3}
f(t,x) = \int_0^t ds \int_0^L dy\, G_L(t-s;x,y) \varphi(s,y).
\eeq
This $\varphi$ is unique. For two such functions $f_1$ and $f_2$ (associated with $\varphi_1$ and $\varphi_2$, respectively),
\beq
\label{5.1.4(*2)}
     \langle f_1, f_2 \rangle_{\cH} = \int_0^T ds \int_0^L dy\, \varphi_1(s,y) \varphi_2(s,y)
     =\langle\varphi_1,\varphi_2\rangle_{L^2([0,T]\times[0,L])}.
     \eeq
\end{lemma}


\begin{proof}
By Definition \ref{ch2'-s3.4.4-d2}, $f\in\hac$ if and only if there is $X\in H$ such that $f(t,x)=E[X v(t,x)]$. By Lemma \ref{ch2'-s3.4.4-l1}, there is $\varphi\in L^2([0,T]\times[0,L])$ such that $X=W(\varphi)$. Therefore, by \eqref{ch2'-s3.4.4.2} with $u_0 \equiv 0$ there,
\beqn
f(t,x)= E[W(\varphi) v(t,x)] = \int_0^t ds \int_0^L dy\ G_L(t-s;x,y) \varphi(s,y).
\eeqn

In order to check uniqueness, suppose that we also have
\beqn
f(t,x) = \int_0^t ds \int_0^L dy\, G_L(t - s; x, y) \psi(s, y),
\eeqn
for some $\psi\in L^2([0,T]\times [0,L])$.
 Then
 \beqn
 f(t,x) = E[W(\varphi) v(t,x)] = E[W(\psi) v(t,x)],
 \eeqn
  so $E[W(\varphi - \psi))v(t,x)] = 0$, for all $(t,x) \in [0, T] \times [0, L]$, that is, $W(\varphi - \psi)$ is orthogonal to $H$. Since $W(\varphi - \psi)$ belongs to $H$ by Lemma \ref{ch2'-s3.4.4-l1}, we conclude that $W(\varphi - \psi) = 0$, or, equivalently, that $\Vert \varphi - \psi \Vert_{L^2([0,T] \times [0, L])} = 0$. This establishes the desired uniqueness.

If $g(t,x) = E[W(\psi) v(t,x)]$, then by Definition \ref{ch2'-s3.4.4-d2},
\beqn
    \langle f, g \rangle_\cH = E[W(\varphi) W(\psi)] = \int_0^T ds \int_0^L dy\ \varphi(s,y) \psi(s,y),
    \eeqn
proving \eqref{5.1.4(*2)} and completing the proof of the lemma.
\end{proof}
\medskip

\noindent{\em Germ-field Markov property of Gaussian random fields}
\medskip

We now state without proof a result on the germ-field Markov property of the Gaussian random field $v$. This is a particular case of a general result on centred Gaussian random fields  (see \cite[Theorem 3.3]{pitt-1971}, or \cite[Theorem 5.1]{kuensch}). In this statement, for a function $f: [0,T]\times[0,L] \rightarrow \IR$, \ ${\text{supp}} f$ denotes the closure of $\{(t,x): f(t,x)\ne 0\}$.

\begin{prop}
\label{ch2'-s3.4.4-p2}
The Gaussian random field $v$ has the germ-field Markov property if and only if its RKHS, denoted by $\hac$, satisfies the following two conditions:
\begin{description}
\item{(a)} if $f_1, f_2\in\hac$ are such that ${\text{supp}} f_1\cap {\text{supp}} f_2=\emptyset$, then $\langle f_1,f_2\rangle_\hac = 0$;
\item{(b)} if $f\in\hac$
is of the form $f=f_1+f_2$, where $f_i: [0,T]\times[0,L] \rightarrow \IR$, $i=1,2$, and  ${\text{supp}} f_1\cap {\text{supp}} f_2=\emptyset$, then $f_i\in \hac$, $i=1,2$.
\end{description}
\end{prop}
\medskip

For its further use in the proof of Proposition \ref{ch2'-s3.4.4-p1}, we define $\cC=\cC([0,T]\times [0,L])$. When $f\in \cC$ and $\varphi \in L^2([0,T] \times [0, L])$ are related by \eqref{ch2'-s3.4.4.3}, we call  $(f, \varphi)$ an {\em $\hac$-couple}. 
We prove in Proposition \ref{app3-s4-new-p1} that $(f, \varphi)$ is an {\em $\hac$-couple} if and only if the following condition holds:
\smallskip

\noindent $({\bf P})$ For all $t \in [0, T]$, for all $\psi \in C^{1,2}([0,t] \times [0, L])$ such that $\psi(\cdot, 0) =\psi(\cdot, L)=0$, we have
\begin{align}
\label{L1-(*1)}
   \int_0^L dx f(t, x) \psi(t,x) &= \int_0^t ds \int_0^L dx f(s, x) \cL^* \psi(s, x)\notag\\
   &\qquad  +  \int_0^t ds \int_0^L dx\, \psi(s,x) \varphi(s, x),
   \end{align}
where $\cL^*= \frac{\partial^2}{\partial x^2} + \frac{\partial}{\partial s}$ is the (opposite of the) adjoint\label{rd06_13p1} of the heat operator $\cL= \frac{\partial}{\partial t}- \frac{\partial^2}{\partial x^2}$.
\bigskip

\noindent{\em Proof of Proposition \ref{ch2'-s3.4.4-p1}.} It suffices to check conditions (a) and (b) of Proposition \ref{ch2'-s3.4.4-p2}.

Consider first $f \in \hac$ and $\varphi \in L^2([0,T] \times [0, L])$ such that $(f, \varphi)$ is an $\hac$-couple. Let $D = {\text{supp}} f$. We are going to show that $\varphi = 0$ a.e.~on $D^c \cap ([0,T] \times [0, L])$. Indeed, consider a smooth function $\psi$ with support in $D^c \cap ([0,T] \times [0, L])$  such that $\psi(\cdot, 0)=\psi(\cdot, L)=0$. Since ${\text{supp}}\, \cL^* \psi \subset {\text{supp}}\, \psi$ and ${\text{supp}}\,\psi \cap D = \emptyset$, by Property $({\bf P})$ above,
\beqn \int_0^T ds \int_0^L dx\, \psi(s,x) \varphi(s, x) = 0.
\eeqn
Since the set of such $\psi$ is dense in $L^2(D^c \cap ([0,T] \times [0, L]))$, we conclude that $\varphi = 0$ a.e.~on $D^c \cap ([0,T] \times [0, L])$.

For (a), let $f_1, f_2\in\hac$ be such that ${\text{supp}} f_1\cap {\text{supp}} f_2=\emptyset$. Then there is an open set $\mathcal{O}\subset [0,T]\times[0,L]$ such that ${\text{supp}} f_1\subset \mathcal{O}$ and ${\text{supp}} f_2\subset \bar{\mathcal{O}}^c$.

 For $i= 1, 2$, consider $\varphi_i \in L^2([0,T] \times [0, L])$ such that $(f_i, \varphi_i)$ is an $\hac$-couple.  Then $\varphi_ 1 = 0$ a.s.~on $({\text{supp}} f_1)^c$ and $\varphi_2 = 0$ a.s.~on $({\text{supp}} f_2)^c$, that is, $\varphi_ 1 \varphi_ 2 = 0 $ a.s.~on $({\text{supp}} f_1 \cap {\text{supp}} f_2)^c = [0,T] \times [0, L]$. By Lemma \ref{ch2'-s3.4.4-l2},
 \beqn
   \langle f_1, f_2 \rangle_{\hac} = \langle \varphi_1, \varphi_2 \rangle_{L^2([0,T] \times [0, L])} = 0,
   \eeqn
and (a) is proved.

Turning to (b), suppose that  $f\in\hac$ and $f=f_1+f_2$, where ${\text{supp}} f_1\cap {\text{supp}} f_2=\emptyset$.
 Let $D_1 := {\text{supp}} f_1$, $D_2 := {\text{supp}} f_2$, and $\mathcal{O}$ be an open set such that $D_1 \subset \mathcal{O}$ and $D_2 \in \bar {\mathcal{O}}^c$. Then $f_2\big|_{\bar{\mathcal{O}}}=0$, $f_1\big|_{\mathcal{O}^c}=0$, so
\beq
\label{ch2'-s3.4.4.5}
f\big|_{\mathcal{\bar O}} = f_1\big|_{\mathcal{\bar O}}\qquad {\text{and}}\qquad f\vert_{\mathcal{O}^c} = f_2\vert_{\mathcal{O}^c}.
\eeq
In particular, $f_1\big|_{\mathcal{\bar O}}$ is continuous. Since $f_1\big|_{\mathcal{O}\setminus D_1}=0$ and $f_1\big|_{\mathcal{O}^c}=0$, $f_1$ is in fact continuous on $[0,T]\times[0,L]$, and the same is true for $f_2$, that is, $f_1, f_2\in\cC$.

Further, consider any $(t,x)$ where $f(t,x)=0$. If $(t,x)\notin D_1\cup D_2$, then $f_1(t,x)=f_2(t,x)=0$.  If
$(t,x)\in D_1$, then $f_2(t,x)=0$, therefore $f_1(t,x)=f(t,x)-f_2(t,x)$ is also equal to $0$. Similarly, if
$(t,x)\in D_2$, then we also have $f_1(t,x)=f_2(t,x)=0$. We conclude from this that
${\text{supp}} f = D_1 \cup D_2$.

It remains to find, for $i = 1, 2$, a function $\varphi_i \in L^2([0,T] \times [0, L])$ such that $(f_i, \varphi_i$) is an $\hac$-couple. Let $\varphi \in L^2([0,T] \times [0, L])$ be such that $(f, \varphi)$ is an $\hac$-couple. We are going to check Property $({\bf P})$ of Proposition \ref{app3-s4-new-p1} for $f_1$ and $\varphi 1_{D_1}$ (the proof for $f_2$ and $\varphi 1_{D_2}$ is similar).

    Let $\phi \in \cC^{\infty}([0,T] \times [0, L])$ be such that $0 \leq \phi \leq 1$, $\phi \equiv 1$ on $D_1$, $\phi \equiv 0$ on $\mathcal{O}^c$ (see e.g. \cite[Section 21.3, p.~188]{gasquet} for a particular construction of such a $\phi$).

   Let $\psi \in \cC^{1, 2}([0,T] \times [0, L])$ be such that $\psi(\cdot, 0)=\psi(\cdot, L)=0$.  Then for $t \in ]0, T]$,
   \begin{align}
   \label{afegit (*1)}
   \int_0^L dx f_1(t,x) \psi(t,x) &= \int_0^L dx\, 1_{D_1}(t,x) f_1(t,x) \psi(t,x)\notag\\
      &= \int_0^L dx\, 1_{D_1}(t,x) f(t,x) \psi(t,x) \phi(t,x)\notag\\
      &= \int_0^L dx\,  f(t,x) \psi(t,x) \phi(t,x),
      \end{align}  
because $ f_1 = f$ on $\cO \supset D_1$, $\phi = 1$ on $D_1$, $f=0$ on  $D_1^c \cap \mathcal{O}$ and $\phi = 0$ on $D_1^c \cap \mathcal{O}^c$. Since $\psi \phi \in \cC^{1, 2}([0, T] \times [0, L])$, we can apply Property $({\bf P})$ for $f$ to see that this is equal to
\begin{align*}
   & \int_0^t ds \int_0^L dx\, f(s,x) \cL^*(\psi\phi)(s,x) + \int_0^t ds \int_0^L dx\, \psi(s,x) \phi(s,x) \varphi(s,x) \\
   &\qquad= \int_0^t ds \int_0^L dx\, 1_{D_1}(s,x) f(s,x) \cL^* (\psi\phi)(s,x)\\
   & \qquad\qquad+ \int_0^t ds \int_0^L dx\, 1_{\mathcal{O}}(s,x) \psi(s,x) \phi(s,x) \varphi(s,x).
   \end{align*}
Recall that $f = f_1$ and $\phi \equiv 1$ on $D_1$, so we can remove $\phi$ in the first integral, to obtain from \eqref{afegit (*1)} that
\begin{align*}
     \int_0^L f_1(t,x) \psi(t,x)\, dx &=  \int_0^t ds \int_0^L dx\, f_1(s, x) \cL^*\psi(s,x)\\
     & \qquad+ \int_0^t ds \int_0^L dx\, 1_{\mathcal{O}}(s, x) \psi(s, x) \phi(s, x) \varphi(s, x).
     \end{align*}
Replace $\mathcal{O}$ by a decreasing sequence of open sets $\mathcal{O}_n \supset D_1$, such that $D_2 \subset \bar \cO_n^c$, with intersection equal to $D_1$, and $\phi$ by $\phi_n$. Observe that $1_{\mathcal{O}_n}\psi\phi_n$ converges to $1_{D_1} \psi$ in $L^2([0,t] \times [0, L])$, therefore
\begin{align*}
   \int_0^L f_1(t,x) \psi(t,x) &=  \int_0^t ds \int_0^L dx\, f_1(s, x) \cL^*\psi(s,x)\\
   & \qquad + \int_0^t ds \int_0^L dx\, \psi(s, x) 1_{D_1}(s, x) \varphi(s, x).
   \end{align*}
We conclude from Property $({\bf P})$ that $(f_1, 1_{D_1} \varphi)$ are an $\hac$-couple.
\hfill\qed
\medskip

\noindent{\em The nonlinear case}
\medskip

We now address the case where $b\not\equiv 0$.

\begin{thm}
\label{ch2'-s3.4.4-t1}
Suppose that the function $b$ in \eqref{ch2'-s3.4.4.1} satisfies the assumptions of Theorem \ref{ch2'-s3.4.1-t1}, and the function $u_0$ in
\eqref{ch2'-s3.4.4.1} is Borel and bounded. Then any weak solution in law $u$ of \eqref{ch2'-s3.4.4.1} has the germ-field Markov property.
\end{thm}

\begin{proof}
Let $W$ be a space-time white noise defined on some stochastic basis $(\Omega,{\tf}, P, ({\tf}_t,\, t \in [0, T]))$.  Let $v$ be as in  \eqref{ch2'-s3.4.4.2} (we do not assume that $u_0\equiv 0$). Define a measure $\tilde P$ by
\begin{align*}
\frac{d\tilde P}{dP} &= \exp\left(\int_0^T \int_0^L b(t,x,v(t,x))\, W(dt,dx)\right.\\
&\left.\qquad\qquad \qquad-\frac{1}{2} \int_0^T dt \int_0^L dx\, b^2(t,x,v(t,x))\right).
\end{align*}
In the proof of Theorem \ref{ch2'-s3.4.1-t1}, we have seen that $\tilde P$ is a probability measure and that $(\Theta, \tilde W, v)$ is a weak solution in law of  \eqref{ch2'-s3.4.4.1}, where $\Theta$ is the filtered probability space $(\Omega, \cF, \tilde P, (\cF_t,\, t \in [0, T]))$. Therefore, it suffices to show that $v$ has the germ-field Markov property under $\tilde P$.

Let $A$ be an open subset of $[0,T]\times [0,L]$ with smooth boundary. Set $\cG_1 = \tf^v(A)$,  $\cG_2 = \mathcal{G}^v(A^c)$.
By Remark \ref{ch2'-s3.4.4-r0}, it suffices to prove that $\cG_1$ and  $\cG_2$ are conditionally independent given $\mathcal{G}^v(\partial A)$.

Let $X$ be a nonnegative $\cG_1$-measurable random variable
and denote $J := \frac{d\tilde P}{dP}.$ Recall that conditional expectations
 relative to $\tilde P$ are given by
\beq
\label{ch2'-s3.4.4.6}
E_{\tilde P}\left[X\,\big\mid\, \cG_2\right] = \frac{E_P\left[X J\,\big\mid\, \cG_2\right]}{E_P\left[J\,\big|\,\cG_2\right]}
\eeq
(this is Bayes' rule as given in \cite[Lemma 3.5.3, p.~193]{ks}).

We notice that $J=J_1J_2$, with
\begin{align*}
J_1&=\exp\left(\int_A b(t,x,v(t,x))\, W(dt,dx)-\frac{1}{2} \int_A b^2(t,x,v(t,x))\, dt dx\right),\\
J_2&=\exp\left(\int_{A^c} b(t,x,v(t,x))\, W(dt,dx)-\frac{1}{2} \int_{A^c} b^2(t,x,v(t,x))\, dt dx\right).
\end{align*}
Observe that because $A$ has a smooth boundary, the $dt dx$-measure of $\partial A$ is $0$, so the integrals in the definition of $J_2$ can be taken over the open set $\bar A^c$.

Apply Lemma \ref{5.1.20a} to $A$ and to $\bar A^c$, respectively, to deduce that $J_1$ is $\cG_1$-measurable, and $J_2$ is measurable with respect to $\tf^v(\bar A^c)$, a $\sigma$-field included in $\cG_2$. Therefore, by \eqref{ch2'-s3.4.4.6}

\beqn
E_{\tilde P}\left[X\, \big|\, \cG_2\,\right] = \frac{E_P\left[XJ_1J_2\,\big|\, \cG_2\right]}{E_P\left[J_1J_2\, \big|\, \cG_2\right]}
=  \frac{E_P\left[XJ_1\, \big|\, \cG_2\right]}{E_P\left[J_1 \, \big|\, \cG_2\right]},
\eeqn
and by the germ-field Markov property of $v$ under $P$ (Proposition \ref{ch2'-s3.4.4-p1}), the right-hand side is equal to
\beqn
\frac{E_P\left[XJ_1\, \big|\, \mathcal{G}^v(\partial A)\right]}{E_P\left[J_1\, \big|\, \mathcal{G}^v(\partial A)\right]}.
\eeqn
Hence, $E_{\tilde P}\left[X\big|\cG_2\right]$ is $\mathcal{G}^v(\partial A)$-measurable. Since
$\mathcal{G}^v(\partial A)\subset \mathcal{G}^v(A^c)$, we deduce that
\beqn
E_{\tilde P}\left[X\, \big|\, \cG_2\right] = E_{\tilde P}\left[X\, \big|\, \mathcal{G}^v(\partial A)\right],
\eeqn
that is, $v$ has the germ-field Markov property under $\tilde P$.
\end{proof}
\begin{remark}
\label{rem-before 5.1.21}
(a)\, If, instead of vanishing Dirichlet boundary conditions, we consider vanishing Neumann boundary conditions, then Proposition \ref{ch2'-s3.4.4-p1}, Lemma \ref{ch2'-s3.4.4-l2} and Theorem \ref{ch2'-s3.4.4-t1} remain valid, with essentially the same proof. This is a consequence of  Remark \ref{app3-s4-new-r1}.

 (b)\, Suppose that instead of \eqref{ch2'-s3.4.4.1}, we consider a nonlinear stochastic heat equation with {\em multiplicative noise.} It is an {\em open problem} to determine whether Theorem \ref{ch2'-s3.4.4-t1} remains true in this situation.

 (c)\, It is known that the solution to \eqref{ch2'-s3.4.4.1} does not satisfy the less-studied sharp Markov property. In fact, the sharp Markov property is only known to hold (for a reasonably large class of sets) for the (reduced) wave equation in spatial dimension $1$ or equivalently, the Brownian sheet (\cite{dalang-walsh-1992-1}, \cite{dalang-walsh-1992-2}), and the Whittle field and certain Bessel fields indexed by $\R^2$ (\cite{pitt-robeva-2003}).
 \end{remark}
 \medskip
 
 \noindent{\em Two technical lemmas}
 \medskip

We end this section with two technical lemmas that have been used in the proofs of Lemma \ref{ch2'-s3.4.4-l1} and Theorem \ref{ch2'-s3.4.4-t1}.
\begin{lemma}
\label{5.1.16a}
Let $v$ be as in \eqref{ch2'-s3.4.4.2}, let $\cL$ be the heat operator on $]0,T[ \times ]0,L[$ and let $ \cL^* = \frac{\partial}{\partial t} + \frac{\partial^2}{\partial x^2}$ be (the opposite of) its adjoint.
 \begin{description}
\item{(a)} Let $\varphi \in \cC_0^\infty(]0,T[ \times ]0,L[)$. Then
\beq
\label{5.1.4(*1)}
   \int_0^T dt \int_0^L dx\, \cL^\star \varphi (t,x) v(t,x) = \int_0^T \int_0^L \varphi (t,x)\, W(dt,dx).      
\eeq
 \item{(b)} Let $A \subset\, ]0,T[ \times ]0,L[$ be an open set. Let $h \in L^2(A)$, and extend $h$ to $]0,T[ \times ]0,L[$ by setting $h = 0$ on $A^c$. Then the random variable
    $(h \cdot W)_T = \int_0^T \int_0^L h(t,x)\, W(dt,dx)$
is $\F^v(A)$-measurable.
\end{description}
\end{lemma}

\begin{proof}
(a)\  On the left-hand side, we replace $v$ by the expression in \eqref{ch2'-s3.4.4.2} and apply the stochastic Fubini's theorem (Theorem \ref{ch1'-tfubini}), whose assumptions are clearly satisfied. We obtain
\begin{align*}
&\int_0^T dt \int_0^L dx\, \cL^*\varphi(t,x) \int_0^L G_L(t; x, y) u_0(y)\\
&\qquad +
   \int_0^T \int_0^L \left(\int_s^T dt \int_0^L dx\, \cL^\star \varphi (t,x) G_L(t-s; x,y)\right) W(ds,dy).
\end{align*}
Since the first term on the right-hand side of \eqref{ch2'-s3.4.4.2}, which we denote $I_0(t,x)$, satisfies $\cL I_0(t,x) = 0$, we use the boundary conditions at $x = 0$, $x = L$, $t = 0$, $t = T$ and integration by parts to replace  $\cL^*\varphi(t,x)  I_0(t,x)$  by $- \varphi(t,x)  \cL I_0(t,x) = 0$. For the second term, applying \eqref{truc},  
we see that the inner integral is equal to $\varphi (s,y)$, proving (a). 

   (b)\ If $\varphi \in \cC_0^\infty(]0,T[ \times ]0,L[)$ with closed support contained in $A$, then the closed support of $\cL^\star \varphi$ is also contained in $A$, so (a) implies that $(\varphi \cdot W)_T$ is $\F^v(A)$-measurable. Indeed, this is the case for the left-hand side  of \eqref{5.1.4(*1)} since $v$ has continuous sample paths, and therefore, the integral is a Riemann integral, which is an $L^2(\Omega)$-limit of linear combinations of the $v(t,x)$. The conclusion for general $h \in L^2(A)$ follows from the fact that $\cC_0^\infty(A)$ is dense in $L^2(A)$. This completes the proof.
   \end{proof}

\begin{lemma}
\label{5.1.20a}
Let $v$ be as in \eqref{ch2'-s3.4.4.2} and
 let $A$ be an open subset of $[0,T] \times [0,L]$ with smooth boundary. Then
    $\int_A b(t,x,v(t,x)) W(dt,dx)$
is $\cF^v(A)$-measurable.  
\end{lemma}
\begin{proof}
 We decompose $A$ into a countable union of closed rectangles with disjoint interiors but possibly overlapping boundaries. Then the integral over $A$ is the sum of the integrals over the rectangles, with convergence in $L^2(\Omega)$, so it suffices to prove the statement for each rectangle separately.

Let $R = [t_0, t_1] \times [a_1, a_2] \subset A$ be one of these rectangles. Let $h \in L^2(R)$, and extend $h$ to $R^c$ by setting $h = 0$ there. According to Lemma \ref{5.1.16a} (b) (applied to the interior of $R$), $(h \cdot W)_T$ is $\cF^v(R)$-measurable.

   Let $(v_i,\, i \geq 1)$ be an orthonormal basis of $V_R := L^2([a_1,a_2])$. For $t \in [t_0, t_1]$, let $\cG_t := \cF^v([t_0, t] \times [a_1, a_2])$, so that $\cG_t \subset \cF_t$. We complete $\cG_t$ with all $\cG_{t_1}$ null sets and make the filtration $(\cG_t,\,  t \in [t_0, t_1])$ right-continuous. Then for $i \geq 1$ and $t \in [t_0, t_1]$, by setting $h(s,x) := 1_{[t_0,t] \times [a_1,a_2]} (s,x)v_i(x)$, we see that the random variable
   \beqn
   U_{i,t} := W_t(v_i) - W_{t_0}(v_i) = (h\cdot W)_t
   \eeqn
is $\cG_t$-measurable. Further, $(U_{i,t},\, t \in [t_0, t_1],\, i\geq 1)$, is a sequence of $(\cG_t)$-Brownian motions with continuous sample paths, independent of $\cF_{t_0}$.

   Let $G$ be a jointly measurable and $(\cF_t)$-adapted stochastic process that satisfies \eqref{weakint} and which vanishes outside of $R$. Using the local property in space of the stochastic integral (Proposition \ref{moreproper}), we see that
   \begin{align*}
   Z := \int_{[t_0,t_1] \times [a_1,a_2]} G(t,x)\, W(dt,dx)
      = \sum_{i=1}^\infty \int_{t_0}^{t_1}  \langle G(t,*), v_i \rangle_{V_R}\, dU_{i,t}.
\end{align*}
   Consider the optional $\sigma$-field $\cO$ on $[t_0,t_1] \times \Omega$ associated with $(\cG_t)$. Suppose that in addition, $(x,t,\omega) \mapsto G(t,x,\omega)$ from $[a_1,a_2]\times [t_0,t_1] \times \Omega$ into $\R$ is $\cB_{[a_1,a_2]} \times \cO$-measurable.  Then for $i \geq 1$, $(Y_t(i) := \langle G(t,*), v_i \rangle_{V_R},\, t\in[t_0,t_1])$, defines an $\cO$-measurable process, and
     $Z_i := \int_{t_0}^{t_1} Y_t(i)  dU_{i,t}$
is an $L^2(\Omega)$-limit of stochastic integrals of simple processes. Each simple process $H= (H_t)$ is the finite sum of terms of the form $H_\ell 1_{]s_1^\ell,s_2^\ell]}(t)$, where $H_\ell$ is a $\cG_{s_1^\ell}$-measurable bounded random variable. In particular, the stochastic integral $(H \cdot U_{i,\cdot} )_{t_1}$ is $\cG_{t_1}$-measurable, hence $\cF^v(A)$-measurable. This implies that $Z_i,$ hence also $Z$, is $\cF^v(A)$-measurable.

   Now, for $(t,x) \in [0,T] \times [0,L]$, let $G(t,x) := b(t,x,v(t,x)) 1_R(t,x)$. Since $b$ is jointly measurable and $(x,t,\omega) \mapsto v(t,x,\omega)$
    defined on $[a_1, a_2] \times [t_0, t_1] \times \Omega$ is
  $\cB_{[a_1,a_2]} \times \cO$-measurable (since it is continuous and adapted to $(\cG_t))$, we see that $(x,t,\omega) \mapsto G(t,x,\omega)$ is also $\cB_{[a_1,a_2]} \times \cO$-measurable.
  For $i \geq 1$, let $(Y_t(i) := \langle G(t,*), v_i \rangle_{V_R},\, t \in [t_0, t_1])$. We conclude from the previous paragraph that
\beqn
    Z_i := \int_{t_0}^{t_1} Y_t(i) \, dU_{i,t} \quad\text{ and }\quad 
    Z := \int_{[t_0,t_1] \times [a_1,a_2]} b(t,x,v(t,x))\, W(dt,dx)
    \eeqn
are $\cF^v(A)$-measurable. This proves the lemma.
\end{proof}

\section{Asymptotic bounds on moments}
\label{ch2'-section5.2}

In this section, we consider the SPDE \eqref{ch1'-s5.0} with the formulation \eqref{ch1'-s5.1}, and the notion of solution given in Definition \ref{def4.1.1}.
The objective is to establish exponential upper bounds on $L^p$-moments of the random field solution $(u(t,x),\, (t,x)\in\re_+\times D)$ that are uniform in $x$ and highlight the explicit dependence on $t$ and $p$.
This is an important ingredient in the study of long-term physical phenomena, such as intermitency (see the notes in Section \ref{notes-ch5}).

\subsection{Main results}
\label{ss5.3-1}

In the formulation of these results, we will use the notion of random field solution $(u(t,x),\, (t,x)\in\re_+\times D)$ to \eqref{ch1'-s5.0} {\em for all time}, which means that Theorem  \ref{ch1'-s5.t1}
on existence and uniqueness of a random field solution to \eqref{ch1'-s5.0}   holds for any $T>0$. To ensure the existence of this object, a slight modification of the assumptions of Section \ref{ch1'-ss5.1} is required, as follows:
\begin{enumerate}
\item  (${\bf H_{I,\infty}}$)\label{rdHIinfinity} and (${\bf H_{L,\infty}}$)\label{rdHLinfinity} denote the assumptions on the initial conditions, and the coefficients $\sigma$ and $b$, respectively. They are the analogue of (${\bf H_I}$) and (${\bf H_L}$), with $[0,T]$ there replaced by $\re_+$.
\item (${\bf H_{\Gamma,\infty}}$).\label{rdHGammainfinity} This is a global-in-time version of $({\bf H_\Gamma})$ formulated as follows:
\begin{description}
\item{$(i^\prime)$} The fundamental solution/Green's function $\Gamma(t,x;s,y)$
 is a jointly measurable mapping from $\{(t,x;s,y)\in\re_+\times D\times\re_+\times D: 0\le s<t < \infty\}$ into $\IR$.
\item{$(ii^\prime)$} There is a Borel function $H: \re_+\times D^2 \longrightarrow \IR_+$ such that
\beqn
\vert\Gamma(t,x;s,y)\vert \le H(t-s,x,y), \qquad 0\le s<t<\infty, \ x,y\in D.
\eeqn
\item{$(iii_a^\prime)$} If in \eqref{ch1'-s5.0}$\  \sigma \not\equiv 0$, then for any $T\ge 0$,
\beqn
\int_0^T ds \ \sup_{x\in D}\int_D dy\, H^2(s,x,y) < \infty.
\eeqn
\item{$(iii_b^\prime)$} If in \eqref{ch1'-s5.0}$\  b \not\equiv 0$, then for any $T\ge 0$,
\beqn
\int_0^T ds \ \sup_{x\in D}\int_D dy\, H(s,x,y) < \infty.
\eeqn
\end{description}
\end{enumerate}

An easy adaptation of the proof of Theorem \ref{ch1'-s5.t1} gives the following.
\begin{thm}
\label{ch1'-s5.t1-global}
Under $(\bf H_{\Gamma,\infty})$, $(\bf H_{I,\infty})$ and  $(\bf H_{L,\infty})$, there exists a random field solution $$u=\left(u(t,x),\, (t,x)\in\re_+\times D\right)$$ to \eqref{ch1'-s5.0} (with $(t,x)\in\re_+\times D$ there). In addition, for any $T>0$ and any $p\ge 2$,
\beqn
\sup_{(t,x)\in[0,T]\times D} E\left(\vert u(t,x)\vert^p\right)<\infty,
\eeqn
and the solution $u$ is unique (in the sense of versions) among random field solutions that satisfy this property with $p=2$.
\end{thm}

Throughout this section, we will also use the following (new) hypothesis
on uniform estimates of integrals of the fundamental solution/Green's function.
\begin{description}
\item{3.} Assumption $({\bf H_{\Gamma-{\text{sup}},\infty}})$.\label{rdHGsinfinity}
Let $H$ be the function considered in $({\bf H_{\Gamma,\infty}})$. There
 exists a finite constant $C$ and real numbers $a_1, a_2>0$ such that, for any $\beta>0$ large enough,
\beq
\label{ch5-s5.5.2-global}
\sup_{x\in D}\int_0^\infty ds\ e^{-\beta s} \int_D dy\ H(s,x,y) \le C\beta^{-a_1}
\eeq
and
\beq
\label{ch5-s5.5.4-global}
\sup_{x\in D}\int_0^\infty ds\  e^{- 2\beta s} \int_D dy\ H^2(s,x,y) \le C^2 \beta^{-2a_2}.
 \eeq
\end{description}

Given a random field $(X(t,x),\, (t,x)\in\IR_+\times D)$ and real numbers $\beta>0$, $p\in[2,\infty[$, set
\beq
\label{ch5-s5.5.6-infty}
\mathcal{N}_{\beta,p,\infty}(X) = \sup_{(t,x)\in\IR_+\times D} e^{-\beta t} \Vert X(t,x)\Vert_{L^p(\Omega)}.
\eeq

For $g \in \{\sigma, b\}$, define
\begin{align}
\label{ch5-s5.5.7-bis}
c(g) & = \sup_{(s,y,\omega)\in \IR_+\times D\times \Omega}|g(s,y,0;\omega)|,\notag\\
L(g) & =
\sup_{\substack{(s,y,\omega)\in \IR_+\times D\times \Omega \\ z_1,\, z_2\in D,\, z_1\ne z_2}}\frac{|g(s,y,z_1;\omega)- g(s,y,z_2;\omega)|}{|z_1-z_2|}.
\end{align}
By $({\bf H_{L,\infty}})$, $c(g)$ and $L(g)$ are finite, and
\beq
\sup_{(s,y,\omega)\in \re_+\times D\times \Omega}|g(s,y,z;\omega)|\le c(g) + L(g)\, |z|.
\label{ch5-s5.5.8}
\eeq

\begin{thm}
\label{ch5-s5.5-t1}
Assume $({\bf H_{\Gamma,\infty}})$, $({\bf H_{I,\infty}})$, $({\bf H_{L,\infty}})$ and   $({\bf H_{\Gamma-{\text{sup}},\infty}})$. Then there exists a constant $K \in \R_+$ such that for all $p\in[2,\infty[$ and for all $t\in \re_+$,
\beq
\label{ch5-s5.5.5-bis}
\sup_{x\in D}E\left[|u(t,x)|^p\right] \le \left(2\Vert I_0\Vert_{\infty}+ 1\right)^p\
\exp\left(\left(K^{\frac{1}{a_1 \wedge a_2}} p^{\frac{1}{2a_2}+1}\right)t\right),
\eeq
 where $\Vert I_0\Vert_{\infty}= \sup_{(t,x)\in\re_+\times D}|I_0(t,x)|$.
\end{thm}
\begin{proof}
Fix $p \geq 2$ and take $\beta>0$ large enough so that \eqref{ch5-s5.5.2-global} and \eqref{ch5-s5.5.4-global} hold. By Assumption $({\bf H_{I.\infty}})$,
\begin{align}
\label{ch5-s5.5.9}
\mathcal{N}_{\beta,p,\infty}(I_0) &= \sup_{(t,x)\in\re_+\times D} e^{-\beta t}|I_0(t,x)|
\le \sup_{(t,x)\in\re_+\times D} |I_0(t,x)|\notag\\
& = \Vert I_0\Vert_{\infty}<\infty.
\end{align}

Set
\beqn
\mathcal{J}(t,x) = \int_0^t ds \int_{D} dy\ \Gamma(t,x;s,y) b(s,y,u(s,y)).
\eeqn
 Applying Minkowski's inequality and using the condition (ii')
of Assumption $({\bf H_{\Gamma,\infty}})$ and \eqref{ch5-s5.5.8} with $g:=b$, we obtain
\begin{align*}
\Vert \mathcal{J}(t,x)\Vert_{L^p(\Omega)} &\le \int_0^t ds \int_D dy \ H(t-s,x,y)\Vert b(s,y,u(s,y))\Vert_{L^p(\Omega)}\\
&\le \int_0^t ds \int_D dy \ H(t-s,x,y) \left( c(b) + L(b)\Vert u(s,y)\Vert_{L^p(\Omega)}\right).
\end{align*}
This implies
\beqn
\mathcal{N}_{\beta,p,\infty}(\mathcal{J}) = \sup_{(t,x)\in\re_+\times D} e^{-\beta t}\Vert \mathcal{J}(t,x)\Vert_{L^p(\Omega)} \le T_1^{(1)} + T_2^{(1)},
\eeqn
where
\begin{align*}
T_1^{(1)}&= c(b)  \sup_{(t,x)\in\re_+\times D} e^{-\beta t} \int_0^t ds \int_D dy \, H(t-s,x,y),\\
T_2^{(1)} &= L(b)  \sup_{(t,x)\in\re_+\times D} e^{-\beta t} \int_0^t ds \int_D dy \, H(t-s,x,y)\Vert u(s,y)\Vert_{L^p(\Omega)}.
\end{align*}

By \eqref{ch5-s5.5.2-global},
\begin{align*}
T_1^{(1)} & \le c(b)  \sup_{(t,x)\in\re_+\times D} \int_0^t ds\ e^{-\beta (t-s)}\int_D dy \, H(t-s,x,y)\\
&\le C\, c(b) \beta^{-a_1},
\end{align*}
and
\begin{align*}
T_2^{(1)} &\le L(b) \sup_{(t,x)\in\re_+\times D} \int_0^t ds\ e^{-\beta(t-s)} \int_D dy\, H(t-s,x,y)\ e^{-\beta s} \Vert u(s,y)\Vert_{L^p(\Omega)}\\
&\le L(b)\ \mathcal{N}_{\beta,p,\infty}(u)\ \sup_{(t,x)\in\re_+\times D} \int_0^t ds\ e^{-\beta(t-s)} \int_D dy\, H(t-s,x,y)\\
& \le C\, L(b)\ \mathcal{N}_{\beta,p,\infty}(u) \beta^{-a_1}.
\end{align*}
In conclusion,
\beq
\label{ch5-s5.5.10}
\mathcal{N}_{\beta,p,\infty}(\mathcal{J}) \le C\beta^{-a_1} \left(c(b) + L(b)\ \mathcal{N}_{\beta,p,\infty}(u) \right).
\eeq

Let
\beqn
\mathcal{I}(t,x) = \int_0^t \int_D \Gamma(t-s;x,y)\, \sigma(s,y,u(s,y))\, W(ds,dy).
\eeqn
Apply first Burkholder's inequality \eqref{ch1'-s4.4.mod} and then Minkowski's inequality. Using \eqref{ch5-s5.5.8} with $g:=\sigma$ we have
\begin{align*}
\Vert \mathcal{I}(t,x)\Vert_{L^p(\Omega)}^2& \le 4p \left\Vert \int_0^t ds \int_D dy\, H^2(t-s,x,y) \sigma^2(s,y,u(s,y))\right\Vert_{L^{\frac{p}{2}}(\Omega)}\\
&\le 8p  \int_0^tds\int_D dy \, H^2(t-s,x,y) \left(c(\sigma)^2 + L(\sigma)^2 \Vert u(s,y\Vert_{L^p(\Omega)}^2\right).
\end{align*}
Therefore,
\beqn
\mathcal{N}_{\beta,p,\infty}(\mathcal{I}) = \sup_{(t,x)\in\re_+\times D} e^{-\beta t}\Vert \mathcal{I}(t,x)\Vert_{L^p(\Omega)}
 \le \left(8p\left(T_1^{(2)} + T_2^{(2)}\right)\right)^\half,
\eeqn
where
\begin{align*}
T_1^{(2)} & = c(\sigma)^2 \sup_{(t,x)\in\re_+\times D} e^{-2\beta t} \int_0^t ds \int_D dy\, H^2(t-s,x,y),\\
 T_2^{(2)} & = L(\sigma)^2 \sup_{(t,x)\in\re_+\times D}e^{-2\beta t}\int_0^t ds \int_D dy\, H^2(t-s,x,y) \Vert u(s,y\Vert_{L^p(\Omega)}^2.
 \end{align*}

 Proceeding as for $T_1^{(1)}$ and $T_2^{(1)}$, but using \eqref{ch5-s5.5.4-global}, we obtain
 \begin{align*}
 T_1^{(2)} &\le C^2 c(\sigma)^2 \beta^{-2 a_2},\\
 T_2^{(2)} &\le  L(\sigma)^2 \left(\mathcal{N}_{\beta,p,\infty}(u)\right)^2  \\
 &\qquad\quad \times \sup_{(t,x)\in\re_+\times D} \int_0^t ds\, e^{-2\beta (t-s)} \int_D dy\, H^2(t-s,x,y)\\
  &\le C^2 L(\sigma)^2 \left(\mathcal{N}_{\beta,p,\infty}(u)\right)^2  \beta^{-2a_2}.
  \end{align*}
  This yields
  \beq
  \label{ch5-s5.5.11}
  \mathcal{N}_{\beta,p,\infty}(\mathcal{I}) \le \sqrt{8p}\ C \beta^{-a_2} \left(c(\sigma) + L(\sigma) 
   \mathcal{N}_{\beta,p,\infty}(u) \right).
  \eeq

  Set
  \beqn
  L(b,\sigma) = \max\left(c(b), c(\sigma), L(b), L(\sigma)\right).
  \eeqn
 Recalling that $u(t,x) = I_0(t,x)+\mathcal{J}(t,x) + \mathcal{I}(t,x)$, from \eqref{ch5-s5.5.9}, \eqref{ch5-s5.5.10} and \eqref{ch5-s5.5.11}, we see that for some finite constants $C_{1}$, $C_{2}$,
  \begin{align}
  \label{ch5-s5.5.12}
&\mathcal{N}_{\beta,p,\infty}(u)\notag\\
&\quad \le \Vert I_0\Vert_{\infty} + C_{1} L(b,\sigma) \left[\frac{1}{\beta^{a_1}} + \frac{\sqrt p }{\beta^{a_2}}\right]
+ C_{2} L(b,\sigma)\left[\frac{1}{\beta^{a_1}} + \frac{\sqrt{p}}{\beta^{a_2}}\right] \mathcal{N}_{\beta,p,\infty}(u)\notag\\
&\quad\leq \Vert I_0\Vert_{\infty} +\bar L(b,\sigma) \max\left(\frac{1}{\beta^{a_1}}, \frac{\sqrt p }{\beta^{a_2}}\right) \left(1 + \mathcal{N}_{\beta,p,\infty}(u)  \right),  
\end{align}
where  $\bar L(b,\sigma)= \max(C_{1},C_{2}) L(b,\sigma)$. Observe that this inequality holds for any $p\in[2,\infty[$ and large enough $\beta > 0$, and that the constant $\bar L(b,\sigma)$ does not depend on $\beta$ or $p$.

Choose a constant $K$ large enough, depending on $\bar L(b,\sigma)$, such that, by defining $\beta:=K^{\frac{1}{a_1}\vee\frac{1}{a_2}} p^{\frac{1}{2a_2}}$, \eqref{ch5-s5.5.2-global} and \eqref{ch5-s5.5.4-global} hold and we have
 \beq
 \label{ch5-s5.5.13}
\bar L(b,\sigma)\max\left(\frac{1}{\beta^{a_1}}, \frac{\sqrt{p}}{\beta^{a_2}}\right)\le \half.
\eeq
With this choice of $\beta$, inequality  \eqref{ch5-s5.5.12} implies that
\beq
 \label{ch5-s5.5.14}
 \mathcal{N}_{\beta,p,\infty}(u) \le 2\Vert I_0\Vert_{\infty} + 1.
 \eeq
 Using the definition of  $\mathcal{N}_{\beta,p,\infty}(u)$, this implies that for some constant $K$ (depending in particular on $\bar L(b,\sigma)$, which in turn depends on the constants
 $c(b)$, $c(\sigma)$, $L(b)$ and  $L(\sigma)$), we have

 \beq
  \label{ch5-s5.5.15}
 \sup_{x\in D} E\left(|u(t,x)|^p\right) \le  \left(2\Vert I_0\Vert_{\infty} + 1\right)^p \exp\left(\left(K^{\frac{1}{a_1\wedge a_2}}p^{\frac{1}{2a_2}+1}\right) t\right),
 \eeq
 for all $t\in\re_+$. This proves \eqref{ch5-s5.5.5-bis}.
\smallskip
\end{proof}
\medskip

\begin{remark}
\label{s-5.3-r1}
Let $T>0$ be fixed. Assume $({\bf H_\Gamma})$, $({\bf H_I})$, $({\bf H_L})$ and furthermore, that conditions \eqref{ch5-s5.5.2-global} and \eqref{ch5-s5.5.4-global} hold with
$\re_+$, $\infty$ and $C$ there replaced by $[0,T]$, $T$ and $C_T$, respectively. Then the same approach as in the proof of Theorem \ref{ch5-s5.5-t1} shows that the solution $u=(u(t,x),\, (t,x)\in[0,T]\times D)$ to \eqref{ch1'-s5.1} satisfies the following property:

 There exists a positive and finite constant $K_T$  (depending on $T$), such that, for all $p\in[2,\infty[$ and $t\in [0,T]$,
\beq
\label{ch5-s5.5.5}
\sup_{x\in D}E\left[|u(t,x)|^p\right] \le \left(2\Vert I_0\Vert_{T,\infty}+ 1\right)^p\,
\exp\left(\left(K_T^{\frac{1}{a_1\wedge a_2}} p^{\frac{1}{2a_2}+1}\right) t\right),
\eeq
where $\Vert I_0\Vert_{T,\infty}= \sup_{(t,x)\in[0,T]\times D}|I_0(t,x)|$.

The assumptions of Theorem \ref{ch5-s5.5-t1} give a version of \eqref{ch5-s5.5.5} that is global in time and with a constant $K$ that does not depend on $T$.
\end{remark}

Property \eqref{ch5-s5.5.5-bis} is related to the notion of {\em Lyapounov exponent,} which is important in stability theory of stochastic evolutions.
 Define the {\em $p$-th moment upper Lyapounov exponent}\index{moment Lyapounov exponent}\index{Lyapounov exponent!moment} of the random field $(u(t,x),\, (t,x)\in\re_+\times D)$ by
\beq
\label{u-Lyapounov}
\bar\gamma_p(x) = \limsup_{t\to\infty}\frac{1}{t}\log E\left[u(t,x)|^p\right].
\eeq
If \eqref{ch5-s5.5.5-bis} holds, then $\sup_{x\in D} \bar\gamma_p(x) =O(p^{\frac{1}{2a_2}+1})$, as $p\to \infty$. In particular,  $\sup_{x\in D} \bar\gamma_p(x)<\infty$ for any $p\in [2,\infty[$. This notion is further discussed in Section \ref{rd1+1anderson}.

\subsection{Examples}
\label{ss5.3-ex1}

The following examples of SPDEs have been considered in Chapter \ref{ch1'-s5}. Since we plan to give statements on the solution for all time, we assume that $({\bf H_{I,\infty}})$ and $({\bf H_{L,\infty}})$ are satisfied. From the results of Section \ref{ch1'-s6}, we already know that the corresponding fundamental solutions satisfy $({\bf H_{\Gamma,\infty}})$. Here, we verify the validity of $({\bf H_{\Gamma-{\text{sup}},\infty}})$, and find the values of the parameters $a_1$ and $a_2$.
 \medskip

\noindent 1. {\em The stochastic heat equation on $\IR$ or $]0,L[$}
\smallskip

In this example, $D=\re$ or $D=\, ]0,L[$ with vanishing Dirichlet or Neumann boundary conditions. In
\eqref{ch1'-s6.1(*1)} and \eqref{ch1'-s6.1(*2)}, we have seen that in the first two cases, if we take
\beqn
H(s,x,y) = \frac{1}{\sqrt{4\pi s}} \exp\left(-\frac{(x-y)^2}{4s}\right)1_{]0,\infty[} (s),\quad x,y\in\re,
\eeqn
then $({\bf H_{\Gamma,\infty}})$ is satisfied, and $\int_0^t ds \int_{\IR} dy\, H(s,x,y) = t$.

 Computing the integrals, we see that
\beq
\label{ex-new(*3)}
\int_0^\infty ds\, e^{-\beta s} \int_\re dy\, H(s,x,y) =  \int_0^{\infty} e^{-\beta s} ds
=\beta^{-1}.
\eeq

As in the calculations that lead to \eqref{heatcauchy-11'}, we see that
\beq
\label{ex-new(*4)}
\int_0^\infty ds\ e^{-2\beta s}\int_D dy\, H^2(s,x,y) \le \int_0^\infty \frac{e^{-2\beta s}}{\sqrt{8\pi s}}\, ds
\le C \beta^{-\frac{1}{2}},
\eeq
where we have used the change of variables $r=2\beta s$ and the definition of the Euler Gamma function $\Gamma_{\text E}$ (see \eqref{Euler-gamma}).
We deduce that $({\bf H_{\Gamma-{\text{sup}},\infty}})$ holds with $a_1=1$ and $a_2 = \frac{1}{4}$.

 In the case of Neumann boundary conditions, we let $H(s, x, y)$ be the function on the right-hand side of \eqref{neumann-bad} to check that $({\bf H_{\Gamma,\infty}})$ is satisfied. In order to bound the expressions in \eqref{ch5-s5.5.2-global} and \eqref{ch5-s5.5.4-global}, we split the integrals from $0$ to $\infty$ into integrals from $0$ to $1$ and $1$ to $\infty$. With the change of variables $r = 2 \beta s$, we see that for large $\beta$, the first integral dominates and is bounded above by the expressions in \eqref{ex-new(*3)}  and \eqref{ex-new(*4)}. We deduce that $({\bf H_{\Gamma-sup,\infty}})$ holds with $a_1 = 1$ and $a_2 = \frac{1}{4}$.
 \smallskip

According to Theorem \ref{ch5-s5.5-t1}, we have the following.
\begin{cor}
\label{5.3-cor1}
Let $u = (u(t,x),\, (t,x) \in \R_+ \times D)$ be the solution for all time of one of the nonlinear stochastic heat equations considered in Theorem \ref{recap}. Then under assumptions $({\bf H_{I,\infty}})$ and $({\bf H_{L,\infty}})$, there are constants $C, \bar K \in \R_+$ such that for all $t \in \R_+$,
\beqn
 \sup_{x\in D} E\left(|u(t,x)|^p\right) \le C \exp(\bar K p^3 t).
 \eeqn
\end{cor}
 \medskip

 \noindent 2. {\em The stochastic wave equation on $\IR$, $\re_+$ or $]0,L[$}
 \smallskip

In the study of this example, we will use the equality
\beq
\sup_{t\in\IR_+}\left(\int_0^t
 s e^{-\beta s}\, ds\right) = \beta^{-2}.  \label{ch5-s5.5.18}
\eeq
We consider the stochastic wave equation on $\re$, $\re_+$ or $]0,L[$  (with Dirichlet boundary conditions in the last two cases), as in Section \ref{ch1'-ss6.2}. When $D = \re$, we can take
\beqn
H(s,x,y) = \frac{1}{2}\ 1_{\{|x-y|\le s\}}, \quad x,y\in \re,\ t\in\re_+.
\eeqn
We have also seen in Section \ref{ch1'-ss6.2} that $({\bf H_{\Gamma,\infty}})$ is satisfied, and
\beqn
\int_{\IR} dy\, H(s,x,y)= 2  \int_{\IR} dy\, H^2(s,x,y) = s .
\eeqn
Applying  \eqref{ch5-s5.5.18}, we see that $({\bf H_{\Gamma-{\text{sup}},\infty}})$
is satisfied with $a_1 = 2$ and  $a_2 = 1$.

When $D = \re_+$, we can check using the formula for the function $H$ satisfying $({\bf H_\Gamma})$ in Section \ref{ch1'-ss6.2} that $H$ can be taken such that for all $s > 0$,
 $ \int_{\re_+} dy\, H(s,x,y) \leq s$.
Therefore, in this case as well, $({\bf H_{\Gamma-{\text{sup}},\infty}})$ is satisfied with $a_1 = 2$ and $a_2 = 1$.
  When $D =\, ]0,L[$, we have seen in Section \ref{ch1'-ss6.2} that the function $H$ can be taken such that for all $s >0$ and $x\in[0,L]$,
  $ \int_0^L dy\, H(s,x,y) \leq C$   and  $ \int_0^L dy\, H^2(s,x,y)  \leq C'$.
Therefore, $({\bf H_{\Gamma-{\text{sup}},\infty}})$  is satisfied with $a_1 = 1$ and $a_2 = \half$. Theorem \ref{ch5-s5.5-t1} implies the following.

\begin{cor}
\label{5.3-c2}
Let $u = (u(t,x),\, (t,x) \in \re_+ \times D)$ be the solution for all time of one of the  stochastic wave equations considered in Theorem \ref{recapwave}. Then under assumptions $({\bf H_{I,\infty}})$ and $({\bf H_{L,\infty}})$, there are constants $C, \bar K \in \re_+$ such that for all $t \in \re_+$ :
\begin{enumerate}
  \item If $D = \re$ or $\re_+$, then
  \beqn
 \sup_{x\in D} E\left(|u(t,x)|^p\right) \le C \exp(\bar K p^{\frac{3}{2}} t).
 \eeqn
  \item If $D =\, ]0,L[$ , then
\beqn
 \sup_{x\in D} E\left(|u(t,x)|^p\right) \le C \exp(\bar K p^{2} t).
 \eeqn
 \end{enumerate}
 \end{cor}
\medskip

 \noindent 3. {\em The fractional stochastic heat equation on $\IR$}
\smallskip

We consider the fractional stochastic heat equation as in Section \ref{ch1'-ss7.3}, in which we can take
\beqn
 H(s,x,y)= \null_\delta G_a(s,x-y),\quad x,y\in\re,\ s>0,
 \eeqn
 with $\null_\delta G_a(t,z)$ given in \eqref{ch1'-ss7.3.3} ($a \in\, ]1, 2[$ and $\vert \delta \vert \leq 2-a$). Therefore, $({\bf H_{\Gamma, \infty}})$ is satisfied (see Lemma \ref{check-HGamma}).

Because of  \eqref{ch1'-ss7.3.4}, we see that the condition \eqref{ch5-s5.5.2-global} holds with $a_1=1$. For the verification of \eqref{ch5-s5.5.4-global}, we use the semigroup property of $\null_\delta G_a(t,z)$ given in \eqref{ch1'-ss7.3.5}, along with \eqref{ch1'-ss7.3.6}, to deduce that
\begin{align*}
\int_0^t ds\, e^{-2\beta s} \int_{\IR} dy\, H^2(s,x,y)
&= \int_0^t ds\, e^{-2\beta s} \int_{\IR} dy\, \null_\delta G_a(s,x-y)^2\\
&=\int_0^t ds\, e^{-2\beta s} \null_\delta G_a(2s,0)
=  C \int_0^t ds\, e^{-2\beta s} s^{-\frac{1}{a}}\\
& = C \beta^{-(1-\frac{1}{a})}\int_0^{2\beta t} e^{-s}\ s^{-\frac{1}{a}}\, dr.
 \end{align*}
 Using the definition of the Euler Gamma function $\Gamma_E$ (see \eqref{Euler-gamma}), we obtain \eqref{ch5-s5.5.4-global}, with $a_2= \frac{a-1}{2a}$. Thus, applying Theorem
 \ref{ch5-s5.5-t1} we obtain the following.
 \begin{cor}
 \label{5.3-c3}
 Let $u = (u(t,x),\, (t,x) \in \re_+ \times \re)$ be the solution for all time of \eqref{ch1'-ss7.3.9} considered in Theorem \ref{ch1'-ss7.3-t1}. Then under assumptions $({\bf H_{I, \infty}})$ and $({\bf H_{L, \infty}})$, there are constants $C, \bar K \in \re_+$ such that for all $t \in \re_+$,
\beqn
 \sup_{x\in D} E\left(|u(t,x)|^p\right) \le C \exp(\bar K p^{\frac{2a-1}{a-1}} t).
 \eeqn
 \end{cor}


\section{Some comparison theorems for the stochastic heat equation}
\label{ch5-added-0}

Comparison theorems\index{comparison theorem}\index{theorem!comparison} for PDEs and SPDEs refer to monotonicity properties of the solution with respect of some of their defining elements, such as the initial value or the coefficients. Comparison theorems relative to the initial value can be used to study the positivity of the solution
(see \cite{mueller-1991}, \cite{kotolenez-1992}), while comparison theorems relative to the drift coefficient can be used as a tool for implementing certain variational methods and establishing
 the existence of solutions (see e.g. \cite{gp2}, \cite{donati-martin-pardoux-1993}). In this section, we present a pathwise comparison theorem relative to both the initial value and the drift coefficient. First, we study the case of a nonlinear stochastic heat equation on a bounded interval and prove an extension of \cite[Theorem 2.1]{donati-martin-pardoux-1993}: see Theorem \ref{ch5-added-0-t1}. Then we address the same problem for the stochastic heat equation on $\re$, using a method different from that of \cite{shiga-1994}: we deduce the result on $\R$ from Theorems \ref{rd08_09t1} and \ref{ch5-added-0-t1}.

\subsection{The case $D = \,  ]0, L[$}
 
For $i = 1, 2$, let $\sigma(t,x,z)$ and $b_i(t, x, z)$ be two functions satisfying $({\bf H_L})$ and let $u_{0, i}$ be a function satisfying $({\bf H_I})$. We consider the stochastic heat equation on $D =\, ]0, L[$ with vanishing  Dirichlet (or Neumann: see Remark \ref{ch5-added-0-r1}) boundary conditions and initial condition $u_{0, i,L}$, as in Section \ref{ch1'-ss6.1}:
\beq
\label{ch5-added-0(*1)}
 \left(\frac{\partial}{\partial t}- \frac{\partial^{2}}{\partial x^{2}}\right)u_{i,L}(t,x) = \sigma(t,x,u_{i,L}(t,x))\dot W(t, x)
 +b_i(t,x,u_{i,L}(t,x)),
\eeq
 $(t, x) \in\, ]0,\infty[ \times\, D$, along with the boundary and initial conditions
 \begin{align*}
 \begin{cases}
 u_{i, L}(0, x) = u_{0,i,L}(x), & {\text{if}}\ x \in D,\\
u_{i, L}(t, 0) = u_{i, L}(t, L) =0, &{\text{if}}\ t \in \, ]0,\infty[.
\end{cases}
\end{align*}

 In the proof of the next theorem, we use the notations introduced in Section \ref{ch1'-s60}, and we denote $u_{1,L}$ and $u_{2,L}$  the solutions given by Theorem \ref{recap}.

\begin{thm}
\label{ch5-added-0-t1}
Suppose that for each $z\in \re$, $b_1(\cdot, \ast, z) \leq b_2(\cdot, \ast, z)$ $dsdxdP$-a.e.~and $u_{0,1,L} \leq u_{0,2,L}$ a.e. Then a.s., for all $(t,x) \in [0, T] \times D$,
    $u_{1,L}(t,x) \leq u_{2,L}(t,x)$.
\end{thm}

\begin{proof}
For simplicity of notation, we remove the index $L$ from $u_{0,i,L}$ and $u_{i,L}$.

The proof will consist of several steps. We begin with an approximation as in Section \ref{ch1'-s60}.
\medskip

\noindent{\em Step 1.~Approximation by finite-dimensional projections.} Consider the CONS 
\beqn
   \left(e_k(x) = \sqrt{\frac{2}{L}}\, \sin\left(\frac{k\pi}{L} x\right),\, k\ge 1\right)
\eeqn 
of $V=L^2(D)$. Let $W^n$ be the noise defined in \eqref{proj-noise} and let $\bar u_{n, i}$ be the approximation by finite-dimensional projections that satisfies the stochastic heat equation \eqref{ch1'-s60.2} with $\sigma$, $b_i$, and $I_0$ associated to $u_{0,i}$. Since the assumptions of Theorem \ref{ch1'-s60-t2} are satisfied, and because $u_i$ and $\bar u_{n, i}$ have continuous sample paths, in order to prove the theorem, it suffices by \eqref{ch1'-s60.4} to show that for all $(t,x) \in [0, T] \times D$,
\beq
\label{ch5-added-0(*2)} 
    \bar u_{n, 1}(t,x) \leq \bar u_{n, 2}(t,x)\quad     a.s.
\eeq

We consider first \eqref{ch5-added-0(*1)} without the index $i$, and denote the solution $u(t,x)$. Its approximation by finite-dimensional projections is denoted by $\bar u$ (where we omit the index $n$, which will remain fixed for the remainder of the proof). The process $\bar u$ is a solution to the SPDE
\beq
\label{ch5-added-0(*3)}
\left(\frac{\partial}{\partial t}- \frac{\partial^2}{\partial x^2}\right) \bar u(t,x)  = \sigma(t, x, \bar u(t,x)) \sum_{j = 1}^n e_j(x) \dot W_t^j + b(t,x,\bar u(t,x)),
\eeq
$(t, x) \in\, ]0,\infty[\, \times D$, where $W_t^j=W_t(e_j)$, along with the initial and boundary conditions
 \begin{align*}
 \begin{cases}
 \bar u(0, x) = u_{0}(x), & {\text{if}}\ x \in D,\\
\bar u(t, 0) = \bar u(t, L) =0, & {\text{if}}\ t \in \, ]0,\infty[.
\end{cases}
\end{align*}
According to \eqref{ch1'-s60.2}, $(\bar u(t,x),\, (t,x)\in \R_+ \times D)$ satisfies the integral equation 
\begin{align}
\label{ch1'-s60.2-ch6}
\bar u(t,x) &= I_0(t,x) + \sum_{j=1}^n \int_0^t  \langle G_L(t-s;x,\ast)\, \sigma(s,\ast,\bar u(s,\ast)), e_j\rangle_V\, dW_s(e_j)\notag\\
&\qquad +  \int_0^t ds \int_D dy\, G_L(t-s;x,y)\, b(s,y,\bar u(s,y)),
\end{align}
where $I_0(t,x) = \int_D G_L(t;x-y) u_0(y)\,dy$ and $G_L(t;x,y)$ is defined in \eqref{ch1'.600}.

\medskip

\noindent{\em Step 2.~Further approximation by a diffusion process.} For $m \in \N^*$, let $V_m$ be the linear subspace of $V$ spanned by $(e_j,\, j=  1,\dots, m)$ and let $\Pi_{V_m}$ denote the orthogonal projection from $V$ onto $V_m$. We are going to further approximate $\bar u$ by a sequence $(v_m(t,x),\, m \in \N^*)$ which, for $m \in \N^*$, is defined by
\beq
\label{ch5-added-0(*3a)}
    v_m(t,x) := \Pi_{V_m}\bar u(t,\ast)(x) = \sum_{k = 1}^m \, \langle \bar u(t,\ast), e_k\rangle_V \, e_k(x)
    \eeq
(this approximation is close to the Galerkin approximation of $\bar u$). For $k=1,\dots, m$, set  $a_k(t):=\langle \bar u(t,\ast), e_k\rangle_V$,  $\lambda_k = - (k \pi /L)^2$, 
\beqn
   \sigma_{j, k}(t) = \langle \sigma(t, \ast, \bar u(t,\ast)) e_j, e_k \rangle_V \quad
\text{and}\quad
   b_k(t) = \langle b(t,*, \bar u(t, *)), e_k \rangle_V.
\eeqn
Take the inner product with $e_k$ on both sides of \eqref{ch1'-s60.2-ch6} to see that
\beq\label{coeff-proj}
     a_k(t) = e^{\lambda_kt}\, a_k(0) + \sum_{j=1}^n \int_0^t e^{\lambda_k(t - s)}  \sigma_{j,k}(s)\, dW_s^j
     + \int_0^t e^{\lambda_k(t - s)} b_k(s)\, ds.
  \eeq
Therefore, $(a_k(t))$ is an Ornstein-Uhlenbeck process that is the solution of the SDE
\beqn
\begin{cases}
      d a_k(t) = \lambda_k\, a_k(t)\, dt + \sum_{j = 1}^n \sigma_{j, k}(t)\, dW_t^j + b_k(t)\, dt,\\
    a_k(0) = \langle u_0, e_k \rangle_V.
    \end{cases}
    \eeqn
Using \eqref{ch5-added-0(*3a)} and \eqref{coeff-proj}, we obtain the following expression for $v_m(t,x)$: 
 \begin{align}
\label{ch5-added-0(*4)-pre}
    v_m(t,x) &= \sum_{k = 1}^m a_{k}(t)  e_k(x)\notag\\
   &= \Pi_{V_m}I_0(t,\ast)(x) 
   + \int_0^t \sum_{k = 1}^m e^{\lambda_k(t - s)}  \sum_{j=1}^n \sigma_{j, k}(s) e_k(x)\, dW_s^j \notag\\
   &\qquad+ \int_0^t ds \sum_{k = 1}^m e^{\lambda_k(t - s)} b_k(s) e_k(x) .
   \end{align}
In particular, using the fact that $\frac{\partial^2}{\partial x^2} e_k = \lambda_k\, e_k$,  we obtain that for fixed $x\in D$, $t\mapsto v_m(t,x)$ is a diffusion process such that
\beq
\label{ch5-added-0(*3a)bis}
   d v_m(t,x) =  \frac{\partial^2}{\partial x^2} v_m(t,x)\,  dt +  \sum_{j = 1}^n  \sum_{k = 1}^m  e_k(x)\sigma_{j,k}(t)\, d W_t^j
    + \sum_{k = 1}^m b_k(t) e_k(x)\, dt,   
   \eeq
and the quadratic variation of $v_m(\cdot,x)$ is
\beq
\label{ch5-added-0(*3b)}
   \langle v_m(\cdot, x) \rangle_t = \int_0^t ds\, \sum_{j = 1}^n     (\Pi_{V_m}(\sigma(s,\ast,\bar u(s,*)) e_j)(x))^2   .      
\eeq
\medskip

\noindent{\em Step 3.~An approximate stochastic heat equation.} For $m\in\N^*$, let
\beqn
G_m(t,x,y) = \Pi_{V_m} G_L(t;x,\ast) = \sum_{k=1}^m e^{\lambda_kt} e_{k}(x) e_{k}(\ast).
\eeqn
 By the definitions of $\sigma_{j, k}$ and $b_k$, the sum of the two last terms in \eqref{ch5-added-0(*4)-pre} is equal to
 \begin{align*}
    &\int_0^t \sum_{j=1}^n   \langle \sigma(s, *, \bar u(t,*)) e_j, \sum_{k = 1}^m e^{\lambda_k(t - s)} e_k(x) e_k \rangle_V\, dW_s^j\\
    &\qquad+ \int_0^t \, \langle b(s, *, \bar u(s,*)), \sum_{k = 1}^m e^{\lambda_k(t - s)} e_k(x) e_k \rangle_V 
    \end{align*}
  which, by the definition of $G_m$ is
    \begin{align*}
    &\sum_{j=1}^n \int_0^t  \langle  \sigma(t, *, \bar u(t,*)) e_j ,  G_m(t-s, x, *) \rangle_V\, dW_s^j\\
    &\qquad + \int_0^t  ds\, \langle b(s, *, \bar u(s,*)), G_m(t-s, x, *) \rangle_V .
    \end{align*}
 We deduce that
    \begin{align}
      v_m(t,x) &= \Pi_{V_m} I_0(t, \ast)(x)\notag \\
    &\qquad+\sum_{j=1}^n \int_0^t  \langle  G_m(t-s; x, *)\sigma(s, *, \bar u(s,*)), e_j \rangle_V\, dW_s^j\notag\\
    &\qquad +  \int_0^t  ds \int_D dy\, G_m(t-s, x, y) b(s, y, \bar u(s,y)).    
   \label{ch5-added-0(*4)}
\end{align}
This is similar to the equation \eqref{ch1'-s60.2-ch6} for $\bar u$, except that the initial condition and $G_L$ there are replaced here by $\Pi_{V_m}I_0(t,\ast)(x)$ and $G_m$, respectively.
\medskip

\noindent{\em Step 4.~Convergence of the $v_m$ to $\bar u$.} Notice that for all $t\ge 0$ and $x\in D$,
\beq
   \label{ch5-added-0(*5)}
    \lim_{m \to \infty} \int_0^t \Vert G_m(r, x, *) - G(r, x,*) \Vert_{V}^2\, dr = 0.      
   \eeq
Indeed, $G_m(r, x, *) = \Pi_{V_m} G(r, x,*)$, so \eqref{ch5-added-0(*5)} follows from the fact that $\Vert G(r, x,*) \Vert_{L^2([0, t] \times D)}  < \infty$.
   This implies that for each $x\in D$ and each $j$,
   \begin{align}
     &E\left[\left(\int_0^t  \langle  \sigma(s, *, \bar u(s,\ast))  (G_m(t-s, x, \ast) - G(t-s,x,\ast)), e_j \rangle_V\, dW_s^j \right)^2\right]\notag\\
      &\quad= E\left[\int_0^t \left\langle  \sigma(s, \ast, \bar u(s,\ast)) (G_m(t-s, x, \ast) - G(t-s,x,\ast)), e_j \right\rangle_V^2\,  ds\right]\notag\\
      &\quad\le  E\Big[\int_0^t ds\,
   \Vert\sigma(s, \ast, \bar u(s,\ast)) (G_m(t-s; x, \ast) - G(t-s;x,\ast))\Vert^2_V\notag\\
   &\quad=  \int_0^t ds  \int_D dy\,  (G_m(t-s, x, y) - G(t-s,x,y))^2 \notag \\
   &\qquad\qquad\qquad\qquad\qquad \times E[\sigma^2(s, y, \bar u(s,y))]\notag\\
   &\quad\leq C \int_0^t \Vert G_m(r, x, *) - G(r, x,*) \Vert_{V}^2\, dr
    \to 0,    
\label{ch5-added-0(*6)}
\end{align}
as $m \to \infty$, where we have used the Cauchy-Schwarz inequality, \eqref{ch1'-s60.3} (with $p=2$) and \eqref{ch5-added-0(*5)}.

 With similar arguments, we see that as  $m \to \infty$,
\beq
\label{ch5-added-0(*7)}
   E\left[\left(\int_0^t ds  \int_D dy\, (G_m(t-s, x, y) - G(t-s, x, y)) b(s, y, \bar u(s,y))\right)^2\right] \to 0,    
   \eeq
and because $u_0\in V$, for a.a.~$x\in D$, we have 
   \beqn
   \lim_{m\to\infty}\left\vert \Pi_{V_m} I_0(t,\ast)(x)-I_0(t,x)\right\vert = 0.
   \eeqn
Together with \eqref{ch5-added-0(*4)}, \eqref{ch5-added-0(*6)} and  \eqref{ch5-added-0(*7)},
we conclude that for all $t \in [0, T]$ and a.a.~$x\in D$,
\beq
\label{ch5-added-0(*7bis)}
    \lim_{m \to \infty} E\left[(v_m(t,x) - \bar u(t,x))^2\right] = 0.   
 \eeq
 \smallskip

\noindent{\em Step 5.~Smoothing the positive part.} We now put the index $i$ back into \eqref{ch5-added-0(*3)}, \eqref{ch5-added-0(*3a)} and the other variables. Let
  $w_m(t,x) = v_{m, 1}(t,x) - v_{m, 2}(t,x)$ and $w(t,x) = \bar u_1(t,x) - \bar u_2(t,x)$.
    Following \cite[Section 2]{donati-martin-pardoux-1993}, for $p \in \N^*$, define $\psi_p : \R \to \R$ by
    \beqn
    \psi_p(v) =
    \begin{cases}
    0&\  {\text{if}}\ v \leq 0,\\
                      2p v &\  {\text{if}}\ v \in [0, 1/p],\\
                      2& \   {\text{if}}\   v \geq 1/p,
                     \end{cases}
                     \eeqn
and $\varphi_p: \R \to \R$ by
\beqn
    \varphi_p(v) = 1_{\R_+}(v) \int_0^v dx \int_0^x dy\, \psi_p(y).
    \eeqn
Then $\varphi_p \in C^2(\R)$ and for any $v\in\re$,
\beqn
     0 \leq \varphi_p'(v) \leq 2 v^+,\qquad   0 \leq \varphi_p''(v) \leq 2 1_{\R_+}(v),
\eeqn
and
\beq
\label{ch5-added-0(*9)}
\varphi_p(v) \uparrow (v^+)^2,\qquad \varphi_p'(v) \uparrow 2 v^+,\qquad \varphi_p''(v) \uparrow 2 1_{\R_+}(v)
\eeq
as $p \to \infty$.     

    Consider the random variables
    \begin{align*}
     \Phi_{p,m}(t) &= \int_D \varphi_p(w_m(t,x))\, dx,\quad  \Phi_p(t) = \int_D \varphi_p(w(t,x))\, dx\\
       \Phi(t) &= \int_D (w^+(t,x))^2\, dx.
     \end{align*}
Observe that for $t \in [0, T]$, $\Phi(t)<\infty$ by \eqref{ch1'-s60.3},
\beq
\label{ch5-added-0(*9a)}
\Phi_p(t) \leq \Phi(t)  \quad {\text{and}}\quad    \lim_{m \to \infty} \Phi_{p, m}(t) = \Phi_p(t)    \quad {\text {a.s.}}
\eeq
Since
\beqn
\lim_{p \to \infty} \int_D  \varphi_p(w(t,x))\, dx =  \int_D (w^+(t,x))^2\, dx
\eeqn
 by \eqref{ch5-added-0(*9)} and monotone convergence, we have
\beq
\label{ch5-added-0(*9b)}
\lim_{p\, \to \infty} \Phi_p(t) = \Phi(t).
    \eeq
 \medskip

\noindent{\em Step 6.~Comparing the diffusion approximations.} We are going to show that for all $t \in [0, T]$, $\Phi(t) = 0$ a.s. This will imply that $w(t,x) \leq 0$ a.s., that is, \eqref{ch5-added-0(*2)} holds, and this will complete the proof of Theorem \ref{ch5-added-0-t1}.

    In the following, we omit for simplicity the variables $t$ and $x$ in $b(t,x,z)$ and $\sigma(t,x,z)$.

     Apply the standard It\^o's formula to obtain for each $x \in D$,
    \begin{align*}
     \varphi_p(w_m(t,x)) &= \varphi_p\left(w_m(0,x)\right)\\
     &\quad \quad+ \int_0^t \varphi_p'(w_m(s,x))\, (dv_{m, 1}(s,x)) - dv_{m, 2}(s,x))\\
     &\quad \quad+ \half\int_0^t \varphi_p''(w_m(s,x))\, d\langle v_{m, 1}(\cdot, x) - v_{m, 2}(\cdot, x)\rangle_s.
     \end{align*}
Using  \eqref{ch5-added-0(*3a)bis} and \eqref{ch5-added-0(*3b)}, we see that this is equal to
     \begin{align*}
   &\varphi_p\left(w_m(0,x)\right) + \int_0^t \varphi_p'(w_m(s,x))\,  \frac{\partial^2}{\partial x^2} w_m(t,x) \, ds\\
   &\qquad + \int_0^t \varphi_p'(w_m(s,x)) \sum_{j = 1}^n  \left( \Pi_{V_m} [ ( \sigma(\bar u_1(s,*)) - \sigma(\bar u_2(s, *)) ) e_j ](x) \right)\, dW_s^j\\
    &\qquad + \int_0^t \varphi_p'(w_m(s,x))\, \Pi_{V_m}[b_1(\bar u_1(s,*)) - b_2(\bar u_2(s,*))](x)\, ds \\
    &\qquad + \half \int_0^t \varphi_p''(w_m(s,x)) \sum_{j = 1}^n\left(\Pi_{V_m}[(\sigma(\bar u_1(s,*)) - \sigma(\bar u_2(s, *))) e_j](x)\right)^2\, ds.
    \end{align*}
Integrate over $x$ to see that
\begin{align*}
    \Phi_{m,p}(t) &= \int_D \varphi_p(w_m(0,x))\, dx
    + \int_D dx \int_0^t  \varphi_p'(w_m(s,x))\,  \frac{\partial^2}{\partial x^2} w_m(t,x)\,  ds\notag\\
     &\quad+ \int_D  dx \int_0^t \varphi_p'(w_m(s,x))\\
     &\qquad\quad\times \sum_{j = 1}^n  (\Pi_{V_m}[(\sigma(\bar u_1(s,*)) - \sigma(\bar u_2(s, *))) e_j](x))\, dW_s^j \notag\\
     &\quad +
\int_D dx \int_0^t  \varphi_p'(w_m(s,x)) (\Pi_{V_m}[b_1(\bar u_1(s,*)) - b_2(\bar u_2(s,*))](x))\, ds\notag\\
&\quad + \half \int_D dx \int_0^t \varphi_p''(w_m(s,x))\notag\\
&\qquad\quad\times \sum_{j = 1}^n (\Pi_{V_m}[(\sigma(\bar u_1(s,*)) - \sigma(\bar u_2(s, *))) e_j](x))^2 ds.
\end{align*}
Integrating once by parts in the second term on the right-hand side of this equality and using the boundary condition  $\varphi_p'(w_m(s,0)) = \varphi_p'(w_m(s,L))=0$, we deduce that
\begin{align}
\label{ch5-added-0(*9c)}
 \Phi_{m,p}(t) &= \int_D \varphi_p(w_m(0,x)\, dx\notag\\
 &\quad -\int_0^t \left\langle \varphi_p''(w_m(s,*)) \frac{\partial}{\partial  x} w_m(s,*),  \frac{\partial}{\partial x} w_m(s,*) \right\rangle_V\, ds\notag\\
&\qquad + \sum_{j = 1}^n \int_0^t \left\langle \varphi_p'(w_m(s,*)), \Pi_{V_m}[(\sigma(\bar u_1(s,*) - \sigma(\bar u_2(s,*))) e_j] \right\rangle_V\,  dW_s^j \notag\\
&\qquad + \int_0^t  \left\langle \varphi_p'(w_m(s,*)), \Pi_{V_m}[b_1(\bar u_1(s,*))  - b_2(\bar u_2(s,*))] \right\rangle_V\, ds\notag\\
&\qquad +  \half \sum_{j=1}^n \int_0^t \left\langle \varphi_p''(w_m(s,*))
\Pi_{V_m}[(\sigma(\bar u_1(s,*)) - \sigma(\bar u_2(s,*))) e_j],\right.\notag\\
&\left.\qquad\qquad\qquad \Pi_{V_m}[(\sigma(\bar u_1(s,*)) - \sigma(\bar u_2(s,*))) e_j] \right\rangle_V\,  ds.
\end{align}
Taking expectations in both sides of this equality yields
\begin{align}
\label{ch5-added-0(*9d)}
  & E[\Phi_{m,p}(t)] \notag\\
  &\quad=  \int_D E\left[\varphi_p\left(w_m(0,x)\right)\right] dx\notag\\
   &\qquad - \int_0^t  E\left[\left\langle \varphi_p''(w_m(s,*)) \frac{\partial}{\partial x} w_m(s,*),   \frac{\partial}{\partial x} w_m(s,*) \right\rangle_V\right]  ds \notag\\
   &\qquad +  \int_0^t  E\left[\left\langle \varphi_p'(w_m(s,*)), \Pi_{V_m}[b_1(\bar u_1(s,*)) - b_2(\bar u_2(s,*))] \right\rangle_V\right] ds\notag\\
   &\qquad + \half \sum_{j=1}^n \int_0^t E\left[\left\langle \varphi_p''(w_m(s,*))
    \Pi_{V_m}[(\sigma(\bar u_1(s,*)) - \sigma(\bar u_2(s,*))] e_j),\right.\right.\notag\\
   &\left.\left.\qquad\qquad\qquad \Pi_{V_m}[(\sigma(\bar u_1(s,*)) - \sigma(\bar u_2(s,*))) e_j] \right\rangle_V\right]  ds.   
   \end{align}
Observe that the second integral on the right-hand side is nonnegative (in fact, even without the expectation, this integral is nonnegative). Therefore, it can be removed to turn the equality into an inequality.  By doing so and then taking the limit as $m \to \infty$, and using the continuity of $\varphi_p''$, $\varphi_p'$ and $\varphi_p$, we  obtain
\begin{align*}
E[\Phi_p(t)]& \leq \int_D E[\varphi_p(w(0,x))]\, dx\\
          &\qquad   + \int_0^t E\left[\left\langle \varphi_p'(w(s,*)),  b_1(\bar u_1(s,*)) - b_2(\bar u_2(s, *)) \right\rangle_V\right] ds\\
           &\qquad    + \half\sum_{j=1}^n \int_0^t E\left[\left\langle \varphi_p''(w(s, *)) (\sigma(\bar u_1(s, *)) - \sigma(\bar u_2(s, *))) e_j,\right.\right.\\
                   &\left.\left.\qquad\qquad\qquad(\sigma(\bar u_1(s, *)) - \sigma(\bar u_2(s, *))) e_j  \right\rangle_V\right] ds.
\end{align*}

In the second integral on the right-hand side, we write
\begin{align*}
   b_1(\bar u_1(s, *)) - b_2(\bar u_2(s, *)) &= b_1(\bar u_1(s, *)) - b_1(\bar u_2(s, *))\\
   &\qquad  + b_1(\bar u_2(s, *)) - b_2(\bar u_2(s, *))\\
  &\leq b_1(\bar u_1(s, *)) - b_1(\bar u_2(s, *))
  \end{align*}
since $b_1 - b_2 \leq 0$ by hypothesis.  Since $\varphi_p'(w(s,*)) \geq 0$, we substitute this expression into the second integral, and use the Lipschitz properties of $b_1$ and $\sigma$ to deduce that
\begin{align}
\label{mars(*2)}
E[\Phi_p(t)] &\leq \int_D E[\varphi_p(w(0,x))]\, dx \notag\\
      &\quad+ C \int_0^t ds\, E\left[\int_D dx\,  \varphi_p'(w(s,x)) \vert \bar u_1(s,x) - \bar u_2(s,x) \vert\right] \notag\\
      &\quad+ \half C^2 \sum_{j=1}^n \int_0^t ds\ \int_D dx\, E\left[ \varphi_p''(w(s, x)) (\bar u_1(s,x) - \bar u_2(s,x))^2 e_j^2(x)\right].       
 \end{align}
We now let $p \to \infty$ and use \eqref{ch5-added-0(*9b)} and \eqref{ch5-added-0(*9)} to obtain
      \begin{align*}
    E[\Phi(t)] &\leq  \int_D E[(w^+(0,x))^2]\, dx\\
       &\quad +  2 C \int_0^t ds E\left[\int_D dx\,  w^+(s,x)\, \vert \bar u_1(s,x) - \bar u_2(s,x) \vert\right]\\
       &\quad+ C^2 \sum_{j=1}^n \int_0^t ds \int_D dx\, E\left[1_{\R_+}(w(s, x))(\bar u_1(s,x) - \bar u_2(s,x))^2 e_j^2(x)\right].
       \end{align*}
The first integral on the right-hand side vanishes by the hypothesis $u_{0, 1} \leq u_{0, 2}$. Because of the factors $w^+(s,x)$ and $1_{\R_+}(w(s,x))$, we can replace the absolute values by positive parts, to obtain
\begin{align*}
   E[\Phi(t)] \leq
    &\ 2 C\int_0^t E\left[\int_D dx\, w^+(s,x) (\bar u_1(s,x) - \bar u_2(s,x))^+\right] ds\\
    &\quad + C^2\sum_{j=1}^n \int_0^t ds \int_D dx\, E\left[1_{\R_+}(w(s,x))\right.\\
     &\left.\qquad \qquad\qquad \times ((\bar u_1(s,x) - \bar u_2(s,x))^+)^2\right] e_j^2(x).
    \end{align*}
Since  $e_j^2 \leq 2/L $, we have
\begin{align}
\label{ch5-added-0(*9e)}
     E[\Phi(t)] &\leq  2 C \int_0^t E\left[\int_D dx\, (w^+(s,x))^2\right] ds\notag\\
   &\qquad+ \frac{2}{L} C^2 n \int_0^t ds \int_D dx\, E\left[(w^+(s,x))^2\right] \notag\\
 & = \left(2C + \frac{2}{L} C^2 n\right) \int_0^t ds\, E[\Phi(s)].   
   \end{align}
Apply the classical Gronwall's lemma (Lemma \ref{lemC.1.1}) to deduce that for all $t \in [0, T]$, $E[\Phi(t)] = 0$. Since $w$ has continuous sample paths, we deduce that a.s., for all $(t,x) \in [0,T] \times D$, $w^+(t,x) = 0$, that is, $\bar u_{1}(t,x) \leq \bar u_{2}(t,x)$. This completes the proof of \eqref{ch5-added-0(*2)} and of Theorem \ref{ch5-added-0-t1}.
 \end{proof}

\begin{remark}
\label{ch5-added-0-r1} (a)\  The same result applies to the stochastic heat equation with vanishing Neumann boundary conditions, with the same proof. Indeed, these boundary conditions also make the boundary terms vanish in \eqref{ch5-added-0(*9c)}.

    (b)\  The term
    \beqn
    -\int_0^t  E\left[\left\langle \varphi_p''(w_m(s,\ast))\frac{\partial}{\partial x} w_m(s,\ast),  \frac{\partial}{\partial x} w_m(t,\ast) \right\rangle_V \right]\,  ds
    \eeqn
     in \eqref{ch5-added-0(*9d)} can be moved to the left-hand side and included further along the calculation. Since the final bound \eqref{ch5-added-0(*9e)} does not depend on $m$, this argument can be used to show that $\bar u_{n,i}$ takes values in $L^2([0,T], H_0^1(D))$, that is, for a.a. $ t\in[0,T]$, $x \mapsto \bar u_{n,i}(t,x)$ is absolutely continuous with a derivative in $L^2(D)$ (see, for example, \cite[Section 2.4]{pardoux79}).

    \end{remark}
 
\subsection{The case $D = \R$}\label{rd05_06ss1}
    
Comparison theorems for the stochastic heat equation on $\R$ are also available (see e.g. \cite{mueller-1991}, \cite{shiga-1994}, \cite{joseph-khosh-mueller-2017}, \cite{chen-kim-2016}). They seem all to involve a discretization of the noise, of the Laplacian, and time and/or space. Here, we use a different approach via Theorem \ref{rd08_09t1}, which discusses an approximation of the stochastic heat equation on $\R$ by the stochastic heat equation on the interval $[-L,L]$ with Dirichlet boundary conditions (as $L\to\infty$). In this way, we easily deduce from Theorem \ref{ch5-added-0-t1} a comparison theorem for the stochastic heat equation on $\R$ (see Theorem \ref{rd08_13t1} below).

For $i = 1, 2$, let $\sigma(t,x,z)$ and $b_i(t, x, z)$ be functions satisfying $({\bf H_L})$ and let $u_{0, i}$ be a function satisfying $({\bf H_I})$ with $D = \R$. We consider the stochastic heat equation on $\R$ with initial condition $u_{0, i}$, as in Section \ref{ch1'-ss6.1}:
\beq
\label{rdch5-added-0(*1)}
 \left(\frac{\partial}{\partial t}- \frac{\partial^{2}}{\partial x^{2}}\right)u_i(t,x) = \sigma(t,x,u_i(t,x))\dot W(t, x)
 +b_i(t,x,u_i(t,x)),
\eeq
 $(t, x) \in\, ]0,\infty[ \times \R$, along with the initial condition
 \begin{align*}
 u_i(0, x) = u_{0,i}(x), & \ x \in \R.
\end{align*}

 In the proof of the next theorem, we use the notations introduced in Section \ref{ch1'-s60}, and we denote by $u_1$ and $u_2$  the continuous versions of the solutions given by Theorems \ref{recap} and \ref{cor-recap}.

\begin{thm}
\label{rd08_13t1}
Suppose that for each $z\in \re$, $b_1(\cdot, \ast, z) \leq b_2(\cdot, \ast, z)$ $dsdxdP$-a.e.~and $u_{0, 1} \leq u_{0, 2}$ a.e. Then a.s., for all $(t,x) \in [0, T] \times \R$,
    $u_1(t,x) \leq u_2(t,x)$.
\end{thm}

\begin{proof}
We let $u_{i,L}$ denote the solution to the stochastic heat equation on $[-L,L]$ with vanishing Dirichlet boundary conditions, the same coefficients $\sigma$, $b_i$ and the same initial condition $u_{0,i}(x)$, $x \in [-L,L]$, as for $u_i$. By Theorem \ref{ch5-added-0-t1}, a.s., for all $(t, x) \in [0, T] \times [-L, L]$, $u_{1,L}(t, x) \leq u_{2,L}(t, x)$. By Theorem \ref{rd08_09t1}, for $i=1, 2$ and $(t, x) \in [0, T] \times \R$,
\beqn
   \lim_{L \to \infty} u_{i, L}(t, x) =  u_i(t, x) \quad\text{ in }  L^p(\Omega).
\eeqn
Therefore, for all $(t, x) \in [0, T] \times \R$, $u_{1}(t, x) \leq u_{2}(t, x)$ a.s. Because $(u_i(t,x))$ has continuous sample paths, the theorem is proved.
\end{proof}

\noindent{\em Application: nonnegativity in the parabolic Anderson model}
\medskip

    We end this section with an application of Theorem \ref{ch5-added-0-t1} to
   a class of equations that includes
    the {\em parabolic Anderson models} on $]0, L[$ and on $\R$ (see Section \ref{ch1-1.2} for a similar SPDE on $\rek$, and Section \ref{rd1+1anderson}).

  Consider the stochastic heat equation \eqref{ch5-added-0(*1)} with vanishing Dirichlet or Neumann boundary conditions (respectively the stochastic heat equation \eqref{rdch5-added-0(*1)}),
 and suppose, in addition to the hypotheses of Theorem \ref{ch5-added-0-t1} (respectively Theorem \ref{rd08_13t1}), that the functions $\sigma$ and $b_1$ are such that $\sigma(\cdot, *, 0) \equiv b_1(\cdot, *, 0) \equiv 0$. Then the solution to equation \eqref{ch5-added-0(*1)} with initial condition $u_{0,1, L} \equiv 0$ is $u_{1, L}(\cdot, *) = 0$ (respectively, to \eqref{rdch5-added-0(*1)} with initial condition $u_{0, 1} \equiv 0$ is $u_1(\cdot, *) = 0$). By Theorem \ref{ch5-added-0-t1}, for any nonnegative initial condition $u_{0,2, L}$, the solution $u_{2, L}$ of \eqref{ch5-added-0(*1)} will satisfy $u_{2, L}(\cdot, *) \geq u_{1, L}(\cdot,*) = 0$. In particular, $u_{2, L}(\cdot, *)$ will remain nonnegative for all time. Similarly, by Theorem \ref{rd08_13t1}, for any nonnegative initial condition $u_{0, 2} \equiv 0$, the solution $u_2$ of \eqref{rdch5-added-0(*1)} will satisfy $u_2(\cdot, *) \geq u_1(\cdot, *) = 0$.
 
 These conclusions are valid in particular for the parabolic Anderson models on $D = \, ]0, L[$ and on $D = \R$ ($\sigma(t, x, u(t,x)) = \rho\, u(t,x)$, $\rho \in \R \setminus \{0\}$, and  $b_i \equiv 0$, $i = 1, 2$).

\section{Polarity of points in high dimensions}
\label{ch2'-s7}

In this section, we consider a $d$-dimensional random field denoted by $u=(u(t,x),\, (t,x)\in\re_+\times \rek)$, where $u(t,x) = (u_1(t,x),\dots, u_d(t,x))$. A relevant notion in probabilistic potential theory is that of {\em polar set}. It is defined as follows.

 Let $I=I_1\times I_2$ be a compact subset of $\re_+\times \rek$, and $A\in\B(\red)$. We denote by $u(I)$ the random
 subset of $\re^d$ consisting of the positions visited by $u$ restricted  to $I$, that is,
 \beqn
 u(I)= \{u(t,x)\in\red: (t,x)\in I_1\times I_2\}.
 \eeqn
 Then $A$ is a {\em polar set}\index{polar set}\index{set!polar} for $u$ restricted to $I$ if
 \beq
 \label{ch2'-s7.1}
 P\{u(I)\cap A\ne\emptyset\}=0,
 \eeq
and is {\em non-polar} if $P\{u(I)\cap A\ne\emptyset\}\ne 0$.

For a random field $u$ restricted to $I_1 \times I_2$ and a set $A$, the property of polarity is related to the regularity of the sample paths of $u$, the geometric measure-theoretic properties of the set $A$ (such as, for example, Hausdorff measure), and the dimensions of $I_1$ and $I_2$. An interesting question is to characterize polarity for classes of random fields. For random field solutions of SPDEs, this problem has been addressed in several papers (see e.g. \cite{m-t2002}, \cite{dn2004}, \cite{dkn07}, \cite{dkn09}, \cite{dss10}, \cite{dss15}, \cite{ssv18}, \cite{dalang-pu-2021}, \cite{hss2020}).

In this section, as an introduction to the topic, we consider the particular case where $A$ is a singleton, except in Section \ref{ch2'-s7.s2}. We will give
 sufficient conditions for a point to be polar for an anisotropic $d$-dimensional random field, considering separately the Gaussian and non-Gaussian cases. The results will be applied to the study of polarity of points for solutions to systems of some classes of  SPDEs.

  \subsection[Sufficient conditions for polarity of points:\texorpdfstring{\\}{} the Gaussian case]{Sufficient conditions for polarity of points: the Gaussian case}
  \label{Ch5-ss7.1}
 We fix the subset of parameters $I=[t_0,t_1]\times[x_1,x_2]$, where
 $[x_1, x_2] = \prod_{\ell=1}^k [x_{1,\ell}, x_{2,\ell}]$, $0\le t_0<t_1$,
$x_1 = (x_{1,1}, \dots, x_{1,k})$, $ x_2 = (x_{2,1}, \dots, x_{2,k})$, with $x_{1,\ell} < x_{2,\ell}$ ($\ell=1,\ldots k$).
For $\varepsilon_0>0$, we denote by $I^{(\varepsilon_0)}$ the closed $\varepsilon_0$-neighbourhood of $I$, that is,
 \beqn
 I^{(\varepsilon_0)} = \{(s,y)\in\re_+\times \rek:\ d((s,y), I)\le \varepsilon_0\},
 \eeqn
 where $d$ here denotes the Euclidean distance.

 \begin{prop}
\label{ch2'-s7-p1}
Assume that $u=(u_1,\dots,u_d)$ is a Gaussian centred $d$-dimensional random field with independent identically distributed components and continuous sample paths a.s. Let $\varepsilon_0>0$.
Suppose that the following conditions are satisfied:
\begin{description}
\item{(i)} There exist $\alpha_\ell\in\, ]0,1]$, $\ell=0,1,\ldots,k$, and a constant $C>0$ such that, for any $(t,x),\, (s,y)\in I^{(\varepsilon_0)}$,
\beq
\label{ch2'-s7.2}
\left\Vert u(t,x)-u(s,y)\right\Vert_{L^2(\Omega)} \le C \left(|s-t|^{\alpha_0}+\sum_{\ell=1}^k|x_\ell-y_\ell|^{\alpha_\ell}\right).
\eeq
\item{(ii)}
$\inf_{(t,x)\in I^{(\varepsilon_0)}} {\rm{Var}}(u_1(t,x))>0$.
\end{description}
Set $Q=\sum_{\ell=0}^k\frac{1}{\alpha_\ell}$. If $d>Q$, then points are polar for $u$ restricted to $I$.
\end{prop}

\begin{remark}
\label{ch2'-s7-r1}
For a random field $u$ satisfying condition (i) of Proposition \ref{ch2'-s7-p1}, the version of Kolmogorov's continuity criterion given in Theorem \ref{app1-3-t1} implies the following:

For any $p> Q$, $\alpha\in \left]\tfrac{Q}{p}, 1\right[$, there exists a constant $C:=C(\alpha, p)$ such that,
\beq
\label{ch2'-s7.200}
E\left[\sup_{(t,x)\ne (s,y)} \frac{|u(t,x)-u(s,y)|^p}{\left(|s-t|^{\alpha_0}+\sum_{\ell=1}^k |x_\ell-y_\ell|^{\alpha_\ell}\right)^{p\alpha- Q}}\right] \le C,
\eeq
where the supremum  is over points $(t,x), (s,y)$ in the compact set $I^{(\varepsilon_0)}$ (see \eqref{ch1'-s7.20}). This fact will be used in the proof of Proposition \ref{ch2'-s7-p1}.

\end{remark}

\noindent{\em Proof of Proposition \ref{ch2'-s7-p1}.}\
For any $j_0\in \N$, $ j_1, \ldots, j_k\in \Z$, $\varepsilon>0$, set $j=(j_0, j_1, \ldots, j_k)$,
\beqn
R_j^\varepsilon =\prod_{\ell=0}^k\left[j_\ell\varepsilon^{\frac{1}{\alpha_\ell}},(j_\ell+1)\varepsilon^{\frac{1}{\alpha_\ell}}\right]\quad {\text{and}}\quad
y_j^\varepsilon = (j_0 \varepsilon^{\frac{1}{\alpha_0}}, \ldots, j_k \varepsilon^{\frac{1}{\alpha_k}}).
\eeqn
Notice that the Lebesgue measure of $R_j^\varepsilon$ is $\ep^Q$.

 We want to prove that for any $z\in\red$,
 \beq
\label{ch2'-s7.3}
P\left\{u(I)\cap \{z\}\ne\emptyset\right\} = 0.
\eeq
 For this, observe that
\beqn
P\left\{u(I)\cap \{z\}\ne\emptyset\right\} = P\{u(t,x)= z,\  \text{for some}\ (t,x)\in I\}.
\eeqn
Clearly, for any $\varepsilon>0$,
\begin{align*}
&P\{u(t,x)= z, \ \text{for some}\ (t,x)\in I\}\\
&\qquad\qquad\le P\{|u(t,x)-z\vert < \varepsilon,\ \text{for some}\ (t,x)\in I\}.
\end{align*}
By a covering argument, the last term is bounded above by
\begin{align*}
&\sum_{j: R_j^\varepsilon\cap I\ne \emptyset} P\left\{|u(t,x)-z\vert < \varepsilon, \text{for some}\ (t,x)\in R_j^\varepsilon\right\}\\
&\qquad \qquad = \sum_{j: R_j^\varepsilon\cap I\ne \emptyset} P\left\{\inf_{(t,x)\in R_j^\varepsilon} |u(t,x)-z|<\varepsilon\right\},
\end{align*}
since the sample paths of $u$ are continuous a.s.
Observe that for the values of $j$ in the sum, if $\varepsilon>0$ is small enough, then $R_j^\varepsilon\subset I^{(\ep_0)}$.
Notice also that the cardinality of the set $J_\varepsilon:=\{j: R_j^\varepsilon\cap I\ne \emptyset\}$ is bounded by a constant times $\varepsilon^{-Q}$.

Our aim is to prove that for all $\eta \in\, ]0,1[$, there exists $C > 0$ such that for all $\ep > 0$ small enough and for all $j\in J_\varepsilon$,
\beq
\label{ch2'-s7.4}
P\left\{\inf_{(t,x)\in R_j^\varepsilon} |u(t,x)-z|<\varepsilon\right\} \le C \varepsilon^{\eta d}.
\eeq
With this, we obtain
\beqn
P\left\{u(I)\cap \{z\}\ne\emptyset\right\} \le \sum_{j\in J_\varepsilon} P\left\{\inf_{(t,x)\in R_j^\varepsilon} |u(t,x)-z|<\varepsilon\right\} \le C \varepsilon^{\eta d-Q}.
\eeqn
We are assuming $d>Q$ and therefore, there exists $\eta\in\, ]0,1[$ such that $d\eta-Q>0$. With this choice of $\eta$ and letting $\varepsilon\downarrow 0$ in the last inequalities, we obtain \eqref{ch2'-s7.3}
and thus, $\{z\}$ is polar for $u$ restricted to $I$.

In order to establish \eqref{ch2'-s7.4}, since the component random fields of $u$ are independent and identically distributed, it suffices to prove that for $\eta \in\, ]0,1[$ and $z_1 \in \R$, there is $C = C(\eta, z_1) > 0$ such that for all small enough $\ep > 0$  and for all $j\in J_{\varepsilon}$,
\beq
\label{ch2'-s7.4-bis}
P\left\{\inf_{(t,x)\in R_j^\varepsilon} |u_1(t,x)-z_1|<\varepsilon\right\} \le C \varepsilon^\eta.
\eeq

The remainder of the proof is devoted to establishing \eqref{ch2'-s7.4-bis}.
Set
\beqn
c_j^\varepsilon(t,x) = \frac{E[u_1(t,x)u_1(y_j^\varepsilon)]}{\text{Var} (u_1(y_j^\varepsilon))}.
\eeqn
Since $u_1$ is Gaussian and centred,
\beq
\label{fce}
E[u_1(t,x)\big| u_1(y_j^\varepsilon)] = c_j^\varepsilon(t,x) u_1(y_j^\varepsilon).
\eeq

In the sequel, we assume that $\varepsilon>0$ is such that $\cup_{j\in J_{\varepsilon}}R_j^\varepsilon\subset I^{(\varepsilon_0)}$. For its further use, we check that there is a constant $C>0$ such that for all $j\in J_{\varepsilon}$,
\beq
\label{ch2'-s7.5}
\sup_{(t,x)\in R_j^\varepsilon}\vert c_j^\varepsilon(t,x) -1\vert \le C \varepsilon.
\eeq

Indeed, by applying the Cauchy-Schwarz inequality and using condition (ii), we have
\begin{align*}
\vert c_j^\varepsilon(t,x) -1\vert & = \frac{\left\vert E\left[u_1(y_j^\varepsilon)(u_1(t,x)-u_1(y_j^\varepsilon))\right]\right\vert}{\text{Var} (u_1(y_j^\varepsilon))}\notag\\
&\le \frac{\left\Vert u_1(y_j^\varepsilon)\right\Vert_{L^2(\Omega)}}{\text{Var} (u_1(y_j^\varepsilon))}
\left\Vert u_1(t,x)-u_1(y_j^\varepsilon)\right\Vert_{L^2(\Omega)}\notag\\
&\le \tilde C\left\Vert u_1(t,x)-u_1(y_j^\varepsilon)\right\Vert_{L^2(\Omega)} . 
\end{align*}
By condition (i), $\left\Vert u_1(t,x)-u_1(y_j^\varepsilon)\right\Vert_{L^2(\Omega)} \le C \varepsilon$,
for any $(t,x)\in R_j^\varepsilon$, and \eqref{ch2'-s7.5} is proved. 

Observe that, for $\varepsilon$ small enough, \eqref{ch2'-s7.5} implies that
\beq
\label{ch2'-s7.6}
 \inf_{(t,x)\in R_j^\varepsilon}c_j^\varepsilon(t,x) \ge \frac{1}{2}
\eeq
and
\beq
\label{5.4.1(*1)}
 \left\vert \frac{1}{c_j^\ep} - 1 \right\vert = \frac{\vert c_j^\ep -1 \vert}{c_j^\ep} \leq 2 C \ep.   
 \eeq

We continue with the proof of \eqref{ch2'-s7.4-bis}.
Set
\begin{align*}
Y_j^\varepsilon &= \inf_{(t,x)\in R_j^\varepsilon} \vert E[u_1(t,x)\big| u_1(y_j^\varepsilon)]-z_1\vert, \\
Z_j^\varepsilon &= \sup_{(t,x)\in R_j^\varepsilon} \vert u_1(t,x) - E[u_1(t,x)\big| u_1(y_j^\varepsilon)]\vert.
\end{align*}
These two Gaussian random variables are independent, and we have
\beqn
P\left\{\inf_{(t,x)\in R_j^\varepsilon} |u_1(t,x)-z_1|<\varepsilon\right\} \le P\left\{Y_j^\varepsilon\le \varepsilon + Z_j^\varepsilon\right\}.
\eeqn
Since $u_1$ and $- u_1$ have the same law, we can assume that $z_1 \geq 0$. Set
\beqn
Z_j^{\varepsilon,1}= \sup_{(t,x)\in R_j^\varepsilon}|u_1(t,x)-u_1(y_j^\varepsilon)|, \qquad  Z_j^{\varepsilon,2}= |
u_1(y_j^\varepsilon)| \sup_{(t,x)\in R_j^\varepsilon} |1-c_j^\varepsilon(t,x)|.
\eeqn
By the triangle inequality, $Z_j^\varepsilon \le Z_j^{\varepsilon,1} + Z_j^{\varepsilon,2}$.

Consider $\eta \in\, ]0,1[$,
$p > Q$ and $\alpha \in\, ]Q/p, 1[$ such that $\alpha > \eta$. Then
\begin{align}
\label{5.4.1(*2)}
     P\{ Y_j^\ep \leq \ep + Z_j^\ep \} &\leq P\{ Y_j^\ep \leq \ep + Z_j^{\ep,1} + Z_j^{\ep,2} \}\notag\\
    &\leq P\{ Y_j^\ep \leq \ep + 2 \ep^\eta \} + P\{ Z_j^{\ep,1} > \ep^\eta \} + P\{ Z_j^{\ep,2} > \ep^\eta \}.
    \end{align}      
We bound each term separately. For the second term, by \eqref{ch2'-s7.200} and the definition of $R_j^\ep$, using Chebychev's inequality, we obtain
\beqn
   P\{ Z_j^{\ep,1} > \ep^\eta \} \leq C(\alpha,p)\, \ep^{p(\alpha - \eta) - Q}.
   \eeqn
Since $\alpha - \eta > 0$, we can choose $p$ large enough so that $p(\alpha - \eta) - Q \geq 1 > \eta$, to get
\beq
\label{5.4.1(*3)}
   P\{ Z_j^{\ep,1} > \ep^\eta \} \leq C(\alpha,p) \ep^\eta.      
   \eeq
For the third term on the right-hand side of \eqref{5.4.1(*2)}, we use \eqref{ch2'-s7.5} to see that
\beqn
   P\{ Z_j^{\ep,2} > \ep^\eta \}  \leq P\{ C \ep\,  \vert u_1(y_j^\ep) \vert > \ep^\eta \}
      = P\{\vert u_1(y_j^\ep) \vert > \ep^{\eta - 1} / C \}.
      \eeqn
Since the variance of $u_1(t,x)$ is a continuous function of $(t,x)$, it is bounded above over $I^{(\ep_0)}$ by some number $\sigma_0$, therefore by \eqref{rdlemC.2.2-before-1} with $a = \ep^{\eta - 1} / C$ there, for $\ep$ small enough,
\beq
\label{5.4.1(*4)}
   P\{ Z_j^{\ep,2} > \ep^\eta \} \leq 
  \exp\left( - \frac{\ep^{2(\eta - 1)}}{2 C^2 \sigma_0^2} \right) \leq K'(\eta)\, \ep^\eta,  
   \eeq
since $\eta - 1 < 0$.
   Finally, we consider the first term on the right-hand side of \eqref{5.4.1(*2)}. Set $r = 3 \ep^\eta$. We will show that
   \beq
   \label{5.4.1(*5)}
     P\{ Y_j^\ep \leq r \} \leq K \ep^\eta.    
\eeq

Indeed, by the definition of $Y_j^\varepsilon$ and by \eqref{fce} we see that the constraint $Y_j^\varepsilon \le r$ implies that for some $(t, x) \in R_j^\ep$,
\beqn
\frac{z_1}{c_j^\varepsilon (t,x)} - \frac{r}{c_j^\varepsilon (t,x)} \le u_1(y_j^\varepsilon) \le \frac{z_1}{c_j^\varepsilon (t,x)} + \frac{r}{c_j^\varepsilon (t,x)}.
\eeqn
By \eqref{5.4.1(*1)} and \eqref{ch2'-s7.6}, we obtain less stringent constraints on $u_1(y_j^\ep)$ if we only require that
    $(1 - 2C \ep) z_1 - 2r \leq u_1(y_j^\ep) \leq (1 + 2C \ep) z_1 + 2r,$
which no longer depends on $(t,x)\in R_j^\ep$. The length of this interval is $4C \ep z_1 + 4r = 4C \ep z_1 + 12 \ep^\eta \leq K'(z_1)\, \ep^\eta$ for small enough $\ep > 0$. Since the density of $u_1(y_j^\ep)$ is bounded uniformly in $j\in J_{\varepsilon}$, we conclude that
\beqn
    P\{ Y_j^\ep \leq \ep + 2 \ep^\eta \} \leq P\{ Y_j^\ep \leq 3 \ep^\eta \} \leq \tilde{K}'(z_1)\, \ep^\eta,
    \eeqn
proving \eqref{5.4.1(*5)}.

Putting together the bounds \eqref{5.4.1(*3)}--\eqref{5.4.1(*5)}, we obtain from \eqref{5.4.1(*2)} that for all
$\varepsilon>0$ small enough and $j\in J_{\varepsilon}$,
\beqn
P\left\{\inf_{(t,x)\in R_j^\varepsilon} |u_1(t,x)-z_1|<\varepsilon\right\} \le  P\left\{Y_j^\varepsilon\le \varepsilon + Z_j^\varepsilon\right\}
\le K{''}\, \varepsilon^\eta,
\eeqn
where $K'' = K''(\eta, z_1)$. This is \eqref{ch2'-s7.4-bis}. The proof of the proposition is complete.
\qed


\subsection[Examples: polarity of points for solutions to \texorpdfstring{\\}{}linear SPDEs]{Examples: polarity of points for solutions to linear SPDEs}
\label{ch2'-s7.e1}

In this section, we apply Proposition \ref{ch2'-s7-p1} to the study of polarity of points for Gaussian random fields that are solutions to the linear SPDEs considered in Chapter \ref{chapter1'}. More specifically, let $\mathcal{L}$ be a partial differential operator on $\re_+\times \rek$ and  $D\subset \rek$ be a bounded or unbounded domain with smooth boundary. Consider the system of linear SPDEs
\beqn
\mathcal{L}u_i(t,x) = \sum_{j=1}^d \sigma_{i,j} \dot W^j(t,x), \quad (t,x)\in\, ]0,\infty[ \times D, \ i=1,\ldots,d,
\eeqn
with vanishing initial conditions and, if $D$ is bounded, also with vanishing Dirichlet or Neumann boundary conditions. In this equation, $\sigma=(\sigma_{i,j})_{1\le i,j\le d}$ is a non-singular deterministic matrix and $\dot W^j$, $j=1,\ldots,d$, are independent copies of a space-time white noise. %

Since $\sigma$ is invertible, the random field $v:=\sigma^{-1} u$ satisfies the system of uncoupled SPDEs
\beq
\label{ch2'-s7.14}
\mathcal{L}v_i(t,x) = \dot W^i(t,x), \quad (t,x)\in\, ]0,\infty[\times D, \ i=1,\ldots,d,
\eeq
with vanishing initial and (if necessary) vanishing boundary conditions. Observe that a point $z\in \red$ is polar for $u$ if and only if $\sigma^{-1}(z)$ is polar for $v$.

Denote by $\Gamma(t,x;s,y)$ the fundamental solution or the Green's function on $\re_+\times D$ associated to $\mathcal{L}$ and suppose that Assumption \ref{ch4-a1-1'} holds. Then, in agreement with Definition \ref{ch4-d1-1'}, the random field solution to the system \eqref{ch2'-s7.14} is the process
$v = (v(t,x), (t,x) \in \re_+ \times D)$, where $v(t,x)$ is the $\red$-valued random variable $v(t,x) = (v_1(t,x), \dots,v_d(t,x))$, and
\beq
\label{ch2'-s7.15}
v_i(t,x) = \int_0^t \int_D \Gamma(t,x;s,y)\, W^i(ds,dy), \quad  i=1,\ldots,d.
\eeq
Notice that these components define independent and identically distributed random fields.

The next theorem gives sufficient conditions for polarity of points for the random field
$(v(t,x),\, (t,x)\in\re_+\times D)$, when $\mathcal{L}$ is the heat operator $\frac{\partial}{\partial t}-\frac{\partial^2}{\partial x^2}$ or the wave operator $\frac{\partial^2}{\partial t^2}-\frac{\partial^2}{\partial x^2}$.
\begin{thm}
\label{ch2'-s7-theorem3}
Consider the following cases:
\begin{enumerate}
\item {\em Linear stochastic heat equation on $D = \re$ or $D =\, ]0,L[$}.

Let $D_1 = \re$, $D_2 = D_3 =\, ]0,L[$. Let $v^{(j)} = (v^{(j)}(t,x),\, (t,x) \in \re_+ \times (D_j \cup \partial D_j))$, $j = 1,2,3$, be the $d$-dimensional random field with independent components, where the $i-$th component $v^{(j)}_i$, $i = 1, \dots, d$, is the solution to the linear stochastic heat equation with vanishing initial conditions \eqref{4.20-1'} if $j = 1$, \eqref{ch1'.HD} if $j = 2$, and  \eqref{1'.14} if $j= 3$ (with vanishing Dirichlet/Neumann boundary conditions in the last two cases), with $\dot W$ there replaced by $\dot W^i$.

If $d>6$, then points are polar for the following processes: $v^{(1)}$ restricted to $]0,\infty[ \times \re$,  $v^{(2)}$ restricted to $]0,\infty[ \times ]0,L[$, and $v^{(3)}$ restricted to $]0,\infty[ \times [0,L]$.
\item {\em Linear stochastic wave equation on $D=\re$, $D=\, ]0,\infty[$ or  $D=\, ]0,L[$}.

Let
$v^{(4)}=\left(v^{(4)}(t,x),\, (t,x)\in\re_+\times D\right)$
 be the $d$-dimensional random field  with independent components, where the $i$-th component $v_i^{(4)}$, $ i=1,\ldots,d$, is the solution to \eqref{wave-1'}--\eqref{3.4(*3)}, with $I_0\equiv 0$ and $W$ replaced by $W^i$.

If $d>4$, then points are polar for $v^{(4)}$ restricted to $]0,\infty[ \times D$.
\end{enumerate}

\end{thm}
The proof of this theorem relies on the following local statement (local in space-time).

\begin{prop}
\label{ch2'-s7-p3}
Let $v^{(j)}$, $j = 1,\dots, 4$, be as in Theorem \ref{ch2'-s7-theorem3}. Define
Define the sets of indices:
\begin{description}
\item{(1)} $I^{(1)}=[t_0,T]\times I_2$, where $t_0\in\, ]0,T[$ and $I_2$ is a compact interval of $\re$.
\item{(2)} $I^{(2)}=[t_0,T]\times [t_1,L-t_1]$, with $t_0\in]0,T[$, $L>0$ and $t_1\in\,]0,L/2[$.
\item{(3)} $I^{(3)}=[t_0,T]\times [0,L]$, with $t_0>0$, $L>0$.
\item{(4)} $I^{(4)}= [t_0,T]\times I_2$, where in the three instances of domains $D=\re$, $D=\, ]0,\infty[$, $D=\, ]0,L[$, the set $I_2$ is a compact interval in $D$.
\end{description}

If $d>6$, then points are polar for $v^{(j)}$ restricted to $I^{(j)}$, $j=1,2,3$. If $d>4$, then points are polar for $v^{(4)}$ restricted to $I^{(4)}$.
\end{prop}

\begin{proof} In each case, we will check that the assumptions (i) and (ii) of Proposition \ref{ch2'-s7-p1} hold and exhibit the values of $\alpha_0$ and $\alpha_1$ in \eqref{ch2'-s7.2}. Observe that in the examples under consideration, the dimension $k$ there equals one.
\medskip

 (1)  From Proposition \ref{ch1'-p0} we see that  assumption (i) holds with $\alpha_0= \frac{1}{4}$, $\alpha_1=\frac{1}{2}$. Thus, $Q = \frac{1}{\alpha_0} + \frac{1}{\alpha_1} = 6$.
Furthermore, recalling  \eqref{heatcauchy-11'}, we see that
\beqn
{\text{Var}} \left(v^{(1)}_1(t,x)\right) = \int_0^t dr \int_\re dz \, \Gamma^2(t-s, x-y) = \left(\frac{t}{2\pi}\right)^{\frac{1}{2}},\quad (t,x)\in\re_+\times \re.
\eeqn
Thus, for $\varepsilon_0>0$ small enough, we have
\beqn
\inf_{(t,x)\in (I^{(1)})^{(\varepsilon_0)}} {\text{Var}} \left(v^{(1)}_1(t,x)\right) > 0,
\eeqn
and therefore assumption (ii) is satisfied. This establishes the assertion of the proposition for $j=1$.
\medskip

 (2)  Assumption (i) holds with $\alpha_0= \frac{1}{4}$, $\alpha_1=\frac{1}{2}$ (see \eqref{4-p1-1'.1-p2}).  Thus, $Q=6$. For $\varepsilon_0$ small enough, $(t,x)\in (I^{(2)})^{(\varepsilon_0)}\mapsto {\text{Var}} \left(v^{(2)}_1(t,x)\right)$ is a continuous strictly positive function.
Indeed,
\beqn
{\rm Var}\left(v_1^{(2)}(t,x)\right) = \int_0^t ds \int_0^L G^2_L(s;x,y)\, dy >0
\eeqn
by Proposition \ref{ch1'-pPD} (ii). Furthermore, the continuity is a consequence of \eqref{4-p1-1'.0}.
Therefore, assumption (ii) holds. This completes the proof for $j=2$.
\medskip

 (3) Assumption (i) holds with $\alpha_0= \frac{1}{4}$, $\alpha_1=\frac{1}{2}$, as follows from \eqref{dia 5}.
Hence, as in the preceding two cases, $Q=6$. Using the expression \eqref{1'.N} (with $u_0\equiv 0$) for $v_1^{(3)}$, we see that for $\varepsilon_0$ small enough, $(t,x)\in (I^{(3)})^{(\varepsilon_0)}\mapsto {\text{Var}} \left(v^{(3)}_1(t,x)\right)$ is a continuous strictly positive function.
In fact, for $(t,x)\in[t_0,T]\times [0,L]$,
\beqn
{\text{Var}} \left(v^{(3)}_1(t,x)\right) \ge \int_0^t ds \int_0^L dy\ \Gamma^2(s, x-y) >0
\eeqn
by \eqref{neumann-bad}, and the continuity is a consequence of \eqref{dia 5}. This implies assumption (ii) and completes the proof for $j=3$.
\medskip

 (4) From item 1. of Theorem  \ref{ch4-pwave-1'}, we see that for the three instances of domains considered here, assumption (i) of Proposition \ref{ch2'-s7-p1} holds with $\alpha_0=\alpha_1=\frac{1}{2}$. Therefore, $Q=4$.

In the case $D=\re$, the validity of Assumption (ii) for $(I^{(4)})^{(\ep_0)}$ (with $\varepsilon_0$ small enough) is ensured by the identity
\beqn
\label{var1}
{\text{Var}} \left(v^{(4)}_1(t,x)\right) = \frac{1}{4} \int_{\re_+} dr \int_{\re} dz\, 1_{D(t,x)}(r,z) = \frac{1}{4} t^2 \ge \frac{1}{4}t_0^2,
\eeqn
which follows from \eqref{anotherone} and \eqref{gamma-cone-norm}.

If $D=\, ]0,\infty[$,
\beq
\label{var2}
{\text{Var}} \left(v^{(4)}_1(t,x)\right) = \frac{1}{4} \int_0^t dr \int_{\re_+} dz\, 1_{E(t,x)}(r,z),
\eeq
where $E(t,x)$ is defined in \eqref{Ewave}. Fix $(t,x)\in [t_0,T]\times [x_2,y_2]$, with $0<x_2<y_2$. The area of the quadrilateral $E(t,x)$ is $t^2$ if $t \leq x$ and is $ t^2 - (t-x)^2 = 2tx - x^2 \geq x^2$ if $t > x$. In either case, it is bounded below by $\min(t_0^2, x_2^2) > 0$. Hence,
${\text{Var}} \left(v^{(4)}_1(t,x)\right)\ge \tfrac{1}{4} \min(t_0x_2,t_0^2)$ and thus, assumption (ii) holds  (with $\varepsilon_0$ sufficiently small there).

Using this result and the relation between the Green's function corresponding to $D=\, ]0,\infty[$ and $D=\, ]0,L[$ (see \eqref{wave-bc1-100'}) we obtain the validity of assumption (ii) for $(I^{(4)})^{(\ep_0)}$ (with $\ep_0$ small enough) in the case $D=\, ]0,L[$. This completes the proof for $j=4$.
\end{proof}
\medskip

\noindent{\em Proof of Theorem \ref{ch2'-s7-theorem3}}.
The set $]0,\infty[\times \re$ can be decomposed into a countable union of sets such as $I^{(1)}$ in
Proposition \ref{ch2'-s7-p3} (1). Since if $d>6$, points are polar for $v^{(1)}$ restricted to each set of this covering,  we deduce that points are polar for $v^{(1)}$ restricted to $]0,\infty[\times \re$. The same argument applies to the remaining cases.
\qed
\medskip

The fact that points are polar for a process $u = (u_1,\dots, u_d)$ in {\em large enough} dimensions $d$ is natural, because in this case, there is lots of freedom of movement for the process $u$, and it is unlikely that any particular point will be visited by $u$. Proposition \ref{ch2'-s7-p1} and Theorem \ref{ch2'-s7-theorem3} make this precise by providing the condition $d > Q$ for polarity of points.
\medskip

\noindent{\em Non-polarity of points in low dimensions}
\medskip

It is natural to ask whether or not points are non-polar for $u$ if $d \leq Q$. This question has been considered by many authors. A classical result, going back to Kakutani and Dvoretzky, Erd\"os and Kakutani (\cite{kakutani-1944}, \cite{kakutani-1945}, \cite{dvoretzky-1950}) states that if $B=(B_t,\, t\in \R_+)$ is a $d$-dimensional Brownian motion (so $\alpha_0 = \half$ and $Q = 2$) and $A\subset \red$ is a compact set, then
 \beq
 \label{hittingb-1}
 P\{B(\re_+^*) \cap A\ne \emptyset\}>0 \quad \Longleftrightarrow\quad {\rm{Cap}}_{d-2}(A) >0,
 \eeq
 where ${\rm{Cap}}_{d-2}(A)$ denotes the Bessel-Riesz capacity of the set $A$ (see \cite[p. 376]{khosh2002} for a definition of this notion). As a consequence, since for any $z\in \red$,
 \beq
 {\rm{Cap}}_{d-2}(\{z\}) =
 \begin{cases}
 1&\ {\text{if}}\ d<2,\\
 0&\ {\text{if}}\  d\ge 2,
 \end{cases}
 \eeq
   points are polar for $(B_t,\, t \in \R_+^*)$ if and only if $d\ge 2$.
  In fact, \eqref{hittingb-1} may be obtained from the following statement:
  Fix $R>0$ and let $A\subset \red$ be a compact set included in $B_R(0)$ (the Euclidean ball centered at $0$ with radius $R$). Let $I=[a,b]$, $0<a<b$. Then
  there exists a constant $C$, depending on $d, R, I$, such that
  \beq
  \label{hitting-1-bis}
  \frac{1}{C}\ {\rm{Cap}}_{d-2}(A)\le P\{B(I) \cap A\ne \emptyset\}\le C\ {\rm{Cap}}_{d-2}(A).
  \eeq

  In \cite[Theorem 1.1]{khosh-shi-99}, these estimates are extended to an $\red$-valued Brownian sheet  $W=(W_{t_1, \ldots, t_k},\, (t_1,\ldots,t_k)\in\rek_+)$ (each of the coordinates of $W$ is a Brownian sheet as defined in Section \ref{ch1-2.3}, and these coordinate processes are independent) in the following way.

  Let $I=[a_1,b_1]\times \cdots\times[a_k,b_k]$, where $0<a_\ell<b_\ell<\infty$, $\ell=1,\ldots,k$. Fix $R>0$. There exists $0<C<\infty$, depending on $k,d,R, I$, such that for all compact sets $A\subset B_R(0)\subset \red$,
  \beq
  \label{hitting-k}
  \frac{1}{C}\ {\rm{Cap}}_{d-2k}(A)\le P\{W(I) \cap A\ne \emptyset\}\le C\ {\rm{Cap}}_{d-2k}(A).
  \eeq
In the same way as for $d$-dimensional Brownian motion, we deduce that points are polar
for $(W_{t_1,\dots, t_k},\, (t_1,\dots, t_k) \in\, ]0,\infty[^k)$ if and only $d \geq 2k$ (note that for the Brownian sheet, $\alpha_\ell = \half$, $\ell = 1, \dots, k$, and $Q = 2k$).

In the context of SPDEs, except for a particular nonlinear system of stochastic wave equations in spatial dimension $1$ studied in \cite{dn2004}, optimal estimates such as  \eqref{hitting-k} are not known. Typical statements provide lower bounds in terms of capacity, as in \eqref{hitting-1-bis} and \eqref{hitting-k}, while Hausdorff measure (see e.g. \cite{rogers} for a definition) replaces capacity in the upper bounds. We refer  to \cite{dkn07}, \cite{dkn09}, \cite{dss10}, \cite{dss15}, \cite{ssv18}, \cite{dalang-pu-2020}, \cite{dalang-pu-2021} for a collection of such results.

Comparing the values of the capacity and the Hausdorff measure of a singleton, ${\rm{Cap}}_{d-Q}(\{z\})$ and  $\mathcal{H}_{d-Q}(\{z\})$, respectively, we find that
\beqn
{\rm{Cap}}_{d-Q}(\{z\}) = \begin{cases}
1,\ {\text{if}}\ d<Q,\\
0, \ {\text{if}}\ d=Q,\\
0,\ {\text{if}}\ d>Q.\\
\end{cases}
\quad \mathcal{H}_{d-Q}(\{z\}) = \begin{cases}
\infty,\, {\text{if}}\ d<Q,\\
1, \ \ {\text{if}}\ d=Q,\\
0,\ \ {\text{if}}\ d>Q.\\
\end{cases}
\eeqn
Therefore, an inequality such as \eqref{hitting-k}, with ${\rm{Cap}}_{d-Q}(A)$ on the left-hand side and $\mathcal{H}_{d-Q}(A)$ on the right-hand side instead of ${\rm{Cap}}_{d-Q}(A)$, implies only that points are polar when $d>Q$ and non-polar when $d<Q$, and polarity or non-polarity of points at the critical dimension $d = Q$ (assuming that $Q \in \N^*$) remains undecided.

To end this digression, we mention that  for Gaussian random fields, a general approach for establishing polarity of points in the critical dimension is developed in \cite{d-m-x-2017}. In particular, their method applies to examples of systems of stochastic heat and wave equations considered in Theorem \ref{ch2'-s7-theorem3}, showing that points are polar in the critical dimensions $d=6$ and $d=4$, respectively.
For the stochastic heat equation only, this result was already established in \cite{m-t2002} using a different method.

\subsection{Polarity for SPDEs with nonlinear drift}
\label{ch2'-s7.s2}

In this section, we consider a generalization of the systems of SPDEs \eqref{ch2'-s7.14} discussed in Section \ref{ch2'-s7.e1}, allowing now an additional nonlinear drift term.
With this aim, we will consider a random vector of independent copies of the solution to \eqref{ch2'-s3.4.1.1}.

Let $D\subset \IR^k$ be a bounded or unbounded domain with smooth boundary. Set
\begin{align}
\label{ch2'-s7.170}
 u^i(t,x) &= I_0^i(t,x) + \int_0^t \int_D \Gamma(t,x;s,y)\,  W^i(ds,dy)\notag\\
&\qquad +  \int_0^t \int_D \Gamma(t,x;s,y) b^i(s,y,u^i(s,y))\, ds dy,
\end{align}
$i=1,\ldots,d$, $(t,x)\in \re_+\times D$, where the $W^i$ are independent space-time white noises.

We assume that  $\Gamma$ and all the $I_0^i$ and $b^i$, $i=1,\ldots,d$, satisfy the conditions of Theorem \ref{ch2'-s3.4.1-t1}.
An (easy) extension of Theorem \ref{ch2'-s3.4.1-t1} to systems of uncoupled equations yields the existence of a weak solution to \eqref{ch2'-s7.170}.

Consider the $d$-dimensional Gaussian vector
\beq
\label{ch2'-s7.20}
v^i(t,x) = I_0^i(t,x) + \int_0^t \int_D \Gamma(t,x;s,y)\,  W^i(ds,dy),\quad i=1,\ldots,d,
\eeq
$(t,x)\in \re_+\times D$, and let $\tilde P$ be the measure defined by
\begin{align}
\label{ch2'-s7.19}
\frac{d\tilde P}{dP}& = \exp\left(-\int_0^T\int_D b(s,y,u(s,y))\cdot W(ds,dy)\right. \notag\\
&\left.\quad\qquad \qquad -\frac{1}{2} \int_0^T\int_D  \vert  b(s,y,u(s,y))\vert^2\ ds dy\right),
\end{align}
where  $b = (b^i, i= 1, \dots, d)$, $W = (W^i, i= 1,\dots, d)$, and we have used the notation ``$\cdot$'' to recall the Euclidean inner product in $\R^d$.
From Proposition \ref{ch2'-s3.4.2-t1} and the independence of the $W^i$, we deduce that $E\left[\frac{d\tilde P}{dP}\right]=1$, therefore $\tilde P$ is a probability measure and, for any $T>0$, restricted to $[0,T]\times D$, the law of $u=(u^i,\, i=1,\ldots,d)$ under $\tilde P$ is the same as that of $v=(v^i, \, i=1,\ldots,d)$ under $P$.
As a consequence, we obtain in the next proposition that sets are polar for the random field $u$ if and only if they are polar for the Gaussian random field $v$.

\begin{prop}
\label{ch2'-s7-p5}
Let $u=(u^i(t,x),\, (t,x)\in [0,T]\times \re,\ i=1,\ldots,d)$ and $v=(v^i(t,x),\, (t,x)\in [0,T]\times \re,\ i=1,\ldots,d)$ be the random fields defined by \eqref{ch2'-s7.170} and \eqref{ch2'-s7.20}, respectively. Suppose that $\Gamma$, $I_0^i$ and $b^i$, $i=1,\ldots,d$, satisfy
the conditions of Theorem \ref{ch2'-s3.4.1-t1}.
 Let $I = [t_0, t_1] \times [x_1, x_2]$, where $0 \leq t_1 < t_2$ and $x_1 < x_2$. Then  $A\in \mathcal{B}(\red)$ is polar for $u$ restricted to $I$ if and only if it is polar for $v$ restricted to $I$.
\end{prop}

\begin{proof}
Fix $A\in \mathcal{B}(\red)$.
Observe that since $\tf$ is complete, $\{u(I)\cap A \neq \emptyset\}\in \tf$, because it is the projection onto $\Omega$ of the measurable set $\{(t,x,\omega)\in I\times \Omega:\ u(t,x,\omega)\in A\}$ (see e.g. \cite[Théorème, p. 252]{dm1}).
Observe that $A$ is polar for $u$ restricted to $I$ if and only if
\begin{align*}
\label{ch2'-s7.23}
0=P\{u(I)\cap A \neq \emptyset\} &= E_P\left[1_{\{u(I)\cap A \neq \emptyset\}}\right] = E_{\tilde P}\left[1_{\{u(I)\cap A \neq \emptyset\}}\frac{dP}{d\tilde P}\right].
\end{align*}
Since $\frac{dP}{d\tilde P} > 0$ a.s. by \eqref{ch2'-s7.19}, we deduce that this is equivalent to
\beqn
   0 = \tilde P \{ u(I) \cap A \neq \emptyset \} = P \{ v(I) \cap A \neq \emptyset \},
   \eeqn
where we have used the fact that, restricted to $I$, the law of $u$ under $\tilde P$ is the same as the law of $v$ under $P$.

 We deduce that $A$ is polar for $u$ restricted to $I$ if and only if it is polar for $v$ restricted to $I$.
\end{proof}

\begin{remark}
\label{5.4.3-r1}
Applying Proposition \ref{ch2'-s7-p5}, we see that the results on polarity for the examples of systems of SPDEs considered in Theorem \ref{ch2'-s7-theorem3}
 also hold for  the corresponding equations with a non vanishing drift, as in  \eqref{ch2'-s7.170}, provided that for $i=1,\ldots,d$,
 $I_0^i\equiv 0$ and $b^i$ satisfies the assumptions of Proposition \ref{ch2'-s7-p5}.
\end{remark}


\subsection[Sufficient conditions for polarity of points:\texorpdfstring{\\}{} the general case]{Sufficient conditions for polarity of points: the general case}
\label{ch2'-s7.s5.4.4}

We have seen that Proposition \ref{ch2'-s7-p1} provides a useful approach to the study of polarity of points for Gaussian processes and, in particular, for random field solutions to systems of linear SPDEs (see Sections \ref{ch2'-s7.e1} and \ref{ch2'-s7.s2}). For SPDEs with multiplicative noise, the random field solution is not Gaussian, and the study of polarity of points requires more sophisticated tools.
In this section, we give a brief account of results concerning polarity of points for the random field solutions to systems of nonlinear stochastic heat and wave equations, such as those studied in Chapter \ref{ch1'-s5}, to which we refer for the setting and assumptions. More precisely,
let $\mathcal{L}$ be a partial differential operator on $\re_+\times \rek$ and  $D\subset \rek$ be a bounded or unbounded domain with smooth boundary. We consider the system of nonlinear SPDEs
\beqn
\mathcal{L}u_i(t,x) = \sum_{j=1}^d \sigma_{i,j} (t,x,u(t,x))\, \dot W^j(t,x) + b_i(t,x,u(t,x)),
\eeqn
$(t,x)\in\, ]0,\infty[\times D$, $i=1,\ldots,d$,
with vanishing initial conditions and, if $D$ has boundaries, also with vanishing boundary conditions, and where $\dot W^j, j=1,\ldots,d$, are independent copies of a space-time white noise. %


As in Section \ref{Ch5-ss7.1}, we fix the subset of parameters $I=[t_0,t_1]\times[x_1,x_2]$, where
 $[x_1, x_2] = \prod_{\ell=1}^k [x_{1,\ell}, x_{2,\ell}]$, $0< t_0<t_1$,
$x_1 = (x_{1,1}, \dots, x_{1,k})$, $ x_2 = (x_{2,1}, \dots, x_{2,k})$, with $x_{1,\ell} < x_{2,\ell}$ ($\ell=1,\ldots k$). For $\varepsilon_0>0$, we denote by $I^{(\varepsilon_0)}$ the closed $\varepsilon_0$-neighbourhood of $I$.

We begin with a
proposition that provides sufficient conditions for polarity of points for random fields that are not necessarily Gaussian.

\begin{prop}
\label{ch2'-s7-p2}
Let $u=(u(t,x),\, (t,x)\in\re_+\times \rek)$ be a $d$-dimensional centred random field.

Suppose that for some $\ep_0 > 0$ and $\alpha_\ell \in\, ]0,1]$, $\ell = 0,1, \ldots, k$, the following
conditions are satisfied:
\begin{description}
\item{(i)}  Let $Q := \sum_{\ell=0}^k \tfrac{1}{\alpha_\ell}$.
Assume that for all $p > Q$, there is a constant $C = C_p$ such that for all $(t,x), (s,y) \in I^{(\ep_0)}$,
\beq
\label{ch2'-s7.8}
\left\Vert u(t,x)-u(s,y)\right\Vert_{L^p(\Omega)} \le C \left(|s-t|^{\alpha_0}+\sum_{\ell=1}^k|x_\ell-y_\ell|^{\alpha_\ell}\right).
\eeq
\item{(ii)} Fix $V \subset \re^d$.
   For any $(t,x) \in I^{(\ep_0)}$, $u(t,x)$ has a density $p_{t,x}$ and there is $C > 0$ such that
   \beq
   \label{ch2'-s7.8-bound}
       \sup_{z \in V^{(\ep_0)}}\, \sup_{(t,x) \in I^{(\ep_0)}} p_{t,x}(z) \leq C.
       \eeq
\end{description}
If $d>Q$, then all points $z \in V$ are polar for $u$ restricted to $I$.
\end{prop}

\begin{proof}
Fix $z\in V$ and assume $d>Q$.
We use the notations in the proof of Proposition \ref{ch2'-s7-p1} and, with the same arguments as there (in particular, up to and just after \eqref{ch2'-s7.4}) we see that if
for some $\gamma \in\, ]Q, d[$, there is $C > 0$ such that for all small enough $\ep > 0$ and all
$j\in J_{\varepsilon }= \{j: R_j^\ep \cap I \neq \emptyset\}$,
\beq
\label{ch2'-s7.10}
P\left\{\inf_{(t,x)\in R_j^\varepsilon} |u(t,x)-z|<\varepsilon\right\} \le C \varepsilon^\gamma,
\eeq
 then $\{z\}$ is polar for $u$ restricted to $I$.

Let $0<\varepsilon<\varepsilon_0$ be such that
$\cup_{j \in J_\varepsilon} R_j^{\varepsilon} \subset I^{(\varepsilon_0)}$
 and define
\beqn
\bar Y_j^\varepsilon = |u(y_j^\varepsilon)-z|, \qquad  \bar Z_j^\varepsilon = \sup_{(t,x)\in R_j^\varepsilon}|u(t,x)- u(y_j^\varepsilon)|.
\eeqn
By the reverse triangle inequality,
\begin{align}
\label{ch2'-s7.133}
P\left\{\inf_{(t,x)\in R_j^\varepsilon} |u(t,x)-z|<\varepsilon\right\} & \le P\left\{\bar Y_j^\varepsilon \le \varepsilon + \bar Z_j^\varepsilon\right\}\notag\\
&\le P\left\{\bar Z_j^\varepsilon\ge\frac{1}{2} \bar Y_j^\varepsilon\right\} + P\left\{\bar Y_j^\varepsilon\le 2\varepsilon\right\}.
\end{align}
By the definition of $\bar Y_j^\varepsilon$ and because of of Assumption (ii), for any $r\in\, ]0,\ep_0[$, we have
\beq
\label{ch2'-s7.134}
P\left\{\bar Y_j^\varepsilon\le r\right\}= P\left\{u(y_j^\varepsilon) \in B_{r}(z)\right\} \le C \ r^d.
\eeq


Using Assumption (i) and by Kolmogorov's continuity criterion Theorem \ref{ch1'-s7-t2}, we obtain for $p>Q$ and $\alpha\in\, ]\tfrac{Q}{p},1[$,
\beq
\label{ch2'-s7.12}
E\left[\vert \bar Z_j^\varepsilon\vert^p\right] \le C(\alpha,p)\, \varepsilon^{p\alpha-Q }
\eeq
(see \eqref{ch1'-s7.20}).

 Fix $\gamma\in\, ]Q, d[$
 and $\alpha \in\, ]\frac{Q}{p}, 1[$ such that  $\alpha > \frac{\gamma}{d}$. Next, we prove that
\beq
\label{ch2'-s7.13}
P\left\{\bar Z_j^\varepsilon \ge \frac{1}{2} \bar Y_j^\varepsilon\right\} \le C \varepsilon^{\gamma},
\eeq
where $C = C(p,\gamma,\alpha,d)$.
 Indeed, consider the decomposition
\beqn
P\left\{\bar Z_j^\varepsilon \ge \frac{1}{2} \bar Y_j^\varepsilon\right\} \le P\left\{\bar Y_j^\varepsilon\le \varepsilon^{\frac{\gamma}{d}}\right\} +
P\left\{\bar Z_j^\varepsilon \ge \frac{1}{2} \varepsilon^{\frac{\gamma}{d}}\right\}.
\eeqn
From \eqref{ch2'-s7.134}, for $\ep>0$ small enough, we have
\beq
\label{ch2'-s7.140}
P\left\{\bar Y_j^\varepsilon\le \varepsilon^{\frac{\gamma}{d}}\right\}\le C\, \varepsilon^\gamma
\eeq
and, by Chebychev's inequality along with \eqref{ch2'-s7.12},
\beqn
P\left\{\bar Z_j^\varepsilon \ge \frac{1}{2} \varepsilon^{\frac{\gamma}{d}}\right\}\le C(\alpha,p) \, \varepsilon^{p(\alpha-\tfrac{\gamma}{d})-Q}.
\eeqn

The above estimates yield
\beqn
P\left\{\bar Z_j^\varepsilon \ge \frac{1}{2} \bar Y_j^\varepsilon\right\} \le C(\alpha,p)\ \varepsilon^{\gamma\wedge (p(\alpha-\tfrac{\gamma}{d})-Q)}.
\eeqn
Since $\alpha>\frac{\gamma}{d}$, we can choose $p$ large enough so that $\gamma\wedge (p(\alpha-\tfrac{\gamma}{d})-Q) = \gamma$.
Thus \eqref{ch2'-s7.13} is proved.
Together with \eqref{ch2'-s7.134} and \eqref{ch2'-s7.133}, we obtain \eqref{ch2'-s7.10} for all sufficiently small $\ep > 0$ and $j\in J_{\varepsilon}$.

The proof of the proposition is complete.
\end{proof}

We devote the remainder of this section to a discussion of some applications of Proposition \ref{ch2'-s7-p2}.

Let $I$ and $\ep_0$ be such that $I^{(\ep_0)}\subset [0,T]\times D$. Assume that the hypotheses of Theorem \ref{ch1'-ss5.2-t1} hold with $I\times \tilde D$ there replaced by  $I^{(\ep_0)}$. The conclusion (a) of this theorem tells us that the random field solution to \eqref{ch1'-s5.1} restricted to $I^{(\ep_0)}$ satisfies the condition (i) of Proposition \ref
{ch2'-s7-p2} (see \eqref{ch1'-ss5.2-t1.1} for the specific values of $\alpha_\ell$, $\ell=0,\ldots, k$).

The methodology for the analysis of condition (ii) of Proposition \ref{ch2'-s7-p2} relies on Malliavin calculus (\cite{malliavin-1997}, \cite{nualart-1995}, \cite{sanz-sole-2005}, \cite{watanabe-1984}). When applied to systems of SPDEs such as those considered in Chapter \ref{ch1'-s5}, some
 additional hypotheses are required,  including smoothness of the coefficient functions $\sigma$ and $b$, and a non-degeneracy condition of elliptic or hypoelliptic type relative to $\sigma$ and $\mathcal{L}$.
There are two problems to address. First, the very existence of the density and second, the uniform bound \eqref{ch2'-s7.8-bound}. For the former, the probabilistic approach to H\"ormander's theorem on hypoellipticity of partial differential operators given in \cite{malliavin-1978} provides a suitable method (see also \cite{nualart-1995} and \cite{sanz-sole-2005}). For the latter, one can use Watanabe's formula  for the density (\cite{watanabe-1984}, \cite[Corollary 3.2.1, p. 161]{nualart-1998}, \cite[Proposition 5.2, p. 63]{sanz-sole-2005}). In many examples, this procedure leads to a verification of condition (ii) of Proposition \ref{ch2'-s7-p2}. Without aiming to be exhaustive, we present a small sample of references related to the examples considered in this book. For systems of stochastic heat equations: \cite{bally-pardoux1998}, \cite{dkn09}; for systems of stochastic wave equations: \cite{dn2004}, \cite{millet-SanzSole-1999}, \cite{dss15}; for the stochastic fractional heat equation: \cite{Pu-thesis}, \cite{dalang-pu-2021}.

Motivated by the examples considered in Theorem \ref{ch2'-s7-theorem3}, it is natural to ask whether the critical dimensions for polarity, that we have found for various examples of linear SPDEs, are preserved for the corresponding nonlinear SPDEs. This is in fact the case (see \cite{dn2004} and \cite{dkn09}).

We do not go into the details of checking condition (ii) of Proposition \ref{ch2'-s7-p2}, but just focus on condition (i) there.
We are assuming  that the initial and boundary conditions vanish. Hence, we see from Theorem \ref {cor-recap} (b) that for the three instances of systems of nonlinear stochastic heat equations in Case 1.~of Theorem \ref{ch2'-s7-theorem3}, assumption (i)  holds with $\alpha_0= \frac{1}{4}$, $\alpha_1= \frac{1}{2}$. Thus $Q=6$.
For systems of stochastic wave equations (Case 2.~in Theorem \ref{ch2'-s7-theorem3}), appealing to Theorem \ref{cor-recapwave}, we deduce the validity of  assumption (i) with $\alpha_0=\alpha_1= \frac{1}{2}$. Thus $Q=4$.

Finally, in the critical dimension $d = Q$ (assuming that  $Q \in \N^* $), except for the nonlinear stochastic wave equation studied in \cite{dn2004}, the issue of polarity of specific points in $\R^d$ remains undecided. For systems of nonlinear stochastic heat equations with multiplicative space-time white noise, a step in this direction is provided by \cite{d-m-x-2020}, in which polarity of ``almost all'' points in $\R^d$ is proved for $d = Q = 6$.


\section{SPDEs with rough initial conditions}
\label{rdrough}

In many situations, it is natural to consider initial conditions which are ``rougher" than those allowed by Assumption ($\bf H_I$), such as unbounded functions or even measures. For instance, in the parabolic Anderson model\index{Anderson model}\index{model!Anderson} (see Section \ref{rd1+1anderson}), the initial condition $u_0(x) = \delta_0(x)$, $x \in \R$, where $\delta_0$ denotes the Dirac delta function, is considered particularly interesting (see e.g.~\cite{konig2016}): it is called the ``narrow wedge condition" and arises for instance in the study of fluctuations in models within the KPZ universality class (see e.g.~\cite{ACQ}). For this initial condition, the solution to the homogeneous heat equation is $I_0(t, x) = \Gamma(t, x)$, where $\Gamma$ is the heat kernel \eqref{heatcauchy-1'}. In particular, assumptions ($\bf H_I$) and (${\bf H_{I,\infty}}$) {\em are not} satisfied. In Theorem \ref{rd08_20t1} below, we extend Theorem \ref{ch1'-s5.t1} to include this more general situation.

Let $\cL$ be a partial differential operator on $\R_+ \times \R$ with associated fundamental solution $\Gamma(t, x; s, y)$ that satisfies Hypothesis (${\bf H_{\Gamma,\infty}}$) (see Section \ref{ch2'-section5.2}). On a complete probability space $(\Omega, \tf, P)$, we consider a space-time white noise $W$ along with a right-continuous complete filtration $(\tf_s,\, s\in\re_+)$, as in the beginning of Sections \ref{ch4-section0} and \ref{ch1'-s7}. Consider the SPDE
\beq
\label{rd08_21e1}
    \cL u(t, x) = \sigma(t, x, u(t, x)) \, \dot W(t, x),
\eeq
$(t,x)\in\re_+^*\times \re$,
subject to suitable initial conditions.

We now give the definition of random field solution to \eqref{rd08_21e1}. It is a slight modification of Definition \ref{def4.1.1} adapted to the context of (possibly) measure-valued initial conditions.

\begin{def1}
\label{rd08_21d1}
A {\em random field solution}\index{solution!random field}\index{random field!solution} to \eqref{rd08_21e1} is a jointly measurable and adapted real-valued random field
\beqn
    u= \left(u(t,x),\, (t,x)\in \R_+^*\times \R\right)
\eeqn
such that, for all $(t,x)\in \R_+^* \times \R$, 
 \beq\label{rd01_22e1}
   E\left[\Vert \Gamma(t, x; \cdot, *) \sigma(\cdot, *, u(\cdot, *)) \Vert_{L^2([0, t] \times \R)}^2\right] < \infty
\eeq
and
\begin{align}\nonumber
     u(t, x) &= I_0(t, x) \\
     &\qquad + \int_0^t \int_\R \Gamma(t, x; s, y)\, \sigma(s, y, u(s, y))\, W(ds, dy)\quad a.s.,
\label{rd08_21e4}
\end{align}
where $I_0$ is the solution to the homogeneous PDE $\cL I_0 = 0$ with the same initial conditions as in \eqref{rd08_21e1}.
\end{def1}
 
 We introduce some notions and notations. Let $J: \re_+^* \times \re \to \re_+$ be a nonnegative function such that
\beq
\label{conv-series-1}
\int_{\re_+}ds \int_{\re} dy\, J(s,y) < \infty.
\eeq
 Set $J^{\star 1}(t,x) = J(t,x)$, and for $n\ge 1$,
\begin{align}\nonumber
   J^{\star(n+1)}(t,x) &= [J^{\star n} \star J](t, x) \\
   & = \int_0^t ds \int_{\re} dy\, J^{\star n}(t-s,x-y)\, J(s,y),\qquad n\ge 1,
\label{space_time_conv}
\end{align}
where for functions $f,g: \re_+\times \re \to \re_+$, the space-time convolution $f \star g$ is defined in \eqref{rd08_14e8}.  

\subsection{Existence, uniqueness and moments}
\label{ss7.1.1}
Instead of Assumption ($\bf H_{\Gamma}$) of Section \ref{ch4-section0}, we will use the following weaker assumption ($\bf H_{\Gamma, I}$) below, that links the fundamental solution $\Gamma$ and the deterministic function $I_0$. We do not specify how $I_0$ is to be computed, since this is a deterministic problem in PDEs (see for instance Example \ref{rd08_21ex1}). However, Assumption ($\bf H_{\Gamma, I}$) and Theorem \ref{rd08_20t1} below encompass all the cases discussed for instance in \cite{chen}--\cite{chendalang2015}.
\bigskip

\noindent ($\bf H_{\Gamma, I}$) {\em Assumptions on the fundamental solution and the initial conditions}\label{rdHGammaI}
\begin{description}
\item{(i)} The mapping $(t,x;s,y)\mapsto \Gamma(t,x;s,y)$ from $\{(t,x;s,y)\in \R_+^* \times \R \times \R_+\times \R : 0 \leq s < t \}$ into $\IR$ is jointly measurable.
\item{(ii)} There is a Borel function $H: \R_+^* \times \R \to \IR_+$ such that
\beqn
\vert\Gamma(t,x;s,y)\vert \le H(t-s,x-y), \qquad 0\le s<t, \quad x,y\in \R.
\eeqn
\item{(iii)} For all $T > 0$,
\beq
\label{rd08_22e6}
\int_0^T ds \int_\R dy\,  H^2(s,y) < \infty.
\eeq
\item{(iv)} For all compact sets $K \subset \R_+^* \times \R$,
\beq
\label{rd08_21e2}
    \sup_{(t, x) \in K} \left[H^2 \star I_0^2 \right](t, x) < \infty.
\eeq
\end{description}
Notice that assumptions (i) - (ii) are already present in ($\bf H_{\Gamma}$), while (iii) and (iv) are not.

\begin{thm}
\label{rd08_20t1}
Consider the SPDE \eqref{rd08_21e1} and suppose that the function $(t,x)\mapsto I_0(t,x)$ from $\R_+^* \times \R$ into $\R$ is $\cB_{\R_+^*\times \R}$-measurable and that $\sigma$ satisfies $({\bf H_{L,\infty}})$ (see Section \ref{ch2'-section5.2}). Assume that the fundamental solution $\Gamma$ and the initial conditions are such that $(\bf H_{\Gamma, I})$ holds. Let  
 \beq
 \label{defofj}
 J(t, x) =2\, \bar c^2\, H^2(t, x),\quad (t,x)\in\re_+^*\times \re,
 \eeq
where $H(t, x)$ comes from $(\bf H_{\Gamma, I})$ and $\bar c$ comes from $({\bf H_{L,\infty}})$. 

Fix an even integer $p \geq 2$ and assume, in addition to $(\bf H_{\Gamma, I})$, that for all $T > 0$ and $ c > 0$,
\beq\label{rd12_19_e3}
     \sum_{\ell = 1}^\infty  \, \sup_{(t, x) \in [1/T, T] \times [-T, T]}\left(c^\ell \, [J^{\star \ell} \star (\bfone + I_0^2)](t, x)\right)^{\frac{1}{p}} < \infty,     
\eeq     
where $\bfone$ is the constant function equal to $1$.
Then there exists a random field solution $(u(t, x),\, (t, x) \in \R_+^* \times \R)$ to \eqref{rd08_21e1}, and there is a finite constant $C = C_{p}$  such that for all $(t, x) \in \R_+^* \times \R$,
\beq
\label{rd08_20e15}
   \Vert u(t, x)\Vert_{L^p(\Omega)}^2 \leq 2 I_0^2(t, x) + \left[\left(\sum_{\ell=1}^\infty C^\ell\, J^{\star \ell} \right) \star (1 + 2 I_0^2) \right](t, x).
\eeq
In addition, for all $(t, x) \in \R_+^* \times \R$,
\beq\label{rd01_31e4}
    \Vert u(t, x) \Vert_{L^2(\Omega)}^2 \leq  I_0^2(t, x) + \left[\left(\sum_{\ell=1}^\infty J^{\star \ell}\right) \star (\bfone + I_0^2)\right](t, x),
\eeq
and the right-hand sides of \eqref{rd08_20e15} and \eqref{rd01_31e4} are finite for a.a.~$(t, x) \in \R_+^* \times \R$. If one of these right-hand sides is infinite for some $(t, x) \in \R_+^* \times \R$, then $I_0(t, x) = \pm\infty$ and a.s., $u_0(t, x) = I_0(t, x) \in \{ \pm\infty\}$.
This solution is unique (in the sense of versions) among random field solutions that satisfy \eqref{rd01_31e4}.
\end{thm}

\begin{remark}
\label{rd08_23r1}
(a) We will see in the proof of this theorem that conditions $(\bf H_{\Gamma, I})$ (iii) and (iv) are precisely the conditions needed for the very first step of the Picard iteration scheme to be well-defined.

(b) For the stochastic heat, wave and fractional heat equations on $\R$, in which the fundamental solutions are of the form $\Gamma(t-s, x-y)$ with $\Gamma  \geq 0$, it will be convenient to choose $H:= \Gamma$ (see the proofs of Propositions \ref{rough-heat}--\ref{rough-fractional-heat}).

(c) While condition \eqref{rd12_19_e3} may seem to be somewhat difficult to check, we will see in Proposition \ref{rd01_07p1} below that there are simpler sufficient conditions that are much easier to verify and that are satisfied in many interesting examples.

(d) The situation $I_0(t_0, x_0) = \pm\infty$ for some $(t_0, x_0) \in \R_+^* \times \R$ 
occurs in certain examples (see for instance Proposition \ref{rd08_25p2} (a)). 
\end{remark}

\medskip

\noindent{\em Towards the proof of Theorem \ref{rd08_20t1}}
\medskip

In order to prove Theorem \ref{rd08_20t1}, we need to be able to bound $L^p(\Omega)$-moments of stochastic integrals without using norms such as those in \eqref{bzp}. A result in this direction is provided by the following lemma (\cite[Lemma 3.4]{foondun-khoshnevisan-2009}, \cite[Lemma 2.4]{ConusKhosh2012}), which uses Burkholder's inequality to bound $L^p(\Omega)$-moments of certain stochastic integrals when $p$ is an even integer, and provides a bound that is sharper than that for arbitrary $p >0$.

\begin{lemma}
\label{rd08_22l1}
We assume that the conditions (i)--(iii) of Assumption $({\bf H_{\Gamma,I}})$ are satisfied. Let $(Z(t, x),\, (t,x)\in\re_+\times \re)$ be a jointly measurable and adapted random field such that, for any $(t, x) \in \R_+^* \times \R$, 
\beqn
    \int_0^t ds \int_\R dy \, \Gamma^2(t, x; s, y) \,  \Vert Z(s, y) \Vert_{L^2(\Omega)}^2  < \infty.  
\eeqn
Define
\beqn
   \cI(t, x) = \int_0^t \int_\R \Gamma(t, x; s, y)\, Z(s, y)\, W(ds, dy),
\eeqn
let $a(t, x)$ be a deterministic real-valued function and let $(w(t, x))$ be a random field such that, for all $(t, x) \in \R_+^* \times \R$,
\beqn
     w(t, x) = a(t, x) + \cI(t, x).
\eeqn

   (a) For all $(t, x) \in \R_+^* \times \R$,
\beq
\label{rd08_22e12}
  \Vert w(t, x) \Vert_{L^2(\Omega)}^2 \leq a^2(t, x) + [H^2 \star \Vert Z(\cdot, *)\Vert_{L^2(\Omega)}^2 ](t, x).
\eeq

  (b) For all even integers $p \geq 2$, let $\tilde C_p$ be the constant that appears in Burkholder's inequality \eqref{ch1'-s4.4.mod}. Then for all $(t, x) \in \R_+^* \times \R$,
\beq
\label{rd08_22e10}
   \Vert \cI(t, x) \Vert_{L^p(\Omega)}^2 \leq (\tilde C_p)^{\frac{2}{p}}\, \left[H^2 \star \Vert Z(\cdot, *) \Vert_{L^p(\Omega)}^2 \right] (t, x)
\eeq
and
\beq
\label{rd08_22e13}
    \Vert w(t, x) \Vert_{L^p(\Omega)}^2 \leq 2\, a^2(t, x) + 2\, (\tilde C_p)^{\frac{2}{p}} \left[H^2 \star \Vert Z(\cdot, *) \Vert_{L^p(\Omega)}^2\right](t, x) .
\eeq
\end{lemma}

\begin{proof}
(a) For $p = 2$, by the Itô isometry, 
\begin{align}\nonumber
   \Vert w(t, x) \Vert_{L^2(\Omega)}^2 &= a^2(t, x) + \Vert \cI(t, x) \Vert_{L^2(\Omega)}^2 \\  \nonumber
      & =  a^2(t, x) + \left\Vert \Gamma(t, x; \cdot, *) \, \Vert Z(\cdot, *)\Vert_{L^2(\Omega)} \right\Vert_{L^2([0, t] \times \R)}^2 \\
   &\leq a^2(t, x) + \left[H^2 \star \Vert Z(\cdot, *) \Vert_{L^2(\Omega)}^2 \right](t, x),
\label{rd01_20e1}
\end{align}
which proves \eqref{rd08_22e12}.

   (b) If $p > 2$, then
\begin{align}\nonumber
   \Vert w(t, x) \Vert_{L^p(\Omega)}^2 &\leq \left(\vert a(t, x) \vert +  \Vert \cI(t, x) \Vert_{L^p(\Omega)}\right)^2 \\
    &\leq 2\, a^2(t, x) + 2\, \Vert \cI(t, x) \Vert_{L^p(\Omega)}^2.
\label{rd08_22e11}
\end{align}
If $p = 2 m$, for some $m \in \N^*$, then for fixed $(t, x) \in \R_+^* \times \R$, we apply Burkholder's inequality \eqref{ch1'-s4.4.mod} to obtain
\begin{align*}\nonumber
   \Vert \cI(t, x) \Vert_{L^p(\Omega)}^p &\leq  \tilde C_p\, E\left[\left(\int_0^t ds \int_\R dy\, \Gamma^2(t, x; s, y)\, Z^2(s, y) \right)^{\frac{p}{2}}\right] \\ \nonumber
   &= \tilde C_p\, E\left[\prod_{n=1}^{m} \int_0^t ds_n \int_\R dy_n\, \Gamma^2(t, x; s_n, y_n)\, Z^2(s_n, y_n) \right] \\ \nonumber
   &= \tilde C_p\, \int_0^t ds_1 \int_\R dy_1\, \cdots \int_0^t ds_m \int_\R dy_m\, \left( \prod_{n=1}^{m}\, \Gamma^2(t, x; s_n, y_n)\right)\\
    &\qquad\qquad\qquad\qquad\times E\left[Z^2(s_1, y_1)\, \cdots\, Z^2(s_m, y_m)\right].
\end{align*}
By Hölder's inequality, the expectation is bounded above by
\beqn
    \prod_{n=1}^{m} \left(E\left[\vert Z(s_n, y_n)\vert^p\right]\right)^{\frac{2}{p}} = \prod_{n=1}^{m} \Vert Z(s_n, y_n) \Vert_{L^p(\Omega)}^{2}.
\eeqn
Therefore,
\begin{align}\nonumber
 \Vert \cI(t, x) \Vert_{L^p(\Omega)}^p &\leq \tilde C_p\, \int_0^t ds_1 \int_\R dy_1 \, \cdots \int_0^t ds_m \int_\R dy_m\, \left[ \prod_{n=1}^{m}\, H^2(t-s_n, x-y_n)\right]\\ \nonumber
 &\qquad\qquad\qquad\qquad\times \left[\prod_{n=1}^{m}\, \Vert Z(s_n, y_n) \Vert_{L^p(\Omega)}^{2} \right] \\ \nonumber
    & = \tilde C_p\, \prod_{n=1}^{m}\, \left[\int_0^t ds_n \int_\R dy_n\, H^2(t-s_n, x-y_n)\, \Vert Z(s_n, y_n) \Vert_{L^p(\Omega)}^{2}\right] \\
    &= \tilde C_p\, \left[\int_0^t ds \int_\R dy\, H^2(t-s, x-y)\, \Vert Z(s, y) \Vert_{L^p(\Omega)}^{2}\right]^{\frac{p}{2}}. 
\label{rd08_22e9}
\end{align}
This proves \eqref{rd08_22e10}.

By \eqref{rd08_22e11} and \eqref{rd08_22e10},
\begin{align}\nonumber
   \Vert w(t, x) \Vert_{L^p(\Omega)}^2 & \leq \left(\vert a(t, x) \vert + \Vert \cI(t,x)\Vert_{L^p(\Omega)} \right)^2 \\
   &\leq 2\, a^2(t, x)+ 2\, (\tilde C_p)^{\frac{2}{p}} \left[H^2 \star \Vert Z(\cdot, *) \Vert_{L^p(\Omega)}^2\right] (t, x).
\label{rd01_20e2}
\end{align}
This proves \eqref{rd08_22e13}.
\end{proof}

\noindent{\em Proof of Theorem \ref{rd08_20t1}.}
In order to establish this theorem, we look back to the proof of Theorem \ref{ch1'-s5.t1} (to which we refer concerning measurability issues) and the Picard iteration scheme $(u^n(t, x),\, (t, x) \in \R_+^* \times \R)$, $n \geq 0$, defined in \eqref{rde4.1.2}: recall that for $(t, x) \in \R_+^* \times \R$, $u^0(t, x) = I_0(t, x)$, and for $n \geq 1$, if $\vert I_0(t, x) \vert < + \infty$, then
\begin{align}\nonumber
   u^n(t, x) &= I_0(t, x) \\
   &\qquad + \int_0^t \int_\R  \Gamma(t, x; s, y)\, \sigma(s, y, u^{n-1}(s, y))\, W(ds, dy),
\label{rd05_01e1}
\end{align}
and if $\vert I_0(t, x) \vert = + \infty$, then we set $u^n(t, x) = I_0(t, x)$ ($\in \{\pm \infty\}$).

 A uniform bound such as \eqref{ch1'-s5.2a} is in general no longer possible (for instance, in the parabolic Anderson model with Dirac delta initial condition, \eqref{ch1'-s5.2a} already fails for $n=0$), and the quantity $M_n(t)$ defined in \eqref{rd08_20e8} may also be infinite here. However, it will be possible to use the space-time Gronwall-type Lemma \ref{rd08_10l2}.

As in \eqref{rd08_10e7} and \eqref{rd08_10e8} from Section \ref{app3-s1}, for $(t, x) \in \R_+^* \times \R$,  we set $\cK_0(t,x) = 0$,
\beq
\label{def-kas}
\cK_n(t,x) =\sum_{\ell=1}^nJ^{\star \ell}(t, x) \quad \text{and} \quad \cK(t,x) =\sum_{\ell=1}^\infty J^{\star \ell}(t, x).
\eeq
We will refer to $\cK(t,x)$ as the 
{\em geometric space-time-convolution series} associated with $J$ (see Section \ref{rd08_26s1}).

The first step consists in proving, by induction on $n$, that for all $n \in \N^*$ and $(t, x) \in \R_+^* \times \R$,
\beq\label{rd04_30e1}
   E\left[ \int_0^t ds \int_\R dy\, \left(\Gamma(t, x; s, y)\, \sigma(s, y, u^{n-1}(s, y))\right)^2 \right] < \infty
\eeq
and
\beq
\label{rd08_22e2}
\Vert u^{n}(t,x)\Vert_{L^2(\Omega)}^2 \leq I_0^2(t, x) + [\cK_{n} \star (\bfone + I_0^2)](t, x),
\eeq
with $\cK_{n}$ given in \eqref{def-kas}.

Indeed, note that $u^0(t, x) = I_0(t, x)$ and that in the first step of the Picard iteration scheme, one defines
\beqn
    u^1(t, x) = I_0(t, x) + \cI^1(t, x),
\eeqn
where
\beq
\label{rd08_21e11}
   \cI^1(t, x) = \int_0^t \int_\R \Gamma(t, x; s, y)\, \sigma(s, y, I_0(s, y)) \, W(ds, dy).
\eeq
Let $\bar c$ be the constant in (${\bf H_{L,\infty}}$). By ($\bf H_{\Gamma, I}$) (ii) and (${\bf H_{L,\infty}}$),
\begin{align} 
\label{rd08_22e8}
   \Vert\cI^1(t, x)\Vert_{L^2(\Omega)}^2&=\int_0^t ds \int_\R dy\, \Gamma^2(t, x; s, y)\, \sigma^2(s, y, I_0(s, y))\notag \\ 
   &\qquad \leq 2\, \bar c^2 \int_0^t ds \int_\R dy\, H^2(t-s, x-y)\, (1 + I_0^2(s, y))\notag \\
   &\qquad = 2\, \bar c^2\, [H^2 \star \bfone](t, x) + 2\, \bar c^2\, \left[H^2 \star I_0^2 \right](t, x).
\end{align}
The first term is finite by Hypothesis ($\bf H_{\Gamma, I}$) (iii) because $[H^2 \star \bfone](t, x) = \Vert H \Vert_{L^2([0, t] \times \R)}^2$, therefore condition \eqref{rd08_21e2} is precisely the condition needed for the second term to be finite and therefore, for \eqref{rd04_30e1} with $n=1$ to hold and for the stochastic integral in \eqref{rd08_21e11} to be well-defined. From \eqref{rd08_22e8}, we obtain
\begin{align}
\label{rd08_23e1}
 \Vert u^1(t,x)\Vert_{L^2(\Omega)}^2 &= I_0^2(t, x) + \Vert\cI^1(t, x)\Vert_{L^2(\Omega)}^2 \notag\\
 &\leq I_0^2(t,x) + [J \star (\bfone + I_0^2)](t, x).
\end{align} 
This is \eqref{rd08_22e2} for $n=1$. 

Assume that  \eqref{rd08_22e2} holds with $n$ there replaced by $n-1$ ($n \geq 2$). From \eqref{rd01_20e1} and using $({\bf H_{L,\infty}})$, 
we deduce  that
\begin{align}\nonumber
 & \Vert u^{n}(t,x)\Vert_{L^2(\Omega)}^2   \\
  &\qquad = I_0^2(t,x)+ \left\Vert \Gamma(t, x; \cdot, *)\, \Vert \sigma(\cdot, *, u^{n-1}(\cdot, *))\Vert_{L^2(\Omega)}\right\Vert_{L^2([0, t] \times \R)}^2 .
  \label{rd05_01e2}
 \end{align}
 Now
 \begin{align}\nonumber
  &\left\Vert \Gamma(t, x; \cdot, *)\, \Vert \sigma(\cdot, *, u^{n-1}(\cdot, *))\Vert_{L^2(\Omega)}\right\Vert_{L^2([0, t] \times \R)}^2 \\ \nonumber
     &\qquad \leq \bar c^2 \left\Vert \Gamma(t, x; \cdot, *) \, \Vert \bfone + \vert u^{n-1}(\cdot, *)\vert \Vert_{L^2(\Omega)} \right\Vert_{L^2([0, t] \times \R)}^2 \\  \nonumber 
     &\qquad \leq \bar c^2 \left\Vert \Gamma(t, x; \cdot, *) \, \left(\bfone + \Vert u^{n-1}(\cdot, *) \Vert_{L^2(\Omega)} \right)\right\Vert_{L^2([0, t] \times \R)}^2 \\ \nonumber
     &\qquad  \leq 2 \bar c^2 \Big[ \Vert \Gamma(t, x; \cdot, *) \Vert_{L^2([0, t] \times \R)}^2 \\ \nonumber
     &\qquad  \qquad\qquad \qquad\qquad + \left\Vert \Gamma(t, x; \cdot, *)\, \Vert u^{n-1}(\cdot, *) \Vert_{L^2(\Omega)}  \right\Vert_{L^2([0, t] \times \R)}^2 \Big] \\ 
   &\qquad  \leq 2 \bar c^2 \left[H^2 \star \left(\bfone + \Vert u^{n-1}(\cdot, *)\Vert_{L^2(\Omega)}^2\right)\right](t, x).
 \label{rd01_07e1}
\end{align}
Using \eqref{rd08_22e2} (with $n$ there replaced by $n-1$), we see that
\begin{align}\nonumber
 & \left\Vert \Gamma(t, x; \cdot, *)\, \Vert \sigma(\cdot, *, u^{n-1}(\cdot, *))\Vert_{L^2(\Omega)}\right\Vert_{L^2([0, t] \times \R)}^2 \\
  &\qquad\qquad \leq   \left[J \star \left(\bfone + I_0^2(t, x) + \cK_{n-1} \star (\bfone + I_0^2) \right)\right](t, x)\notag \\
   &\qquad\qquad= \left[(J + J \star \cK_{n-1}) \star \left(\bfone + I_0^2\right) \right](t, x)\notag \\
   &\qquad\qquad= +  \left[\cK_n \star \left(\bfone + I_0^2\right)  \right](t, x)
\label{l2norm}
\end{align}
by \eqref{rd08_12e1}. Therefore, by \eqref{rd05_01e2}, \eqref{rd08_22e2} holds for all $n \in \N^*$.

We are assuming \eqref{rd12_19_e3}. Thus, we deduce from \eqref{l2norm} that for all $n \in \N^*$,
\beqn
     E\left[\left\Vert \Gamma(t, x; \cdot, *) \sigma(\cdot, *, u^{n-1}(\cdot, *))  \right\Vert_{L^2([0, t] \times \R)}^2 \right]  < \infty.
\eeqn 
This shows that the Picard iteration scheme $(u^n(t, x),\, n \in \N)$  is well-defined. In addition, using \eqref{rd05_01e2}, we see that \eqref{rd08_22e2} holds for all $n \in \N^*$.

For the proof of \eqref{rd08_20e15}, we need an extension of  \eqref{rd08_22e2} to the $L^p$-norm, where $p\geq 2$ is the even integer that appears in \eqref{rd12_19_e3}. Let $p=2m$, where $m\in\N^*$, let $\tilde C_p$ be the constant from Burkholder's inequality \eqref{ch1'-s4.4.mod}, and define 
\beqn
J_{0, p}(t, x) = 4\, \bar c^2 \, (\tilde C_p)^{\frac{2}{p}}\, H^2(t, x) = 2 \, (\tilde C_p)^{\frac{2}{p}}\, J(t, x),
\eeqn
 so that
\beq
\label{rd08_25e2}
   J_{0, p}^{\star \ell}(t, x) = \left(2 \, (\tilde C_p)^{\frac{2}{p}}\right)^\ell\, J^{\star \ell}(t, x).
\eeq

Set $\cK_{0, p}(t,x) = 0$ and, for $n \geq 1$,
\beqn
   \cK_{n, p}(t, x) = \sum_{\ell =1}^n\, J_{0, p}^{\star \ell}(t, x)
\quad \text{and} \quad
   \cK_p(t, x) = \sum_{\ell=1}^\infty\, J_{0, p}^{\star \ell}(t, x).
\eeqn
We will prove by induction on $n$ that  for all $n \in \N^*$,
\beq
\label{rd08_22e1}
    \Vert u^{n}(t, x) \Vert_{L^p(\Omega)}^2 \leq 2 I_0^2(t, x) + [\cK_{n, p} \star (\bfone + 2\, I_0^2)](t, x).
\eeq
Indeed, by \eqref{rd08_22e13} and $({\bf H_{L,\infty}})$,
\begin{equation}\label{rd05_02e1}
   \Vert u^1(t, x) \Vert_{L^p(\Omega)}^2 \leq 2\, I_0^2(t, x) + 2\, (\tilde C_p)^{\frac{2}{p}}\, [H^2 \star \sigma^2(\cdot, *, I_0(\cdot, *))](t, x)
\end{equation}   
and
\begin{align}    
 & 2\, (\tilde C_p)^{\frac{2}{p}}\, [H^2 \star \sigma^2(\cdot, *, I_0(\cdot, *))](t, x)   \notag\\ 
   &\qquad\qquad  \leq 2\, (\tilde C_p)^{\frac{2}{p}}\, [H^2 \star (2 \bar c^2\, \bfone + 2 \bar c^2 \,I_0^2)](t, x)\notag \\
   &\qquad\qquad=  4 \, \bar c^2 \, (\tilde C_p)^{\frac{2}{p}}\, [H^2 \star (\bfone + I_0^2)](t, x)\notag\\
   &\qquad\qquad= [J_{0, p} \star (\bfone + I_0^2)](t, x)\notag \\
   &\qquad\qquad\leq [J_{0, p} \star (\bfone + 2\, I_0^2)](t, x).
\label{rd08_23e3}   
\end{align}
By \eqref{rd05_02e1}, this gives \eqref{rd08_22e1} for $n=1$. 

Assume that \eqref{rd08_22e1} holds with $n$ there replaced by $n-1$ ($n\geq 2$). From \eqref{rd05_01e1} and \eqref{rd08_22e13}, we have
\begin{align}\nonumber
   &\Vert u^{n}(t,x)\Vert_{L^p(\Omega)}^2 \\
   &\qquad\leq  2\, I_0^2(t,x) + 2 (\tilde C_p)^{\frac{2}{p}}\, \left[H^2 \star \Vert \sigma(\cdot, *, u^{n-1}(\cdot, *))\Vert_{L^p(\Omega)}^2\right](t, x).
\label{rd05_02e2}   
\end{align}
Using $({\bf H_{L,\infty}})$, we see that
\begin{align*}
 &2 (\tilde C_p)^{\frac{2}{p}}\, \left[H^2 \star \Vert \sigma(\cdot, *, u^{n-1}(\cdot, *))\Vert_{L^p(\Omega)}^2\right](t, x) \\
   &\qquad\qquad \leq  2 (\tilde C_p)^{\frac{2}{p}}\, \left[H^2 \star (2 \bar c^2)\left( \bfone + \Vert u^{n-1}(\cdot, *) \Vert _{L^p(\Omega)}^2\right) \right](t, x) \\
      &\qquad\qquad \leq  \left[J_{0, p} \star \left(\bfone + \Vert u^{n-1}(\cdot, *)\Vert_{L^p(\Omega)}^2\right)\right](t, x).
\end{align*}
By \eqref{rd08_22e1} (with $n$ there replaced by $n-1$),
\begin{align*}
   &2 (\tilde C_p)^{\frac{2}{p}}\, \left[H^2 \star \Vert \sigma(\cdot, *, u^{n-1}(\cdot, *))\Vert_{L^p(\Omega)}^2\right](t, x)\\
    &\qquad \qquad\leq  \left[J_{0, p} \star \left(\bfone + 2\, I_0^2(t, x) + \cK_{n-1, p} \star (\bfone + 2\, I_0^2) \right)\right](t, x) \\
   &\qquad\qquad=  \left[(J_{0, p} + J_{0, p} \star \cK_{n-1,p}) \star \left(\bfone + 2\, I_0^2\right) \right](t, x) \\
   &\qquad\qquad=  \left[\cK_{n, p} \star \left(\bfone + 2\, I_0^2\right)  \right](t, x)
\end{align*}
by \eqref{rd08_12e1}. Therefore, by \eqref{rd05_02e2},  \eqref{rd08_22e1} holds for all $n \in \N^*$.

Next, we show the convergence of the Picard iteration scheme \eqref{rd05_01e1}. Fix $(t, x) \in \R_+^* \times \R$. If $\vert I_0(t, x) \vert = + \infty$, then we define $u(t, x) = u^n(t, x) = I_0(t, x) \in \{\pm \infty \}$ a.s., for all $n \in \N$.

 Suppose that that $\vert I_0(t, x) \vert < + \infty$. Recall that $u^0(t,x) = I_0(t,x)$, and for $n \in \N$, define
\beq\label{rd12_19_e7}
    F_n(t, x) = \Vert u^{n+1}(t, x) - u^{n}(t, x)\Vert_{L^p(\Omega)}^2.
\eeq
Use \eqref{rd08_22e10} to see that
\begin{align}\nonumber
    F_0(t, x) &= \Vert u^1(t, x) - u^0(t, x)\Vert_{L^p(\Omega)}^2 =  \Vert \cI_1(t, x))\Vert_{L^p(\Omega)}^2 \\
    &\leq 2\, \bar c^2 (\tilde C_p)^{\frac{2}{p}}\, [H^2 \star (\bfone +  I_0^2)](t, x) < \infty
\label{rd08_25e1}
\end{align}
by \eqref{rd08_21e2} in $({\bf H_{\Gamma, I}})$. In fact, \eqref{rd08_21e2} implies that 
$$\sup_{(t,x)\in K} F_0(t, x)<\infty$$ 
for any compact set
$K\subset \re_+^* \times \re$. 

Let $C$ be the Lipschitz constant in $({\bf H_{L,\infty}})$. Using \eqref{rd08_22e10} and the Lipschitz condition ($\bf H_{L,\infty}$), we obtain
\begin{align}
\label{fn}
  F_n(t, x) &\leq (\tilde C_p)^{\frac{2}{p}}\, \left[H^2 \star \Vert \sigma(\cdot, *, u^{n}(\cdot, *)) - \sigma(\cdot, *, u^{n-1}(\cdot, *))\Vert_{L^p(\Omega)}^2\right](t, x)\notag \\
   & \leq C^2\, (\tilde C_p)^{\frac{2}{p}}\, \left[H^2 \star \Vert u^{n} - u^{n-1} \Vert_{L^p(\Omega)}^2\right](t, x)\notag\\
   & = C^2\, (\tilde C_p)^{\frac{2}{p}}\, \left[H^2 \star F_{n-1}\right](t, x)\notag\\
   &= \left[\bar J \star F_{n-1}\right](t, x),
\end{align}
where $\bar J(t,x)= c\, J(t,x)$, with
\beqn
    c =  \frac{C^2(\tilde C_p)^{\frac{2}{p}}}{2\,\bar c^2}.
\eeqn
From \eqref{fn}, we see that the inequality \eqref{rd08_10e9} of Lemma \ref{rd08_10l2} holds with $f_n(t,x):=F_n(t,x)$, $J(t,x):=\bar J(t,x)$, $a\equiv 0$ and $z_0=0$. Therefore, by \eqref{rd12_19_e5},
\beqn
    F_n(t, x) \leq [\bar J^{\star n} \star F_0](t, x).
\eeqn
Consequenly,
\beqn
    \sum_{n=1}^\infty (F_n(t,x))^{1/2} \le \sum_{n=1}^\infty  \left([\bar J^{\star n} \star F_0](t,x)\right)^{1/2}.
 \eeqn
 Notice from \eqref{rd12_19_e7} and \eqref{rd08_25e1} that 
 \begin{align*}
     F_0 &= \Vert \cI^1(t, x) \Vert_{L^p(\Omega)}^2  
       \leq (\tilde C_p)^{2/p} [J \star (\bfone + 2\, I_0^2)](t, x) < \infty
 \end{align*}
by \eqref{rd08_22e6} and \eqref{rd08_21e2}. Therefore,
\beqn
   \sum_{n=1}^\infty (F_n(t, x))^\half \leq (\tilde C_p)^{1/p} \sum_{n=1}^\infty  \left(c^n [J^{\star (n+1)} \star (\bfone + 2\, I_0^2)](t, x)\right)^{\half}.
   \eeqn
In particular, for $(t, x)$ fixed, the sequence $(u_n(t, x))$ converges in $L^p(\Omega)$ and this convergence is in fact locally uniform, because of condition \eqref{rd12_19_e3}.

Denote by $u(t,x)$ the $L^p(\Omega)$-limit of the sequence $(u_n(t, x))$. We check that the process $(u(t,x),\, (t,x)\in\re_+^* \times \re)$ satisfies \eqref{rd01_22e1} and \eqref{rd08_21e4}. Indeed, for \eqref{rd01_22e1}, we use \eqref{rd08_22e2} to see that
\beq\label{rd01_31e2}
    \Vert u^{n}(t,x)\Vert_{L^2(\Omega)}^2 \leq I_0^2(t, x) + [\cK \star (\bfone + I_0^2)](t, x),
\eeq
and therefore, the same property holds for $u(t, x)$. In particular, \eqref{rd01_31e4} holds. 

It follows that
\begin{align*}
    &E\left[ \Vert \Gamma(t, x; \cdot, *)\, \sigma(\cdot, *, u(\cdot, *) \Vert_{L^2([0, t] \times \R)}^2 \right] \\
     &\qquad \leq 2\, \bar c^2\, \left[\Vert \Gamma(t, x; \cdot, *) \Vert_{L^2([0, t] \times \R)}^2 + \left\Vert \Gamma(t, x; \cdot, *) \Vert u(\cdot, *)\Vert_{L^2(\Omega)}) \right\Vert_{L^2([0, t] \times \R)}^2 \right]\\
     &\qquad  \leq 2\, \bar c^2\, \Vert \Gamma(t, x; \cdot, *) \Vert_{L^2([0, t] \times \R)}^2 + 2\, \bar c^2\, \left[H^2 \star\Vert u(\cdot, *)\Vert_{L^2(\Omega)}^2 \right](t, x)] .
\end{align*}
The first term is finite by ($\bf H_{\Gamma, I}$). As for the second term, we can bound it above (using \eqref{rd01_31e2} for $u(t, x)$) by
\beqn
    [J \star (I_0^2 + \cK \star (\bfone + I_0^2))](t, x) = [J \star I_0^2](t, x) + [J \star \cK \star (\bfone + I_0^2)](t, x).
\eeqn
Again, the first term is finite by \eqref{rd08_21e2}. By \eqref{rd08_12e2}, $J \star \cK \leq \cK$, therefore the second term is bounded above by
\begin{align}\nonumber
     [\cK \star (\bfone + I_0^2)](t, x) &= \left[ \left(\sum_{\ell =1}^\infty J^{\star \ell} \right) \star (\bfone + I_0^2) \right](t, x) \\
        &= \sum_{\ell =1}^\infty \, [J^{\star \ell} \star (\bfone + I_0^2) ](t, x) < \infty
\label{rd01_31e3}
\end{align}
by condition \eqref{rd12_19_e3} (we use here the following elementary fact: for $p \geq 1$ and for a sequence $(a_n)$ of nonnegative real numbers, if $\sum_{\ell = 1}^\infty (a_\ell)^{1/p} < \infty$, then $\sum_{\ell = 1}^\infty a_\ell < \infty$). This completes the proof of \eqref{rd01_22e1}.

Inequality \eqref{rd01_31e3} also shows that the right-hand side of \eqref{rd01_31e4} is finite when $\vert I_0(t, x) \vert < +\infty$.

We now prove \eqref{rd08_21e4}. Set
\begin{align*}
\cI^n(t,x) &= \int_0^t \int_{\re} \Gamma(t,x;s,y) \sigma(s,y,u^n(s,y))\, W(ds,dy),\\
{\cI}(t,x)&= \int_0^t \int_{\re} \Gamma(t,x;s,y) \sigma(s,y,u(s,y))\, W(ds,dy)
\end{align*}
(which is well-defined by \eqref{rd01_22e1}) and 
\beqn
G_n(t,x) =\Vert\cI^{n}(t,x) -\cI(t,x)\Vert^2_{L^2(\Omega)}.
\eeqn
Using the It\^o  isometry, we see that there is a finite constant $C$ such that
\begin{align} \nonumber
G_n(t, x)  &\leq C^2 \left[H^2 \star \Vert u^n- u\Vert^2_{L^2(\Omega)}\right](t,x) \\
   &= C^2 \int_0^t ds \int_\R dy\, H^2(t-s, x-y)\, \Vert u^n(t, x) - u(t, x) \Vert_{L^2(\Omega)}^2.
\label{verification}   
\end{align}
The integrand converges to $0$ as $n \to \infty$, because $u^n(t, x) \longrightarrow u(t, x)$ in $L^2(\Omega)$. In order to use the dominated convergence theorem, we observe using \eqref{rd01_31e2} that 
\beqn
  \Vert u^n(t, x) - u(t, x) \Vert_{L^2(\Omega)}^2 \leq 2 \left(I_0^2(t, x) + [\cK \star (\bfone + I_0^2)](t, x)\right).
\eeqn
Clearly, $[H^2 \star I_0^2](t, x) < \infty$ by \eqref{rd08_21e2}, and $[H^2 \star \cK \star (\bfone + I_0^2)](t, x) < \infty$ by the arguments that led to \eqref{rd01_31e3}.

From the dominated convergence theorem, we conclude that $\lim_{n\to\infty} G_n(t,x)=0$. In particular, $\cI^n(t, x) \longrightarrow \cI(t, x)$ in $L^2(\Omega)$. Since $u^{n+1}(t,x) = I_0(t,x) + \cI^n(t,x)$, it follows that $u(t,x) = I_0(t,x) + \cI(t,x)$, therefore \eqref{rd08_21e4} holds.

   We now prove \eqref{rd08_20e15}. Letting $n \to \infty$ in \eqref{rd08_22e1}, we obtain
\beqn
   \Vert u(t, x) \Vert_{L^p(\Omega)}^2 \leq 2 I_0^2(t, x) + [\cK_{p} \star (\bfone + 2\, I_0^2)](t, x).
\eeqn
The right-hand side is equal to
\beq\label{rd05_05e1}
   2 I_0^2(t, x) + \left[ \left(\sum_{\ell = 1}^\infty c^\ell\, J^{\star \ell} \right) \star (\bfone + 2\, I_0^2) \right](t, x),
\eeq
for some constant $c$ (that depends on $p$). 
This establishes \eqref{rd08_20e15}.

If $(t, x) \in \R_+^* \times \R$ is such that $\vert I_0(t, x) \vert < + \infty$, then the expression in \eqref{rd05_05e1} is finite by \eqref{rd12_19_e3}.

Suppose that $\vert I_0(t, x)\vert = \infty$ for some $(t, x) \in \R_+^* \times \R$. By \eqref{rd08_21e4}, $u(t, x) = I_0(t, x) + \cI(t, x)$ a.s., and by \eqref{rd01_22e1}, $\vert \cI(t, x) \vert < \infty$ a.s. Therefore, $\vert u(t, x) \vert = \infty$ a.s., and a.s., $u(t, x) = I_0(t, x)\ (= \pm\infty)$.

We finish the proof of the theorem by establishing uniqueness of the solution. Let $(u(t,x),\, (t,x)\in\re_+^* \times \re)$ and $(\bar u(t,x),\, (t,x)\in\re_+^* \times \re)$ be two random field solutions to \eqref{rd08_21e1} that satisfy \eqref{rd01_31e4}. 

Fix $(t,x)\in\re_+^* \times \re$ and suppose first that $\vert u(t, x) \vert  < \infty$ a.s. Then $\vert I_0(t, x) \vert < \infty$ by \eqref{rd08_21e4} and \eqref{rd01_22e1}, and therefore, $\vert \bar u(t, x) \vert < \infty$ a.s.

By the It\^o isometry and ($\bf H_{L,\infty}$), 
\begin{align*}
   &\Vert u(t,x)-\bar{u}(t,x)\Vert_{L^2(\Omega)}^2 \\
    & \qquad = \left\Vert \Gamma(t, x; \cdot, *)\, \Vert \sigma(\cdot, *, u(\cdot, *))- \sigma(\cdot, *, \bar u(\cdot, *)) \Vert_{L^2(\Omega)} \right\Vert_{L^2([0, t] \times \R)}^2  \\
&\qquad \leq 2 \bar c^2 \left[H^2 \star \Vert u- \bar u\Vert_{L^2(\Omega)}^2\right](t, x).
\end{align*}
Let $w = u- \bar{u}$. By Lemma \ref{rd08_10l2} (d) (with $a \equiv 0$ and $z_0 = 0$), we will be able to conclude that $\Vert w(t, x) \Vert_{L^2(\Omega)}^2 = 0$ if we show that
\beqn
   \lim_{n \to \infty} \left[J^{\star n} \star \Vert w \Vert_{L^2(\Omega)}^2 \right](t, x) = 0.
\eeqn
Clearly, 
\beqn
    \left[J^{\star n} \star \Vert w \Vert_{L^2(\Omega)}\right]^2(t, x) \leq 2 \left(\left[J^{\star n} \star \Vert u \Vert_{L^2(\Omega)}\right](t, x) + [J^{\star n} \star \Vert \bar u \Vert_{L^2(\Omega)}](t, x)\right).
\eeqn
By \eqref{rd01_31e4}, 
\begin{align*}
   [J^{\star n} \star \Vert u \Vert_{L^2(\Omega)}^2](t, x) &\leq [J^{\star n} \star I_0^2](t, x) + \sum_{\ell = 1}^\infty [J^{\star (n +\ell)} \star (\bfone + I_0^2)](t, x) \\
   &\leq \sum_{\ell = n}^\infty [J^{\star (\ell)} \star (\bfone + I_0^2)](t, x) \longrightarrow 0
\end{align*}
as $n \to \infty$ by \eqref{rd01_31e3}, and the same is true with $u$ replaced by $\bar u$. We conclude that $u(t, x) = \bar u(t, x)$ a.s.

Suppose that $P\{\vert u(t, x) \vert = + \infty \} > 0$. By \eqref{rd08_21e4} and \eqref{rd01_22e1}, $\vert I_0(t, x) \vert = + \infty$ and $u(t, x) = I_0(t, x)$ a.s. Therefore, $\vert \bar u(t, x) \vert = + \infty$ and $\bar u(t, x) = I_0(t, x) = u(t, x) $. This completes the proof of Theorem \ref{rd08_20t1}.
\qed
\bigskip

\noindent{\em Some sufficient conditions for Property \eqref{rd12_19_e3}}
\medskip

\begin{prop}\label{rd01_07p1}
Assume \eqref{rd08_22e6} and \eqref{rd08_21e2}. Let $J(t, x)$ be as defined in \eqref{defofj}.

(a) {\em A sufficient condition for \eqref{rd12_19_e3} that applies to the stochastic heat and wave equations.} Suppose that for all $ \ell \in \N^*$, there is a nonnegative nondecreasing function $ t \mapsto B_\ell(t)$ such that for all $(t, x) \in \R_+^* \times \R$, 
\beq\label{rd08_21e3}
    J^{\star \ell}(t, x) \leq B_\ell(t)\, J(t, x),
\eeq
and for all $c > 0$ and $p \geq 1$,
\beq\label{rd08_21e5}
   \sum_{\ell=1}^\infty \, \left[c^\ell \, B_\ell(t)\right]^{1/p} < \infty.
\eeq
Then \eqref{rd12_19_e3} holds.
   
(b) {\em A sufficient condition for \eqref{rd12_19_e3} that applies to the  stochastic heat and fractional heat equations.} Suppose that there are $\gamma: \R_+^* \to \R_+$ and $\tilde \Gamma: \R_+^* \times \R \to \R_+$ such that for all $(t, x) \in \R_+^* \times \R$, 
\beq\label{rd01_31e1}
    \int_0^t \gamma(s)\, ds < \infty
\eeq
and
\beq\label{rd01_03e1}
   J(t, x) \leq \gamma(t)\, \tilde\Gamma(t, x),
\eeq
where $\tilde \Gamma$ satisfies the semigroup property \eqref{semig-heat}, as well as $(\bf H_{\Gamma, I})$ (with the same $H$ as $\Gamma$), and for some $\ell_0 \in \N^*$ and all $\ell \geq \ell_0$, $\gamma^{\smallstar \ell}$ (the $\ell$-fold convolution in time of $\gamma$) is nondecreasing. Then \eqref{rd12_19_e3} holds.
\end{prop}

\begin{proof}
   (a) In this case, for $(t, x) \in [1/T, T] \times [-T, T]$,  by \eqref{rd08_21e3},
\begin{align*}
     [J^{\star \ell} \star (\bfone + I_0^2)](t, x) &= \int_{0}^{t} ds \int_{\R} dy \, (1 + I_0^2(t-s, x-y)) J^{\star \ell} (s, y) \\
        &\leq  \int_{0}^{t} ds \, B_\ell(s)\,  \int_\R dy\, (1 + I_0^2(t-s, x-y)) J(s, y) \\
        &\leq B_\ell(T)\, [J * (\bfone + I_0^2)](t, x),
\end{align*}
so for any $c > 0$, 
\begin{align*}
      &\sum_{\ell = 1}^\infty  \, \sup_{(t, x) \in [1/T, T] \times [-T, T]}\left(c^\ell\, [J^{\star \ell} \star (\bfone + I_0^2)](t, x)\right)^{\frac{1}{p}}  \\
       &\qquad \leq \left[ \sup_{(t, x) \in [1/T, T] \times [-T, T]} \left([J * (\bfone + I_0^2)](t, x)\right)^{\frac{1}{p}} \right] \sum_{\ell = 1}^\infty  \, \left( c^\ell\, B_\ell(T)  \right)^{\frac{1}{p}}  \\
       &\qquad < \infty,
\end{align*}
by \eqref{rd08_21e5}, \eqref{rd08_22e6} and \eqref{rd08_21e2}. This establishes \eqref{rd12_19_e3}.

(b)  We first show that for $\ell \in \N^*$,
\beq\label{rd08_24e3}
    J^{\star \ell}(t, x) \leq \gamma^{\smallstar \ell}(t)\, \tilde\Gamma(t, x)
\eeq
(the convolution in time $f \smallstar g$ is defined in Subsection \ref{rd04_16ss1}). Indeed, for $\ell=1$, \eqref{rd08_24e3}  is just the assumption \eqref{rd01_03e1} on $J$. Suppose by induction that \eqref{rd08_24e3} holds for $\ell-1$ ($\ell \geq 2$). We show that \eqref{rd08_24e3} holds for $\ell$ as follows:
\begin{align*}
    J^{\star \ell}(t,x) &= [J \star J^{\star (\ell-1)}](t, x) = \int_0^t ds \int_\R dy\, J(t-s, x-y)\, J^{\star (\ell-1)}(s, y) \\
     &\leq \int_0^t ds \, \gamma(t-s)\, \gamma^{\smallstar (\ell-1)}(s) \int_\R dy\, \tilde\Gamma(t-s, x-y)\, \tilde \Gamma(s, y) \\
    &= \int_0^t ds \, \gamma(t-s)\, \gamma^{\smallstar (\ell-1)}(s)\, \tilde\Gamma(t, x) \\
    &= [\gamma \smallstar \gamma^{\smallstar (\ell-1)}](t)\, \tilde\Gamma(t, x) =  \gamma^{\smallstar \ell}(t)\, \tilde\Gamma(t, x).
\end{align*}
Therefore, \eqref{rd08_24e3} holds for all $\ell \in \N^*$. 

Using the fact that $\gamma^{\smallstar \ell}$ is nondecreasing when $\ell \geq \ell_0$, we deduce that
\begin{align*}
      &\sum_{\ell = \ell_0}^\infty  \, \sup_{(t, x) \in [1/T, T] \times [-T, T]} \left(c^\ell\, [J^{\star \ell} \star (\bfone + I_0^2)](t, x)\right)^{\frac{1}{p}}  \\
     &\qquad\leq  \sum_{\ell = \ell_0}^\infty  \, \sup_{(t, x) \in [1/T, T] \times [-T, T]} \left(c^\ell\, \int_0^t ds\, \gamma^{\smallstar \ell}(s) \int_\R dy\,  \tilde\Gamma(s, y)  [\bfone + I_0^2 ](t-s, x-y)\right)^{\frac{1}{p}} \\
     &\qquad\leq \sum_{\ell = \ell_0}^\infty    \left(c^\ell\,  \gamma^{\smallstar \ell}(t) \sup_{(t, x) \in [1/T, T] \times [-T, T]}  [\tilde \Gamma \star  (\bfone + I_0^2)](t, x) \right)^{\frac{1}{p}}  \\
     &\qquad  = \left( \sup_{(t, x) \in [1/T, T] \times [-T, T]}  [\tilde \Gamma \star  (\bfone + I_0^2)]](t, x) \right)^{\frac{1}{p}}   \sum_{\ell = \ell_0}^\infty  \left( (c\, \gamma)^{\smallstar \ell} (t) \right)^{\frac{1}{p}}.   
\end{align*}
The first factor is finite by \eqref{rd08_22e6} and \eqref{rd08_21e2}. Notice that
\begin{align*}
    (c\, \gamma)^{\smallstar \ell} (t) & 
     = \int_0^t (c\, \gamma)^{\smallstar (\ell - \ell_0)} (t-s)\, (c\, \gamma)^{\smallstar \ell_0} (s)\, ds \\
    &\leq (c\, \gamma)^{\smallstar \ell_0} (t)  \int_0^t (c\, \gamma)^{\smallstar (\ell - \ell_0)} (s) \, ds.
\end{align*}
By \eqref{rd01_31e1}, we can apply Lemma \ref{A3-l00} (c) with $J$ there replaced by $c \gamma$, to conclude that
\beqn
    \sum_{\ell = \ell_0}^\infty  \left( (c\, \gamma)^{\smallstar \ell} (t) \right)^{\frac{1}{p}} \leq (c\, \gamma)^{\smallstar \ell_0} (T)  \sum_{\ell = \ell_0}^\infty  \left( \int_0^t (c\, \gamma)^{\smallstar (\ell - \ell_0)} (s) \, ds \right)^{\frac{1}{p}} < \infty. 
\eeqn
Therefore, \eqref{rd12_19_e3} holds.
\end{proof}

\begin{remark}\label{rd01_07r1}

Property \eqref{rd01_03e1} will in particular be satisfied whenever the fundamental solution takes the form $\Gamma(t-s, x-y)$, determines a convolution semigroup,
\beq
\label{rd08_26e4}
     \Gamma \geq 0 \quad\text{and}\quad \sup_{x \in \R}\, \Gamma(t, x) =: \gamma(t) < \infty, \text{ for all } t>0.
\eeq
Indeed, in this case, $\Gamma^2(t, x) \leq \gamma(t)\, \Gamma(t, x)$.

We will see in Subsection \ref{rd01_07ss1} that this situation occurs when $\Gamma$ is the fundamental solution of the fractional stochastic heat equation. \end{remark}
\medskip

\subsection{The stochastic heat, wave, and fractional heat equations}\label{rd01_07ss1}
\medskip

In this section, we give examples of SPDEs to which Theorem \ref{rd08_20t1} applies.
\smallskip

\noindent{\em The stochastic heat equation on $\re$}
\vskip 12pt

Let $\cL_{\nu} = \frac{\partial}{\partial t} - \frac{\nu}{2}\frac{\partial^2}{\partial x^2}$, where $\nu\in\re_+^*$. Notice that for $\nu=2$, $\cL_\nu$ is the heat operator considered in several sections of this book. The fundamental solution of $\cL_\nu=0$ is the function $\Gamma_\nu(t,x;s,y):= \Gamma_\nu(t-s,x-y)$,
where
\beq
\label{fsnu}
\Gamma_\nu(t,x) = \frac{1}{\sqrt{2\pi\nu t}}\exp\left(-\frac{x^2}{2\nu t}\right).
\eeq
For $\nu=2$, we will use the notation $\Gamma$ instead of $\Gamma_2$ (in agreement with formula \eqref{heatcauchy-1'}).

\begin{prop}
\label{rough-heat}
Consider the stochastic heat equation on $\re$, that is, \eqref{rd08_21e1} with $\cL=\cL_2$. Assume that the initial condition $u(0,\ast)$ belongs to the set
$\cM_H(\R)$ of signed Borel measures over $\R$ such that
\beq
\label{rd_08_25e1}
     \int_\R e^{-\alpha x^2} |\mu|(d x) < \infty, \quad \text{for all } \alpha >0,
\eeq
where $|\mu|$ denotes the total variation measure of $\mu$. Then Assumption $(\bf{H_{\Gamma,I}})$ and conditions \eqref{rd08_21e3}--\eqref{rd01_03e1} are satisfied.
\end{prop}
\begin{proof}
Let $\mu = u(0, *) \in \cM_H(\R)$, and let $H(t,x):=\Gamma(t,x)$, $(t,x)\in\re_+\times \re$. It is then clear that conditions (i)--(iii) of $(\bf{H_{\Gamma,I}})$ are satisfied. Hence, it remains to prove \eqref{rd08_21e2}.

Recall from the Jordan decomposition\index{Jordan decomposition}\index{decomposition!Jordan} that $\mu=\mu_+-\mu_-$, where $\mu_\pm$ are two nonnegative Borel measures with disjoint support, and $|\mu|:= \mu_+ + \mu_-$.
Observe also that condition \eqref{rd_08_25e1} is equivalent to 
\beqn
    [\Gamma_\nu(t, *) * |\mu|]( x)  < \infty, \qquad\text{for all } (t,x)\in\re_+^*\times \re \quad\text{and } \nu\in\re_+^*.
 \eeqn
 This implies that 
 \beq\label{rd05_06e5}
    I_0(t, x) := [\Gamma(t, *) * \mu](x) 
\eeq
is well-defined and
\beq
     \vert I_0(t, x) \vert \leq \tilde I_0(t, x) := \int_{\re} \Gamma(t,x-y)\, \vert \mu \vert(dy) < \infty,\quad (t,x)\in\re_+^*\times \re.
\label{rd01_29e1}
\eeq

Notice that 
\beq
\label{form-square}
   \tilde I_0^2(t,x) = \int_{\re} \vert \mu\vert (dz_1) \int_\R \vert \mu\vert (dz_2) \, \Gamma(t,x-z_1) \Gamma(t,x-z_2)
\eeq
 and 
 \beqn
 H^2(t,x) = h(t)\, \Gamma_1(t,x), \qquad\text{where } h(t) := \frac{1}{2\sqrt{2\pi t}}.
 \eeqn
Therefore, 
\begin{align*}
[H^2\star I_0^2](t,x) &\leq [H^2\star \tilde I_0^2](t,x) \\
&=\int_0^t ds\, h(t-s)\int_{\re} dy\, \Gamma_1(t-s,x-y)\\
&\qquad\times\int_{\re} \vert\mu\vert(dz_1)\int_\re \vert\mu\vert(dz_2)\,\Gamma(s,y-z_1) \Gamma(s,y-z_2).
\end{align*}
Apply \eqref{rough-heat-l1-1} with $\nu=2$, $s=t$, $x:= y - z_1$, $y:= y-z_2$, to deduce that   
\beqn
\Gamma(s,y-z_1) \Gamma(s,y-z_2)=\Gamma\left(\frac{s}{2}, \frac{(y-z_1)+(y-z_2)}{2}\right) \Gamma(2s,z_2-z_1).
\eeqn
Thus,
\begin{align*}
[H^2\star \tilde I_0^2](t,x)&= \int_0^t ds\, h(t-s)\int_{\re} \vert\mu\vert(dz_1) \int_\re \vert\mu\vert(dz_2)\,\Gamma(2s,z_2-z_1)\\
&\qquad\times \int_{\re} dy\, \Gamma(s/2, y-(z_1+z_2)/2)\Gamma_1(t-s,x-y).
\end{align*}
Let $\bar z= (z_1+z_2)/2$. Because $\Gamma\left(\frac{s}{2}, * \right) = \Gamma_1(s,*)$, the semigroup property of $\Gamma$ (see \eqref{semig-heat}) yields
\beqn
\int_{\re} dy\, \Gamma(s/2, y- \bar z)\Gamma_1(t-s,x-y)= \Gamma_1(t,x-\bar z).
\eeqn
Therefore, 
\begin{align*}
[H^2\star \tilde I_0^2](t,x) &=  \int_0^t ds\, h(t-s)\\
&\qquad\times\int_{\re} \vert\mu\vert(dz_1) \int_{\re} \vert\mu\vert(dz_2)\,\Gamma_1(t,x-\bar z)\, \Gamma(2s,z_2-z_1) \\
&= \int_0^t ds\, h(t-s)\\
&\qquad\times\int_{\re} \vert\mu\vert(dz_1) \int_{\re} \vert\mu\vert(dz_2)\,\Gamma_1(t,x-\bar z) \, \Gamma_1(4s,z_2-z_1).
\end{align*}
Using the inequality \eqref{rough-heat-l1-20} with $a=2$ there, we have
\begin{align*}
[H^2\star \tilde I_0^2](t,x) &\le \int_0^t ds\, \frac{4(t\vee s)}{\sqrt{ts}}\, h(t-s)\\
   &\qquad\times\int_{\re} \vert\mu\vert(dz_1) \int_{\re} \vert\mu\vert(dz_2)\,\Gamma_1(4(t \vee s), x-z_1)\Gamma_1(4(t \vee s), x-z_2)\\
   &= 4\sqrt{t} \int_0^tds\, \frac{h(t-s)}{\sqrt{s}} \\
   &\qquad\times\int_{\re} \vert\mu\vert(dz_1) \int_{\re} \vert\mu\vert(dz_2)\,\Gamma_1(4t, x-z_1)\Gamma_1(4t, x-z_2).
\end{align*}
Since $\Gamma_1(4t,z) = \Gamma(2t,z)$, we obtain
\beqn
[H^2\star I_0^2](t,x) \le 2 \sqrt{t}\, \tilde I_0^2(2t,x)\int_0^tds\, \frac{h(t-s)}{\sqrt{s}}.
\eeqn 
Up to a multiplicative positive constant, 
$$\int_0^tds\, \frac{h(t-s)}{\sqrt{s}} = \int_0^t (t-s)^{-\half}\, s^{-\half}\, ds = t \, B_E\left(\half,\half\right),$$ 
where $B_E(\half,\half)$ is the Euler Beta function (see \eqref{rd08_12e9}). In particular, \eqref{rd08_21e2} is satisfied. This establishes $(\bf{H_{\Gamma,I}})$.

We now establish \eqref{rd08_21e3}--\eqref{rd01_03e1}.  As in \eqref{defofj}, we set
\beq\label{rd01_27e1}
J(t,x)=2\bar c^2\,  \Gamma^2(t,x) = 2\bar c^{2}\frac{1}{2\sqrt{2\pi t}}\, \Gamma_1(t,x).
\eeq
Apply \eqref{rd08_20e12} 
 with $J^h$ there equal to $J$ to see that
\beq\label{rd05_03e1}
J^{\star \ell}(t,x) = \sqrt \pi \left(\sqrt 2\bar c \right)^{2(\ell-1)} \frac{t^{(\ell-1)/2}}{\Gamma_E(\ell/2)}\, J(t,x),
\eeq
where $\Gamma_E$ is the Euler Gamma function. Hence,  \eqref{rd08_21e3} holds with
\beqn
B_\ell(t) = \sqrt \pi \left(\sqrt 2\bar c \right)^{2(\ell-1)} \frac{t^{(\ell-1)/2}}{\Gamma_E(\ell /2)}
\eeqn
and \eqref{rd08_21e5} is satisfied. Further, let 
\beqn
   \gamma(t) = 2\bar c^{2}\frac{1}{2\sqrt{2\pi t}} \qquad\text{and}\qquad \tilde \Gamma(t, x) = \Gamma_1(t, x).
\eeqn
Equality \eqref{rd01_27e1} becomes 
\beqn
   J(t, x) = \gamma(t) \tilde \Gamma(t, x).
\eeqn
Then \eqref{rd01_31e1} and \eqref{rd01_03e1} hold, since we deduce from \eqref{rd01_27e1} and \eqref{rd05_03e1} that $\gamma^{\smallstar \ell} (t) = B_\ell(t) \, \bar c^2 (2 \pi t)^{-1/2}$ and this function is nondecreasing for $\ell \geq 2$. 
\end{proof}
\medskip

\noindent{\em The stochastic wave equation on $\re$}
\medskip

Consider the SPDE \eqref{rd08_21e1} with $\cL = \frac{\partial^2}{\partial t^2} - \frac{\partial^2}{\partial x^2}$ (the one-dimensional wave operator on $\re$) with given deterministic initial conditions
\beqn
u(0,\ast) = u_0, \quad \frac{\partial u}{\partial t}(0,\ast) = v_0,
\eeqn
satisfying certain conditions (see Proposition \ref{rough-wave} below). We recall that the fundamental solution associated to $\cL$ is
$\Gamma(t,x;s,y):=\Gamma(t-s,x-y)$, with
\beq
\label{ex-wave-rough}
\Gamma(t,x) = \frac{1}{2} 1_{\re_+}(t)\, 1_{[-t,t]}(x),
\eeq
(see \eqref{wfs}).

\begin{prop}
\label{rough-wave}
If $(u_0, v_0) \in L_{\text{loc}}^2(\R) \times \cM(\R)$, where $\cM(\R)$ denotes the set of locally finite signed Borel measures on $\R$, 
then Assumption $(\bf{H_{\Gamma,I}})$ and  conditions \eqref{rd08_21e3} and \eqref{rd08_21e5} are satisfied.
\end{prop}

\begin{proof}
Take $H(t,x) = \Gamma(t,x)$ (defined in \eqref{ex-wave-rough}).
We check that condition $(\bf{H_{\Gamma,I}})$ is satisfied.  For this, we notice from \eqref{p93.1} that the solution to the associated homogeneous SPDE ($\sigma\equiv 0$) with initial conditions $(u_0, v_0)$ is 
\beq\label{rd01_29e2}
   I_0(t, x) = \half [u_0(x + t) + u_0(x - t)] + [v_0 * \Gamma(t, *)](t, x),
\eeq
and that under the hypotheses on $(u_0, v_0)$, this expression is well-defined. 
For any $(t,x)\in\re_+^* \times \re$,
\begin{align*}
\left( \Gamma^2\star I_0^2\right)(t,x) &= \frac{1}{4} \int_{\re_+} ds \int_{\re} dy\ 1_{[-(t-s),\, t-s]}(x-y) I_0^2(s,y) \\
   & \le \frac{1}{2} \left[T_1(t,x)+T_2(t,x)\right],
\end{align*}
where
\begin{align*}
T_1(t,x) &= \frac{1}{4}  \int_0^t ds \int_{x-(t-s)}^{x+(t-s)} dy\ [u_0(y+s) + u_0(y-s)]^2,\\
T_2(t,x)&= \int_0^t  ds \int_{x-(t-s)}^{x+(t-s)} dy\, \left(\int_{\re} 1_{[y-s,y+s]}(z)\, |v_0|(dz)\right)^2 . 
\end{align*}
By the change of variables $z=y+s$ and $\bar z = y-s$, we easily see that 
\begin{align}\nonumber
T_1(t,x) &\leq \frac{1}{8}  \int_{x-t}^{x+t} dz \int_{x-t}^{x+t}  d\bar z\, (u_0(z) + u_0(\bar z))^2 \\
    & \le  t \int_{x-t}^{x+t} dz\, u_0^2(z).
    \label{bound-t1}
\end{align}
Since $|v_0|$ is a nonnegative measure, we have 
$$
   \int_{\re} 1_{[y-s,y+s]}(z)\, |v_0|(dz)\le |v_0|([x-t,x+t]), \qquad\text{if } \vert x - y\vert \leq t - s.  
$$
Therefore
\begin{align}\nonumber
T_2(x,y) &\le  |v_0|^2 ([x-t,x+t]) \int_0^t  ds \int_{x-(t-s)}^{x+(t-s)} dy \\
   &= \frac{t^2}{2}\, | v_0|^2 ([x-t,x+t]).
\label{bound-t2}
\end{align}

From \eqref{bound-t1} and \eqref{bound-t2}, we deduce that for some $c < \infty$,
\begin{align*}
\left[H^2\star  I_0^2\right](t,x) &= \left[ \Gamma^2\star I_0^2\right](t,x) \\
    & \le c \left( t \int_{x-t}^{x+t} dz\, u_0^2(z) +  t^2\, |v_0|^2 ([x-t,x+t])\right).
\end{align*}
The right-hand side of this inequality is uniformly bounded over any compact set $K\subset \re_+^*\times \re$. Hence,
 $(\bf{H_{\Gamma,I}})$ holds.
 
We now check \eqref{rd08_21e3} and \eqref{rd08_21e5}.  Let $\rho:= \sqrt{2} \, \bar c$. Then 
\beqn
J(t,x) = \frac{\rho^2}{4}1_{\re_+}(t)\, 1_{[-t,t]}(x).
\eeqn
The equality \eqref{rd08_24e1} with $\nu=1$ (and $J^{\mathrm w} =J$) tells us that
\beqn
J^{\star \ell}(t,x) = \frac{\rho^{2(\ell-1)}(t^2-x^2)^{\ell-1}}{2^{3(\ell-1)}((\ell-1)!)^2}\, J(t,x),
\eeqn
 and therefore $J^{\star \ell}(t,x) \le B_\ell(t) J(t,x)$, with 
 \beqn
 B_\ell(t) =  \frac{\rho^{2(\ell-1)}(t^2-x^2)^{\ell-1}}{2^{3(\ell-1)}((\ell-1)!)^2}.
 \eeqn
From this, \eqref{rd08_21e3} and \eqref{rd08_21e5} follow. 

\end{proof}
\smallskip

\noindent{\em The fractional stochastic heat equation}
\medskip

Consider the SPDE \eqref{rd08_21e1} with the partial differential operator
\beqn
\cL = \frac{\partial}{\partial t} - \null_xD_\delta^a,
\eeqn
where $\null_xD_\delta^a$, $a\in\,]1,2[$, $|\delta| < 2-a$, is the Riesz-Feller fractional derivative defined in \eqref{ch1'-ss7.3.1}. 
We assume that the initial value of the SPDE is  a measure $\mu$ satisfying certain conditions (see Proposition \ref{rough-fractional-heat} below).

From Section \ref{ch1'-ss7.3}, we recall that the fundamental solution associated with $\cL$ is the function $\Gamma(t,x;s,y)= \null_\delta G_a(t-s,x-y)$, where $\null_\delta G_a(t,x)$ is defined in \eqref{fs-fractional}, and that the solution to the homogeneous PDE  with $\mu$ as initial condition is, formally,
 \beq
 \label{initial-frac-heat}
    I_0(t, x) = [\null_\delta G_a(t,*) * \mu](t, x) = \int_{\re} \mu(dy)\, \null_\delta G_a(t,x-y) .
\eeq
 

\begin{prop}
\label{rough-fractional-heat}
Let  $\cM_a(\R)$ denote the set of signed Borel measures $\mu$ such that
\beq
\label{rd08_26e1}
   A_a :=  \sup_{y\in\R} \int_\R\, \frac{|\mu|(d x)}{1+|x-y|^{1+a}} < \infty.
\eeq
If $a\in\,]1,2[$ and $|\delta| < 2-a$, and if the initial condition $\mu$ belongs to  $\cM_a(\R)$, then Assumption $(\bf{H_{\Gamma,I}})$ as well as conditions \eqref{rd01_31e1} and \eqref{rd01_03e1} are satisfied. 
\end{prop}

\begin{proof}
Let $H(t,x)  =  \laplacef(t, x)$.  
Using \eqref{ch1'-ss7.3.6}, let
\beq
\label{ineq-gamma}
    \gamma(t) = \sup_{x \in \R} \Gamma(t, x) = C_0\, t^{-\frac{1}{a}} , \quad \text{where } C_0=\sup_{x\in\re}\null_\delta G_a(1,x).
\eeq
As in Remark \ref{rd01_07r1}, we have 
\beqn
    H^2(t,x) \leq  \gamma(t)\, \null_\delta G_a(t,x).
\eeqn

Assumptions $(\bf{H_{\Gamma,I}})$ (i)--(iii) are clearly satisfied. For the proof of $(\bf{H_{\Gamma,I}})$ (iv), we first notice that according to Lemma \ref{rd01_21l1}, 
for any $t>0$, 
\beq
\label{bound-fs}
\null_\delta G_a(t,x)\le K_{a}\,  t^{-\frac{1}{a}} \, \frac{(t\vee 1)^{1+\frac{1}{a}}}{1+|x|^{1+a}},
\eeq
where $K_{a} \in\R_+^*$. Therefore, for $\mu\in \cM_a(\R)$,
\beq
\label{bound-i0}
\int_{\re} |\mu|(dy)\, \null_\delta G_a(t,x-y) \le K_{a}\, t^{-\frac{1}{a}}(t\vee 1)^{1+\frac{1}{a}} \int_\R \frac{|\mu|(d y)}{1+|x-y|^{1+a}} < \infty
\eeq
by \eqref{rd08_26e1}, which implies that
\beq
\label{bound-i1}
\int_{\re} \mu_\pm(dy)\, \null_\delta G_a(t,x-y) \le K_{a}\,  A_a\, t^{-\frac{1}{a}}(t\vee 1)^{1+\frac{1}{a}} < \infty, 
\eeq
where $\mu_+$ and $\mu_-$ are the nonnegative measures in the Jordan decomposition $\mu = \mu_+-\mu_-$. This implies that \eqref{initial-frac-heat} is well-defined. 

In the calculation of $\left[H^2 \star I_0^2\right](t,x)$, we will use \eqref{bound-i1} and the inequality $(a-b)^2 \le a^2+b^2$, valid for any $a,b\in\re_+$.  Notice that
\begin{align*}
\left[H^2 \star I_0^2\right](t,x)&\leq \int_0^t ds\, \gamma(t-s)  \int_{\re}dy\, \null_\delta G_a(t-s,x-y) \, I^2_0(s,y)\\
&\le \int_0^t ds\, \gamma(t-s)  \int_{\re}dy\, \null_\delta G_a(t-s,x-y) \\
&\qquad \times \left[\left(\int_{\re} \mu_+(dz)\, \null_\delta G_a(s,y-z)\right)^2 + \left(\int_{\re} \mu_-(dz)\, \null_\delta G_a(s,y-z)\right)^2\right].
\end{align*}
Next, we will find a suitable upper bound for 
\beqn 
T^\pm(t,x) := \int_0^t ds\, \gamma(t-s) \int_{\re}dy\, \null_\delta G_a(t-s,x-y) 
\left(\int_{\re} \mu_\pm(dz)\, \null_\delta G_a(s,y-z)\right)^2.
\eeqn 
For this, we use 
\eqref{bound-i1} 
in one factor of the square in the definition of $T^\pm(t,x)$ to see that,
\begin{align*}
T^\pm(t,x) &\le  K_{a} \, A_a \, (t\vee 1)^{1+\frac{1}{a}} \int_0^t ds\, \gamma(t-s)\, s^{-\frac{1}{a}} \int_{\re} dy\, \null_\delta G_a(t-s,x-y) \\
&\qquad\qquad\times\left(\int_{\re} \vert \mu \vert(dz)\,\null_\delta G_a(s,y-z)\right).
\end{align*}
Applying Fubini's theorem and the semigroup property \eqref{ch1'-ss7.3.5}, we have 
\begin{align*}
T^\pm(t,x) &\le  K_{a}\,  A_a\, (t\vee 1)^{1+\frac{1}{a}}\, \int_0^t ds\, \gamma(t-s)\, s^{-\frac{1}{a}} \\
&\qquad\qquad \times\int_{\re} \vert \mu \vert (dz)\,\null_\delta G_a(t,x-z). 
\end{align*}


Summarising, we have proved that
\begin{align*}
\left[H^2 \star I_0^2\right](t,x) &\le 2 K_{a}\, A_a\, (t\vee 1)^{1+\frac{1}{a}}\, \int_0^t ds\, \gamma(t-s)\, s^{-\frac{1}{a}} \\
&\qquad\qquad \times  \int_{\re} |\mu|(dz)\,\null_\delta G_a(t,x-z)  \\
&= 2 K_{a}\, A_a\, C_0\, (t\vee 1)^{2+\frac{1}{a}}\, B_E\left(\frac{1}{a^*},\frac{1}{a^*}\right)\int_{\re} |\mu|(dz)\,\null_\delta G_a(t,x-z),
\end{align*}
by \eqref{ineq-gamma}, where $\frac{1}{a^*} + \frac{1}{a}  = 1$ and $B_E$ denotes the Euler Beta function (see \eqref{rd08_12e9}). This implies $({\bf H_{\Gamma,I}})$ (iv).
  

Observe that
\beqn
    J(t,x) = 2\bar c^2 H^2(t,x) \leq  2\bar c^2 \gamma(t)\, \null_\delta G_a(t,x),
\eeqn
therefore, \eqref{rd01_31e1} and \eqref{rd01_03e1} hold. Using \eqref{rd08_26L1}, the Laplace transform of $\gamma(t)$ is a multiple of $z^{-1 + 1/a}$. Therefore, for $\ell \geq 1$, $\gamma^{\smallstar \ell}(t)$ is a multiple of the inverse Laplace transform of $z^{\ell(-1 + 1/a)}$, which is a multiple of $t^{\ell(1 - 1/a) -1}$, and this function is nondecreasing for  $\ell \geq 2$.
\end{proof}

\medskip

\begin{remark}
\label{rd08_26r1}
Because of standard inequalities on the fundamental solutions of parabolic operators (such as \cite[Theorem VI.1, p.~183, Theorem VI.6, p.~193]{ez}), Theorem \ref{rd08_20t1} also holds for more general parabolic SPDEs than the stochastic heat equation considered in Proposition \ref{rough-heat}, for instance when $\cL$ is the partial differential operator
\beqn
\cL = \frac{\partial}{\partial t} - a(t,x)\, \frac{\partial^2}{\partial x^2} +  b(t,x)\, \frac{\partial}{\partial x} + c(t,x)
 \eeqn
on $\re_+ \times \re$, where $a$, $b$ and $c$ are deterministic functions satisfying suitable conditions.
\end{remark}

\subsection{Examples and counterexamples}\label{rd06_23ss1}

In the three examples considered in Propositions \ref{rough-heat}--\ref{rough-fractional-heat}, the Dirac delta function  $\delta_0$ is an eligible rough initial condition. We write explicitly this claim as follows.
\begin{examp}
\label{rd08_21ex1}
\begin{description}
\item{(a)} {\em The stochastic heat equation on $\R$.}
Let $u_0(x) = \delta_0(x)$, $x\in \R$, so that $I_0(t, x) = \Gamma_\nu(t, x)$ $($defined in \eqref{fsnu}$)$. Then $\mu:=\delta_0$ belongs to $\cM_H(\re)$, and the conclusions of Theorem \ref{rd08_20t1} hold.
\item{(b)} {\em The stochastic wave equation on $\R$.} Consider the intial conditions $u_0 \equiv 0$ and $v_0(x) = \delta_0(x)$, $x \in \R$, so that $I_0(t, x) = \Gamma_\nu(t, x)$ $($defined in \eqref{ex-wave-rough}$)$. Then $(u_0, v_0) \in L_{\text{loc}}^2(\R) \times \cM(\R)$ and the conclusions of Theorem \ref{rd08_20t1} hold.
\item{(c)} {\em The fractional stochastic heat equation on $\R$.} 
Consider the initial condition $u_0(x) = \delta_0(x)$, $x\in \R$, so that $I_0(t, x) = \null_\delta G_a(t, x)$ $($defined in \eqref{fs-fractional}$)$, with $a \in \, ]1, 2[$ and $\vert \delta \vert < 2 - a$. Then \eqref{rd08_26e1} is satisfied and the conclusions of Theorem \ref{rd08_20t1} hold.
\end{description}
\end{examp}

\medskip

We end this section by giving examples in which there is no random field solution of \eqref{rd08_21e1}, or in which the solution takes the value $\infty$ at certain space-time points, or in which a solution exists only in certain regions of space-time.

\begin{prop}
\label{rd08_25p1}
Consider the heat operator $\cL = \frac{\partial}{\partial t} - \frac{\partial^2}{\partial x^2}$
and the associated stochastic heat equation \eqref{rd08_21e1} with coefficient $\sigma(t,x,u(t,x)) = \rho\, u(t,x)$, $\rho\in\re\setminus\{0\}$.  Suppose that the initial condition is $\mu = \delta_0'$ (the derivative of the Dirac delta function at zero). Then \eqref{rd08_21e1} does not have a random field solution (in the sense of Definition \ref{rd08_21d1}).
\end{prop}

\begin{proof}
Suppose that $(u(t, x))$ is a solution of \eqref{rd08_21e1} in the sense of Definition \ref{rd08_21d1}. We follow the proof of \cite[Proposition 2.11]{chendalang2015-2} which consists in checking that for $(u(t,x))$ defined as in \eqref{rd08_21e4}, condition \eqref{rd01_22e1} is not satisfied, which is a contradiction. 

Indeed, by definition, for $(t, x) \in \R_+^* \times \R$,
\beq
\label{initial-c}
I_0(t,x) = \int_{\re} \Gamma(t,x-y)\, \delta_0'(dy) = \frac{\partial}{\partial x} \Gamma(t,x) = -\frac{x}{2 t} \Gamma(t,x).
\eeq
From \eqref{rd08_21e4}, 
we deduce that  $\Vert u(t,x)\Vert_{L^2(\Omega)}\ge\vert I_0(t,x)\vert$ and that
\begin{equation}
\label{main-eq}
\Vert u(t,x)\Vert^2_{L^2(\Omega)}\ge \rho^2 \left[ \Gamma^2 \star \Vert u \Vert^2_{L^2(\Omega)}\right](t,x) 
\ge \rho^2 \left[ \Gamma^2 \star  I_0^2\right](t,x).
\end{equation}

Recall that $\Gamma^2(t,x) = \frac{1}{2\sqrt{2\pi t}}\Gamma_1(t,x)$, with $\Gamma_1(t,x)$ defined by \eqref{fsnu} (with $\nu=1$ there). Therefore, by \eqref{initial-c}, 
\begin{align*}
\left[ \Gamma^2 \star I_0^2\right](t,x)& = \int_0^t ds \int_{\re} dy\, \Gamma^2(t-s,x-y) I_0^2 (s,y)\\
&= \frac{1}{32 \pi} \int_0^t \frac{ds}{s^2\sqrt{s(t-s)}} \int_{\re} dy\, y^2 \, \Gamma_1(t-s,x-y) \Gamma_1(s,y).
\end{align*}
Apply \eqref{rough-heat-l1-1} to see that
\beqn
\Gamma_1(t-s,x-y) \Gamma_1(s,y) = \Gamma_1\left(\frac{(t-s)s}{t}, y-\frac{sx}{t}\right) \Gamma_1(t,x).
\eeqn
Hence,
\beqn
\left[ \Gamma^2 \star I_0^2\right](t,x)= \frac{1}{32 \pi} \Gamma_1(t,x) \int_0^t \frac{ds}{s^2\sqrt{s(t-s)}}\, T(x,s,y),
\eeqn
where
\beqn
T(x,s,y) = \int_{\re} dy\, y^2 \Gamma_1\left(\frac{(t-s)s}{t},\, y-\frac{sx}{t}\right) = E\left[\left(Z+\frac{sx}{t}\right)^2\right] = s-\frac{s^2}{t}+\frac{s^2x^2}{t^2},
\eeqn
and $Z$ is a  ${\rm N}(0, (t-s)s/t)$ random variable. 
This implies that
\beqn
\int_0^t \frac{ds}{s^2\sqrt{s(t-s)}}\, T(x,s,y) = \int_0^t \frac{ds}{s\,\sqrt{s(t-s)}} -\left(\frac{1}{t} - \frac{x^2}{t^2}\right) B_E\left(\frac{1}{2}, \frac{1}{2}\right).
\eeqn
Since 
$$
   \int_0^t \frac{ds}{s\sqrt{s(t-s)}} = +\infty,
$$ 
\eqref{main-eq} implies that for all $(t, x) \in \R_+ \times \R$,
\beqn
    \Vert u(t, x) \Vert_{L^2(\Omega)} \geq \rho^2 \left[ \Gamma^2 \star I_0^2\right](t,x)= +\infty, 
\eeqn
which contradicts \eqref{rd01_22e1}.
\end{proof}

\begin{prop}
\label{rd08_25p2}
Consider the wave operator $\cL = \frac{\partial^2}{\partial t^2}- \frac{\partial^2}{\partial x^2}$, as in Proposition \ref{rough-wave}, and the associated stochastic wave equation \eqref{rd08_21e1}. 

\noindent{(a)}\ Suppose that the initial conditions are $u_0(x) = \vert x \vert^{-1/4}$ and $v_0 \equiv 0$. Then $u_0 \in L_{\text{loc}}^2(\R)$ and
\beqn
   I_0(t, x) = \half \left(\frac{1}{\vert x + t \vert^{1/4}} + \frac{1}{\vert x - t \vert^{1/4}} \right).
\eeqn
The function $I_0(t, x)$ equals $+ \infty$ on the characteristic lines $x=\pm t$ through the origin $(0, 0)$. Nevertheless, 
in agreement with Proposition \ref{rough-wave} and Theorem \ref{rd08_20t1}, a random field solution $(u(t,x),\, (t,x)\in \R_+^*\times\R)$ to \eqref{rd08_21e1} in the sense of Definition \ref{rd08_21d1} does exist. This solution is such that for all $t > 0$, $u(t, -t) = u(t, t) = + \infty$ a.s., and for all $(t,x)\in \R_+^*\times\R$, the stochastic integral
\beqn
   \int_0^t \int_\R 1_{[0,\, t-s]}(\vert x-y\vert)\, \sigma(s,y,u(s,y)) \, W(ds, dy)
\eeqn
is well-defined. 
\smallskip

\noindent{(b)}\ Let $u_0(x)= |x|^{-\half}$ and $\mu \equiv 0$. Then $u_0 \not\in
L_{\text{loc}}^2(\R)$ and therefore, neither Proposition \ref{rough-wave} nor Theorem \ref{rd08_20t1} apply.
The function
\beqn
   I_0^2(t,x) = \frac{1}{4} \left(\frac{1}{|x+ t|^{1/2}}+\frac{1}{|x- t|^{1/2}}\right)^2,
\eeqn
is not locally integrable over domains that intersect the
characteristic lines $x=\pm t$, so there is no random field solution $(u(t, x))$ to \eqref{rd08_21e1} on $\R_+^* \times \R$. However, in agreement with Theorem \ref{rd08_20t1}, a random field solution to \eqref{rd08_21e1} does exist in the domain $\{(t, x):  t < |x|\}$, which is the union of two ``infinite triangles.''
\end{prop}

\begin{proof}
The proof is straightforward and is left to the reader.
\end{proof}

\section{The $(1 + 1)$-Anderson model}
\label{rd1+1anderson}

In this section, we study the particular case of the SPDEs considered in Section \ref{rdrough} in which the function $\sigma(t,x,z)$ is a fixed linear function of $z$. More specifically, we provide an introduction to the $(1 + 1)$-Anderson model
\beq
\label{rd08_14e1}
    \cL u(t, x) = \rho\, u(t, x) \, \dot W(t, x), \qquad  (t, x) \in \, ]0, \infty[ \times \R,
\eeq
where $\rho \in \R \setminus \{ 0 \}$, with given initial conditions.
The qualification $(1 + 1)$ refers to the fact that the couple $(t, x)$ consists of two one-dimensional real variables. 
We will assume that  the hypotheses of Theorem \ref{rd08_20t1} are satisfied and, moreover, that the fundamental solution associated with $\cL$ is $\Gamma(t,x;s,y):= \Gamma(t-s,x-y)$ with $\Gamma \geq 0$. Hence, we take $H(t-s,x-y) = \Gamma(t-s,x-y)$ in ($\bf H_{\Gamma, I}$) and we set 
\beq\label{rd01_27e2}
J(t,x) = \rho^2\, \Gamma^2(t,x)
\eeq
(comparing with \eqref{defofj}, there is no factor $2$ here. Indeed, since $\sigma(t, x, z) = \rho\, z$, no factor $2$ would appear in \eqref{rd01_07e1}).

\subsection{An expression for the second moment} 

Because of the particular form of the right-hand side of the SPDE \eqref{rd08_14e1}, it is possible to give exact formulas for the second moments of its solution, as shown in the next theorem.

\begin{thm}
\label{rd08_14p1}
The hypotheses are as in Theorem \ref{rd08_20t1}.
Let $(u(t, x),\, (t, x) \in \R_+^* \times \R)$ be the solution of the SPDE \eqref{rd08_14e1}, so that 
for all $(t, x) \in \R_+^* \times \R$,
\beq
\label{rd08_14e9}
    u(t, x) = I_0(t, x) +  \rho \int_0^t \int_\R \Gamma(t-s, x - y)\, u(s, y)\, W(ds, dy) ,
\eeq
where $(I_0(t, x))$ is the solution to the homogeneous PDE $\cL u = 0$ with the same initial conditions as for \eqref{rd08_14e1}.
Then
\beq
\label{rd08_20e4}
     E[u^2(t, x)] = I_0^2(t, x) + [\cK \star I_0^2](t, x),
\eeq
where, as in \eqref{def-kas},
\beq\label{rd02_04e1}
   \cK(t, x) = \sum_{\ell = 1}^\infty J^{\star \ell} (t, x)
\eeq
 and $J$ is from \eqref{rd01_27e2}.

In particular, if $I_0(t, x) \equiv c$ for some $c \in \R$, then
\beq
\label{rd08_20e2}
   E[u^2(t, x)] = c^2 \, [\bfone + \cK \star \bfone](t, x).
\eeq
 If $I_0(t, x) = c\, \Gamma(t, x)$, for all $(t, x) \in \R_+^* \times \R$, then
\beq
\label{rd08_20e3}
    E[u^2(t, x)] = (c/\rho)^2\, \cK(t, x).
\eeq
\end{thm}

\begin{proof}
Let $f(t, x) = E[u^2(t, x)]$ By \eqref{rd08_14e9} and the Itô isometry \eqref{ch1'-s4.2},
 \beq\label{rd05_07e2}
     f(t, x) = I_0^2(t, x) + [J \star f](t, x).
 \eeq

If $I_0^2(t, x) = + \infty$, then $E[u^2(t, x)] = + \infty$ by Theorem \ref{rd08_20t1}, therefore \eqref{rd08_20e4} holds because all three terms there are nonnegative. If $I_0^2(t, x) < + \infty$, then $\lim_{n \to \infty} [J^{\star n} \star f](t, x) = 0$ as observed just after \eqref{rd05_05e1}, and as indicated just after \eqref{rd01_31e4}, $[\cK \star I_0^2](t, x) < +\infty$. We can therefore apply Lemma \ref{rd08_10l2} (d) with $a(t,x) = I_0^2(t,x)$ and $z_0 = 0$ there, to see that \eqref{rd08_20e4} holds.

Formula \eqref{rd08_20e2} is an immediate consequence of \eqref{rd08_20e4}. If $I_0(t, x) = c\, \Gamma(t, x)$, then $I_0^2(t, x) = (c/\rho)^2\, J(t, x)$. Therefore, \eqref{rd08_20e3} follows from \eqref{rd08_20e4} and \eqref{rd08_12e2}.
 \end{proof}

Formulas such as \eqref{rd08_20e4}, \eqref{rd08_20e2} and \eqref{rd08_20e3} are particularly useful when there is sufficient information about the function $\cK$. In fact, appealing to the results of Section \ref{rd08_26s1}, it turns out that when the partial differential operator $\cL$ is either the heat or the wave operator, the kernel $\cK$ can be computed explicitly (see \eqref{rd08_12e8} and \eqref{rd08_12e10} respectively). Hence, in these two cases, Proposition \ref{rd08_14p1} leads to exact expressions for the second moments of the solution $u(t,x)$. When $\cL$ is the fractional heat operator, Proposition \ref{rd08_12p1} (c) gives upper and lower bounds for $\cK$ (see \eqref{rd08_17e3} and \eqref{rd08_17e2} respectively), and therefore second moments of the solution can be estimated both from above and from below. These observations will be used in Subsection \ref{rd01_27ss2} below. The next proposition gives these explicit formulas.

\begin{prop}\label{rd05_05p1}
We consider the SPDE \eqref{rd08_14e1} and its solution $(u(t, x),\, (t, x) \in \R_+^* \times \R)$ given by Theorem \ref{rd08_20t1}. Fix $\nu \in \R_+^*$ and $c \in \R$. Let $\delta_0$ be the Dirac delta function.

(a) {\em Heat equation, constant initial data.} Let $\cL := \cL_{\nu} = \frac{\partial}{\partial t} - \frac{\nu}{2}\frac{\partial^2}{\partial x^2}$. If the initial data is $u_0(x) \equiv c$,  then
\beq\label{rd05_06e3}
     E[u^2(t, x)] = 2c^2 \exp\left(\frac{\rho^4\, t}{4\nu}\right) \Phi\left(\sqrt{\frac{\rho^4\, t}{2 \nu}}\, \right),
\eeq
where $\Phi$ is the standard Normal probability distribution function.

(b) {\em Heat equation, Dirac initial data.} Let $\cL = \cL_{\nu} = \frac{\partial}{\partial t} - \frac{\nu}{2}\frac{\partial^2}{\partial x^2}$. If the initial data is $u_0(x) = c\, \delta_0(x)$, then 
\beq\label{rd05_06e1}
     E[u^2(t, x)] = c^2 \left(\frac{1}{\sqrt{4 \pi\, \nu\, t}} + \frac{\rho^2}{2\nu}\, e^{\rho^4\, t/(4\nu)} \Phi\left(\sqrt{\frac{\rho^4\, t}{2 \nu}}\, \right) \right) \Gamma_{\nu/2}(t, x),
\eeq
where $\Gamma_{\nu/2}$ is defined in \eqref{fsnu}.

(c) {\em Wave equation, constant initial position, vanishing initial velocity.} Let $\cL = \frac{\partial^2}{\partial t^2} - \nu^2\, \frac{\partial^2}{\partial x^2}$. If the initial position is $u_0(x) \equiv c$ and the initial velocity vanishes, then
\beq\label{rd05_06e4}
     E[u^2(t, x)] = c^2 \cosh\left(t \, \sqrt{\frac{\rho^2}{2 \nu}}\, \right).
\eeq

(d) {\em Wave equation, vanishing initial position, Dirac initial velocity.} Let $\cL = \frac{\partial^2}{\partial t^2} - \nu^2\, \frac{\partial^2}{\partial x^2}$. If the initial position vanishes and the initial velocity is $v_0(x) = c\, \delta_0(x)$, then 
\beq\label{rd05_06e2}
     E[u^2(t, x)] =  \frac{c^2}{4 \nu^2} \, I_0\left(\sqrt{\frac{\rho^2}{2\, \nu}\, \frac{\left( (\nu\, t)^2 - x^2\right)}{\nu^2}}\, \right)\,
      1_{[0,\, \nu\, t]}(\vert x \vert),
\eeq
where  $I_0$ is the modified Bessel function of the first kind of order $0$ given in \eqref{rd08_26e5}.
\end{prop}

\begin{proof}
(a) In this case, we use \eqref{rd08_12e8} to see that 
\beq\label{rd05_05e2}
   \cK(t, x) = \cK^{\mbox{\scriptsize \rm h}}_\nu(t, x; \rho)  = B_\nu(t; \rho) \, \Gamma_{\nu/2}(t, x),
\eeq
where $\Gamma_{\nu/2}$ is defined in \eqref{fsnu} and
\beqn
   B_\nu(t; \rho) = \frac{\rho^2}{\sqrt{4 \pi\, \nu\, t}} + \frac{\rho^4}{2\nu}\, e^{\rho^4\, t/(4\nu)} \Phi\left(\sqrt{\frac{\rho^4\, t}{2 \nu}} \right).
\eeqn
Therefore,
\beqn
   [\cK \star \bfone](t, x) = \int_0^t ds \int_\R dy\, B_\nu(s; \rho) \, \Gamma_{\nu/2}(s, y) = \int_0^t B_\nu(s; \rho)\, ds.
\eeqn 
Integrating by parts the second term in the formula for $B_\nu(s; \rho)$, we find after simplification that
\beqn
   [\cK \star \bfone](t, x) = 2\, e^{\rho^4\, t/(4\nu)} \Phi\left(\sqrt{\frac{\rho^4\, t}{2 \nu}}\right) - 1.
\eeqn
Together with \eqref{rd08_20e2}, this proves (a). 

   (b) This follows immediately from \eqref{rd08_20e3} and \eqref{rd05_05e2}. 
   
   (c) In this case, we use \eqref{rd08_12e10} to see that 
\beq\label{rd05_05e4}
     \cK(t, x) = \cKw_\nu(t, x) = \frac{\rho^2}{4\, \nu^2}\, I_0\left(\sqrt{\frac{\rho^2}{2\, \nu}\, \frac{\left( (\nu\, t)^2 - x^2\right)}{\nu^2}} \right)\,
      1_{[0,\, \nu\, t]}(\vert x \vert),
\eeq
where  $I_0$ is given in \eqref{rd08_26e5}. By Lemma \ref{rd08_19l1}, 
\begin{align*}
   [\cK \star \bfone](t, x) &= \int_0^t ds \int_\R dy \, \cK(t, x) \\ &= \sqrt{\frac{\rho^2}{2 \nu}} \int_0^t ds\, \sinh\left(s\, \sqrt{\frac{\rho^2}{2 \nu}}\,\right) 
     = \cosh\left(t \, \sqrt{\frac{\rho^2}{2 \nu}} \right) - 1. 
\end{align*}
Together with \eqref{rd08_20e2}, this proves (c). 

(d)  This follows immediately from \eqref{rd08_20e3} and \eqref{rd05_05e4}. 
\end{proof}

\begin{remark}
Further explicit formulas are given in \cite{chendalang2015-2} and \cite{chen-dalang-2015}.
\end{remark}

\subsection{Moment Lyapounov exponents}\label{rd01_27ss1}
As discussed in Section \ref{ch1-1.2}, equation \eqref{rd08_14e1} is particularly interesting because it exhibits intermittency properties. In order to give a mathematical formulation of this notion, define the {\it upper and lower moment Lyapounov exponents}\index{moment Lyapounov exponent}\index{Lyapounov exponent!moment} of order $p>0$ by
\begin{align}\nonumber
 \bar\gamma_p(x)& := \limsup_{t\to\infty}\frac{1}{t}\log E\left[|u(t,x)|^p\right],\\
\underline{\gamma}_p(x)& :=\liminf_{t\to \infty}\frac{1}{t}\log E\left[|u(t,x)|^p\right].
\label{rd08_14e2}
\end{align}
Because we are assuming that $\Gamma(t, x; s, y) = \Gamma(t-s, x-y)$, when the $I_0$ is a constant (which is typically the case when the initial data is constant), then  these two exponents do not depend on $x$. In this case, following Bertini and Cancrini
\cite{bertini-cancrini-1995} and Carmona and Molchanov \cite{carmona-molchanov-1994}, we say that the solution $u$ of \eqref{rd08_14e1} is {\it intermittent}\index{intermittent} if
$\gamma_p:=\underline{\gamma}_p =\bar\gamma_p$ for all $p\in\N$ and if the following
strict inequalities are satisfied:
\begin{align}
\label{rd08_14e3}
     \gamma_1 < \frac{\gamma_2}{2} < \frac{\gamma_3}{3}< \cdots
\end{align}
Carmona and Molchanov \cite{carmona-molchanov-1994} pointed out the following facts.

\begin{lemma}
\label{rd08_13l1}
(a) The inequalities \eqref{rd08_14e3} are always satisfied if ``$<$" is replaced by ``$\leq$":
\beq
\label{rd08_14e4a}
     \gamma_1 \leq \frac{\gamma_2}{2} \leq \frac{\gamma_3}{3} \leq\cdots 
\eeq

(b) The function $p \mapsto \gamma_p$ from $\R_+$ to $\R_+$ is convex:
\beq
\label{rd08_14e4}
    \gamma_p \leq \half (\gamma_{p-h} + \gamma_{p+h})\qquad\text{if } p \geq 0 \text{ and } p \pm h \geq 0.
\eeq
\end{lemma}

\begin{proof}
(a) For $ p \geq 1$, $\Vert u(t, x) \Vert_{L^p(\Omega)} \leq \Vert u(t, x) \Vert_{L^{p+1}(\Omega)}$, that is,
\beqn
    \left(E\left[\vert u(t, x) \vert^p\right]\right)^{\frac{1}{p}} \leq \left(E\left[\vert u(t, x) \vert^{p+1}\right]\right)^{\frac{1}{p+1}},
\eeqn
and this implies \eqref{rd08_14e4a}.

   (b) Clearly,
\beqn
   \vert u(t, x) \vert^p = \vert u(t, x) \vert^{\frac{p+h}{2}}\, \vert u(t, x) \vert^{\frac{p-h}{2}},
\eeqn
and therefore, by the Cauchy-Schwartz inequality,
\beqn
   \left(E\left[\vert u(t, x) \vert^p\right]\right)^{2} \leq E\left[\vert u(t, x) \vert^{p-h}\right] \, E\left[\vert u(t, x) \vert^{p+h}\right].
\eeqn
Taking the logarithm of both sides gives \eqref{rd08_14e4}.
\end{proof}

A consequence of Lemma \ref{rd08_13l1} is the following observation.

\begin{lemma}
If $\gamma_1 = 0$ (which is the case for instance in the parabolic Anderson model: see Remark \ref{rd08_18r1}) and $\gamma_\ell > 0$ for some $\ell \geq 2$, then
\beq
\label{rd08_14e5}
      \frac{\gamma_\ell}{\ell} < \frac{\gamma_{\ell+1}}{\ell+1} < \frac{\gamma_{\ell+2}}{\ell+2} < \cdots
\eeq
\end{lemma}
\begin{proof}
Fix $\ell \geq 2$. Consider $(\ell, \gamma_\ell)$ as a point in the nonnegative quadrant. Then $\gamma_\ell/\ell$ is the slope of the line that passes through the origin and $(\ell, \gamma_\ell)$. Let $\ell$ be the smallest integer for which $\gamma_\ell > 0$. If $\gamma_\ell / \ell = \gamma_{\ell+1}/(\ell+1)$, then the point $(\ell+1, \gamma_{\ell+1})$ is on the line that passes through the origin and the point $(\ell, \gamma_\ell)$. But the convexity property \eqref{rd08_14e4} and the fact that $\gamma_{\ell-1} = 0$ imply that for $m > \ell$, the point $(m, \gamma_m)$ must be above the line that passes through $(\ell-1, 0)$ and $(\ell, \gamma_\ell)$, therefore, since $\ell - 1 >0$, strictly above the line that passes through the origin and $(\ell, \gamma_\ell)$. This contradiction proves \eqref{rd08_14e5}.
\end{proof}

 Property \eqref{rd08_14e5} led Carmona and Molchanov to give the following definition \cite[Definition III.1.1, p.~55]{carmona-molchanov-1994}.

\begin{def1}
\label{rd08_14d1}
Let $\ell$ be the smallest integer for which $\gamma_\ell>0$. If $\ell<\infty$, then we
say that the solution $u(t,x)$ exhibits {\em (asymptotic) intermittency of
order $\ell$,} and if $\ell=2$, then it exhibits {\em full intermittency}.
\end{def1}

For the relationship between ``full intermittency" and the physical notions of ``strong localization" and ``intermittency", see the discussion in Section \ref{ch1-1.2}. Because of Definition \ref{rd08_14d1}, much research has evolved around estimates and formulas for the moments $E[\vert u(t, x) \vert^\ell]]$, with an emphasis on the simplest case $\ell=2$. We have already seen some bounds on $E[\vert u(t, x) \vert^\ell]$ in Corollary \ref{5.3-cor1}.

\begin{remark}
\label{rd08_18r1}
In the case of the stochastic heat equation, if the initial condition is nonnegative, then the solution $u(t, x)$ is also nonnegative (see Subsection \ref{rd05_06ss1}). In particular, $E[\vert u(t, x)\vert] = E[u(t,x)] = I_0(t, x)$, so if $u_0(x) \equiv c \in \R_+^*$, then $I_0(t, x) \equiv c$ and $\gamma_1 = 0$.
\end{remark}


\subsection{The second moment Lyapounov exponent}\label{rd01_27ss2}
\medskip

The second moment Lyapounov exponent 
\beqn
     \gamma_2 = \lim_{t \to \infty}\, \frac{1}{t}\, \log E\left[u^2(t,x)\right],
\eeqn
when the limit exists and does not depend on $x$, gives us the growth rate of $E[u^2(t,x)]$ as $t \to \infty$: for large $t$, $E[u^2(t,x)]$ and $e^{\gamma_2\, t}$ have similar magnitudes.

In the next theorem, we compute $\gamma_2$ explicitly in the case of the stochastic heat and wave equations, for two types of initial data, and we give upper and lower bounds on $\gamma_2$ in the case of the fractional stochastic heat equation. 
\begin{thm}
\label{rd08_18t1}
We consider the SPDE \eqref{rd08_14e1} and its solution $(u(t, x),\, (t, x) \in \R_+^* \times \R)$ given by Theorem \ref{rd08_20t1}. 

(a) In the case of the heat operator 
\beqn
    \cL_\nu = \frac{\partial}{\partial t} - \frac{\nu}{2}\, \frac{\partial^2}{\partial x^2}, \qquad \nu > 0,
\eeqn
consider the following two cases: 
    
   $(1)$  {\em Dirac delta initial data.} $u_0(x) = c\, \delta_0(x)$, for all $x \in \R$, $c \neq 0$; 
    
    $(2)$  {\em Constant initial data.}  $u_0(x) = x_0$, for all $x\in\R$, $x_0 \neq 0$.
    
\noindent In both cases,
\beqn
     \gamma_2 = \frac{\rho^4}{4\nu}.
\eeqn

(b) In the case of the wave operator 
$$
    \cL_\nu = \frac{\partial^2}{\partial t^2} - \nu^2\, \frac{\partial^2}{\partial x^2}, \qquad \nu > 0,
$$ 
consider the following two cases:  

    $(1)$ {\em Dirac delta initial velocity.} The initial position is $u_0 \equiv 0$ and the initial velocity is $v_0(x) = c\, \delta_0(x)$, $x \in \R$, $c \neq 0$;
       
    $(2)$ {\em Constant initial data.}  For some $x_0 \in \R^*$, the initial position is $u_0(x) = x_0$, for all $x \in \R$, and the initial velocity is $v_0 \equiv 0$.

\noindent In both cases, 
\beqn
     \gamma_2 = \frac{\vert \rho \vert}{\sqrt{2\, \nu}}. 
\eeqn

(c) In the setting of the fractional stochastic heat equation on $\R$ driven by space-time white noise, with the operator 
$$
    \cL_\nu = \frac{\partial}{\partial t} - \nu\, \laplacef, \qquad \nu > 0,\qquad a \in \, ]1, 2[,\  \vert\delta\vert < 2 -a,
$$ 
consider the following two cases:  

    $(1)$ {\em Dirac delta initial data.} $u_0(x) = c\, \delta_0(x)$, for all $x \in \R$, $c \neq 0$;
      
    $(2)$ {\em Constant initial data.}  $u_0(x) = x_0$, for all $x\in\R$, $x_0 \neq 0$.

\noindent In both cases, there are positive constants $C_1 = C_1(a, \delta)$ and $C_2 = C_2(a, \delta)$ such that for all $\rho \neq 0$, 
\beqn
  C_1\, \frac{\rho^{2 a^*}}{\nu^{1/(a-1)}} \leq \gamma_2 \leq C_2\, \frac{\rho^{2 a^*}}{\nu^{1/(a-1)}}.
\eeqn
Here, $a^*$ is the dual of $a$: $\frac{1}{a} + \frac{1}{a^*} = 1$.
\end{thm}

\begin{proof}

For Case (1) in part (a), \eqref{rd05_06e1} tells us that
\beq\label{rd01_29e3}
   E\left[u^2(t, x)\right] = c^2 \left(\frac{1}{\sqrt{4 \pi\, \nu\, t}} + \frac{\rho^2}{2\nu}\, e^{\rho^4\, t/(4\nu)} \Phi\left(\sqrt{\frac{\rho^4\, t}{2 \nu}} \right) \right) \Gamma_{\nu/2}(t, x),
\eeq
so since $c \neq 0$,
\beqn
  \gamma_2 = \lim_{t \to \infty} \frac{1}{t}\, \log E\left[u^2(t,x)\right]  = \lim_{t \to \infty} \frac{1}{t} \log\left(e^{\rho^4\, t/(4\nu)} \right) = \frac{\rho^4}{4 \nu}\, .
\eeqn

For Case (1) in part (b), \eqref{rd05_06e2} tells us that
\beq\label{rd01_29e4}
   E\left[u^2(t, x)\right]  = \frac{c^2}{4 \nu^2} \, I_0\left(\sqrt{\frac{\rho^2}{2\, \nu}\, \frac{\left( (\nu\, t)^2 - x^2\right)}{\nu^2}}\, \right)\,
      1_{[0,\, \nu\, t]}(\vert x \vert),
\eeq
so since $c \neq 0$, we use \eqref{rd08_17e5} to see that
\beqn
  \gamma_2 = \lim_{t \to \infty} \frac{1}{t} \log\left(\exp\left(\sqrt{\frac{\rho^2}{2\, \nu}\, \frac{\left( (\nu\, t)^2 - x^2\right)}{\nu^2}} \right) \right) =  \frac{\vert \rho \vert}{\sqrt{2\, \nu}}.
\eeqn

For Case (1) in part (c), \eqref{rd08_17e3} and \eqref{rd08_17e2} give us upper and lower bounds for $\cK(t, x)$ of the form 
\beq\label{rd01_29e5}
    c_i\, \rho^2\,  E_{1/a^*,1/a^*}\left(c_i\, \rho^2\, \nu^{-1/a}\, t^{1/a^*}\right)\, h_i(t,x), \qquad i = 1, 2,
\eeq
where $c_1 > 0$, $c_2 > 0$, $E_{1/a^*,1/a^*}(\cdot)$ is the two-parameter Mittag-Lefler function (see \eqref{E:Mittag-Leffler}), $h_1(t, x) = (\nu\, t)^{-1/a}\, \laplacef(\nu\, t, x)$ and $h_2(t, x) = g_a(\nu\, t, x)$ from \eqref{rd08_17e6}. In view of the asymptotic property of the Mittag-Leffler function given in \eqref{rd08_17e4}, 
\beqn
   \lim_{t \to \infty} \frac{1}{t} \log\left(\exp\left(\left[c_i\, \rho^2\, \nu^{-1/a}\,  t^{1/a^*} \right]^{1/(1/a^*)} \right) \right) = c_i^{a^*} \rho^{2 a^*} \nu^{-a^*/a} = c_i^{a^*} \frac{\rho^{2 a^*}}{\nu^{1/(a-1)}}.
\eeqn
This proves case (1) of parts (a), (b) and (c). 

For case (2) of part (a), 
\eqref{rd05_06e3} tells us that
\begin{align*}
    E\left[u^2(t, x)\right] &= 2x_0^2 \exp\left(\frac{\rho^4\, t}{4\nu}\right) \Phi\left(\sqrt{\frac{\rho^4\, t}{2 \nu}}\, \right) .
\end{align*}
Therefore, since $x_0 \neq 0$,
\beqn
 \gamma_2 = \lim_{t \to \infty} \frac{1}{t}\, \log E\left[u^2(t,x)\right]  = \lim_{t \to \infty} \frac{1}{t}\, \log \exp\left(\frac{\rho^4\, t}{4\nu}\right) = \frac{\rho^4}{4\, \nu}.
\eeqn
This completes the proof of case (2) in part (a).

For case (2) of part (b), \eqref{rd05_06e4} tells us that
\beqn
   E\left[u^2(t,x)\right] = x_0^2\, \cosh\left(t\, \sqrt{\frac{\rho^2}{2 \nu}} \right).
\eeqn
Therefore, since $x_0 \neq 0$,
\beqn
    \gamma_2 = \lim_{t \to \infty}\, \frac{1}{t} \, \log E\left[u^2(t,x)\right] = \lim_{t \to \infty}\, \frac{1}{t} \, \log  \cosh\left(t\, \sqrt{\frac{\rho^2}{2 \nu}} \right) = \frac{\vert \rho \vert}{\sqrt{2 \nu}}\, .
\eeqn
This completes the proof of case (2) in part (b).

  For case (2) of part (c), we notice that $I_0(t, x) \equiv x_0$, so we use \eqref{rd08_20e2} in Proposition \ref{rd08_14p1} to see that 
\begin{align*}
    E[u^2(t, x)] &= x_0^2 \, [\bfone + \cK \star \bfone](t, x) = x_0^2 + x_0^2 \int_0^t ds \int_\re dy\, \cK(s, y).
\end{align*}
We get upper and lower bounds on $\gamma_2$ by replacing $\cK(s, y)$ in this formula by \eqref{rd01_29e5}. For $i = 1, 2$, let
\beqn
   \gamma_{2, i} := \lim_{t \to \infty} \frac{1}{t} \log\left(c_i\,  \rho^2 \int_0^t ds\,  E_{1/a^*,1/a^*}\left(c_i\, \rho^2\, \nu^{-1/a}\,s^{1/a^*}\right) \int_\R dy \, h_i(s,y)\right).
\eeqn
Notice that $\int_\R dy \, h_i(s,y)$ is equal to $(\nu\, s)^{-1/a}$ for $i=1$ and to a constant $k_2/\nu$ that does not depend on $s$ when $i=2$. We use \eqref{rd08_17e4} with $\alpha = \beta = 1/a^* \in\, ]0, \half[\, \subset\, ]0, 2[$ to obtain the asymptotic formula
\beqn
     E_{1/a^*, 1/a^*}(s) 
         = a^* \, s^{a^* - 1}\, \exp\left(s^{a^*}\right) + O\left(s^{-2}\right),\qquad\text{as } s \longrightarrow \infty.
\eeqn
Therefore, for any $t_0 > 0$,
\begin{align*}
   \gamma_{2, 1} &= \lim_{t \to \infty}  \frac{1}{t} \log\Big[\int_{t_0}^t ds\, s^{-1/a} \, \Big[(c_1\, \rho^2\, \nu^{-1/a}\, s^{1/a^*})^{a^* - 1} \exp\big((c_1\, \rho^2\, \nu^{-1/a}\, s^{1/a^*})^{a^*}\big) \\
    &\qquad\qquad\qquad\qquad\qquad\qquad + O\big((s^{1/a^*})^{-2}\big) \Big]  \Big] \\
      &= \lim_{t \to \infty} \frac{1}{t} \log\left[\int_{t_0}^t ds\, \exp\big((c_1\, \rho^2\, \nu^{-1/a})^{a^*} \, s\big) \right].
\end{align*}
Evaluating the integral, we see that $\gamma_{2, 1} = c_1^{a^*}\, \rho^{2 a^*}/\nu^{1/(a-1)}$. Similarly, for any $t_0 > 0$,
\begin{align*}
   \gamma_{2, 2} &= \lim_{t \to \infty}  \frac{1}{t} \log\Big[\int_{t_0}^t ds\, \Big[(c_2\, \rho^2\, \nu^{-1/a}\, s^{1/a^*})^{a^* - 1} \exp\big((c_2\, \rho^2\, \nu^{-1/a}\, s^{1/a^*})^{a^*}\big) \\
    &\qquad\qquad\qquad\qquad\qquad\qquad + O\big((s^{1/a^*})^{-2}\big) \Big]  \Big] \\
      &= \lim_{t \to \infty} \frac{1}{t} \log\left[\int_{t_0}^t ds\, s^{1 - 1/a^*}\, \exp\big((c_2\, \rho^2\, \nu^{-1/a})^{a^*} \, s\big) \right].
\end{align*}
Bounding the factor $s^{1 - 1/a^*}$ above and below by $t^{1 - 1/a^*}$ and $t_0^{1 - 1/a^*}$, respectively, we see that $\gamma_{2, 2} = c_2^{a^*}\, \rho^{2 a^*}/\nu^{1/(a-1)}$. This completes the proof of case (2) in part (c).
\end{proof}

\begin{remark}\label{rd08_19r1}
When, in the SPDE \eqref{rd08_14e1}, we replace the space-time white noise $\dot W$ by a smooth function $f(t, x)$, then the solution $u(t, x)$ can be expressed via a Feynman-Kac-type formula (see for instance \cite[Chap.~4, Theorem 4.2, p.~268]{ks} for the heat equation and \cite[Theorem 3.2]{DMT2008} for the wave equation). This type of formula is also available for certain SPDEs with noises that are smoother than space-time white noise (see for instance \cite{hln}). These formulas can often be used to determine the growth rate of $E[u^2(t, x)]$ as $t \to \infty$ (see e.g.~\cite{dalang-mueller2009}).
\end{remark}

\subsection{Propagation speed of intermittency}\label{rd06_24ss1}

The conclusions of Theorem \ref{rd08_18t1} can be interpreted as follows. For fixed $x \in \R$, as $t \to + \infty$, $E[u^2(t, x)] \sim e^{\gamma_2\, t} \to + \infty$. On the other hand, in the case of the stochastic heat or wave equations considered in Theorem \ref{rd08_18t1} with the initial conditions involving $\delta_0$, by Proposition \ref{rd05_05p1}, for fixed $t > 0$, $E[u^2(t, x)] \sim \Gamma_\nu^2(t, x) \to 0$ as $x \to + \infty$.
\medskip

\noindent{\em Initial data equal to $\delta_0$}
\medskip

The considerations above suggest that when the initial conditions have compact support, then for $\alpha \geq 0$ and $t \to +\infty$, $E[u^2(t, \alpha\, t)]$ will grow exponentially when $\alpha$ is small and will converge to $0$ when $\alpha$ is large. In particular, an observer moving away from the origin at speed $\alpha$ will observe the intermittency phenomenon when moving at slow enough speed, but not when moving at large enough speed. In the case of the heat equation, even though the initial condition propagates at infinite speed, we will see that the intermittency phenomenon only propagates at finite speed (hence the title of this subsection). 

We illustrate this behavior for the two situations considered in cases (1) of Theorem \ref{rd08_18t1} (a) and (b).

\begin{prop} \label{rd05_14p1}
We consider the SPDE \eqref{rd08_14e1} and its solution $(u(t, x),\, (t, x) \in \R_+^* \times \R)$ given in Theorem \ref{rd08_20t1}. 

 (a) {\em Stochastic heat operator on $\R$ with Dirac delta initial condition.} In the setting of Theorem \ref{rd08_18t1} (a) with initial data $u_0(x) = \delta_0(x)$, let
\beqn
    \alphah = \frac{\rho^2}{2}.
\eeqn
For $0\leq \alpha < \alphah$,  
\beqn
    \lim_{t \to \infty} \frac{1}{t} \log E\left[u^2(t, \alpha\, t)\right] = \frac{\rho^4 - 4 \, \alpha^2}{4\, \nu} > 0 .
\eeqn
For $\alpha > \alphah$, $\lim_{t \to \infty}  E\left[u^2(t, \alpha\, t)\right] = 0$.
 
(b) {\em Stochastic wave operator on $\R$ with Dirac delta initial velocity and vanishing initial position}. In the setting of Theorem \ref{rd08_18t1} (b) with initial position  $u_0(x) \equiv 0$ and initial velocity $v_0(x) = \delta_0(x)$, let
\beqn
    \alphaw = \nu.
\eeqn
For $0 \leq \alpha < \alphaw$,
\beqn
    \lim_{t \to \infty} \frac{1}{t} \log E\left[u^2(t, \alpha\, t)\right] = \vert \rho \vert\, \sqrt{\frac{\nu^2 -  \alpha^2}{2\, \nu^3}} > 0 .
\eeqn
For $\alpha > \alphaw$, $\lim_{t \to \infty}  E\left[u^2(t, \alpha\, t)\right] = 0$.
\end{prop}

\begin{proof}
(a) By \eqref{rd05_06e1}, for $\alpha \geq 0$ and as $t \to + \infty$,
\beqn
    E[u^2(t, \alpha\, t)] \sim \exp\left(\frac{\rho^4\, t}{4\, \nu} \right) \Gamma_{\nu/2}(t, \alpha\, t) \sim \exp\left(\left(\frac{\rho^4}{4\, \nu} - \frac{\alpha^2}{\nu} \right) t \right).
\eeqn
Since the exponent is positive when $\alpha < \rho^2/2 = \alphah$, and negative when $\alpha > \alphah$, the conclusion in (a) follows.

(b) By \eqref{rd05_06e2}, when $\alpha > \nu$, $E[u^2(t, \alpha\, t)] = 0$, and when $0 \leq \alpha < \nu$,
\beqn
   E[u^2(t, \alpha\, t)] =  \frac{c^2}{4 \nu^2} \, I_0\left(\sqrt{\frac{\rho^2}{2\, \nu}\, \frac{\left( (\nu\, t)^2 - (\alpha \,t)^2\right)}{\nu^2}}\, \right) \sim \exp\left(\vert \rho \vert\, t\, \sqrt{\frac{\nu^2 - \alpha^2}{2\, \nu^3}} \right)
\eeqn
by \eqref{rd08_17e5}. This establishes (b).
\end{proof}





\medskip

\noindent{\em Initial data with exponential decay at $\pm \infty$}
\medskip

If the initial condition has exponential decay at $\pm \infty$, say $u_0(x) =  e^{- \beta\, \vert x \vert}$, or $v_0(x) = f_\beta(x)$, $x \in \R$, in the stochastic heat or wave equations, respectively, then somewhat surprisingly, the speed of propagation of intermittency changes unless $\beta$ is large enough. This indicates that even an exponentially decaying perturbation of the initial condition may have an everlasting macroscopic effect on the properties of the solution. We illustrate this statement for the two situations considered in Proposition \ref{rd05_14p1}.

\begin{thm} \label{rd05_14p2}

 (a) {\em Stochastic heat equation.} In the setting of Theorem \ref{rd08_18t1} (a) with initial data $u_0(x) = e^{- \beta\, \vert x \vert}$, $x \in \R$, $\beta > 0$, let
\beqn
    \alphah_\beta = \left \{ \begin{array}{ll}
       \frac{\beta\, \nu}{2} + \frac{\rho^4}{8\, \beta\, \nu} & \text{if } 0 < \beta <  \frac{\rho^2}{2 \nu} ,\\
          \frac{\rho^2}{2} & \text{if } \beta \geq  \frac{\rho^2}{2 \nu}.
    \end{array} \right.
\eeqn
For $0< \alpha < \alphah_\beta$, 
\beq\label{rd06_19e12}
    \liminf_{t \to \infty} \frac{1}{t} \log E\left[u^2(t, \alpha\, t)\right] > 0,
\eeq
and for $\alpha > \alphah_\beta$, 
\beq\label{rd06_19e7}
   \lim_{t \to \infty}  E\left[u^2(t, \alpha\, t)\right] = 0.
\eeq
 
(b) {\em Stochastic wave equation.} In the setting of Theorem \ref{rd08_18t1} (b) with initial position  $u_0(x) \equiv 0$ and initial velocity $v_0(x) = e^{- \beta\, \vert x \vert}$, $x \in \R$, $\beta > 0$, let
\beq\label{rd06_21e1}
    \alphaw_\beta =  \nu \left(1 + \frac{\rho^2}{8\, \beta^2\, \nu^3} \right)^{\half}.
\eeq
For $0 < \alpha < \alphaw_\beta$, 
\beq\label{rd06_20e3}
    \liminf_{t \to \infty} \frac{1}{t} \log E\left[u^2(t, \alpha\, t)\right] > 0,
\eeq
and for $\alpha > \alphaw_\beta$, 
\beq\label{rd06_20e4}
   \lim_{t \to \infty}  E\left[u^2(t, \alpha\, t)\right] = 0.
\eeq
\end{thm}

\begin{proof}  (a) In the first part of the proof, $\Gamma_\nu(t, x)$ denotes the heat kernel given in \eqref{fsnu} and we write $\cK(t, x)$ instead of $\cKh_\nu(t, x)$, given in \eqref{rd05_05e2}. Observe that with the initial condition $u_0(x) = e^{- \beta \vert x \vert}$, $x \in \R$, 
\beq\label{rd06_19e3}
    I_0(t, x) = \int_\R \Gamma_\nu(t, x - z)\, e^{- \beta \vert z \vert}\, dz = \left(e^{-\beta\, \vert * \vert} * \Gamma_\nu(t, *) \right)(x).
\eeq
By a direct calculation, one finds that
\beq\label{rd06_19e8}
    I_0(t, x) = e^{\beta^2 \nu t/2} H_{\nu t, \beta}(x),
\eeq
where, for $a, b \in \R_+^*$ and $x \in \R$,
\beq\label{rd06_19e9}
    H_{a, b}(x) = e^{b x}\, \Phi\left(\frac{-a b - x}{\sqrt{a}} \right) + e^{- b x}\, \Phi\left(\frac{-a b + x}{\sqrt{a}}  \right),
\eeq
and $\Phi$ denotes the standard Normal probability distribution function. We will use this formula in the case $0 < \alpha < \alphah_\beta$ below.

   By \eqref{rd08_20e4}, $E[u^2(t, \alpha t)] = I_0^2(t, \alpha t) + [\cK \star I_0^2](t, \alpha t)$. Noticing that 
\begin{align*}
   I_0(t, \alpha t) &=  \int_\R \frac{1}{\sqrt{2 \pi \nu t}} \exp\left(- \frac{(\alpha t - z)^2}{2 \nu t}\right) \, e^{- \beta \vert z \vert}\, dz  \\
    &=  \int_\R \frac{1}{\sqrt{2 \pi \nu t}} \exp\left(- \frac{\alpha^2}{2 \nu }\,  t \, \left(1 - \frac{z}{\alpha t}\right)^2\right) \, e^{- \beta \vert z \vert}\, dz ,
\end{align*}
we see using dominate convergence that $\lim_{t \to \infty} I_0^2(t, \alpha t) = 0$, since $\alpha > 0$. Therefore, we examine $[\cK \star I_0^2](t, \alpha t)$.

\smallskip

\noindent{\em The case $\alpha > \alphah_\beta$.}  Fix $\gamma \in \, ]0, \beta[$. We claim that for all $(t, x) \in \R_+ \times \R$,
\beq\label{rd06_19e4}
    I_0(t, x) \leq \left\{\begin{array}{ll}
       \frac{1}{\sqrt{2 \pi \nu t}}  \exp\left(- \frac{x^2}{2 \nu\, t}\right), & \text{if } \vert x \vert < \gamma \nu\, t, \\[6pt]
        \frac{1}{\sqrt{2 \pi \nu t}} \exp\left(- (\gamma\, \vert x \vert  - \half \gamma^2\, \nu\, t) \right), &\text{if } \vert x \vert > \gamma \nu\, t.
         \end{array}
         \right.
\eeq
Indeed, by \eqref{rd06_19e3},
\begin{align*}
   I_0(t, x) \leq \sup_{y \in \R} \left(\Gamma_\nu(t, x-y)\, e^{- \gamma \, \vert y \vert}\right) \int_\R e^{(\gamma - \beta)\, \vert y \vert}\, dy.
\end{align*}
The integral converges because $\gamma < \beta$, and in order to obtain the supremum over $y$, we notice that it is equal to
\beqn
   \frac{1}{\sqrt{2 \pi \nu t}} \exp\left( - \inf_{y \in \R} \left(\frac{(x - y)^2}{2 \nu t} +\gamma\, \vert y \vert \right) \right),
\eeqn
and that
\beqn
  \inf_{y \in \R} \left(\frac{(x - y)^2}{2 \nu t} +\gamma\, \vert y \vert \right) = \left\{\begin{array}{ll}
         \frac{x^2}{2 \nu\, t}, & \text{if } \vert x \vert < \gamma \nu\, t, \\[6pt]
         \gamma\, \vert x \vert  - \half \gamma^2\, \nu\, t,  &\text{if } \vert x \vert > \gamma \nu\, t.
         \end{array}
         \right.\eeqn
This proves \eqref{rd06_19e4}.

Applying \eqref{rd06_19e2} below, we see that
\begin{align*}
  [\cK \star I_0^2](t, \alpha\, t) \leq \rho^2 \, \sqrt{\frac{\pi t}{\nu}}\, I_0^2(2t, \alpha\, t) \left(1 + 2 \exp\left(\frac{\rho^4\, t}{4 \nu} \right) \right).
\end{align*}
Therefore, for $\alpha > \beta \nu$ and $\gamma \in\, ]0, \beta[$, by \eqref{rd06_19e4},
\begin{align}\label{rd06_19e5}
   E\left[u^2(t, \alpha\, t)\right] &\lesssim 
   \exp\left(- 2\gamma \alpha\, t + \gamma^2 \nu\, t + \frac{\rho^4\, t}{4 \nu} \right), 
\end{align}
where the relation $f(t) \lesssim g(t)$ means here that $\limsup_{t \to \infty} \tfrac{1}{t}\log f(t) \leq \liminf_{t \to \infty} \tfrac{1}{t}\log g(t)$.
The exponent is negative when 
\beq\label{rd06_19e6}
\alpha > \frac{\gamma \nu}{2} + \frac{\rho^4}{8 \gamma \nu} =: \varphi(\gamma).
\eeq
We notice that the function $\gamma \mapsto \varphi(\gamma)$ is decreasing for $\gamma < \rho^2/{2 \nu}$ and increasing when $\gamma > \rho^2/{2 \nu}$, with minimum value $\rho^2/{2}$. 

Therefore, if $\beta \geq \rho^2/{2 \nu}$ and $\alpha > \rho^2/{2} = \alphah_\beta$, then there is $\gamma \in \, ]0, \beta[$ for which $\alpha > \varphi(\gamma)$ (when $\beta > \rho^2/{2 \nu}$, we can simply take $\gamma = \rho^2/{2 \nu}$). From \eqref{rd06_19e5} and \eqref{rd06_19e6}, we conclude that \eqref{rd06_19e7} holds.

On the other hand, if $0 < \beta < \rho^2/{2 \nu}$, then for $\alpha > \varphi(\beta) = \alphah_\beta$, we also conclude from \eqref{rd06_19e5} and \eqref{rd06_19e6} that \eqref{rd06_19e7} holds.

We have proved that for all $\beta > 0$, $\lim_{t \to \infty}  E\left[u^2(t, \alpha\, t)\right] = 0$ when $\alpha > \alphah_\beta$.
\medskip

\noindent{\em The case $0 < \alpha < \alphah_\beta$.}  By \eqref{rd08_20e4},
\beqn
   E[u^2(t, \alpha\, t)] \geq [\cK \star I_0^2](t, \alpha\, t).
\eeqn
   
We are going to apply \eqref{rd06_21e4} and \eqref{rd06_21e5} below, recalling the property $\Phi(y) \sim e^{-y^2/2}/(\vert y \vert \sqrt{2\pi})$ as $y \to -\infty$ (see \cite[7.12.1]{olver}). Let $c \in \, ]0, 1[$. By \eqref{rd06_21e4}, if $\alpha \geq \beta \nu$, then
\begin{align*}
   [\cK \star I_0^2](t, \alpha\, t) &\gtrsim e^{-2\, \alpha \beta\, t + \beta^2 \, \nu\, t} \, \exp\left(- \beta^2 \nu\, (1-c)\, t \right) \, 
       e^{\rho^4\, t/(4\nu)}  .
\end{align*}
while by \eqref{rd06_21e5}, if $\alpha < \beta \nu$, then
\begin{align*}
   [\cK \star I_0^2](t, \alpha\, t) &\gtrsim e^{2 \alpha \beta\, t} \, \exp\left(- \half\left(\alpha \sqrt{\frac{2}{\nu}} + \beta \sqrt{2 \nu}\right)^2 t \right) \\ 
    &\qquad\qquad\times \exp\left(- \beta^2 \nu\, (1-c)\, t) \right)\,  e^{\rho^4 t/(4\nu)}  \\
     &= \exp\left[\left(-\tfrac{\alpha^2}{\nu} + (c-1)\beta^2 \nu + \tfrac{\rho^4}{4 \nu} \right) t \right].
\end{align*}
We conclude that
\beq\label{rd06_22e1}
   \frac{1}{t}\, \log E[u^2(t, \alpha\, t] \gtrsim \left\{\begin{array}{ll}       
      c \beta^2 \nu + \frac{\rho^4}{4 \nu} - 2 \alpha \beta,&\text{if } \alpha \geq \beta\nu \\
        (c-1)\, \beta^2
 \nu + \frac{\rho^4}{4 \nu}  - \frac{\alpha^2}{\nu}, & \text{if } \alpha < \beta\nu.\\    
      \end{array}\right.
\eeq

Suppose that 
\beqn
   \beta <  \frac{\rho^2}{2 \nu}\qquad \text{and}\qquad 0 < \alpha < \alphah_\beta = \frac{\beta \nu}{2}+ \frac{\rho^4}{8 \beta \nu}.
\eeqn
Then for some $c \in \,]0, 1[$ close to $1$, we have $\beta < \frac{\rho^2}{2 \nu \sqrt{2-c}}$ and $\alpha <\frac{c\,  \beta \nu}{2}+
\frac{\rho^4}{8 \beta \nu}$. These two
 inequalities are respectively equivalent to 
\beqn
   \beta \nu < \frac{c\, \beta \nu}{2}+ \frac{\rho^4}{8 \beta \nu}  \qquad \text{and}\qquad  c\beta^2\nu+\frac{\rho^4}{4\nu} -2 \beta\alpha >0.
\eeqn
In particular, the right-hand side of \eqref{rd06_22e1} is positive for $\beta \nu \leq \alpha <\frac{c\: \nu \beta}{2}+ \frac{\rho^4}{8\nu \beta}$, therefore \eqref{rd06_19e12} holds in this case.
For $0 < \alpha < \beta \nu$, since we are in the case where $\beta \nu < \rho^2/2$, we have $\alpha < \rho^2/2$. In particular, for $c$ close enough to $1$, 
\beq\label{rd06_25e1}
   \alpha< \sqrt{\tfrac{\rho^4}{4}+(c-1)\beta^2\nu^2}.
\eeq
This inequality is equivalent to
\beq\label{rd06_25e2}
   (c-1)\beta^2\nu+ \frac{\rho^4}{4\nu} -\frac{\alpha^2}{\nu}>0.
\eeq
In particular, the right-hand side of \eqref{rd06_22e1} is positive in this case also. Therefore, when $\beta < \rho^2/2$, we have proved \eqref{rd06_19e12} for all $\alpha \in \, ]0, \alphah_\beta[$. 

Now suppose that 
\beqn 
   \beta\geq  \frac{\rho^2}{2\nu} \qquad \text{and}\qquad 0 < \alpha < \alphah_\beta = \frac{\rho^2}{2}.
\eeqn
Then for all $c \in \,]0, 1[$, $\beta >  \frac{\rho^2}{2\nu\sqrt{2-c}}$, and there is some $c \in \,]0, 1[$ such that \eqref{rd06_25e1} holds. 
These two inequalities are respectively equivalent to 
$$
   \sqrt{\tfrac{\rho^4}{4}+(c-1)\beta^2\nu^2} \leq \beta \nu 
$$ 
and \eqref{rd06_25e2}.
In particular, $\alpha < \beta \nu$ in this case. By \eqref{rd06_22e1}, \eqref{rd06_19e12} holds in this case.
This shows that for all $\beta > 0$, \eqref{rd06_19e12} holds when $0 < \alpha < \alphah_\beta$.
The proof of part (a) of Theorem \ref{rd05_14p2} is complete.

(b) In this second part of the proof, $\Gamma_\nu(t, x)$ denotes the wave kernel given in \eqref{rd08_21e6b} and we write $\cK(t, x)$ instead of $\cKw_\nu(t, x)$ given in \eqref{rd05_05e4}. In the case of vanishing initial condition and initial velocity $v_0(x) = e^{- \beta\, \vert x \vert}$, a direct calculation using d'Alembert's formula \eqref{p93.1} gives
\beq\label{rd06_23e1}
    I_0(t, x) = \left\{\begin{array}{ll}
         \frac{1}{\beta \nu} \, e^{- \beta \vert x \vert} \sinh(\beta\, \nu\, t), & \text{if } \vert x \vert \geq \nu\, t, \\
         \frac{1}{\beta \nu} \, \left(1 - e^{- \beta \vert x \vert} \cosh(\beta\, x) \right) , & \text{if } \vert x \vert \leq \nu\, t
        \end{array}\right.
\eeq
(see \cite[Lemma 4.4.5]{chen} for details). Fix $\gamma \in \, ]0, \beta[$. We claim that for all $(t, x) \in \R_+ \times \R$,
\beq\label{rd06_20e1}
   I_0^2(t, x) \leq C_{\beta, \gamma}\, e^{2\gamma  \nu\, t -  2\gamma\, \vert x \vert }, \quad \text{with } C_{\beta, \gamma} = \left(\frac{1}{2 \nu} \int_\R e^{(\gamma - \beta)\, \vert y \vert }\, dy\right)^2 < \infty.
\eeq
Indeed, using the triangle inequality,
\begin{align*}
   e^{\gamma \vert x \vert} \left(e^{-\beta \vert * \vert} * \Gamma_\nu(t, *) \right)(x) &= \frac{1}{2 \nu} \int_{x - \nu\, t}^{x + \nu\, t} dy\, e^{- \beta\, \vert y \vert}\, e^{\gamma\, \vert x \vert} \\
   &\leq \frac{1}{2 \nu} \int_{x - \nu\, t}^{x + \nu\, t} dy\, e^{(\gamma - \beta)\, \vert y \vert} e^{\gamma\, \vert x- y \vert} \\
   &\leq e^{\gamma\, \nu\, t} \, \frac{1}{2 \nu} \int_\R dy\, e^{(\gamma - \beta)\, \vert y \vert} .
\end{align*}
This proves \eqref{rd06_20e1}.

By \eqref{rd08_20e4}, \eqref{rd06_20e1} and \eqref{rd06_20e2} below, 
\beqn
  E[u^2(t, x)]  = I_0^2(t, x) + [\cK \star I_0^2](t, x)  \leq C_{\beta, \gamma}\, e^{2\gamma \nu\, t -  2\gamma\, \vert x \vert } + C'\, t\, e^{\sigma_\gamma \, \nu\, t -  2\gamma\, \vert x \vert },
\eeqn
where $\sigma_\gamma = \sqrt{4 \gamma^2 + \rho^2/(2\nu^3)}$. 
Since $\sigma_\gamma > 2 \gamma$,
\beqn
    E[u^2(t, \alpha t)] \lesssim \exp(\sigma_\gamma\, \nu - 2 \alpha\, \gamma)\, t).
\eeqn
The factor in parentheses is negative when 
\beqn
    \alpha > \frac{\sigma_\gamma \, \nu}{2\gamma} = \alphaw_\gamma = \nu \left(1 + \frac{\rho^2}{8 \gamma^2\, \nu^3} \right)^\half .
\eeqn
For $\alpha > \alphaw_\beta$, there is $\gamma \in \, ]0, \beta[$ such that $\alpha > \alphaw_\gamma$,
and in this case, \eqref{rd06_20e4} holds.

In order to prove \eqref{rd06_20e3} for $0 < \alpha < \alphaw_\beta$, we use \eqref{rd08_20e4} to write
\beqn
    E[u^2(t, \alpha\, t)] \geq [\cK \star I_0^2](t, \alpha\, t).
\eeqn

Suppose first that $0 < \alpha \leq \nu$. By \eqref{rd06_21e2} below,  for all $b \in \, ]0, 1[$, using \eqref{rd08_17e5}, we see that
\beqn
    E[u^2(t, \alpha\, t)] \gtrsim \exp\left(\left( 2  \beta\nu (1-b) + b \, \vert \rho \vert\, \nu\, \sqrt{\tfrac{1 - a^2}{2 \nu^3}}\right) t\right).
\eeqn
Since the right-hand side is positive, we see that \eqref{rd06_20e3} holds when $0 < \alpha \leq \nu$.

We now consider the case $\nu < \alpha < \alphaw_\beta$. Fix $a, b \in\,]0, 1[$. By \eqref{rd06_23e2},
\beqn
    E[u^2(t, \alpha\, t)] \gtrsim \exp\left( f(\alpha)\, t\right),
\eeqn
where
\beqn
   f(\alpha) := - 2 \alpha\beta + 2 a b \beta \nu + 2(1-b) \beta \nu + b \rho \nu\, \sqrt{\tfrac{1-a^2}{2 \nu^3}}.
\eeqn
We notice that $f(\alpha) > 0$ when $\alpha < g(a)$, where
\beqn
    g(a) = a b \nu + (1-b) \nu + \frac{b \rho \nu}{\beta} \sqrt{\frac{1 - a^2}{8 \nu^3}}.
\eeqn
The function $a \mapsto g(a)$, $a \in \, ]0, 1[$, attains its maximum at 
\beq\label{rd06_21e3}
    a = \sqrt{\frac{8  \beta^2 \nu^3}{\rho^2 + 8\beta^2  \nu^3 }},
\eeq
and it maximum value is 
\beqn
    g_{\text{max}} = (1-b) \nu + b \nu \, \sqrt {1 + \tfrac{\rho^2}{8 \beta^2 \nu^3}} = (1-b) \nu + b \nu\, \alphaw_\beta.
\eeqn
For $\alpha < \alphaw_\beta$, there is $b \in \, ]0, 1[$ such that $\alpha < g_{\text{max}}$. In this case, choosing $a$ as in \eqref{rd06_21e3}, we see that $f(\alpha) > 0$. This completes the proof of \eqref{rd06_20e3} and of Theorem \ref{rd05_14p2}.
\end{proof}

\begin{remark}
(1) Using the comparison theorem (Theorem \ref{rd08_13t1}), the conclusions of Theorem \ref{rd05_14p2} (a) can be extended to other initial conditions with exponential decay; see also \cite[Theorem 2.12]{chendalang2015-2}.
 
 (2) There is no comparison theorem in the case of the stochastic wave equation, but the results of Theorem \ref{rd08_13t1} (b) do remain valid for a more general class of initial conditions: see \cite[Theorem 2.9]{chen-dalang-2015}. It is interesting to note that in the setting of Theorem \ref{rd08_13t1} (b), an observer who moves at a speed $\alpha$ that is greater that the speed $\nu$ of waves may still observe the intermittency phenomenon, since $\alphah_\beta > \nu$.
\end{remark}

\medskip

\noindent{\em Propagation of intermittency in the fractional heat equation}
\medskip

In the fractional stochastic heat equation, the phenomenon of propagation of intermittency is substantially different than in the stochastic heat and wave equations, as was noticed in \cite{chendalang2015}. Indeed, for any $\alpha > 0$, it turns out that $E[u^2(t, \alpha t)]$ grows exponentially as $t \to +\infty$. This is because for fixed $t>0$, the fundamental solution $\Gamma_\nu(t, x) = \null_\delta G_a(\nu\, t,x)$ has a polynomial rate of decay as $x \to \infty$ (see Lemma \ref{rd01_21l1}), whereas for the heat equation, this rate of decay is exponential, and for the wave equation, the fundamental solution actually vanishes for $\vert x \vert > \nu\, t$.

The kernel $\cK(t, x)$ that appears in \eqref{rd02_04e1} has the same product form 
\beqn
    B(t; \rho, \nu)\, \Gamma(t, x; \nu)
\eeqn
in the stochastic heat equation as the bounds in the case of the fractional stochastic heat equation: see \eqref{rd08_12e8}, \eqref{rd08_17e3} and \eqref{rd08_17e2}. When evaluated at $(t, \alpha\, t)$ the contribution of the second factor is $\exp(- \alpha t/\nu)$ for the heat equation, and merely $\alpha^{-a-1}\, t^{-a}$  for the fractional stochastic heat equation.

This suggests that in the fractional case, we should consider time-space positions of the form $(t, e^{\alpha\, t})$, with $\alpha \geq 0$. In this case, one expects that for small $\alpha \geq 0$, we should have
\beq\label{rd05_14e1}
    \liminf_{t \to \infty} \frac{1}{t} \log E\left[u^2(t, e^{\alpha\, t})\right] > 0,
\eeq
while for large $\alpha$, we should have
\beq\label{rd05_14e2}
    \limsup_{t \to \infty}  E\left[u^2(t, e^{\alpha\, t})\right] = 0.
\eeq
Results of this kind have been obtained in \cite{chendalang2015}, but the existence of a critical value $\alphaf$ such that \eqref{rd05_14e1} holds for $0 \leq \alpha < \alphaf$ and \eqref{rd05_14e2} holds for $\alpha > \alphaf$ appears to still be an open problem.
\medskip

\noindent{\em Upper and lower bounds on $\cK \star I_0^2$}
\medskip

In the next lemma, we establish upper and lower bounds on $[\cKh_\nu \star I_0^2](t, x)$, where $\cKh_\nu$ is given in \eqref{rd05_05e2}; for the lower bound, we assume that $I_0(t, x)$ comes from \eqref{rd06_19e8}. We also establish similar bounds on $[\cKw_\nu \star I_0^2](t, x)$,  where $\cKw_\nu$ is given in \eqref{rd05_05e4} and $I_0(t, x)$ is as in \eqref{rd06_23e1}. These bounds were used in the proof of Theorem \ref{rd05_14p2} (a) and (b), respectively.

\begin{lemma}
\label{rd06_19l1} 

(a) {\em Stochastic heat equation.}  Let $\cKh_\nu$ be as in \eqref{rd05_05e2} and assume that $I_0(t, x)$ is given by \eqref{rd05_06e5} (with $\mu$ there replaced by $u_0$), for some $u_0 \in \cM_H(\R)$. For $(t, x) \in \R_+ \times \R$,
\beq
   [\cKh_\nu \star I_0^2](t, x) \leq \rho^2 \, \sqrt{\frac{\pi t}{\nu}}\, I_0^2(2t, x) \left(1 + 2 \exp\left(\frac{\rho^4\, t}{4 \nu} \right) \right).
\label{rd06_19e2}
\eeq
Now let $I_0(t, x)$ be as in  \eqref{rd06_19e3}. Then for any $c \in \, ]0, 1[$, if $x \geq \beta \nu t$, then
\begin{align}\nonumber
   [\cKh_\nu \star I_0^2](t, x) &\geq  e^{-2 \beta x + \beta^2 \, \nu\, t} \, \Phi^2\left(- \beta\, \sqrt{\nu\, (1-c)\, t}\,\right) \\ 
      &\qquad\qquad\qquad\times  \left(e^{\rho^4 t/(4\nu)} -e^{c\, \rho^4 t/(4\nu)} \right),
   \label{rd06_21e4}
\end{align}
and if $0 < x < \beta \nu t$, then
\begin{align}\nonumber
   [\cK \star I_0^2](t, x) &\geq 2 e^{2\beta\, x } \, \Phi\left(-\left[ x\, \sqrt{\frac{2}{\nu\, t}} + \beta\, \sqrt{2\, \nu\, t} \right] \right)\, \Phi^2\left(- \beta\, \sqrt{\nu\, (1-c)\, t}\,\right)  \\
      &\qquad\qquad\qquad\times \left(e^{\rho^4 t/(4\nu)} -e^{c\, \rho^4 t/(4\nu)} \right).
\label{rd06_21e5}
\end{align}
In the above, $\Phi$ denotes the standard Normal probability distribution function.

(b) {\em Stochastic wave equation.} Let $\cKw_\nu$ be as in \eqref{rd05_05e4} and let $I_0(t, x)$ be as in \eqref{rd06_23e1}. Fix $\gamma \in \, ]0, \beta[$ and define $C_{\beta, \gamma}$ as in \eqref{rd06_20e1}. For $(t, x) \in \R_+ \times \R$,
\beq\label{rd06_20e2}
   [\cKw_\nu \star I_0^2](t, x) \leq C_{\beta, \gamma}\, \frac{\rho^2\, t}{2(\sigma_\gamma - 2 \gamma)}\, e^{\sigma_\gamma\,  \nu\, t -  2 \gamma\, \vert x \vert },\qquad\text{where } \sigma_\gamma = \sqrt{4 \gamma^2 + \tfrac{\rho^2}{2\nu^3}}\, ;
\eeq
for any $a, b \in \, ]0, 1[$, if $\vert x \vert < \nu\, t$, then
\begin{align}\nonumber
    [\cKw_\nu \star I_0^2](t, x) &\geq   \frac{\rho^2}{4 \nu}\, I_0\left( \rho\, b \nu \, t \,\sqrt{\tfrac{1 - a^2}{2 \nu^3}}\right) \\  \nonumber
    &\qquad\times\big[ (8 \beta^4 \nu^3)^{-1}\, 
       \sinh(2 \beta \nu (1-b)\, t)  \\
     & \qquad \qquad\qquad - \frac{(1-b)^2}{\nu \beta^2}\, t^2\big].
  \label{rd06_21e2}
 \end{align}       
 and if $\vert x \vert \geq \nu\, t$, then
 \begin{align}\nonumber
       [\cKw_\nu \star I_0^2](t, x) &\geq  \frac{\rho^2}{8 \nu^3 \beta^3}\, e^{-2 \beta\, \vert x \vert \ + 2 a \nu b \beta\, t}  I_0\left(\rho\, b \nu\, t\, \sqrt{\tfrac{1 - a^2}{2 \nu^3}}\,  \right)  \\
       &\qquad \times \left((4 \beta \nu)^{-1} \sinh(2 (1-b) \beta \nu\, t) - \tfrac{1}{2}\, (1-b)\, t\right),
 \label{rd06_23e2}
 \end{align}
 where $I_0(x)$ denotes the modified Bessel function of the first kind of order $0$ $($see \eqref{rd08_26e5}$)$.
\end{lemma}

\begin{proof}
(a) In this part of the proof, we write $\cK(t, x)$ instead of $\cKh_\nu(t, x)$, and we first prove \eqref{rd06_19e2}.  By \eqref{rd05_06e5},
\begin{align*}
   I_0^2(t, x) = \int_\R u_0(dz_1) \int_\R u_0(dz_2)\, \Gamma_\nu(t, x - z_1)\, \Gamma_\nu(t, x - z_2).
\end{align*}
Using Lemma \ref{rough-heat-l1} (i), we see that
\begin{align}
   I_0^2(t, x) = \int_\R u_0(dz_1) \int_\R u_0(dz_2)\,  \Gamma_\nu\left(\tfrac{t}{2}, y - \tfrac{z_1 + z_2}{2}\right) \Gamma_\nu(2t, z_1 - z_2).
\label{rd06_18e1}
\end{align}

Notice that
\begin{align*}
[\cK \star I_0^2](t, x) = \int_0^t ds \int_\R dy\, I_0^2(s, y)\, \cK(t-s, x-y).
\end{align*}
Using \eqref{rd05_05e2} and \eqref{rd06_18e1}, we see that
\begin{align*}
[\cK \star I_0^2](t, x) &= \int_0^t ds \, B_\nu(t-s; \rho) \int_\R u_0(dz_1) \int_\R u_0(dz_2)\, \Gamma_\nu(2s, z_1 - z_2) \\
     &\qquad\qquad \times \int_\R dy\, \Gamma_{\nu/2}(t-s, x-y) \, \Gamma_\nu\left(\tfrac{s}{2}, y - \tfrac{z_1 + z_2}{2}\right).
\end{align*}
Since $\Gamma_{\nu/2}(t-s, x-y) = \Gamma_{\nu}((t-s)/2, x-y)$, we can use the semigroup property \eqref{semig-heat} to see that
\begin{align*}
[\cK \star I_0^2](t, x) &= \int_0^t ds \, B_\nu(t-s; \rho) \int_\R u_0(dz_1) \int_\R u_0(dz_2) \\
     &\qquad\qquad \times \Gamma_\nu(2s, z_1 - z_2) \, \Gamma_\nu\left(\tfrac{t}{2}, x - \tfrac{z_1 + z_2}{2}\right).
\end{align*}
We use Lemma \ref{rough-heat-l1} (ii) (with $a=4$ there) to bound the product of the two $\Gamma_\nu(\cdots)$ and obtain
\begin{align}\nonumber
[\cK \star I_0^2](t, x) &\leq \int_0^t ds \, B_\nu(t-s; \rho) \int_\R u_0(dz_1) \int_\R u_0(dz_2)\, \frac{2 t}{\sqrt{st}} \\ \nonumber
     &\qquad\qquad \times \Gamma_\nu(2t, x - z_1)\, \Gamma_\nu\left(2t, x - z_2\right) \\
     &=  2 \sqrt{t}\, I_0^2(2 t, x) \int_0^t ds \, \frac{B_\nu(t-s; \rho)}{\sqrt{s}}.
\label{rd06_19e1}
\end{align}
Recalling the expression for $B_\nu(t-s; \rho)$ as the sum of two terms given just after \eqref{rd05_05e2}, we integrate the first term using the formula \eqref{rd08_12e9} for the Beta integral, and we integrate the second term using the formula
\beqn
   \int_0^t \frac{e^{a(t-s)}}{\sqrt{s}}\, ds = \sqrt{\frac{\pi}{a}}\, e^{a\, t}\, \text{erf}\left(\sqrt{a\, t} \right),
\eeqn
which is obtained via the change of variables $s = v^2/a$, and in which $\text{erf}(x) = 2 \Phi(\sqrt{2}\, x) - 1 \leq 1$. This gives
\beqn
  \int_0^t ds \, \frac{B_\nu(t-s; \rho)}{\sqrt{s}} \leq \rho^2\left(\sqrt{\frac{\pi}{4 \nu}} + \sqrt{\frac{\pi}{\nu}}\, \exp\left(\frac{\rho^4\, t}{4 \nu} \right) \right).
\eeqn
Combining this with \eqref{rd06_19e1}, we obtain \eqref{rd06_19e2}.

We now prove \eqref{rd06_21e4} and \eqref{rd06_21e5}. Let $H_{a, b}(x)$ be the function defined in \eqref{rd06_19e9}. 
By \eqref{rd06_19e8} and the first inequality in \eqref{rd06_19e11} below, 
\begin{align*}
   I_0^2(t, x) &\geq e^{\beta^2 \nu\, t}\, \Phi^2\left(- \beta\, \sqrt{\nu\, t} \right) H_{\nu\, t,\, 2 \beta}(x) \\
   & = e^{- \beta^2 \nu\, t}\,  \Phi^2\left(- \beta\, \sqrt{\nu\, t} \right) (e^{-2 \beta\, \vert * \vert} * \Gamma_\nu(t, *))(x),
\end{align*}
where we have also used \eqref{rd06_19e8}, with $2 \beta$ instead of $\beta$. Observe from \eqref{rd05_05e2} that 
\beqn
   \cK(t, x) \geq \frac{\rho^4}{2 \nu} \exp\left(\rho^4 t/(4\nu)\right) \Gamma_{\nu/2}(t, x).
\eeqn 
Therefore,
\begin{align*}
    [\cK \star I_0^2](t, x) &\geq \int_0^t ds\, e^{- \beta^2 \nu (t-s)} \, \Phi^2\left(- \beta\, \sqrt{\nu\, (t-s)}\, \right) \frac{\rho^4}{2 \nu} \\
       &\qquad \times \exp\left(\frac{\rho^4\, s}{4 \nu} \right) \int_\R dy\, \left(e^{-2 \beta\, \vert * \vert} * \Gamma_\nu(t-s, *)\right)(x-y)\, \Gamma_{\nu/2}(s, y) \\
       &= \int_0^t ds\, e^{- \beta^2 \nu (t-s)}\,  \Phi^2\left(- \beta\, \sqrt{\nu\, (t-s)} \, \right) \frac{\rho^4}{2 \nu} \\
        &\qquad \times \exp\left(\frac{\rho^4\, s}{4 \nu} \right) \left(e^{-2 \beta\, \vert * \vert} * \Gamma_\nu\left(t-\tfrac{s}{2}, *\right)\right)(x),
\end{align*}
where we have used the semigroup property \eqref{semig-heat}. By \eqref{rd06_19e3} and Lemma \ref{rd06_19l2} (b) below, for $s \in [0, t]$,
\begin{align*}
   \left(e^{-2 \beta\, \vert * \vert} * \Gamma_\nu\left(t-\tfrac{s}{2}, *\right)\right)(x) &= e^{2\beta^2 \nu \, (t - s/2)} H_{\nu(t-s/2),\, 2\beta}(x) \\
   &\geq e^{2\beta^2  \nu\,(t - s/2)} H_{\nu\, t/2,\, 2\beta}(x).
\end{align*}
It follows that for $c \in [0, 1[$,
\begin{align}\nonumber 
[\cK \star I_0^2](t, x)&\geq H_{\nu\, t/2,\, 2\beta}(x) \, e^{\beta^2 \nu\, t}\, \frac{\rho^4}{2 \nu}  \int_0^t ds\, \Phi^2\left(- \beta\, \sqrt{\nu\, (t-s)}\, \right) e^{\rho^4 s/(4\nu)}\\ \nonumber
&\geq \frac{\rho^4}{2 \nu}  \,  H_{\nu\, t/2,\, 2\beta}(x) \, e^{\beta^2 \nu\, t}\, \Phi^2\left(- \beta\, \sqrt{\nu\, (1-c)\, t}\,\right)  \int_{ct}^t ds\, e^{\rho^4 s/(4\nu)} \\ \nonumber 
&= 2\, H_{\nu\, t/2,\, 2\beta}(x) \, e^{\beta^2 \nu\, t}\, \Phi^2\left(- \beta\, \sqrt{\nu\, (1-c)\, t}\,\right) \\
   &\qquad\qquad\times  \left(e^{\rho^4 t/(4\nu)} -e^{c\, \rho^4 t/(4\nu)} \right).
\label{rd06_19e12a}
\end{align}

Retaining only the second term in \eqref{rd06_19e9}, we see that 
\beqn
   H_{\nu\, t/2,\, 2\beta}(x) \geq e^{-2\beta\, x} \, \Phi\left(\left[x\sqrt{\tfrac{2}{\nu\, t}} - \beta\sqrt{2\, \nu\, t} \right] \right).
\eeqn
If $x \geq \beta \nu\, t$, then $x\, \sqrt{\frac{2}{\nu\, t}} - \beta\sqrt{2\, \nu\, t} \geq 0$, so
\beqn
   H_{\nu\, t/2,\, 2\beta}(x) \geq \half \, e^{-2\beta\, x}.
\eeqn
Together with \eqref{rd06_19e12a}, this establishes \eqref{rd06_21e4}.

If $0 < x < \beta \nu\, t$, then retaining only the first term in \eqref{rd06_19e9}, we see that
\beqn
   H_{\nu\, t/2,\, 2\beta}(x) \geq  e^{2\beta\, x } \, \Phi\left(-\left[ x\, \sqrt{\frac{2}{\nu\, t}} + \beta\, \sqrt{2\, \nu\, t} \right] \right).
\eeqn
Together with \eqref{rd06_19e12a}, this establishes \eqref{rd06_21e5} and completes the proof of part (a).

(b) In this part of the proof, we write $\cK(t, x)$ instead of $\cKw_\nu(t, x)$, defined in \eqref{rd05_05e4}. We first prove \eqref{rd06_20e2}. Fix $\gamma \in \, ]0, \beta[$. 
Looking at the formula \eqref{rd05_05e4}, and because $I_0(z) \leq \cosh(z) \leq e^{\vert z \vert}$, for all $z \in \R$ (see \cite[(10.32.1)]{olver}), we use \eqref{rd06_20e1} to see that
\begin{align}\nonumber
   [\cK \star I_0^2](t, x) &\leq C_{\beta, \gamma}\, \frac{\rho^2}{4 \nu} \int_0^t ds\, e^{2\gamma \nu(t-s)} \int_{- \nu\, t}^{\nu\, t} dy \\
      &\qquad\qquad\times \exp\left(- 2\gamma\, \vert x - y \vert + \sqrt{\frac{\rho^2(\nu^2 \, s^2 - y^2)}{2 \nu^3}} \right).
 \label{rd06_24e1}
\end{align}
Use that fact that $\vert x-y \vert \geq \vert x \vert - \vert y \vert$ to obtain
\begin{align*}
   [\cK \star I_0^2](t, x) &\leq  C_{\beta, \gamma}\, \frac{\rho^2}{2 \nu}\, e^{- 2\gamma\, \vert x \vert} \int_0^t ds\, e^{2\gamma \nu(t-s)} \\
    &\qquad \quad \times \int_{0}^{\nu\, t} dy\, \exp\left(2\gamma\, \vert y \vert + \sqrt{\frac{\rho^2(\nu^2 \, s^2 - y^2)}{2 \nu^3}} \right).
\end{align*}
The function $\psi(y) := 2 \gamma\, y + \left(\rho^2 (\nu^2\, s^2 - y^2) / (2 \nu^3) \right)^{1/2}$ achieves its maximum value $\sigma_\gamma\, \nu\, s$ at $y = 2\sigma_\gamma^{-1}\gamma \nu\, s \in [0, \nu\, s]$, where $\sigma_\gamma$ is defined in \eqref{rd06_20e2}, therefore,
\begin{align*}
   [\cK \star I_0^2](t, x) &\leq C_{\beta, \gamma}\, \frac{\rho^2 \, t}{2}\, e^{- 2\gamma\, \vert x \vert} \int_0^t ds\, e^{2 \gamma \nu\, (t-s) + \sigma_\gamma \nu\, s} \\
   &\leq C_{\beta, \gamma}\, \frac{\rho^2\, t}{2 \nu (\sigma_\gamma - 2 \gamma)}\, e^{- 2 \gamma\, \vert x \vert + \sigma_\gamma\, \nu\, t}.
\end{align*}
This proves \eqref{rd06_20e2}.

We now prove \eqref{rd06_21e2}.  Because $ x \mapsto [\cK \star I_0^2](t, x)$ is an even function (since $\cK(t, *)$ and $I_0(t, *)$ are even), we suppose that $x \in [0,  \nu\, t[$. By \eqref{rd05_05e4}, 
\beqn
   [\cK \star I_0^2](t, x) = \frac{\rho^2}{4\nu} \int_0^t ds \int_{- \nu s}^{\nu s} I_0^2(t-s, x-y)\, I_0\left(\sqrt{\frac{\rho^2\left((\nu\, s)^2 - y^2 \right)}{2 \nu^3}}\, \right).
\eeqn
Let $a, b \in \, ]0, 1[$. We restrict the domain of integration to see that
\begin{align}\nonumber 
   [\cK \star I_0^2](t, x) &\geq  \frac{\rho^2}{4 \nu} \int_{bt}^t ds \int_{- a \nu s}^{a \nu s} dy\, I_0^2(t-s, x-y) \, I_0\left(\sqrt{\frac{\rho^2(1 - a^2)}{2 \nu^3}}\, b \nu \, s\right) \\ \nonumber
      &\geq \frac{\rho^2}{4 \nu} \, I_0\left(\sqrt{\frac{\rho^2(1 - a^2)}{2 \nu^3}}\, b \nu \, t \right) \\
       &\qquad \qquad \times \int_{bt}^t ds  \int_{- a b \nu t}^{a b \nu t} dy\, I_0^2(t-s, x-y).
 \label{rd06_23e3}
\end {align}

We assume that $a$ and $b$ are close enough to $1$ so that
\beqn
   x + \nu (1-b)\, t < ab \nu\, t.
\eeqn
In this case, we claim that $bt \leq s \leq t$ and $\vert y - x \vert \leq \nu(t-s)$ imply that $\vert y \vert \leq ab \nu\, t$. Indeed,
\beqn
   \vert y \vert \leq \vert y-x \vert + \vert x \vert \leq \nu(t-s) + x \leq \nu (1-b) \, t + x < ab \nu\, t
\eeqn
as claimed. In particular,
\beqn
    \int_{bt}^t ds  \int_{- a b \nu t}^{a b \nu t} dy\, I_0^2(t-s, x-y) \geq  \int_{bt}^t ds  \int_{\vert y-x \vert \leq \nu (t-s)} dy\, I_0^2(t-s, x-y).
\eeqn
By \eqref{rd06_23e1}, in this region of integration,
\beqn
   I_0(t-s, x-y) = - (\beta \nu)^{-1} \left( e^{-\beta \nu\, (t-s)} \cosh(\beta(x-y)) - 1 \right) .
\eeqn
Using the inequality $(a + b)^2 \geq \tfrac{a^2}{2} - b^2$ and $\cosh^2(x) =\half(\cosh(2 x) + 1) \geq \half\cosh(2 x)$, we see that
\beqn
   I_0^2(t-s, x-y) \geq (4 \beta^2 \nu^2)^{-1} e^{-2 \beta \nu\, (t-s)} \, \cosh(2\beta(x-y)) - (\beta \nu)^{-2}.
\eeqn
Therefore,
\begin{align*}
  & \int_{bt}^t ds  \int_{\vert y-x \vert < \nu(t-s)} dy\, I_0^2(t-s, x-y) \\
   &\qquad \geq (8 \beta^4 \nu^3)^{-1} \, (1 - e^{-2 b \beta \nu\, t})\, \sinh(2  \beta \nu (1-b) \, t) \\
     & \qquad \qquad\quad 
        - \frac{(1-b)^2}{\nu \beta^2}\, t^2 .
\end{align*}
Together with \eqref{rd06_23e3}, this implies \eqref{rd06_21e2}. 

   We now establish \eqref{rd06_23e2}. Suppose that $\vert x \vert \geq \nu\, t$. Then for all $(s,y)$ such that $s \in [0, t]$ and $\vert y \vert \leq \nu\, s$, we have $\vert x-y \vert \geq \nu\, (t-s)$. By \eqref{rd06_23e1} and \eqref{rd06_23e3},  we see that
\begin{align*}
    [\cK \star I_0^2](t, x) &\geq \frac{\rho^2}{4 \nu^3 \beta^2} \, I_0\left(\sqrt{\frac{\rho^2(1 - a^2)}{2 \nu^3}}\, b \nu \, t\right)   \int_{b t}^t ds\, \sinh^2(\beta \nu\, (t-s))\\
       &\qquad\qquad\qquad \times \int_{- a \nu s}^{a \nu s} dy\, e^{-2\beta\, \vert x-y\vert}.
\end{align*}
Use the inequality $\vert x - y \vert \leq \vert x \vert + \vert y \vert$ to see that
\begin{align*}
    [\cK \star I_0^2](t, x) &\geq \frac{\rho^2 e^{-2\beta\, \vert x \vert + 2 a b \beta \nu\, t}}{8 \nu^3 \beta^3}\, I_0\left(\sqrt{\frac{\rho^2(1 - a^2)}{2 \nu^3}}\, b \nu \, t\right) \\
      &\qquad\qquad\qquad \times \left(\frac{\sinh(2 (1-b) \beta \nu\, t)}{4 \beta \nu} - \half (1-b)\, t \right) .
\end{align*}
This establishes \eqref{rd06_23e2} and completes the proof of Lemma \ref{rd06_19l1}.
\end{proof}

The next lemma was used in the proof of Lemma \ref{rd06_19l1} (a) .

\begin{lemma}\label{rd06_19l2}  Let $H_{a, b}(x)$ be the function defined in \eqref{rd06_19e9}. Then:

(a) for $a, b \in \R_+^*$ and $x \in \R$,
\beq\label{rd06_19e11}
     \Phi(b\, \sqrt{a})\, H_{a, 2b}^{1/2}(x) \leq H_{a, b}(x) \leq e^{-b\, \vert x \vert};
\eeq
(b) for $b>0$ and $x\in \R$ fixed, the map $a \mapsto H_{a, b}(x)$ is decreasing.
\end{lemma}

\begin{proof}
(a) For the second inequality in \eqref{rd06_19e11}, since $x \mapsto H_{a, b}(x)$ is an even function, we assume that $x \geq 0$. We must show that
\beq\label{rd06_19e10}
    e^{b x}\, \Phi\left(\frac{-a b - x}{\sqrt{a}} \right) + e^{- b x}\, \Phi\left(\frac{-a b + x}{\sqrt{a}}  \right) \leq e^{- bx} ,
\eeq
or, equivalently because  $\Phi\left(\frac{-a b + x}{\sqrt{a}} \right)  = 1 - \Phi\left(\frac{ a b - x}{\sqrt{a}} \right) $, that
\beqn
    \kappa(x) := e^{- b x}\, \Phi\left(\frac{ a b - x}{\sqrt{a}} \right)  - e^{b x}\, \Phi\left(\frac{-a b - x}{\sqrt{a}} \right) \geq 0.
\eeqn
Notice that $\kappa(0) = \Phi(b\, \sqrt{a}) -  \Phi(- b\, \sqrt{a}) >  0$ because $b > 0$, and by l'H\^opital's rule, $\lim_{x \to + \infty} \kappa(x) = 0$. By direct calculation, we find that
\beqn
   \kappa'(x) = - b\left( e^{- b x}\, \Phi\left(\frac{ a b - x}{\sqrt{a}} \right)  + e^{b x}\, \Phi\left(\frac{-a b - x}{\sqrt{a}} \right)\right) \leq 0,
\eeqn
so we conclude that $\kappa(x) \geq 0$, which establishes the second equality in \eqref{rd06_19e11}.

For the first inequality in \eqref{rd06_19e11}, observe that
\begin{align*}
   H_{a, b}^2(x) &= \left( e^{b x}\, \Phi\left(\frac{-a b - x}{\sqrt{a}} \right) + e^{- b x}\, \Phi\left(\frac{-a b + x}{\sqrt{a}}  \right)\right)^2 \\
      &\geq e^{- 2 b x}\, \Phi^2\left(\frac{-a b + x}{\sqrt{a}}  \right) \geq e^{- 2 b x}\, \Phi^2\left(b\, \sqrt{a} \right).
\end{align*}
Using the just established second inequality in \eqref{rd06_19e11}, we see that $e^{- 2 b x} \geq H_{a, 2b}(x)$. This completes the proof of \eqref{rd06_19e11}.

   (b) Notice that
\beqn
   \frac{\partial }{\partial a} e^{\pm b x}\Phi\left(\frac{-a b \mp x}{\sqrt{a}}\right)
       = \frac{-a b \pm x}{2 a^{3/2} \sqrt{2\pi}}\, \exp\left(-\frac{a^2 b^2+x^2}{2a}\right).
\eeqn
Adding the two terms that come when writing $\frac{\partial }{\partial a} H_{a, b}(x)$ gives
\beqn
   \frac{\partial H_{a,b}(x)}{\partial a}= \frac{-b}{\sqrt{2\pi a}} \exp\left(-\frac{a^2 b^2+x^2}{2a}\right) < 0.
\eeqn
This proves (b).
\end{proof}

\section{Notes on Chapter \ref{ch2'}}
\label{notes-ch5}

The contents of Section \ref{ch2'-s1} is inspired by the unpublished preprint \cite{gp}. We note however that the notion of weak solution considered in that reference does not coincide with our Definition \ref{ch2'-s3.4.1-d1}. A theorem on uniqueness in law comparable to our Theorem \ref{ch2'-s3.4.2-t2} (along with Proposition \ref{ch2'-s3.4.2-p2}) can be found in \cite[Theorem I.0.2,  p.~251]{lr}.

 Various questions on absolute continuity of the laws of solutions to linear stochastic heat equations on $\re$ and on $[0,L]$ have been discussed in \cite{m-t2002}. Section \ref{ch2'-s3.4.3} is devoted to an illustrative example. Comparing the proof of Theorem \ref{ch5-t-5.1.11} with that of  \cite{m-t2002}, we have simplified the arguments and taken into account a problem mentioned and solved in the unpublished manuscript [C. Mueller and R. Tribe,  A Correction to ``Hitting Properties of the Random String" (EJP 7 (2002), Paper 10, pages 1-29), Nov. 18, 2004].

The results of Section \ref{ch2'-s3.4.4} on the germ-field Markov property for a stochastic heat equation appear in \cite{nualartpardoux1}. Using the non-anticipating version of the Girsanov theorem in \cite{kusuoka-1982}, the germ-field Markov property has also been proved for elliptic SPDEs in \cite{donati-martin-1992}. For an expository account see \cite[Section 4.2]{nualart-1995}.

The main theorem of Section \ref{ch2'-section5.2} gives a precise upper bound on the rate of growth of the $L^p$-moment, as a function of $p$ and time, of solutions to a large class of SPDEs. The growth order $p^3$ for the Lyapounov exponent corresponding to the stochastic heat equation can be found for instance in \cite{foondun-khoshnevisan-2009} and \cite{khosh}, while the order $p^{3/2}$ for the stochastic wave equation appears in \cite{conus-joseph-khoshnevisan-2013}. For the fractional stochastic heat equation, the order $p^{(2a-1)/(a-1)}$ was obtained in \cite {chendalang2015}. Recently, sharp lower bounds which match upper bounds for all moments and for a large class of SPDEs have been established in \cite{hu-wang-2024}.

Among the most important applications of these types of estimates are properties of moment Lyapounov exponents and the study of intermitency phenomena. Intermitency for solutions of SPDEs is a well-developed research area with origins in
\cite{carmona-molchanov-1994} and \cite{bertini-cancrini-1995} for the parabolic Anderson model (see Section \ref{ch1-1.2}). We also refer to \cite{foondun-khoshnevisan-2009}, \cite{chendalang2015}, \cite{chen-dalang-2015}, \cite{chendalang2015-2}.

Exponential $L^p$-bounds are also applied to the study of global solutions to SPDEs with super-linear coefficients. This question has been addressed in \cite{D-K-Z016} for the stochastic heat equation and in \cite{millet-SanzSole-2021} for the stochastic wave equation.

Section \ref{ch5-added-0} aims at introducing comparison theorems for parabolic SPDEs, which have been an important tool that was already present in the early years in the development of the theory (see e.g. \cite{yamada-1973}). Theorem \ref{ch5-added-0-t1} is a more general version of Theorem 2.1 in \cite{donati-martin-pardoux-1993} and Theorem \ref{rd08_13t1} appears in \cite{mueller-1991}, \cite{shiga-1994} and \cite{chen-kim-2016} with various proofs. More sophisticated versions of sample path comparison theorems have been developed by Le Chen and coauthors (see \cite{chen-kim-2016}, for instance). An important related problem concerns {\em moment comparisons}, in which one seeks conditions under which the moments of the solutions are a monotone function of the coefficient $\sigma$: see \cite{cox-fleisch-greven-1996}, \cite{joseph-khosh-mueller-2017}, \cite{f-j-s-2018} and \cite{chen-kim-arxiv}.

Section \ref{ch2'-s7} contains some elements of probabilistic potential theory for random fields. The discussion is restricted to questions of polarity and mainly to systems of linear SPDEs. In the proofs of Propositions \ref{ch2'-s7-p1} and  \ref{ch2'-s7-p2},  we present (and slightly improve) basic ideas that are also used in the proofs of criteria for upper and lower bounds on hitting probabilities for Gaussian and non-Gaussian random fields (see e.g.~\cite{dkn07}, \cite{dss10}, \cite{hss2020}). The development of probabilistic potential theory for SPDEs initiated with the paper
\cite{dn2004} that extended to systems of one-dimensional nonlinear stochastic wave equations the results in  \cite{khosh-shi-99} relative to the Wiener sheet. This was followed by the more applicable approach of \cite{dkn07} and \cite{dkn09} for systems of linear and non-linear stochastic heat equations in spatial dimension one.   The method developed in {\cite{dn2004} and \cite{dkn09} for obtaining lower bounds on hitting probabilities for nonlinear SPDEs relies on Malliavin calculus. This method has been successfully used in several other examples (\cite{dkn2013}, \cite{dss15}, \cite{dalang-pu-2020}). See also \cite{dalang2018} for an overview.

Section \ref{rdrough} originates with \cite{chen}. We have however modified the notion of solution given in Definition \ref{rdHGammaI} and the set of assumptions used in Theorem \ref{rd08_20t1}. The proof of this theorem is adapted from \cite{chen} (see also \cite{chendalang2015-2}). Proposition \ref{rd01_07p1} unifies various results that appear in \cite{chen}, \cite{chendalang2015-2}, \cite{chen-dalang-2015} and \cite{chendalang2015}. The the material in Subsections \ref{rd01_07ss1} and \ref{rd06_23ss1} and in Section \ref{rd1+1anderson}, are also mostly taken (and adapted) from these references (with the exception of Subsection \ref{rd01_27ss1}). The study of the propagation speed of intermitency discussed in Section \ref{rd06_24ss1} was initiated the paper \cite{ConusKhosh2012} (see also \cite{conus-joseph-khoshnevisan-2014}), and further developed in \cite{chen} and in the papers \cite{chendalang2015-2}--\cite{chendalang2015}.

\appendix



\chapter[Some elements of stochastic processes and stochastic \texorpdfstring{\\}{}analysis]{Some elements of stochastic processes and stochastic analysis}
\label{app1}

\pagestyle{myheadings}
\markboth{R.C.~Dalang and M.~Sanz-Sol\'e}{Elements of stochastic processes}

In this chapter, we collect some fundamental notions and results from the theory of stochastic processes, as well as basic facts on the Itô stochastic integral. The first section is devoted to measurability issues. Section \ref{app1-1} deals with distribution-valued stochastic processes, providing the background for the study of space-time white noise. Section \ref{app1-3} addresses sample path  regularity: we present a version of Kolmogorov's continuity theorem for anisotropic random fields. The last two sections are devoted to properties of the Itô integral that are used in particular in Chapter \ref{ch2}: joint measurability of the stochastic integral when the integrand depends on a parameter, and a stochastic Fubini's theorem, respectively.

\section{Stochastic processes and measurability}
\label{rdsecA.4}

A stochastic process is a mathematical model for random evolution. A continuous-time stochastic process is a collection of random variables indexed by $\re_+$, which represents the set of times. However, modeling complex phenomena may require more general ``time sets." This motivates the following definition.
\begin{def1}
\label{rdsecA.4-1}
Let $\T$ be a set and $(S,\mathcal{S})$ a measure space. A stochastic process indexed by $\T$\index{stochastic process}\index{process!stochastic} and taking values in $(S,\mathcal{S})$ is a family $Z=(Z_t,\, t\in \T)$ of measurable mappings $Z_t$ from a probability space $(\Omega,\mathcal{F}, P)$ into $(S,\mathcal{S})$.
\end{def1}
The set $\T$ is called the set of {\em indices} or the {\em index set}.\index{index set} For Brownian motion, $\T=\re_+$; for the Brownian sheet, $\T=\IR_+^2$; other common index sets are $\T=\IN^k$ and $\T=\Z^k$. The set $\T$ can also be a set of functions: for example, in Definition \ref{ch1-d4}, $\T$ is the set $\mathcal{S}(\rek)$ of Schwartz test functions.

In the case where $\T$ is a subset of $\re_+\times \rek$, we often use the term {\em random field}\index{random field}\index{field!random} instead of stochastic process.
The random field solutions to stochastic partial differential equations considered in this book belong to this class of processes.

The measurable space $(S,\mathcal{S})$ is called the {\em state space}.\index{state space} In many examples, $(S,\mathcal{S})= (\red, \mathcal{B}_{\red})$. However, in the framework of SPDEs, it is also natural to consider $Z_t$ as an evolution in time taking values in a space of functions: a Hilbert space, the space of $\alpha$-H\"older continuous functions ${\ca}^\alpha(\red)$, a fractional Sobolev space, etc.

For every $\omega\in\Omega$, the mapping $\T\ni t\mapsto Z_t(\omega)\in S$ is called a {\em trajectory} or a {\em sample path} of the process $Z=(Z_t,\, t\in \T$). This is a deterministic mapping.

Given two stochastic processes $Z$ and $Y$ as above, defined on the same probability space  $(\Omega,\mathcal{F}, P)$, we say that one is a {\em modification}\index{modification} (or a {\em version})\index{version} of the other if, for any $t\in \T$, we have $P\{ Z_t= Y_t\}=1$. We say that  $Z$ and $Y$ are {\em indistinguishable}\index{indistinguishable} if
\beqn
P\{ Z_t = Y_t,\  {\text{for all}}\  t\in\T\} = 1,
\eeqn
that is, if almost all of their sample paths are identical. Obviously, if $Z$ and $Y$ are  indistinguishable, then $Z$ is a modification of $Y$. This implication is strict; however, in some particular cases and under certain conditions, the converse holds. Indeed, suppose that the space of indices $\T$ and the state space $S$ are topological spaces, and $\T$ is separable; if almost all sample paths of $Z$ and $Y$ are continuous and if $Z$ is a modification of $Y$, then $Z$ and $Y$ are indistinguishable.

Let $S^\T$ be the set of all functions from $\T$ into $S$. Observe that all sample paths of a stochastic process $Z$ indexed by $\T$ are  elements of $S^{\T}$. For each  $t\in \T$, we can define the {\em coordinate map}\index{coordinate map} $\pi_t$
 from $S^\T$ into $S$  by $\pi_t(f) = f(t)$.
We denote by $\mathcal{S}^\T$ the smallest $\sigma$-field on $S^\T$ for which all the coordinate maps $\pi_t$ are measurable. It coincides with the $\sigma$-field generated by the cylindrical sets\index{cylindrical!set} (also called measurable rectangles) $\prod_{t\in{\T}}A_t$, where $A_t\in\mathcal{S}$ for all $t$ and, except for a finite set $\{t_1, \ldots,t_n\}$, $A_t=S$. We will refer to $\mathcal{S}^\T$ as the {\em product} $\sigma$-field.\index{product $\sigma$-field}\index{sigma@$\sigma$-field!product} A stochastic process $Z=(Z_t,\, t\in \T)$ defines a measurable mapping
\beqn
\mathcal{Z} : (\Omega,\tf) \longrightarrow (S^\T, \mathcal{S}^\T),
\eeqn
by setting $\mathcal{Z}(\omega)(t) = Z_t(\omega)$.
\medskip

\noindent{\em The law of a stochastic process}
\medskip

The {\em law of the stochastic process}\index{law!of a stochastic process}\index{stochastic process!law of a} $Z$, denoted by $P_Z$, is the probability measure on $\mathcal{S}^\T$ which is the
image of $P$ by $\mathcal{Z}$, that is, $P_Z= P\circ {\mathcal{Z}}^{-1}$. On a measurable rectangle as above, we have
\begin{align}
\label{A.1-finite-dim}
P_Z\left(\prod_{t\in{\T}}A_t\right) &= P\left\{Z_{t_1}\in A_{t_1}, \ldots, Z_{t_n}\in A_{t_n}\right\}\notag\\
&=:\mu_{(t_1,\dots,t_n)}\{A_{t_1}\times \cdots\times A_{t_n}\}.
\end{align}

The collection of probability measures $\mu_{(t_1.\dots,t_n)}$ on the right-hand side of \eqref{A.1-finite-dim}, for all $n\in\N^*$, and all $(t_1, \ldots,t_n)\in\T$, is called the family of {\em finite-dimensional distributions}\index{finite-dimensional!distributions}\index{distributions!finite-dimensional} of the process.

Consider the particular case where $S$ is a complete separable metric space and $\mathcal{S}$ is the $\sigma$-field of Borel subsets of $S$. Then the law of $Z$ is determined by its finite-dimensional distributions \eqref{A.1-finite-dim}. This follows from a version of  Kolmogorov's theorem on extension of measures. We refer to \cite[Section III-3, Th\'eor\`eme, p. 78 and Corollaire p. 79]{neveu-fr}, or \cite[(3.2) Theorem, p. 34]{ry} (without a proof). 
\medskip

\noindent{\em Canonical representation of a process}
\medskip

The canonical process\index{canonical!process}\index{process!canonical} associated to $Z$ is the stochastic process $(\pi_t,\, t\in\T)$ defined on the probability space $(S^\T, \mathcal{S}^\T, P_Z)$.

In this book, we often consider stochastic processes $Z$ indexed by a topological space $\mathbb{T}$ and that possess a continuous version, that is, with trajectories in $\mathcal{C}:=\mathcal{C}(\T; S)$, the space of continuous functions from $\T$ into $S$.
 In this case, if $A \in \mathcal{S}^{\T}$ and $A \supset \mathcal{C}$, then $P_Z(A) = 1$. Although $\mathcal{C}\notin \mathcal{S}^\T$ in general, it is possible to construct a suitable representation of the {\em canonical process} carrying the regularity properties of the process $Z$.
 Indeed, for $t \in \T$, define $\tilde\pi_t: \mathcal{C} \to S$ by $\tilde\pi_t(f) = f(t)$, so that $(\tilde\pi_t,\, t \in \T)$ is the set of coordinate functions from $\mathcal{C}$ into $S$. Let $\mathcal{S}_{\mathcal{C}}^\T$ be the $\sigma$-field on $\mathcal{C}$ for which all the coordinate maps $\tilde\pi_t$ are measurable. One can prove that, for any $\tilde A\in \mathcal{S}_{\mathcal{C}}^\T$, there is $A\in \mathcal{S}^\T$, such that $\tilde A = A\cap \mathcal{C}$ (\cite[p. 35]{ry}). Then one can define a probability measure $\tilde P_Z$ on $\mathcal{S}_{\mathcal{C}}^\T$ by
 \beqn
 \tilde P_Z(\tilde A) = P_Z(A).
\eeqn
The process $(\tilde \pi_t,\, t\in\T)$ defined on the probability space $(\mathcal{C}, \mathcal{S}_{\mathcal{C}}^\T, \tilde P)$ is called a {\em canonical representation}\index{canonical!representation}\index{representation!canonical} of the process $Z$ on $(\mathcal{C}, \mathcal{S}_{\mathcal{C}}^\T, \tilde P_Z)$.

When $Z$ is a $d$-dimensional Brownian motion, $\mathcal{C}$ is the space of continuous functions $\cC(\re_+; \red)$, and the probability measure $\tilde P_Z$ is the Wiener measure.
The probability space $(\mathcal{C}, \mathcal{S}_{\mathcal{C}}^\T, \tilde P_Z)$ is called the Wiener space. We refer to \cite[(3.3) Proposition, p. 35]{ry} for more details (see also \cite[Theorem 2.6]{bass97}).
\medskip

\noindent{\em Joint measurability and separability}
\medskip

Let $\mathcal{T}$ be a $\sigma$-field of subsets of $\T$. The stochastic process $Z$ of Definition \ref{rdsecA.4-1} is {\em jointly measurable}\index{jointly!measurable}\index{measurable!jointly} if the map
$\T\times \Omega\ni(t,\omega)\mapsto Z_t(\omega)\in S$ is measurable with respect to the product $\sigma$-field $\mathcal{T}\times \mathcal{F}$. Using notations of measure theory, this property is expressed in the form
\beqn
Z:(\T\times \Omega, \mathcal{T}\times \mathcal{F}) \rightarrow (S, \mathcal{S}).
\eeqn
If $\T$ is a metric space, a natural and frequent choice is $\mathcal{T}=\mathcal{B}_{\T}$, the Borel $\sigma$-field of $\T$ (generated by the open sets of $\T$).
\smallskip

Assume that $\T$ is a separable metric space; a stochastic process $(Z_t,\, t\in \T)$ is {\em separable}\index{separable!stochastic process}\index{stochastic process!separable} if there is a countable set $\tilde\T\subset \T$ and a $P$-null set $N\in\tf$ such that, for all $\omega\notin N$,  
$\{(t,Z(t,\omega)),\, t\in\tilde\T\}$ is dense in $\{(t,Z(t,\omega)),\, t\in\T\}$ (see \cite[p. 154]{dm1} and also \cite[Définition III-4-2, p. 82]{neveu-fr}).
Every stochastic process $Z$ indexed by $\T$ has a separable version\index{separable!version} (see \cite[Theorem 1, p. 162]{cohn}).
\medskip

\noindent{\em Filtrations and related $\sigma$-fields}
\medskip

In the sequel, we will consider the particular case $(\T, \mathcal{T}) =(\re_+, \mathcal{B}_{\re_+})$ (or $(\T, \mathcal{T}) =([0,T], \mathcal{B}_{[0,T]})$, $T>0$) and $(S, \mathcal{S}) := (\mathcal{E}, \mathcal{B}_{\mathcal{E}})$,
where $\mathcal{E}$ is a metric space and $\mathcal{B}_{\mathcal{E}}$ is the Borel $\sigma$-field of $\mathcal{E}$ (for example,
$(S, \mathcal{S}) = (\red, \mathcal{B}_{\red})$).

 A family $(\tf_t,\, t\in\re_+)$ of sub-$\sigma$-fields of $\tf$ is a {\em filtration}\index{filtration}  if it is increasing, that is, $\tf_s\subset \tf_t$
for any $0\le s<t<\infty$. The {\em natural filtration}\index{filtration!natural}\index{natural filtration} associated with the process $Z$ is defined by
\beqn
\tf_t = \sigma(Z_r,\, 0\le r\le t), \quad t\in\re_+,
\eeqn
where the right-hand side denotes the $\sigma$-field generated by the random variables $Z_r,\ 0\le r\le t$.

A filtration  $(\tf_t,\, t\in\re_+)$ is {\em right-continuous}\index{filtration!right-continuous} if for all $t\in \IR_+$, $\cap_{s>t}\, \tf_s = \tf_t$. It is called {\em complete} if $\tf_0$ contains all $P$-null sets of $\tf$ and therefore, for every $t>0$, $\tf_t$ also contains all $P$-null sets of $\tf$.
To a filtration $(\tf_t^0,\, t \in \re_+)$, we can associate a complete and right-continuous filtration $(\cF_t,\, t \in \re_+)$ by setting $\tf_t = (\cap_{s:\, s > t}\, \tf_s) \vee \mathcal{N}$, where $\mathcal{N}$ is the $\sigma$-field generated by all $P$-null sets.

 The stochastic process $Z$ is {\em adapted}\index{adapted!stochastic process}\index{stochastic process!adapted}\index{process!adapted stochastic} to the filtration  $(\tf_t,\, t\in\re_+)$ if for any $t\in\re_+$, the random variable $Z_t$ is $\tf_t$-measurable, that is, the mapping $Z_t: (\Omega, \tf_t) \rightarrow (\mathcal{E}, \mathcal{B}_{\mathcal{E}})$ is measurable.

A fundamental example of adapted process is the so called $(\cF_t)${\em-standard Brownian motion},\index{standard Brownian motion}\index{Brownian!motion, standard} defined as follows.
Let $(\Omega,\cF,P)$ be a probability space equipped with a complete and right-continuous filtration  $(\cF_t,\, t \in \IR_+)$. An  $(\cF_t)${\em-standard Brownian motion} is a real-valued continuous adapted process $(B_t,\, t\in\re_+)$ such that $B_0= 0$ a.s., the process $(B_t-B_s,\, t\ge s)$ is independent of $\tf_s$, and the increment $B_t-B_s$ is normally distributed with mean zero and variance $t-s$.

In the next definitions, we consider a probability space $(\Omega, \tf, P)$ endowed with a right-continuous complete filtration $(\tf_t,\, t\in\re_+)$. We introduce notions of measurability that are stronger than {\em joint measurability} (defined above).

The stochastic process $Z$ is called {\em progressively measurable}\index{progressively measurable}\index{measurable!progressively} if for each $t\in\re_+$ and $A\in\mathcal{B}_{\mathcal{E}}$, the set
$\{(s,\omega)\in [0,t]\times \Omega: Z_s(\omega)\in A\}$ belongs to $\mathcal{B}_{[0,t]}\times \tf_t$, in other words, for any $t\in\re_+$,
$Z: ([0,t]\times \Omega, \mathcal{B}_{[0,t]}\times \tf_t) \rightarrow (\mathcal{E}, \mathcal{B}_{\mathcal{E}})$.

It is clear that a progressively measurable process is jointly measurable and adapted. Every adapted process $Z$ as above with left- or right-continuous sample paths is progressively measurable (see \cite[Proposition 4.8, p.44]{ry}).

A subset $M$ of $\re_+\times \Omega$ is {\em progressive}\index{progressive!set}\index{set!progressive} if the stochastic process $1_M$ is progressively measurable. The $\sigma$-field consisting of all progressive sets is called the {\em progressive $\sigma$-field}.\index{progressive!$\sigma$-field}\index{sigma@$\sigma$-field!progressive} It is denoted by ${\text{Prog}}$.\label{rdProg}  A process $Z$ is progressively measurable if and only if it is measurable with respect to the progressive $\sigma$-field.

The $\sigma$-field generated by the set of adapted processes $Z$ which are left continuous is called the  {\em predictable $\sigma$-field}.\index{predictable!$\sigma$-field}\index{sigma@$\sigma$-field!predictable} It is denoted by $\mathcal{P}$.\label{rdpred} For such $Z$, by its very definition,
$Z: (\re_+\times \Omega,\ \mathcal{P})\rightarrow  (\mathcal{E}, \mathcal{B}_{\mathcal{E}})$. A process $Z$ is predictable if and only if it is measurable with respect to $\mathcal{P}$.

Consider the $\sigma$-field on $\IR_+ \times \Omega$ generated by the sets of the form $\{0\}\times F_0$ and $]s,t]\times F$, where $F_0 \in \cF_0$ and $F \in \cF_s$ for $s< t$ in $\IR_+$, called {\em predictable rectangles}.
 According to \cite[Chapter IV, Proposition 5.1]{ry}, this $\sigma$-field coincides with $\mathcal{P}$ and also with the $\sigma$-field generated by the set of adapted and continuous processes $X$.

The $\sigma$-field generated by the set of adapted processes $Z$ which are right continuous with left limits (referred to as {\em c\`ad-l\`ag}, an abbreviation of the French phrase ``continu \`a droite avec limites à gauche") is called the  {\em optional $\sigma$-field}.\index{optional!$\sigma$-field}\index{sigma@$\sigma$-field!optional} It is denoted by $\mathcal{O}$.\label{rdOptional} A process $Z$ is optional if and only if it is measurable with respect to $\mathcal{O}$.

From the above definitions, we see that
\beq
\label{app1-rdsecA.4-1}
\mathcal{P} \subset \mathcal{O} \subset {\text{Prog}} \subset ~\mathcal{B}_{\re_+} \times \tf.
\eeq
The inclusions $\mathcal{P} \subset \mathcal{O}$ and $\mathcal{O} \subset {\text{Prog}}$ are in general strict (see \cite[Chapter IV, p. 172]{ry}).
There are $\mathcal{B}_{\re_+} \times \tf$-measurable  adapted processes that are not progressively measurable (see an example in \cite[Chapter 3, p. 62]{chung-williams}).
\medskip

\noindent{\em Stopping times}
\medskip

Next, we recall two notions of random times. Let $(\Omega, \tf, P)$ be a probability space equipped with a filtration $(\tf_t,\, t\in\re_+)$.
A random variable
$\tau: \Omega\rightarrow [0,\infty]$ is a {\em stopping time}\index{stopping time}\index{time!stopping} if for any $t\in\re_+$, the event $\{\tau\le t\}$ belongs to $\tf_t$. The random time $\tau$ is an {\em optional stopping time}\index{stopping time!optional}\index{optional!stopping time} (or just {\em optional time}) if  for any $t\in\re_+$, the event $\{\tau < t\}$ belongs to $\tf_t$. Every stopping time is optional; if the filtration is right-continuous, then the converse is also true. 

A stopping time $\tau$ is said to be {\em predictable}\index{predictable!stopping time}\index{stopping time!predictable} if there exists an increasing sequence $(\tau_n,\, n\ge 1$ of stopping times such that, a.s.,
\begin{description}
\item{(i)}\ $\lim_{n\to\infty}\tau_n = \tau$;
\item{(ii)}\ on the event $\{\tau>0\}$, we have $\tau_n<\tau$.
\end{description}

Given a stopping time $\tau$ with respect to the filtration $(\tf_t,\, t\in\re_+)$, we define
the $\sigma$-field $\tf_\tau$ that consists of all sets $A \in \tf$
satisfying $A\cap \{\tau\le t\}\in\tf_t$, for all $t\in\re_+$. These sets are called {\em events determined prior to}\index{sigma@$\sigma$-field!of events prior to} $\tau$.

Let $Z$ be a jointly measurable stochastic process which is either positive or bounded. The {\em optional projection}\index{optional!projection}\index{projection!optional} of $Z$ is the unique (up to indistinguishability) optional process $Y$ such that
\beqn
E\left(Z_\tau\ 1_{\{\tau<\infty\}}\big|\tf_\tau
\right) = Y_\tau\, 1_{\{\tau<\infty\}}\ a.s.,\ {\text{for any stopping time}}\ \tau.
\eeqn
Existence and uniqueness of such a process $Y$ is proved for instance in \cite[Chapter IV, Theorem 5.6]{ry} (see also  \cite[Theorem 3.6]{chung-williams}).


Random times can be used to define random intervals, which are subsets of $\re_+\times \Omega$. For instance, for $\tau_1$ and $\tau_2$ satisfying $\tau_1\le \tau_2$,
$]\tau_1,\tau_2] = \{(r,\omega): \tau_1(\omega)<r\le \tau_2(\omega)\}$. Similarly, we can define $]\tau_1,\tau_2[$, $[\tau_1,\tau_2]$ and $[\tau_1,\tau_2[$.

The $\sigma$-fields $\cP$ and $\mathcal{O}$ defined above admit a description in terms of random intervals, as follows.
The $\sigma$-field of predictable sets $\cP$ coincides with the $\sigma$-field generated by the random intervals $]\tau,\infty[$, where $\tau$ is a predictable stopping time, while
the optional $\sigma$-field $\mathcal{O}$ coincides with
the $\sigma$-field generated by the stochastic intervals of the form $[\tau,\infty[$, where $\tau$ is a stopping time.
For the proof of these results, we refer to \cite[Sections 2.3 and 3.2]{chung-williams}.
\medskip


\noindent{\em Existence of measurable versions}
\smallskip

The following statement is taken from \cite[Theorem 3, and Remark p.164]{cohn}. It refers to a stochastic process $(Z_t,\, t\in \mathbb{T})$ defined on a (not necessarily complete) probability space $(\Omega, \tf, P)$, where $\mathbb{T}$ is a separable metric space and the $Z_t:\Omega \rightarrow K$ are random variables taking values in a compact metric space $K$ (for example, $K = [0,\infty]$).

Let $M(\Omega,K)$ denote the set of measurable mappings from $(\Omega, \tf)$ to $(K, \mathcal{B}_K)$, in which   mappings that are equal $P$-a.s. are identified. If $Y: (\Omega, \tf) \to (K, \mathcal{B}_K)$ is measurable, then $\bar Y$ denotes the element in  $M(\Omega,K)$ obtained by identifying measurable maps which are $P$-a.s.~equal to $Y$. Denoting by $d$ the metric on $K$, we endow the space
$M(\Omega,K)$ with the distance defined by
\beqn
\rho(\bar Y, \bar Z) = E\left[d(Y, Z)\right],
\eeqn
that corresponds to the topology of convergence in probability.
\begin{thm}
\label{rdsecA.4-t1}
Let $Z:=(Z_t,\, t\in \mathbb{T})$ be a stochastic process as described above. The following conditions are equivalent.
\begin{description}
\item{(i)}  $Z$ has a jointly measurable modification $\tilde Z: (\mathbb{T}\times \Omega, \mathcal{B}_{\mathbb{T}}\times \tf)\to (K, \mathcal{B}_K)$.
\item {(ii)} $Z$ has a separable jointly measurable modification $\tilde Z: (\mathbb{T}\times \Omega, \mathcal{B}_{\mathbb{T}}\times \tf)\to (K, \mathcal{B}_K)$.
\item {(iii)} The map from $\mathbb{T}$ to $M(\Omega,K)$ taking $t$ to $\bar Z_t$ is Borel measurable and has a separable range (if $\T$ is complete, then the requirement that this map has a separable range can be omitted).
\end{description}
If $\mathbb{T}=\re_+$ and $Z$ is adapted to a filtration $(\tf_t,\, t\in \re_+)$, then conditions (i)-(iii) are equivalent to
\begin{description}
\item {(iv)} $Z$ has a separable progressively measurable modification.
\end{description}
\end{thm}
As mentioned in \cite[Remark, p.164]{cohn}, the above theorem extends previous versions from \cite[p.~61]{doob}, \cite[p.~157]{g-s}, \cite[p.~115]{p-r}, \cite[p.~91]{neveu}  and \cite{ch-d}.


The following is a special case of Theorem \ref{rdsecA.4-t1} above.
\begin{thm}
\label{rdsecA.4-t2}
(\cite[Theorem 2]{cohn}
Suppose that the stochastic process  $Z:=(Z_t,\, t\in \mathbb{T})$ is continuous in probability. Then it has a separable jointly measurable modification.
\end{thm}
Indeed, the hypothesis of this theorem implies the validity of Theorem \ref{rdsecA.4-t1} (iii).
\medskip



\section{Distribution-valued stochastic processes}
\label{app1-1}

In this section, we first provide the background needed for the study of stochastic processes indexed by $\mathcal{S}(\R^k)$ that are random linear functionals in the sense of Definition \ref{ch1-d4}. The paradigmatic example is white noise. This is a prelude to the main goal of this section, which is the proof of Theorem \ref{ch1-t1}
concerning the existence of distribution-valued versions of these processes.

\subsection{The space $\mathcal{S}(\R^k)$ as a nuclear space}
\label{app1-s1}

The proof of Theorem \ref{ch1-t1} relies on the structure of the space $\cs(\rek)$ as a nuclear space. In this section, we first introduce some basic ingredients that explain this notion. For a more complete discussion, we refer the reader to \cite[Chapter 1]{ito84}, \cite[Chapter 4]{walsh}, and \cite[p. 143]{rs}.

\medskip

\noindent{\em Hilbertian norms and seminorms}
\medskip

Let $V$ be a vector space. A seminorm $p: V \to \IR_+$ is called Hilbertian,\index{Hilbertian}\index{seminorm!Hilbertian} or an $H$-seminorm,\index{seminorm!$H$-}\index{H@$H$-seminorm} if
\beqn
p^2(x+y) + p^2(x-y) = 2(p^2(x)+p^2(y)), \qquad x,y\in V.
\eeqn
Setting
\beqn
\langle x,y\rangle = \frac{1}{4} (p^2(x+y) - p^2(x-y)),
\eeqn
we define a pre-inner product (see e.g.~\cite[Theorem 1, Section I.5, p.~39]{yosida}) and we will write $\Vert x\Vert = p(x)$.
Suppose that $(V,p)$ is separable.
Consider the quotient space $V/N_p$, where $N_p= \{x\in V: p(x)=0\}$. By an abuse of notation, we write $V:=V/N_p$.
The space $(V,\Vert\cdot\Vert)$ is a separable pre-Hilbert (or inner-product) space.

Assume that we have two Hilbertian norms $\Vert\cdot\Vert_1$ and $\Vert\cdot\Vert_2$ on $V$  which are separable, meaning that the spaces $(V,\Vert\cdot\Vert_i)$, $i=1,2$, are separable. We say that
$\Vert\cdot\Vert_1$ is HS{\em -weaker than}\index{HS-weaker}\label{rdHS-weaker} $\Vert\cdot\Vert_2$ (HS stands for Hilbert-Schmidt), and write
\beqn
\Vert\cdot\Vert_1\le _{\rm{HS}}\Vert\cdot\Vert_2,
\eeqn
 if
\beq
\label{a1}
\sup\left\{\, \sum_{k\ge 1} \Vert e_k\Vert_1^2\right\} < \infty,
\eeq
where the supremum is over all $\Vert\cdot\Vert_2$-complete orthonormal systems (CONS) $(e_k, \, k\ge 1)$.  Clearly,
$\Vert\cdot\Vert_1\le _{{\rm HS}}\Vert\cdot\Vert_2$ is equivalent to
 the property that the inclusion mapping from
 $(E,\Vert\cdot\Vert_2)$ into $(E,\Vert\cdot\Vert_1)$ is a Hilbert-Schmidt operator.

On the space $\cs(\rek)$\label{rdS} of $\cC^\infty(\rek)$ rapidly decreasing functions (also called Schwartz test  functions), we introduce the topology $\tau$ defined by the
family of seminorms \beqn
\Vert \varphi\Vert_{m,\ell} = \sup_{x\in \rek}(1+|x|^2)^\ell |D^{m}\varphi(x)|, \qquad \ell\in\IN,  \ m\in \IN^k .
\eeqn
These are not $H$-norms. However, as we will see further in this section, $\tau$ is also given by a family of $H$-norms.

Consider the sequence $(\bar H_j)$ of Hermite polynomials\index{Hermite polynomial}\index{polynomial!Hermite} on $\re$:
\beq
\label{hermite1/2}
\bar H_j(x)= (-1)^j \exp(x^2)\frac{d^j}{dx^j}\left(\exp(-x^2)\right),\qquad j\in\IN,
\eeq
and set
\beqn
h_j(x)= \left(\pi^\half 2^j j!\right)^{-\half}
 \exp\left(-\frac{x^2}{2}\right)\bar H_j(x),\qquad j\in\IN.
\eeqn
It is well known that $(h_j,\, j\in\IN)$ is a CONS of the Hilbert space $L^2(\re)$ consisting of elements of $\cs(\re)$ (see e.g.~\cite[Lemma 3, p.~142]{rs}). From this, we can construct a CONS of  $L^2(\rek)$.

 Indeed, let $j=(j_1,\ldots,j_k)\in \IN^k$ be a multi-index, and write
$|j|= j_1+\cdots+j_k$. For any $x=(x_1,\ldots,x_k)\in\rek$, we define
\beqn
h_j(x) = \prod_{\ell=1}^k h_{j_\ell}(x_\ell).
\eeqn
The family $(h_j,\, j\in \IN^k)$ is a CONS of $L^2(\rek)$.

For every $\varphi\in L^2(\rek)$, $n\in \Z$, define
\beq
\label{app1.20}
\Vert \varphi\Vert _n^2 = \sum_{j\in \IN^k} (2|j|+k)^{2n}\, \langle\varphi,h_j\rangle_{L^2(\rek)}^2,
\eeq
where $\langle\cdot,\cdot\rangle_{L^2(\rek)}$ denotes the standard inner product in $L^2(\rek)$.

Observe that $n\to \Vert \cdot\Vert _n$ is increasing. Clearly, for $n\in \Z_-$, $\Vert \varphi\Vert _n<\infty$, and
for any $\varphi\in \cs(\rek)$, $\Vert \varphi\Vert _n<\infty$ for all $n\in \Z$ (see \cite[p.~142-143]{rs} and \cite[Section 1.3, p.~6-7]{ito84}). Moreover, $\Vert \cdot\Vert _n$  is an $H$-norm for each $n\in \Z$.

The space $(\cs(\rek), \Vert \cdot\Vert _n)$ is a pre-Hilbert space and
\beqn
h_j^n = (2|j|+k)^{-n}\, h_j, \qquad j\in \IN^k,
\eeqn
form a CONS in $(\cs(\rek), \Vert\cdot\Vert_n)$. The completion of $(\cs(\rek), \Vert\cdot\Vert_n)$ is a separable Hilbert space denoted by $(\cs_n(\rek), \Vert\cdot\Vert_n)$\label{rdSn} and
$(h_j^n,\, j\in \IN^k)$ is also a CONS in $(\cs_n(\rek), \Vert\cdot\Vert_n)$.

The Schwartz topology $\tau$ on $\cs(\rek)$ is equivalent to the Hilbertian topology determined by the sequence of norms $(\Vert\cdot\Vert_n,\, n\in\N)$ (see \cite[Theorem V.13 p.~143]{rs} and \cite[Section 1.3, p.~6-7]{ito84}). A neighbourhood basis of $0$ is given by
\beq
\label{basistau}
\{\varphi\in \cs(\rek): \Vert \varphi\Vert_n <\ep\},\qquad   n\in\N,\ \ep>0.
\eeq
Hence, a sequence $(\varphi^{(\ell)},\, \ell\in\N)$ in $\cs(\rek)$ converges to
$\varphi\in \cs(\rek)$ if and only if
\beqn
\lim_{\ell\to\infty} \Vert\varphi^{(\ell)}-\varphi\Vert_n = 0,\qquad \text{for all } n\in \N.
\eeqn

Observe that
\beqn
\sum_{j\in\IN^k}\Vert h_j^{n}\Vert_m^2 = \sum_{j\in\IN^k} (2|j|+k)^{-2(n-m)},
\eeqn
and that the last series converges if and only if $\frac{k}{2}+m<n$. Therefore, in this case,
\beq
\label{app1.21}
\Vert \cdot\Vert_m\le_{\rm HS} \Vert \cdot\Vert_n.
\eeq
Because of these properties, the space $(\cs(\rek), \tau)$ belongs to the class of functional spaces termed {\em nuclear}\index{nuclear space}\index{space!nuclear} (see \cite[Section 1.2]{ito84} for more details).

Consider now $\csp(\rek)$,\label{rdS'} the dual of $(\cs(\rek), \tau)$, called the space of {\em tempered distributions} (also called Schwartz distributions).\index{tempered distribution}\index{Schwartz distribution}\index{distribution!Schwartz}\index{distribution!tempered} This is the space of linear functionals $\alpha$ on $\cs(\rek)$ such that there are $C \in \IR_+$ and $n\in \IN$ satisfying
\beq
\label{A.2(*A)}
\vert \alpha(\varphi) \vert \leq C \Vert \varphi\Vert_n, \quad {\text{for all}}\ \varphi\in \cs(\rek).
\eeq
 For $\alpha \in \csp(\R^k)$ and $n \in \N$ such that \eqref{A.2(*A)} holds, define
\beqn
\Vert\alpha\Vert_n^\prime = \sup\{|\alpha(\varphi)|: \Vert \varphi\Vert_n\le 1\}.
\eeqn
Then for any CONS $(e_{n,\ell},\, \ell\in\N)$ of $(\cs_n(\rek), \Vert\cdot\Vert_n)$,
\beq\label{rdeA.1.5}
\Vert\alpha\Vert_n^\prime = \left(\sum_{\ell\in\N}\, [\alpha(e_{n,\ell})]^2\right)^{\frac{1}{2}}< \infty.
\eeq
Indeed, fix $\varphi\in\cs(\rek)$ with $\Vert \varphi\Vert_n\le 1$, and consider the expansion
\beqn
\varphi = \sum_{\ell\in\N}\, \varphi_\ell\, e_{n,\ell}
\eeqn
 relative to $(e_{n,\ell},\, \ell\in\N)$. Using the Cauchy-Schwarz inequality, we obtain
\begin{align*}
\vert\alpha(\varphi)\vert &\le \Vert\varphi\Vert_n \left(\sum_{\ell\in\N}\, [\alpha(e_{n,\ell})]^2\right)^{\frac{1}{2}}
 \le \left(\sum_{\ell\in\N}[\alpha(e_{n,\ell})]^2\right)^{\frac{1}{2}}.
\end{align*}
For the converse inequality, we take
$$
   \varphi_\ell = \alpha(e_{n,\ell}) \left(\sum_{m\in\N}\, [\alpha(e_{n,m})]^2\right)^{-1/2}
$$
to get $\Vert \varphi \Vert _n = 1$ and
$$
   \alpha(\varphi) = \left(\sum_{\ell\in\N}[\alpha(e_{n,\ell})]^2\right)^{\frac{1}{2}},
$$
and this proves \eqref{rdeA.1.5}.

It follows from \eqref{rdeA.1.5} that $\Vert\cdot\Vert_n^\prime$ is an $H$-norm on the space
\beqn
\csp_n(\rek) = \{\alpha \in 	\csp(\rek): \Vert\alpha\Vert_n^\prime<\infty\},
\eeqn
and $(\csp_n(\rek), \Vert\cdot\Vert_n^\prime)$\label{rdSn'} is a separable Hilbert space.

For any $n\in \N$, we have the inclusions
\beqn
\cs(\rek)\subset \cs_n(\rek)\subset \cs_0(\rek)\cong L^2(\rek)\cong \csp_0(\rek)\subset \csp_n(\rek)\subset\csp(\rek).
\eeqn


\subsection{Versions with values in $\mathcal{S'}(\IR^{k})$}
\label{app1-s3}

In this section, we give a proof of Theorem \ref{ch1-t1} on existence of $\mathcal{S}^{\prime}(\rek)$-valued versions of random linear functionals. We refer to Definition \ref{ch1-d4} for these notions.
\medskip

\begin{thm} 
\label{app1-ch1-t1}
Fix $m \in \N$. Let $(X(\varphi),\, \varphi \in \mathcal{S}(\IR^{k}))$ be a random linear functional that is continuous in probability for $\Vert \cdot \Vert_m$ (that is, $\Vert \varphi_j - \varphi \Vert_m \to 0$ implies $X(\varphi_j) \to X(\varphi)$ in probability). Then $X$ has a version with values in $\mathcal{S'}(\IR^{k})$. In fact, for $n > \frac{k}{2} + m$, $X$ has a version with values in $\csp_n(\R^k)$.
  \end{thm}
\begin{proof}
Recall that $\Vert Y \Vert_{L^0(\Omega)} := E(\vert Y \vert \wedge 1)$ is a metric that corresponds to the topology of convergence in probability.

 Fix $\ep > 0$. By assumption, there is $\delta > 0$ such that $E(\vert X(\varphi) \vert \wedge 1) \leq \ep$ whenever $\Vert \varphi \Vert_m \leq \delta$.

Fix $n\in \N$ satisfying $\frac{k}{2}+m<n$, which by \eqref{app1.21} implies  $\Vert \cdot\Vert_m\le_{\rm HS} \Vert \cdot\Vert_n$.
\smallskip

 \noindent{\em Step 1.}
 Let $(e_j,\, j\ge 1)$ be a CONS in $(\mathcal{S}(\rek),\Vert\cdot\Vert_n)$. As a preliminary ingredient for the proof of the theorem, we verify that
 the series $\sum_{j\ge 1} X(e_j)^2$ converges a.s.

First, we check that for all $\varphi\in\cs(\rek)$,
 \beq
 \label{a1.30}
 {\rm Re}\ E\left[e^{iX(\varphi)}\right] \ge 1-2\ep - 2\ep \frac{\Vert\varphi\Vert_m^2}{\delta^2}.
 \eeq
 For this, we will use the property
  \beqn
 {\rm Re}\ E\left[e^{iX(\varphi)}\right] = E\left[\cos X(\varphi)\right] \ge 1-\frac{1}{2} E\left[|X(\varphi)|^2 \wedge 4\right],
 \eeqn
 along with the inequalities
 \beqn
 |z|^2\wedge 4 \le 4(|z|^2\wedge 1)\le 4(|z|\wedge 1).
 \eeqn
 Assume first that $\Vert\varphi\Vert_m\le \delta$. Then by our choice of $\delta$,
 \beqn
 E\left[|X(\varphi)|^2 \wedge 4\right] \le 4 E\left[|X(\varphi)| \wedge 1\right] \le 4\ep.
 \eeqn
 Suppose next that $\Vert\varphi\Vert_m > \delta$. Because of the linearity of $X$,
 \beqn
 X(\varphi) = \frac{\Vert\varphi\Vert_m}{\delta}\, X\left(\frac{\delta\varphi}{\Vert\varphi\Vert_m}\right),
 \eeqn
 and then
  \begin{align*}
 E\left[|X(\varphi)|^2 \wedge 4\right] & \le   \frac{\Vert\varphi\Vert_m^2}{\delta^2}\,E\left[\left\vert X\left(\frac{\delta\varphi}{\Vert\varphi\Vert_m}\right)\right\vert^2\wedge 4\right]\\
 & \le 4\ep\,\frac{\Vert\varphi\Vert_m^2}{\delta^2}.
 \end{align*}
 Hence, \eqref{a1.30} is proved.

We continue the proof by introducing a sequence $(Y_j,\, j\ge 1)$ of i.i.d.~Gaussian ${\rm N}(0,\sigma^2)$ random variables, independent of $X$, and defining the $\mathcal{S}(\rek)$-valued random variables
\beqn
\tilde \varphi_N = \sum_{j=1}^N Y_j\, e_j, \qquad N\ge 1.
\eeqn
Given the random variables $(X(e_j),\, j=1,\ldots,N)$, $X(\tilde \varphi_N)= \sum_{j=1}^N Y_j\, X(e_j)$ has the conditional distribution
${\rm N}(0,\sigma^2\sum_{j=1}^N X^2(e_j))$. Therefore, by the properties of conditional expectation, we have
\begin{align}
\label{a1.31}
 {\rm Re}\ E\left[e^{iX(\tilde\varphi_N)}\right] & =  {\rm Re}\ E\left[E\left[e^{i \sum_{j=1}^N Y_j\, X(e_j)}\, \big\vert\, \sigma(X(e_j),\, j=1,\ldots,N)\right]\right] \notag\\
 & = E\left[e^{-\frac{\sigma^2}{2} \sum_{j=1}^NX^2(e_j)}\right].
 \end{align}
By applying \eqref{a1.30},
\begin{align}\nonumber
 {\rm Re}\ E\left[e^{iX(\tilde\varphi_N)}\right] &=
 E\left[{\rm Re}\ E\left[ e^{iX(\tilde\varphi_N)}\, |\, Y_j,\, j=1,\dots,N \right] \right]\\
 &\ge E\left[1-2\ep-2\ep \frac{\Vert\tilde\varphi_N\Vert_m^2}{\delta^2}\right].
 \label{a1.32}
\end{align}
From \eqref{a1.31} and  \eqref{a1.32}, and because $\Vert\cdot\Vert_m$ is Hilbertian
and the $Y_j$ are independent and centred,  we obtain
\begin{align*}
 E\left[e^{-\frac{\sigma^2}{2} \sum_{j=1}^NX^2(e_j)}\right] &\ge E\left[1-2\ep-2\ep \frac{\Vert\tilde\varphi_N\Vert_m^2}{\delta^2}\right]\\
 &= 1 - 2 \ep - 2 \ep \delta^{-2} E\left[\left\langle \sum_{j=1}^N Y_j e_j, \sum_{k=1}^N Y_k e_k \right\rangle_m\, \right] \\
 &= 1 - 2 \ep - 2 \ep \delta^{-2} \sum_{j=1}^N \sum_{k=1}^N\, \langle e_j, e_k \rangle_m E[Y_j Y_k] \\
 &= 1-2\ep-2\ep\delta^{-2}\sigma^2\sum_{j=1}^N\Vert e_j\Vert_m^2.
\end{align*}
Since $\Vert\cdot\Vert_m\le _{\rm HS}\Vert\cdot\Vert_n$ because $n>m+\frac{k}{2}$, the series $\sum_{j=1}^\infty\Vert e_j\Vert_m^2$
converges.
Letting $N \to \infty$ above and using dominated convergence, we obtain
\begin{align*}
P\left\{\sum_{j=1}^\infty X^2(e_j)<\infty\right\}&\ge
E\left[e^{-\frac{\sigma^2}{2} \sum_{j=1}^\infty X^2(e_j)}1_{ \sum_{j=1}^\infty X^2(e_j)<\infty}\right] \\
 &=E\left[e^{-\frac{\sigma^2}{2} \sum_{j=1}^\infty X^2(e_j)}\right] \\
 &\ge 1-2\ep-2\ep\delta^{-2}\sigma^2\sum_{j=1}^\infty\Vert e_j\Vert_m^2.
 \end{align*}
 Finally, by letting $\sigma\to 0$, we have
 \beq
 \label{finals1}
 P\left\{\sum_{j=1}^\infty X^2(e_j)<\infty\right\} \ge 1-2\ep.
 \eeq
 Since $\ep>0$ is arbitrary, we conclude that $\sum_{j \geq 1}  X(e_j)^2$ converges a.s. This completes the proof of Step 1.
 \medskip

 \noindent{\em Step 2.} We now construct the version of $X$ with values in $\csp(\rek)$ (and even in $\csp_n(\rek)$).
 \smallskip

 Let $\Omega_0 = \left\{\sum_{j=1}^\infty X^2(e_j)<\infty\right\}$. We have just proved that $P(\Omega_0)=1$. For $\varphi\in \cs(\rek)$, define
 \beq
 \label{a1.33}
 \tilde X(\varphi, \omega)=
 \begin{cases}
 \sum_{j=1}^\infty \langle\varphi,e_j\rangle_nX(e_j, \omega), & \omega\in \Omega_0,\\
 0, & \omega\notin \Omega_0.
 \end{cases}
 \eeq
 By the Cauchy-Schwartz inequality and Parseval's identity, we see that the series in  \eqref{a1.33} converges absolutely.
 Further,
 \beqn
    \vert \tilde X(\varphi, \omega) \vert \leq \Vert \varphi\Vert_n\left(\sum_{j=1}^\infty X^2(e_j, \omega)\right)^\half,
    \eeqn
so for $\omega \in \Omega_0$, the linear functional $\varphi \mapsto \tilde X(\varphi,\omega)$ belongs to $\mathcal{S}^\prime(\rek)$.      
This yields (b) in Definition \ref{ch1-d4}.
  Moreover, by the property \eqref{rdeA.1.5} of the norm $\Vert \cdot\Vert_{n}^\prime$, on the set $\Omega_0$,
 \beqn
 \Vert \tilde X\Vert_{n}^\prime = \sum_{j=1}^\infty \tilde X^2(e_j) = \sum_{j=1}^\infty X^2(e_j) < \infty.
 \eeqn
 Therefore, for $\omega \in \Omega_0$, the linear functional $\varphi \mapsto \tilde X(\varphi,\omega)$ belongs in fact to $\csp_n(\rek)$.

As for condition (a) in Definition \ref{ch1-d4}, notice that by taking $\varphi = e_j$, we obtain $\tilde X(e_j) = X(e_j)$ a.s. Now fix $\varphi \in \mathcal{S}(\R^k)$ and set $\varphi_N
= \sum_{j=1}^N \langle\varphi,e_j\rangle_n e_j$. Because $X$ is a random linear functional, $\tilde X(\varphi_N) = X(\varphi_N)$ a.s., where the ``a.s." depends on $\varphi$ and $N$.
 Clearly,
 \beq
 \label{A.2-(*1)}
 \lim_{N\to\infty}\Vert\varphi-\varphi_N\Vert_n = 0.
 \eeq
 Since $\Vert \cdot\Vert_m\le \Vert \cdot\Vert_n$, we also have
 \beq
 \label{A.2-(*1)(*2)}
 \lim_{N\to\infty}\Vert\varphi-\varphi_N\Vert_m = 0.
 \eeq
 For $\omega \in \Omega_0$, by \eqref{A.2-(*1)},
 \beqn
    \tilde X(\varphi, \omega) = \lim_{N \to \infty} \tilde X(\varphi_N, \omega),
    \eeqn
that is,
\beqn
\tilde X(\varphi) = \lim_{N \to \infty} \tilde X(\varphi_N)\quad {\text{a.s.}},
\eeqn
and by \eqref{A.2-(*1)(*2)},
\beqn
      X(\varphi) = \lim_{N \to \infty} X(\varphi_N)
      \eeqn
  in probability, therefore $\tilde X(\varphi) = X(\varphi)$ a.s.

   This completes the proof of Theorem \ref{app1-ch1-t1}.
 \end{proof}
 \begin{cor} {\em (Theorem \ref{ch1-t1} of Chapter \ref{ch1})}\
\label{app1-ch1-c1}
 Let $(X(\varphi),\, \varphi \in \mathcal{S}(\IR^{k}))$ be a random linear functional which is continuous in $L^p(\Omega)$, for some $p\ge 1$ $($that is, $\varphi_{n} \to \varphi$ in $\mathcal{S}(\IR^{k})$ implies  $X(\varphi_{n}) \to X(\varphi)$ in $L^p(\Omega))$. Then $X$ has a version with values in $\mathcal{S'}(\IR^{k})$.
  \end{cor}
  \begin{proof}
  First, as in \cite[p. 331]{walsh}, we notice that if $\varphi \mapsto X(\varphi)$ from $(\mathcal{S}(\R^k), \tau)$ into $L^p(\Omega)$ is continuous, then  $\varphi \mapsto X(\varphi)$ from $(\mathcal{S}_m(\rek), \Vert \cdot \Vert_m)$ into $L^p(\Omega)$ is continuous for some $m \in \N$ (the converse is clearly true).

Indeed, assuming the continuity of $\varphi \mapsto X(\varphi)$ from $(\mathcal{S}(\rek), \tau)$ into $L^p(\Omega)$, there is a $\tau$-neighborhood $U$ of $0\in \mathcal{S}(\rek)$ such that $\varphi \in U$ implies $\Vert X(\varphi) \Vert_{L^p(\Omega)} \leq 1$. Thus, there are $m \in \N $ and  $\delta_0 > 0$  such that the member $\{\varphi \in \mathcal{S}(\R^k): \Vert \varphi \Vert_m < \delta_0\}$ of the neighborhood basis of $0$ is contained in $U$. By linearity, for any $\ep > 0$, if $\Vert \varphi \Vert_m < \ep \delta_0$, then $\Vert X(\varphi) \Vert_{L^p(\Omega)} \leq \ep$, that is, $\varphi \mapsto X(\varphi)$ from $(\mathcal{S}_m(\rek), \Vert \cdot \Vert_m)$ into $L^p(\Omega)$ is continuous.

Since convergence in $L^p(\Omega)$ implies convergence in probability, the conclusion follows from Theorem \ref{app1-ch1-t1}.
   \end{proof}

 \medskip

 \begin{remark}
 \label{cont-probab}
Corollary \ref{app1-ch1-c1} (or Theorem \ref{ch1-t1}) also holds if the assumption of continuity in $L^p(\Omega)$, for some $p\ge 1$, of the random linear functional $(X(\varphi),\, \varphi\in \mathcal{S}(\rek))$ is replaced by {\em continuity in probability} (see \cite[Corollary 4.2, p.~332]{walsh}). A more general topological setting is discussed in \cite[[Theorem 2.3.2, p.~24]{ito84} which, in particular, provides a proof of \cite[Corollary 4.2, p.~332]{walsh}.
\end{remark}


\section{Regularity of sample paths}
\label{app1-3}


A fundamental result for the study of the regularity of sample paths of stochastic processes is Kolmogorov's continuity criterion.\index{Kolmogorov!continuity criterion}\index{continuity criterion!Kolmogorov's}\index{criterion!Kolmogorov's continuity}\index{theorem!Kolmogorov's continuity} An elegant version can be found for example in \cite[(2.1) Theorem, page 26]{ry}, which states the following:

 Let $Z= (Z(x),\, x\in [0,1]^k)$ be a random field with values in a separable Banach space, for which there exist constants $p, \varepsilon>0$ and $C<\infty$ such that, for any $x,y\in [0,1]^k$,
\beq
\label{k1}
E\left(\Vert Z(x) - Z(y)\Vert^p\right) \le C |x-y|^{k+\varepsilon}
\eeq
(where $\Vert\cdot\Vert$ denotes the norm in the Banach space).
Then $Z$ has a version $\tilde Z$ such that for every $\alpha\in[0,\varepsilon/p[$,
\beqn
E\left[\left(\sup_{x\ne y} \frac{\Vert\tilde Z(x) - \tilde Z(y)\Vert}{|x-y|^\alpha}\right)^p\right] < \infty.
\eeqn
In particular, the sample paths of $\tilde Z$ are H\"older-continuous with exponent $\alpha$, a.s.

The validity of this result extends without further work to random fields $Z$ consisting of random vectors with values in a complete separable metric space. For $k=1$, we refer to \cite[Theorem 2.9, p. 24]{legall-2013}.

In the above condition \eqref{k1}, all the components of the indices $x\in[0,1]^k$ play the same role, a property that is satisfied for instance when $Z$ is {\em isotropic.}\index{isotropic} However, in many cases,
random field solutions to SPDEs are anisotropic\index{anisotropic} random fields, that is, stochastic processes whose behaviour differs in two or more components of the index set, and the index set may not be a hypercube.
 This motivates various extensions of Kolmogorov's continuity criterion.
In this section, we present one such generalisation that is suitable for the situations addressed in this book. First, we will consider arbitrary random fields $u=(u(t,x),\, (t,x)\in \IR_+ \times D)$, and then we will particularize the results to Gaussian processes.


\subsection{Anisotropic Kolmogorov's continuity criterion}
\label{app1-3-ss1}
In the next statement, $u=(u(t,x),\, (t,x)\in \IR_+ \times D)$ is a real-valued random field, and $D$ is a  domain of $\rek$.
\begin{thm}
\label{ch1'-s7-t2}
Let $I \subset \R_+$ and $D\subset \rek$ be non-empty and bounded sets.
Fix $\alpha_i\in\, ]0,1]$, $i=0,\ldots, k$, and for any $(t,x), (s,y) \in I\times D$, $x=(x_1,\dots,x_k)$, $y=(y_1,\dots,y_k)$, define
\beq
\label{norm-aniso}
{\bf\Delta}(t,x;s,y) = |t-s|^{\alpha_0}+ \sum_{i=1}^k |x_i-y_i|^{\alpha_i} 
\eeq
and $Q = \sum_{i=0}^k \frac{1}{\alpha_i}$.

Suppose that for some constant $K<\infty$ and for some $p> Q$, for all $(t,x),\, (s,y)\in I\times D$,
\beq
\label{ch1'-s7.18}
E\left[|u(t,x)-u(s,y)|^p\right] \le K \left({\bf\Delta}(t,x;s,y)\right)^p.
\eeq
Then $(u(t,x),\, (t,x)\in I\times D)$ has a continuous version $\tilde u=(\tilde u(t,x),\, (t,x)\in I\times D)$, which extends continuously to $\bar I\times \bar D$.
Further, for $\alpha\in\,]Q/p,1[$, there is a constant $a(I,D,\alpha,p,Q)<\infty$, which is an increasing function of $I$ and $D$, and a non-negative random variable $Y$ such that
\beq
\label{ch1'-s7.19bis}
E[Y^p] \le K a(I,D,\alpha,p,Q) <\infty,
\eeq
and for all $(t,x),\, (s,y)\in \bar I\times \bar D$,
\beq
\label{ch1'-s7.19}
|\tilde u(t,x)-\tilde u(s,y)| \le Y \left({\bf\Delta}(t,x;s,y)\right)^{\alpha-\frac{Q}{p}}.
\eeq

Therefore, the paths of $\tilde u$ are jointly Hölder continuous on $\bar I \times \bar D$
 with exponents $(\beta_0, \beta_1, \dots, \beta_k)$  a.s., provided $\beta_i < \alpha_i (1 - \frac{Q}{p})$, $i = 0, \dots, k$. In addition,
\beq
\label{ch1'-s7.20}
E\left[\left(\sup_{(t,x)\ne(s,y)}
\frac{|\tilde u(t,x)-\tilde u(s,y)|}
{\left({\bf\Delta}(t,x;s,y)\right)^{\alpha-\frac{Q}{p}
}}
\right)^p\right] \le K a(I,D,\alpha,p,Q).
\eeq
If $E[\vert u(t_0, x_0) \vert] < \infty$ for some $(t_0, x_0) \in I \times D$, then
\beqn
E\left[\sup_{(t, x) \in I \times D} \vert \tilde u(t,x) \vert^p\right] < \infty.
\eeqn

Moreover, if $0\in \bar I$ and $E\left[\sup_{x\in D}|\tilde u(0,x)|^p\right] \le C_1$, then there exists a finite non-negative constant $c_2(I,D,\alpha,p,Q)$ such that for all $t\in I$,
\beq
\label{ch1'-s7.20bis}
E\left[\sup_{(s,x)\in [0,t]\times D}|\tilde u(s,x)|^p\right] \le 2^{p-1}C_1 + K c_2(I,D,\alpha,p,Q) t^{p\alpha_0(\alpha-Q/p)}.
\eeq
\end{thm}

\begin{proof}
In order to simplify the notation, we replace $\R_+ \times \R^k$ by $\R^{k'}$, with $k' = 1+k$, and identify $(t,x) \in \R_+ \times \R^k$ with $(t,x_1,\dots x_k) \in \R^{k'}$, and $I \times D$ with a bounded subset $D^\prime$ of $\R^{k'}$. Then we remove the primes, so we are dealing with a bounded subset $D$ of $\R^k$, a process $(u(x),\, x \in D)$, and constants $K < \infty$ and $p > Q$ satisfying
     \beq
     \label{kolmogorov(*1)}
     E[\vert u(x) - u(y) \vert^p] \leq K\, ({\bf\Delta}(x,y))^p,
     \eeq
where
$
   {\bf\Delta}(x,y) := \sum_{i=1}^k \vert x_i - y_i \vert^{\alpha_i}
 $
    and $Q := \sum_{i=1}^k \frac{1}{\alpha_i}$.
 We assume further (by translation and scaling) that $D \subset [0,1[^k$.

   Let $\rho(x,y)$ be the metric defined by
   \beqn
   \rho(x,y) = \max(|x_1-y_1|^{\alpha_1}, \dots, |x_k-y_k|^{\alpha_k}).
   \eeqn
Notice that
\beqn
\frac{1}{k} {\bf\Delta}(x,y) \le \rho(x,y)\le {\bf\Delta}(x,y).
\eeqn

 Fix $\alpha \in\, ]Q/p, 1[$. Choose $m \in \N$ large enough so that
    $p - k / m > Q$    and     $k / m < p (1-\alpha)$.
Observe that $1 + 2^{m \alpha_i^{-1}} \leq 2^{m \alpha_i^{-1} + 1}$, $i=1,\dots, k$.

 Consider the subset of $[0,1]$ consisting of elements $s \in [0,1]$ of the form $s = s_{\ell,1,i} := (\ell 2^{- m \alpha_i^{-1}}) \wedge 1$,  $\ell = 0, \dots, \lceil2^{m \alpha_i^{-1}}\rceil$, where, for $s\in \R$, $\lceil s\rceil$
 denotes the {\em ceiling function}, that is, the smallest integer greater than or equal to $s$.

   Let $D_{1,i}$ denote the set $\{s_{\ell,1,i}: \ell \in \N,\  0 \leq \ell < \lceil 2^{m \alpha_i^{-1}}\rceil \}$. All of the intervals $[s_{\ell,1,i},\, s_{\ell + 1,1,i} \wedge 1[$ have equal length $2^{- m \alpha_i^{-1}}$, except for the last one if $2^{\alpha_i^{-1}}$ is not an integer, and then it has even shorter length.   Note that  ${\rm card}\,D_{1,i} \leq 1 + 2^{m \alpha_i^{-1}} \leq 2^{m \alpha_i^{-1} + 1}$.

   For $n \geq 2$, define inductively the set $D_{n,i} \subset [0, 1[$ consisting of elements of the form $s_\ell := t_1 + \ell  2^{-n m \alpha_i^{-1}}$, $\ell\in\N$, with the requirement that $s_\ell < t_2$, where $t_1 < t_2$ are consecutive elements of $D_{n-1,i}$. Since $t_2 - t_1 \leq 2^{- (n-1) m \alpha_i^{-1}}$, for each $t_1 \in D_{n-1,i}$, there are at most $1 + 2^{m \alpha_i^{-1}} \leq 2^{m \alpha_i^{-1} + 1}$ elements of this kind.
 Assuming by induction that $\text{card}\, D_{n-1, i} \leq 2^{(n-1)(m \alpha_i^{-1} + 1)}$, we deduce that
  \beqn
  \text{card}\,D_{n,i} \leq 2^{n (m\alpha_i^{-1} + 1)} = 2^{n m(\alpha_i^{-1} + 1/m)}.
  \eeqn

Intervals whose endpoints are consecutive elements of $D_{n,i}$ have length $2^{- n m \alpha_i^{-1}}$ (except for those intervals with the right endpoint in $D_{n-1,i}$, which are even shorter).
Then each $D_{n,i}$ partitions $[0,1]$ and is a refinement of $D_{n-1,i}$. There are at most $2^{n m (\alpha_i^{-1} + 1/m)}$  elements of $D_{n,i}$, and consecutive elements are equally spaced except possibly if one is in $D_{n-1,i}$ and the other is immediately to its left (if $2^{\alpha_i^{-1}}$ is an integer, these intervals have the same length as the others).

   Let $D_n := D_{n,1} \times \cdots \times D_{n,k}$. Then ${\rm card}\,D_n \leq 2^{nm(Q + k / m)}$. For $x \in D_n$, let
  $I_n(x) := [x_1, y_1 [\, \times \cdots \times [x_k, y_k [$,
where $x_i < y_i $ are consecutive elements of $D_{n,i}$. Then $y_i - x_i \leq  2^{- n m \alpha_i^{-1}}$, $\vert I_n(x) \vert \leq 2^{-n m Q}$ and $\rho(x, y) \leq 2^{-n m}$.  In addition, for distinct points $x,y \in D_n$, the boxes $I_n(x)$ and $I_n(y)$ are disjoint.
Therefore, the boxes $(I_n(x),\, x \in D_n)$ form a partition of $[0,1[^k$.

  Define $\tilde{D}_n
   = \{x \in D_n:  I_n(x) \cap D \neq \emptyset\}$.
   Then $D \subset \cup_{x \in \tilde D_n} I_n(x)$.
    For $n = 1,2,\dots,$ and $x \in \tilde{D}_n$, we denote by $\theta_n(x)$ an element of $I_n(x) \cap D$.
   This element is arbitrary unless there is $1 \leq k < n$ and $y \in \tilde D_k$ such that $\theta_k(y) \in I_n(x)$, in which case we set $\theta_n(x) = \theta_k(y)$.


 For its further use, we point out the following facts regarding $\tilde D_n$:
 \begin{description}
 \item{(i)}  ${\rm card}\,\tilde D_n\leq 2^{n m(Q + k / m)}$, and for each $x\in\tilde D_n$, there are at most $2^{m Q+k}$ possible values for $y\in \tilde{D}_{n+1} \cap I_n(x)$, and then $\rho(\theta_n(x), \theta_{n+1}(y)) \leq 2^{-n m}$.
\item{(ii)} For each $x, y \in \tilde{D}_n$ such that $\rho(x,y) \leq 2^{-n m}$, and for each $i$, $x_i - y_i $ can only take three values (usually $\{-2^{-nm},\, 0,\, 2^{-nm}\}$), so there are at most $3^k$ possible values for
$y$ and then $\rho(\theta_n(x), \theta_n(y)) \leq  2^{1-n m}$.
\end{description}

  Let $S_n = \{\theta_n(x): x \in \tilde{D}_n \}$. Notice that $S_n$ is finite and $S_n \subset S_k$ for all $k > n$. Define $S := \cup_{n \geq 1} S_n$. It is easy to check that $S$ is dense in $D$. Indeed, for each $x \in D$, there is $n\ge 1$ and $y \in \tilde{D}_n$ such that $x \in I_n(y).$ Then $\theta_n(y) \in I_n(y) \cap D \subset S$ and $\rho(x, \theta_n(y)) \leq 2^{-n m}$.
\smallskip

Define the random variables
\begin{align*}
    Y_n &:= \sup \left\{ \vert u(\theta_n(x)) - u(\theta_{n+1}(y))\vert: x \in \tilde{D}_n,\ y \in \tilde{D}_{n+1} \cap I_n(x) \right\},\\
 Z_n &:= \sup \left\{ \vert u(\theta_n(x)) - u(\theta_n(y))\vert: x,\, y \in \tilde{D}_n,\ \rho(x,y) \leq 2^{-n m} \right\}.
 \end{align*}
By \eqref{kolmogorov(*1)}, $E[|u(x)-u(y)|^p] \le K k^p\rho(x,y)^p$, and using (i) and (ii) above respectively, we see that
\begin{align*}
   E[Y_n^p] &\leq \sum_{x \in \tilde{D}_n,\, y \in \tilde{D}_{n+1} \cap I_n(x)} E[\vert u(\theta_n(x)) - u(\theta_{n+1}(y)) \vert^p] \\
   \\
   &\leq 2^{n m(Q + k / m)} 2^{m Q+k} k^{p} K 2^{-p n m} = \bar{c}_0 K 2^{-n m (p - Q - k / m)},
   \end{align*}
with $\bar c_0= 2^{mQ+k}k^p$, and
\begin{align*}
   E[Z_n^p] &\leq \sum_{x, y \in \tilde{D}_n:\, \rho(x,y) \leq 2^{-n m}} E[\vert u(\theta_n(x)) - u(\theta_n(y)) \vert^p] \\
   \\
   &\leq 2^{n m(Q + k / m)} 3^k k^{p} K 2^{p (1 - n m)} = \bar{c}_1 K 2^{-n m (p - Q - k / m)},
\end{align*}
with $\bar c_1= (2k)^p 3^k$, respectively.
Define
\beqn
  \tilde Y_n = 2^{n m (\alpha - Q/p)} Y_n  \quad   {\text{and}}\quad   \tilde{Z}_n = 2^{n m (\alpha - Q/p)} Z_n.
  \eeqn
Then
\beq
\label{kolmogorov(*2)}
     E[\tilde{Y}_n^p] \leq   \bar{c}_0 K 2^{-n m [p (1 -\alpha) - k / m]}\quad   {\text{and}}\quad    E[\tilde{Z}_n^p] \leq \bar c_1 K 2^{-n m [p (1-\alpha) - k / m]} .
\eeq

For distinct $x, y \in S$, choose an integer $n \geq 0$ such that
    $2^{-(n+1)m} < \rho(x,y) \leq 2^{-n m}$,
and choose $\ell \geq n$ large enough so that $x \in S_\ell$ and $y \in S_\ell$. For $j= n, ..., \ell$, choose $x(j),\, y(j) \in \tilde{D}_j$ such that $x \in I_j(x(j))$ and $y \in I_j(y(j))$. By construction, $\theta_\ell(x(\ell)) = x$, $\theta_\ell(y(\ell)) = y$, $x(j) \in \tilde{D}_j$ and
if $j < \ell$, then
$x(j+1) \in \tilde{D}_{j+1} \cap I_j(x(j))$ (because $x \in I_{j+1}(x(j+1)) \cap I_j(x(j))$, so $I_{j+1}(x(j+1)) \subset I_j(x(j)))$. Similar properties hold for $y$. Furthermore, $\rho(x(n), y(n)) \leq 2^{-n m}$. By the triangle inequality,
\begin{align*}
   \vert u(x) - u(y) \vert
      &\leq \left\vert u(\theta_n(x)) - u(\theta_n(y)) \right\vert\\
      &\qquad + \sum_{j=n}^{\ell -1} \left(\vert u(\theta_{j+1}(x(j+1))) - u(\theta_j(x(j))) \vert\right.\\
      &\left.\qquad\qquad\qquad + \vert u(\theta_{j+1}(y(j+1))) - u(\theta_j(y(j))) \vert \right)\\
      &\leq Z_n + 2 \sum_{j=n}^{\ell -1} Y_j\\
      &= 2^{-n m (\alpha - Q/p)}\\
      &\qquad \times\left[2^{n m (\alpha - Q/p} Z_n + 2 \sum_{j=n}^{\ell -1} (2^{(n - j) m (\alpha - Q/p)} 2^{j m (\alpha - Q/p)} Y_j)\right]\\
   & =  2^{-n m (\alpha - Q/p)} \left[\tilde Z_n + 2 \sum_{j=n}^{\ell -1} (2^{-(j - n) m (\alpha - Q/p)}\tilde Y_j)\right]\\
     &\leq (2^m \rho(x,y))^{\alpha - Q/p} \left[\tilde Z_n + 2 (1 - 2^{- m (\alpha - Q/p)})^{-1} \sup_{j \geq n} \tilde Y_j\right].
\end{align*}
In particular, for all distinct $x, y \in S$,
\beq
\label{crucial-S}
    \frac{\vert u(x) - u(y) \vert}{({\bf\Delta}(x,y))^{\alpha - Q/p}}
      \leq \frac{\vert u(x) - u(y) \vert}{(\rho(x,y))^{\alpha - Q/p}} \leq Y,
      \eeq
where
\beqn
  Y = 2^{m(\alpha - Q/p)} \left[ \sup_{n \geq 0} \tilde Z_n + 2 (1 - 2^{- m (\alpha - Q/p)})^{-1} \sup_{j \geq 0} \tilde Y_j \right].
  \eeqn
Notice that
\beqn
   Y^p \leq 2^{mp (\alpha - Q/p)+p-1} \left[ \sup_{n \geq 0} \tilde Z_n^p + 2^p (1 - 2^{- (\alpha - Q/p)})^{-p} \sup_{j \geq 0} \tilde Y_j^p\right]
   \eeqn
and, replacing the suprema by sums, taking the expectation and using \eqref{kolmogorov(*2)}, this gives
\begin{align}
\label{crucial-S-1}
   E[Y^p] &\leq 2^{mp (\alpha - Q/p)+p-1} \left[\sum_{n=0}^\infty E[\tilde{Z}_n^p] + 2^p (1 - 2^{- (\alpha - Q/p)})^{-p} \sum_{j=0}^\infty  E[\tilde{Y}_j^p] \right] \notag\\
     &\leq 2^{mp (\alpha - Q/p)+p-1}\left [\sum_{n=0}^\infty \bar{c}_1 K 2^{-n m [p (1-\alpha) - k / m]}\right.\notag\\
     &\left.\qquad + 2^p (1 - 2^{- (\alpha - Q/p)})^{-p} \sum_{j=0}^\infty \bar{c}_0 K 2^{-j m [p (1-\alpha) - k / m]}\right]\notag\\
     & =:K c_3(p,\alpha,Q,k) <\infty,
     \end{align}
since $k / m < p (1-\alpha)$.

Summarizing, from \eqref{crucial-S} and \eqref{crucial-S-1}, we deduce that on $S$, a.s.,  \eqref{ch1'-s7.19},  \eqref{ch1'-s7.19bis} and \eqref{ch1'-s7.20} hold with $\tilde u$ there replaced by $u$ and the constant $a(I,D,\alpha,p,Q)$ by $c_3(p,\alpha,Q,k)$ from \eqref{crucial-S-1}.


Let $U$ be the event ``the sample paths of $u$ are uniformly continuous on $S$''. For $x \in \bar S = \bar D$, define
  \beqn
  \tilde u(x):= \left\{ \begin{array}{ll} \lim_{y\in S:\, y\to x} u(y) &\text{if } x \in U,\\
         0 & \text{if } x \in U^c.
          \end{array}\right.
  \eeqn
Since $P(U)=1$ by \eqref{crucial-S}, the limit is well-defined a.s.
We now check that on $D$, $\tilde u$ is a version of $u$. Fix $x \in D$ and let $(y_n,\, n \geq 1)$ be a sequence of elements of $S$ that converges to $x$.
By \eqref{kolmogorov(*1)},
$x \mapsto u(x)$ is continuous in probability, therefore $u(y_n) \to u(x)$ in probability, and by definition of $\tilde u$, $u(y_n) \to \tilde u(x)$ a.s.,
therefore $u(x) = \tilde u(x)$ a.s. and so $\tilde u$
  is a continuous version of $u$, extended to $\bar D$.

 From  \eqref{crucial-S}, we obtain \eqref{ch1'-s7.19} and
 \beqn
\sup_{x\ne y}
\frac{|\tilde u(x)-\tilde u(y)|}
{\left({\bf\Delta}(x,y)\right)^{\alpha-\frac{Q}{p}}}
 \le Y,
\eeqn
so \eqref{ch1'-s7.20} follows from \eqref{crucial-S-1}.

Since $\alpha \in\, ]Q/p, 1[$ can be taken arbitrarily close to $1$, we obtain the statement concerning the H\"older exponents $\beta_i \in\, ]0, \alpha_i (1 - \frac{Q}{p})[$ for $\tilde u$.

Next, we argue that the constant $a$ in \eqref{ch1'-s7.19bis} is an increasing function of the domain. Indeed, recall that when $D\subset [0,1[^k$, we have $a(I,D,\alpha,p,Q) = c_3(p,\alpha,Q,k)$ (see \eqref{crucial-S-1}) and we notice that $c_3$ does not depend on the process $u$ nor on the constant $K$. In particular, once $D$ has been scaled into a set $D_1$ that fits into $[0,1[^k$, this constant $a$ does not depend on $D_1$. In order to see how it depends on $D$, translate $D$ so that it fits into $\R_+^k$. For $r > 0$, set $\phi_r(x_1,\dots, x_k) = (r^{\alpha_1^{-1}} x_1,\dots, r^{\alpha_k^{-1}} x_k)$, and choose $r > 0$ such that $D_1 = \{\phi_r(x) : x \in D \} \subset [0,1[^k$. For $y \in D_1$, define $v(y) := u(\phi_r^{-1}(y))$. In \eqref{kolmogorov(*1)}, write $K_u$ instead of $K$. Then from \eqref{kolmogorov(*1)}, we see that for any $x,y\in D_1$,
\beqn
   E[\vert v(x) - v(y) \vert^p] \leq K_u r^{-p} ({\bf \Delta}(x,y))^p.
   \eeqn
    Hence, from \eqref{crucial-S} we obtain for $x, y \in D$,
\begin{align*}
   \vert u(x) - u(y) \vert &= \vert v((\phi_r(x)) -  v(\phi_r(y)) \vert \leq Y_v ({\bf \Delta}(\phi_r(x),\phi_r(y)))^{\alpha - \frac{Q}{p}}\\
   & = Y_v  r^{\alpha - \frac{Q}{p}} ({\bf \Delta}(x,y))^{\alpha - \frac{Q}{p}}.
   \end{align*}
Therefore, we can set $Y_u = Y_v\,  r^{\alpha - \frac{Q}{p}}$, and then
   $E[Y_u^p] = r^{p \alpha - Q} E[Y_v^p] \leq r^{p (\alpha - 1) - Q} K_u\, a$.
So for $u$, we need to set $a_u = r^{p (\alpha - 1) - Q} a$, and since $p (\alpha - 1) - Q < -Q < 0$, the constant $a$ is an increasing function of the domain $D$.
\smallskip

We end the proof by checking \eqref{ch1'-s7.20bis}, so we return to the notations of the theorem. Consider \eqref{ch1'-s7.19} for the particular choice $(s,y)=(0,x)$ and remember that
$\alpha-\frac{Q}{p}
>0$. Then
\beqn
|\tilde u(s,x)|\le Ys^{\alpha_0\left(\alpha -Q/p
\right)} + |\tilde u(0,x)|.
\eeqn
This implies
\begin{align*}
&E\left[\sup_{(s,x)\in [0,t]\times D} |\tilde u(s,x)|^p\right] \notag\\
&\qquad \le 2^{p-1}\left(t^{p\alpha_0\left(\alpha -Q/p \right)} E[Y^p] + E\left[\sup_{x\in D}|\tilde u(0,x)|^p \right]\right) \notag\\
&\qquad \le 2^{p-1}K\, a(I,D,\alpha,p,Q)\, t^{p\alpha_0(\alpha-Q/p)} + 2^{p-1}C_1\notag\\
&\qquad := Kc_2(I,D,\alpha,p,Q)\, t^{p\alpha_0(\alpha-Q/p)} + 2^{p-1}C_1,
\label{ch1'-s7.b}
\end{align*}
where, in the last inequality, we have used \eqref{ch1'-s7.19bis}. This proves \eqref{ch1'-s7.20bis} and completes the proof of the theorem.
\end{proof}
\begin{remark}
\label{ch1'-s7-r100}

\noindent(a)\  Theorem \ref{ch1'-s7-t2} remains valid  if the random variables $u(t,x)$ take values in a separable Banach space $(\mathbb{B}, \Vert \cdot \Vert)$ or a complete separable metric space $(\mathbb{M}, {\rm d})$. In the statement and proof, it suffices to replace $\vert u - v \vert$ by either $\Vert u - v \Vert$ or ${\rm d}(u, v)$.   
\vskip 0.3 cm

\noindent(b)\ In applications, we often encounter situations where the estimate \eqref{ch1'-s7.18} holds for any $p>Q$ (the constant $K$ usually depends on $p$).
In this case, taking, in \eqref{ch1'-s7.19}, $\alpha$ close enough to $1$ and $p$ large enough, we deduce that the sample paths of $\tilde u=(\tilde u(t,x),\, (t,x)\in \bar I\times \bar J)$ are jointly Hölder continuous in $(t,x)=(t,(x_1,\ldots,x_k))$ with exponents $\beta_0\in\, ]0,\alpha_0[$ in the variable $t$, and $\beta_j\in\, ]0,\alpha_j[$ in the variable $x_j$,  $j=1,\ldots,k$.
\vskip 0.3 cm

\noindent(c)\ Notice that Theorem \ref{ch1'-s7-t2} covers the classical Kolmogorov continuity criterion \eqref{k1} for any $p > 0$.
 Indeed, for $k\ge 2$, we identify $[0,1]^k$ with $I \times J$, where $I = [0,1]$ and $J = [0,1]^{k-1}$, we let the generic element of $I \times J$ be $x = (x_1,\dots, x_k)$ instead of $(t, x_1,\dots, x_{k-1})$, and we shift the indices of the $\alpha_i$ from $\{0,\dots, k-1\}$ to $\{1,\dots, k\}$. We now write the right-hand side of \eqref{k1}  as $C (\vert x - y \vert^{(k+\ep)/p})^p$.
 If $(k+\ep)/p \leq 1$, then in order to apply Theorem \ref{ch1'-s7-t2}, for $i = 1,\dots, k$, we should set $\alpha_i=(k + \ep)/p$, $Q = k p / (k+\ep)$, so $p > Q$, and the common Hölder exponent given by Theorem \ref{ch1'-s7-t2} is $\beta_i < [(k + \ep)/p] [1 - Q/p] = \ep/p$.
  Of course, if $\nu := (k+\ep)/p > 1$, then \eqref{k1} implies that $\Vert Z(x) - Z(y) \Vert_{L^p(\Omega)} \leq C \vert x - y \vert^\nu$, and since $\nu > 1$, this implies that $\Vert Z(x) - Z(y) \Vert_{L^p(\Omega)} = 0$, that is, a continuous version $\tilde Z$ of $Z$ is the constant process $\tilde Z(x) := Z(x_0)$, for all $x \in [0,1]^k$, and any $x_0 \in [0,1]^k$.
\vskip 0.3 cm

\noindent(d)\  In the same vein, if \eqref{norm-aniso} is satisfied with $\alpha_{i_0} > 1$ for some $i_0 \in \{1,\dots,k\}$,  then the continuous version $\tilde u$ of $u$ satisfies the following: a.s., for any $x,\, y \in J$ with the same $i_0$-th coordinate, $\tilde u(t,x ) = \tilde u(t,y)$, that is, $\tilde u$ is a constant function of the $i_0$-th coordinate.
\vskip 0.3 cm
\end{remark}


\subsection{Regularity of Gaussian random fields}
\label{app1-3-ss2}
It is well-known that for $p > 0$, $L^p$-moments of Gaussian random variables, are determined by the $L^2$-moments. Hence, when the random field $u=\left(u(t,x),\, (t,x)\in\re_+\times D\right)$ is Gaussian, the assumption \eqref{ch1'-s7.18} of Theorem \ref{ch1'-s7-t2} can be reformulated in terms of second  moments of increments. We make this statement precise in the next theorem. The last part of the section addresses the question of sharpness of the Hölder exponents.
\begin{thm}
\label{app1-3-t1}
Let $I$, $D$, $\alpha_i$, ${\bf\Delta}$ and $Q$ be as in Theorem \ref{ch1'-s7-t2}. Suppose that $u=\left(u(t,x),\, (t,x)\in\re_+\times D\right)$ is a Gaussian random field (not necessarily centred) and for some constant $K_0<\infty$ and for all $(t,x),\, (s,y)\in I\times D$,
\beq
\label{aap1-3.1}
E\left[\left(u(t,x) - u(s,y)\right)^2\right] \le K_0 \left({\bf\Delta}(t,x;s,y)\right)^2.
\eeq
Then the inequality \eqref{ch1'-s7.18} holds for all $p>0$, with $K$ there replaced by
\beq
\label{aap1-3.100}
K_p =  2^p \left(1+\left(\frac{2^p}{\pi}\right)^\half \Gamma_E\left(\frac{p+1}{2}\right)\right) K_0^{\frac{p}{2}},
\eeq
where $\Gamma_E$ is the Gamma Euler function (see \eqref{Euler-gamma}). Therefore, all the conclusions of Theorem \ref{ch1'-s7-t2} hold (with $K$ there replaced by $K_p$). In addition, for any choice of $\beta_i\in\, ]0,\alpha_i[$, $i=0,\ldots,k$, the continuous version $\tilde u$ of $u$ is jointly H\"older-continuous on $\bar I\times \bar D$ with exponents $(\beta_0,\beta_1, \ldots, \beta_k)$.
\end{thm}
\begin{proof}
Since for all $(t,x),\, (s,y)\in I\times D$, the random variable $u(t,x) - u(s,y)$ is Gaussian, by appealing to Claim 3 of Lemma \ref{A3-l1}, we deduce from \eqref{aap1-3.1} that
\beqn
E\left[\left|u(t,x) - u(s,y)\right|^p\right] \le 2^{p}(1+c_p)\ K_0^{\frac{p}{2}} \left({\bf\Delta}(t,x;s,y)\right)^p,
\eeqn
with $c_p = \left(\frac{2^p}{\pi}\right)^ \half \Gamma_E\left(\frac{p+1}{2}\right)$.
Hence, \eqref{ch1'-s7.18} holds with $K$ given by \eqref{aap1-3.100}.

Fix $\beta_i\in\ ]0,\alpha_i[$, $i=0, \ldots, k$. Choose $p>0$ large enough and $\alpha\in\, ]\tfrac{Q}{p},1[$ close enough to $1$ so that $\beta_i<\alpha_i\left(\alpha- Q/p \right)$, for $i=0, \ldots, k$. From \eqref{ch1'-s7.19}, we conclude that $\tilde u$ is jointly H\"older-continuous on $\bar I\times \bar D$ with exponents $(\beta_0,\beta_1, \ldots, \beta_k).$
\end{proof}

The next theorem provides a sufficient condition on $L^2$-increments of a centred Gaussian process that ensures an upper bound on H\"older exponents.
\begin{thm}
\label{app1-3-t2}
Let $J\subset \re$ be a closed interval with positive length. Let $v=(v(x),\, x\in J)$ be a separable centred Gaussian process. Suppose that there is $c_0>0$ and $\alpha\in\, ]0,1]$ such that for all $x,y\in J$,
\beq
\label{app1-3.5}
E\left[\left( v(x) - v(y)\right)^2\right] \ge c_0|x-y|^{2\alpha}.
\eeq
Then $\alpha$ is an upper bound on possible H\"older exponents for $v$, that is,
for $\beta\in\, ]\alpha,1]$, a.s., the sample paths of $v$ are not H\"older-continuous with exponent $\beta$.
 \end{thm}
 
 \begin{proof}
 Suppose by contradiction that with positive probability, the sample paths of $v$ are H\"older-continuous with exponent $\beta\in\, ]\alpha,1]$. Consider the random variable $C$
 \beqn
 C=\sup_{x,\, y\, \in J,\, x\ne y}\frac{|v(x)-v(y)|}{|x-y|^\beta}
 \eeqn
 Then $P\{ C < \infty \} > 0$. Since $C$ is the supremum of the absolute values of a separable Gaussian process (indexed by $J \times J$), the zero-one law for $C$ (see  \cite{landau-shepp-1970}) implies that $P\{ C < \infty \} = 1$.

 By a classical result on centred Gaussian processes (see e.g.~\cite[Theorem 3.2 and Lemma 3.1]{adler}), this implies that
 \beqn
 K:= E\left[\sup_{x,\, \, y\in J,\, x\ne y}\frac{|v(x)-v(y)|}{|x-y|^\beta}\right] < \infty.
 \eeqn
 It follows that for all $x,y\in J$,
 \beqn
 E\left[|v(x) - v(y)|\right] \le K\ |x-y|^\beta,
 \eeqn
 and since the $L^2$-norm of a centred Gaussian random variable is proportional to the square of the $L^1$-norm (see Lemma \ref{A3-l1}), we deduce that
 \beqn
 E\left[\left(v(x)-v(y)\right)^2\right] \le K^2\ |x-y|^{2\beta}.
 \eeqn
 Since $\beta>\alpha$ and $|x-y|$ can be arbitrarily small, this contradicts \eqref{app1-3.5} and proves the theorem.
\end{proof}
\medskip

\subsection{The Garsia-Rodemich-Rumsey lemma}
\label{kolmogorov-GGR}

In this subsection, we present a general formulation of the classical Garsia, or Garsia-Rodemich-Rumsey, lemma of \cite{Garsia:72}, \cite{g-r-r}. This is an important real-variable result that has many implications in stochastic analysis and functional analysis. In particular, it provides an alternate approach to Kolmogorov's continuity theorem Theorem  \ref{ch1'-s7-t2}, and it is also useful in other topics of stochastic analysis, such as estimates of probability density functions (see for instance \cite{florit-nualart}, \cite{Pu-thesis}, \cite{dalang-pu-2020a}).  This lemma also implies various known Poincar\'e inequalities and Besov--Morrey--Sobolev embedding theorems in metric spaces.
Some of these are proved in \cite{Kassmann} using results of \cite{Arnold-Imkeller}. We also refer  to  \cite{BuckleyKoskelaI}, \cite{BuckleyKoskelaII} for work on Sobolev embedding theory.

Our formulation of the Garsia-Rodemich-Rumsey lemma is given in Lemma \ref{pr:Garsia} below, and as an application, we prove in Corollary \ref{ch1'-s7-t2-bis} below a version of Kolmogorov's continuity criterion Theorem \ref{ch1'-s7-t2}. We notice that in comparison with Theorem \ref{ch1'-s7-t2}, the domain of the space variable $x$ is more restrictive (see also Remark \ref{remgarsia}).
\medskip

\noindent{\em The lemma}
\medskip

Recall that $\Psi:\R\to\R_+$ is a \emph{strong Young function}\index{strong!Young function}\index{Young function!strong}
if it is even and convex on $\R$, and strictly increasing on
$\R_+$. Its {\em inverse} is $\Psi^{-1}: \R_+ \to \R_+$.

\begin{lemma}\label{pr:Garsia}
 Let $(S\,,\rho)$ be a metric space and $\mu$ a Radon measure on $S$.
 We assume that $0 < \mu(B_\rho(s, r)) < \infty$, for all $s \in S$ and $r > 0$, where $B_\rho(s,r)$ denotes the open ball of S in the metric $\rho$ of radius $r$.
 Let $\Psi:\R\to\R_+$ be a strong Young function
    with $\Psi(0)=0$ and $\Psi(\infty)=\infty$. Suppose that $p: \re_+
    \to \re_+$ is continuous and strictly increasing, with
    $p(0)=0$. 
  For any $\mu$-measurable function $f : S \to \re$ that is integrable over any open ball in $S$, define 
  \begin{equation}
        \mathcal{C} := \int_S\ \int_S 1_{\{x \neq y\}}\Psi\left(
        \frac{f(x)-f(y)}{p(\rho(x\,,y))}\right)\, \mu(dx)\, \mu(dy).
    \end{equation}
    \begin{description}
    \item{(a)} Suppose that $f: S\to\R$ is continuous.
    Then, for all $s,t\in S$,
    \begin{align}
    \label{eq:garsia2}
        |f(t)-f(s)|
        & \le 5 \int_0^{2 \rho(s,t)}\left[
            \Psi^{-1}\left(\frac{\mathcal{C}}{\left[ \mu(B_{\rho}(s\,,u/4)\right]^2}
            \right)\right.\notag\\
            &\left.\qquad\qquad \qquad+ \Psi^{-1}\left(\frac{\mathcal{C}}{
            \left[ \mu(B_{\rho}(t\,,u/4)\right]^2} \right)
            \right] \, dp(u).
    \end{align}
    \item{(b)}
		Suppose, instead of continuity,
		 that for $\mu$-almost all $x \in S$,
		\begin{equation}\label{cond:leb}
			\lim_{r\to 0^+} \frac{1}{\mu(B_\rho(x\,,r))}
			\int_{B_\rho(x\,,r)} f\,d\mu = f(x).
		\end{equation}
		Then there is a $\mu$-null set $N$ such that
		\eqref{eq:garsia2} holds for all $s,t \in S\setminus N$.
\end{description}
   \end{lemma}

\noindent{\em Proof of Lemma \ref{pr:Garsia}}.\
    Throughout, we choose and fix $s,t\in S$.
    We assume without loss of generality
    that $\mathcal{C}<\infty$ (otherwise, there is nothing to prove
    because $\Psi^{-1}(\infty)=\infty$).

 (a)   For any non-empty open ball $B$ in $S$, define
     \begin{equation}
        \bar{f}_B := \frac{1}{\mu(B)} \int_B f\,d\mu.
    \end{equation}
    Observe that if $A$ and $B$ are non-empty open balls in $S$,
     then
    for all $\alpha>0$,
    \begin{equation}\begin{split}
    	\Psi\left( \frac{\bar{f}_A-\bar{f}_B}{\alpha}\right)
		&=\Psi\left(\frac{1}{\mu(A)\, \mu(B)}\int_A\mu(dx)\int_B\mu(dy)\
			\frac{f(x)-f(y)}{\alpha}\right)\\
		&\le \frac{1}{\mu(A)\, \mu(B)}\int_A\mu(dx)\int_B\mu(dy)
			\Psi\left( \frac{f(x)-f(y)}{\alpha}\right),
    \end{split}\end{equation}
  by Jensen's inequality.  
    It follows from this that if $\alpha\ge p(\rho(x\,,y))$
    for all $x\in A$ and $y\in B$, then
    $\Psi((\bar{f}_A-\bar{f}_B)/\alpha)$ is bounded above
    by $\mathcal{C}/(\mu(A)\, \mu(B))$. Thus, we are led to
    the basic inequality
    \beqn
    \left| \bar{f}_A-\bar{f}_B\right| \le \Psi^{-1}\left(\frac{\mathcal{C}}{\mu(A)\, \mu(B)} \right)\alpha,
    \eeqn
    and taking the infimum over $\alpha\ge p(\rho(x\,,y))$, we have
		 \begin{equation}\label{eq:new:star}
    	\left| \bar{f}_A-\bar{f}_B\right| \le \Psi^{-1}\left(\frac{%
		\mathcal{C}}{\mu(A)\, \mu(B)} \right)\sup_{x\in A,\, y\in B}
		p(\rho(x\,,y)).
    \end{equation}

    Fix $s,t\in S$ with $s\ne t$ (otherwise, \eqref{eq:garsia2} is obvious) and
let $r_0:=\frac12\rho(s\,,t)$.
 Define $r_n$ by
    $p(2r_n) := 2^{-n} p(2r_0)$ for all $n\ge 1$.
    Notice that as $n \to \infty$,  both $r_n$
    and $p(2r_n)$ decrease to $0$.

For $n\ge 0$, define $A_n:= B_{\rho}(s\,,r_n)$,
    $B_n:=B_\rho(t\,,r_n)$ and apply \eqref{eq:new:star} to find that
    \begin{equation}
        \left| \bar{f}_{A_n} - \bar{f}_{A_{n-1}} \right|
        \le \Psi^{-1}\left( \frac{\mathcal{C}}{\left[ \mu(A_n) \right]^2}\right)
        p(2r_{n-1}).
    \end{equation}
    Because $p(2r_n)-p(2r_{n+1})= \frac14 p(2r_{n-1})$,
    \begin{equation}
    \label{aquatresetze}
        \left| \bar{f}_{A_n} - \bar{f}_{A_{n-1}} \right|
        \le 4\Psi^{-1}\left( \frac{\mathcal{C}}{\left[ \mu(A_n) \right]^2}\right)
        \left[ p(2r_n)-p(2r_{n+1})\right].
    \end{equation}
    Note that $\cap_{n=1}^\infty A_n=\{s\}$, whence
    \beq
    \label{A.3-kc}\lim_{n\to\infty}\bar{f}_{A_n}= f(s)
    \eeq
     by continuity. Therefore,
    we can take the sum over all $n\ge 1$ in \eqref{aquatresetze} 
    to find that
    \begin{equation}\label{eq:new:1}
			\left| f(s) - \bar{f}_{A_0} \right|
			\le 4\int_0^{2r_0} \Psi^{-1}\left(
			\frac{\mathcal{C}}{\left[ \mu(B_{\rho}(s\,,u/2))\right]^2
			}\right)\, dp(u).
    \end{equation}
    Similarly,
    \begin{equation}\label{eq:new:2}
			\left| f(t) - \bar{f}_{B_0} \right|
			\le 4\int_0^{2r_0} \Psi^{-1}\left(
			\frac{\mathcal{C}}{\left[ \mu(B_{\rho}(t\,,u/2))\right]^2
			}\right)\, dp(u).
    \end{equation}
    A new application of \eqref{eq:new:star} reveals that
    $| \bar{f}_{A_0} - \bar{f}_{B_0} |$ is at most
    \begin{align}
    \label{eq:new:3}
    	&\Psi^{-1}\left( \frac{\mathcal{C}}{\mu(A_0)\,
			\mu(B_0)}\right) p(4r_0) \notag\\
		&\quad\le \int_0^{4r_0}\Psi^{-1}\left( \frac{\mathcal{C}}{\mu(B_\rho(s\,,u/4))\,
			\mu(B_\rho(t\,,u/4))}\right) dp(u)\notag\\
		&\quad\le \int_0^{4r_0}\Psi^{-1}\left( \frac{\mathcal{C}}{
			\left[\mu(B_\rho(s\,,u/4))\right]^2}\right) dp(u)\notag\\
			&\qquad\qquad\qquad+ \int_0^{4r_0}\Psi^{-1}\left( \frac{\mathcal{C}}{
			\left[\mu(B_\rho(t\,,u/4))\right]^2}\right) dp(u).
    \end{align}
    Since $r_0=\frac12\rho(s\,,t)$, equations \eqref{eq:new:1},
    \eqref{eq:new:2}, and \eqref{eq:new:3} together, along with the triangle inequality, imply \eqref{eq:garsia2}.
    \smallskip

    (b) If, instead of continuity, we assume that  \eqref{cond:leb} holds, then \eqref{A.3-kc} remains valid for $\mu$-a.a.~$s$, and the remainder of the proof is unchanged.
\qed
\medskip

\noindent{\em Application to the Kolmogorov continuity theorem}
\medskip

We recall that $u=(u(t,x),\, (t,x)\in \IR_+ \times D)$ is a random field, and $D\subset \re^k$ is a bounded or unbounded domain. Without making use of the results of Subsections \ref{app1-3-ss1} and \ref{app1-3-ss2}, we establish the following corollary, which is a special case of Theorem \ref{ch1'-s7-t2}, via an alternate proof that only uses Lemma \ref{pr:Garsia} and some technical facts.
\begin{cor}
\label{ch1'-s7-t2-bis}
Let $I$ be a bounded open interval of $\re_+$ and $J\subset D$ a Cartesian product of bounded open intervals of $\re$, all with positive length.
Let ${\bf\Delta}(t,x;s,y)$ be as defined in \eqref{norm-aniso}.

Assume that the assumptions of Theorem \ref{ch1'-s7-t2} are satisfied with $D$ there replaced by $J$. Then all the conclusions of Theorem \ref{ch1'-s7-t2} hold with $D$ (respectively $\bar D$) replaced by $J$ (respectively $\bar J$).
\end{cor}

For the proof of the Corollary \ref{ch1'-s7-t2-bis}, we will use the following technical result.
\begin{lemma}
\label{ch1'-s7-l3}
Let $b\in \re$, $\alpha_i \in \,]0,1]$, $i=0,\dots,k$. Then
\beqn
\mathcal{I} := \int_0^1 dt \int_{-1}^1 dx_1\, \cdots \int_{-1}^1 dx_k  \left(t^{\alpha_0}+|x_1|^{\alpha_1} + \cdots + |x_k|^{\alpha_k} \right)^b <\infty
\eeqn
if and only if $b+Q
>0$.
\end{lemma}
\begin{proof}
Use the change of variables $r_0=t^{\alpha_0}$, $r_1= x_1^{\alpha_1},\dots, r_k= x_k^{\alpha_k}$ to see that $\mathcal{I}$ is equal to
\beq
\label{ch1'-s7.21}
c_k \int_0^1 dr_0\, r_0^{\frac{1}{\alpha_0}-1} \int_0^1 dr_1\, r_1^{\frac{1}{\alpha_1}-1}  \cdots
\int_0^1 dr_k\, r_k^{\frac{1}{\alpha_k}-1}   (r_0+r_1+\cdots +r_k)^b,
\eeq
with $c_k= 2^{k} (\prod_{i=0}^k\alpha_i)^{-1}$.

Since $\frac{1}{\alpha_i}-1\ge 0$, $i=0,\dots,k$, this is bounded above by
\beqn
c_k \int_0^1 dr_0\, \cdots \int_0^1 dr_k\, (r_0+r_1+\cdots +r_k)^{b+
\sum_{i=0}^k\frac{1}{\alpha_i}-k-1}.
\eeqn
Pass to polar coordinates relative to $(r_0,\dots,r_k)$ to bound this by
\beqn
\tilde c_k\int_0^{k+1}\ d\rho\  \rho^{b+Q
-1},
\eeqn
and this is finite provided that $b+Q
>0$.

For the converse implication, we assume $b+Q\le 0$ and prove that $\mathcal{I}=\infty$. Indeed, 
 the multiple integral in \eqref{ch1'-s7.21} is bounded below, up to a multiplicative constant, by
\beqn
 \int_0^1 dr_0\, r_0^{\frac{1}{\alpha_0}-1}\int_0^{r_0} dr_1\, r_1^{\frac{1}{\alpha_1}-1} \cdots \int_0^{r_{k-1}} dr_k\, r_k^{\frac{1}{\alpha_k}-1}  (r_0+r_1+\cdots +r_k)^b.
\eeqn
Observe that for $(r_0,r_1, \ldots, r_k)$ in the domain of integration, we have $r_0+\cdots+r_k\le (k+1)r_0$ and since $b<0$, the integral above is bounded 
 below by
\begin{align*}
  (k+1)^b\int_0^{1} dr_0 \, r_0^{b+\frac{1}{\alpha_0}-1}\! &\int_0^{r_0}\!\! dr_1 \, r_1^{\frac{1}{\alpha_1}-1}\! \cdots  \int_0^{r_{k-1}}\!\! \! dr_{k}\,  r_{k}^{\frac{1}{\alpha_{k}}-1}\\
 &\quad= c \int_0^1 dr_0 \, r_0^{b+Q-1} = \infty,
 \end{align*}
 since $b+Q\le 0$.
\end{proof}
\medskip
\medskip

\noindent{\em Proof of Corollary \ref{ch1'-s7-t2-bis}.}
Fix numbers $p>Q$ such that \eqref{ch1'-s7.18} holds and  $\alpha\in\,]Q/p,1[$.
We observe that \eqref{ch1'-s7.18} implies that $(t,x)\mapsto u(t,x)$ is continuous in probability and therefore, by  Theorem \ref{rdsecA.4-t2}, has a measurable version that we continue to denote by $u$. We note that by \eqref{ch1'-s7.18}, for any fixed $(t_0,x_0)\in I\times J$,
$u(\cdot,\ast) - u(t_0,x_0)\in L^p_{\text{loc}}(dt,dx)$, a.s. Since this shifted process has the same increments as $u$, we can and do assume that $u\in L^p_{\text{loc}}(dt,dx)$, a.s.

We will apply Lemma \ref{pr:Garsia} with $S = I\times J$,
\beqn
\rho(t,x;s,y)= \max(|t-s|^{\alpha_0}, |x_1-y_1|^{\alpha_1}, \ldots, |x_k-y_k|^{\alpha_k}),
\eeqn
so that
\beq
\label{A.3-ne}
\frac{1}{k+1} {\bf\Delta}(t,x;s,y)\le \rho(t,x;s,y)\le {\bf\Delta}(t,x;s,y),
\eeq
$\mu(dt,dx) = dt\, dx_1\cdots dx_k$, $\Psi(w)=|w|^p$ and $p(z)= z^{\alpha+Q/p}$.
With these choices, notice that the hypotheses of Lemma \ref{pr:Garsia} are satisfied. In addition, since $\{(t,x; s,y): (t,x) = (s,y)\}$ has $dtdxdsdy$-measure $0$,
we consider the random variable
\beqn
\mathcal{C}:= \int_I dt \int_J dx \int_I ds \int_J dy\ \frac{|u(t,x)-u(s,y)|^p}{\left(\rho(t,x;s,y)\right)^{p\alpha+Q
}},
\eeqn
and check that $E(\mathcal{C}) < \infty$.

Indeed, assumption \eqref{ch1'-s7.18} and \eqref{A.3-ne} imply that
\beqn
E[\mathcal{C}] \le K (k+1)^p \int_I dt \int_J dx \int_I ds \int_J dy \left(\rho(t,x;s,y)\right)^{p(1-\alpha)-Q}.
\eeqn
Let
\begin{align*}
a(I,J,\alpha,p, Q):&= (k+1)^p \\
&\quad\times\int_I dt \int_J dx \int_I ds \int_J dy \left(\rho(t,x;s,y)\right)^{p(1-\alpha)-Q}.
\end{align*}
Since $p(1-\alpha)>0$, and $I$ and $J$ are bounded, Lemma \ref{ch1'-s7-l3} and \eqref{A.3-ne} above imply that $a(I,J,\alpha,p, Q)$ is finite. Thus,
\beq
\label{ch1'-s7.a}
E[\mathcal{C}] \le K a(I,J,\alpha,p, Q) < \infty,
\eeq
 and clearly, $a(I,J,\alpha,p, Q)$ is increasing  in the first two arguments.

 Since $u\in L^p_{\rm{loc}}(dt dx_1\cdots dx_k)$ a.s. and because $p>1$, a well-known theorem of Jessen, Marcinkiewicz and Zygmund
 (see \cite[Theorem 2.2.1, p. 58]{khosh2002}) implies that a.s., for $\mu$-almost all $(t,x)\in I\times J$ (recall that $\mu$ is Lebesgue measure),
 \begin{align*}
 &\lim_{\delta_i\downarrow 0, i=0,\dots, k} \, \frac{1}{2^{k+1}\delta_0\cdots \delta_k}
 \int_{t-\delta_0}^{t+\delta_0} ds \int_{x_1-\delta_1}^{x_1+\delta_1} dy_1\ \cdots \int_{x_k-\delta_k}^{x_k+\delta_k}  u(s,y_1,\ldots,y_k)\\
  &\qquad\qquad = u(t,x_1,\ldots,x_k).
 \end{align*}
 In particular, since balls in the metric $\rho$ are Cartesian products, condition \eqref{cond:leb} holds. We deduce from Proposition \ref{pr:Garsia} that for a.a.~$\omega\in \Omega$, there is a set $D(\omega)\subset I\times J$ with full Lebesgue measure such that for all $(t,x), (s,y)\in D(\omega)$,
 \begin{align*}
 |u(t,x)(\omega) - u(s,y)(\omega)|
 & \le 10 \left(\alpha + \frac{Q}{p}\right) \sup_{(r,z)\in I\times J} \int_0^{2{\bf\Delta}(t,x;s,y)} dw\ w^{\alpha-1+Q/p}\\
 &\qquad\qquad\qquad\qquad \times \Psi^{-1}\left(\frac{\mathcal{C}}{\left[\mu\left(B_\rho((r,z),w/4)\right)\right]^2}\right).
 \end{align*}
 One checks immediately that since $I\times J$ is a box and $\mu$ is Lebesgue measure on $\R \times \R^k$, there is $c>0$ such that for all $w\in [0,2{\bf\Delta}(t,x;s,y)]$ and  $(r,z)\in I\times J$,
 \begin{equation}\label{rd04_10e1}
 \mu\left(B_\rho((r,z),w/4)\right)\ge c\, w^Q.
 \end{equation}
 Therefore,
 \begin{align}
 \label{A.3-*3}
 & |u(t,x)(\omega) - u(s,y)(\omega)| \le 10 \left(\alpha + \frac{Q}{p}\right) c^{-2/p} \int_0^{2{\bf\Delta}(t,x;s,y)} dw\, w^{\alpha-1+Q/p}\frac{\mathcal{C}^{1/p}}{w^{2Q/p}} \notag\\
  &\qquad = 10 \left(\alpha + \frac{Q}{p}\right)  c^{-2/p} \left(\alpha - \frac{Q}{p}\right)^{-1} \mathcal{C}^{1/p} \, [2 {\bf\Delta}(t,x;s,y)]^{\alpha-Q/p} \notag\\
  &\qquad =c_0\   \mathcal{C}^{1/p} [{\bf\Delta}(t,x;s,y)]^{\alpha-Q/p},
  \end{align}
  where $c_0 = 10(\alpha+Q/p) c^{-2/p} 2^{\alpha-Q/p}/(\alpha-Q/p)$ and we have used the hypothesis $\alpha > Q/p$.

  For $(t,x)\in\bar D = \bar S$, define
  \beqn
  \tilde u(t,x)(\omega):= \limsup_{(s,y)\in D(\omega): (s,y)\to (t,x)} u(s,y)(\omega).
  \eeqn
  Since $u(\cdot,\ast)(\omega)$ is uniformly continuous on $D(\omega)$ by \eqref{A.3-*3}, the $\limsup$ is in fact a limit and $\tilde u(\cdot,\ast)(\omega)$ is continuous on $\bar S$. 
  In addition, by \eqref{ch1'-s7.18}, $(s,y) \mapsto u(s,y)$ is uniformly continuous in $L^p(\Omega)$, therefore 
 $u(s,y)$ converges in $L^p(\Omega)$ to $u(t,x)$ as $(s,y) \longrightarrow (t,x)$ in $\bar S$. Therefore, for all $(t,x)\in S$, $\tilde u(t,x) = u(t,x)$ a.s., that is, $\tilde u$ is a continuous version of $u$, extended to $\bar S$. By
 \eqref{A.3-*3}, \eqref{ch1'-s7.19} holds, and
 \beqn
  \left(\sup_{(t,x)\ne(s,y)}
\frac{|\tilde u(t,x)-\tilde u(s,y)|}
{\left({\bf\Delta}(t,x;s,y)\right)^{\alpha-\frac{Q}{p}}}
\right)^p \le c_0^p\, \mathcal{C},
\eeqn
so \eqref{ch1'-s7.20} follows from \eqref{ch1'-s7.a}. Setting  $Y=\mathcal{C}^{1/p}$, property \eqref{ch1'-s7.19bis} is just \eqref{ch1'-s7.a}.
Since $\alpha \in ]Q/p, 1[$ can be taken arbitrarily close to $1$, we obtain the statement concerning the H\"older exponents $\beta_i \in ]0, \alpha_i (1 - Q/p)[$ for $\tilde u$.

Finally, we prove \eqref{ch1'-s7.20bis}. Consider \eqref{ch1'-s7.19} for the particular choice $(s,y)=(0,x)$ and recall that
$\alpha - Q / p >0$. Then
\beqn
|\tilde u(t,x)|\le c_0\ Yt^{\alpha_0\left(\alpha -Q/p
\right)} + |\tilde u(0,x)|.
\eeqn
This implies that
\begin{align}
&E\left[\sup_{(s,x)\in [0,t]\times J} |\tilde u(s,x)|^p\right] \notag\\
&\qquad \le 2^{p-1}\left(c_0^p\, t^{p\alpha_0\left(\alpha -Q/p \right)}E[Y^p] + E\left[\sup_{x\in J}|\tilde u(0,x)|^p\right]\right) \notag\\
&\qquad \le c_2(I,J,\alpha,p,Q)\, t^{p\alpha_0(\alpha-Q/p)} K+2^{p-1}C_1,
\label{ch1'-s7.b2}
\end{align}
where in the last inequality, we have used \eqref{ch1'-s7.19bis}. This proves \eqref{ch1'-s7.20bis} and completes the proof of the corollary.
\qed

\begin{remark}\label{remgarsia}
In our proof of Corollary \ref{ch1'-s7-t2-bis}, we have used the lower bound \eqref{rd04_10e1}. This lower bound clearly holds when $I \times J$ is a box. However, if we are interested in more general bounded open sets $J$ and want to establish the Kolmogorov continuity criterion as a corollary of the Garsia-Rodemich-Rumsey inequality, then it must be noted that $B_\rho((r, z), w/4)$ is not a ball in $\R_+ \times \R^k$, but a ball in the set $\cO := I \times D$. This is why Da Prato and Zabczyk, who work with the Euclidean metric $\rho$, the Lebesgue measure $\mu$ and $Q = k+1$, and who let $B_\rho((r, z), w)$ denote a Euclidean ball in $\R_+ \times \R^k$, introduce the parameter
\beqn
   \kappa(\cO) := \inf_{w\in\,]0,\, {\rm diam}(\cO)[}\ \inf_{(r,\, x) \in \cO} \ \frac{\mu(\cO \cap B_\rho((r, x), w))}{w^{Q}},
\eeqn
where ${\rm diam}(\cO) = \sup_{(t,\, x),\, (s,\, y) \in \cO} \rho((t, x),\, (s, y))$ (see \cite[(B.1.1), p.~311]{dz-1996}), and the geometric condition 
\begin{equation}\label{rd04_26e1}
   \kappa(\cO) > 0
\end{equation} 
(see \cite[Theorem B.1.5, p.~316]{dz-1996}; this condition is not mentioned in \cite[Theorem C.6]{khosh}). In the case $ k = 1$ and $\cO \subset \R_+ \times \R$, for instance, this geometric condition is satisfied for bounded domains with polygonal space-time boundaries, but in general not for domains with cusps. If, in Corollary \ref{ch1'-s7-t2-bis}, we assume that $D$ is bounded and $\cO := I\times D$ satisfies \eqref{rd04_26e1}, then the conclusion and the proof of Corollary \ref{ch1'-s7-t2-bis} remain valid with $J$ there replaced by $D$.
\end{remark}

\section{Measurability of the It\^o integral}
\label{app1-4}

In this section, we discuss the joint measurability of the stochastic integral with respect to a Brownian motion when the integrand process depends on a parameter. The main result is Theorem \ref{rdthm9.4.1} below. It is applied in particular in Section \ref{ch2'new-s6}, where a similar measurability question is addressed for the stochastic integral with respect to space-time white noise.

Throughout this section, $(\Omega, \cF, P)$ is a complete probability space that we equip with a complete and right-continuous filtration $(\cF_t,\, t \in \R_+)$,
$(X,\cX)$ is a measure space, $(B_t,\, t\in\re_+)$ is an $(\cF_t)${\em-standard Brownian motion},
and $\cO$ is the optional $\sigma$-field defined in Section \ref{rdsecA.4}.

 \begin{thm} \label{rdthm9.4.1} 
 Fix $T>0$ and let $Z:X \times [0,T] \times \Omega \to \IR$ be a function satisfying the following conditions:
\begin{enumerate}
\item $Z$ is  $\cX \times \cB_{[0,T]}\times \cF$-measurable.
\item For all $(x,s) \in X\times [0,T]$, the map $\omega \mapsto Z(x,s,\omega)$ is $\cF_s$-measurable.
\end{enumerate}
Furthermore, suppose that for all $x \in X$,
\beq
\label{rde9.4.1}
   \int_0^T Z^2(x,s) \, ds < +\infty,\qquad a.s.
\eeq
Then there is an $\cX\times \cO$-measurable function $Y : X \times [0,T] \times \Omega \to \R$ such that, for all $x \in X$, $(t,\omega) \mapsto Y(x,t,\omega)$ is a.s.~continuous
and the process $Y(x, \cdot)$ is indistinguishable from the stochastic integral process $Z(x, \cdot)\cdot B.$ Further, for fixed  $t \in [0, T]$, the map $(x, \omega) \mapsto Y(x,t, \omega)$ is $\cX \times \cF_t$-measurable.
 \end{thm}

 The proof of this theorem is based on the next two lemmas.

 \begin{lemma}
\label{rdlem9.4.2} 
Fix $T>0$ and let $Z:X \times [0,T] \times \Omega \to \IR$ be a function satisfying conditions 1.~and 2.~of Theorem \ref{rdthm9.4.1}.
Then:

(a)\  There is an $\cX\times \cO$--measurable function $K: X \times [0,T] \times \Omega \to \IR$ such that for all $(x,s) \in X \times [0,T]$, $Z(x,s) = K(x,s)$ a.s. In particular, for each $x \in X$,
\beq
\label{rde9.4.3}
   E\left[\int_0^T (Z(x,s) - K(x,s))^2\, ds \right] = 0.
\eeq

(b)\ Suppose in addition that for all $x \in X$, $Z$ satisfies \eqref{rde9.4.1}.
Then for all $x\in X$, the  processes $Z(x,\cdot)\cdot B$ and $K(x,\cdot)\cdot B$ are indistinguishable.
\end{lemma}

\begin{proof}
(a)\  By composing $Z$ with a bijection from $\re$ into $]0,1[$, we can assume that $Z$ is positive and bounded. By \cite[Proposition 3]{strickeryor}, there is an $\cX \times \cO$--measurable function $K: X \times [0,T] \times \Omega \to \IR$ such that, for each $x \in X$, $K(x,\cdot)$ is a version of the optional projection of $Z(x,\cdot)$ (the definition of this notion is recalled in Section \ref{rdsecA.4}). 
Fix $(x,s)\in X \times [0,T]$. Since $Z(x,s)$ is $\cF_s$--measurable, by the definition of optional projection,
$$
    K(x,s) = E[ Z(x,s) \mid \cF_s] = Z(x,s) \qquad a.s.
$$
This proves the first part of statement (a) in the lemma.

Since $Z$ and $K$ are jointly measurable in $(x,s,\omega)$, for each $x\in X$,
$$
   A(x) = \{(s,\omega) \in [0,T]\times \Omega: Z(x,s,\omega) \neq K(x,s,\omega) \}
$$
belongs to $\cB_{[0,T]}\times \cF_T$. By Fubini's theorem,
$(ds \times dP)(A(x)) = 0$,
that is,
\beqn
   E\left[\int_0^T (Z(x,s) - K(x,s))^2\, ds \right] = 0.
\eeqn
Therefore, \eqref{rde9.4.3} holds.
\smallskip

(b)\ Notice that if \eqref{rde9.4.1}
 holds, then $K$ also satisfies \eqref{rde9.4.1}, so that the stochastic integral processes $Z(x,\cdot) \cdot B$ and $K(x,\cdot) \cdot B$ are well-defined continuous and adapted processes. By \eqref{rde9.4.3} and the construction of the stochastic integral with respect to Brownian motion, for each $t \in [0,T]$, the random variables $(Z(x,\cdot) \cdot B)_t$ and $(K(x,\cdot) \cdot B)_t$ are equal, a.s.  Since both processes are continuous, they are indistinguishable.
\end{proof}

\begin{remark}
\label{after-rdprop2.6.2}
Let  $G = (G(t,x),\, (t,x) \in [0,T] \times D)$ be a jointly measurable and adapted process (see conditions (1) and (2) at the beginning of Section \ref{ch2new-s2}). From Lemma \ref{rdlem9.4.2} (a), we deduce that $G$ has an optional version $\tilde G$. Here, {\em optional} means that $(x,t,\omega) \mapsto \tilde G(t,x,\omega)$ is $\B_D \times \cO$-measurable. Indeed, Let $Z := G$ and $X := D$ in Lemma \ref{rdlem9.4.2} (a). Condition 1.~of Theorem \ref{rdthm9.4.1} is joint measurability and condition 2.~there is weaker than {\em adapted} (that is, for fixed $s$, $(x,\omega) \mapsto G(s,x,\omega)$ is $\B_D \times \F_s$-measurable). The $\B_D \times \cO$-measurable function $K$ given by Lemma \ref{rdlem9.4.2} is an optional version of $G$.
\end{remark}

\begin{lemma}
\label{2ofsectionA.4} 
We assume that $Z$ satisfies assumptions 1.~and 2.~of Theorem \ref{rdthm9.4.1}, and that for all
$x \in X$,
\beq
\label{*1ofsectionA.4}
  E\left[\int_0^T Z^2(x,s) ds\right] < \infty.
  \eeq
Then there exists an $\cX \times \cO$-measurable function $Y : X \times [0,T] \times \Omega \to \R$ such that for all $x \in X$,  the mapping $(t,\omega) \mapsto Y(x,t,\omega)$ is a.s.~continuous,
 $Y(x)$ and the stochastic integral process $Z(x,\cdot) \cdot B$ are indistinguishable, and for fixed $t \in [0, T]$, $(x, \omega) \mapsto Y(x,t, \omega)$ is $\cX \times \cF_t$-measurable.
 \end{lemma}

 \begin{proof}
 Let $K$ be defined as in Lemma \ref{rdlem9.4.2} (a). By \eqref{rde9.4.3}, $K$ satisfies \eqref{*1ofsectionA.4}. By \cite[Proposition 5]{strickeryor} (notice that $K$ satisfies the hypothesis of this proposition), there exists an $\cX \times \cO$-measurable function $Y : X \times [0,T] \times \Omega \to \R$ such that, for all
 $x \in X$, $Y(x)$ is indistinguishable from the continuous adapted integral process $K(x,\cdot) \cdot B$,
and for fixed $t \in [0, T]$, since  $\cO\vert_{[0, t] \times \Omega} \subset \cB_{[0, t]} \times \cF_t,$ the function  $(x, \omega) \mapsto Y(x,t, \omega)$ is $\cX \times \cF_t$-measurable.   Since by Lemma \ref{rdlem9.4.2} (b), $K(x, \cdot) \cdot B$ is indistinguishable from $Z(x,\cdot) \cdot B$, this is the process $Y$ of the assertion.
 \end{proof}

 \noindent{\em Proof of Theorem \ref{rdthm9.4.1}}.\
 Let $K(x,s)$ be as given in  Lemma \ref{rdlem9.4.2}. From \eqref{rde9.4.3}, we see that $K$ also satisfies
 \eqref{rde9.4.1}.

 For $N\in\N$ and  $x \in X$, let
 \beqn
   \tau_N(x) = \inf\left\{r \in [0,T]: \int_0^r K^2(x,s)\, ds \geq N\right\} \wedge T.
   \eeqn
Then  $(x,\omega) \mapsto \tau_N(x,\omega)$ is $\cX \times \cF_T$-measurable, and $\tau_N(x)\uparrow T$ a.s.~as $N\to\infty$. Moreover, for any $ x \in X$,  $\tau_N(x)$ is a stopping time. Indeed, for $ t \in [0,T[$,
 $\{\tau_N(x) \leq t\} = \{ \int_0^t K^2(x,s)\, ds \geq N\}$, and since $K(x,\cdot)$ is optional (thus, progressively measurable), this event is $\F_t$-measurable.

We observe that for $x \in X$ fixed, because of \eqref{rde9.4.3}, $\tau_N$ can also be expressed in terms of $Z$:
\beqn
  \tau_N(x) = \inf\left\{r \in [0,T]: \int_0^r Z^2(x,s)\, ds \geq N\right\} \wedge T,\qquad    a.s.,
  \eeqn
  where the null set depends on $x$.

  Define $K_N(x) = 1_{[0,\tau_N(x)]}(s) K(x,s)$, $Z_N(x) = 1_{[0,\tau_N(x)]}(s) Z(x,s)$. Then $Z_N(x)$ satisfies the conditions 1.~and 2.~of Theorem \ref{rdthm9.4.1}, and
  \beqn
  E\left[\int_0^T Z_N^2(x,s)\, ds\right] + E\left[\int_0^T K_N^2(x,s)\, ds\right] \le 2N.
  \eeqn

  Applying Lemma \ref{2ofsectionA.4} to $Z_N$, we see that there exists  an $ \cX \times \cO$-measurable function $Y_N : X \times [0,T] \times \Omega \to \R$ such that for all $x \in X$, the process $(t,\omega) \mapsto Y_N(x,t,\omega)$ has continuous sample paths a.s., $Y_N(x)$ and  $Z_N(x,\cdot) \cdot B$ are indistinguishable, and for fixed $t \in [0, T]$, $(x, \omega) \mapsto Y(x,t, \omega)$ is $\cX \times \cF_t$-measurable.

 On the event $\{\tau_N(x) = T\}$, a.s., the sample paths of the continuous processes $Z_N(x,\cdot) \cdot B$ and $Z(x,\cdot) \cdot B$ are identical. Thus, a.s, the sample paths of the continuous processes $Y_N(x)$ and $Z(x,\cdot) \cdot B$ are identical. Let $Y(x,t,\omega) = \limsup_{N\to\infty}Y_N(x,t,\omega)$. Then $Y$ is  $\cX \times \cO$-measurable;
  it is also adapted and, for all $x\in X$,  a.s.~with continuous sample paths, because of the stationary convergence of
 $Y_N(x,\cdot)$ to $Z(x,\cdot)\cdot B$ on $\{\tau_N(x) = T\}$.
 Furthermore, for all $x\in X$, $Y$ is indistinguishable from $Z(x,\cdot) \cdot B$, since a.s., the trajectories are identical on the event $\{\tau_N(x) = T\}$, and a.s., $\{\tau_N(x) = T\}\uparrow_{N\to\infty} \Omega$. The proof of Theorem \ref{rdthm9.4.1}  is complete.
  \qed
\bigskip


We end this section with a restatement, in the notations used in this book, of a result from \cite{strickeryor}.

\begin{lemma}\label{rdlem9.4.2a}
Let $\cG$ be a complete sub-$\sigma$-field of $\cF$
and let $(\mathbb{B}, \Vert \cdot \Vert)$ be a separable Banach space.
 Consider a sequence $(Y_n,\, n\in \N)$ of $\cX \times \cG$-measurable functions from $X\times \Omega$ to $\mathbb{B}$. Suppose that for all $x \in X$, the sequence $(Y_n(x),\, n\in\N)$ converges in probability on $\Omega$. Then there exists an $\cX \times \cG$-measurable function $Y:X\times \Omega \to \IR$ such that, for all $x\in X$, $Y(x) = \lim_{n\to \infty} Y_n(x)$ in probability.
\end{lemma}

For real-valued functions, this lemma is Proposition 1 of \cite{strickeryor}. The proof for functions with values in a separable Banach space is identical. This lemma is used in the proof of Theorem \ref{ch1'-s5.t1}. 
\medskip


\section{Stochastic Fubini's theorem for Brownian motion}
\label{app1-2}

As in classical (deterministic) calculus, the stochastic Fubini's theorem is a fundamental tool in stochastic analysis. There are many versions of this result that depend on the type of stochastic integral and the integrator process.
Section \ref{ch2'new-s4} contains a Fubini's theorem for stochastic integrals with respect to space-time white noise. Its proof relies on a specific Fubini's theorem in the simpler setting where the integrator is Brownian motion. We formulate this statement here.

Throughout this section,  $(X,\cX, \mu)$ is a measure space such that the measure $\mu$ is $\sigma$-finite, and $B=(B_t,\, t\in \re_+)$ is an $\left(\tf_t\right)$-standard Brownian motion.
\begin{thm}
\label{app1-2.t1}
Let $g: X\times [0,T]\times \Omega \to \re$ be jointly measurable and such that for each $t \in [0,T]$, $(x,\omega) \mapsto g(x,t,\omega)$ is $\cX \times \cF_t$-measurable. Suppose that
\beq
\label{condition-fubini}
\int_X \mu(dx) \left(\int_0^T dt\, g^2(x, t)\right)^\half <\infty,\qquad  {\text{a.s.}}
\eeq
Then:
\begin{description}
\item{(a)} There exists $X_0\in\mathcal{X}$ with $\mu(X\setminus X_0)=0$ such that for all $x\in X_0$, $\int_0^T g^2(x, t)\, dt < \infty$ a.s., and there exists a jointly measurable in $(x,t,\omega)$ map
$\Psi: X\times [0,T]\times \Omega \to \re$, such that for all $x\in X_0$, $\Psi(x,\cdot)$ and
 the stochastic integral process $g(x, \cdot) \cdot B$ are indistinguishable, and, in addition,
\beqn
\sup_{t\in[0,T]} \int_X |\Psi(x, t)|\, \mu(dx) < \infty.
\eeqn
In particular, for any $x\in X_0$, $\Psi(x,\cdot)$ has continuous sample paths a.s.
\item{(b)} For almost all $(t,\omega)\in[0,T]\times \Omega$, the function $x\mapsto g(x, t, \omega)$ is $\mu$-integrable and
\beqn
\int_0^T dt \left(\int_X \mu(dx)\, g(x, t)\right) ^2 < \infty,\quad {\text{a.s.}}
\eeqn
\item{(c)} Almost surely, for all $t\in[0,T]$,
\beqn
\int_X \mu(dx) \left(\int_0^t g(x,s)\, dB_s\right) = \int_0^t \left(\int_X\mu(dx)\,  g(x, s)\right) dB_s
\eeqn
$($by definition, the left-hand side is $\int_{X_0} \mu(dx)\, \Psi(x,t))$.
 In particular, a.s., the left-hand side is continuous in $t$.
 \end{description}
\end{thm}

\begin{proof}
This theorem is a slightly modified version of
\cite[Theorem 2.2]{V2012}, where stochastic integrals are with respect to a continuous semimartingale (instead of a Brownian motion), and, for any $x\in X$,  the process $g(x,\cdot)$ is progressively measurable (see also \cite[Lemma 2.6]{krylov-2011} for a related formulation of the theorem).

We give some details on the proof of (a). For the existence of the $\mu(dx)$-null set $X_0$, we notice that the condition \eqref{condition-fubini} implies the following: There exists a $dP$-null set $F_0\in\cF$ such that for all $\omega\notin F_0$,
\beqn
\int_X \mu(dx) \left(\int_0^T dt\, g^2(x, t,\omega)\right)^\half <\infty.
\eeqn
 Therefore, for $\omega\notin F_0$, there is a $\mu(dx)$-null set $X_1(\omega)$ such that for $x\notin X_1(\omega)$,
$\Vert g(x,\cdot,\omega)\Vert_{L^2([0,T])} < \infty$. Since
\beqn
\{(x,\omega)\in X\times \Omega: \Vert g(x,\cdot,\omega)\Vert_{L^2([0,T])}=\infty\} \in \cX \times \cF,
\eeqn
the discussion above implies that this is a $\mu(dx) dP$-null set. Hence, by Fubini's theorem, there exists a $\mu(dx)$-null set $X\setminus X_0\in \cX$ such that, for all $x\in X_0$, $\Vert g(x,\cdot)\Vert_{L^2([0,T])} < \infty$ a.s.

The existence of the function $\Psi$ with the required properties follows from Theorem \ref{rdthm9.4.1}, replacing $X$ and $Z$ there by $X_0$ and $g$.

The remainder of the proof is as in \cite[Theorem 2.2]{V2012}.
\end{proof}

\section{Notes on Appendix \ref{app1}}
\label{notes-app1}

Section \ref{rdsecA.4} presents notions and results that are part of the so-called {\em General Theory of Stochastic Processes} and that are used throughout this book. A comprehensive account of this theory is \cite{dm1}. Here, we have  mainly made use of \cite{chung-williams}, \cite{ks} and \cite{ry}.

The proof of Theorem \ref{app1-ch1-t1} and Corollary \ref{app1-ch1-c1} follows \cite[Theorem 2.3.2, p.~24]{ito84},  as adapted by Walsh in \cite[Theorem 4.1 and Corollary 4.2, p.~332]{walsh} to a family of nuclear spaces that includes $\cs(\rek)$ as a special case. As mentioned in \cite{ito84}, the idea to use the inequality \eqref{a1.30} is borrowed from
 \cite{sazonov}, and the idea to use integration with respect to a Gaussian measure to obtain \eqref{finals1} is taken from \cite{yamazaki}.

Kolmogorov's continuity criterion is extensively used throughout the book, in particular,  in Section \ref{ch1'-s7} (see the proofs of \eqref{ch1'-s7.l1-2000} and \eqref{ch1'-s7.l1-201} in Lemmas \ref{ch1'-s7-l1} and \ref{ch1'-s7-l2}, respectively).
The version given in Theorem \ref{ch1'-s7-t2} of Section \ref{app1-3}  is a variation on that of H.~Kunita \cite[Theorem 1.4.1, p.~31]{kunita}. It is based on the so called {\em chaining argument}. A particularly clear proof in the isotropic case is given by G.~Lowther (https://almostsuremath.com/2020/10/20/the-kolmogorov-continuity-theorem/). Our proof follows this presentation while introducing the ingredients needed to deal with the anisotropy.

A different type of proof of Kolmogorov's continuity criterion makes use of one of the many variations on the Garsia-Rodemich-Rumsey lemma \cite{g-r-r}. This type of proof can be found for instance in \cite[Appendix B, p.~312, Theorem B.1.1]{dz-1996}, \cite[Corollary 2.1.5]{sv-bis}, \cite[Theorem 1.1, p. 271]{walsh}, \cite[Appendix C]{khosh}, \cite[Theorem A.1, p. 573]{friz-victoir} and \cite[Lemma 2.4] {scheutzow-2018}. Our Proposition \ref{pr:Garsia} is yet another version of the Garsia-Rodemich-Rumsey lemma, in the spirit of  \cite{dkn07} and \cite{khosh}, which is similar to \cite[Lemma 2.4, p.~5]{scheutzow-2018}. We use this lemma in the proof of the version of Kolmogorov's continuity criterion given in Corollary \ref{ch1'-s7-t2-bis}, and this lemma is also useful in other contexts (see for instance \cite{florit-nualart} and \cite{dalang-pu-2020}). However, this approach does not seem to allow for index sets $I\times D$ that are as general as those considered in Theorem \ref{ch1'-s7-t2} (see Remark \ref{remgarsia}).

 The statement on measurability of the Itô integral that is proved in Theorem \ref{rdthm9.4.1} is used in Chapter \ref{ch2}. Its proof is an adaptation of a similar result given in \cite[Proposition 5]{strickeryor} in a slightly different setting. The stochastic Fubini's theorem in Section \ref{app1-2} is adapted from \cite{V2012}.



\chapter[Properties of fundamental solutions and Green's \texorpdfstring{\\}{}functions]{Properties of fundamental solutions and Green's functions}
\label{app2}
\pagestyle{myheadings}
\markboth{R.C.~Dalang and M.~Sanz-Sol\'e}{Fundamental solutions and Green's functions}

This chapter contains numerous properties of the fundamental solutions and Green's functions of the classical and fractional heat equations
on the real line and on bounded intervals, as well as for the fundamental solutions and Green's functions of the wave equation on $\R$, $\R_+$ and on a bounded interval. These are mainly integrability properties in the form of $L^p$-estimates ($p>0$), and upper and lower bounds on increments, mostly, but not only, in the $L^2$-norm. These results are extensively used in Chapters \ref{chapter1'} and \ref{ch1'-s5} for the study of random field solutions to the corresponding SPDEs and the properties of their sample paths.

For the heat equation, we also present bounds on the difference between its fundamental solution on $\R$ and its Green's function on $[-L, L]$ (see Section \ref{app2-4}), which are used in Sec. \ref{new-5.5-candil}. Finally, in Section \ref{rd08_26s1}, we give explicit values of some geometric space-time convolution series that are needed in Sections \ref{rdrough} and \ref{rd1+1anderson}.

\section{Heat kernel on $\rek$}
\label{app2-1}

The fundamental solution to the heat equation in $\rek$ is the function
\beqn
\Gamma(t,x; s,y)=\Gamma(t-s, x-y),
\eeqn
where
\beq
\label{app2.0}
\Gamma(r,z) = \frac{1}{(4\pi r)^\frac{k}{2}}\, \exp\left(-\frac{|z|^2}{4r}\right)1_{]0,\infty[}(r), \qquad z\in\rek,
\eeq
is the heat kernel.

We have seen in \eqref{heatcauchy-11'} that if $k=1$, then
\beq
\label{app2.00}
\int_0^t dr \int_{\re} dz\,  \Gamma^2(r,z) = \left(\frac{t}{2\pi}\right)^{\half},
\eeq
therefore, $\Vert\Gamma\Vert_{L^2(\re_+\times \re)}= +\infty$.
\medskip

\noindent{\em Increments in space and in time, $L^2$-norms, upper and lower bounds}
\medskip

The next lemma provides in particular estimates on $L^2(\re_+ \times \re)$-norms of increments in time and space of the heat kernel when the spatial dimension is $k=1$.

\begin{lemma}
\label{ch1'-l0}
Let $k=1$.

\begin{enumerate}
\item For all $h\in\R$,
\beq
 \label{1'.3}
\int_0^\infty dr\int_\IR dz \left[\Gamma(r,z)-\Gamma(r,z+h)\right]^2
 = \half |h|.
 \eeq
 \item For $\barh\ge0$,
  \beq
 \label{1'.4}
\int_0^\infty dr\ \int_\IR d z \left[\Gamma(r+\barh,z)-\Gamma(r,z)\right]^2
= \frac{\sqrt{2}-1}{(2\pi)^{\half}} \sqrt{\barh}.
\eeq
\item Fix $t > 0$. Then
\beq
\label{1'.4-bis}
 \lim_{\bar h \downarrow 0}\frac{1}{\sqrt{\barh}}  \int_0^t d r\int_\IR d z \left[\Gamma(r+\barh,z)-\Gamma(r,z)\right]^2 = \frac{\sqrt{2}-1}{(2\pi)^{\half}}.
\eeq
\item
Fix $C>0$. There is $c_0>0$ (given in \eqref{App-B(*1)} below) such that for all $t\ge 0$ and $h\in\re$ satisfying $|h|\le C\sqrt{t}$,
\beq
\label{1'.50}
 c_0 |h|\le \int_0^t dr\int_\IR dz \left[\Gamma(r,z)-\Gamma(r,z+h)\right]^2.
 \eeq
 \end{enumerate}

 As a consequence of 1. and 2. above, we deduce that for all $s,t \in \IR_+$ and $x,y\in\re$,
\begin{align}
\label{1'.500}
&\int_{0}^\infty dr \int_{\IR}  dz \left(\Gamma(t-r,x-z)-\Gamma(s-r,y-z)\right)^2\notag\\
& \qquad\qquad\qquad\qquad \qquad \le \left[\pi^{-\frac{1}{4}}|t-s|^{\frac{1}{4}} + 2^{-\half}|x-y|^\half \right]^2.
\end{align}
\end{lemma}
\begin{proof}
For the proof of  \eqref{1'.3}, we develop the square of the integrand and apply the identity
\beq
\label{1'.6}
\int_{\re} dz\, \Gamma(s,x-z) \Gamma(r,z) = \Gamma(s+r,x),
\eeq
valid for any $s,r>0$, which is the semigroup property \eqref{semig-heat} of the heat kernel.
We obtain
\begin{align}
\label{app2.001}
& \int_0^\infty dr\int_\IR dz \left[\Gamma(r,z)-\Gamma(r,z+h)\right]^2\notag\\
&\quad =  \int_0^\infty dr \int_\IR dz\, [\Gamma^2(r,z) +
\Gamma^2(r,z+h)-2\Gamma(r,z)\Gamma(r,z+h)]\notag\\
& \quad = \int_0^\infty dr\, [2 \Gamma(2r,0)-2\Gamma(2r,h)]\notag\\
&\quad  = 2 \int_0^\infty \frac{dr}{\sqrt{8\pi r}}\left(1-\exp\left(-\frac{h^2}{8r}\right)\right).
\end{align}
With the change of variables $w = \frac{|h|}{2\sqrt{2r}}$, the last right-hand side is equal to
 \beqn
\frac{1}{2} \frac{|h|}{\sqrt \pi} \int_{0}^\infty \frac{dw}{w^2}\left(1-\exp(-w^2)\right).
 \eeqn
Using integration by parts and the expression of the Gaussian density, one can easily check that
\beqn
\int_{0}^\infty \frac{dw}{w^2}\left(1-\exp(-w^2)\right) = 2\int_0^\infty \exp(-w^2)\, dw = \sqrt \pi.
\eeqn
Thus,
\beqn
\frac{1}{2}\frac{|h|}{\sqrt \pi} \int_{0}^\infty \frac{dw}{w^2}\left(1-\exp(-w^2)\right) =
\frac{|h|}{2},
\eeqn
which establishes \eqref{1'.3}.

Next, we prove  \eqref{1'.4}.  Fix $t\in \re_+$. By using again \eqref{1'.6}, we obtain
\begin{align}\nonumber
& \int_0^t dr\int_\IR d z \left[\Gamma(r+\barh,z)-\Gamma(r,z)\right]^2\\ \nonumber
&\qquad =  \int_0^t dr \left[\Gamma(2(\barh+r),0)+\Gamma(2r,0) - 2\Gamma(2r+\barh,0)\right]\\ \nonumber
&\qquad = \frac{1}{2\sqrt{2\pi}} \int_0^t dr \left(\frac{1}{\sqrt{\barh+r}}+\frac{1}{\sqrt{r}}-\frac{2\sqrt{2}}{\sqrt{2r+\barh}}\right)\\
&\qquad = \frac{1}{\sqrt{2\pi}}\left(\sqrt{\barh+t} + \sqrt{t} + (\sqrt{2}-1)\sqrt{\barh} - \sqrt{2}\sqrt{2t+\barh}\right).
\label{rdB.1.10}
\end{align}
Observe that $\sqrt{\barh+t} + \sqrt{t} - \sqrt{2}\sqrt{2t+\barh}\le 0$ and moreover,
\beqn
\lim_{t\to\infty} \left(\sqrt{\barh+t} + \sqrt{t} - \sqrt{2}\sqrt{2t+\barh}\right)=0.
\eeqn
Thus,
\beq
\label{uptoinfty}
\int_0^\infty dr\int_\IR d z \left[\Gamma(r+\barh,z)-\Gamma(r,z)\right]^2= \frac{\sqrt 2-1}{\sqrt{2\pi}} \sqrt{\barh}.
\eeq
Consequenly \eqref{1'.4} holds.
\smallskip

As for the equality \eqref{1'.4-bis}, it follows from \eqref{rdB.1.10} and the property
$$\lim_{\barh\to 0} \frac{\sqrt{\barh+t} + \sqrt{t} - \sqrt{2}\sqrt{2t+\barh}}{\sqrt{\barh}}=0,$$
which can be checked using l'Hospital's rule.


For \eqref{1'.50}, when $h=0$, the inequality is clear, so we assume $h\neq 0$.
Replacing the upper limit $\infty$ by $t>0$ in the calculations that led to \eqref{app2.001}, we obtain
\begin{align*}
\int_0^t dr\int_\IR dz \left[\Gamma(r,z)-\Gamma(r,z+h)\right]^2& =
2 \int_0^t \frac{dr}{\sqrt{8\pi r}}\left(1-\exp\left(-\frac{h^2}{8r}\right)\right)\\
& = \half\frac{|h|}{\sqrt \pi} \int_{\frac{|h|}{2\sqrt {2 t}}}^\infty \frac{dw}{w^2}\left(1-\exp(-w^2)\right),
\end{align*}
where we have again used the change of variable $w = \frac{|h|}{2\sqrt{2r}}$.

Divide the last expression by $|h|$. Assuming that $|h|\le C\sqrt t$, this is now bounded
from below by
\beq
\label{App-B(*1)}
c_0 := \frac{1}{2\sqrt{\pi}}\int_{\frac{C}{2\sqrt 2}}^\infty \frac{dw}{w^2}\left(1-\exp(-w^2)\right).
\eeq
This yields \eqref{1'.50}.


It remains to prove \eqref{1'.500}. Assume that $s \leq t$ and apply the triangle inequality to see that
\beqn
\left[\int_0^\infty dr \int_{\IR}  dz \left(\Gamma(t-r,x-z)-\Gamma(s-r,y-z)\right)^2\right]^\half \le T_1+T_2,
\eeqn
with
\begin{align*}
T_1^2&:= \int_0^\infty dr \int_{\IR}  dz \left(\Gamma(t-r,x-z)-\Gamma(s-r,x-z)\right)^2,\\
T_2^2&:= \int_0^\infty dr \int_{\IR}  dz \left(\Gamma(s-r,x-z)-\Gamma(s-r,y-z)\right)^2.
\end{align*}
Clearly,
\begin{align*}
T_1^2 & =  \int_0^s dr \int_{\IR}  dz \left(\Gamma(t-r,x-z)-\Gamma(s-r,x-z)\right)^2\\
&\qquad\quad  + \int_s^t dr \int_{\IR}  dz\ \Gamma^2(t-r,x-z).
\end{align*}
By \eqref{1'.4} and \eqref{app2.00},
\beqn
T_1^2  \le \frac{\sqrt 2-1}{(2\pi)^{\half}} \sqrt{t-s} + \frac{\sqrt{t-s}}{(2\pi)^{\half}} = \frac{1}{\pi^{\half}}\sqrt{t-s}.
\eeqn
By \eqref{1'.3},
\beqn
T_2^2 \le \half |x-y|,
\eeqn
and therefore,
\begin{align*}
&\left[\int_0^\infty dr \int_{\IR}  dz \left(\Gamma(t-r,x-z)-\Gamma(s-r,y-z)\right)^2\right]^\half\\
&\qquad\qquad \qquad\qquad\qquad\le \pi^{-\frac{1}{4}} |t-s|^{\frac{1}{4}} + 2^{-\half} |x-y|^{\half},
\end{align*}
which is \eqref{1'.500}.

This completes the proof of Lemma \ref{ch1'-l0}.
\end{proof}
\smallskip

In the remaining lemmas of this section, we examine the $L^p$-norms of the function $(s,y)\mapsto \Gamma(s,y)$, we provide upper bounds for $L^p$-norms of its increments in space, and upper and lower bounds for $L^p$-norms of its increments in time and space. 
\medskip

\noindent{\em $L^p$-norm of $\Gamma$}
\medskip

The next lemma identifies those $p>0$ for which $\Gamma \in L^p([0,T] \times \rek)$, with $k\in \IN^\ast$.

\begin{lemma}
\label{app2-l2}
Let $\Gamma$ be as in \eqref{app2.0}. Then for any $T>0$,
\beqn
\int_0^T dr \int_{\rek} dz\, \left(\Gamma(r,z)\right)^p < \infty\quad  \Longleftrightarrow \quad 0<p<1+\frac{2}{k},
\eeqn
and in this case,
\beq
\label{app2.0a-pre}
\int_0^T dr \int_{\rek} dz\, \left(\Gamma(r,z)\right)^p = C_{k,p} T^{1+\frac{k}{2}(1-p)}.
\eeq
\end{lemma}
\begin{proof}
We have
\beqn
\int_0^T dr \int_{\rek} dz \left(\Gamma(r,z)\right)^p= c_{k,p}\int_0^T dr\, r^{-\frac{kp}{2}} \int_{\rek} dz\, \exp\left(-p\frac{|z|^2}{4r}\right).
\eeqn
The $dz$-integral is equal to $+\infty$ if $p\le 0$. For $p>0$, the double integral is equal to
\begin{align}
\label{app2.0a}
&\tilde c_{k,p}\int_0^T dr\, r^{-\frac{kp}{2}+\frac{k}{2}} \int_{\rek} dz\ \left(\frac{r}{p}\right)^{-\frac{k}{2}}
\exp\left(-p\frac{|z|^2}{4r}\right)\nonumber\\
&\qquad = \bar c_{k,p}\int_0^T dr\, r^{-\frac{kp}{2}+\frac{k}{2}},
\end{align}
because $\int_{\rek} dz \left(\tfrac{p}{4\pi r}\right)^{\tfrac{k}{2}}\exp\left(-p\frac{|z|^2}{4r}\right)=1$.
The last integral in \eqref{app2.0a} is finite if and only if $1-\frac{kp}{2}+\frac{k}{2}>0$, that is, $p<1+\frac{2}{k}$, in which case we have \eqref{app2.0a-pre}.
\end{proof}
\smallskip

\noindent{\em $L^p$-norms, increments in space, upper bounds}
\medskip

The next lemma is Lemma A2 in \cite{SaintLoubert98}.

\begin{lemma}
\label{app2-l1}
Fix $T\in\re_+$ and let $p\in\, ]0,1+\tfrac{2}{k}[$. Set $p_c = 1+\frac{1}{k+1}$. Then  there exists a constant $C = C_{T,k,p} >0$ such that, for all $h\in\rek$,
\beq
\label{app2.1a}
\int_0^T dr \int_{\rek}  dz\, |\Gamma(r,z)-\Gamma(r,z-h)|^p \le C\min(\varphi_p(h),1),
\eeq
where
\beq
\label{app2.1}
   \varphi_p(h) =
\begin{cases}
|h|^p, & \mbox{if } p\in\, ]0,p_c[\, ,\\
|h|^{p_c} \log\left(2+\frac{1}{|h|}\right), & \mbox{if } p=p_c,\\
|h|^{2+k(1-p)}, & \mbox{if } p\in\, ]p_c, 1+\frac{2}{k}[.
\end{cases}
\eeq
In the case $p \in\, ]p_c, 1 + \frac{2}{k}[$, $T$ in \eqref{app2.1a} can be replaced by $+\infty$.
\end{lemma}

\begin{proof}
Set
\beqn
I_T(h): = \int_0^T dr \int_{\rek}  dz\, |\Gamma(r,z)-\Gamma(r,z-h)|^p.
\eeqn
By Minkowski's inequality and Lemma \ref{app2-l2},
\beq
\label{bound-1}
\sup_{h\in \rek}I_T(h) \leq 2^p \Vert \Gamma \Vert^p_{L^p([0,T]\times \rek)}<\infty.
\eeq
 Applying the change of variables $r=u \vert h\vert^2$, $z=w \vert h\vert $, and setting $e_0 = \frac{h}{\vert h \vert}$, we have
\beqn
I_T(h)= C_0\, \vert h \vert^{2 + k(1-p)}\, J_T(h),
\eeqn
where
\beqn
  J_T(h) = \int_0^{\frac{T}{\vert h \vert^2}} du\, u^{-\frac{p k}{2}} \int_{\rek} dw \left\vert \exp\left(-\frac{\vert w\vert^2}{4u}\right) - \exp\left(-\frac{\vert w-e_0\vert^2}{4u}\right) \right\vert^p.
\eeqn
Clearly, $h \mapsto J_T(h)$ is a decreasing function of $\vert h \vert$ and, since $I_T(h)$ is finite,
$J_T(h)$ is also finite.
 In particular, $J_T(h) < J_T(\sqrt T)<\infty$ when $|h|\geq \sqrt{T}$. This implies $J_T(h)\le C$ when $\vert h \vert\ge \sqrt T$.

For $\vert h \vert < \sqrt T$, for some constant $C<\infty$, we can write
\begin{align}
\label{jaja}
J_T(h) &= C\left(1+ \int_1^{\frac{T}{\vert h \vert^2}} du\, u^{-\frac{p k}{2}}\right.\notag\\
&\left. \qquad \quad\times\int_{\rek} dw  \left\vert \exp\left(-\frac{\vert w \vert^2}{4u}\right) - \exp\left(-\frac{\vert w-e_0\vert^2}{4u}\right)\right
\vert^p\right).
\end{align}
Assume without loss of generality that $e_0=(1,0,\dots,0)$. Then the integral in \eqref{jaja} is equal to
\begin{align*}
   & \int_1^{\frac{T}{\vert h \vert^2}} du\, u^{-\frac{p k}{2}} \int_{\re}
   dw_1\, \left\vert \exp\left(-\frac{w_1^2}{4u} \right) - \exp\left(-\frac{(w_1 - 1)^2}{4u} \right) \right\vert^{p} \\
    & \qquad \qquad\qquad \qquad\qquad \times \prod_{j=2}^k \int_{\re} dw_j\, \exp\left(-\frac{w_j^2 p}{4u} \right)\\
    &\qquad = C \int_1^{\frac{T}{\vert h \vert^2}} du\, u^{-\frac{p k}{2} + \frac{k}{2}-\frac{1}{2}} \int_{\re} dw_1\, \left\vert \exp\left(-\frac{w_1^2}{4u} \right) - \exp\left(-\frac{(w_1 - 1)^2}{4u}\right) \right\vert^p.
\end{align*}
Observe that $|w_1-1|\le |w_1|$ if and only if $w_1\ge \frac{1}{2}$. We use this fact to decompose the $dw_1$-integral on the right-hand side into the sum
\beqn
   J_1+J_2,
\eeqn
where
\begin{align*}
J_1 &= \int_{\frac{1}{2}}^{+\infty} dw_1  \left(-\exp\left(-\frac{w_1^2}{4u}\right) + \exp\left(-\frac{(w_1-1)^2}{4u}\right)\right)^p, \\
J_2 &= \int_{-\infty}^{\frac{1}{2}} dw_1  \left(\exp\left(-\frac{w_1^2}{4u}\right) - \exp\left(-\frac{(w_1-1)^2}{4u}\right)\right)^p.
\end{align*}
However, applying the change of variable $v=1-w_1$, we see that $J_2 = J_1$. Moreover,
\begin{align*}
J_1&= \int_{\frac{1}{2}}^{+\infty} dw_1   \left(1-\exp\left(\frac{(w_1-1)^2-w_1^2}{4u}\right)\right)^p \exp\left(-\frac{(w_1-1)^2}{4u}p\right)\\
&= \int_{-\frac{1}{2}}^{+\infty} dy\,  \exp\left(-\frac{y^2 p}{4u}\right) \left(1-\exp\left(\frac{y^2-(y+1)^2}{4u}\right)\right)^p\\
&\le \int_{-\frac{1}{2}}^{+\infty} dy\,  \exp\left(-\frac{y^2p}{4u}\right) \left(\frac{2y+1}{4u}\right)^p,
\end{align*}
where, in the second equality, we have applied the change of variables $y=w_1-1$, and for the third line, the inequality
$1-e^{-x}\le x$, valid for any $x\in\re_+$.

Apply the inequality $(a+b)^{p}\le 2^p(a^p+b^p)$ to obtain
\beqn
J_1 \le C \left(\int_0^{+\infty} dy\,  \exp\left(-\frac{y^2p}{4u}\right) u^{-p}\left(y^p+1\right)
 + \int_{-\frac{1}{2}}^{0} dy\, u^{-p} \exp\left(-\frac{y^2p}{4u}\right) \right).
\eeqn

From these computations, we deduce that
\begin{align*}
I_T(h)& = C |h|^{2+k(1-p)}\left\{1+\int_1^{\frac{T}{\vert h \vert^2}} du\, u^{-\frac{p k}{2} +\frac{k-1}{2}} \right.\\
&\left.\quad\times\left(\int_0^{+\infty} dy\,  u^{-p}\exp\left(-\frac{y^2p}{4u}\right) y^p + \int_{-\frac{1}{2}}^{+\infty} dy\, u^{-p} \exp\left(-\frac{y^2p}{4u}\right) \right)\right\}.
\end{align*}
In the last expression, the second integral with respect to the variable $y$ can be bounded above by the integral of the same function on the whole of $\re$. By doing so, and using the property of the Gaussian density, we obtain an upper bound of the form $C u^{\frac{1}{2}-p}$, for some positive constant $C$. As for the integral
\beqn
J:=\int_0^{+\infty} dy\,  \exp\left(-\frac{y^2p}{4u}\right) \frac{y^p}{u^p},
\eeqn
we apply the change of variables $z:= \frac{y^2p}{4u}$ to obtain
\beqn
J= \frac{2^p}{p^{(p+1)/2}}u^{\frac{1-p}{2}}\left(\int_0^{+\infty} dz\, e^{-z} z^{\frac{p-1}{2}}\right).
\eeqn
The $dz$-integral is  the evaluation of the Euler Gamma function $\Gamma_E(\frac{p+1}{2})$ (see \ref{Euler-gamma}). Hence $J\le C u^{\frac{1-p}{2}}$, with $C>0$ (depending on $p$).
Consequently,
\beq
\label{app2.2}
I_T(h) \le  C |h|^{2+k(1-p)}\left(1+ \int_1^{\frac{T}{\vert h\vert^2}} du\, u^{\frac{k-p - p k}{2}} + \int_1^{\frac{T}{\vert h\vert^2}} du\, u^{\frac{k - 2p -p k}{2}}\right).
\eeq

\noindent Both of these integrals can be evaluated explicitly, and since $k - 2 p - p k \leq k - p - p k$, it turns out that for $|h|<\sqrt T$, $k\ge 1$ and $p\in\,]0,1+\frac{2}{k}[$, the dominating term is always
\beqn
   \int_1^{\frac{T}{\vert h\vert^2}} du\, u^{(k - p -p k)/2}\le C |h|^{-(k-p-pk+2)},
\eeqn
when $k - p - pk + 2 \neq 0$, that is, $p \neq p_c$. Since $k - p - pk + 2 > 0$ if and only if $p < p_c$, this yields the following when $\vert h \vert < \sqrt{T}$:
\medskip

\noindent For $0< p < p_c$,
\beqn
   I_T(h) \leq C \, \vert  h \vert^{2+k(1-p)}\, \vert  h \vert^{-k+p + p k -2} = C \, \vert  h \vert^p.
\eeqn
For $p > p_c$,
\beqn
   I_T(h) \leq C \, \vert h \vert^{2+k(1-p)}.
\eeqn
Notice that in this case, we obtain the same upper bound if we replace $T$ by $+\infty$ in \eqref{app2.2}.

For $p = p_c$,
\beqn
   I_T(h) \leq C\, |h|^{p_c} \left(1+\log\left(\frac{T}{\vert  h \vert^2}\right)\right) \le  K \, |h|^{p_c} \log\left(2+\frac{1}{h}\right),
\eeqn
for $K$ large enough. These inequalities, along with the comments above \eqref{jaja}, end the proof of the Lemma.
\end{proof}

If $k=1$, the value of the critical exponent in the preceding lemma is $p_c = \tfrac{3}{2}$. Hence, $\varphi_p(h) = |h|$ for $p=2$.
Therefore, in this particular case, the upper bound \eqref{app2.1a} follows from \eqref{1'.3} in Lemma \ref{ch1'-l0}.
\medskip

\noindent{\em Increments in time, $L^p$-norms, upper and lower bounds}

 \begin{lemma}
\label{app2-l3}
The following estimates hold.
\begin{description}
\item{(a)} If $p\in \left]1,1+\frac{2}{k}\right[$, then there is a constant $C_{k,p}<\infty$ such that, for all
$t\ge 0$ and $h>0$,
\beq
\label{app2.5}
\int_0^\infty ds \int_{\IR^k} dz\, |\Gamma(t+h-s,z) - \Gamma(t-s,z)|^p \le C_{k,p}h^{1+\frac{k}{2}(1-p)}.
\eeq
If $p\in\, ]0,1[$ and $T>0$, then there is a constant $c_{k,p,T}<\infty$ such that, for all $t\ge 0$ and $h\in[0,T]$,
\beq
\label{app2.6}
\int_0^\infty ds \int_{\IR^k} dz\, |\Gamma(t+h-s,z) - \Gamma(t-s,z)|^p \le c_{k,p,T}h^p.
\eeq
If $p=1$, then there is a constant $C_{k,p} < \infty$ such that, for all $t\ge 0$ and $h\ge 0$,
\begin{align}
\label{app2.20}
&\int_0^\infty ds \int_{\rek} dz \, |\Gamma(t+h-s,z) - \Gamma(t-s,z)| \notag\\
&\qquad\qquad \qquad\le C_{k,p}\, h\, \left(1 + \log\left(\frac{t}{h}\right) 1_{\{h<t \}} \right).
\end{align}
\item{(b)} If $p\in\left]0,1+\frac{2}{k}\right[$, then there is a constant $c_{k,p}>0$ such that, for $t>0$ and $0<h<t$,
\beq
\label{app2.7}
c_{k,p} h^{1+\frac{k}{2}(1-p)}\le \int_0^t ds \int_{\IR^k} dz\, |\Gamma(t+h-s,z) - \Gamma(t-s,z)|^p
\eeq
and
\beq
\label{app2.8}
c_{k,p} h^{1+\frac{k}{2}(1-p)} = \int_t^{t+h}ds \int_{\IR^k} dz\, |\Gamma(t+h-s,z)|^p.
\eeq
\end{description}
\end{lemma}
\begin{proof}
\noindent(a)
We only consider the case $t>0$. Indeed, by \eqref{app2.0a-pre}, the three inequalities are clearly satisfied when $t=0$.

Let
\begin{align}
\label{defis}
I_1 &= \int_0^t ds \int_{\IR^k} dz\, |\Gamma(t+h-s,z) - \Gamma(t-s,z)|^p,\notag\\
I_2&=\int_t^{t+h} ds \int_{\IR^k} dz\, |\Gamma(t+h-s,z)|^p.
\end{align}
By a change of variables and \eqref{app2.0a-pre}, we obtain
\beq
\label{app2.10}
I_2 = \int_0^h ds \int_{\rek} dz\, |\Gamma(s,z)|^p = \bar c_{k,p} h^{1+\frac{k}{2}(1-p)},
\eeq
for any $p\in \left]0,1+\frac{2}{k}\right[$.

We now proceed to study $I_1$, starting with
the following calculation that is partly in \cite[p. 320]{walsh} when $k=1$. Pass to polar coordinates in $z$ to see that
\begin{align*}
I_1&= c_{k}\int_0^t ds \int_0^\infty dr\, r^{k-1}\notag\\
&\qquad \times \left\vert (s+h)^{-\frac{k}{2}} \exp\left(-\frac{r^2}{4(s+h)}\right)
- s^{-\frac{k}{2}}\exp\left(-\frac{r^2}{4s}\right)\right\vert^p.
\end{align*}
Use the change of variables  $s=hv$, $r= h^{1/2}u$, to get
\begin{align}
\label{app2.11}
I_1
&=  c_{k}\, h^{1+\frac{k}{2}(1-p)} \int_0^{t/h} dv \int_0^\infty  du\, u^{k-1}\notag\\
&\qquad\times\left\vert (v+1)^{-\frac{k}{2}} \exp\left( -\frac{u^2}{4(v+1)}\right) - v^{-\frac{k}{2}} \exp\left( -\frac{u^2}{4v}\right)\right\vert^p.
\end{align}
\medskip

\noindent{\em Case $p\in\left]1, 1+\frac{2}{k}\right[$.}
We will prove that
\beq
\label{app2.9}
I_1
 =  \int_0^t ds \int_{\IR^k} dz\, |\Gamma(s+h,z) - \Gamma(s,z)|^p
\le C_{k,p}h^{1+\frac{k}{2}(1-p)}.
\eeq
Along with \eqref{app2.10}, this will imply \eqref{app2.5}.

The double integral in \eqref{app2.11} is bounded above by
\beq
\label{app2.12}
\int_0^\infty  dv \int_0^\infty  du\, u^{k-1}\left\vert (v+1)^{-\frac{k}{2}} \exp\left( -\frac{u^2}{4(v+1)}\right) - v^{-\frac{k}{2}} \exp\left( -\frac{u^2}{4v}\right)\right\vert^p,
\eeq
which, as we now show, is a finite constant that depends on $k$ and $p$.

Indeed, first we restrict the domain of the $v$-variable to $[0,1]$ and we see that
\begin{align}
\label{app2.14}
&\int_0^1  dv \int_0^\infty  du\, u^{k-1}\left\vert (v+1)^{-\frac{k}{2}} \exp\left( -\frac{u^2}{4(v+1)}\right) - v^{-\frac{k}{2}} \exp\left( -\frac{u^2}{4v}\right)\right\vert^p\notag\\
&\qquad \le c_p(J_1 + J_2),
\end{align}
where
\begin{align*}
J_1 &= \int_0^1  dv \int_0^\infty  du\, u^{k-1} (v+1)^{-\frac{kp}{2}}\exp\left( -\frac{p u^2}{4(v+1)}\right),\\
J_2  &= \int_0^1  dv \int_0^\infty  du\, u^{k-1} v^{-\frac{kp}{2}} \exp\left( -\frac{pu^2}{4v}\right).
\end{align*}
For $v\in[0,1]$, we have $(v+1)^{-\frac{kp}{2}}\le 1$ and $\exp\left( -\frac{p u^2}{4(v+1)}\right)\le \exp\left(-\frac{p u^2}{8}\right)$.
Therefore,
\begin{align*}
J_1&\le  \int_0^1  dv \int_0^\infty  du\, u^{k-1}\exp\left(-\frac{p u^2}{8}\right)
 =  \int_0^\infty  du\ u^{k-1}\exp\left(-\frac{p u^2}{8}\right)\\
 &= \tilde c_{p,k}\, \Gamma_E\left(\frac{k}{2}\right) <\infty,
\end{align*}
where $\Gamma_E$ is the Euler gamma function (see \eqref{Euler-gamma}), and we have applied the change of variables $u\mapsto \frac{p u^2}{8}$.

For $J_2$, we apply the change of variables $(v,u)\mapsto (v,\frac{pu^2}{4v})$ to obtain
\begin{align*}
J_2 &= c_{p} \left(\int_0^1 dv\, v^{\frac{k}{2} (1-p)} \right) \left(\int_0^\infty
 dx\,  x^{\frac{k}{2}-1} e^{-x}\right) = \tilde c_{p,k} \Gamma_E\left(\frac{k}{2}\right),
 \end{align*}
 where, in the last equality, we have used that $p<1+\frac{2}{k}$ to obtain that the first integral factor is finite.

In order to conclude that the integral in \eqref{app2.12} is finite, it remains to check that
\beq
\label{app2.15}
\int_1^\infty dv \int_0^\infty du\, u^{k-1}\left\vert (v+1)^{-\frac{k}{2}} \exp\left( -\frac{u^2}{4(v+1)}\right) - v^{-\frac{k}{2}} \exp\left( -\frac{u^2}{4v}\right)\right\vert^p
\eeq
is finite.
For this, we fix $u > 0$ and define
\beqn
f(v) = v^{-\frac{k}{2}}\exp\left(-\frac{u^2}{4v}\right).
\eeqn
Then
\beqn
f^\prime (v)= -\frac{k}{2} v^{-\frac{k}{2}-1} \exp\left(-\frac{u^2}{4v}\right)+v^{-\frac{k}{2}}\exp\left(-\frac{u^2}{4v}\right)\frac{u^2}{4v^2},
\eeqn
and
\begin{align*}
|f^\prime (v)|&\leq c_k \left(v^{-\frac{k}{2}-1} +v^{-\frac{k}{2}-2}u^2\right)\exp\left(-\frac{u^2}{4v}\right)\notag\\
& \leq \bar c_k v^{-(\frac{k}{2}+1)} (1+\tfrac{u^2}{v})\exp\left(-\frac{u^2}{2v}\right).
\end{align*}
Therefore, applying the intermediate value theorem, \eqref{app2.15} is bounded above by
\beqn
c_k \int_1^\infty dv \int_0^\infty du\, u^{k-1} \,v^{-\left(\frac{k}{2}+1\right)p}\left(1+ \tfrac {u^2}{v}\right)^p \exp\left(-\frac{pu^2}{2(v+1)}\right).
\eeqn

We split this term into the sum $c_{k,p}\, (T_1+T_2)$, where
\begin{align*}
T_1&= \int_1^\infty dv \int_0^\infty du\, u^{k-1}\ v^{-\left(\frac{k}{2}+1\right)p} \exp\left(-\frac{pu^2}{2(v+1)}\right),\\
T_2&= \int_1^\infty dv \int_0^\infty du\, u^{k-1}\ v^{-\left(\frac{k}{2}+1\right )p} \left(\frac {u^2}{v}\right)^p \exp\left(-\frac{pu^2}{2(v+1)}\right),
\end{align*}
and prove that both terms $T_1$ and $T_2$ are finite.

Indeed, apply the change of variables $(v,u) \mapsto \left(v,\frac{u}{\sqrt{v+1}}\right)$
 to see that $T_1$ is bounded above by
\begin{align}
\label{app2.16}
&c\int_1^\infty dv\, (v+1)^{\frac{1}{2}+\frac{k-1}{2}} v^{-\left(\frac{k}{2}+1\right)p} \int_0^\infty dx\,  x^{k-1}
\exp\left(-\frac{p x^2}{2}\right) \notag\\
& \qquad \le C_{k,p}\, \Gamma_E \left(\frac{k}{2}\right) \int_1^\infty dv\, (v+1)^{\frac{k}{2}} v^{-\left(\frac{k}{2} + 1 \right)p}\notag\\
&\qquad \le \bar{C}_{k,p} \int_1^\infty dv\,  v^{\frac{k}{2}-(\frac{k}{2}+1)p},
\end{align}
where we use the inequality $v+1\le 2v$ for $v\ge 1$.
The last integral converges provided $\tfrac{k}{2}-\left(\tfrac{k}{2}+1\right)p+1<0$, that is, for any $p>1$.
Hence, $T_1$ is finite.

Next, we apply the change of variables $(v,u)\mapsto \left(v,\frac{pu^2}{2(v+1)}\right)$ to deduce that
\begin{align*}
T_2 &= \int_1^\infty dv \int_0^\infty dx\, v^{-\left(\frac{k}{2}+2\right)p}\ x^{k-1+2p}\exp\left(-\frac{px^2}{2(v+1)}\right)\\
&\le  c_p  \int_1^\infty dv\,  v^{-\left(\frac{k}{2}+2\right)p} (v+1)^{\frac{k+2p}{2}}\int_0^\infty dx\, x^{\frac{k-2+2p}{2}} e^{-x}\\
&\le c_{p} \Gamma_E\left(\frac{k}{2}+p\right) \int_1^\infty dv\, v^{\frac{k}{2} -\left(\frac{k}{2}+1\right)p}.
\end{align*}
 Since we are assuming $p>1$, the integral converges and therefore $T_2$ is finite.
 We conclude that the integral \eqref{app2.12} is finite and this  finishes the proof of
 \eqref{app2.9}.
Claim \eqref{app2.5} is proved.
\medskip

\noindent{\em Case $p=1$.}
For the term $I_1$, we begin with \eqref{app2.11}, but we replace \eqref{app2.12} by
\beqn
\int_0^{\tfrac{t}{h}}  dv \int_0^\infty  du\, u^{k-1}\left\vert (v+1)^{-\frac{k}{2}} \exp\left( -\frac{u^2}{4(v+1)}\right) - v^{-\frac{k}{2}} \exp\left( -\frac{u^2}{4v}\right)\right\vert^p.
\eeqn

If $\tfrac{t}{h}\le 1$, then the calculations that follow \eqref{app2.14} show that this integral is bounded by a constant.
If, on the contrary, $\tfrac{t}{h} > 1$, then up to a multiplicative constant $c_{k,p}$, we proceed as in the previous case to obtain, as in \eqref{app2.16}, the following upper bound (which takes into account the term $T_2)$:
\beq
\label{app2.170}
1+\int_1^{\tfrac{t}{h}} dv\, v^{\frac{k}{2}-(\frac{k}{2}+1)p} =
1+\int_1^{\tfrac{t}{h}} dv\, v^{-p} = 1+\log\left(\frac{t}{h}\right),
\eeq
because $p=1$ here. Together with \eqref{app2.10}, this yields the claim \eqref{app2.20}.
\medskip

\noindent{\em Case $p\in\, ]0,1[$.}
For the term $I_2$, we again use \eqref{app2.10}. As above, for the term $I_1$, we go to \eqref{app2.11}. Then we proceed as in \eqref{app2.12}--\eqref{app2.16}, but with the upper bound $\infty$ in the $dv$-integral replaced by $t/h$ (assuming that $h<t$), and we obtain the following upper bound for $T_1$:
\beq
\label{app2.13a}
c_{k,p} \int_1^{t/h} dv\,  v^{\frac{k}{2}(1-p)-p}
\le \bar c_{k,p}  \left(\frac{t}{h}\right)^{\frac{k}{2} - (\frac{k}{2} + 1)p+1}.
\eeq
Multiplying by the factor $h^{1+\frac{k}{2}(1-p)}$, which appears in \eqref{app2.11}, gives the bound
$c_{k,p,T} h^p$ for $T_1$, as claimed.
For $T_2$, we obtain
    $c \int_1^{t/h} dv\, v^{(1-p)k/2-p}$.
 This is the same bound as in \eqref{app2.13a},
at least when $h < t$. When $T \geq h \geq t$, the exponent of $h$ is $1 + \frac{k}{2} (1 - p) > p$ since $p \in\, ]0, 1[$, so we adjust the constant $c_{k,p,T}$, which gives \eqref{app2.6}.
\medskip

\noindent(b)  For the proof of \eqref{app2.7}, it suffices to bound $I_1$ from below,
 and we use the calculation of $I_1$ up to  \eqref{app2.11}. Let $0<h<t$. The double integral in \eqref{app2.7} is bounded below by
\beqn
\int_0^1dv \int_0^\infty  du\, u^{k-1}
\left\vert (v+1)^{-\frac{k}{2}} \exp\left( -\frac{u^2}{4(v+1)}\right) - v^{-\frac{k}{2}} \exp\left( -\frac{u^2}{4v}\right)\right\vert^p,
\eeqn
which is a constant that depends on $k$ and $p$. This proves \eqref{app2.7}.
\smallskip

The claim \eqref{app2.8} follows from \eqref{app2.10}.
\end{proof}
\medskip

\noindent{\em The case $k\geq 1$ and $p=1$: a lower bound on time increments}
\medskip

In the case $p=1$, a better bound than the one given in part (b) of Lemma \ref{app2-l3} is the following.
\begin{lemma}
\label{app2-l4}
There is a constant $c_k>0$ such that, for $t>0$ and $0<h<t$,
\begin{align}
\label{app2.18}
\int_0^t ds \int_{\rek} dz\, \vert \Gamma(t+h-s,z)-\Gamma(t-s,z)\vert \ge c_k\, h \log \left(\frac{t}{h}\right).
\end{align}
\end{lemma}

\begin{proof} Clearly, the left-hand side of \eqref{app2.18} is equal to
\[
    \int_0^t ds \int_{\rek} dz\, \vert \Gamma(s+h,z)-\Gamma(s,z)\vert .
\]
We now use arguments of
\cite[Lemma A3]{SaintLoubert98}.
Do the change of variables $s=h u$, $z=y\sqrt{h}$, and apply the scaling property $\Gamma(hu,y\sqrt h) = h^{-\frac{k}{2}}\Gamma(u,y)$  to see that
\beqn
\int_0^t ds \int_{\rek} dz\, \vert \Gamma(s+h,z)-\Gamma(s,z)\vert=h\int_0^{t/h}du \int_{\rek} dy\, \vert\Gamma(1+u,y)-\Gamma(u,y)\vert.
\eeqn
Next, we identify those $y$ for which $\Gamma(1+u,y)\ge \Gamma(u,y)$, as follows:
\begin{align*}
\Gamma(1+u,y)\ge \Gamma(u,y) &\Longleftrightarrow (1+u)^{-\frac{k}{2}} \exp\left(-\frac{|y|^2}{4(1+u)}\right) \ge u^{-\frac{k}{2}}
\exp\left(-\frac{|y|^2}{4u}\right)\\
&\Longleftrightarrow \left(1+\frac{1}{u}\right)^{-\frac{k}{2}}\ge \exp\left(\frac{|y|^2}{4}\left(\frac{1}{1+u}-\frac{1}{u}\right)\right)\\
&\Longleftrightarrow -\frac{k}{2}\ln\left(1+\frac{1}{u}\right)\ge -\frac{|y|^2}{4}\frac{1}{u(1+u)}\\
&\Longleftrightarrow |y|^2 \ge 2k u(1+u)\ln \left(1+\frac{1}{u}\right)\\
&\Longleftrightarrow |y| \ge z_0(u),
\end{align*}
where
\beqn
z_0(u) = \left(2k u(1+u)\ln \left(1+\frac{1}{u}\right)\right)^{\frac{1}{2}}.
\eeqn
Notice that $z_0(u)>0$ for $u>0$.

Define
\beqn
I(u) = \int_{\rek} dy\, |\Gamma(1+u,y)-\Gamma(u,y)|.
\eeqn
Then
\begin{align*}
I(u)
& =\int_{|y|\le z_0(u)} dy\, (\Gamma(u,y)-\Gamma(1+u,y))\\
&\qquad +\int_{|y|>z_0(u)} dy\, (\Gamma(1+u,y)-\Gamma(u,y))\\
&=\int_{|y|\le z_0(u)} dy\, \Gamma(u,y)-\int_{|y|>z_0(u)} dy\ \Gamma(u,y)\\
&\qquad -\int_{|y|\le z_0(u)} dy\, \Gamma(1+u,y)+\int_{|y|>z_0(u)} dy\ \Gamma(1+u,y).
\end{align*}
In the first two integrals, apply the change of variables $y=x\sqrt{u}$, and in the last two integrals, the change of variables $y= x\sqrt{u+1}$, to get
\begin{align*}
I(u)&= \int_{|x|\le \frac{z_0(u)}{\sqrt{u}}} dx\, \Gamma(1,x) - \int_{|x| > \frac{z_0(u)}{\sqrt{u}}}dx\, \Gamma(1,x)\\
& \quad -\int_{|x|\le \frac{z_0(u)}{\sqrt{u+1}}} dx\, \Gamma(1,x) + \int_{|x| > \frac{z_0(u)}{\sqrt{u+1}}}dx\, \Gamma(1,x)\\
& = 2 \int_{\frac{z_0(u)}{\sqrt{u+1}}\le |x|\le \frac{z_0(u)}{\sqrt{u}}} dx \, \Gamma(1,x).
\end{align*}
Passing to polar coordinates, we conclude that
\begin{align}
\label{app2.19}
&\int_0^t ds \int_{\rek} dz\, \vert \Gamma(s+h,z)-\Gamma(s,z)\vert \\
 &\qquad =  h\int_0^{t/h} du\, I(u)\notag\\
&\qquad = c_k h \int_0^{t/h} du \int_{\frac{z_0(u)}{\sqrt{u+1}}}^{\frac{z_0(u)}{\sqrt{u}}}dr\, r^{k-1} \exp\left(-\frac{r^2}{4}\right).
\end{align}
From this equality, we can proceed to a lower bound.

For $h<t$, the integral $\int_0^{t/h} du\, I(u)$ is equal to
$I_1+I_2$,
where
\beqn
I_1
   = \int_0^1 du\, I(u), \quad I_2=\int_1^{t/h} du\, I(u).
   \eeqn
Since $I_1\ge 0$, we only need a lower bound on $I_2$.

Observe that
\beqn
\frac{z_0(u)}{\sqrt{u}} = \left(2k(1+u)\ln\left(1+\frac{1}{u}\right)\right)^{\frac{1}{2}},
\eeqn
and
\beq
\label{app2.21}
\lim_{u\to\infty} (1+u)\ln\left(1+\frac{1}{u}\right)=1.
\eeq
Therefore, there is $c_0<\infty$ such that $\frac{z_0(u)}{\sqrt{u}}\le c_0$ for $u\ge 1$. In particular, for $u\ge 1$,
\begin{align*}
&I(u)\\
&\quad\ge \exp\left(-\frac{c_0^2}{4} \right) \int_{\frac{z_0(u)}{\sqrt{u+1}}}^{\frac{z_0(u)}{\sqrt{u}}}dr\, r^{k-1}\\
&\quad\ge  \exp\left(-\frac{c_0^2}{4} \right) z_0(u) \left(\frac{1}{\sqrt u}-\frac{1}{\sqrt{u+1}}\right)\left(\frac{z_0(u)}{\sqrt{u+1}}\right)^{k-1}\\
&\quad= \exp\left(-\frac{c_0^2}{4} \right) \left(2ku(1+u )\ln\left(1+\frac{1}{u}\right)\right)^{\frac{1}{2}}
\frac{\sqrt{u+1}-\sqrt u}{\sqrt u\sqrt{u+1}}\left(\frac{z_0(u)}{\sqrt{u+1}}\right)^{k-1}\\
&\quad=\exp\left(-\frac{c_0^2}{4} \right) \left(2k\ln\left(1+\frac{1}{u}\right)\right)^{\frac{1}{2}}\frac{1}{\sqrt u+\sqrt{u+1}}\left(\frac{z_0(u)}{\sqrt{u+1}}\right)^{k-1}.
\end{align*}
Now
\beqn
\frac{z_0(u)}{\sqrt{u+1}}=\left(2ku\ln\left(1+\frac{1}{u}\right)\right)^{\frac{1}{2}},
\eeqn
and by the same argument as in \eqref{app2.21},
we see that there are $\tilde c_k > 0$ and $c_1 > 0$ such that for all $u \geq 1$, $\frac{z_0(u)}{\sqrt{u+1}} \geq \tilde c_k$ and $\ln(1 + \frac{1}{u}) \geq c_1 \frac{1}{u}$, and therefore for $u\ge 1$, $\ln(1 + \frac{1}{u}) \geq \tilde c_k \frac{1}{u}$ for some constant $\tilde c_k$, and
\beqn
I(u)\ge c_k \left(\ln\left(1+\frac{1}{u}\right)\right)^{\frac{1}{2}} \frac{1}{\sqrt u+\sqrt{u+1}} \ge \tilde C_k\, \frac{1}{u}.
\eeqn

Thus, going back to \eqref{app2.19}, we have
\begin{align}
\label{app2.22}
\int_0^t ds \int_{\rek} dz\, \vert \Gamma(s+h,z)-\Gamma(s,z)\vert
&\ge C_k\, h \int_1^{t/h} \frac{du}{u}\notag\\
& = C_k\, h \ln\left(\frac{t}{h}\right),
\end{align}
for some constant $C_k>0$.
This ends the proof of the lemma.
\end{proof}

\section{Heat kernel with Dirichlet boundary conditions}
\label{app2-2}


Recall from \eqref{ch1'.600} that the Green's function of the heat kernel with Dirichlet boundary conditions is given by
\beq
\label{B.2Green-D}
G_L(t; x,y) = \sum_{n=1}^\infty e^{-\frac{\pi^2}{L^2}n^2 t} e_{n,L}(x)\ e_{n,L}(y), \quad t>0,\quad  x,y\in[0,L],
\eeq
where
\beqn
   e_{n,L}(x) = \sqrt{\frac{2}{L}}\, \sin\left(\frac{n\pi}{L}x\right),\qquad n\ge 1.
\eeqn
In this section, we study integrated squared increments in time and in space of $G_L$.
\medskip

\noindent{\em Increments in space and in time, $L^2$-norms, upper bounds}

\begin{lemma}
\label{ch1'-l2}
The Green's function $G_L(t;x,y)$ given in \eqref{B.2Green-D} satisfies  the following properties.
\begin{enumerate}
\item For any $x,y\in [0,L]$,
\begin{align}
\label{1'.10}
T_{1,L}(x,y):&=\int_0^\infty dr \int_0^L dz\,  [G_L(r; x,z) - G_L(r; y,z)]^2\notag\\
   &=\frac{1}{2}\left(|x-y| - \frac{1}{L} |x-y|^2 \right).
\end{align}
\item For any $h>0$ and $x\in [0,L]$,
\begin{align}
\label{1'.11}
T_{2,L}(h;x):&= \int_0^\infty dr \int_0^L dz\,  [G_L(r+h;x,z) - G_L(r;x,z)]^2\notag\\
&\le \frac{3}{\pi} h^{\half}.
\end{align}
\item For any $t\ge 0$ and $x\in [0,L]$,
\beq
\label{1'.11bis}
T_{3,L}(t; x):=\int_0^t dr \int_0^L dz\, G_L^2(r;x,z)\le \left(\frac{t}{2\pi}\right)^{\half}.
\eeq
\end{enumerate}
As a consequence, there exists $C>0$ such that, for any $t>0$, $0\le s\le t$ and $x,y\in [0,L]$,
\begin{align}\nonumber
& \int_0^t dr\ \int_0^L dz\, [G_L(t-r;x,z)-G_L(s-r;y,z)]^2 \\
   &\qquad\qquad\qquad\qquad \le C \left(|t-s|^\frac{1}{4} + |x-y|^\half \right)^2.
   \label{1'1100}
\end{align}
\end{lemma}

\begin{proof}
Because of the scaling property of the Green's function given in \eqref{scalingD}, by applying
the change of variables $r\mapsto \frac{r}{L^2}$, $z\mapsto \frac{z}{L}$, we see that
\begin{align}
\label{scaling-1}
&\int_0^\infty dr \int_0^L dz \,  [G_L(r; x,z) - G_L(r; y,z)]^2 \notag\\
&\qquad = \frac{1}{L^2} \int_0^\infty dr \int_0^L dz \left[G_1\left(\tfrac{r}{L^2}; \tfrac{x}{L},\tfrac{z}{L}\right) - G_1\left(\tfrac{r}{L^2}; \tfrac{y}{L},\tfrac{z}{L}\right)\right]^2\notag\\
&\qquad= L \int_0^\infty dr \int_0^1 dz \left[G_1\left(r; \tfrac{x}{L},z\right) - G_1\left(r; \tfrac{y}{L},z\right)\right]^2
\end{align}
and
\begin{align}
\label{scaling-2}
&\int_0^\infty dr \int_0^L dz\,  [G_L(r+h;x,z) - G_L(r;x,z)]^2\notag\\
&\qquad =\frac{1}{L^2} \int_0^\infty dr \int_0^L dz \left[G_1\left(\tfrac{r+h}{L^2}; \tfrac{x}{L},\tfrac{z}{L}\right) - G_1\left(\tfrac{r}{L^2}; \tfrac{x}{L},\tfrac{z}{L}\right)\right]^2\notag\\
&\qquad = L\int_0^\infty dr \int_0^1 dz\left[G_1\left(r+\tfrac{h}{L^2};\tfrac{x}{L},z\right) - G_1\left(r;\tfrac{x}{L},z\right)\right]^2.
\end{align}
These two identities show that for the proofs of \eqref{1'.10} and \eqref{1'.11}, it suffices to consider the case $L=1$.


We begin with \eqref{1'.10} for $L=1$.
Using  \eqref{B.2Green-D}, we have
\begin{align*}
T_{1,1}(x,y) & =\int_0^\infty dr\ \int_0^1 dz\, [G_1(r;x,z)-G_1(r;y,z)]^2\\
& = \int_0^\infty dr\ \int_0^1 dz\, \left(\sum_{n=1}^\infty e^{-\pi^2 n^2 r} [e_{n,1}(x)-e_{n,1}(y)]e_{n,1}(z)\right)^2\\
& = \sum_{n=1}^\infty [e_{n,1}(x)-e_{n,1}(y)]^2 \int_0^\infty dr\,  e^{-2\pi^2 n^2 r},
\end{align*}
where, in the last identity, we have used the orthogonality in $L^2([0,1])$ of the sequence $(e_{n,1},\, n\ge 1)$.

Since $e_{n,1}(x) = \sqrt 2\sin\left(n\pi x\right)$ and the integral is equal to $\frac{1}{2\pi^2n^2}$, we see using Lemma \ref{lemB.5.2} that
\begin{align}
\label{proces4}
T_{1,1}(x,y) & = \sum_{n=1}^\infty\, \frac{[\sin(n\pi x)-\sin(n\pi y)]^2}{\pi^2\, n^2}\notag\\
&= \frac{1}{2}\left(\left\vert x-y \right\vert - \left\vert x-y \right\vert^2 \right).
\end{align}
This proves \eqref{1'.10} for $L=1$. Using \eqref{scaling-1}, we obtain \eqref{1'.10} for all $L>0$.

Turning to \eqref{1'.11}, and using again \eqref{B.2Green-D} with $L=1$, we have
\begin{align}
\label{1'.1200}
T_{2,1}(h;x)& = \int_0^\infty dr \int_0^1 dz\, [G_1(r+h; x,z)-G_1(r;x,z)]^2\notag\\
& = \int_0^\infty dr \int_0^1 dz \left[\sum_{n=1}^\infty\left(e^{-\pi^2n^2(r+h)} - e^{-\pi^2n^2r}\right) e_{n,1}(x)e_{n,1}(z)\right]^2\notag\\
& = \int_0^\infty dr\  \sum_{n=1}^\infty\, (e_{n,1}(x))^2\left(e^{-\pi^2n^2(r+h)} - e^{-\pi^2n^2r}\right)^2\notag\\
& \le 2 \int_0^\infty dr\,  \sum_{n=1}^\infty\left(e^{-\pi^2n^2(r+h)} - e^{-\pi^2n^2r}\right)^2.
\end{align}

The integrand in the last expression, is equal to
\beqn
\sum_{n=1}^\infty e^{-2\pi^2 n^2r} \left(1-e^{-\pi^2 n^2h}\right)^2 .
\eeqn
By computing the $dr$-integral of this expression and using the inequality $1-e^{-x} \le \min(1,x)$, valid for any $x\ge 0$, we have
\begin{align}
T_{2,1}(h;x) &\le \frac{1}{\pi^2} \sum_{n=1}^\infty \frac{\left(1-e^{-\pi^2 n^2h}\right)^2}{n^2}\label{T2}\notag\\
& \le \frac{1}{\pi^2} \sum_{n=1}^\infty\, \min(n^{-2}, \pi^4n^2h^2).
\end{align}
Assume first that $h\ge \frac{1}{\pi^2}$. Then for all $n\ge 1$, $\min(n^{-2}, \pi^4n^2h^2)= n^{-2}$ and
\beqn
T_{2,1}(h,x) \le \frac{1}{\pi^2} \sum_{n=1}^\infty \frac{1}{n^2} = \frac{1}{\pi^2}\frac{\pi^2}{6}\le \frac{\pi}{6} h^\half.
\eeqn
Suppose next that $0<h<\frac{1}{\pi^2}$. Then
\begin{align*}
\frac{1}{\pi^2} \sum_{n=1}^\infty\, \min(n^{-2}, \pi^4 n^2h^2) &=\frac{1}{\pi^2} \sum_{n=1}^{\left[\frac{1}{\pi}h^{-\half}\right]}\pi^4n^2 h^2 + \sum_{n=\left[\frac{1}{\pi}h^{-\half}\right]+1}^\infty n^{-2}\\
& \le \frac{1}{\pi^2}\left[\pi h^\half + 2\int_{\frac{1}{\pi}h^{-\half}}^\infty z^{-2}\, dz\right]\\
& = \frac{1}{\pi^2}\left[\pi h^\half + 2\pi h^\half\right]\\
&\le \frac3{\pi}h^\half.
 \end{align*}

 Since $\max\left(\frac{\pi}{6}, \frac{3}{\pi}\right) = \frac{3}{\pi}$, \eqref{1'.11} with $L=1$ is proved. Using \eqref{scaling-2}, we obtain \eqref{1'.11} for all $L>0$.
\medskip

Using Property (ii) in Proposition \ref{ch1'-pPD} and \eqref{app2.00}, we have
\beqn
T_{3,L}(t; x)\le \int_0^t dr \int_0^L dz\, \Gamma^2(r,x-z) \le \left(\frac{t}{2\pi}\right)^{\half},
\eeqn
with $\Gamma$ defined in \eqref{app2.0}, proving \eqref{1'.11bis}.

Finally, \eqref{1'1100} follows from \eqref{1'.10}-\eqref{1'.11bis} by the triangle inequality.
\end{proof}

\begin{remark} If one is only interested in an upper bound for $T_{1,L}(x,y)$, then one can write
\beqn
[\sin(n\pi x)-\sin(n\pi y)]^2 \leq \min\left(4,n^2 \pi^2 \left\vert x-y \right\vert^2 \right)
\eeqn
in \eqref{proces4}, and then study the resulting series as in \eqref{T2}, to obtain the inequality $T_{1,1}(x,y) \leq \tfrac{6}{\pi}\, \vert x - y\vert$. This approach is used in \cite[Lemma A.3]{D-K-Z016}.

\end{remark}
\medskip
\noindent{\em $L^2$-norms, lower bounds}

\begin{lemma}
\label{ch1'-l200} The Green's function given in \eqref{B.2Green-D} satisfies the following.
Let $c_0=\frac{\sqrt{2\pi}-1}{8\pi}$.
Fix $\alpha \in\, ]0, \half[$. There exists $c(\alpha) > 0$ such that, for all $x \in [\alpha L, (1-\alpha) L]$ and every $0 \leq h \leq \min\left(\frac{c_0}{c(\alpha)} L^2, x^2, (L-x)^2\right)$,
\beq
\label{1'.112}
\int_0^h dr \int_{x-\sqrt{h}}^{x+\sqrt{h}} dz\, G_L^2(r;x,z) \ge c_0  h^{\half}.
\eeq
In particular, for any $x \in [\alpha L, (1-\alpha) L]$ and $0 \leq s \leq t$ satisfying $t - s \leq \min\left(\frac{c_0}{c(\alpha)} L^2, x^2, (L-x)^2\right)$,
\beq
\label{1'.1111}
\int_s^t dr \int_{x-\sqrt{t-s}}^{x+\sqrt{t-s}} dz\,  G_L^2(t-r;x,z) \ge c_0 (t-s)^{\half}.
\eeq
\end{lemma}

\begin{proof}
Clearly, \eqref{1'.112} implies \eqref{1'.1111}. By the scaling property \eqref{scalingD}, the left-hand side of \eqref{1'.112} is equal to
\beqn
    L  \int_0^{\frac{h}{L^2}} ds \int_{\frac{x - \sqrt{h}}{L}}^{\frac{x + \sqrt{h}}{L}} dy\, G_1(s, \frac{x}{L}, y).
    \eeqn
The hypotheses of this lemma with $L = 1$ are satisfied by $\frac{x}{L}$ and $\frac{h}{L^2}$. Assuming that \eqref{1'.112} holds when $L = 1$, we conclude that the left-hand side of \eqref{1'.112}  is bounded below by
    $L c_0 (\frac{h}{L^2})^\half = c_0 h^\half$,
which is the right-hand side of \eqref{1'.112}. It remains to prove that \eqref{1'.112}  holds when $L = 1$.

  Assume that $L = 1$ and the hypotheses of the lemma are satisfied with $L = 1$. Apply Lemma \ref{ch1'-ss2.1-l40} to write
\beqn
G_1(r;x,z)=\frac{1}{\sqrt{4 \pi r}}\exp \left(-\frac{{(x-z)}^{2}}{4r}\right)+H_1(r;x,z),\quad z\in[0,1],
\eeqn
with $H_1$ defined in \eqref{HH-previous}. It follows that for $x \in [\alpha, (1-\alpha)]$ and $0 \leq h \leq \min(x^2, (1-x)^2)$, since $x \pm \sqrt{h} \in [0,1]$,
\beq
\label{proces3}
\int_0^h dr \int_{x-\sqrt{h}}^{x+\sqrt{h}} dz\, G_1^2(r;x,z) \ge \half I_1(h,x)-I_2(h,x),
\eeq
with
\begin{align*}
I_1(h,x)&= \int_0^h dr \int_{x-\sqrt{h}}^{x+\sqrt{h}} dz\, \frac{1}{4 \pi r}\exp \left(-\frac{{(x-z)}^{2}}{2r}\right),\\
I_2(h,x)&= \int_0^h dr \int_{x-\sqrt{h}}^{x+\sqrt{h}} dz\, H_1^2(r;x,z).
\end{align*}
Next, we will find a lower bound for $I_1(h,x)$ and an upper bound for $I_2(h,x)$.

According to \eqref{rdlemC.2.2-before-1}, for a random variable $Z$ with distribution ${\text N}(0,\sigma^2)$, for any $a\ge 0$,
\beqn
P(-a\le Z\le a) \ge 1-\sqrt{\frac{2}{\pi}}\frac{\sigma}{a} \exp\left(-\frac{a^2}{2\sigma^2}\right)\ge 1-\sqrt{\frac{2}{\pi}}\frac{\sigma}{a}.
\eeqn
It follows that
\begin{align}
\half I_1(h,x)
&=\frac{1}{2\sqrt{2\pi}}\int_0^h dr \ \frac{1}{2\sqrt r} \int_{x-\sqrt h}^{x+\sqrt{h}} dz\ \frac{1}{\sqrt{2 \pi r}}\exp \left(-\frac{{(x-z)}^{2}}{2r}\right)\notag\\
&\ge \frac{1}{2\sqrt{2\pi}}\int_0^h dr \ \frac{1}{2\sqrt r} \left(1-\left(\frac{2r}{\pi h}\right)^{\half}\right)\notag\\
& = \frac{1}{2\sqrt{2\pi}}\left(\sqrt h-\sqrt{\frac{h}{2\pi}}\right)=  \frac{1}{2\sqrt{2\pi}}\left(1- \frac{1}{\sqrt{2\pi}}\right)\sqrt h \notag \\
 &= 2c_0 \sqrt h.
 \label{1'.113}
\end{align}

Fix $\alpha\in \left]0,\frac{1}{2}\right[$. By Lemma \ref{ch1'-ss2.1-l40}, $H_1(r, x, z)$ is bounded over $[0,1] \times [\alpha, 1 - \alpha] \times [0,L]$, so
\beqn
I_2(h,x) \le c(\alpha) h^{\frac{3}{2}}, \qquad x\in[\alpha,1-\alpha].
\eeqn
Therefore, by \eqref{proces3}, for $x \in [\alpha, (1-\alpha)]$ and $0 \leq h \leq \min\left(\frac{c_0}{c(\alpha)}, x^2, (1-x)^2\right)$,
\begin{align*}
\int_0^h dr \int_{x-\sqrt h}^{x+\sqrt h} dz\, G_1^2(r; x,z) &\ge 2c_0 h^\half - c(\alpha)h^{\frac{3}{2}}\\
 & \ge c_0 h^\half\left(2-\frac{c(\alpha)}{c_0}h\right)\\
 & \ge c_0 h^\half,
 \end{align*}
 since $h\le\frac{c_0}{c(\alpha)}$.
 This completes the proof of the lemma. 
 \end{proof}

\noindent{\em Increments in time, $L^2$-norms, lower bounds}
\medskip

The next lemma concerns lower bounds on integrated squared increments in time of the Green's function.
\begin{lemma}
\label{ch1'.1201bis}
(a) For all $ x \in\, ]0,L[$, $t >0$, and $h \in [0,\frac{1}{9 \pi^2} \min(x^2, (L - x)^2)]$,
\beqn
 \int_0^t dr \int_0^L dy\,  [G_L(r+h;x,y) - G_L(r;x,y)]^2\ge  \frac{1 - e^{-2t}}{20\, \pi}\, (1-e^{-1})^2 \, \sqrt{h}.
\eeqn

(b) For all $ x \in\, ]0,L[$ and $t >0$, there is $c(t,x,L) > 0$ such that, for all $h \in \left[0, \frac{L^2}{36 \pi^2}\right]$,
\beqn
 \int_0^t dr \int_0^L dy\,  [G_L(r+h;x,y) - G_L(r;x,y)]^2 \geq c(t,x,L) \sqrt{h}.
\eeqn
\end{lemma}

\begin{proof}We will prove the lemma when $L=\pi$. By the scaling property of the Green's function given in \eqref{scalingD}, the conclusions will extend to any $L>0$.
Define
$$
   A(t,x,h) := \int_0^t dr \int_0^\pi dy\, [G_\pi(r+h;x,y) - G_\pi(r;x,y)]^2.
$$
Using \eqref{B.2Green-D} with $L=\pi$, we see that
  \begin{align*}
   A(t,x,h) &= \int_0^t dr \int_0^\pi dy\, \left[\sum_{n=1}^\infty (e^{-n^2 (r+h)} - e^{-n^2 r}) e_{n,\pi}(x) e_{n,\pi}(y) \right]^2 \\
	  & = \int_0^t dr\, \sum_{n=1}^\infty\, (e^{-n^2 (r+h)} - e^{-n^2 r})^2 e_{n,\pi}^2(x) \\
		&= \sum_{n=1}^\infty e_{n,\pi}^2(x) (1 - e^{-n^2 h})^2 \int_0^t dr\, e^{-2 n^2 r} \\
		&= \sum_{n=1}^\infty e_{n,\pi}^2(x) (1 - e^{-n^2 h})^2\  \frac{1- e^{-2 n^2 t}}{2 n^2} \\
		&\geq \frac{1 - e^{-2t}}{2} \sum_{n=1}^\infty e_{n,\pi}^2(x)\ \frac{(1 - e^{-n^2 h})^2}{n^2} \\
		&=: \hat A(t,x,h).
\end{align*}
Recall that $e_{n,\pi}(x) = \sqrt{2/\pi} \sin(nx)$. We use the inequality $1 - e^{-n^2 h} \geq 1-e^{-1}$ for $n^2 h \geq 1$, that is, $n \geq h^{-1/2}$. Therefore,
$$
   \hat A(t,x,h) \geq \frac{1 - e^{-2t}}{\pi}\, (1-e^{-1})^2 \sum_{n\geq h^{-1/2}} \frac{\sin^2(nx)}{n^2} .
$$

   For $x \in \, ]0,\pi/2]$, define
$$
   B(x,h) = \left\{n \in \IN: n \geq h^{-1/2}\mbox{ and } (nx) \mbox{ mod } 2\pi \in \,\left]\frac{\pi}{4}, \frac{3\pi}{4}\right]\, \cup \, \left]\frac{5\pi}{4}, \frac{7\pi}{4}\right] \right\}.
$$
For those $n$ that belong to $B(x,h)$, $\sin^2(nx) \geq \frac12$. Therefore,
\begin{equation}\label{BoundHeatA}
   \hat A(t,x,h) \geq \frac{1 - e^{-2t}}{2\pi} (1-e^{-1})^2 \sum_{n \in B(x,h)} \frac{1}{n^2}.
\end{equation}
Observe that
$$
   B(x,h) = \bigcup_{k=1}^\infty I_{x,h,k},
$$
where the $I_{x,h,k}$ are intervals of consecutive integers, ordered so that $I_{x,h,k}$ precedes $I_{x,h,k+1}$, that is,
$$
   \ell_1 \in I_{x,h,k}, \ \ell_2 \in I_{x,h,k+1}\quad \Longrightarrow \quad \ell_1 \leq \ell_2.
$$
Denote $J_{x,h,k}$ the interval of consecutive integers between $I_{x,h,k}$ and $I_{x,h,k+1}$, so that $J_{x,h,k}$ precedes $J_{x,h,k+1}$ and
$$
   J_{x,h,0} \cup \left( \bigcup_{k=1}^\infty (I_{x,h,k} \cup J_{x,h,k}) \right) = \{n\in \IN: n \geq h^{-1/2}\},
$$
where $J_{x,h,0}$ accounts for the integers between $h^{-1/2}$ and $\min I_{x,h,1}$. Note that $J_{x,h,0}$ may be empty and also that the cardinals of the sets $J_{x,h,0}$, $I_{x,h,1}$ and $J_{x,h,1}$ are bounded by $\tfrac{\pi}{2x}$. Moreover,
for $x \in \, ]0,\pi/2]$, except possibly for $k=1$, we have $\vert \mbox{card}(I_{x,h,k}) - \mbox{card}(J_{x,h,k})\vert \leq 2$, and since $I_{x,h,k}$ precedes $J_{x,h,k}$,
$$
   \sum_{n \in I_{x,h,k}} \frac{1}{n^2} \geq \frac{1}{4} \sum_{n \in J_{x,h,k}} \frac{1}{n^2}.
$$
Therefore,
\begin{align*}
   \sum_{n\geq h^{-1/2}} \frac{1}{n^2} &= \sum_{n \in J_{x,h,0}}  \frac{1}{n^2} +\sum_{n \in I_{x,h,1}}  \frac{1}{n^2} +\sum_{n \in J_{x,h,1}}  \frac{1}{n^2}\\
   &\qquad +  \sum_{k=2}^\infty \left[\sum_{n \in I_{x,h,k}} \frac{1}{n^2} + \sum_{n \in J_{x,h,k}} \frac{1}{n^2} \right] \\
	   & \leq 3 \frac{\pi}{2x}\, \frac{1}{(h^{-1/2})^2} + 5 \sum_{k=2}^\infty \sum_{n \in I_{x,h,k}} \frac{1}{n^2} \\
		 & \leq \frac{3 \pi h}{2x} + 5 \sum_{n \in B(x,h)} \frac{1}{n^2}.
\end{align*}
It follows that
\begin{align*}
  \sum_{n \in B(x,h)} \frac{1}{n^2} &\geq \frac{1}{5} \sum_{n\geq h^{-1/2}} \frac{1}{n^2} -  \frac{3 \pi h}{10\, x}
	 \geq \frac{1}{5}\int_{h^{-1/2}}^\infty \frac{dz}{z^2} -  \frac{3 \pi h}{10\, x}\\
	& = \frac{1}{5} \sqrt{h} -  \frac{3 \pi h}{10\, x}
	 = \frac{1}{5} \sqrt{h} \left(1 - \frac{3 \pi \sqrt{h}}{2\, x}\right) \\
	& \geq \frac{1}{10} \sqrt{h}
\end{align*}
provided $h \leq x^2/(9 \pi^2)$. Putting this together with \eqref{BoundHeatA}, we obtain
\begin{align*}
   \hat A(t,x,h) &\geq \frac{1 - e^{-2t}}{20\pi}\, (1-e^{-1})^2  \sqrt{h}.
\end{align*}
This proves (a) for $x \in\, ]0, \pi/2]$.

   For $x \in [\pi/2, \pi[$, we use the fact that $G_\pi(t,x,y) = G_\pi(t,\pi - x ,\pi - y)$, so $A(t,x,h)= A(t,\pi-x,h)$, to get back to the case just treated. This proves (a).

 (b)  For $x \in\, ]0, L[$, set
  $\delta_{x} := \frac{1}{9 \pi^2}\min(x^2,(L- x)^2)$ and $c_0(t) = (1 - e^{-2t}) (1-e^{-1})^2/(20 \pi)$. Since $h \mapsto \hat A(t,x,h)$ is non-decreasing, for $\frac{L^2}{36 \pi^2} \geq h \geq \delta_x$, it follows from (a) that
$$
	\hat A(t,x,h) \geq \hat A(t,x,\delta_{x}) \geq c_0(t) \sqrt{\delta_{x}} = c_0(t) \frac{\sqrt{\delta_{x}}}{\sqrt{h}} \sqrt{h} \geq c_0(t) \frac{6\pi}{ L} \sqrt{\delta_x} \sqrt{h}.
$$
Since $\sqrt{\delta_x} \frac{6\pi}{ L} \leq 1$, we can set $c(t,x,L) = c_0(t) \frac{6\pi}{L} \sqrt{\delta_x} \leq c_0(t)$, so from (a), we get for all $h \in \left[0, \frac{L^2}{36 \pi^2}\right]$,
$$
\hat A(t,x,h) \geq c(t,x,L) \sqrt{h}.
$$
This proves (b).
\end{proof}

\noindent{\em Increments in space, $L^2$-norms, lower bounds}

\begin{lemma}
\label{ch1'-1202}
For every $t\in\re_+$ and $x,y\in[0,L]$,
\begin{align}
\notag
J_L(t;x,y) &:= \int_0^t dr \int_0^L dz\, [G_L(r;x,z)-G_L(r;y,z)]^2 \\
 &\ge \frac{1-e^{-2\frac{\pi^2}{L^2} t}}{2} \left( |x-y| - \frac{1}{L} |x-y|^2 \right).
 \label{1'.dlb0}
\end{align}
\end{lemma}
\begin{proof}

Let $t\ge 0$. Using the formula \eqref{B.2Green-D} for the Green's function and Parseval's identity, for any $r>0$ and $x,y\in[0,L]$, we have
\begin{align*}
&\int_0^L dz\, [G_L(r;x,z)-G_L(r;y,z)]^2\\
&\qquad= \frac{2}{L} \sum_{n=1}^\infty\ e^{-2\frac{\pi^2}{L^2}n^2 r}\left(\sin\left(\frac{n\pi}{L}x\right) - \sin\left(\frac{n\pi}{L}y\right)\right)^2.
\end{align*}
Consequently,
\begin{align*}
J_L(t;x,y) &= \frac{2}{L} \sum_{n=1}^\infty\left(\sin\left(\frac{n\pi}{L}x\right) - \sin\left(\frac{n\pi}{L}y\right)\right)^2 \int_0^t dr\,  e^{-2\frac{\pi^2}{L^2}n^2 r}\\
&= L  \sum_{n=1}^\infty  \left(\sin\left(\frac{n\pi}{L}x\right) - \sin\left(\frac{n\pi}{L}y\right)\right)^2\, \frac{1-e^{-2\frac{\pi^2}{L^2}n^2 t}}{\pi^2n^2}\\
&\ge L  \left(1-e^{-2\frac{\pi^2}{L^2}t}\right)S_1,
\end{align*}
where
\beqn
S_1= \sum_{n=1}^\infty\, \frac{[\sin(n\pi x/L)-\sin(n\pi y/L)]^2}{\pi^2\, n^2} = \half\left(|\tfrac{x}{L}-\tfrac{y}{L}|-|\tfrac{x}{L}-\tfrac{y}{L}|^2\right)
\eeqn
by Lemma \ref{lemB.5.2}. This proves \eqref{1'.dlb0}.
\end{proof}
\smallskip

\noindent{\em Increments in space and in time, $L^p$-norms, upper bounds}
\smallskip

\begin{lemma}
\label{app2-r1}
(a)\  Fix $T>0$. For any $p \in\, ]0, 3[$, there is $C_{p,T} < \infty$ such that for all $t \in [0, T]$ and for all $x, y \in [0, L]$,
\begin{align}
\label{app2.3}
&\int_0^t dr \int_0^L dz\, |G_L(r; x,z)-G_L(r;y,z)|^{p}\notag\\
&\qquad\qquad \le C_{p,T} \times
\begin{cases}
           \vert x - y \vert^{p}&            {\text{if}}\ p \in\, ]0, \frac{3}{2}[,\\
           \vert x - y \vert^{\frac{3}{2}} \log(2 + \frac{1}{\vert x - y \vert})&  {\text{if}}\ p =  \frac{3}{2},\\
         \vert x - y \vert^{3 - p}            &  {\text{if}}\ p \in\, ]\frac{3}{2}, 3[.
         \end{cases}
         \end{align}  

(b)\ For any $p \in\, ]0, 3[$, there exists $C_{p, T} < \infty$ such that for all $s, t \in [0, T]$ and $x \in [0, L]$,
\begin{align}
\label{app2.3(*1)}
   &\int_0^T dr \int_0^L dz \vert  G_L(t-r; x, z) - G_L(s-r; x, z) \vert^p\notag\\
    &\qquad\qquad\leq C_{p, T}
    \begin{cases}
           \vert t - s \vert^p&      {\text{if}}\  p \in\, ]0, 1[,\\
       \vert t - s \vert (1 + \log(\frac{t}{ \vert t - s \vert}) 1_{\{\vert t - s \vert < t\}})&     {\text{if}}\  p = 1,\\
       \vert t - s \vert^{\half (3 - p)}&     {\text{if}}\  p \in\, ]1, 3[.    
       \end{cases}
       \end{align}
\end{lemma}
\begin{proof}
(a) We will use the decomposition in Remark \ref{abans-ch1'-tDL-*1}:
\beq
\label{app2.3(*2)}
   G_L(r; x, z) = \Gamma(r, x - z) - \Gamma(r, x + z) - \Gamma(r, x + z - 2L) + \tilde H_1(r; x, z).
   \eeq
This implies that up to a multiplicative constant,  the left-hand side of \eqref{app2.3} is bounded above by the sum of four terms:
\begin{align*}
   A_1 &= \int_0^T dr \int_0^L dz\, \vert \Gamma(r, x - z) - \Gamma(r, y - z)  \vert^p,\\
       A_2 & = \int_0^T dr \int_0^L dz\, \vert \Gamma(r, x + z) - \Gamma(r, y + z) \vert^p,\\
   A_3 & = \int_0^T dr \int_0^L dz\, \vert \Gamma(r, x + z - 2L) - \Gamma(r, y + z - 2L) \vert^p,\\
      A_4 &=  \int_0^T dr \int_0^L dz\, \vert \tilde H_1(r; x, z) - \tilde H_1(r; y, z) \vert^p .
      \end{align*}
For $A_1$, we  replace the bounds in the $dz$-integral by $\pm \infty$, then apply Lemma \ref{app2-l1} with $k= 1$, $p \in ]0, 3[$, to obtain the desired bound for $A_1$ in each of the three cases of \eqref{app2.3}. For $A_2$ and $A_3$, we do the same, noting that once the $dz$-integrals are over $\R$, we can replace $z$ and $z - 2L$ by $-z$, so we get the same bound as for $ A_1$. Finally, for $A_4$, we use Remark \eqref{abans-ch1'-tDL-*1} to bound $A_4$ by $\tilde C T L \, \vert x - y \vert^p$, which is of smaller order (when  $\vert x - y \vert$ is small). This proves (a).

(b)\ We use again the decomposition \eqref{app2.3(*2)} to bound
  \beqn
  \int_0^T dr \int_0^L dz\, \vert  G_L(t-r; x, z) - G_L(s-r; x, z) \vert^p
  \eeqn
(up to multiplicative constant) by the sum of four terms:
\begin{align*}
  B_1 &=  \int_0^\infty dr \int_0^L dz\, \vert \Gamma(t - r, x - z) - \Gamma(s - r, x - z)  \vert^p,\\
  B_2 &=  \int_0^\infty dr \int_0^L dz\, \vert \Gamma(t - r, x + z) - \Gamma(s - r, x + z) \vert^p,\\
  B_3 &=  \int_0^\infty dr \int_0^L dz\, \vert \Gamma(t - r, x + z - 2L) - \Gamma(s - r, x + z - 2L) \vert^p,\\
   B_4 &=  \int_0^T dr \int_0^L dz\, \vert \tilde H_1(t - r; x, z) - \tilde H_1(s - r; x, z) \vert^p.
   \end{align*}
For the first three terms, we replace the bounds in the $dz$-integral by $\pm \infty,$ then apply Lemma \ref{app2-l3} with $k= 1$, to obtain the desired bound for $B_i$, $i = 1, 2, 3$, in each of the three cases of \eqref{app2.3(*1)}. For $B_4$, we use Remark \eqref{abans-ch1'-tDL-*1} to bound the integral by $\tilde C T L\, \vert t - s \vert^p$, which is of smaller order (when  $\vert t - s \vert$ is small). This proves (b).
\end{proof}
\begin{remark}
\label{app2-r11}
Uniform estimates of $L^1$-norms of increments of $G_L$ in space and in time are proved in \cite{D-K-Z016}, Lemma A.3 and Lemma A.4, respectively. More specifically:
\begin{enumerate}
\item There is a finite constant $C>0$ such that, for all $x,y\in [0,L]$ and for all $t>0$,
\begin{align}
\label{nonameD}
&\int_0^t dr \int_0^L dz\, |G_L(t-r;x,z)-G_L(t-r;y,z)|\notag\\
&\qquad\qquad \le C |x-y|\log\left(e\vee |x-y|^{-1}\right).
\end{align}
\item There is a finite constant $C>0$ such that for all $h\ge 0$,
\begin{align}
\label{nonameD-bis}
&\sup_{t>0} \sup_{x\in[0,L]} \int_0^t dr \int_0^L dz\, |G_L(t+h-r;y,z)-G_L(t-r;y,z)|\notag\\
&\qquad\qquad \le C h^{\frac{1}{2}}.
\end{align}
\end{enumerate}
\end{remark}

\section{Heat kernel with Neumann boundary conditions}
\label{app2-3}

The content of this section is similar to that of Section \ref{app2-2}. However, the results here refer to the Green's function of the heat operator with Neumann boundary conditions,
which behaves quite differently at the boundary than the Green's function with Dirichlet boundary conditions. We recall from \eqref{1'.400} the formula
\beq
\label{B.2Green-N}
 G_L(t;x,y)= \sum^{\infty}_{n=0} e^{-\frac{\pi^2}{L^2}n^{2} t}g_{n,L} (x) g_{n,L} (y),\quad t>0, \ x,y\in[0,L],
\eeq
where
\beqn
g_{0,L}(x)= \frac{1}{\sqrt{L}}, \qquad g_{n,L}(x)= \sqrt{\frac{2}{L}} \cos\left(\frac{n\pi}{L}x\right), \qquad n\ge 1.
\eeqn
\smallskip

\noindent{\em Increments in space and in time, $L^2$-norms, upper bounds}
\medskip

The next lemma is the analogue of Lemma \ref{ch1'-l2} in the case of vanishing Neumann boundary conditions.

\begin{lemma}
\label{ch1'-l3}
The Green's function defined in \eqref{B.2Green-N} satisfies the following properties.
\begin{enumerate}
\item For any  $x,y\in [0,L]$,
\beq
\label{1'.17}
\int_0^\infty dr \int_0^L dz\,  [G_L(r; x,z) - G_L(r; y,z)]^2 = \frac{1}{2}\,|x-y|.
\eeq
\item For any $x\in [0,L]$ and $h\ge0$,
\beq
\label{1'.18}
\int_0^\infty dr \int_0^L dz\, [G_L(r+h; x,z)-G_L(r;x,z)]^2 \le \frac{3}{\pi} h^{\half}.
\eeq
\item Fix $T>0$. There is a finite constant $C=C(T,L)$ such that, for all $t\in[0,T]$ and $x\in[0,L]$,
\beq
\label{1'.19}
\int_0^t dr \int_0^L dz\, G_L^2(r;x,z)\le C\, t^\half.
\eeq
\end{enumerate}
As a consequence, there is a finite constant $C=C(T,L)$ such that, for any $s,t\in[0,T]$ and $x, y\in [0,L]$,
\begin{align}\nonumber
  &\int_0^t dr\ \int_0^L dz\, [G_L(t-r;x,z)-G_L(s-r;y,z)]^2 \\
   &\qquad\qquad\qquad \le C \left(|t-s|^\frac{1}{4} + |x-y|^\half \right)^2.
\label{1'11000}
\end{align}
\end{lemma}
\begin{proof}
The Green's function \eqref{B.2Green-N} satisfies the scaling property \eqref{scalingN}, which is the same as for the Green's function corresponding to Dirichlet boundary conditions. Hence, the equalities \eqref{scaling-1} and \eqref{scaling-2} hold and therefore, without loss of generality, we may take $L=1$.

We start with the proof of \eqref{1'.17} (with $L=1$).  By the expression \eqref{B.2Green-N}, using the fact that the sequence $(g_{n,1},\, n\ge 1)$ is orthonormal in $L^2([0,1])$, we see that for $x,y\in[0,L]$,
\begin{align*}
&\int_0^\infty dr \int_0^1 dz [G_1(r; x,z) - G_1(r; y,z)]^2 \\
&\quad = \int_0^\infty dr \int_0^1 dz \left(\sum_{n=1}^\infty e^{-\pi^2n^2r} [g_{n,1}(x) - g_{n,1}(y)] g_{n,1}(z)\right)^2\\
&\quad = \sum_{n=1}^\infty\, [g_{n,1}(x) - g_{n,1}(y)]^2 \int_0^\infty dr \, e^{-2\pi^2n^2r}\\
&\quad = \sum_{n=1}^\infty\, \frac{[g_{n,1}(x) - g_{n,1}(y)]^2}{2\pi^2 n^2}
 = \sum_{n=1}^\infty\,  \frac{\left(\cos(n\pi x)-\cos(n\pi y)\right)^2}{\pi^2\, n^2}\\
&\quad = \frac{1}{2}\, |x-y|,
\end{align*}
where in the last equality, we have applied Lemma \ref{lemB.5.2}.
This proves \eqref{1'.17} for all $L>0$.

We now prove \eqref{1'.18}.
Using formula \eqref{B.2Green-N}, we see that
\begin{align*}
&\int_0^\infty dr \int_0^1 dz\, [G_1(r+h; x,z)-G_1(r;x,z)]^2\\
&\qquad \qquad\le 2\int_0^\infty dr \, \sum_{n=1}^\infty \left(e^{-\pi^2n^2(r+h)}-e^{-\pi^2n^2 r}\right)^2.
\end{align*}
We have already seen this expression in \eqref{1'.1200} and have bounded it by $\frac{3}{\pi} h^\half$. This ends the proof of  \eqref{1'.18}.

The upper bound \eqref{1'.19} comes from the upper bound on $G_L$ given in \eqref{compneumann}, which compares $G_L$ with a heat kernel, and  \eqref{app2.00}.

The upper bound \eqref{1'11000} follows from \eqref{1'.17}--\eqref{1'.19} by applying the triangle inequality.
\end{proof}

\smallskip

\noindent{\em Increments in space and in time, $L^2$-norms, lower bounds}
\medskip

Next, we present some results concerning lower bounds on integrated squared increments of the Green's function with Neumann boundary conditions. Notice that these bounds apply up to and including the boundary points $0$ and $L$.

\begin{lemma}
\label{ch1-l4}
\begin{enumerate}
\item For all $t\in\re_+$ and $x,y\in[0,L]$,
\beq
\label{1'.2000}
   \frac{1-e^{-2\pi^2 t}}{2}\, |x-y| \le \int_0^t dr \int_0^L dz\ [G_L(r;x,z)-G_L(r;y,z)]^2.
\eeq
\item
For all $t\in\re_+$ and $x\in [0,L]$,
\beq
\label{1'.20}
\frac{1}{\sqrt{2\pi}}\,  t^{\half} \le \int_0^t dr \int_0^L dz\, G_L^2(r;x,z).
\eeq
\item Fix $t>0$.  For all $x \in [0,L]$ and $h\in[0,1]$,
\beq
\label{1'.21}
C h^{\half} \le \int_0^t dr \int_0^L dz\, [G_L(r+h;x,z)-G_L(r;x,z)]^2.
\eeq
with $C=\frac{1}{20\pi}(1-e^{-1})^2 \left(1-\exp(-\tfrac{2\pi^2}{L^2}t)\right)$.
\end{enumerate}
\end{lemma}
\begin{proof}
We begin by checking \eqref{1'.2000}.
As in Lemma \ref{ch1'-l3}, we will prove the statement for $L=1$ and then extend its validity to all $L>0$ by applying the scaling properties \eqref{scaling-1} and \eqref{scaling-2}.

Let $L=1$ and $x,y\in[0,1]$. Set
\beqn
 I(t;x,y)= \int_0^t dr \int_0^1 dz\, [G_1(r;x,z)-G_1(r;y,z)]^2.
\eeqn
Using \eqref{B.2Green-N},
\begin{align}
\label{proces5}
I(t;x,y)& = \int_0^t dr \int_0^1 dz \left(\sum_{n=0}^\infty e^{-\pi^2n^2 r}\, [g_{n,1}(x)-g_{n,1}(y)]\, g_{n,1}(z)\right)^2\notag\\
&= \int_0^t dr\ \sum_{n=1}^\infty e^{-2\pi^2n^2 r}[g_{n,1}(x)-g_{n,1}(y)]^2\notag\\
&=\sum_{n=1}^\infty\, [g_{n,1}(x)-g_{n,1}(y)]^2\, \frac{1-e^{-2\pi^2n^2 t}}{2\pi^2\, n^2}\notag\\
&\ge \left(1-e^{-2\pi^2 t}\right)\sum_{n=1}^\infty\, \frac{[\cos(n\pi x)-\cos(n\pi y)]^2}{\pi^2\, n^2}.
\end{align}
From Lemma \ref{lemB.5.2}, we conclude that \eqref{1'.2000} holds for $L=1$ and therefore also for arbitrary $L>0$.
\medskip

We now prove \eqref{1'.20}.
 Using the properties of the Green's function stated in Proposition \ref{ch1'-pPN}, we have
\begin{align*}
\int_0^t dr \int_0^1 dz\, G_1^2(r;x,z) &= \int_0^t dr\, G_1(2r; x,x)\ge \int_0^t dr\, \Gamma(2r,0)\\
& = \int_0^t \frac{dr}{\sqrt{8\pi r}} =
\left(\frac{t}{2\pi}\right)^\half.
\end{align*}
This proves \eqref{1'.20} for $L=1$, hence for all $L>0$.
\medskip

Next we prove \eqref{1'.21}. We will first consider the case $L= \pi$ to make the calculations more transparent;  then, by the scaling property of the Green's function, we will obtain \eqref{1'.21} for any $L>0$.
Define
$$
   A(t,x,h) := \int_0^t dr \int_0^\pi dz\, [G_\pi(r+h;x,z) - G_\pi(r;x,z)]^2.
$$
Since \eqref{1'.21} clearly holds for $h=0$, we assume that $h>0$. Then
\begin{align*}
   A(t,x,h) &= \int_0^t dr \int_0^\pi dz\, \left[\sum_{n=1}^\infty\, (e^{-n^2 (r+h)} - e^{-n^2 r})\, g_{n,\pi}(x)\, g_{n,\pi}(z) \right]^2 \\
	  & = \int_0^t dr\, \sum_{n=1}^\infty\, (e^{-n^2 (r+h)} - e^{-n^2 r})^2 g_{n,\pi}^2(x) \\
		&= \sum_{n=1}^\infty\, g_{n,\pi}^2(x) (1 - e^{-n^2 h})^2 \int_0^t dr\, e^{-2 n^2 r} \\
		&= \sum_{n=1}^\infty\, g_{n,\pi}^2(x)\, (1 - e^{-n^2 h})^2\, \frac{1- e^{-2 n^2 t}}{2 n^2} \\
		&\geq \frac{1 - e^{-2t}}{2}\, \sum_{n=1}^\infty\, g_{n,\pi}^2(x)\, \frac{(1 - e^{-n^2 h})^2}{n^2}
\end{align*}
Recall that $g_{n,\pi}(x) = \sqrt{2/\pi} \cos(nx)$ for $n\ge 1$. We use the inequality $1 - e^{-n^2 h} \geq 1-e^{-1}$ for $n^2 h \geq 1$, that is, $n \geq h^{-1/2}$, to deduce that
\beq\label{rdeB3.12}
    A(t,x,h) \geq \frac{1 - e^{-2t}}{\pi}\, (1-e^{-1})^2 \sum_{n\geq h^{-1/2}} \frac{\cos^2(nx)}{n^2}
    \eeq
(here, we use $h\in\, ]0,1]$).

Assume first that  $x \in \, ]0,\pi/2]$ and define
\begin{align*}
   &B(x,h)\\
   &\qquad  = \left\{n \in \IN: n \geq h^{-1/2}\mbox{ and } (nx) \mbox{ mod } 2\pi \in\, ]0,\tfrac{\pi}{4}]\,\cup\,] \tfrac{3\pi}{4}, \tfrac{5\pi}{4}]\, \cup \, ] \tfrac{7\pi}{4},2\pi] \right\}.
\end{align*}
For the $n$ that belong to $B(x,h)$, $\cos^2(nx) \geq \frac12$. Therefore,
\begin{equation}\label{eq1}
   A(t,x,h) \geq \frac{1 - e^{-2t}}{2\pi}\, (1-e^{-1})^2 \sum_{n \in B(x,h)} \frac{1}{n^2}.
\end{equation}
Observe that
$$
   B(x,h) = \bigcup_{k=1}^\infty I_{x,h,k},
$$
where the $I_{x,h,k}$ are intervals of consecutive integers, ordered so that $I_{x,h,k}$ precedes $I_{x,h,k+1}$, that is,
$$
   \ell_1 \in I_{x,h,k}, \ \ell_2 \in I_{x,h,k+1}\quad \Longrightarrow \quad \ell_1 \leq \ell_2.
$$
Denote $J_{x,h,k}$ the interval of consecutive integers between $I_{x,h,k}$ and $I_{x,h,k+1}$, so that $J_{x,h,k}$ precedes $J_{x,h,k+1}$ and
$$
   J_{x,h,0} \cup \left( \bigcup_{k=1}^\infty (I_{x,h,k} \cup J_{x,h,k}) \right) = \{n\in \IN: n \geq h^{-1/2}\},
$$
where $J_{x,h,0}$ accounts for the integers between $h^{-1/2}$ and $\min I_{x,h,1}$. Note that $J_{x,h,0}$ may be empty and also that the cardinals of the sets $J_{x,h,0}$, $I_{x,h,1}$ and $J_{x,h,1}$ are bounded by $\tfrac{\pi}{2x}$. Moreover, for $x \in \, ]0,\pi/2]$, except possibly for $k=1$, we have $\vert \mbox{card}(I_{x,h,k}) - \mbox{card}(J_{x,h,k})\vert \leq 2$  and since $I_{x,h,k}$ precedes $J_{x,h,k}$,
$$
   \sum_{n \in I_{x,h,k}} \frac{1}{n^2} \geq \frac{1}{4} \sum_{n \in J_{x,h,k}} \frac{1}{n^2}.
$$
Therefore,
\begin{align*}
   \sum_{n\geq h^{-1/2}} \frac{1}{n^2} &= \sum_{n \in J_{x,h,0}}  \frac{1}{n^2} +\sum_{n \in I_{x,h,1}}  \frac{1}{n^2}\\
   & +\sum_{n \in J_{x,h,1}}  \frac{1}{n^2} +  \sum_{k=2}^\infty \left[\sum_{n \in I_{x,h,k}} \frac{1}{n^2} + \sum_{n \in J_{x,h,k}} \frac{1}{n^2} \right] \\
	   & \leq 3 \frac{\pi}{2x}\, \frac{1}{(h^{-1/2})^2} + 5\, \sum_{k=2}^\infty\, \sum_{n \in I_{x,h,k}} \frac{1}{n^2} \\
		 & \leq \frac{3 \pi h}{2x} + 5 \sum_{n \in B(x,h)} \frac{1}{n^2}.
\end{align*}
It follows that
\begin{align*}
  \sum_{n \in B(x,h)} \frac{1}{n^2} &\geq \frac{1}{5} \sum_{n\geq h^{-1/2}} \frac{1}{n^2} -  \frac{3 \pi h}{10\, x}\\
	& \geq \frac{1}{5}\int_{h^{-1/2}}^\infty \frac{dz}{z^2} -  \frac{3 \pi h}{10\, x}
	 = \frac{1}{5} \sqrt{h} -  \frac{3 \pi h}{10\, x} \\
	& = \frac{1}{5} \sqrt{h} \left(1 - \frac{3 \pi \sqrt{h}}{2\, x}\right)
	 \geq \frac{1}{10} \sqrt{h}
\end{align*}
provided $h \leq x^2/(9 \pi^2)$. Putting this together with \eqref{eq1}, we obtain
\begin{align*}
  A(t,x,h) &\geq \frac{1 - e^{-2t}}{20\pi}\, (1-e^{-1})^2\,  \sqrt{h}.
\end{align*}
This proves \eqref{1'.21} for $x \in\, ]0, \pi/2]$.

For $x=0$, from \eqref{rdeB3.12}, we see that
\begin{align}
\label{lowA}
   A(t,x,h) & \geq \frac{1 - e^{-2t}}{\pi}\, (1-e^{-1})^2 \sum_{n \geq h^{-1/2}} \frac{1}{n^2}\notag \\
    &\geq \frac{1 - e^{-2t}}{2\pi}\, (1-e^{-1})^2 \int_{h^{-1/2}}^\infty \frac{dz}{z^2}\notag\\
    &= \frac{1 - e^{-2t}}{2\pi}\, (1-e^{-1})^2\, \sqrt{h}.
\end{align}
This proves \eqref{1'.21} for $x=0$.

For $x \in\, ]\pi/2, \pi]$, we use that $G_\pi(t,x,z) = G_\pi (t,\pi - x ,\pi - z)$, since $\cos(nx) = \cos(n(\pi-x))$. So, $A(s,x,h)= A(s,\pi-x,h)$ and we go back to the case just treated.

This proves  \eqref{1'.21} for $L=\pi$, hence for all $L>0$ by the scaling property of $G_L$.
\end{proof}



\begin{remark}
\label{app2-r111}
The estimates on $L^p$-increments (respectively, $L^1$-increments) of the Green's function of the heat operator with Dirichlet boundary conditions quoted in Lemma \ref{app2-r1} (respectively, Remark \ref{app2-r11}) are also satisfied in the case of Neumann boundary conditions, with the same proof (respectively see \cite[Lemmas A.3 and A.4]{dkn2013}, where the proof, with slight modifications, also applies to Neumann boundary conditions).
\end{remark}

\section{Comparing heat kernels on $[-L, L]$ and on $\R$}
\label{rdtechnical}

In the following lemma, $\Gamma_L$ is the Green's function for the heat equation on $[-L, L]$ with Dirichlet boundary conditions, that is,
\[
    \Gamma_L(t;x,y) = G_{2L}(t;x+L,y+L), \qquad t > 0, \  x, y \in [-L, L],
\]
where $G_{2L}$ is defined in \eqref{ch1.600} (with $L$ there replaced by $2L$).

\begin{lemma}
\label{rd08_10l3}
Fix $T > 0$. For $t \in \, ]0, T]$ and $x, y \in [-L,L]$, let
\beqn
   H_L(t; x, y) = \Gamma(t, x-y) - \Gamma_L(t; x, y).
\eeqn
Then $H_L(t; x, y) \geq 0$ and
\beq
\label{rd08_10e14}
   \int_{-L}^L dy\, H_L(t; x, y) \leq \half \left[\exp\left(-\tfrac{(L-x)^2}{4t} \right) + \exp\left(-\tfrac{(L+x)^2}{4t} \right) \right],
\eeq

\beq
\label{rd08_10e16}
   \int_0^t ds \int_{-L}^L dy \,  H_L(s; x, y)  \leq \frac{t}{2} \left[\exp\left(-\tfrac{(L-x)^2}{4t} \right) + \exp\left(-\tfrac{(L+x)^2}{4t} \right) \right],
\eeq
and
\beq
\label{rd08_10e15}
   \int_0^t ds \int_{-L}^L dy \,  H_L^2(s; x, y)  \leq \frac{\sqrt{t}}{\sqrt{4\pi}} \left[\exp\left(-\tfrac{(L-x)^2}{4t} \right) + \exp\left(-\tfrac{(L+x)^2}{4t} \right) \right]^2.
\eeq

\end{lemma}

\begin{proof}
The nonnegativity of $H_L$ is mentioned in Proposition \ref{ch1'-pPD} (ii). Let $(t, x) \in\, ]0, T] \times [-L,L]$. 
With the change of variables $z = y + L$, we see that
\beq \label{rd01_06e1}
   \int_{-L}^L dy\, H_L(t; x, y) = \int_0^{2 L} dz\, \Gamma(t, x + L - z) - \int_0^{2 L} dz\, \Gamma_L(t; x, z - L).
\eeq
Notice that $\Gamma_L(t; x, z - L) = G_{2 L}(t; x+L, z)$. By Proposition \ref{ch1-equivGD-bis}, the right-hand side of \eqref{rd01_06e1} is equal to 
\beqn
    \int_0^{2 L}  dz\, \Gamma(t, x + L - z) - \int_{-\infty}^{\infty} dz\, (1_{[0, 2 L]})^{o, p}(z) \, \Gamma(t, x + L - z),
\eeqn
where $(1_{[0, 2 L]})^{o, p} = v^{o, p}$ as defined in \eqref{vop} with $ v = 1_{[0, 2 L]}$. This difference is equal to
\begin{align*}
  & - \int_{\R \setminus [0, 2 L]} dz\,  (1_{[0, 2 L]})^{o, p}(z))\, \Gamma(t, x + L - z) \\
   &\qquad\qquad\leq \int_{\R \setminus [0, 2 L]}  dz\,  \Gamma(t, x + L - z)  \\
    &\qquad\qquad= \int_{x + L}^\infty \Gamma (t, y) dy + \int_{- \infty}^{x - L} \Gamma(t, y) dy.
\end{align*}
Use inequality \eqref{rdlemC.2.2-before-1} to see that this is bounded above by
\beqn
    \half \left[\exp\left(-\tfrac{(L-x)^2}{4t} \right) +  \exp\left(-\tfrac{(L+x)^2}{4t} \right)\right].
\eeqn
This proves \eqref{rd08_10e14}.

For \eqref{rd08_10e16}, it suffices to notice that for $s \leq t$, $\exp(-z^2/(4s)) \leq \exp(-z^2/(4t))$ and to integrate \eqref{rd08_10e14} from $0$ to $t$.

 For \eqref{rd08_10e15}, we first prove that if $x, y \in [-L, L]$, then for $s > 0$,
\beq
\label{rd08_12e6}
    H_L(s; x, y) \leq  \Gamma(s, x+L) + \Gamma(s, x-L).
\eeq
Indeed,
\begin{align*}
   H_L(s; x, y) &= \Gamma(s, x+ y + 2L) + \Gamma(s, x+ y - 2L) \\
    &\qquad - \sum_{m=1}^\infty (\Gamma(s, x - y + 4mL)- \Gamma(s, x+ y + (4m+ 2) L))  \\
    &\qquad - \sum_{m=1}^{\infty} (\Gamma(s, x - y - 4mL) - \Gamma(s, x+ y - (4m + 2) L)),
\end{align*}
and both series are nonnegative. Indeed, for $m \in \N^*$ and $x, y \in [-L, L]$, $0\le x - y + 4mL \leq x+ y + (4m+ 2) L$, showing that the first series is nonnegative.
Furthermore, $0 \leq 4mL + y-x\le (4m+2)L - (x+y)$. 
Therefore, the second series is nonnegative as well.
Since $0 \leq x+L \le x+y+2L$ and $x+y-2L \leq x - L \leq 0$, we obtain \eqref{rd08_12e6}.


It follows that
\begin{align*}
\int_{-L}^L dy \,  H_L^2(s; x, y) &\leq \sup_{z \in [-L, L]} H_L(s; x, z) \int_{-L}^L dy \,  H_L(s; x, y) \\
   &\leq  (\Gamma(s, x-L) + \Gamma(s, x+L)) \\
   &\qquad \times \half \left[\exp\left(-\tfrac{(L-x)^2}{4s} \right) + \exp\left(-\tfrac{(L+x)^2}{4s} \right) \right] \\
   &= \half (4 \pi s)^{-\half} \left[\exp\left(-\tfrac{(L-x)^2}{4s} \right) + \exp\left(-\tfrac{(L+x)^2}{4s} \right) \right]^2,
\end{align*}
where, in the second inequality, we have used \eqref{rd08_10e14}. Replacing $s$ by $t$ in the exponentials and integrating from $0$ to $t$, we obtain \eqref{rd08_10e15}.
\end{proof}

\section{Fractional heat kernel}
\label{app2-4}

In this section, we consider the fundamental solution of the fractional heat kernel. We distinguish two different cases: (1) the equation is given in spatial dimension $1$ and the operator contains a Riesz-Feller fractional derivative; and (2) the equation is given in spatial dimension $k\geq 1$ and the operator contains a fractional Laplacian. The goal here is to establish estimates on integrated squared increments in time and in space of the corresponding fundamental solutions.

\subsection{Spatial dimension $1$}\label{rd04_24ss1}

Let $\mathcal{L} = \tfrac{\partial}{\partial t} - \null_xD_\delta^a$ be the fractional heat operator on $\R_+ \times \R$, where $a\in\,]1,2[$, $|\delta|\le 2-a$ and $\null_xD_\delta^a$ is the Riesz-Feller fractional derivative defined in \eqref{ch1'-ss7.3.1} (the case $a=2$ corresponds to the heat operator and is discussed in Appendix \ref{app2-1}). The associated fundamental solution is the function $\laplacef$ given in \eqref{fs-fractional}. The estimates in this subsection are used in the proof of \eqref{ch1'-ss7.3.102}.
\medskip

\noindent{\em Increments in space and in time, $L^2$-norms, upper bounds}

 \begin{prop}
 \label{ch1'-ss7.3-p1}
 Fix $a\in\, ]1,2[$ and $|\delta|\le 2-a$. There are three constants $C_i$, $i=1,2,3$,
 such that:
 \begin{enumerate}
 \item For all $t>0$ and $x,y\in\re$,
 \beq
 \label{ch1'-ss7.3.15}
 \int_0^t dr \int_\re dz\, \left[\null_\delta G_a(t-r,x-z) - \null_\delta G_a(t-r,y-z)\right]^2 \le C_1\, |x-y|^{a-1}.
 \eeq
 \item For all $0\le s\le t$ and $x\in\re$,
 \beq
 \label{ch1'-ss7.3.16}
 \int_0^s dr \int_{\re} dz\, \left[\null_\delta G_a(t-r,x-z) - \null_\delta G_a(s-r,x-z)\right]^2 \le C_2\, (t-s)^{1-1/a},
 \eeq
 and
 \beq
 \label{ch1'-ss7.3.17}
 \int_s^t dr \int_{\re} dz\, \left[\null_\delta G_a(t-r,x-z)\right]^2 \le C_3\, (t-s)^{1-1/a}.
 \eeq
 \end{enumerate}
 \end{prop}
 \begin{proof}
1.  We recall from \eqref{fs-fractional} and \eqref{ch1'-ss7.3.3} that  for $t>0$,
 \beq
 \label{ch1'-ss7.3.18}
 \tf \null_\delta G_a(t,\ast)(\xi)= \exp\left( t\ \null_\delta\psi_a(\xi)\right)= \exp\left(-t|\xi|^a e^{-i\pi\delta\ {\rm sgn}(\xi)/2}\right),
 \eeq
 and we set
  \beqn
  \beta:= \pi\delta\ {\rm sgn}(\xi)/2 \in\left]-\frac{\pi}{2},\frac{\pi}{2}\right[.
  \eeqn

 Consider the integral on the left-hand side of \eqref{ch1'-ss7.3.15}. By Plancherel's theorem and \eqref{ch1'-ss7.3.18}, it is equal to
 \begin{align*}
 & \frac{1}{2\pi}\int_0^t dr \int_{\re} d\xi \left\vert e^{-i\xi x-(t-r)|\xi|^a e^{-i\beta}}-e^{-i\xi y-(t-r)|\xi|^a e^{-i\beta}}\right\vert^2\\
  &\qquad= \frac{1}{2\pi}\int_0^t dr \int_{\re} d\xi\, \exp\left(-2(t-r)|\xi|^a \cos\beta\right)\left\vert e^{-i\xi x}-e^{-i\xi y}\right\vert^2\\
  &\qquad=\frac{1}{\pi}\int_0^t dr \int_{\re} d\xi\,  \exp\left(-2(t-r)|\xi|^a \cos\beta\right)[1-\cos(\xi(x-y))] \\
  &\qquad= \frac{1}{\pi}\int_{\re} d\xi\,  \frac{1-\exp\left(-2t|\xi|^a \cos\beta\right)}{2|\xi|^a\cos\beta}\, [1-\cos(\xi(x-y))],
  \end{align*}
where the last term is obtained after integrating over $r$. With the change of variables $\xi = u/|x-y|$, we see that this is equal to
\begin{align*}
&\frac{1}{\pi} |x-y|^{a-1} \int_\re \frac{1-\exp\left(-2t\, |u|^a (\cos\beta)/|x-y|^a\right)}{2\, |u|^a \cos\beta} \, (1-\cos u)\, du\\
&\qquad \le C_1\, |x-y|^{a-1},
 \end{align*}
 where
 \beqn
 C_1 := \int_{\re} \frac{1-\cos u}{2\pi\, |u|^a \cos\beta}\, du < \infty,
 \eeqn
 since $a>1$. This proves \eqref{ch1'-ss7.3.15}.
 \medskip

 2. In several places, we will make use of the formula
   \beq
   \label{ch1'-ss7.3.19}
  \int_\re d\xi\, e^{-z|\xi|^a} = 2\int_0^\infty d\xi\, e^{-z|\xi|^a} = 2z^{-1/a}\, \Gamma(1+1/a),
  \eeq
  valid for any $z\in \mathbb{C}$ with $\mathrm{Re}(z)>0$ (see e.g.~\cite[Equation 5.9.1, p.~139]{olver}).

  Let $a^\ast$ be defined by $1/a + 1/a^{\ast} = 1$.
  As in part 1., after applying Plancherel's theorem, the definition of $\null_\delta\psi_a(\xi)$ and \eqref{ch1'-ss7.3.19} with
  $z=2(t-r)\cos\beta$, we have (since $\cos \beta$  does not depend on $\xi$),
  \begin{align}
  \label{ch1'-ss7.3.40}
  \int_s^t dr \int_\re dz\, \left[\null_\delta G_a(t-r,x-z)\right] ^2& = \frac{1}{2\pi}\int_s^t dr \int_\re d\xi\, e^{-2(t-r)|\xi|^a \cos\beta}\notag\\
  & = \frac{\Gamma(1+1/a)}{2^{1/a}\pi\cos^{1/a}(\beta)} \int_s^t \ \frac{dr}{(t-r)^{1/a}}\notag\\
  & = \frac{a^\ast\Gamma(1+1/a)}{2^{1/a}\pi\cos^{1/a}(\beta)}\, (t-s)^{1/a^\ast},
  \end{align}
  proving \eqref{ch1'-ss7.3.17} with
   \beqn
 C_3 := \frac{a^\ast \Gamma(1+1/a)}{2^{1/a} \pi \cos^{1/a}(\pi\delta/2)}.
 \eeqn

Next, we prove \eqref{ch1'-ss7.3.16}. After applying Plancherel's theorem, the left-hand side of
 \eqref{ch1'-ss7.3.16} is equal to
 \begin{align*}
 I:&=\frac{1}{2\pi} \int_0^s dr \int_\re d\xi \left\vert \exp(-i\xi x - (t-r)|\xi|^a e^{i\beta})\right.\\
 &\left.\qquad \qquad \qquad -\exp(-i\xi x - (s-r)|\xi|^a e^{i\beta})\right\vert^2\\
 &=\frac{1}{2\pi} \int_0^s dr \int_\re d\xi \left\vert \exp(-(t-r)|\xi|^a e^{i\beta})\right.\\
 &\left.\qquad \qquad \qquad  -\exp(- (s-r)|\xi|^a e^{i\beta})\right\vert^2.
  \end{align*}
  To simplify the notation, we set
  \beqn
   A_{r,t}:= (t-r)\,|\xi|^a\cos\beta, \qquad B_{r,t}:=(t-r)\,|\xi|^a\sin\beta.
  \eeqn
  Then one can easily check that
  \begin{align}
  \label{ch1'-ss7.3.42}
  &\left\vert \exp(-(t-r)\, |\xi|^a e^{i\beta}) - \exp(- (s-r)\, |\xi|^a e^{i\beta})\right\vert^2\notag\\
  &\qquad = \left\vert e^{-A_{r,t}}\cos(B_{r,t}) - ie^{-A_{r,t}}\sin(B_{r,t})\right.\notag\\
  &\left. \qquad \qquad -  e^{-A_{r,s}}\cos(B_{r,s})+ i e^{-A_{r,s}}\sin(B_{r,s})\right\vert^2\notag\\
  &\qquad = e^{-2A_{r,t}}+e^{-2A_{r,s}}-2e^{-(A_{r,t}+A_{r,s})}\cos(B_{r,t}-B_{r,s}).
   \end{align}
  Apply \eqref{ch1'-ss7.3.19} with $z= 2(t-r)\cos\beta$ and $z= 2(s-r)\cos\beta$, respectively, to obtain
  \begin{align}
  \label{ch1'-ss7.3.41}
  \int_\re d\xi\, e^{-2A_{r,t}}& = \frac{2^{1/a^\ast}\Gamma(1+1/a)}{\cos^{1/a}(\beta)}\, \frac{1}{(t-r)^{1/a}},\notag\\
  \int_\re d\xi\, e^{-2A_{r,s}}& = \frac{2^{1/a^\ast}\Gamma(1+1/a)}{\cos^{1/a}(\beta)}\, \frac{1}{(s-r)^{1/a}}.
  \end{align}
  For the third term in \eqref{ch1'-ss7.3.42}, notice that
  \begin{align*}
  &e^{-(A_{r,t}+A_{r,s})}\cos(B_{r,t}-B_{r,s})\\
  &\qquad = \exp\left(-\left(\frac{t+s}{2}-r\right)2|\xi|^a \cos\beta\right) \cos\left((t-s)|\xi|^a
\sin\beta\right)\\
&\qquad = \mathrm{Re}\left[\exp\left[-\left(\left(\frac{t+s}{2}-r\right) 2\cos\beta + i(t-s)\sin\beta\right)|\xi|^a\right]\right].
  \end{align*}
 Apply  \eqref{ch1'-ss7.3.19} with $z=\left(\frac{t+s}{2}-r\right) 2\cos\beta + i(t-s)\sin\beta$ to obtain
 \begin{align*}
 &\int_\re d\xi\, \exp\left[-\left(\left(\frac{t+s}{2}-r\right) 2\cos\beta + i(t-s)\sin\beta\right)|\xi|^a\right]\\\
 & \quad = 2\Gamma(1+1/a)\left[\left(\left(\frac{t+s}{2}-r\right) 2\cos\beta + i(t-s)\sin\beta\right)\right]^{-1/a}.
 \end{align*}
 From Lemma \ref{ch1'-ss7.3-l1} below with
 \beqn
 c=1/a,\quad  b=\left(\frac{t+s}{2}-r\right) 2\cos\beta, \quad x=(t-s)^2\sin^2(\beta),
 \eeqn
 we have
 \begin{align*}
  &\mathrm{Re}\left[\left(\frac{t+s}{2}-r\right) 2\cos\beta + i(t-s)\sin\beta\right]^{-1/a}\\
 &\qquad \ge \frac{1}{2^{1/a}\cos^{1/a}(\beta)}\,\frac{1}
 {\left((t+s)/2-r\right)^{1/a}}\\
 &\qquad \qquad\qquad- \frac{(a+1)\sin^2(\beta)}{2a^2\left(2\cos\beta\right)^{2+1/a}}\, \frac{(t-s)^2}{\left((t+s)/2-r\right)^{2+1/a}}.
 \end{align*}
 Therefore,
 \begin{align*}
 &2\int_\re d\xi\, e^{-(A_{r,t}+A_{r,s})}\cos(B_{r,t}-B_{r,s})\\
 &\qquad \ge 2^{1+1/a^\ast}\frac{\Gamma(1+1/a)}{\cos^{1/a}(\beta)}
\frac{1}{\left((t+s)/2-r\right)^{1/a}}\\
 &\qquad\qquad -\frac{2\, \Gamma(1+1/a)\, (a+1)\sin^2(\beta)}{a^2\left(2\cos\beta\right)^{2+1/a}}\frac{(t-s)^2}{\left((t+s)/2-r\right)^{2+1/a}}.
 \end{align*}
 Integrating over $r$ and then applying Lemma \ref{ch1'-ss7.3-l2} below, we obtain
 \begin{align*}
 I&\le \frac{\Gamma(1+1/a)}{2^{1/a}\pi\cos^{1/a}(\beta)}\int_0^s dr\, \left(\frac{1}{(t-r)^{1/a}} + \frac{1}{(s-r)^{1/a}}-\frac{2}{\left[(t+s)/2-r\right]^{1/a}}\right)\\
  &\qquad\qquad + \frac{\Gamma(1+1/a)\, (a+1)\sin^2(\beta)}{\pi a^2\left(2\cos\beta\right)^{2+1/a}}\int_0^s dr\, \frac{(t-s)^2}{\left((t+s)/2-r\right)^{2+1/a}}\\
 &\le C_2(t-s)^{1/a^\ast},
  \end{align*}
 with
  \beqn
 C_2 := \left( 2^{1/a}-1\right)C_3 + \frac{\Gamma(1+1/a)\sin^2(\pi\delta/2)}{2\pi a[\cos(\pi\delta/2)]^{2+1/a}}.
 \eeqn
 This proves \eqref{ch1'-ss7.3.16} and ends the proof of the proposition.
 \end{proof}
\medskip

\noindent{\em Upper and lower tail bounds at $\pm \infty$}
\medskip

\begin{lemma}
\label{rd01_21l1}
$(1)$  
 Assume that $a\in \;]1,2[$ and $|\delta| < 2-a$. There is a constant $K_{a} < \infty$ such that for all $x \in \R$,
\begin{align}
\label{E:G-bd}
   0 \leq   \laplacef(1,x) & \le \frac{K_{a}}{1+|x|^{1+a}}\, .
\end{align}
Moreover, for all $T\ge t>0$ and $x\in\R$,
\begin{align}\label{E:Gtx-bd}
 0 \leq   \laplacef(t,x) \le t^{-\frac{1}{a}} \frac{K_{a}}{1+|t^{-1/a}\, x|^{1+a}}
     \le  K_{a}\, t^{-\frac{1}{a}}\, \frac{(T\vee 1)^{1+\frac{1}{a}}}{1+|x|^{1+a}}\:.
\end{align}

$(2)$ 
Assume that  $a\in \;]1,2[$ and $|\delta|< 2-a$. There is a constant $C_{a, \delta} > 0$ such that
for all $t>0$ and $x\in\R$,
\begin{align}\label{E:GlowBd}
   \laplacef(t,x)\ge \frac{C_{a, \delta}\, t}{\left(t^{2/a}+x^2 \right)^{\frac{a}{2}+\frac{1}{2}}} \, .
\end{align}
\end{lemma}

\begin{proof}
(1)  
Set 
 $$
    g_a(t, x) = \frac{t}{\pi} \,  \left(t^{2/a} + x^2\right)^{-(a+1)/2},\quad (t,x) \in \R_+ \times \R,\quad a > 0.
 $$
Since $1<a<2$ and $|\delta|< 2-a$,  $x \mapsto \laplacef(1,x)$ has tails at $\pm\infty$ with
polynomial decay of the same rate as $|x|^{-1-a}$: see \cite[p.~143]{zolotarev1986}. The same is true of $g_a(1,x)$. Since both functions $x \mapsto \laplacef(1,x)$ and $x \mapsto g_a(t, x)$ are continuous, \eqref{E:G-bd} holds. 

    By the scaling property \eqref{ch1'-ss7.3.6} and \eqref{E:G-bd},
\begin{align*}
    \laplacef(t,x) = t^{-\frac{1}{a}}\, \laplacef(1, t^{-\frac{1}{a}}\, x) &\le t^{-\frac{1}{a}}\, \frac{K_{a}}{1+|t^{-1/a}x|^{1+a}} = t^{-\frac{1}{a}}\, \frac{K_{a} \, t^{1+\frac{1}{a}}}{t^{1+\frac{1}{a}}+|x|^{1+a}}.
\end{align*}
Since the function $t\mapsto \frac{t}{t+z}$ is monotone increasing on $\R_+$, the above quantity is less than
\begin{align*}
    K_{a} \, t^{-\frac{1}{a}}\frac{(T\vee 1)^{1+\frac{1}{a}}}{(T\vee 1)^{1+\frac{1}{a}}+|x|^{1+a}}
\le  K_{a} \,  t^{-\frac{1}{a}}\, \frac{(T\vee 1)^{1+\frac{1}{a}}}{1+|x|^{1+a}}\:.
\end{align*}
This proves \eqref{E:Gtx-bd}.
 
 (2) 
 By the scaling property of both $\laplacef$ and $g_a(t,x)$,
\beqn
 \inf_{(t,x)\in\R_+^*\times\R}\;
   \frac{\laplacef(t,x)}{g_a(t,x)} =  \inf_{y\in\R}\; \frac{\laplacef(1,y)}{g_a(1,y)}\;.
\eeqn
Let $f(y)=\laplacef(1,y)/g_a(1,y)$.
Since $y \mapsto \laplacef(t,y)$ is unimodal (see \cite[Lemma 5.10.1]{Lukacs70}), we conclude that $f(y)>0$ for all $y\in\R$, and that
$\lim_{y\rightarrow\pm \infty} f(y)>0$. Therefore, $\inf_{y\in\R} f(y)>0$, and this completes the proof of \eqref{E:GlowBd}.
\end{proof}
\medskip

\noindent{\em Technical results}
\medskip

 The next two technical results have been used in the Proposition \ref{ch1'-ss7.3-p1}.
 \begin{lemma}
 \label{ch1'-ss7.3-l1}
Let $b>0$ and $c\in\, ]0,1]$. Then for all $x\ge 0$,
\beqn
\mathrm{Re}\left[(b\pm i\sqrt x)^{-c}\right] \ge \frac{1}{b^c}-\frac{c(1+c)}{2}\, \frac{x}{b^{2+c}}.
\eeqn
\end{lemma}
\begin{proof}
Let $\theta=\arctan\left(\sqrt x/b\right) \in[0,\pi/2[$. Let $f(x):=\mathrm{Re}\left[(b\pm i\sqrt x)^{-c}\right]$. Clearly
\beqn
f(x) = \left(b^2+x\right)^{-c/2}\cos(c\theta).
\eeqn
By the intermediate value theorem, for some $\eta\in\, ]0,x[$,
\beq
\label{ch1'-ss7.3.45}
\left(b^2+x\right)^{-c/2}= b^{-c}-\frac{1}{2}cx(b^2+\eta)^{-1-c/2} \ge b^{-c}-\frac{1}{2}cxb^{-2-c}.
\eeq
Since $\cos\theta \ge 1-\theta^2/2$ and $\arctan(y)\le y$ for $y\ge 0$, we have
\beq
\label{ch1'-ss7.3.46}
\cos(c\theta)\ge 1-\frac{c^2\theta^2}{2}\ge 1-\frac{c^2x}{2b^2}.
\eeq
The lower bounds \eqref{ch1'-ss7.3.45} and \eqref{ch1'-ss7.3.46} yield the lemma.
\end{proof}

\begin{lemma}
\label{ch1'-ss7.3-l2}
For all $0\le s\le t$ and $a\in\, ]1,2]$, we have
\begin{align*}
&\int_0^s dr \left(\frac{1}{(t-r)^{1/a}} + \frac{1}{(s-r)^{1/a}} - \frac{2}{((t+s)/2-r)^{1/a}}\right)\\
&\qquad \le a^\ast(2^{1/a}-1)(t-s)^{1/a^\ast},
\end{align*}
and
\beqn
\int_0^s dr\, \frac{(t-s)^2}{((t+s)/2-r))^{2+1/a}} \le \frac{a}{a+1} 2^{1+1/a}(t-s)^{1/a^\ast},
\eeqn
where $a^\ast$ is defined by  $1/a + 1/a^\ast=1$.
\end{lemma}

\begin{proof}
By evaluating the antiderivatives, we see that
\begin{align*}
&\frac{1}{a^\ast} \int_0^s dr \left(\frac{1}{(t-r)^{1/a}} + \frac{1}{(s-r)^{1/a}} - \frac{2}{((t+s)/2-r))^{1/a}}\right) \\
&\qquad = s^{1/a^\ast} + t^{1/a^\ast} - (t-s)^{1/a^\ast} + 2^{1/a} (t-s)^{1/a^\ast} - 2^{1/a}(t+s)^{1/a^\ast}.
\end{align*}
We will prove that
\beqn
(t-s)^{-1/a^\ast}\left[s^{1/a^\ast} + t^{1/a^\ast} - (t-s)^{1/a^\ast} + 2^{1/a} (t-s)^{1/a^\ast} - 2^{1/a}(t+s)^{1/a^\ast}\right]
\eeqn
is bounded for  $0\le s < t$, or, equivalently, that
\begin{align*}
g(r):&= \frac{r^{1/a^\ast} +1- (1-r)^{1/a^\ast} + 2^{1/a} (1-r)^{1/a^\ast} - 2^{1/a}(1+r)^{1/a^\ast}}{(1-r)^{1/a^\ast}}\\
&= 2^{1/a^*} - 1 + \frac{1 + r^{1/a^*} - 2^{1/a} (1 + r)^{1/a^*}}{(1 - r)^{1/a^*}}
\end{align*}
is bounded over $[0,1[$.

By applying L'Hospital's rule once, we find that $\lim_{r\uparrow 1} g(r) = 2^{1/a}-1$. Moreover, $g(0)=0$. Consequently,
$\sup_{r\in[0,1[} g(r)<\infty$.

Differentiating the function $g$ yields
\beqn
g^\prime(r) = \frac{(1+r)^{1/a} + (1+1/r)^{1/a} - 2^{1+1/a}}{a^\ast (1-r)^{2-1/a}(1+r)^{1/a}},
\eeqn
and observe that for $r\in[0,1]$,
\begin{align*}
(1+r)^{1/a} + (1+1/r)^{1/a}&\ge 2 \left[(1+r)(1+1/r)\right]^{1/(2a)}\\
& = 2\left(\sqrt r + \frac{1}{\sqrt r}\right)^{1/a} \ge 2^{1+1/a}.
\end{align*}
Thus, $g^\prime(r)\ge 0$ for $r\in[0,1[$, and $\sup_{r\in[0,1[} g(r) = \lim_{r\uparrow 1} g(r) = 2^{1/a}-1$. This proves the first inequality with the constant
$a^\ast(2^{1/a}-1)$.

We now prove the second inequality. By evaluating the antiderivative, we see that
\begin{align*}
&\int_0^s dr\, \frac{(t-s)^2}{((t+s)/2-r)^{2+1/a}}\\
&\qquad =  \frac{a}{a+1}\, 2^{1+1/a}\, \frac{(t+s)^{1+1/a} - (t-s)^{1+1/a}}{(t+s)^{1+1/a}}\, (t-s)^{1/a^\ast}\\
&\qquad \le  \frac{a}{a+1}\, 2^{1+1/a}\, (t-s)^{1/a^\ast}.
\end{align*}
This completes the proof.
\end{proof}

\subsection{Spatial dimension $k \geq 1$}

Let $\mathcal{L} = \tfrac{\partial}{\partial t} + (-\Delta)^{a/2}$ be the fractional heat operator on $\R_+ \times \R^k$, $a > k \geq 1$, where $(-\Delta)^{a/2}$ is the fractional Laplacian defined in \eqref{rd02_08e1}. The associated fundamental solution is given in \eqref{fs-frac}. The estimates in this subsection are used in the proofs of Propositions \ref{ch3-sec3.5-p1} and \ref{ch3-sec3.5-p2}.
\medskip

\noindent{\em  Increments in space and in time, $L^2$-norms, upper bounds}

\begin{lemma}
\label{ch3-sec3.5-lA}
 Let $a > k \geq 1$, $T > 0$ and $L > 0$. There are four constants $C_i= C_i(a,k,T,L) < \infty$, $i = 1,\dots, 4$, such that:
   \smallskip

   (a)\ For all $t \in\, ]0, T]$and $x, y\in [-L, L]^k$,
   \begin{align}
   \label{ch3-sec3.5(*3)}
  &\int_0^t dr \int_{\R^k} dz\, \left[G_a(t-r, x-z) - G_a(t-r, y-z)\right]^2\notag\\
  &\qquad\qquad \leq C_1\, \vert x-y \vert^{(a - k) \wedge 2} \left[1 + 1_{\{a = 2 + k\}} \log\left(\frac{2L}{\vert x-y \vert}\right)\right].
  \end{align}   

   (b)\ For all $0 \leq s \leq t \leq T$ and all $x \in \R^k$,
   \beq
    \label{ch3-sec3.5(*4)}
    \int_0^s dr \int_{\R^k} d z\, \left[G_a(t-r, x-z) - G_a(s-r, x-z)\right]^2 \leq C_2\, (t-s)^{1-k/a},   
\eeq
and
\beq
    \label{ch3-sec3.5(*5)}
   \int_s^t dr \int_{\R^k} d z\, \left[G_a(t-r, x-z)\right]^2 = C_3\, (t-s)^{1-k/a} .   
   \eeq

      (c)\ As a consequence, for all $0 \leq s \leq t$ and $x, y \in [-L, L]^k$,
   \begin{align*}
   &\int_0^t dr \int_{\R^k} dz\, \left[G_a(t-r, x-z) - G_a(s-r, y-z)\right]^2\\
   &\qquad\qquad \leq C_4 \left[\vert t - s \vert^{1-k/a} + \vert x-y \vert^{(a - k) \wedge 2} \left(1 +  1_{\{a = 2 + k\}} \log\left(\frac{2L}{\vert x-y \vert}\right)\right)\right].
   \end{align*}
   \end{lemma}
\begin{proof} (a) Fix $t > 0$ and $x, y \in \R^k$. Apply Plancherel's theorem to see that the left-hand side of \eqref{ch3-sec3.5(*3)} is equal to the constant $1/(2\pi)^k$ multiplied by
\begin{align}
  \label{ch3-sec3.5(*6)}
 &\int_0^t dr \int_{\R^k} d\xi\,  \vert \exp(- i \xi \cdot x - r \vert \xi \vert^a) - \exp(- i \xi \cdot y - r \vert \xi \vert^a) \vert^2\notag\\
  &\qquad=  \int_0^t dr \int_{\R^k} d\xi\, \exp(- 2 r \vert \xi \vert^a) \vert \exp(- i \xi \cdot x) - \exp(- i \xi \cdot y) \vert^2\notag\\
  &\qquad =  \int_{\R^k} d\xi\, \frac{1- \exp(- 2 t \vert \xi \vert^a)}{\vert \xi \vert^a}\, (1 - \cos(\xi\cdot (x - y))),     
  \end{align}
where the last term is obtained by performing the  $dr $-integral and computing the square modulus.

    Let $h =  \vert x - y \vert$ and notice that $h\le 2L$. We will use several times the four inequalities (i) $1 - e^{-s} \leq s$, (ii) $1 - e^{-s} \leq 1$, (iii) $1 - \cos(s) \leq s^2$ and (iv) $1 - \cos(s) \leq 2$, valid for $s \in \R$. We now distinguish two cases.
\medskip

\noindent{\em Case 1}.\  $k < a < k + 2$. With the change of variables $\eta = h \xi$, we see that \eqref{ch3-sec3.5(*6)} is equal, up to a multiplicative constant, to
\beq
  \label{ch3-sec3.5(*7)}
    h^{a - k} \int_{\R^k} d\eta\, \frac{1 - \exp(- 2 t \vert \eta \vert^a / h^a)}{\vert \eta \vert^a}\, (1 - \cos(\eta \cdot e_0)),    
    \eeq
where $e_0$ is an arbitrary unit vector in $\R^k$. We write the integral in \eqref{ch3-sec3.5(*7)} as the sum $I_1+I_2$, where
\beqn
 I_i = \int_{F_i} d\eta\, \frac{1 - \exp(- 2 t \vert \eta \vert^a / h^a)}{\vert \eta \vert^a}\, (1 - \cos(\eta \cdot e_0)),
 \eeqn
 $i=1,2$, $F_1=\{\vert \eta \vert \leq 2L \}$ and
$F_2 = \{\vert \eta \vert > 2L \}$.

   For $I_1$, we use the inequalities  (ii)  and (iii), then pass to polar coordinates, to see that
   \beqn
   I_1 \leq c_k \int_0^{ 2L} d\rho\, \rho^{k-1} \rho^{-a} \rho^2 = c_k \int_0^{2L} d\rho\, \rho^{k + 1 - a}  =  \tilde c_{L,k,a},
   \eeqn
because $k + 2 - a > 0$.

   For $I_2$, we use the inequalities (ii)  and (iv), then pass to polar coordinates, to see that
   \beqn
   I_2 \leq c_k \int_{2L}^\infty d\rho\, \rho^{k-1}\rho^{- a} = \tilde c_{L,k,a} < \infty,
   \eeqn
because $k - a < 0$. Taking into account the factor $h^{a - k}$ which appears in \eqref{ch3-sec3.5(*7)}, this establishes \eqref{ch3-sec3.5(*3)} in this Case.
\smallskip

\noindent{\em Case 2.}\  $a \geq 2 + k$. Set $z = x - y$. By multiplying and dividing by $h^2$, we see that \eqref{ch3-sec3.5(*6)} is equal, up to a multiplicative constant, to
\beq
 \label{ch3-sec3.5(*8)}
     h^2 \int_{\R^k} d\xi\, \frac{1- \exp(- 2 t \vert \xi \vert^a)}{h^2 \vert \xi \vert^a}\,  (1 - \cos(\xi\cdot z)),    
     \eeq
and we decompose this integral into the sum $I_3 + I_4 + I_5$, where
\beqn
   I_i = \int_{F_i} d\xi\, \frac{1- \exp(- 2 t \vert \xi \vert^a)}{h^2 \vert \xi \vert^a}\,  (1 - \cos(\xi\cdot z)),
   \eeqn
   $i=3, 4, 5$, $F_3 = \{\vert \xi\vert \leq 1 \}$,  $F_4 = \{1 < \vert \xi \vert \leq 2L h^{-1} \}$,  $F_5 = \{2L h^{-1} < \vert \xi \vert \}$


  For $I_3$, we use the inequalities (i) and (iii), and pass to polar coordinates to see that
  \beqn
   I_3 \leq c_{k,T} \int_0^{1} d\rho\, \rho^{k-1} h^{-2} (\rho h)^2 = \tilde c_{k,T}.
   \eeqn

  For $I_4$, we use the inequalities (ii)  and (iii), and pass to polar coordinates to see that for $a> k + 2$,
  \begin{align*}
   I_4 &\leq c_k  \int_{1}^{2L h^{-1}} d\rho\, \rho^{k-1} (h^2 \rho^a)^{-1} (\rho h)^2 = c_k  \int_{1}^{2L h^{-1}} d\rho\, \rho^{k-a+1} \\
   &= \tilde c_{k,L} h^{a - k - 2} \leq \tilde c_{a, k, L}
   \end{align*}
for $h \leq  2L$, since the exponent of $h$ is positive; for $a = k+2$, we get
\beqn
    I_4 \leq c_k \log (2L h^{-1}).
    \eeqn

  For $I_5$, we use the inequalities (ii)  and (iv), and pass to polar coordinates to see that
  \beqn
   I_5 \leq  c_k  \int_{2L h^{-1} }^{\infty} d\rho\, \rho^{k-1}  (h^2 \rho^a)^{-1} = \tilde c_{k, L} h^{a - k - 2} \leq \hat c_{k, L}
   \eeqn
for $h \leq  2L$, because $a \geq k + 2$. Taking into account the factor $h^2$ which appears in \eqref{ch3-sec3.5(*8)}, this establishes \eqref{ch3-sec3.5(*3)} in this Case and completes the proof of (a).
\smallskip

    (b)\  Fix $x = y$, $0 \leq s \leq t$ and set $h = t - s$. Apply Plancherel's theorem to see that up to the multiplicative constant $1/(2\pi)^k$, the left-hand side of \eqref{ch3-sec3.5(*4)} is equal to
        \beqn
    A_1 =  \int_0^s dr \int_{\R^k} d\xi\, \vert \exp(- i \xi \cdot x - (t-r) \vert \xi \vert^a) - \exp(-i \xi \cdot x - (s - r) \vert \xi \vert^a) \vert^2,
    \eeqn
and the left-hand side of \eqref{ch3-sec3.5(*5)} is equal to
\beqn
   A_2 = \int_s^t  \int_{\R^k} d\xi\, \exp(- 2 (t-r) \vert \xi \vert^a)) =
    \int_{\R^k} d\xi\,  \frac{1 - \exp(- 2 h \vert \xi \vert^a)}{2 \vert \xi \vert^a}.
   \eeqn

   For $A_1$, we write
   \begin{align*}
   A_1 &=  \int_0^s dr \int_{\R^k} d\xi\, \exp(- 2 (s-r) \vert \xi \vert^a) (1 - \exp(- h \vert \xi \vert^a))^2\\
       & = \int_{\R^k} d\xi\, \frac{1 - \exp(- 2 s \vert \xi \vert^a)}{2 \vert \xi \vert^a} \, (1 - \exp(- h \vert \xi \vert^a))^2,
       \end{align*}
and we decompose $A_1$ into the sum $I_6+I_7$, where
\beqn
     I_i = \int_{F_i} d\xi\, \frac{1 - \exp(- 2 s \vert \xi \vert^a)}{2 \vert \xi \vert^a}\, (1 - \exp(- h \vert \xi \vert^a))^2,
     \eeqn
     $i=6,7$, $F_6 = \{\vert \xi \vert \leq h^{-1/a} \}$, $F_7 = \{\vert \xi \vert > h^{-1/a} \}$.

For $I_6$, we use the inequalities (ii) and (i), then pass to polar coordinates to see that
\beqn
    I_6 \leq  c_{k,T} \int_0^{h^{-1/a}} d\rho\, \rho^{k-1} \rho^{-a} (h \rho^a)^2
      = c_{a, k} h^2 (h^{-1/a})^{k + a} = c h^{1 - k/a}.
      \eeqn
For $I_7$, we use the inequality (ii), then pass to polar coordinates to see that
\beqn
    I_7 \leq  \int_{h^{-1/a}}^\infty d\rho\, \rho^{k-1} \rho^{-a} = c_{a,k} (h^{-1/a})^{k - a} = \tilde c_{a,k} h^{1 - k/a}.
    \eeqn
Therefore, $A_1 \leq \hat C_{a,k} h^{1 - k/a}$. This proves \eqref{ch3-sec3.5(*4)}.

For $A_2$, we use the change of variables $\eta = h^{1/a} \xi$  to see that
\beq
\label{ch3-sec3.5(*9)}
   A_2 = h^{1 - k/a} \int_{\R^k} d\eta\,  \frac{1 - \exp(- 2 \vert \eta \vert^a)}{2 \vert \eta \vert^a} = C_{a,k}\, h^{1 - k/a},  
   \eeq
since the integral converges because $a > k$. This proves \eqref{ch3-sec3.5(*5)}
 and completes the proof of (b).

  Conclusion (c) follows directly from (a), (b) and the triangle inequality.
    \end{proof}
\smallskip

\noindent{\em Increments in space, $L^2$-norms,  lower bounds}

\begin{lemma}\label{ch3-sec3.5-lB}  
Fix $a \in\, ]k, k+2[$ and $t_0 > 0$. There is $c = c_{a, t_0, k} > 0$ such that, for all $t\ge t_0$ and for all $x, y \in \R^k$ with $\vert x - y \vert \leq 1$,
\beq
\label{ch3-sec3.5(*10)}
  \int_0^t dr \int_{\R^k} dz\, \left[G_a(t-r, x-z) - G_a(t-r, y-z)\right]^2 \geq c\, \vert x - y \vert^{a - k} .      
  \eeq
  \end{lemma}

\begin{proof} Fix $a \in\, ]k, k+2[$,  $t > 0$, and let $h = \vert x - y \vert$. We have seen in \eqref{ch3-sec3.5(*7)} that the left-hand side of \eqref{ch3-sec3.5(*10)} is, up to a multiplicative constant, equal to
\beqn
     h^{a - k} \int_{\R^k} d\eta\, \frac{1 - \exp(- 2 t \vert \eta \vert^a / h^a)}{2 \vert \eta \vert^a} \, (1 - \cos(\eta \cdot e_0)),
     \eeqn
where $e_0$ is an arbitrary unit vector in $\R^k$. For $h \in\, ]0, 1]$, this is bounded below by
\beqn
  h^{a - k} \int_{\vert \eta \vert \geq 1} d\eta\, \frac{1 - \exp(- 2 t_0 \vert \eta \vert^a)}{2 \vert \eta \vert^a}\, (1 - \cos(\eta \cdot e_0))
   = c\, h^{a - k}
   \eeqn
with $c > 0$, since the integral is positive.
    This completes the proof.
    \end{proof}

\section{Wave kernel on $\re$}
\label{app2-5}

For any $t\in\re_+$ and $x\in\re$, let
\beq
\label{lightcone-app2}
D(t,x) =\{(s,y)\in\ [0,t]\times \re: \ |x-y|\le t-s\}
\eeq
(see Figure \ref{fig1-chapter3}). The fundamental solution to the wave equation on $\re$ is
\beq
\label{fs-waver}
\Gamma(t,x;s,y)=\Gamma(t-s,x-y) := \half 1_{D(t,x)}(s,y) 1_{\{t>0\}}.
\eeq
The aim of this section is to study integrated squared increments in time and in
space of $\Gamma$.

Throughout this section, we will use the notation $|\cdot |$ for the Lebesgue measure of a set and, for $(t,x), (s,y)\in \ [0,T]\times\R$, we set
\beqn
I(t,x;s,y):= |D(t,x)|-2|D(t,x)\cap D(s,y)|+|D(s,y)|.
\eeqn

\noindent{\em Increments in space-time, $L^2$-norms, upper bounds}

\begin{lemma}\label{P4:G-1'}
Fix $T>0$. For all $(t,x)$, $\left(s,y\right)\in \ [0,T]\times\R$,
\beq
\label{1'.w1}
\int_{\R_+}  dr\int_\R  dz
\left(1_{D(t,x)}(r,z) - 1_{D(s,y)}(r,z)\right)^2
\le 2T\left(\left| x-y\right| + \left|t-s\right| \right),
\eeq
and the constant $2T$ is optimal.

Equivalently,
\beq
\label{1'.w1-bis}
\int_{\R_+}  dr\int_\R  dz
\left(\Gamma(t-r,x-z) -  \Gamma(s-r,y-z)\right)^2
\le \frac{T}{2}\left(\left| x-y\right| + \left|t-s\right| \right)
\eeq
(and the constant $\frac{T}{2}$ is optimal).
\end{lemma}
\begin{proof}

By developing the square on the left-hand side of \eqref{1'.w1}, we see that it is equal to
\begin{align*}
\label{def-inc}
&\int_{\R_+} dr\int_\R  dz
\left(1_{D(t,x)}(r,z)- 1_{D(s,y)}(r,z)\right)^2\notag\\
&\qquad = |D(t,x)|-2|D(t,x)\cap D(s,y)|+|D(s,y)| = I(t,x;s,y).
\end{align*}
Clearly, for any $(t,x)\in [0,T]\times\R$, $|D(t,x)|=t^2$.

Without loss of generality we may and will assume $x-t \le y-s$, since the roles of $(t,x), (s,y)$ in \eqref{1'.w1} are symmetric.
We divide the proof into two cases.
\medskip

\noindent{\em Case 1:} $D(t,x)\cap D(s,y)=\emptyset$.
This happens if $x+t < y-s$, that is, if $t+s\le y - x$.
In this case
\begin{align*}
I(t,x;s,y)=t^2+s^2\leq T(t+s)
 \le  T\, |x-y|,
\end{align*}
giving \eqref{1'.w1} with $T$ instead of $2T$ there.
\medskip

\noindent{\em Case 2:} $D(t,x)\cap D(s,y)\ne\emptyset$. This happens if  $x+t > y-s$, which leads to two possibilities.
\smallskip

(i) $x+t \le y+s$. One can easily compute the area of $D(t,x)\cap D(s,y)$ and obtain
\beq
\label{app2-5(*1)}
|D(t,x)\cap D(s,y)| = \frac{(x-y+t+s)^2}{4},
\eeq
therefore,
\begin{align*}
   I(t,x;s,y) &= t^2-\frac{(x-y)^2}{2} - \frac{(t+s)^2}{2} - (x-y)(t+s) + s^2\\
      &\le t^2-\frac{(t+s)^2}{2} - (x-y)(t+s) + s^2.
\end{align*}
Since $t+s\ge 2(t\wedge s)$ and $t+s\le 2T$, the last expression is bounded from above by
\begin{align*}
t^2 - 2(t\wedge s)^2 + 2T |x-y| + s^2& = |t^2-s^2| +2T\, |x-y|\\
& \le 2T(|t-s| + |x-y|).
\end{align*}

(ii) $x+t > y+s$. In this case $D(t,x)\cap D(s,y) = D(s,y)$. Therefore,
\beqn
I(t,x;s,y) =t^2-s^2 \le 2 T\, |t - s|.
\eeqn
The last inequality is sharp as $t\uparrow T$ and $s\uparrow T$, giving the optimality of the constant $2T$ in \eqref{1'.w1}.
The proof of \eqref{1'.w1} is complete.
Because of \eqref{fs-waver},
the inequality \eqref{1'.w1-bis} follows.
\end{proof}


		\begin{figure}
		\centering
    \begin{subfigure}[b]{0.3\textwidth}
        \begin{tikzpicture}
\draw [<->] (2.5,0) -- (0,0) -- (0,5.5);
\node[below] at (2.5,0) {$r$};
\node[above] at (0,5.5) {$z$};
\draw (0,0.5)--(0.75,1.25) -- (0,2);
\draw (0,2.5)--(1.25,3.75) -- (0,5);
\node[right] at (0.75,1.25) {$(t,x)$};
\node[right] at (1.25,3.75) {$( s , y )$};
\node[left] at (0,5) {$ y + s $};
\node[left] at (0,2.5) {$ y - s $};
\node[left] at (0,2) {$x+t$};
\node[left] at (0,.5) {$x-t$};
\end{tikzpicture}
				
        \caption{Case 1}
    \end{subfigure}
    \begin{subfigure}[b]{0.3\textwidth}
        \begin{tikzpicture}
\draw [<->] (2.5,0) -- (0,0) -- (0,5.5);
\node[below] at (2.5,0) {$r$};
\node[above] at (0,5.5) {$z$};
\draw (0,0.5)--(1.75,2.25) -- (0,4);
\draw (0,2.5)--(1.25,3.75) -- (0,5);
\node[right] at (1.75,2.25) {$(t,x)$};
\node[right] at (1.25,3.75) {$( s , y )$};
\node[left] at (0,5) {$ y + s $};
\node[left] at (0,2.5) {$ y - s $};
\node[left] at (0,4) {$x+t$};
\node[left] at (0,.5) {$x-t$};
\draw [fill=green] (0,2.5)--(0.75,3.25)--(0,4) --(0,2.5);
\end{tikzpicture}
        \caption{Case 2 (i)}
    \end{subfigure}
    \begin{subfigure}[b]{0.3\textwidth}
        \begin{tikzpicture}
\draw [<->] (2.5,0) -- (0,0) -- (0,5.5);
\node[below] at (2.5,0) {$r$};
\node[above] at (0,5.5) {$z$};
\draw (0,0.5)--(2.25,2.75) -- (0,5);
\draw (0,3)--(0.5,3.5) -- (0,4);
\node[right] at (2.25,2.75) {$(t,x)$};
\node[below right] at (0.3,3.5)  {$( s , y )$};
\node[left] at (0,4) {$ y + s $};
\node[left] at (0,3) {$ y - s $};
\node[left] at (0,5) {$x+t$};
\node[left] at (0,.5) {$x-t$};
\draw [fill=green] (0,3)--(0.5,3.5) -- (0,4) --(0,3);
\end{tikzpicture}
        \caption{Case 2 (ii)}
    \end{subfigure}
    \caption{The three cases in the proof of Lemma \ref{ch1'.pW1-lemma}}
    \label{fig:ones_recta_holder1}
\end{figure}
\medskip

\noindent{\em Increments in space-time, $L^2$-norms,  lower bounds}

\begin{lemma}
\label{ch1'.pW1-lemma}
Fix $T>0$, $M>0$ and  $t_0\in\, ]0,T]$. There exists a constant $C=C(M,t_0)>0$ such that, for all $(t,x), (s,y)\in [t_0,T]\times [-M,M]$,
\beq
\label{1'.W100}
\int_{\R_+}  dr\int_\R  dz \left(1_{D(t,x)}(r,z) - 1_{D(s,y)}(r,z)\right)^2
\ge C \left(|t-s| + |x-y|\right).
\eeq
Equivalently,
\beq
\label{1'.W100-bis}
\int_{\R_+}  dr\int_\R  dz\, (\Gamma(t-r,x-z) - \Gamma(s-r,y-z))^2 \ge \frac{C}{4}  \left(|t-s| + |x-y|\right),
\eeq
with $C$ as in \eqref{1'.W100}.
\end{lemma}
\begin{proof}
As in Lemma \ref{P4:G-1'}, we will assume $x-t\le y-s$ and will consider the cases in the proof of that lemma.
\medskip

\noindent{\em Case 1:}\ $D(t,x)\cap D(s,y)=\emptyset$ (see Figure \ref{fig:ones_recta_holder1}(a)).
In this case
\begin{align}
\label{c1}
I(t, x; s, y) =t^2+s^2 \ge 2t_0\, |t-s|+\frac{t_0^2}{M}\, |x-y|.
\end{align}
Indeed, if  $t_0 \le s\le t$, then
\begin{align*}
 t^2 + s^2 &= (s + (t-s))^2 + s^2 = 2 s^2 + 2 s (t-s) + (t-s)^2\\
 & \geq 2 t_0^2 + 2 t_0 (t-s)\\
       &\geq \frac{2 t_0^2}{2M}\,  \vert x - y \vert  +  2 t_0\, \vert t - s \vert,
\end{align*}
while if  $t_0\le t\le s$, then we simply reverse the roles of $s$ and $t$ in the previous argument.
\medskip

\noindent{\em Case 2 (i):}\,  $D(t,x)\cap D(s,y)\ne\emptyset$ and $x+t \le y+s$ (see Figure \ref{fig:ones_recta_holder1} (b)).
In this case, we necessarily have
\beq
\label{B.5 (*0)}
     y - s \leq x + t,
     \eeq
and also $x \leq y$. Indeed, we are assuming $x-t\le y-s$, hence if $s \geq t$, this implies $x \leq y + t - s \leq y$; on the other hand, if $s \leq t$, then the inequality that defines this case gives $x \leq y + s- t \leq y$.

 Recall from \eqref{app2-5(*1)} that the area of $T_g := D(t,x) \cap D(s,y)$ is
\beqn
    |T_g| = \left(\frac{x - y + t + s}{2}\right)^2.
    \eeqn

Set
\beqn
I(t,x;s,y) = |D(t,x)|-2|D(t,x)\cap D(s,y)| + |D(s,y)|.
\eeqn
We distinguish two subcases:
\medskip

\noindent{$(i_1)$}\
$s \ge t$. In this case,
\begin{align*}
    I(t,x ; s, y)  &= |D(t,x)| - |T_g| + |D(s,x)| - |T_g |
   \geq  |D(s,x)| - |T_g| \\
   & = s^2 - \left(\frac{x - y + t + s}{2}\right)^2 \\
   &= \left(s + \frac{x - y + t + s}{2}\right)\left(s - \frac{x - y + t + s}{2}\right)\\
   &= \frac{1}{4} (3s + t + x - y) (s - t + y - x).
   \end{align*}
By \eqref{B.5 (*0)}, the first parenthesis is no less than $3 s + t - s - t = 2 s \geq 2 t_0$, and the second parenthesis is $\vert s - t \vert + \vert y - x \vert$. This gives the lower bound
\beq
\label{B.5(*2)}
  I(t,x ; s, y)  \geq  \frac{t_0}{2} \left(\vert s - t \vert + \vert y - x \vert\right).
\eeq
\smallskip

\noindent{$(i_2)$}\ $s < t$. In this case,
\beqn
      I(t,x ; s, y) \geq |D(t,x)| - |T_g |,
      \eeqn
and by the same type of calculation as above, we get the same lower bound as in \eqref{B.5(*2)}.
\medskip

\noindent{\em Case 2 (ii):}\ $D(t,x)\cap D(s,y)\ne\emptyset$ and $x+t > y+s$ (see Figure \ref{fig:ones_recta_holder1} (c)).
In this case, $(s,y) \in D(t,x)$, so $\vert x - y \vert \leq t-s$ and $D(t,x) \cap D(s,y) = D(s,y)$. Therefore,
\begin{align}
\label{B.5(*3)}
     I(t,x; s,y) &= |D(t,x)| - |D(s,x)| = t^2 - s^2 \notag\\
     &= (t+s)(t-s) \geq 2 t_0 (t-s)
        \geq t_0 (\vert s - t \vert + \vert x - y \vert).
        \end{align}

The inequalities \eqref{c1}, \eqref{B.5(*2)} and \eqref{B.5(*3)} establish the lower bound \eqref{1'.W100}, which is equivalent to
\eqref{1'.W100-bis} by \eqref{fs-waver}.
\end{proof}


\section{Wave kernel on $\re_+$}
\label{app2-6}

For any $t\in\re_+$ and $x\in\re_+$, let
\beqn
E(t,x) =\{(s,y)\in[0,t]\times \re_+: \ |x-t+s|\le y\le x+t-s\}
\eeqn
(see Figure \ref{fig3.2}).
The Green's function of the wave operator on $\re_+ \times \R_+$ with Dirichlet boundary conditions is
\beq
\label{fs-waver-positive} \Gamma(t,x;s,y):= G(t-s;x,y) = \half 1_{E(t,x)}(s,y).
\eeq
In this section, we study integrated squared increments in time and in
space of $\Gamma$.
\medskip

\noindent{\em Increments in space-time, $L^2$-norms,  upper bounds}

\begin{lemma}
\label{rdprop_wave_R_ub}
Fix $T > 0$. For all $(t,x), (s,y)\in [0,T]\times \IR_+$,
\beq
\label{rdt3}
 \int_{\IR_+} dr \int_{\IR_+} dz \left[1_{E(t,x)}(r,z) - 1_{E(s,y)}(r,z) \right]^2 \leq 2T\, (\vert x - y\vert + \vert t - s\vert)
\eeq
(and the constant $2T$ is optimal).

Equivalently,
\beq\label{rdt3-bis}
   \int_{\IR_+} dr \int_{\IR_+} dz \left[G(t-r;x,z) - G(s-r;y,z) \right]^2 \leq \frac{T}{2}\, (\vert x - y\vert + \vert t - s\vert)
\eeq
(and the constant $\frac{T}{2}$ is optimal).
\end{lemma}
\begin{proof}
Throughout the proof, we will often use the following elementary facts:
\begin{description}
\item{(i)}\ If $(t,x)$ is such that $x-t\ge 0$ then $E(t,x) = D(t,x)$ (defined in \eqref{lightcone-app2}).
\item{(ii)}\ If $x-t<0$ then $|E(t,x)|=|D(t,x)|-(t-x)^2$, and this obviously implies $|E(t,x)|\le |D(t,x)|$.
\end{description}

We will also use the formula
\begin{align}
\label{B.6(*0)}
&\int_{\R_+} dr\int_\R  dz
\left(1_{E(t,x)}(r,z)- 1_{E(s,y)}(r,z)\right)^2\notag\\
&\qquad\qquad = |E(t,x)|-2|E(t,x)\cap E(s,y)|+|E(s,y)|
\end{align}
and, without loss of generality, we may and will assume that
$x-t\le y-s$.
\medskip

\noindent{\em Case 1:} $x-t\ge 0$.
Then we also have $y-s\ge 0$,  and from fact  (i), we are in the setting of Lemma \ref{P4:G-1'}. This yields \eqref{rdt3} in this case, as well as the optimality of the constant $2T$.

In the remainder of this proof, we will assume $x-t < 0$ (and $x-t\le y-s$).
\medskip

\noindent{\em Case 2:}  $E(t,x)\cap E(s,y)=D(t,x)\cap D(s,y)$.
Using fact (ii), we see that
\begin{align*}
&|E(t,x)|-2|E(t,x)\cap E(s,y)|+|E(s,y)|\\
&\qquad\qquad \le |D(t,x)|-2|D(t,x)\cap D(s,y)|+|D(s,y)|,
\end{align*}
and again Lemma \ref{P4:G-1'} implies \eqref{rdt3} in this case.

\begin{figure}[H]
\centering
\subcaptionbox{}{
\begin{tikzpicture}[scale=1]
\draw [<->] (4,0) -- (0,0) -- (0,5.5);
\node[below] at (4,0) {$r$};
\node[above] at (0,5.5) {$z$};
\draw (0,0.5)--(.5,0)--(1.75,1.25) -- (0,3);
\draw (0,1) -- (2,3)--(0,5);
\node[right] at (2,3) {$( s , y )$};
\node[right] at (1.75,1.25) {$(t,x)$};
\node[left] at (0,3) {$x+t$};
\node[left] at (0,0.35) {$t-x$};
\node[left] at (0,5) {$ y + s $};
\node[left] at (0,1) {$ y - s $};
\draw [fill=green] (0,1)--(1,2)--(0,3) --(0,1);
\end{tikzpicture}
}%
\hfill 
\subcaptionbox{}{
\begin{tikzpicture}[scale=1]
\draw [<->] (3.25,0) -- (0,0) -- (0,5.5);
\node[below] at (3.25,0) {$r$};
\node[above] at (0,5.5) {$z$};
\draw [fill=green](0,1)--(1,2) -- (0,3);
\draw (0,0.5)--(0.5,0) -- (2.75,2.25)--(0,5);
\node[right] at (2.75,2.25) {$(t,x)$};
\node[right] at (0.75,2.2) {$( s , y )$};
\node[left] at (0,3) {$ y + s $};
\node[left] at (0,1) {$ y - s $};
\node[left] at (0,5) {$x+t$};
\node[left] at (0,0.35) {$t-x$};
\end{tikzpicture}
}%
\caption{$E(t,x)\cap E(s,y)=D(t,x)\cap D(s,y)$ if $|x-t|\le y-s$}\label{fig:ones_semi1}
\end{figure}

\noindent{\em Case 3:}  $E(t,x)\cap E(s,y)\neq D(t,x)\cap D(s,y)$.
The condition $t-x>y-s$ is necessarily satisfied. Indeed, otherwise,  the possible configurations are those illustrated in
Figures \ref{fig:ones_semi1} (a) and \ref{fig:ones_semi1} (b), which fall into Case 2.


We now split the study of this case into three subcases.
\medskip

\noindent{\em Case 3.1:}  $x+t\leq  y + s $.
Looking at Figure \ref{fig:ones_semi2}, we see that in both situations $y - s \geq 0$ and $y - s < 0$, we have
\beq
\label{B.6(*1)}
     \vert E(t,x) \cap E(s,y) \vert = \vert D(t,x) \cap D(s,y) \vert - \frac{1}{4}(t - x - (y-s))^2,
     \eeq
therefore by facts (i) and (ii), the expression in \eqref{B.6(*0)} is bounded above by
\begin{align*}
      &\vert D(t,x) \vert  - (t - x)^2 - 2 \left(\vert D(t,x) \cap D(s,y) \vert  - \frac{1}{4}(t - x - (y-s))^2 \right)\\
      &\qquad\qquad  + \vert D(s,y) \vert \\
     &\qquad \leq \vert D(t,x) \vert - 2 \vert D(t,x) \cap D(s,y) \vert + \vert D(s,y) \vert,
     \end{align*}
     because
\begin{align*}
   &\half (t-x - (y-s))^2 - (t-x)^2 \\
     &\qquad\qquad = \half\, \vert D((t - x - y + s)/2, (t-x + y - s)/2) \vert - \vert D(t-x, 0) \vert \\
    &\qquad\qquad \leq 0,
\end{align*}
since the second triangle contains the first (in fact, $(t-x + y - s)/2 - 0 = t - x - (t-x-y+s)/2$).
We now apply Lemma \ref{P4:G-1'} to obtain \eqref{rdt3}.
\begin{figure}[H]
\centering
\subcaptionbox{}{\begin{tikzpicture}[scale=1]
\draw [<->] (4,0) -- (0,0) -- (0,5.5);
\node[below] at (4,0) {$r$};
\node[above] at (0,5.5) {$z$};
\draw (0,2)--(2,0)--(3,1)--(0,4);
\draw (0,1)--(2,3) --(0,5);
\node[right] at (2,3) {$( s , y )$};
\node[right] at (3,1) {$(t,x)$};
\node[left] at (0,4) {$x+t$};
\node[left] at (0,2) {$t-x$};
\node[left] at (0,5) {$ y + s $};
\node[left] at (0,1) {$ y - s $};
\draw [fill=green] (0,4)--(1.5,2.5)--(0.5,1.5) --(0,2)--(0,4);
\end{tikzpicture}}%
\hfill
\subcaptionbox{}{
\begin{tikzpicture}[scale=0.7]
\draw [<->] (4,0) -- (0,0) -- (0,5.5);
\draw (0,0)--(0,-3.5);
\node[below] at (4,0) {$r$};
\node[above] at (0,5.5) {$z$};
\draw (0,3)--(3,0) --(3.5,0.5)-- (0,4);
\draw (0,0.5)--(0.5,0) -- (2.75,2.25)--(0,5);
\node[right] at (2.75,2.25){$( s , y )$};
\node[right] at (3.5,0.5) {$(t,x)$};
\node[left] at (0,4) {$x+t$};
\node[left] at (0,3) {$t-x$};
\node[left] at (0,-3) {$x-t$};
\node[left] at (0,5) {$ y + s $};
\node[left] at (0,0.5) {$ s - y $};
\node[left] at (0,-0.5) {$ y - s $};
\draw [dashed] (3,0)--(0,-3);
\draw [dashed] (0.5,0)--(0,-0.5);
\draw[fill=green] (0,4)--(2.25,1.75)--(1.75,1.25)--(0,3)--(0,4);
\end{tikzpicture}
}%
\caption{Case $E(t,x)\cap E(s,y)\neq D(t,x)\cap D(s,y)$, with $x+t\le y+s$}\label{fig:ones_semi2}
\end{figure}

\begin{figure}[H]
\centering
\subcaptionbox{}{
\begin{tikzpicture}[scale=1]
\draw [<->] (4,0) -- (0,0) -- (0,5.5);
\node[below] at (4,0) {$r$};
\node[above] at (0,5.5) {$z$};
\draw (0,1)--(1,2) -- (0,3);
\draw (0,2)--(2,0) -- (3.5,1.5)--(0,5);
\node[right] at (3.5,1.5) {$(t,x)$};
\node[right] at (0.8,2) {$( s , y )$};
\node[left] at (0,3) {$ y + s $};
\node[left] at (0,1) {$ y - s $};
\node[left] at (0,5) {$x+t$};
\node[left] at (0,2) {$t-x$};
\draw [fill=green] (0,3)--(1,2)--(0.5,1.5) --(0,2)--(0,3);
\end{tikzpicture}
}%
\hfill 
\subcaptionbox{}{
\begin{tikzpicture}[scale=0.7]
\draw [<->] (5,0) -- (0,0) -- (0,5.5);
\draw (0,0)--(0,-3.5);
\node[below] at (5,0) {$r$};
\node[above] at (0,5.5) {$z$};
\draw (0,3)--(3,0) --(4,1)-- (0,5);
\draw (0,1)--(1,0) -- (2.5,1.5)--(0,4);
\node[right] at (4,1) {$(t,x)$};
\node[below right] at (2,1.6) {$( s , y )$};
\node[left] at (0,4) {$ y + s $};
\node[left] at (0,1) {$ s - y $};
\node[left] at (0,-1) {$ y - s $};
\node[left] at (0,5) {$x+t$};
\node[left] at (0,3) {$t-x$};
\node[left] at (0,-3) {$x-t$};
\draw [dashed] (1,0)--(0,-1);
\draw [dashed] (3,0)--(0,-3);
\draw[fill=green] (0,4)--(2.5,1.5)--(2,1)--(0,3)--(0,4);
\end{tikzpicture}
}%
\caption{Case $E(t,x)\cap E(s,y)\neq D(t,x)\cap D(s,y)$, with $x+t > y+s$}\label{fig:ones_semi4}
\end{figure}
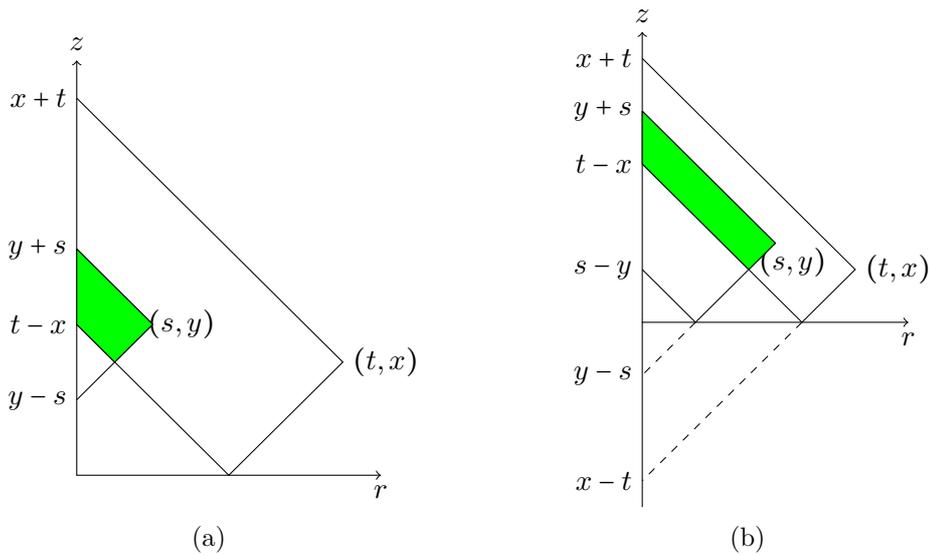

\medskip

\noindent{\em Case 3.2:} $x + t > y + s$ and $t - x \leq  y + s$. Looking at Figure \ref{fig:ones_semi4}, we see that  in both situations $y - s \geq 0$ and $y - s < 0$, \eqref{B.6(*1)} holds and
\beqn
   \vert E(t,x) \vert = \vert D(t,x) \vert - (t - x)^2,
   \eeqn
while
\beqn
    \vert E(s,y) \vert = \vert D(s,y) \vert - 1_{\{y - s < 0\}} (s - y)^2 .
    \eeqn
Therefore, the expression in \eqref{B.6(*0)} is equal to
\begin{align*}
   &\vert D(t,x) \vert - (t - x)^2 - 2 \left(\vert D(t,x) \cap D(s,y) \vert  - \frac{1}{4}(t - x - (y-s))^2 \right)\\
   &\qquad \qquad\qquad + \vert D(s,y) \vert - 1_{\{y - s < 0\}} (s - y)^2\\
   &\qquad\qquad\leq \vert D(t,x) \vert - 2 \vert D(t,x) \cap D(s,y) \vert + \vert D(s,y) \vert,
   \end{align*}
as one easily checks by the same arguments as in Case 3.1. We again apply Lemma \ref{P4:G-1'} to obtain \eqref{rdt3}.
\medskip


\noindent{\em Case 3.3:} $x + t > y + s$ and $t - x > y + s$. Then $E(t,x) \cap E(s,y) = \emptyset$, and therefore, the expression in \eqref{B.6(*0)} is equal to
     $\vert E(t,x) \vert + \vert E(s,y) \vert$.

Recalling that $\vert D(t,x) \vert = t^2$, for all $(t,x) \in \R_+ \times \R_+$, this is equal to
     $ t^2 - (t - x)^2 + s^2 - 1_{\{y - s < 0\}} (s - y)^2$.
For $y - s < 0$, this is equal to
\beqn
    2 t x + 2 s y - x^2 - y^2 \leq 2t\, (x+y) \leq 2t\, \vert t - s \vert
    \eeqn
(dropping the two negative terms, and using that  $s < t - x - y \leq t$ yields $x + y < t - s=|t-s|$);  for $y - s \geq 0$, this is equal to
\beqn
    2 t x - x^2 + s^2 \leq 2tx +y s \leq 2 t x + 2 t y \leq 2t\, (x + y) \leq 2t\, \vert t - s \vert
    \eeqn
(dropping the $- x^2$, using $s \leq y$ and then $s<t-x-y\le t\le 2t$ because  $x + y < t - s$). This ends the proof of the lemma.
\end{proof}
\medskip

\noindent{\em Increments in space-time, $L^2$-norms,  lower bounds}

\begin{lemma}
\label{rdprop_wave_R+_lb}
Fix $0<t_0<T$ and $M > t_0$. There is a constant $c_{t_0,M,T}>0$ such that, for all $(t,x),\, (s,y) \in [t_0,T]\times [t_0,M]$,
\beq
\label{rdt5-lemma}
c_{t_0,M,T}\, (\vert x - y\vert + \vert t - s\vert) \leq \int_{\IR_+} dr \int_{\IR_+} dz \left[ 1_{E(t,x)}(r,z) -1_{E(s,y)}(r,z) \right]^2.
\eeq
Equivalently,
\beq
\label{rdt5-lemma-bis}
\frac{c_{t_0,M,T}}{4}\, (\vert x - y\vert + \vert t - s\vert) \leq \int_{\IR_+} dr \int_{\IR_+} dz \left[ G(t-r;x,z) -G(s-r;y,z) \right]^2.
\eeq

\end{lemma}
\begin{proof}
Consider $s, t \in [t_0, T]$, $x, y \in [t_0, M]$. Without loss of generality, we will assume that
 \beq
 \label{B.6(**1)}
 x - t \leq y - s.
 \eeq

Define
\beqn
J(t,x;s,y) = |E(t,x)|-2|E(t,x)\cap E( s , y )|+|E( s , y )|,
\eeqn
which is the right-hand side of \eqref{rdt5-lemma} (see \eqref{B.6(*0)}).
We now consider several different cases.
\medskip

\noindent{\em Case A}. Suppose that $\vert t - s \vert + \vert x - y \vert \geq t_0/2$. Then $E(t,x) \neq E(s,y)$, and even $1_{E(t,x)} \ne 1_{E(s,y)}$ on a set of positive Lebesgue  measure, so $J(t,x; s,y) > 0$. Since $J$ is a continuous function of $(t,x; s,y)$, over this compact domain, it is bounded below by a constant $c > 0$. Therefore.
\beqn
      J(t,x ; s,y)  \geq c \geq \frac{c}{2T}\, \vert t - s \vert + \frac{c}{2M}\, \vert x - y \vert,
      \eeqn
proving \eqref{rdt5-lemma} in this case.
\medskip

\noindent{\em Case B}. Suppose that $\vert t - s \vert + \vert x - y \vert < t_0/2$. We split this case into two subcases.
\medskip

\noindent{\em Case 1:} $x - t \geq 0$. In this case, by \eqref{B.6(**1)}, we also have $y - s \geq 0$, therefore $E(t,x) = D(t,x)$ and $E(t,y) = D(t,y)$, and \eqref{rdt5-lemma} follows from Lemma \ref{ch1'.pW1-lemma}.
\medskip

 \noindent{\em Case 2:} $x - t < 0$. This case is further divided into subcases.
 \medskip

  \noindent{\em Case 2.1:} $x - t < 0$ and  $y - s < 0.$
  \medskip

   \noindent{\em Case 2.1a:} $x - t < 0$, $y - s < 0$ and (($s \leq t$ and $y \leq x$) or ($s > t$ and $y \geq x$). These two configurations are illustrated in Figure \ref{fig:ones_semi3} (a) and (b), respectively.

   In this case,  $\vert t - s \vert + \vert x - y \vert$ is the vertical distance between the lines $z = -r + t + x$ and $z = -r + s + y$ (in the plane $(r,z))$. The line that is above the other depends on whether $s > t$ or not.
   When $s \leq t$ and $y \leq x$, we have
   \beqn
     \frac{t_0}{2} \geq  \vert t - s \vert + \vert x - y \vert = t - s + x - y.
     \eeqn
Since both terms are nonnegative, we obtain $t - s \leq t_0/2$, that is, $s - t \geq - t_0/2$. Since $x > 0$ and $y \geq t_0$,
     $ t - x \leq s + y - t_0 \leq t + x - t_0$,
because the first inequality is equivalent to $ x + y + s - t \geq t_0$ and is implied by $x + y \geq 3 t_0/2$, which is true since $x \geq t_0$ and $y \geq t_0$. Therefore, the parallelogram $P_1$ with vertices $(t_0/2, s+y-t_0/2), (t_0/2, t + x-t_0/2), (t_0, t + x - t_0), (t_0, s + y - t_0)$ is entirely contained in   $E(t,x) \setminus E(s,y)$. It follows that
\beqn
   J(t,x ; s,y) \geq |P_1| = \frac{t_0}{2}\, (t + x - (s + y)) =  \frac{t_0}{2}\, ( \vert t - s \vert + \vert x - y \vert),
   \eeqn
proving \eqref{rdt5-lemma} in this situation.

    When $s > t$ and $y \geq x$, then $s + y > t + x$ and $s + y - t_0 > t + x - t_0 \geq t > t - x \geq s - y$, where the second inequality holds since $x \geq t_0$ and the last because of \eqref{B.6(**1)}. Therefore, $P_1$ is entirely contained in $E(s,y) \setminus E(t,x)$.
      It follows that
      \beqn
     J(t,x ; s,y) \geq  |P_1| = \frac{t_0}{2}\, (s + y - (t + x)) = \frac{t_0}{2} \, ( \vert t - s \vert + \vert x - y \vert),
     \eeqn
proving \eqref{rdt5-lemma} in this case.
\medskip

\begin{figure}[H]
\centering
\subcaptionbox{}{
\begin{tikzpicture}[scale=1]
\draw [<->] (5,0) -- (0,0) -- (0,5.5);
\draw (0,0)--(0,-3.5);
\node[below] at (5,0) {$r$};
\node[above] at (0,5.5) {$z$};
\draw (0,1.8)--(1.8,0) --(3.4,1.6)-- (0,5);
\draw (0,1)--(1,0) -- (2.5,1.5)--(0,4);
\draw[fill=green] (0,4)--(2.5,1.5)--(1.4,0.4)--(0,1.8)--(0,4);
\node[right] at (3.4,1.6) {$(t,x)$};
\node[left] at (2.5,1.5) {$( s , y )$};
\node[left] at (0,4) {$ y + s $};
\node[left] at (0,1) {$ s - y $};
\node[left] at (0,-1) {$ y - s $};
\node[left] at (0,5) {$x+t$};
\node[left] at (0,1.8) {$t-x$};
\node[left] at (0,-1.8) {$x-t$};
\draw [dashed] (1,0)--(0,-1);
\draw [dashed] (1.8,0)--(0,-1.8);
\end{tikzpicture}
}%
\hfill
\subcaptionbox{}{
\begin{tikzpicture}[scale=1]
\draw [<->] (4,0) -- (0,0) -- (0,5.5);
\draw (0,0)--(0,-3.5);
\node[below] at (4,0) {$r$};
\node[above] at (0,5.5) {$z$};
\draw (0,2)--(2,0) --(2.5,0.5)-- (0,3);
\draw (0,1)--(1,0) -- (3,2)--(0,5);
\node[right] at (3,2){$(s,y)$};
\node[right] at (2.5,0.5) {$(t,x)$};
\node[left] at (0,3) {$x+t$};
\node[left] at (0,2) {$t-x$};
\node[left] at (0,-2) {$x-t $};
\node[left] at (0,5) {$y+s$};
\node[left] at (0,1) {$s-y$};
\node[left] at (0,-1) {$y-s$};
\draw [dashed] (2,0)--(0,-2);
\draw [dashed] (1,0)--(0,-1);
\draw[fill=green] (0,3)--(2,1)--(1.5,.5)--(0,2)--(0,3);
\end{tikzpicture}
}%

\hspace{35pt}
    \caption{The two situations of Case 2.1a in Lemma \ref{rdprop_wave_R+_lb}}\label{fig:ones_semi3}
\end{figure}
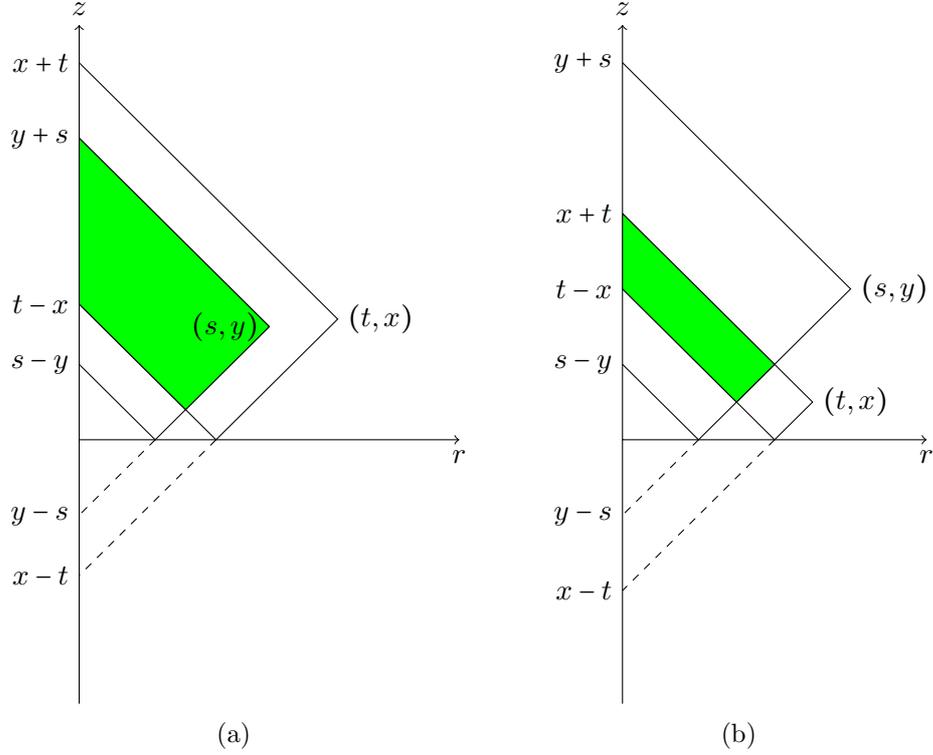

 \noindent{\em Case 2.1b:} $x - t < 0$, $y - s < 0$, $s \leq t$ and $y > x$. This is the situation shown in Figure \ref{fig:ones_semi4} (b). In this case, $\vert t - s \vert + \vert x - y \vert$ is the horizontal distance between the lines $z = r + x - t$ and $z = r+y-s$, and  the second line is above the first.
 Notice that $t - x - t_0/2 < t_0/2 + s - y$, because this is equivalent to the condition that defines Case B. Since the horizontal line with equation $z = t_0/2$ intersects the line $z = r + y - s$ at the point $(r = t_0/2 + s - y, z = t_0/2)$, the parallelogram
$P_2$ with vertices $(t_0/2 + s-y, t_0/2), (t_0/2 + t - x, t_0/2), (t_0 +t - x, t_0), (t_0 + s - y, t_0)$ is entirely contained in  $E(t,x) \setminus E(s,y)$. Therefore,
\beqn
         J(t,x ; s,y) \geq |P_2| = \frac{t_0}{2}\, (t - x - (s - y)) = \frac{t_0}{2}\, ( \vert t - s \vert + \vert x - y \vert),
         \eeqn
proving \eqref{rdt5-lemma} in this case.
\medskip

 \noindent{\em Case 2.2:} $ x - t < 0$, $y - s \geq 0$.
 \medskip

 \noindent{\em Case 2.2a:} $x - t < 0$, $y - s \geq 0$ and (($s \leq t$ and $y \leq x$) or ($s > t $ and $y \geq x$).
 The situation with $s \leq t$ is  shown in Figure \ref{fig:ones_semi1} (b), and the situation with $s > t$ is shown in Figure \ref{fig:ones_semi1} (a). Notice that $t - x < y + s$, since otherwise $t -s \geq y + x \geq 2 t_0$, and this is not possible in Case B. We now use the same argument as in Case 2.1a, which only depends on the situation to the upper left of the two points (no reflection) to see that \eqref{rdt5-lemma} holds. The parallelogram $P_1$ is entirely contained in $E(t,x) \setminus E(s,y)$ if $s \leq t$, or in $E(s,y) \setminus E(t,x)$ if $s > t$.
 \medskip

\noindent{\em Case 2.2b:} $x - t < 0$, $y - s \geq 0$ and $s \leq t$ and $y > x$. In this case, since condition B is satisfied,
\beq
\label{B.6**}
         \frac{t_0}{2} \geq  \vert t - s \vert + \vert x - y \vert = t - s + y - x = (t - x) + (y - s).
         \eeq
Since both terms are positive,
$0 \leq y - s \leq t_0/2$, $0 \leq t - x \leq t_0/2$, and $t_0 - (y - s) \geq t_0/2 - (y - s) \geq 0$ and $t_0/2 + s - y \geq t - x - t_0/2$ by \eqref{B.6**}.
 Therefore, the parallelogram $P_2$ with vertices $(t_0/2 - (y - s), t_0/2), (t_0/2 + t - x, t_0/2), (t_0 + t - x, t_0), (t_0 - (y - s), t_0)$ is entirely contained in $E(t,x) \setminus E(s,y)$.
Hence,
\beqn
      J(t,x ; s,y) \geq  |P_2| =\frac{ t_0}{2}\, (t - x + (y - s)) = \frac{t_0}{2}\, ( \vert t - s \vert + \vert x - y \vert),
      \eeqn
proving \eqref{rdt5-lemma} in this final case.
\end{proof}


\section{Wave kernel on a bounded interval}
\label{app2-7}

Fix $L>0$. For $x,y \in [0,L]$ and $s,t\in \IR_+$ with $s \leq t$, let
\begin{align}
\label{Green-wave-interval}
\Gamma(t,x;s,y)& := G_L(t-s;x,y) \notag\\
&= \sum_{n=1}^\infty \frac{2}{n\pi} \sin\left(\frac{n\pi x}{L}\right)
\sin\left(\frac{n\pi y}{L}\right) \sin\left(\frac{n\pi (t-s)}{L}\right)
\end{align}
be the Green's function of the wave operator (with Dirichlet boundary conditions) on $[0,L]$, introduced in \eqref{wave-bc2-1'}. For $s>t$, we set $G_L(t-s;x,y) =0$. We observe that
\beq
\label{sym}
G_L(t;x,y) = G_L(t;L-x,L-y),
\eeq
because  $\sin(n \pi - x) = \pm \sin x$, for any $n \in \Z$.

In this section, we study integrated squared increments in time and in space of $G_L$.
\medskip

\noindent{\em Increments in space and in time, $L^2$-norms,  upper bounds}

\begin{lemma}
\label{app2-5-l-w1}
\begin{enumerate}
\item
   For any $t\in\IR_+$ and $x,y\in[0,L]$,
\beq
\label{app2.30}
\int_0^t dr \int_0^L dz \left(G_L(t-r;x,z) - G_L(t-r;y,z)\right)^2
 \le t\ |x-y|.
\eeq
\item
   For any $x\in[0,L]$ and $s,t\in \IR_+$ with $s\le t$,
\beq
\label{rd_app2.31}
   \int_s^t dr\int_0^L dz \, G_L^2(t-r;x,z) \leq x \left(1-\tfrac{x}{L}\right) (t-s) \le L(t-s).
\eeq
For all $t\ge s\ge 0$ and $x\in[0,L]$,
\beq
\label{app2.31}
\int_0^s dr \int_0^L dz \left(G_L(t-r;x,z) - G_L(s-r;x,z)\right)^2 \le 2s(t-s).
\eeq
\end{enumerate}
As a consequence of 1.~and 2., for any $T>0$, there is a constant $C=C(T,L)$ such that for all $s,t \in [0,T]$ and $x,y \in [0,L]$,
\begin{align}\nonumber
   &\int_0^T dr \int_0^L dz \left(G_L(t-r;x,z) - G_L(s-r;y,z)\right)^2 \\
   &\qquad\qquad\qquad  \leq C\left(|t-s|^\half + |x-y|^\half\right)^2.
\label{rdeB.5.5a}
\end{align}
\end{lemma}

\begin{proof}
Since $G_L(t;x,y) = G_1\left(\frac{t}{L}; \frac{x}{L}, \frac{y}{L}\right)$, the following identities hold:
\begin{align}
&\int_0^t dr \int_0^L dz \left(G_L(t-r;x,z) - G_L(t-r;y,z)\right)^2 \notag\\
&\qquad = L^2 \int_0^{\tfrac{t}{L}} dr \int_0^1 dv \left(G_1(\tfrac{t}{L}-r;\tfrac{x}{L},v) - G_1(\tfrac{t}{L}-r;\tfrac{y}{L},v)\right)^2, \label{scaling-W-1}\\
&\int_s^t dr\int_0^L dz \, G_L^2(t-r;x,z) = L^2 \int_{\tfrac{s}{L}}^{\tfrac{t}{L}} dr \int_0^1 dz \, G_1^2\left(\tfrac{t}{L}-r;\tfrac{x}{L},z \right), \label{scaling-W-2}\\
 &\int_0^s dr \int_0^L dz \left(G_L(t-r;x,z) - G_L(s-r;x,z)\right)^2 \notag\\
 &\qquad = L^2 \int_0^{\tfrac{s}{L}} dr \int_0^1 dv \left(G_1(\tfrac{t}{L}-r;\tfrac{x}{L},v) - G_1(\tfrac{s}{L}-r;\tfrac{x}{L},v)\right)^2. \label{scaling-W-3}
\end{align}
Hence, we will first prove the lemma for $L=1$ and then these formulas will imply  \eqref{app2.30}--\eqref{app2.31} for any $L>0$.

1. Since $\left(\sqrt{2} \sin (n\pi z),\, n\ge1\right)$ is a CONS of $L^2([0,1])$, using the definition \eqref{Green-wave-interval} with $L=1$, for $x,y \in [0,1]$, we have
\begin{align*}
&\int_0^t dr \int_0^1dz \left(G_1(t-r;x,z) - G_1(t-r;y,z)\right)^2\\
&\quad =  \int_0^t dr \int_0^1dz \left[\sum_{n=1}^\infty\, \frac{2}{n\pi}\, (\sin(n\pi x)-\sin(n\pi y))
\sin(n\pi z) \sin(n\pi (t-r))\right]^2\\
&\quad = \int_0^t dr\, \sum_{n=1}^\infty\, \frac{2}{\pi^2 n^2}\,(\sin(n\pi x)-\sin(n\pi y))^2\sin^2(n\pi (t-r))\\
&\quad \le 2t\, \sum_{n=1}^\infty\, \frac{(\sin(n\pi x)-\sin(n\pi y))^2}{\pi^2\, n^2}\\
&\quad = t \left( |x-y| - |x-y|^2\right) \\
&\quad \le t\, |x-y|,
\end{align*}
where the last equality follows from Lemma \ref{lemB.5.2}. This proves \eqref{app2.30} for $L=1$, hence for all $L>0$.
\smallskip

2. For $x\in [0,1]$,
\begin{align*}
&\int_s^t dr \int_0^1dz\, G_1^2(t-r;x,z) \\
&\qquad = \int_s^t dr \int_0^1dz \left[\sum_{n=1}^\infty\, \frac{2}{\pi n}\, \sin(n\pi x) \sin(n\pi z) \sin(n\pi (t-r))\right]^2\\
&\qquad =  2 \int_s^t dr\, \sum_{n=1}^\infty\, \frac{1}{\pi^2\, n^2}\, \sin^2(n\pi x)\sin^2(n\pi (t-r))\\
&\qquad \le 2 (t-s)\, \sum_{n=1}^\infty\, \frac{\sin^2(n\pi x)}{\pi^2\, n^2} \\
&\qquad = x(1-x) (t-s),
\end{align*}
where the last equality follows from Lemma \ref{lemB.5.2}.  This proves \eqref{rd_app2.31} for $L=1$, hence for all $L>0$.
\smallskip

We now prove \eqref{app2.31} for $L=1$. Applying the trigonometric identity $\sin a - \sin b = 2\cos\left(\frac{a+b}{2}\right)\sin\left(\frac{a-b}{2}\right)$, we see that
\begin{align*}
&\int_0^s dr \int_0^1dz \left(G_1(t-r;x,z) - G_1(s-r;x,z)\right)^2 \\
&\qquad = \int_0^s dr \int_0^1 dz\, \Big[\sum_{n=1}^\infty \frac{2}{\pi n}\, \sin(n\pi x) \sin(n\pi z)\\
&\qquad\qquad\qquad\qquad \qquad\times\left[\sin(n\pi (t-r)) - \sin(n\pi (s-r))\right]\Big]^2\\
&\qquad = 8 \int_0^s dr\, \sum_{n=1}^\infty\, \frac{1}{\pi^2\, n^2}\, \sin^2(n\pi x)\cos^2\left(\frac{n\pi(t+s-2r)}{2}\right)
\sin^2\left(\frac{n\pi (t-s)}{2}\right)\\
&\qquad \le 8s \sum_{n=1}^\infty\, \frac{\sin^2(n\pi (t-s)/2)}{\pi^2\, n^2}\\
&\qquad \le 2s (t-s),
\end{align*}
as follows from the inequality \eqref{rem-lemB.5.2} with $s$ there replaced by $(t-s)/2$.
This proves \eqref{app2.31} for $L=1$, hence for all $L>0$.

\smallskip

It remains to prove \eqref{rdeB.5.5a}. Assume that $s \leq t$ and observe that
\beqn
\left[\int_0^T dr \int_0^L  dz \left(G_L(t-r;x,z)-G_L(s-r;y,z)\right)^2\right]^\half \le T_1+T_2,
\eeqn
with
\begin{align*}
T_1^2&:= \int_0^t dr \int_0^L  dz \left(G_L(t-r;x,z)-G_L(s-r;x,z)\right)^2,\\
T_2^2&:= \int_0^s dr \int_0^L  dz \left(G_L(s-r;x,z)-G_L(s-r;y,z)\right)^2.
\end{align*}
Clearly,
\begin{align*}
T_1^2 &=  \int_0^s dr \int_{\IR}  dz \left(G_L(t-r;x,z)-G_L(s-r;x,z)\right)^2\\
&\qquad + \int_s^t dr \int_{\IR}  dz\, G_L^2(t-r;x,z).
\end{align*}
By \eqref{app2.31}, \eqref{rd_app2.31} and \eqref{app2.30},
\beqn
T_1^2  \le (2T+L)\, (t-s) \qquad\text{and}\qquad T_2^2\le T\, |x-y|.
\eeqn
Therefore,
\begin{align*}
&\left[\int_0^T dr \int_0^L  dz \left(G_L(t-r;x,z)-G_L(s-r;y,z)\right)^2\right]^\half\\
&\qquad\qquad \le c(T,L)(\left( |t-s|^{\half} +  |x-y|^{\half}\right),
\end{align*}
which is equivalent to \eqref{rdeB.5.5a}.
This completes the proof of the lemma.
\end{proof}

The next aim is to establish a lower bound for the $L^2$-norm of increments of $G_L(t;x,y)$. To enable the use of the computations done in Section \ref{app2-6}, it is more convenient to consider the other expression of the Green's function:
\begin{align}
\label{green-bis}
G_L(t; x,y) &= \frac{1}{2} \sum_{m=-\infty}^\infty
 \left[ 1_{\{|x-2mL-y|\le t\}} - 1_{\{|x-2mL+y|\le t\}} \right]\notag\\
 &= \frac{1}{2} \left[1_{F_1(t,x)}(y) - 1_{F_2(t,x)}(y)\right],
 \end{align}
$t>0$, $x,y\in\, ]0,L[$ (see \eqref{wave-bc1-1'} and \eqref{wave-bc1-100'} in Chapter \ref{chapter1'}, as well as Figure \ref{fig3.3}, which is useful as illustration).

\medskip

\noindent{\em Increments in space-time, $L^2$-norms,  lower bounds}

\begin{lemma}
\label{app2-5-l-w2}
Fix $T, L>0$ and $0<t_0<\min\left(T,\tfrac{L}{2}\right)$. Then there is a constant $c_{t_0,L,T} >0$ such that for all $(t,x),(s,y) \in [t_0,T] \times [t_0,L - t_0]$,
\beq\label{rdt7-lemma}
   c_{t_0,L,T} (\vert t-s\vert + \vert x-y\vert) \leq \int_{\IR_+} dr \int_0^L dz\, (G_L(t-r;x,z) - G_L(s-r;y,z))^2.
\eeq
\end{lemma}

\begin{proof}
For $(t,x), (s,y) \in [t_0, T] \times [t_0, L - t_0]$,  denote $J(t,x ; s,y)$ the right-hand side of \eqref{rdt7-lemma}.
As in the proof of Lemma \ref{rdprop_wave_R+_lb}, we only need to establish \eqref{rdt7-lemma} in the case where
\beq
\label{B.7(*1)}
     \vert t - s \vert + \vert x - y \vert \leq \frac{t_0}{2}.
     \eeq
Consider the set
\beqn
     S = ([t_0/2, t_0] \times [0, L]) \cup ([0, T] \times [t_0/2, t_0]) \cup ([0, T] \times [L - t_0, L- t_0/2]).
     \eeqn
     By distinguishing mainly the same cases as in Lemma \ref{rdprop_wave_R+_lb}, and some additional but similar cases, one can check
 that for all points $(t,x),\, (s,y) \in [t_0, T] \times [t_0, L - t_0]$ satisfying \eqref{B.7(*1)}, there is a parallelogram $P_0$ with vertices in $S$ and area $\frac{t_0}{2}(\vert t - s \vert + \vert x - y \vert)$ (entirely contained in the region $[0,T] \times [0,L]$), and where
\beqn
     [ (1_{F_1(t,x)}(r,z) - 1_{F_2(t,x)}(r,z)) - (1_{F_1(s,y)}(r,z) - 1_{F_2(s,y)}(r,z))]^2 = 1.
     \eeqn
When $s$ and $t$ are close to $t_0$, this parallelogram can be taken in the strip $[t_0/2, t_0] \times [0, L]$; otherwise, it may be found at the upper or lower corners of the right-most rectangle in Figure \ref{fig3.3}. This implies that $J(t,x ; s,y) \geq \frac{t_0}{2} (\vert t - s \vert + \vert x - y \vert)$ in this case. This completes the proof of the lemma.
\end{proof}

\section{Geometric space-time convolution series}
\label{rd08_26s1}


We use the term ``geometric" in relation with geometric series of real numbers $\sum_{\ell =1}^\infty x^\ell$. In this section, we deal with series in which the general term is a space-time convolution power, as has been introduced in \eqref{def-kas}:
\beq
\label{rd08_26e2}
   \cK(t, x) =  \sum_{\ell = 1}^\infty \, J^{\star \ell}(t, x).
\eeq
We call these objects {\em geometric space-time-convolution series}.\index{convolution!series}\index{series!convolution} Here, we are mainly interested in the case where $J(t, x) = \rho^2\, \Gamma^2(t, x)$, and $\Gamma$ is the fundamental solution of a PDE. 


In several cases, it is possible to give explicit estimates and sometimes even explicit formulas for $\cK(t, x)$, either by computing the space-time convolutions and the series in \eqref{rd08_26e2} explicitly, or by using the Fourier-Laplace transform operator $\hcF$ (see \eqref{rd08_26e3} for the definition) and Lemma \ref{rd08_11l1} to sum an ordinary geometric series:
\beq
\label{rd08_17e1}
   \hcF[\cK](z, \xi) = \sum_{\ell=1}^\infty\, (\hcF [J](z, \xi))^\ell = \frac{\hcF [J](z, \xi)}{1 - \hcF [J](t, x)}.
\eeq
This requires determining $\hcF [J]$ and the inverse Fourier-Laplace transform of the right-hand side of \eqref{rd08_17e1}.

In the next proposition, we compute explicitly the value of $\cK(t, x)$ when $J$ is (up to a positive multiplicative constant) the square of the fundamental solution of the heat and wave equations, and provide upper and lower bounds when $J$ is the square of the fundamental solution of the fractional heat equation, as needed in Sections \ref{rdrough} and \ref{rd1+1anderson}.
\medskip

\noindent{\em Some special functions}
\medskip

In order to state the proposition, we introduce some notation. We will use $\Phi(x)$ to denote the standard Normal distribution function:\label{rd05_07e1}
$$
   \Phi(x) = \tfrac{1}{\sqrt{2 \pi}}\int_{-\infty}^x  e^{-\frac{y^2}{2}}\, dy,\qquad x \in \R.
$$
For the wave equation, we will need the function
\begin{align}
\label{rd08_26e5}
I_0(x) = \sum_{k=0}^{\infty} \frac{x^{2k}}{2^{2k}\, (k!)^2} \, ,
\end{align}
which is the modified Bessel function of the first kind of order $0$ (see \eqref{rd08_26L10}). It is well-known that
\begin{align}\label{rd08_17e5}
I_0(x) \sim \frac{e^x}{\sqrt{2\pi x}}\quad
\text{ as } x \to \infty
\end{align}
(see \cite[10.30.4]{olver}).

For the fractional heat equation, we will need the  {\em two-parameter Mittag-Leffler
function}\index{Mittag-Leffler function}\index{function!Mittag-Leffler} \cite[Section 1.2]{podlubny1999}:
\begin{align}\label{E:Mittag-Leffler}
E_{\alpha,\beta}(z) := \sum_{k=0}^{\infty} \frac{z^k}{\Gamma_E(\alpha k+\beta)},
\qquad z \in \bC,\, \alpha > 0, \ \beta \in \bC,
\end{align}
where $\Gamma_E(z)$ is the Euler Gamma function \eqref{Euler-gamma}, extended as a meromorphic function to $\bC \setminus (-\N^*)$. It is well-known that for $\alpha \in\, ]0, 2[$ and $\beta > 0$,  as $x \to + \infty$ with $x \in \R$, 
\beq
\label{rd08_17e4}
   E_{\alpha,\beta}(x) = \frac{1}{\alpha}\, x^{(1-\beta)/\alpha}\, \exp\left(x^{1/\alpha}\right)
- \frac{x^{-1}}{\Gamma(\beta-\alpha)} + O\left(x^{-2}\right) 
\eeq
(see \cite[Theorem 1.3 p.~32]{podlubny1999}. Note that when $\alpha = \beta$, $\Gamma(\alpha - \beta) = \Gamma(0)$ and $1/\Gamma(0)$ is interpreted as $0$).  We will also use the function
\beq
\label{rd08_17e6}
    g_a(t, x) = \frac{t}{\pi} \,  \left(t^{2/a} + x^2\right)^{-(a+1)/2},\quad (t,x) \in \R_+ \times \R,\ a > 0.
\eeq

Finally, let $a^*$ be the dual of $a$: $1/a+1/a^*=1$.
\medskip

\noindent{\em Three explicit formulas/bounds}
\medskip

\begin{prop}
\label{rd08_12p1}
Let $\rho \in \R_+^*$ be fixed. 
\medskip

\noindent {\em (a) Case where $\Gamma$ is the heat kernel on $\R$.}
For  $\nu > 0$, as in \eqref{fsnu}, let 
\beqn
   \Gamma_\nu(t,x) = \frac{1}{\sqrt{2\pi\nu t}}\exp\left(-\frac{x^2}{2\nu t}\right), \qquad (t, x) \in \R_+^* \times \R.
\eeqn
This is the fundamental solution associated to the heat operator $\cL_\nu = \frac{\partial}{\partial t} - \frac{\nu}{2}\, \frac{\partial^2}{\partial x^2}$. When $\nu = 2$, this is \eqref{heatcauchy-1'}.

Let $\Jh_\nu(t, x; \rho) =\rho^2\, \Gamma_\nu^2(t, x)$. For $\ell \in \N^*$,
\beq
\label{rd08_20e12}
     (\Jh_\nu)^{\star \ell}(t, x; \rho) = \frac{\rho^{2\ell}}{(4\, \nu)^{\ell/2}}\, \frac{t^{(\ell-2)/2}}{\Gamma_E(\ell/2)}\, \Gamma_{\nu/2}(t, x). 
\eeq
Define $\cK^{\mbox{\scriptsize \rm h}}_\nu(t, x; \rho)$ as in \eqref{rd08_26e2} with $J(t, x)$ there replaced by $\Jh_\nu(t, x; \rho)$. Then
\beq
\label{rd08_12e8}
    \cK^{\mbox{\scriptsize \rm h}}_\nu(t, x; \rho) = B_\nu(t; \rho)\, \Gamma_{\nu/2}(t, x),
\eeq
where
\beqn
   B_\nu(t; \rho) = \frac{\rho^2}{\sqrt{4 \pi\, \nu\, t}} + \frac{\rho^4}{2\nu}\, e^{\rho^4\, t/(4\nu)} \Phi\left(\sqrt{\frac{\rho^4\, t}{2 \nu}} \right).
\eeqn
\smallskip

\noindent{\em  (b)  Case where $\Gamma$ is the fundamental solution of the wave equation on $\R$.}
For  $\nu > 0$ and $(t, x) \in \R_+^* \times \R$, let
\beq
\label{rd08_21e6b}
   \Gamma_\nu(t, x) = \frac{1}{2 \nu}\, 1_{\R_+}(t)\, 1_{[0,\, \nu\, t]}(\vert x \vert).
\eeq
For $\nu > 0$, $\Gamma_\nu$ is the fundamental solution associated to the wave operator $\cL_\nu = \frac{\partial^2}{\partial t^2} - \nu^2 \frac{\partial^2}{\partial x^2}$. When $\nu=1$, this is \eqref{wfs}. 

Define
 $\Jw_\nu(t, x; \rho) = \rho^2\, \Gamma_\nu^2(t, x) = \frac{\rho^2}{2\, \nu}\,  \Gamma_\nu(t, x)$. For $\ell \in \N^*$,
\beq
\label{rd08_24e1}
     (\Jw_\nu)^{\star \ell}(t, x; \rho) = \frac{\rho^{2 (\ell-1)}\, \left((\nu\, t)^2 - x^2 \right)^{\ell-1}}{(2\, \nu)^{3 (\ell-1)}\, ((\ell-1) !)^2} \, \Jw(t, x; \rho). 
\eeq
Define $\cK^{\mbox{\scriptsize \rm w}}_\nu(t, x; \rho)$ as in \eqref{rd08_26e2} with $J(t, x)$ there replaced by $\Jw_\nu(t, x; \rho)$. Then
\begin{align}
\label{rd08_12e10}
      \cK^{\mbox{\scriptsize\rm w}}_\nu(t, x; \rho) &= \frac{\rho^2}{4\, \nu^2}\, I_0\left(\sqrt{\frac{\rho^2}{2\, \nu}\, \frac{\left[(\nu\, t)^2 - x^2\right]}{\nu^2}} \right)\, 1_{[0,\, \nu\, t]}(\vert x \vert) \, 1_{\R_+}(t),
\end{align}
where $I_0$ is the modified Bessel function of the first kind of order $0$ given in \eqref{rd08_26e5}.
\smallskip

\noindent {\em  (c) Case where $\Gamma$ is the fractional heat kernel on $\R$.}
Let $a \in\, ]1, 2[$, $\vert \delta\vert \leq 2 - a$ and let
\beq
\label{rd08_21e7}
     \Gamma_\nu(t, x) = \null_\delta G_a(\nu\, t,x),\quad (t, x) \in \R_+^* \times \R.
\eeq
For $\nu > 0$, $\Gamma_\nu$ is the fundamental solution associated to the fractional heat operator $\cL_\nu = \tfrac{\partial}{\partial t} -\nu\,  \null_x D_\delta^a$ (see Section \ref{ch1'-ss7.3}). When $\nu = 1$, this is \eqref{ch1'-ss7.3.3}.

Define $\Jf_\nu(t, x; \rho) =  \rho^2 \, \Gamma_\nu^2(t,x)$ and let 
\beqn
    \tilde {\Jf_\nu}(t, x; \rho) = C_0\,\rho^2\, (\nu\, t)^{-1/a}\, \Gamma_\nu(t,x), \qquad\text{where }C _0= \sup_{x\in \re}\, \null_\delta G_a(1,x),
\eeqn 
so that $\Jf_\nu(t, x) \leq \tilde {\Jf_\nu}(t, x)$. 

For $\ell \in \N^*$,
\beq
\label{rd08_24e2}
     (\Jf_\nu)^{\star \ell}(t, x; \rho) \leq \frac{ \left(C_0\, \rho^{2}\, \nu^{-1/a}\, \Gamma_E\left(1/a^*\right)\right)^\ell}  {\Gamma_E(\ell/a^*)}\,t^{\frac{\ell}{a^*}-1}\, \null_\delta G_a(\nu\, t,x), 
\eeq
with $\tfrac{1}{a}+\tfrac{1}{a^*}=1$.  Define $\cK^{\mbox{\scriptsize \rm f}}_\nu(t, x; \rho)$ as in \eqref{rd08_26e2} with $J(t, x)$ there replaced by $\Jf_\nu(t, x; \rho)$. There are three constants
$c_i = c_i(a, \delta) > 0$, $i=1, 2, 3$, such that
\beq
\label{rd08_17e3}
  \cK^{\mbox{\scriptsize\rm f}}_\nu(t, x; \rho)
   \leq c_1\, \rho^2\, (\nu\, t)^{-1/a}\, E_{1/a^*,1/a^*}\left(c_1 \rho^2\, \nu^{-1/a}\, t^{1/a^*}\right)\, \null_\delta G_a(\nu\, t,x),
\eeq
where, for $\alpha,\beta>0$, $E_{\alpha,\beta}$ 
is given in \eqref{E:Mittag-Leffler}. 
If $a \in\, ]1, 2[$ and $\vert \delta \vert < 2 - a$ (note the strict inequality), then
\beq
\label{rd08_17e2}
  \cK^{\mbox{\scriptsize\rm f}}_\nu(t, x; \rho) \geq c_2\, \rho^2 \, E_{1/a^*,1/a^*}\left(c_3 \, \rho^2\, \nu^{-1/a}\, t^{1/a^*} \right)\, g_a^2\left(\nu\, t,x\right),
\eeq
where $g_a$ is defined in \eqref{rd08_17e6}.
\end{prop}

\begin{proof} In the proofs of (a), (b) and (c), in order to simplify the notation, we will replace $\Jh_\nu$, $\Jw_\nu$ and $\Jf_\nu$ respectively, by $J$, and $\cK^{\mbox{\scriptsize \rm h}}_\nu$, $\cK^{\mbox{\scriptsize \rm w}}_\nu$ and $\cK^{\mbox{\scriptsize\rm f}}_\nu$ respectively, by $\cK$.

(a) Computing $\Gamma_\nu^2(t,x)$, we obtain
\beqn
   J(t, x; \rho) = \gamma(t) \, \Gamma_{\nu/2}(t, x), \qquad \text{with } \gamma(t) =  \frac{\rho^2}{\sqrt{4 \pi\, \nu\, t}}\, .
   \eeqn

For fixed $t > 0$, the Fourier transform of $x \mapsto J(t, x; \rho)$ is
\beqn
  \cF[J(t,\ast; \rho)](\xi) = \cF\left[\gamma(t)\, \Gamma_{\nu/2}(t, *)\right](\xi) = \gamma(t)\, \exp(- \nu\, t\,  \xi^2/4),
\eeqn
and by \eqref{rd08_26L3} and the fact that $\Gamma_E(1/2) = \sqrt{\pi}$ \cite[5.4.6]{olver}, the Laplace transform of this expression is
\beqn
  \hcF [J(\cdot, *, \rho)](z, \xi) =\underline{\cL}\left[\gamma(\cdot) \exp\left( -(\nu\xi^2 \cdot)/4\right)\right]
(z)  
 =  \frac{\rho^2}{\sqrt{4\nu\, z + \nu^2\, \xi^2}},\quad \text{Re}(z) > 0.
\eeqn
By Lemma \ref{rd08_11l1},
\beqn
   \hcF[J^{\star \ell}(\cdot, *, \rho)](z, \xi) = \left(\frac{\rho^{2}}{\sqrt{4\nu\, z + \nu^2\, \xi^2}}\right)^{\ell}
\eeqn
and by \eqref{rd08_26e2},
\beqn
    \hcF[\cK(\cdot, *, \rho)](z, \xi) = \frac{\rho^2}{\sqrt{4\nu\, z + \nu^2 \, \xi^2} - \rho^2}.
\eeqn
In order to obtain the inverse Laplace transforms of these expressions, we use \eqref{rd08_26L9} to see that
\beqn
    \cF[J^{\star \ell}(t, *; \rho)](\xi) = \frac{\rho^{2\ell}}{(4\, \nu)^{\ell/2} }\, \frac{t^{-1 +\ell/2}}{\Gamma_E(\ell/2)}\, \exp\left(- \tfrac{\nu\, t \, \xi^2}{4} \right).
\eeqn
Since $ \cF^{-1}\left(\exp\left(-\frac{\nu t}{4}\, \ast^2 \right)\right)(x) = \frac{1}{\sqrt{\pi\nu t}}\exp\left(-\frac{x^2}{\nu t}\right)$, we obtain \eqref{rd08_20e12}.

 Further, applying 
 \eqref{rd01_22L1} and using the elementary property $\cLL (e^{- i \alpha\, \cdot} \, f(\cdot) )(z) = \cLL f(z +\alpha)$, we see that 
\beqn
   \cF[\cK(t, *; \rho)](\xi) = \exp\left(- \frac{\nu\, t \, \vert \xi \vert^2}{4} \right) \left(\frac{\rho^2}{\sqrt{4 \pi\, \nu\, t}} + \frac{\rho^4}{2 \nu}\, e^{\frac{\rho^4\, t}{4 \nu}} \Phi\left(\sqrt{\frac{\rho^4\, t}{2 \nu}} \right)\right),
\eeqn
and, as before, taking the inverse Fourier transform of this quantity, we obtain \eqref{rd08_12e8}.

    (b)  Observe that $J(t, x; \rho) = \frac{\rho^2}{2\, \nu}\, \Gamma(t, x)$. Using the fact that for $a < b$, $\F [1_{[a, b]}(*)](\xi) = i(e^{-i \xi b} - e^{-i \xi a})/\xi$, we see that for $t > 0$,
\beqn
   \cF [J(t, *; \rho)](\xi) = \bar\rho^2\, \frac{\sin(\nu\,t \, \xi )}{\nu\, \xi},
\eeqn
with $\bar \rho=\rho/\sqrt{2\, \nu}$. By \eqref{rd08_26L7},
\beqn
   \hcF [J(\cdot, *; \rho)](z, \xi) = \cLL\left[\bar\rho^2\, \frac{\sin(\nu\, t \, \vert \xi \vert)}{\nu\, \vert \xi \vert}\right](z) = \frac{\bar\rho^2}{z^2 + \nu^2\, \xi^2},
\eeqn
\beq
\label{rd08_24e9}
   \hcF [J^{\star \ell}(\cdot, *; \rho)](z, \xi) = \frac{\bar\rho^{2\ell}}{\left(z^2 + \nu^2\, \xi^2\right)^\ell}
\eeq
and
\beqn
   \hcF [\cK(\cdot, *; \rho)](z, \xi) = \frac{\bar\rho^2}{z^2 + \nu^2\, \xi^2 - \bar\rho^2}.
\eeqn

By 
\eqref{rd08_26L7a} and \eqref{rd08_26L7},
\beqn
   \cF [\cK(\cdot, *; \rho)](\xi) = \left\{\begin{array}{ll}
      \bar \rho^2  \frac{\sinh\left(t \, \sqrt{\bar \rho^2 - \nu^2\, \xi^2}\right)}{\sqrt{\bar \rho^2 - \nu^2\, \xi^2}} &\qquad\text{if } \nu^2\, \xi^2 < \bar\rho^2, \\[10pt]
 \bar \rho^2   \frac{\sin\left(t \, \sqrt{\nu^2\, \xi^2 - \bar \rho^2}\right)}{\sqrt{\nu^2\, \xi^2 - \bar \rho^2}}& \qquad\text{if }  \nu^2\, \xi^2 > \bar \rho^2.
    \end{array}\right.
\eeqn
From \eqref{rd08_26L11}, 
noting that $\cF f$ is twice the Fourier-cosine transform of $f$ when $f$ is an even function, we obtain
\beqn
    \cK (t, x; \rho) = \frac{\rho^2}{4\, \nu^2}\, I_0\left(\sqrt{\frac{\rho^2}{2\, \nu}\, \frac{\left( (\nu\, t)^2 - x^2\right)}{\nu^2}}\right)\, 1_{[0,\, \nu\, t]}(\vert x \vert).
\eeqn
This establishes \eqref{rd08_12e10}.

We do not have a ready expression for the inverse Laplace transform of $\hcF[J^{\star \ell}]$ in \eqref{rd08_24e9}. 
Nevertheless, we are able to prove the equality \eqref{rd08_24e1} using the formula just obtained for $\cK(t, x)$. Indeed, recalling the Taylor expansion \eqref{rd08_26e5} of $I_0$, we see that
\begin{align}
    \cK(t, x; \rho) &= \frac{\rho^2}{4\, \nu^2}\, 1_{[0,\, \nu\, t]}(\vert x \vert)\, \sum_{\ell = 1}^\infty\, \frac{1}{2^{2 (\ell- 1)}\, ((\ell - 1)!)^2}\, \left(\frac{\rho^2}{2\, \nu}\, \frac{\left( (\nu\, t)^2 - x^2\right)}{\nu^2} \right)^{\ell -1} ,
\label{rd04_29e1}
\end{align}
and by definition,
\begin{equation}\label{rd04_29e2}
    \cK (t, x; \rho) =  \sum_{\ell = 1}^\infty J^{\star \ell}(t, x) = \sum_{\ell = 1}^\infty \left(\frac{\rho^2}{2\, \nu}\right)^\ell \Gamma_\nu^{\star \ell}(t, x).
\end{equation}
By identifying the terms that contain $\rho^{2 \ell}$ of the two power series (in the variable $\rho$) \eqref{rd04_29e1} and \eqref{rd04_29e2}, we see that
\begin{equation*}
      \left(\frac{1}{2\, \nu}\right)^\ell \Gamma_\nu^{\star \ell}(t, x) = \frac{1}{4\, \nu^2}\, 1_{[0,\, \nu\, t]}(\vert x \vert)\, \frac{1}{2^{2 (\ell- 1)}\, ((\ell - 1)!)^2}\, \left(\frac{1}{2\, \nu}\, \frac{ (\nu\, t)^2 - x^2}{\nu^2} \right)^{\ell -1}.
\end{equation*}
Equivalently,
\begin{align*}
   J^{\star \ell}(t, x; \rho) &= \frac{\rho^2}{2\, \nu}\, \Gamma_\nu(t, x)\, \frac{1}{2^{2 (\ell- 1)}\, ((\ell - 1)!)^2}\, \left(\frac{\rho^2}{2\, \nu}\, \frac{ (\nu\, t)^2 - x^2}{\nu^2} \right)^{\ell -1} \\
   &= J(t, x; \rho)\, \frac{\rho^{2(\ell-1)}\left((\nu\, t)^2 - x^2 \right)^{\ell - 1}}{(2\, \nu)^{3(\ell -1)} ((\ell -1)!)^2},
\end{align*}
which is \eqref{rd08_24e1}.

(c)  In this case, we carry out the proof in two steps, beginning with the case where $\nu = 1$. 
\medskip

\noindent{\em Step 1. The case $\nu = 1$. } Here, we do not have an explicit formula for $\cK$ because we do not have a good expression for $(\null_\delta G_a(t,x))^2$. However, by the scaling property \eqref{ch1'-ss7.3.6},
\beqn
   0 \leq \null_\delta G_a(t,x) = t^{-\frac{1}{a}}\, \null_\delta G_a(1, t^{-\frac{1}{a}}\, x) \leq t^{-\frac{1}{a}}\, \sup_{x \in \R} \null_\delta G_a(1,  x) \leq \gamma(t),
\eeqn
where
\beq
\label{rd08_24e8}
   \gamma(t) =  C_0\, t^{-\frac{1}{a}} \quad\text{and}\quad C_0 = \sup_{x \in \R} \null_\delta G_a(1,  x).
\eeq
Note that $C_0 < \infty$ by \eqref{ch1'-ss7.3.3}. In particular,
\beq
\label{rd08_17e7}
    (\null_\delta G_a(t,x))^2 \leq \gamma(t)\, \null_\delta G_a(t,x).
\eeq
Define
\beqn
    \tilde J(t, x; \rho) = \rho^2\gamma(t)\, \null_\delta G_a(t,x) \quad\text{and} \quad \tilde \cK(t, x; \rho) = \sum_{\ell=1}^\infty \tilde J^{\star \ell}(t, x; \rho).
\eeqn
Then $J(t, x; \rho) = \rho^2 \, (\laplacef (t, x))^2 \leq \tilde J(t, x; \rho)$ and $\cK(t, x; \rho) \leq \tilde \cK(t, x; \rho)$. By \eqref{fs-fractional},
\beqn
   \cF[\tilde J(t, *; \rho)](\xi) = \rho^2\,\gamma(t)\, \exp\left(-t\, |\xi|^a \exp\left(-i\pi\delta\  \text{sgn}(\xi)/2\right) \right) 1_{\R_+}(t).
\eeqn
By \eqref{rd08_26L3}, 
\beqn
   \cLL[\cF[\tilde J(\cdot, *; \rho)](\xi)](z) = C_0\, \rho^2\, \Gamma_E\left(\tfrac{1}{a^*}\right) \left(z + |\xi|^a \exp\left(-i\pi\delta\  \text{sgn}(\xi)/2\right) \right)^{-\frac{1}{a^*}},
\eeqn
where $a^\ast$ is defined by $\tfrac{1}{a}+\frac{1}{a^*} = 1$. 
It follows that
\beqn
   \hcF[\tilde J^{\star \ell}](z, \xi) = \left(C_0\, \rho^2\, \Gamma_E\left(\tfrac{1}{a^*}\right) \right)^\ell \left(z + |\xi|^a \exp\left(-i\pi\delta\  \text{sgn}(\xi)/2\right)\right)^{-\frac{\ell}{a^*}}
\eeqn
and
\beq
\label{rd08_21e9}
    \hcF[\tilde \cK](z, \xi) = \frac{C_0\, \rho^2\, \Gamma_E\left(1/a^*\right)}{(z + |\xi|^a \exp\left(-i\pi\delta\  \text{sgn}(\xi)/2\right))^{\frac{1}{a^*}} - C_0\, \rho^2\, \Gamma_E\left(1/a^*\right)}.
\eeq
Use \eqref{rd08_26L9} 
to see that 
\begin{align}
\label{rd08_21e8}
   \cF[\tilde J^{\star \ell}(t, *; \rho)](\xi) &=
   \frac{\left(C_0\, \rho^{2} \, \Gamma_E\left(1/a^* \right) \right)^\ell}{\Gamma_E(\ell/a^*)}\, t^{\frac{\ell}{a^*}-1}\notag\\
   &\qquad\times \exp\left(- t\, |\xi|^a \exp\left(-i\pi\delta\  \text{sgn}(\xi)/2\right)\right),
\end{align}
and therefore, by \eqref{fs-fractional}, 
\beq
\label{rd08_21e10}
   J^{\star \ell}(t, x; \rho) \leq \tilde J^{\star \ell}(t, x; \rho) =  \frac{ \left(C_0\, \rho^{2}\,\Gamma_E\left(1/a^*\right) \right)^\ell}{\Gamma_E(\ell/a^*)}\,t^{\frac{\ell}{a^*}-1} \null_\delta G_a(t,x).
\eeq
This is \eqref{rd08_24e2} when $\nu = 1$. 

It follows that
\begin{align*}
   \tilde \cK(t, x; \rho) &= \null_\delta G_a(t,x) \, \sum_{\ell = 1}^\infty\, \frac{\left(C_0\, \rho^{2}\, \Gamma_E\left(1/a^*\right)\right)^\ell \, t^{\frac{\ell}{a^*}-1}}{\Gamma_E(\ell/a^*)} \\
       &= C_0\, \rho^{2}\, \Gamma_E\left(1/a^*\right) \, \null_\delta G_a(t,x) \, t^{-1/a}\, \sum_{k = 0}^\infty\, \, \frac{\left(C_0\, \rho^{2}\, \Gamma_E\left(1/a^*\right)\right)^k \, t^{k/a^*}}{\Gamma_E\left(\frac{k}{a^*} + \frac{1}{a^*}\right)} \\
    &= c_1\, \rho^{2}  \, \null_\delta G_a(t,x) \, t^{-1/a}\,  E_{1/a^*,1/a^*}\left(c_1\, \rho^2\, t^{1/a^*}\right)
\end{align*}
by \eqref{E:Mittag-Leffler}. This proves \eqref{rd08_17e3} when $\nu = 1$.

For the lower bound \eqref{rd08_17e2}, notice that for all $(t, x) \in \R_+^* \times \R$,
\beq
\label{rd01_23e1}
   \null_\delta G_a(t,x) \geq C_{a, \delta}\, \pi\, g_a(t, x),
\eeq
where
\beqn
    C_{a, \delta} := \inf_{(t, x) \in \R_+^* \times \R} \frac{\null_\delta G_a(t,x)}{\pi\, g_a(t, x)} > 0
\eeqn
by Lemma \ref{rd01_21l1} (2). 
From \eqref{rd08_26e2} and \eqref{rd01_23e1}, we see that 
\begin{align}\nonumber
\cK(t,x; \rho) &= \sum_{\ell=1}^\infty
\left(\rho^2\laplacef^2\right)^{\star \ell}(t,x) \ge \sum_{\ell=1}^\infty \left(\rho^2\, C_{a, \delta}^2\, \pi^2\,
g_a^2\right)^{\star \ell} (t,x) \\
   &=  \sum_{k=0}^\infty \left(C\, \rho^2\, g_a^2\right)^{\star (k+1)}  (t,x).
\label{rdE:gaIndu}
\end{align}

We now bound space-time convolutions of $g_a^2$ with itself. We claim that for all $\ell \in \N^*$,
\begin{align}\label{E:gaIndu}
 \left(\rho^2 g_a^2\right)^{\star \ell}(t,x)\ge \frac{\rho^{2\ell}
 \, \Theta_a^{\ell - 1}\, (\Gamma_E(1/a^*))^{\ell}}{\Gamma_E(\ell/a^*)}\, t^{(\ell - 1)/a^*}\, g_a^2(t,x), 
\end{align}
where	
\beqn
   \Theta_a := C_{a+1/2}^2\, \pi^{-3}\, 2^{-2(a+3+1/a)}
\eeqn
and $C_{a+1/2}$ is defined in \eqref{rd01_23e2} (with $\nu = a+1/2$ there). 

   Indeed, the case where $\ell=1$ is clear. Consider $\ell \geq 2$ and assume by induction that \eqref{E:gaIndu} holds for $\ell-1$.
By the induction hypothesis and Lemma \ref{rd01_24l1} (3) below,
\begin{align}\nonumber
& \left(\rho^2 g_a^2\right)^{\star \ell}(t,x) = [\left(\rho^2 g_a^2\right)^{\star (\ell - 1)} \star \left(\rho^2 g_a^2\right)](t, x) \\ \nonumber
&\qquad \ge \frac{\rho^{2 \ell} \Theta_a^{\ell-2} (\Gamma_E(1/a^*))^{\ell - 1}}{\Gamma_E((\ell - 1)/a^*)} \\ \nonumber
  &\qquad \qquad \times  \int_0^t ds \, (t-s)^{(\ell - 2)/a^*}
   [g_a^2(t-s, *) \ast g_a^2(s, *)](x)\\
&\qquad\ge K_\ell \,  t^{-1/a} \int_0^t d s \: g_a^2(t-s,x) \left[(t-s)s\right]^{-\frac{1}{a}} (t-s)^{\frac{\ell-2}{a^*}+\frac{2}{a}},
\label{rd11}
\end{align}
where
$$
   K_\ell = \frac{\rho^{2 \ell} \Theta_a^{\ell-2} \, \Gamma_E(1/a^*)^{\ell - 1}  \, C_{a+1/2}^2}
{\Gamma_E((\ell - 1)/a^*)\, \pi^3\, 2^{2a+3}} = \frac{\rho^{2 \ell} \Theta_a^{\ell-1} \, \Gamma_E(1/a^*)^{\ell - 1}\, 2^{3 + 2/a}}
{\Gamma_E((\ell - 1)/a^*)}.
$$
Notice that $t-s \geq t/2$ for $0\le s\le t/2$, so we do the change of variable $r = t-s$ and we apply Lemma \ref{rd01_24l1} (4) below 
to see that
\begin{align}\label{rd11a}
   \left(\rho^2g_a^2\right)^{\star \ell}(t,x) \geq
   \frac{K_\ell}{2^{2+2/a}}\, t^{1/a}g_a^2(t,x) \int_0^{t/2} \, d s \, [s(t-s)]^{-1/a}(t-s)^{\frac{\ell - 2}{a^*}}.
\end{align}
For $0 \leq s \leq t/2$, we have $t-s \geq s$, so we replace the last factor $t-s$ by $s$ to see that
\beqn
\left(\rho^2g_a^2\right)^{\star \ell}(t,x) \geq
  \frac{K_\ell}{2^{2+2/a}}\, t^{1/a}g_a^2(t,x)  \int_0^{t/2} \, d s\, [s(t-s)]^{-1/a}s^{\frac{\ell - 2}{a^*}}.
\eeqn
Use the change of variables $s \mapsto t-s$ in \eqref{rd11a} and add this to this last integral to see that
\beqn
   \left(\rho^2 g_a^2\right)^{\star \ell}(t,x) \geq
     \frac{K_\ell}{2^{2+2/a}}\, t^{1/a}g_a^2(t,x)
\frac{1}{2}\int_0^{t} \, d s\, [s(t-s)]^{-1/a}s^{\frac{\ell - 2}{a^*}}.
\eeqn
Then use Euler's Beta integral \eqref{rd08_12e9} to see that
\begin{align*}
   \left(\rho^2 g_a^2\right)^{\star \ell}(t,x) &\geq  \frac{K_\ell}{2^{3+2/a}}\, t^{1/a}\, g_a^2(t,x) \frac{\Gamma_E(1/a^*)\, \Gamma_E\left(1 +\frac{\ell - 2}{a^*} - \frac{1}{a}\right) }{\Gamma_E\left(1 + \frac{\ell - 1}{a^*} - \frac{1}{a}\right)}\, t^{\frac{\ell - 1}{a^*} - \frac{1}{a}} \\
    & = \frac{\rho^{2 \ell} \, \Theta_a^{\ell - 1} \, \Gamma(1/a^*)^{\ell} }{\Gamma_E(\ell / a^*)}\,  t^{(\ell - 1)/a^*} \, g_a^2(t,x).
\end{align*}
This proves \eqref{E:gaIndu}.

   Therefore, by \eqref{rdE:gaIndu}, \eqref{E:gaIndu} and \eqref{E:Mittag-Leffler},
\begin{align*}
   \cK(t,x; \rho)&\geq \sum_{k=0}^{\infty}  \left(C\, \rho^2\,  g_a^2\right)^{\star (k+1)}(t,x)\\
   &\geq c_2 \, \rho^2\, g_a^2(t, x) \sum_{k=0}^\infty \frac{\left(c_3\, \rho^2\, t^{1/a^*} \right)^k}{\Gamma_E\left(\frac{k}{a^*} + \frac{1}{a^*} \right)} \\
&= c_2\, \rho^2\, g_a^2(t,x)\, E_{1/a^*,1/a^*}\left(c_3 \, \rho^2\,  t^{1/a^*}\right).
\end{align*}
This proves the statement \eqref{rd08_17e2} in Proposition \ref{rd08_12p1} in the case where $\nu = 1$.
\medskip

\noindent{\em Step 2. The case $\nu > 0$. } Recall that $\cL_\nu = \tfrac{\partial}{\partial t} -\nu\,  \null_x D_\delta^a$. Observe that if $f_1$ and $g$ are two functions such that $\cL_1 f_1(t, x) = g(t, x)$, and if we define $f_\nu(t, x) = f_1(\nu\, t, x)$, then
\beqn
    \cL_\nu f_\nu(t, x) = \nu\, \cL_1 f_1(\nu\, t, x) = \nu\, g(\nu\, t, x).
\eeqn
When $g(t, x) = \delta_0(t, x)$, where $\delta_0$ is the Dirac delta function in $\R^2$, the elementary property $\nu\, \delta_0(\nu\, t, x) = \delta_0(t, x)$ show us that the fundamental solutions $\Gamma_\nu(t, x)$ and $\Gamma_1(t, x)$ associated respectively to $\cL_\nu$ and $\cL_1$ satisfy the relation
\beqn
    \Gamma_\nu(t, x) = \Gamma_1(\nu\, t, x).
\eeqn
Therefore,
\beqn
   \Jf_\nu(t, x; \rho) = \rho^2\, \Gamma_\nu^2(t, x) = \rho^2\, \Gamma_1^2(\nu\, t, x) = \Jf_1(\nu\, t, x; \rho).
\eeqn
Observe by induction that for $\ell \geq 1$,
\beq\label{rd05_12e1}
    (\Jf_\nu)^{\star \ell}(t, x; \rho) = \nu^{1 - \ell}\, (\Jf_1)^{\star \ell} (\nu\, t, x ; \rho),
\eeq
therefore, by \eqref{rd08_24e2} for $\nu = 1$,
\begin{align*}
    (\Jf_\nu)^{\star \ell}(t, x; \rho) &\leq \left(C_0\, \rho^2\, \Gamma_E(1/a^*)\right)^\ell\, \nu^{1-\ell}\, (\nu\, t)^{\frac{\ell}{a^*} - 1} \, \null_\delta G_a(\nu\, t,x) \, (\Gamma_E(\ell/a^*))^{-1}\\
       &= \left(C_0\, \rho^2\,\nu^{-1/a}\, \Gamma_E(1/a^*)\right)^\ell\  t^{\frac{\ell}{a^*} - 1}\, \null_\delta G_a(\nu\, t,x)\, (\Gamma_E(\ell/a^*))^{-1}.
\end{align*}
This establishes \eqref{rd08_24e2} for all $\nu > 0$. 

By definition and using \eqref{rd05_12e1}, we see that
\beqn
     \cKf_\nu(t, x; \rho) = \sum_{\ell = 1}^\infty (\Jf_\nu)^{\star \ell} (t, x; \rho) = \nu \sum_{\ell = 1}^\infty \nu^{-\ell}\, (\Jf_1)^{\star \ell}(\nu\, t, x; \rho).
\eeqn
From the definition of $\Jf_1(t, x; \rho)$, we have $\nu^{-1}\, \Jf_1(\nu\, t, x; \rho) = \Jf_1(\nu\, t, x; \rho\, \nu^{-1/2})$, so
\beqn
    \cKf_\nu(t, x; \rho) = \nu \sum_{\ell = 1}^\infty (\Jf_1)^{\star \ell}(\nu\, t, x; \rho\, \nu^{-1/2}) = \nu\, \cKf_1(\nu\, t, x; \rho\, \nu^{-1/2}).
\eeqn
Together with \eqref{rd08_17e3} and \eqref{rd08_17e2} for $\nu = 1$, this yields \eqref{rd08_17e3} and \eqref{rd08_17e2} for all $\nu > 0$. 

We note in passing that the same method could have been used in parts (a) and (b). This completes the proof of Proposition \ref{rd08_12p1}.
\end{proof}
\medskip

\noindent{\em A different method for obtaining $\cK^{\mbox{\scriptsize\rm w}}(t, x) $}
\medskip

In the case of the wave equation, formula \eqref{rd08_12e10} can also be obtained in quite a different way. Indeed, let $(u(t, x))$ be the solution of \eqref{rd08_14e1} with given initial conditions, and let $(J_0(t, x))$ be the solution of the homogeneous PDE $\cL u = 0$ with the same initial conditions. Set $f(t, x) = E[u^2(t, x)]$. Using the notations in Proposition \ref{rd08_12p1} (b), we obtain from \eqref{rd05_07e2} that
\beqn
    f(t, x) = J_0^2(t, x) + \rho^2 \int_0^t ds \int_\R dy\, \Gamma_\nu^2(t-s, x-y)\, f(s, y).
\eeqn
Let $g(t, x) := f(t, x) - J_0^2(t, x)$ and $b := \frac{\rho^2}{2 \nu} $. Since $\Gamma_\nu^2(t, x) = \frac{1}{2 \nu}\, \Gamma_\nu(t, x)$, we can write this equation
\beqn 
    g(t, x) = b \left[\Gamma_\nu \star \left(g + J_0^2(\cdot, *)\right)\right](t, x).
\eeqn
This is the integral formulation of the inhomogeneous wave equation
\beq\label{rd05_08e1}
   \cL_\nu g(t, x) = b \left(g(t, x) + J_0^2(t, x)\right), \qquad (t, x) \in \R_+^* \times \R,
\eeq
with vanishing initial conditions
\beq\label{rd05_07e3}
     g(0, x) = 0,\qquad \frac{\partial g}{\partial t}(0, x) = 0,\qquad x \in \R.
\eeq

Define an operator $\tilde \cL$ on smooth functions $h$ by 
\beqn
   \tilde \cL h(t, x) = \cL_\nu\, h(t, x) - b\, h(t, x),
\eeqn 
so that \eqref{rd05_08e1} becomes
\beq\label{rd05_08e2}
    \tilde \cL_\nu\, g(t, x) = b\, J_0^2(t, x).
\eeq
The fundamental solution associated with the operator $\tilde \cL$ is 
\beqn
    \tilde \cK(t, x) = \frac{1}{2 \nu} I_0\left(\sqrt{\frac{b}{\nu^2} \, \left((\nu t)^2 - x^2 \right)}\, \right)\,  1_{[0,\, \nu\, t]}(\vert x \vert)\, 1_{\R_+}(t),
    \eeqn
where $I_0$ is the modified Bessel function of the first kind of order $0$ (see \cite[4.1.3-1, p.~287]{polyanin}). Taking into account the vanishing initial conditions \eqref{rd05_07e3}, we conclude from \eqref{rd05_08e2} that
\beqn
    g(t, x) = b \left[\tilde \cK \star J_0^2\right](t, x),
\eeqn
or, equivalently, that
\beqn
   E[u^2(t, x)] = J_0^2(t, x) + \frac{\rho^2}{2 \nu} \left[\tilde \cK \star J_0^2\right](t, x).
\eeqn
Identifying with formula \eqref{rd08_20e4}, we see that $\cK^{\mbox{\scriptsize\rm w}}(t, x) = \frac{\rho^2}{2 \nu}\, \tilde \cK(t, x)$, which is the equality \eqref{rd08_12e10}.
\medskip

\noindent{\em Various explicit calculations}
\medskip

We now state and prove four lemmas that were used in some of the previous proofs. 
The first two are needed for the proof of Proposition \ref{rd08_12p1} (c).

\begin{lemma} \label{rd12_20l1} 
Let $b > 0$, $\nu > \half$ and for $x \in \R$, define $f_{b, \nu}(x) = (b^2 + x^2)^{-\nu - \half}$. Then for all $b > 0$ and $\xi \in \R$,
\beqn
     \cF[f_{b, \nu}](\xi) = \int_\R dx \, e^{-i \xi x}\, f_{b, \nu}(x) \geq C_\nu\, b^{- 2 \nu} \exp(- b\, \vert \xi \vert),
\eeqn
where 
\beq\label{rd01_23e2}
    C_\nu = \frac{\Gamma_E(\nu)\, \Gamma_E\left(\half\right)}{2\, \Gamma_E\left(\nu + \half\right)} .
\eeq
\end{lemma}

\begin{proof} Note that $f := f_{b, \nu}$ is an even function, so its Fourier transform is real-valued instead of complex-valued, which means that bounding $\cF f$ from below makes sense. By \eqref{rd08_26L12}, 
for ${\mathrm{Re}}(b) > 0$ and $\nu > -\half$,
\beqn
   \cF[f](\xi) = \left(\frac{\vert \xi \vert}{b} \right)^\nu \frac{\sqrt{\pi}}{2^\nu\, \Gamma_E\left(\nu + \half\right)}\, K_\nu(b\, \vert \xi \vert),
\eeqn
where $K_\nu$ is the modified Bessel function of the second kind given in \eqref{rd01_23L1}. It suffices therefore to show that the function
\beqn
     (b, \xi) \mapsto \left(\frac{\vert \xi \vert}{b} \right)^\nu K_\nu(b\, \vert \xi \vert)\, b^{2 \nu} \exp(b\, \vert \xi \vert),\qquad (b, \xi) \in \R_+^* \times \R,
\eeqn
is bounded away from zero. By choosing $u = b\, \vert \xi \vert$, we reduce this problem to bounding the function
\beqn
     u \mapsto h(u) := u^\nu e^u K_\nu(u), \qquad    u \in \R_+.
\eeqn
By the differential formula $K_\nu'(u) = - K_{\nu - 1}(u) - \frac{\nu}{z} K_\nu(z)$ (see \cite[10.29 (ii), p.~252]{olver}), we see that
\beqn
     h'(u) = e^u u^\nu (K_\nu(u) - K_{\nu - 1}(u)).
\eeqn
By the integral representation
\beqn
      K_\nu(z) = \int_0^\infty e^{- z \cosh(t)} \cosh(\nu t)\, dt, \qquad \mathrm{Arg}(z) < \frac{\pi}{2}
\eeqn
(see \cite[10.32.9, p. 252]{olver}), it follows that
\begin{align*}
         h'(u) &= e^u u^\nu \int_0^\infty e^{- u \cosh(t)} [\cosh(\nu t) - \cosh((\nu - 1) t))]\, dt \\
                 & = e^u u^\nu \half \int_0^\infty e^{- u \cosh(t)} (e^{\nu t} - e^{- (\nu - 1) t}) (1 - e^{-t})\, dt \\
                 &\geq 0
\end{align*}
since $\nu \geq \half$. It follows that
\beqn
     \inf_{u \in \R_+} h(u) = \lim_{u\downarrow 0} h(u) = 2^{\nu - 1} \Gamma_E(\nu),
\eeqn
where we have used the property $K_\nu(u) \sim \half \Gamma_E(\nu) (u/2)^{- \nu}$ as $u \downarrow 0$ (see \cite[10.30.2, p.~252]{olver}). Since $\sqrt{\pi} = \Gamma_E(\half)$ (see \cite[5.4.6]{olver}), we can take
\beqn
    C_\nu = \frac{\sqrt{\pi}}{2^\nu\, \Gamma_E\left(\nu + \half\right)}\, \inf_{u \in \R_+} h(u) = \frac{\sqrt{\pi}}{2^\nu\, \Gamma_E\left(\nu + \half\right)}\, 2^{\nu - 1} \Gamma_E(\nu)= \frac{\Gamma_E\left(\half\right)\, \Gamma_E(\nu)}{2\, \Gamma_E\left(\nu + \half\right)}.
\eeqn
This completes the proof of Lemma \ref{rd12_20l1}.   
\end{proof}

\begin{lemma}\label{rd01_24l1} 
The following properties hold:

(1) For $a>0$ and $t>0$, $g_a(t,x-y) \ge \pi\, 2^{-(a+1)} t^{1/a}g_a(t,x)g_a(t,y)$.
	
(2) For $t >0$ and $\xi \in \R$, $\cF[g_a^2(t, *)] (\xi) \geq C_{a+1/2}\, \pi^{-2}\, t^{-1/a} \exp(- t^{1/a} \vert \xi \vert)$, where the constant $C_{a+1/2}$ is defined in \eqref{rd01_23e2} (with $\nu = a + 1/2$ there).
	
(3) For all $a \geq 1$, $t \ge s > 0$ and $x\in\R$, we have
\beqn
\left[g_a^2\left(t-s,*\right) * g_a^2\left(s, * \right)\right](x) \ge \frac{C_{a+1/2}^2}{\pi^3\: 2^{2a+3}}(ts)^{-1/a} (t-s)^{1/a}g_a^2(t-s,x).
\eeqn

(4)  For $t \geq s \geq t/2>0$, $g_a(s,x) \geq (t/s)^{1/a}\, 2^{-1-1/a} \, g_a(t,x)$.
\end{lemma}

\begin{proof} 
(1) Because $$1+(u-v)^2\le 1+2u^2+2v^2\le (1+2u^2)(1+2v^2),$$ $g_a(t,x-y)$ is bounded from below by
$\pi \, t^{1/a}\, g_a\left(t,\sqrt{2} \, x\right) g_a\left(t,\sqrt{2} \, y\right)$. The conclusion follows from the inequality $g_a\left(t,\sqrt{2} \, x\right) \ge 2^{-(a+1)/2}g_a(t,x)$.

   (2) Setting $\nu = a + 1/2$ and $b = t^{1/a}$, this inequality is a special case of Lemma \ref{rd12_20l1}.

 (3) By statement (1) (with $t$ replaced by $t-s$),
\begin{align*}
  \left[g_a^2\left(t-s, *\right) * g_a^2\left(s, *\right)\right](x) &\ge \pi^2\, 2^{-2(a+1)} (t-s)^{2/a}g_a^2\left(t-s,x\right) \\
   &\qquad\qquad \times \int_\R dy\, g_a^2\left(t-s,y\right) g_a^2\left(s,y\right).
\end{align*}
Using Plancherel's identity and statement (2), we see that
\begin{align*}
  \int_\R dy\, g_a^2\left(t-s,y\right) g_a^2\left(s,y\right)
    &\geq \frac{1}{2\pi} \int_\R d \xi\, \frac{C_{a+1/2}^2}{\pi^4} ((t-s)s)^{-1/a} \exp\{-|\xi |((t-s)^{1/a}+s^{1/a}) \} \\
  &= \frac{C_{a+1/2}^2}{2\pi^5} ((t-s)s)^{-1/a} \frac{2}{(t-s)^{1/a}+ s^{1/a}} \\
	&\geq \frac{C_{a+1/2}^2}{2\pi^5} ((t-s)s)^{-1/a}\, t^{-1/a}.
\end{align*}
This proves (3).

   (4) Notice that for $t \geq s \geq t/2 > 0$,
\begin{align*}
   g_a(s,x)&=\frac{s^{-1/a}}{\pi} \left(1+\frac{x^2}{s^{2/a}}\right)^{-(a+1)/2}
   \ge \frac{s^{-1/a}}{\pi} \left(1+\frac{x^2}{(t/2)^{2/a}}\right)^{-(a+1)/2}\\
	 &\geq \frac{s^{-1/a}}{\pi} \, (t/2)^{1+1/a} \left((t/2)^{2/a} + x^2\right)^{-(a+1)/2} \\
	 & = \frac{s^{-1/a}}{\pi 2^{1+1/a}}\, t^{1+1/a} \left((t/2)^{2/a} + x^2\right)^{-(a+1)/2} \\
        &\ge \frac{s^{-1/a}}{2^{1+1/a}}\,  t^{1/a}\, g_a(t,x).
\end{align*}
This completes the proof of Lemma \ref{rd01_24l1}.
\end{proof}

The next two lemmas were used in the proofs of Proposition \ref{rd05_05p1} and Theorem \ref{rd08_18t1}.

\begin{lemma}
\label{rd08_19l1}


Let $J(t, x) = \rho^2\, \Gamma_\nu^2(t, x)$, where $\Gamma_\nu$ is defined in \eqref{rd08_21e6b}. Let $\cK = \cK^{\mbox{\scriptsize\rm w}}$, where $\cK^{\mbox{\scriptsize\rm w}}$ is defined in \eqref{rd08_12e10}.  For $t \geq 0$,
\beq
\label{rd08_19e1w}
     \int_\R dx\, \cK(t, x) = \frac{\rho}{2 \nu}\, \sinh\left(\frac{\rho\, t}{\sqrt{2 \nu}} \right).
\eeq
\end{lemma}

\begin{proof}
Because the right-hand side of \eqref{rd08_19e1w} is an even function of $\rho$, we assume that $\rho > 0$.
Use the change of variable $y = \sqrt{\rho^2\, [(\nu\, t)^2 - x^2]/(2\nu^3)}$, so that $x = \sqrt{\nu^2\, t^2  - 2\, \nu^3\, y^2/ \rho^2}$, to see that 
\beqn
    \int_\R dx\, \cK(t, x) = \frac{\rho}{\sqrt{2 \nu}} \int_0^{\rho \, t / \sqrt{2 \nu}} dy\,  
     \frac{y\, I_0(y)}{\sqrt{\rho^2 t^2 /(2\nu) - y^2}} \, .
\eeqn
In \cite[(6), p.~365]{EMOT}, we find the formula
\beqn
     \int_0^{\alpha} x^{p+1} (\alpha^2 - x^2)^{\sigma - 1} I_p(x)\, dx = 2^{\sigma-1}\, \alpha^{p+\sigma} \, \Gamma_E(\sigma)\, I_{p+\sigma}(\alpha),
\eeqn
valid for $\text{Re}(p) > -1$, $\text{Re}(\sigma) > 0$, where $I_p$ denotes the modified Bessel function of the first kind of order $p$ (see \eqref{rd08_26L10}).  
We apply this formula with $p = 0$, $\sigma = \half$ and $\alpha = \rho\, t / \sqrt{2 	\nu}$, to obtain
\begin{align*}
     \int_\R dy\, \cK(t, y) &= \frac{\rho}{\sqrt{2 \nu}} \int_0^\alpha  \frac{y}{\sqrt{\alpha^2 - y^2}}\, I_0(y)\, dy \\
      &= \frac{\rho}{\sqrt{2 \nu}}\, \frac{\sqrt{\alpha}}{\sqrt{2}}\, \Gamma_E\left(\half\right) I_{1/2}(\alpha) = \frac{\rho}{\sqrt{2 \nu}}\, \sinh(\alpha),
\end{align*}
where we have used the fact that $\Gamma_E(1/2) = \sqrt{\pi}$ (see \cite[5.4.6]{olver}) and
\beqn
   I_{1/2}(x) = \sqrt{\frac{2}{\pi\, x}} \, \sinh(x)
\eeqn
(see \cite[10.39.1]{olver}). Now replace $\alpha$ by $\rho\, t / \sqrt{2 \nu}$ to obtain \eqref{rd08_19e1w}.
\end{proof}

\begin{lemma}
The following two equalities hold:
\beq
\label{rd08_19e2}
[s \smallstar \sinh(\alpha\, s)](t) =  \alpha^{-2}\, (\sinh(\alpha\, t) - \alpha \, t),
\eeq
and
\beq
\label{rd08_19e3}
     [s^2 \smallstar \sinh(\alpha\, s)](t) = \alpha^{-3}\, (2 \cosh(\alpha\, t) - \alpha^2\, t^2 - 2).
\eeq
\end{lemma}

\begin{proof}
From \eqref{rd08_26L1} and \eqref{rd08_26L4}, 
we see that $\cLL(s)[z] = z^{-2}$ and $\cLL[\sinh(\alpha\, s)](z) = \alpha (z^2 - \alpha^2)^{-1}$, so
\beqn
   \cLL[s \smallstar \sinh(\alpha\, s)](z) = \alpha\, z^{-2}\, (z^2 - \alpha^2)^{-1}.
\eeqn
From \eqref{rd08_26L8}, 
we conclude that \eqref{rd08_19e2} holds.

    From \eqref{rd08_26L1}, 
we see that $\cLL(s^2)[z] = 2\, z^{-3}$, therefore
\beq\label{rd05_08e3}
   \cLL[s^2 \smallstar \sinh(\alpha\, s)](z) = 2 \alpha\, z^{-3}\, (z^2 - \alpha^2)^{-1}.
\eeq
From \eqref{rd08_26L5}, 
we see that $\cLL[2\cosh(\alpha\, t)](z) = 2z (z^2 - \alpha^2)^{-1}$, and by \eqref{rd08_26L2} and \eqref{rd08_26L1}, 
that $\cLL[\alpha^2\, t^2 + 2](z) = 2\, \alpha^2\, z^{-3} + 2\, z^{-1}$. Therefore, the Laplace transform of the righthand side of \eqref{rd08_19e3} is
\beqn
   \frac{1}{\alpha^3} \left(\frac{2 z}{z^2 - \alpha^2}  - \frac{2 \alpha^2}{z^3} - \frac{2}{z} \right) = \frac{2 \alpha}{z^3\, (z^2 - \alpha^2)}.
\eeqn
which matches the right-hand side of \eqref{rd05_08e3}. This establishes \eqref{rd08_19e3}.
\end{proof}

\section{Notes on Appendix \ref{app2}}
\label{AppB-B.8}
This appendix gathers analytic results that are mostly scattered in numerous articles throughout the literature. We have attempted a compilation of the more frequently used upper and lower bounds on $L^2$- and $L^p$-norms of increments of fundamental solutions and Green's functions of the heat, fractional heat and wave equations on the line and on $\rek$. 

The identies proved in Lemma \ref{ch1'-l0} can also be checked using the Fourier transform of the heat kernel and Plancherel's theorem (see e.g.~\cite{chen} and \cite{chendalang2015-2}). In the form of inequalities with non explicit constants, similar estimates can be found in many papers throughout the literature on parabolic SPDEs. However, to the best of our knowledge, inequality \eqref{1'.500} with its explicit and optimal constants, is new.

For $k=1$, some estimates related to those appearing in Lemmas \ref{app2-l2} and \ref{app2-l1} can be found in \cite[pp.~319-320]{walsh}.
Lemma \ref{app2-l1} is an extension of \cite[Lemme A2]{SaintLoubert98}, and the case $p=1$ of Lemma \ref{app2-l3} (a) appears in \cite[Lemme A3]{SaintLoubert98}. The identities in \eqref{1'.17} and \eqref{1'.3}, and the fact that both have the same right-hand side, are particularly useful in the study of the stationary pinned strings of Subsections \ref{rd03_06ss1} and \ref{ch1'-ss3.3.2}, respectively. The statements of Lemma \ref{rd08_10l3} are proved in \cite[Section 2.3]{candil2022}.

The results of Subsection \ref{rd04_24ss1} can be found in \cite{chendalang2015}. Certain related calculations also appear in \cite{dd2005}. The lower bounds in Sections \ref{app2-5} to \ref{app2-7} do not seem to appear elsewhere.

Section \ref{rd08_26s1} contains the analytical results needed for Sections \ref{rdrough} and \ref{rd1+1anderson}, and in particular, for Subsection \ref{rd01_27ss2}. Explicit values for space-time convolution series of the heat and wave kernels appeared first in \cite{chen}, \cite{chendalang2015-2} and \cite{chen-dalang-2015} (see also \cite{borodin-corwin2014}), and the upper and lower bounds in the case of the fractional heat kernel appeared in \cite{chendalang2015}. We have mostly followed these references, with an emphasis on unification and simplification. The method for computing $\cK^{\mbox{\scriptsize\rm w}}(t, x) $ described in \eqref{rd05_08e1}--\eqref{rd05_08e2} was pointed out to us by Le Chen.



\chapter{Miscellaneous results and formulas}
\label{app3}
\pagestyle{myheadings}
\markboth{R.C.~Dalang and M.~Sanz-Sol\'e}{Miscellanea}

This chapter gathers the proofs of various results and formulas that have been used throughout the book.

\section{Gronwall's lemmas}
\label{app3-s1}

Gronwall's lemmas provide estimates on real-valued functions that satisfy certain differential or integral inequalities or identities. They are instrumental in the study of evolution systems, in particular when applying fixed point arguments. In this section, we first recall a classical version of this lemma (Lemma \ref{lemC.1.1} below), then we prove a version that is well-suited to the study of SPDEs (see Lemma \ref{A3-l01}). This version contains contains a convolution integral in the time variable. Then we present a space-time Gronwall-type lemma (Lemma \ref{rd08_10l2} below), which is used in Sections \ref{new-5.5-candil} and \ref{rd1+1anderson}.
 For an extensive compilation of Gronwall's lemmas, we refer for instance to the monographs \cite{B-S-1992}, \cite{dragomir} and \cite{pachpatte}.
\medskip

Fix $T>0$ and let $J: [0,T]\rightarrow \IR_+$ be a nonnegative Borel function such that
\beq
\label{A3.00}
\int_0^T J(s)\, ds < \infty.
\eeq

\subsection{The classical Gronwall's Lemma \index{Gronwall's lemma}\index{lemma!Gronwall}}}

\begin{lemma}
\label{lemC.1.1}
Let $z,f:[0,T] \to \IR_+$ be nonnegative Borel functions. Assume that
\beq\label{rdeC.1.2a}
   \int_0^T f(s)\, J(s)\, ds < \infty
\eeq
and that, for all $t\in [0,T]$,
\beq\label{rdeC.1.3}
   f(t) \leq z(t) + \int_0^t f(s)\, J(s)\, ds.
\eeq
Then for all $t\in[0,T]$,
\beq\label{rdeC.1.4}
   f(t) \leq z(t) + \int_0^t z(s)\, J(s) \exp\left(\int_s^t J(r)\, dr \right) ds.
\eeq
In particular,
\beq\label{rdeC.1.5}
   f(t) \leq \left(\sup_{s\in[0,t]}z(s)\right)\exp\left(\int_0^t J(s)\, ds \right).
\eeq
\end{lemma}

\begin{proof}
We will show by induction that for all $n\geq 0$ and $t \in [0,T]$,
\begin{align}\nonumber
   f(t) &\leq z(t) + \int_0^t z(s)\, J(s) \sum_{k=0}^{n-1}\, \frac{1}{k!} \left(\int_s^t J(r)\, dr  \right)^k\, ds \\
   &\qquad\qquad + \int_0^t  f(s)\, J(s) \frac{1}{n!} \left( \int_s^t J(r)\, dr \right)^n\, ds,
   \label{rdeC.1.6}
\end{align}
where, by convention, if $n=0$ then the second term on the right-hand side is equal to zero.
Indeed, for $n=0$, \eqref{rdeC.1.6} reduces to \eqref{rdeC.1.3}. Assuming that \eqref{rdeC.1.6} holds for some $n\geq 0$, we show \eqref{rdeC.1.6} for $n+1$.
 Indeed, applying \eqref{rdeC.1.3}, we bound $f(s)$ from above in the last term of  \eqref{rdeC.1.6}. This yields
\begin{align}\nonumber
   f(t) &\leq z(t) + \int_0^t z(s)\, J(s) \sum_{k=0}^n\, \frac{1}{k!} \left(\int_s^t J(r)\, dr  \right)^k\, ds\\
   &\qquad + \int_0^t  \left(\int_0^s dv\, f(v)\, J(v) \right) J(s) \frac{1}{n!} \left( \int_s^t J(r)\, dr \right)^n\, ds.
   \label{rdeC.1.7}
\end{align}
Apply Fubini's Theorem to see that the last term is equal to
\beqn
   \int_0^t dv\, f(v)\, J(v) \int_v^t ds\, J(s)\, \frac{1}{n!} \left(\int_s^t dr\, J(r) \right)^n.
\eeqn
Since the $ds$-integral is equal to
\beqn
  \frac{1}{(n+1)!} \left(\int_v^t dr\, J(r) \right)^{n+1},
\eeqn
we obtain \eqref{rdeC.1.6} for $n+1$.

We now let $n \to \infty$ in \eqref{rdeC.1.6}. By monotone convergence, the second term on the right-hand side of  \eqref{rdeC.1.6} converges to the second term in \eqref{rdeC.1.4}, while the third term on the right-hand side of  \eqref{rdeC.1.6} converges to $0$ by \eqref{rdeC.1.2a} and dominated convergence. This proves \eqref{rdeC.1.4}.

From \eqref{rdeC.1.4}, we deduce that
\begin{align*}
   f(t) &\leq z(t) + \left(\sup_{s\in[0,t]} z(s)\right) \int_0^t  J(s)\, \exp\left(\int_s^t J(r)\, dr \right) ds\\
   &= z(t) + \left(\sup_{s\in[0,t]} z(s)\right) \left(\exp\left(\int_0^t J(r)\, dr \right) -1 \right) ,
\end{align*}
which implies \eqref{rdeC.1.5}.
\end{proof}

\subsection{A Gronwall-type lemma with convolution}\label{rd04_16ss1}\index{Gronwall-type lemma}\index{lemma!Gronwall-type}

Many times in this book, we need a version of Lemma \ref{lemC.1.1}, in which the function $J(s)$ in \eqref{rdeC.1.3} is replaced by $J(t-s)$. We begin by introducing some notation and a preliminary lemma.

Let $J: [0,T]\longrightarrow \re_+$ satisfy \eqref{A3.00}. Define
$$
   F(t):=\int_0^t J(s)\, ds, \qquad t\in[0,T].
$$
Throughout this section, we let the symbol ``$\smallstar$'' denote convolution\index{convolution!in time} in the time-variable: for two integrable functions $f,g:[0,T] \to \IR$ and $t\in [0,T]$, $f \smallstar g(t) = \int_0^t f(s) g(t-s)\, ds$ (this is well-defined for a.a.~$t$).

For $n\geq 1$, let $J^{\smallstar n}$ denote the $n$-th convolution power of $J$, that is, $J^{\smallstar1} = J$ and
\beqn
  J^{\smallstar (n+1)}(t) = \int_0^t J^{\smallstar n}(t-s)\, J(s)\, ds,\qquad t\in[0,T].
\eeqn
Further, define $u_0(t) \equiv 0$, and for $n\geq 1$,
\begin{align*}
   u_n(t) &= \sum_{\ell=1}^n J^{\smallstar\ell}(t), \qquad U_n(t) = \int_0^t u_n(s)\, ds,\\
   u(t) &= \sum_{\ell=1}^\infty J^{\smallstar\ell}(t), \qquad U(t) = \int_0^t u(s)\, ds.
   \label{rdeC.1.6}
\end{align*}
Notice that the function $U$ is the {\em renewal function} of \cite[Chapter 5]{asmussen}. The next lemma give some properties of the functions just defined.

\begin{lemma}
\label{A3-l00}

(a) For $t \in [0,T]$,
\beqn
   J(t) + J\smallstar u_n(t) = u_{n+1}(t),
\eeqn
and
\beq\label{rdeC.1.2}
   F(t) + F\smallstar u_n(t) = U_{n+1}(t) = F(t) + J\smallstar U_n(t).
\eeq

 (b) For all $t\in[0,T]$,
\beq
\label{A3.000}
   U_n(t) \leq U(t)\leq U(T) < \infty.
\eeq
   
 (c) For all $p> 0$ and $t \in [0, T]$,
\beq\label{rd05_02e3}
    \sum_{\ell = 1}^\infty \left(\int_0^T J^{\smallstar \ell} (s) \, ds \right)^{\frac{1}{p}} < \infty.
\eeq
\end{lemma}
\begin{proof}
(a) The equality $J+J\smallstar u_n = u_{n+1}$ follows immediately from the definitions. Observe that using elementary properties of the convolution operator and derivatives, we have
\beqn
   (F+F\smallstar u_n)' = F' + F' \smallstar  u_n = J + J\smallstar u_n = u_{n+1} = U_{n+1}',
\eeqn
therefore, the first equality in \eqref{rdeC.1.2} follows from the fact that $F(0) + F\smallstar u_n(0) = 0 = U_{n+1}(0)$. For the second equality in \eqref{rdeC.1.2}, observe that $F \smallstar  u_n = F \smallstar  U_n' = F' \smallstar  U_n = J \smallstar  U_n$.

(b) Consider the truncated Laplace transform\index{truncated Laplace transform}\index{Laplace!transform, truncated} $\hat F(z) := \int_0^T e^{-z s}\, J(s)\, ds$. By \eqref{A3.00}, $\hat F(z) < \infty$ for all $z\in\re$, and $\lim_{z\to\infty} \hat F(z) = 0$ by dominated convergence. Fix $\delta\in\, ]0,1[$ and choose $z_0$ with $\hat F(z_0) <\delta$. Since the truncated Laplace transform of a convolution product of nonnegative functions is bounded above by the product of the truncated Laplace transforms,
\begin{align}
\label{tobeused}
  \int_0^T J^{\smallstar \ell}(s)\, ds \le e^{z_0 T}\int_0^T e^{-z_0 s}\, J^{\smallstar \ell}(s)\, ds \le e^{
  z_0 T}\left[\hat F(z_0)\right]^\ell .
\end{align}
Therefore, for $t\in[0,T]$,
\beqn
U(t) \le U(T) = \sum_{\ell =1}^\infty \int_0^T J^{\smallstar \ell}(s)\, ds \le e^{z_0 T}\, \sum_{\ell=1}^\infty \delta^\ell < \infty.
\eeqn

(c) Fix $p > 0$. By \eqref{tobeused},
\beqn
 \sum_{\ell = 1}^\infty \left(\int_0^T J^{\smallstar \ell} (s) \, ds \right)^{\frac{1}{p}} \leq e^{z_0 T/p}\, \sum_{\ell = 1}^\infty\, \delta^{\ell/p } < \infty.
\eeqn
This proves the lemma.
\end{proof}

\begin{lemma}
\label{A3-l01}
Let $z_0\in\IR_+$, $[0,T]\ni t \mapsto z(t)\in \IR_+$ be a nonnegative Borel function. Let $J$ be as in \eqref{A3.00}, $u_n$, $U_n$, $u$ and $U$ be as above. Consider a sequence $(f_n,\, n\ge 0)$ of non-negative Borel functions defined on $[0,T]$. Assume that for $n\geq 1$ and $t\in[0,T]$,
\beq
\label{A3.03}
f_n(t) \le z(t) + \int_0^t \left(z_0 + f_{n-1}(s)\right) J(t-s)\, ds.
\eeq
\begin{description}
\item{(a)} For all $n\geq 1$ and  $t\in[0,T]$, we have
\begin{align}\nonumber
   f_n(t) &\le z(t)+ \int_0^t z(s)\, u_{n-1}(t-s)\, ds + z_0\ U_n(t) \\
    &\qquad\qquad + \int_0^t f_0(s)\, J^{\smallstar n}(t-s)\, ds.
\label{A3.04} 
\end{align}
In particular,
\begin{align}\nonumber
f_n(t) &\le z(t)+ \int_0^t z(s)\, u(t-s)\, ds + z_0\, U(t) \\
   &\qquad\qquad + \left(\sup_{s\in[0,t]} f_0(s)\right) \int_0^t J^{\smallstar n}(s)\, ds.
\label{rdA3.04}
\end{align}

\item{(b)} If $z_0 = 0$, $z\equiv 0$ and $\sup_{s\in[0,T]} f_0(s)\, ds < \infty$, then for all $p>0$, there is a constant $C_{T,p}<\infty$ such that for all $t\in[0,T]$,
\beq
\label{rdA3.04bis}
\sum_{n=1}^\infty\, \sup_{s\in[0,t]}\left(f_n(t)\right)^{1/p} \leq\ C_{T,p} \sup_{s\in[0,t]}\left(f_0(s)\right)^{1/p}  < \infty.
\eeq
\item{(c)} If the sequence $(f_n,\, n\ge 0)$ is constant, that is, $f_n\equiv f$ for some non-negative function $f$, and if $f$ is bounded, then for all $t\in [0,T]$,
\beq
\label{A3.04bis}
f(t)\le z(t) + \int_0^t z(s)\, u(t-s)\, ds + z_0\, U(t) .
\eeq
\end{description}
\end{lemma}

\begin{proof}
(a) For $n=1$, \eqref{A3.04} and \eqref{A3.03} are the same. Assume that $n\geq 2$ and, by induction, that \eqref{A3.04} holds for $n-1$, that is, for all $t \in [0,T]$,
$$
  f_{n-1}(t)\le z(t)+ z\smallstar  u_{n-2}(t) + z_0\ U_{n-1}(t) + f_0\smallstar  J^{\smallstar (n-1)}(t).
$$
By \eqref{A3.03} and the induction hypothesis,
\begin{align*}
   f_n(t) &\leq z(t) + \int_0^t [z_0 + z(s)+ z\smallstar  u_{n-2}(s) + z_0\ U_{n-1}(s)\\
    &\qquad\qquad\qquad + f_0\smallstar  J^{\smallstar (n-1)}(s)]\, J(t-s)\, ds \\
   &=  z(t) + \left(z_0\, F(t) + z_0\, U_{n-1}\smallstar J(t)\right)  \\
   &\qquad\qquad\qquad+ \left(z\smallstar J(t) + z\smallstar u_{n-2}\smallstar J(t)\right) + f_0\smallstar J^{\smallstar n}(t).
\end{align*}
Using Lemma \ref{A3-l00} (a), this is equal to
$$
   z(t) + z_0\, U_{n}(t) + z\smallstar u_{n-1}(t) + f_0 \smallstar  J^{\smallstar n}(t),
$$
which is the right-hand side of \eqref{A3.04}.

The inequality \eqref{rdA3.04} follows from \eqref{A3.04} and the fact that $u_{n-1} \leq u$ and $U_n \leq U$.

(b) When $z_0 = 0$ and $z \equiv 0$, \eqref{rdA3.04} implies that
$$
  \sup_{s \in [0, t]} f_n(s) \le \sup_{s \in [0, t]} f_0(s) \int_0^t J^{\smallstar n}(s)\, ds, 
$$
consequently, the conclusion follows from Lemma \ref{A3-l00} (c).


(c) In the case where $f_n\equiv f$ for some non-negative and bounded function $f$ and all $n\in\IN$, the inequality \eqref{rdA3.04} becomes
\beqn
  f(t)\le z(t)+ \int_0^t z(s)\, u(t-s)\, ds + z_0\, U(t)  + \left(\sup_{s\in[0,t]} f(s)\right) \int_0^t J^{\smallstar n}(s)\, ds.
\eeqn
Letting $n\to\infty$, we obtain \eqref{A3.04bis} since $\lim_{n\to\infty} \int_0^t J^{\smallstar n}(s)\, ds=0$ by Lemma \ref{A3-l00}(b).
 \end{proof}

\begin{remark}
In some cases, the kernel $u(t)$ defined before Lemma \ref{A3-l00} and that appears in Lemma \ref{A3-l01} can be determined explicitly or can be suitably bounded.
\smallskip

(a) If $J(t) = \rho\, e^{-a t}$ for $t \in \R^+$, with $\rho > 0$ and $a > 0$, then $u(t) = \rho\, e^{(\rho - a)t}$. Indeed, for $z \in \R$, the Laplace transform of $J$ is $\hat J(z) = \rho/(z + a)$, so the Laplace transform $\hat u$ of $u$ is $\hat u(z) = \sum_{\ell = 1}^\infty  (\hat J(z))^\ell = \rho/(z + a - \rho)$. This is the Laplace transform of $\rho\, e^{(\rho - a)t}$, which is therefore $u(t)$.
\smallskip

(b) \cite[Lemma A.2]{chen-hu-nualart-2021}  If $J(t) = \rho\, t^{-a}$ for $t \in [0, T]$, with $a \in \, ]0, 1[$ and $\rho > 0$, then 
\[
   u(t) = \rho\, t^{-a}\, \Gamma_E(1 - a) \, E_{1-a, 1-a}\left(\rho\,  \Gamma_E(1-a)\, t^{1-a}\right), \qquad t \in [0, T],
\]
where $\Gamma_E$ is the Euler Gamma function \eqref{Euler-gamma} and $E_{\alpha, \beta}$ denotes the two-parameter Mittag-Leffler function \eqref{E:Mittag-Leffler}.  In particular, if $z_0 = 0$ 
and $z(t) \equiv b$ for some $b \geq 0$, then there are nonnegative constants $C_a$ and $\gamma_a$ such that
\begin{align*}
   f(t) &\leq b \rho\, \Gamma_E(1-a)\, t^{1-a} \, E_{1-a, 2-a}\left(\rho\, \Gamma_E(1-a)\, t^{1-a}\right) \\
   &\leq b C_a\, \exp\left(\gamma_a\, \rho^{1/(a-1)}\, t \right).
\end{align*}
\smallskip

(c) \cite[Lemma 7.1.2, p.~189]{henry-1981} If, in \eqref{A3.03}, we replace $J(t-s)$ by $J(t, s) := \rho\, (t-s)^{\beta - 1}\, s^{\gamma-1}$, for $0 \leq s \leq t \leq T$, with  $\rho \geq 0$, $ \beta > 0$, $\gamma > 0$ and $\beta + \gamma > 1$, and $z(t) \equiv b$, for some $b \geq 0$, then \eqref{A3.04bis} can be replaced by 
\[
    f(t) \leq b\, E_{\beta,\gamma}\left(\left[\rho\, \Gamma_E(\beta) \right]^{1/(\beta + \gamma - 1)} t \right), \qquad t \in [0, T].
\]

\end{remark}




\subsection{A space-time Gronwall-type lemma}\index{Gronwall-type lemma!space-time}\index{space-time!Gronwall-type lemma} 

Fix $k \in \N^*$, $T > 0$ and let $J:\,  ]0, T] \times \rek \to \R_+$ be a nonnegative Borel function such that
\beq
\label{rd08_10e6}
     \int_0^T ds \int_{\rek} dy\, J(s, y) < \infty.
\eeq
We define the ``space-time" convolution\index{convolution!space-time}\index{space-time!convolution} $f \star g$ of two nonnegative functions $f, g: \, ]0, T] \times \rek \to \R_+$ by
\beq
\label{rd08_14e8}
    [f\star g](t, x) = \int_0^t ds \int_{\rek} dy \, f(t-s, x-y)\, g(s, y).
\eeq
One easily checks that if $f, g \in L^1(]0, T] \times \rek)$, then $f \star g \in L^1(]0, T] \times \rek)$. In particular, since $f$ and $g$  are nonnegative, $[f \star g](t, x) < \infty$ for a.a.~$(t, x)$.

To this kind of convolution, it is natural to associate a Fourier-Laplace operator\index{Fourier-Laplace operator}\index{operator!Fourier-Laplace} $\hcFT$, which corresponds to taking the Fourier transform in the space-variable and, since the $f$ and $g$ are only defined for $t \in\, ]0, T]$, the ``truncated" Laplace transform\index{truncated Laplace transform}\index{Laplace!transform, truncated} in the time-variable:
\beq
\label{rd08_26e3}
     \hcFT [f](z, \xi) = \int_0^T dt\, e^{-zt} \int_{\rek} dx\, e^{-i \xi \cdot x}\, f(t, x), \quad z \geq 0,\ \xi \in \rek,
\eeq
so that $\hcF := {\underline{\mathcal L}}_\infty {\mathcal F}$ is the (non-truncated) Fourier-Laplace transform. Further, $\hcFT[f] = \hcF[f 1_{]0, T]}(\cdot)]$.

This operation is well-defined provided $f \in L^1([0, T] \times {\rek})$. Notice that it is associative and commutative, and
\beqn
    \vert \hcFT f(z, \xi) \vert \leq \hcFT f(z, 0), \quad \xi \in {\rek}.
\eeqn
In addition, it has the following property.

\begin{lemma}
\label{rd08_11l1}
Consider $f, g \in L^1([0, T] \times {\rek})$. Then
\beqn
   \hcF [f\star g](z, \xi) = \hcF[f](z,\xi)\, \hcF[g](z, \xi),  \quad (z, \xi) \in \R_+ \times {\rek}.
\eeqn
If $f\geq 0$ and $g \geq 0$, then
\beqn
   \hcFT [f\star g](z, 0) \leq \hcFT[f](z,0)\, \hcF[g](z, 0), \quad z \in \R_+.
\eeqn
\end{lemma}

\begin{proof}
Let ``$\smallstar$" denote convolution in the time-variable. By definition,
\begin{align*}
   \hcFT [f\star g](z, \xi) &= \int_0^T dt\, e^{-zt} \int_{\rek} dx\, e^{-i \xi \cdot x}\, [f\star g](t, x) \\
   & = \int_0^T dt\, e^{-zt} \int_{\rek} dx\, e^{-i \xi \cdot x} \int_0^t ds \int_{\rek} dy \, f(t-s, x-y)\, g(s, y) \\
   &= \int_0^T dt\, e^{-zt} \int_0^t ds \int_{\rek} dx\, e^{-i \xi \cdot x} \int_{\rek} dy \, f(t-s, x-y)\, g(s, y) \\
   &= \int_0^T dt\, e^{-zt} \int_0^t ds\, \cF f(t-s, *)(\xi)\, \cF g(s, *)(\xi) \\
   & = \int_0^T dt\, e^{-zt} \, (\cF [f(\cdot,*)](\xi)) \smallstar (\cF [g(\cdot,*)](\xi))(t).
\end{align*}
If $T = \infty$, then this is equal to
\begin{align*}
   \cLL (\cF [f(\cdot,*)](\xi) \smallstar \cF [g(\cdot,*)](\xi))(z)
   &= \cLL (\cF [f(\cdot,*)( \xi)](z)\, \cLL(\cF [g(\cdot,*)]( \xi))(z) \\
  &= \hcF [f](z, \xi)\, \hcF [g](z, \xi),
\end{align*}
while for $T \leq \infty$, if $\xi = 0$, $f \geq 0$ and $g \geq 0$, then $\cF f(t-s, *)(0) \geq 0$ and $\cF g(s, *)(0) \geq 0$, so this is a nonnegative number that is equal to
\begin{align*}
   & \int_0^T dt\, e^{-z t} \int_0^t ds \, \cF f(t-s, *)(0)\, \cF g(s, *)(0) \\
  &\qquad = \int_0^T ds\, \int_s^T dt\, e^{-z t} \, \cF f(t-s, *)(0)\, \cF g(s, *)(0)\\
  &\qquad = \int_0^T ds\, \int_0^{T-s} dr\,  e^{-z(s+r)} \, \cF f(r, *)(0)\, \cF g(s, *)(0) \\
   &\qquad \leq \hcFT[f](z, 0)\, \hcFT[g](z, 0).
\end{align*}
This proves the lemma.
\end{proof}
\medskip

Let $J^{\star 1} = J$ and for $n \geq 1$, define
\beqn
  J^{\star (n+1)}(t, x) = \int_0^t ds \int_{\rek} dy \, J^{\star n}(t-s, x-y)\, J(s, y), \qquad t\in\, ]0,T],\ x \in {\rek},
\eeqn
set $\cK_0(t,x) = 0$ and for $n \geq 1$,
\beq
\label{rd08_10e7}
   \cK_n(t, x) = \sum_{\ell =1}^n J^{\star \ell}(t, x)
\eeq
and
\beq
\label{rd08_10e8}
   \cK(t, x) = \sum_{\ell=1}^\infty J^{\star \ell}(t, x).
\eeq

\begin{lemma}
\label{rd08_10l1}
(a) For all $(t, x) \in\, ]0, T] \times {\rek}$ and $n \geq 1$,
\beq
\label{rd08_12e1}
     J(t, x) +  (\cK_{n-1} \star J)(t, x) = \cK_{n}(t, x),
\eeq
\beq
\label{rd08_12e2}
     J(t, x) +  (\cK \star J)(t, x) = \cK(t, x),
\eeq
and
\beq
\label{rd08_12e3}
   \int_0^t ds \int_{\rek} dx\ \cK(s, x) < \infty.
\eeq

(b) Suppose that there are nonnegative nondecreasing functions $t \mapsto B_\ell(t)$, $\ell \in \N$, and a nonnegative function $(t, x) \mapsto J_0(t, x)$ such that for all $(t, x) \in\, ]0, T] \times {\rek}$,
\beq
\label{rd08_22e3}
   J^{\star \ell}(t, x) \leq B_\ell(t)\, J_0(t, x).
\eeq
Then for all $(t, x) \in\, ]0, T] \times {\rek}$,
\beq
\label{rd08_22e5}
     \cK_n(t, x) \leq  J_0(t, x) \, \sum_{\ell = 1}^n \, B_\ell(t)
\eeq 
and
\beq
\label{rd08_22e4}
     \cK(t, x) \leq  J_0(t, x) \, \sum_{\ell = 1}^\infty \, B_\ell(t).
\eeq
\end{lemma}

\begin{proof}
(a) The first two equalities follow immediately from the definition. For the proof of \eqref{rd08_12e3}, use \eqref{rd08_10e6} and  dominated convergence to see that
\beqn
    \lim_{z \to \infty} \int_0^T dt\, e^{-zt}  \int_{\rek} dx\, J(t, x) = 0.
\eeqn
Therefore,
$
    \lim_{z \to \infty} \hcFT[J] (z, 0) = 0.
$
Fix $\delta\in\, ]0,1[$ and choose $z_0 \geq 0$ such that 
$\hcFT[J](z_0, 0) <\delta$. By Lemma \ref{rd08_11l1},
\begin{align*}
   \int_0^t ds \int_{\rek} dx\, J^{\star \ell}(s, x) &\leq e^{z_0 T} \int_0^T ds\, e^{-z_0 s} \int_{\rek} dx\, J^{\star \ell}(s, x) \\
   &= e^{z_0 T}\, \hcFT[J^{\star \ell}] (z_0, 0) \leq e^{z_0 T} \left(\hcFT[J](z_0, 0) \right)^\ell .
\end{align*}
Therefore,
\begin{align}
\nonumber
   \int_0^t ds \int_{\rek} dx\ \cK(s, x) &= \sum_{\ell=1}^\infty \int_0^t ds \int_{\rek} dx\, J^{\star \ell}(t, x) \\
   & \leq e^{z_0 T}\ \sum_{\ell=1}^\infty \delta^\ell < \infty.
   \label{rd08_11e1}
\end{align}

(b) These two inequalities are an immediate consequence of the monotonicity of the $B_\ell$ and the definitions of $\cK_n$ and $\cK$.

This proves the lemma.
\end{proof}

There are several cases where the kernel $\cK$ can be computed explicitly: see Section \ref{rd08_26s1}. Some additional cases where bounds on $\cK$ can be given are discussed in \cite{joq}.

\begin{lemma}
\label{rd08_10l2}
Let $z_0 \in \R_+$ and let $a: \, ]0, T] \times {\rek} \to \R_+$ be a nonnegative Borel function. Let $J$ be as in \eqref{rd08_10e6} and $ \cK_n$ and $\cK$ as in \eqref{rd08_10e7} and \eqref{rd08_10e8}. Consider a sequence $(f_n,\, n\ge 0)$ of non-negative Borel functions defined on $]0,T] \times {\rek}$. Assume that for $n\geq 1$ and $(t, x) \in\, ]0,T] \times {\rek}$,
\beq
\label{rd08_10e9}
    f_n(t, x) \le a(t, x) + \int_0^t ds \int_{\rek} dy\, J(t-s, x - y)\, (z_0 + f_{n-1}(s, y)) .
\eeq
Then the following properties hold:
\begin{description}
\item{(a)} For all $n\geq 1$ and  $(t, x) \in\, ]0, T] \times {\rek}$, we have 
\begin{align}\nonumber
   f_n(t, x) &\le a(t, x) +  [\cK_{n-1} \star (a + z_0\, \bfone)](t, x) \\
    &\qquad\qquad + [J^{\star n} \star (f_0+z_0\, \bfone)](t, x). 
\label{rd08_10e10}
\end{align}
where $\bfone$ denotes the constant function $\bfone(t, x) \equiv 1$. In particular,
\begin{align}\nonumber
f_n(t, x) &\le a(t, x) + [\cK \star (a + z_0\, \bfone)](t, x) \\
  &\qquad\qquad
   + [J^{\star n} \star (f_0+z_0\, \bfone)](t, x) 
\label{rd08_10e11}
\end{align}
(this inequality is vacuous if $[\cK \star a](t, x)$ is infinite).
\item{(b)}  If $a\equiv 0$ and $z_0 = 0$, then \eqref{rd08_10e11} becomes
\beq\label{rd12_19_e5}
   f_n(t, x) \leq [J^{\star n} \star f_0](t, x).
\eeq


\item{(c)} If, for all $(t, x) \in \, ]0,T] \times {\rek}$, \eqref{rd08_10e9} holds with ``$=$" instead of ``$\leq$", then for all $(t, x) \in \, ]0,T] \times {\rek}$, \eqref{rd08_10e10}  also holds with equality. If for some $(t, x) \in \, ]0,T] \times {\rek}$, $\lim_{n \to \infty} [J^{\star n} \star f_0](t, x) = 0$ and $[\cK \star a](t, x) < \infty$, then $f(t, x) = \lim_{n \to \infty} f_n(t, x)$ exists and
\beq
\label{rd08_12e4}
     f(t, x) = a(t, x) + [\cK \star (a + z_0\, \bfone)](t, x).
\eeq
\item{(d)}
If the sequence $(f_n,\, n\ge 0)$ is constant, that is, $f_n\equiv f$ for some non-negative function $f$, then for all $(t, x) \in \, ]0,T] \times {\rek}$ such that $\lim_{n \to \infty} [J^{\star n} \star f](t, x) = 0$ and $[\cK \star a](t, x) < \infty$,
\beq
\label{rd08_10e13}
   f(t, x)\le a(t, x) + [\cK \star (a + z_0\, \bfone)](t, x), 
\eeq
with equality if \eqref{rd08_10e9} holds with equality.
\end{description}
\end{lemma}

\begin{proof}
(a) For $n=1$, \eqref{rd08_10e10} is a rewriting of \eqref{rd08_10e9}. Assume that $n\geq 2$ and, by induction, that \eqref{rd08_10e10} holds for $n-1$, that is, for all $(t, x) \in\, ]0,T] \times {\rek}$,
$$
  f_{n-1}(t,x)\le a(t, x)+ [\cK_{n-2} \star (a + z_0\, \bfone)](t, x)  + [J^{\star (n-1)} \star (f_0+z_0\, \bfone)](t, x).
$$
By \eqref{rd08_10e9}, the induction hypothesis and \eqref{rd08_12e1},
\begin{align*}
   f_n(t, x)) & \leq a(t, x) + [J \star (z_0\, \bfone + f_{n-1})](t, x) \\
    &\leq a(t, x) + [J \star (z_0\, \bfone + a + \cK_{n-2} \star (a + z_0\, \bfone)\\
    &\qquad\qquad + J^{\star (n-1)} \star (f_0+z_0\, \bfone))](t, x)\\
   &=  a(t, x) + [(J + J \star \cK_{n-2}) \star a](t, x) + z_0\, (J + J \star \cK_{n-2}) \star \bfone \\
    &\qquad\qquad  + [J^{\star n} \star (f_0+z_0\, \bfone)](t, x)\\
   &= a(t, x) + [\cK_{n-1}\star a](t, x) + z_0\, [\cK_{n-1} \star \bfone](t, x)\\
   &\qquad\qquad     + [J^{\star n} \star (f_0+z_0\, \bfone)](t, x),
\end{align*}
which is \eqref{rd08_10e10} for $n$. Inequality \eqref{rd08_10e11} follows since $\cK_n \leq \cK$.

(b) If $a \equiv 0$ and $z_0 = 0$, then \eqref{rd12_19_e5} is a rewrite of \eqref{rd08_10e11}. 

(c) If \eqref{rd08_10e9} holds with ``$=$", then the two inequalities in the calculations for part (a) become equalities, so \eqref{rd08_10e10} holds with equality and can be written
\begin{align*}
   f_n(t, x) &= a(t, x) +  [\cK_{n-1} \star (a + z_0\ \bfone)](t, x) + [J^{\star n} \star (f_0+z_0\, \bfone)](t, x) \\
       &= a(t, x) +  [\cK_{n-1} \star (a + z_0\, \bfone)](t, x)  + [J^{\star n} \star f_0] (t, x)\\
       &\qquad \qquad + z_0 \int_0^t ds \int_{\rek} dy\,  J^{\star n}(s, y).
\end{align*}

Suppose that for some $(t, x) \in \, ]0,T] \times {\rek}$, $\lim_{n \to \infty} [J^{\star n} \star f_0](t, x) = 0$ and $[\cK \star a](t, x) < \infty$. The integral converges to $0$ by \eqref{rd08_12e3} and $[J^{\star n} \star f_0](t, x)$ converges to $0$ by assumption. Since $\cK_n \uparrow \cK < \infty$ a.e.~by \eqref{rd08_12e3} and $[\cK_{n-1} \star a](t, x) \leq [\cK \star a](t, x) < \infty$ by assumption, we see that $f(t) = \lim_{n \to \infty} f_n(t)$ exists and \eqref{rd08_12e4} holds by monotone convergence.

(d) If $f_n\equiv f$, then $f = \lim_{n \to \infty} f_n$ and \eqref{rd08_10e13} follows from
\eqref{rd08_10e10} and the arguments in part (c).
\end{proof}


 \section{Integral representation of weak solutions to PDEs}
 \label{app3-s4-new}

In the theory of PDEs, there are two main notions of solutions: {\em classical solutions} as have been described in Section \ref{ch1'-s1}, for which an integral representation is provided via the fundamental solution or the Green's function associated with the PDE operator $\cL$; and {\em weak solutions} for which an equation involving the evaluation of the solution on suitable test-functions is given. 

In the next proposition, we study the relationship between these two notions in the particular case of a non-homogeneous heat equation on $[0,L]$ with vanishing initial condition and vanishing Dirichlet boundary conditions.

More precisely, consider the inhomogeneous linear PDE
\beq
\label{app3-s4-new-1}
\begin{cases}
\frac{\partial u}{\partial t}(t,x)-\frac{\partial^2u}{\partial x^2}(t,x) = \varphi(t,x), & (t,x)\in\,]0,T[\times ]0,L[,\\
u(0,x)= 0, & x\in[0,L],\\
u(t,0)=u(t,L)=0, & t\in\,]0,T].
\end{cases}
\eeq
We assume that $\varphi \in  L^2([0,T]\times [0,L])$, and let $\cC:=\cC([0,T]\times [0,L])$.


Let $G_L$ be the Green's function given in \eqref{ch1'.600}, $f\in\cC$, $\varphi\in L^2([0,T]\times [0,L])$ and consider the following property:
\medskip

\noindent $(\bf{R})$\ For all
$(t,x)\in[0,T]\times [0,L]$,
\beq
\label{app3-s4-new-2}
f(t,x) = \int_0^t ds\int_0^L dy\, G_L(t-s;x,y)\, \varphi(s,y).
\eeq
Notice that if \eqref{app3-s4-new-2} holds, then $f\in\cC$ by \eqref{1'1100}.
Moreover, let $\hac$ be the reproducing kernel Hilbert space of the Gaussian random field
 \beqn
\left(v(t,x)= \int_0^t\int_0^L G_L(t-s;x,y)\, W(ds,dy),\ (t,x)\in[0,T]\times [0,L]\right),
\eeqn
where $\dot W$ is a space-time white noise. By Lemma \ref{ch2'-s3.4.4-l2}, for every $f \in \hac$, there is $\varphi \in L^2([0, T] \times [0, L])$ such that property $({\bf R})$ holds.

In the next proposition,  $\cL= \frac{\partial}{\partial t}- \frac{\partial^2}{\partial x^2}$ and $\cL^\star$ denotes (the opposite of) its adjoint:\label{rd06_13p2}
  $ \cL^\star= \frac{\partial^2}{\partial x^2} + \frac{\partial}{\partial t}$.

\begin{prop}
\label{app3-s4-new-p1}
Given $f\in\cC$ and $\varphi\in L^2([0,T] \times [0, L])$, Property $({\bf R})$ is equivalent to the following property:
\medskip

\noindent $({\bf P})$ For all $t \in [0, T]$, for all $\psi \in C^{1,2}([0,t] \times [0, L])$ such that $\psi(\cdot, 0) =\psi(\cdot, L)=0$, we have
\begin{align}
\label{L1-(*1)-AppC}
   \int_0^L dx\, f(t, x)\, \psi(t,x) &= \int_0^t ds \int_0^L dx\, f(s, x)\, \cL^\star \psi(s, x)\notag\\
   &\qquad  +  \int_0^t ds \int_0^L dx\, \psi(s,x)\, \varphi(s, x).
   \end{align}
   \end{prop}
   
   \begin{remark}
   \label{app3-s4-new-r1}
   Property $({\bf P})$ is the statement that $f$ is a weak solution to \eqref{app3-s4-new-1}, when we use as test functions the set of $\psi \in C^{1,2}([0,t] \times [0, L])$ such that $\psi(\cdot, 0) =\psi(\cdot, L)=0$.
 Property $({\bf R})$ is the statement that the solution $f$ to \eqref{app3-s4-new-1} has the integral representation \eqref{app3-s4-new-2} in terms of the Green's function of $\cL$. Proposition \ref{app3-s4-new-p1} states that for $f \in \cC$ and $\varphi \in L^2([0, T] \times [0, L]$, these two notions of solution are equivalent.
  \end{remark}

\noindent{\em Proof of Proposition \ref{app3-s4-new-p1}}.\
First, we assume Property $({\bf R})$ and prove that Property $({\bf P})$  holds.
Fix $t > 0$, $h \in \cC([0, t] \times [0,L])$ and for $(s, y) \in [0, t] \times [0, L]$, define
   \beqn
   G(h)(s,y) =
   \begin{cases}
    \int_0^L dx\,  h(s,x) G_L(s; x, y),&  {\text{if}}\  s > 0,\\
    h(0,y),&  {\text{if}}\  s = 0.
  \end{cases}
  \eeqn
We will use the following fact:
\medskip

\noindent{\em Fact 1}. For $h \in \cC^{1,2}([0,t] \times [0, L])$ such that $h(s, 0)=h(s, L)=0$ for all $s \in [0, t]$,
\beqn
     G h(s, y) = h(0, y) + \int_0^s dr\, G(\cL^\star h)(r, y).
\eeqn
Indeed, fix $\ep \in\, ]0, s[$. Observe that
\begin{align*}
     \int_\ep^s dr\, G(\cL^\star h) (r, y) &= \int_\ep^s dr  \int_0^L dx\,  \cL^\star h(r, x)  G_L(r; x, y)\\
        &=  \int_\ep^s dr \int_0^L dx \left(\frac{\partial^2}{\partial x^2} h(r, x) + \frac{\partial}{\partial r} h(r, x)\right) G_L(r; x, y).
        \end{align*}
Integrating twice by parts with respect to $x$ in the first term of the last expression, and using the fact that the product terms vanish because of the boundary conditions, we see that this term is equal to
\[
    \int_\ep^s dr \int_0^L dx\, h(r, x) \frac{\partial^2}{\partial x^2} G_L(r; x, y).
\] 
Because $\cL G_L(\cdot\, ; *, y) = 0$, we obtain,
\begin{align*}
   \int_\ep^s dr\, G(\cL^\star h)(r, y) &= \int_\ep^s dr \int_0^L dx \left(h(r, x) \frac{\partial}{\partial r}  G_L(r; x, y)\right.\\
   &\left.\qquad\qquad\qquad\qquad + \frac{\partial}{\partial r}  h(r, x) G_L(r; x, y)\right)\\
   &= \int_0^L dx\, \big[ h(r, x) G_L(r; x, y) \big]_\ep^s \\
   &= G(h)(s, y) - G(h)(\ep, y).
   \end{align*}
Since $h$ is continuous, we can let $\ep \downarrow 0$ to obtain Fact 1.
\medskip

We now establish Property $({\bf P})$.
Observe that
for $f$ and $\varphi$ satisfying \eqref{app3-s4-new-2}, and for $\psi$ as in Property $({\bf P})$,
\begin{align*}
  &\int_0^L dx\, f(t,x)\, \psi(t,x) - \int_0^t du \int_0^L dx\, f(u, x)\, \cL^\star \psi(u, x)\\
  &\qquad\qquad  - \int_0^t ds \int_0^L dy\, \psi(s, y)\, \varphi(s, y)\\
   &\quad=  \int_0^L dx\, \psi(t,x) \int_0^t ds \int_0^L dy\, G_L(t - s; x, y)\, \varphi(s, y)\\
   &\qquad\qquad - \int_0^t du \int_0^L dx\, \cL^\star \psi(u, x) \int_0^u ds \int_0^L dy\, G_L(u - s; x, y)\, \varphi(s, y)\\
   &\qquad\qquad  - \int_0^t ds \int_0^L dy\, \psi(s, y)\, \varphi(s, y).
   \end{align*}
 Apply Fubini's theorem in the first and second terms and integrate first in $x$ to see that this is equal to
 \begin{align}
 \label{L1-(*3)}
   &\int_0^t ds \int_0^L dy\, \varphi(s, y) \Big[G(\psi(s + \cdot, *))(t - s, y)\notag\\
   &\qquad\qquad\qquad  - \psi(s,y) -\int_s^t du\, G(\cL^\star\psi(s+\cdot,\ast))(u-s,y)\Big].
   \end{align}
The $du$-integral is equal to
\beqn
\int_0^{t - s} dv\,  G(\cL^\star\psi(s + \cdot, *))(v, y) = G(\psi(s + \cdot, *))(t-s, y) - \psi(s, y)
\eeqn
by Fact 1, so the expression in \eqref{L1-(*3)} is equal to 0. This shows that Property $({\bf P})$ holds.

Next, we assume that Property $({\bf P})$ holds and we prove Property $({\bf R})$.
   Fix $h \in \cC_0^\infty(]0,L[)$ and $ t\in\, ]0, T]$. For $s \in [0, t[$, define
   \beq
   \label{L1-(*1a)}
        \psi(s, y) = G(t - s, h, y) := \int_0^L dz\,  h(z) G_L(t - s; z, y),
        \eeq     
and $\psi(t, y) = h(y) $.

We have the following:
\medskip

\noindent{\em Fact 2}. $\psi \in \cC^{1,2}([0,t] \times [0, L])$\ and for $(s, x) \in\, ]0, t[ \times ]0, L[$, $\cL^\star \psi(s,x) = 0$, $\psi(s,0) =\psi(s, L)=0$.
\medskip

Indeed, since the Green's function $G_L(t-s;z,y)$ is symmetric in the space variables, $\psi$ solves the backwards heat equation with terminal condition $h$ at time $t$. This solution satisfies the conditions stated in Fact 2 (see e.g.~\cite[Theorem 7, Section 7.1.3, p. 367]{evans}, \cite{ei70}, or \cite[Chapter 1]{friedman64}. This ends the proof of Fact 2.
\medskip

 By Property $({\bf P})$, and because the first term on the right-hand side of \eqref{L1-(*1)-AppC} vanishes by Fact 2,
 \beqn
  \int_0^L dy\, f(t, y)\, \psi(t, y)
   = \int_0^t ds \int_0^L dy\, \psi(s,y) \varphi(s, y).
   \eeqn
Equivalently,
\beq
\label{L1-(*2)}
   \int_0^L dy\, f(t,y)\, h(y) = \int_0^t ds \int_0^L dy\, \psi(s,y)\, \varphi(s, y).
   \eeq     
Fix $y \in\, ]0, L[$ and let $h_0$ be a nonnegative function with compact support contained in $]-1, 1[$ such that $\int_{-1}^1 dy\, h_0(y) = 1$, and for $n \geq 1$, set  $h_n(y) = n h_0(n y)$. The following approximation
holds:
\medskip

\noindent{\em Fact 3}. For $x \in ]0, L[$, the sequence
\[
    \left(\int_0^L dz\,  h_n(x - z)  G_L(t - \cdot\, ; z, *),\ n\in\N\right)
\] 
converges to $G_L(t- \cdot\, ; x, *)$ in $L^2([0,t] \times [0, L])$.
\medskip

Indeed, for $n$ large enough, the difference of the two quantities is
\beqn
\int_0^L dz\,  h_n(x - z)\,  (G_L(t - \cdot\, ; z, *) - G_L(t- \cdot\, ; x, *)).
\eeqn
By Minkowski's inequality,
\begin{align*}
  & \left\Vert \int_0^L dz\,  h_n(x - z)  (G_L(t - \cdot\, ; z, *) - G_L(t- \cdot\, ; x, *)) \right\Vert_{L^2([0,t] \times [0, L])}\\
  & \qquad\leq \int_0^L dz\,  h_n(x - z)  \left\Vert G_L(t - \cdot\, ; z, *) - G_L(t- \cdot\, ; x, *) \right\Vert_{L^2([0,t] \times [0, L])}\\
   &\qquad \leq  \int_0^L dz\,  h_n(x - z) \, \vert z - x \vert^\half \leq \int_{-1}^1 dy\, h_0(y)\left(\frac{\vert y \vert}{n}\right)^\half,
   \end{align*}
where in the second inequality, we have used \eqref{1'1100}. This converges to $0$ as $n \to \infty$ and therefore, Fact 3 is proved.

   For $x\in\, ]0,L[$, replace $h(*)$ by $h_n(x - *)$ in \eqref{L1-(*1a)}, yielding a function $\psi_n(s,y) = \int_0^L dz\,  h_n(x - z) G_L(s; z, y)$, and in \eqref{L1-(*2)}, replace $h(*)$ by $h_n(x - *)$ and $\psi$ by $\psi_n$. Since $f$ is continuous, the left-hand side of \eqref{L1-(*2)} converges, and the right-hand side of \eqref{L1-(*2)} converges by Fact 3, yielding
   \beqn
    f(t,x) = \int_0^t ds \int_0^L dy\, G_L(t-s, x, y)\,  \varphi(s,y).
    \eeqn
    Since both $f$ and the right-hand side of this equality are continuous functions, this equality extends to
    $(t, x) \in [0, T] \times [0, L]$ by continuity.
This is property \eqref{app3-s4-new-2}, and Property ({\bf R}) is proved.
\hfill\qed

\begin{remark}
\label{app3-s4-new-r1-neumann}
The same statement and proof is valid if, in \eqref{app3-s4-new-1}, we replace the vanishing Dirichlet boundary conditions by vanishing Neumann boundary conditions. The only change is that in Property $({\bf P})$ (respectively Fact 1, Fact 2), we should replace the vanishing Dirichlet boundary conditions for $\psi$ (respectively $h$, $\psi$) by vanishing Neumann boundary conditions.
\end{remark}
 \section{Fundamental solutions via Fourier transforms and Fourier series}
 \label{AppC-New-0}
  In this section, with show (informally) how to find the fundamental solution of a partial differential operator using Fourier transforms or Fourier series.
 This is used for instance in \eqref{fs-fractional} and in Section \ref{ch1-s3}.
 
 \subsection{On the whole space}
 
Consider a partial differential operator $\mathcal{L}$ acting on smooth functions $f:\re_+\times \rek \rightarrow \re$. Assume that $\mathcal{L}$ takes the form
$\mathcal{L} = \frac{\partial}{\partial t} - \cA$, where $\cA$ is a partial differential operator acting on smooth functions $\tilde f:\rek\longrightarrow \re$. Suppose that for any smooth and integrable function $\tilde f$ such that  $\cA\tilde f$ is also integrable,
\beq
\label{pde0}
\tf(\cA\tilde f)(\xi) = l(\xi)\tf \tilde f(\xi),
\eeq
for some $l:\rek\rightarrow \re$, where $\tf$ denotes the Fourier transform (see \eqref{defi-fourier}). The function $l$ in \eqref{pde0} is called the Fourier multiplier\index{multiplier!Fourier}\index{Fourier!multiplier} of $\cA$.

We put ourselves in a ``regular setting'', meaning that all the performed operations make sense, and formulate the following statement:

Consider the PDE
\beq
\label{pde1}
\begin{cases}
\left(\frac{\partial}{\partial t} - \cA\right) u(t,x) = g(t,x),& (t,x)\in\, ]0,T]\times \rek, \\
u(0,x)= u_0(x),& x\in \rek.
\end{cases}
\eeq
Then for $(t,x)\in[0,T]\times \re$,
\begin{align}
\label{pde2}
u(t,x) &= \int_{\rek} G(t,x-y) u_0(y)\ dy + \int_0^t ds \int_{\rek}dy\ G(t-s,x-y) g(s,y)\notag\\
&= (G(t,\ast)\ast u_0)(x) + \int_0^t ds \left(G(t-s,\ast)\ast g(s,\ast)\right)(x),
\end{align}
where
\beq
\label{pde3}
G(t,x) = \tf^{-1}\left[\exp(t \, l(\ast))\right](x).
\eeq

Indeed, this can be justified as follows. Take the Fourier transform in the $x$-variable on both sides of \eqref{pde1} to obtain
\beq
\label{pde4}
\begin{cases}
\frac{\partial}{\partial t} \tf u(t,\ast)(\xi)- l(\xi)\, \tf u(t,\ast)(\xi) = \tf g(t,\ast)(\xi),& (t,\xi)\in\, ]0,T]\times \rek,\\
\tf u(0,\ast)(\xi)= \tf u_0(\ast)(\xi),& \xi\in \rek.
\end{cases}
\eeq
Notice that for fixed $\xi\in\rek$, this is a first order linear inhomogeneous ODE for the function $t\mapsto\tf u(t,\ast)(\xi)$. Therefore, we can solve \eqref{pde4}
using for instance the method of variation of constants. This yields
\beqn
\tf u(t,\ast)(\xi) = \tf u_0(\xi)\exp\left[t\, l(\xi)\right] + \int_0^t \exp\left[(t-s)\, l(\xi)\right]\tf g(s,\ast)(\xi)\, ds.
\eeqn
By applying the inverse Fourier transform in the $\xi$-variable, this identity yields
\begin{align*}
u(t,x) &= \left[u_0\ast\tf^{-1}\left[\exp(t\, l(\ast))\right]\right](x)\\
&\qquad + \int_0^t ds\, \tf^{-1}\left[\exp((t-s)\, l(\ast))\ast g(s,\ast)\right](x).
\end{align*}
With the definition \eqref{pde3}, this is formula \eqref{pde2}.

Notice that
\beqn
\lim_{t\downarrow 0} \tf\left[G(t,\cdot)\right](\xi) = \lim_{t\downarrow 0} \exp(t\, l(\xi)) = 1,
\eeqn
and consequently,
\beq
\label{pde5}
\lim_{t\downarrow 0} G(t,x) = \delta_{0}(x).
\eeq
\bigskip

\noindent{\em Examples}
\medskip

\noindent{\em1.} Let $\cA=\Delta$ be the Laplacian operator on $\rek$. In this case, $l(\xi) = -|\xi|^2$, $\xi\in\rek$, and thus
\beqn
G(t,x) = \tf^{-1}\left[\exp((-t\, |\cdot|^2)\right](x).
\eeqn
For $t>0$, the function $\varphi(\xi) = \exp\left(-t\, |\xi|^2\right)$ is the characteristic function of the $k$-dimensional Gaussian probability distribution ${\rm N} _k(0,2t {\rm Id}_k)$. Therefore,
\beqn
G(t,x) = 
      (4\pi t)^{-\frac{k}{2}} \exp\left(-\frac{|x|^2}{4t}\right)\,1_{]0,\infty[}(t),\qquad (t, x) \in \R_+ \times \rek.
\eeqn
\medskip

\noindent{\em 2.} Let $\cA = \null_xD_\delta^a$ be the Riesz-Feller fractional derivative on $\R$ defined by \eqref{ch1'-ss7.3.1} or \eqref{ch1'-ss7.3.2}, with $a\in\, ]1,2[$, $|\delta|\le 2-a$.

From \eqref{ch1'-ss7.3.1}, it follows that \eqref{pde0} holds with
\beqn
l(\xi) = -|\xi|^a \exp\left(-i\pi\delta\ {\rm{sgn}}(\xi)/2\right) =\null_\delta\psi_a(\xi).
\eeqn
From \eqref{pde3}, which defines a function when $t>0$ (because $\exp[t\, l(\ast)]$ is integrable), we obtain
\beq
 \label{fs-fractional.4}
 \null_\delta G_a(t,x) = \tf^{-1}\left[ \exp\left(t\, \null_\delta\psi_a(\cdot)\right)\right](x)\,1_{]0,\infty[}(t),\qquad (t, x) \in \R_+ \times \R,
 \eeq
which is \eqref{fs-fractional}.

\subsection{On a bounded domain}
Consider a bounded domain $D\subset\rek$ with smooth boundary $\partial D$, and the PDE 
\beq
\label{pde22}
\begin{cases}
\left(\frac{\partial}{\partial t} - \cA\right) u(t,x) = g(t,x),& (t,x)\in\, ]0,T]\times D, \\
u(0,x)= u_0(x),& x\in D,\\
u(t,x)=0, & x\in\partial D.
\end{cases}
\eeq
Assume that there exists a CONS $(e_n,\, n\ge 1)$ of $L^2(D)$ consisting of eigenvectors of the operator $\cA$ with corresponding eigenvalues $(\lambda_n,\, n\ge 1)$. For $n \geq 1$, let
\beqn
a_n(t) = \langle u(t,\ast), e_n\rangle,\quad g_n(t) = \langle g(t,*), e_n\rangle,\quad  a_n(0) = \langle u_0,e_n\rangle,
\eeqn
where $\langle \cdot, \cdot \rangle$ denotes the inner product in $L^2(D)$.

Formally, for fixed $t \in [0, T]$, we can expand $u(t, *)$ in the CONS $(e_n,\, n\ge 1)$, so that
\beq
\label{projections-1}
u(t,x) = \sum_{n=1}^\infty a_n(t)\, e_n(x).
\eeq
By taking the inner product with $e_n$ in \eqref{pde22}, we see that the functions $a_n(t)$, $n\ge 1$, must satisfy the ODEs
\beqn
\begin{cases}
a_n'(t) - \lambda_n a_n(t) = g_n(t),\quad  t>0,\\
a_n(0)= \langle u_0,e_n\rangle,
\end{cases}
\eeqn
whose solution is
\beqn
a_n(t) = \left[a_n(0)+\int_0^t ds\, g_n(s) e^{-\lambda_n s}\right] e^{\lambda_nt}.
\eeqn
Turning back to \eqref{projections-1}, we see that
\begin{align*}
u(t,x) &=  \sum_{n=1}^\infty \left[a_n(0)\, e^{\lambda_n\, t}  + \int_0^t e^{\lambda_n(t-s)} g_n(s) \, ds \right]  e_n(x) \\
    &= \int_D dy\, G(t;x,y)\, u_0(y) + \int_0^t ds \int_D dy\, G(t-s;x,y)\, g(s,y).
\end{align*}
with
\beqn
G(t;x,y) = \sum_{n=1}^\infty e^{\lambda_n t} e_n(x) e_n(y),\qquad  t>0,\ x,y\in D.
\eeqn
This is the Green's function associated to the operator $\cL$ on $D$ with vanishing Dirichlet boundary conditions.

In Section \ref{ch1-s3} we used this approach in the case where $\cA$ is the Laplacian operator on $]0,L[$, $e_n(x):=e_{n,L}(x) = \sqrt{\frac{2}{L}} \sin\left(\frac{n\pi}{L}x\right)$, $n\ge 1$ and $\lambda_n= - \frac{\pi^2}{L^2} n^2$. 

\section{A stability result for partial differential equations}
\label{appC-New-1}

In this section, we give two auxiliary results that are used in Section \ref{ch5-ss-in-rev}. These concern stability of solutions to deterministic equations with respect to the driving coefficients.

Recall that $\bD$ is the set of continuous functions $f:[0,L] \to \re$ such that $f(0)=f(1)=0$, and $V= L^2([0,1])$.
Let $G_{a,b}(t; x, y)$ be the Green's function defined in \eqref{s5.3-new(*1)green}, $g \in \bD$, $\zeta : [0, T] \times [0, 1] \to \R$ be a  continuous function, and $F : \bD \to V$ be a Lipschitz function.
In order to simplify the notation, we will write $G(t; x, y)$ instead of $G_{a,b}(t; x, y)$.

Consider the (deterministic) equation
\begin{align}
\label{nbr-(*15b)}
      \kappa(t, x) &=  \int_0^t  ds\, \langle G(t-s; x, *), F(\kappa(s, *)) \rangle_V \notag\\
      &\qquad\qquad+\langle G(t; x, *), g \rangle_V  + \zeta(t,x),\qquad t \geq 0,\ x \in [0, 1].    
      \end{align}
 We will prove a result  about stability of the solution $\kappa(t, x)$ to equation \eqref{nbr-(*15b)}
 with respect to $\zeta$, the coefficient $F$ and the initial condition $g$ (see Proposition \ref{nbr-*3aa}).
 At several points, the inequality
      \beq
      \label{bound-for-G}
      \Vert G(t;x,\ast)\Vert_V \le (2\pi t)^{-\frac{1}{4}},
      \eeq which follows from Lemma \ref{ch1'-pPD} (ii) and \eqref{heatcauchy-11'}, will be used
      \begin{lemma}
      \label{nbr-*3ab}
      The equation \eqref{nbr-(*15b)} has a solution $\kappa \in \cC([0, T], \bD)$ and this solution is unique.
      \end{lemma}

\begin{proof}
Define $H_1, H: \cC([0, T]; \bD) \to  \cC([0, T]; \bD)$ by
\begin{align}
\label{nbr-(*15ac)}
      H_1(\kappa)(t,x) &= \int_0^t  ds\, \langle G(t-s, x, *), F(\kappa(s, *)) \rangle_V, \notag\\
      H(\kappa)(t,x) &= \langle G(t, x, *), g \rangle_V + H_1(\kappa)(t,x) + \zeta(t,x).       
      \end{align}
Notice that for two functions $\kappa^{(1)}, \kappa^{(2)} \in \cC([0, T], \bD)$,
\begin{align}
\label{nbr-(*15ab)}
 & \Vert H(\kappa^{(1)})(t, *) - H(\kappa^{(2)})(t, *)  \Vert_\bD\notag\\
  &\qquad=  \Vert H_1(\kappa^{(1)})(t, *) - H_1(\kappa^{(2)})(t, *)  \Vert_\bD\notag\\
  & \qquad= \sup_{x \in [0, 1]} \left\vert \int_0^t  ds\, \langle G(t-s, x, *), F(\kappa^{(1)}(s, *))
- F(\kappa^{(2)}(s, *)) \rangle_V \right \vert\notag\\
  &\qquad\leq \sup_{x \in [0, 1]} \int_0^t ds\,  \Vert G(t-s, x, *) \Vert_V\,  \Vert F(\kappa^{(1)})(s, *) - F(\kappa^{(2)})(s, *)\Vert_V\notag\\
   &\qquad\leq C \int_0^t ds\, (t - s)^{- \frac{1}{4}}\, \Vert \kappa^{(1)}(s, *) - \kappa^{(2)})(s, *)\Vert_\bD,     
  \end{align}
where in the last inequality, we have used \eqref{bound-for-G} and the Lipschitz property of $F$.

Consider a sequence $(\kappa_n,\, n \in \N)$ defined by $\kappa_0 = \zeta$ and, for $n \geq 1$,
    $\kappa_n(t,x) = H(\kappa_{n-1})(t,x)$. Set
   $ f_n(t) = \Vert \kappa_n(t, *) - \kappa_{n-1}(t, *) \Vert_\bD$. From
\eqref{nbr-(*15ab)}, we deduce that
\beq
      f_n(t) \leq  C\int_0^t ds\, (t - s)^{- \frac{1}{4}}\, f_{n-1}(s),\quad t\in[0,T],
      \eeq
and from Lemma \ref{A3-l01}, we conclude that
\beqn
   \sum_{n=1}^\infty\, \sup_{t \in [0, T]} \Vert \kappa_n(t, *) - \kappa_{n-1}(t, *) \Vert_\bD < \infty.
   \eeqn
Therefore, there is $\kappa \in \cC([0, T]; \bD)$ such that
\beqn
    \lim_{n \to \infty} \sup_{(t,x) \in [0, T] \times [0, 1]} \vert \kappa_n(t, x) - \kappa(t,x) \vert = 0.
    \eeqn
By the continuity property of $F$, we see that $\kappa$ satisfies \eqref{nbr-(*15ac)}. Uniqueness follows immediately from \eqref{nbr-(*15ab)} and Lemma \ref{A3-l01}.
  \end{proof}


 We recall that a sequence of functions $F_n: \bD \to V$, $n\ge 1$, is {\em uniformly Lipschitz continuous}\index{uniformly Lipschitz continuous}\index{Lipschitz!continuous, uniformly} if there is $K < \infty$ such that, for all $w_1, w_2 \in \bD$ and all $n \ge 1$,
\beq
\label{nbr-(*14abb)}
    \Vert F_n(w_1) - F_n(w_2) \Vert_{V} \leq K \Vert w_1 - w_2 \Vert_{\bD}.    
 \eeq

  \begin{prop}
  \label{nbr-*3aa}
  Let $\zeta_n: [0, T] \times [0, 1]$, $n\ge 1$, be continuous functions such that
  \beqn
     \sup_{t \in [0, T]}  \Vert \zeta_n(t,*) - \zeta(t,*) \Vert_\bD \to 0 \quad{\text as}\ n \to \infty.
     \eeqn
Let $g_n$, $g \in \bD$ be such that $\lim_{n \to \infty} \Vert g_n - g \Vert_{\bD} = 0$. For a sequence of uniformly Lipschitz functions $F_n: \bD \to V$, $n\ge 1$, consider the equations
\begin{align}\nonumber
    \kappa_n(t, x) &= \langle G(t, x, *), g_n \rangle_V \\
    &\qquad + \int_0^t  ds\, \langle G(t-s, x, *), F_n(\kappa_n(s, *)) \rangle_V + \zeta_n(t,x),    
\label{nbr-(*15a)}
\end{align}
$(t,x)\in[0,T]\times [0,1]$, $n\ge 1$, and the equation \eqref{nbr-(*15b)}. If for all $w \in \bD, F_n(w) \to F(w)$ in $V$, then $\kappa_n(t,*) \to \kappa(t,*)$ in $\bD$, uniformly on $[0, T]$.
\end{prop}

\begin{proof}
 From \eqref{nbr-(*15a)} and \eqref{nbr-(*15b)}, we have
 \begin{align*}
    \Vert \kappa_n(t, *) - \kappa(t, *) \Vert_\bD&
    =  \sup_{x \in [0, 1]} \Big[\Big\vert \langle G(t, x, *), g_n - g\rangle_V\\
     &\quad+ \int_0^t  ds\, \langle G(t-s, x, *), F_n(\kappa(s, *)) - F(\kappa(s, *)) \rangle_V\\
    &\quad+  \int_0^t  ds\, \langle G(t-s, x, *), F_n(\kappa_n(s, *)) - F_n(\kappa(s,*)) \rangle_V\\
    &\quad+ \zeta_n(t,x) - \zeta(t,x)\Big\vert\Big]\\
    &\leq I_{n,0}(t) + I_{n, 1}(t)  +  I_{n, 2}(t) +  I_{n, 3}(t),
    \end{align*}
where
\begin{align*}
      I_{n,0}(t) &= \sup_{x \in [0, 1]} \vert \langle G(t, x, *), g_n - g \rangle_V \vert,\\
      I_{n, 1}(t) &= \sup_{x \in [0, 1]}  \int_0^t  ds\, \langle G(t-s, x, *), \vert F_n(\kappa(s, *)) - F(\kappa(s, *)) \vert  \rangle_V,\\
      I_{n, 2}(t) &=  \sup_{x \in [0, 1]} \int_0^t  ds\, \langle G(t-s, x, *), \vert F_n(\kappa_n(s, *)) - F_n(\kappa(s, *)) \vert  \rangle_V,\\
     I_{n, 3}(t) &= \Vert \zeta_n(t,*) - \zeta(t,*) \Vert_\bD.
     \end{align*}
Clearly,
\beqn
\label{t01}
     I_{n,0}(t) \leq \sup_{x \in [0, 1]} \Vert G(t, x, *) \Vert_{L^1([0, 1])}\, \Vert g_n - g \Vert_\bD,
          \eeqn
and this implies
\beqn
\label{tn0}
\lim_{n\to\infty}\, \sup_{t\in[0,T]} I_{n,0}(t) = 0.
\eeqn
Similarly,
\begin{align*}
      I_{n, 1}(t)  &\leq \sup_{x \in [0, 1]} \int_0^t  ds\, \Vert G(t-s, x, *) \Vert_V\,  \Vert F_n(\kappa(s, *)) - F(\kappa(s, *)) \Vert_{V}\\
      &
   \leq c \int_0^t  ds\, (t - s)^{-1/4}\, \Vert F_n(\kappa(s, *)) - F(\kappa(s, *)) \Vert_V.
   \end{align*}
 By hypothesis, the integrand in the last expression tends to zero as $n\to\infty$. Moreover, the (uniform) Lipschitz property of $(F_n)$ and $F$ yields
 \beqn
 \sup_{s\in[0,T]}   \Vert F_n(\kappa(s, *)) - F(\kappa(s, *)) \Vert_V \le C\left(1 + \Vert \kappa\Vert_{\cC([0,T]\times [0,1])}\right)
\eeqn
and by Lemma \ref{nbr-*3ab}, this is bounded by a finite constant. Thus, by dominated convergence
\beqn
\label{tn1}
\lim_{n\to\infty}\, \sup_{t\in[0,T]} I_{n,1}(t) = 0.
\eeqn

By hypothesis \eqref{nbr-(*14abb)},
\begin{align*}
   I_{n, 2}(t) &\leq \int_0^t  ds\, \Vert G(t-s, x, *) \Vert_V \  \Vert F_n(\kappa_n(s, *)) - F_n(\kappa(s, *)) \Vert_V\\
   &
   \leq \tilde C \int_0^t  ds\, (t - s)^{-1/4}\, \Vert \kappa_n(s, *) - \kappa(s, *) \Vert_\bD.
   \end{align*}

   Let $\eta_n(t) = I_{n,0}(t) + I_{n,1}(t) + \Vert \zeta_n(t,*) - \zeta(t,*) \Vert_\bD$.
From the above discussion, we deduce that
\beqn
   \Vert \kappa_n(t, *) - \kappa(t, *) \Vert_\bD \leq \eta_n(t) + \tilde C\int_0^t  ds\, (t - s)^{-1/4}\, \Vert \kappa_n(s, *) - \kappa(s, *) \Vert_\bD,
   \eeqn
   with $\lim_{n\to\infty} \sup_{t\in[0,T]}\eta_n(t) =0$.
Applying Lemma \ref{A3-l01}, we conclude the proof.
\end{proof}
\medskip

Proposition \ref{nbr-*3aa} concerns solutions to the deterministic integral equation \eqref{nbr-(*15b)}, but similar arguments apply to ODEs and to PDEs. We give a direct proof for systems of ODEs driven by a linear noise, such as a multidimensional Brownian motion (\eqref{eqforx} for instance), since this property is used in the proof of Proposition \ref{nbr-1}.


 \begin{lemma}
 \label{aproxode}
 Let $A$ be an $n\times n$ matrix of real numbers, $F: \re^n \to \re^n$ be a Lipschitz continuous function and $\xi: \re_+\to \re^n$ be a continuous function. Consider the solution $(x_t)$ of the system of ODEs
 \beq
 \label{ode1}
 \frac{dx_t}{dt} = A x_t + F(x_t)  + \xi_t,\qquad t>0,
 \eeq
 with initial condition $x_0 \in \re$.
 Let $F^{(\ell)}: \re^n \to \re^n$, $\ell\ge 1$, be a sequence of uniformly Lipschitz continuous functions, that is, there exists $K<\infty$ such that for any $x,y\in\re^n$ and $\ell\ge 1$,
$\left\vert F^{(\ell)}(x)-F^{(\ell)}(y)\right\vert\le K|x-y|$. Let $(x_t^{(\ell)})$ be the solution of
\beq
 \label{ode2}
 \frac{dx_t^{(\ell)}}{dt} = A x_t^{(\ell)} + F^{(\ell)}(x_t^{(\ell)})  + \xi_t,\qquad t>0,
\eeq
 with the same initial condition $x_0$. Further, assume that
for any $x\in\re^n$, $F^{(\ell)}(x)\to F(x)$ as $\ell\to\infty$. Then
 \beqn
  \lim_{\ell\to\infty}\left\vert x_t^{(\ell)}- x_t\right\vert =0, \quad {\text{uniformly in}}\ t\in[0,T].
  \eeqn
  \end{lemma}
  
  \begin{proof}
  Integrating \eqref{ode1} and \eqref{ode2} from $0$ to $t$, and using the triangle inequality and Jensen's inequality, we see that
  \begin{align}
  \label{edo-eq1}
  \left\vert x_t^{(\ell)}- x_t\right\vert &\le C\left[\int_0^t \left\vert x_s^{(\ell)}- x_s\right\vert\, ds\right.\notag \\
  &\left.+ \int_0^t \left\vert F^{(\ell)}(x_s^{(\ell)})-F^{(\ell)}(x_s)\right\vert\, ds +  \int_0^t \left\vert F^{(\ell)}(x_s)-F(x_s)\right\vert\, ds\right].
  \end{align}
By assumption, for any fixed $s\in [0,T]$, $\lim_{\ell\to\infty} \left\vert F^{(\ell)}(x_s)-F(x_s)\right\vert =0$. Moreover, the Lipschitz continuity properties of $F^{(\ell)}$ and $F$ imply that
  \beqn
  \sup_{s\in[0,T]} \left\vert F^{(\ell)}(x_s)-F(x_s)\right\vert \le  C \left(1+\sup_{s\in[0,T]}|x_s|\right),
  \eeqn
  where $C$ does not depend on $\ell$, and the right-hand side of this inequality is finite.

  By dominated convergence and using again the uniform Lipschitz property of the $F^{(\ell)}$, we deduce from \eqref{edo-eq1} that
  \beqn
   \left\vert x_t^{(\ell)}- x_t\right\vert \le C_T \int_0^t \left\vert x_s^{(\ell)}- x_s\right\vert\, ds + \varepsilon_\ell,
   \eeqn
   where $\lim_{\ell\to\infty}\varepsilon_{\ell}=0$. The conclusion follows from Gronwall's lemma \ref{lemC.1.1} applied to $f(t):=\sup_{r\in[0,t]}\left\vert x_r^{(\ell)}- x_r\right\vert$, $z(t)\equiv \varepsilon_\ell$ and $J(t)\equiv 1$.
    \end{proof}



\section{Two properties of conditional expectations}
The next two statements give properties concerning conditional expectations.
The first one is used in the proof of Proposition \ref{nbr-2aa}.

\begin{lemma}
\label{ch5-added-l1}
Let $X,Y$ be two random variables defined on a probability space $(\Omega,\cF,P)$. Let $\mathcal{G}$ be a sub $\sigma$-field of $\cF$. Assume that $X$ is independent of $\mathcal{G}$ and $Y$ is $\mathcal{G}$-measurable. Let $h$ be a Borel function defined on $\re^2$ such that $h(X,Y)\in L^1(\Omega)$. Then
\beq
\label{ch5-added-1-20}
E[h(X,Y)\, | \,\mathcal{G}] = E\left[h(X,y)\right]\big\vert_{y=Y},\quad \text{ a.s.}
\eeq
\end{lemma}

\begin{proof}
Take first $h(x,y)= 1_A(x)1_B(y)$, with $A, B\in\cB_{\re}$. By elementary properties of conditional expectation, we
have
\begin{align*}
E[h(X,Y)|\mathcal{G}] &=E[1_A(X)1_B(Y)|\mathcal{G}] =1_B(Y)\, E[1_A(X)|\mathcal{G}] \\
&=1_B(Y)\, E[1_A(X)].
\end{align*}
Moreover,
\beqn
E[h(X,y)] = E[1_A(X)1_B(y)] =1_B(y)\, E[1_A(X)].
\eeqn
Therefore \eqref{ch5-added-1-20} holds in this case.
Since \eqref{ch5-added-1-20} is linear in $h$, the general case follows by taking linear combinations of functions of the form just considered and using the monotone class theorem (see e.g. \cite{dm1}).
\end{proof}

The next lemma is a property of conditional probabilities under two different measures. It is used in Theorem \ref{nbr-t-14}.

\begin{lemma}
\label{nbr-l-11}
Let $(\Omega, \cF, P)$ be a probability space and $Q$ a probability measure on $(\Omega, \cF)$. Let $(S, \cs)$ and $(U, \U)$ be two measure spaces, and $X$ (respectively $Y$) be a random variable with values in $S$ (respectively $U$). Assume that the law on $(S \times U, \cs \times \U)$ of $(X, Y)$ under $P$ is absolutely continuous with respect to the law on $(S \times U, \cs \times \U)$ of $(X, Y)$ under $Q$, with density $\varphi(x, y)$. Suppose also that there is a regular conditional probability distribution $Q_{X \mid Y = y}$ for $X$ given $Y = y$ (under $Q$), denoted $(y, A) \mapsto Q_y(A)$. Then there is a regular conditional probability distribution $P_{X \mid Y = y}$ for $X$ given $Y = y$ (under $P$), denoted $(y, A) \mapsto P_y(A)$, given by
\beq
\label{nbr-(*21)}
      P_y(A) =
\begin{cases}
      \frac{1}{c(y)} \int_S  1_A(x) \varphi(x, y)\, Q_y(dx),&     {\text {if}}\  c(y) > 0,\\
                      1,&        {\text {otherwise}},
 \end{cases}
   \eeq
where
\beq
\label{nbr-(*22)}
     c(y) = \int_S \varphi(x, y)\, Q_y(dx).    
\eeq
\end{lemma}

\begin{proof}
Let $(\Omega_1, \cF_1, P_1)$ be a probability space, $P_2$ be a probability measure on $(\Omega_1, \cF_1)$, $\cG$ be a sub-$\sigma$-field of $\cF$ and $X$ be a random variable. Suppose that $P_1$ is absolutely continuous with respect to $P_2$. We will use the generic formula
\beq
\label{AppC-(*E1)}
     E_{P_1}[X \mid \cG] = \frac{E_{P_2}\left[X \frac{d P_1}{d P_2} \mid \cG\right]}{E_{P_2}\left[\frac{d P_1}{d P_2} \mid \cG\right]},    
     \eeq
(see \eqref{ch2'-s3.4.4.6}) applied to the laws of $(X, Y)$ under $P$ and under $Q$. Notice that the denominator in \eqref{AppC-(*E1)} is positive $P_1$-a.s.

     Let $f: S \to \R$ and $g: S \times U \to \R$ be bounded and measurable. By definition,
 \beqn
      E_P[g(X, Y)] = E_Q[g(X, Y)\, \varphi(X, Y)]
\eeqn
and
\beq\label{AppC-(*E2)}
    E_Q[g(X, Y) \mid Y] = g_Q(Y),   
\eeq
where $g_Q(y) = \int_S g(x, y)\, Q_y(dx)$.

By \eqref{AppC-(*E1)},
\beqn
     E_P[f(X) \mid Y] = \frac{E_Q[ f(X)\, \varphi(X, Y) \mid Y]}{E_Q[\varphi(X, Y) \mid Y]}.
     \eeqn
By \eqref{AppC-(*E2)}, the numerator is equal to $(f \varphi)_Q(Y)$, where
\beqn
     (f \varphi)_Q(y) = \int_S f(x) \varphi(x, y)\, Q_y(dx),
     \eeqn
and the denominator is equal to $\varphi_Q(Y)$, where
\beqn
     \varphi_Q(y) = \int_S \varphi(x, y)\, Q_y(dx).
     \eeqn
Setting $c(y) = 1$ on the set $G = \{y \in U: \varphi_Q(y) = 0 \}$, which satisfies $P\{Y \in G\} = 0$, we obtain \eqref{nbr-(*21)}.
\end{proof}



\section{Facts concerning the Gaussian law}
\label{app3-s3}
In this section, we recall some classical properties of the Gaussian law and Gaussian random variables that are used throughout this book.
\medskip

\noindent {\em Moments of Gaussian random variables}
\smallskip

Recall that the Euler Gamma function\index{Euler!Gamma function}\index{Gamma function!Euler}\index{function!Euler Gamma} is defined for any $p\in\, ]0,\infty[$ by
\beq
\label{Euler-gamma}
\Gamma_E(p) = \int_0^\infty e^{-x}\ x^{p-1}\, dx.
\eeq
For properties of this function, including the basic fact that for $n \in \N^*$, $\Gamma_E(n+1) = n!$, see \cite{olver}. 

\begin{lemma}
\label{A3-l1}
Fix $p\in\, ]-1,\infty[$ and $\sigma\in\, ]0,\infty[$.
\begin{enumerate}
\item Let $f_{\sigma^2}(x) = \frac{1}{\sqrt{2\pi\sigma^2}}\ e^{-\frac{x^2}{2\sigma^2}}$, $x\in \re$,
 be the probability density function of a ${\rm N}(0,\sigma^2)$ random variable $Z$. Then
\beq
\label{A3.4.0}
\int_{\re} |x|^p\, f_{\sigma^2}(x)\, dx = c_p\, \sigma^p,
\eeq
where $c_p=\left(\frac{2^p}{\pi}\right)^{\half} \Gamma_E\left(\frac{p+1}{2}\right)$.
Equivalently,
\beq
\label{A3.4.0-bis}
E[|Z|^p] = c_p \left(E[Z^2]\right)^{\frac{p}{2}}.
\eeq

\item Let $\Gamma(t,x) = \frac{1}{\sqrt{4\pi t}} \exp\left(-\frac{x^2}{4t}\right)\ 1_{]0,\infty[}(t)$ be the fundamental solution to the heat equation on $\re$ $($see \eqref{heatcauchy-1'}$)$. Then
\beq
\label{A3.4}
\int_{\re} |x|^p\, \Gamma(t,x)\, dx = \frac{2^p}{\sqrt \pi} \, \Gamma_E\left(\frac{p+1}{2}\right)\, t^{\frac{p}{2}}.
\eeq
\item Let $Y$ be a ${\rm N}(\mu,\sigma^2)$ random variable ($\mu\in \IR$). For $p>0$,
\beq
\label{A3.4.0-tris}
E[|Y|^p] \le 2^{p}(1+c_p)\left(E[Y^2]\right)^{\frac{p}{2}},
\eeq
with $c_p$ defined in part 1.
\end{enumerate}
\end{lemma}
\begin{proof}
For the proof of 1., we observe that since the integrand on the left-hand side of \eqref{A3.4.0} is an even function, we have
\beqn
\int_{\re} |x|^p\ f_{\sigma^2}(x)\, dx = \left(\frac{2}{\pi\sigma^2}\right)^{\half} \int_0^\infty x^p\, \exp\left(-\frac{x^2}{2\sigma^2}\right) dx.
\eeqn
Apply the change of variable $w:=\frac{x^2}{2\sigma^2}$ to see that
\begin{align*}
\int_0^\infty x^p \exp\left(-\frac{x^2}{2\sigma^2}\right) dx &= 2^{\frac{p-1}{2}}\, \sigma^{p+1}\int_0^\infty e^{-w}\, w^{\frac{p-1}{2}}\, dw\\
& = 2^{\frac{p-1}{2}}\, \sigma^{p+1}\, \Gamma_E\left(\frac{p+1}{2}\right),
\end{align*}
proving \eqref{A3.4.0}.

Formula \eqref{A3.4} is a particular case of \eqref{A3.4.0}, since $\Gamma(t,x) = f_{2t}(x)$.
For the proof of \eqref{A3.4.0-tris}, we write $Y=\mu + Z$, where $Z$ is a ${\rm N}(0,\sigma^2)$ random variable. From \eqref{A3.4.0-bis}, we deduce that for $p>0$,
 \beq
 \label{A3.4.0-quator}
 E[|Y|^p] \le 2^{p}\left(|\mu|^p + E[|Z|^p]\right) = 2^{p}\left(|\mu|^p +c_p \left(E[Z^2]\right)^{\frac{p}{2}}\right).
 \eeq
 Observe that $E[Y^2] = \mu^2 + E[Z^2]$ and therefore $|\mu|^p$ and $(E[Z^2])^{p/2}$ are bounded by
 $(E[Y^2])^{p/2}$. Applying these facts, we see that the right-hand side of \eqref{A3.4.0-quator} is bounded by
 $2^{p}(1 + c_p)(E[Y^2])^{p/2}$, proving \eqref{A3.4.0-tris}.
 \end{proof}

 \noindent{\em Tail estimates}
 \smallskip

 \begin{lemma}
 \label{rdlemC.2.2-before}
 For any  $a>0$, we have the following:
\begin{equation}
\label{rdlemC.2.2-before-1}
 \int_{[-a,a]^c} f_{\sigma^2}(x)\, dx \le 
      \min\left(1, \sqrt{\frac{2}{\pi}}\frac{\sigma}{a}\right) e^{-\frac{a^2}{2\sigma^2}}
\end{equation}
and
\begin{equation}
 \int_{[-a,a]^c}  f_{\sigma^2}^2(x)\, dx \le 
      \frac{1}{2 \sigma\, \sqrt{\pi}} \, e^{-\frac{a^2}{\sigma^2}}.
\label{rdlemC.2.2-before-2}
\end{equation}
\end{lemma}

\begin{proof}
 For any $\gamma>0$, we have
\begin{align}\nonumber
 \int_a^\infty e^{-\gamma x^2}\, dx &= \int_0^\infty e^{-\gamma (z+a)^2}\, dz = e^{-\gamma a^2} \int_0^\infty e^{-2\gamma a z}\, e^{-\gamma z^2}\, dz \\
    & \leq e^{-\gamma a^2} \int_0^\infty e^{-\gamma z^2}\, dz = \half\, \sqrt{\frac{\pi}{\gamma}} \, e^{-\gamma a^2}.
\label{rdlemC.2.2-before-3}
 \end{align}
 Since $f_{\sigma^2}(x)$ is an even function, for any $\eta>0$,
 \begin{align*}
 \int_{[-a,a]^c} f_{\sigma^2}^\eta(x)\, dx = 2 \int_a^\infty f_{\sigma^2}^\eta(x)\, dx 
 = \frac{2}{(2\pi\sigma^2)^{\frac{\eta}{2}}} \int_a^\infty e^{-\frac{\eta}{2\sigma^2}x^2}\, dx.
 \end{align*}
 Taking $\eta = 1$ (respectively, $\eta=2$) and applying \eqref{rdlemC.2.2-before-3} with $\gamma = \frac{1}{2\sigma^2}$ (respectively, $\gamma = \frac{1}{\sigma^2}$),
 we obtain 
 \beq
 \label{part-of 1}
  \int_{[-a,a]^c} f_{\sigma^2}(x)\, dx \le e^{-\frac{a^2}{2\sigma^2}}
 \eeq
 (respectively, \eqref{rdlemC.2.2-before-2}).
 
  To complete the proof of \eqref{rdlemC.2.2-before-1}, we take $\sigma =1$ and write $f(x)$ instead of $f_{1}(x)$. Then
 \begin{align*}
 \int_{[-a,a]^c} f(x)\, dx &=\frac{2}{\sqrt{2\pi}} \int_a^\infty
 dx\, e^{-\frac{x^2}{2}} \le \frac{2}{\sqrt{2\pi}} \int_a^\infty dx\, \frac{x}{a}\, e^{-\frac{x^2}{2}}\\
 &= \sqrt{\frac{2}{\pi}}\frac{1}{a} e^{-\frac{a^2}{2}}.
 \end{align*}
 Since 
 \beqn
  \int_{[-a,a]^c} f_{\sigma^2}(x)\, dx = \int_{[-\frac{a}{\sigma},\frac{a}{\sigma}]^c} f(x)\, dx, 
 \eeqn
 we have
 \beq
 \label{part-of 2}
  \int_{[-a,a]^c} f_{\sigma^2}(x)\, dx \le \sqrt{\frac{2}{\pi}}\frac{\sigma}{a} e^{-\frac{a^2}{2\sigma^2}}.
 \eeq
 Together with \eqref{part-of 1}, this yields \eqref{rdlemC.2.2-before-1}.
 \end{proof}
 \medskip

\noindent{\em Products of Normal density functions}
 \medskip

For $\nu\in\re_+^*$, let $\Gamma_\nu(t,x)$ be defined as in \eqref{fsnu}. The following properties are used in the proof of Proposition \ref{rough-heat}.

\begin{lemma}
\label{rough-heat-l1}
\begin{description}
\item{(i)} For all $s, t\in \re_+^*$ and $x, y\in \re$,
\beq
\label{rough-heat-l1-1}
\Gamma_\nu(t,x)\, \Gamma_\nu(s,y) = \Gamma_\nu\left(\frac{ts}{t+s}, \frac{sx+ty}{t+s}\right)  \Gamma_\nu(t+s, x-y).
\eeq 
\item{(ii)} Fix $a \geq 2$. For all $s, t \in\re_+^*$ and $x, y_1, y_2\in \re$, let $\bar y=\frac{y_1+y_2}{2}$, $\underline{y} = y_1-y_2$. Then
\begin{align}\nonumber
 &\Gamma_1(t,x-\bar y)\, \Gamma_1(s,\underline{y}) \\
  &\qquad\qquad \le \frac{(at)\vee s}{\sqrt{ts}}\, \Gamma_1((a t)\vee s, x-y_1)\, \Gamma_1((a t)\vee s, x-y_2).
\label{rough-heat-l1-20}
\end{align}
\end{description}
\end{lemma}

\begin{proof} 
(i) If we write $z = x-y$, then  identity \eqref{rough-heat-l1-1} can be written
\beq\label{rd04_24e1}
   \Gamma_\nu(t,x)\, \Gamma_\nu(s,z-z) = \Gamma_\nu\left(\frac{ts}{t+s}, x - \frac{t}{t+s}\, z\right)  \Gamma_\nu(t+s, z).
\eeq
This identity expresses the fact that the joint density of $(B_{2t},B_{2(t+s)})$, where $(B_r)$ is a standard Brownian motion, evaluated at $(x, z)$, is equal to the product of the marginal density at time $2(t+s)$, evaluated at $z$, and the conditional density of  $B_{2t}$ given $B_{2(t+s)} = z$, evaluated at $x$. This conditional density is the law at time $2t$ of a Brownian bridge that goes from $0$ to $z$ during the time interval $[0, 2(t+s)]$, and this conditional density is that of a 
\beqn
    {\rm N}\left(\frac{t}{t+s}\, z, \, \frac{ts}{t+s}\right)
\eeqn
random variable.  Checking the identity \eqref{rd04_24e1} by direct calculation is a simple exercise. 

  (ii) Observe that
\beqn
(x-\bar y)^2 = \frac{[(x-y_1)+(x-y_2)]^2}{4}
\eeqn
and
\beqn
   (y_1-y_2)^2 = [(x - y_1) - (x - y_2)]^2.
\eeqn
Therefore, 
\beq\label{c.6.3-2a}
     \left[(x-y_1)+(x-y_2)\right]^2 + (y_2-y_1)^2  = 2\left[ (x-y_1)^2+(x-y_2)^2 \right].
\eeq

Notice that
\begin{align*}
   &\Gamma_1\left(t,x-\bar{y}\right)  \Gamma_1\left(s,\underline{y}\right) \\
   &\qquad\qquad =  
   \frac{1}{2\pi\sqrt{ts}} \exp\left(-\half \left(\frac{[(x-y_1)+(x-y_2)]^2}{4t}+\frac{(y_1-y_2)^2}{s}\right)\right).
\end{align*}
The sum of the two fractions inside the exponential can be written in two ways:
\begin{align}\nonumber
&\frac{[(x-y_1)+(x-y_2)]^2 + (y_1-y_2)^2}{4t} + \left(\frac{1}{s}-\frac{1}{4 t}\right)(y_1-y_2)^2\\
 &\qquad\qquad\  = \frac{(x-y_1)^2+(x-y_2)^2}{2t}+ \left(\frac{1}{s}-\frac{1}{4 t}\right)(y_1-y_2)^2,
 \label{rd04_24e2}
\end{align}
where we have used \eqref{c.6.3-2a}, and
\begin{align}\nonumber
&\frac{[(x-y_1)+(x-y_2)]^2 + (y_1-y_2)^2}{s} + \left(\frac{1}{4t}-\frac{1}{s}\right)[x-y_1)+(x-y_2)]^2\\  
&\qquad\qquad = \frac{(x-y_1)^2+(x-y_2)^2}{s} + \left(\frac{1}{4t}-\frac{1}{s}\right)[x-y_1)+(x-y_2)]^2 
 \label{rd04_24e3}
\end{align}
If $s \leq 4t$ (respectively $s > 4t$), then we use \eqref{rd04_24e2} (respectively \eqref{rd04_24e3}) to see that
\begin{align*}
    \frac{[(x-y_1)+(x-y_2)]^2}{4t}+\frac{(y_1-y_2)^2}{s} &\geq \frac{(x-y_1)^2+(x-y_2)^2}{((2t) \vee s)} \\
     &\geq \frac{(x-y_1)^2+(x-y_2)^2}{((at) \vee s)},
\end{align*}
because $a \geq 2$, and this implies \eqref{rough-heat-l1-20}.
\end{proof}
\medskip

 \noindent{\em Derivatives of the heat kernel}\index{kernel!heat}\index{heat!kernel}
 \medskip

 Upper and lower bounds on derivatives of the fundamental solutions and Green's functions of parabolic Cauchy and boundary value problems are crucial in the classical theory of PDEs. They can be found for instance in \cite[Chapter 9, Theorem 8]{friedman64} and \cite{ez}. We recall here an upper bound in the particular case of the fundamental solution to the heat equation on $\R^k$.

 \begin{lemma}\label{rdlemC.2.2} Let
 \beqn
 \Gamma(t,x;s,y) = \frac{1}{(4\pi(t-s))^{\frac{k}{2}}}\exp\left(- \frac{|x-y|^2}{4(t-s)}\right),\quad 0\le s<t\le T, \ x\in\re^k.
 \eeqn
 Then, for any  $n_1,n_2,n_3\in\mathbb{N}$, there exist positive constants $C = C(T,n_1,n_2,n_3)$, $c = c(T,n_1,n_2,n_3)$ such that,
 \beq
 \label{rdlemC.2.2-before-5}
 \left\vert \frac{\partial ^{n_1+n_2+n_3}\Gamma(t,x;s,y) }{\partial t^{n_1}\partial x^{n_2}\partial y^{n_3}}\right\vert \le C (t-s)^{-(k+2n_1+n_2+n_3)/2} \exp\left(-c\ \frac{|x-y|^2}{t-s}\right).
 \eeq
 \end{lemma}
 \medskip


\section{A property of the Euler Gamma function}
\label{app3-s4}
 \medskip

The next lemma gives the proof of \eqref{f1} in Section \ref{ch1'-ss7.3}.
\begin{lemma}
\label{A3-l0}
Let $a\in\, ]1,2[$ and $q\in\mathbb{C}$ with $\mathrm{Re}(q)\ge 0$. Then
\beq
\label{A3.0}
\int_0^\infty\frac{e^{-qz}-1+qz}{z^{1+a}}\, dz = q^a\, \Gamma_E(-a),
\eeq
where $\Gamma_E$ is the Euler Gamma function defined in \eqref{Euler-gamma}.
\end{lemma}
\begin{proof}
First, we observe that
\beq
\label{A3.1}
\int_0^\infty\frac{\left\vert e^{-qz}-1+qz\right\vert}{z^{1+a}}\, dz <\infty,
\eeq
because, for $z\to\infty$, the numerator is bounded by $Cz$, and for $z\downarrow 0$, the numerator is bounded by $cz^2$.

Clearly
\beq
\label{A3.2}
\int_0^\infty\frac{e^{-qz}-1+qz}{z^{1+a}}\ dz = \int_0^\infty dz\ z^{-1-a} \int_0^z dy\, (z-y)\, q^2\, e^{-qy}.
\eeq
We want to apply Fubini's theorem. However, if $\mathrm{Re}(q) = 0$, $q\ne 0$, then $|e^{-qy}|= 1$ and
\beqn
\int_0^\infty dz\,  z^{-1-a}  \int_0^z dy\, (z-y)\,|q|^2\, |e^{-qy}| = \frac{|q|^2 }{2} \int_0^\infty dz\,  z^{1-a} = \infty,
\eeqn
while if $q_1 = \mathrm{Re}(q) >0$, then
\begin{align*}
&\int_0^\infty dz\, z^{-1-a}  \int_0^z dy\, (z-y)\,|q|^2\, e^{-q_1 y}\\
&\qquad\qquad = \frac{|q|^2}{q_1^2} \int_0^\infty dz\, z^{-1-a} \int_0^z dy\, (z-y)\, q_1^2\, e^{-q_1 y}.
\end{align*}

Applying \eqref{A3.2} with $q_1$ instead of $q$, this is equal to
\beqn
\frac{|q|^2}{q_1^2}\int_0^\infty dz\, \frac{e^{-q_1 z} -1+q_1z}{ z^{1+a}} < \infty,
\eeqn
by \eqref{A3.1} above (since the numerator is positive).

Assume for the moment that $\mathrm{Re}(q) >0$. Applying Fubini's theorem in \eqref{A3.2}, we see that
\begin{align*}
\int_0^\infty dz\,  \frac{e^{-q z} -1+qz}{z^{1+a}}& = \int_0^\infty dy \int_y^\infty
 dz\ \left(z^{-a}-yz^{-1-a}\right) q^2\, e^{-qy}\\
 & = \int_0^\infty dy\, \left(\frac{-y^{1-a}}{1-a} - \frac{y^{1-a}}{a}\right)q^2\, e^{-qy}\\
  & = \frac{q^2}{a(a-1)}\int_0^\infty e^{-qy}\, y^{1-a}\ dy.
  \end{align*}
  Use a change of variables $x=qy$ and \eqref{Euler-gamma} to see that this is equal to
\beqn
  \frac{q^2}{a(a-1)}\, \Gamma_E(2-a)\, \frac{1}{q^{2-a}} = q^a\, \Gamma_E(-a).
\eeqn
 We conclude that
  \beq
  \label{A3.3}
  \int_0^\infty\frac{e^{-qz}-1+qz}{z^{1+a}}\, dz  = q^a\, \Gamma(-a), \qquad  \mathrm{Re}(q) >0.
 \eeq
 In order to extend this identity to the case $\mathrm{Re}(q) = 0$, we fix $q\in\mathbb{C}$ with $\mathrm{Re}(q) = 0$: $q=i q_2$, $q_2\in\re$. We let $q_n=\frac{1}{n}+iq_2$, so that
 \eqref{A3.3} holds for $q_n$. Let
 \beqn
 f(z) = \begin{cases}
 (1+q_2^2)z^{1-a}, &{\text{if}}\ 0<z<1,\\
 (3+|q_2|z) z^{-1-a}, &{\text{if}}\ z\ge 1.
 \end{cases}
 \eeqn
 Then for all $n\in \mathbb{N}^\ast$, $z^{-1-a}\left\vert e^{-q_nz}-1+q_n z\right\vert \le f(z)$ and $\int_0^\infty f(z) \, dz < \infty$, since $a\in\, ]1,2[$.  Passing to the limit $n\to \infty$ in \eqref{A3.3} for $q_n$ and using the dominated convergence theorem, we obtain \eqref{A3.0} for $q$.
\end{proof}
\medskip


\noindent{\em Beta integrals and functions} 
\medskip

In Section \ref{rdrough}, we use the following property.
For $a, b \in \bC$ with ${\mathrm Re}(a) > 0$, and ${\mathrm Re}(b) > 0$,
\beq
\label{rd08_12e9}
   \int_s^t (r-s)^{a-1}\, (t - r)^{b-1} \, dr = B_E(a, b)\,  (t-s)^{a + b - 1},
\eeq
where
\beq\label{rd08_12e9a}
   B_E(a, b) = \frac{\Gamma_E(a)\, \Gamma_E(b)}{\Gamma_E(a+b)}\, 
\eeq 
is the {\em Euler Beta function.}\index{Euler!Beta function}\index{function!Euler Beta}\index{Beta function!Euler}

\section{Explicit formulas for some trigonometric series}
\label{rd08_26s2}
The next lemma is used several times in Sections \ref{stwn-s2} and \ref{ch5-added-1-s2} and in Appendix \ref{app2}. For instance, the identities \eqref{rdeB.5.5} and \eqref{rdeB.5.5} are used in Section \ref{stwn-s2}, \eqref{rdeB.5.5}--\eqref{rdeB.5.6-quatre} are used in the proof of Proposition \ref{rem6.1.15-*1}, the identity \eqref{rdeB.5.4} is applied in the proofs of Lemmas \ref{ch1'-l2}, \ref{ch1'-1202} and \ref{app2-5-l-w1}; the equality \eqref{rdeB.5.5} is used in Lemmas \ref{ch1'-l3} and \ref{ch1-l4}, and \eqref{rdeB.5.6} is used in the proof of Lemma \ref{app2-5-l-w1}.

\begin{lemma}
\label{lemB.5.2}
 Let $x,y \in [0,1]$. Then
\beq\label{rdeB.5.4}
   \sum_{n=1}^\infty\, \frac{[\sin(n\pi x)-\sin(n\pi y)]^2}{\pi^2\, n^2} = \frac{1}{2} \left(\vert x - y\vert - \vert x-y\vert^2 \right),
\eeq
\beq\label{rdeB.5.5}
   \sum_{n=1}^\infty\, \frac{[\cos(n\pi x)-\cos(n\pi y)]^2}{\pi^2\, n^2} = \frac{1}{2} \vert x - y\vert 
\eeq
and
\beq\label{rdeB.5.6}
\sum_{n=1}^\infty\, \frac{\sin^2(n\pi x)}{\pi^2\, n^2} = \frac{1}{2}\, x (1-x).
\eeq
Consequently,
\beq
\label{rdeB.5.6-bis}
\sum_{n=1}^\infty\,
\frac{\cos^2(n\pi x)}{\pi^2\, n^2}= \frac{1}{6} - \frac{x}{2}(1-x),
\eeq
\beq
\label{rdeB.5.6-tris}
\sum_{n=1}^\infty\, \frac{1-\cos(n\pi x)}{\pi^2n^2}= \frac{x}{2}\left(1-\frac{x}{2}\right),
\eeq
\beq
\label{rdeB.5.6-quatre}
\sum_{n=1}^\infty\frac{(1-\cos(n\pi x))(1-\cos(n\pi y))}{\pi^2n^2}= \half(x\wedge y)
\eeq
and
\beq
\label{rdeB.5.6-quatre-bis}
\sum_{n=1}^\infty \frac{\sin(n\pi x)\, \sin(n\pi y)}{\pi^2 n^2} = \frac{1}{2}(x\wedge y)(1-x\vee y).
\eeq
\end{lemma}

\begin{proof} We first recall that the sequence of functions
\beqn
 \varphi_n(z)= \tfrac{1}{\sqrt 2}e^{in\pi z},\quad z\in [-1,1],\ n\in\mathbb{Z},
 \eeqn
 is an orthonormal basis of $L^2([-1,1],\mathbb{C})$.

Assume without loss of generality that $0\leq x < y \leq 1$. For \eqref{rdeB.5.4}, consider the function $f\in L^2([-1,1],\mathbb{C})$ defined by $f(z)= 1_{[x,y]}(z) + 1_{[-y,-x]}(z)$. Using its Fourier expansion in terms of the orthonormal basis $(\varphi_n(z),\, n\in \mathbb{Z})$ and because $f$ is even, we see that
\begin{align*}
 2(y-x) &= \Vert f\Vert_{L^2([-1,1],\mathbb{C})}^2 = |\langle f,\varphi_0\rangle|^2 + 2 \sum_{n=1}^\infty\, |\langle f,\varphi_n\rangle|^2\\
 & =2(y-x)^2 + 4\sum_{n=1}^\infty\, \frac{[\sin(n\pi y)-\sin(n\pi x)]^2}{n^2\pi^2}.
\end{align*}
This yields \eqref{rdeB.5.4}. Formula \eqref{rdeB.5.6} is \eqref{rdeB.5.4} with $y=0$.

For \eqref{rdeB.5.5}, consider the odd function $f\in L^2([-1,1],\mathbb{C})$ defined by $f(z)= 1_{[x,y]}(z) - 1_{[-y,-x]}(z)$. Using its Fourier expansion in terms of the orthonormal basis $(\varphi_n(z),\, n\in \mathbb{Z})$, we see that
\begin{align*}
\Vert f\Vert_{L^2([-1,1],\mathbb{C})}^2 & = 2|x-y|
= 2\sum_{n=1}^\infty\, |\langle f,\varphi_n\rangle|^2\\
& = 4\sum_{n=1}^\infty\, \frac{[\cos(n\pi y)-\cos(n\pi x)]^2}{n^2 \pi^2}.
\end{align*}

Since $\cos^2 z = 1-\sin^2z$, formula \eqref{rdeB.5.6-bis} is a consequence of \eqref{rdeB.5.6} and the identity $\sum_{n=1}^\infty n^{-2} = \pi^2/6$.

For the proof of \eqref{rdeB.5.6-tris}, we use the formula $2(1-\cos z) = (1-\cos z)^2 + \sin^2z$ and apply \eqref{rdeB.5.5} and \eqref{rdeB.5.6}.

For the proof of \eqref{rdeB.5.6-quatre}, we use the identity
\beqn
2(1-a)(1-b) = (1-a)^2 + (1-b)^2 -(a-b)^2,\quad a,b\in\re,
\eeqn
with $a:= \cos (n\pi x)$ and $b:= \cos(n\pi y)$, to write the series on the left-hand side of \eqref{rdeB.5.6-quatre} as $\half$ times
\beqn
  \sum_{n=1}^\infty\, \frac{[1-\cos(n\pi x)]^2}{\pi^2\, n^2} +   \sum_{n=1}^\infty\, \frac{[1-\cos(n\pi y)]^2}{\pi^2\, n^2}
  -   \sum_{n=1}^\infty\, \frac{[\cos(n\pi x)-\cos(n\pi y)]^2}{\pi^2\, n^2}
  \eeqn
which, by \eqref{rdeB.5.5}, is equal to $\frac{1}{4}(x+y-|x-y|)$. This yields  \eqref{rdeB.5.6-quatre}.

Finally, for \eqref{rdeB.5.6-quatre-bis}, we develop the square in \eqref{rdeB.5.4} and then use \eqref{rdeB.5.6}.
\end{proof}

\begin{remark}
\label{app3-rem-lemB5.2}
In the proof of \eqref{app2.31} in Lemma \ref{app2-5-l-w1}, the following extension of \eqref{rdeB.5.6} has been used:

For all $s\geq 0$,
\beq
\label{rem-lemB.5.2}
          \sum_{n=1}^\infty\, \frac{\sin^2(n \pi s)}{n^2 \pi^2} = \half (s - [s]) (1 - (s - [s])) \leq \half s,    \eeq
where $[s]$ denotes the integer part of $s$.

Indeed, the left-hand side is a periodic function of $s$ with period $1$, so we can replace $s$ by $s - [s]$ without changing its value. Since $s - [s] \in [0,1]$, we apply \eqref{rdeB.5.6}to get the equality in \eqref{rem-lemB.5.2}, and the inequality is trivial.
\end{remark}
\smallskip

 \section{Inherited regularity of periodic extensions}
 \label{app3-s2}

In this section, we establish two properties that are used in Sections \ref{ch4-ss2.2-1'} and \ref{ch4-ss2.3-1'}, respectively.

 \begin{lemma}
 \label{A3-l2}
 For $\eta\in\, ]0,1]$ fixed, we consider the spaces of Hölder continuous functions $\mathcal{C}^\eta ([0,L])$ and $\mathcal{C}_0^\eta ([0,L])$ defined in Sections \ref{ch1'-ss3.3.1} and \ref{ch1'-ss3.4.3}, respectively.
 \smallskip

 1. Let  $v\in \mathcal{C}_0^\eta ([0,L])$ and consider the odd extension $v^o$ on $[-L,L]$ and the $2L$-periodic odd extension $v^{o,p}$ defined in \eqref{vop}. Then
 \beq
 \label{A3.5}
 \Vert v^{o,p}\Vert_{\mathcal{C}^\eta(\R)} \le 2 \Vert v\Vert_{\mathcal{C}_0^\eta([0,L])}.
 \eeq

 2. Let $v\in \mathcal{C}^\eta ([0,L])$ and consider the even extension $v^e$ on $[-L,L]$ and the $2L$-periodic even extension $v^{e,p}$ defined in \eqref{vep}. Then
 \beq
 \label{A3.5bis}
 \Vert v^{e,p}\Vert_{\mathcal{C}^\eta(\R)} =  \Vert v\Vert_{\mathcal{C}^\eta([0,L])}.
 \eeq
 \end{lemma}
 \begin{proof}
 1. Assume first that $\Vert v\Vert_{\mathcal{C}^\eta_0([0,L])}=1$. In this case, we show that for $x,y\in\, ]-L,L]$,
 \beq
 \label{A3.6}
 |v^o(x)- v^o(y)|\le 2 |x-y|^\eta.
 \eeq
 Indeed, if $x,y\in[0,L]$ or $x,y\in\, ]-L,0[$, then this inequality holds by hypothesis (with $2$ replaced by $1$). If $x\in\, ]-L,0[$ and $y\in[0,L]$, then
 \beq
 \label{A3.7}
 |v^o(x) - v^o(y)| = |-v(-x) - v(y)| = |v(-x) + v(y)|\le |v(-x)| + |v(y)|,
 \eeq
 and since $v(0)=0$, this is bounded above by $|x|^\eta + |y|^{\eta} \le 2|x-y|^\eta$. Therefore, \eqref{A3.6} is proved.
 From this, it follows that $\Vert v^{o}\Vert_{\mathcal{C}^\eta([-L,L])}\le 2$.

 Next, we consider the case $x, y \in\IR$ with $|x-y|\le 2L$. Without loss of generality, we assume that $x \le y$. Either both $x$ and $y$ are in the same interval of the form $[kL, (k+1)L]$, $k\in\mathbb{Z}$, or are in two adjacent such intervals, or $x \in [kL, (k+1)L[$ and $y \in [(k+2)L,(k+3)L[$.

 In the first two cases, we deduce from \eqref{A3.6} that
 $$
    \left\vert v^{o,p}(x)-v^{o,p}(y)\right\vert \le 2 |x-y|^\eta.
 $$
In the third case, by $2L$-periodicity, there is $z \in [kL, (k+1)L[$ such that $v^{o,p}(z) = v^{o,p}(y)$. Hence,
\begin{align*}
   \left\vert v^{o,p}(x)-v^{o,p}(y)\right\vert = \left\vert v^{o,p}(x)-v^{o,p}(z )\right\vert \le |x-z|^\eta \leq L^\eta
   \leq |x-y|^\eta .
\end{align*}

Therefore,
 \beq
 \label{A3.8}
 \sup_{\substack{x,y\in\IR,\, x\ne y\\ |x-y|\le 2L}} \frac{\left\vert v^{o,p}(x)-v^{o,p}(y)\right\vert}{|x-y|^\eta}\le 2.
 \eeq
 For $x, y\in\IR$ with $|x-y|>2L$, we have, by $2L$-periodicity and \eqref{A3.6},
 \begin{align*}
 \left\vert v^{o,p}(x)-v^{o,p}(y)\right\vert &\le \sup_{z_1, z_2\in [-L,L]} |v^o(z_1) - v^o(z_2)| \\ 
 &\le 2\ (2L)^\eta \le 2 |x-y|^\eta.
 \end{align*}
 We conclude that when $\Vert v\Vert_{\mathcal{C}^\eta_0([0,L])}=1$,
 \beq
 \label{A3.9}
 \Vert v^{o,p}\Vert_{\mathcal{C}^\eta(\R)} \le 2.
 \eeq
 Finally, let $v\in \mathcal{C}_0^\eta ([0,L])$ be arbitrary but $v\ne 0$, so  $\Vert v\Vert_{\mathcal{C}_0^\eta([0,L])}\ne 0$. Set $\bar v = v/\Vert v\Vert_{\mathcal{C}_0^\eta ([0,L])}$. Then,
 $\Vert\bar v\Vert_ {\mathcal{C}_0^\eta ([0,L])}=1$, so by \eqref{A3.9},  $\Vert \bar v^{o,p}\Vert_{\mathcal{C}^\eta(\R)} \le 2$. This is equivalent to \eqref{A3.5}.
 \smallskip

 2. The proof is similar to that of statement 1., with minor changes. Indeed, consider first the case where $\Vert v\Vert_{\mathcal{C}^\eta([0,L])} = 1$. For $x\in\, ]-L,0[$ and $y\in[0,L]$, the calculations in \eqref{A3.7} are replaced by
 \beqn
 |v^e(x) - v^e(y)| = |v(-x) - v(y)| \le |x+y|^\eta \le |x-y|^\eta,
 \eeqn
 so, for $x, y\in\, ]-L,L]$, instead of \eqref{A3.6}, we find that
 \beqn
 |v^e(x)- v^e(y)|\le |x-y|^\eta,
 \eeqn
 and, instead of \eqref{A3.8}, that
\beqn
   \sup_{\substack{x,y\in\IR,\, x\ne y\\ |x-y|\le 2L}} \frac{\left\vert v^{e,p}(x)-v^{e,p}(y)\right\vert}{|x-y|^\eta}\le 1.
 \eeqn
Thus, instead of \eqref{A3.9}, we get $\Vert v^{e,p}\Vert_{\mathcal{C}^\eta(\IR)}\le 1$ when $\Vert v\Vert_{\mathcal{C}^\eta([0,L])} =1$. This gives the conclusion
\beqn
   \Vert v^{e,p}\Vert_{\mathcal{C}^\eta(\R)} \le \Vert v\Vert_{\mathcal{C}^\eta([0,L])}.
\eeqn
Since the converse inequality is trivial, this yields \eqref{A3.5bis}.
 \end{proof}


\section{Some Laplace and Fourier transforms}
\label{app3LF}


We gather in this section various formulas concerning Laplace and Fourier transforms that are used in particular in Section \ref{rd08_26s1}. We denote the Laplace transform\index{Laplace!transform}\index{transform!Laplace} of a function $f: \R_+ \to \R$  by
\beqn
   F(z) = \cLL f(z) = \int_0^\infty e^{-zt}\, f(t)\, dt.
\eeqn
Given $f$, formulas for $F$ can often be found in standard tables such as \cite{EMOT,oberhettinger1973,oberhettinger1990}. These formulas sometimes make use of certain {\em special functions,} whose definitions we now recall.

\bigskip

\noindent{\em Some special functions}
\medskip

\noindent{\em Bessel function of the first kind\index{Bessel function!of the first kind}\index{function!Bessel, of the first kind} of order $\nu$} (see \cite[p.~414]{oberhettinger1973}):
\beq
\label{rd08_26L10a}
   J_\nu(z) = \left(\frac{z}{2}\right)^\nu \sum_{n=0}^\infty \frac{(-1)^n}{n!\, \Gamma_E(\nu+n+1)} \left(\frac{z}{2} \right)^{2n} .
\eeq

\noindent{\em  Modified Bessel function of the first kind\index{Bessel function!modified, of the first kind}\index{modified Bessel function!of the first kind}\index{function!modified Bessel, of the first kind} of order $\nu$} (see \cite[p.~414]{oberhettinger1973}):
\beq
\label{rd08_26L10}
    I_\nu(z) = e^{-i \pi \nu / 2} J_\nu(e^{i \pi / 2}z) = \left(\frac{z}{2}\right)^\nu \sum_{n=0}^\infty \frac{1}{n!\, \Gamma_E(\nu+n+1)} \left(\frac{z}{2} \right)^{2n} .
    \eeq

\noindent{\em Modified Bessel function of the second kind\index{Bessel function!modified, of the second kind}\index{modified Bessel function!of the second kind}\index{function!modified Bessel, of the second kind} of order $\nu$} (see \cite[p.~414]{oberhettinger1973}):
\beq
\label{rd01_23L1}
   K_\nu(x) = \frac{\pi}{2}\, \frac{I_{-\nu}(z) - I_{\nu}(z)}{\sin(\nu\, \pi)}\, .
\eeq

\noindent {\em Legendre functions\index{Legendre function}\index{function!Legendre}} (see \cite[p.~413]{oberhettinger1973}):
\begin{align}\nonumber
    P_\nu^\mu(z) &= [\Gamma_E(1-\mu)]^{-1}\, ((z + 1)/(z-1))^{\mu/2}\\
          &\qquad\qquad \times\  \null_2 F_1(-\nu, \nu+1; 1-\mu; (1-z)/2 ),
\label{rd02_05L1}
\end{align}
where $\null_2 F_1$ is a {\em generalized hypergeometric function} (see \cite[p.~423]{oberhettinger1973}); 
\bigskip

\noindent {\em Error function\index{error function}\index{function!error}} (see \cite[(p.~387]{EMOT}):
\beq\label{rd02_05e1}
   \text{Erfc}(t) = 2\, \pi^{-1/2} \int_{t}^\infty e^{-x^2}\, dx = \sqrt{\pi}\, (1 - \Phi(t \, \sqrt{2})),
\eeq
where $\Phi(\cdot)$ denotes the standard Normal probability distribution function.
\bigskip

\noindent{\em Some Laplace transforms}
\vskip 12pt

We now give various couples $(f, F)$, where $F$ is the Laplace transform of $f$:

\beq
\label{rd08_26L2}
     f(t) = 1, \qquad\qquad F(z) = z^{-1}
\eeq
(see \cite[2.2, p.~12]{oberhettinger1973});
\smallskip

\beq
\label{rd08_26L1}
    f(t) = t^{\nu-1},\quad \nu > 0, \qquad F(z) = \Gamma_E(\nu) \, z^{- \nu}
\eeq
(see \cite[3.3, p.~21]{oberhettinger1973});

\beq
\label{rd08_26L3}
   f(t) = t^{\nu - 1}\, e^{-\alpha t},\  \text{Re}(\nu) > 0,\qquad F(z) = \Gamma_E(\nu)\, (z + \alpha)^{-\nu}
\eeq
(see \cite[5.3, p.~37]{oberhettinger1973});
\smallskip

\beq
\label{rd08_26L4}
    f(t) = \sinh(\alpha t), \qquad\qquad F(z) = \alpha(z^2 - \alpha^2)^{-1}
\eeq
(see \cite[9.1 p.~84]{oberhettinger1973});
\smallskip

\beq
\label{rd08_26L5}
   f(t) = \cosh(\alpha t), \qquad\qquad F(z) = z (z^2 - \alpha^2)^{-1}
\eeq
(see \cite[9.2, p.~84]{oberhettinger1973});

\begin{align}\nonumber
     &f(t) = t^\mu K_\nu(a t), \ \mathrm{Re} > -1,\ a > 0, \\  \nonumber
   &   \qquad\qquad\quad  F(z) = \left(\frac{\pi}{2a} \right)^\half \Gamma_E(\mu+\nu+1) \Gamma_E(\mu+\nu-1) \\
      &\qquad\qquad\qquad\qquad\ \times (a^2 - z^2)^{-(\mu+1)/2} P_{\nu - \half}^{-\mu - \half}(z/a),\text{ for } z \in\,]-a, a[
 \label{rd08_26L6}
\end{align}
(see \cite[15.16, p.~151]{oberhettinger1973});
\smallskip

\beq
\label{rd08_26L7}
     f(t) = a^{-1} \sin(at),\qquad\qquad  F(z) = (z^2 + a^2)^{-1} ,
\eeq
\beq
\label{rd08_26L7a}
     f(t) = a^{-1} \sinh(at),\qquad\qquad  F(z) = (z^2 - a^2)^{-1}
\eeq
(see \cite[2.31 \& 2.32, p.~219]{oberhettinger1973});
\smallskip
\beq
\label{rd08_26L8}
   f(t) =  a^{-3} [\sinh(at) - at], \qquad\qquad  F(z) = z^{-2} (z^2 - a^2)^{-1}
\eeq
(see \cite[2.38, p.~220]{oberhettinger1973});
\smallskip
\beq
\label{rd08_26L9}
   f(t) = [\Gamma_E(\nu)]^{-1} t^{\nu - 1} e^{\mp a t}, \ \mathrm{Re} (\nu) > 0,\qquad F(z) = (z\pm a)^{-\nu}
\eeq
(see \cite[4.2, p.~237]{oberhettinger1973});
\smallskip
\beq
\label{rd01_22L1}
   f(t) = \pi^{-\half}\, t^{-\half} - a\, e^{a^2 t} \,\text{Erfc}(a t^\half), \qquad F(z) = (z^\half + a)^{-1}
\eeq
(here,  see \cite[(4), p.~233]{EMOT}).
 \bigskip

\noindent{\em Some Fourier-cosine transforms}
\bigskip

Given an integrable function $f: \re_+ \to \re$, we denote its Fourier-cosine transform by
\beqn
   \cF_c f(\xi) = \int_0^\infty f(x) \, \cos(\xi\, x)\, dx.
\eeqn
Given $f$, we have the following formulas for $\cF_c$:
\begin{align}
\nonumber
   &f(x) = I_0\left(b (a^2 - x^2)^\half\right)\, 1_{x < a}, \\
     & \qquad\qquad\qquad\qquad \cF_c f(\xi) = \left\{\begin{array}{ll} 
      (b^2 - \xi^2)^{-\half} \sinh\left(a(b^2 - \xi^2)^\half\right)
      &\text{if } \xi < b,\\ 
      (\xi^2 - b^2)^{-\half} \sin\left(a(\xi^2 - b^2)^\half\right)
      &\text{if } \xi > b
     \end{array} \right.
\label{rd08_26L11}
\end{align}
(see \cite[18.31, p.~95]{oberhettinger1990}).
\begin{align}\nonumber
    f(x) &= (a^2 + x^2)^{-\nu - \half}, \quad \mathrm{Re}(\nu) > -\tfrac{1}{2}, \\
     &   \qquad\qquad\qquad\qquad\qquad \cF_c f(\xi) =  \tfrac{\pi^\half  }{\Gamma_E\left(\tfrac{1}{2}+\nu\right)}\, \left(\tfrac{\xi}{2a} \right)^\nu \, K_\nu(a \xi)
\label{rd08_26L12}
\end{align}
(see \cite[2.7, p.~8]{oberhettinger1990}).

\section{Notes on Appendix \ref{app3}}
\label{app3-notes}

The version of Gronwall's lemma presented in Lemma \ref{A3-l01} concerns a Gronwall-Volterra-type inequality that appears often in the study of SPDEs. A partial version  for $J(s) = s^a$ ($a > -1$) appeared in \cite {henry-1981}, and later in \cite[Lemma 3.3, p.~316]{walsh}. However, instances of inequality \eqref{A3.03} were already present in D.A.~Dawson's study of stochastic evolution equations \cite{dawson72} and were resolved there via Laplace transform methods. The statement given in Lemma \ref{A3-l01} can mostly be found in \cite[Lemma 15]{dalang}, where a probabilistic proof is given. The renewal function and the reference \cite{asmussen} appears in the SPDE literature in \cite{khosh} and \cite{chong-2017}. The proof given in Section \ref{app3-s1} uses some ideas presented in \cite{asmussen}. The space-time Gronwall-type Lemma \ref{rd08_10l2} unifies certain ideas from \cite{chen} and \cite{candil2022}. It is used in the proof of Theorem \ref{rd08_09t1} and several times in Chapter \ref{ch2'}. This type of result is also used in \cite{joq} in their study of continuity in law of solutions to SPDEs driven by a continuous family of Gaussian noises.

Statements similar to Proposition \ref{app3-s4-new-p1} for a variety of examples of SPDEs are mentioned in many papers (see for example \cite{walsh} and \cite{gy98}). For a class of non-autonomous SPDEs, the equivalence between a weak formulation of the solution and a mild formulation has been addressed in \cite{ss-v2003}.

In Section \ref{AppC-New-0}, we sketch a general approach, via Fourier transforms and Fourier series, for obtaining the fundamental solution or the Green's function of a partial differential operator. 

Stability for PDEs, SDEs and SPDEs is an important and vast topic. The aim of Section \ref{appC-New-1} is to give a specific result tailored to its direct application in the proof of Proposition \ref{nbr-1} and Lemma \ref{nbr-*3a}. Proposition \ref{nbr-*3aa} and its proof are mostly taken from \cite{zabczyk-1989} (see also \cite{H-S-V-W-2007}). Lemma \ref{nbr-l-11} is quoted from \cite{H-S-V-W-2007}. 

The basic properties in Section \ref{app3-s3} concerning Gaussian random variables and the heat kernel are used throughout the book. In particular, the properties of derivatives of the heat kernel are applied in Chapter \ref{chapter1'} to the study of regularity of the difference of solutions to the stochastic heat equations on the line and on an interval (see Theorem \ref{ch1'-tDL}). Lemma \ref{rough-heat-l1} contains two technical results that are used in the proof of Proposition \ref{rough-heat}. They can be found in \cite[Lemmas A.4 and A.5]{chendalang2015-2}.

Lemma \ref{A3-l0} is quoted from \cite{chendalang2015} (see also \cite[5.9.5, p.~140]{olver}. Lemma \ref{lemB.5.2} provides explicit formulas for several trigonometric series; their proofs are exercices in Fourier expansions. They are used mainly in Sections \ref{stwn-s2} and \ref{ch5-added-1-s2}, and in Appendix \ref{app2}. Section \ref{app3-s2} contains technical results on periodic extensions of functions, and are used in Sections \ref{ch4-ss2.2-1'} and \ref{ch4-ss2.3-1'}.
Finally, Section \ref{app3LF} assembles formulas for Laplace and Fourier transforms of several functions. These are mostly used in Appendix \ref{rd08_26s1}, and  can be found for instance in \cite{EMOT}, \cite{oberhettinger1973} and \cite{oberhettinger1990}.

\backmatter
\addcontentsline{toc}{chapter}{Bibliography}


\bibliographystyle{amsalpha}


\chapter{List of Notations}
\label{notations}

\pagestyle{myheadings}
\markboth{R.C.~Dalang and M.~Sanz-Sol\'e}{Notations}

$\IN=\{0,1,2,\ldots\}$ is the set of natural numbers.\\
$\IN^* = \IN \setminus \{0\}$.\\
$\mathbb{Z}$\ is the set of integers.\\
$\re$ is the set of real numbers.\\
$\IR_+$ is the interval $[0,\infty[$.\\
$\R_+^*$ is the interval $]0, \infty[$.\\
$\IR_+^k$\ is the orthant $[0,\infty[^{\,k}$.\\
$\mathbb{C}$ is the field of complex numbers. When $i\in\mathbb{C}$, then $i^2 = -1$.\\
$n!$ denotes the factorial of $n \in \N$ ($n! = n(n-1) \cdots 2 \cdot 1$ when $n\in \N^*$, $0!$ = 1).\\
$x\wedge y$ is the notation for $\min(x,y)$, when $x, y\in\re$.\\
$x\vee y$ is the notation for $\max(x,y)$, when $x, y\in\re$.\\
$\text{sgn}(x)$ denotes the sign of $x \in \R$.\\
$t \to t_0$ indicates that the real variable $t$ converges to $t_0$.\\
$t \uparrow t_0$ indicates that the real variable $t$ converges to $t_0$ from the left.\\
$t \downarrow t_0$ indicates that the real variable $t$ converges to $t_0$ from the right.\\
$x_n\uparrow x$ denotes convergence of a non-decreasing sequence $(x_n,\, n \in \N)$ to $x$.\\
$x_n\downarrow x$ denotes convergence of a non-increasing sequence $(x_n,\, n \in \N)$ to $x$.\\
$\lfloor x \rfloor$ denotes the integer part of $x \in \R$. \\
$\lceil x \rceil$ denotes the ceiling function of $x \in \R$ \\
${\mathrm{Re}}(z)$ denotes the real part of $z \in \bC$. \\
${\mathrm{Im}}(z)$ denotes the imaginary part of $z \in \bC$. \\
$\bar z$ denotes the complex conjugate of $z$, when $z \in \mathbb{C}$.
\medskip

\noindent
$[a,b]$, $]a,b[$ denote a closed (respectively open) interval of $\re$ when $a,b \in \re$ and $a \leq b.$\\
$[a, b[$ denotes an interval that is closed on the left and open on the right.
$]a, b]$ denotes an interval that is open on the left and closed on the right. \\
$D$ is a domain of $\IR^k$, that is, a non-empty open connected subset of $\IR^k$.\\
$\bar D$ is the closure of the domain $D$.\\
$\partial D$ is the boundary of the domain $D$.\\
$A \cup B$ denotes the union of the two sets $A$ and $B$.\\
$A \cap B$ denotes the intersection of the two sets $A$ and $B$.\\
$A^c$ often denotes the complement of the set $A$.\\
$\emptyset$ denotes the empty set.\\
$:=$ indicates that the value on the right-hand side is assigned to the symbol on the left-hand side. \\
$=:$ indicates that the value on the left-hand side is assigned to the symbol on the right-hand side. \\

\medskip

\noindent
$\mathcal{C}(D)$ space of real-valued continuous functions defined on $D\subset \rek$.\\
$\ca(\IR_+\times D)$\ space of real-valued continuous functions defined on $\IR_+\times D$, where $D\subset \IR^k$.\\
$\ca(\IR_+; D)$\ space of $D$-valued continuous functions defined on $\IR_+$.\\
${\ca}_0^\infty(D)$\ space of real-valued $\mathcal{C}^\infty$ functions defined on $D\subset\IR^k$ and with compact support.\\
 ${\ca}^\eta(D)$\ space of H\"older continuous functions defined on $D\subset\IR^k$ with exponent $\eta\in\, ]0,1]$, \pageref{rd06_14l1}\\
 $\cC_0^\eta([0, L])$\ space of Hölder continuous functions $f$ defined on $[0, L]$ with exponent $\eta \in\, ]0, 1]$ such that $f(0) = f(L) = 0$, \pageref{rd06_14l2}\\
 $\cC^{\eta_1, \eta_2}(A \times B)$ space of jointly Hölder continuous functions defined on $A \times B$, where $A \subset \R$ and $B \subset \R^k$, with exponents $(\eta_1,\eta_2) \in\, ]0, 1]^2$, \pageref{rdjholder}\\
 $\mathbb{D}$ space of continuous functions $f$ on $[0, L]$ with $f(0) = f(L) = 0$.\\
$\mathcal{D}(\IR^k)$\ this is another notation for ${\ca}_0^\infty(\IR^k)$, used in the theory of distributions, \pageref{rdD}\\
$\mathcal{D}^\prime(\IR^k)$\ space of distributions, dual of $\mathcal{D}(\IR^k)$, \pageref{rdD'} \\
$H^n(D)$\  the Sobolev space $W^{n,2}(D)$, $D$ a domain in $\re^k$, $n\in \mathbb{Z}$.\\
$H^n_0(D)$\  the Sobolev space $W^{n,2}_0(D)$, $D$ a domain in $\re^k$, $n\in \mathbb{Z}$.\\
$L^p([0,T];V)$ for $p \geq 1$, denotes the space of equivalence classes (with respect to Lebesgue measure) of functions defined on $[0,T]$ with values in $V$ (a Hilbert space) such that $\int_0^T \Vert f(t)\Vert^p_V\, dt < \infty$.\\
$L^p(\R^k,\nu)$ for $p \geq 1$, denotes the  space of equivalence classes (with respect to a $\sigma$-finite  measure $\nu$) of real-valued functions defined on $\R^k$ such that $\int_{\R^k}|f(x)|^p\, \nu(dx) < \infty$.\\
$L^p_{\rm{loc}}(dt dx)$ for $p \geq 1$, denotes the space of equivalence classes  (with respect to a $\sigma$-finite  measure $dt dx$) of real-valued functions $f$ defined on $\re_+\times \rek$ such that $|f|^p$ is locally integrable.\\
$\mathcal{S}(\IR^k)$\  space of $\mathcal{C}^\infty$ functions with rapid decrease, also called Schwartz test functions, \pageref{rdS}\\
$\mathcal{S}^{\prime}(\IR^k)$\  space of tempered distributions, also termed Schwartz distributions, \pageref{rdS'}\\
$\mathcal{S}_n(\R^k)$ space of functions with $n$ generalized square-integrable derivatives, \pageref{rdSn}\\
 $\mathcal{S}_n^\prime(\R^k)$ is a subset of $\mathcal{S}^\prime(\R^k)$, dual of $\mathcal{S}_n(\R^k)$, \pageref{rdSn'}\\
$\delta_0$ often denotes the Dirac delta function.\\
$\delta_n^m$ denotes the Kroneker symbol: $\delta_n^m = 1$ if $n=m$ and $\delta_n^m = 0$ otherwise,  \pageref{rd05_04i1} \\
$e_{n, L}$  elements of a complete orthonormal system (CONS) in $L^2([0, L])$, $L > 0$.\\  
   $\leq_{HS}$  an ordering of Hilbertian norms, \pageref{rdHS-weaker}\\ 
$\langle \varphi, \psi \rangle_V$ denotes the inner product of two elements $\varphi$ and $\psi$ of the Hilbert space $V$ \\
$\Pi_{V_n}$ orthogonal projection onto the subspace $V_n$ of a Hilbert space $V$, \pageref{rdorthproj}\\
 $\Pi_{V_n^\perp}$ orthogonal projection onto the orthogonal complement of the subspace $V_n$ of a Hilbert space $V$, \pageref{rdorthperp} \\

\medskip

\noindent
$\mathcal{B}_{A}$\ Borel $\sigma$-field of $A\subset\IR^k$.\\
$\mathcal{B}^f_{A}$\  Borel subsets of $A$ with finite Lebesgue measure (or finite measure for some reference measure $\nu$ on $A$), where $A  \subset \R^k$.\\
$\cB_{\cC}$ Borel $\sigma$-field on the complete separable metric space $\cC$.\\
$\mathcal{G}_1\times \mathcal{G}_2$ denotes the product $\sigma$-field of two $\sigma$-fields $\mathcal{G}_1$ and $\mathcal{G}_2$.\\
 $\mathbb{H}^2$ space of continuous martingales, \pageref{rdH2}\\
$\text{Prog}$ is the $\sigma$-field of progressively mesurable sets, \pageref{rdProg}\\
$\mathcal{P}$ is the $\sigma$-field of predictable sets, \pageref{rdpred}\\
 $\mathcal{O}$ is the $\sigma$-field of optional sets, \pageref{rdOptional}\\
\medskip

\noindent
 $({\bf H_\Gamma})$ set of assumptions on the fundamental solution/Green's function, \pageref{rdHGamma}\\ 
$({\bf H_{\Gamma,I}})$ set of assumptions on the fundamental solution and initial conditions, \pageref{rdHGammaI}. \\
$({\bf H_{\Gamma,\infty}})$ global-in-time version of $({\bf H_\Gamma})$, \pageref{rdHGammainfinity}\\   
$({\bf H_{\Gamma-{\text{sup}},\infty}})$ strengthening of $({\bf H_\Gamma,\infty})$, \pageref{rdHGsinfinity}\\
 $({\bf H_I})$  assumptions on the initial condition, \pageref{rdHI}\\ 
$({\bf H_{I,\infty}})$ global-in-time version of $({\bf H_I})$, \pageref{rdHIinfinity}\\   
 $({\bf H_L})$  assumption on the coefficients of an SPDE, \pageref{rdHL}\\ 
$({\bf H_{L,\infty}})$ global-in-time version of $({\bf H_L})$), \pageref{rdHLinfinity}\\   
(${\bf H_w}$) an assumption used when discussing weak uniqueness in law, \pageref{Hw} \\
$({\bf h _\Gamma})$ set of assumptions on the Green's function, \pageref{rdhgamma}\\ 
$({\bf h_I})$ assumptions on the initial condition, \pageref{rdhI}\\ 
$({\bf h_L})$ assumption on the coefficients of an SPDE, \pageref{rdhL}\\ 
\medskip

\noindent
$|\cdot|$ denotes the absolute value of a real number, the Euclidean norm in $\R^k$ or the modulus of a complex number.\\
${\bf \Delta}_i(t,x; s,y)$ usually a metric on $\R_+ \times \R^k$, \pageref{rdDeltai1}, \pageref{rdDeltai2}\\  
$\Vert \cdot\Vert_{L^p(\Omega)}$ denotes the $L^p(\Omega)$-norm.\\
$\Vert \cdot \Vert_{T, \infty, p}$ a norm on random fields, \pageref{rdTinftyp}\\ 
 $\Vert \cdot \Vert_{t, \infty}$ a norm on random fields, \pageref{rdnormtinfty}\\ 
\medskip

\noindent
$|A|$ denotes de Lebesgue measure of a measurable set $A\subset\rek$.\\
 a.a.~is the abbreviation for ``almost all''. In general, it refers to a measure.\\
 a.e.~is the abbreviation for ``almost everywhere''. In general, it refers to a measure.\\
 a.s.~is the abbreviation for ``almost surely''. In general, it refers to a probability measure $P$.\\
 i.i.d.~is the abbreviation for ``independent and identically distributed".\\
 $X\stackrel{\rm d}{=}Y$ means equality in law for the random variables $X$ and $Y$.\\
 $X_n\stackrel{w}\longrightarrow X$ denotes convergence in law (also called ``in distribution'') of a sequence of random variables $(X_n,\, n\ge 1)$ to the random variable $X$.\\
  $u(\cdot, *)$  designates the function $(t,x) \mapsto u(t,x)$.\\
 $u(\cdot, x)$ for a function $(t,x) \mapsto u(t,x)$ designates the partial function $t \mapsto u(t,x)$ ($x$ fixed).\\
 $u(t, *)$ for a function $(t,x) \mapsto u(t,x)$, designates the partial function $x \mapsto u(t,x)$ ($t$ fixed).\\
 \medskip

\noindent
 $f \ast g$  denotes convolution in the space variable of two functions $f(x)$ and $g(x)$, \pageref{4.4-1'} \\
$f \star g$  denotes the space-time convolution of two functions $f(t, x)$ and $g(t, x)$, \pageref{rd08_14e8}\\
$f \smallstar g$  denotes convolution in the time variable of two functions $f(t)$ and $g(t)$, \pageref{rd04_16ss1}\\
$f(t) \lesssim g(t)$ indicates that as $t \to \infty$, $f$ is less than or equal to $g$, \pageref{rd06_19e5}\\ 
$\Delta$ is the Laplacian operator in $\IR^k$, \pageref{rdlaplacian}\\
$(-\Delta)^{a/2}$ is the fractional Laplacian operator in $\IR^k$, \pageref{rdflaplacian}\\
$\cL$ is a space-time partial differential operator, \pageref{ch1'-s1}\\
$\cL^*$, $\cL^{\star}$ denote the adjoint of $\cL$ (or its opposite), \pageref{rd06_13p1}, \pageref{rd06_13p2} \\
$\cF$ denotes the Fourier transform, \pageref{defi-fourier} \\
$\hat F$ denotes the Laplace transform, \pageref{app3LF} \\
$\cLL$ denotes the Laplace transform, \pageref{app3LF} \\
$\hcF$ denotes the Fourier-Laplace transform, \pageref{rd08_26e3} \\
$\hcFT$ denotes the truncated Fourier-Laplace transform, \pageref{rd08_26e3} \\
\medskip

\noindent
$\bfone$ denotes the constant function always equal to $1$, \pageref{rd12_19_e3}\\
$B_E$ denotes the Eular Beta function,  \pageref{rd08_12e9} \\
$E_{\alpha, \beta}$ denotes the  Mittag-Leffler function,  \pageref{E:Mittag-Leffler} \\  
Erfc$(t)$ denotes the Error function, \pageref{rd02_05e1} \\ 
$\Phi$ denotes the standard Normal probability distribution function, \pageref{rd05_07e1}\\
$\Gamma_E$ denotes the Eular Gamma function,  \pageref{Euler-gamma} \\
$g_a(t, x)$: a function whose tails are comparable to those of the fractional heat kernel, \pageref{rd08_17e6} \\
$I_0$ denotes the modified Bessel function of the first kind of order $0$, \pageref{rd08_26e5} \\
$I_\nu$ denotes the modified Bessel function of the first kind of order $\nu$, \pageref{rd08_26L10} \\
$J_\nu$ : denotes the Bessel function of the first kind of order $\nu$, \pageref{rd08_26L10a} \\ 
$P_\nu^\mu(z)$ denotes Legendre functions, \pageref{rd02_05L1} \\
 
\vfill
\eject

\pagestyle{myheadings}
\markboth{R.C.~Dalang and M.~Sanz-Sol\'e}{Index}

\newpage
   \addcontentsline{toc}{chapter}{Index} 
    \printindex 
\end{document}